\documentclass{book}
\pdfoutput=1
\usepackage{graphicx}
\usepackage{makeidx}   

\makeindex

\def\item#1{\par\noindent{\rm#1}\hangindent\parindent}  

\newcommand{\mychapter}[2]{\vspace*{0.5in}%
 \centerline{\Large #1.}\centerline{\vrule width 3cm height .8pt}%
 \par\medskip\par\smallskip\centerline{\Large #2}\bigskip%
 \thispagestyle{empty}\markboth{Kapitel #1. #2}{}%
 \addcontentsline{toc}{chapter}{\protect\numberline{#1.~#2}}%
 \noindent}

\newcommand{\mysection}[1]{\goodbreak\bigskip\vspace{0pt plus \bigskipamount}%
\centerline{\bf#1}\medskip\vspace{0pt plus 4pt}\markright{\protect#1}%
\addcontentsline{toc}{section}{\quad\protect{#1}}}

\newcommand{\mysectionten}[1]{\goodbreak\bigskip\vspace{0pt plus \bigskipamount}%
\centerline{\bf#1}\medskip\vspace{0pt plus 4pt}\markright{\protect#1}%
\addcontentsline{toc}{section}{\quad\hspace{-5pt}\protect{#1}}}

\let\sm=\smallskip
\renewcommand\smallskip{\par\sm}
\let\mm=\medskip
\renewcommand\medskip{\par\mm}
\let\bm=\bigskip
\renewcommand\bigskip{\par\bm}

\makeatletter
\def\eqalign#1{%
\null\,\vcenter{\openup\jot\m@th
  \ialign{\strut\hfil$\displaystyle{##}$&$\displaystyle{{}##}$\hfil
      \crcr#1\crcr}}\,}
\makeatother

\def\pol{{\rm pol}}  
\def\nleq{\not\le}  
\def\anff{\vbox to 1.0pt{\hbox{''}\vss}}  
\let\glqq=\anff
\def\grqq{\hbox{``}}  
\def\Alt{{\rm Alt}}
 \def\ann{{\rm ann}}
 \def\Ant{{\rm Ant}}
 \def\Aut{{\rm Aut}}
 \def\Beta{{\rm B}}    
 \def\bt{\bigtriangleup}
 \def\cart{{\rm cart}}
 \def\Char{{\rm Char}}
 \def\cont{{\rm cont}}
 \def\det{{\rm det}}
 \def\E{{\rm E}}
 \def\emph#1{{\it#1\/}}
 \def\End{{\rm End}}
 \def\Fin{{\rm Fin}}
 \def\GaL{{\rm G}^*{\rm L}}
 \def\GammaL{\Gamma {\rm L}}
 \def\GF{{\rm GF}}
 \def\GL{{\rm GL}}
 \def\ggT{{\rm ggT}}
 \def\Hom{{\rm Hom}}
 \def\I{\mathrel{\hbox{\rm I}}}      
 \def\Inn{{\rm Inn}}
 \def\Inz{(\Pi,\Gamma,\hbox{\rm I})}    
 \def\Kern{{\rm Kern}}
 \def\kgV{{\rm kgV}}
 \def\Ko{{\rm KoRg}}
 \def\La{{\rm L}}
 \def\lhak{\>{\vrule width 5pt height .8pt depth 0pt}
       {\vrule width .8pt height 4.5pt depth 0pt}\>}
 \def\Lin{{\rm Lin}}
 
 \def\mod{{\rm \ mod\ }}
 \def\notI{\mathrel{\hbox{\rm\rlap/I}}}      
 \def\PGaL{{\rm PG}^*{\rm L}}
 \def\PGammaL{{\rm P}\Gamma{\rm L}}
 \def\PGL{{\rm PGL}}
 \def\PGU{{\rm PGU}}
 \def\PSL{{\rm PSL}}
 \def\PSp{{\rm PSp}}
 \def\PSU{{\rm PSU}}
 \def\PTU{{\rm PTU}}
 \def\R{{\bf R}}
 \def\rad{{\rm rad}}
 \def\rhak{\>{\vrule width .8pt height 4.5pt depth 0pt}
       {\vrule width 5pt height .8pt depth 0pt}\>}
 \def\Rg{{\rm Rg}}
 \def\sgn{{\rm sgn}}
 \def\SL{{\rm SL}}
 \def\Sp{{\rm Sp}}
 \def\Spur{{\rm Spur}}
 \def\SU{{\rm SU}}
 \def\T{{\rm T}}
 \def\Tau{{\rm T}}
 \def\TU{{\rm TU}}
 \def\Un{{\rm U}}
 \def\U{{\rm Un}}
 \def\UR{{\rm UR}}

\begin{document}

\frontmatter

 \thispagestyle{empty}

 \noindent\vglue 1in

 \centerline{\huge Projektive Geometrie}

 \bigskip\bigskip

 \centerline{\Large von}
 
 \bigskip\bigskip

 \centerline{\Large Heinz L\"uneburg}

 \vfill
 \centerline{Aus dem Nachlass des Autors herausgegeben}

 \centerline{von Theo Grundh\"ofer und Karl Strambach}

 \vfill

\chapter*{}
\thispagestyle{empty}
\leftline{\LARGE Vorwort der Herausgeber}

\bigskip
\noindent
Heinz L\"uneburg ist am 19.\ Januar 2009 pl\"otzlich verstorben. \"Uber sein Buchprojekt
\anff Projektive Geometrie`` war auf seiner homepage zu lesen:
\smallskip

\leftskip=\parindent
\rightskip=\parindent
\noindent
Dieses Buch, im WS 1967/68 begonnen, sollte mein erstes Buch werden. Dass es dies nicht
wurde, lag an meinem Wechsel nach Kai\-sers\-lau\-tern. Vier Kapitel waren damals mit meiner
Rei\-se\-schreib\-ma\-schi\-ne aufgeschrieben. Hier hatten wir aber nur Studenten im ersten Semester,
so dass ich keine Vor\-le\-sung\-en \"uber den Gegenstand halten konnte, wobei jedoch die
\anff Einf\"uhrung in die Algebra`` von den Vorarbeiten f\"ur jenes Buch profitierte.
Es blieb also liegen. Es blieb liegen, bis ich die Vorlesungen \"uber Fibonacci begann. Diese
Vorlesungen waren einst\"undig. Die eine Stunde bedurfte jedoch einer Woche an Vorbereitung.
Also griff ich in die Konserve und hielt neben der Vorlesung \"uber Fibonacci noch eine
vierst\"undige Vorlesung \"uber projektive Geometrie. ...
Das Buch liegt also jetzt da, 525 Seiten im ehemaligen B.I.-Format, und es fehlt noch ein
Kapitel. Wann ich dieses schreiben kann, wei\ss{} ich nicht. Vielleicht publiziere
ich es ohne das fehlende Kapitel \"uber orthogonale Gruppen.

\leftskip=0pt
\rightskip=0pt

\smallskip
\noindent
Leider enth\"alt der Nachlass von Heinz L\"uneburg nur einen rudiment\"aren Anfang
dieses Kapitels \"uber orthogonale Gruppen, n\"amlich die Konstruktion von Cliffordalgebren, und
keine Spuren von einem weiteren geplanten Kapitel \anff An\-wen\-dung\-en``. Mit der freundlichen
Zustimmung von Frau Karin L\"uneburg machen wir die im Nachlass
vorhandenen Kapitel der \"Offentlichkeit zug\"anglich.

Die Projektive Geometrie von Heinz L\"uneburg spiegelt die zeitlichen Ver\-\"an\-de\-rung\-en wider,
welchen auch die Mathematik unterworfen ist. Das Buch beginnt mit einem verbandstheoretischen
Aufbau, gem\"a\ss\ der grundlagentheoretischen Perspektive der 1960er Jahre. In den Kapiteln
IIII und V wird die L\"uneburgsche Sicht auf die endliche Geometrie und auf Polarit\"aten
deutlich. Die Kapitel VI und VII zeigen, dass Heinz L\"uneburg sich f\"ur die Burausche Welt
der synthetischen algebraischen Geometrie erw\"armen konnte und zu dieser Theorie moderne,
einsichtige Beweise geliefert hat.

Ein Vortrag von Karl Strambach zur Erinnerung an Heinz L\"uneburg ist als Anhang beigef\"ugt.

\bigskip\medskip
\noindent{W\"urzburg, im Juni 2011}\hfill
{Theo Grundh\"ofer und Karl Strambach}

\chapter*{}
\markboth{Inhaltsverzeichnis}{Inhaltsverzeichnis}
\leftline{\LARGE Inhaltsverzeichnis}

\contentsline {chapter}{\numberline {I.\nobreakspace {}Die Grundlagen und ein bisschen mehr}}{1}
\contentsline {section}{\hskip 1em\relax {1. Projektive R\"aume}}{1}
\contentsline {section}{\hskip 1em\relax {2. Projektive Verb\"ande}}{8}
\contentsline {section}{\hskip 1em\relax {3. Der Basissatz}}{16}
\contentsline {section}{\hskip 1em\relax {4. Vollst\"andig reduzible Moduln}}{19}
\contentsline {section}{\hskip 1em\relax {5. Der duale Verband}}{30}
\contentsline {section}{\hskip 1em\relax {6. Homogene, vollst\"andig reduzible Ringe}}{44}
\contentsline {section}{\hskip 1em\relax {7. Endliche projektive R\"aume}}{60}
\contentsline {section}{\hskip 1em\relax {8. Kollineationen und Korrelationen}}{65}
\contentsline {section}{\hskip 1em\relax {9. Der Satz von Desargues}}{75}
\contentsline {section}{\hskip 1em\relax \hspace {-5pt}{10. Der Satz von Pappos}}{77}
\contentsline {chapter}{\numberline {II.\nobreakspace {}Die Strukturs\"atze}}{81}
\contentsline {section}{\hskip 1em\relax {1. Zentralkollineationen}}{81}
\contentsline {section}{\hskip 1em\relax {2. Der Kern von E(H)}}{91}
\contentsline {section}{\hskip 1em\relax {3. Pappossche Geometrien}}{94}
\contentsline {section}{\hskip 1em\relax {4. Der Satz von Hessenberg}}{97}
\contentsline {section}{\hskip 1em\relax {5. Weniger Bekanntes aus der linearen Algebra}}{104}
\contentsline {section}{\hskip 1em\relax {6. Der erste Struktursatz}}{108}
\contentsline {section}{\hskip 1em\relax {7. Der zweite Struktursatz}}{110}
\contentsline {section}{\hskip 1em\relax {8. Der dritte Struktursatz}}{111}
\contentsline {section}{\hskip 1em\relax {9. Quaternionenschiefk\"orper}}{116}
\contentsline {section}{\hskip 1em\relax \hspace {-5pt}{10. Projektive R\"aume mit Clifford-Parallelismus}}{130}
\contentsline {chapter}{\numberline {III.\nobreakspace {}Gruppen von Kollineationen}}{147}
\contentsline {section}{\hskip 1em\relax {1. Erste \"Ubersicht}}{147}
\contentsline {section}{\hskip 1em\relax {2. Die Einfachheit der kleinen projektiven Gruppe}}{155}
\contentsline {section}{\hskip 1em\relax {3. Determinanten}}{161}
\contentsline {section}{\hskip 1em\relax {4. Ausnahmeisomorphismen}}{171}
\contentsline {section}{\hskip 1em\relax {5. Quasiperspektivit\"aten}}{183}
\contentsline {section}{\hskip 1em\relax {6. Zentralisatoren von Involutionen}}{197}
\contentsline {section}{\hskip 1em\relax {7. Die hessesche Gruppe}}{204}
\contentsline {section}{\hskip 1em\relax {8. Isomorphismen der gro\ss en projektiven Gruppe}}{212}
\contentsline {chapter}{\numberline {IIII.\nobreakspace {}Endliche projektive Geometrien}}{225}
\contentsline {section}{\hskip 1em\relax {1. Endliche Inzidenzstrukturen}}{225}
\contentsline {section}{\hskip 1em\relax {2. Inzidenzmatrizen}}{229}
\contentsline {section}{\hskip 1em\relax {3. Kollineationen von projektiven Blockpl\"anen}}{232}
\contentsline {section}{\hskip 1em\relax {4. Korrelationen von projektiven Blockpl\"anen}}{235}
\contentsline {section}{\hskip 1em\relax {5. Taktische Zerlegungen}}{238}
\contentsline {section}{\hskip 1em\relax {6. Endliche desarguessche projektive Ebenen}}{241}
\contentsline {section}{\hskip 1em\relax {7. Ebenen mit vielen Perspektivit\"aten}}{251}
\contentsline {section}{\hskip 1em\relax {8. Einiges \"uber Permutationsgruppen}}{255}
\contentsline {section}{\hskip 1em\relax {9. Geradenhomogene affine Ebenen}}{258}
\contentsline {section}{\hskip 1em\relax \hspace {-5pt}{10. Endliche projektive R\"aume}}{267}
\contentsline {section}{\hskip 1em\relax \hspace {-5pt}{11. Ein Satz von N. Ito}}{270}
\contentsline {chapter}{\numberline {V.\nobreakspace {}Polarit\"aten}}{276}
\contentsline {section}{\hskip 1em\relax {1. Darstellung von Polarit\"aten}}{276}
\contentsline {section}{\hskip 1em\relax {2. Zentralisatoren von Polarit\"aten}}{281}
\contentsline {section}{\hskip 1em\relax {3. Symplektische Polarit\"aten und ihre Zentralisatoren}}{286}
\contentsline {section}{\hskip 1em\relax {4. Polarit\"aten bei Charakteristik 2}}{295}
\contentsline {section}{\hskip 1em\relax {5. Quadratische Formen}}{302}
\contentsline {section}{\hskip 1em\relax {6. Die wittsche Zerlegung}}{305}
\contentsline {section}{\hskip 1em\relax {7. Der Satz von Witt}}{311}
\contentsline {section}{\hskip 1em\relax {8. Unit\"are Gruppen}}{315}
\contentsline {section}{\hskip 1em\relax {9. Endliche unit\"are Gruppen}}{325}
\contentsline {section}{\hskip 1em\relax \hspace {-5pt}{10. Die speziellen unit\"aren Gruppen}}{336}
\contentsline {chapter}{\numberline {VI.\nobreakspace {}Segresche Mannigfaltigkeiten}}{345}
\contentsline {section}{\hskip 1em\relax {1. Tensorprodukte}}{345}
\contentsline {section}{\hskip 1em\relax {2. Homomorphismen projektiver R\"aume}}{354}
\contentsline {section}{\hskip 1em\relax {3. Segresche Mannigfaltigkeiten}}{366}
\contentsline {section}{\hskip 1em\relax {4. Geometrische Erzeugung Segrescher Mannigfaltigkeiten}}{376}
\contentsline {chapter}{\numberline {VII.\nobreakspace {}Gra\ss {}mannsche Mannigfaltigkeiten}}{384}
\contentsline {section}{\hskip 1em\relax {1. Die Gra\ss mannalgebra eines Moduls}}{384}
\contentsline {section}{\hskip 1em\relax {2. Dualit\"at in der Gra\ss mannalgebra}}{392}
\contentsline {section}{\hskip 1em\relax {3. Innere Produkte}}{400}
\contentsline {section}{\hskip 1em\relax {4. Zerlegbare Vektoren}}{405}
\contentsline {section}{\hskip 1em\relax {5. Doppelverh\"altnisse}}{410}
\contentsline {section}{\hskip 1em\relax {6. Gra\ss mannsche Mannigfaltigkeiten I}}{413}
\contentsline {section}{\hskip 1em\relax {7. Gra\ss mannsche Mannigfaltigkeiten II}}{421}
\contentsline {chapter}{\numberline {Anhang: Der Aufbruch der Geometrie}}{433}
\contentsline {chapter}{\numberline {Literatur}}{445}
\contentsline {chapter}{\numberline {Index}}{449}

\mainmatter


 \mychapter{I}{Die Grundlagen und ein bisschen mehr}

 \noindent
 Projektive Geometrien sind zun\"achst Inzidenzstrukturen,
 bestehend aus Punkten und Geraden mit gewissen Eigenschaften. Wie
 sich herausstellt, verbirgt sich dahinter eine reiche Struktur,
 n\"amlich der Verband der Unterr\"aume dieser Geometrien. In diesem
 Kapitel wird es nun vor allem darum gehen, diese Verb\"ande vom
 verbandstheoretischen Standpunkt aus in den Griff zu bekommen,
 soweit dies ohne Verwendung von Methoden der linearen Algebra, die
 hinter all dem steckt, ohne Zwang m\"oglich ist. Insbesondere werden
 wir den Basissatz beweisen und zeigen, dass eine projektive
 Geometrie, deren Dimension nicht gerade $2$ ist, stets
 desarguessch ist. Dies impliziert wiederum, dass ein projektiver
 Verband, dessen Dimension mindestens $3$ ist, dem Unterraumverband
 eines geeigneten Vektorraums isomorph ist. Da die Dimension dieses
 Vektorraumes um $1$ gr\"o\ss er ist als die Dimension der zugeh\"origen
 Geometrie, werden wir das Wort Dimension nur informell benutzen
 und stattdessen vom Rang einer projektiven Geometrie und vom Rang
 eines Vektorraumes reden.
 \par
       Nicht zu den Grundlagen dessen, was weiter folgt, geh\"oren die
 Ausf\"uhrungen \"uber vollst\"andig reduzible Ringe. Sie sind hier
 auf\-ge\-nom\-men, weil man an ihnen demonstrieren kann, wie n\"utzlich
 der Begriff des projektiven Verbandes ist. Hinzu kommt, dass die
 projektive Deutung der einschl\"agigen S\"atze Einsichten vermittelt,
 die man in Algebrab\"uchern\index{Algebrab\"ucher}{} vergeblich sucht.

 \mysection{1. Projektive R\"aume}

 \noindent
 Es sei $\Pi$ eine Menge, deren Elemente wir
 \emph{Punkte}, und $\Gamma$ eine Menge, deren Elemente wir
 \emph{Geraden} nennen. Ferner sei I eine Teilmenge des
 cartesischen Produktes $\Pi \times \Gamma$ von $\Pi$ und $\Gamma$.
 Ist $(P,G) \in {\rm I}$, so sagen wir, dass $P$ mit $G$
 \emph{inzidiere}, andernfalls, dass $P$ und $G$ \emph{nicht
 inzidierten}. Statt $P$ inzidiere mit $G$ werden wir auch andere
 Redewendungen wie $P$ \emph{liege auf} $G$ oder $G$ \emph{gehe
 durch} $P$ oder $P$ \emph{sei enthalten in} $G$ oder \"ahnliche
 verwenden. Statt $(P,G) \in {\rm I}$ bzw. $(P,G) \notin {\rm I}$ werden wir
 meist $P \I G$ bzw. $P \notI G$ schreiben. Das Tripel $\Inz$ hei\ss t
 \emph{Inzidenzstruktur}.\index{Inzidenzstruktur}{}
 \par
      Es sei $\Inz$ eine Inzidenzstruktur und $\Phi \subseteq \Pi$. Die Punkte
 von $\Phi$ hei\ss en \emph{kollinear},\index{kollinear}{} falls es ein $G \in
 \Gamma$ gibt mit $X \, \I \, G$ f\"ur alle $X \in \Phi$. Sind $P$
 und $Q$ zwei verschiedene kollineare Punkte und gibt es nur eine
 Gerade durch $P$ und $Q$, so bezeichnen wir diese mit $P + Q$.
 Diese Bezeichnung soll an die Addition von Unterr\"aumen eines
 Vektorraumes erinnern, da die so bezeichnete Verkn\"upfung am Ende
 nichts anderes als diese sein wird.
 \par
       Aus den un\"ubersehbar vielen Inzidenzstrukturen sondern wir die
 projektiven Geometrien\index{projektive Geometrie}{} durch die folgende
 Definition aus. Die Inzidenzstruktur $\Sigma := \Inz$ hei\ss t \emph{projektive
 Geometrie} oder auch \emph{projektiver Raum},\index{projektiver Raum}{} falls
 $\Sigma$ den folgenden Bedingungen gen\"ugt.\\
 {(P1)} Sind $P$ und $Q$ zwei verschiedene Punkte von $\Sigma$, so gibt es
 genau eine Gerade $G$ von $\Sigma$ mit $P$, $Q \I G$.\\
 {(P2)} Sind $P$, $Q$ und $R$ drei nicht kollineare Punkte von $\Sigma$ und
 sind $D$ und $E$ zwei weitere, verschiedene Punkte mit $D \I P + Q$ und
 $E \I Q + R$, so gibt es einen Punkte $F$ mit $F \I R + P$ und
 $F \I D + E$.\\
 {(P3)} Jede Gerade von $\Sigma$ tr\"agt wenigstens zwei Punkte.
 \vadjust{
 \vglue 4mm
 \centerline{\includegraphics{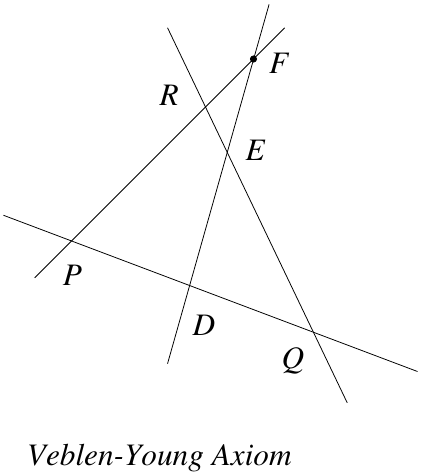}}
 \vglue 3mm}
 \smallskip
       Axiom (P$2$) wird h\"aufig
 \emph{Veblen-Young Axiom}\index{Veblen-Young Axiom}{} genannt. Man
 kann es salopp auch so formulieren: Trifft eine Gerade zwei Seiten
 eines Dreiecks in verschiedenen Punkten, so trifft sie auch die
 dritte Seite.
 \par
       Beispiele von projektiven Geometrien\index{projektive Geometrie}{} sind
 leicht zu verschaffen.  Die einfachsten Beispiele sind die folgenden: Ist $M$
 eine Menge und bezeichnet $P_{2}(M)$ die Menge der Teilmengen mit genau zwei
 Elementen von $M$, so ist $(M, P_{2}(M),\in)$ eine projektive Geometrie. So
 banal diese Beispiele auch erscheinen m\"ogen, so werden sie uns doch sp\"ater
 beim Beweise des Satzes von Wielandt\index{Wielandt, H.}{} \"uber die
 Automorphismengruppe der alternierenden Gruppe ein sehr n\"utzliches Werkzeug
 sein.
 \par
       Projektive Geometrien, deren Geraden mehr als zwei Punkte
 tragen, erh\"alt man folgenderma\ss en: Es sei $V$ ein Vektorraum \"uber
 dem K\"orper $K$. Bezeichnet man mit $\UR_{i}(V)$ die Menge der
 Unterr\"aume des Ranges $i$ von $V$, so ist $(\UR_{1}(V),
 \UR_{2}(V), \subseteq)$ eine projektive Geometrie. (Der Leser
 erinnere sich, dass Rang hier das ist, was anderenorts meist
 Dimension genannt wird.) Dies zu beweisen, sei dem Leser als
 \"Ubungsaufgabe \"uberlassen. Ein Tipp sei jedoch gegeben. Sind $pK$
 und $qK$ zwei verschiedene Punkte, so sind $pK$, $qK$ und $(p + q)K$ drei
 verschieden Punkte auf der Geraden $pK + qK$.
 \par
       Es seien $\Inz$ und $(\Pi', \Gamma', \I')$ zwei Inzidenzstrukturen. Ferner
 sei $\sigma$ eine Bijektion von $\Pi$ auf $\Pi'$ und $\tau$ eine solche von
 $\Gamma$ auf $\Gamma'$. Das Paar $(\sigma, \tau)$ hei\ss t
 \emph{Isomorphismus}\index{Isomorphismus}{} von $\Inz$ auf
 $(\Pi', \Gamma', \I')$, wenn f\"ur alle $P \in \Pi$ und alle $G \in
 \Gamma$ genau dann $P \I G$ gilt, wenn $P^{\sigma} \I' G^{\tau}$ ist. Gibt es
 einen solchen Isomorphismus, so nennen wir
 die beiden Inzidenzstrukturen \emph{isomorph}. Isomorphismen einer
 Inzidenzstruktur auf sich hei\ss en
 \emph{Automorphismen}\index{Automorphismus}{} oder
 \emph{Kollineationen}.\index{Kollineation}{}
 \medskip\noindent
 {\bf 1.1. Satz.} {\it Es sei $\Sigma := \Inz$ ein projektiver Raum.
 Ist $G \in \Gamma$, so setzen wir $G^{\tau} := \{ P \mid P \in \Pi, P
 \I G\}$. Dann ist $(id_{\Pi}, \tau)$ ein Isomorphismus von
 $\Sigma$ auf $(\Pi, \Gamma^{\tau}, \in)$.}
 \smallskip
       Beweis. Das Einzige, was zu beweisen ist, ist die Injektivit\"at von
 $\tau$. Diese folgt aber unmittelbar aus der Tatsache, dass jede
 Gerade mindestens zwei Punkte tr\"agt und zwei verschiedene Punkte
 nur mit einer Geraden inzidieren. \medskip
 \par
       Dieser Satz zeigt, dass man die Geraden\index{Gerade}{} eines projektiven
 Rau\-mes mit den Mengen\index{Menge von Punkten}{} der jeweils auf ihnen
 liegenden Punkte identifizieren kann. Dies ist immer wieder einmal bequem und
 wird, meist ohne es explizit zu sagen, dann auch getan.
 \par
       Es sei $\Sigma := \Inz$ ein projektiver Raum und $U \subseteq \Pi$. Die
 Menge $U$ hei\ss t \emph{Unterraum}\index{Unterraum}{} von $\Sigma$, falls mit
 zwei verschiedenen Punkten stets auch ihre ganze Verbindungsgerade in $U$ liegt.
 (Hier haben wir zum ersten Male Geraden mit Punktmengen identifiziert.)
 Bei\-spie\-le von Unterr\"aumen sind die leere Menge $\emptyset$, die Mengen,
 die nur aus einem Punkt bestehen, die Geraden und $\Pi$. Wir bezeichnen die
 Menge aller Unterr\"aume von $\Sigma$ mit $L(\Sigma)$. Sind $X$, $Y \in
 L(\Sigma)$ und ist $X$ in $Y$ enthalten, so bezeichnen wir diesen Sachverhalt
 mit $X \leq Y$. Wir benutzen in diesem Falle also das Zeichen $\leq$ anstelle
 von $\subseteq$. Die Relation $\leq$ ist \emph{reflexiv}, \emph{antisymmetrisch}
 und \emph{transitiv}, dh. sie ist eine \emph{Teilordnung}\index{Teilordnung}{}
 auf $L(\Sigma)$. Diese Teilordnung gilt es nun zu studieren.
 \par
       Damit der Leser ein Gef\"uhl daf\"ur bekomme, was
 Unterr\"aume\index{Unterraum}{} sind,
 beweise er das Folgende: Es sei $V$ ein $K$-Vektorraum. Ferner sei
 $X$ eine Teilmenge von $\UR_{1}(V)$. Genau dann ist $X$ ein
 Teilraum von $(\UR_{1}(V), \UR_{2}(V), \subseteq)$, wenn es einen
 Teilraum $W$ von $V$ gibt mit $X = \UR_{1}(W)$.
 \par
       Ist $M$ eine Menge, so ist $L (M, P_{2}(M), \in)$ nichts Anderes
 als die Potenzmenge\index{Potenzmenge}{} von $M$.
 \par
       Es ist im Weiteren bequem, die Unterr\"aume, die nur aus einem Punkt
 bestehen, mit diesem Punkt zu identifizieren, so dass sich also
 die Inzidenz des Punktes $P$ mit der Geraden $G$ durch $P \leq G$
 beschreiben l\"a\ss t. Damit sind der Begriff des Punktes und der
 Begriff der Geraden unter den Begriff des Unterraumes subsumiert,
 so dass die Ausnahmestellung dieser Objekte beseitigt ist.
 \par
       Es sei weiterhin $\Sigma$ eine projektive Geometrie. Ist $M \subseteq
 L(\Sigma)$, so ist, wie unmittelbar aus der Definition des Unterraumes folgt,
 auch $\bigcap_{X \in M} X \in L (\Sigma)$. Wie \"ublich definiert man daher
 den von $M$ \emph{auf\-ge\-spann\-ten Unterraum}\index{aufgespannter Unterraum}{}
 $\sum_{X \in M} X$ von $\Sigma$ als den Schnitt \"uber alle Unterr\"aume $Y$ von
 $\Sigma$, f\"ur die $X \leq Y$ f\"ur alle $X \in M$ gilt. Es gilt dann die
 folgende banale, aber n\"utzliche Aussage.
 \medskip\noindent
 {\bf 1.2. Satz.} {\it Es sei $\Sigma$ eine projektive Geometrie. Es
 sei ferner $M \subseteq L(\Sigma)$ und $Y \in L(\Sigma)$.
 Dann gilt: Ist $Y \leq X$ f\"ur alle $X \in M$, so ist
 $Y \leq \bigcap_{X \in M}X$. Ist $X \leq Y$ f\"ur alle $X \in M$, so ist
 $\sum_{X \in M} X \leq Y$.}
 \medskip
       Dieser Satz besagt, dass $\bigcap_{X \in M}X$ der gr\"o\ss te, in allen
 $X \in M$ enthaltene Unterraum von $\Sigma$ ist, w\"ahrend $\sum_{X \in M}X$ der
 kleinste, alle $X \in M$ umfassende Unterraum von $\Sigma$ ist.
 \par
       Die Definition des Operators $\sum$, der im \"Ubrigen das ist, was man
 einen \emph{H\"ullenoperator}\index{hullen@H\"ullenoperator}{} nennt, ist zwar elegant,
 man h\"atte aber auch gerne eine interne Beschreibung des von einer Menge von
 Unterr\"aumen aufgespannten Unterraums. Diese Beschreibung wird in
 K\"urze gegeben werden. Dazu definieren wir zun\"achst eine bin\"are
 Verkn\"upfung $\odot$ auf der Potenzmenge der Punktmenge $\Pi$ eines
 projektiven Raumes.
 \smallskip\noindent
 (1) Es ist $S \odot \emptyset = \emptyset \odot S = S$ f\"ur alle
 $S \subseteq \Pi$.
 \par\noindent
 (2) Ist $S$ eine einelementige Teilmenge von $\Pi$, so setzen wir
 $S \odot S = S$.
 \par\noindent
 (3) Sind $S$ und $T$ Teilmengen von $\Pi$ und gibt es zwei verschiedene Punkte
 $P$ und $Q$ mit $P \in S$ und $Q \in T$, so sei $S \odot T$ die Menge der
 Punkte, die auf den Geraden liegen, die zwei verschiedene Punkte $A$ und $B$ mit
 $A \in S$ und $\Beta \in T$ verbinden.
 \smallskip
       Die wesentlichen Eigenschaften dieser Operation sind im folgenden
 Satz aufgelistet.
 \medskip\noindent
 {\bf 1.3. Satz.} {\it Es sei $\Sigma := \Inz$ eine projektive
 Geometrie. Sind $S$, $T$ und $U$ drei Teilmengen von $\Pi$, so gilt:
 \item{(a)} Ist $S \subseteq T$, so ist $S \odot U \subseteq \Tau \odot U$.
 Insbesondere ist $S \subseteq S \odot S$.
 \item{(b)} Es ist $S \odot T = T \odot S$.
 \item{(c)} Ist $R \in L(\Sigma)$ und gilt $S$, $T \subseteq R$, so ist
 $S \odot T \subseteq R$.
 \item{(d)} Es ist $(S \odot T) \odot U = S \odot (T \odot U)$.
 \item{(e)} Genau dann gilt $S \odot S = S$, wenn $S \in L(\Sigma)$ ist.
 \item{(f)} Sind $S$, $T \in L(\Sigma)$, so ist $S \odot T = S + T$.}
 \smallskip
       Beweis. (a) Die erste Aussage von (a) folgt unmittelbar aus der
 Definition von $\odot$. Wendet man dies nun statt auf $S$, $T$ und
 $U$ auf $\emptyset$, $S$ und $S$ an, so folgt
 $$ S = \emptyset \odot S \subseteq S \odot S. $$
 \par
       (b) und (c) folgen unmittelbar aus der Definition von $\odot$ und
 der definierenden Eigenschaft von Unterr\"aumen mit zwei Punkten
 ihre Verbindungsgerade zu enthalten.
 \par
       Um (d) zu beweisen, nehmen wir zun\"achst an, dass alle drei Mengen
 nur aus einem Punkt bestehen, so dass wir sie gem\"a\ss\ unserer
 Konvention mit diesen Punkten identifizieren. Liegen $S$, $T$ und
 $U$ auf der Geraden $G$, so folgt mit (c), dass $(S \odot T) \odot U$ und
 $S \odot (T \odot U)$ in $G$ liegen. Ist $S \odot T = U$,
 so folgt, dass $S = T = U$ ist. In diesem Falle gilt also (d).
 Enth\"alt $S \odot T$ einen Punkt, der von $U$ verschieden ist, so
 ist $S \neq U$ oder $T \neq U$. Ist $S \neq U$, so folgt
 $$ G = S \odot U \subseteq (S \odot T) \odot U \subseteq G $$
 und
 $$ G = S \odot U \subseteq S \odot (T \odot U) \subseteq G, $$
 so dass (d) auch in diesem Falle gilt. Ist $T \neq U$, so folgt
 (d) ent\-spre\-chend.
 \par
       Es bleibt der Fall zu betrachten, dass $S$, $T$ und $U$ nicht kollinear
 sind. Dann sind sie insbesondere auch paarweise verschieden. Wir zeigen
 zun\"achst, dass $S \odot (T \odot U) \subseteq (S \odot T) \odot U$ gilt. Dazu
 sei $A \in S \odot (T \odot U)$. Auf Grund von (a) d\"urfen wir annehmen, dass
 $A \notin S \odot T$, $U$ gilt. Wegen $A \in S \odot (T \odot U)$ gibt es ein $B
 \in T \odot U$ mit $A \leq S + B$. W\"are $B = T$, so w\"are $A \leq S + T$ im
 Widerspruch zu $A \notin S \odot T$. Also ist $B$ ein von $T$ verschiedener
 Punkt in $T \odot U = T + U$, so dass insbesondere $T + U = T + B$ ist. Hieraus
 folgt weiter, dass $B$, $S$ und $T$ nicht kollinear sind. Weil $A \neq U$ ist,
 gibt es nach (P2) daher einen Punkt $C$ auf $S + T$ mit $C \leq A + U$. Nun ist
 $S + (T \cap T) + U = T$, so dass $C \neq U$ ist. Hiermit folgt zusammen mit
 $S + T = S \odot T$, dass
 $$ A \leq A + U = C + U = C \odot U \subseteq (S \odot T) \odot U $$
 ist. Damit haben wir gezeigt, dass $S \odot (T \odot U) \subseteq
 (S \odot T) \odot U$ gilt. Mit (b) folgt nun
 $$ (S \odot T) \odot U = U \odot (T \odot S) \subseteq (U \odot T) \odot S
			= S \odot (T \odot U).  $$
 Also ist in diesem Falle tats\"achlich $S \odot (T \odot U) =
 (S \odot T) \odot U$.
 \par
      $S$, $T$ und $U$ seien nun beliebige Punktmengen. Ist eine von ihnen
 leer, so ist (d) sicherlich erf\"ullt. Wir d\"urfen daher annehmen,
 dass keine von ihnen leer ist. Ist $P \in (S \odot T) \odot U$, so
 gibt es also Punkte $A$, $B$ und $C$ mit $A \in S$, $B \in S$, $B \in T$
 und $C \in U$ sowie $P \in (A \odot B) \odot C$. Mit dem bereits
 Bewiesenen und (a) folgt
 $$ P \in A \odot (B \odot C) \subseteq S \odot (T \odot U), $$
 so dass $( S \odot T) \odot U \subseteq S \odot (T \odot U)$ ist.
 Hieraus folgt unter Benutzung von (b), dass
 $$ S \odot (T \odot U) = (U \odot T) \odot S \subseteq U \odot (T \odot S)
                        = (S \odot T) \odot U $$
 ist. Damit ist (d) bewiesen.
 \par
       Genau dann ist $S \in L(\Sigma)$, wenn $S$ mit zwei verschiedenen
 Punkten stets auch ihre Verbindungsgerade enth\"alt. Dies bedeutet,
 dass $S$ genau dann in $L(\Sigma)$ liegt, wenn $S \odot S
 \subseteq S$ ist. Aus (a) folgt aber, dass stets $S \subseteq S \odot S$ gilt.
 Somit ist $S$ genau dann ein Unterraum von
 $\Sigma$, wenn $S \odot S = S$ ist. Dies beweist (e).
 \par
       Es seien $S$ und $T$ Unterr\"aume von $\Sigma$. Nach (c) gilt dann
 $S \odot T \subseteq S + T$. Mittels (b), (d) und (e) folgt
 $$ S \odot T \odot S \odot T = S \odot S \odot T \odot T = S \odot T, $$
 so dass $S \odot T$ nach (e) ein Unterraum ist. Wegen $S$, $T
 \subseteq S \odot T$ gilt nach 1.2 daher $S + T \subseteq S \odot
 T$, so dass auch (f) richtig ist. Damit ist alles bewiesen.
 \medskip
       Es sei $\Delta$ eine Menge und $\prec$ sei eine bin\"are Relation
 auf $\Delta$, die reflexiv und transitiv sei. Wir nennen $\Delta$
 bez\"uglich $\prec$ \emph{gerichtet},\index{gerichtete Menge}{} wenn es zu
 $\alpha$, $\beta \in \Delta$ stets ein $\gamma \in \Delta$ gibt mit
 $\alpha \prec \gamma$ und $\beta \prec \gamma$. Prominentestes, nicht triviales
 Beispiel f\"ur diese Situation ist die Menge $\Fin(M)$ der
 endlichen Teilmengen einer Menge $M$, die bez\"uglich der Inklusion
 gerichtet ist, da die Vereinigung zweier endlicher Mengen wieder
 endlich ist.
 \par
       Ist $\Sigma$ eine projektive Geometrie und ist $M \subseteq
 L(\Sigma)$ bez\"uglich der auf $L(\Sigma)$ definierten Relation
 $\leq$ gerichtet, so nennen wir $M$
 \emph{aufsteigendes System}\index{aufsteigendes System}{}
 von Unterr\"aumen von $\Sigma$. Es gilt nun
 \medskip\noindent
 {\bf 1.4. Satz.} {\it Ist $\Sigma$ eine projektive Geometrie und ist
 $M$ ein aufsteigendes System von Teilr\"aumen von $\Sigma$, so ist
 $$ \sum_{X \in M} X = \bigcup_{X \in M} X. $$}
 \par
       Beweis. Es seien $P$ und $Q$ zwei Punkte aus $\bigcup_{X \in M}X$.
 Es gibt dann $Y$, $Z \in M$ mit $P \leq Y$ und $Q \leq Z$. Weil $M$
 gerichtet ist, gibt es ein $U \in M$ mit $Y$, $Z \leq U$. Es folgt
 $$ P + Q \leq U \leq \bigcup_{X \in M} X. $$
 Folglich ist $\bigcup_{X \in M} X$ ein Teilraum von $\Sigma$. Da
 dieser Teilraum in $\sum_{X \in M} X$ liegt, folgt mit 1.2 die
 Gleichheit dieser beiden R\"aume.
 \medskip
       Es sei $M$ eine Teilmenge der Potenzmenge $P(\Pi)$ von $\Pi$,
 wobei $\Pi$ wieder die Punktmenge einer projektiven Geometrie sei.
 Ist $M$ endlich, so ist klar, da $\odot$ assoziativ und kommutativ
 ist, was wir unter
 $$ \bigodot_{X \in M} X. $$
 zu verstehen haben. Wegen 1.3 (a) gilt
 $$ \bigodot_{X \in M}X = \bigcup_{N \subseteq M} \bigodot_{X \in N}X.
                      $$
 Ist $M$ eine nicht notwendig endliche Teilmenge von $P(\Pi)$, so
 setzen wir den obigen Sachverhalt benutzend
 $$ \bigodot_{X \in M}X := \bigcup_{N \in \Fin(M)}\bigodot_{X \in N}X.
                      $$
 \par
       Damit sind wir nun in der Lage, die versprochene interne
 Be\-schrei\-bung der Summe von Unterr\"aumen zu geben.
 \medskip\noindent
 {\bf 1.5. Satz.} {\it Ist $\Sigma$ eine projektive Geometrie und ist
 $M$ eine Teilmenge von $L(\Sigma)$, so ist
 $$ \sum_{X \in M} X = \bigodot_{X \in M}X. $$}
 \par
       Beweis. Induktion nach $|M|$ zeigt unter Verwendung von
 1.3 (f), dass der Satz richtig ist, falls $M$ endlich ist. Daher
 ist in jedem Falle
 $$ \textstyle \bigl\{\bigodot_{X \in N} X\> \big | \> N \in \Fin(M)\bigr\} $$
 ein aufsteigendes System von Teilr\"aumen von $\Sigma$, so dass
 $\bigodot_{X \in M}M$ nach Satz 1.4 ein Teilraum ist. Mittels Satz
 1.2 folgt schlie\ss lich die Behauptung.
 \medskip
       Als n\"achstes beweisen wir die G\"ultigkeit eines eingeschr\"ankten
 Distributivgesetzes in $L(\Sigma)$. Dieses hat Dedekind\index{Dedekind, R.}{}
 als Erster betrachtet, so dass es heute seinen Namen
 tr\"agt (Dedekind 1897).
 \medskip\noindent
 {\bf 1.6. Dedekindsches Modulargesetz.}\index{dedekindsches Modulargesetz}{}
 {\it Es sei $\Sigma$ eine
 projektive Geometrie. Sind $S$, $T$, $U \in L(\Sigma)$ und ist $T \leq
 S$, so ist $S \cap (T + U) = T + (S \cap U)$.}
 \smallskip
       Beweis. Die Gleichung ist richtig, falls $T$ oder $U$ leer ist.
 Wir d\"urfen daher annehmen, dass dies nicht der Fall ist.
 \par
       Es ist ferner trivial, dass $T + (S \cap U) \leq S \cap (T + U)$
 ist. Es sei also $A$ ein Punkt von $S \cap (T + U)$. Ist $A \leq
 T$, so ist nichts zu beweisen. Wir d\"urfen also $A \not\leq T$
 annehmen. Nun ist $A \leq T + U$ und nach 1.3(f) ist $T + U = T
 \odot U$. Weil $T$ und $U$ nicht leer sind, gibt es dann zwei
 Punkte $B$ und $C$ mit $B \leq T$, $C \leq U$ und $A \leq B + C$.
 Wegen $A \not\leq T$ ist $A \neq B$, woraus folgt, dass $B + C = B +
 A$ ist. Somit ist
 $$ C \leq A + B \leq S + T = S, $$
 da ja $T \leq S$ gilt. Also ist $C \leq S \cap U$. Hieraus folgt
 schlie\ss lich
 $$ A \leq B + C \leq T + (S \cap U), $$
 was zu beweisen war.
 \medskip
       Eine weitere wichtige, wenn auch banal zu beweisende Eigenschaft
 einer projektiven Geometrie wird im n\"achsten Satz formuliert.
 \medskip\noindent
 {\bf 1.7. Satz.} {\it Es sei $M$ ein aufsteigendes System von Teilr\"aumen des
 projektiven Raumes $\Sigma$. Ist dann $Y \in L(\Sigma)$, so ist
 $Y \cap \sum_{X \in M}X = \sum_{X \in M}(Y \cap X)$.}
 \smallskip
       Beweis. Mittels Satz 1.4 folgt, da ja auch $\{Y \cap X \mid X \in M\}$
 ein aufsteigendes System von Teilr\"aumen ist,
 $$ Y \cap \sum_{X \in M}X = Y \cap \bigcup_{X \in M}X
                           = \bigcup_{X \in M}(Y \cap X)
                           = \sum_{X \in M} (Y \cap X). $$
 \par\noindent
 {\bf 1.8. Satz.} {\it Es sei $\Sigma$ eine projektive Geometrie. Sind $X$,
 $Y \in L(\Sigma)$ und gilt $X \cap Y = \emptyset$, so gibt es ein
 $Z \in L(\Sigma)$ mit $Y \leq Z$, $X \cap Z = \emptyset$ und $X + Z = \Pi$,
 wobei $\Pi$ wieder die Menge aller Punkte von $\Sigma$ bezeichne.}
 \smallskip
       Beweis. Es sei $N := \{U \mid U \in L(\Sigma), Y \leq U, X \cap U =
 \emptyset\}$. Dann ist $N$ nicht leer, da $Y$ zu $N$ geh\"ort. Es
 sei $M \subseteq N$ ein aufsteigendes System von Unterr\"aumen von
 $\Sigma$. Nach 1.4 ist dann $V := \bigcup_{U \in M}U$ ein Teilraum
 von $\Sigma$, der nat\"urlich  $Y$ enth\"alt. Weil $V$ die
 men\-gen\-the\-o\-re\-ti\-sche Vereinigung der $U$ aus $M$ ist, folgt, dass
 der Schnitt von $V$ mit $X$ leer ist. Somit gilt $V \in N$. Auf
 Grund des zornschen Lemmas gibt es daher einen maximalen Teilraum
 $Z$ in $N$. Es sei $P \in \Pi$. Wir m\"ussen zeigen, dass $P$ in
 $X + Z$ liegt. Dazu d\"urfen wir annehmen, dass $P$ weder zu $X$
 noch zu $Z$ geh\"ort. Dann ist insbesondere $Y \leq Z < Z + P$, so
 dass die Maximalit\"at von $Z$ impliziert, dass $X \cap (Z + P) \neq \emptyset$
 ist. Es gibt also einen Punkt $Q$ mit $Q \leq X \cap (Z + P)$. Weil $P$ nicht in
 $X$ liegt, ist $Q$ von $P$ verschieden. Ferner ist $Z + P = Z \odot P$. Es gibt
 daher einen Punkt $R$ mit $R \leq Z$ und $Q \leq R + P$. Wegen $Q \leq X$ und
 $X \cap Z = \emptyset$ ist $Q \neq R$. Also ist $R + P = Q + R$. Daher ist $P
 \leq Q + R \leq X + Z$. Dies zeigt, dass in der Tat $X + Z = \Pi$
 ist.
 \medskip\noindent
 {\bf 1.9. Korollar.} {\it Ist $X \in L(\Sigma)$, so gibt es ein $Z
 \in L(\Sigma)$ mit $X \cap Z = \emptyset$ und $X + Z = \Pi$.}
 \smallskip
       Beweis. Dies folgt mit $Y := \emptyset$ aus 1.8.
 \medskip
       Ist $\Sigma$ eine projektive Geometrie, sind $X$, $Y$, $Z \in L(\Sigma)$
 und gilt $X \cap Y = \emptyset$ und $X + Y = Z$, so nennen wir $Y$
 \emph{Komplement}\index{Komplement}{} von $X$ in $Z$. Diesen Sachverhalt
 beschreiben wir, wie in der Algebra \"ublich, durch $X \oplus Y = Z$.
 \medskip\noindent
 {\bf 1.10. Satz.} {\it Es sei $\Sigma$ eine projektive Geometrie und
 $\Pi$ ihre Punktmenge. Ist $\Pi = X \oplus Y$ und ist $Y$ das
 einzige Komplement von $X$, so ist $\Pi = X \cup Y$.}
 \smallskip
       Beweis. Es sei $P \in \Pi$ und $P \notin X$. Nach 1.8 gibt es ein
 Komplement $Z$ von $X$ mit $P \leq Z$. Weil $Y$ das einzige
 Komplement von $X$ ist, ist $Z = Y$. Damit ist alles bewiesen.
 \medskip\noindent
 {\bf 1.11. Korollar.} {\it Es sei $\Sigma$ eine projektive Geometrie.
 Tr\"agt jede Ge\-rad\-e von $\Sigma$ wenigstens drei Punkte, so sind
 $\emptyset$ und $\Pi$ die einzigen Unterr\"aume von $\Sigma$, die
 genau ein Komplement besitzen.}
 \smallskip
       Beweis. Es sei $\Pi = X \oplus Y$ und $X$ und $Y$ seien beide
 nicht leer. Ist dann $P$ ein Punkt auf $X$ und $Q$ ein Punkt auf
 $Y$, so enth\"alt die Gerade $P + Q$ noch einen dritten Punkt $R$.
 Wegen $X \cap Y = \emptyset$ ist $P + Q$ weder in $X$ noch in $Y$
 enthalten. Folglich liegt $R$ weder in $X$ noch in $Y$. Nach 1.10
 ist $Y$ daher nicht das einzige Komplement von $X$.
 \medskip
       Es sei $\leq$ eine Teilordnung auf der Menge $L$ und $\leq'$ eine
 solche auf der Menge $L'$. Die Bijektion $\sigma$ von $L$ auf $L'$
 hei\ss t \emph{Isomorphismus}\index{Isomorphismus}{} von $(L, \leq)$ auf
 $(L', \leq')$, wenn f\"ur alle $x$, $y \in L$ genau dann $x \leq y$ gilt, wenn
 $x^{\sigma} \leq' y^{\sigma}$ ist.
 \par
       Ist $V$ ein $K$-Vektorraum, so bezeichnen wir mit $L(V)$ die
 Menge seiner Unterr\"aume. Ferner setzen wir
 $$ \Sigma(V) := \bigl(UR_{1}(V), UR_{2}(V), \subseteq\bigr). $$
 Es gilt dann der folgende Satz, dessen Beweis dem Leser als
 \"U\-bungs\-auf\-ga\-be \"uberlassen bleibe.
 \medskip\noindent
 {\bf 1.12. Satz.} {\it Es sei $V$ ein $K$-Vektorraum. Setze $X^{\sigma}
 := \UR_{1}(X)$ f\"ur alle $X \in L(V)$. Dann ist $\sigma$ ein Isomorphismus
 von $(L(V), \subseteq)$ auf $(L(\Sigma(V)), \leq).$}
 \medskip
       In Kapitel II werden wir sehen, dass die projektiven Geometrien,
 deren Rang mindestens 4 ist, alle von dieser Art sind.

 \mysection{2. Projektive Verb\"ande}
 
\noindent
 Als N\"achstes geht es darum, die Menge der Unterr\"aume einer projektiven
 Geometrie verbandstheoretisch zu charakterisieren.
 \par
       Es sei $L$ eine Menge und $\leq$ sei eine
 Teilordnung,\index{Teilordnung}{} dh. eine reflexive, antisymmetrische und
 transitive Relation auf $L$. Gibt es in $L$ ein Element $\Pi$ mit $X \leq \Pi$
 f\"ur alle $X \in L$, so nennen wir $\Pi$ das \emph{gr\"o\ss te
 Element}\index{grosstes@gr\"o\ss tes Element}{} von $L$. Hat das Element $0 \in L$ die
 Eigenschaft, dass $0 \leq X$ gilt f\"ur alle $X \in L$, so nennen wir $0$ das
 \emph{kleinste Element}\index{kleinstes Element}{} von $L$. Aus der
 Antisymmetrie der Relation $\leq$ folgt, dass eine teilweise geordnete Menge
 h\"ochstens ein gr\"o\ss tes und h\"ochstens ein kleinstes Element besitzt.
 \par
       Ist $M  \subseteq L$ und ist $R$ eine Element von $M$, so dass $X \in M$
 und $R \leq X$ impliziert, dass $R = X$ ist, so hei\ss t
 $R$ \emph{maximales Element}\index{maximales Element}{} von $M$. Entsprechend
 werden \emph{minimale Elemente}\index{minimales Element}{} definiert. Ist
 $S \in L$ und ist $X \leq S$ f\"ur alle $X \in M$, so hei\ss t $S$
 \emph{obere Schranke}\index{obere Schranke}{} von $M$. \emph{Untere
 Schranken}\index{untere Schranke}{} werden entsprechend definiert. Hat
 schlie\ss lich die Menge der oberen Schranken von $M$ ein kleinstes Element, so
 hei\ss t dieses Element \emph{obere Grenze}\index{obere Grenze}{} von $M$.
 Besitzt
 $M$ eine obere Grenze, so bezeichnen wir diese mit $\sum_{X \in M}X$. Besitzt
 die Menge aller unteren Schranken von $M$ ein gr\"o\ss tes Element, so hei\ss t
 dieses Element \emph{untere Grenze}\index{untere Grenze}{} von $M$. Wir
 bezeichnen sie mit $\bigcap_{X \in M}X$. Obere und untere Grenzen sind, falls
 sie existieren, eindeutig bestimmt.
 \par
       Ist $L$ eine teilweise geordnete Menge und besitzt jede nicht leere
 endliche Teilmenge von $L$ eine untere und eine obere Grenze, so hei\ss t $L$
 \emph{Verband}.\index{Verband}{} Hat jede Teilmenge von $L$ eine obere und eine
 untere Grenze, so hei\ss t $L$ \emph{vollst\"andiger
 Verband}.\index{vollst\"andiger Verband}{} Ein vollst\"andiger Verband besitzt
 stets ein gr\"o\ss tes und ein kleinstes Element, n\"amlich $\Pi
 := \sum_{X \in L}X$ und $0 := \bigcap_{X \in L}X$.
 \par
       Bei der Definition des vollst\"andigen Verbandes haben wir des
 Gu\-ten zuviel getan, wie der folgende Satz zeigt.
 \medskip\noindent
 {\bf 2.1. Satz.} {\it Ist $L$ eine teilweise geordnete Menge, so sind
 die folgenden Aussagen \"aquivalent:
 \item{(a)} $L$ ist ein vollst\"andiger Verband.
 \item{(b)} Jede Teilmenge von $L$ hat eine untere Grenze.
 \item{(c)} Jede Teilmenge von $L$ hat eine obere Grenze.}
 \smallskip
       Beweis. (b) und (c) folgen nat\"urlich aus (a).
 \par
       Es gelte (b) und es sei $M \subseteq L$. Mit $M^{*}$ bezeichnen
 wir die Menge der oberen Schranken von $M$ und mit $M^{**}$
 bezeichnen wir die Menge der unteren Schranken von $M^{*}$. Dann
 ist nat\"urlich $M \subseteq M^{**}$. Nach Voraussetzung besitzt
 $M^{*}$ eine untere Grenze $U$. \emph{Per definitionem} ist $U$
 das gr\"o\ss te Element von $M^{**}$, so dass insbesondere $X \leq
 U$ gilt f\"ur alle $X \in M$. Folglich ist $U \in M^{*}$. Wegen $U
 \leq Y$ f\"ur alle $Y \in M^{*}$ ist $U$ das kleinste Element in
 $M^{*}$, so dass $U$ eine obere Grenze von $M$ ist. Damit ist
 gezeigt, dass $L$ ein vollst\"andiger Verband ist.
 \par
       Ganz entsprechend zeigt man, dass (a) auch aus (c) folgt.
 \medskip
       Es sei $L$ ein Verband mit kleinstem Element $0$. Genau die minimalen
 Elemente von $L - \{0\}$ hei\ss en \emph{Atome}\index{Atom}{} von $L$.
 Solche braucht es nicht zu geben. Gibt es jedoch zu jedem $X \in L
 - \{0\}$ ein Atom $A$ mit $A \leq X$, so hei\ss t der Verband $L$
 \emph{atomar}.\index{atomar}{}
 \par
       Es sei $L$ ein Verband mit kleinstem Element $0$ und gr\"o\ss tem Element
 $\Pi$. Sind $A$, $B \in L$ und ist $A + B = \Pi$ sowie $A \cap B = 0$, so
 hei\ss t
 $B$ \emph{Komplement}\index{Komplement}{} von $A$. Diesen Sachverhalt bezeichnen
 wir wie schon zuvor mit $\Pi = A \oplus B$. Hat jedes Element von $L$ ein
 Komplement, so hei\ss t $L$ \emph{komplement\"arer}\index{komplement\"ar}{}
 Verband. Der Verband $L$ hei\ss t
 \emph{ir\-re\-du\-zi\-bel},\index{irreduzibler Verband}{} falls
 $0$ und $\Pi$ die einzigen Elemente von $L$ sind, die genau ein Komplement haben.
 \par
       Der Verband $L$ hei\ss t \emph{modular},\index{modular}{} falls aus
 $A$, $B$, $C \in L$ und $B \leq A$ folgt, dass $A \cap (B + C) = B + (A \cap C)$
 folgt.
 \par
       Es sei $L$ ein Verband und $M$ sei eine Teilmenge von $L$. Ist $M$
 bez\"uglich der auf $L$ gegebenen Teilordnung gerichtet, so nennen
 wir $M$ ein {\it aufsteigendes System\/}.\index{aufsteigendes System}{} Der
 Verband $L$ hei\ss t nach \emph{oben stetig},\index{nach oben stetig}{} falls
 alle aufsteigenden Systeme von $L$ eine obere Grenze haben und falls f\"ur alle
 aufsteigenden Systeme $M$ und alle $Y \in L$ gilt, dass
 $Y \cap \sum_{X \in M}X = \sum_{X \in M}(Y \cap X)$ ist.
 \par
       Einen vollst\"andigen, atomaren, modularen, komplement\"aren und nach oben
 stetigen Verband nennen wir \emph{projektiv}.\index{projektiver Verband}{} Die
 Ergebnisse des ersten Abschnitts lassen sich nun zusammenfassen zu
 \medskip\noindent
 {\bf 2.2. Satz.} {\it Ist $\Sigma$ eine projektive Geometrie, so ist
 $L(\Sigma)$ bez\"uglich der Inklusion als Teilordnung ein projektiver
 Verband.}\index{projektive Geometrie}{}
 \medskip
       Das Hauptziel dieses Abschnittes ist nun zu zeigen, dass jeder projektive
 Verband zum Unterraumverband einer projektiven Geometrie isomorph ist.
 \par
       Es sei $L$ ein Verband und $C$, $D \in L$. Ferner sei $C \leq D$. Wir
 setzen
 $$ D/C := \{X \mid X \in L,\ C \leq X \leq D\} $$
 und nennen $D/C$ den \emph{Quotienten}\index{Quotient}{} von $D$ nach $C$.
 Offenbar ist $D/C$ bez\"uglich der in $L$ definierten Teilordnung ein
 Verband, m. a. W., $D/C$ ist ein \emph{Teilverband}\index{Teilverband}{} von
 $L$.
 \medskip\noindent
 {\bf 2.3. Transformationsregel.} {\it Es sei $L$ ein modularer Verband
 und $A$ und $B$ seien Elemente von $L$. Setze
 $$ X^\sigma := X \cap B\ \hbox{\it f\"ur}\ X \in (A + B)/A $$
 und
 $$ Y^\tau := Y + A\ \hbox{\it f\"ur}\ Y \in B /(A \cap B). $$
 Dann ist $\sigma$ ein Isomorphismus von
 $(A + B)/A$ auf $B/(A \cap B)$ und $\tau$ ist ein Isomorphismus
 von $B/(A \cap B)$ auf $(A + B)/A$. \"Uberdies sind $\sigma$ und
 $\tau$ invers zueinander.}\index{Transformationsregel}{}
 \smallskip
       Beweis. Aus $A \leq X$ folgt $A \cap B \leq X \cap B \leq B$, so
 dass $\sigma$ eine Abbildung von $(A + B)/A$ in $B/(A \cap B)$
 ist. Ist $A \cap B \leq Y \leq B$, so folgt $A \leq Y + A \leq A +
 B$, so dass $\tau$ eine Abbildung von $B/(A \cap B)$ in $(A + B)/A$ ist.
 \medskip
       Ist $A \leq X \leq A + B$, so folgt auf Grund der Modularit\"at von
 $L$, dass
 $$ X^{\sigma \tau} = (X \cap B) + A = X \cap (A + B) = X$$
 ist. Folglich ist $\sigma \tau$ die Identit\"at auf $(A + B)/A$. Ist
 $A \cap B \leq Y \leq B$, so folgt wiederum wegen der
 Modularit\"at von $L$, dass
 $$ Y^{\tau \sigma} = (Y + A) \cap B = Y + (A \cap B) = Y $$
 ist. Also ist $\tau \sigma$ die Identit\"at auf $B/(A \cap B)$. Folglich sind
 $\sigma$ und $\tau$ zueinander inverse Bijektionen.
 \par
       Schlie\ss lich folgt aus $A \leq X_{1} \leq X_{2} \leq A + B$,
 dass $X^{\sigma}_{1} = X_{1} \cap B \leq X_{2} \cap B =
 X^{\sigma}_{2}$ ist, und aus $A \cap B \leq Y_{1} \leq Y_{2} \leq
 B$ folgt, dass auch $Y^{\tau}_{1} = Y_{1} + a \leq Y_{2} + A =
 Y^{\tau}_{2}$ gilt. Damit ist alles bewiesen.
 \medskip\noindent
 {\bf 2.4. Korollar.} {\it Ist $P$ ein Atom des modularen Verbandes $L$, dessen
 kleinstes Element wieder mit $0$ bezeichnet sei, und ist $A \in L$, so ist
 entweder $P \leq A$ oder $(A + P)/A$ und $P/0$ sind isomorph.}
 \smallskip
       Beweis. Nach 2.3 sind $(A + P)/A$ und $P/(A \cap P)$ isomorph. Nun ist
 $0 \leq A \cap P \leq P$ und daher, da $P$ ein Atom ist,
 entweder $A \cap P = 0$ oder $A \cap P = P$, dh., $P \leq A$.
 \medskip\noindent
 {\bf 2.5. Korollar.} {\it Es sei $L$ ein modularer Verband mit kleinstem Element
 $0$. Sind $P$, $Q$ und $R$ Atome von $L$, ist $R \neq P$ und $R \leq P + Q$, so
 ist $P + Q = P + R$.}
 \smallskip
       Beweis. Es ist $P \neq Q$, da sonst $P = Q = R$ w\"are. Somit ist
 $(P + Q)/P$ nach 2.4 zu $Q/0$ isomorph. Da der Quotient $Q/0$ nur aus den
 Elementen $0$ und $Q$ besteht, enth\"alt der Quotient $(P + Q)/P$ nur die beiden
 Elemente $P$ und $P + Q$. Hieraus folgt, dass $P + R = P + Q$ ist, da ja
 $P < P + R \leq P + Q$ ist.
 \medskip\noindent
 {\bf 2.6. Satz.} {\it Es sei $L$ ein modularer Verband mit kleinstem
 Element $0$. Ferner seien $P$ und $Q$ zwei verschiedene Atome von
 $L$. Schlie\ss lich sei $A \in L$ und $Q \leq P + A$. Ist $Q \leq
 A$, so setzen wir $R := Q$. Ist $Q \not\leq A$, so setzen wir $R
 := A \cap (P + Q)$. Dann ist $R$ ein Atom und es gilt $R \leq A$
 und $Q \leq P + R$.}
 \smallskip
       Beweis. Ist $Q \leq A$, so ist nichts zu beweisen. Es sei also $Q \not\leq
 A$. Dann ist auch $P \not\leq A$, da andernfalls $Q \leq P + A = A$ w\"are. Es
 folgt $P \cap A = 0$ und weiter
 $$ P/0 = P/(P \cap A) \cong (A + P)/A = (A + P + Q)/A \cong (P + Q)/R. $$
 Somit ist $(P + Q)/R = \{R, P + Q\}$ und $R \neq P + Q$. Andererseits ist
 $$ (R + P)/R \cong P/(P \cap R) = P/0, $$
 da ja $P \cap R \leq P \cap A = 0$ ist. Also ist $(R + P) = \{R, P + R\}$
 und $R \neq P + R$. Nun ist aber $R + P \leq P + Q + P = P + Q$ und folglich
 $(R + P)/R \subseteq (P + Q)/R$. Hieraus folgt zusammen mit dem bereits
 Bewiesenen, dass $\{R, P + R\} = \{R, P + Q\}$ und damit dass $P + Q = P + R$
 ist. Damit ist zun\"achst gezeigt, dass $Q \leq P + R$ ist.
 \par
       Es bleibt zu zeigen, dass $R$ ein Atom ist. Dies folgt nun aus
 $$ R/0 = R/(P \cap R) \cong (R + P)/P = (Q + P)/P \cong Q/(P \cap Q) = Q/0, $$
 wobei wir erst hier die Voraussetzung $P \neq Q$ benutzt haben.
 \medskip
       Der n\"achste Satz ist der verbandstheoretische Hintergrund des
 Veblen-Young Axioms.\index{Veblen-Young Axiom}{}
 \medskip\noindent
 {\bf 2.7. Satz.} {\it Es sei $L$ ein modularer Verband mit $0$ und $P$,
 $Q$ und $R$ seien drei verschiedene Atome von $L$ mit $R \not\leq P + Q$. Sind
 dann $S$ und $T$ zwei verschiedene Atome mit $S \leq P
 + Q$ und $T \leq Q + R$, so gibt es ein Atom $U$ mit $U \leq S +
 T$ und $U \leq R + P$.}
 \smallskip
       Beweis. Ist $S = P$, so tut's $U := P$. Es sei also $S \neq P$. Nun ist
 $S \leq P + Q$. Nach 2.5 ist daher $P + Q = P + S$. Also ist $Q \leq P + S$.
 Hieraus folgt, dass $T \leq Q + R \leq S + P + R$ ist. Aus 2.6. folgt
 schlie\ss lich, dass es ein Atom $U$ gibt mit $U \leq P + R$ und $T \leq S + U$.
 Weil $T \neq S$ ist, folgt wiederum aus 2.5, dass $S + U = S + T$ ist. Also ist
 auch $U \leq S + T$. Damit ist alles bewiesen.
 \medskip
       F\"ur sp\"atere Verwendung beweisen wir den n\"achsten Satz, der auf
 der Menge der Atome eines modularen Verbandes eine
 \"A\-qui\-va\-lenz\-re\-la\-tion beschreibt. Bei seinem Beweis wird der gerade
 bewiesene Satz benutzt.
 \medskip\noindent
 {\bf 2.8. Satz.} {\it Es sei $L$ ein modularer Verband mit kleinstem Element
 $0$. Sind $P$ und $Q$ Atome von $L$, so setzen wir $P \sim Q$ genau dann, wenn
 entweder $P = Q$ ist oder wenn $P + Q$ wenigstens drei Atome umfasst. Dann ist
 $\sim$ eine \"Aquivalenzrelation auf der Menge der Atome von $L$.}
 \smallskip
       Beweis. Die Relation $\sim$ ist offenkundig reflexiv und
 symmetrisch. Um die Transitivit\"at zu beweisen, seien $P$, $Q$ und
 $R$ drei Atome von $L$ und es gelte $P \sim Q$ und $Q \sim R$.
 Sind zwei der drei Atome gleich, so ist nichts zu beweisen. Wir
 d\"urfen daher annehmen, dass sie paarweise verschieden sind. Ist
 $R \leq P + Q$, so folgt $P + R = P + Q$. Hieraus folgt, dass $P +
 R$ mindestens drei Atome enth\"alt, so dass $P \sim R$ gilt. Es sei
 schlie\ss lich $R \not\leq P + Q$. Wegen $P \sim Q$ und $Q \sim R$
 gibt es ein von $P$ und $Q$ verschiedenes Atom $A$ mit $A \leq P +
 Q$ und ein von $Q$ und $R$ verschiedenes Atom  $B$ mit $B \leq Q +
 R$. Nach 2.7 gibt es ein Atom $C$ mit $C \leq A + B$ und $C \leq R
 + P$. Es bleibt zu zeigen, dass $C$ von $P$ und $R$ verschieden
 ist. W\"are $C = P$, so folgte $C + A = P + A \leq P + Q$ und
 damit $C + A = P + Q$, da ja $P \neq A$ ist. Aus $C \leq A + B$
 folgte weiter $C + A = A + B$ und daher
 $$ B \leq (P + Q) \cap (Q + R) = Q. $$
 Dieser Widerspruch zeigt, dass $C \neq P$ ist.
 Ganz entsprechend zeigt man, dass auch $C \neq R$ gilt. Damit ist
 alles bewiesen.
 \medskip
       Es sei $L$ ein Verband mit kleinstem Element $0$. Der Verband $L$
 hei\ss t \emph{relativ atomar},\index{relativ atomar}{} wenn es zu $A$,
 $B \in L$ mit $B < A$ stets ein Atom $P$ gibt mit $B < B + P \leq A$.
 \medskip\noindent
 {\bf 2.9. Satz.} {\it Es sei $L$ ein Verband mit kleinstem Element $0$. Ist $L$
 relativ atomar, so ist jedes Element von $L$ die obere
 Grenze\index{obere Grenze}{} der in ihm enthaltenen Atome.}
 \smallskip
       Beweis. Es sei $A$ ein von $0$ verschiedenes Element von $L$ und $S$ sei
 die Menge der in $A$ enthaltenen Atome. Dann ist $A$ eine obere Schranke von
 $S$. Es sei $B$ eine weitere obere Schranke von $S$. Dann ist auch $B \cap A$
 eine obere Schranke von $S$. W\"are $B \cap A < A$, so g\"abe es ein Atom $P$ mit
 $B \cap A < (B \cap A) + P \leq A$. Es folgte $P \in S$ und damit der
 Widerspruch $P \leq B \cap A$. Also ist $A = A \cap B \leq B$. Somit ist $A$ die
 kleinste obere Schranke von $S$, was zu beweisen war.
 \medskip
       Der Verband $L$ hei\ss t
 \emph{relativ komplement\"ar},\index{relativ komplement\"ar}{} wenn jeder
 Quotient zweier Elemente von $L$ ein komplement\"arer Verband ist, dh., wenn es
 zu $U$, $V$, $W \in L$ mit $U \leq V \leq W$ stets ein $X \in L$ gibt mit
 $V \cap X = U$ und $V + X = W$.
 \medskip\noindent
 {\bf 2.10. Satz.} {\it Es sei $L$ ein komplement\"arer und modularer Verband mit
 gr\"o\ss tem Element $\Pi$ und kleinstem Element $0$. Ferner seien $U$, $V$ und
 $W$ drei Elemente von $L$ mit $U \leq V \leq W$ und $Y$ sei ein Komplement von
 $V$ in $\Pi$. Setze $X := U + (Y \cap W)$. Dann ist $X \cap V = U$ und
 $X + V = W$. Mit anderen Worten, jeder komplement\"are, modulare Verband ist
 relativ komplement\"ar.}
 \smallskip
       Beweis. Es ist $V + X = V + U + (Y \cap W) = V + (Y \cap W)$. Wegen der
 Modularit\"at von $L$ ist daher
 $$ V + X = W \cap (V + Y) = W. $$
 Andererseits ist, wiederum auf Grund der Modularit\"at von $L$,
 $$\eqalign{
   V \cap X &= V \cap (U + (Y \cap W)) \cr
	    &= U + \bigl(V \cap Y \cap W\bigr) = U + (0 \cap W) = U. \cr} $$
 Damit ist alles bewiesen.
 \medskip\noindent
 {\bf 2.11. Satz.} {\it In einem atomaren, komplement\"aren und modularen Verband
 ist jedes Element die obere Grenze der in ihm enthaltenen Atome.}
 \smallskip
       Beweis. Nach 2.9 gen\"ugt es zu zeigen, dass jeder solche Verband relativ
 atomar ist. Dies folgt aber unmittelbar aus der in Satz 2.10 etablierten
 relativen Komplementarit\"at eines komplement\"aren mo\-du\-la\-ren Verbandes.
 \medskip
       An dieser Stelle ist eine methodische Bemerkung angebracht. Ist $\Sigma$
 eine projektive Geometrie, so ist $L(\Sigma)$ auf Grund seiner Definition
 nat\"urlich relativ atomar. Man ben\"otigt also nicht die relative
 Komplementarit\"at dieses Verbandes, um dies festzustellen. Dies w\"are mit
 Kanonen nach Spatzen geschossen, da man die Komplementarit\"at von $L(\Sigma)$
 nur \"uber das zornsche Lemma erh\"alt.
 \medskip\noindent
 {\bf 2.12. Satz von der endlichen Abh\"angigkeit.} {\it Es sei $L$ ein
 relativ atomarer, vollst\"andiger Verband. Genau dann ist $L$ nach
 oben stetig, wenn gilt: Ist $P$ ein Atom und $S$ eine Menge von
 Atomen von $L$, ist ferner $P \leq \sum_{Q\in S}Q$, so
 gibt es endliche viele Punkte $Q_1$, \dots, $Q_{t} \; in \; S$ mit
 $P \leq \sum_{i := 1}^{t} Q_{i}$.}\index{endliche Abh\"angigkeit}{}
 \smallskip
       Beweis. Es sei $L$ nach oben stetig. Ferner sei $P$ ein Atom und $S$ eine
 Menge von Atomen von $L$ und es gelte $P \leq \sum_{Q \in S}Q$. Das
 System $M := \{\sum_{Q \in \Phi} Q \mid \Phi \in \Fin (S)\}$ ist
 aufsteigend. Ferner ist $\sum_{X \in M} X = \sum_{Q \in S} Q$.
 Folglich ist
 $$ P = P \cap \sum_{Q \in S} Q = P \cap \sum_{X \in M} X
	    = \sum_{X \in M} (P \cap X). $$
 Hieraus folgt die Existenz eines $X \in M$ mit $P \cap X \neq 0$. Weil $P$ ein
 Atom ist, folgt weiter $P \leq X$. Weil $X$ von einer endlichen Teilmenge von
 $S$ erzeugt wird, gibt es also endlich viele $Q_1$, \dots, $Q_t \in S$ mit
 $P \leq \sum_{i := 1}^{t} Q_i$.
 \par
       Umgekehrt gelte: Wann immer $P$ ein Atom und $S$ eine Menge von
 Atomen von $L$ ist, so dass $P \leq \sum_{Q \in S} Q$
 ist, so gibt es endlich viele $Q_1$, \dots, $Q_t \in S$ mit $P \leq
 \sum_{i := 1}^{t} Q_i$. Unter dieser Annahme m\"ussen wir nun
 zeigen, dass $L$ nach oben stetig ist.
 \par
       Ist $X \in L$, so bezeichnen wir mit $S(X)$ die Menge der in $X$
 ent\-hal\-te\-nen Atome. Nach 2.9 ist dann $X = \sum_{Q \in
 S(X)} Q$. Ferner gilt $S(X) \cap S(Y) = S(X \cap Y)$.
 \par
       Es sei nun $M$ ein aufsteigendes System von $L$. Ferner sei $B \in L$.
 Dann ist
 $$ \sum_{X \in M} (B \cap X) \leq B \cap \sum_{X \in M} X. $$
 Es sei $P$ ein Atom von $B \cap \sum_{X \in M} X$. Dann
 ist insbesondere
 $$ P \leq \sum_{X \in M} \sum_{Q \in S(X)} Q.$$
 Es gibt also $Q_1$, \dots, $Q_t \in \bigcup_{X \in M} S(X)$ mit $P \leq
 \sum_{i := 1}^{t} Q_i$. Weil $M$ gerichtet ist, gibt es
 ein $Y \in M$ mit $Q_1$, \dots, $Q_t \in Y$. Es folgt
 $$ P \leq B \cap Y \leq B \cap \sum_{X \in M} X. $$
 Weil in einem relativ atomaren Verband jedes Element die obere
 Grenze seiner Atome ist, folgt schlie\ss lich
 $$ B \cap \sum_{X \in M} X \leq \sum_{X \in M} (B \cap X). $$
 Damit ist Satz 2.12 bewiesen. Er wird uns im n\"achsten Abschnitt
 gute Dienste leisten.
 \medskip
       Nun zeigen wir, dass jeder projektive Verband\index{projektiver Verband}{}
 zum Un\-ter\-raum\-ver\-band einer projektiven
 Geometrie\index{projektive Geometrie}{} isomorph ist.
 \medskip\noindent
 {\bf 2.13. Satz.} {\it Es sei $L$ ein projektiver Verband. Wir definieren eine
 In\-zi\-denz\-stru\-ktur $\Sigma$ wie folgt. Punkte von $\Sigma$ sind die Atome von $L$.
 Geraden von $\Sigma$ sind die Elemente der Form $P + Q$ von $L$, wobei $P$ und
 $Q$ zwei verschiedene Atome von $L$ sind. Dann ist $\Sigma$ eine projektive
 Geometrie und die Verb\"ande $L$ und $L(\Sigma)$ sind isomorph.}
 \smallskip
       Beweis. Durch zwei verschiedene Punkte von $\Sigma$ geht stets eine
 Ge\-ra\-de, und aus 2.5 folgt, dass es auch nur eine Gerade durch zwei
 verschiedene
 Punkte gibt. Also gilt (P1). Nach 2.7 gilt auch (P2). Schlie\ss\-lich gilt (P3)
 auf Grund der Definition einer Geraden.
 \par
       Um die Isomorphie von $L$ und $L(\Sigma)$ zu beweisen, probiere
 man das N\"achst\-lie\-gen\-de. Funktionierte dies nicht, w\"are man in
 Schwie\-rig\-kei\-ten. F\"ur $X \in L$ setzen wir also $X^\sigma$ gleich
 der Menge der in $X$ liegenden Atome und f\"ur $Y \in L(\Sigma)$
 setzen wir $Y^\tau := \sum_{Q \in Y} Q$, wobei $\sum$ die
 obere Grenze in $L$ bezeichne. Dann ist $\sigma$ eine Abbildung
 von $L$ in $L(\Sigma)$ und $\tau$ eine Abbildung von
 $L(\Sigma)$ in $L$. Ist nun $X \in L$, so ist $X$ die obere
 Grenze der in $X$ enthaltenen Atome. Daher ist $X^{\sigma \tau} =
 X$, so dass $\sigma \tau = id_L$ ist.
 \par
       Es sei $Y \in L(\Sigma)$ und $P$ ein Atom von $L$ mit $P \leq Y^\tau$.
 Wir zeigen, dass $P$ in $Y$ liegt. Nach dem Satz von der endlichen
 Abh\"angigkeit gibt es $Q_1$, \dots, $Q_t \in Y$ mit
 $P \leq \sum_{i := 1}^{t} Q_i$, wobei $\sum$ sich wieder auf $L$ bezieht. Ist
 $t = 1$, so ist $P = Q_1 \in Y$. Es sei also $t > 1$. Dann ist
 $P \leq Q_1 + \sum_{i := 2}^t Q_i$. Nach Satz 2.6 gibt es einen Punkt
 $R$ mit $P \leq Q_1 + R$ und $R \leq \sum_{i:= 2}^{t}Q_i$. Nach
 Induktionsannahme ist $R \in Y$. Dann ist aber $Q_1 + R$ eine Gerade von $Y$, so
 dass auch $P \in Y$ gilt. Somit ist $Y^{\tau\sigma} = Y$, so dass
 $\tau \sigma = id_{L(\Sigma)}$ ist. Folglich sind $\sigma$ und $\tau$
 zueinander inverse Bijektionen. Dass beide inklusionstreu sind, ist banal.
 \medskip
       Die Frage, wann ein projektiver Verband irreduzibel ist, ist nun leicht zu
 beantworten.
 \medskip\noindent
 {\bf 2.14. Satz.} {\it Es sei $L$ ein projektiver Verband. Genau dann
 ist $L$ irreduzibel, wenn auf jeder Geraden von $L$ wenigstens
 drei Punkte liegen. Ist $L$ irreduzibel, sind $X$, $Y \in L$ und
 gilt $Y \leq X$, so ist auch $X/Y$ irreduzibel.}
 \smallskip
       Beweis. Nach 2.13 d\"urfen wir annehmen, dass $L = L(\Sigma)$
 ist, wobei $\Sigma$ wie in 2.13 definiert sei. Liegen nun auf
 jeder Geraden von $\Sigma$ drei Punkte, so ist $L$ nach 1.11 irreduzibel.
 \par
       Es sei $L$ irreduzibel. Ferner sei $U$ eine \"Aquivalenzklasse der
 in 2.8 erkl\"arten \"Aquivalenzrelation $\sim$. Dann ist $U \in
 L$. Es sei $V$ die Menge der nicht in $U$ liegenden Punkte von
 $\sigma$. Dann ist auch $V \in L$, da $V$ Vereinigung von
 \"Aquivalenzklassen von $\sim$ ist. Dann ist $U \cup V = \Pi = U
 \oplus V$. Hieraus folgt aber, dass $V$ das einzige Komplement von
 $U$ ist. Da $U$ als \"Aquivalenzklasse nicht leer ist, ist $U =
 \Pi$, da $L$ als irreduzibel vorausgesetzt war. Somit ist $V = 0$
 und alle Geraden von $\Sigma$ sind in $U$ enthalten und tragen
 daher alle mindestens drei Punkte.
 \par
       Es sei $L$ irreduzibel, es seien $X$, $Y \in L$ und es gelte $Y \leq X$.
 Es
 sei $G \in X/Y$ eine Gerade von $X/Y$. Weil $L$ relativ komplement\"ar ist, gibt
 es ein $H$ mit $G = Y \oplus H$. Mittels der Transformationsregel folgt, dass
 $H/0$ zu $G/Y$ isomorph ist.  Also ist $H$ eine Gerade, tr\"agt daher drei
 verschiedene Punkte $P_1$, $P_2$ und $P_3$. Setze $Q_{i} := P_i + Y$. Dann sind
 $Q_1$, $Q_2$, $Q_3$ drei verschieden Punkte von $X/Y$, die alle auf $G$
 liegen. Folglich ist $X/Y$ irreduzibel.
 \medskip
       Es sei $(L_i \mid i \in I)$ eine Familie von Verb\"anden. Wir versehen
 das cartesische Produkt $C:= \cart_{i \in I} L_i$ der $L_i$ mit einer bin\"aren
 Relation $\leq$, indem wir f\"ur alle $F$, $G \in C$ genau dann $F \leq G$
 setzen,
 wenn $F_i \leq G_i$ f\"ur alle $i \in I$ gilt. Es ist schnell verifiziert, dass
 $(C, \leq)$ ebenfalls ein Verband ist. Sind alle $L_i$ projektiv, so ist auch
 $C$ projektiv.
 \medskip\noindent
 {\bf 2.15. Satz.} {\it Es sei $L$ ein projektiver Verband. Auf Grund von 2.13
 d\"urfen wir jedes Element von $L$ mit der Menge der auf ihm liegenden Punkte
 identifizieren. Es sei $\Pi$ die Menge aller Punkte von $L$ und $0$ bezeichne
 die leere Menge. Es sei weiter $\Pi/{\sim}$ die Menge der \"Aquivalenzklassen
 der
 in 2.8 definierten \"Aquivalenzrelation $\sim$. Ist dann $U \in \Pi/\sim$, so
 ist $U \in L$ und der Quotient $U/0$ ist mit der von $L$ ererbten Teilordnung
 ein irreduzibler projektiver Verband. Ist nun $X \in L$, so definieren wir
 $$ \sigma(X) \in \cart_{U \in \Pi/{\sim}} U/0 $$
 durch $\sigma(X)_U := X \cap U$. Dann ist $\sigma$ ein Isomorphismus von $L$ auf
 $\cart_{U \in \Pi/\sim} U/0$.}
 \smallskip
       Beweis. Nat\"urlich ist $\sigma$ eine Abbildung von $L$ in das fragliche
 cartesische Produkt, welches wir abk\"urzend mit $C$ bezeichnen. Es sei
 $\sigma(X) = \sigma(Y)$. Dann ist
 $$\eqalign{
   X &= X \cap \bigcup_{U \in \Pi/\sim} U = \bigcup_{U\in \Pi/\sim} X \cap U \cr
     &= \bigcup_{U \in \Pi/\sim} Y \cap U = Y \cap \bigcup_{U \in \Pi/\sim}
                U = Y. \cr} $$
 Dies zeigt, dass $\sigma$ injektiv ist. Um zu zeigen, dass $\sigma$ auch
 surjektiv ist, sei $F \in C$. Wir setzen $X := \cup_{U \in \Pi / \sim} F_U$.
 Sind $P$ und $Q$ zwei verschiedene Punkte von $X$, so liegt auch $P + Q$ ganz in
 $X$: Dies ist sicherlich richtig, wenn $P$ und $Q$ in ein und demselben $F_U$
 liegen. Liegen sie aber in verschiedenen \"Aquivalenzklassen, so liegen auf
 $P + Q$ nur die beiden Punkte $P$ und $Q$. Somit gilt $X \in L$. Da nun
 $\sigma(X) = F_U$ ist, ist $\sigma$ auch surjektiv.
 \par
        Die Inklusionstreue von $\sigma$ ist banal.
 \medskip
        Dieser Satz zeigt, dass man alle projektiven Verb\"ande kennt,
 wenn man nur die irreduziblen unter ihnen kennt.

 \mysection{3. Der Basissatz}
 
\noindent
 Die Vorgehensweise in diesem Abschnitt mag dem ein oder anderen Leser
 um\-st\"and\-lich erscheinen. Sie erkl\"art sich daraus, dass ich immer versuche,
 ohne das Auswahlaxiom auszukommen.
 \par
       Es sei $L$ ein projektiver Verband, dessen Punktmenge wir wieder
 mit $\Pi$ bezeichnen. Wir nennen $X \in L$ \emph{endlich
 erzeugt},\index{endlich erzeugt}{}
 wenn es ein $\Phi \in \Fin(\Pi)$ gibt mit $X = \sum_{P
 \in \Phi}P$. Ist $X$ endlich erzeugt, so setzen wir
 $$ \textstyle \Rg_L(X) := \min \bigl\{|\Phi| \>\big|\> \Phi \in \Fin(\Pi), X =
                        \sum_{P \in \Phi}P\bigr\}$$
 und lesen $\Rg_L(X)$ als \emph{Rang}\index{Rang}{} von $X$. Ist $\Pi$ endlich
 erzeugt, so
 nennen wir $\Rg_L(\Pi)$ auch Rang von $L$. Der Vollst\"andigkeit halber sei
 erw\"ahnt, dass im Falle der Endlichkeit von $\Rg_L(X)$ die Zahl
 $\Rg_L(X) - 1$ die \emph{Dimension}\index{Dimension}{} von $X$ ist.
 \medskip\noindent
 {\bf 3.1. Satz.} {\it Es sei $L$ ein projektiver Verband. Ferner seien $X$,
 $Y \in L$ und es gelte $Y \leq X$. Ist $X$ endlich erzeugt, so ist auch $Y$
 endlich erzeugt und es gilt $\Rg_L(Y) \leq \Rg_L(X)$. \"Uberdies gilt in
 diesem Falle genau dann $\Rg_L(Y) = \Rg_L(X)$, wenn $Y = X$ ist.}
 \smallskip
       Beweis. Es sei $n := \Rg_L(X)$ und $\Phi$ sei eine Menge von $n$
 Punkten, die $X$ erzeugt. Ist $n=0$, so ist $X = 0$ und dann auch
 $Y = 0$, so dass der Satz in diesem Falle gilt. Es sei also $n >
 0$. Ist $Y = X$, so ist nichts zu beweisen. Es sei also $Y < X$.
 Es gibt dann ein $P \in \Phi$ mit $P \not\leq Y$. Nach der
 Transformationsregel gilt daher $Y/0 \cong (Y + P)/P$. Ist $Q \in
 \Phi - \{P\}$, so ist $Q + P$ ein Atom des Quotienten $X / P$ und
 $X$ ist als Element von $X / P$ das Erzeugnis dieser Atome. Daher
 ist $\Rg_{X/P}(X) \leq n - 1$. Nach Induktionsannahme ist folglich
 $Y + P$ als Element von $X / P$ endlich erzeugt und es gilt
 $\Rg_{X/P}(Y + P) \leq n - 1$. Wegen der Isomorphie von $Y / 0$
 und $(Y + P)/P$ ist also auch $Y$ endlich erzeugt und die
 Ungleichung $\Rg_L(Y) \leq n - 1$ erf\"ullt.
 \medskip
       Statt zu sagen, dass $X$ endlich erzeugt\index{endlich erzeugt}{} sei,
 werden wir in Zukunft auch
 sagen, dass $X$ \emph{endlichen Ranges} sei.\index{endlichen Ranges}{}
 \par
       Ist $\Phi$ eine Menge von Punkten eines projektiven Verbandes $L$ und ist
 $P$ ein Punkt von $L$, so hei\ss t $P$ \emph{abh\"angig}\index{Abh\"angigkeit}{}
 von $\Phi$, falls $P \leq \sum_{Q \in \Phi}Q$ gilt. Hat
 die Menge $\Phi$ von Punkten von $L$ die Eigenschaft, dass keiner
 ihrer Punkte $P$ von $\Phi - \{P\}$ abh\"angt, so nennen wir
 $\Phi$ \emph{unabh\"angig}.\index{Unabh\"angigkeit}{}
 \par
       Ist $X$ ein Teilraum endlichen Ranges eines projektiven Verbandes
 und sind $P_1$, \dots, $P_{\Rg_L(X)}$ Punkte, die $X$ erzeugen, so
 ist die aus diesen Punkten gebildete Menge unabh\"angig.
 \medskip\noindent
 {\bf 3.2. Satz.} {\it Es sei $L$ ein projektiver Verband und $X$ sei ein
 Teilraum endlichen Ranges von $L$. Ist $\Phi$ eine unabh\"angige Teilmenge von
 Punkten von $X$, so ist $\Phi$ endlich und es gilt $|\Phi| \leq \Rg_L(X)$.}
 \smallskip
       Beweis. Der von $\Phi$ erzeugte Unterraum von $L$ ist nach 3.1
 endlich erzeugt und sein Rang ist h\"ochstens gleich dem Rang von
 $X$. Wir d\"urfen daher annehmen, dass $X = \sum_{Q \in
 \Phi}Q$ ist. Setze $n := \Rg_L(X)$. Ist $n = 1$, so ist die
 Aussage des Satzes offenkundig. Es sei also $n > 1$ und $\{P_1,
 \dots, P_n\}$ sei eine Menge von Punkten, die $X$ erzeugt. Setze
 $Y := \sum_{i:=1}^{n-1} P_i$. Dann ist
 $$ Y < X = \sum_{Q \in \Phi}Q. $$
 Es gibt also ein $Q \in \Phi$ mit $Q \not\leq Y$. Nach der Transformationsregel
 gilt
 $$ X/Y = (Y + P_n)/Y \cong P_n/(Y \cap P_n) = P_n/0. $$
 Andererseits sind auch $(Y + Q)/Y$ und $Q/0$ isomorph. Weil $P_n$ und $Q$
 Punkte sind, sind aber auch $P_n/0$ und $Q/0$ isomorph. Somit sind
 $X/Y$ und $(Y + Q)/Y$ isomorph. Weil der zweite Quotient im ersten
 enthalten ist und beide Quotienten nur je zwei Elemente enthalten,
 gilt schlie\ss lich $X = Y + Q$.
 \par
        Setze $Z := \sum _{R \in \Phi - \{Q\}}R$. Dann ist $Z + Q
 = X$ und $Z \cap Q = 0$. Also ist
 $$Y/0 \cong (Y + Q)/Q = (Z + Q) / Q \cong Z/0. $$
 Hieraus folgt $\Rg_L(Z) = \Rg_L(Y) \leq n - 1$. Nach Induktionsannahme ist daher
 $|\Phi - \{Q\}|\leq n - 1$ und folglich $|\Phi| \leq n$.
 \medskip\noindent
 {\bf 3.3. Korollar.} {\it Es sei $L$ ein projektiver Verband und $\Phi$ sei eine
 endliche unabh\"angige Menge von Punkten. Dann ist}
 $$\textstyle |\Phi| = \Rg_L \bigl(\sum_{P \in \Phi}P \bigr). $$
 \smallskip
       Beweis. Es sei $X := \sum_{P \in \Phi}P$. Dann gilt
 \emph{per definitionem} die Ungleichung $\Rg_L(X) \leq |\Phi|$.
 Andererseits gilt nach 3.2 auch $|\Phi| \leq \Rg_L(X)$. Damit ist
 das Korollar bewiesen.
 \smallskip
       Die Punkte eines projektiven Verbandes sind genau die Un\-ter\-r\"au\-me
 des Ranges 1 und die Geraden sind die Unterr\"aume des Ranges 2. Die
 Unterr\"aume des Ranges 3 nennen wir \emph{Ebenen}.\index{Ebene}{} Ferner nennen
 wir die Komplemente von Punkten \emph{Hyperebenen}\index{Hyperebene}{} und
 gelegentlich auch \emph{Ko-Atome}.\index{Ko-Atom}{}
 \medskip\noindent
 {\bf 3.4. Satz.} {\it Es sei $L$ ein projektiver Verband und $U(L)$ sei die
 Menge der unabh\"angigen Mengen von Punkten von $L$. Dann gilt:
 \item{a)} Ist $\Phi$ eine Menge von Punkten von $L$, so gilt genau dann
 $\Phi \in U(L)$, wenn $\Fin(\Phi) \subseteq U(L)$ gilt.
 \item{b)} Sind $\Phi$, $\Psi \in U(L)$, sind $\Phi$ und $\Psi$ beide endlich und
 gilt $|\Psi| = |\Phi| + 1$, so gibt es ein $P \in \Psi - \Phi$ mit
 $\Phi \cup \{P\} \in U(L)$.\par}
 \smallskip
       Beweis. a) Ist $\Phi \in U(L)$, so ist nat\"urlich $\Fin(\Psi)
 \subseteq U(L)$. Es sei $\Phi$ abh\"angig. Es gibt dann ein $P \in
 \Psi$, so dass $P$ von $\Phi - \{P\}$ abh\"angt. Nach dem Satz von
 der endlichen Abh\"angigkeit gibt es eine endliche Teilmenge
 $\Psi$ von $\Phi - \{P\}$, so dass $P$ von $\Psi$ abh\"angt. Es
 folgt, dass $\Psi \cup \{P\}$ eine endliche abh\"angige Teilmenge
 von $\Phi$ ist, so dass $\Fin(\Phi) \not\subseteq U(L)$ gilt.
 \par
       b) Setze $X := \sum_{Q \in \Phi}Q$. Nach 3.2 gibt es dann
 wegen $|\Psi| > |\Phi|$ ein $P \in \Psi$ mit $P \not\leq X$. Setze
 $Y := X + P$. Dann ist $Y$ endlich erzeugt und \"uberdies $X < Y$.
 Nach 3.3 und 3.1 ist daher
 $$ |\Phi| = \Rg_L(X) < \Rg_L(Y) \leq \big|\Phi \cup \{P\}\big|
	   = |\Phi| + 1. $$
 Es folgt $\Rg_L(Y) = \mid \Phi \cup \{P\} \mid$, so dass $\Phi \cup \{P\}$ in
 der Tat unabh\"angig ist.
 \medskip
       Satz 3.4 besagt, dass $U(L)$ ein Beispiel f\"ur das ist, was man
 \emph{Un\-abh\"ang\-ig\-keits\-stru\-ktur}\index{Unabh\"angigkeitsstruktur}{} nennt. Die
 Eigenschaft b) wird gew\"ohnlich \emph{stei\-nitz\-scher
 Austauschsatz}\index{steinitzscher Austauschsatz}{} genannt.
 \par
       Es sei $X$ ein Teilraum des projektiven Verbandes $L$. Ist $\Phi$
 eine Menge von Punkten von $X$, so hei\ss t $\Phi$ eine
 \emph{Basis}\index{Basis}{} von $X$, wenn $\Phi$ unabh\"angig ist und $X$
 erzeugt. Ist $X$ endlichen Ranges und ist $\Phi$ eine Basis von
 $X$, so ist, wie wir oben gesehen haben, $|\Phi| = \Rg_L(X)$. Ferner gilt:
 \medskip\noindent
 {\bf 3.5. Basissatz.} {\it Es sei $L$ ein projektiver Verband. Ist $X \in L$ und
 ist $\Psi$ eine unabh\"angige Menge von Punkten von $X$, so gibt es eine Basis
 $\Phi$ von $X$ mit $\Psi \subseteq \Phi$. Insbesondere hat jedes $X \in L$
 eine Basis.}
 \smallskip
       Beweis. Ist $X$ endlichen Ranges, so gibt es \emph{per definitionem} eine
 Basis $\Omega$ von $X$. Weil $ |\Psi| \leq |\Omega|$ gilt, folgt die Existenz
 von $\Psi$ mittels Induktion aus dem Steinitzschen Austauschsatz.
 \par
       Ist $X$ nicht endlichen Ranges, so erschlie\ss e man die Existenz
 von $\Phi$ mittels des zornschen Lemmas.
 \par
       Die letzte Aussage dieses Satzes folgt mit $\Psi := \emptyset$ aus
 dem bereits Bewiesenen.
 \medskip
       Weil $U(L)$ eine Unabh\"angigkeitsstruktur ist, folgt, dass sich zwei
 Basen eines Elements $X \in L$ stets bijektiv aufeinander abbilden lassen. Wir
 werden diesen Beweis hier nicht durchf\"uhren. F\"ur einen solchen sei der
 interessierte Leser auf L\"uneburg 1989 verwiesen. Die allen Basen von $X$
 gemeinsame Kardinalzahl hei\ss t \emph{Rang}\index{Rang}{} von $X$. Sie werde
 ebenfalls mit $\Rg_L(X)$ bezeichnet.
 \medskip\noindent
 {\bf 3.6. Satz.} {\it Es sei $L$ ein projektiver Verband und $X$ und $Y$ seien
 Elemente von $L$. Ferner sei $\Phi$ eine Basis von $X$ und $\Psi$ eine Basis von
 $Y$. Ist $X \cap Y = 0$, so ist $\Phi \cap \Psi = \emptyset$ und
 $\Phi \cup \Psi$ ist eine Basis von $X + Y$.}
 \smallskip
       Beweis. Die Aussage \"uber den Schnitt von $\Phi$ mit $\Psi$ ist
 trivial. Ferner ist klar, dass $X + Y$ von $\Phi \cup \Psi$
 erzeugt wird. Es ist zu zeigen, dass $\Phi \cup \Psi$ unabh\"angig
 ist. Um dies zu zeigen, setzen wir zun\"achst $\Omega := \Phi \cup
 \Psi$. W\"are $\Omega$ abh\"angig, so g\"abe es ein $P \in \Omega$,
 so dass $P$ von $\Omega - \{P\}$ abhinge. Wir k\"onnten oBdA
 annehmen, dass $P \in \Phi$ g\"alte. Setzte man dann $Z :=
 \sum_{Q \in \Phi - \{P\}}Q$, so folgte $X + Y = Z + Y$
 und weiter, weil $L$ ja modular ist,
 $$ X = X \cup (Z + Y) = Z + (X \cap Y) = Z $$
 im Widerspruch zur Unabh\"angigkeit von $\Phi$. Damit ist alles gezeigt.
 \medskip\noindent
 {\bf 3.7. Rangformel.}\index{Rangformel}{} {\it Es sei $L$ eine projektiver
 Verband. Sind $X$ und $Y$ Elemente endlichen Ranges von $L$, so ist auch $X + Y$
 endlichen Ranges und es gilt}
 $$ \Rg_L(X) + \Rg_L(Y) = \Rg_L(X + Y) + \Rg_L(X \cap Y). $$
 \par
       Beweis. Weil $X$ und $Y$ endlich erzeugt sind, ist es auch $X + Y$. Es sei
 $W$ ein Komplement von $X \cap Y$ in $Y$. Nach 3.6 ist dann
 $$ \Rg_L(Y) = \Rg_L(X + Y) + \Rg_L(W). $$
 Nun ist
 $$ X +Y = X + (X \cap Y) + W = X + W $$
 und
 $$ X \cap W = X \cap Y \cap W = 0. $$
 Also ist $\Rg_L(X + W) = \Rg_L(X) + \Rg_L(W)$. Somit ist
 $$\eqalign{
   \Rg_L(X) + \Rg_L(Y) &= \Rg_L(X) + \Rg_L(X \cap Y) + \Rg_L(W) \cr
	&= \Rg_L(X + Y) + \Rg_L(X \cap Y). \cr} $$
 \par
       Dieser Beweis funktioniert nat\"urlich auch f\"ur unendliche
 Kar\-di\-nal\-zah\-len, doch in diesem Falle taugt die Rangformel nicht viel.

 \mysection{4. Vollst\"andig reduzible Moduln}

 \noindent
 Der Geometrie\index{Geometrie}{} wird nachgesagt, dass sie von jedem Fortschritt
 in der Mathematik profitiere, dass sie aber selbst nur wenig bis gar
 nichts zum Fortschritt der Mathematik beitrage. Dass man die
 Geometrie jedoch manchmal dazu heranziehen kann, Dinge in anderen
 Teilen der Mathematik besser zu verstehen, m\"ochte ich hier am
 Beispiel der vollst\"andig reduziblen Moduln demonstrieren.
 \par
       Es sei zun\"achst $R$ ein Ring, wobei wir nicht voraussetzen, dass $R$
 eine $1$ habe. Es sei ferner $M$ ein $R$-Rechtsmodul. Mit $L_R(M)$ bezeichnen
 wir die Menge aller Teilmoduln von $M$ und f\"ur die Einschr\"ankung der
 Inklusionsrelation auf $L_R(M)$ benutzen wir das Symbol $\leq$. Mit diesen
 Verabredungen gilt nun der folgende Satz.
 \medskip\noindent
 {\bf 4.1. Satz.} {\it Ist $M$ ein $R$-Rechtsmodul, so ist $(L_R(M), \leq)$ ein
 modularer, voll\-st\"an\-di\-ger und nach oben stetiger Verband. Der Schnitt von
 Teilmoduln im mengentheoretischen Sinne ist gleich ihrem Schnitt im
 verbandstheoretischen Sinne, und ihre Summe im mo\-dul\-the\-o\-re\-ti\-schen
 Sinne stimmt \"uberein mit ihrer Summe im
 ver\-bands\-theo\-re\-ti\-schen Sinne.}
 \smallskip
       Der simple Beweis dieses Satzes sei dem Leser als \"Ubungsaufgabe
 \"uberlassen.
 \medskip
       Das Rechtsideal $I$ des Ringes $R$ hei\ss t
 \emph{maximal},\index{maximales Rechtsideal}{} wenn
 $$ \big|L_R(R/I)\big| = 2 $$
 ist. Das Rechtsideal $I$ hei\ss t
 \emph{regul\"ar},\index{regul\"ares Rechtsideal}{} wenn es ein $a \in R$ gibt
 mit $x - ax \in I$ f\"ur alle
 $x \in R$.  Hat $R$ eine Eins, so ist jedes Rechtsideal regul\"ar. Ist $I$ ein
 regul\"ares Rechtsideal, und ist $J$ ein Rechtsideal mit $I \subseteq J$, so ist
 auch $J$ regul\"ar.
 \par
       Der $R$-Modul $M$ hei\ss t \emph{irreduzibel},\index{irreduzibler Modul}{}
 falls $M R \neq \{0\}$ ist und $L_R(M)$ genau zwei Elemente enth\"alt. Ist $M$
 irreduzibel und ist $0 \neq u \in M$, so ist $M = uR$. Ist
 n\"amlich $J := \{v \mid v \in M, vR = \{0\}\}$, so ist $J$ ein
 Teilmodul von $M$, der wegen $MR \neq \{0\}$ und der Irreduzibilit\"at von $M$
 gleich $\{0\}$ ist. Wegen $u \notin J$ ist also $uR = M$.
 \medskip\noindent
 {\bf 4.2. Satz.} {\it Es sei $R$ ein Ring und $M$ sei ein irreduzibler
 $R$-Modul. Ist $0 \neq u \in M$, so definieren wir den Epimorphismus $\sigma$
 des $R$-Rechtsmoduls $R$ auf $uR = M$ durch $\sigma(r) := ur$. Dann ist
 $\Kern(\sigma)$ ein maximales, regul\"ares Rechtsideal von $R$. Ist umgekehrt
 $I$ ein maximales, regul\"ares Rechtsideal von $R$, so ist $R/I$ ein
 irreduzibler $R$-Modul.}
 \smallskip
       Beweis. Weil $M$ irreduzibel ist, ist $\Kern(\sigma)$ ein
 maximales Rechts\-i\-de\-al von $R$. Nun ist $u \in M = uR$. Es gibt
 also ein $a \in R$ mit $u = ua$. Es folgt $u(r - ar) = 0$ f\"ur alle
 $r \in R$. Somit ist $r - ar \in \Kern(\sigma)$ f\"ur alle $r \in
 R$, so dass $\Kern(\sigma)$ regul\"ar ist.
 \par
       Es sei jetzt umgekehrt $I$ ein maximales, regul\"ares Ideal von $R$.
 Dann ist $|L_R(R/I)| = 2$. Es sei nun $a \in R$ mit $r - ar \in I$
 f\"ur alle $r \in R$. Dann ist $A \notin I$, da andernfalls $r \in
 I$ w\"are f\"ur alle $r \in R$. Wegen $r + I = ar + I = (a + I)r$ ist
 daher $(R/I)R=R/I$, so dass $R/I$ irreduzibel ist.
 \medskip
       Zwei Eigenschaften fehlen dem Verband $L_R(M)$, um projektiv zu
 sein, die Komplementarit\"at und die Atomarit\"at, wobei Atome in
 diesem Zusammenhang die irreduziblen Teilmoduln des $R$-Moduls $M$
 sind. Eine sehr n\"utzliche Charakterisierung dieser Moduln
 liefert der n\"achste Satz.
 \medskip\noindent
 {\bf 4.3. Satz.} {\it Es sei $R$ ein Ring und $M$ sei ein $R$-Rechtsmodul. Genau
 dann ist $L_R(M)$ projektiv, wenn $M$ Summe von irreduziblen Teil\-mo\-duln
 ist.}
 \smallskip
       Beweis. Es sei $M$ Summe von irreduziblen Teilmoduln. Ferner  sei
 $X \in L_R(M)$. Mittels des zornschen Lemmas erhalten wir ein $Y
 \in L_R(M)$ maximal bez\"uglich der Eigenschaft $X \cap Y=\{0\}$.
 Wir wollen zeigen, dass $X + Y = M$ ist. Dazu nehmen wir an, dies
 sei nicht der Fall. Weil $M$ Summe von irreduziblen Teilmoduln
 ist, gibt es einen irreduziblen Teilmodul $P$ von $M$ mit $P
 \not\leq X + Y$. Es folgt $(X + Y) \cap P = \{0\}$ und dann auch
 $Y \cap P = \{0\}$. Nun ist $Y \subseteq Y + P$. Mittels der
 Modularit\"at von $L_R(M)$ folgt daraus
 $$ (Y + P) \cap (Y + X) = Y + \bigl(X \cap (Y + P)\bigr). $$
 Andererseits ist $Y \subseteq Y + X$ und daher
 $$ (Y + X) \cap (Y + P) = Y + \bigl(P \cap (Y + X)\bigr) = Y. $$
 Aus diesen beiden Gleichungen folgt $Y = Y + (X \cap (Y + P))$, was wiederum
 $X \cap (Y + P) \subseteq Y$ zur Folge hat.  Also gilt
 $$ X \cap (Y + P) \subseteq X \cap Y = \{0\}.$$
 Nun ist $Y \subseteq Y + P$, so dass die Maximalit\"at von $Y$ erzwingt, dass
 $Y = Y + P$ ist. Dies ergibt aber den Widerspruch $P \leq Y \cap P = \{0\}$.
 Also ist doch $M = X \oplus Y$, so dass $L_R(M)$ komplement\"ar ist.
 \par
       Um zu zeigen, dass $L_R(M)$ atomar ist, sei $X$ ein von $\{0\}$
 verschiedener Teilmodul von $M$. Es sei $0 \neq x \in X$. Weil $M$ Summe von
 irreduziblen Teilmoduln ist, gibt es dann irreduzible Teilmoduln $I_1$, \dots,
 $I_n$ von $M$ mit $x \in \Sigma_{k:=1}^{n}I_k$. Es sei $n$ minimal mit dieser
 Eigenschaft. Dann ist
 $$\sum_{k:=1}^{n} I_k = \bigoplus_{k:=1}^{n} I_k.$$
 Zu jedem $k$ gibt es nun ein $i_k \subseteq I_k$ mit
 $x = \Sigma_{k:=1}^n i_k$. Ist nun $r \in R$ und $xr = 0$, so folgt aus
 der Unabh\"angigkeit der $I_k$, dass $i_kr=0$ ist f\"ur alle $k$.
 F\"ur $m \in M$ setze man $O(m) := \{r \mid r \in R, mr = 0\}$. Dann
 folgt aus dem gerade Bewiesenen
 $$ O(x) \subseteq \bigcap_{k:=1}^n O(i_k). $$
 Es gilt sogar die Gleichheit, wie man unmittelbar sieht,
 uns gen\"ugt es aber zu wissen, dass $O(x) \subseteq O(i_1)$ gilt.
 Wegen der Minimalit\"at von $n$ ist $i_1 \neq 0$ und daher $I_1 =
 i_1R$. Definiere den Modulhomomorphismus $\sigma$ von $R$ auf $xR$
 durch $\sigma(r) := xr$. Weil $L_r(M)$ als modularer und
 komplement\"arer Verband auch relativ komplement\"ar ist, gibt es
 einen Teilmodul $P$ von $xR$ mit $xR = \sigma (O(i_1)) \oplus P$.
 Es folgt, dass
 $$ P \cong xR/\sigma\bigl(O(i_1)\bigr) \cong R/O(i_1) \cong I_1 $$
 ist. Also ist $P$ ein Punkt mit $P \leq xR \leq X$, so dass $L_R(M)$ auch
 atomar ist.
 \par
       Ist $L_R(M)$ ein projektiver Verband, so ist jedes Element dieses
 Verbandes obere Grenze der in ihm enthaltenen Punkte.
 Insbesondere gilt das f\"ur $M$, so dass $M$ nach 4.1 die Summe
 irreduzibler Teilmoduln ist. Damit ist der Satz bewiesen.
 \medskip
       Hat $R$ eine Eins, so kann man mehr beweisen, da in diesem Falle
 alle Rechts\-ide\-a\-le regul\"ar sind.
 \medskip\noindent
 {\bf 4.4. Satz.} {\it Es sei $R$ ein Ring mit Eins und $M$ sei ein unit\"arer
 $R$-Modul. Genau dann ist $L_R(M)$ projektiv, wenn $L_R(M)$ komplement\"ar ist.}
 \smallskip
       Beweis. Es ist nat\"urlich nur zu zeigen, dass die Komplementarit\"at die
 Projektivit\"at nach sich zieht.
 \par
       Es sei also $L_R(M)$ komplement\"ar und $X$ sei ein von $\{0\}$
 verschiedener Teilmodul von $M$. Ferner sei $0 \neq x \in X$. Weil
 $M$ unit\"ar ist, ist $1 \not\in O(x)$. Mit Hilfe des zornschen
 Lemmas erschlie\ss t man die Existenz eines maximalen Rechtsideals
 $I$ mit $O(x) \subseteq I$. Weil $I$ regul\"ar ist, ist $R/I$
 nach 4.2 ein ir\-re\-du\-zib\-ler $R$-Rechtsmodul. L\"asst man dieses
 $I$ die Rolle von $O(i_1)$ im Beweise von 4.3 spielen, so sieht
 man, dass es einen Punkt $P$ gibt mit $P \leq xR \leq X$. Also ist
 $L_R(M)$ atomar und damit projektiv.
 \medskip
       Ist $M$ ein $R$-Rechtsmodul, dessen Teilmodulverband projektiv ist, so
 hei\ss t $M$ \emph{vollst\"andig reduzibel}.\index{vollst\"andig reduzibel}{}
 Dieser Name ist nicht sehr gl\"uck\-lich gew\"ahlt, da er \"ublicherweise nur
 beinhaltet, dass jeder Teilmodul von $M$ ein direkter Summand ist.
 Da wir auch an Ringen ohne Eins interessiert sind, m\"ussen wir
 unter diesen Begriff auch die Atomarit\"at von $L_R(M)$ subsumieren.
 \par
       Ist nun $M$ ein vollst\"andig reduzibler $R$-Modul, so ist Satz 2.15 auf
 $L_R(M)$ anwendbar. Die interessante algebraische Interpretation dieses Satzes
 ist Inhalt des n\"achsten Satzes.
 \medskip\noindent
 {\bf 4.5. Satz.} {\it Es sei $M$ ein vollst\"andig reduzibler $R$-Modul und
 $I_R(M)$ bezeichne die Menge seiner irreduziblen Teilmoduln. Sind $P$,
 $Q \in I_R(M)$, so setzen wir $P \equiv Q$ genau dann, wenn $P$ und $Q$
 isomorphe $R$-Moduln sind. Dann ist $\equiv\ =\ \sim$, wobei $\sim$ wie in 2.8
 definiert sei. F\"ur alle $\Phi \in I_R(M)/{\equiv}$ setzen wir
 $$ H_\Phi := \sum_{P \in \Phi}P. $$
 Dann ist
 $$ M = \bigoplus_{\Phi \in I_R(M)/ \equiv} H_{\Phi}. $$}
 \par
       Beweis. Wir zeigen, dass $\equiv \; = \; \cong$ ist. Dazu seien $A$ und
 $B$ zwei verschiedene Elemente aus $I_R(M)$. Dann ist $A \cap B = \{0\}$.
 \par
       Es gelte $A \equiv B$. Es gibt dann einen Isomorphismus $\sigma$
 von $A$ auf $B$. Setze $C:= \{a + a^\sigma \mid a \in A\}$. Ist $0 \neq
 p \in A$, so ist $A = pR$, da $A$ ja irreduzibel ist. Daher ist
 $$ (p + p^\sigma)R = \bigl\{pr + (pr)^\sigma \mid r \in R\bigr\}
	= \{a + a^\sigma \mid a \in A\} = C. $$
 Somit ist auch $C$ irreduzibel. Ferner
 gilt $C \neq A$, $B$ und $C \leq A + B$. Folglich ist $A \sim B$.
 \par
       Es sei umgekehrt $A \sim B$ und $C$ sei ein von $A$ und $B$
 verschiedener, irredu\-zib\-ler Teilmodul in $A + B$. Nun ist $A + B =
 A \oplus B$. Es gibt daher Projektionen $\alpha$ und $\beta$ von
 $C$ in $A$ beziehungsweise $B$. Weil $A$, $B$ und $C$ irreduzibel
 sind und $C$ von $A$ und auch $B$ verschieden ist, sind $\alpha$
 und $\beta$ Isomorphismen. Es folgt, dass $\alpha^{-1} \beta$
 ein Isomorphismus von $A$ auf $B$ ist. Somit gilt auch $A \equiv
 B$, womit die Gleichheit der beiden \"Aquivalenzrelationen
 bewiesen ist. Mittels 2.15 und 4.1 folgt nun die Behauptung des Satzes.
 \medskip
       Die $H_\Phi$ hei\ss en aus offensichtlichem Grund \emph{homogene
 Kom\-po\-nen\-ten}\index{homogene Komponente}{} von $M$.
 \par
       Das Rechtsideal $I$ des Ringes $R$ hei\ss t
 \emph{minimal},\index{minimales Rechtsideal}{} wenn
 $L_R(I)$ genau zwei Elemente enth\"alt. Man beachte, dass ein
 minimales Rechtsideal als Rechtsmodul \"uber $R$ nicht notwendig
 irreduzibel ist, da ja durch\-aus $IR = \{0\}$ sein kann.
 \par
       Der Ring $R$ hei\ss t \emph{vollst\"andig
 reduzibel},\index{vollst\"andig reduzibel}{} falls er als
 Rechtsmodul \"uber sich selbst vollst\"andig reduzibel ist. Die
 Atome eines solchen Ringes sind gerade die minimalen Rechtsideale
 von $R$. Ist n\"amlich $I$ ein minimales Rechtsideal, so gibt es
 ein Atom, also ein weiteres minimales Rechtsideal $P$, mit $P \leq
 I$. Hieraus folgt $P = I$. Es gilt also f\"ur alle minimalen
 Rechtsideale $I$ von $R$, dass $I R = I$ ist.
 \medskip\noindent
 {\bf 4.6. Satz.} {\it Es sei $R$ ein vollst\"andig reduzibler Ring und
 $I_R(R)/{\equiv}$ sei die Menge der \"Aquivalenzklassen isomorpher, minimaler
 Rechtsideale von $R$. Ist $\Phi \in I_R(R)/{\equiv}$, so ist $H_\Phi$ ein
 zweiseitiges Ideal von $R$. Ferner gilt im ringtheoretischen Sinne
 $$ R = \bigoplus_{\Phi \in I_R(R)/{\equiv}} H_\Phi. $$}
 \par
       Beweis. Als Summe von Rechtsidealen ist $H_\Phi$ nat\"urlich auch
 ein Rechtsideal. Es sei nun $I$ ein minimales Rechtsideal in
 $H_\Phi$. Ist $0 \neq r \in I$, so ist $I = rR$. Es sei $s \in R$.
 Wir definieren $\sigma$ durch $(rk)^\sigma := s r k$ f\"ur alle $k
 \in R$. Dann ist $\sigma$ ein Epimorphismus von $I$ auf $s r R$.
 Weil $I$ minimal ist, ist daher entweder $s r R =\{0\}$ oder $s r
 R$ ist ein zu $I$ isomorphes Rechtsideal von $R$. In beiden
 F\"allen ist $s r R \leq H_\Phi$. Somit sind die homogenen
 Komponenten von $R$ zweiseitige Ideale.
 \par
       Sind $\Phi$ und $\Psi$ verschiedene \"Aquivalenzklassen minimaler
 Rechts\-i\-de\-a\-le, so ist $H_\Phi \cap H_\Psi = \{0\}$ und daher $xy =
 0$ f\"ur $x \in H_\Phi$ und $Y \in H_\Psi$. Hieraus folgt, dass
 $R$ auch im ringtheoretischen Sinne die direkte Summe seiner
 homogenen Komponenten ist. Damit ist alles bewiesen.
 \medskip
       Der vollst\"andig reduzible Ring $R$ hei\ss t
 \emph{homogen},\index{homogener Ring}{}
 falls er nur eine homogene Komponente hat. Dies ist
 gleichbedeutend damit, dass alle minimalen Rechts\-ide\-a\-le von $R$
 als $R$-Rechtsmoduln isomorph sind. Ist $R$ einfach, dh., besitzt
 $R$ nur die beiden zweiseitigen Ideale $\{0\}$ und $R$, so ist $R$
 nach dem gerade bewiesenen Satz homogen. Die Umkehrung gilt nicht,
 wie wir noch sehen werden.
 \medskip\noindent
 {\bf 4.7. Satz.} {\it Ist $R$ ein vollst\"andig reduzibler Ring, und ist $M$ ein
 $R$-Rechtsmodul, so ist $M$ genau dann vollst\"andig reduzibel, wenn $MR = M$
 ist. Ist $MR = M$ und ist $R$ homogen, so ist $L_R(M)$ irreduzibel.}
 \smallskip
       Beweis. Es sei $M$ ein $R$-Rechtsmodul mit $MR = M$. Wir zeigen,
 dass $M$ Summe von Atomen ist. Dazu sei $0 \neq y \in M$. Wir
 definieren die Abbildung $\varphi$ von $R$ auf $yR$ durch
 $r^\varphi := yr$. Dann sind die Moduln $R/\Kern(\varphi)$
 und $yR$ isomorph. Weil $L_R(R)$ ein projektiver Verband ist, ist
 auch $L_R(R/\Kern(\varphi))$ ein projektiver Verband. Also
 ist $R/\Kern(\varphi)$ und dann auch $yR$ Summe von Atomen.
 Folglich ist $M$ wegen $M = MR$ Summe von Atomen, so dass $M$ ein
 vollst\"andig reduzibler $R$-Modul ist.
 \par
       Ist $M$ vollst\"andig reduzibel, so ist $M$ Summe von irreduziblen
 Teilmoduln. Ist $P$ ein solcher und ist $0 \neq y \in P$, so ist
 $P = yR$, wie wir wissen. Hieraus folgt, dass $MR = M$ ist.
 \par
       Es bleibt, die Irreduzibilit\"at von $L_R(M)$ zu beweisen, falls $R$
 homogen ist und $MR = M$ gilt. In diesem Falle sind alle minimalen
 Rechtsideale von $R$ isomorph, da $R$ nur eine homogene Komponente
 hat. Da ein Atom von $M$ aber stets zu einem minimalen Rechtsideal
 von $R$ isomorph ist, wie unmittelbar aus 4.2 folgt, sind auch alle Atome von
 $M$ isomorph, woraus sich die Irreduzibilit\"at von $L_R(M)$ ergibt.
 \medskip
       Jeder K\"orper $K$,\index{korper@K\"orper}{} ob kommutativ oder nicht, ist
 nat\"urlich ein voll\-st\"an\-dig reduzibler
 Ring.\index{vollst\"andig reduzibel}{} Daher erhalten wir aufs Neue den
 Satz, dass $L_K(V)$ ein projektiver Verband ist, falls nur $V$ ein
 Rechtsvektorraum \"uber $K$ ist. Da $K$ ein einfacher Ring ist,
 ist $L_K(V)$ auch irreduzibel. Gibt es weitere vollst\"andig
 reduzible Ringe? Die Antwort lautet: Ja.
 \par
       Um dies einzusehen, holen wir zun\"achst etwas weiter aus. Es sei $V$ ein
 Rechtsvektorraum \"uber dem K\"orper $K$. Es sei ferner $\End_K(V)$ der
 Endomorphismenring von $V$, wobei wir das Bild des Vektors $v$ unter dem
 Endomorphismus $\sigma$ mit $\sigma(v)$ bezeichnen. Schlie\ss lich ben\"otigen
 wir noch die Definition des Begriffs \emph{von
 Neumann-Ring}.\index{von Neumann-Ring}{} Ein Ring $R$ werde so genannt, falls es
 zu jedem $r \in R$ ein $s \in R$ gibt mit $rsr = r$. Ist $rsr = r$, so sind die
 Elemente $rs$ und $sr$ \emph{Idempotente},\index{Idempotent}{} dh., es gilt
 $(rs)^2 = rs$ und $(sr)^2 = sr$, wie man unmittelbar sieht. Die Idempotente
 eines Endomorphismenringes hei\ss en auch
 \emph{Projektionen}.\index{Projektion}{} Ist $R$ ein von Neumann-Ring, sind
 $r$, $s \in R$ und gilt $rsr = r$, so ist $rsR = rR$. Denn einmal ist $rs \in
 rR$ und andererseits ist $r = (rs)r \in rsR$. Jedes
 Haupt\-rechts\-ide\-al eines von Neumann-Ringes wird also von einem
 Idempotenten erzeugt. Hieraus folgt weiter, dass f\"ur ein
 Rechtsideal $I$ von $R$, welches nicht gleich $\{ 0 \}$ ist, auch
 $I^2$ von $\{ 0 \}$ verschieden ist, da jedes derartige Ideal ein
 von $0$ verschiedenes, idempotentes Element enth\"alt. Die gleiche
 Aussage gilt nat\"urlich auch f\"ur Linksideale.
 \medskip\noindent
 {\bf 4.8. Satz.} {\it Ist $V$ ein Rechtsvektorraum \"uber dem K\"orper
 $K$, so ist $\End_K(V)$ ein von Neumann-Ring. Genauer: Es sei $Y$
 ein Komplement von $\Kern(\varphi)$, wobei $\varphi$ ein
 Endomorphismus von $V$ sei. Dann ist $\varphi(V) = \varphi(Y)$ und
 die Einschr\"ankung von $\varphi$ auf $Y$ ist ein Monomorphismus.
 Es sei weiter $Z$ ein Komplement von $\varphi(Y)$. Definiere
 $\psi$ durch $\psi(z) := 0$ f\"ur $z \in Z$ und $\psi(\varphi(y))
 :=y$ f\"ur alle $y \in Y$. Dann ist $\varphi \psi \varphi =
 \varphi$. \"Uberdies sind die R\"ange von $\varphi$ und $\psi$
 gleich.}
 \smallskip
      Beweis. Es ist klar, dass die Einschr\"ankung von $\varphi$ auf
 $Y$ injektiv ist, so dass die Definition von $\psi$ korrekt ist.
 Es sei nun $v \in V$. Es gibt dann ein $y \in Y$ und ein $k \in
 \Kern(\varphi)$ mit $v = y + k$. Es folgt
 $$\eqalign{
       \varphi \psi \varphi(v)
	  &= \varphi \psi \varphi (y + k)
	   = \varphi\psi\bigl(\varphi(y) + \varphi(k)\bigr)
	   = \varphi \bigl(\psi \varphi(y)\bigr)                \cr
	  &= \varphi(y) = \varphi(y) + \varphi(k) = \varphi(v), \cr} $$
 so dass in der Tat $\varphi \psi \varphi = \varphi$ gilt.
 \medskip
       Wir ziehen einige Folgerungen aus diesem Satz.
 \medskip\noindent
 {\bf 4.9. Satz.} {\it Es sei $V$ ein $K$-Vektorraum und $J_K(V)$ sei
 die Menge der Endomorphismen endlichen Ranges von $V$. Dann ist
 $J_K(V)$ ein zweiseitiges Ideal von \emph{End}$_K(V)$. Ferner
 gilt, dass $J_K(V)$, f\"ur sich betrachtet, ein von Neumann-Ring
 ist.}
 \smallskip
       Beweis. Es ist eine simple \"Ubungsaufgabe zu zeigen, dass
 $J_K(V)$ ein zweiseitiges Ideal ist. Ist nun $\rho \in J_K(V)$, so
 gibt es nach 4.8 ein $\sigma \in \End_K(V)$ mit $\rho \sigma \rho
 = \rho$, so dass die R\"ange von $\rho$ und $\sigma$ gleich sind.
 Also gilt sogar $\sigma \in J_K(V)$. Damit ist alles bewiesen.
 \medskip\noindent
 {\bf 4.10. Satz.} {\it Es sei $V$ ein $K$-Vektorraum.
 \item{a)} Ist $\varphi \in \End_K(V)$, so gibt es eine Projektion
 $\pi \in J_K(V)$ mit
 $$ \varphi \End_K(V) = \pi \End_K(V). $$
 \item{b)} Ist $\varphi \in J_K(V)$, so gibt es eine Projektion $\pi \in J_K(V)$
 mit}
 $$ \varphi J_K(V) = \pi J_K(V). $$
 \par
       Beweis. Dies folgt aus der Bemerkung, die vor 4.8 gemacht wurde, und der
 Tatsache, dass die beiden Ringe $\End_K(V)$ und $J_K(V)$ von Neumann-Ringe sind.
 \medskip\noindent
 {\bf 4.11. Satz.} {\it Es seien $\pi$ und $\rho$ Projektionen des
 $K$-Vektorraumes $V$. Ist $\pi(V) = \rho(V)$ und hat $\pi(V)$ den
 Rang $1$, so ist $\pi \End_K(V) = \rho \End_K(V)$ und $\pi
 \End_K(V)$ ist ein Rechtsideal von $\End_K(V)$, welches minimal
 ist. Ferner gilt: Ist $I$ ein minimales Ideal von $\End_K(V)$, so
 gibt es eine Projektion $\pi$ des Ranges $1$ mit $I = \pi
 \End_K(V)$.}
 \smallskip
       Beweis. Setze $P := \pi(V)$. Dann ist
 $$ P \oplus \Kern(\pi) = V = P \oplus \Kern(\rho). $$
 Ist $\Kern(\pi) = \Kern(\rho)$, so ist $\pi = \rho$ und die von $\pi$
 und $\rho$ erzeugten Ideale sind gleich. Es sei also $\Kern(\pi)
 \neq \Kern(\rho)$. Dann ist
 $$ V = \Kern(\pi) + \Kern(\rho) $$
 und daher
 $$\eqalign{
   \Kern(\pi)/\bigl(\Kern(\pi) \cap \Kern(\rho)\bigr)
      &\cong \bigl(\Kern(\pi) + \Kern(\rho)\bigr)/\Kern (\rho) \cr
      &= V/\Kern(\rho). \cr} $$
 Es gibt also einen Punkt $A$ auf $\Kern(\pi)$ mit
 $$ \Kern(\pi) = A + \bigl(\Kern(\pi) \cap \Kern(\rho)\bigr). $$
 Die Gerade $A + P$ ist ein Komplement von $\Kern(\pi) \cap \Kern(\rho)$. Ebenso
 folgt die Existenz eines Punktes $B$ auf Kern$(\rho)$ mit
 $$ \Kern(\rho) = B + \bigl(\Kern(\pi) \cap \Kern(\rho)\bigr). $$
 Die Gerade $B + P$ ist ebenfalls ein Komplement von
 $\Kern(\pi) \cap \Kern(\rho)$. Es seien nun $a$, $b$
 und $p$ Vektoren mit $A = aK$, $B = bK$ und $P = pK$. Wir definieren
 dann eine lineare Abbildung $\lambda$ von $V$ in sich durch
 $\lambda (p) := p$, $\lambda(a) := b$ und $\lambda (u) := u$ f\"ur
 alle $u \in \Kern(\pi) \cap \Kern(\rho)$. Dann ist $
 \pi = \rho \lambda$ und, da $\lambda$ offenbar invertierbar ist,
 $\rho = \pi \lambda^{-1}$, so dass in der Tat $\rho \End_K(V) =
 \pi \End_K(V)$ ist.
 \par
       Es sei nun $0 \neq \sigma \in \pi \End_K(V)$. Nach 4.10 d\"urfen
 wir annehmen, dass $\sigma$ eine Projektion ist. Mit dem bereits
 Bewiesenen folgt, dass $\sigma \End_K(V) = \pi \End_K(V)$ ist.
 Folglich ist $\pi \End_K(V)$ ein minimales Ideal.
 \par
       Es sei $0 \neq \pi \in I$. Dann ist $\pi \End_K(V) = I$, da
 $\End_K(V)$ ja eine Eins hat und $I$ minimal ist. Nach 4.10
 d\"urfen wir annehmen, dass $\pi$ eine Projektion ist. Es sei $P
 := \pi(V)$. Ferner sei $P = Q \oplus C$ mit einem Teilraum $Q$ des
 Ranges 1. Ist dann $\sigma$ die Projektion von $V$ auf $Q$, deren
 Kern gleich $C \oplus \Kern(\pi)$ ist, so ist $\sigma = \pi
 \sigma$. Weil $I$ minimal und $\sigma \neq 0$ ist, folgt $\sigma
 \End _K(V) = \pi \End_K(V)$. Es gibt also ein $\gamma \in
 \End_K(V)$ mit $\pi = \sigma \gamma$. Daher ist
 $$ P = \pi(V) = \sigma \gamma (V) \leq \sigma (V) = Q, $$
 so dass in der Tat $\Rg (P) = 1$ ist.
 \medskip\noindent
 {\bf 4.12. Satz.} {\it Es sei $V$ ein Rechtsvektorraum \"uber dem K\"orper $K$.
 Die mi\-ni\-mal\-en Rechtsideale von $\End_K(V)$ sind genau die minimalen Rechtsideale
 von $J_K(V)$, wobei $J_K(V)$ f\"ur sich als Ring betrachtet ist.}
 \smallskip
       Beweis. Es sei $I$ ein minimales Rechtsideal von $J_K(V)$. Ferner sei
 $0 \neq \sigma \in I$. Nach 4.10 b) gibt es eine Projektion $\pi \in J_K(V)$
 mit $\sigma J_K(V) = \pi J_K(V)$. Es folgt $\pi \neq 0$ und weiter $0 \neq \pi
 = \pi^2 \in \pi J_K(V)$. Die Minimalit\"at von $I$ erzwingt daher, dass
 $$ \sigma J_K(V) = \pi J_K(V) = I $$
 ist. Da $\End_K(V)$ eine Eins hat und $J_K(V)$ ein Ideal von $\End_K(V)$ ist,
 gilt
 $$ J_K(V) \End_K(V) = J_K(V). $$
 Also ist
 $$ I \End_K(V) = \sigma J_K(V) \End_K(V) = \sigma J_K(V) = I, $$
 so dass $I$ ein Rechtsideal von $\End_K(V)$ ist. Hieraus folgt schlie\ss lich
 $$ I = \sigma J_K(V) \subseteq \sigma \End_K(V) \subseteq I, $$
 so dass $I$ auch ein minimales Ideal von $\End_K(V)$ ist.
 \par
       Es sei umgekehrt $I$ ein minimales Rechtsideal von $\End_K(V)$. Es
 gibt dann eine Projektion $\pi$ vom Range 1 mit $I = \pi
 \End_K(V)$. Es folgt $\pi \in J_K(V)$ und damit $I \subseteq
 J_K(V)$. Also ist $I$ ein Rechtsideal von $J_K(V)$. Es sei nun
 $I'$ ein von $\{0\}$ verschiedenes Rechtsideal von $J_K(V)$,
 welches in $I$ enthalten ist. Dann enth\"alt $I'$ eine Projektion
 $\sigma$ ungleich Null. Dann ist $\sigma J_K(V)$ ein von $\{0\}$
 verschiedenes Rechtsideal von $\End_K(V)$, welches in $I$
 enthalten ist. Aus der Minimalit\"at von $I$ folgt daher
 $$ I = \sigma J_K(v) \subseteq I' \subseteq I, $$
 so dass $I = I'$ ist. Also ist $I$ auch ein minimales Ideal von $J_K(V)$.
 \medskip\noindent
 {\bf 4.13. Satz.} {\it Ist $V$ ein $K$-Rechtsvektorraum, so ist
 $J_K(V)$ ein einfacher, vollst\"andig reduzibler Ring.}
 \smallskip
       Beweis. Wir zeigen zun\"achst, dass $J_K(V)$ vollst\"andig
 reduzibel ist. Dazu sei $\sigma \in J_K(V)$. Weil $J_K(V)$ ein von
 Neumann-Ring ist, gibt es ein $\rho \in J_K(V)$, so dass $\sigma =
 \sigma \rho \sigma$ ist. Setze $\pi := \sigma \rho$. Dann ist
 $\pi$ eine Projektion und $\pi(V)$ hat endlichen Rang. Es gibt
 also eine Basis $b_1$, \dots, $b_n$ von $\pi(V)$. Wir definieren
 Projektionen $\pi$, \dots, $\pi_n$ durch
 $$ \pi_i(b_j) := \cases{b_i & f\"ur $i = j$ \cr
	        	   0 & f\"ur $i \neq j$ \cr} $$
 und $\pi_i(h) := 0$ f\"ur $h \in \Kern(\pi)$. Dann ist $\pi =
 \sum_{i :=1}^n \pi_i$ und daher
 $$ \sigma = \pi\sigma = \sum_{i := 1}^n \pi_i \sigma. $$
 Hieraus und aus Satz 4.12 folgt, dass $J_K(V)$ Summe von
 minimalen Idealen ist. Weil die minimalen Ideale Projektionen
 enthalten, sind sie auch irreduzible $J_K(V)$-Moduln, so dass
 $J_K(V)$ ein vollst\"andig reduzibler Ring ist.
 \par
       Als N\"achstes zeigen wir, dass $J_K(V)$ homogen ist. Dazu seien
 $I$ und $I'$ minimale Rechtsideale dieses Ringes. Es gibt dann
 zwei Punkte $P$ und $Q$ von $L_K(V)$ mit
 $$ I = \bigl\{\sigma \mid \sigma \in J_K(V), \sigma(V) \leq P\bigr\} $$
 und
 $$ I' = \bigl\{ \sigma \mid \sigma \in J_K(V), \sigma(V) \leq Q\bigr\}.$$
 Ist $P = Q$, so ist $I = I'$ und nichts weiter zu beweisen. Es sei also
 $P \neq Q$. Es sei $C$ ein Komplement von $P + Q$. Ferner sei $P = pK$ und
 $Q = qK$ sowie $R := (p + q)K$.  Setze $H := R + C$. Dann ist $H$ ein Komplement
 von $P$ wie auch von $Q$. Es sei $\pi$ die Projektion von $V$ auf $P$ mit
 $\Kern(\pi) = H$ und $\pi'$ die Projektion von $V$ auf $Q$ mit
 $\Kern(\pi') = H$. Dann ist $I = \pi J_K(V)$ und $I' = \pi'J_K(V)$.
 Definiere $\lambda$ durch $\lambda(p) := q$, $\lambda (q) := p$ und
 $\lambda(c) := 0$ f\"ur alle $c \in C$. Dann ist
 $$ \lambda \pi \lambda(q) = \lambda \pi (p) = \lambda(p) = q = \pi'(q) $$
 und
 $$ \lambda\pi\lambda (p + q) = \lambda\pi(p + q) = 0 = \pi'(p + q)$$
 sowie
 $$ \lambda \pi \lambda (c) = 0 = \pi'(c) $$
 f\"ur alle $c \in C$. Weil $p + Q = p + R$ und daher auch $V = (P + R) \oplus C$
 gilt, folgt $\pi' = \lambda \pi \lambda$. Ferner ist klar, dass
 $\lambda$ in $J_K(V)$ liegt.
 \par
       Wir definieren nun $\varphi$ durch $\varphi (\pi \sigma) :=
 \lambda \pi \sigma$. Dann ist $\varphi$ ein Mo\-dul\-epi\-mor\-phis\-mus von
 $I$ auf $\varphi (I)$. Nun ist $\varphi (\pi \lambda) = \lambda
 \pi \lambda = \pi'$. Weil $I$ und $I'$ minimale Rechtsideale sind,
 folgt hieraus, dass $\varphi (I) = I'$ ist. Daher ist $\varphi$
 ein Isomorphismus von $I$ auf $I'$. Damit ist gezeigt, dass
 $J_K(V)$ homogen ist.
 \par
       Es sei schlie\ss lich $N$ ein von $\{0\}$ verschiedenes, zweiseitiges
 Ideal von $J_K(V)$. Ferner sei $\pi J_K(V)$ ein minimales Rechtsideal, welches
 in $N$ enthalten sei, und $\pi$ sei idempotent. Eine solche Konstellation
 existiert, da $N$ ja ungleich $\{0\}$ ist. Es sei weiter $I$ ein minimales
 Rechtsideal von $J_K(V)$. Weil $J_K(V)$ homogen ist, gibt es einen
 Mo\-dul\-iso\-mor\-phis\-mus $\varphi$ von $\pi J_K(V)$ auf $I$. Setze
 $\iota := \varphi(\pi)$. Dann ist
 $$ 0 \neq \iota = \varphi(\pi) = \varphi(\pi^2) = \iota \pi \in I \cap N, $$
 so dass $I \subseteq N$ gilt. Weil $J_K(V)$ Summe von minimalen Rechtsidealen
 ist, ist also $J_K(V) \subseteq N$, so dass $J_K(V)$, wie behauptet, einfach
 ist.
 \par
       Das zuletzt benutzte Argument wird uns sp\"ater noch einmal gute
 Dienste leisten.
 \medskip
       Mit den Ringen $J_K(V)$ haben wir nun viele Beispiele von
 voll\-st\"an\-dig reduziblen Ringen gewonnen, die in aller Regel keine
 K\"orper sind. Diese Ringe haben \"uberdies keine Eins, wenn $V$
 nicht endlichen Ran\-ges ist. Hat $V$ endlichen Rang, so ist $J_K(V)
 = \End_K(V)$, so dass in diesem Falle $\End_K(V)$ ein einfacher,
 vollst\"andig reduzibler Ring ist. Man kann noch etwas mehr sagen, n\"amlich:
 \medskip\noindent
 {\bf 4.14. Satz.} {\it Es sei $V$ ein $K$-Rechtsvektorraum. Ist $N$ ein von
 $\{0\}$ verschiedenes, zweiseitiges Ideal von $\End_K(V)$, so ist
 $J_K(V) \subseteq N$.}
 \smallskip
       Beweis. Wie wir wissen, enth\"alt $N$ eine von Null verschiedene
 Projektion $\pi$. Es sei $P$ ein Punkt mit $P \leq \pi(V)$ und $\rho$ sei eine
 Projektion von $V$ auf $P$ mit $\Kern(\pi) \subseteq \Kern(\rho)$. Ist
 $v \in V$, so gibt es ein $c \in \pi(V)$ und ein $x \in \Kern(\pi)$ mit
 $v = c + x$. Es folgt
 $$ \rho\pi(v) = \rho\bigl(\pi(c) + \pi (x)\bigr) = \rho(c) = \rho(c) + \rho(x)
		 = \rho(v). $$
 Also ist
 $$ 0 \neq \rho = \rho \pi \in J_K(V) \cap N, $$
 dh., $J_K(V) \cap N$ ist ein von $\{0\}$ verschiedenes, zweiseitiges Ideal von
 $\End_K(V)$. Weil $J_K(V)$ einfach ist, folgt hieraus $J_K(V) \cap N = J_K(V)$,
 so dass in der Tat $J_K(V) \subseteq N$ gilt.
 \medskip
       Es ist nat\"urlich zu fragen, was $L_{J_{K}(V)}(J_K(V))$ mit
 $L_K(V)$ zu tun hat. Man erwartet, dass diese beiden Verb\"ande
 isomorph sind, und die Erwartung tr\"ugt nicht, wie der n\"achste Satz lehrt.
 \medskip\noindent
 {\bf 4.15. Satz.} {\it Es sei $V$ ein $K$-Rechtsvektorraum. Definiert
 man $\Phi$ und $\Psi$ durch
 $$ \Phi(U) := \bigl\{\sigma \mid \sigma \in J_K(V), \sigma(V) \leq U\bigr\} $$
 f\"ur alle $U \in L_K(V)$ und
 $$ \textstyle \Psi(I) := \sum_{\sigma \in I} \sigma(V) $$
 f\"ur alle Rechtsideale $I$ von $J_K(V)$, so ist $\Phi$ ein Isomorphismus von
 $L_K(V)$ auf $L_{J_{K}(V)} (J_K(V))$ und es gilt $\Psi = \Phi^{-1}$.}
 \smallskip
       Beweis. Es ist trivial, dass $\Psi(U)$ f\"ur alle $U \in L_K(V)$ ein
 Rechts\-i\-de\-al von $J_K(V)$, und noch banaler, dass $\Psi(I)$ f\"ur
 alle Rechtsideale $I$ von $J_K(V)$ ein Teilraum von $V$ ist.
 \par
       Ist $U \in L_K(V)$, so gilt
 $$\eqalign{
      \Psi \Phi (U) &= \Psi\bigl(\{\sigma \mid \sigma \in J_K(V),
			                  \sigma(V) \leq U\}\bigr) \cr
	 &= \sum_{\sigma \in J_K(V), \sigma(V) \leq U} \sigma(V)  
	 \ \leq\ U. \cr} $$
 Ist $P$ ein Punkt mit $P \leq U$, so gibt es eine
 Projektion $\pi$ mit $\pi(V) = P$. Es folgt $P \leq \Psi \Phi
 (U)$. Weil $U$ die obere Grenze der in $U$ enthaltenen Punkte ist,
 folgt weiter $\Psi \Phi(U) = U$. Somit ist $\Psi \Phi = id_{L_{K}(V)}$.
 \par
       Es sei $I$ ein Rechtsideal von $J_K(V)$. Dann gilt
 $$
 \Phi \Psi (I) = \Phi \biggl(\sum_{\sigma \in I} \sigma(V) \biggr) 
	       = \bigl\{\tau \mid \tau \in J_K(V), \tau(V) \leq
		      \sum_{\sigma \in I} \sigma(V)\bigr\} 
	       \ \leq\ I. 
 $$
 
       Es sei $\tau$ eine Projektion von $V$ auf einen Punkt von
 $\sum_{\sigma \in I} \sigma (V)$. Wir zeigen, dass $\tau
 \in I$ gilt. Es gibt eine endliche Teilmenge $E$ von $I$ mit
 $$ \tau(V) \leq \sum_{\sigma \in E} \sigma(V), $$
 da $\tau$ ja endlichen Rang hat. Zu jedem $\sigma \in E$ gibt es
 endlich viele minimale Rechtsideale in $I$, deren Summe $\sigma$
 enth\"alt. Hieraus schlie\ss t man, dass es endlich viele
 Projektionen $\pi_1$, \dots, $\pi_n \in I$ gibt, so dass alle
 $\pi_i$ den Rang 1 haben und
 $$ \tau (V) \leq \sum_{i:=1}^n \pi_i(V) $$
 gilt. Ist $\tau(V) = \pi_1(V)$ --- dies ist gewiss dann
 der Fall, wenn $n = 1$ ist ---, so gilt nach 4.12 und 4.11 die
 Gleichung
 $$ \tau J_K(V) = \tau \End_K(V) = \pi \End_K(V) = \pi_1J_K(V), $$
 so dass in diesem Falle $\tau \in I$ gilt. Wir d\"urfen daher annehmen, dass
 $\tau (V) \neq \pi_1(V)$ gilt.
 \par
       Es ist $\tau (V) \leq \pi_1(V) + \sum_{i := 2}^n \pi_i (V)$. Nach
 2.6 gibt es daher eine Projektion $\sigma$ von $V$ auf
 einen Punkt von $\sum_{i := 2}^n \pi_i (V)$ mit $\tau(V)
 \leq \pi_1(V) + \sigma (V)$. Nach der nicht explizit formulierten
 Induktionsannahme ist $\sigma \in I$, so dass wir oBdA annehmen
 d\"urfen, dass $\sigma = \pi_2$ ist. Wegen 4.11 d\"urfen wir
 weiter annehmen, dass $\pi_1(V) \leq \Kern(\pi_2)$ und $\pi_2(V)
 \leq \Kern(\pi_1)$ gilt. Dann ist $\Kern(\pi_1) \cap \Kern(\pi_2)$ ein
 Kom\-ple\-ment von $\pi_1(V) + \pi_2(V)$. Es
 sei $\rho$ die Projektion von $V$ auf $\pi_1(V) + \pi_2(V)$ mit
 $\Kern(\rho) = \Kern(\pi_1) \cap \Kern(\pi_2)$.
 Dann ist $\rho \tau(v) = \tau (v)$ f\"ur alle $v \in V$ und daher
 $\rho \tau = \tau$. Ist nun $v \in V$, so gibt es ein $u_1 \in
 \pi_1(V)$, ein $u_2 \in \pi_2(V)$ und ein $x \in \Kern (\rho)$ mit
 $v = u_1 + u_2 + x$. Es folgt
 $$
    \rho(v) = \rho (u_1 + u_2 + x) = u_1 + u_2          
	    = \pi_1(v) + \pi_2(v) = (\pi_1 + \pi_2)(v), 
 $$
 so dass $\rho = \pi_1 + \pi_2 \in I$ gilt. Es folgt $\tau = \rho
 \tau \in I.$
 \par
       Es gilt also, dass alle minimalen Rechtsideale, die in $\Phi\Psi(I)$
 ent\-hal\-ten sind, bereits in $I$ liegen. Daher gilt auch $\Phi
 \Psi (I) \leq I$, so dass $\Phi \Psi (I) = I$ ist. Damit ist
 gezeigt, dass $\Phi$ eine Bijektion und dass $\Psi$ ihre Inverse
 ist. Da beide Abbildungen offensichtlich inklusionstreu sind, ist
 alles bewiesen.
 \medskip
       Einmalig neugierig geworden, m\"ochte man m\"oglichst alles \"uber
 vollst\"andig reduzible Ringe wissen. Wir werden daher im
 \"uber\-n\"achs\-ten Abschnitt noch einmal auf sie zu sprechen kommen.
 Zuvor jedoch wollen wir die Basis unserer Kenntnisse noch verbreitern.

\mysection{5. Der duale Verband}

 \noindent
 Es sei $(L, \leq)$ ein Verband. Wir definieren $(L^d, \leq^d)$,
 indem wir $L^d := L$ setzen und $\leq^d$ dadurch definieren, dass
 $A \leq^d B$ genau dann gelte, wenn $B \leq A$ gilt. Bezeichnet
 man die obere Grenze zweier Elemente $A$ und $B$ von $L^d$ mit $A
 +^d B$ und ihre untere Grenze mit $A \cap^d B$, so ist $A +^d B =
 A \cap B$ und $A \cap^d B = A + B$. Daher ist auch $(L^d, \leq^d)$
 ein Verband. Er hei\ss t der zu $(L, \leq)$
 \emph{duale Verband}.\index{dualer Verband}{}
 \medskip\noindent
 {\bf 5.1. Satz.} {\it Ist $L$ ein Verband und ist $L^d$ der zu $L$
 duale Verband, so gelten die folgenden Aussagen.
 \item{a)} Ist $L$ modular, so ist auch $L^d$ modular.
 \item{b)} Ist $L$ komplement\"ar, so ist auch $L^d$ komplement\"ar.
 \item{c)} Ist $L$ vollst\"andig, so ist auch $L^d$ vollst\"andig.}
 \smallskip
       Beweis. a) Der Verband $L$ sei modular. Ferner seien $A$, $B$, $C \in
 L^d$ und es gelte $B \leq^d A$. Dann ist also $A \leq B$. Hieraus
 folgt
 $$ A \cap^d ( B +^d C) = A + (B \cap C) = B \cap (A + C)
		  = B +^d (A \cap^d C), $$
 so dass auch $L^d$ modular ist.
 \par
       b) und c) folgen unmittelbar aus der Definition der Komplementarit\"at,
 bzw. der Vollst\"andigkeit.
 \smallskip
 Wichtig ist auch der n\"achste Satz.
 \medskip\noindent
 {\bf 5.2. Satz.} {\it Es sei $L$ ein modularer und komplement\"arer
 Verband. Ist $L$ atomar, so ist auch $L^d$ atomar.}
 \smallskip
       Beweis. Die Atome von $L^d$ sind gerade die Komplemente der Atome
 von $L$. Wir nennen sie, wie schon im Fall der projektiven Verb\"ande,
 \emph{Hyperebenen}\index{Hyperebene}{} von $L$ oder auch
 \emph{Ko-Atome}\index{Ko-Atom}{} von $L$.
 \par
      Wir m\"ussen zeigen, dass jedes vom gr\"o\ss ten Element $\Pi$ von
 $L$ verschiedene Element $A$ in einer Hyperebene von $L$ liegt.
 Weil $L$ komplement\"ar ist, gibt es ein $B \in L$ mit $\Pi = A \oplus B$. Wegen
 $A \neq \Pi$ ist $B \neq 0$. Es gibt also ein Atom $P$ von $L$ mit $P \leq B$.
 Nach 2.10 ist $L$ relativ komplement\"ar. Es gibt daher ein $C$ mit
 $B = P \oplus C$. Setze $K := A + C$. Dann ist $A \leq K$. Ferner gilt
 $$ P + K = P + A + C = A + B = \Pi. $$
 Andererseits ist
 $$P \cap K \leq B \cap K = B \cap (C + A) = C + (B \cap A) = C $$
 und daher $P \cap K \leq P \cap C = 0$. Folglich ist $K$ ein Ko-Atom von $L$,
 welches $A$ umfasst. Damit ist alles bewiesen.
 \medskip\noindent
 {\bf 5.3. Korollar.} {\it Ist $L$ ein modularer, komplement\"arer und atomarer
 Verband, so ist jedes vom gr\"o\ss ten Element verschiedene Element von $L$
 gleich dem Schnitt der es umfassenden Hyperebenen.}
 \smallskip
       Beweis. Dies folgt mittels 5.2, 5.1 und 2.11.
 \medskip
       Der Verband $L$ hei\ss t \emph{noethersch},\index{noetherscher Verband}{}
 wenn jede nicht leere Teilmenge von $L$ ein maximales Element enth\"alt. Er
 hei\ss t \emph{artinsch},\index{artinscher Verband}{} wenn jede seiner nicht
 leeren Teilmengen ein minimales Element enth\"alt.
 \medskip\noindent
 {\bf 5.4. Satz.} {\it Ein modularer, komplement\"arer Verband ist
 genau dann artinsch, wenn er noethersch ist.}
 \smallskip
       Beweis. Wegen 5.1 und der Bemerkung, dass $(L^d)^d = L$ ist,
 gen\"ugt es zu zeigen, dass ein modularer und komplement\"arer
 Verband ar\-tinsch ist, falls er noethersch ist.
 \par
       Es sei also $L$ ein solcher Verband und $M$ sei eine nicht leere
 Teilmenge von $L$. Wir nehmen an, dass $M$ kein minimales Element
 enth\"alt. Dann gibt es zu jedem $X \in M$ ein $Y \in M$ mit $Y <
 X$. Auf Grund des Auswahlaxioms gibt es also eine Abbildung $f$
 von $M$ in sich mit $f(X) < X$ f\"ur alle $X \in M$. Ist nun $Y
 \in M$, so folgt mittels des dedekindschen Rekursionssatzes die
 Existenz einer Abbildung $g$ der nicht negativen ganzen Zahlen in
 $M$ mit $g(0) = Y$ und $g(n + 1) < g(n)$ f\"ur alle $n$. Weil $L$
 nach 2.10 relativ komplement\"ar ist, folgt mittels des
 Auswahlaxioms die Existenz einer Abbildung $c$ der Menge der
 nat\"urlichen Zahlen in $L$ mit $g(n) = g(n + 1) \oplus c(n + 1)$
 f\"ur alle $n \geq 0$. Setze $W(n) := \sum_{i:=1}^n
 c(i)$. Dann ist
 $$ W(n) \leq W (n + 1) $$
 f\"ur alle $n$.  Wir zeigen, dass $W(n) \cap g(n) = 0$ ist. Dies ist sicher
 richtig, falls $n = 1$ ist. Es sei also $n \geq 1$. Dann gilt auf
 Grund der Induktionsannahme und der Modularit\"at von $L$, dass
 $$ W(n + 1) \cap g(n + 1)
    \leq \biggl( \sum_{i:=1}^n c(i) + c(n + 1) \biggr) \cap g(n) = c(n + 1) $$
 ist. Also ist
 $$ W(n + 1) \cap g(n + 1) \leq c(n + 1) \cap g(n + 1) = 0. $$
 Weil $L$ noethersch ist, gibt es nun ein $N$ mit $W(N) =
 W(n)$ f\"ur alle $n \geq N$. Hieraus folgt
 $$ c(N + 1) + W(N) = W(N + 1) = W(N) $$
 und damit $c(N + 1) \leq W(N)$.  Andererseits ist auch $c(N + 1) \leq g(N)$.
 Folglich gilt $c(N + 1) \leq W(N) \cap g(N) = 0$ und damit
 $$ g(N) = c(N +  1) + g(N + 1) = g(N + 1) $$
 im Widerspruch zu $g(N + 1) < g(N)$. Also
 enth\"alt $M$ doch ein minimales Element, so dass $L$ artinsch
 ist. Damit ist Satz 5.4 bewiesen.
 \medskip\noindent
 {\bf 5.5. Satz.} {\it Ein noetherscher, atomarer, komplement\"arer und
 modularer Verband ist projektiv.}
 \smallskip
       Beweis. Es sei $L$ ein solcher Verband. Es ist zu zeigen, dass $L$
 vollst\"andig und nach oben stetig ist.
 \par
       Es sei $M \subseteq L$ und $S$ sei die Menge der unteren Schranken
 von $M$. Dann ist $S$ nicht leer, da $0 \in S$ gilt. Weil $L$
 noethersch ist, enth\"alt $S$ also ein maximales Element $B$. Ist
 $C \in S$, so ist auch $B + C \in S$ und wegen $B \leq B + C$
 daher $B = B + C$. Folglich ist $C \leq B$, so dass $B$ das
 gr\"o\ss te Element von $S$ ist. Somit ist $B = \bigcap_{X \in M} X$,
 so dass $L$ nach 2.1 vollst\"andig ist.
 \par
       Es sei $M$ ein aufsteigendes System von $L$. Dann enth\"alt $M$
 ein maximales Element $G$. Weil $M$ aufsteigend ist, folgt $X \leq
 G$ f\"ur alle $X \in M$. Also ist $\sum_{X \in M} X = G$.
 Es sei nun $Y \in L$. Dann ist
 $$ Y \cap \sum_{X \in M} X = Y \cap G. $$
 Ferner ist $Y \cap X \leq Y \cap G$ f\"ur alle $X \in M$. Somit gilt
 $$ Y \cap G \leq \sum_{X \in M} (Y \cap G) \leq Y \cap G, $$
 so dass $L$ auch nach oben stetig ist.
 \medskip
       Als N\"achstes charakterisieren wir die noetherschen projektiven
 Ver\-b\"an\-de als diejenigen projektiven Verb\"ande, deren Rang
 endlich ist.
 \medskip\noindent
 {\bf 5.6. Satz.} {\it Ein projektiver Verband ist genau dann noethersch, wenn er
 end\-lich\-en Rang hat.}
 \smallskip
       Beweis. Es sei $L$ ein projektiver Verband endlichen Ranges. Ist $M$ eine
 nicht leere Teilmenge, so gilt nach 3.1, dass $\Rg_L (X) \leq \Rg(L)$ ist f\"ur
 alle $X \in M$. Es gibt daher ein $X$
 maximalen Ranges in $M$. Ist nun $Y \in M$ und gilt $X \leq Y$, so
 folgt $\Rg_L(Y) \leq \Rg_L(X)$. Hieraus folgt mit 3.1, dass $X = Y$
 ist, so dass $X$ maximal ist. Also ist $L$ noethersch.
 \par
       Es sei umgekehrt $L$ ein noetherscher projektiver Verband. Ist
 dann $M$ die Menge der endlich erzeugten Teilr\"aume von $L$, so
 enth\"alt $M$ ein maximales Element $\Pi$. Ist $P$ ein Punkt von
 $L$, so ist auch $\Pi + P$ endlich erzeugt und daher $\Pi = \Pi +
 P$, was $P \leq \Pi$ zur Folge hat. Folglich ist $\Pi$ das
 gr\"o\ss te Element von $L$, so dass $L$ endlichen Rang hat.
 \par
       Wir sind nun in der Lage, die Frage zu beantworten, unter wel\-chen
 Be\-ding\-ung\-en der zu einem projektiven Verband duale Verband
 projektiv ist. Auf Grund von 2.15 gen\"ugt es, diese Frage f\"ur
 irreduzible projektive Verb\"ande zu beantworten. Diese lapidare
 Feststellung ist nicht so ganz ohne. Sie besagt n\"amlich, dass
 $L^d$ genau dann projektiv ist, wenn f\"ur jeden irreduziblen
 Bestandteil $I$ der in 2.15 beschriebenen Zerlegung von $L$ gilt,
 dass auch $I^d$ projektiv ist. Um dies zu beweisen, muss man auf
 der Menge der Hyperebenen von $L$ eine \"Aquivalenzrelation $\sim$
 definieren, so dass neben $H \sim H$ f\"ur alle Hyperebenen $H$
 gilt, dass f\"ur zwei verschiedene Hyperebenen $H$ und $H'$ genau
 dann $H \sim H'$ gilt, wenn es eine dritte Hyperebene $H''$ gibt
 mit $H \cap H' \leq H''$. Dass dies eine \"Aquivalenzrelation ist,
 folgt aus 5.1a) und 2.8. Es sei dem Leser \"uberlassen, den
 Beweis in allen Einzelheiten auszuf\"uhren.
 \medskip\noindent
 {\bf 5.7. Satz.} {\it Ist $L$ ein irreduzibler projektiver Verband, so ist $L^d$
 genau dann projektiv, wenn der Rang von $L$ endlich ist.}\index{dualer Verband}{}
 \smallskip
       Beweis. $L$ habe endlichen Rang. Dann ist $L$ nach 5.6 noethersch
 und nach 5.4 dann auch artinsch. Somit ist $L^d$ noethersch. Da
 $L^d$ nach 5.1 modular und komplement\"ar und nach 5.2 auch atomar
 ist, ist $L^d$ nach 5.5 auch projektiv.
 \par
       Es sei $L$ ein projektiver Verband unendlichen Ranges. Ferner sei
 $H$ eine Hyperebene von $L$ und $\Phi$ sei die Menge aller
 unabh\"angigen Punktmengen von $L$, die keinen Punkt von $H$
 enthalten. Dann ist $\Phi$ von endlichem Charakter und enth\"alt
 folglich nach dem Lemma von Teichm\"uller und Tukey ein maximales
 $B$. Setze $U := \sum_{P \in B} P$. W\"are $U \cap H < H$,
 so g\"abe es einen Punkt $Q \leq H$ mit $Q \not\leq U \cap H$. Ist
 $P \in B$, so w\"are $P + Q$ eine Gerade, die weder in $U$ noch in
 $H$ l\"age. Weil $L$ irreduzibel ist, g\"abe es nach 2.14 auf $P +
 Q$ einen dritten Punkt $R$, der dann weder in $U$ noch in $H$
 l\"age. Aus der Maximalit\"at von $B$ folgte die Abh\"angigkeit
 von $B \cup \{R\}$. Hieraus folgte der Widerspruch $R \leq H$.
 Also ist doch $U \cap H = H$. Weil $U$ wenigstens einen Punkt
 enth\"alt, der nicht auf $H$ liegt, ist $U = \Pi$, wobei $\Pi$
 wieder das gr\"o\ss te Element von $L$ bezeichne. Somit ist $B$
 eine Basis von $\Pi$, so dass $B$ insbesondere unendlich ist.
 Folglich gibt es eine abz\"ahlbare Teilmenge $\{P_1, P_2, P_3,
 \dots \} \subseteq B$. Setze $H_i := \sum_{j:=i}^{\infty} P_j$. Dann bilden
 die $H_i$ eine absteigende Kette. Wir zeigen, dass
 $\bigcap_{i:=i}^{\infty} H_i = 0$ ist. Es seien $m$ und
 $n$ nat\"urliche Zahlen. Dann gelten die Gleichungen
 $$\eqalign{
   \Rg_L \biggl(\sum_{k:=1}^{n} P_k \bigg) &= n,  \cr
   \Rg_L \biggl( \sum_{k:=1}^{n + m} P_k \biggr) &= n + m, \cr
   \Rg_L \biggl( \sum_{k:= n + 1}^{n+m} P_k \biggr) &= m. \cr} $$
 Mittels der Rangformel folgt hieraus, dass
 $$ \biggl( \sum_{k:=1}^{n} P_k \biggr) \cap
	\biggl(\sum_{k:= n+1}^{n+m} P_k \biggr) = 0 $$
 ist. Der Satz von der endlichen Abh\"angigkeit,
 der ja eine Folge der Stetigkeit nach oben ist, liefert dann, dass
 $$ \biggl(\sum_{k:=1}^{n} P_k \biggr) \cap H_{n+1} = 0 $$
 ist. Hieraus folgt wiederum, dass
 $$ \bigcap_{i = 1}^{\infty} H_i = 0 $$
 gilt. W\"are n\"amlich $W$ ein Punkt in diesem Schnitt, so l\"age $W$ in $H_1$.
 Es g\"abe dann $P_1$, \dots, $P_n$ mit $W \leq \sum_{k:=1}^{n} P_k$, woraus der
 Widerspruch $W \leq (\sum_{k:=1}^{n} P_k) \cap H_{n+1} = 0$ folgte.
 \par
       W\"are nun $L^d$ projektiv, so w\"are $L^d$ nach oben stetig.
 Folglich w\"are
 $$ H = H + \bigcap_{i:=1}^\infty H_i = \bigcap_{i:=1}^\infty (H+H_i) = \Pi. $$
 Dieser Widerspruch zeigt, dass $L^d$ nicht projektiv ist, falls $L$ nicht
 endlichen Ranges ist.
 \medskip
       Ist $N$ eine maximale Kette aus Elementen des Verbandes $L$, so nennen wir
 $N$ ein \emph{Nest}\index{Nest}{} von $L$. Jedes Nest von $L$ enth\"alt
 nat\"urlich das gr\"o\ss te und das kleinste Element von $L$, falls $L$ solche
 Elemente besitzt. Enth\"alt das Nest $N$ von $L$ nur endlich viele Elemente, so
 nennen wir $|N| - 1$ die {\it L\"ange\/}\index{L\"ange eines Nestes}{} des
 Nestes.
 \medskip\noindent
 {\bf 5.8. Satz.} {\it Ist $L$ ein projektiver Verband, so sind die folgenden
 Bedingungen \"aquivalent:
 \item{a)} Der Rang von $L$ ist endlich.
 \item{b)} Jedes Nest von $L$ hat endliche L\"ange.
 \item{c)} Es gibt ein Nest endlicher L\"ange von $L$.

 \noindent
       Ist $\Rg(L)$ endlich, so ist die L\"ange eines jeden Nestes von $L$ gleich
 $\Rg(L)$.}
 \smallskip
       Beweis. a) impliziert b). Gilt $\Rg(L) \leq 1$, so ist $L$ das
 einzige Nest von $L$, so dass b) und auch der Nachsatz in diesem
 Falle korrekt ist. Es sei also $\Rg(L) \geq 2$ und $N$ sei ein
 Nest von $L$. Es sei $M$ die Menge der von $0$ verschiedenen
 Elemente in $N$. Weil der Rang von $L$ endlich ist, ist $L$
 artinsch, so dass $M$ ein minimales Element $X$ enth\"alt. Ist nun
 $P$ ein Punkt auf $X$, so gilt $0 < P \leq X$. Weil $N$ eine maximale Kette ist,
 ist daher $X = P$. Es ist folglich $\Rg(\Pi/X) = \Rg(L) - 1$.
 Weil $M$ ein Nest des Quotienten $\Pi/X$ ist, ist
 $M$ also endlich und es gilt, dass $\Rg(L) - 1$ die L\"ange von
 $M$ ist. Also ist auch $N$ endlich und $\Rg(L)$ ist die L\"ange
 von $N$. Damit ist b) aus a) hergeleitet und auch der Nachsatz
 vollst\"andig bewiesen.
 \par
       b) impliziert c). Auf Grund des hausdorffschen Ma\-xi\-mum\-prin\-zips
 besitzt $L$ ein Nest. Dieses hat dann nach b) endliche L\"ange.
 \par
       c) impliziert a). Es sei $N$ ein Nest endlicher L\"ange von $L$.
 Enth\"alt $N$ h\"ochstens zwei Elemente, so ist $N = L$ und daher
 $\Rg(L) \leq 1$. Wir d\"urfen daher annehmen, dass $N$ ein von $0$
 und $\Pi$ verschiedenes Element $X$ enthalte. Dann ist
 $$ N_1 := \{Y \mid Y \in N,\ Y \leq X\} $$
 ein Nest von $X/0$ und
 $$ N_2 := \{Z \mid Z \in N,\ X \leq Z\} $$
 ein Nest von $\Pi/X$. Weil diese beiden Nester k\"urzere L\"ange als $N$ haben,
 folgt mittels Induktion, dass $X/0$ und $\Pi/X$ endlichen Rang haben. Daher hat
 auch $L$ endlichen Rang.
 \par
       Damit ist alles bewiesen.
 \medskip\noindent
 {\bf 5.9. Satz.} {\it Ist $L$ ein irreduzibler projektiver Verband
 endlichen Ran\-ges, so ist auch $L^d$ ein irreduzibler Verband
 endlichen Ranges und es gilt $\Rg(L) = \Rg(L^d)$.}
 \smallskip
       Beweis. Der Verband $L^d$ ist nach 5.7 projektiv. Weil $0$ und
 $\Pi$ die einzigen Elemente von $L$ sind, die genau ein Komplement
 haben, haben sie auch als Elemente von $L^d$ genau ein Komplement,
 so dass auch $L^d$ irreduzibel ist. Weil jedes Nest von $L$ auch
 ein Nest von $L^d$ ist, hat $L^d$ nach 5.8 ebenfalls endlichen
 Rang und es gilt $\Rg(L) = \Rg(L^d)$.
 \medskip\noindent
 {\bf 5.10. Satz.} {\it Es sei $L$ ein irreduzibler projektiver Verband endlichen
 Ranges. Ist dann $X \in L$, so gilt}
 $$ \Rg_L(X) + \Rg_{L^d} (X) = \Rg(L). $$
 \par
       Der Beweis sei dem Leser als \"Ubungsaufgabe \"uberlassen. Eine Beweisidee
 findet sich im Beweise von 5.8.
 \medskip
       Ist $L$ ein projektiver Verband, ist $X \in L$ und ist der Rang des
 Quotienten $\Pi/X$ endlich, so setzen wir $\Ko_\Pi (X) := \Rg(\Pi/X)$. Wir
 nennen diese Zahl \emph{Ko-Rang}\index{Ko-Rang}{} von $X$ in $\Pi$.
 \medskip\noindent
 {\bf 5.11. Satz.} {\it Es sei $L$ ein irreduzibler projektiver Verband. Ist
 $X \in L$, so hat $X$ genau dann endlichen Ko-Rang, wenn $X$ Schnitt von
 endlichen vielen Hyperebenen ist. Ist $\Ko_\Pi (X)$ endlich, so ist $X$ Schnitt
 von $\Ko_\Pi (X)$ Hyperebenen, jedoch nicht Schnitt von weniger als
 $\Ko_\Pi (X)$ Hyperebenen.}
 \smallskip
       Beweis. Es sei $X$ Schnitt von endlich vielen Hyperebenen und
 $H_1$, \dots, $H_n$ sei eine minimale Menge von Hyperebenen, deren
 Schnitt $X$ sei. Wegen der Minimalit\"at von $n$ ist dann $H_i
 \cap H_n$ f\"ur $i := 1$, \dots, $n - 1$ eine Hyperebene in $H_n$.
 Ferner ist
 $$ X = \bigcap_{i:=1}^{n-1} (H_i \cap H_n). $$
 Mittels Induktion folgt $\Ko_{H_n} (X) = n - 1$ und damit dann
 $$ \Ko_\Pi (X) = n. $$
 \par
       Es sei umgekehrt der Ko-Rang von $X$ in $\Pi$ endlich. Dann ist
 $\Pi/X$ ein endlich erzeugter projektiver Verband, der nach Satz
 2.14 auch irreduzibel ist, so dass auch $(\Pi/X)^d$ ein endlich
 erzeugter projektiver Verband ist. Daher ist $X$ Schnitt von
 endlich vielen Hyperebenen.
 \par
       Die letzte Aussage des Satzes folgt schlie\ss lich daraus, dass
 der Ko-Rang von $X$ in $\Pi$ gerade der Rang von $(\Pi/X)^d$ ist.
 \medskip\noindent
 {\bf 5.12. Satz.} {\it Es sei $L$ ein irreduzibler projektiver
 Verband. Haben $X$, $Y \in L$ beide endlichen Ko-Rang in $\Pi$, so
 haben auch $X + Y$ und $X \cap Y$ endlichen Ko-Rang in $\Pi$ und
 es gilt}
 $$ \Ko_\Pi (X) + \Ko_\Pi (Y) = \Ko_\Pi (X + Y) + \Ko_\Pi (X \cap Y). $$
 \par
       Beweis. Weil $X$ und $Y$ Schnitte von endlich vielen Hyperebenen sind, ist
 auch $X \cap Y$ Schnitt von endlich vielen Hyperebenen. Die Ko-Rangformel ist
 daher nichts anderes als die Rangformel f\"ur
 $ \bigl(\Pi/(X \cap Y)\bigr)^d$.
 \medskip
       Ist $V$ ein Rechtsvektorraum \"uber dem K\"orper $K$, so ist
 $L_K(V)$ ein projektiver Verband, wie wir wissen. Hat $V$
 endlichen Rang, so hat auch $L_K(V)$ endlichen Rang und
 $L_K(V)^d$ ist ebenfalls ein projektiver Verband des glei\-chen
 Ranges. Es ist richtig, wie wir sehen werden, dass dieser Verband
 zum Teilraumverband des zu $V$ dualen Vektorraumes isomorph ist. Da die
 Dualit\"atstheorie in vielen B\"uchern\index{bucher@B\"ucher}{} der linearen Algebra
 miserabel dargestellt wird, erl\"autern wir sie hier noch einmal.
 Dabei werden dann gleich auch die ent\-sprech\-en\-den Notationen f\"ur sp\"ater
 fixiert. --- Viele B\"ucher\index{bucher@B\"ucher}{} der linearen Algebra handeln
 nur von Vektorr\"aumen \"uber kommutativen K\"orpern und
 unterscheiden nicht zwischen Links- und Rechtsvektorr\"aumen, was
 im Falle von nicht kommutativen K\"orpern zwingend erforderlich
 ist. Ihre Autoren\index{Autoren}{} glauben,
 didaktisch geschickt\index{didaktisch}{} vorzugehen. In
 Wirklichkeit aber erweisen sie ihren Lesern einen B\"arendienst,
 da die Realit\"at viele ihrer Leser irgendwann einholt, und diese
 dann oft erhebliche Verst\"andnisschwierigkeiten haben.
 Hinzu kommt, dass diese Autoren h\"aufig auch noch sagen, sie
 beschr\"ankten sich auf den Fall kommutativer
 K\"orper,\index{kommutativer K\"orper}{} im Falle nicht kommutativer K\"orper
 ginge eh alles genauso. Dabei
 ist z.~B.\ das Transponieren\index{transponieren}{} im Ring aller
 $(n \times n)$-Matrizen \"uber einem K\"orper $K$ genau dann ein
 Antiautomorphismus\index{Antiautomorphismus}{} dieses
 Ringes, wenn $K$ kommutativ ist. Es geht also nicht alles genauso.
 \par
       Es sei $V$ ein $K$-Rechtsvektorraum. Ist $\sigma$ ein
 Homomorphismus von $V$ in einem weiteren $K$-Vektorraum, so
 bezeichnen wir das Bild von $v \in V$ unter $\sigma$ mit $\sigma
 v$. Mit $V^*$ bezeichnen wir die Menge aller linearen Abbildungen
 von $V$ in $K$. Dann ist $V^*$ mit der punktweise definierten
 Addition als Verkn\"upfung eine abelsche Gruppe. Ist $f \in V^*$
 und $k \in K$, so definieren wir $kf$ durch $(kf)v := k(fv)$ f\"ur
 alle $v \in V$. Auf diese Weise wird $V^*$ zu einem
 $K$-Linksvektorraum, dem \emph{Dualraum}\index{Dualraum}{} von $V$.
 \medskip\noindent
 {\bf 5.13. Satz.} {\it Es sei $V$ ein $K$-Rechtsvektorraum und $B$ sei
 eine Basis von $V$. Zu jedem $b \in B$ definieren wir $b^* \in
 V^*$ durch die Vorschrift
 $$ b^*c := \cases{1 & f\"ur $c = b$     \cr
	           0 & f\"ur $c \neq b$. \cr} $$
 Dann ist die Menge
 $$ B^* := \{b^* \mid b \in B\} $$
 linear unabh\"angig.}
 \smallskip
       Beweis. Es sei $E$ eine endliche Teilmenge von $B$ und $k$ sei
 eine Abbildung von $E$ in $K$ mit $0 = \sum_{b \in E} k_b
 b^*$. Ist dann $c \in E$, so folgt
 $$ 0 = 0c = \sum_{b \in E} k_b b^* c = k_c. $$
 Somit sind alle endlichen Teilmengen von $B^*$ linear unabh\"angig und damit
 auch $B^*$.
 \medskip
       Zur Notation sei noch gesagt, dass wir bei Linksvektorr\"aumen und
 Links\-mo\-duln $\Rg^K (X)$, $L^K (V)$, $\End^R (M)$ etc.\ schreiben.
 \medskip\noindent
 {\bf 5.14. Satz.} {\it Es sei $V$ ein Rechtsvektorraum \"uber dem
 K\"orper $K$ und $B$ sei eine Basis von $V$. Genau dann ist $B^*$
 eine Basis von $V^*$, wenn der Rang von $V$ endlich ist. Ist der
 Rang von $V$ endlich, so ist $\Rg_K (V) = \Rg^K (V^*)$.}
 \smallskip
       Beweis. Der Rang von $V$ sei endlich. Ist $f \in V^*$ und $v \in
 V$, so gibt es eine Abbildung $k$ von $B$ in $K$ mit $v =
 \sum_{b \in B}bk_b.$ Es folgt
 $$\eqalign{
       fv &= f\biggl(\sum_{b \in B} bk_b\biggr) = \sum_{b \in B} fbk_b \cr
	  &= \sum_{b \in B} \sum_{c \in B} (fc) c^* bk_b
	   = \sum_{c \in B} \sum_{b \in B} (fc) c^* bk_b   \cr
	  &= \biggl(\sum_{c \in B} (fc) c^*\biggr) \sum_{b \in B} bk_b
	   = \biggl(\sum_{c \in B} (fc)c^* \biggr) v. \cr} $$
 Hieraus folgt, dass $f = \sum_{c \in B} (fc)c^*$. Damit ist gezeigt, dass $B^*$
 ein Erzeugendensystem von
 $V^*$ ist. Weil $B^*$ nach 5.13 auch linear unabh\"angig ist, ist
 $B^*$ in der Tat eine Basis von $V^*$. Ferner folgt, dass die
 R\"ange von $V$ und $V^*$ gleich sind. Die Basis $B^*$ hei\ss{}t die
 zu $B$ \emph{duale Basis}\index{duale Basis}{} von $V^*$.
 \par
       Der Rang von $V$ sei nun nicht endlich. Wir definieren $f \in V^*$
 durch $fb := 1$ f\"ur alle $b \in B$. W\"are $f$ ein Element von
 $\sum_{b \in B} Kb^*$, so g\"abe es eine Abbildung $k$
 von $B$ in $K$ mit endlichem Tr\"ager und $f = \sum_{b \in B} k_b
 b^*$. Weil der Tr\"ager von $f$ endlich, die Menge $B$ aber
 unendlich ist, g\"abe es ein $c \in B$ mit $k_c = 0$. Hieraus
 folgte der Widerspruch
 $$ 1 = fc = \sum_{b \in B}k_b b^*c = k_c = 0, $$
 so dass $B^*$ kein Erzeugendensystem von $V^*$ ist.
 \medskip
       Ist $V$ ein Linksvektorraum\index{Linksvektorraum}{} \"uber dem K\"orper
 $K$, so definiert
 man ganz entsprechend wie im Falle der Rechtsvektorr\"aume den
 Dualraum von $V$. Der Dualraum von $V$ ist dann ein
 Rechtsvektorraum \"uber $K$. Die bislang bewiesenen S\"atze
 behalten nat\"urlich \emph{mutatis mutandis} ihre G\"ultigkeit.
 Ist $V$ ein Rechtsvektorraum, so bezeichnen wir mit $V^{**}$ den
 Dualraum von $V^*$. Der Raum $V^{**}$ ist dann wieder ein
 Rechtsvektorraum \"uber $K$. Man nennt ihn den zu $V$
 \emph{bidualen Raum}.\index{bidualer Raum}{} Der n\"achste Satz beschreibt den
 Zusammenhang zwischen $V$ und seinem Bidualen.
 \medskip\noindent
 {\bf 5.15. Satz.} {\it Es sei $V$ ein Rechtsvektorraum \"uber dem
 K\"orper $K$. Definiert man f\"ur $v \in V$ die Abbildung
 $v^\varphi$ von $V^*$ in $K$ verm\"oge $fv^\varphi := fv$ f\"ur alle
 $f \in V^*$, so ist $\varphi$ ein Monomorphismus von $V$ in
 $V^{**}$. Genau dann ist $\varphi$ surjektiv, wenn der Rang von
 $V$ endlich ist. }
 \smallskip
       Beweis. Routinerechnungen zeigen, dass $\varphi$ ein
 Homomorphismus von $V$ in $V^{**}$ ist. Es bleibt zu zeigen, dass
 $\varphi$ injektiv ist. Dazu sei $0 \neq v \in V$. Es gibt dann
 eine Hyperebene $H$ mit $V = vK \oplus H$. Definiere $f \in V^*$
 durch $fv := 1$ und $fh := 0$ f\"ur alle $h \in H$. Es folgt $1 = fv
 = fv^\varphi$, so dass $v^\varphi \neq 0$ ist. Somit ist
 $\Kern(\varphi) = \{0\}$.
 \par
      Ist $V$ endlichen Ranges, so folgt mit 5.14, dass
 $$ \Rg_K(V) = \Rg^K(V^*) = \Rg_K(V^{**}) $$
 ist. Hieraus folgt mittels 3.1, dass $V^\varphi = V^{**}$ ist. In diesem Falle
 ist $\varphi$ surjektiv, wie behauptet.
 \par
       Es sei $V$ nicht endlichen Ranges und $B$ sei eine Basis von $V$.
 Nach 5.14 und 5.13 ist $U := \sum_{b \in B} Kb^*$ ein
 echter Teilraum von $V^*$. Es sei $f \in V^* - U$. Weil nat\"urlich
 auch $L^K (V^*)$ ein projektiver Verband ist, gibt es  eine
 Hyperebene $H$ mit $f \notin H$ und $U \leq H$. Es gibt daher ein
 $g \in V^**$ mit $fg = 1$ und Kern$(g) = H$. W\"are nun $g =
 v^\varphi$ mit einem $v \in V$, so w\"are
 $$ 0 = b^* g = b^* v^\varphi = b^* v $$
 f\"ur alle $b \in B$. Nun w\"are aber $v = \sum_{c \in B} ck_c$. Wendete man
 hierauf die Abbildung $b^*$ an, so folgte $0 = k_b$, so dass $v = 0$ w\"are im
 Widerspruch zu $g  \neq 0$. Damit ist alles bewiesen.
 \medskip
       Ist der Rang von $V$ nicht endlich, so besteht eine Unsymmetrie zwischen
 $V$ und $V^*$. Um diese Unsymmetrie zu beseitigen, studieren wir nach dem
 Vorgange von Artin\index{Artin, E.}{} eine etwas allgemeinere Situation, in der,
 da sie allgemeiner ist, nicht mehr alle S\"atze gelten, die f\"ur $V$ und $V^*$
 richtig sind. Diesen Nachteil nehmen wir jedoch in Kauf, da er durch den Gewinn
 der Symmetrie fast wett gemacht wird. Diese  Situation wird wie folgt
 beschrieben. Es sei $K$ ein K\"orper und $W$ sei ein Links- sowie $V$ ein
 Rechtsvektorraum \"uber $K$. Ferner sei $f$ eine Abbildung des cartesischen
 Produktes von $W$ mit $V$ in $K$. Das Bild von $(w,v)$ unter $f$ werde mit $wv$
 bezeichnet. Wir sagen, dass $W$ mit $V$
 \emph{gekoppelt}\index{gekoppelte R\"aume}{} sei, bzw.\ dass $(W,V,f)$ ein
 \emph{duales Raumpaar}\index{duales Raumpaar}{} bez.\ des
 \emph{Skalarproduktes}\index{Skalarprodukt}{} $f$ sei,
 falls gilt:
 \smallskip\noindent
 1) Es ist $(w + w')v = wv + w'v$ f\"ur alle $w$, $w' \in W$ und alle $v \in V$.
 \par\noindent
 2) Es ist $w(v + v') = wv + wv'$ f\"ur alle $w \in W$ und alle $v$, $v' \in V$.
 \par\noindent
 3) Es ist $(kw)v = k(wv)$ f\"ur alle $k \in K$, alle $w \in W$ und alle
 $v \in V$.
 \par\noindent
 4) Es ist $w(vk) = (wv)k$ f\"ur alle $k \in K$, alle $w \in W$  und alle
 $v \in V$.
 \par\noindent
 5) Ist $w \in W$ und $wV = \{0\}$, so ist $w = 0$.
 \par\noindent
 6) Ist $v \in V$ und $Wv= \{0\}$, so ist $v = 0$.
 \smallskip
       F\"ur $X \in L_K(V)$ setzen wir
 $$ X^\perp := \bigl\{w \mid w \in W, wX = \{0\}\bigr\} $$
 und f\"ur $Y \in L^K(W)$ setzen wir
 $$y^\perp := \bigl\{v \mid v \in V, Yv = \{0\}\bigr\}.$$
 \par
       Sind $X$, $Y \in L_K(V)$ und gilt $X \leq Y$, so gilt $Y^\perp \leq
 X^\bot$. Ferner gilt $X \leq X^{\bot \top}$. Sind $X$, $Y \in
 L^K(W)$ und gilt $X \leq Y$, so gilt $Y^\top \leq X^\top$ wie auch
 $X \leq X^{\top \bot}$. Dies folgt unmittelbar aus den
 Definitionen.
 \medskip\noindent
 {\bf 5.16. Satz.} {\it Es sei $K$ ein K\"orper, $W$ sei ein
 $K$-Links- und $V$ ein $K$-Rechts\-vek\-tor\-raum. \"Uberdies seien $W$
 und $V$ gekoppelt. Ist $X \in L_K(V)$ und $Y \in L^K(W)$, so gilt:
 \item{a)} Es gibt genau einen Monomorphismus $\varphi$ von $X^\bot$ in $(V/X)^*$ mit
 $$ y^\varphi (v + X) = yv $$
 f\"ur alle $y \in X^\bot$ und alle $v \in V$.
 \item{a$'$)} Es gibt genau einen Monomorphismus $\psi$ von $Y^\top$ in $(W/Y)^*$ mit
 $$ (w + Y) x^\varphi = wx $$
 f\"ur alle $w \in W$ und alle $x \in Y^\top$.
 \item{b)} Ist $w \in W$, so definieren wir $w^\zeta$ durch $w^\zeta := wx$
 f\"ur alle $x \in X$. Dann ist $\zeta$ ein Homomorphismus von $W$ in $X^*$
 und es gilt $\Kern(\zeta)= X^\bot$. Ist der Rang von $X$ endlich, so ist
 $\zeta$ surjektiv.
 \item{b$'$)} Ist $v \in V$, so definieren wir $v^\xi$ durch $yv^\xi
 := yv$ f\"ur alle $y \in Y$. Dann ist $\xi$ ein Homomorphismus von
 $V$ in $Y^*$ und es gilt $\Kern(\xi) = y^\top$. Ist der Rang von
 $Y$ endlich, so ist $\xi$ surjektiv.\par}
 \smallskip
       Beweis. Aus Symmetriegr\"unden gen\"ugt es, a) und b) zu beweisen.
 \par
       a) Ist $y \in X^\bot$ und $v + X = v' + X$, so folgt $yv  = yv'$.
 Definiert man daher, $y^\varphi$ durch $y^\varphi(v + X) := yv$,
 so ist $y^\varphi$ wohldefiniert. Es folgt, dass $\varphi$ eine
 Abbildung von $X^\bot$ in $(V/X)^*$ ist. Triviale Rechnungen
 zeigen, dass $\varphi$ sogar ein Homomorphismus ist. Es sei
 $y^\varphi = 0$. Dann ist $yV = \{0\}$ und daher $y = 0$, so dass
 $\varphi$ ein Monomorphismus ist.
 \par
       b) Ist $w \in W$, so definieren wir $w^\zeta$ durch
 $w^\zeta x := wx$ f\"ur alle $ x \in X$. Triviale Rechnungen
 zeigen, dass $\zeta$ ein Homomorphismus von $W$ in $X^*$ ist.
 Ist $w^\zeta = 0$, so ist $wX = \{0\}$, so dass $w \in X^\bot$
 gilt.
 \par
       Der Rang von $X$ sei endlich. Dann ist nach dem gerade Bewiesenen
 $$ \Rg^K (W/X^\bot) \leq \Rg^K(X^*) = \Rg_K(X). $$
 Nach a$'$) gibt es einen Monomorphismus von $X^{\bot\top}
 $ in $(W/X^\bot)^*$. Weil $X \leq X^{\bot\top}$ gilt, ist also
 $$ \Rg_K(X) \leq \Rg_K(X^{\bot\top}) \leq \Rg_K(W/X^\bot)^* \leq \Rg_K(X).$$
 Somit ist $X = X^{\bot\top}$. Es folgt weiter, dass
 $\Rg^K(W/X^\bot) = \Rg^K(X^*)$ ist. Folglich ist $\zeta$ in
 diesem Falle auch surjektiv.
 \par
       Damit ist alles bewiesen.
 \medskip
       Im Beweise von 5.16 haben wir mehr bewiesen als im Satz
 formuliert. Darauf werden wir gleich zur\"uckkommen. Zuvor zeigen
 wir jedoch, dass die Theorie der gekoppelten R\"aume im Falle von
 R\"aumen endlichen Ranges nichts Neues liefert.
 \medskip\noindent
 {\bf 5.17. Satz.} {\it Es sei $K$ ein K\"orper, $W$ sei ein $K$-Links-
 und $V$ ein $K$-Rechts\-vek\-tor\-raum. \"Uberdies seien $W$ und $V$
 gekoppelt. Genau dann hat $V$ endlichen Rang, wenn $W$ endlichen
 Rang hat. Hat einer der beiden R\"aume endlichen Rang, so haben
 also beide endlichen Rang und es gilt:
 \item{a)} Ist $w \in W$ und definiert man $w^\zeta$ durch $w^\zeta v := wv$
 f\"ur alle $v \in V$, so ist $\zeta$ ein Isomorphismus von $W$ auf $V^*$.
 \item{b)} Ist $v \in V$ und definiert man $v^\xi$ durch $wv^\xi := wv$ f\"ur
 alle $w \in W$, so ist $\xi$ ein Isomorphismus von $V$ auf $W^*$.\par}
 \smallskip
       Beweis. Es sei $V$ endlichen Ranges. Dann folgt mit 5.16 b), dass
 $\zeta$ ein Epimorphismus von $W$ auf $V^*$ ist und dass
 $ \Kern(\zeta) = V^\bot = \{0\} $
 gilt. Also hat auch $W$ endlichen Rang und es gilt a).
 \par
       Hat $W$ endlichen Rang, so folgt mit 5.16 b$'$), dass $V$ endlichen
 Rang hat und dass b) gilt. Damit ist alles bewiesen.
 \medskip\noindent
 {\bf 5.18. Satz.} {\it Es sei $K$ ein K\"orper und $W$ und $V$ seien
 gekoppelte Vektorr\"aume \"uber $K$. Dann gilt:
 \item{a)} Sind $X$, $Y \in L_K(V)$ und gilt $X \leq Y$, so ist
 $Y^\bot \leq X^\bot$.
 \item{b)} Sind $X$, $Y \in L^K(W)$ und gilt $X \leq Y$, so ist
 $Y^\top \leq X^\top$.
 \item{c)} Es ist $X \leq X^{\bot\top}$ f\"ur alle $X \in L_K(V)$. Hat $X$
 endlichen Rang, so gilt $X = X^{\bot\top}$.
 \item{d)} Es ist $X \leq X^{\top\bot}$ f\"ur alle $X \in L^K(W)$. Hat $X$
 endlichen Rang, so gilt $X = X^{\top\bot}$.
 \item{e)} Es ist $X^{\bot\top\bot} = X^\bot$ f\"ur alle $X \in L_K(V)$.
 \item{f)} Es ist $X^{\top\bot\top} = X^\top$ f\"ur alle $X \in L^K(W)$.}
 \smallskip
       Beweis. a), b), c) und d) sind entweder schon als trivial erkannt
 oder bereits bewiesen, n\"amlich beim Beweise von 5.16. Es bleiben
 e) und f) zu beweisen.
 \par
       Es sei $U \in L_K(V)$. Dann ist $U \leq U^{\bot\top}$ nach c) und
 folglich $U^{\bot\top\bot} \leq U^\bot$ nach a). Andererseits ist
 $U^\bot \in L^K(W)$, so dass nach d) gilt, dass $U^\bot \leq
 U^{\bot\top\bot}$ ist. Damit ist e) bewiesen. Die Aussage f)
 beweist sich analog.
 \medskip
       Das Paar $(\bot, \top)$ ist eine
 \emph{Galoisverbindung}\index{Galoisverbindung}{} von $L_K(V)$ mit $L^K(W)$.
 \medskip\noindent
 {\bf 5.19. Satz.} {\it Ist $V$ ein Rechtsvektorraum endlichen Ranges \"uber dem
 K\"orper $K$ und ist $V^*$ der Dualraum von $V$, so sind $V^*$ und $V$ gekoppelt
 und $\bot$ ist ein Isomorphismus des Verbandes $(L_K(V)^d, \leq ^d)$ auf den
 Verband $(L^K(V^*), \leq)$.}
 \smallskip
       Beweis. Nach 5.18 a) und b) ist $\bot$ eine ordnungserhaltende Abbildung
 von $(L_K(V)^d, \leq^d)$ in $(L^K(V^*), \leq)$ und $\top$ ist eine
 ordnungserhaltende Abbildung von $(L^K(V^*), \leq)$ in $(L_K(V)^d, \leq^d)$.
 Nach 5.18 c) ist ferner $\bot\top$ die Identit\"at auf $L_K(V)$. Weil der Rang
 von $V$ endlich ist, ist es auch der Rang von $V^*$. Daher folgt mit 5.18 d),
 dass $\top\bot$ die Identit\"at auf $L^K(V^*)$ ist. Also ist $\bot$ bijektiv,
 wie behauptet.
 \medskip
       Mit diesem Satz haben wir eine Beschreibung von $L_K(V)^d$ gewonnen, die
 diesen Verband der Behandlung mittels Methoden der li\-ne\-a\-ren
 Algebra\index{lineare Algebra}{} zu\-g\"ang\-lich macht.
 Puristen\index{Purist}{} --- und in diesem Zusammenhang sind wir Puristen ---
 werden bem\"angeln, dass der duale Verband mit Hilfe eines
 Linksvektorraumes\index{Linksvektorraum}{} dargestellt wird. Man w\"unscht sich
 auch f\"ur ihn eine Beschreibung durch einen
 Rechtsvektorraum.\index{Rechtsvektorraum}{} Eine solche Beschreibung ist nun
 rasch gefunden.  Es sei $K$ ein K\"orper. Wir definieren auf $K$ eine
 Multiplikation $\circ$ durch $a \circ b := ba$ f\"ur alle $a, b \in K$. Dann ist
 auch $(K, +, \circ)$ ein K\"orper und die Identit\"at ist ein
 Antiisomorphismus\index{Antiisomorphismus}{} von $K$ auf jenen K\"orper,
 den wir mit $K_\circ$ be\-zeich\-nen. Ist nun $V$ ein Linksvektorraum
 \"uber $K$, so wird $V$ zu einen Rechtsvektorraum \"uber $K_\circ$
 durch die Vorschrift $v \circ a:= av$. Offenbar gilt $L^K(V) = L_{K_\circ}(V)$.
 Wendet man dies nun auf die obige Situation an, so erh\"alt man in
 $L_{K_\circ}(V^*)$ eine Beschreibung von $L_K(V)^d$ durch einen Rechtsvektorraum
 \"uber dem zu $K$ antiisomorphen K\"orper $K_\circ$. Ist $K$
 kommutativ,\index{kommutativer K\"orper}{} so ist $K$
 nat\"urlich gleich $K_\circ$. Daher sind in diesem Falle auch
 $L_K(V)$ und $L_K(V)^d$ gleich. Im allgemeinen ist dies jedoch
 nicht der Fall. Mehr dar\"uber im n\"achsten Kapitel.
 \par
       Eine interessante Folgerung aus 5.19 ist der n\"achste Satz.
 \medskip\noindent
 {\bf 5.20. Satz.} {\it Es sei $V$ ein Vektorraum endlichen Ranges
 \"uber dem K\"orper $K$. Ist $U \in L_K(V)$, so ist
 $$\Rg_K(V) = \Rg_K(U) + \Rg^K(U^\bot).$$ Dabei ist $\bot$ als
 Abbildung von $L_K(V)$ in $L^K(V^*)$ aufzufassen.}
 \smallskip
       Beweis. Weil $\bot$ nach 5.19 ein Isomorphismus von $L_K(V)^d$ auf
 $L^K(V^*)$ ist, ist $\Ko_K(U) = \Rg^K(U^\bot)$. Hieraus folgt die Behauptung.
 \medskip\noindent
 {\bf 5.21. Satz.} {\it Es sei $K$ ein K\"orper und $W$ und $V$ seien gekoppelte
 Vektorr\"aume \"uber $K$. Es sei $\Xi$ die Menge aller
 Unterr\"aume der Form $X^\bot$ von $W$, wobei $X$ ein Teilraum
 endlichen Ranges von $V$ ist. Dann gilt:
 \item{a)} Ist $Y \in \Xi$ und ist $Z$ ein Teilraum von $W$ mit $Y \leq Z$, so ist
 $Z \in \Xi$.
 \item{b)} Sind $Y$, $Z \in \Xi$, so ist $Y \cap Z \in \Xi$.\par
 \noindent  Entsprechendes gilt,
 wenn man die Rollen von $W$ und $V$ vertauscht.\par}
 \smallskip
       Beweis. a) Es gibt ein $X$ endlichen Ranges in $L_K(V)$ mit $Y =
 X^\bot$. Wir definieren eine Koppelung, der wir diesmal einen
 Namen, n\"amlich $\pi$, geben, von $W/Y$ mit $X$ durch $(w +
 Y)x:=wx$ f\"ur alle $w \in W$ und alle $x \in X$. Wegen $Y =
 X^\bot$ ist $\pi$ wohldefiniert. Wir zeigen, dass $\pi$ eine
 Koppelung ist. Dazu sind nur die Eigenschaften 5) und 6)
 nachzuweisen, da die \"ubrigen sich von selbst verstehen. Es sei
 $(w + Y) X= \{0\}$. Dann ist $wX = \{0\}$ und daher $w \in X^\bot
 = Y$. Ist $x \in X$ und $(W/Y)x = \{0\}$, so ist $Wx = \{0\}$ und
 daher $x = 0$. Somit ist $\pi$ tats\"achlich eine Koppelung. Die
 durch $\pi$ definierten Abbildungen von Unterr\"aumen des einen
 Raumes auf ihren annullierenden Unterraum des anderen bezeichnen
 wir mit $\bot(\pi)$ bzw. $\top(\pi)$.
 \par
       Es sei nun $Z/Y$ ein Teilraum von $W/Y$. Dann ist
 $$\eqalign{
   (Z/Y)^{\top(\pi)} &= \bigl\{x \mid x \in X,\ (Z/Y)x = \{0\}\bigr\} \cr
		     &= \bigl\{x \mid x \in X,\ Zx = \{0\}\bigr\}     \cr
		     &= Z^\top. \cr} $$
 \par
      Es sei andererseits $U$ ein Teilraum von $X$. Dann ist
 $$\eqalign{
     U^{\bot(\pi)} &= \bigl\{w + Y \mid w \in W,\ (w + Y)U = \{0\}\bigr\} \cr
		   &= \bigl\{w + Y \mid w \in W,\ wU = \{0\}\bigr\}       \cr
		   &= U^\bot/Y. \cr} $$
 \par
       Weil $X$ endlichen Rang hat, sind die Abbildungen $\bot(\pi)$ und
 $\top(\pi)$ Bijektionen. Daher gilt
 $$ Z/Y = (Z/Y)^{\top(\pi)\bot(\pi)} = (Z^\top)^{\bot(\pi)} = Z^{\top\bot}/Y.$$
 Hieraus folgt $$Z = Z^{\top \bot}$$ f\"ur alle Teilr\"aume $Z$ von
 $W$, die $Y$ umfassen. Weil $X$ endlichen Rang hat, gilt aber auch
 $X = X^{\bot \top}$. Somit ist die Einschr\"ankung von $\bot$ auf
 $L_K(X)$ eine Bijektion dieser Menge auf die Menge der Teilr\"aume
 von $W$, die $Y$ enthalten. Daher gilt a).
 \par
      b) Es gibt Unterr\"aume $A$ und $B$ endlichen Ranges von $V$ mit
 $Y = A^\bot$ und $Z = B^\bot$. Nun ist aber
 $$ Y \cap Z = A^\bot \cap B^\bot = (A + B)^\bot. $$
 Weil auch $A + B$ endlichen Rang hat, ist demnach $Y \cap Z \in \Xi$.
 \medskip
       Die n\"achsten vier S\"atze formulieren wir nur f\"ur den Fall
 eines Rechtsvektorraumes und seines Dualraums. Jegliche
 Ver\-all\-ge\-mei\-ne\-rung w\"urde technisch aufwendig und w\"are nur
 schwer zu merken. Falls der Leser allgemeinere Situationen
 antrifft, wird es ihm nicht schwer fallen, die entsprechenden
 S\"atze aus den hier wiedergegebenen herauszupr\"aparieren.
 \par
 Ist $V$ ein Linksvektorraum \"uber $K$ und ist $\sigma$ eine
 lineare Abbildung von $V$ in einen weiteren $K$-Linksvektorraum,
 so bezeichnen wir das Bild von $v \in V$ unter $\sigma$ mit
 $v\sigma$.
 \medskip\noindent
 {\bf 5.22. Satz.} {\it Es seien $V$ und $W$ Rechtsvektorr\"aume \"uber
 dem K\"orper $K$ und $V^*$ und $W^*$ seien ihre Dualr\"aume. Ist
 $\varphi \in \Hom_K(V,W)$, so gibt es genau ein $\varphi^* \in
 \Hom^K(W^*,V^*)$ mit $(f \varphi^*)v = f(\varphi v)$ f\"ur alle $v \in V$ und
 alle $f \in W^*$. Die Abbildung $*$ ist ein
 Monomorphismus von $\Hom_K(V,W)$ in $\Hom^K(W^*, V^*)$. Haben $V$
 und $W$ beide endlichen Rang, so ist $*$ ein Isomorphismus der
 abelschen Gruppe $\Hom_K(V,W)$ auf $\Hom^K(W^*,V^*)$. Ist $V$
 endlichen Ranges und gilt $V = W$, so ist $*$ ein
 Ringisomorphismus von $\End_K(V)$ auf $\End^K(V^*)$.}
 \smallskip
       Der Beweis dieses Satzes sei dem Leser als \"Ubungsaufgabe
 \"uber\-las\-sen.
       Man nennt $\varphi^*$ die zu $\varphi$ \emph{duale
 Abbildung}.\index{duale Abbildung}{}
 \medskip\noindent
 {\bf 5.23. Satz.} {\it Es seien $V$ und $W$ Rechtsvektorr\"aume \"uber
 dem K\"orper $K$. Ist $\varphi \in \Hom_K(V,W)$ und ist $\varphi^*$
 die zu $\varphi$ duale Abbildung, so ist}
 $$ \Kern(\varphi) = (W^*\varphi^*)^\top\ \ \ \hbox{\it und}\ \ \ 
                  \Kern(\varphi^*) = (\varphi V)^\bot. $$
 \par
       Beweis. Es sei $v \in \Kern(\varphi)$. Dann ist
 $$ (f\varphi^*) v= f(\varphi v) = 0 $$
 f\"ur alle $f \in W^*$. Daher ist
 $v \in (W^* \varphi^*)^\top$. Es sei umgekehrt $v \in
 (W^*\varphi^*)^\top$. Dann ist
 $$ f(\varphi v) = (f\varphi^*)v = 0 $$
 f\"ur alle $f \in W^*$. Hieraus folgt $\varphi v = 0$ und damit
 $v \in \Kern(\varphi)$. Also ist $\Kern(\varphi) = (W^*\varphi^*)^\top$.
 \par
       Es sei $f \in \Kern(\varphi^*)$. Dann ist
 $$ f(\varphi v) = (f\varphi^*)v = 0 $$
 f\"ur alle $v \in V$, so dass $f \in (\varphi V)^\bot$ gilt. Ist umgekehrt
 $f \in (\varphi V)^\bot$, so ist
 $$(f \varphi^*)v = f(\varphi v) = 0$$
 f\"ur alle $v \in V$. Daher ist $f\varphi^* = 0$, so dass
 $f \in \Kern(\varphi^*)$ gilt. Demnach ist $\Kern(\varphi^*) =
 (\varphi V)^\bot$. Damit ist alles bewiesen.
 \medskip\noindent
 {\bf 5.24. Satz.} {\it Es seien $V$ und $W$ Rechtsvektorr\"aume \"uber
 dem K\"orper $K$. Ist $\varphi \in \Hom_K(V,W)$ und ist $\varphi^*$
 die zu $\varphi$ duale Abbildung, so gilt:
 \item{a)} Genau dann ist $\varphi$ surjektiv, wenn $\varphi^*$ injektiv ist.
 \item{b)} Genau dann ist $\varphi$ injektiv, wenn $\varphi^*$ surjektiv ist.
 \item{c)} Genau dann ist $\varphi$ bijektiv, wenn $\varphi^*$ bijektiv ist.}
 \smallskip
       Beweis. a) Nach 5.23 ist $\Kern(\varphi^*) = (\varphi V)^\bot$. Ist nun
 $\varphi^*$ injektiv, so ist $(\varphi V)^{\bot \top} = \{0\}^\top = W$. Mit
 5.18 c) folgt hieraus $\varphi V = W$, so dass $\varphi$ surjektiv ist. Ist
 umgekehrt $\varphi$ surjektiv, so ist $\Kern(\varphi^*) = W^\bot = \{0\}$, 
 also $\varphi^*$ injektiv.
 \par
       b) Nach 5.23 ist $\Kern(\varphi) = (W^*\varphi^*)^\top$. Ist
 $\varphi^*$ surjektiv, so ist also $\Kern(\varphi) =
 (V^*)^\top = \{0\}$, so dass $\varphi$ injektiv ist.
 \par
       Es sei $\varphi$ schlie\ss lich injektiv. Ferner sei $B$ eine Basis
 von $V$. Weil $\varphi$ injektiv ist, ist $\{\varphi b \mid b \in B\}$
 eine Menge von linear unabh\"angigen Vektoren von $W$. Es gibt
 somit einen Teilraum $U$ von $W$ mit
 $$ W = U \oplus \bigoplus_{b \in B} (\varphi b) K.$$
 \par
       Ist $g \in V^*$, so definieren wir $f \in W^*$ durch $f(\varphi b) := gb$
 f\"ur alle $b \in B$ und $fu := 0$ f\"ur alle $u \in U$. Es
 folgt, dass $f\varphi^* = g$ ist. Somit ist $\varphi^*$ surjektiv.
 \par
       c) folgt aus a) und b).
 \medskip\noindent
 {\bf 5.25. Satz.} {\it Es seien $V$ und $W$ Rechtsvektorr\"aume \"uber $K$.
 Ferner sei $\varphi \in \Hom_K(V,W)$ und $\varphi^*$ sei die zu $\varphi$ duale
 Abbildung. Genau dann ist $\varphi V$ endlich erzeugt, wenn $W^* \varphi^*$
 endlich erzeugt ist. Ist einer dieser beiden R\"aume endlich erzeugt, so ist
 $\Rg(\varphi) = \Rg(\varphi^*)$. Ist $V$ oder $W$ endlich erzeugt, so sind
 $\varphi V$ und $W^* \varphi^*$ endlich erzeugt.}
 \smallskip
       Beweis. $W$ und $W^*$ sind gekoppelte R\"aume und $\varphi V$ ist
 ein Teilraum von $W$. Ist nun $\varphi V$ endlich erzeugt, so sind
 nach 5.16 b) die R\"aume $(\varphi V)^*$ und $W^*/(\varphi
 V)^\bot$ isomorph. Nach 5.14 sind daher die R\"ange von $\varphi
 V$ und $W^*/(\varphi V)^\bot$ gleich. Nach 5.21 ist $(\varphi
 V)^\bot$ der Kern von $\varphi^*$, so dass $W^* \varphi^*$ endlich
 erzeugt ist und die R\"ange von $\varphi$ und $\varphi^*$ gleich sind.
 \par
       Es sei $W^* \varphi^*$ endlich erzeugt. Nach 5.18 d) ist dann
 $$ W^* \varphi^* = (W^* \varphi^*)^{\top \bot}. $$
 Nach 5.16 a) ist daher $W^* \varphi^*$ zu $(v/(W^* \varphi^*)^\top)^*$
 isomorph. Nach 5.21 ist daher $(V/\Kern(\varphi))^*$ zu $W^* \varphi^*$
 isomorph. Mittels 5.14 folgt, dass $\varphi V$ endlich erzeugt ist. Dann ist
 aber, wie wir gesehen haben, $\Rg(\varphi) = \Rg(\varphi^*)$.
 \par
       Ist $V$ oder $W$ endlich erzeugt, so ist in beiden F\"allen auch
 $\varphi V$ endlich erzeugt.
 \medskip
       Die Situation von 5.25 verallgemeinernd formulieren wir
 \medskip\noindent
 {\bf 5.26. Satz.} {\it Es sei $K$ ein K\"orper und $W$ und $V$ seien gekoppelte
 Vektorr\"aume \"uber $K$. Es sei ferner $R$ ein Ring, $W$ sei ein Rechtsmodul
 und $V$ ein Linksmodul \"uber $R$. Wir setzen zudem voraus, dass die folgenden
 Rechenregeln gelten:
 \item{a)} Es ist $(wr)v=w(rv)$ f\"ur alle $r \in R$, alle $w \in W$ und alle
 $v \in V$.
 \item{b)} Es ist $(kw)r = k(wr)$ f\"ur alle $k \in K$, alle $w \in W$ und alle
 $r \in R$.
 \item{c)} Es ist $(rv)k = r(vk)$ f\"ur alle $r \in R$, alle $v \in
 V$ und alle $k \in K$.

 \noindent Ist $r \in R$, so ist $rV$ genau dann
 endlich erzeugt, wenn $Wr$ endlich erzeugt ist. Ist dies der Fall,
 so gilt $\Rg_K(rV) = \Rg^K(Wr)$.}
 \smallskip
       Beweis. Setze $Y := (rV)^\bot$. Ist dann $y \in Y$ und $v \in V$,
 so folgt $(yr)v = y(rv) = 0$. Wegen $V^\bot = \{0\}$ ist daher $yr
 = 0$. Ist andererseits $w \in W$ und $wr = 0$, so folgt $0 = (wr)v
 = w(rv)$ f\"ur alle $v \in V$, so dass $w \in Y$ gilt. Es folgt,
 dass $Wr$ zu $W/Y$ isomorph ist. Nach Fr\"uherem ist $Wr$ daher
 endlich erzeugt und es gilt $\Rg^K(Wr) = \Rg_K(rV)$. Aus
 Symmetriegr\"unden ist damit alles gezeigt.
 \medskip
       Dieser Satz wird sp\"ater noch eine Rolle spielen.

\bigskip

\mysection{6. Homogene, vollst\"andig reduzible Ringe}%

\smallskip

\noindent
 In diesem Abschnitt greifen wir, wie versprochen, die Frage nach
 der Struktur der homogenen, vollst\"andig reduziblen Ringe wieder
 auf. Kennen wir diese, so kennen wir auf Grund von Satz 4.6 alle
 vollst\"andig reduziblen Ringe, da die ringtheoretische direkte
 Summe von homogenen, vollst\"andig reduziblen Ringen ein
 vollst\"andig re\-du\-zib\-ler Ring ist.
 \par
       Wir definieren das \emph{Jacobson-Radikal}\index{Jacobson-Radikal}{}
 $J(R)$ eines Ringes $R$ als den Schnitt \"uber alle regul\"aren, maximalen
 Rechtsideale von $R$. Ist $M$ ein irreduzibler $R$-Rechtsmodul und ist $0 \neq
 x \in M$, so ist $O(x)$ nach 4.2 ein regul\"ares, maximales
 Rechtsideal. Daher ist $J(R) \subseteq O(x)$, so dass $xj = 0$ ist
 f\"ur alle $j \in J(R)$. Somit gilt $M J (R) = \{0\}$ f\"ur alle
 irreduziblen $R$-Rechtsmoduln $M$. Ist $I$ ein regul\"ares,
 maximales Rechtsideal von $R$, so ist $R/I$ nach 4.2 ein
 irreduzibler $R$-Modul. Nach der gerade gemachten Bemerkung ist
 daher $RJ(R) \subseteq I$. Da dies f\"ur alle regul\"aren,
 maximalen Ideale gilt, folgt $RJ(R) \subseteq J(R)$, so dass
 $J(R)$ ein zweiseitiges Ideal von $R$ ist.
 \medskip\noindent
 {\bf 6.1. Satz.} {\it Ist $R$ ein vollst\"andig reduzibler Ring, so
 ist $J(R) \neq R$ und $J(R)^2 = \{0\}$.}
 \smallskip
       Beweis. Es sei $I$ ein minimales Rechtsideal, welches in $J(R)$
 ent\-hal\-ten ist. Dann ist $I$ ein irreduzibler $R$-Rechtsmodul.
 Daher ist $IR \neq \{0\} = IJ(R)$. Hieraus folgt $R \neq J(R)$.
 Weil $J(R)$ Summe von minimalen Rechtsidealen von $R$ ist, folgt,
 dass auch $J(R)^2 = \{0\}$ gilt.
 \medskip\noindent
 {\bf 6.2. Satz.} {\it Ist $R$ ein vollst\"andig reduzibler Ring und ist $I$ ein
 minimales Rechtsideal von $R$, so gilt genau dann $I \subseteq J(R)$, wenn
 $I^2 = \{0\}$ ist.}
 \smallskip
       Beweis. Ist $I \subseteq J(R)$, so ist $I^2 \subseteq J(R)^2 = \{0\}$.
 \par
       Es sei $I \not\subseteq J(R)$. Es gibt dann ein regul\"ares,
 maximales Rechts\-i\-de\-al $M$ von $R$, welches $I$ nicht enth\"alt. Es
 folgt $R = I \oplus M$ als $R$-Rechtsmodul. Weil $M$ regul\"ar
 ist, gibt es ein $a \in R$ mit $x - ax \in M$ f\"ur alle $x \in
 R$. Es ist $a = i + m$ mit $i \in I$ und $m \in M$. Es folgt, dass
 auch $x - ix \in M$ gilt f\"ur alle $x \in R$. Insbesondere ist
 also
 $$ i - i^2 \in I \cap M = \{0\}. $$
 Also ist $i^2 = i$.  Nun ist aber $i \neq 0$, da andernfalls $R = M$ w\"are.
 Somit ist $I^2 \neq \{0\}$.
 \medskip
       Im Verlaufe des Beweises des letzten Satzes haben wir gezeigt, dass das
 minimale Ideal $I$ ein von Null verschiedenes Idempotent\index{Idempotent}{}
 enth\"alt, wenn $I$ nicht im Jacobson-Radikal des Ringes liegt.
 Allgemeiner gilt:
 \medskip\noindent
 {\bf 6.3. Satz} {\it Es sei $R$ ein Ring und $I$ sei ein minimales
 Rechtsideal von $R$. Ist $I^2 \neq \{0\}$, so gibt es ein von $0$
 verschiedenes Idempotent $i \in I$. Insbesondere ist $I = iR$.}
 \smallskip
       Beweis. Wegen $I^2 \neq \{0\}$ gibt es ein $x \in I$ mit $xI \neq \{0\}$.
 Weil $xI$ ein in $I$ enthaltenes Rechtsideal ist, folgt wegen der Minimalit\"at
 von $I$ daher, dass $xI = I$ ist. Es gibt also ein $i \in I$ mit $xi = x$. Es
 folgt $x(i^2 - i) = 0$. Setze $I' := \{y \mid y \in I, xy = 0\}$. Dann ist $I'$
 ein Rechtsideal, welches in $I$ enthalten ist. Nun ist $x \neq 0$ und daher
 $i \in I'$. Folglich ist $I' = \{0\}$. Wegen $i^2 - i \in I'$ ist also
 $i^2 = i$.
 \par
       Schlie\ss lich ist $i = i^2 \in iR \subseteq I$ und somit $iR = I$
 auf Grund der Minimalit\"at von $I$.
 \medskip\noindent
 {\bf 6.4. Satz.} {\it Es sei $R$ ein homogener, vollst\"andig reduzibler Ring.
 Fer\-ner sei $S$ die Summe \"uber alle minimalen Ideale $I$ von $R$, f\"ur die
 $I^2 \neq \{0\}$ gilt. Dann ist $S$ ein einfacher, vollst\"andig reduzibler Ring
 und es gilt im modultheoretischen Sinne
 $$ R = J(R) \oplus S. $$
 Sind $j$, $j' \in J(R)$ und $s$, $s' \in S$, so ist}
 $$ (j + s)(j' + s') = js' + ss'. $$
 \par
       Beweis. Weil $R$ die Summe seiner minimalen Rechtsideale ist,
 folgt zun\"achst, dass $RJ(R) = \{0\}$ ist. Ferner folgt mit 6.2,
 dass $R$ im modultheoretischen Sinne die direkte Summe von $J(R)$
 und $S$ ist. Weil $S$ ein Rechtsideal ist, ist nat\"urlich $SS \subseteq S$,
 so dass $S$ ein Ring ist.
 \par
       Es sei $I$ ein minimales Rechtsideal von $R$, welches in $S$
 enthalten ist. Ist $0 \neq x \in I$, so gilt
 $$ I = xR \subseteq xJ(R) + xS = xS \subseteq I. $$
 Daher ist $I$ auch ein
 minimales Rechtsideal von $S$. Da $S$ Summe mi\-ni\-ma\-ler Rechtsideale
 von $R$ ist, folgt hieraus schon, dass $S$ vollst\"andig reduzibel ist.
 \par
       Es sei nun $N$ ein von $\{0\}$ verschiedenes, zweiseitiges Ideal
 von $S$. Wegen
 $$ NR \subseteq NJ (R) + NS = NS \subseteq N $$
 ist $N$ ein Rechtsideal von $R$. Es gibt also ein minimales
 Rechtsideal $I$ von $R$, welches in $N$ enthalten ist. Nach der
 Definition von $S$ gilt $I^2 \neq \{0\}$, so dass $I$ nach 6.3 ein
 von Null verschiedenes Idempotent $n$ enth\"alt. Es sei $I'$ ein
 weiteres minimales Rechtsideal von $R$, welches in $S$ enthalten
 ist. Weil $R$ homogen ist, gibt es einen $R$-Isomorphismus
 $\sigma$ von $I$ auf $I'$. Setze $i := \sigma (n)$. Dann ist
 $$0 \neq i = \sigma (n) = \sigma (n^2) = in \in I' \cap N.$$
 (Dieses Argument haben wir schon einmal gefeiert.) Weil
 $I'$, wie wir gesehen haben, auch ein minimales Rechtsideal von
 $S$ ist, folgt $I' \subseteq N$. Benutzt man nun noch einmal, dass
 $S$ die Summe \"uber alle minimalen Rechtsideale von $R$ ist,
 deren Quadrat nicht Null ist, so folgt, dass $N = S$ ist, womit
 die Einfachheit von $S$ nachgewiesen ist.
 \par
       Ist $r \in R$ und $j \in J(R)$, so ist $rj = 0$, wie zuvor schon
 bemerkt. Hieraus folgt schlie\ss lich auch noch die G\"ultigkeit
 der letzten Behauptung des Satzes.
 \medskip
       Ist $S$ ein einfacher, vollst\"andig reduzibler Ring und ist $M$ ein
 $S$-Rechtsmodul mit $MS = M$, so zeigt Satz 6.4, wie man aus $S$ und $M$ einen
 homogenen, vollst\"andig reduziblen Ring $R$ machen kann, so dass $M$ gerade das
 Jacobson-Radikal von $R$ wird.  Man kennt also alle homogenen, vollst\"andig
 reduziblen Ringe, wenn man alle einfachen, vollst\"andig  reduziblen Ringe
 kennt. Um diese besser in den Griff zu kriegen, beweisen wir zun\"achst das
 schursche Lemma,\index{schursches Lemma}{} welches ein \"uberaus n\"utzliches
 Werkzeug ist.
 \medskip\noindent
 {\bf 6.5. Schursches Lemma.} {\it Es sei $R$ ein Ring und $M$ sei ein
 ir\-re\-du\-zib\-ler $R$-Modul. Dann ist $\End_R(M)$ ein K\"orper.}
 \smallskip
       Beweis. Es sei $0 \neq \sigma \in \End_R(M)$. Dann ist $\sigma(M)$
 ein von $\{0\}$ verschiedener Teilmodul von $M$. Somit gilt wegen
 der Irreduzibilit\"at von $M$ die Gleichung $\sigma(M) = M$, so
 dass $\sigma$ surjektiv ist. Ferner gilt $\Kern(\sigma) \neq
 M$ und daher $\Kern(\sigma) = \{0\}$, so dass $\sigma$ auch
 injektiv ist. Folglich ist $\sigma$ ein Automorphismus von $M$, so
 dass $\sigma$ eine Einheit von $\End_R(M)$ ist. Damit ist alles bewiesen.
 \medskip
       Es sei $R$ ein Ring und $M$ sei ein $R$-Rechtsmodul. Man nennt $M$
 {\it treuen\/}\index{treuer Modul}{} $R$-Rechtsmodul, falls aus $r \in R$ und
 $Mr = \{0\}$ folgt, dass $r = 0$ ist. Der Ring $R$ hei\ss t \emph{primitiv},
 \index{primitiver Ring}{}falls er einen treuen, irreduziblen Modul besitzt.
 \medskip\noindent
 {\bf 6.6. Satz.} {\it Es sei $R$ ein primitiver Ring und $M$ sei ein treuer,
 irreduzibler $R$-Rechtsmodul. Setze $K := \End_R(M)$. Dann ist $K$ ein K\"orper
 und $M$ ist ein Linksvektorraum \"uber $K$.  Es bezeichne $E^K(M)$ die Menge der
 Teilr\"aume endlichen Ranges des $K$-Vektorraumes $M$. F\"ur $X \in E^K(M)$ sei
 $O(X) := \{r \mid r \in R,\ Xr = \{0\}\}$. Dann ist $O$ ein Monomorphismus von
 $E^K(V)$ in $L_R(R)$ mit der Eigenschaft, dass f\"ur $X$, $Y \in E^K(V)$
 genau dann $X \leq Y$ gilt, wenn $O(Y) \subseteq O(X)$ ist.}
 \smallskip
       Beweis. Dass $K$ ein K\"orper ist, besagt gerade das schursche Lemma.
 \par
       Aus der Definition von $O$ folgt, dass $X \leq Y$ impliziert, dass
 $O(Y) \subseteq O(X)$ ist. Wir zeigen, dass aus $X < Y$ folgt,
 dass $O(Y) \subset O(X)$ ist. Dazu machen wir Induktion nach dem
 Rang von $X$. Es sei $y \in Y - X$. Ist der Rang von $X$ gleich
 0, so ist $X = \{0\}$. Weil $M$ irreduzibel ist, ist $yR = M$,
 so dass es ein $s \in R$ gibt mit $ys \neq 0$. Also ist $O(Y) \subset R = O(X)$.
 \par
       Der Rang von $X$ sei nun gr\"o\ss er als 0 und $w$ sei ein von
 $0$ verschiedener Vektor in $X$. Es gibt dann einen Teilraum $V$
 von $X$ mit $X = V \oplus Kw$. Wir nehmen nun an, dass $O(X)
 \subseteq O(y)$ gelte, wobei $y$ immer noch das oben gew\"ahlte
 Element ist. Nach Induktionsannahme ist $O(V) \not \subseteq
 O(w)$, da andernfalls $O(V) = O(X)$ w\"are. Daher ist $wO(V)$ ein
 von $\{0\}$ verschiedener Teilmodul von $M$. Folglich gilt $wO(V)
 = M$. Es seien $a$, $b \in O(V)$ und es gelte $wa = wb$. Dann ist
 $$ a - b \in O(w) \cap O(V) = O(X) \subseteq O(y), $$
 so dass $ya = yb$ ist. Folglich wird durch $\kappa (wa) := ya$ ein
 Endomorphismus $\kappa$ des $R$-Rechtsmoduls $wO(V) = M$
 definiert, dh., es gilt $\kappa \in K$. Es folgt $(\kappa w - y)a
 = 0$ f\"ur alle $a \in O(V)$. Nach In\-duk\-ti\-ons\-an\-nah\-me ist daher
 $\kappa w - y \in V$ und daher $y \in V + Kw = X$. Dieser
 Widerspruch zeigt, dass $O(X) \not \subseteq O(y)$ ist. Folglich
 ist $O(Y) \subset O(X)$. Es sei schlie\ss lich $O(Y) \subseteq
 O(X)$. Dann ist
 $$ O(X + Y) = O(X) \cap O(Y) = O(Y). $$
 Wegen $Y \leq X + Y$ ist daher nach dem bereits Bewiesenen $Y = X + Y$, so dass
 $X \leq Y$ gilt. Damit ist alles bewiesen.
 \medskip\noindent
 {\bf 6.7. Dichtesatz.} {\it Es sei $R$ ein primitiver Ring und $M$ sei
 ein irreduzibler Rechtsmodul \"uber $R$. Setze $K := \End_R(M)$.
 Dann ist $K$ ein K\"orper und $M$ ist ein Linksvektorraum \"uber
 $K$. Ist dann $n$ eine nat\"urliche Zahl, sind $b_1, \dots, b_n$
 linear unabh\"angige Elemente des $K$-Vektorraumes $M$ und sind
 $v_1, \dots, v_n$ irgendwelche Elemente von $M$, so gibt es ein
 $r \in R$ mit $b_ir=v_i$ f\"ur $i := 1, \dots, n$.}\index{Dichtesatz}{}
 \smallskip
       Beweis. Nach 6.6 gibt es Elemente $s_i \in R$ mit $b_js_i = 0$
 f\"ur $j \neq i$ und $b_is_i \neq 0$. Es folgt weiter, dass
 $(b_is_i)R =M$ ist. Es gibt also ein $t_i \in R$ mit $b_is_it_i =
 v_i$. Setzt man nun $r := \sum_{j:=1}^{n} s_jt_j$, so ist
 $$ b_ir = \sum_{j:=1}^{n} b_is_jt_j = b_is_it_i = v_i. $$
 Damit ist der Dichtesatz bewiesen.
 \smallskip
       Man nennt einen Ring $R$ \emph{artinsch},\index{artinscher Ring}{} wenn der
 Verband der Rechtsideale von $R$ artinsch ist.
 \medskip\noindent
 {\bf 6.8. Satz.} {\it Es sei $R$ ein primitiver Ring. Genau dann gibt
 es einen Rechtsvektorraum $V$ endlichen Ranges \"uber einem
 K\"orper $K$, so dass $R$ zu $\End_K(V)$ isomorph ist, wenn $R$
 artinsch ist. Insbesondere ent\-h\"alt jeder solche Ring eine Eins.}
 \smallskip
       Beweis. Es sei $V$ ein Rechtsvektorraum endlichen Ranges \"uber
 dem K\"orper $K$. Dann ist $J_K(V) = \End_K(V)$, so dass $\End_K(V)$
 nach 4.15 artinsch ist. Somit ist $R$ artinsch, wenn $R$ zu
 $\End_K(V)$ isomorph ist.
 \par
       Es sei $R$ artinsch. Ferner sei $M$ ein treuer, irreduzibler
 $R$-Rechts\-mo\-dul. Dann ist $M$ ein Linksvektorraum \"uber dem
 K\"orper $K := \End_R(M)$. Setze
 $$ \Psi := \bigl\{O(X) \mid X \in E^K (M)\bigr\}. $$
 Weil $R$ artinsch ist, enth\"alt $\Psi$ ein minimales
 $O(Y)$. Nach 6.6 ist $Y$ dann ein maximales Element von $E^K(M)$.
 Es folgt, dass $Y = M$ ist. Ist $n$ der Rang von $M$ und ist $b_1$,
 \dots, $b_n$ eine Basis von $M$, sind ferner $v_1$, \dots, $v_n$
 irgendwelche Elemente von $M$, so gibt es nach dem Dichtesatz ein
 $r \in R$ mit $b_ir = v_i$ f\"ur alle $i$. Hieraus folgt, da $R$
 auf $M$ ja treu operiert, dass $R$ zu $\End^K(M)$ isomorph ist. Ist
 $V$ der Dualraum des $K$-Linksvektorraumes $M$, so folgt mit 5.20,
 dass $R$ auch zu $\End_K(V)$ isomorph ist.
 \medskip\noindent
 {\bf 6.9. Satz.} {\it Es sei $R$ ein Ring. Dann sind die folgenden
 Bedingungen \"aquivalent:
 \item{a)} $R$ ist als Ring die direkte Summe von endlich vielen
 Endomorphismenringen von Rechtsvektorr\"aumen end\-li\-chen Ranges.
 \item{b)} $R$ ist vollst\"andig reduzibel und hat eine Eins.
 \item{c)} $R$ ist artinsch und f\"ur jedes minimale Rechtsideal $I$ von $R$ gilt
 $I^2 \neq \{0\}$.}
 \smallskip
       Beweis. Endomorphismenringe von Rechtsvektorr\"aumen end\-li\-chen
 Rang\-es haben eine Eins und sind nach 4.13 vollst\"andig reduzibel sind; daher ist b) eine
 Folge von a).
 \par
       Es gelte b). Nach 4.6 ist $R = \bigoplus_{\Phi \in I_R(R)/\equiv}H_\Phi$.
 Weil $R$ eine Eins hat, gibt es endlich viele homogene Komponenten
 $H_1$, \dots, $H_t$ von $R$ mit $1 \in \sum_{k:=1}^{t} H_k$. Es folgt
 $R \subseteq \sum_{k:=1}^{t} H_k$, so dass $R$ nur endlich viele
 homogene Komponenten hat. Ist $H$ eine von ihnen, so ist $H$ die
 Summe von minimalen Rechtsidealen. Mit dem gleichen Argument
 erh\"alt man, dass $H$ Summe von endlich vielen minimalen
 Rechtsidealen ist. Hieraus folgt, dass der Rechtsidealverband von
 $H$ eine projektive Geometrie endlichen Ranges ist. Da $R$ die
 direkte Summe von endlich vielen homogenen Komponenten ist, ist
 $R$ daher artinsch. Es sei schlie\ss lich $I$ ein minimales Ideal
 von $R$. Weil $R$ vollst\"andig reduzibel ist, gibt es ein
 Rechtsideal $J$ mit $R = I \oplus J$. Es gibt daher ein $i \in I$
 und ein $j \in J$ mit $1 = i + j$. Es folgt
 $$ i = 1i = i^2 + ji $$
 und weiter
 $$ ji = i - i^2 \in I \cap J= \{0\}. $$
 Also ist $i^2 = i$. W\"are nun $i = 0$, so w\"are $1 \in J$ und
 daher $R \subseteq J$. Also ist $i \neq 0$, so dass auch $I^2 \neq
 \{0\}$ ist. Damit ist c) aus b) hergeleitet.
 \par
       Es gelte c). Weil $R$ artinsch ist, gibt es zu jedem von $\{0\}$
 verschiedenen Rechtsideal $J$ von $R$ ein minimales Rechtsideal
 $I$ mit $I \subseteq J$. Es gibt also eine Auswahlfunktion $\mu$,
 so dass $\mu(J)$ f\"ur alle solchen $J$ ein minimales, in $J$
 enthaltenes Rechtsideal ist. Setzt man noch $\mu(\{0\}) := \{0\}$,
 so ist $\mu$ auf der Menge aller Rechtsideale von $R$ definiert.
 \par
       Auf Grund unserer Annahme \"uber die minimalen Rechtsideale von
 $R$ gibt es nach 6.3 eine Auswahlfunktion $e$, so dass $e_I$ f\"ur
 alle minimalen Rechtsideale $I$ ein Idempotent ist, welches in $I$
 liegt und von $0$ verschieden ist. Ist $I$ das Nullideal, so
 setzen wir noch $e_I = 0$. Es sei nun $I$ ein minimales
 Rechtsideal von $R$. Dann ist $e_I \notin O(e_I)$. Andererseits ist
 $O(e_I)$ nach 4.2 ein regul\"ares, maximales Rechtsideal. Daher
 gilt
 $$ R = I + O(e_I). $$
 Aus dem dedekindschen Rekursionssatz folgt die Existenz einer Abbildung $g$ der
 Menge der nicht negativen Zahlen in die Menge der Rechtsideale von $R$
 mit $g(0) = R$ und
 $$ g(n + 1) = g(n) \cap O(e_{\mu(g(n))}). $$
 Es folgt
 $$ g(n) = g(n) \cap \bigl(\mu(g(n)) + O(e_{\mu(g(n))})\bigr)
	 = \mu\bigl(g(n)\bigr) + g(n + 1). $$
 Mittels Induktion folgt die G\"ultigkeit von
 $$ R = \sum_{i:=0}^{n} \mu\bigl(g(i)\bigr) + g(n+1) $$
 f\"ur alle $n$. Nun ist aber $g(n+1) \subseteq g(n)$ f\"ur alle $n$. Weil $R$
 artinsch ist, gibt es
 also ein $N$ mit $g(N) = g(N+1)$. Hieraus folgt weiter
 $$ \mu\bigl(g(N)\bigr) \subseteq g(N) \cap O(e_{\mu(g(N))}). $$
 Dies hat
 $$ R = \mu\bigl(g(N)\bigr) + O(e_{\mu(g(N))}) = O(e_{\mu(g(N))}) $$
 zur Folge, was wiederum $g(N) = \{0\}$ nach sich zieht. Also ist $R$
 Summe von minimalen Rechtsidealen und folglich, da deren Quadrat
 ja niemals das Nullideal ist, vollst\"andig reduzibel.
 \par
       Nach 4.6 ist $R$ direkte Summe seiner homogenen Komponenten, von
 denen es nur endlich viele gibt, da $R$ artinsch ist. Mittels 6.2
 und 6.4 folgt weiter, dass die homogenen Komponenten von $R$
 allesamt einfache Ringe sind. Einfache, vollst\"andig reduzible
 Ringe sind aber stets primitiv. Mittels 6.8 folgt dann, dass die
 homogenen Komponenten von $R$ zu Endomorphismenringen von
 Rechtsvektorr\"aumen isomorph sind. Damit ist a) aus c)
 hergeleitet und der Satz bewiesen.
 \medskip
       Die Bedingung c) wird in der Literatur meist so formuliert, dass man sagt,
 $R$ sei ein halbeinfacher, artinscher Ring, wobei
 {\it halbeinfach\/}\index{halbeinfacher Ring}{} bedeutet, dass $J(R) = \{0\}$
 ist. Benutzt man diese Formulierung, so muss man zeigen, dass das
 Ja\-cob\-son-Radikal\index{Jacobson-Radikal}{} eines Ringes $R$ alle nilpotenten
 Rechtsideale von $R$ enth\"alt. Ferner muss man zeigen, dass das
 Jacobson-Radikal eines artinschen Ringes nilpotent ist.\index{artinscher Ring}{}
 \par
       Der soeben bewiesene Satz ist grundlegend f\"ur die
 Darstellungstheorie\index{Darstellungstheorie von Gruppen}{}
 endlicher Gruppen, wie wir gleich belegen werden. Zuvor jedoch m\"ussen wir noch
 den Begriff der \emph{Gruppenalgebra}\index{Gruppenalgebra}{} einer endlichen
 Gruppe einf\"uhren. Dazu sei $K$ ein kommutativer K\"orper und $G$ sei eine
 endliche Gruppe. Wir interpretieren $G$ als Basis eines Vektorraumes $K[G]$
 \"uber $K$, von dem wir auf Grund der Kommutativit\"at von $K$ annehmen
 d\"urfen, dass er zweiseitig ist. Sind dann $x:= \sum_{g
 \in G}k_g g$ und $y:= \sum_{h \in G} l_h h$ zwei Elemente
 von $K[G]$, so definieren wir ihr Produkt $xy$ durch
 $$ xy := \sum_{{g \in G}\atop } \sum_{{a,b \in G} \atop { ab = g}} k_al_bg. $$
 Auf diese Weise wird $K[G]$ zu einer $K$-Algebra, der
 Gruppenalgebra von $G$ \"uber $K$. Es gilt nun der
 \medskip\noindent
 {\bf 6.10. Satz von Maschke.} {\it Ist $K$ ein kommutativer K\"orper,
 ist $G$ eine endliche Gruppe und ist $K[G]$ ihre Gruppenalgebra
 \"uber $K$, so ist $K[G]$ ein artinscher Ring mit Eins. Genau dann
 ist $K[G]$ vollst\"andig reduzibel, wenn die Charakteristik von
 $K$ kein Teiler von $|G|$ ist.}\index{Satz von Maschke}{}
 \smallskip
       Beweis. Die Eins von $G$ ist auch die Eins von $K[G]$. Dass $K[G]$
 artinsch ist, liegt daran, dass $K[G]$ ein Vektorraum endlichen
 Ranges \"uber $K$ ist. Damit ist der banale Teil des Satzes
 bewiesen.
 \par
       Die Charakteristik von $K$ sei ein Teiler von $|G|$. Setze $a:=
 \sum_{g \in G}g$. Weil $G$ eine Basis von $K[G]$ ist, ist
 $a \neq 0$. Ist $x \in G$, so ist $ax = a = xa$, so dass $a$ im
 Zentrum von $K[G]$ liegt. Dies beinhaltet, dass $aK[G]$ ein vom
 Nullideal verschiedenes, zweiseitiges Ideal von $K[G]$ ist. Nun ist
 $$ a^2 = a \sum_{g \in G} g = |G|a = 0. $$
 Daher ist $(a K [G])^2 = \{0\}$. Weil $R$ artinsch ist, enth\"alt
 $aK[G]$ ein minimales Rechtsideal $I$. F\"ur dieses gilt dann $I^2 = \{0\}$,
 so dass $K[G]$ nach 6.9 nicht vollst\"andig reduzibel ist.
 \par
       Es sei $K[G]$ nicht vollst\"andig reduzibel. Dann enth\"alt $K[G]$
 nach 6.9 ein minimales Rechtsideal $I$ mit $I^2 = \{0\}$. Es sei
 $0 \neq i \in I$. Dann ist $i = \sum_{g \in G} k_gg$,
 wobei es wenigstens ein $h \in G$ gibt mit $k_h \neq 0$. Dann ist
 auch $ih^{-1}$ ein von $0$ verschiedenes Element von $I$, so dass
 wir annehmen d\"urfen, dass in der Entwicklung f\"ur $i$ das
 Element $k_1$ von Null verschieden ist.
 \par
       Ist $a\in K[X]$, so ist die Abbildung $\alpha_a$, die durch
 $\alpha_a(x) := xa$ f\"ur alle $x \in K[G]$ definiert wird, eine
 $K$-lineare Abbildung von $K[G]$ in sich. Ist insbesondere $a \in
 G$, so ist $\Spur(\alpha_a) = 0$, wenn $a \neq 1$ gilt, und
 $\Spur(\alpha_a) = |G|$, wenn $a = 1$ ist. Hieraus folgt
 $$ \Spur(\alpha_i)= \sum_{g \in G}k_g \Spur(\alpha_g) = k_1|G|. $$
 Nun ist aber $i^2 = 0$, so dass alle Eigenwerte von $\alpha_i$ gleich $0$ sind.
 Dann ist aber auch die Spur von $\alpha$ gleich $0$,
 da sie die Summe der Eigenwerte ist. Also ist $k_1|G| = 0$. Weil $k_1$ aber
 nicht $0$ ist, ist $|G|$ in $K$ gleich $0$. Folglich ist die Charakteristik
 von $K$ ein Teiler von $|G|$. Damit ist der Satz von Maschke bewiesen.
 \medskip
       Zusammen mit Satz 6.9 besagt der Satz von Maschke, dass der
 Gruppenring\index{Gruppenring}{}
 einer endlichen Gruppe, falls die Charakteristik des zu Grunde liegenden
 kommutativen K\"orpers die Gruppenordnung nicht teilt, die direkte Summe von
 endlich vielen Endomorphismenringen von Vektorr\"aumen endlichen Ranges ist. Es
 ist zu beachten, dass die K\"orper, \"uber denen diese Vektorr\"aume definiert
 sind, zwar alle $K$ enthalten, aber durchaus nicht immer gleich $K$ sind.
 Einzelheiten findet der Leser in den B\"uchern zur Darstellungstheorie endlicher
 Gruppen.
 \par
       Wir wenden uns nun wieder der Untersuchung beliebiger einfacher,
 voll\-st\"an\-dig reduzibler Ringe zu. Dabei folgen wir bislang
 un\-ver\-\"of\-fent\-lich\-ten Notizen von Herrn
 U.\ Dempwolff.\index{Dempwolff, U.}{}
 \medskip\noindent
 {\bf 6.11. Satz.} {\it Es sei $R$ ein einfacher, vollst\"andig
 reduzibler Ring und $M$ sei ein treuer, irreduzibler Rechtsmodul
 \"uber $R$. Setze
 $$ K:= \End_R(M). $$
 Dann ist $K$ ein K\"orper und $M$ ist ein Linksvektorraum \"uber $K$. Weil $R$
 auf $M$ treu operiert, d\"urfen wir $R \subseteq \End^K(M)$ annehmen.
 Wir setzen
 $$ \Xi := \bigl\{U \mid U \in L^K(M)\ \hbox{\it und es gibt ein}\ r \in R\ 
		\hbox{\it mit}\ U = \Kern (r)\bigr\}. $$
 Dann gilt:
 \item{a)} Ist $U \in \Xi$, so ist $\Ko^M(U)$ endlich, dh. es ist $R \subseteq
 J^K(M)$.
 \item{b)} Ist $U \in \Xi$, ist $V \in L^K(M)$ und gilt $U \leq V$, so ist $V \in
 \Xi$.
 \item{c)} Sind $U$, $V \in \Xi$, so ist auch $U \cap V \in \Xi$.
 \item{d)} Es ist $\bigcap_{U \in \Xi}U = \{0\}$.
 \smallskip\noindent
 Ferner gilt: $R = \bigl\{r \mid r \in \End^K(M), \Kern (r) \in \Xi\bigr\}$.}
 \smallskip
       Beweis. Wir beweisen zun\"achst
 \par
       1) Genau dann ist $I$ ein minimales Rechtsideal von $R$,
 wenn es eine Hyperebene $H \in \Xi$ gibt mit $I = O(H)$. Ist $I$
 ein minimales Rechtsideal, so gibt es ein Idempotent $\pi \in I$
 mit $I=\pi R$.
 \par
       Es sei $I$ ein minimales Rechtsideal von $R$. Weil $R$ einfach
 ist, gibt es dann nach 6.4 und 6.3 ein Idempotent $\pi \in I$ mit
 $I = \pi R$. Setze $H:=\{m - m\pi \mid m \in M\}$. Dann ist $H \in
 L^K(M)$. \"Uberdies gilt $H \leq \Kern (\pi)$. Ist andererseits
 $m \in \Kern (\pi)$, so ist $m = m - m\pi \in H$, so dass $H =
 \Kern(i)$ ist. Es seien nun $u$ und $v$ zwei Elemente von
 $M$, so dass $u\pi$ und $v\pi$ linear unabh\"angig sind. Auf Grund
 des Dichtesatzes, der ja f\"ur $R$ gilt, gibt es dann ein $r \in
 R$, so dass $u \pi r \neq 0$ und $v \pi r = 0$ ist. Es folgt $0
 \neq \pi r \in I$ und damit $\pi r R = I$. Es gibt also ein $s \in
 R$ mit $\pi r s = \pi$. Es folgt der Widerspruch $v \pi = v \pi r
 s = 0$. Damit ist gezeigt, dass $H$ eine Hyperebene von $M$ ist.
 Es gilt ferner, dass $I = \pi R \subseteq O(H)$ ist. Es sei $s \in
 O(H)$ und $m \in M$. Dann ist $m = m \pi + m - m \pi$. Es folgt
 $ms = m \pi s$, so dass $s = \pi s \in I$ gilt. Somit ist $I = O(H)$.
 \par
       Es sei $H$ eine Hyperebene, die zu $\Xi$ geh\"ort. Es gibt dann
 ein $\lambda \in R$ mit $\Kern (\lambda) = H$. Es folgt $\lambda R
 \subseteq O(H)$. Weil $\lambda \neq 0$ ist, gibt es ein in
 $\lambda R$ enthaltenes minimales Rechtsideal $I$ von $R$. Wie wir
 bereits gesehen haben, enth\"alt $I$ ein Idempotent $\pi$ mit $\pi
 R = I = O(\Kern (\pi))$. Weil $\Kern(\pi)$ eine Hyperebene ist,
 die $H$ enth\"alt, ist $H = \Kern (\pi)$. Somit ist $O(H) = I
 \subseteq \lambda R \subseteq O(H)$, so dass $O(H)$ in der Tat ein
 minimales Rechtsideal ist.
 \par
       2) Es seien $H_1$, \dots, $H_t$ Hyperebenen, die zu $\Xi$
 geh\"oren. Setze
 $$ V_i := \bigcap_{k:=1, \dots, t,\ k \neq i}H_k $$
 f\"ur $i:= 1$, \dots, $t$. Es gelte ferner
 $$ V_i \not\leq H_i $$
 f\"ur $i := 1$, \dots, $t$. Es gibt dann ein Idempotent $\pi \in R$ mit
 $\bigcap_{i:=1}^{t} H_i = \Kern (\pi)$. Insbesondere ist
 $O(\bigcap_{i:=1}^{t} H_i) = \pi R$.
 \par
       Ist $t = 1$, so ist dies richtig nach 1). Es sei also $t \geq 2$. Wir
 setzen
 $$ U:= \bigcap_{i:=1}^t H_i. $$
 Es gibt dann $a_i \in V_i$ mit $V_i = Ka_i \oplus U$. Auf Grund der Annahme
 \"uber die Hyperebenen ist
 $$ M = U + \sum_{i:=1}^t Ka_i. $$
 Zu jedem $i$ gibt es nach 1) ein Idempotent $\pi_i$ mit $\Kern(\pi_i)
 = H_i$. Weil $V_i \not\leq H_i$ gilt, folgt $a_i \pi_i \neq 0$.
 Daher gelten die Gleichungen
 $$\eqalign{
     a_1O(H_1) &= M,        \cr
     a_2O(H_2) &= M,        \cr
	       &\ \ \vdots      \cr
     a_tO(H_t) &= M.        \cr} $$
 Es gibt also $\lambda_i \in O(H_i)$ mit $a_i \lambda_i = a_i$. Ist $i \neq j$,
 so ist $a_j \in H_i$, so dass $a_j\lambda_i = 0$ ist. Setze
 $\pi := \sum_{i:=1}^t \lambda_i$. Ist nun $m \in M$, so gibt es $k_i \in K$ und
 ein $u \in U$ mit
 $$ m = u + \sum_{i:=1}^t k_i a_i. $$
 Dann folgt $m \lambda_i = k_ia_i$ und weiter
 $$ m\pi = \sum_{i:=1}^t k_ia_i. $$
 Weil die $a_i$ linear unabh\"angig sind, folgt, dass der Kern von
 $\pi$ gleich $U$ ist. Ferner folgt $m\pi^2 = m\pi$ f\"ur alle $m
 \in M$, so dass $\pi$ in der Tat ein Idempotent ist.
 \par
       Um auch die letzte Aussage zu beweisen, bemerken wir zun\"achst,
 dass $\pi R \subseteq O(\bigcap_{i:=1}^t H_i)$ gilt. Es sei $r \in
 O(\bigcap_{i:=1}^t H_i)$. Ist dann wieder $m = u +
 \sum_{i:=1}^t k_ia_i$, so folgt
 $$ mr = \sum_{i:=1}^t k_ia_i r = m \pi r. $$
 Also ist $r = \pi r \in \pi R$. Damit ist 2)
 bewiesen.
 \par
       Als N\"achstes beweisen wir:
 \par
       3) Es seien $H_1$, \dots, $H_t$ Hyperebenen, die zu $\Xi$
 geh\"oren. Setze
 $$ V_i := \bigcap_{k:=1, \dots, t,\ k \neq i} H_k $$
 f\"ur $i := 1$, \dots, $t$. Es gelte ferner
 $$ V_i \not\leq H_i $$
 f\"ur $i := 1$, \dots, $t$. Ist dann $H$ eine Hyperebene von $M$ mit
 $\bigcap_{i:=1}^t H_i \leq H$, so ist $H \in \Xi$.
 \par
       Dies ist richtig, falls $H$ gleich einem der $H_i$ ist. Dies ist
 sicher dann der Fall, wenn $t = 1$ ist. Wir d\"urfen daher
 annehmen, dass $H$ von allen $H_i$ verschieden ist und dass der
 Satz f\"ur $t - 1$ gilt. Wendet man 2.6 auf $L^K (M/\bigcap_{i:=1}^t H_i)^d$ an,
 so folgt die Existenz einer
 Hyperebene $H'$ mit $H_1 \cap H' \leq H$ und $\bigcap_{i:=2}^t H_i
 \leq H'$. Wir d\"urfen daher annehmen, dass $t = 2$ und --- nach
 Induktionsannahme --- dass $H' = H_2$ ist. Weil $H$ von $H_1$ und
 $H_2$ verschieden ist, gibt es ein von Null verschiedenes $p \in
 H$ mit $H = Kp \oplus (H_1 \cap H_2)$. Es folgt $pO(H_1) = M$ und
 $pO(H_2) = M$. Es gibt also ein $r_i \in O(H_i)$ f\"ur $i:= 1$, 2 mit
 $pr_1 = p = -pr_2$. Setze $f:= r_1 + r_2$. Dann ist $H \leq \Kern
 (f)$, wie man unmittelbar sieht. W\"are $f = 0$, so folgte $r_1 =
 -r_2$. Weil weder $r_1$ noch $r_2$ Null sind, folgte mit 1) weiter
 $O(H_1) = r_1R = r_2R = O(H_2)$, da die Ideale $O(H_i)$ ja minimal
 sind. Hieraus folgte aber der Widerspruch $H_1 = H_2$. Also ist $f
 \neq 0$ und daher $\Kern (f) = H$.
 \par
       a) Es sei $r \in R$. Es gibt dann endlich viele minimale Rechtsideale
 $I_1$, \dots, $I_t$ mit $r \in \sum_{k:=1}^t I_k$. Nach 1) gibt es Hyperebenen
 $H_1$, \dots, $H_t$ mit $O(H_k) = I_k$. Ist $m \in \bigcap_{k:=1}^t H_k$, so
 folgt $mr = 0$, so dass
 $$ \bigcap_{k:=1}^t H_k \leq \Kern (r) $$
 gilt. Folglich hat $\Kern (r)$ endlichen Ko-Rang, so dass a) bewiesen ist.
 \par
       b) Wir haben gerade gesehen, dass es Hyperebenen $H_1$, \dots, $H_t
 \in \Xi$ gibt mit $\bigcap_{i:=1}^t H_i \leq U$. Da alle
 Hyperebenen von $M$, die oberhalb des Schnitts der $H_i$
 liegen, nach 3) zu $\Xi$ geh\"oren, und da $V$ Schnitt von
 endlich vielen dieser Hyperebenen ist, geh\"ort auch $V$ nach 2)
 zu $\Xi$.
 \par
       Der Beweis von c) ist analog dem von b).
 \par
       d) Es sei schlie\ss lich $m \in \bigcap_{U \in \Xi} U$.
 Dann ist $mR = \{0\}$ und daher $m = 0$, da $M$ ja irreduzibel
 ist. Damit ist auch d) bewiesen.
 \par
       Um die letzte Aussage zu beweisen, sei $U \in \Xi$. Da $U$ sich
 dann als Schnitt von unabh\"angigen Hyperebenen darstellen
 l\"asst, die alle zu $\Xi$ geh\"oren, gibt es eine Projektion
 $\pi$ mit $O(U) = \pi R$. Es folgt, dass $M = M \pi \oplus U$ ist.
 Es sei nun $a_1$, \dots, $a_n$ eine Basis von $M \pi$. Ferner sei
 $s \in \End^K (M)$ und es gelte $U \leq \Kern (s)$. Auf Grund des
 Dichtesatzes gibt es ein $r \in R$ mit $a_ir = a_is$ f\"ur alle
 $i$. Hieraus folgt $a_i \pi r = a_i s$ f\"ur alle $i$. Da $U$
 sowohl von $\pi r$ als auch von $s$ annulliert wird, gilt also $s =
 \pi r \in R$. Damit ist alles bewiesen.
 \medskip
       Eine unmittelbare Folgerung aus der Aussage 2) des Beweises von
 6.11 notieren wir als
 \medskip\noindent
 {\bf 6.12. Korollar.} {\it Ist $R$ ein einfacher, vollst\"andig
 reduzibler Ring, so gibt es zu jedem endlich erzeugten Rechtsideal
 $I$ von $R$ ein Idempotent $\pi \in I$ mit $I = \pi R$.}
 \medskip
       Die Teilmenge $I$ des Verbandes $L$ hei\ss t \emph{duales
 Ideal}\index{duales Ideal}{} von
 $L$, wenn aus $A, B \in I$ stets auch $A \cap B \in I$ folgt und
 mit $A$ auch alle $X$ mit $A \leq X$ zu $I$ geh\"oren. Die in Satz
 6.11 definierte Menge $\Xi$ ist also ein duales Ideal von $L^K
 (M)$.
 \medskip\noindent
 {\bf 6.13. Satz.} {\it Es sei $K$ ein K\"orper und $M$ sei ein von
 $\{0\}$ verschiedener Linksvektorraum \"uber $K$. Ferner sei $\Xi$
 ein duales Ideal von $L^K (M)$ aus Unterr\"aumen endlichen
 Ko-Ranges und es gelte $\bigcap_{U \in \Xi} u = 0$. Ist dann
 $$ R:= \bigl\{r \mid r \in \End^K(M),\ \Kern (r) \in \Xi\bigr\}, $$
 so ist $R$ einfach und vollst\"andig reduzibel. \"Uberdies ist $M$ ein
 treuer, irreduzibler $R$-Rechtsmodul.}
 \smallskip
       Beweis. Weil $M$ nicht trivial ist, folgt aus $\bigcap_{U \in \Xi} U =
 \{0\}$, dass $\Xi$ nicht leer ist. Weiter folgt, dass es zu
 linear unabh\"angigen $a_1$, \dots, $a_n$ ein $U \in \Xi$ gibt mit
 $(\sum_{i:=1}^n Ka_i) \cap U = \{0\}$. Um dies
 einzusehen, setzen wir 
 $$ X := \sum_{i:=1}^n Ka_i. $$
 Weil $X$ endlichen Rang hat, gibt es dann ein $U \in \Xi$, so dass $X
 \cap U$ von minimalen Rang ist. Ist nun $W \in \Xi$, so ist auch
 $U \cap W \in \Xi$. Hieraus folgt, dass
 $$ X \cap U \cap W = X \cap U $$
 ist. Somit ist $X \cap U \leq W$ f\"ur alle $W \in
 \Xi$. Folglich ist $X \cap U= \{0\}$, wie behauptet. Es gibt nun
 ein Komplement $V$ von $\sum_{i:=1}^n Ka_i$ in $M$ mit $U \leq V$.
 Sind $x_1$, \dots, $x_n$ beliebige Elemente von $M$, so definieren
 wir $r \in \End^K (M)$ durch $a_ir := x_i$ und $vr := 0$ f\"ur alle
 $v \in V$. Es folgt $U \leq V \leq \Kern (r)$, so dass $\Kern (r)
 \in \Xi$ gilt, da $\Xi$ ja ein duales Ideal ist. Folglich ist $r
 \in R$, so dass die Folgerung des Dichtesatzes f\"ur $R$ gilt.
 Hieraus folgt insbesondere, dass $M$ ein irreduzibler Rechtsmodul
 \"uber $R$ ist.
 \par
       Es sei $H$ eine Hyperebene, die zu $\Xi$ geh\"ort. Ferner sei $0
 \neq x \in O(H)$ sowie $y \in O(H)$. Ferner sei $u \in M - H$.
 Setze $v := ux$ und $w := uy$. Weil die Folgerung des Dichtesatzes
 f\"ur $R$ gilt, gibt es ein $r \in R$ mit $vr = w$. Es folgt $u(xr
 - y) = 0$. Wegen $M = Ku \oplus H$ ist daher $xr = y$ und folglich
 $O(H) = xR$, so dass $O(H)$ ein minimales Rechtsideal von $R$ ist,
 f\"ur das \"uberdies $O(H)R = O(H)$ gilt.
 \par
       Es sei $r \in R$. Setze $U := \Kern (r)$. Weil $\Xi$ ein duales
 Ideal ist, gibt es dann Hyperebenen $H_1$, \dots, $H_n \in \Xi$ mit
 $U = \bigcap_{i:=1}^n H_i$. Es folgt
 $$ r \in O(U) = \sum_{i:=1}^n O(H_i), $$
 so dass $R$ Summe von minimalen Rechtsidealen ist. Also ist $R$
 voll\-st\"an\-dig reduzibel.
 \par
       Es seien nun $H$ und $H'$ zwei verschiedene Hyperebenen, die beide
 zu $\Xi$ geh\"oren. Es sei $P$ ein Punkt auf $H$, der nicht auf
 $H'$ liegt, und $Q$ ein Punkt auf $H'$, der nicht auf $H$ liegt.
 Es sei $P = Kp$ und $Q = Kq$. Schlie\ss lich sei $S := K(p+q)$. Es
 gibt dann einen Endomorphismus $r$ von $M$ mit $pr = q$, $qr = p$
 und $\Kern (r) = H \cap H'$. Weil $\Xi$ ein duales Ideal ist, ist
 $\Kern (r) \in \Xi$, so dass $r \in R$ gilt. Es sei weiter $\pi$
 die Projektion von $M$ auf $S$ mit $\Kern(\pi) = H$ und $\pi '$ die
 Projektion von $M$ auf $S$ mit $\Kern(\pi') H'$. Es folgen die
 Gleichungen
 $$\eqalign{
           qr \pi r &= p \pi r = 0 = q \pi',   \cr
    (p + q) r \pi r &= p + q = (p + q) \pi',   \cr
	   ur \pi r &= 0 = u \pi'              \cr} $$
 f\"ur alle $u \in H \cap H'$. Daher ist $r \pi r =
 \pi'$. Definiert man nun die Abbildung $\lambda$ durch $\lambda
 (\pi s) := r \pi s$, so ist $\lambda$ ein Homomorphismus von $O(H)
 = \pi R$ auf $r \pi R$. Wegen $r \pi r = \pi' \in r \pi R \cap
 O(H')$ ist $\lambda$ also ein Isomorphismus von $O(H)$ auf
 $O(H')$. Hieraus folgt insbesondere, dass $R$ homogen ist.
 \par
       Die Ideale der Form $O(H)$ mit einer Hyperebene $H \in \Xi$ sind
 aber alle minimalen Rechtsideale von $R$. Ist n\"amlich $I$ ein
 minimales Rechtsideal und ist $H$ eine Hyperebene in $\xi$, so
 gibt es, da $R$ homogen ist, einen Isomorphismus $\lambda$ von
 $O(H)$ auf $I$. Wie wir bereits gesehen haben, wird $O(H)$ von
 einem Idempotent $\pi$ erzeugt. Setze $r := \lambda(\pi)$. Dann
 ist $r = r\pi \in I \cap rR$. Wegen der Minimalit\"at von $I$ ist
 daher $I \subseteq rR$, so dass $I = rR$ ist. Wegen $r \in O(\Kern
 (r))$ gilt weiter $I \subseteq O(\Kern (r))$. Weil der Rang von
 $\pi$ gleich $1$ und $ r \neq 0$ ist, ist $1$ auch der Rang von
 $r$. Folglich ist $\Kern (r)$ eine Hyperebene. Somit ist $O(\Kern
 (r))$ ein minimales Rechtsideal und folglich $I = O(\Kern (r))$.
 \par
       Weil nun alle minimalen Rechtsideale Idempotente enthalten, ist
 ihr Quadrat niemals das Nullideal. Daher ist das Jacobson-Radikal
 von $R$ gleich dem Null\-ide\-al, so dass $R$ einfach ist. Hieraus
 folgt schlie\ss lich, dass $M$ auch ein treuer $R$-Modul ist.
 \medskip\noindent
 {\bf 6.14. Satz.} {\it Es sei $R$ ein einfacher, vollst\"andig
 reduzibler Ring und $M$ sei ein treuer, irreduzibler
 $R$-Rechtsmodul. Setze $K:= \End_R(M)$. Ist $L$ ein Linksideal
 von $R$, so setzen wir
 $$ \Psi (L) := \sum_{\sigma \in L} M \sigma. $$
 Ist $U \in L^K (M)$, so setzen wir
 $$ \Phi (U) := \{\sigma \mid \sigma \in R,\ M\sigma \leq U\}. $$
 Dann ist $\Psi$ ein Isomorphismus von $L^R (R)$ auf $L^K (M)$ und es gilt
 $\Psi^{-1} = \Phi$.}
 \smallskip
       Beweis. Es ist trivial, dass $\Phi (U)$ ein Linksideal von $R$ und
 dass $\Phi (L)$ ein Teilraum von $M$ ist.
 \par
       Es sei $\Xi$ die Menge der Kerne der Elemente von $R$. Nach 6.11
 ist $\Xi$ dann ein duales Ideal von $L^K(M)$ aus Unterr\"aumen
 endlichen Ko-Ranges, f\"ur das au\ss erdem $\bigcap_{U \in \Xi} U
 = \{0\}$ gilt. Weil jedes $U \in \Xi$ Schnitt von Hyperebenen
 ist, die ebenfalls zu $\Xi$ geh\"oren, gilt sogar
 $$ \bigcap_{H \in \Xi,\ H\mathrm{\ ist\ Hyperebene}} H = \{0\}. $$
 Ferner ist
 $$ R = \bigl\{r \mid r \in \End^K(M),\ \Kern(r) \in \Xi\bigr\}. $$
 Es sei $U \in L^K (M)$. Dann ist
 $$ \Psi \Phi (U) = \Psi\bigl(\{\sigma \mid \sigma \in R,\ 
			M \sigma \leq U\}\bigr)
                 = \sum_{\sigma \in R, M \sigma \leq U} M \sigma \leq U. $$
 Es sei $P$ ein Punkt mit $P \leq U$. Nach der zuvor
 gemachten Bemerkung, dass der Schnitt \"uber die zu $\Xi$
 geh\"orenden Hyperebenen Null ist, gibt es also eine Hyperebene $H
 \in \Xi$ mit $M = P \oplus H$. Es sei $\pi$ die Projektion von $M$
 auf $P$ mit Kern $H$. Dann ist $\pi \in R$ und wegen $M \pi = P
 \leq U$ gilt sogar $\pi \in \Phi (U)$. Folglich ist $P = M \pi
 \leq \Psi \Phi (U)$. Weil $U$ die obere Grenze der in $U$
 liegenden Punkte ist, ist also $U \leq \Psi \Phi (U)$, so dass in
 der Tat $U = \Phi \Psi (U)$ ist.
 \par
       Es sei $L$ ein Linksideal von $R$. Dann ist
 $$ \textstyle\Phi \Psi (L) = \Phi \bigl( \sum_{\sigma \in L} M \sigma \bigr) 
	       = \{ \tau \mid \tau \in R, M \tau \leq \sum_{\sigma \in L}\ 
	       M \sigma\} \supseteq L. $$
 \par
       Es sei $\tau \in \Phi \Psi (L)$. Wir m\"ussen zeigen, dass $\tau
 \in L$ ist. Nun hat $M \tau$ endlichen Rang. Wie wir weiter oben
 schon gesehen haben, gibt es ein $C \in \Xi$ mit $M = M \tau
 \oplus C$. Ist $\pi$ die Projektion von $M$ auf $M \tau$ mit
 $\Kern(\pi) = C$, so ist $\tau \in \Phi \Psi (L)$ und es gilt $\tau \pi =
 \tau$. Es gen\"ugt daher zu zeigen, dass $\pi \in L$ gilt. Ist
 $b_1$, \dots, $b_n$ eine Basis von $M \tau$, so definieren wir
 Projektionen $\rho_i$ durch $b_j \rho_i := 0$, f\"ur $j \neq i$
 und $b_i \rho_i := b_i$ sowie $c \rho_i := 0$ f\"ur alle $c \in
 C$. Weil $\Xi$ ein duales Ideal ist, liegen alle $\rho_i$ in $R$
 und damit in $\Phi \Psi (L)$. Weil $\pi = \sum_{i := 1}^n
 \rho_i$ ist, d\"urfen wir daher annehmen, dass der Rang von $\pi$
 gleich $1$ ist.
 \par
       Es gibt nun $\sigma_1$, \dots, $\sigma_t$ in $L$ mit
 $$ M \pi \leq \sum_{i:=1}^t M \sigma_i. $$
 Es sei $M \pi = Kp$. Ein solches $p$ gibt es, da der Rang von $\pi$ ja Eins
 ist. Es gibt dann $m_i \in M$ mit $p = \sum_{i := 1}^t m_i \sigma_i$.
 Es sei $0 \neq m \in M$. Dann ist $mR = M$. Es gibt daher
 $\lambda_i \in R$ mit $m_i=m \lambda_i$. Es folgt
 $$ p = m \sum_{i := 1}^t \lambda_i \sigma_i, $$
 so dass also $p = m \sigma$ ist mit einem $\sigma \in L$. Es gibt eine
 Projektion $\alpha \in R$ von $M$ auf $K m$. Es folgt 
 $$ p = m \alpha \sigma, $$
 so dass wir annehmen d\"urfen, dass $\sigma$ den Rang $1$ hat.
 \par
       Weil der Rang von $\sigma$ gleich Eins ist, ist also $\sigma$, $\pi
 \in \Psi (Kp)$. Der Ko-Rang von $U := \Kern (\pi) \cap \Kern
 (\sigma)$ sei gleich $2$. Dann geh\"ort $U$ wiederum zu $\Xi$. Es
 gibt nun Elemente $a_1$ und $a_2$ mit $\Kern (\pi) = K a_1 \oplus
 U$ und $\Kern (\sigma) = K a_2 \oplus U$. Es folgt
 $$ M = K a_1 \oplus K a_2 \oplus U. $$
 \par
       Wir definieren $\rho$ durch $a_1 \rho := a_2$, $a_2 \rho = a_1$ und
 $\Kern (\rho) := U$. Es folgt $a_1 \rho \sigma = 0$, so dass
 $\Kern (\rho \sigma) = \Kern (\sigma)$ ist. Es folgt weiter, dass
 $M \rho \sigma = Kp$ ist. Wir d\"urfen daher annehmen, dass der
 Kern von $\sigma$ gleich dem Kern von $\pi$ ist.
 \par
       Definiere $\mu \in R$ durch $m \mu = km$ und $\Kern(\mu) = \Kern(\pi)$,
 wobei $k \in K$ durch $m \pi = kp$ definiert sei. Dann ist
 $$ m \mu \sigma = km \sigma = kp = m \pi $$ 
 und folglich $\pi = \mu \sigma$, da die Kerne der drei Abbildungen
 identisch sind. Also ist $\pi = \mu \sigma \in L$, so dass in der
 Tat $L = \Phi \Psi (L)$ ist. Damit ist der Satz bewiesen.
 \medskip
       Ist $R$ ein einfacher, vollst\"andig reduzibler Ring, so ist er
 auch als Linksmodul \"uber sich vollst\"andig reduzibel, wie der
 gerade bewiesene Satz zeigt. Dennoch besteht scheinbar eine
 Unsymmetrie zwischen dem Verband der Linksideale und dem Verband
 der Rechtsideale von $R$. Dieser Eindruck der Unsymmetrie wird
 noch best\"arkt durch das folgende Beispiel.
 \par
       Es sei $M$ ein Linksvektorraum \"uber dem K\"orper $K$ und $B$ sei
 eine Basis von $M$. Wir nehmen an, dass $B$ nicht endlich sei. Es
 sei $\Xi$ die Menge der Teilr\"aume $U$ von $M$, f\"ur die $B - U$
 endlich ist. Dann ist $\Xi$ ein duales Ideal von $L^K (M)$ mit
 $\Xi \subseteq E^K (M)$ und $\bigcap_{U \in \Xi} U=\{0\}$. Setzt
 man
 $$ R := \bigl\{\sigma \mid \sigma \in J^K (M),\ \Kern (\sigma) \in \Xi\bigr\}, $$
 so ist $R$ nach 6.13 ein vollst\"andig reduzibler, einfacher Ring. Dieser Ring
 ist nun nicht gleich $J^K (M)$. Ist n\"amlich $a \in B$ und definiert man $\rho$
 durch $b \rho := a$ f\"ur alle $b \in B$, so ist $\rho \in J^K (M)$ und es gilt
 $B \cap \Kern (\rho) = \{0\}$. Somit liegt $\rho$ nicht in $R$.
 \par
       Diese scheinbare Unsymmetrie wird aber beseitigt, wenn man ne\-ben einem
 irreduziblen Rechtsmodul \"uber $R$ einen geeigneten ir\-re\-du\-zib\-len
 Linksmodul \"uber $R$ in die Betrachtungen mit einbezieht. Der Rest dieses
 Abschnitts sei der Darstellung dieses Sachverhalts gewidmet.
 \par
       Das Rechtsideal $I$ des Ringes $R$ hei\ss t
 \emph{nilpotent},\index{nilpotent}{} wenn
 es eine nat\"urliche Zahl $n$ gibt, so dass $I^n = \{0\}$ ist.
 Nilpotente Linksideale, bzw.\ nilpotente Ideale werden entsprechend definiert.
 \medskip\noindent
 {\bf 6.15. Satz.} {\it Es sei $R$ ein Ring. Dann sind die folgenden Aussagen
 gleichbedeutend:
 \item{a)} $R$ enth\"alt kein von $\{0\}$ verschiedenes nilpotentes Rechtsideal.
 \item{b)} $R$ enth\"alt kein von $\{0\}$ verschiedenes nilpotentes Ideal.
 \item{c)} $R$ enth\"alt kein von $\{0\}$ verschiedenes nilpotentes Linksideal.}
 \smallskip
       Beweis. Weil zweiseitige Ideale auch Rechtsideale
 sind, ist b) eine Konsequenz von a).
 \par
       Es sei $I$ ein nicht triviales nilpotentes Rechtsideal von $R$.
 Ist dann $RI = \{0\}$, so ist $I$ sogar ein zweiseitiges
 nilpotentes Ideal von  $R$. Es sei also $RI \neq \{0\}$ und $I^n =
 \{0\}$. Dann ist $RI$ ein nicht triviales zweiseitiges Ideal von
 $R$. Wegen
 $$(RI)^n = R(IR)^{n-1} I \subseteq RI^n = \{0\} $$
 ist $RI$ nilpotent. Somit ist a) eine Folge von b)
 \par
       Die \"Aquivalenz von b) mit c) beweist sich entsprechend.
 \medskip\noindent
 {\bf 6.16. Satz.} {\it Es sei $R$ ein Ring ohne von $\{0\}$ verschiedene,
 zweiseitige Ideale. Ist $e$ ein von $0$ verschiedenes Idempotent von $R$, so
 sind die folgenden Aussagen \"aquivalent:
 \item{a)} $eR$ ist ein minimales Rechtsideal von $R$.
 \item{b)} $eRe$ ist ein K\"orper.
 \item{c)} $Re$ ist ein minimales Linksideal von $R$.}
 \smallskip
       Beweis. Es gelte a). Wegen $e = e^3 \in eRe$ besteht $eRe$ nicht
 nur aus der Null. Weiter folgt $(ere)R = eR$, falls nur $ere \neq 0$ ist. Es
 gibt dann ein $s \in R$ mit $eres = e^2 = e$. Hieraus folgt weiter
 $$ e = e^2 = erese = ereese, $$
 so dass $ere$ ein Rechtsinverses hat. Da $e$ eine Rechtseins von $eRe$ ist,
 bilden die von Null verschiedenen Elemente von $eRe$ eine Gruppe, so dass
 $eRe$ ein K\"orper ist.
 \par
       Es gelte b) und $I$ sei ein von $\{0\}$ verschiedenes Rechtsideal,
 welches in $eR$ ent\-hal\-ten ist. Weil $R$ keine nicht trivialen,
 nilpotenten zweiseitigen Ideale enth\"alt, folgt mit 6.15, dass
 $I^2 \neq \{0\}$ ist. Nun ist $I^2 \subseteq IeR$ und daher $IeR
 \neq \{0\}$. Hieraus folgt die Existenz eines $a \in I$ mit $ aeR
 \neq \{0\}$. Wegen $I \subseteq eR$ ist $ea = a$. Es folgt
 $$ eae = ae = ae^2 \in IeR \subseteq I. $$
 Somit ist $eae$ ein von Null verschiedenes Element im Schnitt von $I$ mit $eRe$.
 Weil $eRe$ ein K\"orper ist, gibt es ein $b \in R$ mit $e = eaeebe$.
 Hieraus folgt $e \in I$ und damit $eR \subseteq I$, so dass $eR$
 in der Tat ein minimales Rechtsideal ist.
 \par
       Die \"Aquivalenz von b) mit c) beweist sich entsprechend.
 \medskip\noindent
 {\bf 6.17. Satz.} {\it Es sei $R$ ein einfacher Ring und $e$ sei ein
 Idempotent von $R$, so dass $K := eRe$ ein K\"orper ist. Setze $W := eR$ und
 $V := Re$. Dann gilt:
 \item{a)} Es ist $W \cap V = K$.
 \item{b)} $W$ ist ein Links- und $V$ ein Rechtsvektorraum \"uber $K$.
 \item{c)}\, Die Abbildung, die $w \in W$ und $v \in V$ ihr Produkt $wv$ in $R$
 zuordnet, ist eine Koppelung von $W$ mit $V$ \"uber $K$.
 \item{d)}\, Ist $\varphi \in \End_R (W)$, so ist $\varphi e \in K$ und es gilt
 $\varphi w = (\varphi e) w$ f\"ur alle $w \in W$. Somit ist $K$ zu $\End_R (W)$
 isomorph.
 \item{e)} \,Ist $\varphi \in \End^R (V)$, so ist $\varphi e \in K$ und es gilt
 $v \varphi = v(\varphi e)$ f\"ur alle $v \in V$. Somit ist $K$ zu $\End^R (V)$
 isomorph.
 \item{f)}\, F\"ur $R$ gilt der Dichtesatz sowohl als Unterring von $\End_K (V)$ als
 auch als Unterring von $\End^K (W)$.\par}
 \smallskip
       Beweis. a) Es ist $K = eRe \subseteq  Re \cap eR = W \cap V$. Ist
 andererseits $x \in Re \cap eR$, so ist $x = re = es$ und folglich
 $x = exe \in K$.
 \par
       b) ist trivial.
 \par
       c) Von den nachzuweisenden Bedingungen sind nur die Bedingungen 5)
 und 6) nicht v\"ollig banal. Es sei also $w \in W$ und es gelte
 $wV = \{0\}$. Ferner sei $w = er$. Dann ist also $erse = 0$ f\"ur
 alle $s \in R$. W\"are $er \neq 0$, so folgte $erR = eR$, da $eR$
 wegen der Einfachheit von $R$ nach 6.16 ein minimales Ideal ist.
 Es g\"abe daher ein $t \in R$ mit $ert = e$. Hieraus folgte der
 Widerspruch $0 = erte = e$. Entsprechend zeigt man, dass aus $Wv =
 \{0\}$ folgt, dass $v = 0$ ist.
 \par
       d) Es ist $\varphi e \in W = eR$. Hieraus folgt
 $$ \varphi e = \varphi (e^2) = (\varphi e) e \in eRe = K. $$
 Alles weitere ist banal.
 \par
       e) beweist sich wie d).
 \par
       f) folgt mittels d) bzw. e), da $R$ als einfacher Ring ja auf jedem
 irreduziblen Modul treu operiert, also primitiv ist.
 \medskip\noindent
 {\bf 6.18. Satz.} {\it Es sei $R$ ein einfacher, vollst\"andig
 reduzibler Ring. Dann enth\"alt $R$ ein Idempotent $e$, so dass
 $eR$ ein minimales Rechtsideal ist. Setze $K := eRe$. Dann ist $K$
 ein K\"orper. Setze $W := eR$ und $V := Re$. Dann sind $W$ und $V$
 gem\"a\ss\ der Beschreibung von Satz 6.17 gekoppelte Vektorr\"aume
 \"uber $K$. Setze
 $$ \Xi^W := \bigl\{X^\bot \mid X \in E_K (V)\bigr\} $$
 und
 $$ \Xi_V := \bigl\{Y^\top \mid Y \in E^K (W)\bigr\}. $$
 Dann gilt:
 \item{a)} Der Verband $L^R (R)$ der Linksideale von $R$ ist isomorph zu
 $L^K (W)$.
 \item{b)} Der Verband $L_R (R)$ der Rechtsideale von $R$ ist isomorph zu
 $L_K (V)$.
 \item{c)} Es ist
 $$ R = \bigl\{\sigma \mid \sigma \in J^K (W),\ \Kern(\sigma) \in \Xi^W\bigr\}.$$
 \item{d)} Es ist}
 $$ R = \bigl\{\sigma \mid \sigma \in J_K (V),\ \Kern(\sigma) \in \Xi_V\bigr\}.$$
 \smallskip
       Beweis. a) ist nichts Anderes als die Aussage des Satzes 6.14. Sie besagt
 dar\"uber hinaus, dass $R$ auch als Linksmodul \"uber sich vollst\"andig
 reduzibel ist. Weil $V$ ein irreduzibler Linksmodul \"uber $R$ ist, gilt daher
 auch b).

 Zur Notation $E_K(V)$, $E^K(W)$ siehe 6.6 und die Bemerkung vor 5.14.

       Um c) zu beweisen, sei zun\"achst $r \in R$. Dann hat $rV$
 endlichen Rang. Ist $w \in W$, so gilt $w(rV) = \{0\}$ genau dann,
 wenn $(wr)V = \{0\}$ ist. Dies ist wiederum genau dann der Fall,
 wenn $wr = 0$ ist. Also ist $\Kern (r) = (rV)^\bot$, wobei Kern
 sich hier auf die Wirkung von $r$ auf $W$ bezieht. Es gilt also
 $\Kern (r) \in \Xi^W$. Es sei $Y \in \Xi^W$. Es gibt dann einen
 Teilraum $X \in E_K (V)$ mit $Y = X^\bot$. Es sei $b_1, \ldots,
 b_n$ eine Basis von $X$. Auf Grund des Dichtesatzes, der nach 5.17
 ja gilt, gibt es ein $r \in R$ mit $rb_i = b_i$ f\"ur alle $i$. Es
 ist $X \leq rV$  und daher $\Kern (r) = (rV)^\bot \leq X^\bot$.
 Nach 6.11 gibt es also ein $s \in R$ mit $\Kern (s) = X^\bot = Y$.
 Damit ist gezeigt, dass $\Xi^W$ gerade die Menge der $\Kern (r)$
 mit $r \in R$ ist. Eine nochmalige Anwendung von 6.11 zeigt dann
 die G\"ultigkeit von c).
 \par
       d) gilt aus Symmetriegr\"unden.
 \medskip
       Was in diesem Abschnitt \"uber das Jacobson-Radikal eines Ringes
 gesagt wurde, ist nat\"urlich alles auf den uns allein
 interessierenden Fall vollst\"andig reduzibler Ringe zugeschnitten
 worden. Der Leser, der mehr erfahren m\"ochte, sei auf die
 B\"ucher \"uber Ringtheorie verwiesen.

\mysection{7. Endliche projektive R\"aume}

\noindent
 In diesem Abschnitt wollen wir die Anzahl der Unterr\"aume vom
 Range $i$ eines endlichen projektiven Raumes sowie die Anzahl
 seiner Basen bestimmen. Dies bedarf noch einiger Vorbereitung. Wir
 beweisen zun\"achst den folgenden Satz.
 \medskip\noindent
 {\bf 7.1. Satz.} {\it Es sei $L$ ein irreduzibler projektiver Verband. Sind $A$,
 $B$ und $D$ Elemente von $L$ mit der Eigenschaft, dass $D$ maximal ist
 bez\"uglich der Eigenschaft, mit $A$ und $B$ trivialen Schnitt zu haben, so ist
 $A + D \leq B + D$ oder $B + D \leq A + D$.}
 \smallskip
       Beweis. Es sei $P$ ein Punkt von $A + B + D$, der nicht in $D$ liegen
 m\"oge. Dann ist $A \cap (D + P) \neq 0$ oder $B \cap (D + P) \neq 0$.
 \par
       Es sei etwa $A \cap (D + P) \neq 0$. Es gibt dann einen Punkt $Q
 \leq A \cap (D + P)$. Ist $Q = P$, so ist $P \leq A + D$. Es sei
 also $Q \neq P$. Nach 2.6 gibt es dann einen Punkt $R \leq D$ mit
 $Q \leq P + R$. Wegen $A \cap D = 0$ ist $Q \neq R$. Daher ist $P
 + R = Q + R$ und folglich $P \leq A + D$. Ist $B \cap (D + P) \neq
 0$, so folgt entsprechend $P \leq B + D$. Also liegen 
 die Punkte von $A + B + D$ in $A + D$ oder in $B + D$. 
 \par
       G\"abe es nun einen Punkt $U$ auf $A + D$, der nicht in $B + D$
 l\"age, und einen Punkt $V$ auf $B + D$, der nicht in $A + D$
 l\"age, dann g\"abe es, da $L$ ja irreduzibel ist, einen Punkt $W$
 auf $U + V$, der weder zu $A + D$ noch zu $B + D$ geh\"orte. Dies
 widerspr\"ache aber $R \leq A + B + D$. Also ist $A + D \leq B +
 D$ oder $B + D \leq A + D$.
 \medskip\noindent
 {\bf 7.2. Satz.} {\it Es sei $L$ ein irreduzibler projektiver Verband
 mit dem gr\"o\ss ten Element $\Pi$.
 Sind $A$ und $B$ Elemente endlichen Ranges von $L$ und gilt $\Rg_L
 (A) = \Rg_L (B)$, so gibt es ein $C \in L$ mit}
 $$ \Pi = A \oplus C = B \oplus C. $$
 \smallskip
       Beweis. Weil $A$ und $B$ endlichen Rang haben, hat auch $A + B$
 endlichen Rang, so dass der Quotient $(A + B)/0$ noethersch ist.
 Unter allen Teilr\"aumen von $A + B$, die mit $A$ und $B$
 trivialen Schnitt haben, gibt es also einen maximalen Ranges.
 Es sei $D$ ein solcher. Nach 7.1 ist dann oBdA  $A + D \leq B +
 D$. Nun ist aber
 $$ \Rg_L (A + D) = \Rg_L (A) + \Rg_L (D) = \Rg_L (B) + \Rg_L (D)
		  = \Rg_L (B + D), $$
 so dass $A + D = B + D = A + B$ ist.
 \par
       Ist nun $F$ ein Komplement von $A + B$ in $\Pi$, so ist $C := D + F$ ein
 gemeinsames Komplement von $A$ und $B$ in $\Pi$.
 \medskip\noindent
 {\bf 7.3. Korollar.} {\it Ist $L$ ein irreduzibler projektiver Verband, sind
 $A$, $B \in L$ und gilt $\Rg_L (A) = \Rg_L (B) < \infty$, so sind die Quotienten
 $\Pi /A$ und $\Pi /B$ bzw. $A / 0$ und $ B / 0$ isomorph.}
 \smallskip
       Beweis. Nach 7.2 gibt es ein gemeinsames Komplement $C$ von $A$
 und $B$. Daher gilt nach der Transformationsregel
 $$\eqalign{
     \Pi /A &= (A + C)/A \cong C/(A \cap C) = C/0    \cr
	    &= C/(B \cap C) \cong (B + C)/B = \Pi /B \cr} $$
 und
 $$\eqalign{
        A/0 &= A/(A \cap C) \cong (A + C)/C         \cr
	    &= (B + C)/C \cong B/(B \cap C) = B/0,  \cr} $$
 womit das Korollar bewiesen ist.
 \medskip\noindent
 {\bf 7.4. Satz.} {\it Ist $L$ ein projektiver Verband, so sind die folgenden
 Aussagen \"aquivalent:
 \item{a)} $L$ ist irreduzibel oder $L$ ist isomorph zu $(P(A), \subseteq)$, wobei
 $A$ die Menge der Punkte von $L$ bezeichne.
 \item{b)} Ist $k$ eine nat\"urliche Zahl mit $k \leq \Rg(L)$ und sind $X$ und $Y$
 zwei Teilr\"aume des Ranges $k$ von $L$, so sind die Quotienten $X/0$ und $Y/0$
 isomorph.
 \item{c)} Sind $G$ und $H$ zwei Geraden von $L$, so sind $G/0$ und $H/0$
 isomorph.}
 \smallskip
       Beweis. a) impliziert b): Ist $L$ irreduzibel, so gilt b) nach 7.3. Ist
 $L$ zu $(P(A), \subseteq)$ isomorph, so ist $|X| = k = |Y|$, so dass b) auch in
 diesem Falle gilt.
 \par
       b) impliziert nat\"urlich c).
 \par
 c) impliziert a): Gibt es eine Gerade, die mehr als zwei Punkte tr\"agt, so
 tragen alle Geraden mehr als zwei Punkte, so dass $L$ nach 2.14 irreduzibel ist.
 Gibt es keine solche Gerade, so ist $L$ nach Satz 2.15 zu dem cartesischen
 Produkt $\cart_{P \in A} P/0$ isomorph. In jedem Falle gilt also a).
 \medskip
       Wir setzen im Folgenden stets voraus, dass der projektive Verband
 $L$ mindestens den Rang 2 hat und dass alle seine Geraden
 gleichviele Punkte tragen.
 \par
       Sind $G$ und $H$ zwei Geraden des endlichen projektiven Verbandes $L$, so
 tragen die Geraden $G$ und $H$ nach unserer Annahme gleichviele Punkte. Ist
 $q + 1$ diese Anzahl, so hei\ss t $q$ \emph{Ordnung}\index{Ordnung}{} von $L$.
 Ist $L$ irreduzibel, so ist $q \geq 2$. Im anderen Fall ist $q = 1$.
 \par
       Es sei weiterhin $L$ endlich. Wir bezeichnen mit $N_i (k,L)$ die
 Anzahl der Teilr\"aume des Ranges $i$, die in einem Teilraum des
 Ranges $k$ von $L$ enthalten sind. Mit 7.4 folgt auf Grund unserer
 Annahme \"uber die Gleichm\"achtigkeit der Geraden von $L$, dass
 diese Zahl wirklich nur von $i$ und $k$ abh\"angt, nicht jedoch
 von der speziellen Auswahl des Unterraumes des Ranges $k$.
 Offenbar gilt
 $ N_0 (1,L) = 1 = N_1 (1,L)$.
 Weiter gilt

 \medskip
\noindent
 {\bf 7.5. Satz.} {\it Es sei $L$ ein endlicher projektiver Verband des Ranges
 $r$, dessen Geraden alle gleichviele Punkte tragen. Ist $0 \leq i \le j \leq n
 \leq r$, so ist}
 $$ N_i (n,L) N_{j-i} (n-i, L) = N_i (j, L) N_j (n,L). $$

       Beweis. Wie wir schon festgestellt haben, h\"angen die fraglichen
 Zahlen alle nur von den angegebenen Parametern ab.
       Es sei nun $U$ ein Unterraum des Ranges $n$ von $L$. Ferner sei
 $\I$ die Menge der Paare $(A,B)$ von Unterr\"aumen $A$ und $B$ mit
 $\Rg_L (A) = i$, $\Rg_L (B) = j$ und $A \leq B \leq U$. Die Anzahl
 der $B$ vom Range $j$, die ein gegebenes $A$ vom Range $i$
 umfassen und gleichzeitig in $U$ liegen, ist gleich der Anzahl der
 Teilr\"aume des Ranges $j-i$ des Quotienten $U/A$. Ist $C$ ein
 Komplement von $A$ in $U$, so gilt nach der Transformationsregel
 $U/A \cong C/0$. Daher ist diese Anzahl gleich $N_{j-i} (n - i,
 L)$. Also ist
 $$ |\I| = N_i (n,L) N_{j-i} (n-i, L). $$
 Andererseits ist $N_i (j,L)$ die Anzahl der Teilr\"aume vom Range
 $i$, die in einem gegebenen Teilraum des Ranges $j$ enthalten
 sind. Analoges gilt f\"ur $N_j (n,L)$. Daher ist
 $$ |\I| = N_i (j,L) N_{j} (n, L), $$
 so dass der Satz bewiesen ist.
 \medskip
       Es sei $q$ eine nat\"urliche Zahl. Wir definieren die
 $q$-Analoga\index{qanalog@$q$-Analogon}{} der Fakult\"aten\index{Fakult\"at}{} durch
 $(0, q)! := 1$ und
 $$ (n + 1,q)! := (n,q)! \sum_{i:=0}^n q^i $$
 f\"ur $n \geq 0$. Ist $q = 1$, so ist $(n,q)! = n!$.
 \par
       Sind $k$ und $n$ nicht negative ganze Zahlen und gilt $k \leq n$,
 ist ferner $q$ eine nat\"urliche Zahl, so setzen wir
 $$ {n \choose k,q} := {(n,q)! \over (k,q)!(n - k,q)!}. $$
 Ist $q = 1$, so sind dies die
 Binomialkoeffizienten.\index{Binomialkoeffizient}{} Ist $q > 1$, so hei\ss en
 diese Zahlen \emph{gau\ss sche Zahlen}.\index{gau\ss sche Zahlen}{}
 \par
       Ist $q >1$, so gilt
 $$ {n \choose k,q} = \prod_{i := 1}^{k} {q^{n+1-i} - 1 \over q^i - 1}. $$
 Diese Darstellung ist f\"ur das Auge gef\"alliger.
 \medskip\noindent
 {\bf 7.6. Satz.} {\it Es sei $L$ ein endlicher projektiver Verband,
 dessen Ge\-ra\-den alle gleichviele Punkte tragen. Es sei $q$ die
 Ordnung und $r$ der Rang von $L$. Sind dann $k$ und $n$ nicht
 negative ganze Zahlen mit $k \leq n \leq r$, so gilt
 $$ N_k(n,L) = {n \choose k,q}. $$
 Insbesondere h\"angen die Zahlen $N_k(r,L)$ also nur von $k,r$ und $q$ ab.}
 \smallskip
       Beweis. Nach 7.4 sind auf Grund unserer Annahme alle Unterr\"aume
 gegebenen Ranges von $L$ isomorph. Ist $U$ ein Unterraum des
 Ranges $n$ und ist $P$ ein Punkt von $U$, so ist die Anzahl der
 Geraden durch $P$, die in $U$ liegen, gleich der Anzahl der
 Unterr\"aume des Ranges $1$ in $U/P$. Ist $U = P \oplus C$, so
 folgt mittels der Transformationsregel, dass diese Anzahl gleich
 der Anzahl der Punkte auf $C$, dh., dass sie gleich $N_1(n-1,L)$
 ist. Auf jeder Geraden durch $P$ liegen $q$ Punkte, die von $P$
 verschieden sind, und je zwei dieser Geraden haben nur $P$
 gemeinsam. Da schlie\ss lich jeder Punkt von $U$ mit $P$ durch
 eine Gerade verbunden ist, ist
 $$ N_1 (n,L) = qN_1 (n-1,L) + 1. $$
 Weil $N_1(1,L) = 1$ ist, folgt mittels Induktion, dass
 $$\textstyle N_1 (n,L) = \sum_{i:=0}^{n-1} q^i = {n \choose 1,q} $$
 ist. Setzt man in 7.5 nun $i:=1$ und $j := k$, so folgt
 $$\textstyle( \sum_{i:=0}^{n-1} q^i ) N_{k-1} (n - 1, L) =
          ( \sum_{i:=0}^{k-1} q^i ) N_k(n,L). $$
 Vollst\"andige Induktion liefert nun das gew\"unschte Ergebnis.
 \medskip
       Der soeben bewiesene Satz besagt unter anderem, dass die
 Binomialkoeffizienten\index{Binomialkoeffizient}{} wie auch die gau\ss schen
 Zahlen\index{gau\ss sche Zahlen}{} nat\"urliche Zahlen sind, letztere zumindest
 f\"ur solche $q$, die Ordnung eines endlichen projektiven Raumes sind, was nicht
 f\"ur alle nat\"urlichen Zahlen zutrifft. Es ist jedoch eine einfache
 \"Ubungsaufgabe zu zeigen, dass diese einschr\"ankende Voraussetzung an $q$
 nicht gemacht werden muss. F\"ur die gau\ss schen Zahlen gilt n\"amlich die
 Funktionalgleichung
 $$ {n + 1 \choose k,q} = q {n \choose k,q} + {n \choose k - 1,q}. $$
 \medskip\noindent
 {\bf 7.7. Satz.} {\it Es sei $L$ ein endlicher projektiver Verband des Ranges
 $r$, dessen Geraden alle gleichviele Punkte tragen. Ist $0 \leq k \leq r$, so 
 ist}
 $$ N_k (r,L) = N_{r-k} (r,L). $$
 \par
       Beweis. Weil alle Geraden von $L$ gleichviele Punkte tragen, ist
 $L$ nach 7.4 irreduzibel oder aber zum Potenzmengenverband einer
 geeigneten Menge isomorph. Weil $L$ endlich ist, ist im ersten
 Fall $L^d$ nach Fr\"uherem projektiv. Im zweiten Fall ist $L^d$
 ebenfalls projektiv. Die Komplementbildung vermittelt ja einen
 Isomorphismus des Potenzmengenverbandes auf seinen dualen Verband.
 Wir d\"urfen daher im Folgenden $L^d$ in unserer Argumentation verwenden.
 \par
       Es sei nun $G$ eine Gerade von $L$ und $H$ sei ein Komplement von
 $G$ in $\Pi$. Dann ist $H$ eine Gerade von $L^d$. Aus der
 Isomorphie von $G/0$ und $\Pi/H$ folgt daher, dass $L$ und $L^d$
 die gleiche Ordnung haben. Da diese Verb\"ande nach 5.9 den
 gleichen Rang haben und die Unterr\"aume vom Ko-Rang $r-k$ von $L$
 gerade die Unterr\"aume vom Range $k$ von $L^d$ sind, folgt mit
 7.6
 $$ N_k(r,L) = N_k (r,L^d) = N_{r-k} (r,L).  $$
 \par
       Wir werden sp\"ater sehen, dass $q$ im Falle $r \geq 4$ stets
 Potenz einer Primzahl ist. Im Falle $r = 3$, dh. im Falle, dass
 $L$ eine projektive Ebene ist, ist das eine bislang unbewiesene
 Vermutung. Man wei\ss\ nach einem Satz von Bruck und
 Ryser,\index{Satz von Bruck und Ryser}{} dass
 eine Zahl $q$, die kongruent $1$ oder $2$ modulo $4$ ist und die
 sich nicht als Summe von zwei Quadraten darstellen l\"asst,
 niemals die Ordnung einer endlichen projektiven Ebene ist.
 \medskip\noindent
 {\bf 7.8. Satz.} {\it Ist $L$ ein endlicher, irreduzibler projektiver
 Verband der Ordnung $q$ und des Ranges $r$, so ist die Anzahl
 $A_S(r,q)$ von geordneten $s$-Tupeln unabh\"angiger Punkte gleich
 $$ (q - 1)^{-s} q^{{1 \over 2} s(s-1)} \prod_{i:= r+1-s}^r (q^i - 1). $$
 Insbesondere ist die Anzahl der geordneten Basen von $L$ gleich}
 $$ (q - 1)^{-r} q^{{1 \over 2} r(r-1)} \prod_{i:=1}^r (q^i - 1). $$
 \smallskip
       Beweis. Es ist $A_1(r,q) = N_1(r,q) = (q - 1)^{-1} (q^r - 1)$.
 Es sei bereits bewiesen, dass $A_s(r,q)$ die angegebene Form hat.
 Es sei $(P_1, \dots, P_s)$ ein unabh\"angiges $s$-Tupel von
 Punkten. Dann ist das $(s+1)-$Tupel $(P_1, \dots, P_s, Q)$ genau
 dann unabh\"angig, wenn $Q \not\leq \sum_{i:=1}^{s} P_i$
 gilt. Somit ist die Anzahl dieser $(s+1)$-Tupel gleich
 $$ N_1(r,q) - N_1(s,q). $$
 Folglich ist
 $$A_{s+1} (r,q) = A_s(r,q)(N_1(r,q) - N_1(s,q)). $$
 Hieraus folgt mit einer simplen Rechnung die Behauptung.
 \medskip
       Es sei $L$ ein projektiver Verband des Ranges $r$. Sind $P_1$, \dots,
 $P_{r+1}$ Punkte von $L$ mit der Eigenschaft, dass je $r$ von ihnen unabh\"angig
 sind, also eine Basis von $L$ bilden, so hei\ss t $(P_1, \ldots, P_{r + 1})$
 \emph{Rahmen}\index{Rahmen}{} von $L$.
 \medskip\noindent
 {\bf 7.9. Satz.} {\it Ist $L$ ein endlicher irreduzibler projektiver
 Verband des Ranges $r$ und der Ordnung $q$, so ist die Anzahl der
 Rahmen von $L$ gleich}
 $$ q^{{1 \over 2} r(r-1)} \prod^r_{i:=2} (q^i - 1). $$
 \par
       Beweis. Es sei $P_1$, \dots, $P_r$ eine Basis. Dann ist die Anzahl
 der Rahmen $(P_1, \ldots, P_r, X)$ gleich $(q - 1)^{r-1}$ : Dies
 ist sicherlich richtig, falls $r = 2$ ist, da in diesem Falle
 jedes Punktetripel ein Rahmen ist. Es sei nun $r > 2$. Es ist eine
 einfache \"Ubungsaufgabe zu zeigen, dass $(P_1, \ldots, P_r, X)$
 genau dann ein Rahmen ist, wenn $P_r \neq X$ und $X \not \leq
 \sum_{i:=1}^{r-1} P_i$ ist und wenn $(P_1, \ldots,
 P_{r-1}, X')$ mit $X' := (P_r + X) \cap \sum_{i:=1}^{r-1}
 P_i$ ein Rahmen von $\sum_{i:=1}^{r-1} P_i/0$ ist.
 Hieraus folgt, dass die Anzahl der Rahmen $(P_1, \ldots, P_r, X)$
 gleich $q-1$ mal der Anzahl der Rahmen $(P_1, \ldots, P_{r-1},
 X')$, dh., gleich $(q-1)^{r-1}$ ist. Somit ist die Anzahl der
 Rahmen von $L$ gleich $(q-1)^{r-1} A_r(r,q)$. Hieraus folgt mit
 7.7 die Behauptung des Satzes.

\mysection{8. Kollineationen und Korrelationen}

\noindent
 In Abschnitt 1 hatten wir den Begriff des Isomorphismus einer projektiven
 Geometrie auf eine andere definiert und in Abschnitt 2 den des
 Isomorphismus\index{Isomorphismus}{} zwischen zwei Verb\"anden. Dass es bei
 projektiven Geometrien\index{projektive Geometrie}{} nicht darauf ankommt,
 welchen der beiden
 Isomorphiebegriffe man der Theorie zu Grunde legt, zeigt der folgende Satz.
 \medskip\noindent
 {\bf 8.1. Satz.} {\it Es seien $L$ und $L'$ zwei projektive
 Verb\"ande, es sei $\Sigma$ die projektive Geometrie aus den
 Punkten und Geraden von $L$ und $\Sigma'$ die projektive Geometrie
 aus den Punkten und Geraden von $L'$. Ist nun $\sigma$ ein
 Isomorphismus von $\Sigma$ auf $\Sigma'$, so gibt es genau einen
 Isomorphismus $\tau$ von $L$ auf $L'$ mit $P^\sigma = P^\tau$ und
 $G^\sigma = G^\tau$ f\"ur alle Punkte $P$ und alle Geraden $G$ von
 $\Sigma$. Ist umgekehrt $\tau$ ein Isomorphismus von $L$ auf $L'$
 und definiert man $\sigma$ durch $P^\sigma := P^\tau$ und
 $G^\sigma := G^\tau$ f\"ur alle Punkte $P$ und alle Geraden $G$
 von $\Sigma$, so ist $\sigma$ ein Isomorphismus von $\Sigma$ auf $\Sigma'$.}
 \smallskip
       Beweis. Es sei $U \in L$ und $S(U)$ bezeichne wieder die Menge der
 in $U$ enthaltenen Atome. Dann ist, wie wir wissen, $U = \sum_{P \in S(M)} P$.
 Wir definieren $\tau$ durch $U^\tau := \sum_{P \in S(M)} P^\sigma$. Dann ist
 $\tau$ eine inklusionstreue Abbildung von $L$ in $L'$.
 \par
       Es sei $V \in L'$. Wir setzen
 $$ A := \bigl\{P \mid P \in S(\Pi), P^\sigma \in S(V)\bigr\} \quad
 \mathrm{und}\quad
  U := \sum_{P \in A} P. $$
 Dann ist $A \subseteq S(U)$ und folglich $V \leq U^\tau$. Wir zeigen, dass
 $V = U^\tau$ ist, indem wir zeigen, dass $A = S(U)$ gilt. Dazu sei $Q$ ein Punkt
 von $U$. Es gibt dann endlich viele Punkte $P_1$, \dots, $P_n \in A$ mit
 $Q \leq \sum_{i:=1}^n P_i$. Ist $Q = P_1$, so ist $Q \in A$.
 Wir d\"urfen daher annehmen, dass $Q$ von $P_1$ verschieden ist.
 Wegen $Q \leq P_1 + \sum_{i:=2}^n P_i$ gibt es nach 2.6
 einen Punkt $R$ in $\sum_{i:=2}^n P_i$ mit $Q \leq P_1 +
 R$. Nach Induktionsannahme ist $R \in A$, so dass wir $R = P_2$
 annehmen d\"urfen. Weil $\sigma$ ein Isomorphismus ist, ist
 $Q^\sigma $ ein Punkt auf $P_1^\sigma + P_2^\sigma$. Weil
 Letzteres eine Gerade von $V$ ist, ist $Q^\sigma$ ein Punkt von
 $V$, so dass $Q \in A$ gilt. Also ist $A = S(U)$ und damit $U^\tau
 = V$. Folglich ist $\tau$ surjektiv. Weil $\sigma$ injektiv und
 $L$ relativ atomar ist, ist aber auch $\tau$ injektiv. Somit ist
 $\tau$ eine inklusionstreue Bijektion von $L$ auf $L'$. Wendet man
 die gerade ausgef\"uhrte Argumentation auf $\sigma^{-1}$ an, so
 sieht man, dass auch $\tau^{-1}$ inklusionstreu ist, so dass $\tau$
 in der Tat ein Isomorphismus von $L$ auf $L'$ ist.
 \par
       Die Einzigkeit von $\tau$ ist banal, wie auch die zweite Aussage
 des Satzes.
 \medskip
       Sind $L$ und $L'$ projektive Verb\"ande und ist $\sigma$ ein Isomorphismus
 von $L$ auf $L^{'d}$, so hei\ss t $\sigma$
 \emph{Antiisomorphismus}\index{Antiisomorphismus}{} von $L$ auf $L'$. Nach 8.1
 kann ein Antiisomorphismus $\sigma$ auch dadurch beschrieben werden, dass
 $\sigma$ eine Bijektion der Menge der Punkte von $L$ auf die Menge
 der Hyperebenen von $L'$ ist, so dass drei Punkte $P$, $Q$ und $R$
 von $L$ genau dann kollinear sind, wenn der Schnitt der
 Hyperebenen $P^\sigma$, $Q^\sigma$ und $R^\sigma$ einen Unterraum
 vom Ko-Rang 2 enth\"alt, wenn also die Hyperebenen $P^\sigma$,
 $Q^\sigma$ und $R^\sigma$ so zueinander liegen wie die Bl\"atter eines Buches.
 --- So beschrieb E. Witt\index{Witt, E.}{} diese Situation einmal in
 einem Telefongespr\"ach mit mir.
 \par
       Ein Antiisomorphismus von $L$ auf sich hei\ss t
 \emph{Korrelation}\index{Korrelation}{}, und eine involutorische Korrelation
 hei\ss t \emph{Polarit\"at}.\index{Polarit\"at}{}
 \smallskip
       Es sei $L$ ein projektiver Verband und $\UR_r(L)$ bezeichne die
 Menge der Unterr\"aume vom Range $r$ von $L$. Ein Isomorphismus
 von $L$ auf $L'$ bildet unabh\"angige Punktmengen auf
 unabh\"angige Punktmengen ab. Daher induziert jeder Isomorphismus
 von $L$ auf $L'$ eine bijektive Abbildung von $\UR_r(L)$ auf
 $\UR_r(L')$. Dies ist besonders dann interessant, wenn $r$ endlich
 ist. Ist $\Rg (L) = n$ endlich und ist $\sigma$ ein
 Antiisomorphismus von $L$ auf $L'$, so induziert $\sigma$ eine
 Bijektion von $\UR_r(L)$ auf $\UR_{n-r}(L')$. Wir werden nun die
 Frage beantworten, wie die von Isomorphismen bzw.
 Antiisomorphismen induzierten Abbildungen von $\UR_r(L)$ auf
 $\UR_s(L')$ gekennzeichnet werden k\"onnen.
 \par
       Es sei $r$ eine nat\"urliche Zahl. Sind $X$, $Y \in \UR_r (L)$ und
 ist $\Rg_L(X+Y) = r + a$, so nennen wir $a$ den
 \emph{Abstand}\index{Abstand}{} von $X$ und $Y$. Nach 3.7 gilt
 $$ \Rg_L (X + Y) + \Rg_L (X \cap Y) = \Rg_L (X) + \Rg_L (Y) = 2r. $$
 Also ist $\Rg_L (X \cap Y) = r - a$, falls $a$ der Abstand von $X$ und $Y$ ist.
 Man nennt $X$ und $Y$ {\it benachbart\/},\index{benachbart}{} falls $X$ und $Y$
 den Abstand 1 haben.  Eine Menge $M \subseteq \UR_r (L)$ von paarweise
 benachbarten Unterr\"aumen hei\ss t \emph{maximal},\index{maximal}{} wenn aus
 $M \subseteq M' \subseteq \UR_r(L)$ und der Eigenschaft, dass die Elemente aus
 $M'$ paarweise benachbart sind, $M = M'$ folgt. Diese maximalen
 Mengen lassen sich nun wie folgt kennzeichnen.
 \medskip\noindent
 {\bf 8.2. Satz.} {\it Es sei $L$ ein projektiver Verband. Ferner sei
 $r$ eine na\-t\"ur\-li\-che Zahl mit $2 \leq r + 1 \leq \Rg(L)$. Ist $M
 \subseteq \UR_r (L)$ eine maximale Menge von paarweise
 benachbarten Unterr\"aumen, so gibt es entweder einen Unterraum
 $U$ vom Range $r-1$ mit
 $$ M = \bigl\{X \mid X \in \UR_r (L),\ U \leq X \bigr\} $$
 oder einen Unterraum $W$ vom Range $r + 1$ mit
 $$ M = \bigl\{X \mid X \in \UR_r (L),\ X \leq W \bigr\}. $$
 Die beiden F\"alle schlie\ss en sich gegenseitig aus, falls $\Rg (L) > 2$ ist.}
 \smallskip
       Beweis. Es seien $X$, $Y$, $Z$ drei verschiedene Elemente aus $M$.
 Ferner sei $Z \not \leq X + Y$. Dann ist
 $$ \Rg_L\bigl(Z \cap (X + Y)\bigr) \leq r - 1. $$
 Andererseits ist
 $$ Z \cap X \leq Z \cap (X + Y) $$
 und
 $$ Z \cap Y \leq Z \cap (X + Y). $$
 Weil die R\"ange der links stehenden R\"aume gleich $r - 1$ sind, folgt
 $Z \cap X = Z \cap Y = Z \cap (X + Y)$. Hieraus folgt weiter
 $Z \cap X \cap Y = Z \cap X$. Daher ist
 $$ \Rg_L (Z \cap X \cap Y) = r - 1 = \Rg_L(X \cap Y). $$
 Dies hat $X \cap Y \leq Z$ zur Folge.
 \par
       Es sei $S$ ein viertes Element aus $M$. Ist $S \not\leq X + Y$, so
 ist $X \cap Y \leq S$, wie wir gerade gesehen haben. Es sei also
 $S \leq X + Y$. Dann ist $S + X = X + Y$, da $S$ und $X$, bzw.,
 $X$ und $Y$ ja benachbart sind. Andererseits ist $X + Y \neq X + Z$
 und daher $X = (X + Y) \cap (X + Z)$. Folglich ist $S \not\leq X +
 Z$. Nach dem bereits Bewiesenen ist daher $X \cap Y = X \cap Z \leq S$.
 Damit ist gezeigt, dass entweder alle Elemente von $M$ unterhalb $X + Y$ oder
 dass alle Elemente von $M$ oberhalb von $X \cap Y$ liegen. Aus der Maximalit\"at
 von $M$ folgt damit die erste Behauptung des Satzes.
 \par
       Um die zweite Behauptung zu beweisen, nehmen wir an, dass
 $$ M = \{Z \mid Z \in L, X \cap Y < Z < X + Y\} $$
 gilt. Es sei $P$ ein Punkt, der nicht in $X \cap Y$ liege. Aus der Maximalit\"at
 von $M$ folgt, dass $(X \cap Y) + P \in M$ ist. Dann ist aber
 $$ P \leq (X \cap Y) + P \leq X + Y. $$
 Es folgt, dass $X + Y = \Pi$ ist. Andererseits folgt aus der Maximalit\"at von
 $M$, dass jeder Teilraum des Ko-Ranges 1 in $X + Y$ zu $M$ geh\"ort. Hieraus
 folgt $X \cap Y = 0$, so dass $r-1 = 0$ ist. Somit ist $\Rg (L) = 2$.
 Damit ist dann alles bewiesen.
 \medskip\noindent
 {\bf 8.3. Satz.} {\it Es sei $L$ ein projektiver Verband und $r$ sei
 eine nat\"urliche Zahl. Sind $X, Y \in \UR_r (L)$ und ist $a$ der
 Abstand von $X$ und $Y$, so gibt es $X_0, \dots, X_a \in \UR_r(L)$ mit
 $X = X_0$ und $X_a = Y$, so dass $X_i$ und $X_{i+1}$
 f\"ur $i := 0, \dots, a - 1$ benachbart sind.
 \par
       Sind umgekehrt $Y_0, \dots, Y_n \in \UR_r (L)$, ist $X = Y_0$ und
 $Y_n = Y$ und sind $Y_i$ und $Y_{i+1}$ f\"ur alle in Frage kommenden $i$
 benachbart, so ist $a \leq n$.}
 \smallskip
       Beweis. Dies ist sicherlich richtig, falls $a \leq 1$ ist. Es sei
 also $a \geq 2$. Es sei $P$ ein Punkt von $X$, der nicht in $Y$
 liegt. Dann ist
 $$ \Rg_L\bigl(P + (X \cap Y)\bigr) = 1 + \Rg_L (X \cap Y) = r - a + 1. $$
 Weil $X \cap Y$ den Ko-Rang $a$ hat, gibt es einen Teilraum $Q$ von $Y$ des
 Ranges $a - 1$ mit $X \cap Y \cap Q = \{0\}$. Setze
 $$ Z := P + (X \cap Y) + Q. $$
 Dann ist der Rang von $Z$ gleich $r$. Mit Hilfe des modularen Gesetzes folgt
 ferner
 $$ Z \cap X = P + (X \cap Y) $$
 und
 $$ Z \cap Y = (X\cap Y) + Q. $$
 Somit ist
 $$ \Rg_L (Z \cap X) = r - a + 1 $$
 und
 $$\Rg_L(Z \cap Y) = r - 1. $$
 Hieraus folgt mittels Induktion die erste Behauptung.
 \par
       Die zweite Aussage gilt sicherlich, falls $n \leq 1$. Es sei also
 $n \geq 2$. Ferner sei $b$ der Abstand von $X$ und $Y_{n-1}$. Nach
 Induktionsannahme ist dann $b \leq n - 1$, so dass wir $b < a$
 annehmen d\"urfen. Nun ist
 $$ \Rg_L\bigl((X \cap Y) + (Y_{n-1} \cap Y)\bigr) = 2r - a - 1 - 
		\Rg_L(X \cap Y \cap Y_{n-1}). $$
 Andererseits ist
 $$ \Rg_L\bigl((X \cap Y_{n-1}) + (Y_{n-1} \cap Y)\bigr) = 2r - b - 1 -
		\Rg_L(X \cap Y \cap Y_{n-1}). $$
 Nun ist
 $$ Y_{n-1} \cap Y \leq (X \cap Y) + (Y_{n-1} \cap Y) \leq Y $$
 und
 $$Y_{n-1} \cap Y \leq (X \cap Y_{n-1}) + (Y_{n-1} \cap Y) \leq Y_{n-1}. $$
 Daher gilt
 $$ r - 1 \leq \Rg_L\bigl((X \cap Y) + (Y_{n-1} \cap Y)\bigr) \leq r $$
 und
 $$ r - 1 \leq \Rg_L\bigl((X \cap Y_{n-1}) + (Y_{n-1} \cap Y)\bigr) \leq r. $$
 Wegen $b < a$ ist daher
 $$\Rg_L\bigl((X \cap Y) + (Y_{n-1} \cap Y)\bigr) = r - 1 $$
 und
 $$ \Rg_L\bigl((X \cap Y_{n-1}) + (Y_{n-1} \cap Y)\bigr) = r. $$
 Somit ist $a = b + 1$ und folglich $a \leq n$. Damit ist alles bewiesen.
 \medskip
       Ist $\sigma$ ein Isomorphismus des projektiven Verbandes $L$ auf
 den projektiven Verband $L'$, so induziert $\sigma$ eine Bijektion
 von $\UR_r(L)$ auf $\UR_r(L')$ mit der Eigenschaft, dass $X$, $Y \in
 \UR_r(L)$ genau dann benachbart sind, wenn $X^\sigma$ und $Y^\sigma$ benachbart
 sind. Da zwei verschiedene Punkte,
 bzw. zwei verschiedene Hyperebenen stets benachbart sind, hat im
 Falle $r = 1$ oder $r + 1 = \Rg(L)$ jede Bijektion von $\UR_r(L)$
 auf $\UR_r(L')$ die Eigenschaft, die Nachbarschaft zu erhalten.
 Dies zeigt, dass die Annahme \"uber $r$ im Folgenden, im Wesentlichen von
 W. L. Chow\index{Satz von Chow}{} stammenden Satz wirklich n\"otig
 ist (Chow 1949). Chow setzte von Anfang an voraus, dass $r = s$ sei.
 \medskip\noindent
 {\bf 8.4. Satz.} {\it Sind $L$ und $L'$ projektive Verb\"ande, sind
 $r$ und $s$ nat\"urliche Zahlen und ist $2 < r + 1 < \Rg(L)$, ist
 ferner $\sigma$ eine Bijektion von $\UR_r(L)$ auf $\UR_s(L')$ mit
 der Eigenschaft, dass $X$, $Y \in \UR_r(L)$ genau dann benachbart
 sind, wenn $X^\sigma$ und $Y^\sigma$ benachbart sind, so wird $\sigma$ entweder
 durch einen Isomorphismus von $L$ auf $L'$ induziert, in welchen
 Falle $r = s$ ist, oder aber $\sigma$ wird durch einen
 Antiisomorphismus von $L$ auf $L'$ induziert, in welchem Falle
 $r + s = \Rg(L)$ ist. Die beiden F\"alle schlie\ss en sich
 gegenseitig aus. \"Uberdies gilt, dass der $\sigma$ induzierende
 Isomorphismus bzw. Antiisomorphismus eindeutig bestimmt ist.}
 \smallskip
       Beweis. Es sei $\sigma$ eine Bijektion von $\UR_r(L)$ auf
 $\UR_s(L')$, so dass $X$, $Y \in \UR_r(L)$ genau dann benachbart
 sind, wenn $X^\sigma$ und $Y^\sigma$ es sind.  Dann ist auch
 $\sigma^{-1}$ eine nachbarschaftserhaltende Abbildung von
 $\UR_s(L')$ auf $\UR_r(L)$. Daher bildet $\sigma$ maximale Mengen
 von paarweise benachbarten Unterr\"aumen wieder auf solche ab und
 jede maximale Menge von paarweise benachbarten R\"aumen aus
 $\UR_s(L')$ ist Bild einer ebensolchen aus $\UR_r(L)$.
 \par
       Es sei $U \in \UR_{r-1}$ und $M := \{X \mid X \in \UR_r(L), U \leq X\}$.
 Es gibt dann nach 8.2 genau ein $U' \in \UR_{s-1} (L') \cup \UR_{s+1} (L')$ mit
 $$ M^\sigma = \bigl\{X' \mid X' \in \UR_s(L'),\ U' \leq X'\bigr\} $$
 oder
 $$ M^\sigma = \bigl\{X' \mid X' \in \UR_s(L'),\ X' \leq U'\bigr\}. $$
 Ist $U \in \UR_{r+1} (L)$, so findet man zu $M := \{X \mid X \in \UR_r (L),
 X \leq U\}$ ebenfalls genau ein $U' \in \UR_{s-1} (L') \cup \UR_{s+1} (L')$ mit
 $$ M^\sigma = \bigl\{X' \mid X' \in \UR_s(L'),\ U' \leq X'\bigr\} $$ 
 oder
 $$ M^\sigma = \bigl\{X' \mid X' \in \UR_s(L'),\ X' \leq U'\bigr\}. $$
 In beiden F\"allen setzen wir $U^\sigma := U'$. Dann ist $\sigma$
 nach der zuvor gemachten Bemerkung eine Bijektion von
 $$ \UR_{r-1} (L) \cup \UR_r (L) \cup \UR_{r+1}(L) $$
 auf
 $$ \UR_{s-1} (L') \cup \UR_s (L') \cup \UR_{s+1}(L'). $$
 Wegen $\UR_r (L)^\sigma = \UR_s (L')$ gilt
 $$ \bigl(\UR_{r-1} (L) \cup \UR_{r+1} (L)\bigr)^\sigma =
	    \UR_{s-1} (L') \cup \UR_{s+1} (L') .$$
 Angenommen es seien $X$, $Y \in \UR_{r-1} (L)$ und es gelte
 $X^\sigma \in \UR_{s-1} (L')$ und $Y^\sigma \in \UR_{s+1} (L')$.
 Nach 8.3 gibt es endlich viele $X_i \in \UR_{r-1} (L)$ mit $X =
 X_0$ und $Y = X_n$, so dass $X_i$ und $X_{i+1}$ benachbart sind.
 Es gibt dann ein $i$, so dass $X^\sigma_i \in \UR_{s-1} (L')$ und
 $X^\sigma_{i+1} \in \UR_{s+1} (L')$ gilt. Wir d\"urfen daher
 annehmen, dass $X$ und $Y$ benachbart sind.
 \par
       Es sei $M$ die Menge der Unterr\"aume des Ranges $r$, die $X$
 enthalten. Auf Grund der Definition von $X^\sigma$ ist dann
 $$ M^\sigma = \bigl\{X' \mid X' \in \UR_s (L'),\ X^\sigma \leq X'\bigr\}. $$
 Ist $N$ die Menge der Unterr\"aume des Ranges $r$, die $Y$ enthalten, so folgt
 entsprechend
 $$ N^\sigma = \bigl\{Y' \mid Y' \in \UR_s (L'),\ Y' \leq Y^\sigma\bigr\}.$$
 Ist nun $X^\sigma < Z' < Y^\sigma$, so gibt es ein $Z \in M \cap N$ mit
 $Z^\sigma = Z'$. Nun ist aber $M \cap N = \{X + Y\}$, so dass es zwischen
 $X^\sigma$ und $Y^\sigma$ genau einen Unterraum des Ranges $s$ gibt. Dies
 widerspricht aber der relativen Atomarit\"at des Verbandes $L'$. Also ist doch
 $\Rg_{L'} (X^\sigma) = \Rg_{L'} (Y^\sigma)$, falls nur $X$, $Y \in \UR_{r-1}(L)$
 gilt. Ganz analog folgt auch f\"ur $X$, $Y \in \UR_{r+1} (L)$ die Gleichheit der
 R\"ange von $X^\sigma$ und $Y^\sigma$. Es sind also die folgenden beiden F\"alle
 zu betrachten:
 \smallskip\noindent
 1. $\UR_{r-1}(L)^\sigma = \UR_{s-1}(L')$ und $\UR_{r+1}(L)^\sigma =
 \UR_{s+1} (L').$
 \smallskip\noindent
 2. $\UR_{r-1} (L)^\sigma = \UR_{s+1}(L')$ und $\UR_{r+1}(L)^\sigma =
 \UR_{s-1} (L').$
 \smallskip
       Wegen $2 < r + 1 < \Rg (L)$ enthalten $\UR_{r-1}(L)$ und $\UR_{r+1}(L)$
 je mindestens zwei Elemente. Daher gilt dies auch f\"ur $\UR_{s-1} (L)$ und
 $\UR_{s+1} (L')$, so dass auch $2 < s+1 <\Rg (L')$ gilt.
 \par
       1. Fall. Hier gilt: Sind $X$, $Y \in \UR_{r-1}(L) \cup \UR_r(L) \cup
 \UR_{r+1}(L)$, so ist genau dann $X \leq Y$, wenn $X^\sigma \leq Y^\sigma$ ist.
 Insbesondere folgt, dass zwei Unterr\"aume des Ranges $r - 1$ von $L$ genau dann
 benachbart sind, wenn ihre Bilder es sind.
 \par
       Es sei $r = 2$. Unter Zuhilfenahme von $\sigma^{-1}$ folgt, dass in
 diesem Falle auch $s=2$ ist. Daher bildet $\sigma$ Punkte auf
 Punkte und Geraden auf Geraden unter Erhaltung der Inzidenz ab.
 Nach 8.1 wird $\sigma$ folg\-lich von einem Isomorphismus von $L$
 auf $L'$ induziert. Vollst\"andige Induktion f\"uhrt nun zum Ziel.
 \par
       2. Fall. In diesem Falle gilt: sind $X$, $Y \in \UR_{r-1} (L) \cup
 \UR_r (L) \cup \UR_{r+1} (L)$, so ist genau dann $X \leq Y$, wenn
 $Y^\sigma \leq X^\sigma$ ist. Benachbarte Elemente werden also
 auch hier auf benachbarte Elemente abgebildet. Ist nun $r=2$, so
 erschlie\ss t man mit Hilfe von $\sigma^{-1}$, dass die $X' \in
 \UR_s (L')$ den Ko-Rang 2 haben. Also ist in diesem Falle $r+s =
 \Rg(L')$. Es folgt ferner, dass $\sigma$ kollineare Punkte auf
 \emph{konfluente}\index{konfluente Hyperebenen}{} Hyperebenen abbildet. Nach 8.1
 wird $\sigma$
 folglich von einem Antiisomorphismus induziert. Induktion f\"uhrt
 nun auch hier zum Ziele. Wird $\sigma$ durch einen Isomorphismus
 induziert, so liegt offensichtlich der erste Fall vor und es
 folgt, dass $\sigma$ nicht durch einen Antiisomorphismus induziert
 wird. Wird $\sigma$ durch einen Antiisomorphismus induziert, so
 folgt, dass $\sigma$ nicht durch einen Isomorphismus induziert
 wird. Die beiden F\"alle schlie\ss{}en sich also gegenseitig aus.
 \par
       Es seien nun $\sigma$ und $\tau$ zwei Isomorphismen oder zwei
 Antiisomorphismen von $L$ auf $L'$, die auf $\UR_r(L)$ die gleiche
 Abbildung induzieren. Dann ist $\sigma \tau^{-1}$ ein
 Automorphismus, welcher auf $\UR_r (L)$ die Identit\"at induziert.
 Da jeder Punkt von $L$ sich als Schnitt von Unterr\"aumen des
 Ranges $r$ darstellen l\"asst, folgt, dass $\sigma \tau^{-1}$
 alle Punkte von $L$ festl\"asst. Daher ist $\sigma = \tau$,
 womit alles bewiesen ist.
 \medskip
       Der Satz von Chow sowie 8.1 sind Beispiele daf\"ur, dass eine
 Abbildung eines Teils von $L$ auf einen Teil von $L'$ unter
 gewissen Voraussetzungen durch einen Isomorphismus von $L$ auf
 $L'$ induziert wird. Zwei weitere Beispiele dieser Art werden wir
 noch sehen. Zuvor jedoch wollen wir zeigen, dass der Satz von Chow
 auch in anderen Teilen der Mathematik von Nutzen ist, indem wir
 nach dem Vorgange von H. M\"aurer\index{maurer@M\"aurer, H.}{} einen Satz von
 H. Wielandt\index{Wielandt, H.}{} beweisen.
 \par
       Ist $X$ eine Menge, so bezeichnen wir mit $S_X$ die
 \emph{symmetrische Gruppe}\index{symmetrische Gruppe}{} auf $X$, dh., die Gruppe
 aller Bijektionen von $X$ auf sich. Ferner bezeichnen wir mit $A_X$ die
 \emph{alternierende Gruppe}\index{alternierende Gruppen}{} auf $X$, das ist die
 von allen Zyklen der L\"ange 3 erzeugte Untergruppe von $S_X$. Dann ist $A_X$
 ein Normalteiler von $S_X$.
 \medskip\noindent
 {\bf 8.5. Satz.} {\it Es sei $X$ eine Menge, die nicht genau sechs
 Punkte ent\-hal\-te. Ist dann $\alpha$ ein Automorphismus von $A_X$,
 so gibt es ein $\gamma \in S_X$ mit $\alpha (\xi) = \gamma \xi
 \gamma^{-1}$ f\"ur alle $\xi \in A_X$.}
 \smallskip
       Beweis. Ist $|X| \leq 2$, so ist $A_X = \{1\}$ und daher nichts zu
 beweisen. Ist $|X| = 3$, so ist $A_X$ zyklisch der Ordnung 3 und
 hat somit genau zwei Automorphismen, die beide durch innere
 Automorphismen der $S_X$ induziert werden, da die $S_X$ ja nicht abelsch ist.
 \par 
       Im folgenden enthalte $X$ mindestens vier Elemente. Ist $i$ eine
 nat\"urliche Zahl, so bezeichnen wir mit $Z_i$ die Menge der
 Elemente aus $A_X$, die Produkt von genau $i$ disjunkten Zyklen
 sind, die alle die L\"ange 3 haben. Die $Z_i$ sind
 Konjugiertenklassen von $A_X$, so dass $\alpha$ jedes $Z_i$ auf
 ein $Z_j$ abbildet. Wir zeigen zun\"achst, dass $\alpha (Z_1) = Z_1$ ist.
 \par
       Weil die Elemente der $A_5$ nur die Ordnung 1, 2, 3 oder 5
 haben, haben auch die Elemente aus dem Komplexprodukt $Z_1 Z_1$
 nur diese Ordnungen. Damit $Z_i$ nicht leer ist, muss $X$
 mindestens $3i$ Elemente haben. F\"ur $i \geq 2$ bedeutet dies auf
 Grund unserer Annahme \"uber $X$, dass $X$ in diesem Falle
 mindestens sieben Elemente enth\"alt. Wegen
 $$ (123)(124)(567)(567) = (14)(23)(576) $$
 enth\"alt $Z_2 Z_2$ daher Elemente der Ordnung 6. Ferner ist
 $$ (123)(456)(789)(147)(258)(369) = (159267348). $$
 Hieraus folgt, dass $Z_i Z_i$ f\"ur $i \geq 3$ stets Elemente enth\"alt, deren
 Ordnung gleich $9$ ist. Somit ist $\alpha (Z_1) = Z_1$.
 \par
       Ist $\{a,b,c\}$ eine 3-Teilmenge von $X$, so sind $(abc)$ und
 $(acb)$ die einzigen Elemente aus $Z_1$ mit dieser Menge als
 Tr\"agermenge. Ist $\alpha(abc) = (a'b'c')$, so ist $\alpha (acb)
 = (a'c'b')$. Daher wird durch
 $$ \beta \bigl(\{a,b,c\}\bigr) : = \{a',b',c'\} $$
 eine Abbildung auf der Menge der 3-Teilmengen von $X$ definiert.
 Es seien $\{a,b,c\}$ und $\{a,b,d\}$ zwei benachbarte 3-Teilmengen
 von $X$. --- Man erinnere sich: Die Potenzmenge von $X$ ist mit der
 Inklusion als Teilordnung ein projektiver Verband und der Rang
 eines Teilraumes ist nichts anderes als seine Kardinalit\"at. ---
 Die von $(abc)$ und $(abd)$ erzeugte Untergruppe von $A_X$ ist
 isomorph zur $A_{\{a,b,c,d\}}$. Es seien $\{a,b,c\}$ zwei
 Teilmengen des Abstandes 2. Dann ist die von $(abc)$ und $(ade)$
 erzeugte Untergruppe von $A_X$ isomorph zur $A_{\{a,b,c,d,e\}}$. ---
 Es gen\"ugt zu bemerken, dass $(abc)(ade) = (abcde)$ ist. --- Weil
 disjunkte Dreierzyklen eine abelsche Gruppe der Ordnung 9
 erzeugen, erh\"alt $\beta$ daher die Nachbarschaft. Es gibt also
 nach dem Satz von Chow einen Automorphismus oder
 Antiautomorphismus $\gamma$ von $P(M)$, der $\beta$ induziert.
 Weil $X$ aber nicht genau 6 Elemente enth\"alt, kann $\gamma$
 ebenfalls nach dem Satz von Chow kein Antiautomorphismus sein.
 Somit ist $\gamma$, wenn man nur seine Wirkung auf die Punkte von
 $X$ betrachtet, ein Element von $S_X$. Definiere $\lambda$ durch
 $$ \lambda (\xi) := \gamma \xi \gamma^{-1}. $$
 Dann ist $\alpha^{-1} \lambda$ ein Automorphismus von $A_X$, der alle
 Untergruppen von $A_X$, die von Elementen aus $Z_1$ erzeugt
 werden, invariant l\"asst. Setze $\mu := \alpha^{-1} \lambda$. Dann
 ist also $\mu (\zeta) = \zeta^{\epsilon (\zeta)}$ mit $\epsilon
 (\zeta) = \pm 1$ f\"ur alle $\zeta \in Z_1$.
 \par
       Es ist $(123)(124) = (14)(23)$ und $(132)(124) = (134)$. W\"are
 nun $\mu (123) =(132)$, so folgte $\mu (124) = (142)$. Weil man
 nach Satz 8.3 je zwei 3-Teilmengen von $X$ durch eine Kette
 benachbarter 3-Teilmengen verbinden kann, folgte weiter $\mu(abc)
 = (acb)$ f\"ur alle 3-Teilmengen $\{a, b, c\}$. Dann folgte aber der Widerspruch
 $$\eqalign{
       (143) &= \mu (134) = \mu\bigl((132)(124)\bigr)          \cr
	     &= (132) \mu (124) = (123)(142) = (234). \cr}$$
 Also ist $\mu$ die Identit\"at und daher $\alpha = \lambda$, was zu beweisen war.
 \medskip
       Ohne Beweis sei hier mitgeteilt, dass die $S_X$ im Falle $|X| = 6$
 einen Automorphismus besitzt, der kein innerer Automorphismus ist.
 Weil die $A_X$ eine charakteristische Untergruppe von $S_X$ ist,
 induziert dieser Automorphismus einen Automorphismus in $A_X$.
 Dieser Automorphismus bildet $Z_1$ auf $Z_2$ ab, kann also nicht
 durch einen inneren Automorphismus von $S_X$ induziert werden.
 \medskip\noindent
 {\bf 8.6. Satz.} {\it Es seien $L$ und $L'$ projektive Verb\"ande und
 $r$ und $s$ seien nat\"urliche Zahlen. Es sei $V$ eine Menge von
 Unterr\"aumen von $L$, von denen jeder einen Unterraum des Ranges
 $r+1$ enthalte und die \"uberdies die weitere Eigenschaft habe,
 dass jeder Teilraum von $L$, dessen Rang $r + 1$ ist, sich als
 Schnitt von Elementen aus $V$ darstellen l\"asst.
 Entsprechend sei $V'$ eine Menge von Unterr\"aumen von $L'$, von
 denen jeder einen Teilraum des Ranges $s+1$ enthalte und die
 wei\-ter\-hin die Eigenschaft habe, dass jeder Teilraum von
 $L'$, dessen Rang $s + 1$ ist, sich als Schnitt von Elementen
 aus $V'$ darstellen l\"asst. Ist dann $\rho$ eine Bijektion von
 $\UR_r (L) \cup V$ auf $\UR_s (L') \cup V'$, so dass f\"ur $X \in
 \UR_r (L)$ und $Y \in V$ genau dann $X \leq Y$ gilt, wenn $X^\rho
 \leq Y^\rho$ ist, so gibt es genau einen Isomorphismus $\tau$ von
 $L$ auf $L'$, der $\rho$ induziert.}
 \smallskip
       Beweis. Setze $U:=\UR_r (L)$ und $U' := \UR_s (L')$. Ist $Z \in L$
 und gibt es ein $X \in U$ mit $X \leq Z$, so definieren wir $Z^\sigma$ durch
 $$ Z^\sigma := \sum_{X \in U,\ X \leq Z} X^\rho. $$
 Es ist klar, dass $X^\sigma = X^\rho$ f\"ur alle $X \in U$ gilt.
 Sind $Z$, $Z' \in L$ und sind $Z^\sigma$ und $Z'^\sigma$ definiert,
 so folgt aus $Z \leq Z'$, dass $Z^\sigma \leq Z^{'\sigma}$ ist.
 Ist $Y \in V$, so ist $X^\rho \leq Y^\rho$ f\"ur alle $X \in U$,
 f\"ur die $X \leq Y$ gilt. Daher ist $Y^\sigma \leq Y^\rho$.
 \par
       Es sei $Z \in L$ und der Rang von $Z$ sei $r+1$. Es sei ferner
 $X'$ ein Unterraum des Ranges $s$ von $Z^\sigma$. Auf Grund
 unserer Annahme \"uber $V$ gibt es ein $Y \in V$ mit $Z \leq Y$.
 Es gibt weiterhin ein $X \in U$ mit $X^\rho = X'$. Es folgt
 $$ X^\rho = X^\sigma \leq Z^\sigma \leq Y^\sigma \leq Y^\rho $$
 und damit $X \leq Y$. Nach unserer Annahme ist $Z$
 Schnitt von Elementen aus $V$. Somit gilt $X \leq Z$. Dies
 besagt, dass $\sigma$ die Menge der Teilr\"aume des Ranges $r$ von
 $Z$ bijektiv auf die Menge der Teilr\"aume des Ranges $s$ von
 $Z^\sigma$ abbildet.
 \par
       Es sei $Z'$ ein Teilraum des Ranges $s+1$ von $Z^\sigma$. Setzt
 man $M' := \{Y' \mid Y' \in V', Z' \leq Y'\}$, so gilt nach Voraussetzung
 $$ Z' = \bigcap_{Y' \in M'} Y'. $$
 Es gibt nun zwei verschiedene R\"aume $X$ und $Y$ des Ranges $r$ von
 $L$ mit $Z' = X^\sigma + Y^\sigma$. Es folgt $Z = X+Y$. Ist nun $A
 \in V$, so gilt genau dann $X \leq A$, wenn $X^\rho \leq A^\rho$
 gilt. Entsprechendes gilt f\"ur $Y$. Setzt man $M := \{A \mid A \in V,
 Z \leq A\}$, so ist daher $M^\rho = M'$. Ist $B$ ein Teilraum des
 Ranges $r$ von $Z$, so gilt $B \leq A$ f\"ur alle $A \in M$.
 Hieraus folgt
 $$ B^\rho \leq \bigcap_{Y' \in M'} Y' = Z' $$
 f\"ur alle Teilr\"aume $B$ des Ranges $r$ von $Z$. Daher ist $Z^\sigma \leq Z'$
 und damit $Z^\sigma = Z'$. Somit ist der Rang von $Z^\sigma$ gleich $s+1$.
 \par
       Ist nun $r=1$, so folgt die Behauptung aus 8.1, da man mit Hilfe
 von $\rho^{-1}$ erschlie\ss t, dass in diesem Falle auch $s=1$
 ist. Ist $r > 1$, so folgt die Behauptung mittels des Satzes 8.4,
 da $\rho$ offensichtlich nicht von einem Antiisomorphismus
 induziert wird.
 \medskip
       Ein typisches Beispiel f\"ur die Situation dieses Satzes ist die,
 dass $r = 1$ und $V$ die Menge der Hyperebenen von $L$ ist. Der
 Leser formuliere den zu 8.6 dualen Satz.
 \par
       Es sei $L$ ein projektiver Verband und $A$ sei ein Element von
 $L$, welches vom gr\"o\ss ten Element von $L$ verschieden ist. Mit
 $L_A$ be\-zeich\-nen wir die Menge der $X \in L$ mit $X \not\leq A$.
 Dann hei\ss t $L_A$ der \emph{verm\"oge} $A$ \emph{geschlitzte
 Raum}.\index{geschlitzter Raum}{} Ist $A = 0$, so ist $L_A$ nat\"urlich im
 Wesentlichen dasselbe wie $L$. Ist $A$ eine Hyperebene, so hei\ss t $L_A$ auch
 \emph{affiner Raum}.\index{affiner Raum}{}
 \par
       Wir fragen uns nun, wann eine Abbildung $\sigma$ der Punkte von
 $L_A$ auf die Punkte von $L'_{A'}$ durch einen Isomorphismus
 induziert wird. Sicherlich muss $\sigma$ kollineare Punkte in
 kollineare Punkte \"uberf\"uhren. Diese notwendige Bedingung ist,
 falls auf jeder Geraden von $L$ mindestens vier Punkte liegen,
 auch hinreichend. Liegen auf jeder Gerade von $L$ genau drei
 Punkte und sind $L_A$ und $L'_{A'}$ etwa affine R\"aume, so
 erh\"alt jede Bijektion von der Punktmenge von $L_A$ auf die
 Punktmenge von $L'_{A'}$ die Kollinearit\"at. In diesem Falle
 muss man noch verlangen, dass auch komplanare Punkte wieder in
 komplanare Punkte \"ubergehen.
 \par
       Zwei Elemente $X$ und $Y$ eines projektiven Verbandes hei\ss en
 \emph{windschief},\index{windschief}{} falls $X \cap Y = 0$ ist.
 \medskip\noindent
 {\bf 8.7. Satz.} {\it Es sei $L$ ein projektiver Verband und $G$ und
 $H$ seien zwei Geraden von $L$. Es sei ferner $P$ ein Punkt von
 $G+H$, der weder auf $G$ noch auf $H$ liege. Genau dann sind $G$
 und $H$ windschief, wenn es genau eine Gerade durch $P$ gibt, die
 sowohl mit $G$ als auch mit $H$ einen Punkt gemeinsam hat.}
 \smallskip
       Beweis. Wegen $P \leq G + H$ ist
 $$\eqalign{
         4 - \Rg_L(G \cap H) &= \Rg_L(G + H)          \cr
                       &= \Rg_L(G + P + H + P)  \cr
		       &= \Rg_L(G + P) + \Rg_L(H + P) -
				  \Rg_L\bigl((G + P) \cap (H + P)\bigr) \cr
		       &= 6 - \Rg_L\bigl((G+P) \cap (H + P)\bigr). \cr} $$
 Somit sind $G$ und $H$ genau dann windschief, wenn $(G + P) \cap (H + P)$ eine
 Gerade ist. Da alle Geraden durch $P$, die $G$ und
 $H$ treffen, in $(G + P) \cap (H + P)$ liegen, folgt die
 Behauptung des Satzes.
 \medskip\noindent
 {\bf 8.8. Satz.} {\it Es seien $L$ und $L'$ irreduzible projektive Verb\"ande,
 die verm\"oge $A$ bzw. $A'$ geschlitzt seien. Ist
 $\sigma$ eine bijektive Abbildung der Menge der Punkte von $L_A$
 auf die Menge der Punkte von $L'_{A'}$, die kollineare Punkte auf
 kollineare Punkte sowie nicht kollineare auf nicht kollineare
 abbildet, und die im Falle, dass jede Gerade von $L$ genau drei
 Punkte tr\"agt, noch die weitere Eigenschaft hat, komplanare
 Punkte auf komplanare Punkte und nicht komplanare Punkte auf
 ebensolche abzubilden, so gibt es genau einen Isomorphismus $\tau$
 von $L$ auf $L'$ mit $P^\sigma = P^\tau$ f\"ur alle Punkte $P$ von $L_A$.}
 \smallskip
       Beweis. Ist $\Rg(L) \leq 2$, so ist die Aussage des Satzes banal.
 Es sei also $\Rg(L) > 2$.
 \par
       Da jede Gerade von $L$, die nicht in $A$ liegt, auf Grund der
 Ir\-re\-du\-zi\-bi\-li\-t\"at von $L$ mindestens zwei Punkte tr\"agt, die
 ebenfalls nicht in $A$ liegen, k\"onnen wir auf Grund der Annahme,
 dass $\sigma$ nicht kollineare Punkte auf nicht kollineare Punkte
 und kollineare Punkte auf kol\-li\-ne\-a\-re Punkte abbildet, die
 Abbildung $\sigma$ fortsetzen auf die Menge der Geraden von $L$,
 die nicht in $A$ liegen, indem wir $G^\sigma := P^\sigma +
 Q^\sigma$ setzen, falls $G = P + Q$ ist.
 \par
       Es sei nun $P$ ein Punkt auf $A$ und $\Gamma_P$ sei die Menge der
 Ge\-ra\-den von $L$, die durch $P$ gehen, aber nicht in $A$ liegen.
 Wir setzen $\Gamma'_P := \{G^\sigma \mid G \in \Gamma_P\}$. Wegen $\Rg(L) > 2$
 enth\"alt $\Gamma_P$ und damit $\Gamma'_P$ mindestens zwei Geraden. Es seien nun
 $G$ und $H$ zwei verschiedene Geraden
 aus $\Gamma_P$. Dann ist $G+H$ eine Ebene von $L$. Enth\"alt nun
 jede Gerade von $L$ genau drei Punkte, so ist auf Grund unserer
 Annahme $G^\sigma + H^\sigma$ eine Ebene von $L'$. Enth\"alt jede
 Gerade von $L$ mehr als drei Punkte, so folgt aus 8.7, dass auch
 $G^\sigma + H^\sigma$ eine Ebene von $L'$ ist, da es auf Grund der
 Irreduzibilit\"at von $L$ einen Punkt in $G+H$ gibt, der weder auf
 $G$ noch auf $H$ liegt. Hieraus folgt, dass $G^\sigma \cap
 H^\sigma$ ein Punkt ist, der nicht zu $L'_{A'}$ geh\"ort. Somit
 ist
 $ G^\sigma \cap H^\sigma = A' \cap G^\sigma = A' \cap H^\sigma$.
 Dies zeigt, dass $G^\sigma \cap A'$ ein Punkt ist, der
 nicht von der Auswahl von $G \in \Gamma_P$ abh\"angt.
 \par
       Ist nun $P$ ein Punkt von $L$, so setzen wir $P^\tau := P^\sigma$,
 falls $P$ nicht auf $A$ liegt, und $P^\tau := G^\sigma \cap A'$
 mit $G \in \Gamma_P$, falls $P$ auf $A$ liegt. (Hier haben wir vom
 Auswahlaxiom Gebrauch gemacht.)
 Es ist nun ein Leichtes zu zeigen,
 dass $\tau$ eine Bijektion der Punkte von $L$ auf die Menge der
 Punkte von $L'$ ist, die kollineare Punkte auf kollineare Punkte
 und nicht kollineare Punkte auf nicht kollineare Punkte abbildet.
 Die Eindeutigkeit von $\tau$ ist ebenfalls banal.

\mysection{9. Der Satz von Desargues}

\noindent
 Es sei $L$ ein projektiver Verband. Wir nennen $L$
 \emph{desarguessch}, falls in $L$ der folgende \emph{Schlie\ss ungssatz}, der
 sog. Satz von Desargues gilt:\index{Satz von Desargues}{}

       Sind $G_1$, $G_2$ und $G_3$ drei verschiedene Geraden, die in einer
 Ebene von $L$ liegen, ist $P$ ein Punkt mit $P \leq G_i$ f\"ur
 $i:= 1$, 2, 3 und sind $P_1$, $P_2$, $P_3$, $Q_1$, $Q_2$, $Q_3$ Punkte mit
 $P_i \leq G_i$ und $Q_i \leq G_i$ sowie $P_i$, $Q_i \neq P$ f\"ur
 $i:= 1$, 2, 3, ist schlie\ss lich f\"ur $i \neq j$ auch $P_i + P_j
 \neq Q_i + Q_j$, so sind die Punkte $R_{ij} := (P_i + P_j) \cap
 (Q_i + Q_j)$ kollinear.
 \vadjust{
 \vglue 7mm
 \centerline{\includegraphics{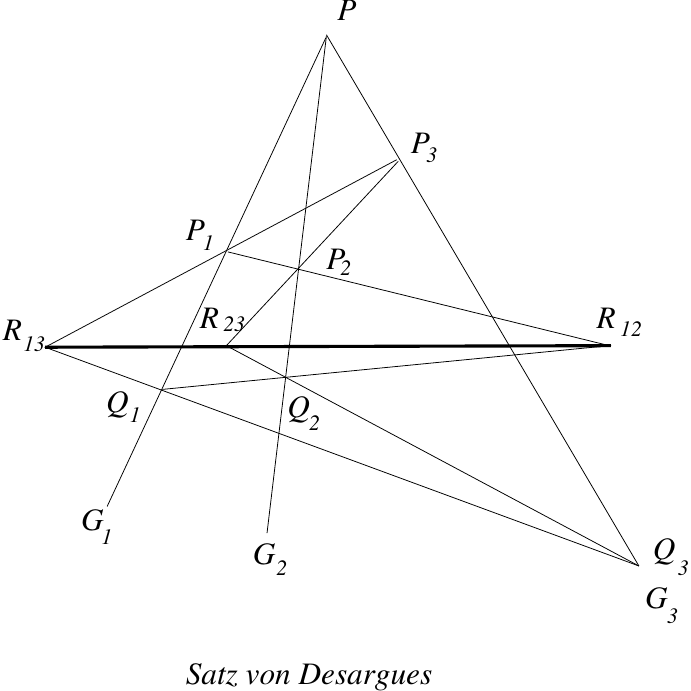}}
 \vglue 3mm}

       Es ist zu bemerken, dass die $R_{ij}$ unter den gemachten
 Voraussetzungen tats\"achlich Punkte sind. Man kann den Satz von
 Desargues grob so formulieren, dass man sagt, dass perspektive
 Dreiecke stets auch axial sind. Dabei ist $P$ als das
 Perspektivit\"atszentrum und die Gerade, auf der die Punkte
 $R_{ij}$ liegen, als die Achse der beiden Dreiecke aufzufassen.
 \par
       Es erhebt sich die Frage, welche projektiven R\"aume desarguessch
 sind. Es ist recht einfach, projektive Ebenen zu konstruieren, die
 es nicht sind, so dass projektive Geometrien vom Rang 3 nicht
 immer desarguessch sind. Beispiele hierf\"ur werden wir im letzten
 Kapitel kennen lernen, die allerdings nicht so einfach zu
 konstruieren sind, wie gerade behauptet. F\"ur die Theorie der
 projektiven Ebenen, die uns hier nur am Rande interessiert, sei
 der Leser auf die B\"ucher Pickert 1955, Hughes \&\ Piper 1973,
 L\"uneburg 1980 verwiesen.
 \par
       F\"ur projektive R\"aume h\"oheren Ranges als 3 gilt jedoch
 \medskip\noindent
 {\bf 9.1. Satz.} {\it Ist $L$ ein irreduzibler projektiver Verband und
 ist $\Rg (L) \geq 4$, so ist $L$ desarguessch.}
 \bigskip\goodbreak
       Beweis. Es seien $G_1$, $G_2$ und $G_3$ drei verschiedene Geraden,
 die in einer Ebene $E$ von $L$ liegen. Ferner sei $P$ ein Punkt,%
 \vadjust{
 \vglue 4mm
 \centerline{\includegraphics{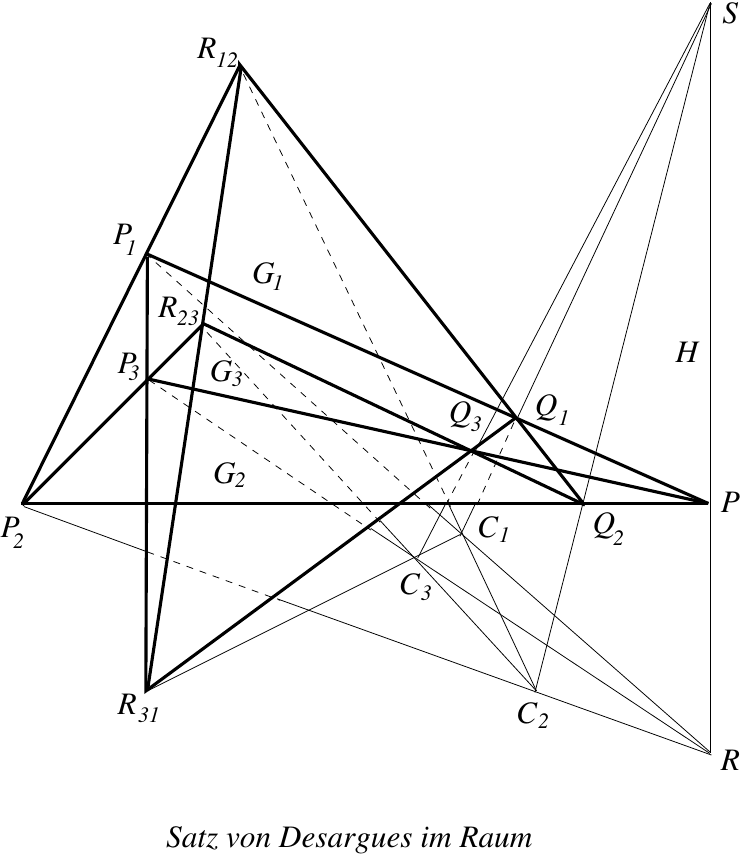}}
 \vglue 4mm}
 der mit die\-sen drei Geraden inzidiere. Schlie\ss lich seien $P_1$,
 $P_2$, $P_3$, $Q_1$, $Q_2$, $Q_3$ Punkte mit $P_i$, $Q_i \leq G_i$ und $P_i$,
 $Q_i \neq P$ sowie $P_i + P_j \neq Q_i + Q_j$, falls nur $i \neq j$
 ist. Da der Rang von $L$ gr\"o\ss er als 3 ist, gibt es eine
 Gerade $H$ mit $P \leq H$ und $H \not\leq E$. Weil $L$ irreduzibel
 ist, gibt es zwei verschiedene Punkte $R$ und $S$ auf $H$, die
 nicht in $E$ liegen. Die Geraden $P_i + R$ und $Q_i + S$ liegen in
 der Ebene $G_i + H$. Wegen $R \neq S \neq P \neq R$ und $H
 \not\leq E$ ist $P_i + R \neq Q_i + S$. Also ist $C_i := (P_i + R)
 \cap (Q_i + S)$ ein Punkt. Weil die Geraden $G_i$ paarweise
 verschieden sind, sind auch die Punkte $C_i$ paarweise verschieden.
 \par
       Gibt es eine Gerade $G$ mit $C_i \leq G$ f\"ur alle $i$, so liegen
 $P_1$, $P_2$ und $P_3$ in der Ebene $R + G$ und $Q_1$, $Q_2$ und $Q_3$
 in der Ebene $S + G$. Hieraus folgt, dass sowohl $P_1$, $P_2$, $P_3$ als
 auch $Q_1$, $Q_2$, $Q_3$ kollinear sind. In diesem Falle ist $R_{12} =
 R_{23} = R_{31}$ und daher nichts zu beweisen. Wir d\"urfen also
 annehmen, dass $C_1 + C_2 + C_3$ eine Ebene ist. Ist $i \neq j$,
 so ist $P_i + P_j + R \neq Q_i + Q_j + S$. Da $C_i$ und $C_j$ in
 beiden Ebenen liegen, ist somit
 $$ C_i + C_j = (P_i + P_j + R) \cap (Q_i + Q_j + S). $$
 Nun ist $R_{ij}$ ein Punkt, der sowohl in
 $P_i + P_j + R$ als auch in $Q_i + Q_j + S$ liegt. Also ist
 $R_{ij} \leq C_i + C_j$ und folglich
 $ R_{ij} \leq E \cap (C_1 + C_2 + C_3)$.
 Weil die $C_i$ nicht in $E$ liegen, ist
 $ \Rg_L \bigl(E \cap (C_1 + C_2 + C_3)\bigr) \leq 2$.
 Somit sind die $R_{ij}$ kollinear, was zu beweisen war.
 \medskip
       Der Satz von Desargues wird im n\"achsten Kapitel eine entscheidende Rolle
 bei der Darstellung projektiver R\"aume mittels Vek\-tor\-r\"au\-men spielen.

\mysectionten{10. Der Satz von Pappos}

\noindent
 Eine wichtige Rolle spielt auch der Satz von Pappos beim Aufbau
 der projektiven Geometrie, wenn auch keine so zentrale wie der
 Satz von Desargues. Er gilt nicht in allen projektiven R\"aumen.
 Er ist vielmehr, wie Hilbert als erster bemerkte, gleichbedeutend
 mit der Kommutativit\"at des der Geometrie zugrunde liegenden
 K\"orpers. Hier formulieren wir zun\"achst den Satz von Pappos und
 bringen ihn mit gewissen anderen r\"aumlichen Tatsachen in
 Zusammenhang.
 \par
       Es seien $G$ und $H$ zwei verschiedene Geraden des projektiven
 Verbandes $L$, die sich im Punkte $S$ schneiden. Ferner seien
 $P_1$, $P_2$ und $P_3$ drei verschiedene Punkte auf $G$, die auch
 von $S$ verschieden seien, und $Q_1$, $Q_2$ und $Q_3$ seien drei
 verschiedene Punkte auf $H$, die ebenfalls von $S$ verschieden
 seien. Ist dann $R_3 := (P_1 + P_2) \cap (Q_1 + Q_2)$, $R_1 := (P_2
 + P_3) \cap (Q_2 + Q_3)$ und $R_2 := (P_3 + P_1) \cap (Q_3 +
 Q_1)$, so sind $R_1$, $R_2$ und $R_3$ kollinear.
 \vadjust{
 \vglue 4mm
 \centerline{\includegraphics{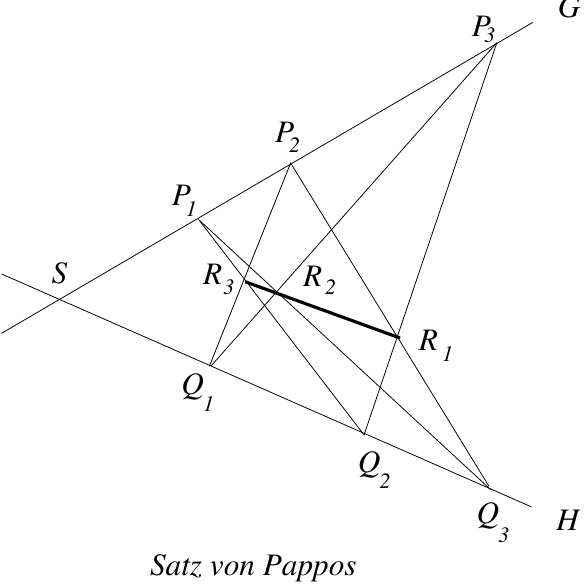}}
 \vglue 3mm}
 \par
 \goodbreak
       Wir folgen nun \"Uberlegungen von 
 Dandelin,\index{Dandelin, G. P.}{} der die G\"ul\-tig\-keit
 des Satzes von Pappos aus der Existenz von gewissen r\"aumlichen
 Konfigurationen erschloss (Germinal Pierre Dandelin, 1794--1847).
 \par
       Es sei $L$ ein irreduzibler projektiver Verband mit $\Rg (L) \geq 4$.
 Sind $G_1$, $G_2$ und $G_3$ paarweise windschiefe Geraden von
 $L$ wie auch $H_1$, $H_2$ und $H_3$ und ist $G_i \cap H_j \neq 0$
 f\"ur $i$, $j := 1$, 2, 3, so nennen wir die Konfiguration dieser sechs
 Geraden {\it Hexagramme mystique\/}.\index{Hexagramme mystique}{}
 \par
       Ist $G_1$, $G_2$, $G_3$, $H_1$, $H_2$, $H_3$ ein Hexagramme mystique, so
 ist $G_i \neq H_j$. Andernfalls h\"atte n\"amlich $G_i$ mit $G_k$
 einen Punkt gemein, was jedoch nicht der Fall ist. Setze $S_i :=
 G_i \cap H_i$ f\"ur $i := 1$, 2, 3. Dies sind dann Punkte, die
 \"uberdies nicht kollinear sind. Sie sind ja sicherlich paarweise
 verschieden, da die $G_i$ windschief sind. Nun ist $S_1 \leq G_1$
 und $S_3 \leq H_3$. W\"are $S_2 \leq S_1 + S_3$, so w\"are $S_2
 \leq G_1 + H_3$. Nun ist $S_2 \leq H_3$, da $H_2$ und $H_3$
 windschief sind. Ferner ist $S_2 \leq G_2$ und $G_2 \cap H_3$ ein
 Punkt von $G_1 + H_3$, der von $S_2$ verschieden ist. Also ist
 $$ G_2 = S_2 + (G_2 \cap H_3) \leq G_1 + H_3. $$
 Daraus folgt, da $G_1 + H_3$ ja eine Ebene ist, dass $G_1 \cap G_2 \neq 0$ ist.
 Dieser Widerspruch zeigt, dass die $S_i$ nicht kollinear sind.
 \par
       Nach diesen Vorbemerkungen sind wir nun in der Lage, den Satz von
 Dandelin zu beweisen.\index{Satz von Dandelin}{}
 \medskip\noindent
 {\bf 10.1. Satz von Dandelin.} {\it Es seien $G$ und $H$ zwei
 verschiedene Ge\-ra\-den eines irreduziblen projektiven Verbandes $L$,
 dessen Rang mindestens $4$ sei. Ferner sei $S := G \cap H$ ein
 Punkt. Weiter seien $P_1$, $P_2$ und $P_3$ drei verschiedene Punkte
 auf $G$ und $Q_1$, $Q_2$ und $Q_3$ drei verschiedene Punkte auf $H$
 und alle diese Punkte seien von $S$ verschieden. Setze
 $$\eqalign{
      R_3 &:= (P_1 + Q_2) \cap (P_2 + Q_1)  \cr
      R_1 &:= (P_2 + Q_3) \cap (P_3 + Q_2)  \cr
      R_2 &:= (P_3 + Q_1) \cap (P_1 + Q_3).  \cr} $$
 Die Punkte $R_1$, $R_2$ und $R_3$ sind sicher dann kollinear, wenn
 es ein Hexagramme mystique $G_1$, $G_2$, $G_3$, $H_1$, $H_2$, $H_3$ gibt mit
 $P_i \leq H_i$, $Q_i \leq G_i$ und $G_i$, $H_i \leq G + H$ f\"ur $i := 1$, 2,
 3.}\index{Satz von Dandelin}{}
 \smallskip
       Beweis. Setze $E_{ij} := G_i + H_j$. Dann ist $E_{ij}$, da $G_i$ und
 $H_j$ verschiedene Geraden sind, die sich in einem Punkte
 schneiden, eine Ebene. Setze ferner $S_i := G_i \cap H_i$. Nun ist
 $$ R_3 \leq P_1 + Q_2 \leq H_1 + G_2 = E_{21} $$
 und
 $$ R_3 \leq P_2 + Q_1 \leq H_2 + G_1 = E_{12}. $$
 \vadjust{
 \vglue 4mm
 \centerline{\includegraphics{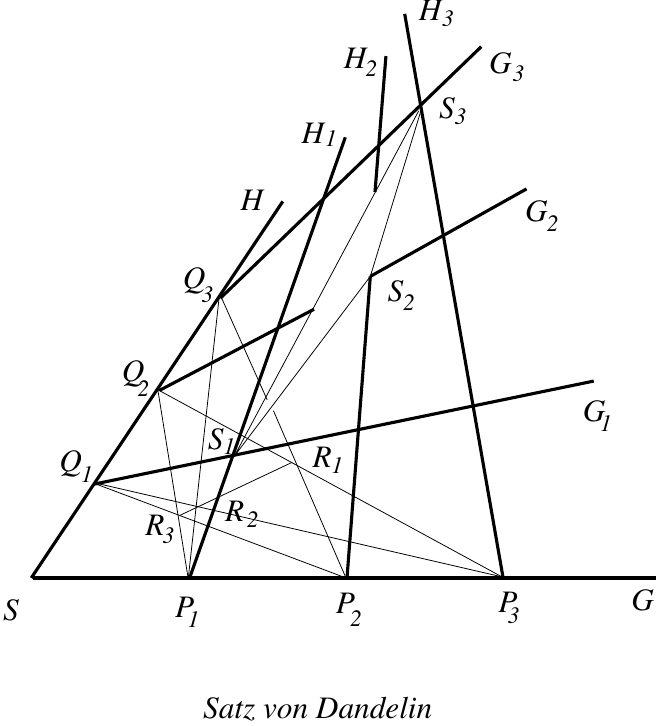}}
 \vglue 1mm}

 \noindent
 Ferner folgt $E_{12} \neq E_{21}$, da $G_1$ und $G_2$ ja windschief sind.
 Au\ss erdem ist $S_1 \leq H_1 \leq E_{21}$ und $S_2 \leq G_2 \leq E_{21}$.
 Ebenso folgt, dass $S_1, S_2 \leq E_{12}$ ist. Weil die $S_i$ nach
 unserer Vorbemerkung nicht kollinear sind, ist also
 $$ S_1 + S_2 = E_{12} \cap E_{21}. $$
 Hieraus folgt
 $$ R_3 \leq E_{12} \cap E_{21} = S_1 + S_2 \leq S_1 + S_2 + S_3. $$
 Entsprechend folgt, dass auch $R_1$ und $R_2$ in $S_1 + S_2 + S_3$
 liegen. W\"aren die $R_i$ nicht kollinear, so w\"are
 $$ G + H = R_1 + R_2 + R_3 = S_1 + S_2 + S_3. $$
 Da die $S_i$ nicht kollinear sind, liegt wenigstens
 einer von ihnen nicht auf $H$. Es sei etwa $S_i \leq H$. Dann ist
 $S_i \neq Q_i$ und daher $G_i = S_i + Q_i \leq G + H$ im
 Widerspruch zu unserer Voraussetzung. Dieser Widerspruch zeigt,
 dass die $R_i$ doch kollinear sind.
 \medskip
       Der Satz von Dandelin gestattet auch eine Umkehrung, die wir nun
 formulieren werden. Dieser Satz ist mit den noch zu entwickelnden analytischen
 Methoden leicht zu beweisen. Ein Beweis an dieser Stelle ist m\"oglich, aber
 knifflig. Der Leser ist herausgefordert, seine Kr\"afte zu erproben.
 \par
       Es seien $G_1$, $G_2$ und $G_3$ drei paarweise windschiefe Geraden
 eines ir\-re\-du\-zib\-len projektiven Verbandes. Ist $H$ eine Gerade und
 gilt $H \cap G_i \neq 0$ f\"ur alle $i$, so hei\ss t $H$ eine
 {\it Transversale\/}\index{Transversale}{} von $G_1$, $G_2$, $G_3$. Eine solche
 gibt es nur dann, wenn der von den $G_i$ aufgespannte Teilraum den Rang 4 hat.
 In diesem Falle gibt es aber sehr viele Transversalen, n\"amlich
 durch jeden Punkt von $G_i$ genau eine (Satz 8.7).
 \medskip\noindent
 {\bf 10.2. Satz.} {\it Es sei $L$ ein irreduzibler projektiver Verband
 mit $\Rg(L) \geq 4$. Genau dann gilt in $L$ der Satz von Pappos,
 wenn f\"ur jedes Hexagramme mystique $G_1$, $G_2$, $G_3$, $H_1$, $H_2$, $H_3$
 gilt: ist $H$ eine Transversale der $G_i$ und $G$ eine
 Transversale der $H_j$, so ist $G \cap H \neq 0$.}
 \smallskip
       Beweis. Nicht triviale \"Ubungsaufgabe.
 \medskip
       Wir beenden das lange erste Kapitel mit einer \"Ubungsaufgabe, die
 man mit Hilfe des Satzes von Pappos in jeder papposschen Ebene
 l\"osen kann, deren Geraden je mindestens vier Punkte tragen (Jackson 1821).\index{Jackson, J.}{}
 $$ {\vbox{\hsize=4 cm\noindent
	   Your aid I want\par\noindent
           nine trees to plant\par\noindent
           in rows just half a score.\par\noindent
	   And let there be\par\noindent
	   in each row three.\par\noindent
	   Solve this, I ask no more.}} $$


 \newpage
       
 \mychapter{II}{Die Strukturs\"atze}

 \noindent
 In diesem Kapitel geht es nun darum, die Strukturs\"atze der projektiven
 Geometrie zu etablieren. Es sind drei S\"atze. Der erste besagt, dass ein
 des\-ar\-gues\-scher, irreduzibler projektiver Verband nichts Anderes ist als der
 Unterraumverband eines geeigneten Vektorraumes. Der zweite besagt, dass
 Isomorphismen zwischen solchen Verb\"anden durch semilineare Abbildungen
 induziert werden, w\"ahrend der letzte besagt, dass eine Korrelation eines
 solchen Verbandes stets durch eine Semibilinearform dargestellt wird. Mit diesen
 S\"atzen hat man dann das Instrumentarium der linearen
 Algebra\index{lineare Algebra}{} an der Hand,
 wenn man tiefer in die Theorie der projektiven Geometrien eindringen will.

\mysection{1. Zentralkollineationen}

\noindent
 Im Folgenden betrachten wir nur irreduzible
 projektive Verb\"ande.\index{Zentralkollineation}{}
 Es sei $L$ ein solcher mit $\Rg (L) \geq 3$. Wir bezeichnen
 wieder, wie bisher, mit $\Pi$ das gr\"o\ss te und mit $0$ das
 kleinste Element von $L$. Ist $\sigma$ eine Kollineation von $L$
 und gibt es eine Hyperebene $H$ von $L$ mit der Eigenschaft
 \smallskip
       (A)\quad Ist $X \leq H$, so ist $X^\sigma = X$,
 \smallskip\noindent
 so nennen wir $\sigma$ {\it axiale Kollineation\/}\index{axiale Kollineation}{}
 und $H$ {\it Achse\/}\index{Achse}{} von
 $\sigma$. Da nach I.2.11 jedes Element von $L$ die obere Grenze der in ihm
 liegenden Punkte ist, ist (A) gleich\-be\-deu\-tend mit
 \smallskip
       (A*)\enspace Ist $P$ ein Punkt auf $H$, so ist $P^\sigma = P$,
 \smallskip
       Ist $\sigma$ eine Kollineation von $L$ und gibt es einen Punkt $P$
 in $L$ mit der Eigenschaft
 \smallskip
       (Z)\quad Ist $P \leq X$, so ist $X^\sigma = X$,
 \smallskip\noindent
 so nennen wir $\sigma$ \emph{zentral}\index{Zentralkollineation}{} und $P$
 hei\ss t {\it Zentrum\/}\index{Zentrum}{} von $\sigma$. Da $P \leq P$ ist, ist
 $P^\sigma = P$. Daher induziert $\sigma$ eine Kollineation in $\Pi / P$, die
 gleich der Identit\"at ist. Weil $\Pi / P$ ein projektiver Verband
 ist, ist (Z) gleichbedeutend mit
 \smallskip
       (Z*)\enspace Ist $G$ eine Gerade durch $P$, so ist $G^\sigma = G$.
 \smallskip
       Es gilt nun der folgende Satz.
 \medskip\noindent
 {\bf 1.1. Satz.} {\it Ist $L$ ein irreduzibler projektiver
 Verband, dessen Rang mindestens $3$ ist, und ist $\sigma$ eine
 Kollineation von $L$, so gilt:
 \item{a)} Ist $\sigma$ axial und hat $\sigma$ zwei verschiedene Achsen,
 oder
 \item{b)} ist $\sigma$ zentral und hat $\sigma$ zwei verschiedene
 Zentren,
 \par\noindent
 so ist $\sigma = 1_L$.}
 \smallskip
       Beweis. a) Es seien $H$ und $K$ zwei verschiedene Achsen von
 $\sigma$. Ferner sei $P$ ein Punkt von $L$. Liegt $P$ auf $H$ oder
 $K$, so ist $P^\sigma = P$. Es liege $P$ weder auf $H$ noch auf
 $K$. Da $\Rg(L) \geq 3$ ist, ist $\Rg_L(H) \geq 2$. Es gibt folglich eine Gerade
 $G$ mit
 $G \leq H$ und $G \not\leq H \cap K$. Wegen $\Ko_H (H \cap K) = 1$
 ist $G \cap H \cap K$ ein Punkt. Weil $L$ irreduzibel ist, gibt es
 folglich zwei Punkte $R$ und $S$ auf $G$, die nicht in $H \cap K$
 liegen. Nun ist $R + S = G \cap H$. Hiermit folgt, dass $R + P \neq S + P$, und
 weiter, dass
 $$ (R + P) \cap (S + P) = P $$
 ist. Ferner sind $U := K \cap (R + P)$ und $V := K \cap (S + P)$
 Punkte, die von $R$ und $S$ verschieden sind, da weder $R$ noch
 $S$ in $K$ liegt. Daher ist $R + P = R + U$ und $S + P = S + V$.
 Es folgt, dass
 $$ P = (R + U) \cap (S + V) $$
 ist. Damit erhalten wir schlie\ss lich
 $$\eqalign{
 P^\sigma &= \bigl((R + U) \cap (S + V)\bigr)^\sigma
 	   = (R^\sigma + U^\sigma) \cap (s^\sigma + V^\sigma) \cr
	  &= (R + U) \cap (S + V) = P. \cr} $$
 Damit ist gezeigt, dass $\sigma$ alle Punkte von $L$ festl\"asst, woraus folgt,
 dass $\sigma = 1_L$ ist.
 \par
       b) Da der zu einem projektiven Verband duale
 Verband\index{dualer Verband}{} nur dann wieder ein projektiver Verband ist, wenn der Rang des
 gegebenen Verbandes endlich ist --- wir reden hier von irreduziblen
 projektiven Verb\"anden ---, k\"onnen wir uns zum Beweise von b)
 nicht auf die Aussage a) berufen, wir m\"ussen sie vielmehr
 getrennt beweisen.
 \par
      Es seien $P$ und $Q$ zwei verschiedene Zentren von $\sigma$. Ferner
 sei $X$ ein Punkt, der nicht auf $P + Q$ liegt. Dann sind $P + X$
 und $Q + X$ zwei verschiedene Geraden, so dass $X = (P + X) \cap
 (Q + X)$ gilt. Weil die Geraden durch $P$ wie auch die Geraden
 durch $Q$ festbleiben, folgt $X^\sigma = X$. Somit l\"asst
 $\sigma$ alle Punkte von $L_{P + Q}$ fest. Mittels I.8.8 folgt
 schlie\ss lich, dass $\sigma = 1_L$ ist.
 \par
       Damit ist alles bewiesen.
 \medskip\noindent
 {\bf 1.2. Satz.} {\it Es sei $L$ ein irreduzibler projektiver
 Verband mit $\Rg(L) \geq 3$ und $\sigma$ sei eine Kollineation
 von $L$. Genau dann ist $\sigma$ axial, wenn $\sigma$ zentral
 ist.}
 \smallskip
       Beweis. Es sei $\sigma$ axial und $H$ sei eine Achse von $\sigma$.
 Ist $P$ ein Punkt von $L$ mit $P^\sigma = P \not\leq H$ und ist
 $G$ eine Gerade durch $P$, so ist $G = P + (G \cap H)$. Hieraus
 folgt
 $$ G^\sigma = \bigl(P + (G \cap H)\bigr)^\sigma = P^\sigma + (G \cap H)^\sigma
      = P + (G \cap H) = G, $$
 so dass $P$ wegen (Z*) ein Zentrum von $\sigma$ ist.
 \par
       Es sei nun $P \neq P^\sigma$ f\"ur alle Punkte $P$, die nicht auf
 $H$ liegen. Es sei $P$ ein solcher Punkt. Dann ist $P + P^\sigma$
 eine Gerade. Ferner gilt $\Pi = P + H = P^\sigma + H$. Hieraus
 folgt mit Hilfe des Modulargesetzes
 $$\eqalign{
     P + P^\sigma &= (P + P^\sigma) \cap (P^\sigma + H)
		   = P^\sigma + \bigl((P + P^\sigma) \cap H\bigr) \cr
		  &= P^\sigma + \bigl((P + P^\sigma) \cap H\bigr)^\sigma
		   = \bigl(P + ((P + P^\sigma) \cap H)\bigr)^\sigma \cr
		  &= \bigl((P + P^\sigma) \cap (P + H)\bigr)^\sigma
		   = (P + P^\sigma)^\sigma, \cr}$$
 so dass $P + P^\sigma$ eine Fixgerade von $\sigma$ ist.
 \par
       Es sei $F$ ein Punkt mit $F \not\leq H$. Dann ist also $F +
 F^\sigma$ eine Fixgerade von $\sigma$. Setze $C := (F + F^\sigma)
 \cap H$. Dann ist $C$ ein Punkt. Wir zeigen, dass $C$ Zentrum von
 $\sigma$ ist. Es sei $G$ eine Gerade durch $C$. Ist $G \leq H$
 oder $G = F + F^\sigma$, so ist $G$ eine Fixgerade von $\sigma$.
 Es sei also $G \not\leq H$ und $G \neq F + F^\sigma$. Es gibt dann
 einen Punkt $P$ auf $G$, der nicht auf $F + F^\sigma$ liegt. Wegen
 $$ (F + F^\sigma) \cap G = C = G \cap H $$
 ist $P \not\leq H$ und daher $P \neq P^\sigma$. Setze $R := (P + F) \cap H$.
 Dann ist $R$ ein Punkt, da ja $P \neq F$ ist. Es folgt
 $ P + F = P + R = F + R$.
 Wegen $R^\sigma = R$ ist daher
 $$\eqalign{
   P + P^\sigma + F + F^\sigma &= P^\sigma + P + F + F^\sigma
	 		        = P^\sigma + P + R + F^\sigma \cr
		 	       &= P^\sigma + P + R^\sigma + F^\sigma
			        = P^\sigma + P + (R + F)^\sigma \cr
			       &= P^\sigma + P + (P + F)^\sigma
			        = P^\sigma + P + F^\sigma. \cr}$$
 Wegen $P \not\leq F + F^\sigma$ ist also
 $$ 3 \leq \Rg_L (P + P^\sigma + F + F^\sigma)
	 = \Rg_L (P^\sigma + P + F^\sigma) \leq 3. $$
 Hieraus folgt, dass $(P + P^\sigma) \cap (F + F^\sigma)$ ein
 Punkt ist. Weil $P + P^\sigma$ und $F + F^\sigma$ Fixgeraden von
 $\sigma$ sind, ist dieser Punkt ein Fixpunkt von $\sigma$, liegt
 also auf $H$. Daher ist
 $ (P + P^\sigma) \cap (F + F^\sigma) = C$.
 Hieraus folgt weiter, dass $P + P^\sigma = P + C = G$ ist. Somit ist
 $G^\sigma = G$, so dass $C$ in der Tat ein Zentrum von $\sigma$ ist.
 \par
       Es sei nun $\sigma$ zentral und $P$ sei ein Zentrum von $\sigma$.
 Ist $H$ eine Hyperebene, die nicht durch $P$ geht und die unter
 $\sigma$ festbleibt, so ist $H$ eine Achse von $\sigma$, wie
 unschwer zu sehen ist. Es gelte also $H^\sigma \neq H$ f\"ur alle
 nicht durch $P$ gehenden Hyperebenen $H$. Es sei $H$ eine solche Hyperebene.
 Dann ist $P \not\leq H \cap H^\sigma$, so dass $A := P + (H \cap H^\sigma)$ eine
 Hyperebene ist. Weil $P \leq A$ ist, ist $A^\sigma = A$. Nun ist
 $$ A \cap H = \bigl(P + (H \cap H^\sigma)\bigr) \cap H
	     = (P \cap H) + (H \cap H^\sigma) = H \cap H^\sigma. $$
 Ebenso folgt, dass $A \cap H^\sigma = H \cap H^\sigma$
 ist. Also ist $A \cap H = A \cap H^\sigma$ und daher
 $$ (A \cap H)^\sigma = A^\sigma \cap H^\sigma = A \cap H^\sigma
		      = A \cap H. $$
 Insbesondere folgt, dass $(H \cap H^\sigma)^\sigma = H \cap H^\sigma$ ist.
 \par
       Es sei nun $K$ ein Komplement von $P$. Dann ist $(K \cap K^\sigma)^\sigma
 = K \cap K^\sigma$. Weil alle Geraden durch $P$
 unter $\sigma$ festbleiben, l\"asst $\sigma$ alle Punkte von $K
 \cap K^\sigma$ fest. Somit ist $K \cap K^\sigma \leq A$ und daher
 $A \cap K = K \cap K^\sigma$. Ist nun $X$ ein von $P$
 verschiedener Punkt von $A$, so gibt es nach I.1.8 ein Komplement
 $K$ von $P$, welches $X$ enth\"alt. Es folgt $X \leq A \cap K = K
 \cap K^\sigma$, so dass $X$ ein Fixpunkt von $\sigma$ ist. Damit
 ist $A$ als Achse von $\sigma$ erkannt und der Satz in allen
 seinen Teilen bewiesen.
 \medskip
       Die zentralen bzw.\ axialen Kollineationen eines projektiven Verbandes
 werden wir im folgenden auch
 {\it Perspektivit\"aten\/}\index{Perspektivit\"aten}{} nennen. Nach dem gerade
 bewiesenen Satz besteht zwischen ihnen ja kein Unterschied. Ist $\sigma$ eine
 Perspektivit\"at mit dem
 Zentrum $P$ und der Achse $H$ und ist $P \leq H$, so nennen wir
 $\sigma$ {\it Elation\/}.\index{Elation}{} Ist $P \not\leq H$, so nennen wir
 $\sigma$ {\it Streckung\/}\index{Streckung}{} oder
 {\it Homologie\/}.\index{Homologie}{} Ist $H$ eine
 Hyperebene, so be\-zeich\-ne $\Delta(H)$ die Menge aller
 Perspektivit\"aten mit der Achse $H$. Offensichtlich ist $\Delta
 (H)$ eine Untergruppe der Gruppe aller Kollineationen von $L$.
 Entsprechend bezeichnen wir f\"ur einen Punkt $P$ mit $\Delta(P)$
 die Menge aller Perspektivit\"aten mit dem Zentrum $P$. Auch dies
 ist eine Untergruppe der Gruppe aller Kollineationen von $L$.
 Schlie\ss lich setzen wir $\Delta(P,H) := \Delta(P) \cap \Delta(H)$. Diese
 Gruppen werden wir nun eingehender untersuchen.
 \medskip\noindent
 {\bf 1.3. Satz.} {\it Ist $P$ ein Punkt und $H$ eine Hyperebene
 des irreduziblen projektiven Verbandes $L$, ist ferner $1 \neq
 \sigma \in \Delta(P,H)$ und ist $X$ ein Element von $L$ mit
 $X^\sigma = X$, so ist $P \leq X$ oder $X \leq H$.}
 \smallskip
       Beweis. Es sei $P \not\leq X$. Ist dann $Q$ ein Punkt von $X$, so
 ist $Q$ Fixpunkt von $\sigma$, da ja $Q = (P + Q) \cap X$ ist.
 W\"are nun $X \not\leq H$, so g\"abe es insbesondere einen Punkt
 $Q$ auf $X$, der nicht in $H$ l\"age. Dieser Punkt w\"are dann
 aber ein zweites Zentrum von $\sigma$, da ja jede Gerade durch $Q$
 die Hyperebene $H$ in einem Punkte tr\"afe und somit unter
 $\sigma$ fix bliebe. Diese widerspr\"ache aber 1.1b).
 \medskip\noindent
 {\bf 1.4. Satz.} {\it Es seien $\sigma$ und $\tau$ zwei Elationen
 des irreduziblen projektiven Verbandes $L$. Haben $\sigma$ und
 $\tau$ die gleiche Achse aber verschiedene Zentren oder das
 gleiche Zentrum aber verschiedene Achsen, so ist $\sigma \tau =
 \tau \sigma$.}
 \smallskip
       Beweis. Es sei $H$ Achse von $\sigma$ und $\tau$. Ferner sei $P$
 ein Zentrum von $\sigma$ und $Q$ ein solches von $\tau$ und es sei
 $P \neq Q$. Es ist $\tau^{-1} \sigma \tau$ eine Elation mit dem
 Zentrum $P^\tau$ und der Achse $H^\tau$. Wegen $P \leq H$ und
 $\tau \in \Delta(Q,H)$ bleiben $P$ und $H$ unter $\tau$ fest.
 Also gilt $\sigma^{-1}, \tau^{-1} \sigma \tau \in \Delta(P,H)$
 und damit $\sigma^{-1} \tau^{-1} \sigma \tau \in \Delta(P,H)$.
 \par
       F\"angt man mit $\sigma^{-1} \tau^{-1} \sigma$ an, so sieht man
 ganz entsprechend, dass auch $\sigma^{-1} \tau^{-1} \sigma \tau
 \in \Delta(Q,H)$ gilt. Mit Satz 1.1 folgt daher
 $$ \sigma^{-1} \tau^{-1} \sigma \tau \in \Delta(P,H) \cap \Delta(Q,H)
	  = \{1\}. $$
 Also ist $\sigma \tau = \tau \sigma$.
 \par
       Die zweite Aussage von 1.4 beweist sich analog.
 \medskip
       Mit $\E(H)$ bezeichnen wir im folgenden die Menge aller Elationen
 mit der Achse $H$ und mit $\E(P)$ die Menge aller Elationen mit dem
 Zentrum $P$. Ferner setzen wir $\E(P,H) := \E(P) \cap \E(H)$. Dies
 ist nat\"urlich nur dann von Interesse, wenn $P \leq H$ ist.
 \medskip\noindent
 {\bf 1.5. Satz.} {\it Ist $L$ ein irreduzibler projektiver
 Verband, ist $H$ eine Hyperebene und $P$ ein Punkt von $L$, so
 sind $\E(H)$ und $\E(P)$ Untergruppen der Kollineationsgruppe von
 $L$.}
 \smallskip
       Beweis. Es seien $\sigma$, $\tau \in \E(H)$. Haben $\sigma$ und
 $\tau$ das gleiche Zentrum, so hat auch $\sigma \tau$ dieses
 Zentrum, so dass das Produkt in $\E(H)$ liegt. Wir nehmen nun an,
 dass $\sigma$ das Zentrum $P$ und $\tau$ das Zentrum $Q$ habe und
 dass $P \neq Q$ sei. Wir d\"urfen weiter annehmen, dass $\sigma$
 und $\tau$ beide von der Identit\"at verschieden sind. Dann ist
 $\sigma \tau \neq 1$. Es sei $R$ das eindeutig bestimmte Zentrum
 der Perspektivit\"at $\sigma \tau$. Dann ist $R^\sigma$ das
 Zentrum von $\sigma^{-1} \sigma \tau \sigma = \tau \sigma$. Nach
 1.4 ist $\tau \sigma = \sigma \tau$, so dass also $R^\sigma = R$
 ist. Nach 1.3 ist daher $R = P$ oder $R \leq H$. Wegen $P \leq H$
 gilt in jedem Falle $R \leq H$. Somit ist $\sigma \tau \in \E(H)$.
 Schlie\ss lich ist klar, dass mit $\sigma$ auch $\sigma^{-1}$ zu
 $\E(H)$ geh\"ort.
 \par
      Ebenso zeigt man, dass auch $\E(P)$ eine Gruppe ist.
 \medskip
       Wir stellen nun die Frage nach der Existenz von Perspektivit\"aten
 und wir werden sehen, dass uns mit dem Satz von Desargues ein
 Mittel in die Hand gegeben ist, solche zu konstruieren. Ist $L$
 ein irreduzibler projektiver Verband und ist $\Rg (L) \geq 4$, so
 gilt nach I.9.1 in $L$ der Satz von Desargues.\index{Satz von Desargues}{} Ist
 $\Rg (L) = 3$, so m\"ussen wir seine G\"ultigkeit voraussetzen. Wir beweisen
 zun\"achst.
 \medskip\noindent
 {\bf 1.6. Satz.} {\it Ist $L$ ein desarguesscher, irreduzibler
 projektiver Verband, ist $H$ eine Hyperebene von $L$ und sind $P$,
 $A$ und $B$ drei verschiedene, kollineare Punkte von $L$ mit $A$,
 $B \not\leq H$, so gibt es genau ein $\sigma \in \Delta(P,H)$ mit
 $A^\sigma = B$.}
 \smallskip
       Beweis. Aus 1.3 folgt, dass es h\"ochstens ein solches $\sigma$
 gibt. Ist n\"amlich $\tau$ ein weiteres Element aus $\Delta(P,H)$
 mit $A^\tau = B$, so ist $\sigma \tau^{-1} \in \Delta(P,H)$.
 Andererseits ist $A$ ein von $P$ verschiedener Fixpunkt von $\sigma
 \tau^{-1}$, der auch nicht in $H$ liegt. Also ist $\sigma
 \tau^{-1} = 1$, so dass $\sigma = \tau$ ist.
 \par
       Wir zeigen nun die Existenz eines solchen $\sigma$. Dies ist sehr
 einfach, falls auf jeder Geraden von $L$ genau drei Punkte liegen.
 Wegen $P \leq A + B$ und $A$, $B \not\leq H$ ist dann $P \leq H$.
 Wir definieren $\sigma$ durch $X^\sigma = X$, falls $X$ ein Punkt
 von $H$ ist. Ist $X$ ein Punkt des geschlitzten Raumes $L_H$, so
 sei $X^\sigma$ der von $X$ und $P$ verschiedene dritte Punkt auf
 $X + P$. Offenbar ist $\sigma^2 = 1$, so dass $\sigma$ eine
 Bijektion der Punkt\-men\-ge von $L$ auf sich ist. Um zu zeigen, dass
 $\sigma$ eine Kollineation ist, gen\"ugt es wegen $\sigma^2 = 1$
 zu zeigen, dass $\sigma$ kollineare Punkte auf kollineare Punkte
 abbildet. Dies ist aber, da auf jeder Geraden von $L$ nur drei
 Punkte liegen, eine unmittelbare Konsequenz des Veblen-Young
 Axioms. Wir d\"urfen im Folgenden also annehmen, dass auf jeder
 Geraden von $L$ mindestens vier Punkte liegen.
 \par
       Wir betrachten den geschlitzten Raum $L_{A+B}$. Es sei $X$ ein
 Punkt dieses Raumes. Ist $X$ ein Punkt von $H$, so setzen wir
 $X^\sigma := X$ und $X^\tau := X$. --- Wir definieren neben
 $\sigma$ also auch gleichzeitig die zu $\sigma$ inverse Abbildung
 $\tau$. --- Liegt $X$ nicht in $H$, so setzen wir
 $$ X^\sigma := \bigl(((X + A) \cap H) + B\bigr) \cap (P + X) $$ 
 bzw.
 $$ X^\tau := \bigl(((X + B) \cap H) + A\bigr) \cap (P + X). $$
 Offensichtlich ist $\sigma \tau = 1 = \tau \sigma$, so
 dass $\sigma$ und $\tau$ zueinander inverse Bijektionen der
 Punktmenge von $L_{A+B}$ auf sich sind. Es gen\"ugt nun im
 Folgenden zu zeigen, dass $\sigma$ kollineare Punkte auf
 kollineare Punkte abbildet, da $\tau$ als Abbildung des gleichen
 Typs dann ebenfalls diese Eigenschaft hat.
 \par
       Es seien nun $X$, $Y$ und $Z$ drei kollineare Punkte von
 $L_{A+B}$.
 Wir d\"urfen oBdA annehmen, dass $Z \leq H$ und $X$, $Y
 \not\leq H$ gilt. Ist $P \leq X + Y$ oder $A \leq X + Y$, so folgt
 unmittelbar aus der Definition von $\sigma$, dass $X^\sigma$,
 $Y^\sigma$ und $Z^\sigma = Z$ kollinear sind.%
 \vadjust{
 \vglue 4mm
 \centerline{\includegraphics{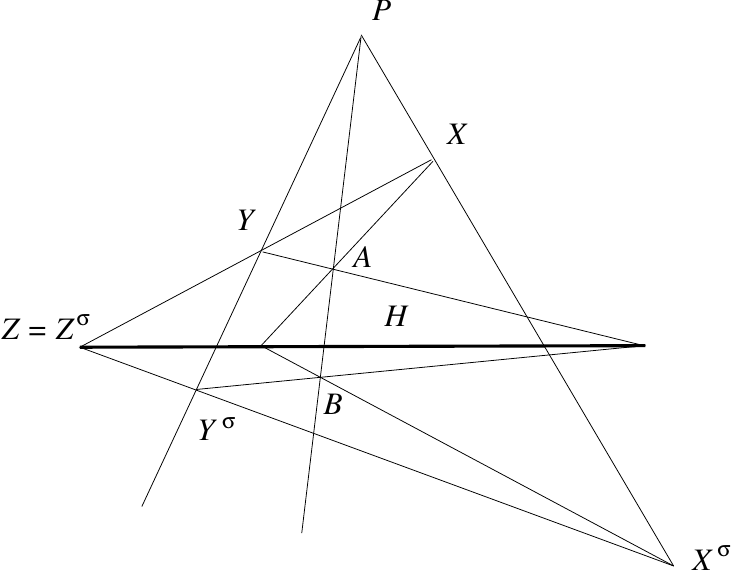}}
 \vglue 3mm}
 Es sei also $P$, $A
 \not\leq X + Y$. Dann ist
 auch $P$, $B \not\leq X^\sigma + Y^\sigma$. Die
 Dreiecke $A$, $X$, $Y$ und $B$, $X^\sigma$,
 $Y^\sigma$ erf\"ullen daher die Voraussetzungen des Satzes von
 Desargues. Somit sind die Punkte
 $$\eqalign{
 (A + Y) &\cap (B + Y^\sigma), \cr
 (A + X) &\cap (B + X^\sigma), \cr
 (X + Y) &\cap (X^\sigma + Y^\sigma) \cr} $$
 kollinear. Da die ersten beiden Punkte in $H$ liegen, liegt auch
 der dritte Punkt in $H$. Folglich ist
 $$ (X + Y) \cap (X^\sigma + Y^\sigma) = (X + Y) \cap H = Z, $$
 so dass die Punkte
 $X^\sigma$, $Y^\sigma$ und $Z = Z^\sigma$ kollinear sind. Damit ist
 gezeigt, dass $\sigma$ kollineare Punkte auf kollineare Punkte
 abbildet. Nach unserer Vorbemerkung bildet $\sigma$ dann auch
 nicht kollineare Punkte auf nicht kollineare Punkte ab. Weil auf
 jeder Geraden von $L$ mindestens vier Punkte liegen, ist daher
 Satz I.8.8 anwendbar, so dass $\sigma$ durch eine Kollineation von
 $L$ induziert wird. Diese geh\"ort offensichtlich zu $\Delta(P,H)$ und bildet
 $A$ auf $B$ ab. Damit ist alles bewiesen.
 \medskip
       Projektive Ebenen\index{projektive Ebene}{} sind nicht immer
 desarguessch. Es ist daher
 angebracht, die Situation in diesem Falle etwas genauer zu ana\-ly\-sie\-ren.
 Die nun folgende Analyse stammt von R. Baer (Baer 1942).\index{Baer, R.}{}
 \par
       Es sei $L$ eine projektive Ebene, dh.\ ein irreduzibler
 projektiver Verband des Ranges 3. Ferner sei $P$ ein Punkt und $G$
 eine Gerade von $L$. Wir sagen, dass in $L$ der $(P,G)$-{\it Satz von
 Desargues\/} erf\"ullt sei, wenn Folgendes gilt:
 \par
      Sind $G_1$, $G_2$, $G_3$ drei verschiedene Geraden von $L$, die
 durch $P$ gehen, und sind $P_i$, $Q_i$ f\"ur $i := 1$, 2, 3
 Punkte auf $G_i$, die von $P$ verschieden sind und auch nicht auf
 $G$ liegen, ist $P_i + P_j \neq Q_i + Q_j$ f\"ur $i \neq j$ und
 liegen die Punkte
 $$ R_{12} := (P_1 + P_2) \cap (Q_1 + Q_2) $$
 und
 $$ R_{13} := (P_1 + P_3) \cap (Q_1 + Q_3) $$
 auf $G$, so liegt auch der Punkt
 $$ R_{23} := (P_2 + P_3) \cap (Q_2 + Q_3) $$
 auf $G$.
 \par
       Ferner sagen wir, dass $L$ eine
 $(P,G)$-{\it transitive\/}\index{transitiv}{} Ebene
 ist, wenn es zu jedem Punktepaar $A,$ $B$ mit $P \neq A \not\leq G$ und
 $P \neq B \not\leq G$ sowie $P + A = P + B$ ein $\sigma \in \Delta(P,G)$ mit
 $A^\sigma = B$ gibt. Der Beweis von 1.6
 liefert auch die G\"ultigkeit einer Schlussrichtung in
 \medskip\noindent
 {\bf 1.7. Satz.} {\it Eine projektive Ebene ist genau dann $(P,G)$-transitiv,
 wenn in ihr der $(P,G)$-Satz von Desargues gilt.}
 \smallskip
       Dies ist wieder eine gute Gelegenheit f\"ur den Leser, seine
 be\-reits erworbenen F\"ahigkeiten zu erproben und zu zeigen, dass
 aus der $(P,G)$-Transitivit\"at der $(P,G)$-Satz von Desargues
 folgt. Dazu w\"ahle er sich ein $\sigma \in \Delta(P,G)$ mit
 $P^\sigma_1 = Q_1$ und verfolge die Wirkung von $\sigma$ auf den
 Rest der Figur.
 \medskip\noindent
 {\bf 1.8. Korollar.} {\it Eine projektive Ebene ist genau dann desarguessch,
 wenn sie $(P,G)$-transitiv ist f\"ur alle Punkt-Geradenpaare $(P,G)$.}
 \medskip
 Weiter gilt
 \medskip\noindent
 {\bf 1.9. Korollar.} {\it Ist $L$ eine desarguessche projektive
 Ebene, so ist auch $L^d$ desarguessch.}
 \smallskip
       Dies folgt aus 1.8 und der Bemerkung, dass $L$ offensichtlich
 genau dann $(P,G)$-transitiv ist, wenn $L^d$ eine $(G,P)$-transitive Ebene ist.
 \medskip
 Ein Element $g$ einer Gruppe $G$ hei\ss t \emph{Involution}, wenn
 $g^2 = 1 \neq g$ ist.
 \medskip\noindent
 {\bf 1.10. Satz.} {\it Es sei $L$ ein irreduzibler projektiver Verband, dessen
 Rang mindestens $3$ sei. Ferner seien $P$ und $Q$ zwei verschiedene Punkte und
 $H$ eine Hyperebene von $L$ mit $P$, $Q \not\leq H$.
 Ist dann $\rho$ eine Involution aus $\Delta(P,H)$ und $\sigma$ eine Involution
 aus $\Delta(Q,H)$, so ist $1 \neq \rho \sigma \in \E((P + Q) \cap H,H)$.}
 \smallskip
       Beweis. Nat\"urlich gilt $\rho \sigma \in \Delta(H)$. Ferner gilt
 $(P + Q)^{\rho \sigma} = P + Q$. Wegen $P + Q \not\leq H$ liegt
 das Zentrum $R$ von $\rho \sigma$ nach 1.3 auf $P + Q$, da ja
 offenbar $\rho \sigma \neq 1$ ist. Es bleibt zu zeigen, dass $R$
 auf $H$ liegt.
 \par
       Wegen $R^{\rho \sigma} = R$ und $\sigma^2 = 1$ ist $R^\rho =
 R^{\rho \sigma^2} = R^\sigma$.  Hieraus folgt
 $$ (R^\rho)^{\rho \sigma} = R^{\rho^2\sigma} = R^\sigma = R^\rho. $$
 Nach 1.3 gilt daher $R^\rho = R$ oder $R^\rho \leq H$. Im letzteren Falle folgt
 $R^\rho = R^{\rho^2} = R$ und damit $R \leq H$. Im ersten Fall
 folgt wieder mit 1.3, dass $R \leq H$ oder $R = P$ ist. W\"are $R = P$ so folgte
 $P = R = R^{\rho \sigma} = P^\sigma$ und mit $P \not\leq H$ weiter der
 Widerspruch $P = Q$. Also ist in jedem Falle $R \leq H$.
 \medskip
 Der n\"achste Satz ist vor allem f\"ur projektive Ebenen interessant, da eine
 ebene projektive Geometrie nicht notwendig desarguessch ist.
 \medskip\noindent
 {\bf 1.11. Satz.} {\it Es sei $L$ ein irreduzibler projektiver
 Verband mit $\Rg(L) \geq 3$ und jede Gerade von $L$ trage
 wenigstens vier Punkte. Ferner sei $H$ eine Hyperebene und $P$ und
 $Q$ seien zwei verschiedene Punkte von $L$, die nicht auf $H$
 liegen. Ist $\Rg(L) = 3$, so werde noch vorausgesetzt, dass $L$
 sowohl $(P,H)$- als auch $(Q,H)$-transitiv sei. Dann ist
 $\Delta(X,H)$ f\"ur alle Punkte $X$ auf $P + Q$ eine
 Untergruppe der von $\Delta(P,H)$ und $\Delta(Q,H)$
 erzeugten Gruppe. Ist $\Rg (L) = 3$, so ist $L$ f\"ur alle diese
 Punkte $(X,H)$-transitiv.}
 \smallskip
       Beweis. Es sei $G$ die von 
 $\Delta(P,H)$ und $\Delta(Q,H)$ erzeugte Gruppe.
 \smallskip
       j) Die Gruppe $G$ operiert scharf zweifach transitiv auf der Menge
 $\Omega$ der von $(P + Q) \cap H$ verschiedenen Punkte auf $P + Q$.
 \smallskip
        Nach 1.6 bzw. auf Grund unserer Annahme operiert $\Delta(P,H)$
 auf $\Omega - \{P\}$ transitiv. Weil $\Omega - \{P\}$ auf Grund
 der Annahme, dass jede Gerade mindestens vier Punkte tr\"agt,
 wenigstens zwei Punkte enth\"alt, enth\"alt $\Delta(Q,H)$
 wenigstens ein von 1 verschiedenes Element, so dass $P$ auf Grund
 von 1.6 bzw. unserer Annahme \"uber $\Delta(Q,H)$ und wegen
 1.1 kein Fixpunkt von $\Delta(Q,H)$ ist. Daher ist $G$ auf
 $\Omega$ zweifach transitiv. Weil $G$ nur aus axialen
 Kollineationen mit der Achse $H$ besteht, folgt mit 1.1, dass das
 einzige Element in $G$, welches zwei Fixpunkte au\ss erhalb $H$
 hat, die Identit\"at ist. Folglich ist $G$ scharf zweifach
 transitiv auf $\Omega$.
 \smallskip
       ij) Es ist $\E((P + Q) \cap H, H) \cap G \neq \{1\}.$
 \smallskip
       Es seien $A$ und $B$ zwei verschiedene Punkte aus $\Omega$. Nach
 j) gibt es ein $\gamma \in G$ mit $A^\gamma = B$ und $B^\gamma = A$. Es folgt
 $A^{\gamma^2} = A$ und $B^{\gamma^2} = B$. Nach 1.1 ist daher
 $\gamma^2 = 1$. Wegen $A \neq B$ ist $\gamma$ daher eine
 Involution. Ist $\gamma \in \E((P + Q) \cap H,H)$, so sind wir
 fertig. Es sei also $\gamma \notin \E((P + Q) \cap H,H)$. Dann
 liegt das Zentrum von $\gamma$ in $\Omega$. Weil dieses Zentrum
 auf Grund von j) kein Fixelement von $G$ ist, gibt es noch eine
 zweite involutorische Streckung in $G$, deren Zentrum vom Zentrum
 von $\gamma$ verschieden ist. Mit 1.10 folgt daher die Behauptung.
 \smallskip
       iij) Es ist $\E((P+Q) \cap H,H) \subseteq G$ und $\E((P + Q) \cap
 H,H)$ operiert transitiv auf $\Omega$.
 \smallskip
       Die Gruppe $\E((P + Q) \cap  H,H) \cap G$ ist nach ij) ein nicht
 tri\-vi\-a\-ler Normalteiler von $G$. Weil $G$ nach j) zweifach transitiv
 operiert,
 operiert dieser Normalteiler auf $\Omega$ transitiv. Mit 1.1 folgt hieraus, dass
 $\E((P + Q) \cap H,H) \cap G = \E((P + Q) \cap H,H)$ ist.
 \par
       Dass schlie\ss lich auch alle $\Delta(X,H)$ f\"ur $X \in \Omega$ in $G$
 enthalten sind und auf $\Omega - \{X\}$ transitiv operieren, folgt aus der
 Transitivit\"at von $G$ auf $\Omega$. Damit ist alles bewiesen.
 \medskip
       Soviel an Bemerkungen zu der Situation im ebenen Fall. Wir wenden
 uns nun wieder der allgemeinen Situation zu.
 \medskip\noindent
 {\bf 1.12. Satz.} {\it Es sei $H$ eine Hyperebene des irreduziblen
 projektiven Verbandes $L$. Sind $\sigma$, $\tau \in \Delta(H)$ und
 ist $\sigma$, $\tau$, $\sigma \tau \neq 1$, ist ferner $P$ das
 Zentrum von $\sigma$, ist $Q$ das Zentrum von $\tau$ und $R$ das
 von $\sigma \tau$, so ist $R \leq P + Q$.}
 \smallskip
       Beweis. Ist $P = Q$, so ist $R = P = Q \leq P + Q$. Es sei also
 $P \neq Q$. Setze $G := P + Q$. Dann ist $G^\sigma = G = G^\tau$. Ist
 $G \not\leq H$, so folgt aus 1.3, dass $R \leq G$ ist. Es sei also
 $G \leq H$. Weil $\sigma$ und $\tau$ Elationen sind, ist auch
 $\sigma \tau$ nach 1.5 eine solche. Folglich ist $R \leq H$. Es
 gibt eine Ebene $E$ --- evtl. das gr\"o\ss te Element von $L$ ---
 mit $G = E \cap H$. Weil $E$ unter $\sigma$ und $\tau$ festbleibt,
 bleibt sie auch unter $\sigma \tau$ fest. Wegen $E \not\leq H$
 folgt mit 1.3, dass $R \leq E$ ist. Also ist $R \leq E \cap H = G$.
 \medskip
       Es sei $L$ ein irreduzibler projektiver Verband, dessen Rang
 mindestens gleich 3 sei. Ferner sei $H$ eine Hyperebene von $L$.
 Ist $0 \neq X \in L$, so bezeichnen wir mit $\Delta(X,H)$ die
 Menge der Perspektivit\"aten, deren Zentrum  in $X$ liegen. Mit
 1.12 folgt dann, dass $\Delta(X)$ eine Untergruppe von
 $\Delta(H)$ ist. Wir setzen noch $\Delta(0,H) := \{1\}$. Es ist
 $\Delta(\Pi,H) = \Delta(H)$ und $\Delta(H,H) = \E(H)$. Ist
 $X \leq H$, so setzen wir $\E(X,H) := \Delta(X,H)$. Dann ist
 $\E(X,H)$ eine Untergruppe von $\E(H)$.
 \par
       Ist $\E(H) \neq \{1\}$, so setzen wir
 $$ \pi(H) := \bigl\{\E(P,H) \mid P \leq H,\ \Rg_L(P) = 1,\ 
      \E(P,H) \neq \{1\}\bigr\}.$$
 Da jedes Element von $\E(H)$ ein Zentrum hat und von 1 verschiedene Elemente
 aber auch nur eines, gilt
 $$ \E(H) = \bigcup_{\Xi \in \pi(H)} \Xi$$
 und
 $$ \Xi \cap \Psi = \{1\}, $$
 falls nur $\Xi$ und $\Psi$ verschiedene Elemente von $\pi(H)$ sind. Dies
 bedeutet, dass $\pi(H)$ eine Partition von $\E(H)$ ist. Dabei nennen wir eine
 Menge $\pi$ von nicht trivialen Untergruppe einer Gruppe $G$
 {\it Partition\/}\index{Partition}{} von $G$, falls $G = \bigcup_{X \in \pi} X$
 ist und sich zwei verschiedene Elemente von $\pi$ stets trivial schneiden. Die
 Elemente von $\pi$ hei\ss en die {\it Komponenten\/}\index{Komponente}{} der
 Partition $\pi$. Die Partition $\pi$ hei\ss t {\it nicht
 trivial\/},\index{nicht triviale Partition}{} falls sie
 mehr als eine Komponente besitzt.
 \par
       Wir werden sehen, dass die Gruppen $\E(H)$ und $\E(P)$ in allen uns
 interessierenden F\"allen abelsch sind. Dies wird aus dem
 folgenden, allgemeineren Satz folgen.
 \medskip\noindent
 {\bf 1.13. Satz.} {\it Es sei $\pi$ eine nicht triviale Partition
 der Gruppe $G$. Genau dann ist $G$ abelsch, wenn alle Komponenten
 von $\pi$ Normalteiler von $G$ sind. Ist $G$ abelsch und enth\"alt
 $G$ ein von $1$ verschiedenes Element endlicher Ordnung, so gibt es
 eine Primzahl $p$, so dass $G$ eine elementarabelsche $p$-Gruppe
 ist.}
 \smallskip
       Beweis. Ist $G$ abelsch, so sind alle Untergruppen von $G$ normal
 in $G$, insbesondere also auch die Komponenten von $\pi$.
 \par
       Es seien alle Komponenten von $\pi$ normal in $G$. Ferner seien $g$ und
 $h$ zwei Elemente aus $G$ und $U$, $V \in \pi$ mit $g \in U$ und $h \in V$. Weil
 $V$ ein Normalteiler ist, gilt $h^{-1}$, $g^{-1}hg \in V$. Also ist
 $$ h^{-1}g^{-1}hg \in V. $$
 Analog folgt $h^{-1}g^{-1}h$, $g \in U$ und damit
 $$ h^{-1}g^{-1}hg \in U. $$
 Ist nun $U \neq V$, so ist $U \cap V = \{1\}$, so dass $gh = hg$ ist. Es sei
 $U = V$. Da $\pi$ nicht
 trivial ist, gibt es eine Komponente $W$ von $\pi$, die von $U$
 verschieden ist. Es sei $1 \neq k \in W$. Da $gh \in U$ ist, folgt
 $ghk = kgh$. Nun ist $kg \notin U$, da sonst $k$ in $U$ l\"age.
 Also ist $kgh = hkg$. Insgesamt ist also
 $$ ghk = kgh = hkg = hgk. $$
 Hieraus folgt, dass $gh = hg$ ist. Damit ist gezeigt, dass $G$ abelsch ist.
 \par
       Nun sei $G$ abelsch. Enth\"alt $G$ ein Element endlicher Ordnung, so
 ent\-h\"alt $G$ auch ein Element von Primzahlordnung $p$. Es sei $g$ ein solches.
 Ferner sei $U \in \pi$ und $g \in U$. Es sei $h \in G - U$. Es gibt dann eine
 Komponente $V$ von $\pi$ mit $h \in V$. Dann ist $hg \notin V$, da $g \notin V$.
 Es gibt also eine von $V$ verschiedene Komponente $W$ mit $hg \in W$. Es folgt,
 da $G$ abelsch ist,
 $$ h^p = h^p g^p = (hg)^p \in V \cap W = \{1\}. $$
 Damit ist gezeigt, dass alle Elemente in $G - U$ die Ordnung $p$ haben. Weil
 $\pi$ nicht trivial ist, ist $G - U \neq \emptyset$, so dass $G$ von $G - U$
 erzeugt wird. Folglich ist $G$ eine elementarabelsche $p$-Gruppe.
 \medskip
       Ist $\kappa$ eine Kollineation von $L$, so ist
 $$ \kappa^{-1}\E(X,H)\kappa = \E(X^\kappa,H^\kappa). $$
 Ist $\kappa \in \Delta(H)$, so ist $X^\kappa = X$ und $H^\kappa = H$.
 Somit gilt
 \medskip\noindent
 {\bf 1.14. Satz.} {\it Ist $L$ ein irreduzibler projektiver
 Verband, ist $H$ eine Hyperebene von $L$ und ist $X \leq H$, so
 ist $\E(X,H)$ normal in $\Delta(H)$ und daher erst recht in
 $\E(H)$. Insbesondere ist $\E(P,H)$ f\"ur jeden Punkt $P$ von $H$
 normal in $\E(H)$.}
 \medskip
       Eine wichtige Folgerung aus diesem Satz ist das n\"achste
 Korollar.
 \medskip\noindent
 {\bf 1.15. Korollar.} {\it Es sei $H$ eine Hyperebene des
 irreduziblen projektiven Verbandes $L$. Ist $\Rg_L(H) \geq 3$, so
 ist $\E(H)$ abelsch. Ist $\Rg_L(H) = 2$, und ist $\pi(H)$ nicht
 trivial, so ist $\E(H)$ abelsch. In beiden F\"allen folgt, dass
 $\E(H)$ eine elementarabelsche $p$-Gruppe ist, wenn $\E(H)$ ein von
 $1$ verschiedenes Element endlicher Ordnung ent\-h\"alt.}
 \medskip
 Ist $L$ ein desarguesscher, irreduzibler projektiver Verband, so
 o\-pe\-riert die Gruppe $\E(H)$ auf der Menge der Punkte von $L_H$
 scharf transitiv, dh., zu je zwei Punkten $P$ und $Q$, die nicht
 in $H$ liegen, gibt es genau ein $\gamma \in \E(H)$ mit $P^\gamma =
 Q$. Ist $L$ endlich und ist $n$ die Anzahl der Punkte von $L_H$,
 so ist also $n = |\E(H)|$. Nach I.7.6 ist, wenn $r$ der Rang und
 $q$ die Ordnung von $L$ ist,
 $$ n = (q - 1)^{-1} (q^r - 1 - q^{r-1} + 1) = q^{r-1}. $$
 Andererseits ist $\E(H)$ endlich und daher nach 1.15 eine elementarabelsche
 $p$-Gruppe, so dass $\E(H)$ eine Potenz von $p$ ist. Folglich ist $q$ eine
 Potenz von $p$. Nennt man eine projektive Ebene $L$ eine
 {\it Translationsebene\/},\index{Translationsebene}{}
 wenn $\E(H)$ auf der Menge der Punkte von $L_H$ transitiv operiert, so liefert
 der eben gef\"uhrte Beweis auch noch, dass die
 Ordnung\index{Ordnung}{} einer endlichen Translationsebene ebenfalls Potenz einer Primzahl
 ist. Es gilt also der im Anschluss an 7.7 in Kapitel I an\-ge\-k\"un\-dig\-te
 \medskip\noindent
 {\bf 1.16. Satz.} {\it Ist $L$ ein endlicher, irreduzibler
 projektiver Verband mit $\Rg(L) \geq 4$ oder ist $L$ eine
 endliche Translationsebene, so ist die Ordnung von $L$ Potenz
 einer Primzahl.}
 \medskip
 Da eine desarguessche projektive Ebene Translationsebene
 bez.\index{endliche desarguessche Ebene}{} jeder ihrer Geraden ist, ist auch die Ordnung einer
 endlichen desarguesschen Ebene Potenz einer Primzahl.

\mysection{2. Der Kern von E(H)}

\noindent
 Es sei $L$ ein irreduzibler projektiver Verband mit $\Rg(L) \geq 3$ und $H$ sei
 eine Hyperebene von $L$. Ist $\Rg(L) = 3$, so
 setzen wir voraus, dass $L$ eine Translationsebene bez.\ $H$ ist.
 Dann ist $\E(H)$ in jedem Falle abelsch. Mit $K(H)$ bezeichnen wir
 die Menge der Endomorphismen $\eta$ von $\E(H)$, die die
 Eigenschaft haben, dass $\E(P,H)^\eta \subseteq \E(P,H)$ ist
 f\"ur alle Punkte $P$ von $H$. Da offensichtlich die Summe, die
 Differenz und das Produkt zweier Elemente aus $K(H)$ wieder in
 $K(H)$ liegen, ist $K(H)$ ein Unterring des Endomorphismenringes
 von $\E(H)$. Man nennt $K(H)$ den {\it Kern\/}\index{Kern}{} von $\E(H)$.
 \medskip\noindent
 {\bf 2.1. Satz.} {\it Es sei $L$ ein irreduzibler projektiver
 Verband mit $\Rg(L) \geq 3$ und $H$ sei eine Hyperebene von $L$.
 Im Falle $\Rg(L) = 3$ sei $L$ eine Translationsebene bez. $H$.
 Dann gilt: $K(H)$ ist ein K\"orper und $\E(H)$ ist ein
 Rechtsvektorraum \"uber $K(H)$. Ist $L$ desarguessch, was f\"ur
 $\Rg(L) > 3$ der Fall ist, so sind die Komponenten von $\pi(H)$
 Unterr\"aume des Ranges $1$ des $K(H)$-Vektorraumes $\E(H)$. Ist $P$
 ein Punkt von $L$, der nicht in $H$ liegt, so ist $\Delta(P,H)$
 zur multiplikativen Gruppe von $K(H)$ isomorph.}
 \smallskip
       Beweis. Es sei $1 \neq \sigma \in \E(H)$ und $\eta \in K(H)$.
 \"Uberdies gelte $\sigma^\eta = 1$. Nun sei $\tau \in \E(H)$ und
 das Zentrum $Q$ von $\tau$ sei vom Zentrum $P$ von $\sigma$
 verschieden. Weil $\sigma$ genau ein Zentrum hat, folgt $\tau
 \sigma \notin \E(Q,H)$. Also ist $\tau \sigma \in \E(R,H)$,
 wobei $R$ ein von $Q$ verschiedener Punkt auf $H$ ist. Hieraus
 folgt
 $$ \tau^\eta = \tau^\eta \sigma^\eta
	      = (\tau\sigma)^\eta \in \E(Q,H) \cap \E(R,H) = \{1\}, $$
 so dass $\tau^\eta = 1$ ist f\"ur alle $\tau \in \E(H) - \E(P,H)$. Da diese
 Menge auf Grund unserer Annahme \"uber $L$ nicht leer ist und $\E(H)$ daher von
 ihr erzeugt wird, gilt $\E(H)^\eta = \{1\}$. Anders ausgedr\"uckt: Jedes von
 Null verschiedene Element in $K(H)$ ist injektiv.
 \par
       Wir zeigen als N\"achstes, dass jedes von Null verschiedene Element aus
 $K(H)$ auch surjektiv ist. Dazu sei $1 \neq \sigma \in \E(P,H)$ und
 $1 \neq \tau \in \E(Q,H)$ sowie $P \neq Q$. Ferner
 sei $X$ ein Punkt, der nicht in $H$ liegt, und $0 \neq \eta \in K(H)$. Dann ist
 $\tau^\eta \neq 1$, wie wir gerade gesehen haben.
 Nun ist $\tau^\eta \in \E(Q,H)$. Daher ist $\sigma \tau^\eta
 \notin \E(P,H)$. Es sei etwa $\sigma \tau^\eta \in \E(R,H)$.
 Nach 1.12 sind die Punkte $P$, $Q$ und $R$ kollinear. Sie sind
 aber auch paarweise verschieden. Die Gerade $X + R$ liegt in der
 Ebene $P + Q + X$, da $R$ ja auf $P + Q$ liegt. Weil diese Ebene durch
 das Zentrum $Q$ von $\tau$ geht, bleibt sie unter $\tau$ fest, so
 dass auch die Gerade $X^\tau + P$ in $P + Q + X$ liegt. Weil die
 Geraden $X + R$ und $X^\tau + P$ nicht in $H$ liegen, sind sie
 verschieden, so dass $(X + R) \cap (X^\tau + P)$ ein Punkt ist, der
 nicht auf $P + Q$ liegt. Somit ist
 $$ Y := \bigl(((X + R) \cap (X^\tau + P)) + Q\bigr)
 \cap (X+P)$$
 ein Punkt auf $P + X$. Weil $Q$ das Zentrum von $\tau$ ist und $P$ auf $H$
 liegt, folgt weiter
 $$ Y^\tau = \bigl(((X + R) \cap (X^\tau + P)) + Q\bigr) \cap (X^\tau + P). $$
 Hieraus folgt mittels des Modulargesetzes, dass
 $$ Y^\tau = (X + R) \cap (X^\tau + P) \leq X + R $$
 ist.
 \par
       Weil $\E(H)$ auf Grund unserer Annahme auf den Punkten des affinen
 Raumes $L_H$ transitiv operiert, gibt es ein $\rho \in \E(P,H)$ mit
 $X^\rho = Y$. Es folgt
 $$ X^{\rho \tau} = Y^\tau \leq X + R. $$
 Dies besagt, dass $\rho \tau \in \E(R,H)$ gilt. Dann ist
 aber auch $\rho^\eta \tau^\eta = (\rho \tau)^\eta \in \E(R,H)$.
 Ferner ist $\sigma \rho^{-\eta} \in E(P, H)$, da $\sigma$,
 $\rho^\eta \in \E(P,H)$ gilt. Also ist
 $$ \sigma \rho^{-\eta} = \sigma \tau^\eta \tau^{-\eta} \rho^{-\eta} \in
       \E(P,H) \cap \E(R,H) = \{1\}. $$
 Folglich ist $\sigma = \rho^\eta$, so dass $\eta$ in der Tat surjektiv ist.
 Damit ist gezeigt, dass die von Null verschiedenen Elemente aus $K(H)$
 Automorphismen sind. Da ihre Inversen offenbar auch zu $K(H)$
 geh\"oren, ist $K(H)$ als K\"orper erkannt. Insbesondere ist
 $\E(H)$ damit ein Rechtsvektorraum \"uber seinem Kern.
 \par
       Es sei $P$ ein Punkt, der nicht in $H$ liegt, und $\delta$ sei ein
 Element aus $\Delta(P,H)$. Dann ist die durch $\tau^{\delta^*}
 := \delta^{-1} \tau \delta$ f\"ur $\tau \in \E(H)$ erkl\"arte
 Abbildung $\delta^*$ nach 1.14 ein Automorphismus von $\E(H)$, der
 sogar in $K(H)$ liegt. Ist $\eta$ ebenfalls in $\Delta(P,H)$
 und ist $\delta^* = \eta^*$, so ist $\delta^{-1} \tau \delta =
 \eta^{-1} \tau \eta$ f\"ur alle $\tau \in \E(H)$. Hieraus folgt,
 dass $\eta \sigma^{-1}$ im Zentralisator von $\E(H)$ liegt, dh.,
 dass $\eta \delta^{-1}$ mit allen Elementen von $\E(H)$
 vertauschbar ist. Weil $\eta \delta^{-1}$ den Fixpunkt $P$ hat und
 $\E(H)$ auf der Menge der Punkte von $L_H$ transitiv operiert,
 l\"asst $\eta \delta^{-1}$ alle Punkte von $L_H$ fest. Folglich
 ist $\eta \delta^{-1} = 1$, so dass die Abbildung $*$ injektiv
 ist. Somit ist $*$ ein Monomorphismus von $\Delta(P,H)$ in die
 multiplikative Gruppe von $K(H)$.
 \par
       Wir zeigen, dass die Abbildung $*$ auch surjektiv ist. Dazu sei
 $0 \neq \eta \in K(H)$. Wir definieren die Abbildung $\sigma$ wie
 folgt. Es sei $Q$ ein Punkt von $L_H$. Es gibt dann genau ein
 $\tau \in \E(H)$ mit $P^\tau = Q$. Wir setzen $Q^\sigma := P^{\tau^\eta}$. Dann
 ist $\sigma$ eine Abbildung der Menge der
 Punkte von $L_H$ in sich. Mittels $\eta^{-1}$ definieren wir auf
 die gleiche Weise eine Abbildung $\sigma'$. Dann ist
 $$ P^{\tau \sigma \sigma'} = (P^{\tau^\eta})^{\sigma'}
                            = P^{\tau^{\eta \eta^{-1}}}
			    = P^\tau = P^{\tau^{\eta^{-1} \eta}}
			    = P^{\tau \sigma' \sigma}, $$
 so dass also $\sigma \sigma' = 1 = \sigma' \sigma$ ist. Somit ist $\sigma$
 bijektiv und $\sigma^{-1} = \sigma'$.
 \par 
       Es seien $Q$, $R$ und $S$ drei kollineare Punkte von $L_H$. Es
 gibt dann ein $\rho \in \E(H)$ mit $Q^\rho = P$. Die Punkte $P$,
 $R^\rho$ und $S^\rho$ sind dann ebenfalls kollinear. Es gibt
 $\lambda$, $\mu \in \E(H)$ mit $P^\lambda = R^\rho$ und $P^\mu = S^\rho$. Wegen
 der Kollinearit\"at von $P$, $R^\rho$ und $S^\rho$
 folgt, dass $\lambda$ und $\mu$ das gleiche Zentrum haben. Dann
 haben aber auch $\lambda^\eta$ und $\mu^\eta$ das gleiche Zentrum.
 Es folgt, dass auch die Punkte $P$, $P^{\lambda^\eta}$ und
 $P^{\mu^\eta}$ kollinear sind. Nun ist aber $Q^{\rho \sigma} = P^\sigma = P$,
 $R^{\rho \sigma} = P^{\lambda^\eta}$ und $S^{\rho \sigma} = P^{\mu^\eta}$.
 Also sind auch die Punkte $Q^{\rho \sigma}$, $R^{\rho \sigma}$ und
 $S^{\rho \sigma}$ kollinear. Wir
 wissen bereits, dass $\sigma^{-1}$ $\rho^{-1}$ $\sigma \in \E(H)$
 ist. Daher sind auch die Punkte $Q^{\rho \sigma \sigma^{-1}
 \rho^{-1} \sigma} = Q^\sigma$, $R^{\rho \sigma \sigma^{-1}
 \rho^{-1} \sigma} = R^\sigma$, $S^{\rho \sigma \sigma^{-1}
 \rho^{-1} \sigma} = S^\sigma$ kollinear. Da $\sigma^{-1}$ eine
 Abbildung gleicher Bauart wie $\sigma$ ist, bildet auch $\sigma^{-1}$ kollineare
 Punkte auf kollineare Punkte ab. Aus I.8.8 folgt daher, falls
 jede Gerade von $L$ wenigstens drei Punkte tr\"agt, dass $\sigma$
 durch genau eine Kollineation von $L$ induziert wird, die, da
 $\sigma$ alle Geraden durch $P$ festl\"asst, in $\Delta(P,H)$ liegt. Enthalten
 alle Geraden von $L$ genau drei Punkte, so
 ist $|\E(Q,H)| = 2$ f\"ur alle Punkte $Q$ von $H$. Dann ist aber
 auch $|K(H)| = 2$, so dass in diesem Falle nichts zu beweisen ist.
 Damit ist gezeigt, dass $*$ ein Isomorphismus von $\Delta(P,H)$
 auf die multiplikative Gruppe von $K(H)$ ist.
 \par

Die Gruppen $\E(Q,H)$ sind Unterr\"aume des $K(H)$-Vektorraumes
 $\E(H)$. Ist $L$ desarguessch, so m\"ussen wir noch zeigen, dass
 sie alle den Rang 1 haben. Es sei also $L$ desarguessch. Sind
 $\sigma$ und $\tau$ von 1 verschiedene Elemente aus $\E(Q,H)$ und
 ist $P$ ein Punkt, der nicht in $H$ liegt, so ist $P \neq P^\sigma$,
 $P^\tau$ und $P + P^\sigma = P + P^\tau$. Es gibt also ein
 $\delta \in \Delta(P,H)$ mit $P^{\sigma \delta} = P^\tau$. Nun
 ist $P^{\delta^{-1}} = P$ und $\delta^{-1} \sigma \delta \in E(H)$. Aus
 $$ P^\tau = P^{\sigma \delta} = P^{\delta^{-1} \sigma \delta} $$
 folgt daher, dass
 $$ \sigma^{\delta^*} = \delta^{-1} \sigma \delta = \tau $$
 ist. Damit ist gezeigt, dass die $\E(Q,H)$ Unterr\"aume des Ranges 1 sind,
 falls $L$ desarguessch ist. Hiermit ist 2.1 in allen seinen Teilen bewiesen.
 \medskip\noindent
 {\bf 2.2. Korollar.} {\it Ist $L$ ein endlicher, desarguesscher
 irreduzibler projektiver Verband der Ordnung $q$ und ist $H$ eine
 Hyperebene von $L$, so ist $K(H) \cong \GF(q)$.}
 \smallskip
       Beweis. Dies folgt sofort aus $|K(H)| = q$ und aus der Tatsache,
 dass es bis auf Isomorphie nur einen endlichen K\"orper mit $q$
 Elementen gibt, n\"amlich das Galoisfeld $\GF(q)$.
 \medskip\noindent
 {\bf 2.3. Korollar.} {\it Es sei $L$ ein endlicher, irreduzibler
 projektiver Verband und $H$ sei eine Hyperebene. Ferner sei $L$
 eine Translationsebene bez.\ $H$, falls $\Rg(L) = 3$ ist. Ist
 dann $P$ ein Punkt von $L_H$, so ist $\Delta(P,H)$ zyklisch.}
 \smallskip
       Beweis. In diesen F\"allen ist $K(H)$ ein endlicher K\"orper, so
 dass die mul\-ti\-pli\-ka\-ti\-ve Gruppe von $K(H)$ nach bekannten S\"atzen
 zyklisch ist.
 \medskip\noindent
 {\bf 2.4. Satz.} {\it Es sei $L$ ein irreduzibler, desarguesscher
 projektiver Verband und $\Sigma$ sei eine Untergruppe von $\E(H)$.
 Genau dann ist $\Sigma$ ein Unterraum des $K(H)$-Vektorraumes
 $\E(H)$, wenn es ein $X \in L$ gibt mit $X \leq H$ und $\Sigma = \E(X,H)$.}
 \smallskip
       Beweis. Aus der Definition von $K(H)$ folgt, dass
 $\E(X,H)$ f\"ur alle $X \in L$ mit $X \leq H$ ein Unterraum des
 $K(H)$-Vektorraumes $\E(H)$ ist.
 \par
       Es sei also umgekehrt $\Sigma$ ein Unterraum des Vektorraumes
 $\E(H)$. Ferner sei $S$ die Menge der Zentren der von 1 verschiedenen Elemente
 aus $\Sigma$. Schlie\ss lich sei $X := \Sigma_{P \in S} P$. Ist
 $S = \emptyset$, so ist $X = 0$ und $\Sigma = \{1\}$, so dass
 $\Sigma = \E(0,H)$ ist. Es sei also $S \neq \emptyset$. Sicherlich ist
 $\Sigma \subseteq \E(X,H)$. Wir m\"ussen also zeigen, dass
 $\E(X,H) \subseteq \Sigma$ ist. Dazu gen\"ugt es zu zeigen, dass $S$ alle Punkte
 von $X$ enth\"alt. Ist n\"amlich $P$ ein Punkt von $X$, der in $S$ liegt, so ist
 $\E(P,H) \cap \Sigma \neq \{1\}$. Hieraus folgt mit 2.1, da $L$ als
 desarguessch vorausgesetzt wurde, dass $\E(P,H) \subseteq \Sigma$
 gilt. Um nun zu zeigen, dass $S$ alle Punkte von $X$ enth\"alt,
 gen\"ugt es nach I.2.12 zu zeigen, dass mit zwei verschiedenen
 Punkten auch alle Punkte ihrer Verbindungsgeraden in $S$ liegen.
 \par
       Es seien also $P$ und $Q$ zwei verschiedene Punkte von $S$. Dann
 ist, wie wir bereits bemerkten, $\E(P,H) \subseteq \Sigma$ und
 $\E(Q,H) \subseteq \Sigma$. Also ist auch
 $$ \E(P,H) E(Q,H) \subseteq \Sigma. $$
 Es sei $R$ ein von $P$ und $Q$ verschiedener Punkt auf $P + Q$. Schlie\ss lich
 sei $U$ ein Punkt von $L_H$ und $1 \neq \tau \in \E(R,H)$. Wir setzen
 $V := (U^\tau + Q) \cap (U + P)$ und $W := (U^\tau + P) \cap (U + Q)$. Es gibt
 dann ein $\rho \in \E(P,H)$ mit $U^\rho = V$ und ein $\sigma \in \E(Q,H)$ mit
 $U^\sigma = W$. Dann ist aber $U^{\rho \sigma} = U^\tau$ und folglich
 $\tau = \rho \sigma$. Somit ist $1 \neq \tau \in \Sigma$, was wiederum
 $R \in S$ nach sich zieht.

\mysection{3. Pappossche Geometrien}

\noindent
 In diesem und dem n\"achsten Abschnitt werden
 wir der Frage nachgehen, in welchen projektiven Geometrien der Satz von
 Pappos\index{Satz von Pappos}{} gilt. Dabei werden wir in diesem Abschnitt nur
 den Fall der de\-sar\-gues\-schen Geometrien betrachten, der ja sicher dann
 vorliegt,
 wenn die betrachteten Geometrien mindestens den Rang 4 haben. Im
 n\"achsten Abschnitt werden wir den Satz von Hessenberg beweisen,
 der besagt, dass eine pappossche projektive Ebene stets auch
 desarguessch ist.
 \par
       Es sei $L$ ein irreduzibler, desarguesscher projektiver Verband,
 in dem der Satz von Pappos gelte. Es sei $H$ eine Hyperebene und
 $P$ ein Punkt von $L$, der auch auf $H$ liegen darf. Weiter seien
 $\delta$ und $\eta$ zwei von 1 verschiedene Elemente aus $\Delta
 (P,H)$ und $G$ und $K$ seien zwei verschiedene Geraden durch
 $P$, die nicht in $H$ liegen. Schlie\ss lich sei $X$ ein von $P$
 und $G \cap H$ verschiedener Punkt auf $G$ und $Y$ ein von $P$ und
 $K \cap H$ verschiedener Punkt auf $K$. Ist $\delta = \eta$, so
 ist $\delta \eta = \eta \delta$. Es sei also $\delta \neq \eta$.
 Dann sind $Y$, $Y^\delta$ und $Y^\eta$ drei verschiedene Punkte
 auf $K$, die alle von $P$ verschieden sind. Ebenso sind
 $X^\delta$, $X^\eta$ und $X^{\delta \eta}$ drei verschiedene
 Punkte auf $G$, die auch von $P$ verschieden sind. Weil wir die
 G\"ultigkeit des Satzes von Pappos vorausgesetzt haben, sind die
 Punkte
 $$\eqalign{
           (X^\delta + Y^\delta) &\cap (X^\eta + Y^\eta),   \cr 
                  (X^\delta + Y) &\cap (X^{\delta \eta} + X^\eta), \cr
    (X^{\delta \eta} + Y^\delta) &\cap (X^\eta + Y)  \cr} $$
 kollinear. Nun ist
 $$ (X^\delta + Y^\delta) \cap H = (X + Y) \cap H = (X^\eta + Y^\eta) \cap H, $$
 so dass der erste dieser Punkte auf $H$ liegt.
 \vadjust{
 \vglue 4mm
 \centerline{\includegraphics{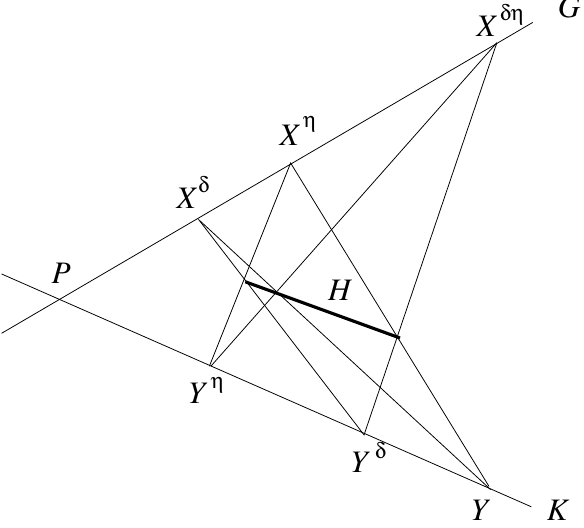}}
 \vglue 3mm}
 Der zweite Punkt liegt aus dem gleichen Grund auf $H$. Somit auch
 der dritte. Nun ist
 $$ X^{\eta \delta} = \bigl(((X^\eta + Y) \cap H) + Y^\delta\bigr) \cap G. $$
 Aus der gerade gemachten Bemerkung
 folgt, dass
 $$ \bigl((X^\eta + Y) \cap H\bigr) + Y^\delta = Y^\delta + X^{\delta \eta} $$
 ist. Also ist
 $$ X^{\eta \delta} = (Y^\delta + X^{\delta \eta}) \cap G = X^{\delta \eta}. $$
 Hieraus folgt schlie\ss lich, dass $\delta \eta = \eta \delta$ ist. Somit
 ist $\Delta(P,H)$ abelsch.
 \par
       Es sei nun umgekehrt $\Delta(P,H)$ abelsch f\"ur alle
 Punkt-Hy\-per\-e\-be\-nen\-paa\-re des desarguesschen projektiven Verbandes
 $L$. (Wie wir wissen, ist $\Delta(P,H)$ sicher dann abelsch,
 wenn $P$ auf $H$ liegt.) Es seien $G$ und $K$ zwei verschiedene
 Geraden durch $P$ und $P_1$, $P_2$ und $P_3$ seien drei
 verschiedene Punkte auf $G$, die alle von $P$ verschieden seien.
 Entsprechend seien $Q_1$, $Q_2$ und $Q_3$ drei verschiedene Punkte
 auf $K$, die ebenfalls von $P$ verschieden seien. Schlie\ss lich sei
 $$ J:= \bigl((P_1 + Q_2) \cap (Q_1 + P_2)\bigr)
		    + \bigl((P_1 + Q_3) \cap (Q_1 + P_3)\bigr). $$
 Dann ist $J$ eine Gerade in der Ebene $G + K$. Es sei $A$ ein Komplement von
 $G + K$. Dann ist $H := J + A$ eine Hyperebene von $L$ mit $H \cap (G + K) = J$.
 Hieraus folgt, dass die Punkte $P_i$ und $Q_j$ nicht in $H$ liegen. Es gibt
 $\delta$, $\eta \in \Delta(P,H)$ mit $Q_3^\delta = Q_2$ und
 $Q_3^\eta = Q_1$. Setze
 $$ X:= \bigl(Q_3 + ((P_1 + Q_2) \cap (Q_1 + P_2))\bigr) \cap G. $$
 Weil die Gerade $G$ durch $P$ geht und der
 Punkt $(P_1 + Q_2) \cap (Q_1 + P_2)$ auf $H$ liegt, folgt mit
 Hilfe des Modulargesetzes
 $$\eqalign{
   X^\delta &= \bigl(Q_3^\delta + ((P_1 + Q_2) \cap (Q_1 + P_2))\bigr) \cap G \cr
            &= \bigl(Q_2 + ((P_1 + Q_2) \cap (Q_1 + P_2))\bigr) \cap G        \cr
	    &= (P_1 + Q_2) \cap (Q_2 + Q_1 + P_2) \cap G            \cr
	    &= (P_1 + Q_2) \cap G                                   \cr
	    &= P_1.                                                 \cr} $$
 Ebenso folgt
 $$\eqalign{
     X^\eta &= \bigl(Q_3^\eta + ((P_1 + Q_2) \cap (Q_1 + P_2))\bigr) \cap G \cr
	    &= \bigl(Q_1 + ((P_1 + Q_2) \cap (Q_1 + P_2))\bigr) \cap G      \cr
	    &= (Q_1 + P_2) \cap (Q_1 + P_1 + Q_2) \cap G          \cr
	    &= (Q_1 + P_2) \cap G                                 \cr
	    &= P_2 \cr} $$
 und
 $$\eqalign{
    P^\eta_1 &= \bigl(Q_3^\eta + ((P_1 + Q_3) \cap (Q_1 + P_3))\bigr) \cap G \cr
             &= \bigl(Q_1 + ((P_1 + Q_3) \cap (Q_1 + P_3))\bigr) \cap G      \cr
	     &= (Q_1 + P_3) \cap (Q_1 + P_1 + Q_3) \cap G            \cr
	     &= (Q_1 + P_3) \cap G                                   \cr
	     &= P_3 \cr} $$
 Es ist also $X^\delta = P_1$, $X^\eta = P_2$ und
 $X^{\delta \eta} = P^\eta_1 = P_3$. Nun ist $\delta \eta = \eta
 \delta$, woraus folgt, dass
 $$ P_2^\delta = X^{\eta \delta} = X^{\delta \eta} = P_3 $$
 ist. Daraus folgt wiederum, dass $(P_2 + Q_3) \cap (Q_2 + P_3)$ ein Punkt von
 $H$ ist. Also ist
 $$ (P_2 + Q_3) \cap (Q_2 + P_3) \leq H \cap (G + K) = J, $$
 so dass in $L$ der Satz von Pappos gilt. Damit ist ein Teil des folgenden
 Satzes bewiesen.
 \medskip\noindent
 {\bf 3.1. Satz.} {\it Ist $L$ ein irreduzibler, desarguesscher
 projektiver Verband mit $\Rg(L) \geq 3$, so sind die folgenden
 Bedingungen \"aquivalent:
 \item{(a)} In $L$ gilt der Satz von Pappos.
 \item{(b)} $\Delta(P,H)$ ist f\"ur alle Punkt-Hyperebenenpaare $(P,H)$ abelsch.
 \item{(c)} Es gibt ein nicht inzidentes Punkt-Hyperebenenpaar
 $(P,H)$, so dass $\Delta(P,H)$ abelsch ist.
 \item{(d)} Es gibt eine Hyperebene $H$, so dass $K(H)$ kommutativ
 ist.}
 \smallskip
       Beweis. Die \"Aquivalenz von (a) und (b) haben wir gerade bewiesen
 und die \"Aquivalenz von (c) und (d) folgt aus Satz 2.1. Bedingung
 (c) ist nat\"urlich eine Folge von (b). Es bleibt zu zeigen, dass
 (b) eine Folge von (c) ist. Um dies zu zeigen, beweisen wir
 zun\"achst, dass die von allen $\E(H)$ erzeugte Gruppe auf der
 Menge der nicht inzidenten Punkt-Hyperebenenpaare transitiv
 operiert.
 \par
       Es seien $H$ und $K$ zwei verschiedene Hyperebenen von $L$. Dann
 ist $H \cap K$ eine Ko-Gerade von $L$. Es gibt also eine Gerade
 $G$ mit $\Pi = G \oplus (H \cap K)$. Weil $G$ mindestens drei
 Punkte tr\"agt, gibt es einen Punkt $P$ auf $G$, der weder auf $H$
 noch auf $K$ liegt. Setze $V := P + (H \cap K)$. Dann ist $V$ eine
 Hyperebene. Weil $L$ desarguessch ist, gibt es ein $\tau \in \E(P,V)$ mit
 $(G \cap H)^\tau = G \cap K$. Es folgt, dass $H^\tau = K$ ist. Die fragliche
 Gruppe ist also auf der Menge der Hyperebenen transitiv. Weil $\E(H)$ auf der
 Menge der Punkte von $L_H$ transitiv operiert, ist die fragliche Gruppe auf der
 Menge der nicht inzidenten Punkt-Hyper\-ebenen\-paare transitiv.
 \par
 Sei nun $(P,H)$ ein nicht inzidentes Punkt-Hyperebenenpaar
 und die Gruppe $\Delta(P,H)$ sei abelsch. Ist $(Q,K)$ ein
 weiteres, nicht inzidentes Punkt-Hyp\-er\-ebe\-nen\-paar, so gibt es also
 eine Kollineation $\gamma$ mit $P^\gamma = Q$ und $H^\gamma = K$.
 Es folgt
 $$ \gamma^{-1}\Delta(P,H)\gamma = \Delta(P^\gamma,H^\gamma) = \Delta(Q,K), $$
 so dass auch $\Delta(Q,K)$ abelsch ist.
 \par
       Da die Gruppen $\E(P,H)$ f\"ur alle inzidenten
 Punkt-Hy\-per\-e\-be\-nen\-paa\-re stets abelsch sind, ist alles bewiesen.
 \medskip
       Unsere Zwischenbemerkung ist wichtig genug, um als Satz formuliert zu
 werden.
 \medskip\noindent
 {\bf 3.2. Satz.} {\it Ist $L$ ein irreduzibler, desarguesscher
 projektiver Verband, so ist die von allen Elationen erzeugte
 Untergruppe der Kol\-li\-ne\-a\-ti\-ons\-grup\-pe von $L$ auf der Menge der
 nicht inzidenten Punkt-Hy\-per\-ebe\-nen\-paa\-re von $L$ transitiv.}
 \medskip
       Aus 3.1 und 3.2 folgt noch das
 \medskip\noindent
 {\bf 3.3. Korollar.} {\it In allen endlichen, desarguesschen
 irreduziblen projektiven Geometrien gilt der Satz von Pappos.}
 \medskip
       Alle Beweise von 3.3, die ich kenne, benutzen den Satz von
 Wedderburn,\index{Satz von Wedderburn}{} dass alle endlichen K\"orper kommutativ
 sind. Einen ge\-o\-met\-ri\-schen Beweis zu finden, scheint also sehr schwierig
 zu sein. Den Beweis von Tecklenburg\index{Tecklenburg, H.}{} 1987 kann ich nicht
 als einen solchen werten, da er nur den einfachen wittschen\index{Witt, E.}{}
 Beweis des
 wedderburnschen Satzes in eine komplizierte geometrische Sprache \"ubersetzt.

\mysection{4. Der Satz von Hessenberg}

\noindent
 Hessenberg\index{Hessenberg, G.}{} zeigt, wenn
 auch mit einigen L\"ucken im Beweis, dass pappossche projektive
 Ebenen\index{pappossche Ebene}{} auch desarguessch sind (Hessenberg 1905). Sein
 Beweis und einige andere sind in Pickerts Buch abgedruckt. Diese Beweise haben
 jedoch niemanden so recht befriedigt, da sie wegen der vielen zu betrachtenden
 Ent\-ar\-tungs\-f\"alle\index{Entartungsf\"alle}{} sehr un\"ubersichtlich
 sind. Entartungsf\"alle in algebraischen
 Situationen\index{algebraische Situation}{} sind meist banal, w\"ahrend sie in geometrischen
 Situationen\index{geometrische Situation}{} h\"aufig sehr
 viel Kopfzerbrechen be\-rei\-ten. Gesucht war also ein Beweis, bei dem
 man schon fr\"uh von algebraischen Methoden Gebrauch machen
 konnte. Vorbild war der Satz von Baer, dass die G\"ultigkeit des
 $(P,G)$-Satzes von Desargues der $(P,G)$-Transitivit\"at der
 fraglichen Ebene \"aquivalent ist. Einen ersten Beweis dieser Art
 fand ich\index{luneburg@L\"uneburg, H.}{} im Jahre 1968, den ich eine Weile als den
 sch\"onsten Beweis f\"ur den hessenbergschen Satz ansah (L\"uneburg 1969b). Dann
 jedoch publizierte A.~Herzer\index{Herzer, A.}{} seinen Beweis dieses
 Satzes (Herzer 1972). Dieser Beweis
 z\"ahlt f\"ur mich zu den mathematischen Juwelen. Ich werde ihn
 also, Altruist\index{Altruist}{} der ich bin, dem Leser nicht vorenthalten.
 \par
       Vom geometrischen Standpunkt her gesehen sind
 Kollineationen\index{Kollineation mit Fixpunkten}{} mit vielen
 Fixpunkten einfacher zu handhaben als solche ohne
 Fixpunkte. Das lassen schon unsere S\"atze \"uber
 Perspektivit\"aten ahnen. Antiautomorphismen eines projektiven
 Verbandes haben keine Fixpunkte. Dennoch gibt es auch hier Punkte,
 die vor anderen aus\-ge\-zeich\-net sind. Bevor wir sie definieren,
 verabreden wir noch, dass Antiautomorphismen eines projektiven
 Verbandes in Zukunft {\it Korrelationen\/}\index{Korrelation}{} oder auch
 {\it Dualit\"aten\/}\index{Dualit\"at}{} genannt werden.
 \par
       Es sei $\kappa$ eine Korrelation des projektiven Verbandes $L$.
 Der Punkt $P$ von $L$ hei\ss t {\it absolut\/},\index{absoluter Punkt}{} falls
 $P \leq P^\kappa$ gilt.  Analog hei\ss t die Hyperebene $H$ von $L$ {\it
 absolut\/},\index{absolute Hyperebene}{} wenn $H^\kappa \leq H$ ist.
 \medskip\noindent
 {\bf 4.1. Satz.} {\it Es sei $\kappa$ eine Korrelation der
 projektiven Ebene $L$. Sind $P$ und $Q$ zwei verschiedene absolute
 Punkte von $\kappa$ und gilt $Q \leq P^\kappa$, so sind $P$ und
 $Q$ die einzigen absoluten Punkte auf $P^\kappa$. \"Uberdies ist
 $P^{\kappa^2} = Q$ und $\kappa^2 \neq 1$.}
 \smallskip
       Beweis. Setze $G:=P^\kappa$. Die Menge der Geraden durch $P$ wird
 von $\kappa$ bijektiv auf die Menge der Punkte von $G$ abgebildet.
 Da $G$ nach Voraussetzung mindestens zwei absolute Punkte tr\"agt,
 gibt es eine von $G$ verschiedene Gerade $H$ durch $P$, so dass $H^\kappa$ ein
 absoluter Punkt auf $G$ ist, so dass also $H^\kappa \leq H^{\kappa^2}$ gilt.
 Hieraus folgt, da auch $K^{-1}$ eine Korrelation ist, dass
 $$ (H^{\kappa^2})^{\kappa^{-1}} \leq (H^{\kappa})^{\kappa^{-1}}, $$
 dh., dass
 $$H^\kappa \leq H$$ 
 ist. Also ist
 $$ H^\kappa = H \cap G = P. $$
 Dies besagt, dass es nur eine von $G$ verschiedene Gerade durch $P$
 gibt, deren Bild unter $\kappa$ ein absoluter Punkt auf $G$ ist,
 n\"amlich die Gerade $P^{\kappa^{-1}}$. Weil $Q$ ein zweiter
 absoluter Punkt auf $G$ ist, ist also $G^\kappa = Q$. Hieraus
 folgt weiter
 $$ P^{\kappa^2} = G^\kappa = Q. $$
 Weil schlie\ss lich $P \neq Q$ ist, ist $\kappa^2 \neq 1$. Damit ist alles
 bewiesen.
 
 \medskip

       Es sei $\kappa$ eine Korrelation der projektiven Ebene $L$ und $G$
 und $H$ seien zwei verschiedene Geraden von $L$. Die Korrelation
 $\kappa$ hei\ss e $(G,H)$-{\it Korrelation\/},\index{Korrelation}{} falls
 die Punkte von $G$ wie auch die Punkte von $H$ absolute Punkte von $\kappa$
 sind. Es gilt nun der folgende Satz.
 \medskip\noindent
 {\bf 4.2. Satz.} {\it Es seien $G$ und $H$ zwei verschiedene Geraden der
 projektiven Ebene $L$ und $\kappa$ sei eine $(G,H)$-Korrelation von $L$. Es
 seien ferner $P$ und $Q$ die beiden
 verschiedenen Punkte von $L$ mit $P^\kappa = G$ und $Q^\kappa = H$.
 Schlie\ss lich sei $S:= G \cap H$ und $V := P + Q$. Dann gilt:
 \item{a)} $P \not\leq G$ und $Q \not\leq H$.
 \item{b)} Ist $X$ eine Gerade durch $P$, so ist $X$ absolut und
 $X^\kappa = X \cap G$. Ist $Y$ eine Gerade durch $Q$, so ist $Y$
 absolut und $Y^\kappa = Y \cap H$.
 \item{c)} Ist $X$ ein Punkt von $L_V$, so ist
 $$ X^\kappa = \bigl((X + P) \cap G\bigr) + \bigl((X + Q) \cap H\bigr). $$
 \item{d)} Es ist $S \leq V$.
 \item{e)} Es ist $G^\kappa = Q$ und $H^\kappa = P$.
 \item{f)} Es ist $S^\kappa = V$ und $V^\kappa = S$.
 \item{g)} Ist $X$ eine Gerade mit $S \not\leq X$, so ist
 $$ X^\kappa = \bigl((X \cap G) + Q\bigr) \cap \bigl((X \cap H) + P\bigr). $$
 \item{h)} Durch $P$, $Q$, $G$ und $H$ wird $\kappa$ eindeutig festgelegt.
 \item{i)} Ist $A$ ein absoluter Punkt von $\kappa$, so ist $A \leq
 G$ oder $A \leq H$.}
 \smallskip
       Beweis. a) Weil $G$ mindestens drei Punkte tr\"agt und alle diese
 Punkte absolut sind, folgt mit 4.1 wegen $P^\kappa = G$, dass $P$
 nicht auf $G$ liegt. Ebenso folgt, dass $Q$ nicht auf $H$ liegt.
 \par
       b) Es sei $X$ eine Gerade durch $P$. Dann ist $X^\kappa \leq
 P^\kappa = G$, so dass $X^\kappa$ nach Voraussetzung absolut ist.
 Also ist $X^\kappa \leq X^{\kappa^2}$. Anwendung von $\kappa^{-1}$
 zeigt, dass $X^\kappa \leq X$ ist. Somit ist $X$ absolut. Wegen
 $X^\kappa \leq G$ und $X^\kappa \leq X$ folgt schlie\ss lich
 $X^\kappa = X \cap G$. Ebenso folgt die Aussage \"uber $Y$.
 \par
       c) Es sei $X$ ein Punkt,
 \vadjust{
 \vglue 4mm
 \centerline{\includegraphics{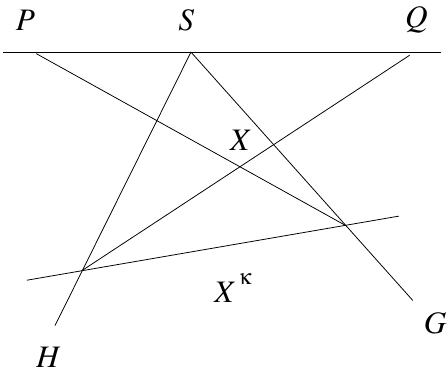}}
 \vglue 3mm}
 der nicht auf $V$ liegt. Wegen $V = P + Q$ ist
 dann $X = (X + P) \cap (X + Q)$. Mit b) folgt daher
 $$ X^\kappa = (X + P)^\kappa + (X + Q)^\kappa
	     = \bigl((X + P) \cap G\bigr) + \bigl((X + Q) \cap H\bigr).$$
 \par
       d) W\"are $S \not\leq V$, so folgte mit c) der Widerspruch
 $$ S^\kappa = \bigl((S + P) \cap G\bigr) + \bigl((S + Q) \cap H\bigr)
	     = S + S = S. $$
 \par
       e) Mit a) folgt $V \neq G$, so dass nach d) gilt, dass $S = G \cap V$ ist.
 Ist nun $X$ ein von $S$ verschiedener Punkt auf $G$, so
 ist also $X \not\leq V$. Nach c) und dem Modulargesetz ist also
 $$ X^\kappa = \bigl((X + P) \cap G\bigr) + \bigl((X + Q) \cap H\bigr)
	     = X + \bigl((X + Q) \cap H\bigr) = X + Q. $$
 Ist nun $Y$ ein weitere von $X$ und $S$ verschiedener Punkt auf $G$, so folgt,
 wie gerade gesehen, $Y^\kappa = Y + Q$ und daher
 $$ G^\kappa = (X + Y)^\kappa = X^\kappa \cap Y^\kappa = (X + Q) \cap (Y + Q)
	     = Q. $$
 ebenso folgt, dass $H^\kappa = P$ ist.
 \par
       f) Es ist
 $$ S^\kappa = (G \cap H)^\kappa = G^\kappa + H^\kappa = Q + P = V. $$
 Entsprechend gilt
 $$ V^\kappa = (P + Q)^\kappa = P^\kappa \cap Q^\kappa = G \cap H = S. $$
 \par
       g) Ist $X$ eine Gerade, die nicht durch $S$ geht, so ist
 $X = (X \cap G) + (X \cap H)$. Daher ist
 $$ X^\kappa = (X \cap G)^\kappa \cap (X \cap H)^\kappa
	     = \bigl((X \cap G) + Q\bigr) \cap \bigl((X \cap H) + P\bigr). $$
 \par
       h) Es sei $\lambda$ eine weitere $(G,H)$-Korrelation mit
 $P^\lambda = G$ und $Q^\lambda = H$. Nach c) gilt dann die
 Gleichung $X^\kappa = X^\lambda$ f\"ur alle Punkte $X$ von $L_V$.
 Daher gilt $X^{\kappa \lambda^{-1}} = X$ f\"ur alle diese Punkte.
 Nach I.8.8 ist daher $\kappa \lambda^{-1} = 1$, so dass in der Tat
 $\kappa = \lambda$ ist. Es sei $A$ ein absoluter Punkt von
 $\kappa$. Ist $A \leq V$, so ist $S = V^\kappa \leq A^\kappa$.
 W\"are $A \neq S$, so folgte $A^\kappa \neq S^\kappa = V$. Weil
 $A$ absolut ist, folgte weiter $A = A^\kappa \cap V = S$. Also ist
 doch $A = S$. Es sei also $A \not\leq V$. Dann ist
 $$ A \leq A^\kappa = \bigl((A + P) \cap G\bigr) + \bigl((A + Q) \cap H\bigr). $$
 Wir
 d\"urfen annehmen, dass $A$ nicht auf $G$ liegt. Dann ist $A$ ein
 von $(A + P) \cap G$ verschiedener Punkt auf $A^\kappa$. Es ist
 also $A^\kappa = A + ((A + P) \cap G)$. Mittels des
 Modulargesetzes folgt weiter
 $$ A^\kappa = (A + P) \cap (A + G) = A + P. $$
 Also gilt $(A + Q) \cap H \leq A + P$. Weil $A$
 nicht auf $V$ liegt, ist $A + P \neq A + Q$. Daher ist
 $$ (A + Q) \cap H = (A + P) \cap (A + Q) = A, $$
 so dass $A$ in der Tat auf $H$ liegt.
 \par
       Die zweite Aussage von i) ist dual zur ersten. Damit ist der Satz
 vollst\"andig bewiesen.
 \medskip
 Hier nun der alles entscheidende Satz von Herzer, der die
 Ver\-bin\-dung des Satzes von Pappos zu den $(G,H)$-Korrelationen herstellt.
 \medskip\noindent
 {\bf 4.3. Satz.} {\it Ist $L$ eine projektive Ebene, so sind die
 beiden folgenden Aussagen \"aquivalent:
 \item{1)} In $L$ gilt der Satz von Pappos.
 \item{2)} Sind $G$ und $H$ zwei verschiedene Geraden und $P$ und $Q$ zwei
 verschiedene Punkte von $L$, gilt ferner $P$, $Q \not\leq G$, $H$ und sind die
 Punkte $P$, $Q$ und $G \cap H$ kollinear, so gibt es eine $(G,H)$-Korrelation
 $\kappa$ von $L$ mit $P^\kappa = G$ und $Q^\kappa = H$.}
 \smallskip
       Beweis. 1) impliziert 2): Es seien $P$ und $Q$ zwei Punkte und $G$
 und $H$ zwei Geraden von $L$, die die Voraussetzungen von 2)
 erf\"ullen. Setze $V := P + Q$ und $S := G \cap H$. Setze ferner
 $$ L_{S,V} := L_V \cap (L^d)_S. $$
 Dann besteht $L_{S, V}$
 also aus den Punkten, die nicht auf $V$ liegen, und den Geraden,
 die nicht durch $S$ gehen. Wir definieren zwei Abbildungen
 $\kappa$ und $\lambda$ von $L_{S, V}$ in sich wie folgt: ist $X$
 ein Punkt von $L_{S, V}$, so setzen wir
 $$ X^\kappa := \bigl((X + P) \cap G\bigr) + \bigl((X + Q) \cap H\bigr) $$
 und ist $X$ eine Gerade von $L_S, V$, so setzen wir
 $$ X^\kappa := \bigl((X + P) \cap G\bigr) + \bigl((X + Q) \cap H\bigr) $$
 und ist $X$ eine Gerade von $L_{S,V}$, so setzen wir
 $$ X^\kappa := \bigl((X \cap G) + Q\bigr) \cap \bigl((X \cap H) + P\bigr). $$
 Ist $Y$ ein Punkt von $L_{S, V}$, so setzen wir ferner
 $$ Y^\lambda := \bigl((Y + Q) \cap  G\bigr) + \bigl((Y + P) \cap H\bigr). $$
 Schlie\ss lich setzen wir
 $$ Y^\kappa := \bigl((Y \cap G) + P\bigr) \cap \bigl((Y \cap H) + Q\bigr),$$
 falls $Y$ eine Gerade von $L_{S, V}$ ist.
 \par
       Ist $X$ ein Punkt von $L_{S, V}$, so folgt mittels des Modulargesetzes
 $X^\kappa \cap G = (X + P) \cap G$ und $X^\kappa \cap H = (X + Q) \cap H$.
 Hieraus folgt weiter
 $$ X^\kappa \cap G) + P = \bigl((X + P) \cap G\bigr) + P = X + P $$
 und
 $$ X^\kappa \cap H) + Q = \bigl((X + Q) \cap H\bigr) + Q = X + Q. $$
 Daher ist
 $$ X^{\kappa \lambda} = \bigl((X^\kappa \cap G) + P\bigr) \cap
			       \bigl((X^\kappa \cap H) + Q\bigr)
                       = (X + P) \cap (X + Q) = X.$$
 \par
       Ist $X$ eine Gerade von $L_{S,V}$, so folgt mittels des Modulargesetzes
 $X^\kappa + P = (X \cap H) + P$ und $X^\kappa + Q = (X \cap G) + Q$. Hieraus
 folgt weiter
 $$ (X^\kappa + P) \cap H = X \cap H $$
 und
 $$ (X^\kappa + Q) \cap G = X \cap G. $$
 Somit ist
 $$X^{\kappa \lambda} = \bigl((X^\kappa + Q) \cap G\bigr) +
			      \bigl((X^\kappa + P) \cap H\bigr)
		      = (X \cap G) + (X \cap H) = X. $$
 Also ist $\kappa \lambda = 1$. Aus Dualit\"atsgr\"unden, --- den Satz von
 Pappos haben wir ja noch nicht benutzt ---, ist dann auch $\lambda \kappa = 1$,
 so dass $\kappa$ eine Bijektion von $L_{S,V}$ auf
 sich und dass $\lambda$ die zu $\kappa$ inverse Abbildung ist.
 \par
       Es sei $X$ ein Punkt und $Y$ eine Gerade von $L_{S, V}$ und es
 gelte $X \leq Y$. Wir zeigen, dass $Y^\kappa \leq X^\kappa$ gilt.
 Es sei zun\"achst $X \leq G$. Dann folgt unter Benutzung des
 Modulargesetzes
 $$\eqalign{
    X^\kappa &= \bigl((X + P) \cap G\bigr) + \bigl((X + Q) \cap H\bigr) \cr
	     &= X + \bigl((X + Q) \cap H\bigr)                \cr
	     &=(X + Q) \cap (X + H) = X + Q. \cr} $$
 Andererseits ist, da ja $Y \cap G = X$ ist,
 $$ Y^\kappa  = \bigl((Y \cap G) + Q\bigr) \cap \bigl((Y \cap H) + P\bigr)
	      = (X + Q) \cap \bigl((Y \cap H) + P\bigr), $$
 so dass in diesem Falle $Y^\kappa \leq X + Q = X^\kappa$ gilt. Liegt $X$ auf
 $H$, so folgt genauso, dass $Y^\kappa \leq X^\kappa$ ist. Die F\"alle, dass $P$
 oder $Q$ auf $Y$ liegen, sind dual zu den behandelten, so dass auch hier
 $Y^\kappa \leq X^\kappa$ gilt, da wir den Satz von Pappos noch immer nicht
 benutzt haben.
 \par
       Wir d\"urfen daher des Weiteren annehmen, dass $X$ weder auf $G$
 noch auf $H$ liegt und dass $Y$ weder durch $P$ noch durch $Q$
 geht. Dann sind die Punkte $P$, $S$ und $Q$ drei verschiedene
 Punkte auf $V$, die von $V \cap Y$ verschieden sind. Ferner sind
 $Y \cap G$, $X$ und $Y \cap H$ drei verschiedene Punkte auf $Y$,
 die ebenfalls von $V \cap Y$ verschieden sind. Weil in $L$ der
 Satz von Pappos gilt, sind daher die Punkte
 $$\eqalign{
    \bigl((Y \cap G) + Q\bigr) &\cap \bigl((Y \cap H) + P\bigr), \cr
                       (X + P) &\cap \bigl((Y \cap G) + S\bigr), \cr
		       (X + Q) &\cap \bigl((Y \cap H) + S\bigr)  \cr} $$
 kollinear. Wegen $(Y \cap G) + S = G$ und $(Y \cap H) + S = H$
 folgt daher, dass der Punkt
 $$Y^\kappa = \bigl((Y \cap G) + Q\bigr) \cap \bigl((Y \cap H) + P\bigr) $$
 auf der Geraden
 $$ \bigl((X + P) \cap G\bigr) + \bigl((X + Q) \cap H\bigr) = X^\kappa $$
 liegt. Aus $X \leq Y$ folgt also stets $Y^\kappa \leq X^\kappa$.
 \par
       Da $\lambda$ von gleicher Bauart ist --- es haben ja nur $P$ und
 $Q$ ihre Rollen vertauscht ---, gilt auch f\"ur $\lambda$, dass
 $Y^\lambda \leq X^\lambda$ eine Folge von $X\leq Y$ ist. Somit ist
 $\kappa$ ein Antiautomorphismus von $L_{S,V}$.
 \par
       Es seien nun $X$ und $Y$ zwei verschiedene Punkte von $L_{S,V}$, die in
 $L$ auf einer Geraden durch $S$ liegen. W\"are nun $Z := X^\kappa \cap Y^\kappa$
 ein Punkt von $L_{S,V}$, so w\"are
 $Z^\lambda$ eine Gerade von $L_{S,V}$ mit $X$, $Y \leq Z^\lambda$.
 Es folgte der Widerspruch $S \leq X + Y = Z^\lambda$. Dies
 impliziert, dass die Punkte einer Geraden durch $S$ von $\kappa$
 auf konfluente Geraden abgebildet werden, deren gemeinsamer
 Schnittpunkt auf $V$ liegt. Bei diesem Schluss haben wir nur
 davon Gebrauch gemacht, dass $\kappa$ ein Antiautomorphismus von
 $L_{S,V}$ ist und dass $\lambda = \kappa^{-1}$ gilt. Aus
 Dualit\"atsgr\"unden bildet $\lambda$ daher konfluente Geraden,
 die sich in einem von $S$ verschiedenen Punkt auf $V$ schneiden,
 auf Punkte ab, die auf einer Geraden durch $S$ liegen. Dies ist
 nun mehr als genug, um mit I.8.8 zu schlie\ss en, dass $\kappa$
 durch eine Korrelation von $L$ induziert wird, die wir ebenfalls
 $\kappa$ nennen.
 \par
       Es sei $X$ ein von $S$ verschiedener Punkt auf $G$. Dann ist
 $$\eqalign{
   X^\kappa &= \bigl((X + P) \cap G\bigr) + \bigl((X + Q) \cap H\bigr)
	     = X + \bigl((X + Q) \cap H\bigr) \cr
	    &= (X + Q) \cap (X + H) = X + Q, \cr} $$
 da ja $X$ nicht auf $H$ liegt. Hieraus folgt, dass $X$ ein
 absoluter Punkt und dass $G^\kappa = Q$ ist. Entsprechend folgt,
 dass die von $S$ verschiedenen Punkte von $H$ absolut sind und
 dass $H^\kappa = P$ ist. Also ist $S^\kappa = (G \cap H)^\kappa = Q + P = V$.
 Somit ist auch $S$ absolut und $\kappa$ als $(G,H)$-Korrelation erkannt.
 \par
       Es sei $Z$ eine von $V$ verschiedene Gerade durch $P$. Setzt man
 $X := Z \cap G$ und $Y := Z \cap H$, so ist
 $$ Z^\kappa = (X + Y)^\kappa = X^\kappa \cap Y^\kappa = (X + Q) \cap (Y + P)
	     = (X \cap Q) \cap Z = X. $$
 Hieraus folgt, dass $P^\kappa = G$ ist. Ebenso folgt $Q^\kappa = H$. Damit ist
 gezeigt, dass 2) eine Folge von 1) ist.
 \par
       Es bleibe dem Leser \"uberlassen zu zeigen, dass 1) eine Folge von 2) ist.
 \medskip
 Nun sind wir endlich in der Lage, den Satz von
 Hessenberg\index{Satz von Hessenberg}{} zu beweisen.
 \medskip\noindent
 {\bf 4.4. Satz von Hessenberg.} {\it Jede pappossche Ebene ist desarguessch.}
 \smallskip
       Beweis. Es sei $(P,H)$ ein nicht inzidentes Punkt-Geradenpaar
 der pap\-pos\-schen projektiven Ebene $L$. Ferner sei $V$ eine Gerade
 durch $P$. Setze $S := V \cap H$ und $G$ sei eine von $V$ und $H$
 verschiedene Gerade durch $S$. Schlie\ss lich seien $X$ und $Y$
 zwei Punkte auf $V$, die von $P$ und $S$ verschieden sind. Nach
 4.3 gibt es dann zwei $(G,H)$-Korrelationen $\kappa$ und
 $\lambda$ mit $P^\kappa = P^\lambda = G$ und $X^\kappa = Y^\lambda
 = H$. Setze $\sigma := \kappa \lambda^{-1}$. Dann ist $P^\sigma = P$ und
 $X^\sigma = Y$. Ist $Z$ ein Punkt auf $H$, so folgt mit 4.2
 $$ Z^\sigma = Z^{\kappa \sigma^{-1}} = (Z + P)^{\lambda^{-1}} = Z. $$
 Folglich ist $\sigma \in \Delta(P,H)$
 und $L$ ist als $(P,H)$-transitiv erkannt. Da dies f\"ur alle
 nicht inzidenten Punkt-Geradenpaare von $L$ gilt, folgt mit 1.11,
 dass $L$ f\"ur alle Punkt-Geradenpaare eine $(P,H)$-transitive
 Ebene ist, seien sie inzident oder nicht. Damit ist der Satz von
 Hessenberg bewiesen.
 \medskip
       Ich denke, dass der Leser mit mir einer Meinung ist, dass dieser
 Beweis des hessenbergschen Satzes ein Juwel ist.

\mysection{5. Weniger Bekanntes aus der linearen Algebra}

\noindent
 Ist $V$ ein Rechtsvektorraum
 \"uber dem K\"orper $K$, so ist
 $L(V)$ mit der Inklusion als Teilordnung ein projektiver Verband.
 Ist $H$ eine Hyperebene dieses Verbandes und $P$ ein Punkt, der
 nicht auf $H$ liegt, so haben wir die Objekte $\Delta(P,H)$,
 $\E(H)$ und $K(H)$, die wir in den ersten beiden Abschnitten dieses
 Kapitels definiert haben. Es geht uns nun darum, f\"ur diese
 Objekte eine algebraische Beschreibung zu finden. Es wird sich
 herausstellen, dass $K(H)$ zu $K$ isomorph ist, --- wie k\"onnte
 es anders sein ---, und dass $\E(H)$ als $K(H)$-Vektorraum zum
 $K$-Vektorraum $H$ isomorph ist. Das hat nach Fr\"uherem dann
 wieder zur Folge, dass $\Delta(P,H)$, falls der Punkt $P$ nicht
 auf der Hyperebene $H$ liegt, zur multiplikativen Gruppe $K^*$
 isomorph ist.
 \par
       Wir beginnen mit der Darstellung von $\E(H)$. Kandidaten zur
 Be\-schrei\-bung von Elationen sind solche linearen Abbildungen von
 $V$ in sich, die eine Hyperebene vektorweise festlassen und ebenso
 den Faktorraum nach dieser Hyperebene. Es wird sich herausstellen,
 dass sich in der Tat jede Elation durch eine solche lineare
 Abbildung darstellen l\"asst.
 \par
       Es sei $V$ ein Vektorraum \"uber dem K\"orper $K$ mit $\Rg_K(V) \geq 2$.
 Ferner sei $H$ eine Hyperebene von $V$. Ist $\tau$ ein Endomorphismus von $V$,
 so nennen wir $\tau$ {\it Transvektion\/}\index{Transvektion}{}
 mit der {\it Achse\/}\index{Achse}{} $H$ von $V$, falls $\tau$ auf $H$ und $V/H$
 die Identit\"at induziert. Mit $\T(H)$ bezeichnen wir die Menge
 aller Transvektionen mit der Achse $H$. Offensichtlich ist $\T(H)$
 eine multiplikativ abgeschlossene Teilmenge des
 Endomorphismenringes von $V$, die die Identit\"at enth\"alt. Es
 sei $\tau \in \T(H)$. Dann gilt zun\"achst
 $$ v^\tau = v + v^\tau - v $$
 f\"ur alle $v \in V$. Weil $\tau$ auf $V/H$ die
 Identit\"at induziert, folgt $v^\tau - v \in H$. Wir definieren $\tau'$ durch
 $$ v^{\tau'} := v - v^\tau + v. $$
 Es folgt $\tau' \in \T(H)$ und $\tau \tau' = 1 = \tau' \tau$, wie eine
 einfache Rechnung best\"atigt. Dies besagt, dass $\T(H)$ sogar eine
 Untergruppe der Einheitengruppe des Endomorphismenringes von $V$ ist.
 \par
       Ist $\tau \in \T(H)$, so induziert $\tau$ in $L(V)$ nat\"urlich
 eine Perspektivit\"at mit der Achse $H$. L\"asst $\tau$ einen
 Punkt au\ss erhalb von $H$ invariant, so induziert $\tau$ auf
 diesem Punkt die Identit\"at, ist also selbst die Identit\"at.
 Folglich induziert jede Transvektion mit der Achse $H$ eine
 Elation mit der Achse $H$ in $L(V)$. Mehr noch: Die Gruppe $\T(H)$
 ist isomorph zu einer Untergruppe von $\E(H)$,  da ja nur die
 Identit\"at aus $\T(H)$ Fixpunkte au\ss erhalb $H$ hat. Dies
 impliziert, dass $\T(H)$ abelsch ist. Wir werden gleich sehen, dass
 jede Elation durch eine Transvektion induziert wird.
 \par
       Ginge man nun didaktisch geschickt vor, um das Folgende zu motivieren, so
 k\"ame man bald in technische Schwierigkeiten, so dass man nach der
 Motivation\index{Motivation}{} noch einmal von vorne beginnen m\"usste.
 Dies zeigten mir jedenfalls meine Versuche. Da ich aber
 wei\ss, wie es wei\-ter\-geht, spare ich mir die Motivation und dem
 Leser das Lesen meiner Schmierzettel.\index{Schmierzettel}{}
 \medskip\noindent
 {\bf 5.1. Satz.} {\it Es sei $V$ ein Rechtsvektorraum \"uber dem K\"orper $K$
 mit $\Rg_K(V) \geq 2$. Ferner sei $P$ ein Punkt und $H$ eine Hyperebene von
 $L(V)$ mit $V = P \oplus H$. Schlie\ss lich sei $P = pK$. Wir definieren die
 Abbildung $\varphi$ von $V$ in $K$ durch
 $$ \varphi (pk + h) := k $$
 f\"ur alle $k \in K$ und alle $h \in H$. Dann ist $\varphi \in V^*$ und
 $H = \Kern(\varphi)$. F\"ur $h \in H$ definieren wir die Abbildung $\tau(h)$
 durch
 $$ v^{\tau (h)} := v + h \varphi (v) $$
 f\"ur alle $v \in V$. Dann ist $\tau$ ein Isomorphismus der Gruppe $H$
 auf die Gruppe $\T(H)$ aller Transvektionen mit der Achse $H$.
 \par
        Ist $\Rg_K(V) \geq 3$, so wird $\E(H)$ von $\tau(H)$ treu induziert. Ist
 $0 \neq h \in H$, so ist $hK$ das Zentrum der von $\tau (h)$ induzierten
 Elation.}
 \smallskip
       Beweis. Es ist banal nachzurechnen, dass $\varphi \in V^*$ und
 $\Kern (\varphi) = H$ ist. Weil $\varphi$ linear ist, ist auch
 $\tau(h)$ linear, so dass $\tau(h)$ ein Endomorphismus von $V$
 ist. Da $\tau (h)$ offensichtlich auf $H$ und $V/H$ die
 Identit\"at induziert, gilt $\tau(h) \in \T(H)$. Weiter gilt
 $$\eqalign{
    v^{\tau(h)\tau(h')} &= \bigl(v + h \varphi(v)\bigr)^{\tau(h')} \cr
                        &= v^{\tau(h')} + h^{\tau(h')}\varphi(v)   \cr
	                &= v + h'\varphi(v) + h\varphi(v)          \cr
			&= v + (h + h') \varphi (v)                \cr
			&= v^{\tau(h + h')}, \cr} $$
 so dass $\tau(h + h') = \tau (h)\tau(h')$ ist. Schlie\ss lich
 folgen aus $\tau(h) = 1$ die Gleichungen
 $$ p = p^{\tau (h)} = p + h $$
 und damit $h = 0$. Dies zeigt, dass $\tau$ ein Monomorphismus von $H$ in
 $\T(H)$ ist.
 \par
       Wir m\"ussen noch zeigen, dass $\tau$ sogar surjektiv ist. Dazu
 sei $Q$ ein von $P$ verschiedenes Komplement von $H$ und $C := (P + Q) \cap H$.
 Dann ist $C$ ein Punkt auf $H$, so dass $P + Q = Q + C$ ist. Es gibt also ein
 $q \in Q$ und ein $c \in C$, so dass $p = q - c$ ist. Es folgt
 $$ p^{\tau (c)} = p + c \varphi (p) = p + c = q. $$
 Weil $p$ nicht Null ist, ist auch $q$ nicht Null.
 Folglich ist $P^{\tau (c)} = Q$. Somit ist $\tau(H)$ --- dies ist
 kein Druckfehler --- auf der Menge der Punkte, die nicht auf $H$
 liegen, transitiv. Weil die einzige Abbildung in $\T(H)$, die
 einen Punkt au\ss erhalb $H$ festl\"asst, die Identit\"at ist,
 folgt hieraus, dass $\tau(H) = \T(H)$ ist.
 \par
       Es sei schlie\ss lich $0 \neq h \in H$ und $X$ sei ein Teilraum
 von $V$, der $hK$ enth\"alt. Ist dann $y \in X$, so ist
 $$ y^{\tau (h)} = y + h \varphi (y) \in X, $$
 so dass $X^{\tau (h)} = X$ ist. Damit ist alles bewiesen.
 \medskip
       Man beachte, dass der Isomorphismus $\tau$ von der Wahl von $p$
 ab\-h\"angt.
 \par
       Es sei wieder $V$ ein Vektorraum \"uber dem K\"orper $K$. Ferner
 sei $P$ ein Punkt und $H$ eine Hyperebene von $V$ mit $V = P \oplus H$. Ist
 $\sigma$ ein Endomorphismus von $V$, der auf $H$
 die Identit\"at induziert, und der weiterhin $P^\sigma = P$
 erf\"ullt, so hei\ss t $\sigma$ {\it Homothetie\/}\index{Homothetie}{} mit dem
 {\it Zentrum\/}\index{Zentrum}{} $P$ und der {\it Achse\/}\index{Achse}{}
 $H$ von $V$. Die Homothetien mit dem Zentrum $P$ und der Achse $H$ bilden eine
 Untergruppe der Einheitengruppe des Endomorphismenringes von $V$,
 die wir mit $\Sigma (P,H)$ bezeichnen.
 \par
       Es sei $\lambda \in \Sigma (P,H)$. Ist $0 \neq p \in P$, so gibt
 es genau ein $a \in K^*$ mit $p^\lambda = pa$. Hieraus folgt
 $$ (pk + h)^\lambda = p^\lambda k + h^\lambda = pak + h $$
 f\"ur alle $k \in K$ und alle $h \in H$. Homothetien lassen sich
 also sehr einfach beschreiben. Diese Beschreibung wird nun im
 n\"achsten Satz zur Definition benutzt. Dabei sei f\"ur den Leser,
 der eine der weniger guten Vorlesungen\index{Vorlesung}{} \"uber lineare Algebra
 geh\"ort hat --- alle be\-trach\-te\-ten K\"orper seien kommutativ ---,
 hier ausdr\"ucklich bemerkt, dass der Koeffizient $a$ zwischen $p$
 und $k$ stehen muss, da nur so die Li\-ne\-a\-ri\-t\"at der gleich zu
 definierenden Abbildung $\delta(a)$ erzwungen wird.
 \medskip\noindent
 {\bf 5.2. Satz.}  {\it Es sei $V$ ein Vektorraum \"uber dem
 K\"orper $K$ und es gelte $\Rg_K(V) \geq 2$. Ferner seien $P$ ein
 Punkt und $H$ eine Hyperebene von $L(V)$ mit $V = P \oplus H$.
 Schlie\ss lich sei $0 \neq p \in P$. F\"ur $a \in K^*$ definieren
 wir die Abbildung $\delta (a)$ von $V$ in sich durch
 $$ (pk + h)^{\delta(a)} := pa^{-1} k + h $$
 f\"ur alle $k \in K$ und alle $h \in H$. Dann ist $\delta$ ein Isomorphismus
 von $K^*$ auf $\Delta(P,H)$.
 \par
       Ist $\Rg_K (V) \geq 3$, so wird $\Delta(P,H)$ von $\Sigma (P,H)$ treu
 induziert.}
 \smallskip
       Beweis. Eine banale Rechnung zeigt, dass $\delta(a)$ linear ist.
 Hieraus folgt dann, dass $\delta(a) \in \Sigma (P,H)$ gilt.
 \par
       Es sei $0 \neq a$, $b \in K^*$. Dann ist
 $$\eqalign{
   (pk + h)^{\delta (ab)} &= p(ab)^{-1} k + h                         \cr
			  &= pb^{-1} a^{-1} k + h                     \cr
			  &= p^{\delta (b)} a^{-1} k + h^{\delta (b)} \cr
			  &= (pa^{-1} k + h)^{\delta (b)}             \cr
			  &= (pk + h)^{\delta (a) \delta (b)}, \cr} $$
 so dass $\delta(ab) = \delta(a)\delta(b)$ ist. Folglich ist $\delta$ ein
 Homomorphismus.
 \par
       Um die Injektivit\"at zu beweisen, beweisen wir etwas mehr,
 n\"amlich, dass $\delta(a)$ genau dann einen von $P$
 verschiedenen Punkt $Q$ von $L(V)_H$ festl\"asst, wenn $a = 1$
 ist. Es sei also $Q$ ein solcher Punkt und es gelte $Q^{\delta(a)} = Q$. Weil
 $P + Q$ eine Gerade ist, ist $C := (P + Q) \cap H$
 ein Punkt, der \"uberdies von $P$ und $Q$ verschieden ist. Es gibt
 daher von $0$ verschiedene Vektoren $p'$, $q'$ und $c'$ mit $p' \in P$,
 $q' \in Q$ und $c' \in C$ und $q' = p' + c'$. Es gibt ferner
 ein von $0$ verschiedenes $k \in K$ mit $p = p'k$. Setze $q :=
 q'k$ und $c:= c'k$. Dann sind auch $p$, $q$ und $c$ von $0$
 verschieden und es gilt $q = p + c$. Wegen $Q^{\delta(a)} = Q$
 gibt es ein $l \in K$ mit $q^{\delta(a)} = ql$. Es folgt
 $$ pl + cl = q^{\delta(a)} = (p + c)^{\delta(a)} = pa^{-1} + c. $$
 Weil $p$ und $c$ linear unabh\"angig sind, folgt hieraus $a^{-1} = l = 1$,
 so dass $a = 1$ ist, und zum Andern, dass $\Sigma (P,H)$ eine Untergruppe von
 $\Delta(P,H)$ treu induziert, falls nur $\Rg_K(V) \geq 3$ ist, da nur in
 diesem Falle $\Delta(P,H)$ definiert ist.
 \par
       Ist $\lambda \in \Sigma (P,H)$, so haben wir schon gesehen, dass
 es ein $a \in K^*$ gibt mit $\delta (a^{-1}) = \lambda$. Somit ist
 $\delta$ auch surjektiv.
 \par
       Es sei $\Rg_K (V) \geq 3$. In diesem Falle m\"ussen wir noch
 zeigen, dass $\Sigma (P,H)$ alle Kollineationen in $\Delta(P,H)$ induziert. Zu
 diesem Zweck seien $P_1$ und $P_2$ zwei von $P$
 verschiedene Komplemente von $H$ mit $P + P_1 = P + P_2$. Setze $C
 := (P + P_1) \cap H$. Dann ist $C$ ein Punkt auf $H$. Es gibt also
 wieder von $0$ verschiedene $p_i \in P_i$ und $c_1$, $c_2 \in C$
 mit $p_i = p + c_i$ f\"ur $i:= 1$, 2. Es gibt ferner ein $a \in
 K^*$ mit $c_2 = c_1 a$. Damit folgt nun
 $$ p_1^{\delta(a)} = pa^{-1} + c_1 = (p + c_2) a^{-1} = p_2 a^{-1}. $$
 Also ist
 $P_1^{\delta (a)} = P_2$. Dies beweist nach nun schon sattsam
 bekannten S\"atzen auch die noch offene letzte Behauptung des
 Satzes.
 \smallskip
 Auf Grund der beiden S\"atze 5.1 und 5.2 d\"urfen wir im folgenden
 die Gruppe $\T(H)$ mit $\E (H)$ und die Gruppe $\Sigma (P,H)$
 mit $\Delta(P,H)$ identifizieren.
 \medskip\noindent
 {\bf 5.3. Satz.} {\it Es sei $V$ ein Vektorraum \"uber dem
 K\"orper $K$ und es gelte $\Rg_K(V) \geq 3$. Ferner sei $P$ ein
 Punkt und $H$ eine Hyperebene von $L(V)$ mit $V = P \oplus H$.
 Wir definieren eine Abbildung $\delta^*$ von $K$ in den
 Endomorphismenring von $\E(H)$ durch $\delta^*(0) := 0$ und
 $$ \rho^{\delta^* (a)} := \delta (a)^{-1} \rho \delta (a) $$
 f\"ur alle $\rho \in \E(H)$ und alle $a \in K^*$. Dann ist
 $\tau(h)^{\delta^*(a)} = \tau(ha)$ f\"ur alle $h \in H$ und alle $a \in K$.
 Insbesondere ist $(\tau,\delta^*)$ ein Isomorphismus des
 $K$-Vektorraumes $H$ auf den $K(H)$-Vektorraum $\E(H)$.}
 \smallskip
       Bevor wir mit dem Beweise des Satzes beginnen, scheint noch ein
 Wort des Kommentars angebracht. Auch gute Vorlesungen\index{Vorlesung}{}
 \"uber lineare Algebra kommen meist nicht auf den Begriff der semilinearen
 Abbildung\index{semilineare Abbildung}{} zu sprechen. In der Geometrie ist man
 aber durch die Sache gezwungen, auch von diesen zu reden, wie sich hier
 zum ersten Male zeigt. Es ist also an der Zeit den Begriff der
 semilinearen Abbildung zu definieren. Dazu sei $V$ ein Vektorraum
 \"uber $K$ und $V'$ ein solcher \"uber $K'$. Ist dann $\sigma$ ein
 Homomorphismus der abelschen Gruppe $V$ in die abelsche Gruppe
 $V'$, ist $\alpha$ ein Isomorphismus des K\"orpers $K$ auf den
 K\"orper $K'$ und gilt
 $$ (vk)^\sigma = v^\sigma k^\sigma $$
 f\"ur alle $v \in V$ und alle $k \in K$, so hei\ss t das Paar
 $(\sigma, \alpha)$ {\it semilineare Abbildung\/} von $V$ in $V'$.
 Es bedarf nun nicht mehr viel an Fantasie, um den Begriff des
 Isomorphismus von Vektorr\"aumen zu definieren.
 \smallskip
       Beweis. Mit 2.1 folgt, dass $\delta^*$ eine surjektive Abbildung
 von $K$ auf $K(H)$ ist. Ferner ist klar, dass
 $\tau(h)^{\delta^*(0)} = \tau(h 0)$ ist. Es sei also $a \in K^*$. Zu $h \in H$
 gibt es dann ein $h' \in H$ mit $\tau(h)^{\delta^*(a)} = \tau(h')$.
 Um $\tau (h')$ zu bestimmen, gen\"ugt es, die Wirkung von $\tau(h')$ auf den
 Punkt $P$ zu bestimmen. Wegen $\varphi (p) = 1$ ist nun
 $$\eqalign{
       p + h' &= p^{\tau(h')} = p^{\delta(a)^{-1}\tau(h) \delta(a)} \cr
	      &= \bigl(pa + h\varphi(pa)\bigr)^{\delta(a)} = p + ha. \cr} $$
 Somit ist $h' = ha$, so dass $\tau(h)^{\delta^*(a)} = \tau(ha)$ ist.
 Routinerechnungen zeigen nun, dass $\delta^*$ ein Homomorphismus von $K$ auf
 $K(H)$ ist. Weil nicht triviale Homomorphismen von K\"orpern stets
 Monomorphismen sind, ist bereits alles bewiesen.

\mysection{6. Der erste Struktursatz}

\noindent
 In diesem Abschnitt werden wir den schon
 angek\"undigten Satz beweisen, dass sich jeder desarguessche projektive Verband
 als Unterraumverband eines Vektorraumes darstellen l\"asst. Eine
 unmittelbare Folgerung aus ihm wird der von Hilbert\index{Hilbert, D.}{}
 stammende Satz sein, dass eine desarguessche projektive Geometrie genau dann
 pappossch ist, wenn der zu Grunde liegende Koordinatenk\"orper
 kommutativ ist (Hilbert 1899). Damit werden wir zusammen mit dem hessenbergschen
 Satz eine befriedigende Beschreibung aller papposschen Geometrien erhalten.
 \medskip\noindent
 {\bf 6.1. Erster Struktursatz.} {\it Ist $L$ ein irreduzibler projektiver
 Verband vom Rang 
 $\Rg (L) \geq 3$ und ist $L$ in Falle $\Rg(L) = 3$ desarguessch, so
 gibt es einen und bis auf Isomorphie auch nur einen Vektorraum $V$ \"uber einem
 K\"orper $K$, so dass $L$ und $L(V)$ isomorph sind.}
 \smallskip
       Beweis. Wir beweisen zuerst die Eindeutigkeitsaussage. Es seien $L$ und
 $L'$ zwei desarguessche projektive Verb\"ande und $\sigma$ sei ein Isomorphismus
 von $L$ auf $L'$. Ferner sei $H$ eine Hyperebene von $L$ und $H' := H^\sigma$.
 Ist $\tau \in \E(H)$ und $\kappa \in K(H)$, so setzen wir
 $\tau^{\sigma^*} := \sigma^{-1} \tau \sigma$ und
 $\kappa^{\sigma^{**}} := (\sigma^*)^{-1} \kappa \sigma^*$. Eine simple Rechnung
 zeigt, dass das so definierte
 Abbildungspaar ein Isomorphismus des $K (H)$-Vektorraumes $\E(H)$ auf den
 $K(H')$-Vektorraum $\E(H')$ ist. Ist nun $L = L(V)$ und $L' = L(V')$, so folgt
 aus dieser Bemerkung, aus $\Rg_K(V) = \Rg (L) = \Rg (L') = \Rg_{K'} (V')$ und
 aus 5.3 sowie aus dem Struktursatz f\"ur
 Vektorr\"aume,\index{Struktursatz f\"ur Vektorr\"aume}{} dass n\"amlich ein Vektorraum
 \"uber dem K\"orper $K$ durch seinen Rang und eben diesen K\"orper
 bis auf Isomorphie eindeutig bestimmt ist, dass $V$ und $V'$
 isomorph sind.
 \par
       Den Beweis der Existenzaussage beginnen wir mit einer
 Vorbemerkung. Es sei $L$ ein desarguesscher projektiver Verband
 und $H$ sei eine Hyperebene von $L$. Die Gruppe $\E(H)$ operiert
 dann scharf trans\-i\-tiv auf der Menge der Punkte von $L_H$. Ist $O$
 ein Punkt von $L_H$, so gibt es zu jedem Punkt $P$ von $L_H$ genau
 ein $\tau_P \in \E (H)$ mit
 $$ O^{\tau_P} = P. $$
 Die Abbildung $\tau$ ist eine Bijektion der Punktmenge von $L_H$ auf
 $\E(H)$. Ist nun $X \leq H$, so gilt offensichtlich $P \leq O + X$ genau dann,
 wenn $\tau_P \in \E(X)$ ist. Bezeichnet $\Pi$
 wieder das gr\"o\ss te Element von $L$, so folgt mittels Satz 2.4,
 dass $\tau$ einen Isomorphismus von $\Pi/O$ auf den Verband der
 Unterr\"aume des $K(H)$-Vektorraumes $\E(H)$ induziert.
 Hieraus folgt schlie\ss lich auf Grund der Transitivit\"at von $E(H)$, dass
 $\tau$ einen Isomorphismus von $L_H$ auf den Verband
 der Rechtsrestklassen nach allen Unterr\"aumen von $\E(H)$
 induziert. (Hieran k\"onnte man nun didaktische Bemerkungen \"uber
 freie Vektoren,\index{freier Vektor}{} den Elementen von $\E(H)$, und
 Ortsvektoren,\index{Ortsvektor}{} der
 Beschreibung der Punkt\-men\-ge von $L_H$ mittels $O$ und $\tau$,
 anschlie\ss en. Diese niemals sauber definierten Begriffe sind eine st\"andig
 flie\ss ende Quelle der Konfusion.\index{Konfusion}{} Der Leser wird
 jedoch gen\"ugend Fantasie besitzen, anhand der vorstehenden
 Bemerkungen eine Trennung der beiden Begriffe vornehmen zu
 k\"onnen, so dass sich weitere Bemerkungen meinerseits
 er\"ubrigen. Sie w\"urden den richtigen Adressaten ja doch nicht
 erreichen.)
 \par
       Mit Hilfe der Vorbemerkung ist es nun ein Leichtes, die
 Ex\-is\-tenz\-aus\-sa\-ge des Satzes zu beweisen. Es sei also $L$ ein
 desarguesscher projektiver Verband und $H$ sei eine Hyperebene von
 $L$. Ferner sei $V$ ein Vektorraum \"uber $K(H)$ mit $\Rg_{K(H)}(V) = \Rg(L)$.
 Einen solchen Vektorraum gibt es stets. (Um dies
 einzusehen, betrachte man die Menge aller Abbildungen einer Basis
 von $L$ in $K(H)$ mit endlichem Tr\"ager. Diese Menge versehen mit
 der punktweise definierten Addition und Skalarmultiplikation ist
 ein solcher.) Es sei weiter $H'$ eine Hyperebene von $V$. Nach 5.3
 sind $K(H')$ und $K(H)$ isomorph. Mittels 2.4 und $\Rg_{K(H)}(V) = \Rg (L)$
 folgt
 $$ \Rg_{K(H)} (\E (H)) = \Rg_\Pi (H) = \Rg_{\Pi '} (H')
			= \Rg_{K(H')}\bigl(\E (H')\bigr). $$
 Hieraus folgt weiter, dass der $K(H)$-Vektorraum $\E(H)$ zu dem 
 $K(H')$-Vek\-tor\-raum $\E(H')$ isomorph ist. Unsere Vorbemerkung sagt dann aber,
 dass auch $L_H$ und $L(V)_{H'}$ isomorph sind. Hieraus folgt schlie\ss\-lich
 mittels I.8.8, dass auch $L$ und $L(V)$ isomorph sind. Damit ist alles bewiesen.
 \medskip
       Wie der Beweis zeigt, gilt auch der folgende, von Andr\'e stammende
 Satz (Andr\'e 1954).
 \medskip\noindent
 {\bf 6.2. Satz.} {\it Ist $L$ eine Translationsebene bez\"uglich der Geraden
 $H$, so ist $L$ genau dann desarguessch, wenn der $K (H)$-Vektorraum $\E (H)$
 den Rang $2$ hat.}
 \medskip
       Ist $L$ ein desarguesscher projektiver Verband, so gibt es also
 einen Vektorraum $V$ \"uber einem K\"orper $K$ mit $L \cong L(V)$. Man nennt $V$
 den $L$
 {\it zu Grunde liegenden Vektorraum\/}\index{zu Grunde liegender Vektorraum}{}
 und $K$ den
 {\it Koordinatenk\"orper\/}\index{Koordinatenk\"orper}{} von $L$. Mit 3.1 und
 3.5 folgt daher der schon angek\"undigte Satz von Hilbert (Hilbert
 1899).\index{Hilbert, D.}{}
 \medskip\noindent
 {\bf 6.3. Satz.} {\it Ein desarguesscher projektiver Verband ist
 genau dann pappossch, wenn sein Koordinatenk\"orper kommutativ ist.}
 \medskip
       Da es K\"orper gibt, die nicht kommutativ sind, gibt es auch
 projektive R\"aume, in denen der Satz von Pappos nicht gilt.
 Beispiele solcher K\"orper werden wir sp\"ater noch kennenlernen.
 \par
 Ist $q$ Potenz einer Primzahl, so gibt es einen und bis auf
 Isomorphie auch nur einen K\"orper mit $q$ Elementen, n\"amlich
 das Galoisfeld $\GF(q)$. Daher gilt auch das
 \medskip\noindent
 {\bf 6.4. Korollar.} {\it Ist $q$ Potenz einer Primzahl und ist
 $n$ eine nat\"urliche Zahl mit $n \geq 4$, so gibt es einen und
 bis auf Isomorphie auch nur einen endlichen irreduziblen
 projektiven Verband der Ordnung $q$ und des Ranges $n$.}
 \medskip
 Dieser Satz ist f\"ur $n = 3$ falsch, da es nicht desarguessche
 projektive Ebenen gibt.

\mysection{7. Der zweite Struktursatz}

\noindent
 Sind $V$ und $V'$ Vektorr\"aume \"uber $K$ und
 $K'$, so induziert jeder Isomorphismus von $V$ auf $V'$ einen Isomorphismus von
 $L(V)$ auf $L(V')$. Die Umkehrung dieses Sachverhaltes ist Inhalt des zweiten
 Struktursatzes, den wir jetzt formulieren und beweisen werden. Dabei sei daran
 erinnert, dass Isomorphismen von Vektorr\"aumen auch se\-mi\-li\-ne\-ar sein
 k\"onnen.\index{Struktursatz}{}
 \medskip\noindent
 {\bf 7.1. Zweiter Struktursatz.} {\it Sind $V$ und $V'$ zwei
 Vektorr\"aume \"uber $K$ bzw. $K'$, ist $\Rg_K(V) \geq 3$ und ist
 $\sigma$ ein Isomorphismus von $L(V)$ auf $L(V')$, so wird
 $\sigma$ durch einen Isomorphismus von $V$ auf $V'$ induziert.}
 \smallskip
       Beweis. Aus dem ersten Struktursatz folgt, dass $V$ und $V'$
 isomorph sind, so dass wir annehmen d\"urfen, dass $V = V'$ und $K = K'$ ist.
 Satz 5.1 lehrt, dass alle Elationen von $L(V)$ sogar
 durch li\-ne\-a\-re Automorphismen von $V$ induziert werden, und die
 Gruppe von Kollineationen, die von allen Elationen erzeugt wird,
 ist nach 3.2 auf der Menge der nicht inzidenten
 Punkt-Hyperebenenpaare von $L(V)$ transitiv. Wir d\"urfen daher
 des Weiteren annehmen, dass es ein nicht inzidentes
 Punkt-Hyperebenenpaar $(P,H)$ gibt mit $P^\sigma = P$ und
 $H^\sigma = H$. Dann ist
 $$ \sigma^{-1} \E(H) \sigma = \E(H^\sigma) = \E(H), $$
 so dass der durch $\sigma$ induzierte
 Automorphismus der Kollineationsgruppe von $L(V)$ einen
 Automorphismus auf $\E(H)$ induziert. Wir definieren die
 Abbildung $\sigma^*$ durch
 $$ \tau (h^{\sigma^*}) := \sigma^{-1}\tau(h)\sigma. $$
 Mit 5.1 folgt, dass $\sigma^*$ ein Automorphismus der abelschen Gruppe $H$ ist.
 \par
       Es sei $a \in K$. Ist $a = 0$, so ist $\tau ((ha)^{\sigma^*}) = 1
 = \tau (h^{\sigma^*} a)$. Es sei also $a \neq 0$. Nach 5.3 ist dann
 $$\eqalign{
    \tau\bigl((ha)^{\sigma^*}\bigr) &= \sigma^{-1}\tau(ha)\sigma \cr
    &= \sigma^{-1}\delta(a)^{-1}\sigma\sigma^{-1}\tau(h) \sigma
                \sigma^{-1}\delta(a)\sigma          \cr
    &= \sigma^{-1}\delta(a)^{-1}\sigma\tau(h^{\sigma^*})\sigma^{-1}\delta(a)
		\sigma. \cr}$$
 \par
 Nun ist
    $$ \sigma^{-1}\Delta(P,H)\sigma = \Delta(P^\sigma,H^\sigma) = \Delta(P,H). $$
 Also induziert $\sigma$ einen Automorphismus in $\Delta(P, H)$. Es gibt daher zu
 jedem $a \in K^*$ ein $a^{\sigma^{**}} \in K^*$ mit
 $$ \sigma^{-1}\delta(a)\sigma = \delta(a^{\sigma^{**}}). $$
 Setzt man noch $0^{\sigma^{**}} := 0$, so ist $\sigma^{**}$ eine
 Bijektion von $K$ auf sich und es gilt
 $$ \tau \bigl((ha)^{\sigma^*}\bigr) = \tau (h^{\sigma^*} a^{\sigma^{**}}) $$
 f\"ur alle $h \in H$ und alle $a \in K$. Hieraus folgt wiederum, da
 $\tau$ injektiv ist, dass
 $$ (ha)^{\sigma^*} = h^{\sigma^*} a^{\sigma^{**}} $$
 f\"ur alle $h \in H$ und alle $k \in K$ gilt.
 Routinerechnungen zeigen nun, dass $(\sigma^*,\sigma^{**})$ ein
 Automorphismus des $K(H)$-Vektorraumes $H$ ist.
 \par
       Es sei $0 \neq p \in P$. Dann ist $V = pK \oplus H$. Wir
 definieren die Abbildung $\rho^*$ von $V$ in sich durch
 $$ (pa + h)^{\rho^*} := pa^{\sigma^{**}} + h^{\sigma^*}. $$
 Dann ist $(\rho^*, \sigma^{**})$ ein Automorphismus von $V$. Nun ist
 $p^{\sigma^*} = p$ und daher
 $$\eqalign{
      p^{\rho^{* -1}\tau (h) \rho^*} &= p^{\tau (h) \rho^*}      \cr
	     &= (p + h)^{\rho^*}                                 \cr
	     &= p^{\tau (h^{\sigma^*})}.   \cr}$$
 Also ist $\rho^{*-1} \tau (h) \rho^* = \tau (h^{\sigma^*})$. Ist $\rho$ die von
 $\rho^*$ in $L(V)$ induzierte Kollineation, so ist also
 $\rho^{-1} \tau (h) \rho = \sigma^{-1} \tau (h) \sigma$. Daher ist
 $\sigma \tau^{-1}$ f\"ur alle $h \in H$ mit $\tau (h)$ vertauschbar. Nun ist
 $P^{\sigma \rho^{-1}} = P$ und folglich
 $$ P^{\tau(h)} = P^{\sigma \rho^{-1}\tau(h)}
                = (P^{\tau (h)})^{\sigma \rho^{-1}} $$
 f\"ur alle $h \in H$.  Weil $\E (H)$ auf der Menge der Punkte von $L(V)_H$
 transitiv operiert, folgt, dass $\sigma \rho^{-1}$ alle Punkte von
 $L(V_H)$ zu Fixpunkten hat, so dass nach I.8.8 die Gleichung
 $\sigma \rho^{-1} = 1$ gilt. Folglich ist $\sigma = \rho$, was wir
 zu beweisen hofften.

\mysection{8. Der dritte Struktursatz}

\noindent
 Der erste in diesem Abschnitt zu beweisende
 Satz verdient eigentlich auch den Namen Struktursatz\index{Struktursatz}{}. Er
 lautet:
 \medskip\noindent
 {\bf 8.1. Satz.} {\it Ist $L$ ein desarguesscher projektiver
 Verband, so ist $L$ genau dann selbstdual, wenn der Rang von $L$
 endlich ist und der $L$ zu Grunde liegende Koordinatenk\"orper
 einen Antiautomorphismus besitzt.}
 \smallskip
       Beweis. Es sei $K$ der $L$ zu Grunde liegende Koordinatenk\"orper.
 \par
       Ist $L$ selbstdual, so folgt mit I.5.7, dass der Rang von $L$
 endlich ist. Die nach I.5.19 gemachte Bemerkung zeigt, dass
 $K_\circ$ der Koordinatenk\"orper von $L^d$ ist. Weil $L$ zu $L^d$
 isomorph ist, folgt aus dem ersten Struktursatz, dass es einen
 Isomorphismus $\alpha$ von $K$ auf $K_\circ$ gibt. Dann ist
 $\alpha$ aber nichts Anderes als ein Antiautomorphismus von $K$.
 \par
       Es sei umgekehrt der Rang von $L$ endlich und $\alpha$ sei ein
 Antiautomorphismus von $K$. Dann ist $\alpha$ ein Isomorphismus
 von $K$ auf $K_\circ$. Ist $V$ der $L$ zu Grunde liegende
 Vektorraum, so ist der $K_\circ$-Vektorraum $V^*$ der $L^d$
 zu Grunde liegende Vektorraum. Da nach I.5.17 die R\"ange von $V$
 und $V^*$ gleich sind, ist der $K$-Vektorraum $V$ zum
 $K_\circ$-Vektorraum $V^*$ isomorph. Folglich sind auch $L$ und
 $L^d$ isomorph, so dass $L$ selbstdual ist.
 \medskip
       Da bei einem kommutativen K\"orper die Identit\"at stets auch ein
 Antiautomorphismus ist, gilt, wie schon in Abschnitt 5 des ersten
 Kapitels, wenn auch in etwas anderer Formulierung, bemerkt, das
 \medskip\noindent
 {\bf 8.2. Korollar.} {\it Ein papposscher projektiver Raum
 endlichen Ranges ist stets selbstdual. Insbesondere sind alle
 endlichen desarguesschen projektiven R\"aume selbstdual.}
 \medskip
 Es gibt K\"orper, die keinen Antiautomorphismus gestatten. Ist
 n\"am\-lich $K$ ein K\"orper des Ranges $n < \infty$ \"uber seinem
 Zentrum, so zeigt die Theorie der einfachen Algebren, dass $K$
 h\"ochstens dann einen Antiautomorphismus besitzt, wenn $n$ eine
 Potenz von 2 ist (siehe etwa Deuring 1968, Satz 11, S.~45 und Satz 2,
 S.~59).
 \medskip
 Es sei $f$ eine Abbildung von $V \times V$ in $K$ und $\alpha$ sei
 ein Antiautomorphismus von $K$. Die Abbildung $f$ hei\ss t
 $\alpha$-\emph{Semibilinearform},\index{Semibilinearform}{} falls gilt:
 \smallskip\noindent
 a) Es ist $f(u + v, w) = f(u,w) + f(v,w)$ und $f(u,v + w) = f(u,v) + f(u,w)$
 f\"ur alle $u$, $v$, $w \in V$.
 \smallskip\noindent
 b) Es ist $f(uk,v) = k^\alpha f(u,v)$ und $f(u, vk) = f(u,v)k$ f\"ur alle
 $u$, $v \in V$ und alle $k \in K$.
 \smallskip
       Ist $\kappa$ eine Korrelation von $L(V)$, so sagen wir, dass
 $\kappa$ durch die $\alpha$-Semi\-bi\-li\-ne\-ar\-form $f$\index{Semibilinearform}{}
 {\it dargestellt\/} werde, falls
 $$ U^\kappa = \bigl\{v \mid v \in V, f(u,v) = 0 \ \hbox{\rm f\"ur alle}\
		  u \in U\bigr\} $$
 f\"ur alle $U \in L(V)$ ist. Es gilt nun:
 \medskip\noindent
 {\bf 8.3. Dritter Struktursatz.} {\it Es sei $V$ ein Vektorraum
 \"uber dem K\"orper $K$ mit $\Rg_K (V) \geq 3$. Ist $\kappa$ eine
 Korrelation von $L(V)$, so wird $\kappa$ durch eine
 $\alpha$-Semibilinearform dargestellt.}
 \smallskip
       Beweis. Es sei $W$ der Dualraum zu $V$. Beide R\"aume haben dann
 nach 8.1 den gleichen, endlichen Rang. F\"ur $X \in L(V)$
 setzen wir
 $$ X^\pi := \big\{w \mid w \in W,\ wX = \{0\}\bigr\}. $$
 Dies ist die Abbildung, die wir im f\"unften Abschnitt des ersten
 Kapitels mit $\perp$ bezeichneten. Nach I.5.19 ist $\pi$ ein
 Isomorphismus von $L(V)^d$ auf $L(W)$, so dass $\kappa \pi$
 ein Isomorphismus von $L(V)$ auf $L(W)$ ist. Nach dem
 zweiten Struktursatz gibt es also einen Isomorphismus $(\rho,\alpha)$ des
 $K$-Vektorraumes $V$ auf den $K_\circ$-Vektorraum $W$
 mit $X^\rho = X^{\kappa \pi}$ f\"ur alle $X \in L(V)$. Weil
 $\alpha$ ein Isomorphismus von $K$ auf $K_\circ$ ist, ist $\alpha$
 ein Antiautomorphismus von $K$. Ferner gilt
 $$ (vk)^\rho = v^\rho \circ k^\alpha = k^\alpha v^\rho $$
 f\"ur alle $v \in V$ und alle $k \in K$. Wir definieren nun $f$ durch
 $$ f(u, v) := u^\rho v $$
 f\"ur alle $u$, $v \in V$. Triviale Rechnungen zeigen, dass $f$ eine
 $\alpha$-Se\-mi\-bi\-li\-ne\-ar\-form ist.
 \par
       Es sei $U \in L(V)$. Dann ist $U^{\kappa \pi} = U^\rho$. Ist
 $v \in U^\kappa$ und $u \in U$, so ist $u^\rho \in U^{\kappa \pi}$ und daher
 $f(u, v) = u^\rho v = 0$. Ist umgekehrt $u^\rho v = f(u,v) = 0$ f\"ur alle
 $u \in U$, so ist $v \in U^{\rho \pi^{-1}} = U^\kappa$. Also ist
 $$ U^\kappa = \bigl\{v \mid v \in V, f(u, v) = 0
     \ \hbox{\rm f\"ur alle}\ u \in U\bigr\}, $$
 was zu beweisen war.
 \medskip
       Nicht alle Semibilinearformen stellen Korrelationen dar und
 verschiedene Semibilinearformen k\"onnen durchaus auch ein und
 dieselbe Korrelation darstellen.
 \par
       Eine $\alpha$-Semibilinearform $f$ hei\ss e {\it nicht ausgeartet\/}
 oder {\it nicht entartet\/},
 falls\index{nicht ausgeartet}\index{nicht entartet} aus der G\"ultigkeit von $f(u,v) = 0$
 f\"ur alle $u \in V$ folgt, dass $v = 0$ ist.
 \medskip\noindent
 {\bf 8.4. Satz.} {\it Es sei $V$ ein $K$-Vektorraum mit $3 \leq
 \Rg_K (V) < \infty$. Die $\alpha$-Semi\-bi\-li\-ne\-ar\-form $f$ auf $V$ induziert genau
 dann eine Korrelation auf $L(V)$, wenn $f$ nicht ausgeartet ist.}
 \smallskip
       Beweis. Die Korrelation $\kappa$ werde von $f$ induziert. Aus
 $V^\kappa = \{0\}$ folgt, dass $f$ nicht ausgeartet ist.
 \par
       Es sei nun $f$ nicht ausgeartet. Mit $W$ bezeichnen wir den
 Dualraum von $V$ aufgefasst als Rechtsvektorraum \"uber
 $K_\circ$. Wir definieren eine Abbildung $\varphi$ von $V$ in $W$
 durch
 $$ v^\varphi x := f(x, v)^{\alpha^{-1}} $$
 f\"ur alle $v$, $x \in V$. Weil $\alpha$ ein Antiautomorphismus ist, ist auch
 $\alpha^{-1}$ ein solcher. Dies hat zur Folge, dass in der Tat
 $v^\varphi \in W$ ist. Man verifiziert dar\"uber hinaus m\"uhelos,
 dass $(\varphi,\alpha^{-1})$ eine semilineare Abbildung von $V$
 in $W$ ist. Weil $f$ nicht ausgeartet ist, ist $\varphi$ injektiv.
 Mittels der Endlichkeit des Ranges von $V$ erschlie\ss en wir
 hieraus, dass $(\varphi,\alpha^{-1})$ sogar ein Isomorphismus von
 $V$ auf $W$ ist.
 \par
       Wir definieren eine Abbildung $\kappa$ von $L(V)$ in sich durch
 $$ U^\kappa := \bigl\{v \mid v \in V, f(u,v) = 0\ \hbox{\rm f\"ur alle}\ 
		  u \in U\bigr\}. $$
 Es sei $X^\kappa = Y^\kappa$. Ist $y \in Y$, jedoch $y \notin X$, so gibt es ein
 $g \in W$ mit $X \subseteq \Kern(g)$ und $gy = 1$. Nach unserer Vorbemerkung
 gibt es ein $v \in V$ mit $g = v^\varphi$. Daher ist
 $$ 0 = gx = v^\varphi x= f(x, v)^{\alpha^{-1}} $$
 f\"ur alle $x \in X$. Hieraus folgt, dass $f(x,v) = 0$ ist f\"ur alle $x \in X$.
 Somit ist $v \in X^\kappa = y^\kappa$ und daher $f(y, v) = 0$. Andererseits ist
 $1 = (gy)^\alpha = f(y,v)$. Dieser Widerspruch zeigt, dass $Y \leq X$ ist.
 Ebenso folgt $X \leq Y$, so dass $\kappa$ als injektiv erkannt ist.
 \par
       Es ist klar, dass $Y^\kappa \leq X^\kappa$ von $X \leq Y$ impliziert wird.
 Es sei nun $Y^\kappa \leq X^\kappa$. Es ist $Y \leq X + Y$ und daher
 $(X + Y)^\kappa \leq Y^\kappa$. Ist $u \in Y^\kappa$, so ist $f(x,u) = 0$ f\"ur
 alle $x \in X$, da ja $Y^\kappa \leq X^\kappa$ ist. Trivialerweise gilt auch
 $f(y,u) = 0$ f\"ur alle $y \in Y$. Daher gilt $f(x+y,u) = 0$ f\"ur alle
 $x \in X$ und alle $y \in Y$. Somit ist $u \in (X + Y)^\kappa$. Damit haben wir
 auch $Y^\kappa \leq (X + Y)^\kappa$. Also ist $Y^\kappa = (X + Y)^\kappa$. Weil
 $\kappa$ injektiv ist, folgt weiter $Y = X + Y$ und damit $X \leq Y$. Folglich
 gilt $X \leq Y$ genau dann, wenn $Y^\kappa \leq X^\kappa$ ist.
 \par 
       Es sei $b_1$, \dots, $b_n$ eine Basis von $V$. Ferner sei
 $u = \sum^{n}_{i:=1} b_i u_i$ und $v = \sum^{n}_{i:=1} b_i v_i$. Dann ist
 $$ f(u,v) = \sum_{i:=1}^n u_i^\alpha \sum_{j:=1}^n f(b_i, b_j) v_j. $$
 Genau dann ist $f(u, v) = 0$ f\"ur alle $u \in V$, wenn
 $$ \sum_{j:=1}^n f(b_i, b_j) v_j =0 $$
 ist f\"ur $i := 1$, \dots, $n$. Da $f$ nicht ausgeartet ist, hat dieses System
 linearer Gleichungen nur die triviale L\"osung. Daher ist der
 Rechts\-spal\-ten\-rang der Matrix $(f(b_i, b_j))$ gleich $n$. Hieraus
 folgt, dass der Linkszeilenrang dieser Matrix auch gleich $n$ ist.
 Definiere $g$ durch
 $$ g(u,v) := f(v,u)^{\alpha^{-1}} $$
 f\"ur alle $u$, $v \in V$. Dann ist $g$ eine $\alpha^{-1}$-Semibilinearform auf
 $V$.  Gilt nun $\sum^n_{j:=1} g(b_i, b_j)v_j = 0$ f\"ur alle $i := 1$, \dots,
 $n$, so ist
 $$ 0 = \sum_{j:=1}^n f(b_j, b_i)^{\alpha^{-1}} v_j
      = \biggl(\sum_{j:=1}^n v_j^\alpha f(b_i,b_j)\biggr)^{\alpha^{-1}} $$
 f\"ur alle $i$. Weil der Linkszeilenrang der Matrix $(f(b_j,b_i))$ gleich $n$
 ist, folgt $v_j = 0$ f\"ur alle $i$. Also ist
 auch $g$ nicht ausgeartet. Definiere $\lambda$ durch
 $$ X^\lambda := i\bigl\{ v \mid v \in V, g(u, v) = 0\ \hbox{\rm f\"ur alle}\ 
		     u \in X\bigr\}.$$
 Dann ist $\lambda$ aus dem gleichen Grunde wie $\kappa$ eine die Inklusion
 um\-keh\-ren\-de, injektive Abbildung von $L(V)$ in sich. Ferner gilt
 $X \leq X^{\kappa \lambda}$ und $X \leq X^{\lambda \kappa}$ f\"ur alle
 $X \in L(V)$. Weil $X$, $X^{\kappa \lambda}$ und $X^{\lambda \kappa}$ den
 gleichen Rang haben, dies folgt aus der Injektivit\"at von $\kappa \lambda$ und
 $\lambda \kappa$, und da dieser Rang endlich ist, folgt $X = X^{\kappa \lambda}
 = X^{\lambda \kappa}$. Folglich ist $\kappa$ bijektiv und
 $\lambda = \kappa^{-1}$.
 \medskip
       Der gerade gef\"uhrte Beweis liefert mehr als im Satz formuliert.
 Was wir mehr bewiesen haben, formulieren wir in den n\"achsten
 beiden Korollaren. --- Der Leser beachte, dass in 8.5 die Rollen der beiden
 Argumente von $f$ gegen\"uber der Definition des Nicht-ausgeartet-seins
 vertauscht sind!
 \medskip\noindent
 {\bf 8.5. Korollar.} {\it Es sei $V$ ein $K$-Vektorraum mit $3 \leq \Rg_K (V) <
 \infty$ und $f$ sei eine $\alpha$-Semibilinearform auf $V$. Genau dann ist $f$
 nicht ausgeartet, wenn aus der G\"ultigkeit von $f(u,v) = 0$ f\"ur alle
 $v \in V$ folgt, dass $u = 0$ ist.}
 \medskip
       Ferner haben wir gesehen, wie die zu einer Korrelation inverse
 Korrelation sich darstellen l\"asst.
 \medskip\noindent
 {\bf 8.6. Korollar.} {\it Es sei $V$ ein $K$-Vektorraum mit
 $3 \leq \Rg_K (V) < \infty$ und $\kappa$ sei eine Korrelation von
 $L(V)$. Wird $\kappa$ durch die $\alpha$-Se\-mi\-bi\-li\-ne\-ar\-form $f$
 dargestellt, so wird $\kappa^{-1}$ durch die
 $\alpha^{-1}$-Semibilinearform $g$ dar\-ge\-stellt, die durch $g(u,v)
 := f(v, u)^{\alpha^{-1}}$ definiert wird.}
 \medskip
       Schlie\ss lich geben wir noch Antwort auf die Frage, wann zwei
 Semibilinearformen die gleiche Korrelation darstellen.
 \medskip\noindent
 {\bf 8.7. Satz.} {\it Es sei $V$ ein Vektorraum \"uber $K$ mit
 $3 \leq \Rg_K (V) < \infty$. Ferner sei $f$ eine nicht ausgeartete
 $\alpha$-Semibilinearform auf $V$.
 \item{(a)} Ist $k \in K^*$ und definiert man $g$ durch $g(u, v) :=
 k^{-1}f(u,v)$ f\"ur alle $u$, $v \in V$, so ist $g$ eine nicht
 ausgeartete $\beta$-Semibilinearform, wobei $\beta$ der durch
 $x^\beta := k^{-1}x^{\alpha}k$ definierte Antiautomorphismus von
 $K$ ist. \"Uberdies induzieren $f$ und $g$ die gleiche
 Korrelation.
 \item{(b)} Ist $g$ eine nicht ausgeartete $\beta$-Semibilinearform
 auf $V$ und stellen $f$ und $g$ die gleiche Korrelation von $L(V)$ dar, so gibt
 es ein $k \in K^*$ mit $kg(u,v) = f(u,v)$ f\"ur alle $u$, $v \in V$.\par}
 \smallskip
       Beweis. (a) Der Leser verifiziert m\"uhelos, dass $g$ eine
 $\beta$-Se\-mi\-bi\-li\-ne\-ar\-form ist. Dass $f$ und $g$ die gleiche
 Korrelation von $L(V)$ darstellen, folgt aus der Bemerkung,
 dass $f(u,v) = 0$ genau dann gilt, wenn $g(u,v) = 0$ ist.
 \par
       (b) Wir definieren zwei Abbildungen $\varphi$ und $\psi$ von $V$ in
 $V^*$ durch $v^\varphi(x) := f(v,x)$ bzw. $v^\psi(x) := g(v,x)$. Weil $f$ und $g$
 nicht ausgeartet sind, sind $\varphi$
 und $\psi$ injektiv. Fasst man $V^*$ als Rechtsvektorraum \"uber
 $K_\circ$ auf, so ist $(\psi,\beta)$ eine semilineare Abbildung
 von $V$ in $V^*$. Daher sind $u^\psi$ und $v^\psi$ genau dann
 linear unabh\"angig, wenn $u$ und $v$ es sind.
 \par
       Weil $f$ und $g$ die gleiche Korrelation darstellen, gilt
 $$ \Kern (v^\varphi) = \Kern (v^\psi) $$
 f\"ur alle $v \in V$. Es gibt also zu jedem $v \in V - \{0\}$ genau ein
 $l(v) \in K^*$ mit $v^\varphi = l(v)v^\psi$. Sind $u$, $v \in V - \{0\}$, so ist
 $$\eqalign{
       l(u + v)u^\psi + l(u + v)v^\psi &= l(u + v)(u + v)^\psi     \cr
				       &=(u + v)^v                 \cr
                                       &=u^\varphi + v^\varphi     \cr
				       &= l(u)u^\psi + l(v)v^\psi. \cr}$$
 Sind nun $u$ und $v$ linear unabh\"angig, so sind es auch $u^\psi$
 und $v^\psi$, wie wir bemerkten. In diesem Falle gilt also $l(u) =
 l(u + v) = l(v)$. Sind $u$ und $v$ linear abh\"angig, so gibt es
 wegen $\Rg_K(V) \geq 2$ ein $w \in V$ mit $u$, $v \in wK$. Nach dem
 bereits Bewiesenen ist dann $l(u) = l(w) = l(v)$. Es gibt also ein
 $k \in K^*$ mit $l(u) = k$ f\"ur alle $v \in V - \{0\}$. Wegen
 $0^\varphi = 0 = 0^\psi$ gilt also $v^\varphi = kv^\psi$ f\"ur
 alle $v \in V$. Hieraus folgt schlie\ss lich
 $$ f(u,v) = u^\varphi v = (ku^\psi)v = k(u^\psi v) = kg(u,v) $$
 f\"ur alle $u$, $v \in V$, so dass auch (b) bewiesen ist.
 \medskip
       Es gibt K\"orper, die einen Antiautomorphismus besitzen, aber keinen
 involutorischen solchen; siehe Morandi {\it et alii\/} 2005.

\mysection{9. Quaternionenschiefk\"orper}

\noindent
 Es sei $K$ ein K\"orper und $Z(K)$ sei sein
 Zentrum. In diesem
 Falle be\-zeich\-net man mit $[K:Z(K)]$ den Rang von $K$ als
 $Z(K)$-Vektorraum. Ist $K \neq Z(K)$, so ist $[K:Z(K)] \geq 4$,
 wie man wei\ss\ oder sich schnell \"uberlegt. Ist $[K:Z(K)] = 4$,
 so hei\ss t $K$ \emph{Quaternionenschiefk\"orper}. Diese lassen
 sich sehr sch\"on geometrisch kennzeichnen, wie wir im n\"achsten
 Abschnitt sehen werden. Zun\"achst wollen wir aber zeigen, dass es
 Quaternionenschiefk\"orper in H\"ulle und F\"ulle gibt. Hierzu
 l\"osen wir eine Serie von Aufgaben aus meinem Algebrabuch
 (L\"uneburg 1973, S.\ 145ff.). Diese Aufgabenserie hat eine
 Geschichte, die ich dem Leser nicht vorenthalten will: Sie spielt
 im akademischen Jahr 1970/71, f\"ur unsere immer noch junge
 Universit\"at\footnote{Sie ist mittlerweile zur Technischen Universit\"at
 verk\"ummert.} das Jahr Eins ihres Bestehens.
 \par
 Wir waren bei der Diskussion der Diplompr\"ufungsordnung, die im \"Ubrigen bis
 heute (WS 1989/90 und Jahre dar\"uber hinaus) nicht abgerissen ist. Eines der
 studentischen Mitglieder des Fach\-be\-reichs\-ra\-tes, Gert Schneider,
 sp\"ater zu makabrer Ber\"uhmtheit gelangt, verlangte, dass die
 Pr\"aambel in der Pr\"u\-fungs\-ord\-nung ersatzlos zu strei\-chen sei. In
 dieser Pr\"aambel stand n\"amlich, 
 dass es Zweck der Pr\"ufung sei fest zu stellen, ob der
 Kandidat nach wissenschaftlichen Grunds\"atzen zu arbeiten gelernt
 habe. Schneider meinte zur Begr\"undung seines Verlangens, es
 garantiere ja niemand, dass die Studenten w\"ahrend ihres Studiums
 auch wirklich die M\"oglichkeit h\"atten, dies zu lernen. Dies
 traf mich zutiefst in meiner Berufsehre. Voller Zorn verfasste
 ich daraufhin eine Serie von drei \"Ubungsbl\"attern mit eben
 jenen Aufgaben, wobei ich auf dem ersten Blatt verschiedene
 M\"oglichkeiten des wissenschaftlichen Arbeitens erl\"auterte,
 darunter die M\"oglichkeit der Ver\-all\-ge\-mei\-ne\-rung bekannter
 Resultate. Im vorliegenden Falle hatte ich n\"amlich die
 fraglichen Resultate f\"ur den Spezialfall des Ringes der ganzen
 Zahlen in der Vorlesung vorgerechnet. Diese Aufgabenreihe war also
 f\"ur die damaligen Studenten eine fr\"uhe --- nicht die erste ---
 Anleitung zum wissenschaftlichen Arbeiten. Rosemarie Rink\index{Rink, R.}{} und
 Manfred Dugas,\index{Dugas, M.}{} dem mathematischen Publikum durch sch\"one
 Arbeiten bekannt, haben damals schon durch die L\"osung dieser Aufgaben ihr
 Talent bewiesen.
 \par
       Im Folgenden bezeichne $\omega$ die Menge der nicht negativen
 ganzen Zahlen.
 \par
       Es sei $R$ ein Hauptidealbereich und $Q(R)$ sei sein
 Quo\-ti\-en\-ten\-k\"or\-per. Wir fassen $Q(R)$ auf als Modul \"uber $R$.
 Ist $p$ ein Primelement, so bezeichnen wir mit $Z(p^\infty)$ den
 Teilmodul des Faktormoduls $Q(R)/R$, der aus allen Elementen der
 Form
 $$ {r \over p^n} + R $$
 mit $r \in R$ und einer nicht negativen ganzen Zahl $n$ besteht. Die
 Bezeichnung $Z(p^\infty)$ erinnert an die aus der Theorie der abelschen Gruppen
 be\-kann\-ten Pr\"ufergruppen\index{prufergruppe@Pr\"ufergruppe}{}, deren
 Verallgemeinerungen sie sind. Wir
 nennen $Z(p^\infty)$ sinn\-ge\-m\"a\ss\ {\it Pr\"ufermodul zum Primelement\/}
 $p$\index{prufermodul@Pr\"ufermodul}{}
 \"uber $R$. Ferner setzen wir $R_p := \End_R(Z(p^\infty))$. Diesen Ring werden
 wir zun\"achst untersuchen.
 \medskip\noindent
 {\bf 9.1. Satz.} {\it Es sei $R$ ein Hauptidealbereich und $Z(p^\infty)$ sei der
 Pr\"u\-fer\-mo\-dul zum Primelement $p$ von $R$.  F\"ur $i \in \omega$ setzen
 wir $U_i := \{a \mid a \in Z(p^\infty),\, p^ia = 0\}$. Dann gilt:
 \item{a)} Es ist $U_i \subseteq U_{i+1}$ f\"ur alle $i \in \omega$.
 \item{b)} Es ist $\bigcup_{i\in \omega} U_i = Z(p^\infty)$.
 \item{c)} Es ist $U_i = ({1 \over p^{i}} + R) R$ f\"ur alle $i$.
 \item{d)} Ist $V$ ein Teilmodul von $Z(p^\infty)$, so ist entweder
 $V = Z(p^\infty)$, oder es gibt ein $i \in \omega$ mit $V = U_i$.\par}
 \smallskip
       Beweis. a), b) und c) sind banal zu beweisen. d) Es sei $V$ ein
 von $Z(p^\infty)$ verschiedener Teilmodul. Es gibt dann eine
 nat\"urliche Zahl $n$ und ein $r \in R$ mit ${r \over p^n} + R
 \not\in V$. Unter all diesen $n$ gibt es ein kleinstes, welches
 wir wiederum mit $n$ bezeichnen. Es folgt, dass $U_{n-1} \subseteq
 V$ ist. Ist andererseits ${x \over p^m} + R \in V$, so d\"urfen wir
 annehmen, dass $p$ kein Teiler von $x$ ist. Weil $p$ ein
 Primelement ist, ist dann $x$ zu $p^m$ teilerfremd, so dass es
 Elemente $u$ und $v$ in $R$ gibt mit $1 = ux + p^m v$. Hieraus
 folgt ${1 \over p^m} + R \in V$. Dies hat $U_m \subseteq V$ zur
 Folge. Also ist $m \leq n - 1$ und daher
 $$ {x \over p^m} + R \in U_m \subseteq U_{n-1}. $$
 Folglich ist $V \subseteq U_{n-1}$.  Damit ist auch d) bewiesen.
 \medskip\noindent
 {\bf 9.2. Satz.} {\it Es sei $R$ ein Hauptidealbereich und $p$ sei
 ein Primelement von $R$. Ist dann $\alpha$ ein Element des
 Endomorphismenringes $R_p$ des Pr\"ufermoduls $Z(p^\infty)$, so
 gilt $\alpha (U_i) \subseteq U_i$ f\"ur alle $i \in \omega$. Dabei
 sei $U_i$ wie in 9.1 definiert.}
 \smallskip
       Beweis. Es sei $x \in U_i$. Dann ist $p^i\alpha(x) =  \alpha(p^i x) = 0$.
 Also ist $\alpha (x) \in U_i$.
 \medskip\noindent
 {\bf 9.3. Satz.} {\it Es sei $R$ ein Hauptidealbereich und $p$ sei
 ein Primelement von $R$. Wir setzen $a_i := {1 \over p^i} + R$
 f\"ur alle $i \in \omega^*$, wobei $\omega^* := \omega - \{0\}$
 gesetzt wurde. Ist dann $\alpha$ ein Element des
 Endomorphismenringes $R_p$ von $Z(p^\infty)$, so gibt es eine
 Abbildung $f$ von $\omega^*$ in $R$ mit $\alpha (a_i) = f_i a_i$
 und $f_{i+1} \equiv f_i$ mod $p^i$ f\"ur alle $i \in \omega^*$.
 Ist umgekehrt $f$ eine Abbildung von $\omega^*$ in $R$ und gilt
 $f_{i+1} \equiv f_i$ mod $p_i$ f\"ur alle $i \in \omega^*$, so
 gibt es genau ein $\alpha \in R_p$ mit $\alpha (a_i) = f_i a_i$.}
 \smallskip
       Beweis. Nach 9.2 ist $\alpha (U_i) \subseteq U_i$ f\"ur alle
 $i \in \omega ^*$. Daher ist die Menge aller $r \in R$ mit $\alpha(a_i) = ra_i$
 nicht leer, so dass es auf Grund des Auswahlaxioms
 eine Abbildung $f$ von $\omega^*$ in $R$ gibt mit $\alpha (a_i) =
 f_ia_i$ f\"ur alle $i$. Es folgt
 $$ pf_ia_{i+1} = f_ia_i = \alpha(a_i) = \alpha (pa_{i+1}) = pf_{i+1} a_{i+1} $$
 und damit
 $$ pf_i \equiv pf_{i+1} \mod p^{i+1}, $$
 dh.,
 $$ f_i \equiv f_{i+1} \mod p^i. $$
 \par
       Es sei umgekehrt $f$ eine Abbildung von $\omega^*$ in $R$ und es
 gelte $f_i \equiv f_{i+1} \mod p^i$ f\"ur alle $i \in \omega^*$.
 Mittels Induktion folgt die G\"ultigkeit von
 $$ f_i \equiv f_{i+k} \mod p^i $$ f\"ur alle $i$, $k \in \omega^*$. Wir
 definieren eine bin\"are Relation $\alpha$ auf $Z(p^\infty)$ durch
 $(u,v) \in \alpha$ genau dann, wenn es ein $x \in R$ und ein $i
 \in \omega^*$ gibt mit $u = xa_i$ und $v = xf_ia_i$. Sind $(u,v)$,
 $(u,w) \in \alpha$, so ist also $u = xa_i$, $v = xf_ia_i$, $u =
 ya_k$ und $w = yf_ka_k$ f\"ur passende $x$, $y$, $i$ und $k$. Es
 sei oBdA $i \leq k$. Dann ist $xp^{k-i} a_k = u = ya_k$ und daher
 $$ y \equiv xp^{k-i} \mod p^k. $$
 Andererseits ist $f_k \equiv f_i$ mod $p^i$ und daher
 $$ f_kp^{k-i} \equiv f_ip^{k-i} \mod p^k. $$
 Folglich ist
 $$ yf_k \equiv xp^{k-i} f_k \equiv xp^{k-i} f_i \mod p^k. $$
 Hiermit folgt
 $$ v \equiv xf_ia_i = xf_ip^{k-i} a_k = yf_ka_k = w. $$
 Ist andererseits $u \in Z(p^\infty)$, so gibt es
 nach 9.1 ein $x \in R$ und ein $i \in \omega^*$ mit $u = xa_i$.
 Daher ist $(u, xf_ia_i) \in \alpha$. Damit ist $\alpha$ als
 Abbildung von $Z(p^\infty)$ in sich erkannt.
 \par
       Es seien $u$, $v \in Z(p^\infty)$. Es gibt dann $x$, $y \in R$ und $i$,
 $k \in \omega^*$ mit $u = xa_i$ und $v = ya_k$. Es sei oBdA $i\leq k$. Dann ist
 $$\eqalign{
       \alpha (u+v) &= \alpha \bigl((xp^{k-i} + y)a_k\bigr) \cr
		    &= (xp^{k-i} +y)f_ka_k        \cr
		    &= xf_ka_i + yf_ka_k. \cr} $$
 Wegen $f_i \equiv f_k \mod p^i$ folgt weiter
 $\alpha(u + v) = \alpha(u) + \alpha(v)$, so dass $\alpha$ additiv ist. Ist
 $r \in R$, so folgt auch noch
 $$ \alpha (ru) = rxf_ia_i = r\alpha (u), $$
 so dass $\alpha$ in der Tat ein Endomorphismus mit den
 verlangten Eigenschaften ist. Die Einzigkeit von $\alpha$ folgt
 aus der Bemerkung, dass die Menge der $a_i$ ein Erzeugendensystem
 von $Z(p^\infty)$ ist.
 \medskip\noindent
 {\bf 9.4. Satz.} {\it Es sei $R$ ein Hauptidealbereich und $p$ sei
 ein Primelement von $R$. Sind $\alpha$, $\beta \in R_p$ und werden
 $\alpha$ und $\beta$ gem\"a\ss\ 9.3 durch die beiden Abbildungen
 $f$ und $g$ von $\omega^*$ in $R$ dargestellt, so werden
 $\alpha+\beta$ und $\alpha \beta$ durch $f + g$ bzw. $fg$
 dargestellt, wobei Summe und Produkt von $f$ mit $g$
 punktweise definiert sei. Insbesondere folgt, dass $R_p$
 kommutativ ist.}
 \smallskip
       Beweis. Banal.
 \medskip\noindent
 {\bf 9.5. Satz.} {\it Es sei $R$ ein Hauptidealbereich und $p$ sei ein
 Primelement von $R$. Wir definieren $\pi \in R_p$ durch $\pi (x) := px$ f\"ur
 alle $x \in Z(p^\infty)$. Ist dann $Q$ ein von $\{0\}$ verschiedenes Ideal von
 $R_p$ und ist $V := \bigcap_{\alpha \in Q} \Kern(\alpha)$, so ist $V = U_n$ mit
 einem $n \in \omega$ und $Q = \pi^n R_p$.}
 \smallskip
       Beweis. Weil $V$ ein Teilmodul von $Z(p^\infty)$ ist, ist nach 9.1
 entweder $V = Z(p^\infty)$ oder $V = U_n$ mit einem passenden $n$. Weil
 $Q \neq \{0\}$ ist, kann der erste Fall nicht eintreten.  Also ist $V = U_n$.
 \par
       Es sei $0 \neq \alpha Q$ und $U_m$ sei der Kern von $\alpha$. Dann
 ist $n \leq m$. Es werde $\alpha$ gem\"a\ss\ 9.3 durch die Abbildung
 $f$ dargestellt. Dann ist $f_m \equiv 0 \mod p^m$ und daher
 $f_{m+i} \equiv 0 \mod p^m$ f\"ur alle $i \in \omega^*$.
 Andererseits ist $f_{m+i}$ f\"ur alle $i \in \omega^*$ nicht durch
 $p^{m+1}$ teilbar, da sonst $a_{m+1}$ wegen
 $f_{m+i} \equiv f_{m+1} \mod p^{m+1}$ im Kern von $\alpha$ l\"age. Insbesondere
 ist also $f_{m+i} \neq 0$ f\"ur alle $i \in \omega^*$. Es gibt
 daher zu jedem $i \in \omega^*$ genau ein $r_i \in R$ mit $f_{m+i}
 = r_ip^m$. Es folgt $r_{i+1}p^m \equiv r_ip^m \mod p^{m+1}$, was
 die Kongruenz $r_{i+1} \equiv r_i \mod p^i$ nach sich zieht. Nach
 9.3 wird also durch $\rho (a_i) := r_ia_i$ ein $\rho \in R_p$
 definiert. Nun ist
 $$ \rho \pi^m (a_i) = \rho (p^ma_i) = 0 = \alpha (a_i), $$
 falls $i \leq m$ ist. Andererseits ist
 $$ \rho \pi_m(a^{m+i}) = \rho(a_i) = r_ia_i
		        = r_ip^m a_{m+i} = \alpha (a_{m+i}) $$
 f\"ur alle $i>0$. Also ist $\rho \pi^m = \alpha$. W\"are nun $\rho$ keine
 Einheit, so w\"are $U_1 \subseteq \Kern(\rho)$, was wegen
 $r_1 \not\equiv 0 \mod p$ nicht der Fall ist. Weil $Q$ ein Ideal ist, ist also
 $\pi^m \in Q$.
 \par
       Es sei nun $k$ das kleinste unter allen $m \in \omega$ mit $\pi^m
 \in Q$. Dann gibt es also zu jedem $\alpha \in Q$ ein $\sigma \in
 R_p$ mit $\alpha  = \sigma \pi^k$. Also ist
 $$ Q \subseteq \pi^k R_p \subseteq Q. $$
 Folglich ist $Q = \pi^k R_p$. Dann ist aber $U_n = \Kern(\pi^k)$ und damit
 $k = n$, was zu beweisen war.
 \medskip
 Der Beweis von 9.5 liefert auch noch den folgenden Satz.
 \medskip\noindent
 {\bf 9.6. Satz.} {\it Es sei $R$ ein Hauptidealbereich und $p$ sei
 ein Primelement von $R$. Es sei weiter $\pi$ der durch $\pi (x) :=
 px$ definierte Endomorphismus von $Z (p^\infty)$. Ist dann $0 \neq
 \alpha \in R_p$, so gibt es eine Einheit $\rho$ in $R_p$ und ein
 $n \in \omega$ mit $\alpha = \rho \pi^n$.}
 \medskip
       Weil Einheiten stets bijektiv sind und $\pi$ offensichtlich
 surjektiv ist, gilt auch
 \medskip\noindent
 {\bf 9.7. Korollar.} {\it Jeder von Null verschiedene
 Endomorphismus von $Z (p^\infty)$ ist surjektiv.}
 \medskip
       Weil das Produkt zweier surjektiver Abbildungen surjektiv ist,
 haben wir ferner
 \medskip\noindent
 {\bf 9.8. Korollar.} {\it Die Ringe $R_p$ sind Integrit\"atsbereiche.}
 \medskip
       Da $R_p$ Integrit\"atsbereich ist, existiert der
 Quotientenk\"orper $Q(R_p)$. Um ein hinreichendes Kriterium daf\"ur zu
 erhalten, dass gewisse dieser K\"orper nicht isomorph sind,
 untersuchen wir nun die Gruppen von Einheitswurzeln in $R_p$.
 Bevor wir dies tun, sei noch darauf hin\-ge\-wie\-sen, dass $R_p$ der
 Ring der ganzen henselschen $p$-adischen Zahlen ist, falls $R$ der
 Ring der ganzen Zahlen ist. In diesem Falle ist $Q(R_p)$ der
 K\"orper der henselschen $p$-adischen
 Zahlen.\index{henselsche $p$-adische Zahlen}{}
 \medskip\noindent
 {\bf 9.9. Satz.} {\it Es sei $R$ ein Hauptidealbereich und $p$ sei ein
 Primelement von $R$. Ist $\sigma \in R_p$ und wird $\sigma$ gem\"a\ss\ 9.3 durch
 $f$ dargestellt, so setzen wir $\psi(\sigma) := f_i + p^i R$. Dann ist $\psi$
 ein Epimorphismus von $R_p$ auf $R/p^iR$. Ferner gilt $\Kern(\psi) = \pi^i R_p$,
 so dass $R_p/\pi^i R_p$ und $R/p^iR$ isomorph sind.}
 \smallskip
       Beweis. Wird $\sigma$ gem\"a\ss\ 9.3 auch durch $g$ dargestellt,
 so ist $g_i \equiv f_i \mod p^i$, so dass $\psi$ wohldefiniert
 ist. Dann ist aber klar, dass $\psi$ ein Homomorphismus ist.
 \par
       Ist $r \in R$, so wird durch $\rho(x) := rx$ f\"ur alle
 $x \in Z(p^\infty)$ ein $\rho \in R_p$ definiert. Wegen $\psi(\rho) = r + p^iR$
 ist $\psi$ auch surjektiv.
 \par
       Es sei $\psi(\sigma) = 0$. Dann ist $f_i \equiv 0 \mod p^i$. Ist
 $\sigma \neq 0$, so gibt es nach 9.6 eine Einheit $\tau$ in $R_p$
 und ein $m \in \omega$ mit $\sigma = \tau \pi^m$. Wird $\tau$
 durch $t$ dargestellt, so folgt $f_i \equiv t_ip^m \mod p^i$.
 Weil $\tau$ eine Einheit ist, ist $t_i$ nicht durch $p$ teilbar.
 Wegen $f_i \equiv 0 \mod p^i$ ist daher $i \leq m$, womit alles
 bewiesen ist.
 \medskip\noindent
 {\bf 9.10. Satz.} {\it Es sei $R$ ein Hauptidealbereich und $p$ sei ein
 Primelement von $R$. Mit $W(R_p)$ bezeichnen wir die Gruppe der
 Einheitswurzeln\index{Gruppe der Einheitswurzeln}{}
 von $R_p$, dh., die Torsionsgruppe\index{Torsionsgruppe}{} der Einheitengruppe
 von $R_p$. Setzt man
 $$ \varphi (\epsilon) := \epsilon + \pi R_p $$
 f\"ur alle $\epsilon \in W(R_p)$, so ist
 $\pi$ ein Epimorphismus von $W(R_p)$ auf die Gruppe der
 Einheitswurzeln $W(R_p/\pi R_p)$ des K\"orpers $R_p/\pi R_p$.}
 \smallskip
       Beweis. Es sei $\epsilon \in W(R_p)$ und es gelte $\epsilon^n = 1$. Dann
 ist $\varphi(\epsilon)^n = \varphi(\epsilon^n) = 1$, so
 dass $\varphi$ ein Homomorphismus in $W (R_p/\pi R_p)$ ist.
 \par
       Um die Surjektivit\"at von $\varphi$ zu beweisen, sei $n$ die
 Ordnung des Elementes $\eta + \pi R_p \in W(R_p/\pi R_p)$. Wendet
 man 9.9 mit $i = 1$ an, so folgt die Existenz eines $f_1 \in R$
 mit $\psi(\eta) = f_1 + pR$. Es folgt, dass die Ordnung von $f_1 + pR$ gleich
 $n$ ist. Insbesondere ist also $f^n_1 \equiv 1 \mod p$. Es sei $i \geq 1$ und es
 gebe $f_1$, \dots, $f_i \in R$ mit
 $f_{j+1} \equiv f_j$ mod $p^j$ und $f^n_j \equiv 1$ mod $p^j$
 f\"ur alle $j \leq i$. Insbesondere ist dann $f_i^n = 1 + kp^i$
 mit einem $k \in R$. Ist die Charakteristik $q$ von $R/pR$
 ungleich $0$, so hat $R/pR$ nur die triviale $q$-te Einheitswurzel
 1, so dass $q$ kein Teiler von $n$ ist. Dann ist aber $p$ kein
 Teiler von $q \cdot 1$, wobei 1 die Eins von $R$ bezeichne. Weil
 auch $f_i$ zu $p$ teilerfremd ist, gibt es ein $v \in R$ mit
 $$ nf_i^{n-1} v \equiv -k \mod p. $$
 Es folgt weiter
 $$ nf^{n-1}_i vp^i + kp^i \equiv 0 \mod  p^{i+1}. $$
 Setze $f_{i+1} := f_i + vp^i$. Dann ist $f_{i+1} \equiv f_i \mod p^i$ und
 $$\eqalign{
    f^n_{i+1} &\equiv f^n_i + nvp^i f^{n-1}_i     \cr
              &\equiv 1 + kp^i + nvp^if^{n-1}_i   \cr
	      &\equiv 1 \mod p^{i+1}. \cr} $$
 Damit ist --- auf die \"ubliche, aber schlampige Weise --- die Existenz einer
 Abbildung $f$ der nat\"urlichen Zahlen in $R$ sichergestellt mit
 $f_{i+1} \equiv f_i \mod p^i$ und $f^n_i \equiv 1 \mod p^i$ f\"ur alle $i$. Nach
 9.3 gibt es nun ein $\epsilon \in R_p$, welches durch $f$ dargestellt wird. Es
 folgt
 $$ \epsilon^n = 1 $$ sowie
 $$ \psi(\epsilon) = f_1 + pR = \psi(\eta) $$
 und damit
 $$ \varphi(\epsilon) = \epsilon + \pi R_p = \eta + \pi R_p, $$
 womit auch die Surjektivit\"at von $\varphi$ nachgewiesen ist.
 \medskip
       Ist $R$ ein Integrit\"atsbereich, so erbt $R$ von seinem
 Quotientenk\"orper $Q(R)$ die Charakteristik,\index{Charakteristik}{} da $R$ ja zu
 einem Teilring von $Q (R)$ isomorph ist.
 \medskip\noindent
 {\bf 9.11. Satz.} {\it Es sei $R$ ein Hauptidealbereich und $p$
 sei ein Primelement von $R$. Ferner sei $\varphi$ der in 9.10
 definierte Epimorphismus von $W(R_p)$ auf $W(R_p/\pi R_p)$. Ist
 $\Kern(\varphi) \neq \{1\}$, so ist die Charakteristik von
 $R$ gleich $0$ und die Charakteristik $q$ von $R_p/\pi R_p$ ist
 gr\"o\ss er als $0$. \"Uberdies ist $p$ ein Teiler von $q \cdot
 1$, wobei $1$ die Eins von $R$ bezeichne.}
 \smallskip
       Beweis. Es sei $1 \neq \alpha \in \Kern(\varphi)$ und $n$
 sei die Ordnung von $\alpha$. Ferner sei $q$ eine $n$ teilende
 Primzahl. Setze $\beta := \alpha^{n/q}$. Dann ist $\beta$ ein
 Element der Ordnung $q$ im Kern von $\varphi$. Es gibt daher ein
 $\rho \in R_p$ mit $\beta = 1 + \rho \pi$. Es folgt
 $$ 1 = \beta^q = 1 + \sum_{i=1}^q {q \choose i} \rho^i \pi^i. $$
 Weil $R_p$ nullteilerfrei ist, folgt hieraus
 $$ q \cdot 1 \equiv 0 \mod \pi. $$
 Daher ist $q$ die Charakteristik von $R_p/\pi R_p$ und es folgt, dass
 $\Kern(\varphi)$ eine $q$-Gruppe ist. W\"are die Charakteristik von
 $R_p$ nicht Null, so w\"are sie gleich $q$. Es folgte
 $$ 1 = \beta^q = (1 + \rho \pi)^q = 1 + \rho^q \pi^q $$
 und damit $\rho^q \pi^q = 0$, so dass $\rho = 0$ w\"are im Widerspruch zu
 $\beta \neq 1$. Damit ist der Satz in allen seinen Teilen
 bewiesen.
 \medskip
       Ist $R$ ein Hauptidealbereich, so erbt $R$, wie schon bemerkt, von
 $Q(R)$ die Charakteristik. Ist diese positiv, so hat auch $R/pR$
 diese Charakteristik, falls nur $p$ ein Primelement von $R$ ist.
 Ist $R/pR$ endlich, so k\"onnen wir $R_p$ auf eine andere Art
 darstellen, die das Isomorphieproblem in diesem Falle l\"ost.
 Genaueres sagt der n\"achste Satz.
 \medskip\noindent
 {\bf 9.12. Satz.} {\it Es sei $R$ ein Hauptidealbereich der
 Charakteristik $q > 0$. Ist $p$ ein Primelement von $R$ und ist
 $R/pR$ algebraisch \"uber seinem Primk\"orper, so ist der
 Endomorphismenring $R_p$ von $Z(p^\infty)$ dem Ring der formalen
 Potenzreihen \"uber $R/pR$ isomorph.}\index{Ring der formalen Potenzreihen}{}
 \smallskip
       Beweis. Weil $R/pR$ zu $R_p/\pi R_p$ isomorph ist, ist auch dieser
 K\"orper algebraisch \"uber seinem Primk\"orper, der zu $\GF(q)$
 isomorph ist. Dies impliziert, dass der K\"orper $R_p/\pi R_p$ die
 Vereinigung der in ihm ent\-hal\-te\-nen endlichen Teilk\"orper ist.
 Hieraus folgt wiederum
 $$ W(R_p/\pi R_p) = (R_p/\pi R_p)^*. $$
 \par
       Setze $K:= W(R_p) \cup \{0\}$. Wir zeigen, dass $K$ ein K\"orper
 ist. Da $K$ offensichtlich multiplikativ abgeschlossen ist, m\"ussen
 wir nur noch zeigen, dass $K$ auch additiv abgeschlossen ist. Dazu
 seien $\alpha$ und $\beta$ zwei von $0$ verschiedene Elemente aus
 $K$. Dann sind $\alpha + \pi R_p$ und $\beta + \pi R_p$ von Null
 verschiedene Elemente von $R_p/\pi R_p$, die \"uberdies in einem
 endlichen Teilk\"orper von $R_p/\pi R_p$ liegen. Dieser habe die
 Ordnung $q^n$. Dann gilt
 $$ (\alpha + \pi R_p)^{q^n-1} = 1 + \pi R_p $$
 und
 $$ (\beta + \pi R_p)^{q^n-1} = 1 + \pi R_p. $$
 Wegen 9.10 und 9.11 ist dann $\alpha ^{q^n-1} = 1$ und $\beta^{q^n-1} = 1$.
 Hieraus folgt weiter $\alpha^{q^n} = \alpha$ und $\beta^{q^n} = \beta$. Weil
 $q$ auch die Charakteristik von $R_p$ ist, folgt schlie\ss lich
 $$ (\alpha + \beta)^{q^n} = \alpha + \beta. $$
 Ist nun $\alpha + \beta = 0$, so ist nichts zu beweisen. Ist dies nicht der
 Fall, so folgt
 $$ (\alpha + \beta)^{q^n -1} = 1 $$
 und damit $\alpha + \beta \in K$. Folglich ist $K$ ein K\"orper. Weil
 $K \cap \pi R_p = \{0\}$ ist, folgt schlie\ss lich, dass $K$ zu $R_p/\pi R_p$
 isomorph ist.
 \par
       Offensichtlich gilt $\bigcap_{n:=1}^{\infty} \pi^n R_p = \{0\}$.  Nimmt
 man diese Ideale als Umgebungsbasis\index{Umgebungsbasis}{} der Null,
 so erh\"alt man eine hausdorffsche Topologie\index{hausdorffsche Topologie}{}
 $T$ auf $R_p$, die $R_p$ zu einem topologischen Ring\index{topologischer Ring}{}
 macht. Das Element $\pi$ ist transzendent
 \"uber $K$. Ist n\"amlich $n$ eine nat\"urliche Zahl, sind $n$ und
 $k_0$, \dots, $k_n \in K$ und gilt $0 = \sum^n_{i:=0} k_i \pi^i$,
 so folgt wegen $K \cap \pi R_p = \{0\}$ zun\"achst $k_0 = 0$ und
 dann mittels Induktion $k_i = 0$ f\"ur $i := 1$, \dots, $n$. Somit
 ist $K[\pi]$ zum Polynomring in einer Unbestimmten \"uber $K$
 isomorph. Mit $K_n[\pi]$ bezeichnen wir die Menge der Polynome in
 $\pi$, deren Grad $n$ nicht \"ubersteigt. Dann ist
 $$ R_p = K_n [\pi] \oplus \pi^{n+1} R_p. $$
 Dies folgt mit einer leichten Induktion \"uber $n$. Zu
 jedem $r \in R_p$ und zu jedem $n \in \omega$ gibt es also genau
 ein $c_n \in K_n [\pi]$ mit
 $$r \equiv c_n \mod \pi^{n+1}. $$
 Es folgt, dass $c$ bez. der Topologie $T$ eine Cauchyfolge ist und dass
 $r = \lim c$ gilt. Dies impliziert, dass $R_p$ ein Teilring des Rings $K[[\pi]]$
 der formalen Potenzreihen \"uber $K$ ist. Wir m\"ussen nun nur noch
 zeigen, dass jede Cauchyfolge \"uber $K[\pi]$ in $R_p$ konvergiert.
 \par
       Um dies zu zeigen, sei $c$ eine Cauchyfolge \"uber $K[\pi]$. Es
 sei $n \in \omega$.  Es gibt dann ein $M \in \omega$ mit
 $$ c_u - c_v \in \pi^n R_p $$
 f\"ur alle $u$ und $v$ mit $u$, $v \geq M$. Setzt
 man wieder $a_i := {1 \over \pi^i} + R$, so ist also
 $$ c_u(a_i) = c_v(a_i) $$
 f\"ur alle $i \leq n$ und alle $u$, $v \geq M$. Es sei $M(n)$ das kleinste unter
 diesen $M$. Dann ist $M(n) \leq M(n+1)$ f\"ur alle $n$. Es gibt ferner
 $r_1$, \dots, $r_n \in R$ mit
 $$ c_u(a_i) = r_ia_i $$
 f\"ur alle $i \geq n$ und alle $u \leq M(n)$. Es gibt also eine Auswahlfunktion
 $F$, die jedem $n$ ein $n$-Tupel $F_n$ von Elementen aus $R$ zuordnet, so dass
 $$ c_u(a_i) = F_{n,i}a_i $$
 gilt f\"ur alle $i \leq n$ und alle $u \geq M(n)$. Definiere $f$ durch
 $f_n := F_{n,n}$. Dann ist $c_u(a_n) = f_na_n$ f\"ur alle $u \geq M(n)$. Wegen
 $M(n+1) \geq M(n)$ ist dann
 $$\eqalign{
   f_{n+1}a_n &= f_{n+1}pa_{n+1}       \cr
	      &= pc_{M(n+1)}(a_{n+1})  \cr
              &= c_{M(n+1)} (a_n)      \cr
	      &= c_{M(n)}(a_n)         \cr
	      &= f_na_n. \cr} $$
 Hieraus folgt
 $$ f_{n+1} \equiv f_n \mod p^n, $$
 so dass es nach 9.3 ein Element $\rho \in R_p$ gibt mit $\rho (a_n) = f_na_n$
 f\"ur alle $n$. Es folgt
 $$ (\rho - c_u)a_n = f_na_n - F_{n,n}a_n = 0 $$
 f\"ur alle $u \geq M(n)$ und damit
 $$ \rho - c_u \in \pi^n R_p $$
 f\"ur alle $u \geq M(n)$, so dass $\lim c = \rho$ gilt. Damit ist alles bewiesen.
 \medskip\noindent
 {\bf 9.13. Korollar} {\it Es sei $R$ Hauptidealbereich mit von
 Null ver\-schie\-de\-ner Charakteristik. Ist $p$ ein Primelement von $R$
 und ist $R/pR$ algebraisch \"uber seinem Primk\"orper, so ist der
 Quotientenk\"orper $Q(R_p)$ von $R_p$ isomorph zum K\"orper der
 formalen Laurentreihen \"uber
 $R/pR$.}\index{korper@K\"orper!der formalen Laurentreihen}{}
 \smallskip
       Beweis. Dies folgt unmittelbar aus 9.12 und der Bemerkung, dass
 der Quotientenk\"orper des Ringes der formalen
 Potenzreihen\index{Ring der formalen Potenzreihen}{} \"uber
 einem K\"orper der K\"orper der formalen Laurentreihen \"uber eben
 die\-sem K\"orper ist.
 \medskip\noindent
 {\bf 9.14. Satz.} {\it Es seien $R$ und $S$ Hauptidealbereiche mit
 von Null verschiedener Charakteristik. Ferner seien $p$ und $q$
 Primelemente von $R$ bzw. $S$. Ist $R/pR$ algebraisch \"uber
 seinem Primk\"orper, so sind die folgenden Aussagen \"aquivalent:
 \item{a)} Die Ringe $R_p$ und $S_q$ sind isomorph.
 \item{b)} Die K\"orper $R_p/\pi R_p$ und $S_q/\rho S_q$ sind
 isomorph. Dabei seien mit $\pi R_p$ und $\rho S_q$ die maximalen
 Ideale der Ringe $R_p$ bzw. $S_q$ bezeichnet.
 \item{c)} Die Ringe $R/pR$ und $S/qS$ sind isomorph.
 \item{d)} Der K\"orper $S_q/\rho S_q$ ist algebraisch \"uber seinem
 Primk\"orper und die Quotientenk\"orper $Q(R_p)$ und $Q(S_q)$ von
 $R_p$ bzw. $S_q$ sind isomorph.\par}
 \smallskip
       Beweis. a) impliziert b). Jeder Isomorphismus von $R_p$ auf $S_q$
 bildet das maximale Ideal $\pi R_p$ von $R_p$ auf das maximale
 Ideal $\rho S_q$ von $S_q$ ab und induziert daher einen
 Isomorphismus von $R_p/\pi R_p$ auf $S_q/\rho S_q$.
 \par
 b) impliziert c). Weil $R/p R$ zu $R_p/\pi E_p$ und $S/qS$ zu
 $S_q/\rho S_q$ isomorph ist, sind $R/pR$ und $S/qS$ isomorph.
 \par
 c) impliziert a). Sind $R/pR$ und $S/qS$ isomorph, so sind es auch
 $R_p/\pi R_p$ und $S_q/\rho S_q$. Weil $R/R_p$ \"uber seinem
 Primk\"orper algebraisch ist, gilt dies daher auch f\"ur die
 zuletzt genannten K\"orper. Weil die Charakteristiken von $R$ und
 $S$ von Null verschieden sind, sind nach 9.12 die Ringe $R_p$ und
 $S_q$ folglich zu den Ringen von formalen Po\-tenz\-rei\-hen \"uber
 $R_p/ \pi R_p$ bzw. \"uber $S/\rho S_q$ isomorph. Daher sind $R_p$
 und $S_q$ isomorph.
 \par
       Isomorphe Integrit\"atsbereiche haben nat\"urlich isomorphe
 Quotientenk\"orper. Wir m\"ussen daher nur noch zeigen, dass a)
 eine Folge von d) ist. Hierzu m\"ussen wir auf den Satz
 zur\"uckgreifen, dass Haupt\-i\-de\-al\-be\-rei\-che --- und $R_p$ und $S_q$
 sind ja solche --- in ihren Quotientenk\"orpern ganz abgeschlossen
 sind. Dies bedeutet, dass die in $Q(R_p)$ liegenden Nullstellen
 eines Polynoms aus $R_p[x]$, dessen Leitkoeffizient 1 ist, in
 $R_p$ liegen. Dies zu beweisen, sei dem Leser als \"Ubungsaufgabe
 \"uberlassen. Es ist eine einfache Folgerung aus dem Satz \"uber
 die eindeutige Primfaktorzerlegung, der ja in Hauptidealbereichen
 gilt. Aus der Tatsache, dass $R_p$ und $S_q$ in ihren
 Quotientenk\"orpern ganz abgeschlossen sind, folgt, dass
 $W(Q(R_p)) = W (R_p)$ und $W(Q(S_q)) = W(S_q)$ ist.
 \par
       Es sei nun $\alpha$ ein Isomorphismus von $Q(R_p)$ auf $Q(S_q)$,
 also $\alpha(W(R_p)) = W(S_q)$. Setzt man $K:= W(R_p)
 \cup \{0\}$ und $L:= W(S_q) \cup \{0\}$, so induziert $\alpha$
 einen Isomorphismus des K\"orpers $K$ auf den K\"orper $L$. Weil
 sowohl $R_p/\pi R_p$ als auch $S_q/\rho S_q$ \"uber ihren
 Primk\"orpern algebraisch sind, sind die K\"orper $R_p/\pi R_p$
 und $S_q/\rho S_q$ isomorph, so dass b) und damit a) eine Folge
 von d) ist. Damit ist alles bewiesen.
 \medskip
 \goodbreak
       Mit den S\"atzen 9.13 und 9.14 erhalten wir daher zu jeder von 0
 verschiedenen Charakteristik \"uberabz\"ahlbar viele Beispiele von
 K\"orpern $Q(R_p)$, da $\GF(r)$ f\"ur jede Primzahl $r$
 \"uberabz\"ahlbar viele, nicht isomorphe, algebraische
 Er\-wei\-te\-rung\-en besitzt (Steinitz 1910).
 \par
       Mit $Z$ bezeichnen wir im folgenden den Ring der ganzen Zahlen.
 \medskip\noindent
 {\bf 9.15. Satz.} {\it Es sei $p$ eine Primzahl im Ring $Z$ der
 ganzen Zahlen. Dann ist $|W(Z_p)| = p - 1$, es sei denn, es ist
 $p = 2$. In diesem Falle gilt $|W(Z_2)| = 2$.}
 \smallskip
       Beweis. Es sei $\varphi$ der in 9.10 definierte Epimorphismus von
 $W(Z_p)$ auf $\GF(p)^*$. Ferner sei $1 \neq \alpha \in \Kern(\varphi)$. Nach 9.11
 ist dann $o (\alpha) = p^n$ mit einer nat\"urlichen Zahl $n$. Wir setzen
 $\beta := \alpha^{p^{n-1}}$. Dann ist $o(\beta) = p$. Es werde $\beta$
 entsprechend 9.3 durch $f$ dargestellt, wobei wir annehmen
 d\"urfen, dass $0 \leq f_i < p^i$ ist f\"ur alle $i$. Es folgt
 $$ f^p_i \equiv 1 \mod p^i $$
 f\"ur alle $i$. Insbesondere folgt $f_1 = 1$. Es sei $m$ die kleinste
 nat\"urliche Zahl, f\"ur die $f_{m+1} \neq 1$ gilt. Dann ist
 $$ f_{m+1} = f_m + vp^m = 1 + vp^m. $$
 Wegen $1< f_{m+1} < p^{m+1}$ folgt $0 < v < p$. Nun ist
 $$ f_{m+2} = f_{m+1} + wp^{m+1} $$
 und daher
 $$ 1 \equiv f^p_{m+2} \equiv f^p_{m+1} + pf^{p-1} wp^{m+1}
      \equiv f^p_{m+1} \mod p^{m+2}. $$
 W\"are $m \geq 2$, so folgte weiter
 $$ 1 \equiv (1 + vp^m)^p \equiv 1 + pvp^m \mod p^{m+2} $$
 und dann
 $$ 0 \equiv vp{m+1} \mod p^{m+2}, $$
 so dass
 $$ 0 \equiv v \mod p $$
 w\"are im Widerspruch zu $0 < v < p$. Also ist
 $m = 1$. Dann ist aber
 $$ 1 \equiv 1 + pvp + {p \choose 2} v^2 p^2 \mod p^3 $$
 und weiter
 $$ 0 \equiv pvp + {p \choose 2} v^2p^2 \mod p^3. $$
 W\"are $p > 2$, so w\"are $p \choose 2$ durch $p$ teilbar, und es folgte wieder
 der Widerspruch $v \equiv 0 \mod p$. Also ist $p = 2$.
 \par
 Es bleibt uns noch zu zeigen, dass $\alpha = -1$ ist. Um dies zu
 zeigen, m\"ussen wir zeigen, dass $n = 1$ ist. Dazu werde $\alpha$
 gem\"a\ss\ 9.3 durch $g$ dargestellt. Weil $\beta^2 = 1 \neq
 \beta$ ist, ist $\beta = -1$, dh., es ist $\alpha^{2^{n-1}} = -1$.
 Daher gilt insbesondere
 $$ g^{2^{n-1}}_2 \equiv -1 \mod 4. $$
 Da $-1$ kein Quadrat modulo 4 ist, folgt hieraus $n = 1$, so dass in der Tat
 $\alpha = -1$ ist. Damit ist alles bewiesen.
 \medskip\noindent
 {\bf 9.16. Satz.} {\it Es seien $p$ und $q$ zwei Primzahlen des
 Ringes $Z$ der ganzen Zahlen. Genau dann sind die
 Quotientenk\"orper $Q(Z_p)$ und $Q(Z_q)$ isomorph, wenn $p = q$
 ist.}
 \smallskip
 Beweis. Ist $p=q$, so ist nichts zu beweisen. Wir nehmen daher an,
 dass die fraglichen K\"orper isomorph seien. Es gilt wieder
 $W(Z_p) = W(Q(Z_p))$ und $W(Z_q) = W(Q(Z_q))$.  Aus der Isomorphie
 von $Q(Z_p)$ und $Q(Z_q)$ folgt nat\"urlich die Isomorphie
 $W(Z_p)$ und $W(Z_q)$, so dass mit 9.15 entweder $p - 1 = q - 1$ und
 damit $p=q$ oder aber, falls wir $q<p$ annehmen, $q=2$ und $p=3$
 folgt.
 \par
       Um zu zeigen, dass $Q(Z_2)$ und $Q(Z_3)$ nicht isomorph sind,
 zeigen wir, dass das Polynom $x^3 + x + 1$ in $Q(Z_3)$ eine
 Nullstelle hat, in $Q(Z_2)$ aber nicht. Letzteres ist sofort zu
 sehen. W\"are n\"amlich $\zeta$ eine Nullstelle dieses Polynoms in
 $Q(Z_2)$, so l\"age $\zeta$ in $Z_2$, da $Z_2$ ja ganz
 abgeschlossen ist. Dann w\"are aber $\zeta + \pi Z_2$ eine
 Nullstelle dieses Polynoms \"uber $\GF(2)$ im Widerspruch zu der
 Tatsache, dass dieses Polynom \"uber $\GF(2)$ nur den Wert 1
 annimmt.
 \par
       Setze $f_1 := 1$. Dann ist $f^3_1 + f_1 + 1 \equiv 0 \mod 3$. Es sei
 $n \geq 1$ und es gebe $f_i \in Z_3$ mit $f^3_i + f_i + 1 \equiv 0 \mod 3$
 f\"ur alle $i \leq  n$, f\"ur die au\ss erdem noch gilt, dass
 $f_{i+1} \equiv f_i \mod 3^i$ ist f\"ur alle $i < n$. Es gibt dann ein $k \in Z$
 mit $f^3_n + f_n + 1 = k3^n$. Setze $f_{n+1} := f_n - k3^n$. Dann gilt
 zum einen $f_{n+1} \equiv f_n \mod 3^n$ und zum Anderen
 $$ f^3_{n+1} = f^3_n - 3f^2_nk3^n + 3f_nk^23^{2n} - k^33^{3n}
	    \equiv f^3_n \mod 3^{n+1}. $$
 Daher ist
 $$ f^3_{n+1} + f_{n+1} + 1 \equiv f^3_n + f_n - k3^n + 1 = 0 \mod 3^{n+1}. $$
 Nach 9.3 folgt die Existenz eines
 $\eta \in Z_3$ mit $\eta (a_n) = f_na_n$. Dann ist aber
 $$ (\eta^3 + \eta + 1) (a_n) = (f^3_n + f_n + 1) a_n = 0, $$
 so dass $\eta$ in der Tat eine Nullstelle von $x^3 + x + 1$ ist. Damit ist
 auch gezeigt, dass auch $Q(Z_2)$ und $Q(Z_3)$ nicht isomorph sind.
 \medskip
       Dem Kenner wird nicht verborgen geblieben sein, dass wir in den Beweisen
 des letzten Satzes und des Satzes 9.10 Spezialf\"alle des hen\-sel\-schen
 Lemmas\index{henselsches Lemma}{} etabliert haben. Dem, der dieses Lemma nicht
 kennt, sei als gute \"Ubungsaufgabe empfohlen, dieses Lemma aus
 den beiden Beweisen herauszupr\"aparieren und zu beweisen. Als
 Hinweis sei erw\"ahnt, dass das henselsche Lemma eine Aussage \"uber
 die Fak\-to\-ri\-sie\-rung von Polynomen macht. Die hier \"uber die
 Existenz von Nullstellen gewisser Polynome gemachten Aussagen
 muss man also als Teilbarkeit durch Polynome (ersten Grades)
 interpretieren.
 \par
       Wir sind nun in der Lage, nicht isomorphe
 Qua\-ter\-ni\-o\-nen\-schief\-k\"or\-per \emph{en masse} zu konstruieren.
 \medskip\noindent
 {\bf 9.17. Satz.} {\it Es sei $R$ ein Hauptidealbereich und $p$
 sei ein Primelement von $R$. Ferner sei $c \in R_p$ und das
 Polynom $x^2 + x - c$ sei aufgefasst als Polynom \"uber $R_p/\pi
 R_p$ irreduzibel. Ist dann $H_p$ der Zentralisator der beiden
 Matrizen
 $$\pmatrix{ 0 & 1 & 0 & 0 \cr
	     c & 1 & 0 & 0 \cr
	     0 & 0 & 0 & 1 \cr
	     0 & 0 & c & 1 \cr}\quad
 \mathit{und}\quad
 \pmatrix{ 0   &  0  & 1 & 0 \cr
	     0   &  0  & 1 &-1 \cr
	     \pi &  0  & 0 & 0 \cr
	     \pi &-\pi & 0 & 0 \cr}$$
 im Ring der $(4 \times 4)$-Matrizen \"uber $Q(R_p)$, so ist
 $H_p$ ein Quaternionen\-schief\-k\"or\-per mit $Z(H_p) \cong Q(R_p)$.
 \"Uberdies besteht $H_p$ gerade aus den Matrizen
 $$\pmatrix{\alpha & \beta & \gamma & \delta \cr
	    c\beta & \alpha + \beta & c\delta  & \gamma + \delta \cr
	    \pi(\gamma + \delta) & -\pi \delta & \alpha + \beta & -\beta \cr
	    -\pi c\delta & \pi \gamma & -c\beta & \alpha \cr} $$
 mit $\alpha$, $\beta$, $\gamma$, $\delta \in Q(R_p)$.}
 \smallskip
 Beweis. F\"ur den Erfahrenen ist es eine langweilige, f\"ur den
 Un\-ge\-\"ub\-ten eine n\"utzliche und dar\"uber hinaus eine gehaltvolle
 \"Ubung in Matrizenmultiplikation nachzuweisen, dass $H_p$ gerade
 aus den zuletzt genannten Matrizen besteht. Man sieht dann, dass
 man diese Matrizen mit ihrer ersten Zeile identifizieren kann.
 Sind $(\alpha, \beta, \gamma, \delta)$ und $(\alpha', \beta',
 \gamma', \delta')$ zwei dieser Zeilen, so ergibt sich die erste
 Zeile des Produktes ihrer zugeh\"origen Matrizen zu
 $$\eqalign{ 
      \alpha'' &= \alpha \alpha' + c\beta \beta' + \pi \gamma \gamma' +
                         \pi \gamma \delta' - \pi c \delta \delta', \cr
       \beta'' &= \alpha \beta' + \alpha' \beta + \beta \beta' -
			 \pi \gamma \delta' + \pi \gamma' \delta,   \cr
      \gamma'' &= \alpha \gamma' + c\beta \delta' + \gamma \alpha' +
			 \gamma \beta' - c \delta \beta',           \cr
      \delta'' &= \alpha \delta' + \beta \gamma' + \beta \delta' -
		         \gamma \beta' + \delta \alpha'. \cr} $$
 Die Elemente $(\zeta,0,0,0)$ liegen nat\"urlich im Zentrum von
 $H_p$, da sie ja Diagonalmatrizen darstellen, deren
 Hauptdiagonaleintr\"age alle gleich $\zeta$ sind. Es sei
 $(\zeta_1, \zeta_2, \zeta_3, \zeta_4)$ ein Element des Zentrums.
 Dann ist
 $$ (\zeta_1,\zeta_2,\zeta_3,\zeta_4)(0,1,0,0) =
        (c\zeta_2,\zeta_1 + \zeta_2,\zeta_3 - c\zeta_4,- \zeta_3) $$
  und
 $$(0,1,0,0)(\zeta_1,\zeta_2,\zeta_3,\zeta_4) =
        (c\zeta_2,\zeta_1 + \zeta_2,c\zeta_4,\zeta_3 + \zeta_4).$$
 Koeffizientenvergleich liefert $\zeta_3 = 2c\zeta_4$ und
 $\zeta_4 = -2\zeta_3$. Hieraus folgt $\zeta_3 = -4c\zeta_3$.
 W\"are $\zeta_3$ nicht Null, so folgte $-c = {1 \over 4}$ und
 damit
 $$ x^2 + x - c = \biggl(x + {1 \over 2}\biggr)^2, $$
 so dass dieses Polynom \"uber $R_p/\pi R_p$ nicht irreduzibel w\"are. Also
 ist $\zeta_3 = 0$ und damit auch $\zeta_4 = 0$.
 \par
       Weiter folgt
 $$ (\zeta_1,\zeta_2,0,0)(0,0,1,0)=(0,0,\zeta_1,\zeta_2) $$
 und
 $$(0,0,1,0)(\zeta_1,\zeta_2,0,0) = (0,0,\zeta_1 + \zeta_2,-\zeta_2), $$
 so dass auch $\zeta_2 = 0$ ist. Also besteht $Z(H_p)$ aus genau den Elementen
 $(\zeta,0,0,0)$.
 \par
       Es ist
 $$
    (\alpha, \beta, \gamma, \delta)(\alpha + \beta, -\beta, -\gamma, -\delta) 
  = \bigl(\alpha^2 + \alpha \beta - c\beta^2
       - \pi(\gamma^2 + \gamma \delta - c\delta^2), 0, 0, 0\bigr). 
 $$
 Hieraus folgt, dass die zu $(\alpha,\beta,\gamma,\delta)$ geh\"orende Matrix
 sicher dann invertierbar ist, wenn
 $$(\alpha^2 + \alpha \beta - c\beta^2 - \pi(\gamma^2 +
        \gamma \delta - c\delta^2) \neq 0$$
 ist. Es sei also
 $$ (\alpha^2 + \alpha \beta - c\beta^2 - \pi(\gamma^2 + \gamma\delta
       - c\delta^2) = 0. $$
 Indem wir diese Gleichung mit dem Quadrat des Hauptnenners der in sie
 eingehenden Koeffizienten multiplizieren, sehen wir, dass wir
 annehmen d\"urfen, dass sie in $R_p$ liegen. Weil $R_p$ ein
 Hauptidealring ist, d\"urfen wir weiter annehmen, dass sie
 teilerfremd sind, falls sie nicht alle Null sind. Nun ist
 $$ \alpha^2 + \alpha \beta - c\beta^2 \equiv 0 \mod \pi. $$
 Wegen der vorausgesetzten Irreduzibilit\"at von $x^2 + x -c$ folgt hieraus,
 dass sowohl $\alpha$ also auch $\beta$ durch
 $\pi$ teilbar ist. Dann ist aber
 $$ \alpha^2 + \alpha \beta - c\beta^2 \equiv 0  \mod \pi^2, $$
 so dass
 $$ \gamma^2 + \gamma \delta - c\delta^2 \equiv 0 \mod \pi, $$
 ist. Hieraus folgt, dass auch $\gamma$ und $\delta$ durch
 $\pi$ teilbar sind, so dass die Teilerfremdheit von $\alpha$,
 $\beta$, $\gamma$ und $\delta$ nicht zu erreichen ist. Also sind alle
 diese Koeffizienten Null, so dass $H_p$ in der Tat ein K\"orper
 ist. Da offensichtlich $[H_p:Z(H_p)] = 4$ ist, ist $H_p$ ein
 Quaternionenschiefk\"orper. Damit ist alles bewiesen.
 \medskip
       In dem gerade bewiesenen Satz steht eine Voraussetzung, deren
 Realisierbarkeit noch nicht gekl\"art ist, n\"amlich die
 Voraussetzung, dass $c$ ein Element sei, f\"ur das $x^2 + x - c$
 \"uber $R_p/\pi R_p$ irreduzibel ist. Um die anstehende Frage zu
 beantworten, betrachten wir einen K\"orper $K$ der Charakteristik
 $p>0$, von dem wir annehmen, dass er algebraisch \"uber seinem
 Primk\"orper sei. Wir setzen
 $$ E(K) := \bigl\{n \mid n \in N, \GF(p^n) \subseteq K\bigr\}. $$
 Dabei ist $\GF(p^n) \subseteq K$ so zu
 verstehen, dass $K$ einen zu $\GF(p^n)$ isomorphen Teilk\"orper
 enthalte. Dann hat $E(K)$ die folgenden Eigenschaften:
 \smallskip\noindent
 a) Es ist $1 \in E(K)$.
 \smallskip\noindent
 b) Ist $n \in E(K)$ und ist $m$ ein Teiler von $n$, so ist $m \in E(K)$.
 \smallskip\noindent
 c) Sind $m, n \in E(K)$, so ist $\kgV (m,n) \in E(K).$
 \smallskip
       Ist umgekehrt $M$ eine Menge von nat\"urlichen Zahlen
 mit den Eigenschaften a), b) und c), ist ferner $p$ eine Primzahl,
 so gibt es bis auf Isomorphie genau eine algebraische Erweiterung
 $K$ von $\GF(p)$ mit $E(K) = M$.
 \par
       Die Mengen $M$, die die Eigenschaften a), b) und c) haben, kann
 man auch folgenderma\ss{}en beschreiben. Ist $\Pi$ die Menge aller
 Prim\-zah\-len, so definieren wir durch
 $$ \alpha(q) := \max\{e \mid e \in \omega,\ \hbox{\rm es gibt ein}\ n \in M\ 
     \hbox{\rm mit}\ n \equiv 0 \mod q^e\} $$
 eine Abbildung von $\Pi$ in $\omega \cup \{\infty\}$. Ist dann $n \in \omega*$,
 so ist genau dann $n \in M$, wenn f\"ur alle $q \in \Pi$ gilt: ist $q^e$ ein
 Teiler von $n$, so ist $e \leq \alpha(q)$. Zu jeder algebraischen
 Erweiterung $K$ von $\GF(p)$ geh\"ort also eine solche Abbildung,
 die wir mit $\alpha_K$ bezeichnen. Ist $\alpha_K = \alpha_L$, so
 sind $K$ und $L$ isomorph.
 \medskip\noindent
 {\bf 9.18. Satz.} {\it Es sei $K$ eine algebraische Erweiterung
 von $\GF(p)$. Ist $\alpha_K(2) < \infty$, so gibt es ein $\gamma \in K$, so dass
 das Polynom $x^2 + x - \gamma$ irreduzibel ist.}
 \smallskip
       Beweis. Setze $n := 2^{\alpha_K (2)}$. Mit b) und der Definition
 von $\alpha_K$ folgt, dass $K$ einen zu $\GF(p^n)$ isomorphen
 Teilk\"orper $L$ enth\"alt. Wir betrachten die durch
 $$ f(\lambda) := \lambda^2 + \lambda $$
 definierte Abbildung $f$ von
 $L$ in sich. Wegen $f(0) = f(-1) = 0$ ist $f$ nicht injektiv und
 folglich wegen der Endlichkeit von $L$ auch nicht surjektiv. Es
 gibt also ein $\gamma \in L$ mit
 $$\lambda^2 + \lambda - \gamma \neq 0 $$
 f\"ur alle $\lambda \in L$. Daher ist $x^2 + x - \gamma$ \"uber  $L$ irreduzibel.
 H\"atte dieses Polynom eine Nullstelle $\kappa \in K$, so w\"are $L[\kappa]$ ein
 zu $\GF(p^{2n})$ isomorpher Teilk\"orper von $K$ und es folgte der
 Widerspruch
 $$ 2n = 2^{\alpha_K(2)+1} \leq 2^{\alpha_K(2)} = n. $$
 Also ist $x^2 + x - \gamma$ auch \"uber $K$ irreduzibel.
 \medskip
       Die S\"atze 9.18 und 9.17 zeigen, dass es
 Quaternionenschiefk\"orper gibt, deren Zentrum zu $Q(Z_p)$, bzw.
 zu $K((x))$, dem K\"orper der formalen Laurentreihen \"uber $K$,
 isomorph ist, wenn nur $K$ eine algebraische Erweiterung eines
 endlichen Primk\"orpers ist mit $\alpha_K(2) < \infty$. Sind die
 Zentren zweier solcher Quaternionenschiefk\"orper nicht isomorph,
 so sind auch die fraglichen Quaternionenschiefk\"orper nicht
 isomorph, so dass die fr\"uheren S\"atze zeigen, dass es
 Quaternionenschiefk\"orper wie Sand am Meer gibt. Ist $H_p$ einer
 dieser K\"orper, so ist die durch
 $$ \alpha(a,b,c,d) := (a + b,-b,-c,-d) $$
 definierte Abbildung $\alpha$ ein Antiautomorphismus von $H_p$.
 \par
       Um zu sehen, dass es sehr viele Quaternionenschiefk\"orper gibt,
 h\"atte man sich nicht so anstrengen m\"ussen. Hierauf machte mich
 T.\ Grundh\"ofer aufmerksam, von dem ich folgende Konstruktion
 lernte. Es sei $Q$ ein Quaternionenschiefk\"orper --- den muss
 man haben ---. Setze $K := Z(Q)$. Schlie\ss lich sei $L$ ein
 kommutativer K\"orper, der $K$, jedoch keine quadratische
 Erweiterung von $K$ enthalte. Dann ist das Tensorprodukt $Q
 \otimes_K L$ ein Quaternionenschiefk\"orper mit Zentrum $L$. Nimmt
 man nun f\"ur $L$ rein transzendente Erweiterungen von $K$, so
 erh\"alt man sogar Quaternionenschiefk\"orper beliebiger
 M\"achtigkeit.
 \par
 Es gibt noch eine dritte interessante M\"oglichkeit,
 Quaternionenschiefk\"orper zu konstruieren. Ist n\"amlich $K$ ein
 angeordneter K\"orper, so ist der Zentralisator der beiden
 Matrizen
 $$\pmatrix{0 & 1 & 0 & 0  \cr
	   -1 & 0 & 0 & 0  \cr
	    0 & 0 & 0 & -1 \cr
	    0 & 0 & 1 & 0  \cr}\quad
 \mathrm{und}\quad
 \pmatrix{0 & 0 & 1 & 0 \cr
	    0 & 0 & 0 & 1 \cr
	   -1 & 0 & 0 & 0 \cr
	    0 & -1 & 0 & 0 \cr} $$
 im Ring der der $(4 \times 4)$-Matrizen \"uber $K$ ein
 Quaternionenschiefk\"or\-per. Dies zu beweisen, sei dem Leser als
 \"Ubungsaufgabe \"uberlassen.

\mysectionten{10. Projektive R\"aume mit Clifford-Parallelismus}

\noindent
 William Kingdon Clifford 
 ist der Entdecker\index{Clifford, W. K.}{} einer bemerkenswerten Eigenschaft des
 dreidimensionalen projektiven Raumes \"uber dem K\"or\-per der
 reellen Zahlen, n\"amlich der Eigenschaft dieses Raumes, zwei
 Par\-al\-le\-li\-t\"ats\-re\-la\-ti\-o\-nen auf der Menge der Geraden zu
 besitzen,
 die beide das euklidische Parallelenpostulat erf\"ullen und
 dar\"uber hinaus eng miteinander verkn\"upft sind. Ohne weitere
 Motivation fangen wir nun an, wie in der Mathematik \"ublich, das
 Pferd beim Schwanz auf\-zu\-z\"au\-men.\index{Clifford-Parallelismus}{}
 \par
       Es sei $G$ eine Gruppe und $\pi$ sei eine Partition von $G$. Aus
 diesen Daten machen wir wie folgt eine Inzidenzstruktur $\pi(G)$.
 Die Punkte von $\pi(G)$ sind die Elemente von $G$ und die Geraden
 die Rechts\-rest\-klas\-sen $X_g$ mit $X \in \pi$ und $g \in G$. Als
 Inzidenz nehmen wir die Relation $\in$. Sind dann $g$ und $h$ zwei
 verschiedene Elemente von $G$, so ist $gh^{-1} \neq 1$, so dass es
 genau ein $X \in \pi$ gibt mit $gh^{-1} \in X$. Es folgt $g \in
 Xh$. Da auch $h \in Xh$ gilt, ist $Xh$ eine Gerade durch die
 beiden Punkte $g$ und $h$. Sie ist aber auch die einzige Gerade,
 die mit diesen beiden Punkten inzidiert. Sind n\"amlich $g$, $h \in
 Yf$, so gilt $gf^{-1}$, $hf^{-1} \in Y$ und daher $gh^{-1} =
 gf^{-1} (hf^{-1}) \in Y$, so dass $X = Y$ ist. dann ist aber $Xh
 \cap Xf \neq 0$ und folglich $Xh = Xf$. Man beachte, dass die
 durch
 $$ x^\gamma(g) := xg $$
 f\"ur alle $x$, $g \in G$ definierte Abbildung $\gamma$ ein Monomorphismus von
 $G$
 in die Kollineationsgruppe von $\pi (G)$ ist.
 \par
       Die Partition $\pi$ der Gruppe $G$ hei\ss e {\it normal\/}, falls
 f\"ur alle $X \in \pi$ und alle $g \in G$ gilt, dass auch $g^{-1}
 Xg \in \pi$ ist. Ist $\pi$ eine normale Partition von $G$ und ist
 $Xg$ eine Gerade von $\pi(G)$, so ist $Xg = g(g^{-1} Xg)$, so
 dass in diesem Falle die Menge der Geraden von $\pi(G)$ auch
 gleich der Menge der Linksrestklassen nach den Komponenten von
 $\pi$ ist.
 \par
       Der Satz von Dandelin,\index{Satz von Dandelin}{} den wir im ersten
 Kapitel bewiesen haben, ist nicht nur sch\"on, er ist auch n\"utzlich, wie wir
 sehen werden.
 \medskip\noindent
 {\bf 10.1. Satz.} {\it Es sei $G$ eine Gruppe und $\pi$ sei eine normale, nicht
 tri\-vi\-a\-le Partition von $G$. Ist $\pi(G)$ eine irreduzible projektive
 Geometrie, so ist $\pi(G)$ pappossch.}\index{papposscher Raum}{}
 \smallskip
       Beweis. Zun\"achst beachten wir, dass es in $\pi(G)$ windschiefe
 Geraden gibt. Ist n\"amlich $X \in \pi$, so ist $X \neq G$, da
 $\pi$ nicht trivial ist. Es gibt also ein $g \in G - X$. Dann sind
 aber $X$ und $Xg$ zwei windschiefe Geraden. Hieraus folgt, dass
 der Rang von $\pi(G)$ mindestens gleich $4$ ist.
 \par
       Sind $X$ und $Y$ zwei verschiedene Komponenten von $\pi$, sind
 $g_1$, $g_2$, $g_3$ drei verschiedene Elemente aus $Y - \{1\}$ und
 $h_1$, $h_2$, $h_3$ drei verschiedene Elemente in $X - \{1\}$, so
 bilden die Geraden $Xg_1$, $Xg_2$, $Xg_3$ und $h_1Y$, $h_2Y$, $h_3Y$ ein
 {\it Hexagramme mystique\/}. Um dies einzusehen, sei zun\"achst
 $Xg_i = Xg_k$. Dann ist
 $$ g_i {g_k}\!^{-1} \in X \cap Y = \{1\} $$
 und daher $i = k$, so dass die Geraden $Xg_1$, $Xg_2$ und
 $Xg_3$ paarweise verschieden sind. Da sie verschieden sind, sind
 sie auch paarweise windschief. Analog sieht man dass die Geraden
 $h_1Y$, $h_2Y$, $h_3Y$ --- dies sind Geraden, da $\pi$ normal ist ---,
 paarweise windschief sind. Wegen
 $$ h_i g_k \in Xg_k \cap h_iY $$
 bilden sie also ein {\it Hexagramme mystique\/}.
 \par
       Es seien nun $X$ und $Y$ zwei verschiedene Geraden von $\pi(G)$,
 die sich im Punkte $s$ schneiden. Dann sind $Xs^{-1}$ und
 $Ys^{-1}$ zwei Geraden, die sich im Punkte $1$ schneiden, so dass
 wir des weiteren annehmen d\"urfen, dass $s=1$ ist. Dann sind aber
 $X$ und $Y$ zwei verschiedene Komponenten von $\pi$. Sind nun
 $g_1$, $g_2$ und $g_3$ drei verschiedene Punkte aus $Y - \{1\}$ und
 $h_1$, $h_2$ und $h_3$ drei verschiedene Punkte aus $X - \{1\}$, so
 ist $Xg_1$, $Xg_2$, $Xg_3$, $h_1Y$, $h_2Y$, $h_3Y$ ein {\it Hexagramme
 mystique\/}, wie wir gerade gesehen haben. Nach dem Satz I.1.10 von
 Dandelin sind daher die Punkte
 $$\eqalign{
    (g_1 + h_2) &\cap (g_2 + h_1) \cr
    (g_2 + h_3) &\cap (g_3 + h_2) \cr
    (g_3 + h_1) &\cap (g_1 + h_3) \cr}$$
 kollinear, so dass in $\pi(G)$, wie behauptet, der Satz von Pappos gilt.
 \medskip
       Bevor wir uns auf die Suche nach Gruppen $G$ begeben, die eine
 normale Partition $\pi$ besitzen, so dass $\pi(G)$ eine
 irreduzible projektive Geometrie ist, beweisen wir einen Satz von
 Karzel,\index{Karzel, H.}{} der die Suche erheblich erleichtert. Der hier
 vorgetragene Beweis stammt von mir (Karzel 1965, L\"uneburg 1969a).
 \medskip\noindent
 {\bf 10.2. Satz.} {\it Es sei $G$ eine Gruppe und $\pi$ sei eine
 normale, nicht triviale Partition von $G$. Ist $\pi(G)$ eine
 irreduzible projektive Geometrie, so ist $G$ eine
 elementarabelsche $2$-Gruppe oder es ist $\Rg(\pi(G)) = 4$.}
 \smallskip
       Beweis. Wir definieren die Bijektion $\rho$ von $G$ auf sich durch
 $x^\rho := x^{-1}$. Ist dann $X_g$ eine Gerade von $\pi(G)$, so
 ist $(X_g)^\rho = g^{-1} X$ wegen der Normalit\"at von $\pi$
 ebenfalls eine Gerade von  $\pi(G)$, so dass $\rho$ eine
 Kollineation von $\pi(G)$ ist. Ist $\rho = 1$, so ist $G$ eine
 elementarabelsche $2$-Gruppe. Es sei also $\rho \neq 1$. Dann ist
 $\rho$ eine involutorische Kollineation von $\pi(G)$, die alle
 Geraden durch den Punkt 1 festl\"asst, so dass $\rho$ eine
 Zentralkollineation mit dem Zentrum 1 ist. Nach Satz 1.2 besitzt
 $\rho$ daher eine Achse $H$. Ist $x \in H$, so gilt offenbar $x^2 = 1$. Ist
 andererseits $y \neq 1 = y^2$, so ist $y$ ein von 1
 verschiedener Fixpunkt von $\rho$ und daher $y \in H$.
 \par
       Es sei $h \in H$ und $I$ sei die Menge der mit $h$ vertauschbaren
 Elemente aus $H$. Wir zeigen, dass
 $$ I = (H \cap Hh) \cup \{h\} $$
 ist. Es sei $j \in I$ und es gelte $j \neq h$. Dann ist $jh \neq 1$ und
 $$ (jh)^2 = j^2 h^2 = 1. $$
 Daher ist nach obiger Bemerkung $jh \in H$. Ferner folgt
 $$ j = (jh)h \in H \cap Hh. $$
 Also ist $I \subseteq (H \cap Hh) \cup \{h\}$. Ist
 umgekehrt $k \in (H \cap Hh) \cup \{h\}$ und ist $k \neq h$, so
 gibt es ein $x \in H$ mit $k = xh$. Es folgt $kh = x$ und weiter
 $$ 1 = x^2 = (kh)^2, $$
 so dass $kh = hk$ ist. Also ist $k \in I$ und obige Gleichung damit bewiesen.
 \par
       Wir betrachten nun zus\"atzlich den Fall, dass $\rho$ eine Elation
 ist. In diesem Fall hat $\rho$ keinen Fixpunkt au\ss erhalb $H$.
 Es gilt also
 $$ H = \{x \mid x \in G, x^2 = 1\}. $$
 Man beachte, dass in diesem Falle jede Untergruppe von $G$, die in $H$
 enthalten ist, abelsch ist.
 \par
       Es sei $1 \neq h \in H$ und $C$ bezeichne den Zentralisator von
 $h$ in $G$, dh., die Menge aller mit $h$ vertauschbaren Elemente
 von $G$. Der Zentralisator $C$ von $h$ ist nat\"urlich eine
 Untergruppe von $G$. Wir zeigen, dass $C$ in $H$ enthalten ist. Es
 sei $U$ die Komponente von $\pi$, die $h$ enth\"alt. Weil $U$ die
 beiden Punkte 1 und $h$ enth\"alt, folgt $U \subseteq H$. Es sei
 $g \in C$ und es gelte $g \notin U$. Es sei $V$ die Komponente von
 $\pi$, die $hg$, und $W$ die Komponente, die $g$ enth\"alt. Dann
 ist $V \neq W$, da sonst $h \in U \cap V = \{1\}$ w\"are. Also ist
 $$ g^2 = h^2g^2 = (hg)^2 \in V \cap W = \{1\}, $$
 so dass $g \in H$ gilt. Dies zeigt, dass $C \subseteq H$ ist. Nach dem, was
 wir weiter oben gesehen haben, ist also $C = (H \cap Hh) \cup \{h\}$. Wegen
 $h = 1h$ ist aber $h \in H \cap Hh$, so dass
 $$ C = H \cap Hh $$
 gilt.
 \par
       Es sei $g \in G - H$. Dann ist $H \neq Hg$, da ja $g = 1g \in Hg$
 ist. W\"are $H$ eine Untergruppe von $G$, so folgte $H \cap Hg = \emptyset$ und
 damit, da $H$ und $Hg$ ja Hyperebenen sind, $\Rg(\pi(G)) = 2$ im Widerspruch zu
 $\Rg(\pi(G)) \geq 4$. Also ist $C \neq H$, so dass $C$ als Schnitt von zwei
 Hyperebenen den $\Ko$-Rang $2$ hat. Es sei $1 \neq d \in C$ und $D$ sei der
 Zentralisator von $d$ in $G$. Dann ist $C \subseteq D$, da $C$
 abelsch ist. Weil auch $D$ ein Unterraum des Ko-Ranges $2$ ist,
 folgt $C = D$. Es gibt nun ein $e \in H - C$. Es sei $E$ der
 Zentralisator dieses Elementes. Dann ist also $C \neq E$ und
 folglich $C \cap E = \{1\}$, wie wir gerade gesehen haben. Damit
 haben wir zwei Unterr\"aume des Ko-Ranges $2$ gefunden, die sich
 in genau einem Punkte schneiden. Hieraus folgt, dass der Rang von
 $\pi(G)$ gleich $4$ ist.
 \par
       Es bleibt der Fall zu betrachten, dass $\rho$ eine Streckung ist.
 Hier gilt
 $$ H = \{x \mid x \in G, x^2 = 1 \neq x\}. $$
 Es sei $h \in H$ und $I$ sei der Zentralisator von $h$ in $H$. Dann ist, wie
 wir gesehen haben
 $$ I = (H \cap Hh) \cup \{h\}. $$
 Ferner gilt $h \notin H \cap Hh$, da andernfalls $1 \in H$ w\"are. Der
 Rang von $H \cap Hh$ ist mindestens 2. Es gibt also ein $b \in H \cap Hh$. Es
 folgt $b = ah$ mit einem $a \in H$. Es folgt $h = ab$
 und damit wieder $ab = ba$. Wir betrachten den Schnitt $Ha \cap Hb$. Dann ist
 $$ 1 = aa= bb \in Ha \cap Hb, $$
 da $a$ und $b$ ja in $H$ liegen. Wegen $h = ab = ba$ ist auch
 $h \in Ha \cap Hb$. Es sei nun $d \in Ha \cap Hb$. Es folgt $da$, $db \in H$.
 Also ist $(da)^2 = 1 = (db)^2$. Hieraus folgt $ada = d^{-1} = bdb$. Und weiter
 $$ hdh = abdba = ad^{-1}a = d, $$
 so dass $d$ und damit $Ha \cap Hb$ im Zentralisator von $h$ in $G$ liegt. Es
 folgt
 $$ Ha \cap Hb \cap H \subseteq I = (H \cap Hh) \cup \{h\}. $$
 Nun ist $Ha \cap Hb \cap H$ ein Unterraum von $\pi(G)$,
 der $h$ enth\"alt. Weil $h$ nicht in $H \cap Hh$ liegt und jede
 Gerade, die $h$ mit einem Punkt von $H \cap Hh$ liegt, folgt
 $$ Ha \cap Hb \cap H = \{h\}. $$
 Dies impliziert wieder, dass der Rang von $\pi(G)$ gleich 4 ist. Damit ist der
 Satz bewiesen.
 \medskip
       Es sei $\pi$ eine normale Partition der Gruppe $G$ und $\pi(G)$ sei eine
 projektive Geometrie. Sind $X$ und $Y$ zwei Geraden von $\pi(G)$, so hei\ss en
 $X$ und $Y$ {\it rechtsparallel\/},\index{rechtsparallel}{} falls es ein
 $U \in \pi$ und $g$, $h \in G$ gibt mit $X = Ug$ und $Y = Uh$. Die Geraden $X$
 und $Y$ hei\ss en {\it linksparallel\/},\index{linksparallel}{} falls
 es ein $U \in \pi$ und $g$, $h \in G$ gibt mit $X = gU$ und $Y = hU$. Diese
 beiden Parallelit\"atsrelationen erf\"ullen
 offensichtlich das euklidische Parallelenpostulat. Sie sind
 Verallgemeinerungen der oben erw\"ahnten cliffordschen
 Parallelit\"atsrelationen.\index{Clifford-Parallelismus}{} Sie werden explizit
 nicht mehr auftauchen.
 \medskip\noindent
 {\bf 10.3. Satz.} {\it Es sei $\pi$ eine normale, nicht triviale
 Partition der Gruppe $G$. Ist $\pi(G)$ eine projektive Geometrie,
 so gilt:
 \item{a)} Ist $K$ eine Hyperebene von $\pi(G)$, so ist $G = K^{-1} K$.
 Insbesondere ist $G = HH$, wobei $H$ wieder Achse der oben definierten
 Zentralkollineation $\rho$ ist.
 \item{b)} Ist $X \in \pi$, so ist $X$ abelsch.
 \item{c)} $G$ ist eine elementarabelsche $2$-Gruppe oder es ist $Z(G) = \{1\}$.}
 \smallskip
       Beweis. a) Es sei $g \in G$. Da der Rang von $\pi(G)$ mindestens
 4 ist, ist der Schnitt der beiden Hyperebenen $K$ und $Kg$ nicht
 leer. Es gibt also $k$, $l \in K$ mit $k = lg$, so dass $g = l^{-1} k$ ist.
 \par
       b) Es sei $1 \neq g \in X$. Nach a) gibt es $u$, $v \in H$ mit $g = uv$.
 Es ist $v \in Xv$ und $u = uvv = gv$. Weil $g \neq 1$ ist,
 ist $u \neq v$, so dass die Gerade $Xv$ in $H$ liegt. Daher ist
 $(xv)^2 = 1$ und damit $vxv = x^{-1}$ f\"ur alle $x \in X$. Hieraus folgt
 $$ vxyv = (xy)^{-1} = y^{-1} x^{-1} = vyvvxv = vyxv $$
 und dann $xy = yx$ f\"ur alle $x$, $y \in X$, so dass $X$ abelsch ist.
 \par
       c) Es sei $u$ eine Involution von $G$ und $v$ sei ein mit $u$
 vertauschbares Element von $G$. Dann liegen $u$ und $v$ entweder
 in der gleichen Komponente von $\pi$ oder aber es ist $v^2 = 1$.
 Sind n\"amlich $V$, $W \in \pi$ und gilt $v \in V$ und $uv \in W$,
 so ist $V \neq W$ und daher
 $$ v^2 = (uv)^2 \in V \cap W = \{1\}. $$
 (Dieser Schluss wurde nun schon so h\"aufig gezogen,
 dass er es ei\-gent\-lich verdient h\"atte, als Hilfssatz formuliert zu werden.)
 \par
       Ist nun $1 \neq g \in Z(G)$, so gibt es eine Involution in $G$,
 die nicht in der gleichen Komponente wie $g$ liegt. Nach der
 Vorbemerkung ist also $g^2 = 1$, so dass $Z(g)$ eine
 elementarabelsche 2-Gruppe ist. W\"are nun $\rho$ nicht die
 Identit\"at, so l\"age $Z(G)-\{1\}$ in der Achse $H$ von $\rho$
 und die Elemente au\ss erhalb $H$ h\"atten nicht die Ordnung 2 im
 Widerspruch zu unserer Vorbemerkung.
 \medskip
       Wir ziehen nun einige interessante Folgerungen aus Satz 10.2.
 Zun\"achst eine geometrische, die wir durch den folgenden Satz
 vorbereiten.
 \medskip\noindent
 {\bf 10.4. Satz.} {\it Es sei $G$ eine Gruppe und $U$ 
 eine Teilmenge von $G$ mit $|U-\{1\}| \geq 2$. Gilt dann $U^{-1} x = U$
 f\"ur alle von $1$ verschiedenen $x \in U$, so ist $U$ eine Untergruppe von $G$.}
 \smallskip
       Beweis. Es sei $H := \{g \mid g \in G,\ Ug = U\}$. Dann ist $H$ eine
 Untergruppe von $G$. Es sei $u \in U$. Es gibt dann ein $v \in U$
 mit $v \neq 1$, $u$. Es folgt $U = U^{-1} v$. Es gibt daher ein $x
 \in U$ mit $u = x^{-1} v$. Wegen $u \neq v$ ist $x \neq 1$ und
 daher $U = U^{-1} x$. Somit ist $U^{-1} = Ux^{-1}$. Es folgt
 $$ U = U^{-1} v = Ux^{-1} v = Uu $$
 und weiter $u \in H$. Es ist also $U \subseteq H$. Weil $H$ eine Gruppe ist,
 folgt $U^{-1} \subseteq H$. Aus der Definition von $H$ folgt schlie\ss lich
 $UU^{-1} \subseteq U$. Dies zeigt, dass $U$ eine Untergruppe von
 $G$ ist.
 \medskip\noindent
 {\bf 10.5. Satz.} {\it Es sei $G$ eine Gruppe und die Elemente von
 $G$ seien gleichzeitig die Punkte einer irreduziblen projektiven
 Geometrie $\Lambda(G)$ mit $\Rg(\Lambda(G)) \geq 3$. F\"ur alle
 $g \in G$ sei die durch $x^{\kappa(g)} := xg$ definierte
 Abbildung $\kappa(g)$ von $G$ in sich eine Kollineation von
 $\Lambda(G)$. Ist dann auch die durch $x^\rho := x^{-1}$
 definierte Abbildung $\rho$ eine Kollineation von $\Lambda(G)$,
 so ist die Menge $\pi$ der Geraden durch den Punkt $1$ --- die
 Geraden als Punktmengen aufgefasst --- eine normale Partition
 von $G$, und $\pi(G)$ ist die Geometrie aus den Punkten und
 Geraden von $\Lambda (G)$. Insbesondere ist $G$ also eine
 elementarabelsche $2$-Gruppe oder aber es ist $\Rg(\Lambda (G)) = 4$.}
 \smallskip
       Beweis. Ist $U \in \pi$, so ist $|U - \{1\}| \geq 2$, da die Geraden einer
 irreduziblen projektiven Geometrie mindestens drei Punkte tragen. Ferner folgt,
 da $\rho$ eine Kollineation ist, dass auch $U^{-1} \in \pi$ gilt. Es sei $1$,
 $x \in U$. Dann ist 1, $x^{-1} \in U^{-1} x$ eine Gerade durch die Punkte $x$
 und 1. Da auch $U$ eine Gerade durch diese beiden Punkte ist, folgt
 $U^{-1}x = U$. Nach 10.4 ist $U$ daher eine Untergruppe von $G$, so dass
 $\pi$ eine Partition von $G$ ist, da jeder von 1 verschiedene
 Punkt von $\Lambda(G)$ auf genau einer Geraden durch 1 liegt.
 Diese Partition ist wegen $\Rg(\Lambda(G)) \geq 3$ nicht
 trivial. Wegen $gx = x^{\rho \kappa(g^{-1}) \rho}$ ist die
 Abbildung $x \rightarrow gx$ ebenfalls eine Kollineation von
 $\Lambda(G)$. Hieraus folgt $g^{-1} Ug \in \pi$ f\"ur alle $U \in
 \pi$. somit ist $\pi$ eine normale Partition von $G$. Weil die
 Gruppe der Kollineationen $\kappa(g)$ auf der Menge der Punkte
 von $\Lambda(G)$ transitiv operiert, folgt schlie\ss lich, dass
 $\pi(G)$ die Geometrie aus den Punkten und Geraden von $\Lambda(G)$ ist. Die
 restlichen Behauptungen folgen nun aus 10.2.
 
 \medskip

      Die weiteren Folgerungen sind algebraischer Natur, so dass die
 Geometrie\index{Geometrie}{} auch einmal anderen Teilen der
 Mathematik\index{Mathematik}{} dient, was
 selten genug vorkommt, wie schon zu Beginn von Abschnitt 4 des
 Kapitels I bemerkt wurde.
 \medskip\noindent
 {\bf 10.6. Satz.} {\it Es sei $M$ ein K\"orper und $K$ sei ein in $Z(M)$
 enthaltener Teilk\"orper von $M$. Hat jedes $a \in M - K$ ein Minimalpolynom des
 Grades $2$ \"uber $K$, so gilt eine der folgenden Aussagen:
 \item{a)} Es ist $K = M$.
 \item{b)} Es ist $[M:K] = 2$.
 \item{c)} Es ist $Z(M) = K$ und $M$ ist ein Quaternionenschiefk\"orper \"uber $K$.
 \item{d)} Der K\"orper $M$ ist kommutativ, die Charakteristik von
 $M$ ist $2$ und $M$ ist eine rein inseparable Erweiterung von $K$
 mit $a^2 \in K$ f\"ur alle $a \in M$.\par}
 \smallskip
       Beweis. Ist $[M:K] \leq 2$, so liegt einer der beiden F\"alle a)
 oder b) vor. Es sei also $[M:K] \geq 3$. Setze $G := M^*/K^*$. Da
 $K$ im Zentrum von $M$ liegt, ist $K^*$ ein Normalteiler von
 $M^*$, so dass $G$ eine Gruppe ist. Die Elemente von $G$ sind die
 Restklassen $mK^*$ mit $m \in M^*$, so dass die durch $\varphi(mK) := mK^*$
 definierte Abbildung $\varphi$ eine Bijektion der Menge der Punkte von
 $\La_K(M)$ auf $G$ ist.
 \par
       Es sei $X$ eine Gerade von $L_K(M)$ durch den Punkt $K$. Ferner
 sei $a \in X - K$. Dann ist $X = K + Ka$. Weil das Minimalpolynom
 von $a$ \"uber $K$ den Grad 2 hat, ist $X$ ein K\"orper. Ist
 umgekehrt $Y$ eine in $M$ liegende, quadratische Erweiterung von
 $K$, so ist $Y$ eine durch $K$ gehende Gerade von $L_K(M)$. Es
 folgt, dass die Menge  der Punkte einer Geraden durch $K$ unter
 $\varphi$ auf die Menge $Y^*/K^*$ abgebildet wird, wobei $Y$ der
 K\"orper ist, der die fragliche Gerade repr\"asentiert. Es folgt,
 dass die Menge $\pi$ der Gruppen $Y^*/K^*$, wobei $Y$ eine in $M$
 liegende quadratische Erweiterung von $K$ ist, eine Partition von
 $G$ ist. Die Partition $\pi$ ist aber auch normal, da innere
 Automorphismen von $M$ quadratische Erweiterungen von $K$ auf ebensolche
 abbilden. Sind nun $aK$ und $bK$ zwei Punkte von $\La_K(M)$, so
 wird durch $x^\sigma := xa^{-1}b$, da $K$ im Zentrum von $M$
 liegt, eine bijektive lineare Abbildung von $M$ auf sich definiert
 mit $(aK)^\sigma = bK$. Dies k\"onnen wir so interpretieren, dass
 $G$ auf der Menge der Punkte von $\La_K(M)$ transitiv operiert.
 Aus all dem folgt nun, dass $\pi(G)$ eine Beschreibung der
 Inzidenzstruktur aus den Punkten und Geraden von $L_K(M)$ ist.
 Mit Satz 10.2 erhalten wir daher, dass $G$ eine elementarabelsche
 2-Gruppe oder dass $\Rg(\pi(G)) = 4$ ist.
 \par 
       Es sei $G$ keine elementarabelsche 2-Gruppe. Dann ist $Z(G) = \{1\}$ nach
 10.3. Hieraus folgt $Z(M) = K$, so dass $M$ ein
 Quaternionenschiefk\"orper \"uber $K$ ist.
 \par
       Es sei $G$ eine elementarabelsche 2-Gruppe. Dann ist $a^2 \in K$
 f\"ur alle $a \in M$. Ist $a \in M - K$, so ist also
 $$ 2a = (a + 1)^2 - a^2 - 1 \in K. $$
 Dies impliziert, dass die Charakteristik von $M$ gleich 2 ist, da wegen
 $2 \in K$ sonst $a \in 2^{-1} K = K$ w\"are. Es seien $a$, $b \in M^*$. Weil $G$
 abelsch ist, ist $b^{-1}ab \in K[a]$. Es gibt also $k$, $l \in K$ mit
 $b^{-1} ab = k + la$. Hieraus folgt $aba^{-1} = kba^{-1} + lb$. Andererseits ist
 nat\"urlich auch $aba^{-1} \in K[b]$. Es folgt $kba^{-1} \in K[b]$. Ist
 $a \in K[b]$, so ist $ab = ba$. Es sei also $a \notin K[b]$. Dann ist $k = 0$
 und somit $b^{-1} ab= la$. hieraus folgt, da $a^2 \in K$ gilt,
 $$ a^2 = b^{-1}a^2b = l^2 a^2. $$
 Da die Charakteristik 2 ist, impliziert dies $l=1$. Also gilt auch in diesem
 Falle $ab = ba$, so dass $M$ als kommutativ erkannt ist. Damit ist alles
 bewiesen.
 \medskip
       Man beachte, dass im Fall d) der Rang von $M$ \"uber $K$ eine
 Potenz von 2 ist, falls er endlich ist.
 \par
       Der Beweis von 10.6 zeigt, wie man sich Gruppen $G$ mit einer
 Partition $\pi$ verschaffen kann, so dass $\pi(G)$ eine projektive
 Geometrie ist. Eine n\"ahere Analyse zeigt, dass dies bereits alle
 Beispiele sind. Diese Analyse werden wir jedoch nicht
 durchf\"uhren. Der an ihr interessierte Leser sei an Karzel 1965 verwiesen.
 \par
       Wir zeigen eine weitere Folgerung aus den bisherigen
 Entwicklungen, die sich auf Quaternionenschiefk\"orper der
 Charakteristik 2 bezieht. Sie ist Spezialfall eines viel
 allgemeineren Satzes. Da wir sie hier jedoch gratis bekommen, ---
 {\it for love\/}, wie es im Englischen so poetisch hei\ss t, ---
 wollen wir sie hier formulieren.
 \medskip\noindent
 {\bf 10.7. Korollar.} {\it Es sei $M$ ein
 Quaternionenschiefk\"orper \"uber $K$ und die Charakteristik von
 $M$ sei $2$. Setze
 $$ H := \{a \mid a \in M, a^2 \in K\}. $$
 Dann ist $H$ eine Hyperebene des $K$-Vektorraumes $M$.
 \item{a)} F\"ur $a \in H - K$ ist $K[a]$ eine inseparable Erweiterung des
 Grades $2$ von $K$.
 \item{b)} F\"ur $a \in M - H$ ist $K[a]$ eine separable Erweiterung
 des Grades $2$ von $K$.
 \par\noindent
       $M$ enth\"alt also sowohl inseparable als
 auch separable Erweiterungen des Grades $2$ von $K$.}
 \smallskip
       Beweis. Setze $G := M^*/K^*$ und $\pi$ sei die im Beweis von 10.6
 beschriebene Partition von $G$. Weil $M$ nicht kommutativ ist, ist
 $G$ nach 10.6 keine elementarabelsche 2-Gruppe. Also ist die
 $g^\rho := g^{-1}$ definierte Abbildung $\rho$ eine involutorische
 Perspektivit\"at von $\pi(G)$. Weil die Charakteristik von $K$
 gleich 2 ist, k\"onnen die Streckungsgruppen von $\pi(G)$, da sie
 zu $K^*$ isomorph sind, keine Involutionen enthalten. Also ist
 $\rho$ eine Elation. Daher gilt f\"ur die Achse $A$ von $\rho$ die
 Gleichung
 $$ A = \{g \mid g \in G,\ g^2 = 1\}. $$
 Der im Satz definierte Unterraum $H$ des $K$-Vektorraumes $M$ ist offenbar der
 $A$ entsprechende Unterraum von $M$.
 \par
       Ist $a \in H - K$, so ist $a^2 \in K$, so dass $K[a]$ eine
 inseparable Erweiterung von $K$ ist. Ist $b \in M - H$, so ist
 $b^2 \notin K$, so dass $K[b]$ eine separable Erweiterung von $K$
 ist. Damit ist alles bewiesen.
 \medskip
       Es sei $K$ ein angeordneter, kommutativer K\"orper. Wir nennen $K$
 {\it reell ab\-ge\-schlos\-sen\/},\index{reell abgeschlossen}{} falls jedes
 Polynom ungeraden Grades
 \"uber $K$ eine Nullstelle in $K$ hat. Prominentestes Beispiel
 eines reell ab\-ge\-schlos\-se\-nen K\"orpers ist der K\"orper der reellen
 Zahlen. Der zweite Beweis, den Gau\ss\ f\"ur den Fundamentalsatz
 der Algebra gab, l\"asst sich auf reell abgeschlossene K\"orper
 \"ubertragen, so dass also
 $$ K[x]/(x^2 + 1)K[x] $$
 algebraisch abgeschlossen ist, falls nur $K$ reell abgeschlossen
 ist. In einem reell abgeschlossenen K\"orper ist jedes positive
 Element ein Quadrat, so dass reell abgeschlossene K\"orper nur
 eine Anordnung besitzen.
 \par
       Ferner haben alle \"uber $K$ irreduziblen Polynome den Grad 1 oder
 2. Dies ist der Schl\"ussel zum Beweise des folgenden Satzes, der
 f\"ur den Fall des K\"orpers der reellen Zahlen von
 Frobenius\index{Frobenius, G. F.}{} stammt (Frobenius 1877).
 \medskip\noindent
 {\bf 10.8. Satz.} {\it Es sei $M$ ein K\"orper und $K$ sei ein
 reell abgeschlossener K\"orper mit $K \subseteq Z(M)$. Ist dann
 $[M:K] < \infty$, so gilt einer der folgenden F\"alle:
 \item{a)} Es ist $M = K$.
 \item{b)} $M$ ist der algebraische Abschluss von $K$.
 \item{c)} $M$ ist der bis auf Isomorphie eindeutig bestimmte
 Quaternionenschiefk\"orper \"uber $K$.\par}
 \smallskip
       Beweis. Da ein angeordneter K\"orper die Charakteristik 0 hat,
 folgt mit 10.6, dass $[M:K] = 1$, 2 oder 4 ist, da irreduzible
 Polynome \"uber $K$ ja den Grad 1 oder 2 haben. Im ersten Fall
 ist $M = K$. Im zweiten Fall folgt aus dem oben Referierten, dass
 $M$ der algebraische Abschlu\ss\ von $K$ ist. Es gelte also $[M:K]
 = 4$. Dann ist $M$ ein Quaternionenschiefk\"orper \"uber $K$. Es
 bleibt zu zeigen, dass $M$ bis aus Isomorphie eindeutig bestimmt
 ist.
 \par
       Es sei $H$ die Vereinigung \"uber alle $aK$ mit $a \notin K$ und
 $a^2 \in K$. Die schon verschiedentlich benutzte Perspektivit\"at
 $\rho$ zeigt, dass $H$ eine Hyperebene des $K$-Vektorraumes $M$
 ist. Ist $a \notin K$ und $a^2 \in K$, so ist das Polynom
 $x^2 - a^2$ \"uber $K$ irreduzibel. Weil $K$ reell abgeschlossen
 ist, folgt hieraus
 $$ a^2 = -k^2 $$
 mit einem positiven $k\in K$. Setzt man $b := k^{-1}a$, so folgt $bK = aK$ und
 $b^2 = -1$. Hieraus folgt
 $$ H = \bigcup_{b \in M,\ b^2 + 1 = 0} bK. $$
 (An dieser Stelle sieht man sehr sch\"on, dass es die Kommutativit\"at ist, die
 erzwingt, dass ein Polynom \"uber einem kommutativen K\"orper h\"ochstens so
 viele Nullstellen\index{Nullstellen eines Polynoms}{} besitzt,
 wie der Grad des Polynoms angibt, da im vorliegenden Falle das
 Polynom $x^2 + 1$ ja unendlich viele Nullstellen hat.)
 \par
       Es sei $i \in H$ und es gelte $b^2 = -1$. Ist $ib = -bi$, so setzen
 wir $j := b$. Es sei $ib \neq -bi$. Wegen
 $$ ib + bi = (i + b)^2 - i^2 - b^2 $$
 ist $ib + bi \in K$. Wir machen nun nach dem Vorgang von
 Fibonacci\index{Fibonacci}{} den falschen Ansatz\index{falscher Ansatz}{}
 $\lambda := ib + bi$ und $\mu := 2$, der uns dann doch zum rechten Ziele
 f\"uhrt. Mit diesem Ansatz folgt
 $$ i(\lambda i + \mu b) + (\lambda i + \mu b)i = -2\lambda + \mu(ib + bi) = 0.$$
 Weil $\lambda$ und $\mu$ nicht Null und $i$ und $b$ linear unabh\"angig sind,
 ist $(\lambda i + \mu b)K$ ein von $iK$ verschiedener Punkt in $H$. Wie wir oben
 gesehen haben, gibt es ein $j \in (\lambda i + \mu b)K$ mit $j^2 = -1$. Es
 folgt, dass $ij = -ji$ ist. Es folgt weiter
 $$ (ij)^2 = ijij = -jiij = (-1)^3 = -1 $$
 und damit $j \in K[i] \cap H = iK$. Dieser Widerspruch zeigt, dass $i$, $j$,
 $ij$ eine Basis von $H$ ist. Folglich ist $1$, $i$, $j$, $ij$ eine Basis von
 $M$. Hieraus folgt alles Weitere.
 \medskip
       Der Beweis des letzten Satzes zeigt auch, dass es \"uber einem
 angeordneten K\"orper, dessen positive Elemente allesamt Quadrate
 sind, bis auf Isomorphie nur einen Quaternionenschiefk\"orper gibt.
 \par
       Nach Satz 8.1 ist $L(V)$ genau dann selbstdual, wenn der Rang
 von $V$ endlich ist und wenn au\ss erdem der $L(V)$ zu Grunde
 liegende Koordinatenk\"orper einen Antiautomorphismus gestattet.
 Daher ist es interessant, auch nicht kommutative K\"orper zu
 kennen, die einen Antiautomorphismus  besitzen. Es ist nun nicht
 schwer zu sehen, dass alle von uns konstruierten
 Quaternionenschiefk\"orper einen involutorischen
 Antiautomorphismus haben. Wir werden zeigen, dass sogar alle
 Quaternionenschiefk\"orper einen solchen besitzen. Dazu beweisen wir zun\"achst:
 \medskip\noindent
 {\bf 10.9. Satz.} {\it Ist $M$ ein Quaternionenschiefk\"orper
 \"uber $K$, so gibt es Elemente $i$, $j$, $k \in M$ und Elemente
 $\gamma$, $\delta \in K$, so dass $1$, $i$, $j$, $k$ eine Basis des
 $K$-Vektorraumes $M$ ist und die folgenden Relationen gelten:}
 $$\eqalign{
     i^2 &= i + \gamma,       \cr
      ij &= k,                \cr
      ik &= k + \gamma j,     \cr
      ji &= j - k,            \cr
     j^2 &= \delta,           \cr
      jk &= \delta - \delta i,\cr
      ki &= -\gamma j,        \cr
      kj &= \delta i,         \cr
     k^2 &= -\gamma \delta. \cr}$$
 \smallskip
       Beweis. Wir beginnen den Beweis mit der folgenden Bemerkung. Sind
 $u$, $v \in M$, so ist $u^2 \in Ku + K$ und $v^2 \in Kv + K$ sowie
 $(u + v)^2 \in K(u + v) + K$. Hieraus folgt, dass
 $$ uv+vu \in Ku + Kv + K $$
 ist.
 \par
       W\"are $x^2 \in K$ f\"ur alle $x \in M$, so w\"are $M^*/K^*$ eine
 elementarabelsche 2-Gruppe. Wie der Beweis von 10.6 zeigt, w\"are
 $M$ ein kommutativer K\"orper im Widerspruch zu unserer
 Voraussetzung. Es gibt also ein $y \in M$ mit $y^2 \notin K$. Es
 gibt somit Elemente $r$, $s \in K$ mit $r \neq 0$ und $y^2 = ry = ry + s$. Setze
 $i := yr^{-1}$ und $\gamma := sr^{-2}$. Dann ist
 $$ i^2 = y^2r^{-1} = (ry + s)r^{-2} = i + \gamma. $$
 \par
       $Ki + K$ ist ein kommutativer Teilk\"orper von $M$. Es gibt
 folglich ein $z \in M$, welches nicht in $Ki + K$ liegt. Nach der
 eingangs gemachten Bemerkung gibt es $r$, $s$, $t \in K$ mit
 $$ iz + zi = ri + sz + t. $$
 Es gibt ferner ein $u \in K$ mit $z^2 - uz \in K$. Sind $l$, $m \in K$, so folgt
 $$\eqalign{ -(li + mz)^2 &+ (l + mr)li + (mu + ls) mz    \cr
	 &= -l^2(i^2 - i) - m^2(z^2 - uz) - ml(iz + zi - ri - sz). \cr} $$
 Dies zeigt, dass
 $$ -(li + mz)^2 + (l + mr)li + (mu + ls)mz \in K $$
 ist. Es gibt auch ein $v \in K$ mit
 $$ (li + mz)^2 - v(li + mz) \in K.$$
 Aus diesen beiden Gleichungen folgt
 $$(l + mr - v)li + (mu + ls - v)mz \in K. $$
 Weil $1$, $i$, $y$ \"uber $K$ linear unabh\"angig sind, folgt
 $$(l + mr - v)l = 0 = (mu + ls - v)m. $$
 Letztere Gleichung gilt f\"ur alle $l$, $m \in K$, da
 $l$ und $m$ ja irgendwelche Elemente von $K$ waren. W\"ahlt man
 insbesondere $l$ und $m$ in $K^*$, so folgt
 $$ l + mr = v = mu + ls. $$
 Dies hat wiederum
 $$ l(1 - s)= m(r - u) $$
 zur Folge. Mit $l = m = 1$ folgt
 $$ 1 - s = r - u. $$
 Weil $K$ das Zentrum von $M$ ist, ist $K$ nicht endlich. Es gibt also ein
 $l \in K^*-\{1\}$. Mit $m = 1$ folgt
 $$ l(1 - s) = r - u = 1 - s, $$
 so dass $s - 1 = 0$ und daher $s = 1$ und $r = u$ ist. Also ist
 $$ iz + zi = ri + z + t $$
 und
 $$z^2 - rz \in K. $$
 In (1) ist $4\gamma + 1 \neq 0$, sonst w\"are n\"amlich
 $$ (2i - 1)^2 = 4i^2 - 4i + 1 = 4(i + \gamma) - 4i + 1 = 0 $$
 und folglich $2i-1=0$. Dies h\"atte $2 \neq 0$ und dann den Widerspruch
 $i \in K$ zur Folge. Wir setzen nun
 $$ j := z - (r + 2t)(4\gamma + 1)^{-1}i + (t - 2r\gamma)(4\gamma + 1)^{-1}. $$
 Unter Benutzung von (2) und (1) erh\"alt man nach einfacher Rechnung
 $$ ij + ji = j. $$
 Dr\"uckt man nun $z$ gem\"a\ss\ (4) durch $i$ und $j$ aus
 und berechnet $z^2 - rz$, so erh\"alt man unter Benutzung von (5)
 und (1), dass die Koeffizienten bei $j$ und $i$ Null sind, dh.,
 man erh\"alt eine Gleichung der Form
 $$ z^2 - z = y^2 + w $$
 mit einem $w \in K$. Setzt man $\delta := z^2-z-w$, so folgt mit (5),
 dass $\delta \in K$ gilt.
 \par
       W\"are $ij \in Ki + Kj + K$, so g\"abe es $u$, $v$, $w \in K$ mit
 $ij = ai + vj + w$. Es folgte
 $$ (i - v)j = ai + w \in Ki+K. $$
 Wegen $j \notin K$ folgte $i - v = 0$ und damit der Widerspruch $i \in K$.
 Dies zeigt, dass $1$, $i$, $j$, $ij$ eine Basis von $M$ \"uber $K$ ist.
 Setze $k := ij$. Es gilt nun, die neun Relationen des Satzes
 herzuleiten. Die erste ist bereits nachgewiesen und die zweite ist
 Inhalt der Definition von $k$. Es ist
 $$ ik = i^2j = (i + \gamma)j = k + \gamma j. $$
 Mit (5) folgt
 $$ ji = j - ij = j - k. $$
 Dass $j^2 = \delta$ ein Element von $K$ ist, wurde dem Leser \"uberlassen
 nachzurechnen. Es ist
 $$ jk = jij = (j - ij)j = j^2 - ij^2 = \delta - \delta i, $$
 da $d$ ja im Zentrum $K$ von $M$ liegt. Weiter ist
 $$ ki = iji = i(j - ij) = ij - (i + \gamma)j) = -\gamma j $$
 und
 $$ kj = ij^2 = \delta i. $$
 Schlie\ss lich ist
 $$ k^2 = ijij = i(j - ij)j = i(\delta - i\delta) = i\delta -
        (i + \gamma)\delta = -\gamma \delta. $$
 Damit ist alles bewiesen.
 \medskip
       Mittels der in 10.9 notierten Relationen verifiziert man, dass
 Folgendes gilt: Ist
 $$ (a + bi + cj + dk)(a' + b'i + c'j + d'k) = a'' + b''i + c''j + d''k, $$
 so ist
 $$\eqalign{
     a'' &= aa' + bb'\gamma + cc'\delta + cd'\delta - dd'\gamma \delta    \cr
     b'' &= ab' + ba' + bb' - cd'\delta + dc'\delta                       \cr
     c'' &= ac' + ca' + bd'\gamma + cb' - db'\gamma                       \cr
     d'' &= ad' + da' + bc' + bd' - cb'. \cr} $$
 \medskip\noindent
 {\bf 10.10. Satz.} {\it Es sei $M$ ein Quaternionenschiefk\"orper
 \"uber $K$ und $1$, $i$, $j$, $k$ sei die in 10.9 beschriebene $K$-Basis
 von $M$. Ist $x = z + ui + vj + wk$ mit $z$, $u$, $v$, $w \in K$ und setzt man
 $$ x^\alpha := z + u - (ui + vj + wk), $$
 so ist $\alpha$ ein involutorischer Antiautomorphismus von $M$.
 \par
       Ist $\Char(K) \neq 2$, so sind die Fixelemente von $\alpha$ genau
 die Elemente von $K$. Ist $\Char(K) = 2$, so sind die Fixelemente
 von $\alpha$ genau die Elemente mit $u = 0$.}
 \smallskip
       Beweis. Dass $\alpha$ involutorisch ist, ist unmittelbar zu sehen.
 Ebenfalls, dass $(x + y)\alpha = x^\alpha + y^\alpha$ ist. Ebenso
 ist die Aussage \"uber die Fixelemente von $\alpha$ wegen der linearen
 Unabh\"angigkeit von $1$, $i$, $j$, $k$ banal. Das einzige, nicht triviale, was
 nachzuweisen ist, ist, dass $(xy)^\alpha = y^\alpha x^\alpha$ ist. Dazu
 berechne man mittels der vor Satz 10.10 aufgelisteten Relationen f\"ur $x$,
 $y \in M$ einmal $(xy)^\alpha$ und zum anderen $y^\alpha x^\alpha$, um
 herauszufinden, dass die beiden Ergebnisse \"ubereinstimmen.
 \medskip
       Ist $M$ ein Quaternionenschiefk\"orper der Charakteristik 2, so
 bilden die Fix\-el\-em\-en\-te des in Satz 10.10 beschriebene
 Antiautomorphismus $\alpha$ einen Unterraum des Ranges 3 des
 $Z(M)$-Vektorraumes $M$. Daher wird $M$ als K\"orper von diesen
 Fixelementen erzeugt. Ist die Charakteristik von $M$ von 2
 verschieden, so machen die Fixelemente von $\alpha$ gerade das
 Zentrum von $M$ aus. Diese Situation ist typisch f\"ur
 Quaternionenschiefk\"orper mit von 2 verschiedener Charakteristik,
 wie der n\"achste, von Dieudonné\index{Dieudonn\'e, J.}{} stammende Satz zeigt.
 \medskip
 \noindent
 {\bf 10.11. Satz.} {\it Es sei $K$ ein K\"orper und $\alpha$ sei
 ein Antiautomorphismus von $K$ mit $\alpha^2 = 1$. Setze
 $$ L := \{x \mid x \in K, x^\alpha = x\}. $$
 Dann gilt einer der folgenden F\"alle.
 \item{a)} Der K\"orper $K$ ist kommutativ, $L$ ist ein
 Teilk\"orper von $K$ und 
 $[K:L] \leq 2$.
 \item{b)} $L$ ist das Zentrum von $K$ und $K$ ist ein
 Qua\-ter\-ni\-o\-nen\-schief\-k\"or\-per \"uber $L$ mit von $2$ verschiedener
 Charakteristik.
 \item{c)} $K$ wird von $L$ erzeugt.}
 \smallskip
       Beweis. Ist $K$ kommutativ, so ist $\alpha$ ein Antiautomorphismus
 von $K$. In diesem Falle ist $L$ ein Teilk\"orper von $K$ und es
 gilt nach be\-kann\-ten S\"atzen, dass der Grad von $K$ \"uber $L$
 h\"ochstens 2 ist, dh., es gilt a).
 \par
       Es sei des Weiteren $K$ nicht kommutativ. Wir nehmen zun\"achst
 an, dass $L \subseteq Z(K)$ gelte. Dann ist $L$ ein Teilk\"orper
 von $Z(K)$. Wegen
 $$ (k + k^\alpha)^\alpha = k^\alpha + k^{\alpha^2} = k + k^\alpha $$
 und
 $$ (k^\alpha k)^\alpha = k^\alpha + k^{\alpha^2} = k^\alpha k $$
 ist
 $$ k + k^\alpha,\ k^\alpha k \in L $$
 f\"ur alle $k\in K$. Nun ist
 $$ k^2 - (k + k^\alpha)k + k^\alpha k = 0 $$
 f\"ur alle $k \in K$. Nach 10.6 ist daher $K$ ein Quaternionenschiefk\"orper
 \"uber $L$, da wir $K$ als nicht kommutativ vorausgesetzt hatten. Es
 bleibt zu zeigen, dass die Charakteristik von $K$ nicht 2 ist.
 \par
       Wir nehmen an, dass die Charakteristik von $K$ gleich 2 sei. Setze
 $$ H:= \{a \mid a \in K,\ a^2 \in K\}. $$
 Dann ist $H$ nach 10.7 eine Hyperebene des $Z(K)$-Vektorraumes $K$. Offenbar
 wird $H$ von $\alpha$ invariant gelassen und $\alpha$ induziert eine
 semilineare Abbildung auf dem $Z(K)$-Vektorraum $H$. Es sei $0 \neq a \in H$.
 Dann ist $a^2 \in Z(K)$ und andererseits auch
 $aa^\alpha \in L \subseteq Z(K)$. Es folgt  $$a^{-1}
 a^\alpha = a^{-2} aa^\alpha \in Z(K).$$ Dies besagt, dass $\alpha$
 alle Punkte von $H$ festl\"asst. Nach Satz III.1.2 gibt es daher
 ein $k \in Z(K)$ mit $h^\alpha = hk$ f\"ur alle $h \in H$.
 Insbesondere ist dann $1 = 1^\alpha = k$. Folglich ist $h^\alpha = h$ f\"ur alle
 $h \in H$ im Widerspruch zu unserer Annahme, dass $L \subseteq Z(K)$ gilt. Also
 gilt auch b).
 \par
       Es bleibt der Fall zu betrachten, dass $L \not\subseteq Z(K)$
 gilt. Es sei $R$ der von $L$ erzeugte Teilring von $K$. Ist $0 \neq r \in R$, so
 ist auch $r^\alpha \in R$. Setze $d := rr^\alpha$. Dann ist $0 \neq d \in L$.
 Folglich ist auch $d^{-1} \in L \subseteq R$. Dann ist aber auch
 $$ r^{-1} = r^\alpha d^{-1} \in R, $$
 so dass $R$ sogar ein K\"orper ist. Wir wollen
 zeigen, dass $R = K$ ist. Dazu nehmen wir an, dass dies nicht der
 Fall sei. Es sei $a \in K - R$. Es sei $u \in R$. Dann ist
 $$ au + u^\alpha a^\alpha = au + (au)^\alpha \in L. $$
 Es folgt
 $$ au - u^\alpha a = au + u^\alpha a^\alpha - u^\alpha(a + a^\alpha) \in R. $$
 Andererseits ist auch
 $$ u^\alpha a - au^\alpha
     = u^\alpha a + a^\alpha u - (a^\alpha + a) u^\alpha \in R.$$
 Somit ist schlie\ss lich $a(u-u^\alpha) \in R$. Weil $a$
 nicht in $R$ liegt, folgt $u^\alpha = u$. Somit ist $R \subseteq L$ und damit
 $R = L$, so dass $L$ ein K\"orper ist. Sind $u$, $v \in L$, so ist auch
 $uv \in L$ und es gilt
 $$ uv = (uv)^\alpha = v^\alpha u^\alpha = vu, $$
 so dass $L$ ein kommutativer K\"orper ist.
 \par
       Durch $D(u) := au - u^\alpha a$ wird, wie schon gesehen, eine
 Abbildung von $R$ in sich definiert. Weil $R = L$ ist, ist also
 $$ D(u) = au - ua. $$
 \par
       Es sei zun\"achst $\Char(K) \neq 2$. Dann ist
 $$ a = \textstyle{1 \over 2}(a + a^\alpha) +
	 \textstyle{1 \over 2}(a - a^\alpha). $$
 Das Element $a - a^\alpha$ ist nicht in $R$. Wir k\"onnen daher $a$ durch
 $a - a^\alpha$ ersetzen. Es ist dann $a^\alpha = -a$ und demzufolge
 $a^2 = -aa^\alpha \in R$. Es ist
 $$ D^2(u) = a(au - ua) - au - ua)a = a^2u -  2aua + ua^2 \in R. $$
 Wegen $a^2 u \in R$ ist daher $aD(u) \in R$. Weil $R$ ein K\"orper ist und $a$
 nicht in $R$ liegt, ist daher $D(u) = 0$, so dass $au = ua$ ist. Dies
 zeigt, dass jedes $u \in R$ mit jedem $a \in K$ vertauschbar ist,
 f\"ur das $a^\alpha = -a$ ist. Nun ist aber
 $$ a = \textstyle{1 \over 2}(a + a^\alpha) +
	\textstyle{1 \over 2}(a - a^\alpha), $$
 so dass $a$ die Summe eines Elementes $b \in R$ und einem Element $c$ mit
 $c^\alpha = -c$ ist. Weil $R$ kommutativ ist, ist daher jedes $u \in R$ mit
 allen $a \in K$ vertauschbar. Also liegt $R$ und damit
 $L$ im Zentrum von $K$ im Gegensatz zu unserer Annahme.
 \par
       Es bleibt der Fall zu betrachten, dass $\Char(K) = 2$ ist. Es sei
 wieder $a \in K-R$. Setze
 $$ V := aR + R. $$
 Dann ist $V$ ein Rechtsvektorraum \"uber $R$. Wir zeigen, dass $V$ ein
 Teilk\"orper von $K$ ist. Dazu ist zu zeigen, dass $V$
 multiplikativ ab\-ge\-schlos\-sen ist und dass das Inverse eines von
 Null verschiedenen Elementes aus $V$ wieder zu $V$ geh\"ort. Es
 ist
 $$ a^2 = a(a + a^\alpha) - aa^\alpha $$
 und $a + a^\alpha$, $aa^\alpha \in R$, so dass $a^2 \in V$ gilt.
 \par
       Als N\"achstes zeigen wir, dass auch $ra \in V$ ist f\"ur alle $r \in R$.
 Dazu d\"urfen wir nat\"urlich annehmen, dass $r \neq 0$
 ist. Von den im Folgenden zu betrachtenden Elementen aus $L$
 nehmen wir an, dass sie von Null verschieden sind. Es sei $u \in L$. Dann ist
 $$ au + ua^\alpha = au + (au)^\alpha \in L. $$
 Hieraus folgt weiter
 $$ au - ua = au + ua^\alpha - u(a + a^\alpha) \in R. $$
 Es sei $n > 1$, es seien $u_1$, \dots, $u_n \in L$ und es gelte
 $$ au_1 \cdots u_{n-1} - u_1 \cdots u_{n-1} a \in V. $$
 Weil die $u_i$ allesamt als von Null
 verschieden angenommen sind, ist $aR = au_1 \cdots u_{n-1}R$. Wir
 k\"onnen und werden daher in dem vorangegangenen Argument $a$
 durch $au_1 \cdots u_{n-1}$ ersetzen. Man erh\"alt dann, dass
 $$ au_1 \cdots u_{n} - u_n u_1 \cdots u_{n-1} a \in V. $$
 liegt. Weil $R$ kommutativ ist, folgt weiter, dass auch
 $$ au_1 \cdots u_n-u_1 \cdots u_na \in V $$
 liegt. Da $R$ als Ring additiv abgeschlossen ist, folgt schlie\ss lich, dass
 $$ ar - ra \in V $$
 f\"ur alle $r\in R$ gilt. Dies impliziert
 wiederum, dass $ra \in V$ ist f\"ur alle $r \in R$. Weiter folgt, dass
 $$\eqalign{
    (ar + s)(ar' + s') &= arar' + ars' + sar' + ss'       \cr
		       &= a^2rr' + at + t' \cr} $$
 ist mit $t$, $t' \in R$. Weil $a^2 \in V$ ist, folgt, dass $V$ multiplikativ
 ab\-ge\-schlos\-sen ist.
 \par
       Es ist $a + a^\alpha \in L = R$. Hieraus folgt, dass
 $$ a^\alpha \in aR + R = V $$
 ist. Es folgt weiter, dass $x^\alpha \in V$ ist f\"ur alle $x \in V$. Ist nun
 $0 \neq x \in V$, so ist $d := xx^\alpha \in L$. Dann ist auch $d^{-1} \in L$
 und folglich $x^{-1} = x^\alpha d^{-1} \in V$. Damit ist gezeigt, dass $V$ ein
 K\"orper ist.
 \par
       Es gilt
 $$ a(a + a^\alpha)a + aa^\alpha a = a^3 = a^2(a + a^\alpha) + a^2 a^\alpha. $$
 Also ist
 $$ aD(a + a^\alpha) = D(aa^\alpha) \in R. $$
 Hieraus folgt $D(a + a^\alpha) = 0$ und dann auch $D(aa^\alpha) = 0$, so dass
 $a$ mit $a + a^\alpha$ und $aa^\alpha$ vertauschbar ist. Setze
 $b := a + a^\alpha$ und $c := aa^\alpha$. Dann ist
 $$(ab^{-1})^2 = a^2b^{-2} = (ab + c)b^{-2} = ab^{-1} + cb^{-2}. $$
 Ersetzt man $a$ durch $ab^{-1}$, so ist also $a + a^\alpha = 1$ und daher
 $$ a^2 = a + aa^\alpha. $$
 \"Uberdies ist $a$ mit $aa^\alpha$ vertauschbar und es gilt nach wie vor
 $V = aR + R$.
 \par
       Es sei $Z$ der Zentralisator von $a$ in $R$. Weil $R$ kommutativ
 ist, liegt $Z$ im Zentrum von $V$. Es sei $x \in R$. Dann ist wegen
 $a^2 = a + aa^\alpha$ und $aa^\alpha \in R$ und weil die Charakteristik ja 2
 ist,
 $$\eqalign{
       D^2(x) &= D(ax - xa)              \cr
	      &= a^2x - axa - axa + xa^2 \cr
              &=(a + aa^\alpha)x + x(a + aa^\alpha) \cr
	      &= D(x) \cr} $$
 f\"ur alle $x \in R$. Also ist
 $$ D\bigl(D(x) + x\bigr) = 0, $$
 so dass $D(x) + x \in Z$ gilt f\"ur alle $x \in R$. Hieraus
 folgt
 $$ D(x)x - xD(x) = \bigl(D(x) + x\bigr)x - x\bigl(D(x) + x\bigr) = 0. $$
 Weiter folgt
 $$ D(x^2) = D(x^2) - 2xD(x) = D(x)x - xD(x) = 0. $$
 Also ist $x^2 \in Z$ f\"ur alle $x \in R$. Es sei $u \in V$. Dann ist
 $u = ar + s$ mit $r$,$s \in R$. Es folgt, da $r^2\in Z$ und $a + a^\alpha = 1$
 ist,
 $$\eqalign{
      uu^\alpha &= (ar + s)(ra^\alpha + s)                    \cr
		&= ar^2 a^\alpha +ars + sra^\alpha + s^2      \cr
		&= aa^\alpha r^2 + D(rs) + rs(a+a^\alpha)+s^2 \cr
                &= aa^\alpha r^2 + D(rs) + rs + s^2. \cr} $$
  Wir haben schon gesehen, dass $aa^\alpha$ mit $a$
 vertauschbar ist. Da $aa^\alpha$ auch in $R$ liegt, ist
 $aa^\alpha$ ein Element von $Z$. Weil auch $r^2$, $s^2$ und $D(rs) + rs$ in $Z$
 liegen, ist $uu^\alpha$ f\"ur alle $u\in V$ ein Element von $Z$. Weiter gilt
 $$ u + u^\alpha = ar + s + ra^\alpha + s
      = ar + ra + r(a + a^\alpha) = D(r) + r, $$
 so dass auch $u+u^\alpha$ f\"ur alle $u \in V$ ein Element von $Z$ ist. Wegen
 $$ u^2 + u(u + u^\alpha) + uu^\alpha = 0 $$
 haben wir f\"ur $V$ und $Z$ die in Satz 10.6 beschriebene Situation.
 \par
       Wir nehmen zun\"achst an, $V$ sei ein Quaternionenschiefk\"orper
 \"uber $Z$. Es ist $uu^\alpha \in Z$ f\"ur alle $u \in V$. Ist $u
 \neq 0$, so folgt $u^{-1} Z^* = u^\alpha Z^*$, so dass die
 Kollineation $uZ^* \rightarrow u^{-1} Z^*$ durch die $Z$-lineare
 Abbildung $\alpha$ induziert wird. Diese Kollineation ist aber
 eine Elation, wie wir schon fr\"uher bemerkten. Daher l\"asst
 $\alpha$ eine Hyperebene $H$ des $Z$-Vektorraumes $V$ punktweise
 fest. Weil der Rang von $H$ gleich 3 ist, gibt es ein $k \in Z$
 mit $h^\alpha = hk$ f\"ur alle $h \in H$. (Dies wieder im Vorgriff auf III.1.2.)
 Wegen $R \subseteq H$
 folgt $k = 1$. Dies impliziert wiederum $H = L = R$. Nun hat $V$
 aber den Rang 2 \"uber $R$ und den Rang 4 \"uber $Z$. Das
 bedeutet, dass $R$ den Rang 2 \"uber $Z$ hat. Damit erhalten wir
 das widerspr\"uchliche Ergebnis $2 = 3$. Also ist $V$ kein
 Quaternionenschiefk\"orper \"uber $Z$. Daher ist $V$ nach 10.6 ein
 kommutativer K\"orper. Dies besagt, dass $R$, da $V$ im Laufe des
 Beweises nicht ge\"andert wurde, wie wir bemerkten, im
 Zentralisator des urspr\"unglich gew\"ahlten Elementes $a$ liegt. Da
 dieses Element ein beliebiges Element aus $K-R$ war und da $R$
 kommutativ ist, liegt $R$ im Zentrum von $K$ entgegen unserer
 Annahme. Damit ist der Satz in allen seinen Teilen bewiesen.


 \newpage
       
 \mychapter{III}{Gruppen von Kollineationen}

 \noindent
 Mit jeder Struktur geht ihre Automorphismengruppe einher und das
 Studium der Automorphismengruppe l\"asst R\"uckschl\"usse auf
 die Ausgangsstruktur zu. Dies haben wir in dem uns hier
 interessierenden Falle der projektiven Geometrien im letzten
 Kapitel bereits gesehen, wo uns das eingehende Studium der Gruppen
 $\E(H)$ und $\Delta(P,H)$ zu den Strukturs\"atzen\index{Struktursatz}{} der
 projektiven Geometrie f\"uhrte. In diesem Kapitel nun werden wir
 die volle Kollineationsgruppe einer projektiven Geometrie und
 gewisse ihrer Untergruppen etwas sys\-te\-ma\-ti\-scher studieren.

\mysection{1. Erste \"Ubersicht}

\noindent
 Es sei $V$ ein Vektorraum \"uber dem K\"orper $K$. Mit $\GammaL(V,K)$
 bezeichnen wir die Gruppe aller Automorphismen des $K$-Vektorraumes $V$,
 dh., die Gruppe aller bijektiven semilinearen Abbildungen von $V$
 auf sich. Mit $\GaL(V,K)$ bezeichnen wir die Gruppe aller
 bijektiven li\-ne\-a\-ren Abbildungen von $V$ auf sich, mit $\GL(V,K)$
 die von allen Homothetien und Transvektionen erzeugte Gruppe und
 mit $\SL(V,K)$ die von allen Transvektionen erzeugte Gruppe. Weil
 Homothetien und Transvektionen lineare Abbildungen sind, gilt
 $$ \SL(V,K) \subseteq \GL(V,K) \subseteq \GaL(V,K) \subseteq \GammaL(V,K). $$
 Da Transvektionen, Homothetien und lineare
 Abbildungen unter inneren Automorphismen von $\GammaL(V,K)$
 wieder in solche \"ubergehen, folgt, dass die Gruppen $\SL(V,K)$ und
 $\GL(V,K)$ Normalteiler von $\GammaL(V,K)$ sind.
 \medskip\noindent
 {\bf 1.1. Satz.} {\it Es sei $V$ ein Vektorraum \"uber dem
 K\"orper $K$. Ist der Rang von $V$ mindestens $2$ und enth\"alt
 $K$ mindestens drei Elemente, so wird $\GL(V,K)$ bereits von den
 Homothetien von $V$ erzeugt.}
 \smallskip
       Beweis. Es sei $\tau$ eine von 1 verschiedene Transvektion von
 $V$. Dann gibt es eine lineare Abbildung $\varphi$ von $V$ in $K$
 und einen Vektor $a \in \Kern(\varphi)$ mit $x^\tau = x + a\varphi(x)$ f\"ur
 alle $x \in V$. Weil $\tau \neq 1$ ist, ist $a \neq 0$
 und $\Kern(\varphi) \neq V$. Es folgt, dass $\varphi$ surjektiv
 ist. Weil $K$ mehr als zwei Elemente enth\"alt, gibt es daher einen
 Vektor $p \in V$ mit $\varphi(p) \neq 0$, 1. Setze $k : = -\varphi(p) + 1$.
 Dann ist $k \neq 0$, so dass $k^{-1}$ existiert. Ferner
 setzen wir $q := p + a$. Weil $p$ wegen $\varphi(p) \neq 0$ nicht
 in $H := \Kern(\varphi)$ liegt, liegt $q$ nicht in $H$, und weil
 $a \neq 0$ ist, folgt, dass $p$ und $q$ linear unabh\"angig sind.
 Somit sind $pK$ und $qK$ zwei verschiedene Komplemente von $H$.
 Wir definieren nun zwei Homothetien $\delta$ und $\eta$ verm\"oge
 $(px + h)^\delta := pkx + h$ bzw. $(qx + h)^\eta := qk^{-1}x + h$ f\"ur
 alle $x \in K$ und alle $h \in H$. Dann ist
 $$ \eqalign{
    (px + h)^{\delta\eta} &= (pkx + h)^\eta = p^\eta kx + h    \cr
	  &= (q - a)^\eta kx + h = (qk^{-1} - a)kx + h         \cr
	  &= qx - akx + h = px + h + a(1 - k)x                 \cr
	  &= px +h + a\varphi(p)x = px + h + a \varphi(px + h) \cr
	  &= (px + h)^\tau. \cr} $$
 Also ist $\delta\eta = \tau$. Damit ist gezeigt, dass alle
 Transvektionen in der von allen Homothetien erzeugten Gruppe liegen.
 \medskip
       Mit $\PSL(V,K)$, $\PGL(V,K)$, $\PGaL(V,K)$ und $\PGammaL(V,K)$
 be\-zeich\-nen wir die von den Gruppen $\SL(V,K)$, $\GL(V,K)$, $\PGaL(V,K)$
 bzw. $\GammaL(V,K)$ in $L_K(V)$ induzierten Kollineationsgruppen. Diese
 Gruppen hei\ss en der Reihe nach die {\it kleine
 projektive Gruppe\/},\index{kleine projektive Gruppe}{}
 {\it projektive Gruppe\/}\index{projektive Gruppe}{} und
 {\it gro\ss e projektive Gruppe\/}\index{gro\ss e projektive Gruppe}{} von
 $L_K(V)$. Nach II.7.1 ist $\PGammaL(V,K)$ die volle Kollineationsgruppe von
 $L_K(V)$, falls $\Rg_K(V) \geq 3$ ist. Die Gruppe $\PGL(V,K)$ ist dann die
 von allen Streckungen und Elationen erzeugte Gruppe und wird nach
 1.1 bereits von allen Streckungen erzeugt, wenn $K$ mindestens
 drei Elemente enth\"alt. Die Gruppe $\PSL(V,K)$ ist die Untergruppe von
 $\PGammaL(V,K)$, die von allen Elationen erzeugt wird.
 \par
       Sind keine Verwechslungen zu bef\"urchten, so lassen wir im
 Folgenden den zweiten Parameter $K$ bei diesen Gruppen weg.
 \par
       Wir fragen nun nach dem Kern $M(V)$\index{Kern}{} des Homomorphismus
 von $\GammaL(V)$ auf $\PGammaL(V)$. Ist $\Rg_K (V) = 1$, so ist
 nat\"urlich $M (V) = \GammaL(V)$. Andernfalls gilt:
 \medskip\noindent
 {\bf 1.2. Satz.} {\it Es sei $V$ ein Vektorraum \"uber dem
 K\"orper $K$ 
 mit $\Rg_K(V) \geq 2$. Dann ist
 $$ M(V) = \bigl\{\mu(k) \mid k \in K^*\bigr\}, $$
 wobei $\mu (k)$ die durch $x^{\mu(k)} := xk$ definierte Abbildung ist.
 Der Zentralisator $C_{\GammaL(V)} (M(V))$ von $M(V)$ in $\GammaL(V)$ ist gleich
 $\GaL(V)$. Ferner gilt}
 $$ M(V) \cap \GaL(V) = \bigl\{\mu(k) \mid k \in Z(K^*)\bigr\}. $$
 \par
       Beweis. Jedes $\mu(k)$ ist eine bijektive semilineare Abbildung,
 die alle Punkte von $L_K(V)$ invariant l\"asst. Somit gilt
 $\mu (k) \in M(V)$.
 \par
       Es sei $\rho \in M(V)$. Dann gibt es zu jedem $v \in V - \{0\}$
 genau ein $k_v \in K^*$ mit $v^\rho = vk_v$. Sind $u$, $v \in V - \{0\}$ und
 ist $u + v \neq 0$, so gilt
 $$ (u + v) k_{u+v} = (u + v)^\rho = u^\rho + v^\rho = uk_u + vk_v. $$
 Sind $u$ und $v$
 linear unabh\"angig, so ist
 $$ k_u = k_{u+v} = k_v. $$
 Sind $u$ und $v$ linear abh\"angig, ist $uK = vK$. Wegen $\Rg_K(V) \geq 2$
 gibt es ein $w \in V - uK$. Dann sind
 $w$ und $u$ wie auch $w$ und $v$ linear unabh\"angig. Nach dem
 bereits Bewiesenen ist daher
 $$ k_u = k_w = k_v.$$ 
 Hieraus folgt $\rho = \mu (k_u)$, so dass die erste Aussage des Satzes
 bewiesen ist.
 \par
       Es sei $\gamma \in C_{\GammaL(V)}(M(V))$. Ferner sei $\alpha$
 der begleitende Automorphismus von $\gamma$. Ist dann $k \in K^*$,
 so folgt
 $$ v^\gamma k = v^{\gamma\mu(k)} = v^{\mu(k)\gamma} = (vk)^\gamma
	       = v^\gamma k^\alpha. $$
 Also ist $\alpha = 1$ und daher $\gamma \in \GaL(V)$. Ist umgekehrt
 $\gamma \in \GaL(V)$, so ist
 $$ v^{\gamma\mu(k)} = v^\gamma k = (vk)^\gamma = v^{\mu(k)\gamma}, $$
 so dass $\gamma \in C_{\GammaL(V)} (M(V))$ ist.
 \par
       Es sei $\mu(k) \in M(V) \cap \GaL(V)$. Dann ist
 $$ vlk = (vl)^{\mu(k)} = v^{\mu (k)} l = vkl $$
 f\"ur alle $v \in V$ und alle $l \in K$. Es folgt $lk = kl$ und damit
 $k \in Z(K^*)$. Ist umgekehrt $k \in Z(K^*)$, so ist $\mu(k)$ linear, wie man
 unmittelbar sieht. Damit ist der Satz in allen seinen Teilen bewiesen.
 \medskip
      Ist $K$ ein K\"orper, so bilden die inneren Automorphismen von $K$
 eine Gruppe $\Inn(K)$, die sogar ein Normalteiler der Gruppe
 $\Aut(K)$ ist. Die Faktorgruppe $\Aut(K)/\Inn(K)$ hei\ss t
 {\it \"au\ss ere Automorphismengruppe\/} von
 $K$.\index{auss@\"au\ss ere Automorphismengruppe}{}
 \par
       Wir nennen eine Kollineation eines projektiven Raumes {\it projektiv\/},
 falls\index{projektive Kollineation}{} sie in der gro\ss en projektiven Gruppe
 enthalten ist.
 \medskip\noindent
 {\bf 1.3. Satz.} {\it Ist $V$ ein Vektorraum \"uber dem K\"orper
 $K$ mit $\Rg_K(V) \geq 2$, so ist
 $$ \PGammaL(V) / \PGaL(V) $$
 zur \"au\ss eren Automorphismengruppe von $K$ isomorph.
 Insbesondere gilt: ist $\gamma \in \GammaL(V)$, so ist die von
 $\gamma$ in $L_K(V)$ induzierte Kollineation genau dann projektiv,
 wenn der begleitende Automorphismus von $\gamma$ ein innerer
 Automorphismus von $K$ ist.}
 \smallskip
       Beweis. Es sei $\gamma \in \GammaL(V)$ und $\alpha$ sei der
 begleitende Automorphismus von $\gamma$.
       Wir nehmen an, dass $\gamma$ eine projektive Kollineation
 induziere. Es gibt dann nach 1.2 ein $\delta \in \GaL(V)$ und ein
 $k \in K^*$ mit $\gamma \delta ^{-1} = \mu(k)$. Weil $\delta$
 nach 1.2  mit $\mu (k)$ vertauschbar ist, ist $\gamma = \delta\mu(k)$. Ist
 nun $l \in K$, so folgt
 $$ v^\delta lk = (vl)^{\delta\mu(k)} = (vl)^\gamma = v^\gamma l^\alpha
     = v^\delta kl^\alpha. $$
 Also ist $l^\alpha = k^{-1} lk$, so dass $\alpha \in \Inn(K)$ gilt.
 \par
       Es sei umgekehrt $\alpha \in \Inn(K)$. Es gibt dann ein $k \in K^*$
 mit $l^\alpha = k^{-1} lk$ f\"ur alle $l \in K$. Setze
 $\delta := \gamma\mu(k^{-1})$. Dann induziert $\delta$ die
 gleiche Kollineation wie $\gamma$. Wegen
 $$ (vl)^\delta = (vl)^\gamma k^{-1} = v^\gamma k^{-1} lkk^{-1} = v^\delta l $$
 ist $\delta$ linear und $\gamma$ somit projektiv.
 \par
       Das bislang Bewiesene zeigt, dass $\GaL(V) M(V)$ die Menge aller
 Elemente aus $\GammaL(V)$ ist, deren Begleitautomorphismen innere
 Automorphismen von $K$ sind. Ist nun $\beta$ die Abbildung von
 $\GammaL(V)$ in $\Aut(K)$, die jedem Element seinen begleitenden
 Automorphismus zuordnet, so ist $\beta$ ein Homomorphismus, der
 sogar surjektiv ist, wie man unschwer nachpr\"uft. Definiert man nun
 $\varphi(\gamma)$ f\"ur $\gamma \in \GammaL(V)$ durch $\varphi (\gamma) :=
 \beta (\gamma) \Inn(K)$, so ist $\varphi$ ein Epimorphismus
 von $\GammaL(V)$ auf die \"au\ss ere Automorphismengruppe von $K$
 und es gilt, wie wir gesehen haben,
 $$ \Kern(\varphi) = \GaL(V) M(V). $$
 Wegen
 $$\eqalign{
 \GammaL(V)/\GaL(V)M(V)
	&\cong \bigl(\GammaL(V)/M(V)\bigr)/\bigl(\GaL(V)M(V)/M(V)\bigr) \cr
        &= \PGammaL(V)/\PGaL(V) \cr} $$
 folgt schlie\ss lich die noch offene Behauptung des Satzes.
 \medskip\noindent
 {\bf 1.4. Satz.} {\it Es sei $V$ ein Vektorraum \"uber $K$. Sind
 $P$ und $Q$ Punkte von $L_K (V)$ und ist $T$ ein Teilraum, der
 weder $P$ noch $Q$ enth\"alt, so gibt es eine Hyperebene $H \in
 L_K(V)$ mit $T \geq H$ und $P$, $Q \not\leq H$.}
 \smallskip
       Beweis. Nach I.1.8 gibt es Hyperebenen $H_1$ und $H_2$ mit $T \leq
 H_i$ und $P + H_1 = V = Q+H_2$. Ist $Q \not\leq H_1$, so ist $H_1$
 eine Hyperebene der gesuchten Art. Ist $P \not\leq H_2$, so ist $H_2$
 eine solche. Wir d\"urfen also annehmen, dass $P \leq H_2$ und $Q
 \leq H_1$ ist. Ist dann $R$ ein von $P$ und $Q$ verschiedener
 Punkt auf $P + Q$, so setzen wir $H := R + (H_1 \cap H_2)$. W\"are $H$
 keine Hyperebene, so folgte $R \leq H_1 \cap H_2$ und weiter
 $$ P \leq P + Q = R + Q \leq H_1. $$
 Dieser Widerspruch zeigt, dass
 $H$ doch eine Hyperebene ist. Mittels des Mo\-du\-lar\-ge\-setz\-es folgt
 $$\eqalign{
       P + H &= P + R + (H_1 \cap H_2) \cr
             &= P + Q + (H_1 \cap H_2) \cr
	     &= P + (H_1 \cap (Q + H_2)) \cr
	     &= P + H_1 = V. \cr} $$
 Ebenso folgt, dass auch $Q + H = V$ ist.
 \medskip\noindent
 {\bf 1.5. Satz.} {\it Es sei $V$ ein Vektorraum \"uber dem
 K\"orper $K$ und $U$ sei ein Unterraum endlichen Ranges von $V$.
 Setze $n := \Rg_K(U)$. Ist $\gamma \in \GaL(V)$, so gibt es
 $\sigma_1$, \dots, $\sigma_n \in \GL(V)$ mit den Eigenschaften:
 \item{(1)} Jedes $\sigma_i$ l\"asst eine Hyperebene vektorweise fest.
 \item{(2)} $\gamma \sigma_1 \dots \sigma_n$ l\"asst $U$ vektorweise fest.}
 \smallskip
       Beweis. Wir machen Induktion nach $n$. Ist $n = 0$, so ist nichts zu
 beweisen. Es sei also $n \geq 1$. Ferner sei $U = uK \oplus W$ mit
 einem $u \neq 0$. Dann ist $\Rg_K(W) = n - 1$, so dass es nach
 Induktionsannahme $\sigma, \dots, \sigma_{n-1} \in \GL(V)$ gibt,
 so dass $(1)$ und $(2)$ gelten, wenn man nur $n$ durch $n - 1$ und
 $U$ durch $W$ ersetzt. Setze $\rho := \gamma \sigma_1 \dots
 \sigma_{n-1}$. Dann gilt $u$, $u^\rho \in W$. Nach 1.4 gibt es eine
 Hyperebene $H$ mit $W \subseteq H$, die weder $u$ noch $u^\rho$
 enth\"alt. Es gibt daher eine lineare Abbildung $\sigma_n$, die
 $H$ vektorweise festl\"asst und die $u^\rho$ auf $u$ abbildet.
 Es folgt, dass $\rho\sigma_n$ den Teilraum $W$ vektorweise
 festl\"asst und dass $u$ ebenfalls von dieser Abbildung
 festgelassen wird. Weil $\rho \sigma_n$ linear ist und $U = uK \oplus W$ gilt,
 bleibt jeder Vektor aus $U$ bei $\rho \sigma_n$
 fest. Damit ist der Satz bewiesen.
 \medskip\noindent
 {\bf 1.6. Satz.} {\it Es sei $V$ ein Vektorraum \"uber dem
 K\"orper $K$. Ist $\sigma \in \GaL(V)$, so gilt genau dann
 $\sigma \in \GL(V)$, wenn es einen Unterraum endlichen Ko-Ranges
 gibt, der von $\sigma$ vektorweise festgelassen wird.}
 \smallskip
       Beweis. Ist $\sigma \in \GL(V)$, so ist $\sigma$ Produkt von
 endlich vielen Trans\-vek\-ti\-o\-nen und Homothetien. Der Schnitt
 der Achsen dieser Elationen und Homothetien ist dann ein Unterraum
 endlichen Ko-Ranges, der von $\sigma$ vektorweise festgelassen wird.
 \par
       Es sei $W := \{v \mid v \in V, v^\sigma = v\}$ und $W$ habe endlichen
 Ko-Rang. Setze $n := \Ko_K(W)$. Ist $n = 0$, so ist nichts zu
 beweisen. Es sei also $n > 0$ und $b_1, \dots, b_n$ seien $n$
 linear unabh\"angige Vektoren, die nicht in $W$ liegen. Dann ist
 auch $b_n^\sigma \not\in W$. Nach 1.4 gibt es also eine Hyperebene
 $H$ die zwar $W$, nicht aber $b_n$ und $b_n^\sigma$ enth\"alt. Es
 gibt daher ein $\tau \in \GL(V)$, welches $H$ vektorweise
 festl\"asst und $b^\sigma_n$ auf $b_n$ abbildet. Dann l\"asst
 aber $\sigma \tau$ den Unterraum $b_nK + W$ vektorweise fest, so
 dass Induktion zum Ziele f\"uhrt.
 \medskip\noindent
 {\bf 1.7. Satz.} {\it Ist $V$ ein Vektorraum \"uber dem K\"orper
 $K$ so sind die folgenden Aussagen \"aquivalent:
 \item{a)} $\Rg_K(V)$ ist endlich.
 \item{b)} Es ist $\GL(V) = \GaL(V)$.
 \item{c)} Es ist $\PGL(V) = \PGaL(V)$.}
 \smallskip
       Beweis. a) impliziert b): Dies folgt mit 1.6.
 \par
       b) impliziert c): Banal.
 \par
       c) impliziert a): Wir zeigen, dass die Verneinung von a)
 die Verneinung von c) nach sich zieht. Ist der Rang von $V$ nicht
 endlich, so gibt es eine abz\"ahlbare Teilmenge
 $ \{b_i \mid i \in \omega\} $
 von linear unabh\"angigen Vektoren in $V$. Es sei
 $$ U := \sum_{i:=0}^\infty b_iK $$
 und $W$ sei ein
 Komplement von $U$ in $V$. Wir definieren $\sigma \in \GaL(V)$
 durch $b^\sigma_i := b_i + b_{i+1}$ und $w^\sigma := w$ f\"ur alle
 $i \in \omega$ und alle $w \in W$. Es sei $vK$ ein Fixpunkt von
 $\sigma$. Es gibt dann ein $n \in \omega$, sowie $k_0, \dots, k_n
 \in K$ und ein $w \in W$ mit
 $$ v=\sum_{i:=0}^n b_ik_i + w. $$
 Weil $vK$ ein Fixpunkt ist, gibt es ein $a \in K^*$ mit
 $$
    \sum_{i:=0}^n b_ik_ia + wa = va = v^\sigma
                                = \sum_{i:=0}^n b_i^\sigma k_i + w 
                               = \sum_{i:=0}^n (b_i + b_{i+1})k_i + w. 
 $$
 Hieraus folgt $0 = k_n$ und $k_ia = k_i + k_{i-1}$ f\"ur $i := 0$, \dots, $n$.
 Dies hat zur Folge, dass alle $k_i$ gleich Null sind. Somit
 gilt $vK \subseteq W$. W\"are nun die von $\sigma$ induzierte
 Kollineation ein Element von $\PGL(V)$, so w\"are diese
 Kollineation ein Produkt von endlich vielen Elationen und
 Streckungen. Der Schnitt ihrer Achsen w\"are ein Teilraum
 endlichen Ko-Ranges in $V$. Dieser Teilraum l\"age nach unserer
 Vorbemerkung aber in $W$, so dass $W$ endlichen Ko-Rang h\"atte,
 was nicht der Fall ist. Damit ist alles bewiesen.
 \medskip
       Als N\"achstes wollen wir die papposschen
 R\"aume\index{papposscher Raum}{} mittels ihrer
 Kol\-li\-ne\-a\-ti\-ons\-grup\-pen charakterisieren. Dazu beweisen wir
 zun\"achst den folgenden Satz.
 \medskip\noindent
 {\bf 1.8. Satz.} {\it Es sei $V$ ein Vektorraum \"uber dem
 K\"orper $K$ und $U$ sei ein Teilraum von $V$ mit $2 \leq \Rg_K(U)<\infty$.
 Setze $n:= \Rg_K(U)$. Ist dann $P_0$, \dots, $P_n$
 ein Rahmen\index{Rahmen}{} von $U$, so gibt es Vektoren $v_0$, \dots, $v_n$ mit
 $P_i = v_iK$ f\"ur $i := 0$, \dots, $n$ und $v_0 = \sum_{i:=1}^n v_i$.}
 \smallskip
       Beweis. Es sei $P_i = u_iK$. Weil $P_1$, \dots, $P_n$ eine Basis von
 $U$ ist, gibt es $k_i \in K$ mit $u_0 = \sum_{i:=1}^n
 u_ik_i$. Weil je $n$ der $n + 1$ Punkte eines Rahmens von $U$ eine
 Basis von $U$ bilden, folgt, dass die $k_i$ allesamt ungleich Null
 sind. Setzt man daher $v_0 := u_0$ und $v_i := u_ik_i$ f\"ur
 $i:=1$, \dots, $n$, so tun's diese $v_i$.
 \medskip\noindent
 {\bf 1.9. Satz.} {\it Es sei $V$ ein Vektorraum \"uber dem
 K\"orper $K$. Ferner sei $n$ eine nat\"urliche Zahl gr\"o\ss er
 als $1$. Sind dann $U$ und $W$ zwei Unterr\"aume des Ranges $n$
 von $V$, ist $P_0$, \dots, $P_n$ ein Rahmen von $U$ und $Q_0$, \dots,
 $Q_n$ ein solcher von $W$, so gibt es ein $\sigma \in \PGL(V)$ mit
 $P_i^\sigma = Q_i$ f\"ur $i := 1, \dots, n$. Ist der Rang von $V$
 endlich, so ist $\PGL (V)$ also transitiv auf der Menge der Rahmen
 von $V$.}
 \smallskip
       Beweis. Nach I.7.2 haben $U$ und $W$ ein gemeinsames Komplement
 $C$ in $V$. Nach 1.8 gibt es ferner $u_0, \dots, u_n$ und $w_0,
 \dots, w_n$ mit $P_i = u_iK$ und $Q_i = w_iK$ f\"ur $i := 0,
 \dots, n$ sowie $u_0 = \sum_{i:=1}^n u_i$ und $w_0 =
 \sum_{i:=1}^n w_i$. Definiert man nun $\rho$ durch
 $c^\rho := c$ f\"ur alle $c \in C$ und $u^\rho_i := v_i$ f\"ur
 $i:=1, \dots, n$, so ist $\rho$ nach 1.6 ein Element in $\GL(V)$.
 Da offenbar $u^\rho_0 = w_0$ ist, leistet die von $\rho$
 induzierte Kollineation $\sigma$ das Verlangte.
 \medskip
       Es sei $V$ ein Vektorraum des Ranges 2 \"uber dem
 K\"orper $K$. Motiviert durch den Satz II.6.3 nennen wir $L(V)$
 pappossch, wenn $K$ kommutativ ist.
 \medskip\noindent
 {\bf 1.10. Satz.} {\it Es sei $V$ ein Vektorraum \"uber dem
 K\"orper $K$ und es gelte $\Rg_K(V) \geq 2$. Dann sind
 \"aquivalent:
 \item{a)} $L_K(V)$ ist pappossch.\index{papposscher Raum}{}
 \item{b)} Ist $U$ ein Teilraum des endlichen Ranges $n \geq 2$ von
 $V$, ist $P_0$, \dots, $P_n$ ein Rahmen von $U$, ist $\sigma \in
 \PGaL(V)$ und gilt $P^\sigma_i = P_i$ f\"ur $i := 0$, \dots, $n$, so
 induziert $\sigma$ auf $U/0$ die Identit\"at.
 \item{c)} Es gibt einen Teilraum $U$ des endlichen Ranges $n \geq 2$ von $V$
 und einen Rahmen\index{Rahmen}{} $P_0$, \dots, $P_n$ von $U$, so dass jedes
 $\sigma \in \PGL(V)$, f\"ur das $P_i^\sigma = P_i$ f\"ur $i := 0$, \dots,
 $n$ gilt, auf $U/0$ die Identit\"at induziert.\par}
 \smallskip
       Beweis. a) impliziert b): Nach II.6.3 ist $K$ kommutativ. Davon
 werden wir gleich Gebrauch machen.
 \par
       Es sei nun $P_i = u_iK$ und es
 gelte $u_0 = \sum_{i:=1}^n u_i$. Es sei $\sigma \in \GaL(V)$ und es gelte
 $P^\sigma_i=P_i$ f\"ur alle $i$. Es gibt dann
 $k_i \in K^*$ mit $u_i^\sigma = u_ik_i$ f\"ur alle $i$. Es folgt
 $$ \sum_{i:=1}^n u_ik_i = \biggl(\sum_{i:=1}^n u_i \biggr)^\sigma
			 = u_0^\sigma = u_0k_0 = \sum^n_{i:=1} u_ik_0. $$
 Wegen der linearen Unabh\"angigkeit von $u_1$, \dots, $u_n$ folgt
 $k_0 = k_1 = \dots = k_n$. Ist nun $v \in U$, so gilt $v = \sum^n_{i:=1} u_ia_i$
 mit gewissen $a_i \in K$. Hieraus
 folgt zusammen mit der Kommutativit\"at von $K$, dass
 $$ v^\sigma = \sum_{i:=1}^n u_ik_0a_i
      = \biggl(\sum^n_{i:=} u_ia_i \biggr) k_0 = vk_0 $$
 ist. Somit gilt b).
 \par
       b) impliziert c): Weil der Rang von $V$ mindestens 2 ist, ist
 c) nat\"urlich eine Folge von b).
 \par 
       c) impliziert a): Es gibt wieder $u_i \in U$ mit $P_i = u_i K$ und
 $u_0 = \sum_{i:=1}^n u_i$. Es sei $k \in K^*$. Dann ist
 $u_1k$, \dots, $u_nk$ eine Basis von $U$. Unter Zuhilfenahme eines
 Komplementes von $U$ folgt die Existenz eines $\sigma \in \GL(V)$
 mit $u^\sigma_i = u_i k$ f\"ur $i := 1$, \dots, $n$. Es folgt, dass
 auch $u^\sigma_0 = u_0k$ ist. Also ist $P^\sigma_i = P_i$ f\"ur
 alle $i$, so dass $\sigma$ nach Voraussetzung auf der Menge der
 Punkte von $L_K(V)$ die Identit\"at induziert. Mit Satz 1.2
 folgt, dass $v^\sigma = vk$ gilt f\"ur alle $v \in U$ und dass
 dar\"uber hinaus $k$ ein Element von $Z(K)$ ist, da $\sigma$ ja in
 $\GL(V)$ liegt. Also ist $K^* \subseteq Z(K)$, so dass a) eine
 Folge von c) ist.
 \medskip\noindent
 {\bf 1.11. Korollar.} {\it Es sei $V$ ein Vektorraum \"uber dem
 K\"orper $K$ mit $2 \leq \Rg_K(V) < \infty$. Genau dann ist
 $L_K(V)$ pappossch, wenn $\PGL(V)$ auf der Menge der Rahmen von
 $L_K(V)$ scharf transitiv operiert.}
 \smallskip
 Beweis. Dies folgt unmittelbar aus 1.9 und 1.10.
 \medskip
       Ist $V$ ein Vektorraum endlichen Ranges \"uber einem kommutativen
 K\"orper und ist $\gamma \in \GL(V)$, so bezeichnen wir mit
 $\det(\gamma)$ die Determinante von $\gamma$.\index{Determinante}{}
 \medskip\noindent
 {\bf 1.12. Satz.} {\it Ist $V$ ein Vektorraum \"uber dem
 kommutativen K\"orper $K$ mit $1 \leq \Rg_K(V) < \infty$, so gilt:
 \item{a)} Die Abbildung $\det$ ist ein Homomorphismus von $\GL(V)$ auf $K^*$.
 \item{b)} $\SL(V)$ ist der Kern von $\det$.}
 \smallskip
       Beweis. a) Nach 1.7 ist $\GL(V) = \GaL(V)$, so dass a) dem Leser
 aus dem Anf\"angerunterricht bekannt sein sollte.
 \par
       b) Es sei $V = P \oplus H$ mit einem Punkt $P$ und einer
 Hyperebene $H$. Dann ist $\SL(V)\Sigma(P,H) \subseteq \GL(V)$,
 wobei $\Sigma(P,H)$ wieder die Gruppe aller Homothetien mit dem
 Zentrum $P$ und der Achse $H$ bezeichne. Da die Gruppe $\SL(V)$ auf
 der Menge der nicht inzidenten Punkt-Hyperebenenpaare von $L(V)$
 transitiv operiert, enth\"alt $\SL(V)\Sigma(P,H)$ alle
 Homothetien. Daher ist $\SL(V)\Sigma(P,H) = \GL(V)$.
 \par
       Es sei $\tau$ eine Transvektion.\index{Transvektion}{} Es gibt dann eine
 lineare Abbildung $\varphi$ von $V$ auf $K$ und ein $a \in \Kern(\varphi)$
 mit $x^\tau = x + a\varphi(x)$ f\"ur alle $x \in V$. Ist $b_1$,
 \dots, $b_n$ eine Basis von $V$ mit $b_2$, \dots, $b_n \in H$, so
 wird $\tau$ bez\"uglich dieser Basis durch eine Dreiecksmatrix mit
 lauter Einsen auf der Hauptdiagonalen dargestellt. Daher ist
 $\det(\tau) = 1$, so dass $\SL(V) \subseteq \Kern(\det)$ gilt. Ist
 nun $\kappa \in \Kern(\det)$, so gibt es ein $\sigma \in \SL(V)$
 und ein $k \in K^*$ mit $\kappa = \sigma\delta(k)$. Es folgt
 $\det(\kappa) = \det(\delta(k))$. W\"ahlt man nun eine Basis $b_1$,
 \dots, $b_n$ mit $b_1 \in P$ und $b_2$, \dots, $b_n \in H$, so sieht
 man, dass $\delta (k)$ durch eine Diagonalmatrix dargestellt wird,
 wo an der ersten Stelle der Diagonalen $k^{-1}$ steht und alle
 \"ubrigen Diagonalelemente gleich $1$ sind. Also ist
 $$ 1 = \det(\kappa) = k^{-1}, $$
 so dass $k = 1$ und damit $\kappa \in \SL(V)$ ist.
 \medskip
       Wir sind nun in der Lage, die Ordnungen\index{Gruppenordnung}{} der
 Gruppen $\GammaL(V)$
 etc.\ auszurechnen, falls $V$ ein endlicher Vektorraum ist. Ist $V$
 ein Vektorraum des Ranges $n$ \"uber $\GF(q)$, so schreiben wir
 statt $\GammaL(V)$, usw., $\GammaL(n,q)$, usw. 
 \medskip\noindent
 {\bf 1.13. Satz.} {\it Sind $n$ und $r$ nat\"urliche Zahlen, ist
 $p$ eine Primzahl und ist $q := p^r$, so gilt:
 \smallskip
 \item{a)} $|\GammaL(n,q)| = rq^{{1 \over 2} n(n-1)} \prod_{i:=1}^n (q^i - 1).$
 \smallskip
 \item{b)} $|\GL(n,q)| = q^{{1 \over 2} n(n-1)} \prod_{i:=1}^n (q^i - 1).$
 \smallskip
 \item{c)} $| SL(n,q)| = q^{{1 \over 2} n(n-1)} \prod_{i:=2}^n (q^i - 1).$
 \smallskip
 \item{d)} $|\PGammaL(n,q)| = rq^{{1 \over 2} n(n-1)} \prod_{i:=2}^n (q^i - 1).$
 \smallskip
 \item{e)} $|\PGL(n,q)| = q^{{1 \over 2} n(n-1)} \prod_{i:=2}^n (q^i - 1).$
 \smallskip
 \item{f)} $|\PSL(n,q)| = \ggT(n,q - 1)^{-1} q^{{1 \over 2} n(n-1)}
		       \prod_{i:=2}^n (q^i - 1).$}
 \medskip
       Beweis. Da alle endlichen K\"orper kommutativ sind, ist die Gruppe
 $\PGL(n,q)$ auf der Menge der Rahmen des zu Grunde liegenden
 projektiven Raumes nach 1.11 scharf transitiv. Mit I.7.9 folgt daher
 $$ \big|\PGL(n,q)\big| = q^{{1 \over 2} n(n-1)} \prod_{i:=2}^n (q^i - 1).$$ 
 \par
       Aus 1.7 und 1.3 folgt, dass $\PGammaL(n,q)/\PGL(n,q)$ zu
 $\Aut(\GF(q))$ isomorph ist. Also ist
 $$ |\PGammaL(n,q)| = \big|\Aut(GF(q))\big|\,\big|\PGL(n,q)\big| 
		    = r\big|\PGL(n,q)\big|.$$
 \par
       Aus 1.2 folgt, dass der Kern $M$ des Homomorphismus der Gruppe $\GL(n,q)$
 auf die Gruppe $\PGL(n,q)$ die Ordnung $q - 1$ hat. Daher ist
 $$ \big|\GL(n,q)\big| = (q - 1)\big|\PGL(n,q)\big|. $$
       Wegen $\GammaL(n,q)/GL(n,q) \cong (\GammaL(n,q)/M)/(\GL(n,q)/M)$ ist
 $$ \big|\GammaL(n,q)\big| = r\big|\GL(n,q)\big|. $$
       Schlie\ss lich folgt aus 1.12
 $$ \big|\GL(n,q)\big| = (q - 1)\big|\SL(n,q)\big|. $$
 \par
       Ist $\sigma$ eine Abbildung aus dem Kern des Homomorphismus von
 $\SL(n,q)$ auf $\PSL(n,q)$, so ist $v^\sigma = vk$ f\"ur alle $v \in V$ und
 einem passenden $k \in K^*$. Es folgt $1 = \det(\sigma) = k^n$. Die Anzahl der
 L\"osungen dieser Gleichungen ist aber gleich
 $\ggT(n,q - 1)$, da die multiplikative Gruppe zyklisch ist. Daher ist
 $$ \big|\PSL(n,q)\big| = \ggT(n,q - 1)^{-1}\big|\SL(n,q)\big|. $$
 Damit ist alles bewiesen.

\mysection{2. Die Einfachheit der kleinen projektiven Gruppe}

\noindent
 In diesem Abschnitt wollen wir zeigen, dass die Gruppen $\PSL(V)$
 bis auf zwei Ausnahmen einfache Gruppen sind, dh., dass sie nur
 die beiden unvermeidbaren Normalteiler haben, die jede Gruppe hat.
 Auf dem Wege dorthin und als Folgerungen daraus werden wir noch
 eine Reihe weiterer interessanter Ergebnisse gewinnen. Zun\"achst
 m\"ussen wir aber unser Vokabular erweitern.
 \par
       Wir haben in Kapitel II h\"aufig den Begriff der Transitivit\"at
 einer Gruppe benutzt, der Begriff der Bahn ist aber noch nicht
 gefallen. Hier seine Definition. Es sei $G$ eine
 Permutationsgruppe auf der Menge $M$. Definiert man auf 
 $M$
 die Relation $\sim$
 durch $x\sim y$ genau dann, wenn es ein $\gamma \in G$ gibt mit
 $x^\gamma = y$, so folgt, dass $\sim$ eine \"Aquivalenzrelation
 ist. Die \"Aquivalenzklassen von $\sim$ hei\ss en {\it Bahnen\/}\index{Bahn}{}
 von $G$ auf $M$. Ist $M$ selbst eine Bahn, so ist das
 gleichbedeutend mit der Transitivit\"at von $G$ auf $M$.
 \par
       Es sei $G$ transitiv auf $M$. Gibt es eine Teilmenge $T$ von $M$,
 die wenigstens zwei Elemente enth\"alt, jedoch von $M$ verschieden
 ist, und gilt, dass f\"ur alle $\gamma \in G$, f\"ur die $T \cap
 T^\gamma \neq \emptyset$ ist, $T = T^\gamma$ ist, so sagen wir,
 dass $G$ auf $M$ {\it imprimitiv\/}\index{imprimitiv}{} operiere und dass $T$
 ein
 {\it Im\-pri\-mi\-ti\-vi\-t\"ats\-ge\-biet\/}\index{Imprimitivit\"atsgebiet}{}
 von $G$ sei. Gibt es kein solches $T$, so operiert $G$ auf $M$ 
 {\it primitiv\/}.\index{primitiv}{}
 \par
       Es sei $N$ ein Normalteiler von $G$ und $T$ sei eine Bahn von $N$.
 Ist $\gamma \in G$, so ist auch $T^\gamma$ eine Bahn von $N$. Aus
 $T \cap T^\gamma \neq \emptyset$ folgt daher, dass $T = T^\gamma$
 ist. Ist $N \neq \{1\}$, so folgt aus der Transitivit\"at von $G$
 auf $M$, dass alle Bahnen von $N$ gleichm\"achtig sind und
 wenigstens zwei Elemente enthalten. Ist $T \neq M$, so ist $T$
 also ein Imprimitivit\"atsgebiet von $G$. Ist $G$ primitiv, so ist
 jeder von $\{1\}$ verschiedene Normalteiler von $G$ transitiv
 auf $M$.
 \par
       Ist $G$ eine Gruppe, so bezeichnen wir mit $G'$ die
 {\it Kommutatorgruppe\/}\index{Kommutatorgruppe}{} von $G$, das ist die von
 allen {\it Kommutatoren\/}\index{Kommutator}{} $a^{-1}b^{-1}ab$ mit $a$,
 $b \in G$ erzeugte
 Untergruppe von $G$. Es folgt, dass $G'$ ein Normalteiler von $G$ ist und
 dass $G/G'$ abelsch ist. Es folgt ferner, dass $G'$ in allen Normalteilern
 $N$ von $G$ enthalten ist, f\"ur die $G/N$ abelsch ist. $G'$ ist also der
 kleinste unter diesen Normalteilern.
 \par
       Die {\it absteigende Kommutatorreihe\/}\index{absteigende Kommutatorreihe}{}
 der Gruppe $G$ ist rekursiv
 de\-fi\-niert durch $G^{(0)} := G$ und $G^{(i+1)} := (G^{(i)})'$. Die
 Gruppe $G$ hei\ss t {\it aufl\"osbar},\index{aufl\"osbar}{} wenn es ein $n$
 gibt mit $G^{(n)} = \{1\}$. Sie hei\ss t
 {\it perfekt\/},\index{perfekte Gruppe}{} wenn $G = G'$ ist.
 \par
       Diese Begriffe werden nun zur Formulierung eines von Iwasawa
 stammenden Satzes benutzt (Iwasawa 1962).\index{Iwasawa, K.}{}
 \medskip\noindent
 {\bf 2.1. Satz.} {\it Ist $G$ eine Permutationsgruppe auf der
 Menge $M$, operiert $G$ auf $M$ primitiv, ist $G$ perfekt und
 enth\"alt der Stabilisator eines Elementes aus $M$ in $G$ einen
 aufl\"osbaren Normalteiler, der zusammen mit seinen Konjugierten
 die Gruppe $G$ erzeugt, so ist $G$ einfach.}
 \smallskip
       Beweis. Es sei $N$ ein nicht trivialer Normalteiler von $G$. Da
 $G$ primitiv ist, operiert $N$ nach obiger Bemerkung auf $M$
 transitiv. Es sei $m \in M$ und $B$ sei der nach Voraussetzung
 existierende Normalteiler von $G_m$, der zusammen mit seinen
 Konjugierten die Gruppe $G$ erzeugt. Es sei $\gamma \in G$. Wegen
 der Transitivit\"at von $N$ auf $M$ gibt es ein $\nu \in N$ mit
 $m^{\gamma\nu} = m$. Weil $B$ in $G_m$ normal ist, folgt
 $(\gamma \nu)^{-1} B(\gamma \nu) = B$. Folglich ist
 $$ \gamma^{-1} B \gamma = \nu B \nu^{-1}. $$
 Somit enth\"alt die Gruppe $NB$ alle zu
 $B$ in $G$ konjugierten Untergruppen. Dies besagt, dass $G = NB$ ist.
 \par
       Es ist also $G = NB^{(0)}$. Es sei $i \leq 0$ und es gelte $G = NB^{(i)}$.
 Dann ist
 $$G/N = (NB^{(i)})/N \cong B^{(i)}/(B^{(i)} \cap N). $$
 Wegen $G'=G$ gibt es kein von
 $\{1\}$ verschiedenes epimorphes Bild von $G/N$, welches abelsch
 ist. Daher ist $B^{(i)} = B^{(i+1)} (B^{(i)} \cap N)$. Hieraus
 folgt weiter
 $$ G = NB^{(i)} = NB^{(i+1)} (B^{(i)} \cap N) = NB^{(i+1)}. $$
 Weil $B$ aufl\"osbar ist, gibt es ein $n$ mit
 $B^{(n)} = \{1\}$. Daher ist $G = N$, womit die Einfachheit von
 $G$ nachgewiesen ist.
 \medskip\noindent
 {\bf 2.2. Satz.} {\it Es sei $V$ ein Vektorraum \"uber dem
 K\"orper $K$. Dann gilt:
 \item{a)} Ist $\Rg_K(V)=2$, so ist $\PSL(V)$ auf der Menge der
 Punkte von $L_K(V)$ zweifach transitiv.
 \item{b)} Ist $\Rg_K(V)\geq 3$, so ist $\PSL(V)$ auf der Menge der
 Tripel $(P,Q,H)$ transitiv, wobei $P$ und $Q$ zwei verschiedene
 Punkte und $H$ eine Hyperebene von $L_K(V)$ ist, die weder $P$
 noch Q enth\"alt. Insbesondere ist $\PSL(V)$ auf der Menge der
 Punkte von $L_K(V)$ zweifach transitiv.\par}
 \smallskip
       Beweis. a) Ist $H$ eine Hyperebene von $V$, was im vorliegenden
 Falle gleichbedeutend damit ist, dass $H$ ein Punkt ist, so
 permutiert die Gruppe $\Tau(H)$ aller Transvektionen mit der Achse
 $H$ die von $H$ verschiedenen Punkte transitiv. Da der Punkt $H$
 kein Fixpunkt von $\PSL(V)$ ist, ist diese Gruppe also auf der
 Menge der Punkte von $L_K(V)$ zweifach transitiv.
 \par
       b) Es seien $(P,Q,H)$ und $(P',Q',H')$ zwei Tripel der verlangten
 Art. Nach II.3.2 operiert $\PSL(V)$ auf der Menge der nicht
 inzidenten Punkt-Hy\-per\-ebe\-nen\-paare transitiv. Wir d\"urfen daher
 annehmen, dass $P = P'$ und $H = H'$ ist. Wir d\"urfen ferner
 annehmen, dass $Q \not\leq Q'$ ist. Dann ist $Q + Q'$ eine Gerade
 von $L_K (V)$.
 \par
       1. Fall: Es ist $P \not \leq Q + Q'$. Wir setzen $R:= (Q+Q') \cap
 H$. Dann ist $R$ ein Punkt, da $Q+Q'$ eine Gerade ist, die nicht
 in $H$ liegt. Nach 1.4 gibt es eine Hyperebene $H^*$ mit $P+R \leq
 H^*$ und $Q$, $Q' \not \leq H^*$. Es gibt also eine Elation $\tau \in
 \E(R,H^*)$ mit $Q^\tau = Q'$. Nun ist $P \leq H^*$ und $R \leq H$.
 Daher ist $P^\tau = P$ und $H^\tau = H$. Damit ist in diesem Falle
 alles bewiesen.
 \par
       2. Fall: Es ist $P \leq Q + Q'$. Wegen $\Rg_K(V) \geq 3$, gibt es
 einen Punkt $S$ mit $S \not\leq Q + Q'$ und $S \not\leq H$. Dann
 erf\"ullen die beiden Tripelpaare $(P,Q,H)$, $(P,S,H)$ und
 $(P,S,H)$, $(P,Q',H)$ die Voraussetzungen des Falles 1, so dass Fall
 2 auf diesen zur\"uckgef\"uhrt ist.
 \par
       Dass $\PSL(V)$ auf der Menge der Punkte von $L_K(V)$ zweifach
 trans\-i\-tiv operiert, folgt schlie\ss lich daraus, dass zwei Punkte
 von $L_K(V)$ nach 1.4 stets ein gemeinsames Komplement besitzen.
 \medskip\noindent
 {\bf 2.3. Satz.} {\it Ist $V$ ein Vektorraum \"uber dem K\"orper
 $K$, so gilt:
 \item{a)} Ist $\Rg_K(V) = 2$, so sind alle von $1$ verschiedenen
 Trans\-vek\-ti\-o\-nen\index{Transvektion}{} in $\GL(V)$ konjugiert.
 \item{b)} Ist $\Rg_K(V) \geq 3$, so sind alle von $1$ verschiedenen
 Transvektionen in $\SL(V)$ konjugiert.\par}
 \smallskip
       Beweis. a) Es seien $\sigma$ und $\tau$ zwei von 1 verschiedene
 Transvektionen. Dann gibt es zwei linear unabh\"angige Vektoren
 $a$ und $b$, so dass
 $$ (ar + bs)^\sigma = a(r + s) + bs, $$
 und zwei linear unabh\"angige Vektoren $c$ und $d$, so dass
 $$ (cr + ds)^\tau = c(r + s) + bs $$
 ist f\"ur alle $r$, $s \in K$. Es gibt
 ferner eine $\gamma \in \GL(V)$ mit $a^\gamma = c$ und $b^\gamma =
 d$. Hiermit folgt
 $$
    (cr + ds)^{\gamma^{-1}\sigma\gamma} = (ar + bs)^{\sigma\gamma}
        = \bigl(a(r + s) + bs\bigr)^\gamma = c(r + s) + ds 
       =(cr + ds)^\gamma. 
 $$
 Also ist $\gamma^{-1}\sigma\gamma = \tau$, womit a) bewiesen
 ist.
 \par
       b) Es seien $\sigma$ und $\tau$ von 1 verschiedene Transvektionen.
 Wie \"ublich stellen wir sie dar durch $x^\sigma = x + a\varphi(x)$ und
 $x^\tau = x + b\psi(x)$. Dann sind $H := \Kern(\varphi)$ und $H' := \Kern(\psi)$
 zwei Hyperebenen. Es seien $P$ und $Q$ Punkte mit $P \not\leq H$ und
 $Q \not\leq H'$. Dann
 ist auch $P^\sigma \not\leq H$ und $Q^\tau \not\leq H'$. Es gibt ein $p \in P$
 mit $\varphi(p) = 1$ und ein $q \in Q$ mit $\psi(q) = 1$.
 Es gilt $p^\sigma = p + a$ und $q^\tau = q + b$.
 \par 
       Nach 2.2 gibt es ein $\gamma \in \SL(V)$ mit $P^\gamma = Q$,
 $P^{\sigma\gamma} = Q^\tau$ und $H^\gamma = H'$. Es folgt
 $p^\gamma = qk$, $p^{\sigma\gamma} = q^\tau m$ und wegen
 $aK = (P + P^\sigma) \cap H$ und $bK = (Q + Q^\tau) \cap H'$ auch
 $a^\gamma = bn$ mit $k$, $m$, $n \in K^*$. Somit ist
 $$
       qk + bn = p^\gamma + a^\gamma = (p + a)^\gamma = p^{\sigma
                                         \gamma} = q^\tau m = (q + b)m 
	       = qm + bm. 
 $$
 Weil $q$ und $b$ linear unabh\"angig sind, folgt $k = m = n$. Ist
 nun $r \in K$ und $h \in H'$, so folgt wegen $h^{\gamma^{-1}} \in H$, dass
 $h^{\gamma^{-1}\sigma \gamma} = h$ ist. Hiermit folgt weiter
 $$\eqalign{
 (qr+h)^{\gamma^{-1}\sigma \gamma} &= q^{\gamma^{-1}\sigma \gamma}r + h
     = \bigl(q^{\gamma^{-1}} + a\varphi(q^{\gamma^{-1}})\bigr)^\gamma r + h  \cr
    &= qr + h + a^\gamma\varphi(pk^{-1})r
     = qr + h + bkk^{-1}r                                                    \cr
    &= qr + h + b\psi(qr + h) = (qr + h)^\tau. \cr} $$
 Also ist $\gamma^{-1}\sigma\gamma = \tau$, so dass $\sigma$ und
 $\tau$, wie behauptet, in $\SL(V)$ konjugiert sind.
 \medskip\noindent
 {\bf 2.4. Satz.} {\it Es sei $V$ ein Vektorraum \"uber dem
 K\"orper $K$ und es gelte
 $ V = u_1K \oplus u_2K \oplus W $
 mit von $0$ verschiedenen Vektoren $u_1$ und $u_2$. Ist dann
 $a \in K^*$, so liegt die durch
 $$ (u_1r + u_2s + w)^\sigma := u_1ar + u_2a^{-1}s + w $$
 definierte Abbildung $\sigma$ in $\SL(V)$.}
 \smallskip
       Beweis. F\"ur $i := 1$, 2 setzen wir $H_i := u_iK + W$. Dann sind
 $H_1$ und $H_2$ Hyperebenen von $L_K(V)$. Wir definieren
 Transvektionen $\tau_1$, $\tau_2$, $\tau_3$ und $\tau_4$ durch
 $$\eqalign{
    (u_2r + h_1)^{\tau_1} &= u_2r + h_1 - u_1r          \cr
    (u_1r + h_2)^{\tau_2} &= u_1r+h_2-u_2(1-a)a^{-1}r   \cr
    (u_2r + h_1)^{\tau_3} &= u_2r + h_1 + u_1ar         \cr
    (u_1r + h_2)^{\tau_4} &= u_1r + h_2 + u_2(1 - a)a^{-2}r, \cr}$$
 wobei die $h_i$ Elemente aus $H_i$ bezeichnen. Eine einfache
 Rechnung zeigt dann, dass $\sigma = \tau_1\tau_2\tau_3\tau_4$ ist.
 \medskip\noindent
 {\bf 2.5. Satz.} {\it Ist $V$ ein Vektorraum \"uber dem K\"orper
 $K$ mit 
 $\Rg_K(V) \geq 2$, so ist $\SL(V)$ perfekt, es
 sei denn, es ist $\Rg_K(V) = 2$ und $|K| = 2$ oder $|K| = 3$.}
 \smallskip
       Beweis. Auf Grund von 2.3 gen\"ugt es zu zeigen, dass die $\SL(V)$
 unter den gemachten Voraussetzungen eine Transvektion enth\"alt,
 die der Kommutator zweier Elemente aus $\SL(V)$ ist.
 \par
       1. Fall: $\Rg_K(V) = 2$. Dann ist $V = uK \oplus vK$. Weil $K$
 mindestens $4$ Elemente enth\"alt, gibt es ein $a \in K$, welches
 von 0, 1 und $-1$ verschieden ist. Nach 2.4 liegt die durch
 $$ (ur + vs)^\sigma := uar + va^{-1}s $$
 definierte Abbildung $\sigma$ in $\SL(V)$. Wir definieren ferner eine
 Transvektion $\tau$ durch
 $$ (ur + vs)^\tau := ur + v(s - r). $$
 Dann ist, wie man leicht nachrechnet,
 $$ (ur + vs)^{\sigma^{-1} \tau^{-1}\sigma \tau} = ur + vs + v(a^{-2} - 1)r. $$
 Dies zeigt, dass $\sigma^{-1}\tau^{-1}\sigma\tau$ eine Transvektion ist, die von
 1 verschieden ist, da ja $a \neq 1$, $-1$ ist.
 \par
       2. Fall:  Es ist $\Rg_K(V) \geq 3$. Es seien $H$ und $H'$ zwei
 verschiedene Hyperebenen von $L_K(V)$ und $\varphi$ seien zwei
 Linearformen auf $V$ mit $H = \Kern(\varphi)$ und $H' = \Kern(\psi)$. Weil der
 Rang von $V$ mindestens gleich 3 ist, gibt
 es ein von Null verschiedenes $b \in H \cap H'$. Es sei ferner $a
 \in H - H'$. Wir definieren die Transvektionen $\sigma$ und $\tau$
 durch $x^\sigma := x + a\varphi(x)$ und $x^\tau := x + b\psi(x)$.
 Eine einfache Rechnung zeigt dann, dass
 $$ x^{\sigma^{-1}\tau^{-1}\sigma \tau} = x + b\psi(a)\varphi(x) $$
 ist. Nun ist $\psi(a) \neq 0$, da ja $a \not\in H' = \Kern(\psi)$
 ist. Also ist $\sigma^{-1}\tau^{-1}\sigma\tau$ eine von 1
 verschiedene Transvektion, so dass der Satz auch im zweiten Falle
 etabliert ist.
 \medskip\noindent
 {\bf 2.6. Satz.} {\it Ist $V$ ein Vektorraum \"uber dem K\"orper
 $K$ mit 
 $\Rg_K(V) \geq 2$, so ist $\PSL(V)$ einfach, es
 sei denn, es ist $\Rg_K(V) = 2$ und $|K| = 2$ oder $|K| = 3$.}
 \smallskip
 Beweis. Es sei $|K| \geq 4$, falls $\Rg_K(V) = 2$ ist. Nach 2.5
 ist $\SL(V)$ dann perfekt, so dass auch $\PSL(V)$ perfekt ist.
 Ferner gilt, dass $\PSL(V)$ auf der Menge der Punkte von $L_K(V)$
 zweifach transitiv operiert, also erst recht primitiv. Ist $P$ ein
 Punkt, so ist die Gruppe $\E(P)$ aller Elationen mit dem Zentrum
 $P$ ein abelscher Normalteiler von $\PSL(V)_P$, der zusammen mit
 seinen Konjugierten die Gruppe $\PSL(V)$ erzeugt. Da abelsche
 Gruppe nat\"urlich aufl\"osbar sind, folgt aus dem Satz 2.1 von
 Iwasawa, dass $\PSL(V)$ einfach ist.
 \par
       Nach 1.13 ist $|\PSL(2,2)| = 2 \cdot 3$ und $|\PSL(2,3)| = 3 \cdot 4$,
 so dass diese Gruppen aufl\"osbar sind. Im Falle der $\PSL(2,2)$
 gibt es n\"amlich $3 \cdot 1$ Elemente der Ordnung 2, die von
 Transvektionen induziert werden. Zusammen mit den 3 Elementen
 einer $3$-Sylowgruppe erh\"alt man alle 6 Elemente der
 $\PSL(2,2)$, so dass die 3-Sylowgruppe von $\PSL(2,2)$ ein
 Normalteiler dieser Gruppe ist. Im Falle der $\PSL(2,3)$ gibt es $4
 \cdot 2 = 8$ Elemente der Ordnung 3, die von Transvektionen
 induziert werden. Zusammen mit den 4 Elementen einer
 2-Sylowgruppe erh\"alt man alle 12 Elemente von $\PSL(2,3)$, so
 dass in diesem Falle die 2-Sylowgruppe normal ist. Somit sind
 diese beiden Gruppen wirkliche Ausnahmen zu Satz 2.6.
 \medskip\noindent
 {\bf 2.7. Satz.} {\it Ist $V$ ein Vektorraum \"uber dem K\"orper
 $K$, ist $\Rg_K(V) \geq 2$ und ist $|K| > 2$, falls $\Rg_K(V) = 2$,
 so gilt:
 \smallskip
  \item{a)} Es ist $\GL(V)' = \SL(V)$.
  \smallskip
  \item{b)} Es ist $\PGL(V)' = \PSL(V)$.}
 \smallskip
       Beweis. Wir beginnen mit einer Vorbemerkung. Es sei $G$ eine
 Gruppe, $S$ ein Normalteiler und $U$ eine Untergruppe von $G$.
 Ferner sei $S$ perfekt und $G = SU$. Da $S$ als Untergruppe von
 $SU'$ diese Gruppe normalisiert und da $S$ und $U'$ von $U$
 normalisiert werden, folgt, dass $SU'$ ein Normalteiler von $SU
 = G$ ist. Daher ist
 $$ G/(SU') = \bigl(SU'U)/(SU'\bigr) \cong U/\bigl(U \cap SU'\bigr). $$
 Wegen $U'\subseteq U \cap SU'$ ist $U/(U \cap S U')$ und damit
 $G/(S U')$ abelsch. Also ist
 $$ G' \subseteq S U' = S'U' \subseteq G' $$
 und folglich $G' = SU'$.
 \par
       Es sei nun $G := \GL(V)$, $S := \SL(V)$ und $U := \sum (P,H)$, wobei $P$
 ein Punkt und $H$ eine Hyperebene von $V$ sei mit $V = P\oplus H$.
 Ist dann $|K| > 3$, falls $\Rg_K(V) = 2$ ist, so folgt mit der
 Vorbemerkung und Satz 2.5, dass
 $$ \GL(V)' = \SL(V)\Sigma(P,H)' $$
 ist. Wir m\"ussen also zeigen, dass $\Sigma(P,H)' \subseteq \SL(V)$ ist. Dazu
 seien $a$, $b \in K^*$. Setze $c := b^{-1}a^{-1}ba$. Wir zeigen, dass
 $\delta(c) \in \SL(V)$ ist.
 \par
       Es sei $P = uK$ und $H = vK \oplus W$ mit $v \neq 0$. Wir definieren
 Abbildungen $\sigma_1$, $\sigma_2$ und $\sigma_3$ durch
 $$\eqalign{
    (ur + vs + w)^{\sigma_1} &:= ubr + vb^{-1}s + w, \cr
    (ur + vs + w)^{\sigma_2} &:= uar + va^{-1}s + w, \cr
    (ur + vs + w)^{\sigma_3} &:= ua^{-1}b^{-1}r + vbas + w \cr} $$
 f\"ur alle $r$, $s \in K$ und alle $w \in W$. Nach 2.4 ist dann
 $\sigma_i \in \SL(V)$ f\"ur alle $i$. Dann ist aber auch
 $\sigma_1\sigma_2\sigma_3 \in \SL(V)$. Eine einfache Rechnung zeigt, dass
 $$ (ur + vs + w)^{\sigma_1\sigma_2\sigma_3} = uc^{-1}r + vs + w
	 = (ur + vs + w)\delta(c) $$
 ist. Folglich ist $\delta (c) = \sigma_1\sigma_2 \sigma_3$, so dass in der Tat
 $$ \Sigma(p,H)' \subseteq \SL(V)$$
 ist. Daher ist $\GL(V)' = \SL(V)$.
 \par
       Es sei nun $|K| = 3$ und $V = uK \oplus vK$. Dann ist $K$ kommutativ.
 Die Gruppe $\sigma (uK,vK)$
 ist zur multiplikativen Gruppe von $K$ isomorph, also
 abelsch. Wegen $\GL(V) = \SL(V)\Sigma(uK,vK)$ ist daher $G' \subseteq \SL(V)$.
 Es sei $\tau$ eine nicht triviale Transvektion
 mit dem Zentrum $vK$. Dann ist $(ur + vs)^\tau = ur + vs + var$ mit
 $a\in K^*$. Ferner ist die durch $(ur + vs)^\sigma := ur - vs$
 definierte Abbildung $\sigma$ ein Element aus $\GL(V)$. Dann ist aber
 $$
    (ur + vs)^{\sigma^{-1}\tau^{-1}\sigma\tau} = ur + vs + 2var
                                                = ur + vs - var     
                                               = (ur + vs)^{\tau^{-1}},  
 $$
 so dass jede Transvektion ein Kommutator in $\GL(V)$ ist. Also ist
 $\SL(V) \subseteq \GL(V)'$, so dass a) bewiesen ist.
 \par
       b) ist eine einfache Folgerung aus a).
 \medskip
       Wir haben gerade gesehen, dass
 $\Sigma(P,H)' \subseteq \SL(V) \cap \Sigma(P,H)$ ist. Wissbegierig wie wir
 sind, stellen wir die Frage, ob vielleicht
 $$ \Sigma(P,H)' = SL(V) \cap \Sigma (P,H) $$
 ist. Im Falle eines kommutativen K\"orpers und
 endlicher Dimension von $V$ ist das ja so, wie wir beim Beweise
 von 1.13 gesehen haben. Diese Gleichheit gilt in der Tat immer,
 wie wir im n\"achsten Abschnitt sehen werden.
 \medskip\noindent
 {\bf 2.8. Satz.} {\it Es sei $V$ ein Vektorraum \"uber dem
 K\"orper $K$ mit 
 $\Rg_K(V) \geq 2$. Ist $V$ nicht der
 Vektorraum vom Range $2$ \"uber $\GF(2)$ oder $\GF(3)$ und ist $N$
 eine von $\{1\}$ verschiedene Untergruppe von $\PGaL(V)$, die von
 $\PSL(V)$ normalisiert wird, so ist $\PSL(V) \subseteq N$.}
 \smallskip
       Beweis. Wir beginnen wieder mit einer Vorbemerkung, die dem, der
 sich mit Permutationsgruppen auskennt, gel\"aufig ist. Es seien
 $M$ und $N$ zwei Normalteiler der Gruppe $G$ mit $M \cap N = \{1\}$. Ist
 $m \in M$ und $n \in N$, so ist $n^{-1}m^{-1}n \in M$ und daher
 $n^{-1}m^{-1}nm \in M$. Andererseits ist $m^{-1}nm \in N$ und daher auch
 $n^{-1}m^{-1}nm \in N$. Also ist
 $$ n^{-1}m^{-1}nm \in M \cap N = \{1\}. $$
 Dies zeigt, dass $M$
 und $N$ sich gegenseitig zentralisieren. Ist \"uberdies $G$ eine
 auf der Menge $\Omega$ primitive Permutationsgruppe, so sind $M$
 und $N$ beide transitiv, falls sie beide von $\{1\}$ verschieden
 sind. Es sei $\alpha \in \Omega$, $m \in M$ und es gelte $\alpha^m = \alpha$.
 Ist dann $n \in N$, so ist
 $$ \alpha^{nm} = \alpha^{mn} = \alpha^n, $$
 so dass auch $\alpha^n$ ein Fixelement
 von $m$ ist. Wegen der Transitivit\"at von $N$ auf $\Omega$ folgt
 daher $m = 1$. Folglich ist $M$ auf $\Omega$ scharf transitiv.
 Ebenso folgt nat\"urlich, dass auch $N$ auf $\Omega$ scharf
 transitiv operiert.
 \par
       Zur\"uck zu unserem Satz. Weil $N$ von $\PSL(V)$ normalisiert wird,
 ist $G := \PSL(V)N$ eine Untergruppe von $\PGaL(V)$ und $\PSL(V)$ und
 $N$ sind Normalteiler dieser Gruppe. Weil $\PSL(V)$ auf der Menge
 der Punkte von $L_K(V)$ nach 2.2 zweifach transitiv operiert,
 operiert sie auf dieser Menge erst recht primitiv. W\"are nun
 $\PSL(V) \not\subseteq N$, so folg\-te $\PSL(V) \cap N = \{1\}$, da
 $\PSL(V)$ unter den gemachten Voraussetzungen ja einfach ist. Nach
 unserer Vorbemerkung folgte dann aber, dass $\PSL(V)$ auf der Menge
 der Punkte von $L_K(V)$ scharf transitiv operierte. Dieser
 Widerspruch zeigt, dass doch $\PSL(V) \subseteq N$ gilt.
 \medskip\noindent
 {\bf 2.9. Korollar.} {\it Es sei $V$ ein Vektorraum \"uber dem
 K\"orper $K$, dessen Rang mindestens $2$ sei. Ferner sei $|K| > 3$,
 falls $\Rg_K(V) = 2$ ist. Ist $N$ eine Untergruppe von $\GaL(V)$, die von
 $\SL(V)$ normalisiert wird, so ist $N$ entweder im
 Zentrum von $\GaL(V)$ enthalten oder aber $N$ enth\"alt $\SL(V)$.}
 \smallskip
       Beweis. Ist $k \in K^*$, so sei $\mu(k)$, wie schon zuvor, die
 durch $v^{\mu(k)} := vk$ definierte Abbildung von $V$ auf sich
 und $M(V)$ bezeichne die Gruppe aller dieser $\mu(k)$. Ist $Z$
 das Zentrum von $\GaL(NV)$, so gilt nach 1.2, dass
 $M(V) \cap \GaL(V) \subseteq Z$ ist. W\"are $M(V) \cap \GaL(V) \neq Z$, so
 induzierte $Z$ einen von $\{1\}$ verschiedenen abelschen
 Normalteiler in $\PGaL(V)$, der nach 2.8 die nicht abelsche Gruppe
 $\PSL(V)$ enthielte. Also ist doch $M(V) \cap \GaL(V) = Z$.
 \par
       Die fragliche Gruppe $N$ sei nicht in $Z$ enthalten. Wegen
 $Z = M(V) \cap \GaL(V)$ induziert $N$ dann eine nicht triviale
 Untergruppe $N'$ in $\PGaL(V)$, die von $\PSL(V)$ normalisiert
 wird. Nach 2.8 ist daher $\PSL(V) \subseteq N'$. Hieraus folgt,
 dass $\SL(V) \subseteq NZ$ ist. Nun ist $(NZ)/N \cong Z/(N \cap
 Z)$, so dass $(NZ)/N$ abelsch ist. Daher ist $(NZ)' \subseteq N$
 und weiter
 $$ \SL(V) = \SL(V)' \subseteq (NZ)' \subseteq N. $$
 Damit ist das Korollar bewiesen.
 \medskip
       Die S\"atze 2.7 und 2.9 geben eine vollst\"andige \"Ubersicht
 \"uber alle Normalteiler von $\GL(V)$. Ist der Rang von $V$ nicht
 endlich, so bleiben noch die Normalteiler\index{Normalteiler}{} von $\GaL(V)$ zu
 bestimmen. Der interessierte Leser findet eine vollst\"andige
 \"Ubersicht \"uber die Normalstruktur von $\GaL(V)$ in Rosenberg 1958.

\mysection{3. Determinanten}

\noindent
 In diesem Abschnitt soll nun das Versprechen eingel\"ost werden, den Beweis
 daf\"ur zu liefern, dass $\Sigma(P,H)' = \SL(V) \cap \Sigma (P,H)$ ist. Dazu
 ben\"otigen wir die von Dieudonn\'e\index{Dieudonn\'e, J.}{}
 eingef\"uhrte Determinantenfunktion\index{Determinante}{} f\"ur
 Endomorphismen\index{Endomorphismus}{} von
 Vektorr\"aumen \"uber nicht notwendig kommutativen K\"orpern. Die
 hier 
 angegebene koordinatenfreie Konstruktion der
 Determinantenfunktion stammt von U.~Demp\-wolff\index{Dempwolff, U.}{}. Es ist zu
 beachten, dass man die Determinantenfunktion ohne Um\-schwei\-fe f\"ur
 jeden, also auch nicht endlich erzeugten Vektorraum $V$ auf
 $\GL(V)$ definieren kann. (Dempwolff 1993, Dieudonn\'e 1943,
 L\"u\-ne\-burg 1993, S.~193ff.)
 \par
       Wir beginnen damit, eine gemeinsame Beschreibung der
 Dilatationen\index{Dilatation}{}
 und Transvektionen\index{Transvektion}{} zu geben. Da wir im Folgenden st\"andig
 Dilatationen und Transvektionen unter dem gleichen Blickwinkel
 betrachten, geben wir der Menge aus Dilatationen und
 Transvektionen eines Vektorraumes den Namen $\Sigma$. Man beachte,
 dass auch $1_V \in \Sigma$ ist.
 \medskip\noindent
 {\bf 3.1. Satz.} {\it Es sei $V$ ein Vektorraum \"uber dem K\"orper $K$ und
 $\Sigma$ sei die Menge der Dilatationen und Transvektionen von $V$.
 \item{a)} Ist $\sigma \in \Sigma$, so gibt es ein $a \in V$ und
 ein $f \in V^*$ mit $f(a) \neq 1$ und
 $$ v^\sigma = v - af(v) $$
 f\"ur alle $v \in V$. Ist $\sigma$ nicht die
 Identit\"at, so ist $a \neq 0$ und auch $f \neq 0$.
 \item{b)} Ist $v^\sigma = v - af(v)$ und $v^\tau = v - a'f'(v)$, und
 ist $\sigma$, $\tau \neq 1_V$, so ist genau dann $\sigma = \tau$,
 wenn es ein $k \in K^*$ gibt mit $a' = ak$ und $f'=k^{-1}f$.\par}

 \smallskip
       Beweis. a) Es sei $\sigma$ eine Transvektion und $H$ sei ihre
 Achse. Es gibt dann {\it per definitionem\/} ein $f \in V^*$ mit
 $\Kern(f) = H$ sowie ein $a \in H$, so dass $v^\tau = v + af(v)$ ist.
 Indem man $a$ durch $-a$ ersetzt, erh\"alt man f\"ur $\tau$ die
 Darstellung
 $$ v^\tau = v - af(v). $$
 \"Uberdies ist $f(a) = 0$ und daher $f(a)\neq 1$.
 \par
       Es sei $\sigma$ eine Dilatation mit der Achse $H$ und dem Zentrum
 $P$. Ferner sei $0 \neq a \in P$. Ist dann $v \in V$, so gibt es
 ein $k \in K$ und ein $h \in H$ mit $v = ak + h$. Es gibt ferner
 ein $b \in K^*$ mit $(ak + h)^\sigma = abk + h$. Definiere $f\in V^*$
 durch
 $$ f(ak + h) := (1 - b)k. $$
 Dann ist
 $$ (ak + h)^\delta = abk + h = ak + h - a(1 - b)k = ak + h - af(ak + h), $$
 m. a. W., es ist
 $$ v^\delta = v - af(v) $$
 f\"ur alle $v \in V$. Weil $\delta$ injektiv ist, ist $f(a) \neq 1$. Die
 restliche Aussage von a) versteht sich von selbst.
 \par
       b) Es sei $\sigma = \tau$. Dann ist $v - af(v) = v - a'f'(v)$ und
 damit $af(v)=a'f'(v)$ f\"ur alle $v\in V$. Weil $\sigma$ nicht die
 Identit\"at ist, ist $f' \neq 0$. Es gibt also ein $v$ mit $f'(v) = 1$. Setze
 $k := f(v)$. Dann ist $k \in K^*$ und es gilt $a' = ak$. Es
 folgt $af(v) = akf'(v)$. Wiederum weil $\sigma$ nicht die
 Identit\"at ist, ist $a \neq 0$. Es folgt $f(v) = kf'(v)$ f\"ur
 alle $v\in V$. Die Umkehrung ist banal. Also gilt auch b).
 \medskip
       Es sei $\sigma \in \Sigma$. Dann gibt es einen Vektor $a \in V$
 und ein $f \in V^*$ mit $x^\sigma = x - af(x)$ f\"ur alle $x \in V$.
 Ist auch $x^\sigma = x - a'f'(x)$, so gibt es nach 3.1 b) ein $k \in K^*$ mit
 $a' = ak$ und $f' = k^{-1}f$. Da $f(a) \neq 1$ ist folgt,
 $$\eqalign{
     1 - f'(a') &= 1 - k^{-1}f(ak) = k^{-1}\bigl(1 - f(a)\bigr)k     \cr
		&= \bigl(1 - f(a)\bigr)\bigl(1 - f(a)\bigr)^{-1}
		       k^{-1}\bigl(1 - f(a)\bigr)k. \cr}$$
 Setzt man nun
 $$ \det(\sigma) := \bigl(1 - f(a)\bigr)(K^*)', $$
 so ist $\det$ also eine wohldefinierte Abbildung von $\Sigma$ in die
 Kommutatorfaktorgruppe $K^*/(K^*)'$.
 \par 
       Ist $\sigma \in \Sigma$ und ist $\tau \in \GL(V)$, ist ferner
 $v^\sigma = v - af(v)$, so folgt
 $$ v^{\tau^{-1}\sigma \tau} = v - a^\tau f(v\tau^{-1}). $$
 Hieraus folgt weiter
 $$ \det(\tau^{-1}\sigma\tau) = \bigl(1 - f(a^{\tau\tau^{-1}})\bigr)(K^*)'
			      = \det(\sigma), $$
 so dass $\det$ unter Konjugation invariant bleibt.
 \par
       Ist $K$ ein K\"orper, so setzen wir im Folgenden
 $$ K_A := K^*/(K^*)'. $$
 Den kanonischen Epimorphismus von $K^*$ auf $K_A$ bezeichnen wir mit $\pi$.
 \medskip\noindent
 {\bf 3.2. Satz.} {\it Es sei $V$ ein Vektorraum \"uber dem
 K\"orper $K$ und $\Sigma$ sei die Menge der Transvektionen und
 Dilatationen von $V$. Sind $\sigma$, $\tau \in \Sigma$ und haben
 $\sigma$ und $\tau$ das gleiche Zentrum oder die gleiche Achse, so
 ist auch $\sigma\tau \in \Sigma$ und es gilt $\det(\sigma\tau)
 = \det(\sigma)\det(\tau)$.}
 \smallskip
       Beweis. Wir betrachten zun\"achst den Fall, dass $\sigma$ und
 $\tau$ das gleiche Zentrum haben. Es gibt dann ein $a \in V$ und
 $f$, $g \in V^*$ mit $v^\sigma = v - af(v)$ und $v^\tau = v - ag(v)$. Es folgt
 $$\eqalign{
      v^{\sigma\tau} &= \bigl(v - af(v)\bigr)^\tau = v^\tau - a^\tau f(v) \cr
	             &= v - ag(v) - \bigl(a - ag(a)\bigr)f(v)             \cr
		     &= v - a\bigl(g(v) + f(v) - g(a)f(v)\bigr). \cr} $$
 Daher ist $\sigma \tau \in \Sigma$ und es gilt
 $$\eqalign{
      \det(\sigma\tau) &= \pi\bigl(1 - g(a) - f(a) + g(a)f(a)\bigr)    \cr
		       &= \pi\bigl(1 - g(a)\bigr)\pi\bigl(1-f(a)\bigr) \cr
		       &= \det(\sigma)\det(\pi), \cr} $$
 da die Multiplikation in $K_A$ ja kommutativ ist.
 \par
       Es bleibt der Fall zu betrachten, dass $\sigma$ und $\pi$ die
 gleiche Achse haben. Dann gibt es $a$, $b \in V$ und ein $f \in V^*$
 mit $v^\sigma = v - af(v)$ und $v^\tau = v - bf(v)$. Es folgt
 $$\eqalign{
      v^{\sigma \tau} &= \bigl(v - af(v)\bigr)^\tau             \cr
		      &= v^\tau -a^\tau f(v)                    \cr
		      &= v - bf(v) - \bigl(a - bf(a)\bigr)f(v)  \cr
		      &= v - \bigl(b + a - bf(a)\bigr)f(v).     \cr}$$
 Also ist $\sigma\tau \in \Sigma$ und
 $$\eqalign{
  \det(\sigma\tau) &= \pi\bigl(1 - f(b + a - bf(a)\bigr)             \cr
		   &= \pi\bigl(1 - f(b) - f(a) + f(b)f(a)\bigr)      \cr
		   &= \pi\bigl(1 - f(b)\bigr)\pi\bigl(1 - f(a)\bigr) \cr
		   &= \det(\sigma)\det(\tau). \cr}$$
 Damit ist der Satz bewiesen.
 \medskip
       Bei diesem Beweis haben wir zweimal von der Kommutativit\"at von
 $K_A$ Gebrauch gemacht. Das ist der Sache nicht inh\"arent, liegt
 vielmehr daran, dass wir die Abbildungen auf die gleiche Seite wie
 die Skalare schreiben, n\"amlich als Exponenten rechts von den
 Vektoren.
 \medskip\noindent
 {\bf 3.3. Satz.} {\it Es sei $V$ ein $K$-Vektorraum und $\Sigma$
 sei die Menge der Transvektionen und Dilatationen von $V$. Sind
 $\sigma_1$, $\sigma_2 \in \Sigma$ und haben $\sigma_1$ und
 $\sigma_2$ verschiedene Zentren $P_1$ und $P_2$ und verschiedene
 Achsen und ist $P$ ein Punkt auf der Geraden $P_1 + P_2$, so gibt
 es $\tau_1$, $\tau_2 \in \Sigma$ mit den folgenden Eigenschaften:
 \item{a)} $\tau_2$ hat Zentrum $P$.
 \item{b)} Es ist $\sigma_1\sigma_2 = \tau_1\tau_2$.
 \item{c)} Es ist $\det(\sigma_1)\det(\sigma_2) = \det(\tau_1)\det(\tau_2)$.}
 \smallskip
       Beweis. Ist $P = P_2$, so tun es $\tau_1 := \sigma_1$ und $\tau_2
 := \sigma_2$. Ist $P = P_1$, so tun es $\tau_1 := \sigma_1\sigma_2\sigma_1^{-1}$
 und $\tau_2 := \sigma_1$.
 \par
       Es sei $P \neq P_1$, $P_2$. Ist $0 \neq a_1 \in P_1$ und $0 \neq a_2
 \in P_2$, so gibt es $f_1$, $f_2 \in V^*$ mit
 $$ v^{\sigma_i} = v - a_if_i(v) $$
 f\"ur $i := 1$, 2 und alle $v \in V$. Setze $H_i := \Kern(f_i)$. Nach 3.2
 h\"angt $H_i$ nur von $\sigma_i$, nicht aber
 von der Wahl von $a_i \in P_i$ ab. Davon werden wir im Folgenden
 Gebrauch machen, indem wir die $a_i$ der Situation entsprechend w\"ahlen.
 \par
       Setze $D := H_1 \cap H_2$ sowie $G := P_1 + P_2$. Nach Voraussetzung
 ist dann
 $$ \Rg_K(G) = 2 = \Ko_K (D). $$
 \par
       Es sei $P_1 \leq H_2$. Hier w\"ahlen wir die $a_i$ so, dass
 $P = (a_1 + a_2) K$ gilt. Hiermit definieren wir $\tau_1$ und $\tau_2$
 durch
 $$ v^{\tau_1} := v - a_1(f_1 - f_2)(v) $$
 bzw.
 $$ v^{\tau_2} := v - (a_1 + a_2)f_2(v). $$
 Dann ist $a_1^{\tau_2} = a_1 = a_1^{\sigma_2}$ und daher
 $$\eqalign{
     v^{\tau_1\tau_2} &= \bigl(v - a_1(f_1(v) - f_2(v)r)\bigr)^{\tau_2}  \cr
                      &= v^{\tau_2} - a_1f_1(v) + a_1f_2(v)  \cr
		      &= v - (a_1 + a_2)f_2(v) - a_1f_1(v) + a_1f_2(v) \cr
		      &= v - a_1f_1(v) - a_2f_2(v)           \cr
		      &= v^{\sigma_1\sigma_2}.               \cr}$$
 Ferner ist
 $$ \det(\tau_1) = \pi\bigl(1 - f_1(a_1) + f_2(a_1)\bigr)
		 = \pi\bigl(1 - f_1(a_1)\bigr)
                 = \det(\sigma_1) $$
 und
 $$ \det(\tau_2) = \pi\bigl(1 - f_2(a_1 + a_2)\bigr)
		 = \pi\bigl(1 - f_2(a_2)\bigr)
		 = \det(\sigma_2) $$
 \par
       Es sei $P_2 \leq H_1$. Wegen $\sigma_1\sigma_2 = \sigma_1
 \sigma_2 \sigma_1^{-1}\sigma_1$ erledigt sich dies mit dem gerade
 behandelten Fall.
 \par 
       Wir d\"urfen daher im Folgenden annehmen, dass $P_1 \not\leq H_2$ und
 $P_2 \not\leq H_1$ gilt. Dann ist $P_1$, $P_2 \not\leq D$ und folglich
 $$ \Rg_K(G \cap D) \leq 1. $$
 \par
       Es sei $C$ ein Komplement von $G \cap D$ in $D$. Weil $\sigma_1$
 und $\sigma_2$ auf $C$ die Identit\"at induzieren, d\"urfen wir
 annehmen, dass $C = \{0\}$ ist. Dann ist also $D \leq G$ und
 $\Rg_K(D) \leq 1$.
 \par
       1. Fall: Es ist $\Rg_K(D) = 1$. Dann ist $\Rg_K(V) = 3$. Wegen
 $P_1$, $P_2 \neq D$ k\"onnen wir die $a_i$ so w\"ahlen, dass
 $$ D = (a_1 + a_2)K $$
 ist. Setze $b_1 := a_1 + a_2$ und $b_2 := a_1$. Dann
 ist $b_1$, $b_2$ eine Basis von $G$. Wegen $P \neq P_1 = b_2K$ gibt
 es also ein $x \in K$ mit
 $$ P = (b_1 + xb_2)K. $$
 Es sei $b_3 \in H_2 - D$. Dann ist $b_3 \not\in G$, da $g\cap H_2 = D$ ist.
 Folglich ist $\{b_1, b_2, b_3\}$ eine Basis von $V$. Ferner gilt,
 man beachte, dass $b_1 \in D = H_1 \cap H_2$ ist,
 $$ b_1^{\sigma_1} = b_1 $$
 und, mit implizit definierten $u$, $v \in K$,
 $$\eqalign{
    b_2^{\sigma_1} &= a_1 - a_1 f_1 (a_1)    \cr
		   &= a_2\bigl(1 + f_1(a_2)\bigr) - (a_1 + a_2)f_1(a_2) \cr
		   &= b_2\bigl(1 - f_1(b_2)\bigr) - b_1f_1(b_2)         \cr
		   &= b_1u + b_2v, \cr} $$
 sowie
 $$ b_3^{\sigma_1} = b_3. $$
 Entsprechend erhalten wir
 $$ b_1^{\sigma_2} = b_1 $$
 und, mit wiederum implizit definierten $y$, $z \in K$,
 $$ b_2^{\sigma_2} = a_2 - a_2f_2(a_2) = b_2\bigl(1-f_2(b_2)\bigr)=b_2y $$
 sowie
 $$ b_3^{\sigma_2} = b_3 - a_2f_2(b_3) = b_3 - b_2f_2(b_3) = b_2z + b_3. $$
 Es folgt
 $$\eqalign{
     b_1^{\sigma_1\sigma_2} &= b_1          \cr
     b_2^{\sigma_1\sigma_2} &= b_1u + b_2yv \cr
     b_3^{\sigma_1\sigma_2} &= b_2z+b_3.    \cr}$$
 Wir setzen $s := 1$, falls $x = 0$, dh., falls $P = b_1K$ ist, und $s := x^{-1}$,
 falls $x \neq 0$ ist. Wir definieren nun Abbildungen
 $\tau_1$ und $\tau_2$ durch
 $$\eqalign{
       b_1^{\tau_1} &:= b_1                       \cr
       b_2^{\tau_1} &:= b_1\bigl(u + s(1 - yv)\bigr) + b_2  \cr
       b_3^{\tau_1} &:= -b_1sz + b_3 \cr} \qquad
 \mathrm{und}\qquad
\eqalign{
       b_1^{\tau_2} &:= b_1    \cr
       b_2^{\tau_2} &:= b_1s(yv - 1) + b_2yv   \cr
       b_3^{\tau_2} &:= b_1sz + b_2z + b_3.    \cr}
 $$
 Einfache Rechnungen zeigen, dass $b_i^{\sigma_1\sigma_2} = b_i^{\tau_1\tau_2}$
 ist f\"ur alle $i$. Somit ist $\sigma_1\sigma_2 = \tau_1\tau_2$.
 \par
       Ist $v = b_1x_1 + b_2x_2 + b_3x_3$, so folgt
 $$ v^{\tau_1} = v - b_1\bigl(szx_3 - (u + s(1 - yv))x_2\bigr) $$
 und
 $$ v^{\tau_2} = v - (b_1s + b_2)\bigl((1 - yv)x_2 - zx_3)\bigr), $$
 so dass $\tau_i \in \Sigma$ gilt.
 Ferner ist $\det(\tau_1) = \pi(1)$ und $\det(\tau_i) = \pi(yv)$. Es
 bleibt die Aussage \"uber das Zentrum von $\tau_2$ zu beweisen.
 \par
       Ist $x \neq 0$, so ist $s = x^{-1}$ und daher
 $$ (b_1s + b_2)K = (b_1 + b_2x)K = P. $$
 Ist $x = 0$, so ist $s = 0$ und $P$ ist
 das Zentrum von $\tau_1$. Dann folgt aber mit $\tau_1\tau_2 =
 \tau_1\tau_2\tau_1^{-1}\tau_1$ und Ersetzen von $\tau_1$ und
 $\tau_2$ durch $\tau_1\tau_2\tau^{-1}$ und $\tau_1$ die Behauptung.
 \par
       2. Fall: Es ist $D = \{0\}$. In diesem Falle w\"ahlen wir die $a_i$
 so, dass $P=(a_1 - a_2)K$ ist. Wir definieren $u$, $v$, $x$, $y \in K$ durch
 $$\eqalign{
       a_1^{\sigma_1} &= a_1\bigl(1 - f_1(a_1)\bigr) = a_1u       \cr
       a_2^{\sigma_1} &= a_2 - a_1f_1(a_2) = a_1v + a_2 \cr}$$
 und
 $$\eqalign{
       a_1^{\sigma_2} &= a_1 - a_2f_2(a_1) = a_1 + a_2x \cr
       a_2^{\sigma_2} &= a_2\bigl(1 - f_2(a_2)\bigr) = a_2y.      \cr} $$
 Es folgt
 $$\eqalign{
     a_1^{\sigma_1\sigma_2} &= a_1u + a_2xu         \cr
     a_2^{\sigma_1\sigma_2} &= a_1v + a_2 (xv + y). \cr} $$
 Es ist $\det(\sigma_1) = \pi(u)$ und $\det(\sigma_2) = \pi(y)$. Wir
 haben nun zwei F\"alle zu unterscheiden.
 \par
       2.1. Fall: Es ist $u = v$. In diesem Falle definieren wir $\tau_1$
 und $\tau_2$ durch
 $$\eqalign{
      a_1^{\tau_1} &:= a_1(1 - xu) + a_2xu           \cr
      a_2^{\tau_1} &:= a_1(1 - xu - y) + a_2(xu + y) \cr} $$
 und
 $$\eqalign{
      a_1^{\tau_2} &:= a_1u         \cr
      a_2^{\tau_2} &:= a_1 u + a_2. \cr} $$
 Eine einfache Rechnung zeigt, dass
 $$ a_1^{\tau_1\tau_2} = a_1u + a_2xu = a_1^{\sigma_1\sigma_2} $$
 und
 $$ a_2^{\tau_1\tau_2} = a_1u + a_2(xu + y) $$
 ist. Wegen $u = v$ ist also auch
 $$ a_2^{\tau_1\tau_2} = a_2^{\sigma_1\sigma_2}, $$
 so dass $\sigma_1\sigma_2 = \tau_1\tau_2$ ist. Ist $v = a_1k_1 + a_2k_2$,
 so folgt
 $$ v_1^\tau = v - (a_1 - a_2)\bigl(xuk_1 + (y + xu - 1)k_2\bigr) $$
 und
 $$ \det(\tau_1) = \pi\bigl(1 - xu - (y + xu - 1)(-1)\bigr) = \pi(y) $$
 sowie
 $$ v^{\tau_2} = v - a_1\bigl((1 - u)k_1 - uk_2\bigr) $$
 und
 $$ \det(\tau_2) = \pi(1 - 1 + u) = \pi(u). $$
 Weil die Multiplikation in $K_A$ kommutativ ist, ist also
 $$ \det(\sigma_1)\det(\sigma_2) = \det(\tau_1)\det(\tau_2). $$
 Ferner ist $P$ das Zentrum von $\tau_1$, woraus mit der schon mehrfach
 angewandten Identit\"at $\tau_1\tau_2 = (\tau_1\tau_2)\tau_1^{-1}\tau_1$ die
 Behauptung in diesem Falle folgt.
 \par
       2.2. Fall: Es ist $u \neq v$. Hier definieren wir $\tau_1$ und
 $\tau_2$ durch
 $$\eqalign{
       a_1^{\tau_1} &:= a_1u + a_2(1 - u)    \cr
       a_2^{\tau_1} &:= a_1v + a_2(1 - v)    \cr} $$
 und
 $$\eqalign{
       a_1^{\tau_2} &:= a_1 + a_2\bigl(x + y(u - v)^{-1}(u - 1)\bigr)   \cr
       a_2^{\tau_2} &:= a_2y(u - v)^{-1}u. \cr}$$
 Dann ist wieder $\sigma_1\sigma_2 = \tau_1\tau_2$. Ferner gilt,
 falls $v = a_1k_1 + a_2k_2$ ist,
 $$ v^{\tau_1} = v - (a_1 - a_2)\bigl((1 - u)k_1 - vk_2\bigr) $$
 und daher
 $$ \det(\tau_1) = \pi\bigl(1 - (1 - u) - v(-1)\bigr) = \pi(u-v) $$ 
 sowie
 $$ v^{\tau_2}
	   = v - a_2\bigl((x + y(u - v)^{-1}(1 - u)k_1
		 - (y(u - v)^{-1}u + 1)k_2)\bigr) $$
 und folglich
 $$ \det(\tau_2) = \pi\bigl(1 + (y(u - v)^{-1}u + 1)(-1)\bigr)
		 = \pi\bigl(y(u - v)^{-1}u\bigr). $$
 Der Punkt $P$ ist wieder Zentrum von $\tau_1$. Aus all diesem
 folgt auch in diesem letzten Fall die Behauptung des Satzes, wobei
 wieder zu beachten ist, dass $K_A$ eine abelsche Gruppe ist.
 \medskip\noindent
 {\bf 3.4. Satz.} {\it Es sei $V$ ein Vektorraum \"uber dem
 K\"orper $K$ und $\Sigma$ sei die Menge der Dilatationen und
 Transvektionen von $V$. Sind $\sigma_1$, \dots, $\sigma_n \in
 \Sigma$ und gilt
 $$ \sigma_1 \cdots \sigma_n = 1, $$
 so ist}
 $$ \det(\sigma_1) \cdots \det(\sigma_n) = 1.$$
 \par
       Beweis. Wir machen Induktion nach $n$. F\"ur $n = 1$ ist der Satz
 richtig. Es sei $n = 2$. Dann ist also $1 = \sigma_1\sigma_2$. Es
 folgt, dass $\sigma_1$ und $\sigma_2$ gleiches Zentrum und gleiche
 Achse haben. Nach 3.2 ist daher
 $$ 1 = \det(1) = \det(\sigma_1)\det(\sigma_2). $$
 Es sei $n \geq 3$ und der Satz gelte f\"ur $n - 1$. Ferner sei
 $$ \sigma_1 \cdots \sigma_n = 1. $$
 Gibt es ein $i$, so dass $\sigma_i$ und
 $\sigma_{i+1}$ gleiches Zentrum oder gleiche Achse haben, so ist
 $\sigma_i \sigma_{i+1} \in \Sigma$ und
 $\det(\sigma_i\sigma_{i+1}) = \det(\sigma_i)\det(\sigma_{i+1})$. Nun ist
 $$ \sigma_1 \cdots (\sigma_i\sigma_{i+1}) \cdots \sigma_n = 1. $$
 Nach Induktionsannahme ist daher
 $$ \det(\sigma_1) \cdots \det(\sigma_n)
   = \det(\sigma_1) \cdots \det(\sigma_i \sigma_{i+1})\cdots \det(\sigma_n) = 1.
                         $$
 Es bleibt der Fall zu betrachten, dass
 $\sigma_i$ und $\sigma_{i+1}$ f\"ur $i := 1$, \dots, $n - 1$
 verschiedene Zentren und verschiedene Achsen haben. Es sei $P_i$
 das Zentrum von $\sigma_i$. Dann ist
 $$ \sigma_n^{-1} = \sigma_1 \cdots \sigma_{n-1} $$
 und $P_n$ ist auch das Zentrum von $\sigma_n^{-1}$. Setze
 $$ H := P_1 + \dots + P_{n-1}. $$
 Weil $\sigma_i$ auf $V/P_i$ die Identit\"at induziert, induziert
 $\sigma_i$ auch auf $V/H$ die Identit\"at. Somit induziert auch
 $\sigma_n^{-1}$ auf $V/H$ die Identit\"at, da $\sigma_n^{-1}$ ja
 das Produkt von $\sigma_1$, \dots, $\sigma_{n-1}$ ist. Es folgt
 $$ P_n = V^{\sigma_n^{-1} -1} \leq H, $$
 da $P_n$ ja auch das Zentrum von $\sigma_n^{-1}$ ist.
 \par
       Es sei $r$ eine nat\"urliche Zahl. Ferner seien $\tau_1$, \dots,
 $\tau_r \in \Sigma$ und $Q_r$ sei das Zentrum von $\tau_r$.
 Schlie\ss lich gelte
 $$ \sigma_n^{-1} = \tau_1 \cdots \tau_r\sigma_{r+1} \cdots \sigma_{n-1} $$
 und
 $$ P_n \leq Q_r + P_{r+1} + \dots + P_{n-1} $$
 sowie
 $$ \det(\sigma_1) \cdots \det(\sigma_{n-1})
        = \det(\tau_1) \cdots \det(\tau_r)\det(\sigma_{r+1}) \cdots
                \det(\sigma_{n-1}). $$
 F\"ur $r = 1$ erh\"alt man diese Situation,
 indem man $\tau_1 = \sigma_1$ und $Q_1 := P_1$ setzt, wie wir gerade
 gesehen haben. Es sei $1 \leq r < n - 1$ und es m\"ogen $\tau_1, \dots,
 \tau_r$ mit den verlangten Eigenschaften geben. Ist $Q_r = P_{r+1}$
 oder ist die Achse von $\tau_r$ gleich der Achse von
 $\sigma_{r+1}$, so schlie\ss t man wie zuvor, dass
 $$ \det(\sigma_1) \cdots \det(\sigma_n) = 1 $$
 ist. Es sei also $Q_r \neq P_{r+1}$ und auch die Achsen von $\tau_r$ und
 $\sigma_{r+1}$ seien verschieden.
 \par 
       Es sei $P_i = a_iK$ f\"ur $i := 1$, \dots, $n$ und $Q_r = b_rK$. Wegen
 $P_n \leq Q_r + P_{r+1} + \dots + P_{n-1}$ gibt es $k_r$, \dots,
 $k_{n-1} \in K$ mit
 $$ a_n = b_rk_r + a_{r+1}k_{r+1} + \dots + a_{n-1}k_{n-1}. $$
 Ist $b_rk_r + a_{r+1}k_{r+1} \neq 0$, so setzen wir
 $$ Q_{r+1} := (b_rk_r + a_{r+1}k_{r+1})K. $$
 Ist $b_rk_r + a_{r+1}k_{r+1} = 0$, so sei $Q_{r+1}$ ein beliebiger Punkt
 auf der Geraden $Q_r + P_{r+1}$. In jedem Falle gilt
 $$ P_n \leq Q_{r+1} + P_{r+2} + \dots + P_{n-1}. $$
 Nach Satz 3.3 gibt es $\tau'_r$, $\tau_{r+1} \in \Sigma$ mit
 $$ \tau_r\sigma_{r+1} = \tau'_r\tau_{r+1}, $$
 so dass $Q_{r+1}$ das Zentrum von $\tau_{r+1}$ ist und \"uberdies
 $$ \det(\tau_r)\det(\sigma_{r+1}) = \det(\tau'_r)\det(\tau_{r+1}) $$
 gilt. Ersetzt man $\tau_r$ durch $\tau'_r$, so hat man die
 Induktion um einen Schritt weitergetrieben oder aber erkannt, dass
 der Satz auch f\"ur $n$ gilt. Wir d\"urfen daher annehmen, dass es
 $\tau_1$, \dots, $\tau_{n-1} \in \Sigma$ gibt mit
 $$ \sigma_n = \tau_1 \cdots \tau_{n-1} $$ 
 und 
 $$ \det(\sigma_1) \cdots \det(\sigma_{n-1})=\det(\tau_1) \cdots
	  \det(\tau_{n-1}), $$
 so dass $P_n$ auch Zentrum von $\tau_{n-1}$ ist. Setze
 $\rho:=\tau_{n-1}\sigma_n$. Nach Satz 3.2 ist $\rho \in \Sigma$
 und es gilt $\det(\rho) = \det(\tau_{n-1})\det(\sigma_n)$. Es folgt
 $$ 1 = \tau_1 \cdots \tau_{n-2} \rho $$
 und daher
 $$ 1 = \det(\tau_1) \cdots \det(\tau_{n-2}) \rho = \det(\sigma_1)
          \cdots \det(\sigma_n). $$
 Damit ist der Satz bewiesen.
 \medskip\noindent
 {\bf 3.5. Satz.} {\it Es sei $V$ ein Vektorraum \"uber dem K\"orper
 $K$ und $\Sigma$ sei die Menge der Dilatationen und Transvektionen
 von $V$. Ist $\gamma \in \GL(V)$, so gibt es $\sigma_1, \dots,
 \sigma_n \in \Sigma$ mit $\gamma = \sigma_1 \cdots \sigma_n$.
 Setzt man
 $$ \det(\gamma) := \det(\sigma_1) \cdots \det(\sigma_n), $$
 so ist $\det$ ein Homomorphismus von $\GL(V)$ auf
 $K^*/(K^*)'$ und es gilt, falls $\Rg_K(V) \geq 2$ ist,}
 $$ \Kern(\det) = \SL(V).$$
 \par
       Beweis. Mit 3.4 folgt, dass $\det$ wohldefiniert ist. Dass $\det$
 dann ein Epimorphismus ist, ist banal. Es bleibt die Aussage
 \"uber den $\Kern$ von $\det$ zu beweisen.
 \par
       Ist $\tau$ eine Transvektion, so ist $\det(\tau) = 1$, wie wir
 wissen. Daher gilt $\SL(V) \subseteq \Kern(\det)$.
 \par
       Es sei $\alpha \in \Kern(\det)$. Ferner sei $V = P \oplus H$ mit
 einem Punkt $P$ und einer Hyperebenen $H$. Dann ist, wie wir
 wissen, $\GL(V) = \SL(V)\Sigma(P,H)$. Es gibt also ein $\sigma \in \SL(V)$ und
 ein $k\in K^*$ mit $\alpha = \sigma\delta (k)$. Es
 folgt
 $$ \pi(1) = \det(\alpha) = \det(\sigma)\det\bigl(\delta(k)\bigr) = \pi(k). $$
 Dies besagt, dass $k \in (K^*)'$, bzw., dass
 $\delta (k) \in \Sigma(P,H)'$ ist. Ist $|K| = 2$, so folgt
 $\delta (k) = 1$ und damit $\alpha = \sigma \in \SL(V)$. Ist $|K| > 2$, so ist
 $\GL(V)' = \SL(V)$ nach 2.7 a). Also gilt $\delta(k)
 \in \SL(V)$ und dann auch $\alpha \in \SL(V)$. Damit ist alles
 bewiesen.
 \medskip
       Nun k\"onnen wir endlich unser Versprechen einl\"osen.
 \medskip\noindent
 {\bf 3.6. Satz.} {\it Ist $V$ ein Vektorraum \"uber dem K\"orper
 $K$ und ist der Rang von $V$ mindestens gleich $2$, ist ferner $P$
 ein Punkt und $H$ eine Hyperebene von $L_K(V)$ mit $V = P \oplus H$,
 so ist}
 $$ \SL(V) \cap \Sigma(P,H) = \Sigma(P,H)'. $$
 \par
       Beweis. Dies ist richtig, falls $|K| = 2$ ist. Es sei also $|K| > 2$.
 Dann ist $\Sigma(P,H)' \subseteq \SL(V) \cap \Sigma(P,H)$, wie
 wir bereits wissen. Es sei $\delta(k) \in \SL(V) \cap
 \Sigma(P,H)$. Es gibt dann Transvektionen $\tau_1$, \dots, $\tau_n$
 mit 
 $$ \delta(k) = \tau_1 \cdots \tau_n. $$
 Es folgt
 $$ \pi (k) = \det\bigl(\delta(k)\bigr) = \pi(1) $$
 ist. Also ist $k \in (K^*)'$ und damit $\delta(k) \in \Sigma(P,H)'$, so dass
 auch dieser Satz bewiesen ist.
 \medskip
       Nach Cohn (1977, S. 117) gibt es nicht kommutative K\"orper, deren
 innere Automorphismengruppe auf $K^* - \{1\}$ transitiv operiert.
 Dies hat zur Folge, dass $(K^*)' = K^*$ ist. Ist $V$ ein Vektorraum
 \"uber einem solchen K\"orper, so ist nach unseren Entwicklungen also
 $\GL(V) = \SL(V)$.
 \par
       Die hier konstruierte Funktion $\det$ stimmt\index{Determinante}{} im
 Falle endlich er\-zeug\-ter Vektorr\"aume mit der in der sonstigen Literatur
 $\det$ genannten Funktion \"uberein, da die Einschr\"ankungen beider
 Funktionen auf die Men\-ge der Transvektionen und Dilatationen
 \"ubereinstimmen (siehe etwa L\"u\-ne\-burg 1993).

\mysection{4. Ausnahmeisomorphismen}

\noindent
 F\"ur diesen Abschnitt wird unterstellt, dass der Leser mit den
 Grundbegriffen der Theorie endlicher Gruppen vertraut ist. Unser
 Hauptziel ist, den folgenden Satz zu beweisen.\index{Ausnahmeisomorphien}{}
 \medskip\noindent
 {\bf 4.1. Satz.} {\it Die folgenden Gruppen sind isomorph:
 \smallskip
 \item{a)} $\PSL(2,2)$ und $S_3$.
 \smallskip
 \item{b)} $\PSL(2,3)$ und $A_4$.
 \smallskip
 \item{c)} $\PSL(2,4)$, $\PSL(2,5)$ und $A_5$ (H\"older 1892).
 \smallskip
 \item{d)} $\PSL(2,7)$ und $\PSL(3,2)$ (H\"older 1892).
 \smallskip
 \item{e)} $\PSL(2,9)$ und $A_6$.
 \smallskip
 \item{f)} $\PSL(4,2)$ und $A_8$. (Jordan 1870/1989, S. 380 ff.)}
 \medskip
       Um a) zu beweisen, hat man nur zu beachten, dass $\PSL(2,2)$ eine
 Permutationsgruppe vom Grade 3 ist, die nach 1.13 die Ordnung
 6 hat.
 \par
       Um b) zu beweisen, bemerkt man, dass $\PSL(2,3)$ eine
 Permutationsgruppe vom Grade 4 ist, die von Transvektionen
 erzeugt wird. Da die Transvektionen in diesem Falle 3-Zyklen
 sind, folgt, dass $\PSL(2,3)$ zu einer Untergruppe der $A_4$
 isomorph ist. Schlie\ss lich folgt aus 1.13, dass
 $|\PSL(2,3)| = 12 = |A_4|$ ist. Dies zeigt, dass $\PSL(2,3)$ und $A_4$
 isomorph sind.\index{alternierende Gruppen}{}
 \par
       Um c) zu beweisen, bemerken wir zun\"achst, dass $\PSL(2,4)$ eine
 einfache Permutationsgruppe vom Grade 5 ist. Da eine einfache
 Permutationsgruppe\index{einfache Permutationsgruppe}{} keine ungerade
 Permutation\index{ungerade Permutation}{} enthalten kann,
 folgt, dass $\PSL(2,4)$ zu einer Untergruppe der $A_5$ isomorph
 ist. Wegen $|\PSL(2,4)| = 60 = |A_5|$ folgt daher, dass $\PSL(2,4)$ und
 $A_5$ isomorph sind. Damit ist eine der beiden Aussagen von c)
 bewiesen. Dar\"uber hinaus ist gezeigt, dass $A_5$\index{alternierende Gruppen}{}
 einfach ist. Nun
 ist $|\PSL(2,5)| = 60$ und $\PSL(2,5)$ ist einfach. Daher ist c) eine
 Konsequenz des folgenden allgemeineren Satzes.
 \medskip\noindent
 {\bf 4.2. Satz.} {\it Es gibt bis auf Isomorphie nur eine einfache
 Gruppe der Ordnung $60$ (H\"older 1892).}
 \smallskip
       Beweis. Es sei $G$ eine einfache Gruppe der Ordnung 60. Enth\"alt
 $G$ eine Untergruppe vom Index $5$, so besitzt $G$ wegen der
 Einfachheit eine treue Darstellung vom Grade 5 und es folgt,
 dass $G$ zur $A_5$ isomorph ist. Es ist also zu zeigen, dass $G$
 eine solche Untergruppe enth\"alt. Dazu zeigen wir zun\"achst,
 dass $G$ keine Untergruppe vom Index 3 enth\"alt. Enthielte $G$
 eine solche Gruppe, so folgte aus der Einfachheit von $G$, dass
 $G$ zu einer Untergruppe der $S_3$ isomorph w\"are, was wegen
 $|G| = 60 > 6 = |S_3|$ nicht sein kann. Es seien $S$ und $T$ zwei
 verschiedene 2-Sylowgruppen von $G$. Ferner sei $1 \neq s \in  S \cap T$. Die
 Ordnung einer 2-Sylowgruppe von $G$ ist gleich 4.
 Daher sind $S$ und $T$ abelsch und es folgt, dass sie beide im
 Zentralisator $C$ von $S$ liegen. Daher gilt $|C| \geq 3 \cdot 4$.
 Weil $G$ einfach ist, ist $C$ von $G$ verschieden. Also ist
 $|C| = 12$ oder 20. Weil $G$ keine Untergruppe vom Index 3
 besitzt, ist also $|C| = 12$, so dass $C$ eine Untergruppe vom
 Index 5 ist. Wir d\"urfen daher annehmen, dass je zwei
 2-Sylowgruppen von $G$ trivialen Schnitt haben.
 \par 
       Es seien wieder $S$ und $T$ zwei verschiedene $2$-Sylowgruppen von
 $G$. Ferner sei $s \in S$ und es gelte $s^{-1}Ts = T$. Dann
 erzeugen $T$ und $s$ eine 2-Gruppe, so dass $s \in T$ folgt. Auf
 Grund unserer Annahme ist daher $s = 1$. Hieraus folgt, dass die
 Anzahl der 2-Sylowgruppen von $G$ kongruent 1 modulo 4 ist.
 Da diese Anzahl aber auch ein Teiler von $|G| = 60$ ist und wegen
 der Einfachheit von $G$ nicht 1 sein kann, folgt, dass $G$ genau
 f\"unf 2-Sylowgruppen hat, so dass der Normalisator einer solchen
 den Index 5 hat. Damit ist 4.2 bewiesen.
 \medskip
       Da $|\PSL(2,7)| = 168 = |\PSL(3,2)|$ ist, folgt d) aus
 \medskip\noindent
 {\bf 4.3. Satz.} {\it Es gibt bis auf Isomorphie nur eine einfache
 Gruppe der Ordnung $168$ (H\"older 1892).}
 \smallskip
       Beweis. Es sei $G$ eine einfache Gruppe der Ordnung
 $168 = 8 \cdot 3 \cdot 7$. Die Anzahl der 7-Sylowgruppen von $G$ ist kongruent
 1 modulo 7 und ein Teiler von 168. \"Uberdies ist sie
 gr\"o\ss er als 1, da $G$ einfach ist. Es folgt, dass $G$ genau
 acht 7-Sylowgruppen enth\"alt. Es sei $N$ der Normalisator einer
 solchen. Dann ist also $|N| = 3 \cdot 7$. Stellt man $G$ dar als
 Permutationsgruppe auf der Menge der Rechtsrestklassen nach $N$,
 so ist diese Darstellung treu, da $G$ einfach ist. Somit hat $G$ eine
 Darstellung als transitive Gruppe von Grade $8$. Der Stabilisator
 von $N$ in dieser Darstellung ist $N$ selber und hat daher die
 Ordnung $21$. Hieraus folgt, dass $G$ sogar zweifach transitiv
 ist.
 \par
       Wir bezeichnen die Rechtsrestklassen von $N$ mit $P_\infty$, $P_0$,
 $P_1$, \dots, $P_6$, wobei wir annehmen, dass $N = P_\infty$ ist. Es
 sei $\sigma$ ein Element der Ordnung 7 in $N$. Indem wir die
 $P_i$ gegebenenfalls umnummerieren, k\"onnen wir erreichen, dass
 $P_i^\sigma = P_{i+1}$ ist, wobei die Indizes modulo $7$ zu
 reduzieren sind. Es sei $\nu \in N$. Dann ist $\nu^{-1} \sigma \nu = \sigma^a$
 mit einer Zahl $a$, die zwischen $0$ und $6$ liegt. Definiere $f$ durch
 $P^\nu_i = P_{f(i)}$. Dann ist
 $$ P_{f(i+1)} = P^\nu_{i+1} = P^{\sigma \nu}_i = P^{\nu \sigma^a}
	       = P_{f(i)+a}.$$
 Also ist $f(i+1) = f(i) + a$. Setzt man $b := f(0)$, so erh\"alt man,
 dass $f(i) = ai + b$ ist. Schlie\ss lich folgt aus $|N| = 21$, dass
 $\nu^3$ in der von $\sigma$ erzeugten Gruppe liegt. Dies
 impliziert wiederum die Kongruenz $a^3 \equiv 1 \mod 7$. Hat
 $\nu$ zwei von $P_\infty$ verschiedene Fixpunkte $P_i$ und $P_j$,
 so ist $i \neq j$ und $i = ai + b$ sowie $j = aj + b$, woraus $a = 1$ und
 $b = 0$ folgt. Also induziert $\nu$ die Identit\"at in
 $\{P_0, \dots, P_\infty \}$, woraus $\nu = 1$ folgt, da $G$ ja treu
 operiert. Damit ist gezeigt, dass $N$ aus allen Abbildungen $\nu$
 mit $P_\infty^\nu = P_\infty$ und $P_i^\nu = P_{ai + b}$ besteht,
 wobei $a$, $b \in \GF(7)$ und $a^3 = 1$ ist. Schlie\ss lich folgt
 noch, dass die Identit\"at die einzige Permutation aus $G$ ist,
 welche drei verschiedene der Punkte $P_i$ festl\"asst.
 \par
       Setze $K := \GF(7)$. Es sei $V = uK \oplus vK$ ein Vektorraum vom Rang
 $2$ \"uber $K$. Setzt man $P_\infty := vK$ und $P_i := (u + vi)K$, so
 wird $G$ zu einer Permutationsgruppe auf der Menge der Punkte von
 $L_K(V)$. Ferner sieht man, dass die Gruppe, die von allen
 linearen Abbildungen $\lambda$ der Form $u^\lambda = ua + vb$ und
 $v^\lambda = va^{-1}$ mit $a$, $b  \in K$ und $a^3 = 1$ auf der Menge
 der Punkte von $L_K(V)$ induziert wird, gleich $N$ ist. Wegen
 $\det(\lambda) = aa^{-1} = 1$ ist daher $N \subseteq \PSL(V)$.
 \par
       Die $3$-Sylowgruppen von $G$ sind gerade die Stabilisatoren zweier
 Punkte aus $L_K(V)$. Weil die Identit\"at die einzige
 Permutation aus $G$ ist, die drei verschiedene Punkte festl\"asst, folgt, dass
 die Anzahl der 3-Sylowgruppen von $G$ gleich ${8 \choose 2} = 28$ ist. Somit
 enth\"alt $G$ genau $2 \cdot 28 = 56$ Elemente der Ordnung 3. Ferner enth\"alt
 $G$ genau $8 \cdot 6 = 48$ Elemente der Ordnung 7.
 \par
       Es seien nun $P$ und $Q$ zwei verschiedene Punkte von $L_K(V)$.
 Dann hat der Stabilisator $H$ der Menge $\{P,Q\}$ die Ordnung 6
 und die 3-Sylowgruppe $S$ von $H$ ist ein Normalteiler von $H$.
 Weil die Identit\"at die einzige Permutation in $G$ ist, die drei
 Fixpunkte hat, zerlegt $S$ die Menge der von $P$ und $Q$
 verschiedenen Punkte in zwei Bahnen der L\"ange $3$. Da
 Involutionen keine Fixpunkte haben, werden diese beiden Bahnen von
 den in $H$ liegenden Involutionen vertauscht, so dass $H$ neben
 der Bahn $\{P,Q\}$ noch eine weitere Bahn der L\"ange $6$ hat.
 \par
       W\"are $H$ zyklisch, so enthielte $H$ genau zwei Elemente der
 Ordnung 6. An Hand der Bahnen von $H$ sieht man, dass $G$ dann
 mindestens $2 \cdot {8 \choose 2} = 56$ Elemente der Ordnung 6
 enthielte. Da $G$, wie wir schon wissen, 56 Elemente der Ordnung
 3 und 48 Elemente der Ordnung 7 enth\"alt, enthielte $G$
 h\"ochstens $168 - 2 \cdot 56 - 48 = 8$ Elemente, deren Ordnung eine
 Potenz von 2 ist. Folglich enthielte $G$ nur eine
 2-Sylowgruppe im Widerspruch zur Einfachheit von $G$. Also ist
 $H$ nicht zyklisch und damit nicht abelsch.
 \par
       Es sei jetzt insbesondere $P = P_0$ und $Q = P_\infty$. Dann induziert
 die durch $u^\rho := 2u$ und $v^\rho := 4v$ definierte Abbildung
 $\rho$ ein Element der Ordnung 3 aus $H$. Die von $\rho$
 erzeugte Gruppe hat die beiden Bahnen $B_1 := \{P_1, P_2, P_4\}$
 und $B_2 := \{P_{-1}, P_{-2}, P_{-4}\}$. Es sei $\lambda$ eine
 Involution aus $H$. Wie wir gesehen haben, ist $B_1^\lambda =
 B_2$. Weil $H$ als nicht abelsche Gruppe drei Involutionen
 enth\"alt, d\"urfen wir annehmen, dass $P_1^\lambda = P_{-1}$ ist.
 Ist $i \neq 0$, so gibt es ein $g(i) \in K^*$ mit $P_i^\lambda =
 P_{g(i)}$. Weil $H$ nicht zyklisch ist, ist $\rho \lambda = \lambda \rho^{-1}$.
 Also ist
 $$ P_{g(2i)} = P^{\rho \lambda}_i = P^{\lambda \rho^{-1}}_i = P_{4g(i)}. $$
 Daher ist $g(2i) = 4g(i)$. Wegen $g(1) = -1$ ist daher $g(2) = -4$ und
 $g(4) = 4g(2) = -2$. Aus $\lambda ^2 = 1$ folgt somit
 $$ \lambda = (P_0P_\infty)(P_1P_6)(P_2P_3)(P_4P_5). $$
 definiert man die Abbildung $\kappa$ durch $u^\kappa := -v$ und $v^\kappa := u$,
 so sieht man, dass $\kappa \in L_K(V)$ die Permutation $\lambda$
 induziert. Wegen $\det(\kappa) = 1$ ist daher $\lambda \in \PSL(V)$.
 Weil $G$ das Erzeugnis von $N$ und $\lambda$ ist, folgt weiter
 $G \subseteq \PSL(V)$. Aus $|G| = 168 = |\PSL(V)|$ folgt schlie\ss lich
 die Behauptung.
 \medskip
       Um e) zu beweisen, nutzen wir aus, dass die Grup\-pe
 $\PSL(2,9)$ eine zur $A_5$ isomorphe Untergruppe enth\"alt. Dass
 dem so ist, folgt aus dem Satz, der besagt, dass $\PSL(2,q)$ genau
 dann eine zur $A_5$ isomorphe Untergruppe enth\"alt, wenn
 $|\PSL(2,q)|$ durch $5$ teilbar ist. Zum Beweise dieses Satzes
 m\"ussen wir etwas weiter ausholen und beweisen zun\"achst den folgenden
 Satz (Moore 1897).
 \medskip\noindent
 {\bf 4.4. Satz.} {\it Es sei $m$ eine nat\"urliche Zahl mit $m \geq 3$. Ist
 $F$ die freie Gruppe in den Erzeugenden $a_1, \dots,
 a_{m-2}$ und ist $N$ der von den Elementen $a_1^3$, $a_i^2$ f\"ur
 $i := 2, \dots, m - 2$, $(a_ia_{i+1})^3$ f\"ur $i := 1, \dots, m - 3$ und
 $(a_ia_j)^2$ f\"ur $1 \leq i \leq j - 2 \leq m - 4$ erzeugte
 Normalteiler von $F$, so ist $F/N$ zur alternierenden Gruppe $A_m$
 vom Grade $m$ isomorph.}\index{alternierende Gruppen}{}
 \smallskip
       Beweis. Dies ist gewiss richtig f\"ur $m = 3$. In diesem Falle ist $F$
 ja die von $a_1$ erzeugte unendliche zyklische Gruppe und $N$ ist ihre
 Untergruppe vom Index 3. Andererseits hat die alternierende Gruppe vom Grade
 3 die Ordnung 3.
 \par
       Es sei also $m > 3$ und der Satz gelte f\"ur
 $m - 1$. Dann ist die von den Elementen $a_1N$, \dots, $a_{m-3}N$
 erzeugte Untergruppe $U$ von $F/N$ nach Induktionsannahme ein epimorphes Bild
 der $A_{m-1}$. Es folgt $|U| \leq {1 \over 2} (m - 1)!$. Es sei $V$ die
 Untergruppe von $F$ mit $V/N = U$. Dann gilt $a_i \in V$ f\"ur
 $i := 1$, \dots, $m - 3$. Wir definieren die Rechtsrestklassen $R_i$ von
 $V$ rekursiv durch $R_{m-1} := V$ und $R_i := R_{i+1}a_i$ f\"ur
 $i := m - 2$, \dots, 1 sowie $R_0 := R_1a_1$ und zeigen, dass dies
 bereits alle Rechtsrestklassen nach $V$ sind.
 \par
       Zun\"achst beachten wir, dass $R_iy = R_i$ ist f\"ur alle $i$
 und alle $y \in N$, da $N$ ein in $V$ enthaltener Normalteiler
 ist. Es ist $R_{m-1} a_1 = R_{m-1}$. Es sei $i \geq 4$ und es
 gelte $R_ia_1 = R_i$. Dann ist $R_{i-1}a_1 = R_ia_{i-1}a_1$. Wegen
 $i - 1 \geq 3$ ist $(a_1a_{i-1})^2 \in N$. Es gibt daher $y$, $z \in N$
 mit $a_{i-1}a_1 = ya_1^{-1}a_{i-1} = za^2_1a_{i-1}$. Hieraus folgt
 $$ R_{i-1}a_1 = R_iza^2_1a_{i-1} = R_{i-1}. $$
 Also gilt $R_ia_1 = R_i$ f\"ur alle $i \geq 3$. Nun ist $R_2a_1 = R_1$,
 $R_1a_1 = R_0$ und wegen $a^3_1 \in N$ ist $R_0a_1 = R_2$. Damit ist
 gezeigt, dass Multiplikation von rechts mit $a_1$ die $R_i$
 untereinander permutiert. Wir beachten ferner, dass die
 Rechtsmultiplikation mit $a_1$ auf den Indizes die Permutation
 $(012)$ bewirkt.
 \par
       Es sei $j > 1$. Ist $i > j + 1$, so ist $a_j$ mit $a_i$ modulo $N$
 vertauschbar, da in diesem Falle ja $a_j^2$, $a_i^2$, $(a_ja_i)^2 \in N$ ist.
 Daher folgt mit einer simplen Induktion
 $$ R_ia_j = R_{i-1}a_ia_j = R_{i-1}a_ja_i = R_{i-1}a_i = R_i $$
 f\"ur $i := j + 2$, \dots, $m - 1$. Es ist
 $$ R_ja_j = R_{j+1} $$
 und
 $$ R_{j+1}a_j = R_j. $$
 \par
       Es sei $j = 2$. Es ist $a_1^3$, $a_2^3$, $(a_1a_2)^3 \in N$. Es folgt
 $$\eqalign{
      R_1a_2 &= R_2a_1a_2 = R_2a_2a_1^{-1}a_2a_1^{-1}   \cr
 	    &= R_3a_1^2 a_2a_1^2 = R_3a_2a_1^2         \cr
 	    &= R_2a_1^2=R_0 \cr} $$
 und
 $$\eqalign{
     R_0a_2 &= R_1a_1a_2=R_1a_2a_1^{-1} a_2a_1^{-1}     \cr
	    &= R_0a_1^{-1} a_2a_1^{-1} = R_1a_2a_1^{-1} 
	    = R_1 \cr} $$
 Multiplikation von rechts mit $a_2$ permutiert also die $R_i$ und
 bewirkt auf den Indizes die Permutation $(01)(23)$.
 \par
       Es sei $j > 2$. Dann ist $a^2_{j-1}$, $a^2_j$, $(a_{j-1}a_j)^3 \in N$.
 Hiermit und mit dem bereits Bewiesenen folgt
 $$\eqalign{
       R_{j-1}a_j &= R_ja_{j-1}a_j = R_ja_ja_{j-1}a_ja_{j-1}       \cr
		  &= R_{j+1}a_{j-1}a_ja_{j-1} = R_{j+1} a_ja_{j-1} \cr
		  &= R_ja_{j-1} = R_{j-1}. \cr}$$
 Mittels Induktion folgt hieraus
 $$ R_ia_j = R_i $$
 f\"ur $i := 2, \dots, j - 1$. Wegen $j>2$ ist $(a_1a_j)^2 \in N$. Hieraus
 folgt, wenn man noch bemerkt, dass $2 \leq j - 1$ ist,
 $$ R_1a_j = R_2a_1a_j = R_2a_ja_1^{-1} = R_2a_1^{-1} = R_0$$
 und
 $$ R_0a_j = R_1a_1a_j = R_1a_ja_1^{-1} = R_0a_1^{-1} = R_1. $$
 Multiplikation von rechts mit $a_j$ permutiert also die $R_i$ und
 bewirkt auf den Indizes die Permutation $(01)(j,j+1)$. Damit ist
 gezeigt, dass der Index von $V$ in $F$ h\"ochstens gleich $m$ ist.
 Also ist auch der Index von $U$ in $F/N$ h\"ochstens gleich $m$.
 Daher ist
 $$ |F/N|\leq {1 \over 2} (m - 1)!m = |A_m|.$$
 \par
       Setze $\alpha_1 := (012)$ und $\alpha_i := (01)(i,i+1)$ f\"ur
 $i:= 2$, \dots, $m - 2$. Dann zeigt eine simple Induktion, dass $A_m$
 von $\{\alpha_i \mid i := 1, \dots, m - 2\}$ erzeugt wird. Ferner gelten die
 Relationen $\alpha^3_1 = 1, \alpha^2_i = 1$ f\"ur $i := 2$, \dots, $m - 2$,
 $(\alpha_i\alpha_{i+1})^3 = 1$ f\"ur $i: = 1$, \dots, $m - 2$ und
 $(\alpha_i\alpha_j)^2 = 1$ f\"ur $1 \leq i \leq j - 2$. Hieraus
 folgt, dass es einen Epimorphismus $\pi$ von $F/N$ auf $A_m$ gibt
 mit $\pi(a_iN) = \alpha_i$ f\"ur alle $i$. Somit enth\"alt $F/N$
 mindestens $|A_m|$ Elemente, so dass $|F/N| = |A_m|$ ist. Dies
 impliziert schlie\ss lich, dass $\pi$ ein Isomorphismus ist. Damit
 ist der Satz bewiesen.\index{alternierende Gruppen}{}
 \medskip\noindent
 {\bf 4.5. Satz.} {\it Ist $G$ eine von $\{1\}$ verschiedene Gruppe
 und hat $G$ Erzeugende $a$ und $b$ mit $a^5 = b^2 = (ab)^3=1$, so ist
 $G$ zur $A_5$ isomorph.}
 \smallskip
       Beweis. Wir zeigen, dass $G$ ein epimorphes Bild der $A_5$ ist.
 Hieraus folgt wegen der Einfachheit der $A_5$ und der
 Voraussetzung, dass $G \neq \{1\}$ ist, die Behauptung \"uber $G$.
 \par
       Wir setzen $\alpha_1 := ba^{-1}$, $\alpha_2 := (ba)b(ba)^{-1}$ und
 $\alpha_3 := \alpha_2aba^3$. Dann ist
 $$ (\alpha_1\alpha_2^{-1}\alpha_3)^2 = a^{3 \cdot 2} = a \quad
 \mathrm{und}\quad
  \alpha_1 a = b, $$
 so dass $G$ auch von $\{\alpha_1, \alpha_2, \alpha_3\}$ erzeugt wird.
 \par 
       Es ist
 $$ \alpha_1^3 = (ba^{-1})^3 = ((ab)^{-1})^3 = 1. $$
 Ferner ist $\alpha^2_2 = 1$, da $\alpha_2$ zu $b$ konjugiert ist.
 \par
       Aus $(ab)^3 = 1$ folgt $bab = a^{-1}ba^{-1}$ und damit dann
 $$ \alpha_3 = baba^{-1}baba^3 = a^{-1}ba^{-3}ba^2.$$
 Daher ist
 $$\eqalign{
    \alpha^2_3 &= a^{-1}ba^{-3}ba^2a^{-1}ba^{-3}ba^2 = a^{-1}ba^2baba^2ba^2  \cr
               &=a^{-1}bababa^2 = a^{-1}a^{-1}b^2a^2 = 1. \cr}$$
 Es ist
 $$\eqalign{
   \alpha_1\alpha_2 &= a^{-1}baba^{-1}b = ba^{-2}ba^{-2}b=ba^{-3}(ab)a^{-2}b \cr
                    &= (a^{-2}b)^{-1}(ab)a^{-2}b \cr}$$
 und daher $(\alpha_1\alpha_2)^3 = 1$.
 \par 
       Es ist $ba^{-1}b = aba$ und somit
 $$ (\alpha_2\alpha_3)^2 = aba^4ba^3 = aba^{-1}ba^3=a^2ba^4. $$
 Hiermit folgt dann
 $$ (\alpha_2\alpha_3)^3 = a^2ba^4aba^3 = 1.$$
 Wegen $ba^4b = ba^{-1}b = aba$ erhalten wir
 $$ \alpha_1 \alpha_3 = ba^{-1}baba^{-1}baba^3 = aba^3ba^2ba^3 $$
 und weiter
 $$\eqalign{
       (\alpha_1\alpha_3)^2 &= aba^3ba^2ba^4ba^3ba^2ba^3
                             = aba^3ba^3ba^4ba^2ba^3             \cr
			    &= aba^3ba^4ba^3ba^3 = aba^4ba^4ba^3 
                            = a^2ba^5ba^3 = 1. \cr}$$
 \par
       Wir haben gezeigt, dass die Gruppe $G$ von $\{\alpha_1, \alpha_2,
 \alpha_3\}$ erzeugt wird und dass die Relationen
 $ \alpha_1^3 = \alpha_2^2 = \alpha^2_3 = 1$
 und
 $(\alpha_1\alpha_2)^3 = (\alpha_2\alpha_3)^3=(\alpha_1\alpha_3)^2 = 1 $
 gelten. Mit 4.4 folgt daher, dass $G$ ein epimorphes Bild der
 $A_5$ ist. Damit ist, wie eingangs bemerkt, der Satz bewiesen.
 \medskip
       Satz 4.5 sagt nichts dar\"uber, ob $A_5$ eine Pr\"asentation der
 in diesem Satz beschriebenen Art hat. Dies ist jedoch leicht
 nachzupr\"ufen und wird au\ss erdem aus dem gleich noch zu
 Beweisenden hervorgehen.
 \medskip\noindent
 {\bf 4.6. Satz.} {\it Es sei $K$ ein kommutativer K\"orper und
 $V 
 $ sei ein Vektorraum vom Rang $2$ \"uber $K$. Ist
 $\sigma \in \SL(V)$, so gilt:
 \item{a)} Ist $\Spur(\sigma)=1$, so ist $\sigma^3=-1$.
 \item{b)} Ist $\Spur(\sigma)=-1$, so ist $\sigma^3=1$.
 \item{c)} Ist $\Spur(\sigma)=0$, so ist $\sigma^2=-1$.}
 \smallskip

 Beweis. Nach dem Satz von Cayley--Hamilton ist $\sigma^2 - \Spur(\sigma)\sigma +1 =0$,
 da ja $\det (\sigma) = 1$. Also gilt c). Wegen $\sigma^3 = -\Spur(\sigma)\sigma^2 -\sigma$
 folgen auch die Behauptungen a) und b).

 \medskip\noindent
 {\bf 4.7. Satz.} {\it Es sei $K$ ein kommutativer K\"orper und
 $\alpha$ sei ein involutorischer Automorphismus von $K$. Ferner
 sei $V = uK \oplus vK$ ein Vektorraum vom Rang $2$ \"uber $K$ und
 $f$ sei die durch
 $$ f(ua + vb,uc + vd) := a^\alpha c + b^\alpha d $$
 definierte $\alpha$-Semibilinearform auf $V$. Bezeichnet
 $\SU(V,f)$ die Gruppe aller linearen Abbildungen $\sigma \in \SL(V,K)$
 mit $f(x^\sigma, y^\sigma) = f(x,y)$ f\"ur alle $x$, $y \in V$, so gilt:
 \item{a)} Die Gruppe $\SU(V,f)$ ist isomorph zur Gruppe aller Matrizen
 $$ \pmatrix{
	               a & b          \cr
	       -b^\alpha & a^\alpha   \cr} $$
 mit $a$, $b \in K$ und $a^{1+\alpha} + b^{1+\alpha} = 1$.
 \item{b)} Gibt es ein $x \in V - \{0\}$ mit $f(x,x) = 0$ und ist $W$
 ein Vektorraum vom Rang $2$ \"uber dem Fixk\"orper
 $$ F := \{k \mid k \in K, k^\alpha = k\} $$
 von $\alpha$, so ist $\SU(V,f)$ zu $\SL(W,F)$ isomorph.\par}
 \smallskip
       Beweis. a) Es sei $\sigma \in \SL(V,K)$. Es gibt dann Elemente
 $a$, $b$, $c$, $d \in K$ mit $u^\sigma = ua + vb$, $v^\sigma = uc + vd$ und
 $ad - bc = 1$. Nun liegt $\sigma$ genau dann in $\SU(V,f)$, wenn
 $$ f(u^\sigma,u^\sigma) = 1 = f(v^\sigma,v^\sigma) $$
 und
 $$ f(u^\sigma,v^\sigma) = 0 $$
 ist. Also liegt $\sigma$ genau dann in $\SU(V,f)$, wenn
 $$\eqalign{
       1 &= a^{1 + \alpha} + b^{1 + \alpha},         \cr
       1 &= c^{1 + \alpha} + d^{1 + \alpha},         \cr
       0 &= a^\alpha c + b^\alpha d,                 \cr
       1 &= ad - bc \cr} $$
 ist. Man verifiziert m\"uhelos, dass diese Gleichungen genau dann erf\"ullt
 sind, wenn $d = a^\alpha$, $c = -b^\alpha$ und
 $a^{1 + \alpha} + b^{1 + \alpha} = 1$ ist. Damit ist a) bewiesen.
 \par
       b) Aus $f(x,x)=0$ und $x \neq 0$ folgt, dass $x \not\in vK$ ist.
 Ist also $x = ur + vs$, so ist $r \neq 0$. Es folgt
 $$ f(x,u) = r^\alpha \neq 0. $$
 Wir d\"urfen daher annehmen, dass $f(x,u)=1$ ist.
 \par
       Es gilt nat\"urlich auch $x \not\in uK$. Daher ist $V = xK \oplus uK$. Es
 sei $y := u + xr$. Dann ist
 $$\eqalign{
        f(y,y) &= f(u + xr,u + xr)     \cr
	       &= f(u,u)+r^\alpha f(x,u) + rf(u,x) + r^{1 + \alpha}f(x,x) \cr
	       &= 1 + r + r^\alpha. \cr} $$
 Die durch $\varphi(r) := r + r^\alpha$ definierte Abbildung
 $\varphi$ ist eine lineare Abbildung des $F$-Vektorraumes $K$ in
 $F$. W\"are $\varphi = 0$, so w\"are $r^\alpha = -r$ f\"ur alle
 $r\in K$. Es folgte $1 = 1^\alpha = -1$ und damit $r^\alpha = r$, so
 dass $\alpha$ die Identit\"at w\"are. Involutionen sind aber
 {\it per definitionem\/} von  der Identit\"at verschieden. Also ist
 $\varphi$ nicht Null und damit surjektiv. Es gibt also ein $r \in K$ mit
 $1 + r + r^\alpha = 0$. Ist nun $y := u + xr$ mit eben diesem $r$, so
 ist $f(y,y)=0$. Ferner ist
 $$ f(x,y) = f(x,u) + f(x,x)r = 1. $$
 Es gibt ein $s\in K$ mit $s^\alpha \neq s$. Setze $k := s - s^\alpha$.
 Dann ist $k \neq 0$ und $k^\alpha = -k$. Setze $z := yk$. Dann ist
 $f(z,z) = 0$ und $f(x,z) = f(x,y)k  = k = -f(z,x)$. Schlie\ss lich gilt
 $V = xK \oplus zK$, da ja $y \not\in xK$ ist.
 \par 
       Es sei nun $\sigma \in \SU(V,f)$ und $x^\sigma = xa + zb$ sowie
 $z^\sigma = xc + zd$. Dann gelten die Gleichungen
 $$\eqalign{
        1 &= ad - bc,       \cr
	0 &= f(x,x) = f(x^\sigma, x^\sigma) = k(a^\alpha b - ab^\alpha), \cr
	0 &= f(z,z) = f(z^\sigma, z^\sigma) = k(c^\alpha d - cd^\alpha), \cr
	k &= f(x,z) = f(x^\sigma, z^\sigma) = k(a^\alpha d - b^\alpha c).\cr}$$
 Weil $k \neq 0$ ist, gelten daher auch die Gleichungen
 $$\eqalign{
        1 &= ad - bc,    \cr
	0 &= a^\alpha b - ab^\alpha, \cr
	0 &= c^\alpha d - cd^\alpha, \cr
	1 &= a^\alpha d - b^\alpha c.\cr}$$
 Aus den beiden letzten Gleichungen sowie der ersten folgt
 $$ d^\alpha = (a^\alpha d^\alpha - b^\alpha c^\alpha)d = (ad-bc)^\alpha d = d $$
 und
 $$ c^\alpha = (a^\alpha d^\alpha - b^\alpha c^\alpha)c
	     = (ad - bc)^\alpha c = c. $$
 Also liegen $c$ und $d$ in $F$. Aus der zweiten und vierten Gleichung folgt
 unter Benutzung der ersten
 $$a = (ad - bc)a^\alpha = a^\alpha $$
 und
 $$ b = (ad - bc)b^\alpha = b^\alpha. $$
 Also liegen auch $a$ und $b$ in $F$. Sind umgekehrt $a$, $b$, $c$, $d \in F$ und
 ist $ad - bc = 1$, so verifiziert man ohne M\"uhe, dass die durch
 $x^\sigma := xa + zb$, $z^\sigma := xc + zd$ definierte Abbildung $\sigma$ in
 $\SU(V,f)$ liegt. Damit ist gezeigt, dass $\SU(V,f)$ zur Gruppe
 aller Matrizen
 $$ \pmatrix{
              a & b \cr
	      c & d \cr} $$
 mit $a$, $b$, $c$, $d \in F$ und $ad - bc = 1$ isomorph ist. Da
 diese Gruppe wiederum mit $\SL(W,F)$ isomorph ist, ist alles bewiesen.
 \medskip\noindent
 {\bf 4.8. Korollar.} {\it Die Gruppe $\SL(2,q)$ ist isomorph zur
 Gruppe aller Matrizen
 $$ \pmatrix{
	         a & b   \cr
	      -b^q & a^q \cr} $$
 mit $a$, $b \in \GF(q^2)$ und $a^{1+q} + b^{1+q} = 1$.}
 \smallskip
       Beweis. Es sei $K := \GF(q^2)$ und $V = uK \oplus vK$ sei ein
 Vektorraum vom Rang $2$ \"uber $K$. Ferner sei $f$ die durch
 $$ f(ua + vb, uc + vd) := a^qc + b^qd $$
 definierte Semibilinearform auf $V$. Das Korollar 4.8 folgt aus 4.7, wenn wir
 nachweisen k\"onnen, dass es einen Vektor $x \in V - \{0\}$ gibt mit
 $f(x,x)=0$.
 \par
       Die multiplikative Gruppe von $K$ ist zyklisch und hat die Ordnung
 $q^2 - 1$. Daher ist die durch $r^\eta := r^{q+1}$ definierte
 Abbildung $\eta$ ein Homomorphismus der multiplikativen Gruppe von
 $K$ auf die multiplikative Gruppe von $\GF(q)$. Es gibt daher ein
 $r \in K$ mit $r^{1+q} = -1$. Setzt man $x := u + vr$, so ist $x \neq 0$ und
 $f(x,x) = 1 + r^{1+q} = 0$, so dass $x$ die gew\"unschten
 Eigenschaften hat.
 \medskip\noindent
 {\bf 4.9. Satz.} {\it Genau dann enth\"alt $\PSL(2,q)$ eine zur
 $A_5$ isomorphe Untergruppe, wenn die Ordnung von $\PSL(2,q)$ durch
 $5$ teilbar ist.}
 \smallskip
       Beweis. Die Notwendigkeit der Bedingung ist trivial. Es sei also 5
 ein Teiler von $|\PSL(2,q)|$. Wir setzen $K := \GF(q)$ und $V = uK \oplus vK$
 sei ein Vektorraum vom Range $2$ \"uber $K$. Wie wir wissen, ist
 $$ \big|\PSL(2,q)\big| = {1 \over \ggT(2,q - 1)} q(q^2 - 1). $$
 Drei F\"alle sind zu unterscheiden.
 \par
       1. Fall: Es ist $q = 5^r$. Betrachtet man die Menge aller Abbildungen
  $\sigma$ der Form $u^\sigma = ua + vb$ und $v^\sigma = uc + vd$ mit
 $a$, $b$, $c$, $d \in \GF(5)$ und $ad - bc = 1$, so ist diese eine zur
 $\SL(2,5)$ isomorphe Untergruppe von $\SL(2,q)$. Hieraus folgt, dass
 $\PSL(2,q)$ eine zur $\PSL(2,5)$ isomorphe Untergruppe enth\"alt. Da
 diese nach 4.1c) zur $A_5$ isomorph ist, gilt in diesem Falle die
 Behauptung.
 \par
       2. Fall: Es ist 5 ein Teiler von $q - 1$. In diesem Falle gibt es
 ein $a \in K$ mit $a^5 = 1 \neq a$. Es folgt $a - a^{-1} \neq 0$.
 Wir setzen $b := (a - a^{-1})^{-1}$ und w\"ahlen $c$, $d \in K$, so
 dass $b^2 + cd = -1$ ist. Wir definieren die Abbildungen $\sigma$ und
 $\tau$ von $V$ in sich durch $u^\sigma := ua$ und $v^\sigma := va^{-1}$ sowie
 $u^\tau := ub + vc$ und $v^\tau := ud - vb$. Dann ist
 zun\"achst $\det(\sigma) = 1 = \det(\tau)$. Ferner ist $\sigma^5 = 1$
 und $\Spur(\tau) = b - b = 0$, so dass nach 4.6b) also $\tau^2 = -1$ ist.
 Schlie\ss lich ist $u^{\sigma \tau}=uba + vca$ und $v^\sigma
 \tau = uda^{-1} - vba^{-1}$ und daher $\Spur(\sigma \tau) = b(a - a^{-1}) = 1$,
 so dass nach 4.6a) die Gleichung $(\sigma \tau)^{-3} = -1$ gilt.
 Sind $\bar{\sigma}$ und $\bar{\tau}$ die von $\sigma$ und $\tau$
 in $\PSL(2,q)$ induzierten Abbildungen, so ist also $\bar{\sigma}
 \neq 1$ und $\bar{\sigma}^5 = \bar{\tau}^2 = (\bar{\sigma}\bar{\tau})^3 = 1$,
 so dass nach 4.5 die von $\bar{\sigma}$ und
 $\bar{\tau}$ erzeugte Untergruppe von $\PSL(2,q)$ zur $A_5$ isomorph ist.
 \par
       3. Fall: Es ist 5 ein Teiler von $q + 1$. Es gibt ein $a \in \GF(q^2)$
 mit $a^5 = 1 \neq a$. Weil 5 ein Teiler von $q + 1$ ist,
 folgt $a^{1+q} = 1$. Da $q - 1$ nicht durch $5$ teilbar ist, ist $a^q \neq a$.
 Wir setzen $b := (a - a^q)^{-1}$. Dann ist $b^q = -b$ und
 folglich $1-b^{1+q} \in K$, so dass es ein $c \in \GF(q^2)$ gibt
 mit $c^{1+q} = 1 - b^{1+q}$. Wir betrachten die Matrizen
 $$ A := \pmatrix{ a & 0   \cr
   		   0 & a^q \cr}\ \ \ {\rm und }\ \ \ 
    B := \pmatrix{ b & c   \cr
        	-c^q & b^q \cr}. $$
 Dann ist $\det(A) = a^{1+q} = 1 = b^{1+q} + c^{1+q} = \det(B)$. Ferner ist
 $a^5 = 1$ und $\Spur(B) = b + b^q = b - b = 0$, so dass $B^2 = 1$ ist.
 Schlie\ss lich ist $\Spur(AB) = b(a - a^q) = 1$, woraus $(AB)^3 = -1$ folgt. Aus
 4.8 und 4.5 folgt daher, dass $\PSL(2,q)$ auch in diesem Falle eine zur
 $A_5$ isomorphe Untergruppe enth\"alt. Damit ist 4.9 bewiesen.
 \medskip
       Um 4.1 e) zu beweisen, beachten wir zun\"achst, dass
 $|\PSL(2,9)| = 360 = |A_6|$ ist. Nach 4.9 enth\"alt $\PSL(2,9)$ eine zur
 $A_5$ isomorphe Untergruppe, so dass sie, da sie einfach ist, eine
 treue Darstellung als Permutationsgruppe vom Grad $6$ besitzt.
 Wiederum wegen der Einfachheit kann sie nur gerade Permutationen
 enthalten, so dass sie zu einer Untergruppe der $A_6$ isomorph
 ist. Aus der Gleich\-heit der Ordnungen folgt daher die Isomorphie
 der Gruppen.\index{alternierende Gruppen}{}
 \par
       Es bleibt zu zeigen, dass die Gruppen $\PSL(4,2)$ und $A_8$
 isomorph sind. Es gibt mehrere M\"oglichkeiten, dies zu beweisen.
 Der sch\"onste Beweis scheint mir der von Moore zu sein,
 der hier wie\-der\-ge\-ge\-ben sei (Moore 1899).
 \par
       Wir betrachten das folgende Schema:
 $$\matrix{
              0 & 1 & 2 & 3 & 4 & 5 & 6 \cr
	      1 & 2 & 3 & 4 & 5 & 6 & 0 \cr
              3 & 4 & 5 & 6 & 0 & 1 & 2 \cr} $$
 Nennt man die Ziffern 0 bis 6
 Punkte und die Ziffernmengen, die aus den Ziffern einer Spalte
 bestehen, Geraden und definiert man die Inzidenz als die
 $\in$-Relation, so ist die so 
 entstehende Inzidenzstruktur
 die projektive Ebene der Ordnung 2. Aus dem zweiten Struktursatz
 und Satz 1.13 folgt, dass $\PSL(3,2)$ die volle
 Kollineationsgruppe dieser Ebene ist. Weil diese Gruppe einfach
 ist, ist sie eine Untergruppe der $A_7$. Mit 1.13 folgt, dass ihr
 Index in der $A_7$ gleich $15$ ist. Hieraus folgt, dass die
 Wirkung der $A_7$ auf der Menge $\{0, \dots, 6\}$ insgesamt 15
 verschiedene Realisierungen der projektiven Ebene der Ordnung 2
 liefert. Diese nennen wir Punkte einer Geometrie $\Sigma$, deren
 Geraden die 3-Teilmengen von $\{0, \dots, 6\}$ sind.
 (Ziffern hei\ss en im Folgenden wieder Ziffern.)
 Wir definieren Inzidenz dadurch, dass f\"ur einen Punkt $P$ und eine
 Gerade $G$ von $\Sigma$ genau dann $P \I G$ gelte, wenn $G$ eine
 Gerade der projektiven Ebene $P$ ist. Wir werden zeigen, dass
 $\Sigma$ aus den Punkten und Geraden des projektiven Raumes des
 Ranges 4 \"uber $\GF(2)$ besteht.
 \par
       1) $\Sigma$ enth\"alt $v = 15$ Punkte und $b = 35$ Geraden und jeder
 Punkt inzidiert mit genau $r = 7$ Geraden.
 \par
       Dies folgt unmittelbar aus der Definition von $\Sigma$.
 \par
       2) $A_7$ operiert als Automorphismengruppe auf $\Sigma$ und ist
 sowohl auf der Menge der Punkte als auch auf der Menge der Geraden
 von $\Sigma$ transitiv.
 \par
       Auch dies ist eine unmittelbare Konsequenz aus der Definition von
 $\Sigma$.
 \par
       3) Jede Gerade von $\Sigma$ tr\"agt genau drei Punkte.
 \par
       Aus 2) folgt, dass jede Gerade mit gleich vielen, etwa $k$ Punkten
 inzidiert. Das beim Beweise von I.7.5 benutzte Verfahren der zweifachen
 Abz\"ahlung\index{zweifache Abz\"ahlung}{} liefert im vorliegenden Falle die
 Gleichung $vr = bk$, dh. die Gleichung $15 \cdot 7 = 35 \cdot k$, so dass in
 der Tat $k = 3$ ist.
 \par
       4) Zwei verschiedene Geraden $G$ und $H$ haben genau dann einen
 Punkt gemeinsam, wenn sie als Tripel genau eine Ziffer gemeinsam
 haben.
 \par
       Liegt der Punkt $Q$ auf den beiden Geraden $G$ und $H$, so sind
 die Tripel $G$ und $H$ Geraden der projektiven Ebene $Q$, so dass
 sie genau eine Ziffer gemeinsam haben. Haben sie umgekehrt genau
 eine Ziffer gemeinsam, so ist $|G \cup H| = 5$. Weil die $A_7$ auf
 $\{0, \dots, 6\}$ f\"unffach transitiv ist, k\"onnen wir annehmen,
 dass $G = \{0,1,3\}$ und $H = \{1,2,4\}$ ist. Dann sind $G$ und $H$
 aber Geraden der Ebene, von der wir ausgegangen sind, so dass sie
 als Geraden von $\Sigma$ einen Punkt gemeinsam haben.
 \par
       5) Zwei verschiedene Geraden haben h\"ochstens einen Punkt ge\-mein\-sam.
 \par
       Es seien $G$ und $H$ zwei Geraden, die einen Punkt $Q$ gemeinsam
 haben. Nach 4) ist dann $|G \cup H|=1$. Weil die
 Kollineationsgruppe $K$ der projektiven Ebene $Q$ auf den Rahmen
 dieser Ebene transitiv ist, ist sie erst recht auf der Menge der
 2-Teilmengen der Geradenmenge transitiv. Da es sieben Geraden
 gibt, hat diese Menge die L\"ange 21. Daher gilt f\"ur den
 Stabilisator $K_{\{G,H\}}$ der Menge $\{G,H\}$, dass
 $$ |K_{\{ G,H\}}| = |K|/21 = 168/21 = 8 $$
 ist. Da die Menge aller
 2-Teilmengen von Tripeln mit genau einer gemeinsamen Ziffer gleich
 $$ {1 \over 2} \cdot 7 \cdot {6 \choose 2} \cdot {4 \choose 2}
 = 7 \cdot 3^2 \cdot 5 $$
 ist, folgt
 $ \big|(A_7)_{\{G,H\}}\big| = 8$,
 so dass $K_{\{G,H\}} = (A_7)_{\{G,H\}}$ ist. Hieraus folgt, dass es au\ss er
 $Q$ keinen weiteren Punkt mehr gibt, der auf $G$ und auf $H$ liegt.
 \par
       6) Durch zwei verschiedene Punkte geht genau eine Gerade.
 \par
       Sind $P$und $Q$ zwei verschiedene Punkte, so bezeichnen wir mit
 $r_{\{P,Q\}}$ die Anzahl der mit $P$ und $Q$ gleichzeitig
 inzidierenden Ge\-ra\-den. Mit 5) folgt, dass $r_{\{P,Q\}} \leq 1$
 ist. Da die Anzahl der 2-Teilmengen von Punkten, die auf einer
 Geraden liegen, gleich 3 ist, liefert zwei\-fa\-che Abz\"ahlung
 $$ 35 \cdot 3 = \sum_{\{P,Q\}} r_{\{P,Q\}} \leq {15 \choose 2} = 15\cdot 7. $$
 Hieraus folgt die Behauptung.
 \par
       Es bleibt, die G\"ultigkeit des Veblen-Young-Axioms nachzuweisen.
 Dazu seien $P$, $Q$ und $R$ drei nicht kollineare Punkte von $\Sigma$. Ferner
 seien $D$ und $E$ zwei verschiedene Punkte mit $D \I P + Q$ und $E \I Q + R$.
 Wir haben zu zeigen, dass die Geraden $D + E$ und $R + P$ einen Punkt gemeinsam
 haben. Wir d\"urfen wieder
 annehmen, dass $P + Q=\{0,1,3\}$ und $Q + R=\{1,2,4\}$ ist. Die Gerade
 $R + P$ hat mit den beiden Geraden $P + Q$ und $Q + R$ jeweils genau
 eine Ziffer gemeinsam. W\"are $1\in R + P$, so w\"are daher
 $R + P = \{1,5,6\}$. Da dies ein Tripel der Ebene $Q$ ist, w\"aren
 $P$, $Q$ und $R$ kollinear. Dieser Widerspruch zeigt, dass $1 \not\in R + P$
 ist. Indem wir nun gegebenenfalls eine der Permutationen
 $(03)(24), (03)(56)$ oder $(24)(56)$ anwenden, k\"onnen wir
 erreichen, dass $0$, $2 \in R + P$ ist. Also ist $R + P = \{0,2,5\}$ oder
 $R + P = \{0,2,6\}$. W\"are $R + P = \{0,2,6\}$, so w\"aren $P$, $Q$ und $R$
 aber kollinear. Also ist $R + P = \{0,2,5\}$.
 \par
       Die Ziffernmenge $D + E$ hat mit $P + Q = \{0,1,3\}$ und
 $Q + R = \{1,2,4\}$ jeweils genau eine Ziffer gemeinsam. Daher ist $D + E$
 gleich einem der neun Tripel $\{1,5,6\}$, $\{0,2,x\}$, $\{0,4,x\}$, $\{3,2,x\}$,
 $\{3,4,x\}$, wobei $x = 5$ oder 6 ist. H\"atten nun $D + E$ und $R + P$
 keinen Punkt gemeinsam, so folgte $D + E = \{0,2,6\}$, $\{0,4,5\}$,
 $\{3,2,5\}$, $\{3,4,6\}$. In jedem Falle folgte $Q \I D + E$. Dann
 w\"are aber entweder $Q = D$ und dann $D + E  = Q + R$ oder $Q = E$ und
 $D + E = P+Q$, was beides nicht der Fall ist. Dieser Widerspruch
 zeigt, dass $D + E$ und $R + Q$ doch einen Punkt gemeinsam haben.
 Damit ist gezeigt, dass $\Sigma$ eine projektive Geometrie ist.
 \par
       Da auf jeder Geraden von $\Sigma$ genau drei Punkte liegen, ist
 die Ordnung von $\Sigma$ gleich 2. Ist $r$ der Rang  der
 Geometrie, so ist also $15 = 2^r - 1$, so dass $r = 4$ ist. Auf Grund
 des ersten Struktursatzes ist $\Sigma$ daher die projektive
 Geometrie vom Range 4 \"uber $\GF(2)$. Von dieser wissen wir nun,
 dass sie eine Kollineationsgruppe besitzt, die zur $A_7$ isomorph
 ist. Aus Satz 1.13 folgt, dass
 $$ |\PGammaL(4,2)| = |\PGL(4,2)| = |\PSL(4,2)| = |A_8| $$
 ist. Ferner enth\"alt $\PSL(4,2)$, wie gerade gesehen, eine zur $A_7$
 isomorphe Untergruppe. Diese hat nat\"urlich den Index 8 in $\PSL(4,2)$,
 so dass $\PSL(4,2)$ auf den Rechtsrestklassen dieser Untergruppe als
 Permutationsgruppe vom Grade 8 operiert. Da die $\PSL(4,2)$ einfach ist, folgt,
 dass sie in der alternierenden Gruppe vom Grade 8 enthalten ist. Weil beide
 Gruppen die gleiche Ordnung haben, sind sie also gleich.
 \medskip
       Damit ist Satz 4.1 vollst\"andig bewiesen. Zum Schluss haben wir
 mehr bewiesen als in diesem Satze formuliert. Es ist also noch ein
 Korollar zu formulieren.
 \medskip\noindent
 {\bf 4.10. Korollar.} {\it Ist $V$ der Vektorraum vom Range $4$
 \"uber $\GF(2)$, so besitzt $\PSL(V)$ eine zur $A_7$ isomorphe
 Kollineationsgruppe, die auf den Punkten von $L_{\GF(2)}(V)$
 zweifach transitiv operiert.}
 \smallskip
       Beweis. Die $A_7$ operiert in obiger Darstellung auf der Menge der
 Geraden transitiv. Wir setzen $\gamma := (456)$ und bezeichnen mit
 $P$ die Ebene, mit der wir starteten. Dann ist $P \neq P^\gamma
 \neq P^{\gamma^2} \neq P$. Folglich permutiert die von $\gamma$
 erzeugte Gruppe die mit der Geraden $\{0,1,2\}$ inzidierenden
 Punkte transitiv. Weil zwei verschiedene Punkte mit genau einer
 Geraden inzidieren, ist $A_7$ daher auf der Menge der 2-Teilmengen der
 Punktmenge transitiv. Weil die Ordnung der $A_7$ ge\-ra\-de ist, folgt
 schlie\ss lich, dass $A_7$ auf der\index{alternierende Gruppen}{}
 Punktmenge der betrachteten Geometrie zweifach transitiv operiert.
 \medskip
       Bei diesen Untersuchungen ist als Nebenergebnis herausgekommen, dass die
 Gruppen $A_5$, $A_6$ und $A_8$ einfach sind. Generell gilt, dass $A_n$ einfach
 ist, wenn nur $n \neq 4$ ist. Dies beweist Camille Jordan in seinem {\it
 Trait\'e\/} (Jordan 1870/1989, S. 66). Mit diesem Buche beginnen auch die
 systematischen Untersuchungen der in diesem Kapitel betrachteten und anderer in
 der Geometrie vorkommenden Gruppen, wobei Jordan als
 K\"orper nur die Galoisfelder $\GF(p)$ betrachtet, wo $p$ eine Primzahl ist.
 Jordan kannte den Begriff der Faktorgruppe noch nicht, der erst durch Otto
 H\"older eingef\"uhrt wurde (H\"older 1889). Bei ihm, Jordan, findet sich
 die Einfachheit der $\PSL(n,p)$ unter dem Satz verborgen, dass ein
 Normalteiler der Gruppe $\GL(n,p)$ entweder im Zentrum von $\GL(n,p)$ enthalten
 ist oder aber die $\SL(n,p)$ enth\"alt.

\mysection{5. Quasiperspektivit\"aten}

\noindent
 In diesem Abschnitt werden wir Kollineationen untersuchen, die
 eine Verallgemeinerung der Perspektivit\"aten darstellen und wie
 diese ge\-wis\-se Unterraumkonfigurationen punktweise festlassen. Da
 sie sich also qua\-si wie Perspektivit\"aten verhalten, werden wir
 sie Quasiperspektivit\"aten nennen. Wir beginnen mit der
 Formulierung eines Sat\-zes, der eine banale, jedoch sehr
 n\"utzliche Folgerung aus Satz 1.2 ist.
 \medskip\noindent
 {\bf 5.1. Satz.} {\it Es sei $V$ ein Vektorraum \"uber dem K\"orper $K$ 
 mit $\Rg_K(V) \geq 2$. Ist $\gamma \in \PGammaL(V)$ und l\"asst $\gamma$ eine
 Gerade von $L_K(V)$ punktweise fest, so ist $\gamma$
 projektiv.}\index{projektive Kollineation}{}
 \smallskip
       Beweis. Es sei $G$ eine Gerade von $L_K(V)$, die von $\gamma$
 punktweise festgelassen wird. Ferner sei $\delta \in \GammaL(V)$
 und $\gamma$ werde von $\delta$ induziert. Nach 1.2 gibt es dann
 ein $k \in K^*$ mit $g^\delta = gk$ f\"ur alle $g \in G$. F\"ur den
 begleitenden Automorphismus $\alpha$ von $\delta$ gilt dann
 $x^\alpha = k^{-1}xk$, so dass $\alpha$ ein innerer
 Automorphismus ist. Nach 1.3 ist $\gamma$ daher projektiv.
 \medskip
       Es sei $V$ ein Vektorraum \"uber dem K\"orper $K$ und $\tau$ sei
 eine Kollineation von $L_K(V)$, die einen Unterraum $U$
 festl\"asst. Dann induziert $\tau$ eine Kollineation $\tau^*$ in
 $U$ und eine Kollineation $\tau^{**}$ in $V/U$. Wir nennen $\tau$
 {\it Quasiperspektivit\"at\/},\index{Quasiperspektivit\"at}{} falls
 $\tau^* = 1$ und $\tau^{**} = 1$
 ist. Ist $\tau$ eine Quasiperspektivit\"at und ist $U = \{0\}$ oder
 $U = V$, so ist $\tau = 1$. Ist $\Rg_K(U) = 1$, so ist $\tau$ eine
 Perspektivit\"at mit dem Zentrum $U$, und ist $\Rg_K(V/U) = 1$, so
 ist $\tau$ eine Perspektivit\"at mit der Achse $U$. Es sei daher
 $\Rg_K(U) \geq 2$ und $\Rg_K(V/U) \geq 2$. Nach 5.1 ist $\tau$
 projektiv. Es gibt daher eine lineare Abbildung $\sigma$, die
 $\tau$ induziert. Wegen $\Rg_K(U) \geq 2$ folgt aus 1.2 die
 Existenz eines von Null verschiedenen Elementes $m$ in $Z(K)$ mit
 $u^\sigma = um$ f\"ur alle $u \in U$. Wir k\"onnen ohne
 Einschr\"ankung der Allgemeinheit annehmen, dass $m = 1$ ist.
 Wegen $\Rg_K(V/U) \geq 2$ folgt fernerhin die Existenz eines von
 Null verschiedenen Elementes $k \in Z(K)$ mit $v^\sigma + U = vk + U$
 f\"ur alle $v \in V$. Definiert man $\alpha$ durch $v^\alpha := v^\sigma - vk$
 f\"ur alle $v \in V$, so ist $\alpha$ eine lineare
 Abbildung von $V$ in $U$. Es sind nun zwei F\"alle zu
 unterscheiden, n\"amlich  die F\"alle $k \neq 1$ und $k = 1$.
 \par
       Ist $k \neq 1$, so nennen wir $\tau$ eine {\it Quasistreckung\/}.
 In\index{Quasistreckung}{} diesem Falle setzen wir
 $$ W := \bigl\{v - v^\alpha (1 - k)^{-1}\mid v\in V\bigr\}. $$
 Dann ist $W$ ein Unterraum von $V$, da $k$
 ja ein Zentrumselement ist. Aus $v^\alpha \in U$ und
 $$ v = v^\alpha (1-k)^{-1} + v-v^\alpha (1-k)^{-1} $$
 folgt $V = U + W$. Ist $u \in U$, so ist $u^\alpha = u^\sigma - uk = u(1 - k)$.
 Wegen $v^\alpha \in U$ ist daher
 $$ v^{\alpha^2} = v^\alpha (1 - k). $$
 Hiermit folgt
 $$ \bigl(v - v^\alpha(1 - k)^{-1}\bigr)^\alpha
			= v^\alpha - v^\alpha (1 - k)(1 - k)^{-1} = 0. $$
 Dies besagt, dass $W$ im
 Kern von $\alpha$ liegt. Ist nun $x \in U \cap W$, so ist also $0 = x^\alpha =
 x(1 - k)$, so dass $x = 0$. Folglich ist $V = U \oplus W$. Ferner ist
 $$\eqalign{
      \bigl(v - v^\alpha(1 - k)^{-1}\bigr)^\sigma
	   &= v^\sigma - v^{\alpha\sigma} (1-k)^{-1}               \cr
	   &= v^\sigma -v^\alpha (1-k)^{-1}                        \cr
	   &= v^\alpha + vk - v^\alpha (1 - k)^{-1}                \cr
           &= \bigl(v^\alpha(1-k)+vk(1-k)-v^\alpha\bigr)(1-k)^{-1} \cr
	   &= \bigl(vk(1 - k)-v^\alpha k\bigr)(1 - k)^{-1}         \cr
	   &= \bigl(v-v^\alpha(1-k)^{-1}\bigr)k. \cr} $$
 Daher ist $w^\sigma = wk$ f\"ur alle $w \in W$. Die Kollineation
 $\tau$ l\"asst also in diesem Falle die beiden windschiefen
 Unterr\"aume $U$ und $W$ punkt\-wei\-se fest, was den Namen
 Quasistreckung rechtfertigt.
 \par
       Im Falle $k = 1$ gilt
 $$ v^{\alpha^2} = v^{\alpha\sigma} - v^\alpha k = v^\alpha-v^\alpha = 0, $$
 da ja $v^\alpha \in U$ ist. Folglich ist $\alpha^2 = 0$. Ferner ist
 $u^\alpha = u^\sigma - u = u - u = 0$ f\"ur alle $u \in U$. Also ist
 $U^\alpha \subseteq U \subseteq \Kern(\alpha)$. Es sei $v^\sigma = va \neq 0$.
 Dann ist $va = v^\sigma = v + v^\alpha$. Hieraus folgt
 $v(a - 1) = v^\alpha \in U$.
 W\"are $a \neq 1$, so w\"are $v \in U$ und daher $v = v^\sigma = va$
 und also doch $a = 1$. Also ist $v = v^\sigma = v + v^\alpha$ und
 damit $v^\alpha = 0$. Somit ist $v \in \Kern(\alpha)$. Dies besagt
 wiederum, dass $\tau$ au\ss erhalb $\Kern(\alpha)$ keinen
 Fixpunkt hat. Deshalb nennen wir $\tau$ eine
 {\it Quasielation\/}.\index{Quasielation}{}
 Ferner nennen wir $\Kern(\alpha)$ {\it Achse\/}\index{Achse}{} und $V^\alpha$
 {\it Zentrum\/}\index{Zentrum}{} von $\tau$.
 \medskip\noindent
 {\bf 5.2. Satz.} {\it Es sei $V$ ein Vektorraum \"uber dem
 K\"orper $K$ mit $\Rg_K(V) \geq 3$. Ferner seien $U$ und
 $W$ Unterr\"aume von $V$ mit $V = U \oplus W$. Schlie\ss lich
 bezeichne $\Lambda (U,W)$ die Gruppe aller Quasistreckungen, die
 $U$ und $W$ punktweise festlassen. 
 \item{a)} Ist $\Rg_K(U) = 1$ oder $\Rg_K(W) = 1$, so ist $\Lambda(U,W)$ zu
 $K^*$ isomorph.
 \item{b)} Ist $\Rg_K(U) \geq 2$ und $\Rg_K(W) \geq 2$, so ist
 $\Lambda(U,W)$ zu $Z(K^*)$ isomorph.
 \\
       In jedem Falle ist $\Lambda (U,W) \subseteq \PGaL(V)$.}
 \smallskip
       Beweis. a) ist nichts Neues. Es ist hier nur noch einmal notiert,
 um den Kontrast zu b) deutlich werden zu lassen.
 \par
       b) Ist $k \in Z(K^*)$ und definiert man $\lambda(k)$ durch
 $$ (u + w)^{\lambda(k)} := u + wk, $$
 so induziert $\lambda(k)$ eine Kollineation aus $\Lambda(U,W)$. Wie wir schon
 gesehen haben, erh\"alt man auf diese Weise alle Abbildungen aus $\Lambda(U,W)$.
 Mit 1.2 folgt, dass $\lambda(k)$ genau dann die
 Identit\"at induziert, wenn $k = 1$ ist. Da $\lambda$
 offensichtlich ein Homomorphismus ist, ist bereits alles gezeigt.
 \medskip
       Es sei wieder $V = U \oplus W$. Die Geometrie $L_K(V)$ werde
 $(U,W)$-{\it transitiv\/}\index{transitiv}{} genannt, falls es zu zwei
 verschiedenen
 Punkten $P$ und $Q$, die weder auf $U$ noch auf $W$ liegen und
 f\"ur die $U \cap (P + Q) \neq \{0\}$ sowie $W \cap (P + Q) \neq \{0\}$ gilt,
 stets ein $\lambda \in \Lambda (U,W)$ gibt mit $P^\lambda = Q$.
 \medskip\noindent
 {\bf 5.3. Satz.} {\it Ist $V$ ein Vektorraum \"uber dem K\"orper
 $K$ mit $\Rg_K(V) \geq 4$, so sind die folgenden Aussagen
 \"aquivalent:
 \item{a)} $L_K(V)$ ist pappossch.\index{papposscher Raum}{}
 \item{b)} Sind $U, W \in L_K(V)$, ist $\Rg_K(U) \geq 2$ und
 $\Rg_K(W) \geq 2$ und ist $V = U \oplus W$, so ist $L_K(V)$ ein
 $(U,W)$-transitiver Raum.
 \item{c)} Es gibt zwei Unterr\"aume $U$, $W \in L_K(V)$ mit $\Rg_K(U) \geq 2$,
 $\Rg_K(W) \geq 2$ und $V = U\oplus W$, so dass $L_K(V)$ ein
 $(U,W)$-transitiver Raum ist.\par}
 \smallskip
       Beweis. a) impliziert b): Gilt a), so ist $K$ kommutativ. Es sei
 $V = U \oplus W$. Ferner seien $P = pK$ und $Q = qK$ Punkte, die weder
 in $U$ noch in $V$ liegen, f\"ur die jedoch $(P + Q)\cap U \neq \{0\}$ und
 $(P + Q)\cap V \neq \{0\}$ gilt. Es gibt dann von $0$
 verschiedene Vektoren $u$, $u' \in U$ und $w$, $w' \in W$ mit $p = u + w$
 und $q = u' + w'$. Weil $(P + Q)\cap U$ und $(P + Q) \cap W$ Punkte sind,
 gibt es von 0 verschiedene Elemente $a$, $b \in K$ mit $u' = ua$ und
 $w' = wb$. Setzt man $c: = ba^{-1}$ und definiert $\sigma$ durch
 $(x + y)^\sigma := x + yc$ f\"ur alle $x \in U$ und alle $y \in W$, so
 induziert $\sigma$ wegen der Kommutativit\"at von $K$ eine
 Kollineation, die zu $\Lambda (U,W)$ geh\"ort. Ferner ist
 $$ p^\sigma = u + vc = (ua + vb)a^{-1} = qa^{-1} $$ 
 und daher $P^\sigma = Q$. Damit ist b) aus a) hergeleitet.
 \par 
       b) impliziert c): Nat\"urlich, denn wir haben ja $\Rg_K(V) \geq 4$
 vorausgesetzt.
 \par
       c) impliziert a): Es sei $0 \neq u \in U$, $0 \neq w \in W$ und
 $0 \neq c \in K$. Dann sind $P := (u + w)K$ und $Q := (u + wc)K$ zwei
 verschiedene Punkte, die weder auf $U$ noch auf $W$ liegen und
 f\"ur die $(P + Q)\cap U = uK$ und $(P + Q)\cap W = wK$ gilt. Nach
 Voraussetzung gibt es ein $a \in Z(K^*)$, so dass die durch
 $(x + y)^\sigma := x + ya$ f\"ur alle $x \in U$ und alle $y \in W$
 definierte Abbildung $\sigma$ den Punkt $P$ auf den Punkt $Q$
 abbildet. Es gibt daher ein $b \in K$ mit
 $$ u + wa = p^\sigma = qb = (u + wc)b. $$
 Weil $u$ und $w$ linear unabh\"angig
 sind, folgt $b = 1$ und $a = c$, so dass $c \in Z(K)$ gilt. Daher ist
 $K$ kommutativ, so dass $L_K(V)$ pappossch ist.
 \medskip\noindent
 {\bf 5.4. Satz.} {\it Es sei $V$ ein Vektorraum \"uber dem
 K\"orper $K$ und es gelte $\Rg_K(V) \geq 3$. Ferner seien $U$ und
 $W$ von $\{0\}$ verschiedene Unterr\"aume von $V$ mit $V = U \oplus
 W$. Schlie\ss lich sei $\Lambda_1$ die Gruppe aller Kollineationen
 aus $\PGaL(V)$, die $U$ als Ganzes und $W$ punktweise festlassen,
 und $\Lambda_2$ die Gruppe aller Kollineationen aus $\PGaL(V)$,
 die $U$ punktweise und $W$ als Ganzes festlassen. Dann gilt:
 \item{a)} Die Gruppen $\Lambda_1$ und $\Lambda_2$ normalisieren
 sich gegenseitig.
 \item{b)} Es ist $\Lambda_1 \cap \Lambda_2 = \Lambda(U,W)$.
 \item{c)} \,$\Lambda_1$ induziert in $L_K(U)$ die Gruppe $\PGaL(U)$
 und $\Lambda_2$ induziert in $L_K(W)$ die Gruppe $\PGaL(W)$.
 \item{d)} \,$\Lambda_1 \Lambda_2$ ist die Untergruppe aller
 Kollineationen aus $\PGaL(V)$, die $U$ und $W$ festlassen.
 \item{e)} Es ist $(\Lambda_1\Lambda_2)/\Lambda(U,W) \cong
 \PGaL(U) \times \PGaL(W)$.
 \item{f)} Es ist $Z(\Lambda_1\Lambda_2) = Z(\Lambda(U,W)).$}
 \smallskip
       Beweis. a) und b) sind banal.
 \par
       c) Wir zeigen, dass $\Lambda_1$ in $L_K(U)$ die Gruppe
 $\PGaL(U)$ induziert. Die zweite Aussage beweist sich analog.
 \par
       Weil $\lambda_1$ den Raum $U$ festl\"asst, induziert $\Lambda_1$
 eine Kollineationsgruppe $L$ in $L_K(U)$, die in $\PGaL(U)$
 enthalten ist, da $\Lambda_1$ nur aus projektiven Kollineationen
 besteht. Es sei $\sigma \in \PGaL(U)$. Es gibt dann ein $\lambda \in \GaL(U)$,
 welches $\sigma$ induziert. Wir definieren die
 Abbildung $\mu \in \GaL(V)$ durch $(u + w)^\mu = u^\lambda + w$ f\"ur
 alle $u \in U$ und alle $w \in W$. Dann induziert $\mu$ eine
 Kollineation aus $\Lambda_1$, die auf $L_K(U)$ mit $\sigma$
 \"ubereinstimmt. Also ist $L = \PGaL(U)$.
 \par
       d) Es sei $\gamma$ eine projektive Kollineation mit $U^\gamma = U$
 und $W^\gamma = W$. Nach c) gibt es ein $\lambda \in \Lambda_1$
 mit $P^\gamma = P^\lambda$ f\"ur alle Punkte $P$ von $LK(U)$.
 Es folgt $\gamma \lambda^{-1} \in \Lambda_2$ und damit $\gamma \in
 \Lambda_2\Lambda_1 = \Lambda_1\Lambda_2$.
 \par 
       e) Aus a) und b) folgt
 $$ (\Lambda_1\Lambda_2)/\Lambda(U,W)
      = \bigl(\Lambda_1/\Lambda(U,W)\bigr)
		  \times \bigl(\Lambda_2/\Lambda(U,W)\bigr), $$
 woraus mit c) die Behauptung folgt.
 \par
       f) Es sei $\zeta \in Z(\Lambda_1\Lambda_2)$. Da 
 $Z(\PGaL(U)) = \{1\}$ und 
 $Z(\PGaL(W)) = \{1\}$ ist, folgt mit
 c), dass $\zeta \in \Lambda(U,W)$ gilt. Somit ist
 $Z(\Lambda_1\Lambda_2) \subseteq Z(\Lambda(U, W))$.
 \par 
       Es sei umgekehrt $\zeta \in Z(\Lambda(U,W))$: Dann gibt es ein $k
 \in Z(K^*)$ mit $(u + w)^\zeta = u + wk$ f\"ur alle $u \in U$ und alle
 $w\in W$. (Hier schon haben wir benutzt, dass $\zeta \in
 Z(\Lambda(U,W))$ gilt.) Ist $\lambda \in \Lambda_1\Lambda_2$, so folgt
 $$
     (u + w)^{\zeta\lambda}
	= (u+wk)^\lambda = u^\lambda + w^\lambda k
	 = (u^\lambda + w^\lambda)^\zeta             
        = (u + w)^{\lambda \zeta}.                 
 $$
 Also ist $\zeta\lambda = \lambda\zeta$. Damit ist alles bewiesen.
 \medskip
       Mit $C_G(U)$ bezeichnen wir den Zentralisator der Teilmenge $U$
 von $G$ in $G$, dh. die Untergruppe aller mit jedem Element
 von $U$ vertauschbaren Elemente von $G$.
 \medskip\noindent
 {\bf 5.5. Satz.} {\it Es sei $V$ ein Vektorraum \"uber dem
 K\"orper $K$ mit $\Rg_K(V) \geq 3$. Es sei $V = U \oplus W$
 mit von $\{0\}$ verschiedenen Unterr\"aumen $U$ und $W$. Ist dann
 $1 \neq \lambda \in Z(\Lambda(U,W))$, so gilt
 \item{a)} Gibt es keine Kollineation in $\PGaL(V)$, die $U$ mit
 $W$ vertauscht, oder ist $\lambda^2 \neq 1$, so ist
 $$ C_{\PGaL(V)}(\lambda) = \Lambda_1\Lambda_2.$$
 \item{b)} Ist $\lambda^2 = 1$ und gibt es eine Kollineation in
 $\PGaL(V)$, die $U$ mit $W$ vertauscht, so ist
 $$ \big|C_{\PGaL(V)}(\lambda) : \Lambda_1\Lambda_2\big| = 2,$$
 und es gilt $Z(C_{\PGaL(V)}(\lambda)) = \{1,\lambda\}$.\par}
 \smallskip

       Beweis. Nach 5.4 f) ist
 $\Lambda_1\Lambda_2 \subseteq C_{\PGaL(V)}(\lambda)$. Es sei $\gamma$ ein
 Element aus dem
 Zentralisator von $\lambda$. Weil $\gamma$ die Unterr\"aume, die
 aus Fixpunkten von $\lambda$ bestehen, unter sich permutiert, ist
 $U^\gamma = U$ und $W^\gamma = W$ oder $U^\gamma = W$ und
 $W^\gamma = U$. Aus 5.4 d) folgt daher, dass
 $|C_{\PGaL(V)}(\lambda):\Lambda_1 \Lambda_2| = 1$ oder 2 ist. Es
 sei nun $\gamma \not\in \Lambda_1 \Lambda_2$. Dann ist also
 $U^\gamma = W$ und $W^\gamma = U$. Es sei $(u + w)^\lambda = u + wk$
 mit einem geeigneten $k \in Z(K^*)$. Dann gibt es ein $r \in K$
 mit
 $$ w^\gamma+u^\gamma k = (u + w)^{\gamma\lambda}
			= (u + w)^{\lambda\gamma}r
			= (u^\gamma + w^\gamma k)r. $$
 Hieraus folgt $k = r$ und $kr = 1$. Also ist $k = 1$ oder $k = -1$, so
 dass a) bewiesen ist.
 \par
       Ist $\lambda^2 = 1$, so ist $(u + w)^\lambda = u - w$, da ja
 $\lambda \neq 1$ ist. Ist $\gamma$ eine Abbildung aus $\GaL(V)$, die $U$ mit $W$
 vertauscht, so ist
 $$\eqalign{
    (u + w)^{\lambda\gamma}
       &= (u - w)^\gamma = u^\gamma - w^\gamma  \cr
       &= -(w^\gamma - u^\gamma) = -(w^\gamma + u^\gamma)^\lambda
	= -(u + w)^{\gamma\lambda}, \cr} $$
 woraus folgt, dass $|C_{\PGaL(V)} (\lambda):\Lambda_1\Lambda_2| = 2$ ist. Aus
 a) folgt schlie\ss\-lich noch, dass
 $$ Z\bigl(C_{\PGaL(V)}(\lambda)\bigr) = \{1, \lambda\} $$
 ist, da $Z(K^*)$ nur eine Involution, n\"amlich $-1$ enth\"alt.
 \medskip
       Sind $U$ und $W$ Unterr\"aume des Vektorraumes $V$, so bezeichnen
 wir mit $E(U,W)$ die Menge aller Quasielationen, deren Zentren in
 $U$ liegen und deren Achsen $W$ enthalten.
 \par
       Sind $X$ und $Y$ zwei Rechtsvektorr\"aume \"uber dem K\"orper $K$,
 so bezeichnen wir mit $\Hom_K(X,Y)$ die Menge aller linearen
 Abbildungen von $X$ in $Y$. Die Menge $\Hom_K(X,Y)$ wird mit der
 punktweise definierten Addition eine abelsche Gruppe.
 \medskip\noindent
 {\bf 5.6. Satz.} {\it Es sei $V$ ein Vektorraum \"uber dem
 K\"orper $K$ mit $\Rg_K(V) \geq 3$. Sind $U$ und $W$
 Unterr\"aume von $V$ mit $U \leq W$, so ist\/ $\E(U,W)$ eine zu\/
 $\Hom_K(V/W,U)$ isomorphe Gruppe.}
 \smallskip
       Beweis. Setze
 $$ H := \bigl\{\alpha \mid \alpha \in \Hom_K(V,U), W \subseteq
		      \Kern(\alpha)\bigr\}. $$
 Dann  ist $H$ eine Untergruppe von
 $\Hom_K(V,U)$, die, wie unmittelbar zu sehen, zu $\Hom_K(V/W,U)$
 isomorph ist. Ist $\alpha \in H$ und definiert man $\tau(\alpha)$ durch
 $$ H := \bigl\{\alpha \mid \alpha \in \Hom_K(V,U),
		  W \subseteq \Kern(\alpha)\bigr\}. $$
 Dann ist $H$ eine Untergruppe von $\Hom_K(V,U)$,
 die wie unmittelbar zu sehen, zu $\Hom_K(V/W,U)$ isomorph ist. Ist
 $\alpha \in H$ und definiert man $\tau(\alpha)$ durch
 $$ x^{\tau(\alpha)} := x + x^\alpha $$
 f\"ur alle $x \in X$ und ist $\tau^*(\alpha)$ die von
 $\tau (\alpha)$ in $L_K(V)$ induzierte
 Kollineation, so ist $\tau^*(\alpha) \in \E(U,W)$. Ferner ist
 $\tau^*$ nach dem, was wir bereits wissen, eine Abbildung von $H$
 auf $\E (U,W)$.
 \par
       Sind $\alpha$, $\beta \in H$, so ist $V^\alpha \leq U \leq W \leq
 \Kern (\beta)$ und daher $\alpha \beta = 0$. Folglich ist
 $$\eqalign{
       x^{\tau(\alpha)\tau(\beta)} &= x + x^\alpha + (x + x^\alpha)^\beta  \cr
	     &= x + x^\alpha + x^\beta = x+x^{\alpha+\beta}                \cr
	     &= x^{\tau(\alpha + \beta)}. \cr}$$
 Dies besagt, dass $\tau^*$ ein Homomorphismus von $H$ auf
 $\E(U,W)$ ist.
 \par 
       Ist $\tau^*(\alpha) = 1$, so gibt es ein $k \in K$ mit
 $xk = x^{\tau(\alpha)} = x + x^\alpha$ f\"ur alle $x \in V$. Nun ist
 $$ v^\alpha \subseteq U \subseteq W \subseteq \Kern(\alpha). $$
 Hieraus folgt, dass $\alpha$ nicht injektiv sein kann, weil das
 dann $V = \{0\}$ zur Folge h\"atte. Somit ist der $\Kern$ von
 $\alpha$ nicht trivial. Es gibt also ein $y \in \Kern(\alpha)$ mit
 $y \neq 0 = y^\alpha$. Dann ist aber $yk = y$ und folglich $k = 1$.
 Hieraus folgt $x^\alpha = 0$ f\"ur alle $x \in V$. Dies zeigt, dass
 $\tau^*$ injektiv ist. Damit ist der Satz bewiesen.
 \medskip
       Zur Bestimmung des Zentralisators einer Quasielation ben\"otigen wir noch
 einige Informationen \"uber die Struktur von Schiefk\"orpern.
 \medskip\noindent
 {\bf 5.7. Satz von Cartan-Brauer-Hua.} {\it Es sei $K$ ein K\"orper und $F$ sei
 ein Teilk\"orper von $K$. Ist dann $x^{-1}Fx = F$ f\"ur alle $x \in K^*$, so ist
 $F \subseteq Z(K)$ oder $F = K$.} \index{Satz von Cartan--Brauer--Hua}{}
 \smallskip
       Beweis. Es sei $k \in K$ und es gebe ein $f \in F$, welches mit $k$ nicht
 vertauschbar sei. Setze $f_1 := kfk^{-1}$. Dann ist also $f_1 \neq f$, so dass
 insbesondere $k \neq -1$ ist. Setze $f_2 := (1 + k)f(1 + k)^{-1}$. Dann ist auch
 $f_2 \neq f$. Nun ist $(1 + k)f = f_2(1 + k)$ und folglich
 $f_2k - kf = f - f_2$. Weiter ist $kf = f_1k$. Also ist
 $$ 0 \neq f - f_2 = f_2k - kf = f_2k - f_1k = (f_2 - f_1)k, $$
 so dass $f_2-f_1 \neq 0$ ist. Daher ist
 $$ k = (f_2 - f_1)^{-1}(f - f_2). $$
 Weil $x^{-1}Fx=F$ f\"ur alle $x\in K^*$ gilt, sind $f_1$, $f_2 \in F$. Daher ist
 $k \in F$. Dies zeigt, dass alle Elemente von $K$, die $F$ nicht
 zentralisieren, in $F$ liegen.
 \par
       Es sei nun $F \not\subseteq Z(K)$. Ist $C$ der Zentralisator von $F$
 in $K$, so ist dann $K - C \neq \emptyset$. Wie wir gerade gesehen
 haben, gilt $K - C \subseteq F$. Weil $K - C$ nicht leer ist, gibt es
 ein $k$ in dieser Menge. Dann ist aber $k + C \subseteq K - C$.
 Folglich gilt $k \in F$. Weil $F$ ein K\"orper ist, ist dann auch
 $C \subseteq F$, so dass, wie behauptet, $K = F$ ist.
 \medskip\noindent
 {\bf 5.8. Korollar.} {\it Es sei $K$ ein K\"orper und $k \in K$.
 Ist $x^{-1}kx$ f\"ur alle $x \in K^*$ mit $k$ vertauschbar, so ist $k \in Z(K)$.}
 \smallskip
       Beweis. Wir d\"urfen annehmen, dass $k \neq 0$ ist. Es sei $C$ der
 Zen\-tra\-li\-sa\-tor von $k$ in $K$ und $F$ sei der von allen $x^{-1}kx$
 erzeugte Teilk\"orper. Dieser ist definiert, da $k$ von $0$
 verschieden ist. Auf Grund unserer Annahme ist $F \subseteq C$.
 Ferner ist $y^{-1}Fy = F$ f\"ur alle $y \in K^*$. Nach dem Satz von
 Cartan-Brauer-Hua ist daher entweder $F \subseteq Z(K)$ oder $F = K$. Nun ist
 aber $k = k^{-1}kk \in F$, so dass in beiden F\"allen
 $k \in Z(K)$ gilt, da aus $F = K$ ja $C = K$ folgt.
 \medskip\noindent
 {\bf 5.9. Satz.} {\it Es sei $K$ ein K\"orper der Charakteristik
 $p > 0$. Ist $S$ eine Untergruppe von $K^*$ und ist $Z$ eine
 Untergruppe von $Z(K^*)$ mit $Z \subseteq S$, so gilt:
 \item{a)} Ist $S/Z$ eine abelsche $p$-Gruppe, so ist $S$ abelsch.
 \item{b)} Ist $S/Z$ ein abelscher $p$-Normalteiler von $K^*/Z$, so
 ist $S \subseteq Z(K^*)$.}
 \smallskip
       Beweis. a) Sind $s$, $t \in S$, so gibt es ein $z \in Z$ mit 
 $t^{-1}st = sz$. Ist $p^a$ die Ordnung von $s$ modulo $Z$, so ist
 $$ s^{p^a} = t^{-1}s^{p^a}t = (t^{-1} st)^{p^a} = s^{p^a} z^{p^a}, $$
 da $Z$ ja im Zentrum von $K^*$ liegt. Also ist
 $z^{p^a} = 1$. Dies hat $z = 1$ zur Folge, da $p$ ja die
 Charakteristik von $K$ ist.
 \par
       b) Ist $k \in K^*$ und $s \in S$, so gibt es ein $t$ in $S$ mit
 $k^{-1}sk = st$. Ist $p^a$ die Ordnung von $s$ modulo $Z$, so ist
 $$ s^{p^a} = k^{-1}s^{p^a} k = (k^{-1} sk)^{p^a} = (st)^{p^a}
                              = s^{p^a} t^{p^a}, $$
 da ja $Z$ im Zentrum von $K$ liegt. Hieraus
 folgt wiederum $t^{p^a} = 1$ und damit $t = 1$, so dass $k^{-1}sk=s$
 ist. Damit ist alles bewiesen.
 \medskip
       Wir nennen {\it Normalreihe\/}\index{Normalreihe}{} einer Gruppe eine
 Kette von Nor\-mal\-tei\-lern dieser Gruppe, wobei sich Kette auf die Inklusion als
 Teilordnung bezieht. Der Leser beachte, dass die Terminologie in der Literatur
 nicht einheitlich ist.
 \medskip\noindent
 {\bf 5.10. Satz.} {\it Es sei $V$ ein Vektorraum \"uber dem
 K\"orper $K$ mit $\Rg_K(V) \geq 3$. Ferner sei $\tau$
 eine von eins verschiedene Quasielation von $L_K(V)$ mit dem
 Zentrum $U$ und der Achse $W$. Setze $H := C_{\PGaL(V)}(\tau)$. Dann
 besitzt $H$ eine Normalreihe 
 $$ \{1\} \subseteq H_2 \subseteq H_1 \subseteq H_0 \subseteq H $$
 mit folgenden Eigenschaften:
 \item{a)} Ist $\Rg_K(V/W) \geq 2$, so ist $H/H_0 \cong \PGaL(V/W)$, ist jedoch
 $\Rg_K(V/W)=1$, so ist $H/H_0 \cong K^*/Z(K^*)$.
 \item{b)} Es ist $H_0/H_1 \cong \GaL(W/U)$.
 \item{c)} Es ist $H_1/H_2 \cong \Hom_K(W/U,U)$.
 \item{d)} Es ist $H_2 \cong \Hom_K(V/W,W)$.\par
 \noindent  
 Ferner gilt:
 \par\noindent
 \item{e)} $Z(H)$ ist zur additiven Gruppe von $Z(K)$ isomorph.
 \item{f)} Es ist $Z(H_1) \cong \Hom_K(V/W,U)$.
 \item{g)} Ist $U \neq W$, so ist $Z(H_1)$ echt in $H_2$ enthalten.
 \item{h)} Es ist $C_H(Z(H_1)) = H_0$.
 \item{i)} Ist die Charakteristik $p$ von $K$ von Null verschieden,
 so ist $H_1$ der gr\"o\ss te aufl\"osbare $p$-Normalteiler von $H$.\par}
 \smallskip
       Beweis. Setze $\alpha := \tau - 1$. Dann ist $U = V^\alpha$ und
 $W = \Kern(\alpha)$. Dabei haben wir $\tau$ mit der linearen
 Abbildung identifiziert, die $\tau$ induziert. Es sei $\gamma \in \GaL(V)$ und
 $\gamma$ induziere in $L_K(V)$ eine Kol\-li\-ne\-a\-ti\-on,
 die mit der Kollineation $\tau$ vertauschbar sei. Es gibt dann ein
 $k \in K^*$ mit $x^{\tau\gamma} = x^{\gamma\tau} k$ f\"ur alle
 $x \in V$. Hieraus folgt
 $$ x^\gamma + x^{\alpha\gamma} = x^\gamma k + x^{\gamma\alpha}k $$
 f\"ur alle $x \in V$. Weil $\tau \neq 1$ ist, ist $U \neq \{0\}$ und daher auch
 $W \neq \{0\}$. Es gibt also ein von $0$ verschiedenes $y \in W$. Wegen
 $W = \Kern(\alpha)$ ist daher
 $$ y^\gamma (1-k) = y^{\gamma \alpha}. $$ W\"are $k \neq 1$, so folgte hieraus
 $$ y^\gamma \in V^\alpha \subseteq \Kern(\alpha) $$
 und damit wiederum $y^\gamma = 0$, was seinerseits den Widerspruch $y = 0$ nach
 sich z\"oge. Also ist doch $k=1$, so dass $\alpha \gamma = \gamma \alpha$ ist.
 Ist umgekehrt $\gamma \in \GaL(V)$ und ist $\gamma$
 mit $\alpha$ vertauschbar, so induziert $\gamma$ ein Element im
 Zentralisator der Quasielation $\tau$. Somit wird $H$ von
 $C_{\GaL(V)}(\alpha)$ induziert, so dass wir zun\"achst den
 Zentralisator von $\alpha$ in $\GaL(V)$ untersuchen. Diesen
 Zentralisator bezeichnen wir im Folgenden mit $C$.
 \par
       Es sei $W = U \oplus S$ und $V = W \oplus T$. Dann ist
 $V = U \oplus S \oplus T$. Wegen $W = \Kern(\alpha)$ ist $T^\alpha = U$ und
 wegen $W \cap T = \{0\}$ gibt es zu jedem $u \in U$ genau ein $t \in T$ mit
 $t^\alpha = u$.
 \par
       Es sei $\gamma \in C$. Dann ist $U^\gamma = U$. Ist $u  \in  U$, so ist
 also
 $$ u^\gamma = t^{\alpha\gamma} = t^{\gamma_{11}\alpha} $$
 mit einem $t \in T$ und einem Endomorphismus $\gamma_{11}$ von $T$. Weil die
 Einschr\"ankung von $\gamma$ auf $U$ in $\GaL(U)$ liegt, ist $\gamma_{11}$
 sogar ein Element aus $\GaL(T)$.
 \par 
       Weil $W$ der Kern von $\alpha$ ist, ist auch $W^\gamma = W$. Ist
 $s \in S$, so gibt es daher ein $\gamma_{21} \in \Hom_K(S,U)$ und
 ein $\gamma_{22} \in \End_K(S)$ mit
 $$ s^\gamma = s^{\gamma_{21}} + s^{\gamma_{22}}. $$
 Ist $s^{\gamma_{22}} = 0$, so ist $s^\gamma \in U$. Weil $\gamma$ bijektiv ist
 und $U$ invariant l\"asst, folgt $s \in U$ und damit $s = 0$, da ja $U \cap S =
 \{0\}$ ist. Dies zeigt, dass $\gamma_{22}$ injektiv ist. Wegen
 $$\eqalign{
      U + S &= W = W^\gamma = U^\gamma + S^\gamma    \cr
	    &\subseteq U + S^{\gamma_{21}} + S^{\gamma_{22}}
	     = U + S^{\gamma_{22}} \subseteq U + S   \cr} $$
 ist $U \oplus S = U + S^{\gamma_{22}}$. Hieraus folgt mittels des
 Modulargesetzes, dass $S = S^{\gamma_{22}}$ ist. Also ist
 $\gamma_{22} \in \GaL(S)$.
 \par 
       Ist schlie\ss lich $t \in T$, so gibt es 
 $\gamma_{31} \in \Hom_K(T,U)$, $\gamma_{32} \in \Hom_K(T,S)$ und
 $\gamma_{33} \in \End_K(T)$
 mit
 $$ t^\gamma = t^{\gamma_{31}} + t^{\gamma_{32}} + t^{\gamma_{33}}. $$
 Wegen $U + S = W = \Kern(\alpha)$ folgt
 $$ t^{\gamma_{11}\alpha} = t^{\alpha\gamma} = t^{\gamma\alpha}
	= t^{\gamma_{31}\alpha} + t^{\gamma_{32}\alpha} + t^{\gamma_{33}\alpha}
	= t^{\gamma_{33}\alpha}. $$
 Weil die Einschr\"ankung von $\alpha$ auf $T$ injektiv ist, folgt hieraus
 die Gleichheit von $\gamma_{11}$ und $\gamma_{33}$.
 \par
       Seien umgekehrt $\gamma_{11} \in \GaL(T)$ und $\gamma_{22} \in \GaL(S)$
 sowie $\gamma_{21} \in \Hom_K(S,U)$, $\gamma_{31} \in \Hom_K(T,U)$ und
 $\gamma_{32} \in \Hom_K(T,S)$. Definieren wir $\gamma$ durch
 $$ (t^\alpha+s+t_1)^\gamma := t^{\gamma_{11}\alpha} + s^{\gamma_{21}}
	+ s^{\gamma_{22}} + t^{\gamma_{31}}_1 + t_1^{\gamma_{32}}
	+ t^{\gamma_{11}}_1, $$
 so ist $\gamma\in C$, wie wir jetzt zeigen werden.
 \par 
       Weil $t^{\gamma_{11}\alpha} + s^{\gamma_{21}} + s^{\gamma_{22}}
 + t_1^{\gamma_{31}} + t_1^{\gamma_{32}}$ im Kern von $\alpha$ liegt, ist
 $$ (t^\alpha+s + t_1)^{\gamma\alpha} = t^{\gamma_{11}\alpha}. $$
 Andererseits ist auch $t^\alpha + s$ ein Element des Kerns von $\alpha$. Daher
 ist
 $$ (t^\alpha + s + t_1)^{\alpha\gamma} = t_1^{\alpha\gamma}. $$
 Benutzt man schlie\ss lich noch einmal die Definition von $\gamma$, so
 folgt $t_1^{\alpha\gamma} = t_1^{\gamma_{11}\alpha}$. Fasst man
 alles zusammen, so erh\"alt man $\alpha \gamma = \gamma \alpha$.
 \par
       Wir zeigen, dass $\gamma$ surjektiv ist. Es ist
 $$ U^\gamma = T^{\alpha\gamma} = T^{\gamma_{11}\alpha} = T^\alpha = U. $$
 Ferner ist
 $$ s^{\gamma_{22}} = s^\gamma - s^{\gamma_{21}}\in S^\gamma + U $$
 und daher $S = S^{\gamma_{22}} \subseteq S^\gamma +U$. Folglich ist
 $$ W = U + S\subseteq U + S^\gamma = U^\gamma + S^\gamma = W^\gamma, $$
 so dass $W=W^\gamma$ ist. Wegen
 $$ t^{\gamma_{11}} = t^\gamma - t^{\gamma_{31}} - t^{\gamma_{32}} \in
                          T^\gamma + W $$
 ist $T = T^{\gamma_{11}} \subseteq T^\gamma + W$ und daher
 $$ V = W + T \subseteq W + T^\gamma = W^\gamma + T^\gamma = V^\gamma. $$
 Dies besagt aber, dass $\gamma$ surjektiv ist.
 \par
       Es sei schlie\ss lich $0 = (t^\alpha + s + t_1)^\gamma$. Weil $V$ die
 direkte Summe von $U$, $S$ und $T$ ist, folgt zun\"achst
 $t_1^{\gamma_{11}} = 0$ und damit $t_1 = 0$. Folglich ist
 $$ t^{\gamma_{11}\alpha} + s^{\gamma_{21}} + s^{\gamma_{22}} = 0. $$
 Dies zieht $s^{\gamma_{22}} = 0$ und damit $s = 0$ nach sich. Dann ist
 aber $t^{\gamma_{11}\alpha} = 0$, was schlie\ss lich $t = 0$
 impliziert. Also ist $\gamma$ auch injektiv. Daher gilt das
 folgende Zwischenresultat:
 \par
       Ist $\gamma \in C$, so gibt es Abbildungen $\gamma_{11} \in \GaL(T)$,
 $\gamma_{22} \in \GaL(S)$, $\gamma_{21} \in \Hom_K(S,U)$,
 $\gamma_{31} \in \Hom_K(T,U)$ und $\gamma_{32} \in \Hom_K(T,S)$ mit
 $$ (t^\alpha + s + t_1)^\gamma = t^{\gamma_{11}\alpha} + s^{\gamma_{21}}
           + s^{\gamma_{22}} + t^{\gamma_{31}}_1 + t^{\gamma_{32}}_1
	   + t^{\gamma_{11}}_1; $$
 sind umgekehrt $\gamma_{ij}$ wie gerade beschrieben gegeben und
 definiert man $\gamma$ durch diese Gleichung, so ist $\gamma \in C$.
 \par 
       Es sei $C_0$ diejenige Untergruppe von $C$, die auf dem Faktorraum
 $V/W$ die Identit\"at induziert. Es ist $V/W = \{t + W \mid t\in T\}$.
 Ist $\gamma \in C_0$, so ist also $t^{\gamma_{11}} + W = t + W$ f\"ur
 alle $t \in T$. Wegen $T \cap W =\{0\}$ folgt, dass $\gamma_{11} = 1$
 ist. Also gilt $\gamma \in C_0$ genau dann, wenn $\gamma_{11} = 1$
 ist. Dies impliziert einmal, dass
 $$ C/C_0 \cong \GaL(T) = \GaL(V/W) $$ 
 ist und dass $C_0$ auch den Unterraum $U$ vektorweise festl\"asst.
 \par 
       Es sei $C_1$ diejenige Untergruppe von $C_0$, die auf $W/U$ die
 Identit\"at induziert. Ist $\gamma \in C_1$, so folgt, da
 $W/U = \{s + U \mid s \in S\}$ ist, dass $s^{\gamma_{22}} + U = s + U$ ist f\"ur
 alle $s \in S$. Weil $S$ und $U$ trivialen Schnitt haben,
 folgt $\gamma_{22} = 1$. Es ist also genau dann $\gamma \in C_1$,
 wenn $\gamma_{11} = 1$ und $\gamma_{22} = 1$ ist. Hieraus folgt
 wiederum
 $$ C_0/C_1 \cong \GaL(S) \cong \GaL(W/U). $$
 \par
       Es sei $C_2$ die Untergruppe von $C_1$, die $W$ vektorweise
 festl\"asst. Dann folgt, dass $C_1/C_2$ zu einer Gruppe von
 Quasitransvektionen von $W$ isomorph ist, deren Zentren in $U$
 liegen und deren Achsen $U$ enthalten. Es sei umgekehrt $\lambda$
 eine Quasitransvektion von $W$, deren Zentrum in $U$ liegt und
 deren Achse $U$ enth\"alt. Dann ist also $w^\lambda = w + w^\beta$
 mit $\beta \in \Hom_K(W,U)$ und $U \subseteq \Kern (\beta)$.
 Definiert man nun $\gamma$ durch
 $$ (t^\alpha + s + t_1)^\gamma := t^\alpha + s^\beta + s + t_1, $$
 so ist $\gamma \in C_1$ und wegen $U \subseteq \Kern(\beta)$ induziert $\gamma$
 in $W$ die Abbildung $\lambda$. Nach 5.6 ist folglich
 $$ C_1/C_2 \cong \Hom_K(W/U,U). $$
 \par
       Weil die Elemente von $C_2$ auf $W$ und $V/W$ die Identit\"at induzieren,
 ist $C_2$ eine Gruppe von Quasitransvektionen von $V$, deren Zentren in $W$
 ent\-hal\-ten sind und deren Achse $W$ enthalten. Ist umgekehrt $\mu$
 eine solche Quasitransvektion, so gibt es ein $\mu_{31} \in \Hom_K(T,U)$ und ein
 $\mu_{32} \in \Hom_K(T,S)$ mit
 $$ t^\mu = t + t^{\mu_{31}} + t^{\mu_{32}}. $$
 Hieraus folgt
 $$ (t^\alpha + s + t_1)^\mu = t^\alpha + s + t^{\mu_{31}} + t^{\mu_{32}} + t, $$
 so dass $\mu \in C_2$ ist. Somit ist
 $$ C_2 \cong \Hom_K(V/W,W). $$
 Wir betrachten nun die von $C$, $C_0$, $C_1$ und $C_2$
 in $L_K(V)$ induzierten Kollineationsgruppen $H$, $H_0$, $H_1$ und
 $H_2$. Weil $C_i$ offensichtlich ein Normalteiler von $C$ ist, ist
 $H_i$ ein Normalteiler von $H$. Auf Grund des Zwischenresultates
 ist der Kern des Homomorphismus von $C$ auf $H$ gerade die Gruppe
 $M$ aller Multiplikationen $\mu (k)$ mit $k \in Z(K^*)$, wobei
 $\mu (k)$, wie schon zuvor, durch $v^{\mu(k)} := vk$ definiert ist.
 Es folgt, dass $MC_i$ das Urbild von $H_i$ in $C$ ist. Ferner ist
 $M \cap C_i = \{1\}$, da $U$ wegen $\tau \neq 1$ mindestens den
 Rang 1 hat und $U$ von $C_i$ vektorweise festgelassen wird.
 Setzt man noch $C_3 := \{1\}$ und $H_3 := \{1\}$, so folgt  f\"ur
 $i := 0$, 1, 2, dass
 $$\eqalign{
       H_i/H_{i+1} &\cong (MC_i)/(MC_{i+1})
                    = (MC_{i+1}C_i)/(MC_{i+1})                  \cr
                   &\cong C_i/(C_i \cap MC_{i+1}) = C_i/C_{i+1} \cr}$$
 ist, da ja wegen $C_{i+1} \subseteq C_i$ auf Grund des
 Modulargesetzes die Glei\-chun\-gen
 $$ C_i \cap MC_{i+1} = (C_i \cap M)C_{i+1} = C_{i+1} $$
 gelten. Damit ist die G\"ultigkeit von b), c) und d) bewiesen.
 \par
       Betrachtet man die Abbildungen $\gamma \in C$ mit $\gamma_{21} = 0$,
 $\gamma_{31} = 0$, $\gamma_{32} = 0$ und $\gamma_{22}=1$, so bilden diese
 eine Untergruppe $D$ von $C$. Es ist $C = DC_0$ und $D \cap C_0 = \{1\}$. Ferner
 sieht man sofort, dass $D \cap MC_0 = Z(D)$ ist. Daher ist
 $$ H/H_0 \cong C/(MC_0) = (DMC_0)/(MC_0)\cong D/(D \cap MC_0) = D/Z(D), $$
 so dass auch a) gilt.
 \par
       Als N\"achstes bestimmen wir $Z(H)$. Es sei $\delta M \in Z(H)$.
 Ist $\gamma \in C$, so gibt es wegen $\Rg_K(V) \geq 3$ nach 1.2
 ein $k_\gamma \in Z(K^*)$ mit $v^{\gamma \delta} = v^{\delta\gamma} k_\gamma$
 f\"ur alle $v \in V$. Insbesondere ist dann
 $$
    t^{\gamma_{11}\delta_{11}\alpha} = t^{\gamma_{11}\alpha\delta}
			      = t^{\alpha \gamma \delta}
			     =t^{\alpha \delta \gamma}k_\gamma
			      = t^{\delta_{11}\gamma_{11}\alpha}k_\gamma     
			     = (t^{\delta_{11} \gamma_{11}} k_\gamma)^\alpha 
			     $$
 und daher
 $$ t^{\gamma_{11}\delta_{11}} = t^{\delta_{11}\gamma_{11}} k_\gamma $$
 f\"ur alle $t \in T$. Weil
 $\gamma_{11}$ alle Werte in $\GaL(T)$ annimmt und das Zentrum der
 Gruppe $\PGaL(T)$ trivial ist, folgt, dass $\delta_{11}$ alle
 Punkte von $L_K(T)$ festl\"asst.
 \par
       Ist $\Rg_K(T) \leq 2$, so gibt es nach 1.2 ein $k \in Z(K^*)$ mit
 $t^{\delta_{11}} = tk$ f\"ur alle $t \in T$. Es sei $\Rg_K(T) = 1$
 und $T = t_0K$. Es gibt dann ein $k \in K^*$ mit
 $t^{\delta_{1}}_0 = t_0k$. Ist nun $m \in K^*$, so gibt es ein
 $\gamma_{11}\in \GaL(T)$ mit $t_0^{\gamma_{11}} = t_0m$. Es gibt
 ferner ein $\gamma \in C$, welches dieses $\gamma_{11}$ induziert.
 Damit folgt
 $$ t_0km = t^{\delta_{11}}_0 m = t_0^{\gamma_{11}\delta_{11}}
          = t_0^{\delta_{11} \gamma_{11}}k_\gamma = t_0^{\gamma_{11}} k_\gamma
	  = t_0mkk_\gamma. $$
 Hieraus folgt die Gleichung $km = mkk_\gamma$. Weil $k_\gamma$ im Zentrum
 von $K$ liegt ist daher
 $$ km^{-1}km = k^2k_\gamma = kk_\gamma k = m^{-1}kmk. $$
 Weil dies f\"ur alle von Null verschiedenen $m$ aus
 $K$ gilt, folgt mit 5.8, dass $k \in Z(K^*)$ ist. Es gilt also in
 jedem Falle, dass $t^\delta_{11} = tk$ f\"ur alle $t \in T$ gilt,
 wobei $k$ ein passendes Elemente aus dem Zentrum von $K$ ist.
 Hiermit folgt weiter
 $$ t^{\gamma_{11}}k = t^{\gamma_{11}\delta_{11}}
		     = t^{\delta_{11} \gamma_{11}}k_\gamma
		     = t^{\gamma_{11}} kk_\gamma $$
 und dann $k=kk_\gamma$, so dass $k_\gamma = 1$ ist. Da dies f\"ur alle
 $\gamma \in C$ gilt, liegt $\delta$ im Zentrum von $C$, so dass also $Z(C)$ das
 Urbild von $Z(H)$ ist.
 \par
       Weil $k$ im Zentrum von $K$ liegt, folgt, dass die durch
 $v^\eta := v^\delta k^{-1}$ erkl\"arte Abbildung $\eta$ in $C$  und damit
 in $C_0$ liegt. Da $\eta$ und $\delta$ die gleiche Kollineation in
 $L_K(V)$ induzieren, folgt, dass $Z(C) \cap C_0$ bei dem
 Homomorphismus von $C$ auf $H$ auf $Z(H)$ abgebildet wird. Weil
 die Einschr\"ankung dieses Homomorphismus auf $C_0$ wegen $M \cap C_0 = \{1\}$
 ein Monomorphismus ist, sind $Z(C) \cap C_0$ und $Z(H)$ isomorph.
 \par
       Es sei $\gamma \in C$ und $\delta \in Z(C) \cap C_0$. Wegen
 $\delta_{11} = 1$ ist dann
 $$ s^{\gamma\delta} = (s^{\gamma_{21}} + s^{\gamma_{22}})^\delta
                     = s^{\gamma_{21}} + s^{\gamma_{22}\delta_{21}}
		       + s^{\gamma_{22}\delta_{22}}. $$
 Andererseits ist
 $$ s^{\delta\gamma} = (s^{\delta_{21}} + s^{\delta_{22}})^\gamma
		     = s^{\delta_{21}\gamma} + s^{\delta_{22}\gamma_{21}}
		       + s^{\delta_{22}\gamma_{22}}. $$
 Wegen $\gamma \delta = \delta \gamma$ folgen die Gleichungen
 $$ s^{\gamma_{21}} + s^{\gamma_{22}\delta_{21}}
	       = s^{\delta_{21}\gamma} + s^{\delta_{22}\gamma_{21}} $$
 und
 $$ s^{\gamma_{22}\delta_{22}} = s^{\delta_{22}\gamma_{22}}. $$
 Da $\gamma_{22}$ aller Werte in $\GaL(S)$ f\"ahig ist, folgt wieder
 die Existenz eines $k \in Z(K^*)$ mit $s^{\delta_{22}} = sk$ f\"ur
 alle $s \in S$. Hieraus folgt mittels der ersten Gleichung
 $$ s^{\gamma_{21}} (1 - k)
      = s^{\delta_{21} \gamma} - s^{\gamma_{22}\delta_{21}}. $$
 Mit Hilfe des Zwischenresultates folgt, dass $\gamma_{21}$ auch dann noch aller
 Werte in $\Hom_K(S,U)$ f\"ahig ist, wenn $\gamma$ auf $U$ die Identit\"at
 induziert und $\gamma_{22} = 1$ ist. Daher gilt
 $$ s^{\gamma_{21}}(1 - k) = 0 $$
 f\"ur alle $s \in S$ und alle $\gamma_{21} \in \Hom_K(S,U)$. Ist $S \neq \{0\}$,
 so folgt $k = 1$. Daher ist in jedem Falle $\delta_{22} = 1$ und 
 $$ 0 = s^{\delta_{21} \gamma}-s^{\gamma_{22} \delta_{21}}. $$
 Hieraus folgt weiter, wiederum auf Grund des Zwischenresultates, dass
 $\delta_{21}=0$ ist, es sei denn, es ist $K = \GF(2)$ und
 $\Rg_K(U) = \Rg_K(S) = 1$ und folglich $\Rg_K(V) = 3$.
 Dieser Fall kann aber nicht eintreten, wie wir gleich sehen werden.
 \par
       Nun ist einerseits
 $$ t^{\gamma\delta} = t^{\gamma_{31}} + t^{\gamma_{32}\delta_{21}}
	 + t^{\gamma_{32}} + t^{\gamma_{11}\delta_31} +
	     t^{\gamma_{11}\delta_{32}} + t^{\gamma{11}} $$
 und andererseits
 $$ t^{\delta\gamma} = t^{\delta_{31}\gamma} + t^{\delta_{32}\gamma_{21}}
	 + t^{\delta_{32}\gamma_{22}} + t^{\gamma_{31}} + t^{\gamma_{32}}
	 + t^{\gamma_{11}}. $$
 Wegen $\gamma\delta = \delta\gamma$ folgen die beiden Gleichungen
 $$ t^{\gamma_{32} \delta_{21}} + t^{\gamma_{11} \delta_{31}}
	  = t^{\delta_{31} \gamma} + t^{\delta_{31}\gamma_{21}} $$
 und
 $$ t^{\gamma_{11} \delta_{32}} = t^{\delta_{32}\gamma_{22}}. $$
 W\"are $\delta_{21} \neq 0$, so w\"are also $\Rg_K(U) = \Rg_K(S) = \Rg_K(T)=1$
 und $K = \GF(2)$. Daher w\"are $\gamma_{11} = 1$ und $\gamma_{22} = 1$ und
 folglich $t^{\gamma_{32}\delta_{21}} = t^{\delta_{32}\gamma_{21}}$, da ja
 dann $t^{\delta_{31} \gamma} = t^{\delta_{31}}$ ist. Mit
 $\gamma_{21} = 0$ und $\gamma_{32} \neq 0$ folgte der Widerspruch
 $\delta_{21} = 0$. Also ist in jedem Falle $\delta_{21} = 0$.
 \par
       W\"are $\delta_{32} \neq 0$, so folgte wiederum, dass $K = \GF(2)$ und
 $$ \Rg_K(U) = \Rg_K(S) = \Rg_K(T) = 1 $$
 w\"are. Dann w\"are $\gamma_{11} = 1$ und folglich
 $$ t^{\delta_{31}} = t^{\delta_{31}} + t^{\delta_{32} \gamma_{21}}, $$
 woraus $\gamma_{21} = 0$ folgte. Das z\"oge aber den Widerspruch
 $$ \{0\} = \Hom_K(S,U) \neq \{0\} $$
 nach sich. Also ist doch $\delta_{32} = 0$.
 \par
       Damit ist gezeigt: Ist $\delta \in Z(C) \cap C_0$, so ist
 $$ (t^\alpha + s + t_1)^\delta = t^\alpha + s + t^{\delta_{31}}_1 + t_1, $$
 wobei $\delta_{31}$ der Bedingung $\gamma_{11}\delta_{31} = \delta_{31}\gamma$
 f\"ur alle $\gamma \in C$ gen\"ugt. Umgekehrt sieht man unmittelbar, dass auch
 jede solche Abbildung $\delta$ in $Z(C) \cap C_0$ liegt.
 \par
       Ist $1 \neq \delta \in Z(C) \cap C_0$, so ist $\delta_{31} \neq 0$.
 Bezeichnet man mit $\beta$ die Einschr\"ankung von $\alpha$
 auf $T$, so ist $\beta$ ein Isomorphismus von $T$ auf $U$. Es
 folgt, dass auch $\delta_{31}\beta^{-1} \neq 0$ ist. Ferner ist
 $$ \gamma_{11} \delta_{31} = \delta_{31} \gamma
	= \delta_{31} \beta^{-1} \alpha \gamma
	= \delta_{31} \beta^{-1}\gamma_{11} \alpha, $$
 so dass also
 $$ \gamma_{11}\delta_{31} \beta^{-1} = \delta_{31}\beta^{-1}\gamma_{11} $$
 ist. Hieraus folgt, dass $\delta_{31}\beta^{-1} \in Z(\GaL(T))$ ist. Dies
 impliziert, auch im Falle $\Rg_K(T) = 1$, dass es ein $k \in Z(K^*)$ gibt mit
 $$ t^{\delta_{31} \beta^{-1}} =tk $$
 dh. mit
 $$ t^{\delta_{31}} = t^\alpha k $$
 f\"ur alle $t \in T$. Damit ist gezeigt, dass
 $$ (t^\alpha + s + t_1)^\delta = t^\alpha + s + t^\alpha_1k + t_1 $$
 ist.
 \par
       Ist umgekehrt $k \in Z(K^*)$ und definiert man $\delta$ durch die
 gerade etablierte Gleichung, so liegt $\delta$ in $Z(C) \cap C_0$.
 Hieraus folgt die Behauptung unter e).
 \par
       Es seien $\gamma$, $\delta \in C_1$ und es gelte $\delta M \in
 Z(H_1)$. Es gibt dann ein $k \in Z(K^*)$ mit $v^{\gamma \delta} =
 v^{\delta \gamma} k$. Hieraus folgt wegen $\gamma_{11} =
 \delta_{11} = 1$, dass
 $$ t^{\alpha\gamma\delta} = t^{\alpha\delta\gamma} k $$
 ist. Dies hat $k = 1$ zur Folge.
 Hieraus folgt weiter, dass $Z(H_1)$ von $Z(C_1)$ induziert wird.
 \par
       Es sei $\gamma \in C_1$ und $\delta  \in Z(C_1)$. Dann ist
 $$ t^{\gamma \delta} = t^{\gamma_{31}} + t^{\gamma_{32}\delta_{21}}
	    + t^{\gamma_{32}} + t^{\delta_{31}} + t^{\delta_{32}} + t $$
 und
 $$ t^{\delta\gamma} = t^{\delta_{31}} + t^{\delta_{32}\gamma_{21}}
	   + t^{\delta_{32}} + t^{\gamma_{31}} + t^{\gamma_{32}} + t. $$
 Somit ist $\gamma_{32}\delta_{21} = \delta_{32}\gamma_{21}$ f\"ur alle
 $\gamma_{21}$. Hieraus folgt $\delta_{21} = 0$ und $\delta_{32} = 0$. Nun sind
 die Elemente aus $C_2$ gerade die Elemente $\eta$ mit $\eta_{11} = 1$,
 $\eta_{21} = 0$. Ist nun $U \neq W$, so gibt es ein $\eta \in C_2$ mit
 $\eta_{32} \neq 0$, so dass in diesem Falle $Z(C_1)$ echt in $C_2$
 enthalten ist. Hieraus folgt die G\"ultigkeit von g).
 \par
       Es sei $\delta \in C_2$. Dann ist $\delta_{11} = 1$, $\delta_{22} = 1$
 und $\delta_{21} = 0$. Ist auch noch $\delta_{32} = 0$, so
 \"uberzeugt man sich leicht, dass $\delta \in Z(C_1)$ ist. Dies
 impliziert f).
 \par
       Es sei $\gamma \in C$ und $v^{\gamma \delta} = v^{\delta \gamma}k_\delta$
 f\"ur alle $\delta \in Z(C_1)$, wobei $k_\delta$ ein von
 $\delta$ abh\"angiges Element aus $Z(K^*)$ sei. Ist $u \in U$, so ist auch
 $u^\gamma \in U$. Daher ist $u^\delta = u$ und $u^{\gamma\delta} = u^\gamma$.
 Also ist
 $$ u^\gamma = u^{\gamma\delta} = u^{\delta \gamma}k_\delta
	     = u^\gamma k_\delta, $$
 so dass $k_\delta=1$ ist, da ja $U \neq \{0\}$ ist. Somit sind $\gamma$
 und $\delta$ miteinander vertauschbar.
 \par
       Es ist
 $$ t^{\gamma \delta} = t^{\gamma_{31}} + t^{\gamma_{32}}
			 + t^{\gamma_{11}\delta_{31}} + t^{\gamma_{11} }$$
 und
 $$ t^{\delta \gamma} = t^{\delta_{31} \gamma} + t^{\gamma_{31}}
			   + t^{\gamma_{32}} + t^{\gamma_{11}}. $$
 Also ist $\gamma_{11}\delta_{31} = \delta_{31}\gamma$ f\"ur alle
 $\delta_{31} \in \Hom_K(T,U)$. Dies impliziert, wie schon
 mindestens zweimal gesehen, dass $t^{\gamma_{11}} = tz$ ist mit
 einem $z \in Z(K^*)$. Nun induzieren $\gamma$ und die durch
 $v^\eta := v^\gamma z^{-1}$ definierte Abbildung $\eta$ die
 gleiche Kollineation in $L_K(V)$. Daher ist $C_H (Z(H_1))
 \subseteq H_0$. Eine triviale Rechnung zeigt schlie\ss lich, dass
 $C_H(Z(H_1)) = H_0$ ist.
 \par
       Aus c) und d) folgt, dass $H_1^{(2)}=\{1\}$ ist. Daher ist $H_1$
 aufl\"osbar. Ist $p := \Char (K) >0$ , so folgt ebenfalls aus c) und
 d), dass $\eta^{p^2} = 1$ ist f\"ur alle $\eta \in H_1$. Somit ist
 $H_1$ ein aufl\"osbarer $p$-Normalteiler von $H$. Es sei $G$ ein
 aufl\"osbarer $p$-Normalteiler von $H$. Dann ist auch $GH_1$ ein
 aufl\"osbarer $p$-Normalteiler, so dass wir annehmen d\"urfen,
 dass $H_1 \subseteq G$ ist.
 \par
       Ist $G \cap H_0 \neq H_1$, so ist $(G \cap H_0)/H_1$ ein nicht
 trivialer $p$-Nor\-mal\-tei\-ler von $H_0/H_1$. Aus b) folgt dann, dass
 $\GaL(W/U)$ einen nicht trivialen $p$-Normalteiler $N$ enth\"alt.
 Hieraus folgt $\Rg_K(W/U) \leq 2$, denn andernfalls enthielte
 $K^*$ einen nicht trivialen $p$-Normalteiler, was wegen $\Char(K) = p > 0$ nicht
 der Fall ist. Nun ist $Z(\GaL(W/U))$ zu $Z(K^*)$
 isomorph, so dass auch das Zentrum von $\GaL(W/U)$ keinen nicht
 trivialen $p$-Normalteiler enth\"alt. Weil $\SL(W/U)$ keine
 $p$-Gruppe ist, folgt $\SL(W/U) \not\subseteq N$. Nach 2.9 ist daher
 $W/U$ der Vektorraum vom Range 2 \"uber $\GF(2)$ oder $GF(3)$. Dies
 kann aber nicht sein, da $\GL(2,2)$ keinen nicht trivialen
 $2$-Normalteiler und $\GL(2,3)$ keinen nicht trivialen
 3-Normalteiler besitzt. Also ist $G \cap H_0 = H_1$.
 \par
       Es sei $G \neq H_1$. Wegen $G \cap H_0 = H_1$ ist dann
 $(GH_0)/H_0$ ein nicht trivialer $p$-Normalteiler von $H/H_0$.
 W\"are $\Rg_K(V/W) \geq 2$, so h\"atte $\PGaL(V/W)$ einen nicht
 trivialen $p$-Normalteiler, was mit Hilfe von 2.8 auf einen
 Widerspruch f\"uhrte. Also ist $\Rg_K(V/W) = 1$ und es folgt, dass
 $K^*/Z(K^*)$ einen nicht trivialen, aufl\"osbaren $p$-Normalteiler
 $N/Z(K^*)$ hat. Es gibt eine nat\"urliche Zahl $i$ mit $N^{(i)}
 \not\subseteq Z(K^*)$ und $N^{(i+1)} \subseteq Z(K^*)$. Weil
 $N^{(i)}$ in $N$ charakteristisch ist, ist $N^{(i)}$ ein
 Normalteiler von $K^*/Z(K^*)$ ist. Nach 5.9 b) ist daher $N^{(i)}
 \subseteq Z(K^*)$. Dieser Widerspruch zeigt, dass $G = H_1$ ist, so
 dass 5.10 vollst\"andig bewiesen ist.
 \medskip
       Wir schlie\ss en diesen Abschnitt mit einem Korollar zu Satz 5.10,
 das zu beweisen dem Leser als \"Ubungsaufgabe \"uberlassen sei.
 \medskip\noindent
 {\bf 5.11. Korollar.} {\it Ist $V$ ein Vektorraum \"uber dem
 K\"orper $K$ mit $\Rg_K(V) \geq 3$, so sind die folgenden
 Aussagen \"aquivalent:
 \item{a)} $L_K(V)$ ist pappossch.
 \item{b)} Ist $(OP,H)$ ein inzidentes Punkt-Hyperebenenpaar von
 $L_K(V)$ und sind $\sigma$, $\tau$ zwei von eins verschiedene
 Elationen aus $\E(P,H)$, so ist $C_{\PGaL(V)}(\sigma) =
 C_{\PGaL(V)}(\tau)$.
 \item{c)} Es gibt ein inzidentes Punk-Hyperebenenpaar $(P,H)$ von
 $L_K(V)$, so dass f\"ur alle von eins verschiedenen Elationen
 $\sigma$ und $\tau$ aus $\E(P,H)$ gilt, dass $C_{\PGaL(V)}
 (\sigma) = C_{\PGaL(V)} (\tau)$ ist.\par}

\mysection{6. Zentralisatoren von Involutionen}

\noindent
 Involutorische Streckungen und involutorische Elationen lassen
 sich an Hand ihrer Zentralisatoren von anderen involutorischen
 Kol\-li\-ne\-a\-ti\-o\-nen unterscheiden. Dies zu zeigen ist 
 Ziel dieses
 Abschnitts. Im n\"achsten Abschnitt 
 machen wir dann davon Gebrauch
 und zeigen, dass die projektive Geometrie $L_K(V)$ durch
 die Gruppe $\PGaL(V)$ bis auf Isomorphie oder Antiisomorphie
 festgelegt ist. 
 \medskip\noindent
 {\bf 6.1. Satz.} {\it Es sei $V$ ein Vektorraum \"uber dem
 K\"orper $K$ mit $\Rg_K(V) \geq 1$. Ferner sei $\alpha$
 ein innerer Automorphismus von $K$. Ist dann $\sigma \in \GammaL(V)$ eine
 semilineare Abbildung mit begleitendem Automorphismus
 $\alpha$, ist $\gamma \in \GaL(V)$, sind $k$, $m \in K$ und gilt
 $v^{\sigma^2} = vk$ sowie $v^{\sigma\gamma} = v^{\gamma\sigma}m$
 f\"ur alle $v \in V$, so ist $m = 1$ oder $m = -1$.}
 \smallskip
       Beweis. Es ist
 $$ v^{\gamma \sigma} mx^\alpha = v^{\sigma\gamma} x^\alpha
	       = (vx)^{\sigma\gamma} = (vx)^{\gamma\sigma} m
	       = v^{\gamma \sigma} x^\alpha m. $$
 Also ist $mx^\alpha = x^\alpha m$ f\"ur alle $x \in K$. Somit liegt $m$ im
 Zentrum von $K$. Weil $\alpha$ ein innerer Automorphismus von $K$ ist, ist daher
 $m^\alpha = m$. Hiermit folgt, dass 
 $$ v^\gamma k = v^{\gamma\sigma^2} = (v^{\sigma\gamma} m^{-1})^\sigma
	       = v^{\sigma\gamma\sigma} m^{-1} = v^{\sigma^2\gamma}m^{-2}
	       = v^\gamma km^{-2} $$
 ist. Somit ist $k = km^{-2}$, woraus $m^2 = 1$ und damit die Behauptung folgt,
 da $m$ ja im Zentrum liegt.
 \medskip\noindent
 {\bf 6.2. Satz.} {\it Es sei $V$ ein Vektorraum \"uber dem
 K\"orper $K$ mit $\Rg_K(V) \geq 3$. Es sei $\sigma$ eine
 Involution aus $\PGaL(V)$, die in $L_K(V)$ keinen Fixpunkt habe,
 und $C$ sei der Zentralisator von $\sigma$ in $\PGaL(V,K)$. Mit
 $\rho$ bezeichnen wir eine lineare Abbildung, die $\sigma$
 induziert, und mit $F$ den Teilring von $\End_Z(V)$, der von $K$
 und $\rho$ erzeugt wird. Dann gilt:
 \item{a)} $F$ ist ein K\"orper und quadratische Erweiterung von $K$.
 \item{b)} $V$ ist auch ein Vektorraum \"uber $F$.
 \item{c)} Ist $\Char (K) = 2$, so ist $C/Z(C) \cong \PGaL(V,F)$
 und $Z(C) \cong Z(F^*)/Z(K^*)$.
 \item{d)} Ist $\Char (K) \neq 2$, so besitzt $C$ einen
 Normalteiler $N$ vom Index $2$ mit $N/Z(N) \cong \PGaL(V,F)$ und
 $Z(N) \cong Z(F^*)/Z(K^*)$. Ferner ist $Z(C) = \{1,\sigma\}$.\par}
 \smallskip
       Beweis. Weil $\sigma$ eine Involution ist, gibt es ein $k \in Z(K^*)$ mit
 $v^{\rho^2} = vk$ f\"ur alle $v \in V$. Weil $\rho$
 eine lineare Abbildung ist und $k$ im Zentrum von $k$ liegt, ist
 $$ F = \{a + b \rho \mid a,b \in K\}. $$
 Es sei $0 \neq v \in V$. Dann ist $vF$ ein Teilmodul des $F$-Moduls $V$ und es
 ist $vF = vK + v^\rho K$. Weil $\sigma$ keinen Fixpunkt hat, ist
 $v^\rho \not\in vK$. Daher ist sogar $vF = vK \oplus v^\rho K$. Dies
 impliziert, dass aus $0 = v(a + b\rho) = va + v^\rho b$ folgt, dass
 $a = b = 0$ ist.
 \par
       Es sei $R$ ein von $\{0\}$ verschiedenes Rechtsideal von $F$. Dann
 ist nach dem gerade Bewiesenen $vR$ ein von Null verschiedener
 $F$-Teil\-mo\-dul von $vF$. Weil $K$ in $F$ enthalten ist, ist $vF$
 auch ein Teilraum des $K$-Vektorraumes $V$. Weil $\sigma$ keinen
 Fixpunkt hat, folgt aus $vR \neq \{0\}$, dass der Rang von $vR$
 als $K$-Vektorraum mindestens $2$ ist. Da andererseits
 $\Rg_K(vF) = 2$ ist, ist $vR = vF$. Es sei nun $f \in F$. Es gibt dann
 ein $r \in R$ mit $vr = vf$. Hieraus folgt $v(r - f)=0$. Dies hat
 $r - f=0$ zur Folge, wie wir oben gesehen haben. Also ist $F = R$, so
 dass $F$ das einzige von Null verschiedene Rechtsideal von $F$
 ist. Dies impliziert, dass $F$ ein K\"orper ist. Damit sind a) und
 b) bewiesen.
 \par
       Es sei $\alpha$ die durch $(a + b\rho)^\alpha := a - b\rho$ definierte
 Abbildung. Tri\-vi\-a\-le Rechnungen zeigen, dass $\alpha$ ein
 Automorphismus von $F$ ist. Ist die Charakteristik von $K$ gleich
 $2$, so ist $\alpha = 1$. Ist die Charakteristik von $K$ von $2$
 verschieden, so ist $\alpha$ involutorisch. Ist nun $\gamma$ ein
 Element aus $\GaL(V,K)$, welches eine Kollineation aus $C$ in
 $L_K(V)$ hervorruft, so ist nach 6.1 entweder $\rho \gamma = \gamma\rho$ oder
 $\rho\gamma = -\gamma\rho$. Im ersten Falle
 ist $(v(a + b\rho))^\gamma = v^\gamma(a + b \rho)$, dh., es ist
 $\gamma \in \GaL(V,F)$, und im zweiten Falle ist
 $(v(a + b\rho))^\gamma = v^\gamma(a + b\rho)^\alpha$, dh., $\gamma$
 ist eine semilineare Abbildung des $F$-Vektorraumes $V$ mit
 begleitendem Automorphismus $\alpha$. Umgekehrt sieht man, dass
 jedes Element aus $\GaL(V,F)$ und jede bijektive semilineare
 Abbildung des $F$-Vektorraumes $V$ mit $\alpha$ als begleitendem
 Automorphismus eine Kollineation in $L_K(V)$ induziert, die in
 $C$ liegt. Aus diesem Bemerkungen folgen die Aussagen c) und d)
 bis auf die Aussage \"uber das Zentrum von $C$ unter d).
 \par
       Es ist klar, dass $Z(C) \subseteq Z(N)$ ist. Es sei also $f \in Z(F^*)$
 und $\gamma$ sei eine bez\"uglich $\alpha$ semilineare
 Abbildung des $F$-Vektorraumes $V$. Ferner induziere die durch
 $v^\zeta := vf$ definierte Abbildung $\zeta$ eine Kollineation aus
 $Z(C)$. Dann gibt es ein $k \in K^*$ mit
 $$ vfk = v^\zeta k = v^{\gamma^{-1} \zeta \gamma}
		    = (v^{\gamma^{-1}} f)^\gamma = vf^\alpha. $$
 Also ist $fk=f^\alpha$. Weil $K$ von
 $\alpha$ elementweise festgelassen wird und weil $\alpha^2=1$ ist,
 ist also $k^2 = 1$ und damit $k = 1$ oder $k = -1$. Im ersten Fall folgt
 $f \in K$, so dass $\zeta$ die Identit\"at induziert. Im zweiten
 Falle folgt $v^\zeta = v^\rho  b$ mit einem $b \in K^*$, so dass
 $\zeta$ die Abbildung $\rho$ induziert. Damit ist 6.2 bewiesen.
 \medskip
       Eine Bemerkung ist angebracht. Hat im vorliegenden Falle $V$
 end\-li\-chen Rang, so ist $\Rg_K(V) = 2 \cdot \Rg_F(V)$, so dass
 $\Rg_K(V)$ gerade ist.
 \medskip\noindent
 {\bf 6.3. Satz.} {\it Es sei $K$ ein K\"orper und $\alpha$ sei ein
 involutorischer Automorphismus von $K$. Ist dann
 $F := \{f \mid f \in K, f^\alpha = f\}$ der Fixk\"orper von $\alpha$, so gilt:
 \smallskip
 \item{a)} Ist $\Char (K) \neq 2$, so gibt es ein $j \in K$ mit
 $j^\alpha = -j \neq 0$. Ist $j^\alpha = -j \neq 0$, so ist $K = F \oplus jF$.
 \smallskip
 \item{b)} Ist $\Char (K) = 2$, so gibt es ein $j \in K$ mit
 $j^\alpha = j+1$. Ist $j^\alpha = j+1$, so ist $K = f \oplus jF$.
 \smallskip\noindent
 In beiden F\"allen ist $[K:F]=2$}.
 \smallskip
       Beweis. a) Weil $\alpha \neq 1$ ist, gibt es ein $x \in K$ mit
 $x^\alpha \neq x$. Setze $j := x - x^\alpha$. Dann ist $j^\alpha = -j \neq 0$.
 \par
       Es sei $j^\alpha = -j \neq 0$. Dann ist $F \cap jF=\{0\}$. Ist $k \in K$,
 so ist
 $$ k = {1 \over 2} (k + k^\alpha) + {1 \over 2}(k - k^\alpha). $$
 Ferner ist ${1 \over 2}(k + k^\alpha) \in F$. Weil
 $j$ von $0$ verschieden ist, gibt es ein $r \in K$ mit
 ${1 \over 2}(k - k^\alpha) = jr$. Dann ist
 $$ -jr = {1 \over 2}(k - k^\alpha)
	= \biggl({1 \over 2}\bigl(k - k^\alpha\bigr)\biggr)^\alpha
        =(jr)^\alpha = -jr^\alpha, $$
 so dass $r^\alpha = r$ ist. Folglich ist $r \in F$,
 so dass $K = F + jF$ ist. Also ist $K = F \oplus jF$.
 \par
       b) Es sei $x^\alpha \neq x$ und $r(x^\alpha + x) = 1$ Wegen
 $x^\alpha + x \in F$ ist dann auch $r \in F$. Setze $j := rx$. Dann ist
 $$ j^\alpha + j = (rx)^\alpha + rx = r(x^\alpha + x) = 1 $$
 und folglich, da die Charakteristik ja $2$ ist, $j^\alpha = j + 1$.
 \par
       Es sei nun $j^\alpha = j + 1$. Wegen $j \neq j + 1$ ist $j^\alpha \neq j$
 und daher $F \cap jF = \{0\}$. Ist $k \in K$, so ist
 $k + k^\alpha \in F$. Ferner ist
 $$ \bigl(k + j(k + k^\alpha)\bigr)^\alpha = k^\alpha + (j + 1)(k + k^\alpha)
					   = k + j(k + k^\alpha). $$
 Also ist auch $k + j(k + k^\alpha) \in F$. Schlie\ss lich ist
 $$ k = k + j(k + k^\alpha) + j(k + k^\alpha) \in F + jF. $$
 Damit ist alles bewiesen.
 \medskip
       Als N\"achstes beweisen wir einen Spezialfall von Hilberts Satz
 90.\index{Hilberts Satz 90}{}
 \medskip\noindent
 {\bf 6.4. Satz.} {\it Es sei $K$ ein K\"orper und $\alpha$ sei ein
 involutorischer Automorphismus von $K$. Ist $a \in K$, so
 sind die folgenden Aussagen \"aquivalent:
 \item{a)} Es ist $aa^\alpha = 1$.
 \item{b)} Es gibt ein $b \in K$ mit $a = b^\alpha b^{-1}$.
 \item{c)} Es gibt ein $c \in K$ mit $a = cc^{-\alpha}$.}
 \smallskip
 Beweis. a) impliziert b). Es sei $aa^\alpha = 1$. W\"are
 $x + a^{-1}x^\alpha = 0$ f\"ur alle $x \in K$, so folgte mit $x = 1$, dass
 $a^{-1} = -1$ w\"are. Es folgte der Widerspruch $x^\alpha = x$ f\"ur alle
 $x \in K$. Es gibt folglich ein $x \in K$, so dass $x + a^{-1} x^\alpha \neq 0$
 ist. Setze $b := x + a^{-1} x^\alpha$. Dann ist
 $$b^\alpha b^{-1} = (x^\alpha + a^{-\alpha} x)b^{-1}
		   = (x^\alpha + ax)b^{-1} = abb^{-1}=a.$$
 \par
       b) impliziert c). Es sei $a = b^\alpha b^{-1}$. Setze
 $c := b^\alpha$. Dann ist $c \neq 0$ und $a = cc^{-\alpha}$.
 \par
       c) impliziert a). Ist n\"amlich $a = cc^{-\alpha}$, so ist
 $aa^\alpha = cc^{-\alpha}c^\alpha c^{-1} =1$. Damit ist alles bewiesen.
 \medskip\noindent
 {\bf 6.5. Satz.} {\it Es sei $V$ ein Vektorraum \"uber dem
 K\"orper $K$ mit $\Rg_K(V) \geq 3$. Ferner sei $\sigma$
 eine Involution aus $\PGaL(V,K)$, die einen Fixpunkt besitzt.

 Wird
 $\sigma$ nicht durch eine Involution aus $\GaL(V,K)$ induziert, so
 gibt es einen inneren Automorphismus $\alpha$ der Ordnung $2$ von
 $K$ und eine bez\"uglich $\alpha$ semilineare Abbildung $\eta$ von
 $V$ in sich mit $\eta^2 = 1$, welche $\sigma$ induziert.

 Ist $F$
 der Fixk\"orper von $\alpha$ und ist $W$ die Menge der Fixvektoren
 von $\eta$, so ist $V$ auch ein $F$-Vektorraum und $W$ ist ein
 Teilraum des $F$-Vektorraumes $V$. Ferner gilt $\Rg_V (W) = \Rg_K(V)$.
 \par 
       Ist $C$ der Zentralisator von $\sigma$ in $\PGaL(V,K)$, so gilt:
 \item{a)} Ist $\Char (K)=2$, so ist
 $$ C/Z(C) \cong \PGaL(W,F)\ \ \ und\ \ \ Z(C) \cong Z(F^*)/Z(K^*). $$
 \item{b)} Ist $\Char (K) \neq 2$, so enth\"alt $C$ einen Normalteiler $N$ vom Index
 $2$ und es ist
 $$ N/Z(N) \cong \PGaL(W,F)\ \ \ sowie\ \ \ Z(N) \cong Z(F^*)/Z(K^*). $$
 Ferner ist $Z(C) = \{1,\sigma\}$.\par}
 \smallskip
       Beweis. Es sei $\rho \in \GaL(V,K)$ und $\rho$ induziere $\sigma$.
 Ferner sei $P = pK$ ein nach Voraussetzung existierender Fixpunkt
 von $\sigma$. Es gibt dann ein $k \in K^*$ mit $p^\rho = pk$ und ein
 $z\in Z(K^*)$ mit $v^{\rho^2} = vz$ f\"ur alle $v \in V$. Hieraus
 folgt, dass $k^2 = z\in Z(K)$ ist. Definiert man $\alpha$ durch
 $x^\alpha := kxk^{-1}$ f\"ur alle $x \in K$, so ist $\alpha$ ein
 innerer Automorphismus von $K$ mit $\alpha^2=1$.
 \par
       Es sei $\eta$ die durch $v^\eta := v^\rho k^{-1}$ definierte
 semilineare Abbildung von $V$ in sich. Dann wird $\sigma$ auch von
 $\eta$ induziert. Ferner ist $v^{\eta^2} = v^{\rho^2}k^{-2} = v$, so
 dass $\eta^2 = 1$ ist. Weil $\sigma$ von $\eta$ induziert wird,
 folgt auf Grund unserer Annahme, dass $\eta$ nicht linear ist.
 Folglich ist $k\not \in Z(K)$ und daher $\alpha \neq 1$. Ist $x \in K$, so ist
 $$ (vx)^\eta = (vx)^\rho k^{-1} = v^\rho xk^{-1} = v^\eta x^\alpha, $$
 so dass $\alpha$ der begleitende Automorphismus von $\eta$ ist.
 \par 
       Die Existenz von $\eta$ und $\alpha$ etabliert, sei $F$ der
 Fixk\"orper von $\alpha$. Dann ist $V$ auch ein Vektorraum \"uber
 $F$ und es gilt $\eta \in \GaL(V,F)$. Wie im Satz notiert, sei $W$
 die Menge der Fixvektoren von $\eta$. Dann ist $W$ ein Teilraum
 des $F$-Vektorraumes $V$.
 \par
       1) Ist $M$ eine Menge von linear unabh\"angigen Vektoren des
 $F$-Vektorraumes $W$, so ist $M$ eine linear unabh\"angige
 Teilmenge des $K$-Vektorraumes $V$.
 \par
       Es seien $m_1$, \dots, $m_n \in M$ und $x_1$, \dots, $x_n \in K$ und
 es gelte $\sum^n_{i:=1} m_ix_i = 0$. Dann ist auch
 $\sum_{i:=1}^n m_ix_iy=0$ f\"ur alle $y \in K$. Hieraus folgt
 $$ 0 = \sum_{i:=1}^n m_i^\eta (x_iy)^\alpha = \sum_{i:=1}^n m_i(x_ia)^\alpha. $$
 Also ist auch
 $$ 0 = \sum_{i:=1}^n m_i\bigl(x_iy + (x_iy)^\alpha\bigr). $$
 Weil $M$ aus $F$-linear unabh\"angigen Vektoren besteht, folgt
 $$ 0 = x_iy + (x_iy)^\alpha $$
 f\"ur alle $i$ und alle $y$. W\"are $x_j \neq 0$, so folgte mit $y := x_j^{-1}$,
 dass $0 = 1 + 1$ w\"are. Mit
 $y:=x_j^{-1}z$ folgte daher $0 = z + z^\alpha$ und dann $z^\alpha = z$.
 Dies widerspr\"ache aber der Tatsache, dass $\alpha$ nicht die
 Identit\"at ist. Damit ist 1) bewiesen.
 \par 
       2) Ist $B$ eine $F$-Basis von $W$, so ist $B$ eine $K$-Basis von $V$.
 \par
       Wegen 1) gen\"ugt es zu zeigen, dass $B$ ein $K$-Erzeugendensystem
 von $V$ ist.
 \par
       Es sei zun\"achst die Charakteristik von $K$ von $2$ verschieden.
 Nach 6.3 gibt es dann ein $j \in K$ mit $j^\alpha = -j$. Es folgt
 $(j^2)^\alpha = j^2$, so dass $j^2 \in F$ ist. Ferner ist $V$ als abelsche
 Gruppe die direkte Summe von $W$ und $W' := \{v \mid v \in V, v^\eta = -v\}$.
 Offensichtlich ist $Wj \subseteq W'$ und $W'j \subseteq W$, so dass
 $Wj^2 \subseteq W$ ist. Weil $W$ ein $F$-Vektorraum ist und $j^2$ in $F^*$
 liegt, ist daher $Wj^2 = W$. Hieraus folgt
 $$ W = Wj^2 \subseteq W'j \subseteq W, $$
 so dass $W'j = W$ ist. Hieraus folgt $V = W + Wj^{-1}$, so dass $B$ in der
 Tat ein $K$-Erzeugendensystem von $V$ ist.
 \par
       Ist die Charakteristik von $K$ gleich $2$, so gibt es nach 6.3 ein
 $j \in K$ mit $j^\alpha = j + 1$. Ist $v \in V$, so ist $v + v^\eta \in W$ und
 $$ \bigl(v + (v + v^\eta)j\bigr)^\eta = v^\eta + (v + v^\eta)(j + 1)
			               = v + (v + v^\eta)j, $$
 so dass $v + (v + v^\eta)j \in W$ gilt. Aus $v = v + (v^\eta)j + (v + v^\eta)j$
 folgt daher, dass $V = W + Wj$ ist. Also ist $B$ auch in diesem Falle
 ein $K$-Erzeugendensystem von $V$.
 \par
       3) Ist $P$ ein Fixpunkt von $\sigma$, so ist $P \cap W$ ein
 $F$-Unterraum des Ranges $1$ von $W$.
 \par
       Es sei $P = pK$. Dann ist $p^\eta = pa$ mit einem $a \in K^*$.
 Ferner ist $p = p^{\eta^2} = paa^\alpha$, so dass $aa^\alpha = 1$.
 Nach 6.4 gibt es daher ein $c \in K^*$ mit $a = cc^{-\alpha}$. Setze
 $q := pc$. Dann ist $0 \neq q \in P$. Ferner ist
 $$q^\eta = (pc)^\eta = p^\eta c^\alpha = pac^\alpha = pc = q. $$
 Also ist $P \cap W \neq \{0\}$. Mit 1) folgt
 $$ \Rg_F(P\cap W) \leq \Rg_K(V), $$
 so dass, wie behauptet, $\Rg_F(P\cap W) = 1$ ist.
 \par
       4) Es sei $L$ die Menge der Unterr\"aume von $L_K(V)$, die von
 Fixpunkten von $\sigma$ aufgespannt werden. Dann erbt $L$ von
 $L_K(V)$ die Teilordnung. Definiert man die Abbildung $\varphi$
 von $L$ in $L_F(W)$ durch $X^\varphi := X \cap W$ f\"ur alle $X \in L$,
 so ist $\varphi$ ein Isomorphismus von $L$ auf $L_F(W)$.
 Insbesondere ist $L$ also ein projektiver Verband.
 \par
       Es seien $X$, $Y \in L$. Ist $X \leq Y$, so ist
 $$ X^\varphi = X \cap W \leq Y \cap W = Y^\varphi. $$
 Es sei umgekehrt $X^\varphi \leq Y^\varphi$. Die Definition von $L$ liefert
 zusammen mit 3), dass
 $$ X = \sum_{x \in X^\varphi} xK \leq \sum_{y \in Y^\varphi} yK=Y $$
 ist. Damit ist $\varphi$ als Monomorphismus von $L$ in $L_F(W)$ erkannt.
 \par 
       Es sei $Y \in L_F(W)$. Setze $X := \sum_{y \in Y} yK$. Wegen
 $(yK)^\eta = yK$ ist $X \in L$. \"Uberdies ist $X^\varphi = X \cap W > Y$. Es
 sei $x \in X \cap W$. Es gibt dann $y_1$, \dots, $y_n \in Y$
 und $k_1$, \dots, $k_n \in K$ mit $x = \sum_{i:=1}^n y_ik_i$. Wir
 d\"urfen annehmen, dass die $y_i$ als Vektoren des
 $K$-Vektorraumes $V$ linear unabh\"angig sind. Dann sind sie aber
 auch $F$-linear unabh\"angig. Somit lassen sie sich zu einer Basis
 von $W$ erg\"anzen. Weil nach 2) jede $F$-Basis von $W$ eine
 $K$-Basis von $V$ ist, folgt aus $x \in W$, dass $k_1$, \dots, $k_n \in F$ ist.
 Also ist $x \in Y$, so dass $X^\varphi = Y$ ist.
 Folglich ist $\varphi$ ein Isomorphismus von $L$ auf $L_F(W)$.
 \par
       Es sei $\gamma \in \GaL(V,K)$ und es gelte $v^{\gamma \eta} =
 v^{\eta \gamma} e$ f\"ur alle $v \in V$. Nach 6.1 ist dann
 entweder $e = 1$ oder $e = -1$. Weil die von $\gamma$ in $L_K(V)$
 induzierte Kollineation im Zentralisator $C$ von $\sigma$ liegt,
 ist $L^\gamma = L$. Daher ist $\varphi^{-1}\gamma\varphi$ eine
 Kollineation von $L_F(W)$. Wegen 2) ist $\Rg_F(W) \geq 3$. Nach
 dem zweiten Struktursatz II.7.1 gibt es daher eine semilineare
 Abbildung $(\delta, \beta)$ des $F$-Vektorraumes $W$ in sich,
 welche in $L_F(W)$ die Kollineation $\varphi^{-1} \gamma\varphi$
 induziert. Hieraus folgt zusammen mit $\Rg_F(W) \geq 3$
 die Existenz eines Elementes $a \in K^*$ mit $w^\delta = w^\gamma a$
 f\"ur alle $w \in W$. Es ist
 $$ w^{\gamma\eta} a^\alpha = (w^\gamma a)^\eta = w^{\delta\eta} = w^\delta
                            = w^{\eta\delta} = w^{\eta\gamma}a
			    = w^{\gamma \eta} ea. $$
 Also ist entweder $a^\alpha = a$ oder $a^\alpha = -a$. Ferner ist, falls
 $f \in F$ ist,
 $$ w^\delta f^\beta = (wf)^\delta = (wf)^\gamma a = w^\gamma fa
		     = x^\delta a^{-1} fa. $$
 Folglich ist $f^\beta = a^{-1}fa$ f\"ur alle $f \in F$. Ist $a^\alpha = a$, so
 ist $a \in F$, so dass wir auf Grund von 1.3 annehmen d\"urfen, dass
 $\beta = 1$ ist. Ist $a^\alpha = -a$ und ist $j \in K$ ein Element,
 f\"ur welches $j^\alpha = -j$ gilt, so ist $a = jf_0$ mit $f_0 \in
 F$ nach 6.3). Dann ist $f^\beta = f_0^{-1} j^{-1} fjf_0$, so dass
 wir wiederum nach 1.3 annehmen d\"urfen, dass $a = j$ ist. (Der
 Leser fragt sich vielleicht, warum so umst\"andlich, wo doch
 $a^\alpha = -a$ gilt. Die Antwort lautet: Man kann $j$ {\it a priori\/}
 w\"ahlen, und darauf kommt es an.) Es ist
 $-j = j^\alpha = kjk^{-1}$. Hieraus folgt, dass
 $$ k \neq -k = j^{-1} kj = k^\beta $$
 ist. Nun ist $F$ aber der Fixk\"orper von $\alpha$, dh., der Zentralisator von
 $k$ in $K$. Daher liegt $k$ im Zentrum von $F$. Somit gilt nach dem
 Vorstehenden, dass $\beta$ kein innerer Automorphismus von $F$ ist. Es sei $H$
 das Urbild von $C$ in $\GaL(V,K)$. Ist dann $\gamma \in H$, so ist
 $v^{\gamma\eta} = v^{\eta \gamma}$ oder
 $v^{\gamma\eta} = -v^{\eta \gamma}$, wie wir schon bemerkten.
 \par
       Es sei nun $\Char(K) = 2$. Dann ist $H = C_{\PGaL(V,K)}(\eta)$.
 Daher ist $W^H = W$. Ist $\epsilon(\gamma)$ die Einschr\"ankung von
 $\gamma \in H$ auf $W$, so ist $\epsilon$ ein Homomorphismus von
 $H$ in $\GaL(W,F)$. Weil jede $F$-Basis von $W$ gleichzeitig eine $K$-Basis von
 $V$ ist, ist $\epsilon$ injektiv. Es sei $\gamma \in \GaL(W,F)$. Ist $B$ eine
 $F$-Basis von $W$ und definiert man die
 lineare Abbildung $\gamma$ von
 $V$ in sich durch $b^\gamma := b^\lambda$ f\"ur alle $b \in B$, so
 ist $\gamma$ bijektiv und es gilt $\epsilon(\gamma) = \lambda$.
 Folglich ist $\epsilon$ ein Isomorphismus von $H$ auf $\GaL(W,F)$.
 Hieraus folgen nun die Behauptungen unter a).
 \par
       Es sei schlie\ss lich $\Char(K) \neq 2$. Ferner sei $B$ wieder eine
 $F$-Basis von $W$ und $j$ sei ein Element von $K$ mit $j^\alpha = -j$. Ist nun
 $v \in V$ und $v = \sum_{b\in B} bx_b$ und definiert
 man $\gamma$ durch $v^\gamma := \sum_{b\in B} bjx_b$, so ist
 $\gamma \in \GaL(V,K)$. \"Uberdies ist
 $$ v^{\gamma \eta} = \sum_{b \in B} b^\eta j^\alpha x^\alpha_b
		    = -\sum_{b \in B} b_jx^\alpha_b
		    = -(\sum_{b \in B} bx_b^\alpha)^\gamma
		    = -v^{\eta\gamma}. $$
 Daher ist $\gamma \in H$. Hieraus folgt, dass $H_0 := C_{\GaL(V,K)}(\eta)$ den
 Index 2 in $H$ hat. Weil $W \cap W^\gamma = \{0\}$ ist, folgt nach unserer
 Vorbemerkung, dass es
 eine bez\"uglich $\beta$ semilineare Abbildung $\delta$ von $W$ in
 sich gibt mit $w^\delta = w^\gamma j$ f\"ur alle $w \in W$. Weil
 $\beta$ kein innerer Automorphismus ist, wird $H_0$ bei dem
 Homomorphismus von $H$ auf $C$ auf einen Normalteiler $N$ vom
 Index 2 in $C$ abgebildet. Ferner folgt ohne M\"uhe, wie schon im
 Falle der Charakteristik 2, dass $N/Z(N) \cong \PGaL(W,F)$ und
 $Z(N) \cong Z(F^*)/Z(K^*)$ ist. Es bleibt zu zeigen, dass $Z(C)$
 nur aus den beiden Elementen 1 und $\sigma$ besteht.
 \par
       Offensichtlich ist $Z(C) \subseteq Z(N)$. Es sei $1 \neq \rho \in Z(C)$
 und $\rho$ werde von $\zeta$ induziert. Es gibt dann ein $f \in Z(F^*)$ mit
 $w^\zeta = wf$ f\"ur alle $w \in W$. Weil $\beta$
 in $Z(F)$ einen Automorphismus der Ordnung 2 induziert und da
 $k^\beta = -k$ ist, gibt es nach 6.3 zwei Elemente $g$, $h \in
 Z(F)$ mit $f = g + kh$ und $g^\beta = g$ sowie $h^\beta = h$. Wegen
 $x^\beta = j^{-1} xj$ hei\ss t das, dass $g$ und $h$ mit $j$
 vertauschbar sind. Weil $g$ und $h$ auch im Zentrum von $F$ liegen
 und nach 6.3 die Gleichung $K = F + jF$ gilt, liegen $g$ und $h$ im
 Zentrum von $K$. Aus $\rho \neq 1$ folgt, dass $h \neq 0$ ist. Wir
 d\"urfen daher annehmen, dass $h = 1$ ist. Wegen $\rho \in Z(C)$ ist
 $v^{\zeta \gamma} = v^{\gamma \zeta} r$ mit einem passenden $r \in Z(K^*)$. Ist
 $b \in B$, so ist also
 $$ bjf = b^\gamma f = (bf)^\gamma = b^{\zeta \gamma} = b^{\gamma \zeta}r
		 = b^\zeta jr = bfjr. $$
 Somit ist, da ja $f = g + k$ ist,
 $$ j(g + k) = jf = fjr = gjr + kjr = jr(g - k). $$
 Weil $g$ in $F$ und $k$ nicht in $F$ liegt, folgt $g = 0$ und $r = -1$. Somit
 ist
 $$v^\zeta = \sum_{b\in B} b^\zeta x_b = \sum_{b\in B} bkx_b
           = \sum_{b \in B} b^\eta x_b^\alpha k = v^\eta k, $$
 dh., es ist $\rho = \sigma$. Damit ist dann alles bewiesen.
 \medskip
       Wir haben nun gute Kenntnisse \"uber die Zentralisatoren von
 Involutionen, die keine Fixpunkte haben, und von Involutionen, die
 zwar Fixpunkte haben, die aber nicht durch involutorische lineare
 Abbildungen induziert werden. F\"ur die noch verbleibenden
 projektiven Involutionen, die von involutorischen linearen
 Abbildungen induziert werden, ist die Arbeit aber auch schon
 erledigt, wie der n\"achste Satz zeigt.
 \medskip\noindent
 {\bf 6.6. Satz.} {\it Es sei $V$ ein Vektorraum \"uber dem
 K\"orper $K$ mit $\Rg_K(V) \geq 3$. Ist $\sigma$ eine
 Kollineation von $L_K(V)$, die von einer involutorischen
 linearen Abbildung von $V$ in sich induziert wird, so gilt:
 \item{a)} Ist $\Char (K) = 2$, so ist $\sigma$ eine Quasielation.
 \item{b)} Ist $\Char (K) \neq 2$, so ist $\sigma$ eine
 Quasistreckung.}
 \smallskip
       Der Beweis sei dem Leser als \"Ubungsaufgabe \"uberlassen.

\mysection{7. Die hessesche Gruppe}

\noindent
 Schlitzt man die projektive Ebene \"uber $\GF(3)$ bez\"uglich einer
 Ge\-ra\-den, so erh\"alt man die
 {\it affine Ebene\/}\index{affine Ebene der Ordnung 3}{} \"uber $\GF(3)$.
 Diese ist eine Inzidenzstruktur aus 9 Punkten und 12 Geraden. Sie
 findet auch bei Anh\"angern der algebraischen Geometrie Interesse,
 da sie sich n\"amlich in vielen projektiven Ebenen ---
 insbesondere in der \"uber dem K\"orper der komplexen Zahlen ---
 als Konfiguration der Wendepunkte einer irreduziblen Kurve dritten
 Grades wiederfindet. Der Stabilisator dieser Konfiguration in der
 projektiven Gruppe ist die nach Ludwig Otto Hesse\index{Hesse, L. O.}{} benannte
 {\it hessesche Gruppe\/}.\index{hessesche Gruppe}{} Eine Untergruppe dieser
 Gruppe wird im n\"achsten Abschnitt ihren Beitrag zur L\"ange der dort
 vorgef\"uhrten Argumentationskette leisten, so dass in diesen Abschnitt ein
 wenig \"uber diese Situation berichtet werden soll.
 \par
       Hesse publizierte seine Ergebnisse im crelleschen Journal (Hesse 1844, 1847,
 Salmon 1850).
 Sie finden sich wiedergegeben in Jordans \anff Trait\'e`` und in Nettos
 \anff Substitutionen\-the\-o\-rie`` (Jordan, S.~302ff, Netto 1882, S.~232ff).
 Heutige Leser haben es bequemer, wenn sie sich an die Ausarbeitung der
 brieskornschen Vorlesungen des WS 1975/76 und des SS 1976 halten (Brieskorn und
 Kn\"orrer o. J., S. 379--388).
 \par
       Die Punkte der affinen Ebene \"uber $\GF(3)$ kann man
 identifizieren mit den Paaren $(x,y)$ mit $x$, $y \in \GF(3)$. Die
 Geraden werden dann durch inhomogene lineare Gleichungen
 beschrieben. Kodiert man das Paar $(x,y)$ durch die Zahl $3x + y$,
 so werden die Punkte der affinen Ebene \"uber $\GF(3)$ durch die
 Ziffern $0$ bis $8$ beschrieben und die Ge\-ra\-den durch die Spalten
 der folgenden Matrix.
 $$ \matrix{
      0\ 1\ 2\ \ &0\ 1\ 2\ \ &0\ 1\ 2\ \ &0\ 3\ 6\cr
      3\ 4\ 5\ \ &4\ 5\ 3\ \ &5\ 3\ 4\ \ &1\ 4\ 7\cr
      6\ 7\ 8\ \ &8\ 6\ 7\ \ &7\ 8\ 6\ \ &2\ 5\ 8\cr}$$
 \par
       Wir gehen nun der Frage nach, unter welchen Bedingungen sich diese
 affine Ebene in eine gegebene projektive Ebene einbetten l\"asst. Diese Frage
 erscheint nach den vorstehenden Bemerkungen vielleicht nicht ganz so
 k\"unstlich, wie sie zun\"achst klingen mag. Hinzukommt, wie schon gesagt, dass
 sie sich im n\"achsten Abschnitt zwangsl\"aufig stellen wird.
 \medskip\noindent
 {\bf 7.1. Satz.} {\it Es sei $V$ ein Vektorraum des Ranges $3$
 \"uber dem K\"orper $K$ und $b_0$, $b_1$, $b_2$ sei eine Basis von
 $V$. Ferner sei $t$ ein von $0$ und $1$ verschiedenes Element aus
 $K$. Wir definieren eine Abbildung $P$ von $\{0,1,2,3,4,5,6,7,8\}$
 in die Punktmenge von $L_K(V)$ durch
 \smallskip\noindent\hskip 25 mm
       $P_0 := b_0K$, \smallskip\noindent\hskip 25 mm
       $P_3 := b_1K$, \smallskip\noindent\hskip 25 mm
       $P_4 := b_2K$, \smallskip\noindent\hskip 25 mm
       $P_7 := (b_0 + b_1 + b_2)K$,   \smallskip\noindent\hskip 25 mm
       $P_6 := (b_0 + b_1t)K$,        \smallskip\noindent\hskip 25 mm
       $P_2 := (b_0 + b_1t + b_2)K$,  \smallskip\noindent\hskip 25 mm
       $P_1 := (b_0t + b_1t + b_2)K$, \smallskip\noindent\hskip 25 mm
       $P_8 := (b_0t + b_2)K$,        \smallskip\noindent\hskip 25 mm
       $P_5 := (b_0(t^2 - t + 1) + b_1t + b_2t)K$. \smallskip\noindent
 Genau dann ist $P$ ein Monomorphismus der affinen Ebene der
 Ordnung $3$ in $L_K(V)$, wenn
 $$t^2 - t + 1 = 0$$
 ist.}
 \smallskip
       Beweis. Weil $t$ von $0$ und $1$ verschieden ist, ist $P$ injektiv.
       Es ist $P_0 = b_0K, P_3 = b_1K$ und $P_6 = (b_0 + b_1t)K$, so dass diese
 drei Punkte kollinear sind.
 \par
       Es ist $P_1 = (b_0t + b_1t + b_2)K, P_4 = b_2K$ und
 $P_7 = (b_0 + b_1 + b_2)K$, so dass diese Punkte wegen $t \neq 1$ kollinear
 sind.
 \par 
       Es ist $P_2 = (b_0 + b_1t + b_2)K, P_5 = (b_0(t^2 - t + 1)
 + b_1t + b_2t)K$ und $P_8 = (b_0t + b_2)K$. Wegen
 $$ (b_0 + b_1t + b_2) + (b_0t + b_2)(t - 1) = b_0(t^2 - t + 1) + b_1t + b_2t $$
 sind auch diese drei Punkte kollinear.
 \par
       Es ist $P_0 = b_0K$, $P_4 = b_2K$ und $P_8=(b_0t + b_2)K$, so dass diese
 Punkte kollinear sind.
 \par
       Es ist $P_1 = (b_0t + b_1t + b_2)K$, $P_5 =
 (b_0(t^2 - t + 1) + b_1t + b_2t)K$ und $P_6 = (b_0 + b_1t)K$. Wegen
 $$ (b_0t + b_1t + b_2)t - (b_0 + b_1t)(t - 1) = b_0(t^2 - t + 1) + b_1t + b_2t$$
 sind diese drei Punkte ebenfalls kollinear.
 \par
       Es ist $P_2 = (b_0 + b_1t + b_2)K$, $P_3 = b_1K$ und
 $P_7 = (b_0 + b_1 + b_2)K$,
 so dass diese drei Punkte wegen $t \neq 1$ kollinear sind.
 \par
       Es ist $P_0 = b_0K$, $P_5 = (b_0(t^2 - t + 1) + b_1t + b_2t)K$ und
 $P_7 = (b_0 + b_1 + b_2)K$, so dass diese Punkte wegen $t \neq 1$ kollinear
 sind.
 \par
       Es ist $P_1 = (b_0t + b_1t + b_2)K$, $P_3 = b_1K$ und
 $P_8 = (b_0t + b_2)K$, so dass auch diese drei Punkte kollinear sind.
 \par
       Es ist $P_2 = (b_0 + b_1t + b_2)K$, $P_4 = b_2K$ und
 $P_6 = (b_0 + b_1t)K$.  Somit sind auch diese drei Punkte kollinear.
 \par
       Bislang ist also noch nichts weiter geschehen. Doch nun finden wir
 die Bedingung, die $t$ erf\"ullen muss, damit $P$ ein
 Monomorphismus ist.
 \par
       Es ist $P_3 = b_1K$, $P_4 = b_2K$ und
 $P_5 = (b_0(t^2 - t + 1) + b_1t + b_2t)K$.
 Damit diese drei Punkte kollinear sind, ist notwendig und
 hinreichend, dass es $\lambda$ und $\mu$ in $K$ gibt mit
 $$ b_0(t^2 - t + 1) + b_1t + b_2t = b_1\lambda + b_2\mu. $$
 Somit sind diese drei Punkte genau dann kollinear, wenn $t^2 - t + 1 = 0$ ist.
 Diese Bedingung ist also notwendig daf\"ur, dass $P$ ein Monomorphismus
 ist. Sie ist aber auch hinreichend. Um dies einzusehen, sei
 $t^2 - t + 1 = 0$. Wir m\"ussen dann nur noch zeigen, dass auch die
 Tripel $P_0$, $P_1$, $P_2$ und $P_6$, $P_7$, $P_8$ kollinear sind.
 \par
       Es ist $P_0 = b_0K$, $P_1 = (b_0t + b_1t + b_2)K$ und
 $P_2 = (b_0 + b_1t + b_2)K$.  Wegen $t^2 - t + 1 = 0$ ist
 $$ b_0 = -(b_0t + b_1t + b_2)t + (b_0 + b_1t + b_2)t, $$
 so dass auch diese Punkte kollinear sind.
 \par
       Schlie\ss lich ist $P_6 = (b_0 + b_1t)K$, $P_7 = (b_0 + b_1 + b_2)K$ und
 $P_8 = (b_0t + b_2)K$. Nun ist aber
 $$ b_0 + b_1t = (b_0 + b_1 + b_2)t - (b_0t + b_2)t, $$
 so dass auch dieses letzte Punktetripel kollinear ist. Damit ist der Satz
 bewiesen.
 \medskip
       Die in Satz 7.1 beschriebene Einbettung der affinen Ebene \"uber $\GF(3)$
 in eine desarguessche projektive Ebene ist im Wesentlichen auch die einzige
 M\"oglichkeit, diese affine Ebene in eine des\-ar\-gues\-sche projektive Ebene
 einzubetten. Dies ist die Aussage des n\"achsten Satzes.
 \medskip\noindent
 {\bf 7.2. Satz.} {\it Es sei $V$ ein Vektorraum des Ranges $3$
 \"uber dem K\"orper $K$. Ist dann $P$ ein Monomorphismus der
 affinen Ebene \"uber $\GF(3)$ in $L_K(V)$, so gibt es eine Basis
 $b_0$, $b_1$, $b_2$ von $V$ sowie ein $t \in K$, so dass $P$ mittels
 dieser Daten wie in Satz 7.1 beschrieben dargestellt ist.
 Insbesondere gibt es genau dann einen Monomorphismus der affinen
 Ebene \"uber $\GF(3)$ in $L_K(V)$, wenn es ein $t \in K$ gibt mit
 $t^2 - t + 1 = 0$.}
 \smallskip
       Beweis. Da die Punkte $P_0$, $P_3$, $P_4$, $P_7$ einen Rahmen
 von $L_K(V)$ bilden, gibt es nach 1.8 Vektoren $b_0$, $b_1$,
 $b_2$ mit $P_0 = b_0K$, $P_3 = b_1K$, $P_4 = b_2K$ und
 $P_7 = (b_0 + b_1 + b_2)K$. Wegen $P_6 \leq P_0 + P_3$ und $P_6 \neq P_3$
 gibt es ein $t \in K$ mit $P_6 = (b_0 + b_1t)K$. Nun ist
 $$\eqalign{
      P_2 &= (P_3 + P_7) \cap (P_4 + P_6),      \cr
      P_1 &= (P_4 + P_7) \cap (P_0 + P_2),      \cr
      P_8 &= (P_0 + P_4) \cap (P_1 + P_3),      \cr
      P_5 &= (P_1 + P_6) \cap (P_2+P_8). \cr}$$
 Hieraus folgt alles weitere.
 \medskip\noindent
 {\bf 7.3. Satz.} {\it Es sei $V$ ein Vektorraum des Ranges $3$
 \"uber dem K\"orper $K$. Es sei $b_1$, $b_2$, $b_3$ eine Basis von
 $V$ und $t$ sei ein Element in $K$ mit $t^2 - t + 1 = 0$. Schlie\ss lich sei $P$
 der in 7.1 beschriebene Monomorphismus der affinen
 Ebene \"uber $\GF(3)$ in $L_K(V)$. Mit $C$ bezeichnen wir den
 Zentralisator von $t$ in $K^*$. Ist dann $\sigma \in \GL(V)$, so
 gilt $P^\sigma_i = P_i$ f\"ur alle $i$ genau dann, wenn es ein
 $\lambda \in C$ gibt mit $b^\sigma_j = b_j\lambda$ f\"ur $j := 0$,
 $1$, $2$. Die Gruppe der projektiven Kollineationen, die alle
 Punkte $P_i$ festlassen, ist isomorph zu $C/Z(K^*)$.}
 \smallskip
       Beweis. Ist $\lambda \in C$ und gilt $b_j^\sigma = b_j\lambda$ f\"ur
 $j := 0$, 1, 2, so ist nat\"urlich $P^\sigma_i = P_i$ f\"ur alle $i$. Es
 sei also $P^\sigma_i = P_i$ f\"ur alle $i$. Weil $\sigma$ den Rahmen
 $P_0$, $P_3$, $P_4$, $P_7$ festl\"asst, gibt es ein $\lambda \in
 K^*$ mit $b_j^\sigma = b_j\lambda$ f\"ur $j := 0$, 1, 2. Mit
 $P_6^\sigma = P_6$ folgt die Existenz eines $\mu \in K^*$ mit
 $$ b_0\lambda + b_1\lambda t = (b_0 + b_1t)^\sigma = (b_0 + b_1t)\mu. $$
 Hieraus folgt $\lambda = \mu$ und weiter $\lambda t = t\lambda$, so dass
 $\lambda \in C$ gilt. Damit ist die erste Aussage des Satzes bewiesen. Die
 zweite Aussage folgt unmittelbar aus der ersten.
 \medskip
       Die Kollineationsgruppe der affinen
 Ebene\index{affine Ebene der Ordnung 3}{} \"uber $\GF(3)$ ist
 schnell bestimmt. Nach I.8.8 l\"asst sich jede Kollineation
 dieser Ebene auf genau eine Weise zu einer Kollineation der
 projektiven Ebene \"uber $\GF(3)$ fortsetzen. Daher ist die
 Kollineationsgruppe der affinen Ebene gleich dem Stabilisator
 einer Geraden in der Kollineationsgruppe der projektiven Ebene.
 Mit 1.13 folgt daher, dass diese Gruppe die Ordnung $3^3 \cdot
 2^4$ hat. Sie enth\"alt den Normalteiler der Ordnung 9, der aus
 den Elationen besteht, deren Achse die entfernte Gerade ist. Diese Elationen
 hei\ss en in diesem Zusammenhang {\it Translationen\/}.\index{Translation}{}
 Die Faktorgruppe nach dem Normalteiler der Translationen ist
 isomorph zur $\GL(2, 3)$, wobei diese Gruppe innerhalb der
 Kollineationsgruppe sich als der Stabilisator eines Punktes der
 affinen Ebene wiederfindet. Wir wollen uns nun \"uberlegen, was
 von diesen Kollineationen bei der Einbettung der affinen Ebene
 \"uber $\GF(3)$ in eine desarguessche projektive Ebene im
 Stabilisator der eingebetteten Ebene in der projektiven Gruppe der
 einbettenden Ebene wiederzufinden ist.
 \medskip\noindent
 {\bf 7.4. Satz.} {\it Es sei $K$ ein K\"orper und $t$ 
 ein
 Element aus $K$ mit $t^2 - t + 1 = 0$. Es sei ferner $V$ ein Vektorraum
 vom Rang $3$ \"uber $K$ und $b_0$, $b_1$, $b_2$ sei eine Basis von
 $V$. Schlie\ss lich sei $P$ der in 7.1 mittels dieser Daten
 definiert Monomorphismus der affinen Ebene der Ordnung $3$ in
 $L_K(V)$. Wir definieren die Abbildungen $\rho$, $\sigma$ und
 $\tau$ in $\GL(V)$ durch
 $$\eqalign{
          b^\rho_0 &:= -b_0             \cr
	  b^\rho_1 &:= b_0t^{-1} + b_1  \cr
	  b_2^\rho &:= b_0t + b_2,      \cr
	b_0^\sigma &:= b_0 + b_1t       \cr
	b_1^\sigma &:= -b_1             \cr
        b_2^\sigma &:= b_1 + b_2,       \cr
	  b_0^\tau &:= b_0+b_2t^{-1}    \cr
	  b_1^\tau &:=b_1+b_2           \cr
	  b_2^\tau &:= -b_2. \cr} $$
 Dann sind $\rho$, $\sigma$ und $\tau$ Involutionen, die in
 $L_K(V)$ involutorische Perspektivit\"aten induzieren. Nennt man
 diese wiederum $\rho$, $\sigma$ und $\tau$, so gilt:
 \par
       Das Zentrum von $\rho$ ist $P_0$, die Achse ist
 $$ \bigl\{b_0\lambda_0 + b_1(2t\lambda_0 - t^2\lambda_2) + b_2\lambda_2 \mid
		   \lambda_0, \lambda_2 \in K\bigr\}$$
 und $\rho$ wirkt auf den $P_i$ durch
 $ (P_1P_2)(P_3P_6)(P_4P_8)(P_5P_7)$.
 \par
       Das Zentrum von $\sigma$ ist $P_3$, die Achse ist
 $$ \bigl\{b_0\lambda_0 + b_1\lambda_1 + b_2(2\lambda_1 - t\lambda_2) \mid
	 \lambda_0, \lambda_1 \in K\bigr\} $$
 und $\sigma$ wirkt auf den $P_i$ durch
 $ (P_0P_6)(P_1P_8)(P_2P_7)(P_4P_5)$.
 \par
       Das Zentrum von $\tau$ ist $P_4$,  die Achse ist
 $$ \bigl\{b_0(2t\lambda_2 - t\lambda_1) + b_1\lambda_1 + b_2\lambda_2 \mid
         \lambda_1, \lambda_2 \in K\bigr\}$$
 und $\tau$ wirkt auf den $P_i$ 
 durch
 $ (P_0P_8)(P_1P_7)(P_2P_6)(P_3P_5)$.
 \par
       Die von $\rho$, $\sigma$ und $\tau$ erzeugte Kollineationsgruppe induziert
 auf der affinen Ebene der Ordnung drei eine Gruppe der Ordnung $9 \cdot 2$
 bestehend aus den Translationen und Streckungen dieser Ebene.}
 \smallskip
       Beweis. Dies zu beweisen ist eine banale Rechenaufgabe.
 \medskip\noindent
 {\bf 7.5. Satz.} {\it Es sei $K$ ein K\"orper und $t$ sei ein
 Element aus $K$ mit $t^2 - t + 1 = 0$. Es sei ferner $V$ ein Vektorraum
 vom Rang $3$ \"uber $K$ und $b_0$, $b_1$, $b_2$ sei eine Basis von
 $V$. Schlie\ss lich sei $P$ der in 7.1 mittels dieser Daten definierte
 Monomorphismus der affinen Ebene der Ordnung $3$ in $L_K(V)$. Wir definieren
 die Abbildungen $\rho$ und $\sigma$ in $\GL(V)$ durch
 $$\eqalign{
        b_0^\rho &:= b_0(t-1)         \cr
	b_1^\rho &:= b_1(t-1)         \cr
	b_2^\rho &:= b_0 + b_1 + b_2, \cr
	b_0^\sigma &:= b_0(t - 1)     \cr
	b_1^\sigma &:= -(b_0 + b_1 + b_2)t \cr
	b_2^\sigma &:=b_2(t-1). \cr} $$
 Dann sind $\rho$ und $\sigma$ Elemente der Ordnung $3$. Die von
 ihnen in $L_K(V)$ induzierten Kollineation lassen die durch $P$
 beschriebene affine Ebene invariant. Bezeichnet man diese
 Kollineationen ebenfalls mit $\rho$ und $\sigma$, so gilt:
 \par
       Die Kollineation $\rho$ hat die Fixpunkte $P_0$, $P_3$ und $P_6$ und
 die Wir\-kung von $\rho$ auf den restlichen Punkten wird durch
 $$ (P_1 P_4 P_7)(P_2P_8P_5) $$
 beschrieben.
 \par
       Die von diesen beiden Kollineationen erzeugte Gruppe induziert auf
 der affinen Ebene der Ordnung $3$ eine zur $\SL(2,3)$ isomorphe
 Gruppe.}
 \smallskip
       Beweis. Auch dies ist banal zu verifizieren.
 \medskip
       Die beiden S\"atze 7.4 und 7.5 implizieren, dass der Stabilisator
 der Menge der Punkte der $P_i$ auf dieser Menge zweifach transitiv operiert.
 \par 
       Geraden einer affinen Ebene, deren Schnittpunkt in der
 zu\-ge\-h\"o\-ri\-gen projektiven Ebene auf der entfernten Gerade liegt,
 hei\ss en {\it parallel\/}.\index{parallel}{} Die Parallelit\"atsrelation ist
 eine \"Aquivalenzrelation. Ihre \"Aqui\-va\-lenz\-klas\-sen hei\ss en
 {\it Parallelenscharen\/}.\index{Parallelenschar}{}
 \medskip\noindent
 {\bf 7.6. Satz.} {\it Es sei $K$ ein K\"orper und $t$ sei ein
 Element aus $K$ mit $t^2 - t + 1 = 0$. Es sei ferner $V$ ein Vektorraum
 vom Rang $3$ \"uber $K$ und $b_0$, $b_1$, $b_2$ sei eine Basis von
 $V$. Schlie\ss lich sei $P$ der in 7.1 mittels dieser Daten
 definierte Monomorphismus der affinen Ebene der Ordnung $3$ in
 $L_K(V)$. Das Bild dieser affinen Ebene unter $P$ werde mit
 $\alpha$ bezeichnet. Dann sind die folgenden Aussagen
 \"aquivalent:
 \item{a)} Es ist $t = -1$.
 \item{b)} Die Charakteristik von $K$ ist gleich $3$.
 \item{c)} Die Geraden einer jeden Parallelschar von $\alpha$ sind
 als Geraden von $L_K(V)$ konfluent.
 \item{d)} Es gibt eine Parallelenschar in $\alpha$, deren Geraden
 als Geraden von $L_K(V)$ konfluent sind.\par}
 \smallskip
       Beweis. a) und b) sind \"aquivalent: Ist $t = -1$, so ist
 $0 = t^2 - t + 1 = 3$, so dass b) eine Folge von a) ist. Ist die
 Charakteristik von $K$ gleich 3, so ist $0 = t^2 - t + 1 = (t + 1)^2$, da
 ja $-1 = 2$ ist. Folglich ist $t = -1$.
 \par
       b) impliziert c): Da die Charakteristik von $K$ gleich 3 ist,
 enth\"alt $L_K(V)$ auch eine Kopie der projektiven Ebene \"uber
 $\GF(3)$. Da $-1$ die einzige L\"osung der Gleichung $t^2 - t + 1 = 0$ ist,
 folgt mit 7.2 die G\"ultigkeit von c).
 \par
       c) ist zu d) \"aquivalent: Aus c) folgt nat\"urlich d). Weil der
 Sta\-bi\-li\-sa\-tor von $\alpha$ in der projektiven Gruppe von $L_K(V)$
 auf der Punktmenge von $\alpha$ zweifach transitiv operiert, ist
 c) auch eine Konsequenz von d).
 \par
       c) impliziert b): Es ist $(b_0 + b_1)K = (P_0 + P_3)\cap (P_4 + P_7)$.
 Nach Voraussetzung ist $(b_0 + b_1)K \leq P_2 + P_5$. Es gibt also $r$,
 $s \in K$ mit $b_0 + b_1 = (b_0 + b_1t + b_2)r + (b_1 + b_2)s$. Hieraus folgen
 die Gleichungen $r = 1$, $tr - s = 1$ und $s + r = 0$. Hieraus folgt $t = 2$ und
 weiter $0 = t^2 - t + 1 = 3$. Damit ist alles bewiesen.
 \medskip\noindent
 {\bf 7.7. Satz.} {\it Es sei $K$ ein K\"orper und $t$ sei ein
 Element aus $K$ mit $t^2 - t + 1 = 0$, sei ferner $V$ ein Vektorraum
 vom Rang $3$ \"uber $K$ und $b_0$, $b_1$, $b_2$ sei eine Basis von
 $V$. Schlie\ss lich sei $P$ der in 7.1 mittels dieser Daten
 definierte Monomorphismus der affinen Ebene der Ordnung $3$ in
 $L_K(V)$. Das Bild dieser affinen Ebene unter $P$ werde mit
 $\alpha$ bezeichnet und $G$ sei der Stabilisator von $\alpha$ in
 $\PGL(V)$. 
 \item{a)} Ist die Charakteristik von $K$ gleich $3$, so induziert
 $G$ die volle Kollineationsgruppe in der affinen Ebene $\alpha$.
 \item{b)} Ist die Charakteristik von $K$ von $3$ verschieden, so
 induziert $G$ in $\alpha$ die Gruppe $T \cdot SL(2,3)$, wobei $T$
 die Gruppe der Translationen von $\alpha$ bezeichnet. Diese Gruppe
 ist vom Index $2$ in der vollen Kollineationsgruppe von $\alpha$.\par}
 \smallskip
       Beweis. a) ist banal.
 \par
       b) Die Gruppe, die von $G$ auf $\alpha$ induziert wird, enth\"alt
 die Gruppe $T \cdot \SL(2,3)$, wie wir fr\"uher schon feststellten.
 Diese Gruppe ist vom Index 2 in der vollen Kollineationsgruppe von
 $\alpha$. Es gen\"ugt also, eine Kollineation von $\alpha$
 aufzuzeigen, die nicht durch eine Kollineation aus $G$ induziert wird.
 \par
       Die affine Ebene $\alpha$ besitzt eine involutorische
 Kollineation, die die Punkte $P_0$ $P_3$ und $P_6$ festl\"asst
 und die $P_4$ auf $P_8$ abbildet. Wir nehmen an, dass diese
 Kollineation durch eine Kollineation $\tau \in G$ induziert werde.
 Es gibt dann ein mit $t$ vertauschbares Element $\lambda \in K^*$
 mit $b_0^\tau = b_0\lambda$ und $b_1^\tau = b_1\lambda$ sowie ein $\mu \in K^*$
 mit $b_2^\tau = (b_0t+b_2)\mu$. Weil $\tau$ involutorisch
 ist, gibt es ein $\nu \in K^*$ mit
 $$ b_2\nu = b^{\tau^2}_2 = b_0(\lambda t\mu + t\mu^2) + b_2\mu^2. $$
 Weil $\lambda$ mit $t$ vertauschbar ist, folgt
 $$ 0 = \lambda t\mu + t\mu^2 = t(\lambda\mu + \mu^2), $$
 so dass $\mu = -\lambda$ ist.  Also ist $b^\tau_2 = -(b_0t + b_2) \lambda$.
 Weiter folgt $P^\tau_7 = P - 2$. Es gibt daher ein $\kappa \in K^*$ mit
 $$ (b_0 + b_1t + b_2) \kappa = (b_0 + b_1 + b_2)^\tau
        = b_0(\lambda - t\lambda) + b_1\lambda - b_2\lambda. $$
 Hieraus folgt $\kappa = \lambda - t\lambda$, $t\kappa = \lambda$ und
 $\kappa = -\lambda$ und dann $t = -1$, so dass nach 7.6 die Charakteristik von
 $K$ gleich 3 ist. Dies widerspricht aber unserer Annahme. Damit ist
 der Satz bewiesen.
 \medskip
       Im Falle der Charakteristik 3 ist nichts au\ss ergew\"ohnlich an
 der bislang betrachteten Situation. Wir werden sie daher im n\"achsten Satz
 ausschlie\ss en. In ihm studieren wir die Einbettung einer affinen Ebene der
 Ordnung 3 in eine pappossche projektive Ebene, da in einer solchen Ebene die
 Verh\"altnisse besonders einfach sind.
 \medskip\noindent
 {\bf 7.8. Satz.} {\it Es sei $K$ ein kommutativer K\"orper mit von
 $3$ verschiedener Charakteristik und $t$ sei ein Element aus $K$ mit
 $t^2 - t + 1 = 0$. Es sei ferner $V$ ein Vektorraum vom Rang $3$ \"uber
 $K$ und $b_0$, $b_1$, $b_2$ sei eine Basis von $V$. Schlie\ss lich
 sei $P$ der in 7.1 mittels dieser Daten definierte Monomorphismus
 der affinen Ebene der Ordnung $3$ in $L_K(V)$. Das Bild dieser
 affinen Ebene unter $P$ werde mit $\alpha$ be\-zeich\-net und $G$ sei
 der Stabilisator von $\alpha$ in $\PGL(V)$. Dann enth\"alt $G$
 einen Normalteiler $T$ der Ordnung $9$ und $G/T$ ist zur $\SL(2,3)$
 isomorph. Dar\"uber hinaus gilt $T = C_{\PGL(V)}(T)$.}
 \smallskip
       Beweis. Weil $K$ kommutativ ist, ist $L_K(V)$ pappossch, so dass
 mit 1.11 folgt, dass $G$ auf $\alpha$ treu operiert. Mit 7.7 b)
 folgt daher die erste Aussage \"uber die Gruppe $G$.  Es seien
 $\rho$ und $\sigma$ die unter diesen Namen in 7.4 definierten
 Abbildungen. Bezeichnet man die von diesen Abbildungen induzierten
 Kollineationen ebenfalls mit $\rho$ bzw. $\sigma$, so folgt $\rho
 \sigma \in T$. Ferner gilt
 $$\eqalign{
       b_0^{\rho \sigma} &= -b_0 - b_1t     \cr
       b_1^{\rho \sigma} &= b_0t^{-1}       \cr
       b_2^{\rho \sigma} &= b_0t + b_1t + b_2. \cr} $$
 Nun ist $(b_0r_0 + b_1r_2 + b_2r_2)K$ genau dann ein Fixpunkt von
 $\rho \sigma$, wenn es ein $\mu \in K^*$ gibt mit
 $$\eqalign{
        r_0\mu &= -r_0 + t^{-1}r_1 + tr_2   \cr
	r_1\mu &= -tr_0 + tr_2              \cr
	r_2\mu &= r_2. \cr}$$
 Es sei zun\"achst $r_2 = 0$. Dann ist $r_0 \neq 0$, da andernfalls
 auch $r_1 = 0$ w\"are. Ferner ist
 $$
       r_0\mu^2 = -r_0\mu + t^{-1}r_1\mu     
	         = -r_0\mu + t^{-1}(-tr_0)   
                 = -r_0(\mu + 1). 
 $$
 Wegen $r_0 \neq 0$ folgt $\mu^2 + \mu + 1 = 0$.
 \par
       Ist umgekehrt $\mu \in K$ und $\mu^2 + \mu + 1 = 0$, so ist
 $(b_0 - b_1t\mu^{-1})K$ ein Fixpunkt von $\rho \sigma$, wie man leicht
 nachpr\"uft.
 \par
       Weil $K$ kommutativ ist, hat das Polynom $x^2 + x + 1$ h\"ochstens
 zwei Nullstellen. Eine davon ist $-t$. Da die Charakteristik von
 $K$ von drei verschieden ist, ist $-t$ eine primitive dritte
 Einheitswurzel, so dass $t^2$ die zweite Nullstelle ist. Somit
 sind $(b_0 + b_1)K$ und $(b_0 - b_1t^{-1})K$ zwei Fixpunkte von $\rho\sigma$.
 \par
       F\"ur jeden weiteren Fixpunkt d\"urfen wir $r_2 = 1$ annehmen. Dann
 ist auch $\mu = 1$ und somit
 $$\eqalign{
           r_0 &= -r_0 + t^{-1}r_1 + t    \cr
	   r_1 &= -tr_0 + t               \cr}$$ 
 Hieraus folgt $r_0 = 1 + t$. Weil die Charakteristik nicht 3 ist, gibt
 es also noch genau einen weiteren Fixpunkt, so dass $\rho\sigma$
 genau drei Fixpunkte hat. Weil die von 1 verschiedenen Elemente
 von $T$ in $G$ alle konjugiert sind, haben sie alle genau drei
 Fixpunkte.
 \par
       Es sei nun $\gamma$ aus dem Zentralisator von $T$. Da $\gamma$ mit
 jedem Element von $T$ vertauschbar ist, permutiert $\gamma$ die
 drei Fixpunkte eines jeden von 1 verschiedenen Elementes von $G$
 unter sich. Hieraus folgt, dass $\alpha$ von $\gamma$ festgelassen
 wird, so dass $\gamma \in G$ gilt. In $G$ ist $T$ aber sein
 eigener Zentralisator. Damit ist alles bewiesen.
 \medskip
       Erst als wir $\mu$ mit $-t$ bzw. $t^2$ identifizierten, haben wir
 von der Kommutativit\"at von $K$ Gebrauch gemacht. Dies zeigt,
 dass es im Nichtkommutativen Kollineationen mit ungewohnten
 Fixpunktkonfigurationen\index{Fixpunktkonfiguration}{} gibt. Dies ist die eine
 Bemerkung, die zu 7.8 zu machen ist. Die andere ist die, dass $G$ im Falle eines
 nicht kommutativen K\"orpers nicht treu auf $\alpha$ operiert.
 Dies hat zur Folge, dass $T$ nicht sein eigener Zentralisator ist.
 \par
       Dieser Abschnitt f\"uhrt den Namen {\it hessesche
 Gruppe\/}\index{hessesche Gruppe}{} im Titel, so dass zum Schlu\ss\ noch einmal
 gesagt sein soll, dass
 sich unter diesem Namen die Gruppe $G$ des Satzes 7.8 verbirgt.

\mysection{8. Isomorphismen der gro\ss en projektiven Gruppe}

\noindent
 Ist $\sigma$ ein Isomorphismus oder ein Antiisomorphismus von
 $L_K(V)$ auf $L_{K'}(V')$ und definiert man $\varphi$ durch
 $\gamma^\varphi := \sigma^{-1}\gamma\sigma$ f\"ur alle $\gamma \in \PGaL(V)$,
 so ist $\varphi$, falls $\Rg_K(V) \geq 3$ ist, ein
 Isomorphismus von $\PGaL(V)$ auf $\PGaL(V')$, da $\sigma$ wegen
 $\Rg_K(V) \geq 3$ nach den Strukturs\"atzen durch eine semilineare
 Abbildung von $V$ auf $V'$, bzw. von $V$ auf $V'{^*}$, aufgefasst
 als Vektorraum \"uber $K_\circ$, induziert wird. Wir werden in
 diesem Abschnitt zeigen, dass auch umgekehrt jeder Isomorphismus
 von $\PGaL(V)$ auf $\PGaL(V')$ durch einen Isomorphismus oder
 Antiisomorphismus von $L_K(V)$ auf $L_{K'}(V')$ induziert
 wird, falls nur $\Rg_K(V) \geq 3$ und $\Rg_{K'}(V') \geq 3$ ist.
 Hat wenigstens eine der beiden projektiven Geometrien den Rang 2,
 so gilt dieser Sachverhalt nicht mehr, wie die
 Ausnahmeisomorphismen\index{Ausnahmeisomorphien}{} zeigen. Man kann jedoch
 mehr zeigen, als wir hier tun, muss sich dann aber anderer
 Beweismethoden\index{Beweismethoden}{} bedienen. Mehr
 zu diesem Thema findet sich in Dieudonn\'e (1955).
 \par
       Es sei daran erinnert, dass wir mit $\E(X,Y)$ die Gruppe aller
 Qua\-si\-e\-la\-ti\-o\-nen bezeichnen, deren Zentren in $X$ liegen und deren
 Achsen $Y$ enthalten. Ist $X = Y$, so schreiben wir abk\"urzend
 $\E(X)$ an Stelle von $\E(X,Y)$.
 \bigskip\noindent
 {\bf 8.1. Satz.} {\it Es sei $V$ ein Vektorraum \"uber dem K\"orper
 $K$ mit $\Rg_K(V) \geq 3$. Ferner sei $G$ eine Gruppe von
 Elationen von $L_K(V)$. Genau dann gilt $G = \E(U)$, wobei $U$
 entweder ein Punkt oder eine Hyperebene von $L_K(V)$ ist, wenn
 $G = C_{\PGaL(V)}(G)$ ist.}
 \smallskip
       Beweis. Es sei $G = \E(U)$ und $U$ sei eine Hyperebene von $L_K(V)$.
 Weil $\E(U)$ abelsch ist, ist $G \subseteq C_{\PGaL(V)}(G)$. Es sei
 $\gamma \in C_{\PGaL(V)}(G)$. Ist $P$
 Zentrum einer von Eins verschiedenen Elation aus $G$, so ist
 $P^\gamma = P$. Folglich ist $\gamma$ eine Perspektivit\"at mit der
 Achse $H$. Weil $\E(U)$ auf der Menge der nicht in $U$ liegenden
 Punkte transitiv ope\-riert, ist $\gamma = 1$, falls $\gamma$ einen
 Fixpunkt au\ss erhalb $U$ hat. Also ist $\gamma \in \E(U)=G$. Dies
 zeigt, dass $G  = C_{\PGaL(V)}(G)$ ist.
 \par
       Ist $U$ ein Punkt und ist $G = \E(U)$, so zeigt man ent\-spre\-chend,
 dass $G = C_{\PGaL(V)}(G)$ ist.
 \par
       Es sei nun umgekehrt $G = C_{\PGaL(V)}(G)$. Dann ist $G$
 insbesondere abelsch. Au\ss erdem ist $G \neq \{1\}$. Es sei $U$
 der von den Zentren der Elemente aus $G - \{1\}$ aufgespannte
 Unterraum und $W$ sei der Schnitt der Achsen der Elationen
 aus $G - \{1\}$. Ist $P$ Zentrum der Elation $\sigma \in G - \{1\}$
 und ist $H$ die Achse einer Elation $\tau \in G - \{1\}$, so folgt
 aus $\sigma \tau = \tau \sigma$, dass $P \leq H$ ist. Also ist $U
 \leq W$. Dies hat $G \subseteq \E (U,W)$ zur Folge. Nach 5.6 ist
 $\E(U,W)$ abelsch, so dass $G = \E(U,W)$ gilt. Es sei $P$ ein Punkt
 auf $W$ und $H$ eine Hyperebene, die $U+P$ enthalte. Dann ist
 $\E(P,H) \subseteq C_{\PGaL(V)}(G) = G$. Hieraus folgt $P \leq U$
 und $W \leq H$. Dies impliziert $U = W$. W\"are nun $U$ weder ein
 Punkt noch eine Hyperebene, so g\"abe es eine Gerade $X$ und eine
 Ko-Gerade $Y$ mit $X \leq U \leq Y$. Ferner g\"abe es ein $\alpha
 \in \Hom_K(V,X)$ mit $V^\alpha=X$ und $\Kern(\alpha) = Y$. Definierte
 man $\tau$ verm\"oge $v^\tau := v + v^\alpha$, so w\"are $\tau \in
 \E(U)=G$. Es folgte der Widerspruch $\Rg_K(X)=1$, da $G$ ja eine
 Gruppe von Elationen ist. Dieser Widerspruch zeigt, dass $U$ eine
 Punkt oder eine Hyperebene ist.
 \medskip\noindent
 {\bf 8.2. Satz.} {\it Es sei $V$ ein Vektorraum \"uber dem
 K\"orper $K$ und der Rang von $V$ sei mindestens $3$. Ferner sei
 $G$ eine Untergruppe von $\PGaL(V)$ und $\sigma$ eine
 involutorische Streckung von $L_K(V)$. Gilt
 \item{a)} Die Gruppe $G$ enth\"alt keine Involution,
 \item{b)} Es ist $G = C_{\PGaL(V)}(G)$,
 \item{c)} Ist $\gamma \in G$, so gibt es ein $\eta \in G$ mit
 $\gamma \sigma = \eta^{-1}\sigma \eta$,

 \noindent so gilt eine der folgenden Aussagen:
 \item{1)} Es gibt eine Hyperebene $H$ von $L_K(V)$ mit $G = \E(H)$.
 \item{2)} Es gibt einen Punkt $P$ von $L_K(V)$ mit $G = \E(P)$.
 \item{3)} Es ist $\Rg_K(V) = 3$, die Gruppe $G$ ist endlich und hat
 die Ordnung 9 und die Zentren der involutorischen Streckungen aus
 $GZ$, wobei $Z := \{1,\sigma\}$ gesetzt wird, sind die Punkte einer
 affinen Unterebene der Ordnung $3$ von $L_K(V)$.\par}
 \smallskip
       Beweis. Nach c) ist $\gamma\sigma$ eine involutorische Streckung.
 Daher ist $(\gamma \sigma)^2 = 1$ und folglich
 $\sigma\gamma\sigma = \gamma^{-1}$. Ist $Z := \{1,\sigma\}$, so ist also
 $M := GZ$ eine Gruppe. Ferner folgt aus $\gamma = \gamma\sigma\sigma$,
 dass $M$ von den Elementen $\gamma \sigma$ mit $\gamma \in G$ erzeugt wird.
 \par
       Es seien $\gamma$ und $\delta$ zwei verschiedene Elemente aus $G$.
 Haben $\gamma\sigma$ und $\delta\sigma$ beide das Zentrum $P$,
 so ist
 $$ \gamma\delta^{-1} = \gamma\sigma\delta\sigma $$
 eine Perspektivit\"at mit dem Zentrum $P$. Weil $\gamma\delta^{-1}$ von $G$
 nach b) zentralisiert wird, folgt, dass $P$ unter $G$ festbleibt. Daher haben
 wegen c) alle involutorischen Streckungen in $M$ das Zentrum $P$. Hieraus folgt,
 da die involutorischen Streckungen ja $M$ erzeugen, dass alle Elemente aus $M$
 Perspektivit\"aten mit dem Zentrum $P$ sind. Es sei $1 \neq \gamma \in G$ und
 $H$ sei die Achse von $\gamma$. Weil $G$ abelsch ist, folgt $H^G = H$. Wegen
 $\sigma^{-1}\gamma\sigma = \sigma\gamma \sigma = \gamma^{-1}$ gilt auch
 $H^\sigma = H$. L\"age $P$ nicht auf $H$, so w\"are $M$ eine Untergruppe von
 $\Delta(P,H)$. Da diese Gruppe zur multiplikativen Gruppe von $K$ isomorph
 w\"are, enthielte $M$ nur eine Involution. Hieraus folgt
 $\gamma\sigma = \sigma$ f\"ur alle $\gamma \in G$, so dass $G = \{1\}$ w\"are.
 Dieser Widerspruch zu b) zeigt, dass $P \leq H$ ist. Somit ist $G$ eine Gruppe
 von Elationen, so dass in diesem Falle mit 8.1 die Aussage 1) folgt.
 \par
       Haben $\gamma \sigma$ und $\delta\sigma$ beide die Achse $H$, so sieht man
 ganz ent\-spre\-chend, dass $G = \E(H)$ ist, dh., dass die Aussage 2) gilt.
 \par
       Wir nehmen nun an, dass f\"ur alle verschiedenen $\gamma$ und
 $\delta$ in $G$ die Kollineationen $\gamma\sigma$ und $\delta\sigma$
 verschiedene Zentren und verschiedene Achsen haben. Mit
 $C_\gamma$ bezeichnen wir das Zentrum und mit $H_\gamma$ die Achse
 von $\gamma \sigma$. Ferner setzen wir
 $$ A := \sum_{\gamma \in G - \{1\}}C_\gamma \qquad
 \mathrm{und}\qquad
  B := \bigcap_{\gamma \in G - \{1\}} H_\gamma. $$
 Dann gilt:
 \smallskip
       (1) Es ist $A = V$ oder $B = \{0\}$.
 \par
       Ist n\"amlich $A \neq V$ und $B \neq \{0\}$, so gibt es einen
 Punkt $P \leq B$ und eine Hyperebene $H$ mit $A \leq H$. Ist dann
 $\rho$ eine von 1 verschiedene Perspektivit\"at mit dem Zentrum
 $P$ und der Achse $H$, so ist $C_\gamma^\rho = C_\gamma$ und
 $H_\gamma^\rho = H_\gamma$ f\"ur alle $\gamma \in G$. Weil es in
 $\Delta(C_\gamma, H_\gamma)$ nur eine involutorische
 Perspektivit\"at gibt, folgt weiter $\rho^{-1} \gamma\sigma\rho = \gamma \sigma$
 f\"ur alle $\gamma \in G$. Weil $M$ von den
 Elementen $\gamma\sigma$ erzeugt wird, folgt schlie\ss lich
 $$ \rho \in C_{\PGaL(V)} (M) \subseteq C_{\PGaL(V)} (G) = G. $$
 Hieraus folgt
 $\rho = \sigma^{-1} \rho \sigma = \rho^{-1} $
 im Widerspruch zu a).
 \smallskip
       (2) Die Gruppe $G$ operiert auf $\{C_\gamma \mid \gamma \in G\}$
 und auf $\{H_\gamma \mid \gamma \in G\}$ regul\"ar.
 \par
       Da $G$ ja abelsch ist und $\sigma\eta = \eta^{-1}\sigma$ f\"ur
 alle $\eta \in G$ gilt, ist
 $$ \eta^{-1}\gamma\sigma\eta = \eta^{-2}\gamma\sigma $$
 f\"ur alle $\gamma$, $\eta \in G$. Also ist
 $C_\gamma^\eta = C_{\eta^{-2}\gamma}$ und $H_\gamma^\eta = H_{\eta^{-2}\gamma}$.
 Weil die Abbildungen $C$ und $H$ injektiv sind, folgt (2) somit aus a).
 \par
       Es sei $1 \neq \gamma \in G$. Dann l\"asst $\gamma$ wegen
 $\gamma = \gamma \sigma \sigma$ die Gerade $X := C_1+C_\gamma$
 hyperebenenweise und die Ko-Gerade $Y := H_1 \cap H_\gamma$
 punktweise fest. F\"ur dieses $\gamma$ gelten die n\"achsten beiden Aussagen.
 \smallskip
       (3) Es gibt ein $\eta \in G$ mit $X^\eta \neq X$ oder
 $Y^\eta \neq Y$.
 \par
       Andernfalls w\"are $X^G = X$ und $Y^G = Y$. Aus c) und $\Rg_K(V) \geq 3$
 folgte dann, dass $A = X\neq V$ und $B = Y\neq \{0\}$ w\"are im
 Widerspruch zu (1).
 \smallskip
       (4) Alle Fixpunkte von $\gamma$, die nicht in $Y$ liegen,
 liegen auf $X$ und alle Fixhyperebenen von $\gamma$, die $X$ nicht
 enthalten, enthalten $Y$.
 \par
       Es sei $F$ ein Fixpunkt von $\gamma$, der nicht in $Y$ liege. Aus
 (2) folgt, dass weder $C_1$ noch $C_\gamma$ ein Fixpunkt von
 $\gamma$ ist, so dass $F$ von diesen beiden Punkten verschieden
 ist. W\"are $F \leq H_1$, so w\"are $H_1 = F + Y$ und daher
 $H_1^\gamma = H_1$ im Widerspruch zu (2). Also ist $F \not\leq H_1$.
 Ebenso folgt $F \not\leq H_\gamma$. Daher ist $F \neq F^\sigma$ und
 $F \neq F^{\gamma \sigma}$. Es folgt weiter $C_1\leq F + F^\sigma$
 und $C_\gamma \leq F + F^{\gamma \sigma}$. Nun ist $F = F^\gamma$ und
 folglich $F^\sigma = F^{\gamma \sigma}$. Also ist
 $$ X = C_1 + C_\gamma \leq F + F^\sigma $$
 und damit $X = F + F^\sigma$, so dass $F \leq X$ gilt.
 \par
 Die zweite Aussage von (4) beweist sich analog.
 \smallskip
       (5) Es ist $\Rg_K(V) = 3$.
 \par
       Es sei $\Rg_K(V)\geq 4$. Dann ist $\Rg_K(Y)\geq 2$. Nach (3) gibt
 es ein $\eta \in G$ mit $X^\eta \neq X$ oder $Y^\eta \neq Y$. Es
 sei $Y^\eta \neq Y$. Die Punkte von $Y^\eta$, die nicht in $Y$
 liegen, spannen $Y^\eta$ auf. \"Uberdies sind diese Punkte
 Fixpunkte von $\gamma$, da $G$ ja abelsch ist. Hieraus und aus (4)
 folgt, dass $Y^\eta \leq X$ ist. Wegen $\Rg_K(X) = 2$ und
 $\Rg_K(Y) \geq 2$ ist daher $Y^\eta = X$. Weil die Punkte von
 $Y^\eta$ Fixpunkte von $\gamma$ sind und $C_1 \leq X = Y^\eta$ ist,
 ist $C_1^\gamma = C_1$. Dies widerspricht aber (2). Also ist
 $Y = Y^\eta$ und daher $X^\eta \neq X$. Weil $G$ abelsch ist,
 l\"asst $\gamma$ die Gerade $X^\eta$ hyperebenenweise fest. Nach
 I.5.3 ist $X^\eta$ die untere Grenze der es umfassenden
 Hyperebenen und es folgt, dass $X^\eta$ die untere Grenze
 derjenigen Hyperebenen ist,  die $X^\eta$ aber nicht $X$
 umfassen. Mit (4) folgt daher, dass $Y \leq X^\eta$ ist. Wegen
 $\Rg_K(X) = 2 \leq \Rg_K(Y)$ folgt $Y = X^\eta$. Hieraus folgt dann
 $H_1^\gamma=H_1$ im Widerspruch zu (2). Also ist doch $\Rg_K(V)=3$.
 \par
       Weil der Rang von $V$ endlich ist, ist der zu $L_K(V)$ duale
 Verband nach I.5.7 ebenfalls ein projektiver Verband. Wir d\"urfen
 daher annehmen, dass es ein $\eta \in G$ gibt mit $Y^\eta \neq Y$.
 Nach (4) ist dann $Y^\eta \leq X$, da $Y^\eta$ ein von $Y$
 verschiedener Fixpunkt von $\gamma$ ist. Weil $C_1$ kein Fixpunkt
 von $\gamma$ ist, ist $Y^\eta \neq C_1$. Weil $X$ eine Fixgerade
 von $\gamma$ ist, $H_1$ aber nicht, ist $X^\eta$ auch von $X \cup H_1$
 verschieden. Folglich gilt auch $Y^\eta \neq Y^{\eta\sigma}$. Damit haben wir in
 $Y^{\eta \sigma}$ einen dritten Fixpunkt von $\gamma$ gefunden.
 \par
       Es sei $X^G = X$. Dann ist $Y \leq X$, da ja $Y^\eta \leq X$ ist.
 Hieraus folgt $A = X$ und mit (1) dann $B = \{0\}$. Es gibt daher ein
 $\delta \in G$ mit $Y \not\leq H_\delta$. Nun ist
 $$ Y = Y \cup X \leq H_1\cup X < X. $$
 Daher ist $Y = H_1\cup X$. Setzt man $Z := H_1\cup H_\delta$, so ist also
 $Z \neq Y$ und folglich $Z \not\leq X$. Ersetzt man in unserer bisherigen
 Argumentation $\gamma$ durch
 $\delta$, so sieht man, dass $A = Z + C_1$ ist. Hiermit folgt der
 Widerspruch $Z \leq A = X$. Also ist $X$ keine Fixgerade von $G$.
 \par
       Da $X$ keine Fixgerade von $G$ ist, hat $\gamma$ nur die drei
 Fixpunkte $Y$, $Y^\eta$ und $Y^{\eta \sigma}$. Diese Fixpunkte
 sind \"uberdies nicht kollinear. Ferner gilt
 $$ Y^{\eta\sigma\eta} = Y^{\sigma\eta^{-1} \eta} = Y. $$
 Hieraus folgt, dass $Y^{\eta^2} =  Y^{\eta \sigma}$ ist, so dass die Fixpunkte
 von $\gamma$ von der von $\eta$ erzeugten Gruppe transitiv permutiert
 werden. Alles, was wir von $\gamma$ sagten, gilt aber auch f\"ur
 $\eta$, so dass auch $\eta$ genau drei Fixpunkte hat. Diese sind
 von den Fixpunkten von $\gamma$ verschieden. Es bezeichne $S$ die
 Menge der Kollineationen aus $G$, die die Fixpunkte von $\gamma$
 einzeln festlassen. Die Gruppe $S$ l\"asst die Menge der
 Fixpunkte von $\eta$ invariant. Da $\gamma$ zu $S$ geh\"ort,
 operiert $S$ auf der Menge der Fixpunkte von $\eta$
 transitiv. Weil die von 1 verschiedenen Elemente aus $G$ aber
 alle nur drei Fixpunkte haben, folgt, dass $S$ die Ordnung $3$
 hat. Dies wiederum hat zur Folge, dass $G$ die Ordnung 9 hat.
 \par
       Weil $G$ die Ordnung 9 hat, enth\"alt $GZ$ genau 9
 involutorische Streckungen. Sind nun $\gamma$, $\delta \in G$, so
 ist $C^{\delta\sigma}_\gamma$ ein von $C_\gamma$ und $C_\delta$
 verschiedenes Zentrum einer involutorischen Streckung aus $GZ$ auf
 $C_\gamma + C_\delta$. Beachte man, dass $\sigma\xi\sigma = \xi^{-1} = \xi^2$
 ist f\"ur alle $\xi \in G$, so folgt
 $$ (\delta\sigma)^{-1} \gamma\sigma\delta\sigma = (\gamma\delta)^2\sigma. $$
 Hieraus folgt $C_\gamma^{\delta\sigma} = C_{(\gamma\delta)^2}$, so dass die
 Zentren $C_\gamma$, $C_\delta$ und $C_{(\gamma \delta)^2}$ stets kollinear sind.
 Dieses Tripel ist eine Bahn der von $\gamma \delta^{-1}$ erzeugten
 Untergruppe von $G$. Denn es ist ja, wie wir fr\"uher gesehen
 haben, $C_\xi^\eta = C_{\eta^{-2}\xi} = C_{\eta \xi}$. Nun haben wir
 aber gezeigt, dass $A = V$ ist. Somit sind nicht alle $C_\xi$
 kollinear. Mit den zuvor gemachten Bemerkungen folgt dann, dass
 jede Gerade von $L_K(V)$, die zwei der Zentren tr\"agt, genau
 drei von ihnen tr\"agt. Hieraus folgt dann aber m\"uhelos, dass
 die $C_\xi$'s die Punkte einer affinen Unterebene von $L_K(V)$ sind.
 \smallskip
       Dass die letzte Situation tats\"achlich vorkommt, haben wir im
 letz\-ten Abschnitt gesehen.
 \medskip\noindent
 {\bf 8.3. Satz.} {\it Es sei $V$ ein Vektorraum \"uber dem K\"orper
 $K$, dessen Rang mindestens $3$ sei. Ferner seien $P$ ein Punkt,
 $H$ eine Hyperebene sowie $U$ und $W$ Unterr\"aume von $L_K(V)$
 und es gelte $V = P \oplus H=U \oplus W$. Ist $\sigma$ eine
 Involution aus $\Lambda(U,W)$ und gilt $P^\sigma = P$ sowie
 $H^\sigma = H$, so ist $\sigma$ mit allen Streckungen aus
 $\Delta(P,H)$ vertauschbar.}
 \smallskip
 Beweis. Weil $P$ ein Fixpunkt von $\sigma$ ist, ist $P \leq U$
 oder $P \leq W$. Wir d\"urfen annehmen, dass $P \leq U$ ist. Wir
 d\"urfen ferner annehmen, dass $u\sigma = -u$ ist f\"ur alle $u \in U$ und dass
 $w^\sigma = w$ ist f\"ur alle $w \in W$. Wegen $H^\sigma=H$ ist
 $$ H = (H \cup U) \oplus (H \cup W). $$
 Wegen $P \leq U$ folgt
 $$ V = P + H = U + (H \cup U) + (H \cup W) = U + (H\cup W). $$
 Hieraus erhalten wir mit Hilfe des Modulargesetzes
 $$ W = W \cup V = W \cup \bigl(U + (H \cup W)\bigr) = H \cup W. $$
 Somit ist $W \leq H$. Ist nun $\lambda \in \Delta (P,H)$, so ist also
 $U^\lambda = U$ und $W^\lambda = W$. Somit ist $\lambda^{-1}\sigma\lambda$ eine
 Involution aus $\Lambda (U,W)$. Da diese Gruppe aber nur eine
 Involution enth\"alt, folgt $\lambda^{-1}\sigma\lambda = \sigma$,
 womit der Satz bewiesen ist.
 \medskip\noindent
 {\bf 8.4. Satz.} {\it Es seien $K$ und $L$ K\"orper der
 Charakteristik $p > 0$. Ferner sei $X$ ein Vektorraum \"uber $L$
 des Ranges mindestens $2$. Ist dann $K^*/Z(K^*)$ zu $\PGaL(X,L)$
 isomorph, so ist $\Rg_L(X) = 2$ und $L$ ist kommutativ.}
 \smallskip
       Beweis. Wir beginnen den Beweis mit einer Vorbemerkung. Es sei $X$
 ein Vektorraum des Ranges 2 \"uber $L$. Ist $x_1$, $x_2$ eine
 Basis von $X$ und definiert man $\tau$ durch $x^\tau_1 := x_1$ und
 $x_2^\tau := x_1 + x_2$, so ist die Ordnung von $\tau$ gleich $p$. Ist
 $\gamma$ definiert durch $x_1^\gamma := x_1a + x_2c$ und
 $x^\gamma_2 := x_1b + x_2d$ und ist $v^{\tau\gamma} = v^{\gamma\tau}$ f\"ur
 alle $v\in X$, so gelten die Gleichungen
 $$\eqalign{
           a &= (a + c)r,  \cr
           c &= cr,        \cr
       a + b &= (b + d)r,  \cr
           d &= dr. \cr} $$
 Da $c$ und $d$ nicht beide gleich Null sein k\"onnen, ist $r = 1$.
 Dann ist aber $c = 0$ und $d = a$, so dass also $x^\gamma_1 = x_1a$ und
 $x^\gamma_2 = x_1b + x_2a$ ist.
 \par
       Sind umgekehrt $a$, $b \in L$, ist $a \neq 0$ und definiert man die
 Abbildung $\gamma$ durch $x^\gamma_1 := x_1a$ und
 $x^\gamma_2 := x_1b + x_2a$, so ist $\tau\gamma = \gamma\tau$. Ist
 nun $H$ der Zentralisator von $\tau$ in $\PGL(X,L)$ und ist $H_2$
 diejenige Untergruppe von $H$, die von allen $\gamma$ der Form
 $x_1^\gamma = x_1$ und $x_2^\gamma = x_1b + x_2$ induziert wird, so ist
 $H_2$ ein abelscher $p$-Normalteiler von $H$ und $H/H_2$ ist
 isomorph zu $L^*/Z(L^*)$. Ferner folgt, wie eine einfache Rechnung
 zeigt, dass $L$ genau dann kommutativ ist, wenn $H_2 \subseteq Z(H)$ gilt.
 \par
       Es sei nun $kZ(K^*)$ ein Element der Ordnung $p$ aus $K^*/Z(K^*)$.
 Ferner sei $x \in K^*$ und es gelte $x^{-1}kZ(K^*)x = kZ(K^*)$. Dann gibt es ein
 $z \in Z(K^*)$ mit $x^{-1}kx = kz$. Wegen $k^p \in Z(K^*)$ gilt daher
 $$ k^p = x^{-1}k^px = k^pz^p, $$
 so dass $z = 1$ ist. Der Zentralisator $C$ von $kZ(K^*)$ in $K^*/Z(K^*)$
 wird also durch den Zentralisator $F^*$ von $k$ in $K^*$ induziert.
 Dann ist aber $F:=F^* \cup \{0\}$ ein Teilk\"orper von $K$.
 \par
       Es sei $\varphi$ ein Isomorphismus von $\PGaL(X,L)$ auf
 $K^*/Z(K^*)$. Ferner sei $\tau$ eine von $1$ verschiedene Elation
 aus $\PGaL(X,L)$ und $H$ sei der Zentralisator von $\tau$ in
 $\PGaL(X,L)$. Mit $H_2$ bezeichnen wir die Gruppe aller Elationen,
 deren Achse gleich der Achse von $\tau$ ist. Dann ist $H_2$ ein
 abelscher $p$-Normalteiler von $H$. Ferner sei $S$ definiert durch
 $S/Z(K*) = H_2^\varphi$. Dann ist $S/Z(K^*)$ ein abelscher $p$-Normalteiler
 von $F^*/Z(K^*)$, wobei $F$ wie im Absatz zuvor bestimmt sei, wobei
 hier $\varphi(\tau)$ die Rolle von $k$ spielt. Es ist dann
 $H^\varphi = F^*/Z(K^*)$. Ferner ist klar, dass $Z(K^*) \subseteq Z(F^*)$ ist.
 Weil $F$ ein K\"orper ist, folgt mit 5.9 b), dass $S \subseteq Z(F^*)$ ist.
 Somit ist $H^\varphi_2 \subseteq Z(H)$.
 Mittels 5.10 folgt, dass der Rang von $X$ nicht gr\"o\ss er als
 2 sein kann, so dass $\Rg_L(X) = 2$ ist. Aus unserer Vorbemerkung
 folgt dann weiter, dass $L$ kommutativ ist, da ja $H_2$ hier wie
 dort die gleiche Bedeutung hat.
 \medskip
       Nun haben wir alle Vorbereitungen getroffen, um den Beweis des
 folgenden Satzes anzugehen, der jedoch immer noch eine ganze Weile
 in Anspruch nehmen wird.
 \medskip\noindent
 {\bf 8.5. Satz.} {\it Es sei $V$ ein Vektorraum \"uber dem
 K\"orper $K$ und $W$ sei ein Vektorraum \"uber dem K\"orper $L$.
 Beide Vektorr\"aume haben mindestens den Rang $3$. Ist $\varphi$
 ein Isomorphismus von $\PGaL(V)$ auf $\PGaL(W)$, so gibt es einen
 Isomorphismus oder einen Antiisomorphismus $\sigma$ von $L_K(V)$
 auf $L_L(W)$, so dass $\gamma^\varphi = \sigma^{-1}\gamma \sigma$
 ist f\"ur alle $\gamma \in \PGaL(V)$.}
 \smallskip
       Beweis. Wir zeigen zun\"achst, dass $K$ genau dann die
 Charakteristik 2 hat, wenn $L$ die Charakteristik 2 hat. Weil
 $\varphi^{-1}$ ein Isomorphismus von $\PGaL(W)$ auf $\PGaL(V)$
 ist, gen\"ugt es dazu zu zeigen, dass $L$ die Charakteristik 2
 hat, wenn $K$ die Charakteristik 2 hat.
 \par
       Es sei $\Char(K) = 2 \neq \Char(L)$. Ferner sei $\tau$ eine von 1
 verschiedene Elation von $L_K(V)$. Mit $H$ bezeichnen wir den
 Zentralisator von $\tau$ in $\PGaL(V)$. Ferner sei
 $$\{1\} \subseteq H_2\subseteq H_1 \subseteq H_0\subseteq H $$
 die in 5.10 beschriebene Normalreihe von $H$.
       Setze $\sigma := \tau^\varphi$. Dann ist $\sigma$ eine Involution
 aus $\PGaL(W)$, da $\tau$ wegen $\Char(K) = 2$ eine Involution ist.
 Ist $C$ der Zentralisator von $\sigma$ in $\PGaL(W)$, so ist
 $H^\varphi = C$.
 \par
       Eine l\"angere Argumentationskette wird zeigen, dass $\sigma$
 durch eine Involution aus $\GaL(W)$ induziert wird, was dann
 seinerseits auf einen Widerspruch f\"uhrt. Wir nehmen daher
 zun\"achst an, $\sigma$ werde nicht durch eine Involution aus
 $\GaL(W)$ induziert. Nach 6.2 und 6.5 enth\"alt $C$ dann einen
 Normalteiler $N$ vom Index 2. Ferner ist nach diesen S\"atzen
 $N/Z(N) \cong \PGaL(X,F)$, wobei $F$ ein K\"orper mit
 $Z(L)\subseteq F$ ist. Operiert $\sigma$ fixpunktfrei, so ist $F$
 eine quadratische Erweiterung von $L$, so dass also $|F| > 3$ ist.
 Ferner ist $\Rg_F(X) = {1 \over 2} \Rg_L(W) \geq {3 \over 2}$, so
 dass $\Rg_F(X) \geq 2$ ist. Hat $\sigma$ einen Fixpunkt, so ist
 $\Rg_F(X) = \Rg_L(W)\geq 3$. Es ist also stets $\Rg_F(X) \geq 2$ und
 $X$ ist nicht der Vektorraum vom Range 2 \"uber $\GF(2)$ oder
 $\GF(3)$. Schlie\ss lich ist $Z(N) \cong Z(F^*)/Z(L^*)$ und
 $Z(C) = \{1,\sigma\}$.
 \par
       Es sei $M$ ein in $H^\varphi_1$ enthaltener Normalteiler von $C$.
 Angenommen, es w\"are $MZ(N) \cap N \neq Z(N)$. Dann enthielte
 $(MZ(N) \cap N)/Z(N)$ nach 2.8 eine zu $\PSL(X,F)$ isomorphe
 Untergruppe, da ja $\Rg_F(X) \geq 2$ und da $X$ nicht der
 Vektorraum vom Range 2 \"uber $\GF(2)$ oder $\GF(3)$ ist. Dies
 widerspricht aber der Tatsache, dass $\PSL(X,F)$ keine 2-Gruppe
 ist. Also ist $MZ(N) \cap N = Z(V)$. Hie\-raus folgt $M \cap N \subseteq Z(N)$.
 Daher ist entweder $M \subseteq Z(N)$ oder
 $MN = C$, da ja $N$ ein Normalteiler vom Index 2 in $C$ ist.
 \par
       Weil $\tau$ eine Elation und weil $\Rg_K(V) \geq 3$ ist, ist $H_1$
 nicht abelsch. Daher ist $H^\varphi_1$ keine Untergruppe von
 $Z(N)$. Nach der gerade ge\-mach\-ten Bemerkung ist daher $C = H_1^\varphi N$.
 \par
       Es sei $Z(H_1^\varphi) \subseteq Z(N)$. Wegen
 $C_H(Z(H_1^\varphi) = H^\varphi_0$ ist dann $N \subseteq H_0^\varphi$. Daher ist
 $$ H^\varphi = C = NH^\varphi_1 \subseteq NH_0^\varphi = H_0^\varphi, $$
 so dass $H = H_0$ ist. Mit 5.10 folgt hieraus, dass $K$ kommutativ ist. Hieraus
 folgt ebenfalls mit 5.10, dass das Zentrum von $H$ zur additiven Gruppe
 von $K$ isomorph ist. Aus $|Z(C)| = 2$ folgt daher, dass $|K| = 2$
 ist. Dies hat wiederum zur Folge, dass die Gruppen $\GaL(A/P)$ und
 $\PGaL(A/P)$ isomorph sind, wenn $A$ die Achse und $P$ das Zentrum
 von $\tau$ bezeichnen. Nun ist $\PGaL(X,F)$ nicht aufl\"osbar, so
 dass $C$ und damit $H$ nicht aufl\"osbar ist. Dann ist aber auch
 $\PGaL(A/P)$ nicht aufl\"osbar. Aus 2.8 folgt daher, dass
 $\PGaL(A/P)$ und damit $H_0/H_1 = H/H_1$ keinen von $\{1\}$
 verschiedenen abelschen Normalteiler enth\"alt. Folglich ist
 $$ H_1Z(N)^{\varphi^{-1}}/H_1 = \{1\}, $$
 so dass $Z(N)^{\varphi^{-1}} \subseteq H_1$ ist. Also ist $Z(N)$ eine
 abelsche 2-Gruppe vom Exponenten 2 oder 4. Nun ist $Z(N) \cong Z(F^*)/Z(L^*)$.
 Ist $f \in Z(F^*)-Z(L^*)$ und $f^2 \in Z(L^*)$, so
 ist, da der Exponent von $Z(N)$ ein Teiler von 4 ist,
 $$ 1 + 4f + 6f^2 + 4f^3 + f^4 = (1 + f)^4 \in Z(L). $$
 Weil die Charakteristik von $L$ von $2$ verschieden ist und wir angenommen
 haben, dass $f^2 \in Z(L)$ ist, folgt weiter
 $$ f(f^2 + 1) \in Z(L). $$
 Weil $f$ aber kein Element von $Z(L)$ ist, folgt $f^2 + 1 = 0$. Dies besagt,
 dass $f$ eine primitive vierte Einheitswurzel von $Z(F)$ ist. Das
 impliziert wiederum, dass $Z(N)$ h\"ochstens eine und damit genau
 eine Involution enth\"alt. Da $Z(N)$ eine abelsche Gruppe vom
 Exponenten 2 oder 4 ist, folgt, dass $Z(N)$ eine zyklische Gruppe der Ordnung
 2 oder 4 ist. Weil $Z(F)$ und $Z(L)$ K\"orper sind,
 ist $|Z(F^*)/Z(L^*)| > 2$, daher ist $Z(N)$ die zyklische Gruppe der
 Ordnung $4$. Wegen $N \cap H_1^\varphi \subseteq Z(N)$ ist $N \cap H_1 = Z(N)$.
 Daher folgt aus
 $$ C/N = (NH_1^\varphi)/N\cong H_1^\varphi/(N \cap H_1^\varphi)
	= H_1^\varphi/Z(N), $$
 dass $|H_1| = 2|Z(N)|=8$ ist. Dies impliziert, dass $V$ endlich ist.
 \par
       Es sei $n$ der Rang von $V$. Dann ist, da $\Rg_K(A) = n - 1$ und
 $\Rg_K(P) = 1$ sowie $|K| = 2$ ist,
 $$ \big|\Hom_K(A/P,P)\big| = 2^{n - 2} $$
 und
 $$ \big|\Hom_K(V/A,A)\big| = 2^{n-1}. $$
 Mit 5.10 folgt hieraus
 $$ 2^3 = |H_1/H_2|\cdot |H_2| = 2^{2n-3} $$
 und damit $n = 3$. Dies impliziert $\Rg_K(A/P) = 1$, was wiederum
 $\GaL(A/P) = \{1\}$ nach sich zieht, so dass $H$ aufl\"osbar ist.
 Dieser Widerspruch zeigt, dass $Z(H_1^\varphi)$ nicht in $Z(N)$
 enthalten ist. Dann ist aber, wie wir gesehen haben,
 $NZ(H_1^\varphi) = C$. Ferner ist $N \cap Z(H_1^\varphi) \subseteq Z(N)$.
 Hieraus folgt
 $$ N \cap Z(H_1^\varphi) \subseteq Z(NH_1^\varphi) = Z(C). $$
 Andererseits ist
 $$ Z(C) \subseteq N \cap Z(H_1^\varphi). $$
 Also ist
 $$ \big|N \cap Z(H_1^\varphi)\big| = 2. $$
 Aus $|C/N| = 2$ und $C = NZ(H_1^\varphi)$ folgt
 daher $|Z(H_1^\varphi)| = 4$. Nach 5.10 ist $Z(H_1)$ zu
 $\Hom_K(V/A,P)$ isomorph. Da hiernach $\Hom_K(V/A,P)$ genau vier
 Elemente enth\"alt, hat auch $K$ genau vier Elemente. Somit ist
 $K$ kommutativ. Nach 5.10 ist daher $H = H_0$ und daher
 $Z(H) = Z(H_1)$. Hiermit folgt der Widerspruch
 $$ 2 = \big|Z(H)\big| = \big|Z(H_0)\big| = 4. $$
 Dieser Widerspruch r\"uhrt daher, dass wir
 annahmen, $\sigma$ w\"urde nicht durch eine Involution aus
 $\GaL(W)$ induziert, so dass diese Annahme also zu verwerfen ist.
 \par
       Die Kollineation $\sigma$ wird also durch eine Involution
 induziert. Ist $H$ aufl\"osbar, so ist auch $C$ aufl\"osbar. Nach
 6.6, 5.10 bzw. 5.4 und 2.6 folgt dann $|K| = 2$ und $|L| = 3$ sowie
 $\Rg_K(V) \leq 4$ und $\Rg_L(W) \leq 4$. Nach 1.7 und 1.13 ist aber
 $$\eqalign{
       \big|\PGaL(3,2)\big| &= 2^3\cdot 3 \cdot 7,   \cr
       \big|\PGaL(4,2)\big| &= 2^6\cdot 3^2 \cdot 5 \cdot 7, \cr
       \big|\PGaL(3,3)\big| &= 3^3\cdot 2^4 \cdot 13,        \cr
       \big|\PGaL(4,3)\big| &= 3^6\cdot 2^8 \cdot 5 \cdot 13,\cr}$$
 so dass $\PGaL(V)$ und $\PGaL(W)$ nicht isomorph sein k\"onnen.
 Also ist $H$ und damit auch $C$ nicht aufl\"osbar.
 \par
       Nach 6.6 gibt es Unterr\"aume $X$ und $Y$ von $W$ mit $W = X \oplus Y$
 und $\sigma \in \Lambda(X,Y)$. Nach 5.2 ist $\Lambda(X,Y)$ zur
 multiplikativen Gruppe von $L$ oder von $Z(L)$ isomorph. Weil
 $\sigma$ involutorisch ist, liegt $\sigma$ in jedem Falle im
 Zentrum von $\Lambda(X,Y)$. Es sei nun $\Lambda_1$ diejenige
 Untergruppe von $\PGaL(W)$, die $X$ als Ganzes und $Y$ punktweise
 festl\"asst, und $\Lambda_2$ sei diejenige Untergruppe von
 $\PGaL(W)$, die $X$ punktweise und $Y$ als Ganzes festl\"asst.
 Nach 5.5 ist dann
 $$ |C:\Lambda_1\Lambda_2| \leq 2. $$
 \par
       Liegt $H_1^\varphi$ nicht in $\Lambda_1\Lambda_2$, so ist
 $C \neq \Lambda_1\Lambda_2$ und daher $|C:\Lambda_1\Lambda_2| = 2$. Es
 folgt $C = \Lambda_1\Lambda_2H_1^\varphi$. Dies impliziert
 seinerseits, dass es ein $\eta \in H_1^\varphi$ gibt mit
 $\eta^{-1}\Lambda_1\eta = \Lambda_2$. Aus $|C:\Lambda_1\Lambda_2| = 2$ und
 $C = \Lambda_1\Lambda_2H_1^\varphi$ folgt, dass
 $$ \big|H_1^\varphi : \bigl((\Lambda_1\Lambda_2) \cap H_1^\varphi\bigr)\big|
	  = 2 $$
 ist. Ferner ist $(\Lambda_1\Lambda_2) \cap H_1^\varphi$ ein 2-Normalteiler von
 $\Lambda_1\Lambda_2$. Hieraus folgt, dass
 $$ \bigl(\Lambda(X,Y)((\Lambda_1\Lambda_2)\cap H_1^\varphi)\bigr)/\Lambda(X,Y)$$
 ein 2-Normalteiler von
 $$ (\Lambda_1\Lambda_2)/\Lambda(X,Y) $$
 ist. Nach 5.4 e) ist $(\Lambda_1\Lambda_2)/\Lambda(X,Y)$ zu
 $$ \PGaL(X) \times \PGaL(Y) $$
 isomorph. Weil $\eta^{-1}\Lambda_1\eta = \Lambda_2$
 ist, ist $X^\eta = Y$, so dass $\Rg_L(X) = \Rg_L(Y) \geq 2$ ist.
 Enthielte $\Lambda_1\Lambda_2/\Lambda(X,Y)$ einen nicht trivialen
 $2$-Normalteiler, so folgt $|L| = 3$ und $\Rg_L(X) = \Rg_L(Y)=2$, so
 dass $C$ und damit $H$ aufl\"osbar w\"are. Dieser Widerspruch
 zeigt, dass
 $$ (\Lambda_1\Lambda_2) \cap H_1\varphi \subseteq \Lambda(X,Y) $$
 gilt. Wegen $\Rg_L(X) = \Rg_L(Y) \geq 2$ ist $\Lambda(X,Y) \cong Z(K^*)$ nach
 5.4 e). Weil der Exponent von $H_1^\varphi$ gleich 4 ist, folgt daher
 $$ \big|(\Lambda_1\Lambda_2) \subseteq H_1^\varphi\big| \leq 4, $$
 so dass $|H_1| = |H_1^\varphi| \leq 8$ ist. Folglich ist wieder $|K| = 2$ und
 $\Rg_K(V) = 3$, was den Widerspruch der Aufl\"osbarkeit von $H$ nach
 sich zieht. Also ist doch $H_1^\varphi \subseteq \Lambda_1\Lambda_2$.
 \par
       Weil der Rang von $V$ mindestens 3 ist, folgt mit 5.10, dass
 $H_2$ mindestens zwei Involutionen enth\"alt. Daher kann
 $H_2^\varphi$ nicht in der Gruppe $\Lambda(U,V)$ enthalten sein,
 da diese Gruppe zur multiplikativen Gruppe von $L$ oder $Z(L)$
 isomorph ist. Hieraus folgt, dass
 $$ \bigl(\Lambda(X,Y)H_1^\varphi\bigr)/\Lambda(X,Y) $$
 ein nicht trivialer 2-Normalteiler von $(\Lambda_1\Lambda_2)/\Lambda(X,Y)$ ist.
 Nun ist
 $$ (\Lambda_1\Lambda_2)/\Lambda(X,Y)\cong \PGaL(X) \times \PGaL(Y). $$
 Weil $C$ nicht aufl\"osbar ist, k\"onnen
 $\PGaL(X)$ und $\PGaL(Y)$ nicht beide aufl\"osbar sein, da sonst
 $L$ nur drei Elemente enthielte und $C$ dann doch aufl\"osbar
 w\"are. Daher kann nur eine der Gruppen $\PGaL(X)$ und $\PGaL(Y)$
 einen nicht trivialen 2-Normalteiler enthalten. Wir d\"urfen
 daher annehmen, dass die Gruppe $\PGaL(X)$ einen zu
 $(\Lambda(X,Y)H_1^\varphi)/H_1^\varphi)/\Lambda(X,Y)$ isomorphen
 2-Normalteiler ent\-h\"alt. Hieraus folgt $\Rg_L(X) = 1$ oder
 $|L| = 3$ und $\Rg_L(X) = 2$. Weil $\PGaL(X)$ nicht trivial ist, kann
 aber der Fall $\Rg_L(X) =  1$ nicht eintreten. Dies impliziert, dass
 der maximale 2-Normalteiler der Gruppe
 $(\Lambda_1\Lambda_2)/\Lambda(X,Y)$ die Ordnung 4 hat. Weil
 $\Lambda(X,Y)$ zur multiplikativen Gruppe von $K$ isomorph ist und
 folglich die Ordnung 3 hat, hat der maximale 2-Normalteiler
 von $\Lambda_1\Lambda_2$ die Ordnung 8. Hieraus folgt, dass
 $|H_1| = 8$ und $\Rg_K(V) = 3$ ist, was wiederum die Aufl\"osbarkeit
 von $H$ impliziert. Dieser Widerspruch zeigt endlich, dass die
 Charakteristik von $L$ doch gleich 2 ist.
 \par
       Es sei weiter
 $\Char(K) = 2$. Dann ist, wie wir jetzt wissen, auch
 $\Char(L) = 2$. Wir zeigen nun, dass jede Elation aus $\PGaL(V)$ auf
 eine Elation aus $\PGaL(W)$ abgebildet wird. Es sei also $\tau$
 eine von eins verschiedene Elation aus $\PGaL(V)$. Dann ist
 $\sigma := \tau^\varphi$ eine Involution aus $\PGaL(W)$. Es sei
 wieder $H$ der Zentralisator von $\tau$ in $\PGaL(V)$ und $C$ der
 Zentralisator von $\sigma$ in $\PGaL(W)$.
 \par
       Wird $\sigma$ nicht durch eine Involution aus $\GaL(W)$ induziert,
 so ist die Faktorgruppe $C/Z(C)$ nach 6.2 bzw. 6.5 isomorph zur
 Gruppe $\PGaL(X,F)$, wobei $F$ ein K\"orper mit $Z(L) \subseteq F$
 ist. Operiert $\sigma$ fixpunktfrei auf $L_L(W)$, so ist $F$
 eine quadratische Erweiterung von $L$, so dass $|F| > 3$ ist. Ferner
 ist in diesem Falle $\Rg_F(X)= {1 \over 2} \Rg_L(W) \geq 2$. Hat
 $\sigma$ einen Fixpunkt, so ist $\Rg_F(X) = \Rg_L(W) \geq 3$. Es ist
 also stets $\Rg_F(X) \geq 2$ und $X$ ist nicht der Vektorraum vom
 Range $2$ \"uber $\GF(2)$ bzw. $\GF(3)$. Mit Hilfe von 2.8 folgt
 daher, dass $C/Z(C)$ keinen von eins verschiedenen
 2-Normalteiler enth\"alt. Hieraus folgt, dass $H_1^\varphi \subseteq Z(C)$ ist.
 Dies widerspricht aber der Tatsache, dass $H_1$ nicht abelsch ist. Dieser
 Widerspruch zeigt, dass $\sigma$ eine Quasielation ist.
 \par
       Weil $\sigma$ eine Quasielation ist, besitzt $C$ eine Normalreihe
 $\{1\} \subseteq C_1 \subseteq C_0 \subseteq C$, wobei $C$ der
 gr\"o\ss te aufl\"osbare 2-Normalteiler von $C$ und $C_0$ der
 Zentralisator von $Z(C_1)$ in $C$ ist. Daher ist $H_1^\varphi =
 C_1$ und $H_0^\varphi = C_0$. Hieraus folgt, dass $H/H_0$ zu
 $C/C_0$ isomorph ist. Ist $K$ kommutativ, so ist $H = H_0$ und damit
 $C = C_0$. Dies besagt, dass $\sigma$ eine Elation ist und dass auch
 $L$ kommutativ ist. Ist $\sigma$ keine Elation, so ist also $K$
 nicht kommutativ. Weil $\varphi^{-1}$ ein Isomorphismus von
 $\PGaL(W)$ auf $\PGaL(V)$ ist, folgt, dass auch $L$ nicht
 kommutativ ist. Ist nun $B$ die Achse von $\sigma$, so ist
 $\Rg_L(W/B) \geq 2$. Ferner ist $C/C_0$ zu $\PGaL(W/B)$ isomorph,
 so dass wegen $H/H_0 \cong K^*/Z(K^*)$ die Gruppen $K^*/Z(K^*)$
 und $\PGaL(W/B)$ isomorph sind. Wegen $\Char(K) = 2 = \Char(L)$ folgt
 daher mit 8.4 der Widerspruch, dass $L$ kommutativ ist. Damit ist
 gezeigt, dass im Falle $\Char(K) = 2$ der Isomorphismus $\varphi$
 Elationen auf Elationen abbildet.
 \par
       Es sei nun $\Char(K) \neq 2$. Wir wollen zeigen, dass $\varphi$
 auch in diesem Falle Elationen auf Elationen abbildet. Dazu zeigen
 wir zun\"achst, dass $\varphi$ involutorische Streckungen auf
 involutorische Streckungen abbildet. Zun\"achst notieren wir, dass
 auch die Charakteristik von $L$ von $2$ verschieden ist, wie wir
 nicht ohne M\"uhe zeigten. Es sei nun $\rho$ eine involutorische
 Streckung von $L_K(V)$ mit der Achse $A$ und dem Zentrum $P$. Es
 sei $C$ der Zentralisator von $\rho$ in $\PGaL(V)$. Es sei
 $\Lambda_1$ die Gruppe, die $P$ punktweise und $A$ als Ganzes
 festl\"asst, und $\Lambda_2$ sei die Gruppe, die $A$ punktweise
 und $P$ als Ganzes festl\"asst. Wegen $\Rg_K(P) = 1$ ist dann
 $\Lambda_2 = \Delta(P,A)$. Ferner ist $\Lambda_2 \subseteq \Lambda_1$, so dass
 insbesondere $\Lambda_1 = \Lambda_1\Lambda_2$
 ist. Mit 5.4 folgt weiter, dass $|C:\Lambda_1| = 1$ oder 2 ist.
 Wegen $\Rg_K(P) = 1 < \Rg_K(A)$ kann der zweite Fall nach 5.4 nicht
 eintreten. Also ist $C = \Lambda_1$. Nochmals 5.4 zitierend sehen
 wir, dass $Z(C) = Z(\Delta(P,A))$ ist.
 \par
       Wir setzen $\sigma := \rho^\varphi$ und 
 bezeichnen mit $D$  den
 Zentralisator von $\sigma$ in $\PGaL(W)$. Wird $\sigma$ nicht
 durch eine Involution aus $\GaL(W)$ induziert, so enth\"alt $D$
 nach 6.2 bzw. 6.5 einen Normalteiler $N$ vom Index $2$ mit
 $N/Z(N) \cong \PGaL(X,F)$, wobei $F$ ein K\"orper mit $Z(L) \subseteq F$ ist.
 Ferner ist, wie wir schon zweimal bemerkten,
 $\Rg_F(X) \geq 2$ und $X$ ist nicht der Vektorraum vom Range 2
 \"uber $\GF(2)$ oder $\GF(3)$. Schlie\ss lich ist $Z(N) \cong Z(F^*)/Z(L^*)$ und
 $Z(D) = \{1,\sigma\}$. W\"are $|K| = 3$ und
 $\Rg_K(A) = 2$, so w\"are $C$ und damit $D$ aufl\"osbar, was nicht
 der Fall ist. Somit enth\"alt $\Lambda_1/\Delta(P,A)$ wegen 2.8
 keinen abelschen Normalteiler ungleich Eins. Hieraus folgt, dass
 $Z(N)^{\varphi^{-1}} \subseteq \Delta(P,A)$ ist. Wegen
 $$ Z(C) = Z(D)^{\varphi^{-1}} \subseteq Z(N)^{\varphi^{-1}} $$
 enth\"alt die multiplikative Gruppe von $K$ einen abelschen
 Normalteiler $A$ ungleich eins mit $Z(K^*) \subseteq A$. Da $A$
 einen kommutativen, normalen Teilk\"orper von $K$ erzeugt, folgt
 nach dem Satz von Cartan-Brauer-Hua, dass $A \subseteq Z(K^*)$
 ist. Also ist $A = Z(K^*)$ und daher $Z(C) = Z(N)^{\varphi^{-1}}$.
 Hieraus folgt $Z(N) = Z(D)$. Dies ergibt den Widerspruch
 $|Z(F^*)/Z(L^*)| = 2$. Somit ist $\sigma$ eine Quasistreckung.
 \par
       Angenommen $\sigma$ ist keine Streckung. Ist $\Rg_K(A) = 2$ und
 $|K| = 3$, so ist $\Rg_K(V) = 3$ und $C$ ist aufl\"osbar.
 Dann ist auch $D$ aufl\"osbar, so dass nach
 5.4 auch $|L| = 3$ sowie $\Rg_L(W)\leq 4$ ist. Ein Vergleich der
 Ordnungen von $\PGaL(V)$ und $\PGaL(W)$ zeigt, dass auch der Rang
 von $W$ gleich 3 ist. Dann ist $\sigma$ aber eine Streckung.
 Also ist $V$ nicht der Vektorraum vom Range $3$ \"uber $\GF(3)$.
 \par
       Es sei $W = X \oplus Y$ und $\sigma \in \Lambda(X,Y)$. Ferner sei
 $Q$ ein Punkt in $X$ und $B$ ein $Y$ umfassendes Komplement von
 $Q$. Ist $\lambda$ die involutorische Streckung aus $\Delta(Q,B)$,
 so ist $X^\lambda = Y$ und $Y^\lambda = Y$. Daher ist
 $\lambda \sigma = \sigma \lambda$. Folglich ist
 $\lambda^{\varphi^{-1}}\rho = \rho\lambda^{\varphi^{-1}}$. Hieraus folgt, dass
 $\varphi^{-1}$ sowohl
 $P$ als auch $A$ festl\"asst. Nach dem bereits Bewiesenen ist
 $\lambda^{\varphi^{-1}}$ eine involutorische Quasistreckung. Nach
 8.3 wird folglich $\Delta(P,A)$ von $\lambda^{\varphi^{-1}}$
 zentralisiert. Daher wird $\lambda$ von $\Delta(P,A)^\varphi$
 zentralisiert, woraus folgt, dass $Q$ ein Fixpunkt von
 $\Delta(P,A)^\varphi$ ist. Weil $Q$ ein beliebiger Punkt in $X$
 war, hei\ss t das, dass $X$ von $\Delta(P,A)^\varphi$ punktweise
 festgelassen wird. Ebenso zeigt man, dass $\Delta(P,A)^\varphi$
 auch $Y$ punktweise festl\"asst. Somit ist $\Delta(P,A)^\varphi$
 eine Untergruppe von $\Lambda(X,Y)$. Weil $\sigma$ keine Streckung
 ist, ist $\Lambda(X,Y)$ nach 5.2 b) abelsch. Folglich ist
 $\Delta(P,A)$ abelsch, so dass nach 5.4 f) die Gleichung
 $Z(C) = \Delta(P,A)$ gilt. Weil $C/\Delta(P,A)$ nicht aufl\"osbar
 ist, folgt aus 2.8, dass $C/\Delta(P,A)$ keinen nicht trivialen
 abelschen Normalteiler enth\"alt. Da andererseits
 $\Lambda(X,Y)^{\varphi^{-1}}/\Delta(P,A)$ ein abelscher
 Normalteiler dieser Gruppe ist, folgt
 $\Lambda(X,Y)^{\varphi^{-1}} = \Delta(P,A)$. Hieraus folgt, dass
 $\varphi^{-1}$ einen Isomorphismus von $D/\Lambda(X,Y)$ auf $C/\Delta(P,A)$
 induziert. Daher enth\"alt die Grup\-pe $\PGaL(A)$ einen Normalteiler
 $N = N_1 \times N_2$ mit $N_1 \cong \PGaL(X)$ und $N_2 \cong \PGaL(Y)$. Aus 2.8
 folgt wiederum, dass $\PSL(A)$ in $N$
 enthalten ist, so dass sowohl $N_1$ als auch $N_2$ die Gruppe
 $\PSL(A)$ normalisieren. Eine nochmalige Anwendung von 2.8 liefert
 daher den Widerspruch $\PSL(A) \subseteq N_1 \cap N_2 = \{1\}$.
 Dieser Widerspruch zeigt, dass $\sigma$ doch eine involutorische
 Streckung ist.
 \par
       Wir zeigen nun, dass $\varphi$ auch im Falle einer von 2
 verschiedenen Charakteristik die Elationen aus $\PGaL(V)$ auf die
 Elationen von $\PGaL(W)$ abbildet. Um dies zu zeigen, sei $H$ ---
 der Buchstabe $H$ steht wieder zur Verf\"ugung --- eine Hyperebene
 von $L_K(V)$ und $G$ sei die Gruppe aller Elationen mit der
 Achse $H$. Schlie\ss lich sei $\sigma$ eine involutorische
 Streckung mit der Achse $H$ und $Z$ sei die von $\sigma$ erzeugte
 Gruppe der Ordnung 2. Nach 8.1 ist die Gruppe $G$ ihr eigener
 Zentralisator in $\PGaL(V)$ und weil die Charakteristik von $K$
 von $2$ verschieden ist, enth\"alt $G$ auch keine Involution.
 Ferner ist $|GZ:G| = 2$, so dass die Kollineationen $\gamma\sigma$
 mit $\gamma \in G$ involutorische Streckungen mit der Achse $H$
 sind. Weil $G$ auf der Punktmenge von $L_K(V)_H$ transitiv
 operiert, gibt es zu jedem $\gamma \in G$ ein $\eta \in G$ mit
 $\gamma \sigma = \eta^{-1} \sigma \eta$. Folglich erf\"ullen $G$
 und $\sigma$ die Voraussetzungen von 8.2. Nun ist $\sigma^\varphi$
 eine involutorische Streckung von $L_L(W)$, wie wir gesehen
 haben. Daher erf\"ullen auch $G^\varphi$ und $\sigma^\varphi$ die
 Voraussetzungen von 8.2. Somit ist $G^\varphi$ eine Gruppe von
 Elationen, es sei denn, es ist $\Rg_L(W) = 3$ und $G$ hat die
 Ordnung 9. In diesem Falle besitzt $L_K(V)_H$ aber nur neun
 Punkte, so dass $V$ der Vektorraum vom Rang $3$ \"uber $\GF(3)$
 ist. Weil die Gruppen $\PGaL(V)$ und $\PGaL(W)$ isomorph sind,
 folgt dann, dass auch $L = \GF(3)$ ist. Weil die Zentren der
 Streckungen aus $(GZ)^\varphi$ die Punkte einer affinen Unterebene
 der Ordnung 3 von $L_L(W)$ ist, folgt mittels 7.2 und 7.6,
 dass die Streckungen aus $(GZ)^\varphi$ doch alle die gleiche
 Achse haben. Also kann dieser Fall nicht eintreten. Dies zeigt,
 dass $\varphi$ Elationen auf Elationen abbildet.
 \par
       Es sei nun $P$ ein Punkt und $H$ eine Hyperebene von $L_K(V)$.
 Offensichtlich gilt genau dann $P \leq H$, wenn $\E(P) \cap \E(H) \neq \{1\}$
 ist. Weil $\varphi$ Elationen auf Elationen abbildet,
 sind $\E(P)^\varphi$ und $\E(H)^\varphi$ Gruppen von Elationen.
 Mit 8.1 folgt die Existenz von Unterr\"aumen $P^\rho$ bzw.
 $H^\rho$ von $W$ mit $\Rg_L(P^\rho) = 1$ bzw. $\Ko_L(P^\rho) = 1$ und
 $\Ko_L(H^\rho) = 1$ bzw. $\Rg_L(H^\rho) = 1$, so dass $\E(P)^\varphi =
 \E(P^\rho)$ und $\E(H)^\varphi = \E(H^\rho)$ ist. Nun sind $\E(P)$
 und $\E(H)$ nicht konjugiert unter $\PGaL(V)$, so dass
 $\E(P^\rho)$ und $\E(H^\rho)$ in verschiedenen
 Konjugiertenklassen von $\PGaL(W)$ liegen. Daher ist entweder
 $$ \Rg_L(P^\rho) = \Ko_L(H^\rho) = 1 $$
 f\"ur alle Punkte $P$ und alle Hyperebenen $H$ von $L_K(V)$ oder es ist
 $$ \Ko_L(P^\rho) = \Rg_L(H^\rho) = 1 $$
 f\"ur alle Punkte $P$ und alle
 Hyperebenen $H$ von $L_K(V)$. In beiden F\"allen ist $\rho$ eine
 Bijektion der Menge der Punkte und der Hyperebenen von $L_K(V)$
 auf die Menge der Punkte und Hyperebenen von $L_L(V)$. Im ersten
 Falle werden dabei die Punktmenge auf die Punktmenge und die
 Hyperebenenmenge auf die Hyperebenenmenge abgebildet, w\"ahrend im
 zweiten Falle die Punktmenge auf die Hyperebenenmenge und die
 Hyperebenenmenge auf die Punktmenge abgebildet werden. Weil genau
 dann $\E(P) \cap \E(H) \neq \{1\}$ gilt, wenn $\E(P^\rho) \cap
 \E(H^\rho) \neq \{1\}$ gilt, ist $\rho$ in beiden F\"allen
 inzidenztreu. Mit I.8.6 und der \"Ubungsaufgabe zu diesem Satz
 folgt, dass $\rho$ durch einen Isomorphismus oder
 Antiisomorphismus $\sigma$ von $L_K(V)$ auf $L_L(W)$ induziert
 wird.
 \par
       Es sei schlie\ss lich $\gamma \in \PGaL(V)$. Dann ist
 $$
    \E(P^{\sigma \gamma^\varphi}) = \gamma^{-\varphi}\E(P^\sigma)\gamma^\varphi
    = \gamma^{-1}\E(P)^\varphi \gamma^\varphi
    =(\gamma^{-1}\E(P)\gamma)^\varphi      
   = \E(P^\gamma)^\varphi = \E(P^{\gamma \sigma}).
 $$
 Hieraus folgt $P^{\sigma\gamma^\varphi} = P^{\gamma \sigma}$ f\"ur
 alle Punkte $P$ von $L_K(V)$. Da ein Isomorphismus bzw. ein
 Antiisomorphismus durch seine Wirkung auf die Punkte bereits
 festgelegt wird, folgt $\sigma\gamma^\varphi = \gamma\sigma$, so dass
 $$ \gamma^\varphi = \sigma^{-1}\gamma\sigma $$
 f\"ur alle $\gamma \in \PGaL(V)$ gilt. Damit ist 8.5 endlich bewiesen.
 \medskip
       Das folgende Korollar ist eine unmittelbare Folgerung aus 8.5.
 \medskip\noindent
 {\bf 8.6. Korollar.} {\it Sei $V$ ein Vektorraum \"uber dem
 K\"orper $K$ mit $\Rg_K(V) \geq 3$.
 \item{a)} Ist $L_K(V)$ nicht selbstdual, so ist $\PGammaL(V)$
 zur Automorphismengruppe von $\PGaL(V)$ isomorph.
 \item{b)} Ist $L_K(V)$ selbstdual, so enth\"alt die
 Automorphismengruppe von $\PGaL(V)$ einen zu $\PGammaL(V)$
 isomorphen Normalteiler vom Index $2$.
 \item{c)} Die Gruppen $\PSL(V)$ und $PGL(V)$ sind charakteristische
 Untergruppen von $\PGaL(V)$, dh., sie sind unter allen
 Automorphismen von $\PGaL(V)$ invariant.\par}
 \medskip
       Der Leser wird vielleicht fragen, wie man auf all dies komme. Diese Frage
 kann ich hier ausnahmsweise einmal beantworten. Geometrisch sind Elationen
 und Streckungen am besten zu fassen, gruppentheoretisch Involutionen.
 Daher stellt man zuerst die Frage, wie man die unterschiedlichen Involutionen
 in den fraglichen Gruppen unterscheiden kann. Da dr\"angen sich ihre
 Zentralisatoren auf. Der Rest ist eine genaue Analyse der Situation gepaart
 mit Routine. Dann entgeht einem auch die hessesche
 Gruppe\index{hessesche Gruppe}{} nicht, die ich nirgends in der Li\-te\-ra\-tur in der vorliegenden
 Allgemeinheit behandelt fand.


 \newpage
       
 \mychapter{IIII}{Endliche projektive Geometrien}

\noindent
 In diesem Kapitel studieren wir zun\"achst kombinatorische Eigenschaften
 endlicher projektiver Geometrien. Wir werden herausfinden, dass jede
 Kollineation einer projektiven Geometrie ebenso viele Fixpunkte wie
 Fixhyperebenen hat und dass die Anzahl der Punktbahnen einer Kollineationsgruppe
 gleich der Anzahl ihrer Hyperebenenbahnen ist. Wir werden auch sehen, dass es
 keine endlichen el\-lip\-ti\-schen Geometrien gibt, da Polarit\"aten endlicher
 projektiver R\"aume stets absolute Punkte haben. Sp\"ater werden wir
 verschiedene Kenn\-zeich\-nun\-gen der endlichen, desarguesschen projektiven
 Ebenen innerhalb aller endlichen projektiven Ebenen geben. Diese
 Kenn\-zeich\-nun\-gen bedienen sich der Kollineationsgruppen dieser Ebenen.
 Au\-\ss er\-dem geben wir noch verschiedene kombinatorische und
 grup\-pen\-the\-o\-re\-ti\-sche Kennzeichnungen der
 endlichen projektiven R\"aume innerhalb der Klasse der projektiven Blockpl\"ane.
 \par
       Dieses Kapitel sieht beim ersten Hinschauen aus wie ein Sammelsurium der
 unterschiedlichsten S\"atze. Dem aufmerksamen Leser wird aber nicht entgehen,
 dass ihm eine Koh\"arenz innewohnt, die sich bis zu dem Satz 11.3 von N. Ito
 erstreckt, dass also nichts in diesem Kapitel \"uberfl\"ussig ist.

\mysection{1. Endliche Inzidenzstrukturen}

\noindent
 Es sei
 $\Pi$ eine Menge, deren Elemente wir Punkte, und $\Beta$ eine Menge,
 deren Elemente wir Bl\"ocke nennen. Ferner sei $\I\ \subseteq \Pi \times \Beta$.
 Das Tripel $T := (\Pi,\Beta,\I)$ hei\ss t {\it
 Inzidenzstruktur\/}\index{Inzidenzstruktur}{} (siehe Kapitel I, Abschnitt 1). Die
 Elemente von $\I$ nennen wir {\it Fahnen\/}\index{Fahne}{} von $T$. Die
 Inzidenzstruktur $(\Pi,\Beta,\I)$ hei\ss t
 {\it endlich\/},\index{endliche Inzidenzstruktur}{} falls $\Pi$ und
 $\Beta$ endlich sind.
 \par
       Ist $T:=(\Pi,\Beta,\I)$ eine Inzidenzstruktur, so definieren wir die
 Inzidenzstruktur $T^d:= (\Pi^d,\Beta^d,\I^d)$ durch $\Pi^d := \Beta$,
 $\Beta^d := \Pi$ und $\I^d := \{(b,P) \mid (P,b) \in\ \I\}$. Die Inzidenzstruktur
 $T^d$ hei\ss t die zu $T$ {\it duale
 Inzidenzstruktur\/}.\index{duale Inzidenzstruktur}{} Offenbar ist $T^{dd} = T$.
 \par
       Ist $(\Pi,\Beta,\I)$ eine Inzidenzstruktur, so definieren wir auf $\I$
 zwei \"Aqui\-va\-lenz\-re\-la\-tio\-nen $\sim$ und $\approx$ durch $(P,b) \sim (Q,c)$ genau
 dann, wenn $P = Q$, bzw. $(P,b) \approx (Q,c)$ genau dann, wenn $b = c$ ist.
 F\"ur $P \in \Pi$ bezeichne $\I_P$ die durch $P$ bestimmte \"Aquivalenzklasse
 von $\sim$ und f\"ur $b \in \Beta$ bezeichne $\I_b$ die durch $b$
 bestimmte \"Aquivalenzklasse von $\approx$. Es gilt dann
 $$ \I\ = \bigcup_{P \in \Pi} \I_P\ = \bigcup_{b \in B} \I_b $$
 und $I_P \cap I_Q = \emptyset$, falls $P \neq Q$, und
 $\I_b \cap \I_c = \emptyset$, falls $b \neq c$ ist. Ist $r_P$ die Anzahl der
 Elemente in $\I_P$, dh. die Anzahl der Bl\"ocke, die mit $P$
 inzidieren, und ist $k_b$ die Anzahl der mit $b$ inzidierenden
 Punkte, so gilt also der folgende Satz.
 \medskip\noindent
 {\bf 1.1. Satz.} {\it Ist $(\Pi,\Beta,\I)$ eine endliche Inzidenzstruktur, so
 ist}
 $$ \sum_{P \in \Pi} r_P = |\I| = \sum_{b \in \Beta} k_b.$$
 
       Die Inzidenzstruktur $T := (\Pi,\Beta,\I)$ hei\ss t {\it linksseitige
 taktische Konfiguration\/},\index{taktische Konfiguration}{} falls
 es eine ganze Zahl $r$ gibt mit $r_P=r$ f\"ur alle $P \in \Pi$, und sie hei\ss t
 {\it rechtsseitige taktische Konfiguration\/},\index{taktische Konfiguration}{}
 falls es eine ganze Zahl $k$ gibt mit $k_b = k$ f\"ur alle
 $b \in \Beta$. Schlie\ss lich hei\ss t $T$ {\it ta\-kt\-ische
 Konfiguration\/},\index{taktische Konfiguration}{} falls $T$ sowohl
 linksseitige als auch rechtsseitige taktische Konfiguration ist.
 \par
       Ist $T := (\Pi,\Beta,\I)$ eine endliche Inzidenzstruktur, so setzen wir
 $v := |\Pi|$ und $b := |\Beta|$. Damit d\"urfen wir den Buchstaben $b$ nicht
 mehr als Namen f\"ur einen Block verwenden. Die Anzahl der Bl\"ocke $b$ zu
 nennen, versteht sich von selbst. Die Anzahl der Punkte $v$ zu nennen, r\"uhrt
 aus der Statistik,\index{Statistik}{} wo man endliche Inzidenzstrukturen beim
 Planen\index{endliche Inzidenzstruktur}{} von Versuchen benutzt und $v$ dann
 f\"ur das englische {\it variety\/} steht. Der Buchstabe $r$ als
 Anzahl der Bl\"ocke durch einen Punkt steht f\"ur {\it replication\/}. Ist $T$
 linksseitige taktische Konfiguration, so hei\ss en $v$, $r$, $k_c$ mit $c \in B$
 Parameter\index{Parameter}{} von $T$, und ist
 $T$ rechtsseitig taktische Konfiguration, so werden die Zahlen $b$,
 $k$, $r_p$ mit $P \in \Pi$ Parameter von $T$ genannt. Ist $T$
 taktische Konfiguration, so hei\ss en $v$, $b$, $k$, $r$ Parameter von
 $T$. Mit diesen Definitionen folgt aus 1.1 unmittelbar
 \medskip\noindent
 {\bf 1.2. Korollar.} {\it Es sei $T := (\Pi,\Beta,\I)$ eine endliche
 Inzidenzstruktur. Dann gilt:
 \item{a)} Ist $T$ eine linksseitige taktische Konfiguration mit den
 Pa\-ra\-me\-tern $v$, $r$, $k_c$ mit $c \in \Beta$, so ist
 $vr = \sum_{c \in \Beta}k_c$.
 \item{b)} Ist $T$ eine rechtsseitige taktische Konfiguration mit den
 Pa\-ra\-me\-tern
 $b$, $k$, $r_p$ mit $P \in \Pi$ so ist $bk = \sum_{P \in \Pi} r_P$.
 \item{c)} Ist $T$ eine taktische Konfiguration mit den Parametern $v$, $b$, $r$,
 $k$, so ist $vr = bk$.\par}
 \medskip
       Es seien $t$, $v$, $k$ und $\lambda$ nat\"urliche Zahlen. Die
 Inzidenzstruktur $(\Pi,\Beta,\I)$ hei\ss t
 $t$-$(v,k,\lambda)$-{\it Blockplan\/},\index{Blockplan}{} falls gilt:
 \smallskip
       (1) Es ist $v = |\Pi|$.
 \smallskip
       (2) Es ist $k_c = k$ f\"ur alle $c \in \Beta$.
 \smallskip
       (3) Je $t$ verschiedene Punkte aus $\Pi$ inzidieren mit genau $\lambda$
 Bl\"ocken aus $\Beta$.
 \smallskip\noindent
       Statt $t$-$(v,k,\lambda)$-Blockplan werden wir h\"aufig nur $t$-Blockplan
 sagen und $2$-Block\-pl\"ane werden wir meist nur Blockpl\"ane nennen.
 \par
       Die 1-Blockpl\"a\-ne sind genau die taktischen Konfigurationen mit den
 Parametern $v$, $b$, $r = \lambda$ und $k$. Auf Grund von $(2)$ sind alle
 $t$-Blockpl\"ane rechtsseitige taktische Konfigurationen. Wir werden sehen, dass
 sie stets sogar taktische Konfigurationen sind. Bevor wir dies beweisen, noch
 eine Definition. Ist $T := (\Pi,\Beta,\I)$ eine In\-zi\-denz\-struk\-tur und ist
 $P$ ein Punkt von $T$, so definieren wir die {\it abgeleitete Struktur\/} $T_P$
 wie folgt:\index{Abgeleitete Struktur}{}
 \smallskip
       (a) Punkte von $T_P$ sind die von $P$ verschiedenen Punkte von $T$.
 \smallskip
       (b) Bl\"ocke von $T_P$ sind die mit $P$ inzidierenden Bl\"ocke von $T$.
 \smallskip
       (c) Inzidenz in $T_P$ ist gleichbedeutend mit Inzidenz in $T$.
 \smallskip\noindent
       Sind $P_1$, \dots, $P_n$ verschiedene Punkte von $T$, so definieren
 wir $T_{P_1, \dots, P_n}$ rekursiv verm\"oge
 $$ T_{P_1, \dots, P_n} := (T_{P_1, \dots, P_{n-1}})_{P_n}. $$
 \medskip\noindent
 {\bf 1.3. Satz.} {\it Ist $T$ ein $t$-$(v,k,\lambda)$-Blockplan und ist
 $1 \leq s \leq t$, so ist $T$ ein $s$-$(v,k,\lambda_s)$-Blockplan, wobei
 $$ \lambda_s = \lambda {(v - s)(v - s - 1) \cdots (v - t + 1) \over
				    (k - s)(k - s - 1) \cdots (k - t + 1)}$$
 ist. Insbesondere ist $T$ auch eine taktische Konfiguration. Ist
 $b$ die Anzahl der Bl\"ocke von $T$, so ist
 $$ b = \lambda {v(v - 1) \cdots (v - t + 1)
				     \over k(k - 1) \cdots (k - t + 1)}.$$}
 \par
       Beweis. $P_1$, \dots, $P_s$ seien $s$ verschiedene Punkte von $T$.
 Wir bezeichnen mit $\lambda(P_1, \dots, P_s)$ die Anzahl der Bl\"ocke
 von $T_{P_1, \dots, P_s}$. Ist $s = t$, so ist $\lambda(P_1, \dots,
 P_s) = \lambda$, so dass in diesem Falle $\lambda(P_1, \dots, P_s)$
 von der Auswahl der Punkte $P_1, \dots, P_s$ unabh\"angig ist. Es
 sei $s < t$ und es sei bereits gezeigt, dass es ein $\lambda_{s+1}$ gibt mit
 $$ \lambda(X_1, \dots, X_{s+1}) = \lambda_{s+1} $$
 f\"ur jede Wahl der $X_1$, \dots, $X_{s+1}$. Dann ist $T_{P_1, \dots, P_s}$
 eine taktische Konfiguration mit den Parametern $v_s = v - s$, $b_s =
 \lambda (P_1, \dots, P_s)$, $k_s = k - s$ und $r_s = \lambda_{s+1}$.
 Nach 1.2c) ist $v_sr_s = b_sk_s$. Somit ist
 $$ \lambda(P_1, \dots, P_s) = {v - s \over k - s} \lambda_{s+1}. $$
 Damit ist gezeigt, dass $T$ auch ein $s$-$(v,k,\lambda_s)$-Blockplan ist. Nun
 ist $r_P=\lambda_1$ f\"ur alle Punkte $P$ von $T$. Daher ist
 $v\lambda_1 = bk$ nach 1.2c) und folglich
 $$ b = \lambda {v(v - 1) \cdots (v - t+1)\over k(k - 1) \cdots (k - t + 1)}. $$
 Damit ist alles bewiesen.
 \medskip
       Es sei $L$ ein projektiver Verband. Mit $\UR_i(L)$ bezeichnen wir,
 wie schon zuvor, die Menge der Unterr\"aume des Ranges $i$ von
 $L$. Ist $L$ ein endlicher projektiver Verband der Ordnung $q$ und
 des Ranges $n$, so bezeichnen wir wieder die Anzahl der Elemente
 in $\UR_i(L)$ mit $N_i(n,q)$. F\"ur $0 \leq i \leq j \leq n$,
 setzen wir
 $$ L_{i,j} := \bigl(\UR_i(L), \UR_j(L), \leq\bigr). $$
 Es gilt dann:
 \medskip\noindent
 {\bf 1.4. Satz.} {\it Es sei $L$ ein endlicher projektiver Verband
 der Ordnung $q$ und des Ranges $n$. Ferner seien $i$ und $j$
 nat\"urliche Zahlen mit $1 \leq i \leq j \leq n$. Dann gilt:
 \item{a)} $L_{i,j}$ ist eine taktische Konfiguration mit den Parametern
 $v = N_i(n,q)$, $b = N_j(n,q)$, $k = N_i(j,q)$ und $r = N_{j-i} (n-i,q)$.
 \item{b)} Ist $1 < j$, so ist $L_{1,j}$ ein $2$-Blockplan mit $v = N_1(n,q)$,
 $b = N_j(n,q)$, $k = N_1(j,q)$ und $\lambda = N_{j-2}(n-2,q)$.\par}
 \smallskip
       Beweis. Dass $v,b$ und $k$ die angegebenen Werte haben, folgt aus
 der Definition der $N_a(b,q)$. Es sei $U \in \UR_i(L)$ und $\Pi$
 sei das gr\"o\ss te Element von $L$. Dann ist die Anzahl der
 Unterr\"aume vom Range $j$, die $U$ umfassen, gleich der Anzahl
 der Unterr\"aume vom Range $j - i$ in dem Quotienten $\Pi/U$.
 Folglich ist diese Anzahl gleich $N_{j-i}(n - i,q)$. Damit ist a)
 bewiesen.
 \par
       Dass $N_{j-2}(n - 2,q)$ die Anzahl der Unterr\"aume vom Rang $j$
 ist, die zwei verschiedene vorgegebene Punkte umfassen, folgt aus
 a), wenn man nur bemerkt, dass die Summe zweier verschiedener
 Punkte stets eine Gerade ist. Somit gilt auch b).
 \medskip\noindent
 {\bf 1.5. Satz.} {\it Ist $L$ eine projektive Ebene der Ordnung
 $q$, so ist $L_{1,2}$ ein $2$-$(q^2+q+1,q+1,1)$-Blockplan. Ist
 umgekehrt $q$ eine nat\"urliche Zahl mit $q \geq 2$ und ist $E$
 ein $2$-$(q^2 + q + 1,q + 1,1)$-Blockplan, so gibt es eine projektive
 Ebene $L$, so dass $E$ und $L_{1,2}$ isomorph sind.}
 \smallskip
       Beweis. Nach I.7.6 ist $N_1(3,q) = q^2 + q + 1$, $N_1(2,q) = q + 1$ und
 $N_0(1,q) = 1$. Daher ist $L_{1,2}$ auf Grund von 1.4 ein
 2-$(q^2 + q + 1, q + 1,1)$-Blockplan.
 \par
       Ist umgekehrt $E$ ein 2-$(q^2 + q + 1,q + 1,1)$-Blockplan mit $q \geq 2$,
 so tr\"agt jeder Block von $E$ mindestens $3$ Punkte und durch
 zwei verschiedene Punkte geht genau ein Block. Nach 1.3 gehen
 durch jeden Punkt von $E$ genau $(q^2 + q)q^{-1} = q + 1$ Bl\"ocke. Es
 seien $g$ und $h$ zwei verschiedene Bl\"ocke. Weil zwei Bl\"ocke
 h\"ochstens einen Punkt gemeinsam haben, gibt es einen Punkt $P$
 auf $g$, der nicht auf $h$ liegt. Auf $h$ liegen $q + 1$ Punkte.
 Diese sind alle mit $P$ durch genau einen Block verbunden. Diese
 Bl\"ocke sind allesamt verschieden voneinander, da $P$ nicht auf
 $h$ liegt. Da $P$ mit genau $q + 1$ Bl\"ocken inzidiert, ist $g$
 einer dieser Bl\"ocke. Folglich haben $g$ und $h$ einen Punkt
 gemeinsam. Da zwei Geraden, wie gerade gesehen, stets einen Punkt
 gemeinsam haben, gilt auch das Veblen-Young Axiom. Es gibt daher
 einen projektiven Verband $L$, so dass $E$ und $L_{1,2}$ isomorph
 sind. Da zwei verschiedene Geraden von $L$ stets einen
 Schnittpunkt haben, ist der Rang von $L$ h\"ochstens gleich $3$.
 Weil nicht alle Punkte von $L$ kollinear sind, es ist ja
 $q^2 + q + 1 > q + 1$, ist der Rang von $L$ aber auch mindestens gleich $3$.
 Somit ist $L$ eine projektive Ebene.
 \medskip
       Es sei $G$ eine Gruppe und $M$ sei eine Menge und jedes $g \in G$ wirke
 als un\"arer Operator auf $M$. Wir nennen $G$
 {\it Operatorgruppe\/}\index{Operatorgruppe}{} auf $M$, wenn folgende
 Bedingungen erf\"ullt sind.
 \smallskip
       (j) Es ist $x^{gh}=(x^g)^h$ f\"ur alle $x \in M$ und alle $g$, $h \in G$.
 \smallskip
       (ij) Es ist $x^1=x$ f\"ur alle $x \in M$.
 \smallskip
       Ist $x \in M$, so hei\ss t $x^G := \{x^g \mid g \in G\}$
 {\it Bahn\/}\index{Bahn}{} von $G$. Zwei verschiedene Bahnen haben leeren
 Durchschnitt. $G$ hei\ss t {\it transitiv\/},\index{transitiv}{} falls $M$ eine
 Bahn ist. Sind $G$ und $M$ endlich, so gilt, wie auch bei Permutationsgruppen,
 die Gleichung $|G| = |x^G||G_x|$, wobei $G_x := \{g \mid g \in G,\ x^g = x\}$
 ist.
 \par
       Sind $G$ und $M$ endlich, so bezeichnen wir mit $\chi(g)$ die Anzahl der
 Fixelemente von $g \in G$. Die Abbildung $\chi$ hei\ss t {\it
 Permutationscharakter\/}\index{Permutationscharakter}{} von $G$. 
 \medskip\noindent
 {\bf 1.6. Satz.} {\it Ist $G$ eine endliche Operatorgruppe auf der
 endlichen Menge $M$ und ist $a$ die Anzahl der Bahnen von $G$, so
 ist}
 $$ a|G| = \sum_{g\in G} \chi(g). $$
 \par
       Beweis. Wir betrachten die Inzidenzstruktur $(M,G,\I)$ mit $x \I g$
 genau dann, wenn $x^g = x$ ist. Dann ist $k_g = \chi(g)$ und
 $r_x = |G_x|$. Nach 1.1 ist daher $\sum_{g \in G}\chi(g) =
 \sum_{x \in M} |G_x|.$ Es seien $M_1$, \dots, $M_a$ die
 Bahnen von $G$ und $x_i$ sei ein Element aus $M_i$. Dann ist
 $$ \sum_{x \in M} |G_x| = \sum_{i:=1}^a \sum_{x\in M_i}
 |G_x|= \sum_{i:=1}^a |M_i||G_{x_i}| = \sum^a_{i:=1} |G| = a|G|.$$
 Damit ist der Satz bewiesen.
 \medskip\noindent
 {\bf 1.7. Satz.} {\it Ist $G$ eine endliche Operatorgruppe auf der
 endlichen Menge $M$, ist $G$ transitiv auf $M$ und ist $b$ die
 Anzahl der Bahnen von $G_x$ mit $x \in M$, so ist}
 $$b|G|=\sum_{g\in G} \chi(g)^2.$$
 \par
       Beweis. Weil $G$ transitiv ist, sind die Gruppen $G_x$ f\"ur $x \in M$
 alle konjugiert. Somit sind $G_x$ und $b$ unabh\"angig von $x \in M$. Wir
 betrachten die Inzidenzstruktur $(M \times M,G,\I)$, wobei genau dann
 $(x, y) \I g$ gilt, wenn $x^g = x$ und $y^g = y$ ist. Es ist dann
 $k_g = \chi(g)^2$ und $r_{(x,y)} = |G_{x,y}|$.  Nach 1.1 ist daher
 $$ \sum_{g\in G} \chi(g)^2 = \sum_{x \in M} \sum_{y\in M} |G_{x,y}|. $$
 Ersetzt man im Beweise von 1.6 die Gruppe $G$ durch die Gruppe $G_x$, so zeigt
 die letzte Zeile dieses Beweises, dass $\sum_{y\in M}|G_{x,y}| = b|G_x|$
 ist. Also ist
 $$ \sum_{g \in G} \chi(g)^2 = b\sum_{x \in M}\mid G_x| = b|M||G_x| = b|G|, $$
 q. o. o.

\mysection{2. Inzidenzmatrizen}

\noindent
 Es sei 
 $T$ eine endliche Inzidenzstruktur, $P_1$, \dots, $P_v$ seien
 ihre Punkte und $c_1$, \dots, $c_b$ ihre Bl\"ocke. Wir definieren
 die $(v \times b)$-Matrix $A$ durch
 $$ A_{ij} := \cases{1, & falls $P_i \I c_j$ \cr
		     0, & falls $P_i \notI c_j$. \cr} $$
 Die Matrix $A$ hei\ss t {\it Inzidenzmatrix\/}\index{Inzidenzmatrix}{} von $T$.
 Mit $E$ bezeichnen wir die $(v \times v)$-Einheitsmatrix und mit $J$ die
 $(v \times v)$-Matrix, deren s\"amtliche Eintr\"age gleich $1$
 sind. Schlie\ss lich bezeichne $A^t$ die zu $A$ transponierte Matrix.
 \medskip\noindent
 {\bf 2.1. Satz.} {\it Ist $A$ Inzidenzmatrix eines
 $2$-$(v,k,\lambda)$-Blockplanes $T$ und ist $r$ die Anzahl der Bl\"ocke
 durch einen Punkt von $T$, so ist $AA^t = (r - \lambda)E + \lambda J$.}
 \smallskip
 Beweis. Der Punkt $P_i$ inzidiert mit genau $r$ Bl\"ocken. Daher
 sind von den Zahlen $A_{ij}$ f\"ur $j := 1$, \dots, $b$ genau $r$
 gleich $1$, w\"ahrend alle \"ubrigen Null sind. Daher ist einmal $A_{ij}^2 =
 A_{ij}$ und weiter
 $$ r = \sum_{j:=1}^b A_{ij} = \sum_{j:=1}^b A^2_{ij}. $$
 Folglich sind die Diagonalelemente von $AA^t$ alle gleich $r$. Ist $i \neq j$,
 so ist genau dann $A_{ik} A_{jk} \neq 0$, wenn $A_{ik} = A_{jk} = 1$ ist, dh.,
 genau dann, wenn $P_i$, $P_j \I c_k$ gilt.
 Weil zwei verschiedene Punkte mit genau $\lambda$ Bl\"ocken
 inzidieren, ist daher
 $$ \sum_{k:=1}^b A_{ik} A_{jk} = \lambda, $$
 falls nur $i \neq j$ ist. Damit ist 2.1 bewiesen.
 \medskip\noindent
 {\bf 2.2. Satz.} {\it Es seien $a$ und $b$ zwei reelle Zahlen und
 es sei $B := (a - b)E + bJ$ eine $(v \times v)$-Matrix. Dann gilt:
 \item{a)} Die Eigenwerte von $B$ sind $a + (v - 1)b$ und $a - b$. Die
 Vielfachheit von $a + b(v - 1)$ ist $1$ und die Vielfachheit von $a - b$
 ist $v - 1$.
 \item{b)} Es ist $\det(B) = (a + (v - 1)b)(a - b)^{v-1}$.}
 \smallskip
       Beweis. Es sei $E_v$ das $v$-Tupel aus lauter Einsen. Ferner sei
 f\"ur $i := 1$, \dots, $v - 1$ das $v$-Tupel $E_i$ erkl\"art durch
 $$ E_{ij} := \cases{ 1, & falls $j = i$ \cr
		     -1, & falls $j = i + 1$ \cr
		      0 & sonst. \cr} $$
 Dann ist
 $$
 BE_v = (a - b)EE_v + bJE_v = (a - b)E_v + bvE_v  
		= \bigl(a + (v - 1)b\bigr)E_v 
 $$
 und f\"ur $i < v$ ist
 $$ BE_i = (a - b)EE_i + bJE_i = (a - b)E_i. $$
 Weil die Charakteristik von $\R$ Null ist, sind die Vektoren $E_1$, \dots,
 $E_v$ linear unabh\"angig. Folglich gilt a).
 \par
       Die Aussage b) ist eine unmittelbare Folgerung aus a).
 \medskip\noindent
 {\bf 2.3. Korollar.} {\it Ist $A$ die Inzidenzmatrix eines
 $2$-$(v,k,\lambda)$-Blockplanes und ist $k < v$, so ist der Rang von $A$ gleich
 $v$.}
 \smallskip
       Beweis. Ist $r$ die Anzahl der Bl\"ocke durch einen Punkte des
 frag\-li\-chen
 Blockplanes, so gilt nach 1.3 die Gleichung $r(k - 1) = \lambda(v - 1)$. Wegen
 $k < v$ ist daher $\lambda < r$. Aus 2.1 und 2.2 folgt somit 
 $$ \det(AA^t) = \bigl(r + (v - 1)\lambda\bigr)(r - \lambda)^{v-1} \neq 0. $$
 Weil $A$ eine $(v \times b)$- und $AA^t$ eine $(v \times v)$-Matrix ist, ist
 daher
 $$ v \geq \Rg(A) \geq \Rg(AA^t) = v, $$
 woraus die Behauptung folgt.
 \medskip
       Wir ziehen sofort Nutzen aus diesem Korollar.
 \medskip\noindent
 {\bf 2.4. Fishersche Ungleichung.}\index{fishersche Ungleichung}{} {\it Ist $b$
 die Anzahl der Bl\"ocke eines
 $2$-$(v,k,\lambda)$-Blockplanes und ist $k < v$, so ist $v \leq b$.}
 \smallskip
       Beweis. Ist $A$ eine Inzidenzmatrix eines solchen Blockplanes, so
 ist $A$ eine $(v \times b)$-Matrix. Nach 2.3 ist daher $b \geq \Rg(A) = v$,
 q. e. d.
 \medskip\noindent
 {\bf 2.5. Satz.} {\it Ist $T$ eine Inzidenzstruktur mit den Eigenschaften
 \item{a)} $T$ besitzt genau $v$ Punkte und genau $v$ Bl\"ocke,
 \item{b)} Jeder Block von $T$ inzidiert mit genau $k$ Punkten,
 \item{c)} zwei verschiedene Bl\"ocke haben genau $\lambda$ Punkte gemeinsam,
 \item{d)} Es ist $\lambda < k$,\par
 \noindent
 so ist $T$ ein $2$-$(v,k,\lambda)$-Blockplan.}
 \smallskip
       Beweis. Es sei $A$ eine Inzidenzmatrix von $T$. Dann folgt mit b) und c),
 dass
 $$ A^tA = (k - \lambda)E + \lambda J $$
 ist. Nach 2.2 b) ist daher
 $$ \det(A^tA) = \bigl(k + (v - 1)\lambda\bigr)(k - \lambda)^{v-1}. $$
 Weil $\lambda < k$ ist, ist daher $\det(A^tA) \neq 0$. Also ist auch
 $\det(A) \neq 0$.  Folglich ist $A$ regul\"ar und $A^{-1}$ existiert. Aus b)
 folgt, dass $JA = kJ$ ist. Somit ist $k^{-1}J = JA^{-1}$. Weiterhin ist
 $$ JA^tA = (k - \lambda)JE + \lambda J^2 = (k - \lambda +\lambda v)J. $$
 Daher ist
 $$ JA^t = (k - \lambda + \lambda v)JA^{-1} = k^{-1}(k - \lambda +\lambda v)J. $$
 Transponieren liefert
 $$ AJ = k^{-1}(k - \lambda + \lambda v)J $$
 ist. Also ist
 $$ JAJ = k^{-1}v(k-\lambda + \lambda v)J. $$
 Andererseits ist
 $$ JAJ = (JA)J = kJ^2 = kvJ. $$
 Folglich ist $k^{-1} v(k-\lambda + \lambda v) = kv$ und daher
 $k = k^{-1}(k - \lambda + \lambda v)$. Somit ist
 $$ AJ = k^{-1}(k - \lambda + \lambda v)J = kJ = JA. $$
 Hieraus folgt
 $$ AA^t = A(A^tA)A^{-1} = A((k - \lambda)E + \lambda J)a^{-1}
                         = (k - \lambda)E + \lambda J). $$
 Dies besagt schlie\ss lich, dass zwei verschiedene Punkte von $T$ mit genau
 $\lambda$ Bl\"ocken inzidieren, q. e. d.
 \medskip
       Wir nennen einen $2$-Blockplan $T$
 {\it projektiv\/},\index{projektiver Blockplan}{} falls die Anzahl seiner
 Bl\"ocke gleich der Anzahl seiner Punkte
 ist und falls kein Block von $T$ mit allen Punkten von $T$ inzidiert.
 \medskip\noindent
 {\bf 2.6. Korollar} {\it Es sei $T$ ein $2$-$(v,k,\lambda)$-Blockplan mit
 $k < v$. Genau dann ist $T$ projektiv, wenn
 $T^d$ ein $2$-Blockplan ist. Ist $T$ projektiv, so ist $T^d$
 ebenfalls ein $2$-$(v,k,\lambda)$-Blockplan.}
 \smallskip
       Beweis. Es sei $T$ ein $2$-$(v,k,\lambda)$ Blockplan. Dann ist $v \leq b$
 auf Grund der fisherschen Ungleichung. Ist $T^d$ ein $2$-Blockplan, so ist
 ebenfalls auf Grund dieser Ungleichung $b \leq v$. Also ist $v = b$ und $T$ ist
 projektiv.
 \par
       Es sei umgekehrt $T$ projektiv. Aus $v = b$ und $vr = bk$ folgt dann
 $r = k$. Folglich erf\"ullt $T^d$ die Voraussetzungen von 2.5. Also
 ist $T^d$ ein $2$-$(v,k,\lambda)$-Blockplan. Damit ist alles
 bewiesen.
 \medskip
       Die Inzidenzstrukturen
 $L_{1, n-1}$, wobei $L$ ein endlicher projektiver Verband des Ranges $n \geq 3$
 ist, sind Beispiele f\"ur projektive Blockpl\"ane.
 \par
       Das n\"achste Korollar folgt unmittelbar aus 2.6 und 2.1.
 \medskip\noindent
 {\bf 2.7. Korollar} {\it Jede Inzidenzmatrix $A$ eines projektiven
 Blockplanes ist normal, dh. es ist $AA^t = A^tA$.}

\mysection{3. Kollineationen von projektiven Blockpl\"anen}

\noindent
 Es sei 
 $Q$ eine $(n \times n)$-Matrix mit $Q_{ij} \in \{0,1\}$. Die Matrix $Q$
 hei\ss t {\it Permutationsmatrix\/},\index{Permutationsmatrix}{} falls
 $QQ^t = E$ ist. Ist $Q$ eine Permutationsmatrix, so ist $Q$ regul\"ar und es
 gilt $Q^{-1} = Q^t$. Folglich ist auch $Q^tQ = E$. Aus $QQ^t = Q^tQ = E$
 folgt, dass in jeder Zeile und in jeder Spalte von $Q$ genau eine $1$ steht. Ist
 $Q$ eine Permutationsmatrix und definieren wir $\pi$ durch $i^\pi := j$ genau
 dann, wenn $Q_{ij} = 1$ ist, so ist $\pi \in S_n$, falls $S_n$ wieder die
 symmetrische Gruppe\index{symmetrische Gruppe}{} vom
 Grade $n$ bezeichnet. Ist umgekehrt $\pi \in S_n$ und definiert
 man $Q(\pi)$ durch
 $$ Q(\pi)_{ij} := \cases{1, & falls $i^\pi = j$,    \cr
	                  0, & falls $i^\pi \neq j$, \cr} $$
 so ist $Q(\pi)$ eine Permutationsmatrix. Die so definierte
 Abbildung $Q$ ist, wie man sich leicht \"uberzeugt, ein
 Monomorphismus von $S_n$ in die Gruppe der regul\"aren $(n \times
 n)$-Matrizen, so dass also die Per\-mu\-ta\-ti\-ons\-ma\-tri\-zen bez\"uglich
 der Matrizenmultiplikation eine zu $S_n$ isomorphe Gruppe bilden.
 Hieraus folgt, dass jede Permutationsmatrix endliche Ordnung hat,
 so dass die Eigenwerte einer
 Permutationsmatrix\index{Eigenwerte einer Permutationsmatrix}{}
 Einheitswurzeln\index{Einheitswurzel}{} sind.
 \medskip\noindent
 {\bf 3.1. Satz.} {\it Es seien $T$ und $T'$ isomorphe Inzidenzstrukturen. $P_1,
 \dots, P_v$ seien die Punkte und $c_1, \dots, c_b$ seien die Bl\"ocke von $T$
 und $A$ sei die zu dieser Nummerierung der Punkte und Bl\"ocke geh\"orende
 Inzidenzmatrix von $T$. Ferner seien $P'_1, \dots, P'_v$ die Punkte und
 $c'_1, \dots, c'_b$ seien die Bl\"ocke von $T'$ und $A'$ sei die zugeh\"orige
 Inzidenzmatrix von $T'$. Ist $\sigma$ ein Isomorphismus von $T$
 auf $T'$, so definieren wir die $(v \times v)$-Permutationsmatrix
 $Q(\sigma)$ durch
 $$ Q(\sigma)_{ij} := \cases{1, & falls $P_i^\sigma = P'_j$,    \cr
                             0, & falls $P_i^\sigma \neq P'_j$, \cr} $$
 und die $(b \times b)$-Permutationsmatrix $R(\sigma)$ durch
 $$ R(\sigma)_{ij} := \cases{1, & falls $c_i^\sigma = c'_j$,    \cr
                             0, & falls $c_i^\sigma \neq c'_j$. \cr} $$
 Es gilt dann:
 \smallskip
 \item{a)} Die Abbildung $\sigma \to (Q(\sigma),R(\sigma))$ ist eine Bijektion der
 Menge der Isomorphismen von $T$ auf $T'$ auf die Menge der Paare $(Q,R)$ von
 $(v \times v)$- bzw.  $(b \times b)$-Permutationsmatrizen mit $PA' = AQ$.
 \smallskip
 \item{b)} Ist $T = T'$ sowie $P_i = P'_i$ und $c_j = c'_j$ f\"ur alle in Frage
 kommenden $i$ und $j$ und sind $\sigma$ und $\tau$ Automorphismen von $T$, so
 ist $Q(\sigma \tau) = Q(\sigma) Q(\tau)$ und $R(\sigma \tau)=R(\sigma)r(\tau)$.\par}
 \smallskip
       Beweis. a) Es sei $\sigma$ ein Isomorphismus von $T$ auf $T'$.
 Dann ist
 $$ \bigl(Q(\sigma) A'\bigr)_{ij} = \sum_{k:=1}^v Q(\sigma)_{ik} A'_{kj}
				  = A'_{fj}, $$
 wobei $f$ durch $P^\sigma_i = P'_f$ bestimmt ist. Andererseits ist
 $$ \bigl(AR(\sigma)\bigr)_{ij} = \sum^b_{k:=1} A_{ik} R(\sigma)_{kj}
				= A_{ig}, $$
 wobei $g$ durch $c_g^\sigma = c'_j$ bestimmt ist. Die Frage ist also, ob
 $A'_{fj} = A_{ig}$ ist. Nun ist $P^\sigma_i = P'_f$ und $c_g^\sigma = c'_j$.
 Weil $\sigma$ ein Isomorphismus ist, gilt daher $P_i \I c_g$ genau dann, wenn
 $P_i^\sigma \I c_g^\sigma$, dh., es ist $A_{ig} = 1$ genau dann, wenn
 $A'_{fj} = 1$ ist. dies zeigt, dass $Q(\sigma)A' = AR(\sigma)$ ist. Also ist
 $\sigma \to (Q(\sigma),R(\sigma))$ eine Abbildung in die
 fragliche Menge, die offensichtlich auch injektiv ist.
 \par
       Um zu zeigen, dass sie auch surjektiv ist, seien $C$ und $D$
 Permutationsmatrizen mit $CA' = AD$. Wir definieren $\sigma$ durch
 $P^\sigma_i := P'_j$, falls $C_{ij} = 1$ ist, und entsprechend
 $c_i^\sigma := c_j$ falls $D_{ij} = 1$ ist. Weil $C$ und $D$
 Permutationsmatrizen sind, ist $\sigma$ eine bijektive Abbildung
 der Menge der Punkte von $T$ auf die Menge der Punkte von $T'$ und
 der Menge der Bl\"ocke von $T$ auf die Menge der Bl\"ocke von
 $T'$. Es sei $P_i^\sigma = P'_f$ und $c_g^\sigma = c'_j$. Dann ist
 $$ A_{ig} = \sum_{k:=1}^b A_{ik} D_{kj} = \sum^v_{k:=1} C_{ik}A'_{kj}
	   = A'_{fj}. $$
 Hieraus folgt, dass $\sigma$ inzidenztreu ist. Es ist klar, dass dann auch
 $(Q(\sigma), R(\sigma)) = (C,D)$ gilt, womit die Surjektivit\"at nachgewiesen
 ist. Damit ist a) bewiesen.
 \par
       Es seien $\sigma$ und $\tau$ Automorphismen von $T$. Dann ist
 $$ \bigl(Q(\sigma)Q(\tau)\bigr)_{ij} = \sum_{k:=1}^v Q(\sigma)_{ik}Q(\tau)_{kj}
			    = Q(\sigma)_{if}Q(\tau)_{fj} = Q(\tau)_{fj}, $$
 wobei $f$ durch $P_i^\sigma = P_f$ bestimmt ist. Nun gilt genau dann die
 Ungleichung $(Q(\sigma)Q(\tau))_{ij} \neq 0$, wenn $P_f^\tau = P_j$, dh. genau
 dann, wenn $P^{\sigma \tau}_i = P_j$ ist. Also ist
 $Q(\sigma) Q(\tau)_{ij} = Q(\sigma \tau)_{ij}$. Ebenso beweist man die
 Multiplikativit\"at von $R$.
 \medskip\noindent
 {\bf 3.2. Satz.} {\it Jeder Automorphismus eines projektiven Blockplanes hat
 ebenso viele Fixpunkte wie Fixbl\"ocke.}
 \smallskip
       Beweis. Es sei $T$ ein projektiver Blockplan und $A$ sei eine
 Inzidenzmatrix von $T$. Ferner sei $\sigma$ ein Automorphismus von $T$. Dann ist
 $\Spur (Q(\sigma))= \sum_{i:=1}^v Q(\sigma)_{ii}$ die Anzahl der Fixpunkte von
 $\sigma$. Entsprechend ist $\Spur(R(\sigma))$ die Anzahl der Fixbl\"ocke von
 $\sigma$. Nun ist $Q(\sigma)A = AR(\sigma)$ und $A$ ist nach 2.3 regul\"ar. Also
 ist $Q(\sigma) = AR(\sigma)A^{-1}$ und daher $\Spur(Q(\sigma)) =
 \Spur(R(\sigma))$.
 \medskip
       Ist $\sigma$ eine Permutation auf der Menge $M$ und $\tau$ eine solche auf
 der Menge $N$, so hei\ss en $\sigma$ und $\tau$ {\it
 \"ahnlich\/},\index{ahnlich@\"ahnliche Permutationen}{} wenn es eine Bijektion $\beta$
 von $M$ auf $N$ gibt mit $\sigma = \beta\tau\beta^{-1}$. Entsprechend wird die
 \"Ahnlichkeit von Permutationsgruppen
 definiert.
 \medskip\noindent
 {\bf 3.3. Satz.} {\it Ist $Z$ eine zyklische Gruppe von Automorphismen des
 projektiven Blockplanes $T$, so ist die Darstellung von $Z$ als
 Per\-mu\-ta\-ti\-ons\-grup\-pe auf der Menge der Punkte von $T$ \"ahnlich zu
 der Darstellung von $Z$ auf der Menge der Bl\"ocke.}
 \smallskip
       Beweis. Es gen\"ugt zu zeigen, dass f\"ur alle $n$ die Anzahl
 $i_n$ der Punktbahnen der L\"ange $n$ gleich der Anzahl $i'_n$ der
 Blockbahnen der L\"ange $n$ ist.
 \par
       Ist $i_n \neq 0$, so gibt es eine Bahn der L\"ange $n$. Es folgt, dass $n$
 Teiler von $|Z|$ ist. Ist $n$ kein Teiler von $|Z|$, so ist $i_n = 0$ und
 nat\"urlich auch $i'_n = 0$. In diesem Falle ist also $i_n = i'_n$.
 \par
       Es sei also $n$ Teiler von $|Z|$. Ferner sei $\zeta$ ein erzeugendes
 Element von $Z$. Dann ist $i_1$ die Anzahl der Fixpunkte und
 $i'_1$ die Anzahl der Fixbl\"ocke von $\zeta$. Nach 3.2 ist
 $i_1 = i'_1$. Es sei also $n > 1$ und es gelte, dass $i_x = i'_x$ ist
 f\"ur alle $x < n$. Die Anzahl der Fixpunkte von $\zeta^n$ ist
 $$ \sum_{x \, {\rm teilt} \, n} xi_x$$ und die Anzahl der
 Fixbl\"ocke von $\zeta^n$ ist
 $$ \sum_{x \, {\rm teilt} \, n} xi'_x. $$
 Nach 3.2 sind diese beiden Zahlen aber
 gleich. Auf Grund unserer In\-duk\-tions\-an\-nahme folgt daher $ni_n = ni'_n$ und
 weiter $i_n = i'_n$. Damit ist das Korollar bewiesen.
 \medskip
       Dass man in 3.3 auf eine Annahme wie die, dass $Z$ zyklisch ist,
 nicht verzichten kann, sieht man am Beispiel der Gruppe $\E(H,H)$.
 Diese hat nur eine Fixhyperebene, n\"amlich $H$, l\"asst jedoch
 alle Punkte von $H$ fest.
 \par
       Der n\"achste Satz, der ebenfalls eine Folgerung aus 3.2 ist, wird
 h\"aufig Satz von Dembowski--Hughes--Parker
 genannt.\index{Satz von Dembowski--Hughes--Parker}{} Er wurde jedoch schon 1941 von
 Richard Brauer\index{Brauer, R.}{} bewiesen (Brauer 1941, Lemma 3),
 also vierzehn Jahre vor den Publikationen
 von Dembowski, Hughes und Parker. Dieser Satz l\"asst sich auf beliebige
 Blockpl\"ane verallgemeinern. Dies wird in Satz 5.5 geschehen.
 \medskip\noindent
 {\bf 3.4. Satz.} {\it Ist $G$ eine Gruppe von Automorphismen eines projektiven
 Blockplanes $T$, so hat $G$ ebenso viele Punkt- wie Blockbahnen. Insbesondere
 ist $G$ genau dann punkttransitiv, wenn
 $G$ blocktransitiv ist.}
 \smallskip
        Beweis. Es sei $a_1$ die Anzahl der Punkt- und $a_2$ die Anzahl
 der Blockbahnen von $G$. Ferner sei $\chi_1$ der Permutationscharakter von $G$
 aufgefasst als Operatorgruppe auf der Menge der Punkte von $T$ und $\chi_2$
 sei der Permutationscharakter von $G$ aufgefasst als Operatorgruppe auf der
 Menge der Bl\"ocke von $T$. Nach 3.2 ist dann $\chi_1(g) = \chi_2(g)$ f\"ur alle
 $g \in G$. Nach 1.6 ist daher
 $$ a_1|G| = \sum_{g\in G} \chi_1(g) = \sum_{g \in G} \chi_2(g) = a_2|G| $$
 und somit $a_1 = a_2$.
 \medskip\noindent
 {\bf 3.5. Satz.} {\it Es sei $G$ eine punkttransitive oder blocktransitive
 Automorphismengruppe eines projektiven Blockplanes $T$. Ist $P$ ein Punkt und
 $c$ ein Block von $T$, so hat $G_P$ ebenso viele Punktbahnen wie $G_c$
 Blockbahnen. Insbesondere Ist $G$ genau dann zweifach transitiv auf der Menge
 der Punkte von $T$, wenn $G$ zweifach transitiv auf der Menge der Bl\"ocke von
 $T$ ist.}
 \smallskip
       Beweis. Ist $G$ punkttransitiv, so ist $G$ nach 3.4 auch
 blocktransitiv und umgekehrt. Ist $b_1$ die Anzahl der Punktbahnen
 von $G_P$ und $b_2$ die Anzahl der Blockbahnen von $G_c$, so ist
 nach 1.7 und 3.2 also
 $$ b_1|G| = \sum_{g \in G} \chi_1 (g)^2 = \sum_{g \in G} \chi_2(g)^2 = b_2|G| $$
 und daher $b_1 = b_2$, q. o. o.

\mysection{4. Korrelationen von projektiven Blockpl\"anen}

\noindent
 Ist 
 $T := (\Pi,\Lambda,\I)$ eine Inzidenzstruktur und ist $\kappa$ ein
 Isomorphismus von $T$ auf $T^d$, so hei\ss t $\kappa$ {\it
 Korrelation\/}\index{Korrelation}{} von $T$. Eine Korrelation $\kappa$ ist also
 eine Bijektion von $\Pi$ auf $\Lambda$ und von $\Lambda$ auf $\Pi$ mit der
 Eigenschaft, dass genau dann $P \I c$ gilt, wenn $P^\kappa \I^d c^\kappa$, dh.
 genau dann, wenn $c^\kappa \I P^\kappa$ gilt. Ist $\kappa$ eine Korrelation, so
 ist $\kappa^2$ ein Automorphismus von $T$. Ist $\kappa^2 = 1$, so hei\ss t
 $\kappa$ {\it Polarit\"at\/}.\index{Polarit\"at}{}
 \par
       Ist $\kappa$ eine Korrelation\index{Korrelation}{} von $T$ und ist $P$
 ein Punkt von
 $T$ mit $P \I P^\kappa$, so hei\ss t $P$ {\it absoluter
 Punkt\/}\index{absoluter Punkt}{} von $\kappa$. Entsprechend hei\ss t der Block
 $c$ {\it absolut\/}, wenn $c^\kappa \I c$ ist.\index{absoluter Block}{}
 \medskip\noindent
 {\bf 4.1. Satz.} {\it Ist $\kappa$ eine Korrelation des projektiven Blockplanes
 $T$, so hat $\kappa$ ebenso viele absolute Punkte wie absolute Bl\"ocke.}
 \smallskip
       Beweis. Es sei $a$ die Anzahl der absoluten Punkte und $A$ die
 Anzahl der absoluten Bl\"ocke von $\kappa$. Ist $P \I P^\kappa$,
 so ist $(P^\kappa)^\kappa \I P^\kappa$. Folglich ist $P^\kappa$
 ein absoluter Block, woraus folgt, dass $a \leq A$ ist. Ist
 umgekehrt $c^\kappa \I c$, so ist $c^\kappa \I (c^\kappa)^\kappa$,
 so dass $c^\kappa$ ein absoluter Punkt ist. Daher ist $A \leq a$.
 Damit ist alles bewiesen.
 \medskip
       Ist $A$ eine Inzidenzmatrix der Inzidenzstruktur $T$, so ist $A^t$
 eine Inzidenzmatrix von $T^d$. Ist $\kappa$ eine Korrelation von
 $T$, so ist nach 3.1 also $Q(\kappa)A^t = AR(\kappa)$. Nummeriert man
 die Bl\"ocke von $T$ so, dass $c_i = P_i^\kappa$ ist, so ist
 $Q(\kappa) = E$ und daher $AR(\kappa) = A^t$. Nun ist genau dann
 $A_{ii} = 1$, wenn $P_i\I c_i = P_i^\kappa$ ist. Somit ist
 $\Spur(A)$\index{Spur}{} die Anzahl der absoluten Punkte von $\kappa$.
 \medskip\noindent
 {\bf 4.2. Satz.} {\it Es sei $T$ ein projektiver
 $2$-$(v,k,\lambda)$-Blockplan. Ferner sei $k - \lambda = ns^2$ mit
 quadratfreiem $n$. Ist dann $a$ die Anzahl der absoluten Punkte
 einer Korrelation von $T$, so ist $a \equiv \lambda \mod ns$.}
 \smallskip
       Beweis. Wie wir bereits bemerkten, gibt es eine Inzidenzmatrix $A$
 von $T$ und eine Permutationsmatrix $R$ mit $a = \Spur(A)$ und
 $AR = A^t$. Nach 2.7 ist $A$ normal. Es gibt folglich eine unit\"are
 Matrix $U$ mit $A = \bar{U}^tDU$, wobei $D$ die Diagonalmatrix aus
 den Eigenwerten $\mu_1$, \dots, $\mu_v$ von $A$ ist. Wegen
 $$ A^t = \bar{A}^t = \bar{U}^t \bar{D}^t U = \bar{U}^t \bar{D} U $$
 sind $\bar{\mu}_1$, \dots, $\bar{\mu}_v$ die Eigenwerte von $A^t$. Weil
 $U$ unit\"ar ist, ist $U\bar{U}^t = E$. Weil die Koeffizienten von
 $A$ reell sind, ist daher
 $$ AA^t = \bar{U}^tDU\bar{U}^t\bar{D}U = \bar{U}^t D\bar{D} U. $$
 Somit sind die $\mu_i \bar{\mu}_i$ die Eigenwerte von $AA^t$. Die Zeilensummen
 von $A$ sind alle gleich $k$, so dass $k$ ein Eigenwert von $A$ ist. Wir
 d\"urfen annehmen, dass $\mu_1 = k$ ist. Dann ist
 $\mu_1\bar{\mu}_1 = k^2 = k + \lambda(v - 1)$. Aus 2.2 a) folgt dann, dass
 $\mu_i\bar{\mu}_i = k - \lambda$ ist f\"ur $i := 2$, \dots, $v$. Setze
 $m := k - \lambda$. Dann ist $\mu_i = \epsilon_i \sqrt{m}$ und
 $\epsilon_i \bar{\epsilon}_i = 1$. Weil $A$ regul\"ar ist, folgt
 $R = A^{-1} A^t$. Hieraus folgt, dass die $\bar{\mu}_i \mu_i^{-1}$
 die Eigenwerte von $R$ sind. Nun ist $\bar{\mu}_1 \mu_1^{-1} = 1$ und
 $\bar{\mu}_i \mu_i^{-1} = \epsilon^{-2}_i$ f\"ur $i: = 2$, \dots, $v$.
 Weil $R$ eine Permutationsmatrix ist, sind die Eigenwerte von $R$
 Einheitswurzeln. Hieraus folgt wiederum, dass die $\epsilon_i$
 Einheitswurzeln sind. Folglich ist ($\sum_{i:=1}^v \epsilon_i)^2$
 eine ganz algebraische Zahl. Nun ist
 $$ (a - k)^2 = \bigl(\Spur(A) - k\bigr)^2 = \biggl(\sum_{i:=2}^v \mu_i\biggr)^2
	       = m\biggl(\sum^v_{i:=2} \epsilon_i\bigg)^2. $$
 Somit ist $(\sum_{i:=2}^v \epsilon_i)^2$ rational
 und als ganz algebraische Zahl sogar ganz rational. Also ist
 $(a - k)^2 \equiv 0 \mod m$ und folglich $a \equiv k \mod ns$.
 Wegen $k = m + \lambda$ folgt hieraus $a \equiv \lambda \mod ns$,
 q. o. o.
 \medskip
       Die Begriffe {\it absoluter Punkt\/} und {\it absolute Hyperebene\/} einer
 Korrelation eines projektiven Verbandes wurden in Abschnitt 4 von Kapitel II
 definiert.
 \medskip\noindent
 {\bf 4.3. Korollar} {\it Es sei $L$ ein endlicher projektiver Verband mit
 $r:= \Rg (L) \geq 3$.
 \item{a)} Ist $\kappa$ eine Korrelation von $L_{1,r-1}$, so hat $\kappa$
 wenigstens einen absoluten Punkt.
 \item{b)} Ist $\kappa$ eine Korrelation von $L$, so hat $\kappa$
 wenigstens einen absoluten Punkt.}
 \smallskip
       Beweis. Da jede Korrelation von $L$ eine Korrelation in
 $L_{1,r-1}$ induziert, gen\"ugt es, a) zu beweisen.
 \par
       Es sei $q$ die Ordnung von $L$. Nach 1.4 b) ist dann
 $k = \sum^{r-2}_{i:=0} q^i$ und $\lambda = \sum^{r-3}_{i:=0} q^i$.
 Daher ist $m = k - \lambda = q^{r-2}$. Ist $m = ns^2$ mit quadratfreiem
 $n$ und ist $p$ ein Primteiler von $q$, so ist $p$ ein Teiler von
 $ns$. Ferner ist $\lambda \equiv 1 \mod p$, so dass nach 4.2 die
 Kongruenz $a \equiv 1 \mod p$ gilt. Hieraus folgt $a \geq 1$, q. e. d.
 \smallskip
       Es sei $V$ ein Vektorraum \"uber dem kommutativen K\"orper $K$.
 Die Abbildung $Q$ von $V$ in $K$ hei\ss t {\it quadratische
 Form\/}\index{quadratische Form}{} auf $V$, falls gilt: 
 \item{(1)} Es ist $Q(vk) = Q(v)k^2$ f\"ur alle $v \in V$ und alle $k \in K$.
 \item{(2)} Die durch $f(u,v) := Q(u + v) - Q(u) - Q(v)$ definierte
 Abbildung $f$ von $V \times V$ in $K$ ist bilinear.
 \medskip
       Als Anwendung von 4.3 beweisen wir den folgenden, rein
 al\-ge\-bra\-i\-schen Satz
 \"uber quadratische Formen \"uber endlichen K\"orpern.
 \medskip\noindent
 {\bf 4.4. Satz.} {\it Es sei $V$ ein Vektorraum \"uber $\GF(q)$ und $Q$ sei eine
 quadratische Form auf $V$. Ist $U \in L(V)$ und ist $3 \leq \Rg_K(U) <
 \infty$, so gibt es ein $u \in U$ mit $u \neq 0$ und $Q(u) = 0$.}
 \smallskip
       Beweis. Wir betrachten zun\"achst den Fall, dass die Charakteristik von
 $\GF(q)$ gleich $2$ ist. Ferner sei $u \in U$ und $Q(u) \neq 0$. Wegen
 $0 = f(0,0) =  Q(0) - 2Q(0) = Q(0)$ ist $u \neq 0$. Definiert man $g$ durch
 $g(x) := f(u,x)$ f\"ur alle $x \in U$, so ist $g$ nach (2) eine lineare
 Abbildung von $U$ in $\GF(q)$. Nun ist $\Rg(U) \geq 3$ und daher
 $\Rg(\Kern(g)) \geq 2$. Es gibt folglich ein $x \in U$ mit $x \notin u\GF(q)$
 und $g(x) = 0$. Es sei $k \in \GF(q)$. Dann ist
 $$ 0 = g(x)k = f(x,u)k = f(x,uk) = Q(x + uk) - Q(x) - Q(uk) $$
 und daher
 $$ Q(x + uk) = Q(x) + Q(u)k^2. $$
 Nun ist $\GF(q) = \GF(q)^2$, da die Charakteristik von $\GF(q)$ ja $2$ ist. Weil
 au\ss erdem $Q(u) \neq 0$ ist, gibt es ein $k \in \GF(q)$ mit $Q(x + uk) = 0$.
 W\"are $x + uk = 0$, so w\"are $x \in u\GF(q)$, was nicht der Fall ist.
 Damit ist 4.4 f\"ur den Fall, dass $\GF(q)$ die Charakteristik 2 hat, bewiesen.
 \par
       Die Charakteristik von $\GF(q)$ sei nun von $2$ verschieden. Wir
 nehmen an, dass $f(u,u) \neq 0$ sei f\"ur alle $u$ mit $0 \neq u \in U$. Die
 Einschr\"ankung von $f$ auf $U$ ist dann eine nicht
 ausgeartete Bilinearform auf $V$ und definiert somit nach II.8.4
 eine Korrelation $\kappa$ in $L(U)$. Wegen $\Rg_K(U) \geq 3$
 gibt es daher nach 4.3 einen absoluten Punkt $wU$ in $L(U)$. Es
 folgt $0 \neq w$ und $f(w,w) = 0$. Es gibt also doch ein $u \in U$
 mit $u \neq 0$ und $f(u,u) = 0$. Es folgt
 $$ 0 = f(u,u) = Q(u + u) - Q(u)- Q(u) = 4Q(u) - 2Q(u) = 2Q(u) $$
 und daher $Q(u) = 0$, da die Charakteristik von $\GF(q)$ ja ungleich $2$ ist.
 Damit ist alles bewiesen.
 \medskip
       Im Falle einer Polarit\"at kann man mehr \"uber die Anzahl der
 absoluten Punkte aussagen. Es gilt n\"amlich:
 \medskip\noindent
 {\bf 4.5. Satz.} {\it Ist $T$ ein projektiver $2$-$(v,k,\lambda)$-Blockplan und
 ist $\pi$ eine Polarit\"at von $T$, so ist die Anzahl der absoluten Punkte von
 $\pi$ gleich $k + s\sqrt{k - \lambda}$, wobei $s$ eine geeignete ganze Zahl
 ist.}
 \smallskip
       Beweis. $P_1$, \dots, $P_v$ seien die Punkte von $T$. Die Bl\"ocke
 von $T$ nummerieren wir so, dass $c_i = P_i^\pi$ ist. F\"ur die zu
 dieser Nummerierung geh\"orende Inzidenzmatrix gilt dann wegen
 $\pi^2 = 1$, dass $A = A^t$ ist. Ferner ist $\Spur(A)$ wieder die
 Anzahl der absoluten Punkte von $\pi$. Wie wir bereits wissen, ist
 $k$ ein Eigenwert von $A$ und die \"ubrigen Eigenwerte von $A$
 sind von der Form $\epsilon_i \sqrt{k - \lambda}$ mit
 $|\epsilon_i| = 1$. Weil $A$ symmetrisch ist, sind die Eigenwerte
 von $A$ alle reell. Daher ist $\epsilon_i \in \{1, -1\}$. Hieraus
 folgt alles weitere, da die Spur von $A$ ja gleich der Summe
 \"uber die Eigenwerte von $A$ ist.
 \medskip\noindent
 {\bf 4.6. Korollar.} {\it Ist $T$ ein projektiver
 $2$-$(v,k,\lambda)$-Blockplan, ist $k - \lambda$ kein Quadrat und ist
 $\pi$ eine Polarit\"at von $T$, so ist die Anzahl der absoluten
 Punkte von $\pi$ gleich $k$.}

\mysection{5. Taktische Zerlegungen}

\noindent
 Es sei 
 $A$ eine $(v \times b)$-Matrix mit Koeffizienten in einem kommutativen
 K\"orper $K$. Die Zeilen seien dabei mit den Elementen aus $\{1, \dots, v\}$ und
 die Spalten mit den Elementen aus $\{1, \dots, b\}$ indiziert. $A$ sei also das,
 was ich anderswo ein rechteckiges Schema nannte. Ferner sei $\{Z_1, \dots,
 Z_t\}$ eine Partition von $\{1, \dots, v\}$ mit nicht leeren $Z_i$
 und $\{S_1, \dots, S_{t'}\}$ seien eine Partition von $\{1, \dots, b\}$ mit
 nicht leeren $S_i$. Mit $A^{ij}$ bezeichnen wir die Einschr\"ankung der
 Abbildung $A$ auf $Z_i \times S_j$. Wir nennen $(Z_1, \dots, Z_t; S_1, \dots,
 S_{t'})$ {\it linksseitige tak\-ti\-sche Zerlegung\/}\index{taktische Zerlegung}{}
 von $A$, falls alle Matrizen $A^{ij}$ konstante Zeilensummen
 $C_{ij}$ haben. Die Matrix $C$ hei\ss t {\it assoziierte Zeilensummenmatrix\/}
 der\index{assoziierte Zeilensummenmatrix}{} linksseitigen taktischen
 Zerlegung $(Z_1, \dots, Z_t; S_1, \dots, S_{t'})$.
 \par
       Wir nennen $(Z_1, \dots, Z_t;S_1, \dots, S_{t'})$ {\it rechtsseitige
 taktische Zer\-le\-gung\/}\index{taktische Zerlegung}{} von $A$,
 falls die Matrizen $A^{ij}$ konstante Spaltensummen $D_{ij}$ haben.
 $D$ hei\ss t dann sinngem\"a\ss\ {\it assoziierte Spaltensummenmatrix\/}
 der\index{assoziierte Spaltensummenmatrix}{} rechtsseitigen taktischen Zerlegung
 $(Z_1, \dots, Z_t;$ $S_1, \dots, S_{t'})$.
 \par
       Schlie\ss lich hei\ss t $(Z_1, \dots, Z_t; S_1, \dots, S_t')$
 {\it taktische Zerlegung\/}\index{taktische Zerlegung}{} der Matrix $A$, falls
 $(Z_1, \dots, Z_t; S_1, \dots, S_{t'})$ sowohl rechtsseitige als auch
 linksseitige taktische Zerlegung der Matrix $A$ ist.
 \par
       Der Leser fragt sich vielleicht, wie man auf solch eine Definition kommt.
 Das ist hier ausnahmsweise einmal zu beantworten. Nichts ist nat\"urlicher, als
 die Struktur der Punkt- und Blockbahnen von Kollineationsgruppen von
 Inzidenzstrukturen zu untersuchen. Nimmt man die Inzidenzstruktur aus einer
 Punkt- und einer Blockbahn, so sieht man unmittelbar, dass sie eine taktische
 Konfiguration ist. Verfolgt man nun dies f\"ur alle m\"oglichen Paarungen in
 einer Inzidenzmatrix der fraglichen Inzidenzstruktur, so erh\"alt man eine
 tak\-ti\-sche Zerlegung der Inzidenzmatrix. Von da ist es dann klar, dass
 man die Eigenschaften \anff konstante Zeilensummen`` und \anff konstante
 Spaltensummen`` noch separiert. Dies ist in ganz groben Z\"ugen auch die
 historische Entwicklung, wie sie sich in den Publikationen widerspiegelt.
 \medskip\noindent
 {\bf 5.1. Satz.} {\it Es sei $A$ eine $(v \times b)$-Matrix mit Koeffizienten in
 einem kommutativen K\"orper. Ferner sei $\{Z_1, \dots, Z_t\}$ eine Partition von
 $\{1, \dots, v\}$ und $\{S_1, \dots, S_{t'}\}$ eine Partition von
 $\{1, \dots, b\}$. Dann gilt:
 \item{a)} Ist $(Z_1, \dots, Z_t;S_1, \dots, S_t')$ eine rechtsseitige taktische
 Zerlegung von $A$, ist $\rho$ der Rang von $A$ und $\rho_D$ der Rang der
 assoziierten Spaltensummenmatrix, so ist
 $$ t \leq \rho_D + v - \rho. $$
 Insbesondere gilt
 $$ t \leq t' + v - \rho. $$
 \item{b)} Ist $(Z_1, \dots, Z_t; S_1, \dots, S_{t'})$ eine linksseitige taktische
 Zerlegung von $A$, ist $\rho$ der Rang von $A$ und $\rho_C$ der Rang der
 assoziierten Zeilensummenmatrix, so ist
 $$ t' \leq \rho_C + b - \rho. $$
 Insbesondere gilt
 $$ t' \leq t + b - \rho. $$\par}
 \par
       Beweis. Ist $(Z_1, \dots, Z_t; S_1, \dots, S_{t'})$ eine
 linksseitige taktische Zerlegung von $A$ und ist $C$ die
 assoziierte Zeilensummentmatrix, so ist $(S_1, \dots, S_{t'}; Z_1,
 \dots, Z_t)$ eine rechtsseitige taktische Zerlegung von $A^t$ und
 $C^t$ ist die assoziierte Spaltensummenmatrix dieser Zerlegung.
 Weil sich die R\"ange beim Transponieren nicht \"andern, gen\"ugt
 es daher, 5.1a) zu beweisen.
 \par
       Da $\rho$ der Rang von $A$ ist, gibt es $\rho$ linear unabh\"angige
 Zeilenvektoren von $A$. Die Indizes der \"ubrigen $v-\rho$ Zeilenvektoren liegen
 in h\"ochstens $v-\rho$ Klassen $Z_1$, \dots, $Z_t$. Es gibt also $t - v + \rho$
 Klassen in $\{Z_1, \dots, Z_t\}$, so dass die Zeilenvektoren mit Indizes in
 der Ver\-ei\-ni\-gung dieser $v - t + \rho$ Klassen linear unabh\"angig sind.
 Wir d\"urfen annehmen, dass dies die Klassen $Z_1$, \dots, $Z_{v-t+\rho}$ sind.
 Es seien $z_1$, \dots, $z_{v-t+\rho}$ die ersten $t - v + \rho$
 Zeilenvektoren der assoziierten Spaltensummenmatrix. Ferner sei
 $\sum_{i:=1}^{t-v+\rho} z_i\lambda_i = 0$ und es sei $a_l$
 der $l$-te Zeilenvektor von $A$. Dann folgt aus der Definition der
 rechtsseitigen taktischen Zerlegung, dass
 $$ \sum_{i:=1}^{t-v+\rho} \sum_{l\in Z_i} a_l \lambda_i = 0 $$
 ist. Weil die Zeilenvektoren $a_l$ mit $l \in \bigcup_{i:=1}^{t-v+\rho} Z_i$
 linear unabh\"angig sind, ist daher $\lambda_1 = \dots = \lambda_{t-v+\rho}=0$.
 Also ist $t - v + \rho \leq \rho_D$ und folglich $t \leq \rho_D + v - \rho$.
 Schlie\ss lich ist $\rho_D \leq t'$, da ja $D$ eine $(t \times t')$-Matrix ist.
 Damit ist 5.1 bewiesen.
 \medskip
       Und nun zur Anwendung dieses Satzes auf endliche Inzidenzstrukturen.
       Es sei $T := (\Pi, \Lambda, \I)$ eine endliche Inzidenzstruktur.
 Ferner sei $\{\Pi_1, \dots, \Pi_t\}$ eine Partition von $\Pi$ und
 $\{\Lambda_1, \dots, \Lambda_{t'}\}$ eine solche von $\Lambda$.
 Wir nennen dann $(\Pi_1, \dots, \Pi_t; \Lambda_1, \dots, \Lambda_{t'})$
 {\it linksseitige taktische Zer\-le\-gung\/}\index{taktische Zerlegung}{} von $T$, falls die
 Inzidenzstrukturen $T^{ij} := (\Pi_i, \Lambda_j, \I_{ij})$ mit
 $$ \I_{ij} := \I \cap \,(\Pi_i \times \Lambda_j) $$
 f\"ur alle $i$ und $j$ linksseitige taktische Konfigurationen sind. Sind die
 $T^{ij}$ allesamt rechtsseitige taktische Konfigurationen, so 
 hei\ss t
 die Zer\-le\-gung $(\Pi_1, \dots, \Lambda_{t'})$ {\it rechtsseitige
 taktische Zerlegung\/}. Eine Zerlegung, die sowohl rechtsseitig als
 auch linksseitig ist, hei\ss t
 {\it taktische Zerlegung\/}
 schlechthin.\index{taktische Zerlegung}{}
 \medskip\noindent
 {\bf 5.2. Satz.} {\it Es sei $T := (\Pi,\Lambda,\I)$ eine endliche
 Inzidenzstruktur. $P_1, \dots, P_v$ seien ihre Punkte und $c_1, \dots, c_b$
 ihre Bl\"ocke und $A$ sei die zu dieser Nummerierung der Punkte und Bl\"ocke
 geh\"orende Inzidenzmatrix. Schlie\ss lich sei $\{\Pi_1, \dots, \Pi_t\}$ eine
 Partition von $\Pi$ und $\{\Lambda_1, \dots, \Lambda_{t'}\}$ eine Partition von
 $\Lambda$. Wir definieren eine Partition
 $ \{Z_1, \dots, Z_t\} $
 von $\{1, \dots, v\}$ sowie eine Partition
 $ \{S_1, \dots, S_{t'}\} $
 von $\{1, \dots, b\}$ durch
 $$ Z_i := \bigl\{i \mid i \in \{1, \dots, v\},\ P_i \in \Pi_i\bigr\} $$
 bzw.
 $$ S_i := \bigl\{i \mid i \in \{1, \dots, b\},\ c_i \in \Lambda_i\bigr\}. $$
 Es gilt dann:
 \smallskip
 \item{a)} Genau dann ist $(\Pi_1, \dots, \Pi_t; \Lambda_1, \dots, \Lambda_{t'})$
 eine linksseitige tak\-ti\-sche Zerlegung von $T$, wenn
 $(Z_1, \dots, Z_t; S_1, \dots, S_{t'})$ eine linksseitige taktische Zerlegung
 von $A$ ist.
 \item{b)} Genau dann ist $(\Pi_1, \dots, \Pi_t; \Lambda_1, \dots, \Lambda_{t'})$
 eine rechtsseitige taktische Zerlegung von $T$,
 wenn $(Z_1, \dots, Z_t; S_1, \dots, S_{t'})$ eine rechts\-sei\-ti\-ge
 taktische Zerlegung von $A$ ist.
 \item{c)} Genau dann ist $(\Pi_1, \dots, \Pi_t; \Lambda_1, \dots, \Lambda_{t'})$
 eine taktische Zerlegung von $T$, wenn $(Z_1, \dots, Z_t; S_1, \dots, S_{t'})$
 eine taktische Zerlegung von $A$ ist.\par}
 \medskip
       Zum Beweise von 5.2 hat man nur zu bemerken, dass die $A^{ij}$
 Inzidenzmatrizen von $T^{ij}$ sind.
 \medskip\goodbreak
 Aus 5.1 a), 5.2 und 2.3 folgt unmittelbar
 \medskip\noindent
 {\bf 5.3. Korollar.} {\it Ist $T$ ein $2$-$(v,k,\lambda)$-Blockplan mit $k < v$
 und besitzt $T$ eine taktische Zerlegung mit $t$ Punkt- und $t'$ Blockklassen
 und ist $b$ die Anzahl der Bl\"ocke von $T$, so ist
 $$ t \leq t' \leq t + b - v. $$
 Ist $T$ projektiv, so ist $t = t'$.}
 \medskip
       Dass taktische Zerlegungen projektiver Blockpl\"ane stets ebenso viele
 Punkt- wie Blockklassen haben, wurde zuerst von P.
 Dembowski\index{Dembowski, P.}{} bewiesen. Es ist eines der ersten Resultate
 \"uber taktische Zerlegungen \"uberhaupt. Die fishersche
 Ungleichung\index{fishersche Ungleichung}{} ist ein Spezialfall von 5.3, da die
 Zerlegung der Punkt- bzw. Blockmenge von $T$ in einelementige Teilmengen eine
 taktische Zerlegung von $T$ ist.
 \par
       Wie oben schon gesagt, war die im folgenden Satz ausgesprochene
 Bemerkung Ausgangspunkt f\"ur alle Untersuchungen \"uber taktische
 Zerlegungen.
 \medskip\noindent
 {\bf 5.4. Satz.} {\it Ist $G$ eine Gruppe von Automorphismen einer endlichen
 Inzidenzstruktur $T := (\Pi, \Lambda, \I)$, sind $\Pi_1$, \dots, $\Pi_t$  die
 Punktbahnen und $\Lambda_1$, \dots, $\Lambda_{t'}$ die Blockbahnen von $G$, so
 ist
 $$ (\Pi_1, \dots, \Pi_t; \Lambda_1, \dots, \Lambda_{t'}) $$
 eine taktische Zerlegung von $T$.}
 \smallskip
       Beweis. Es ist zu zeigen, dass $T^{ij} := (\Pi_i, \Lambda_j, \I_{ij})$
 eine taktische Konfiguration ist. Dies ist aber trivial, da $G$ in $T^{ij}$ eine
 Gruppe von Automorphismen induziert, die sowohl auf $\Pi_i$ als auch auf
 $\Lambda_j$ trans\-i\-tiv operiert.
 \medskip
       Mit 5.4 und 5.3 folgt schlie\ss lich die schon angek\"undigte
 Verallgemeinerung von 3.4.
 \medskip\noindent
 {\bf 5.5. Korollar.} {\it Ist $G$ eine Gruppe von Automorphismen
 eines end\-li\-chen Blockplanes $T$ und ist die Anzahl der Punktbahnen
 von $G$ gleich $t$ und die Anzahl der Blockbahnen gleich $t'$, so
 ist
 $$ t \leq t' \leq t + b - v. $$
 Dabei ist $v$ die Anzahl der Punkte und $b$ die Anzahl der Bl\"ocke von $T$.}
 \medskip
       Das Korollar 3.4 folgt auch aus 5.5, da bei projektiven Blockpl\"anen ja
 $v = b$ und demzufolge $t = t'$ ist.

\mysection{6. Endliche desarguessche projektive Ebenen}

\noindent
 Ist 
 $E$ eine projektive Ebene, so hei\ss t $E$ {\it
 Moufangebene\/},\index{Moufangebene}{} falls $E$ be\-z\"ug\-lich jeder ihrer
 Geraden Translationsebene ist (siehe Kapitel II, Ende von Abschnitt 1). Alle
 desarguesschen Ebenen sind somit auch Mou\-fang\-e\-be\-nen, aber nicht
 alle Moufangebenen sind desarguessch (siehe etwa Pickert 1955, S.
 176--178). F\"ur endliche projektive Ebenen gilt jedoch
 \goodbreak
 \medskip\noindent
 {\bf 6.1. Satz.} {\it Alle endlichen Moufangebenen sind desarguessch.}
 \medskip
       Dieser Satz ist die Grundlage f\"ur viele Charakterisierungen der
 endlichen desarguesschen projektiven Ebenen. Sein Beweis bedarf einiger
 Vorbereitung. Der Leser beachte, dass bei den S\"atzen 6.2 bis 6.14 nicht die
 Endlichkeit der betrachteten Moufangebenen und Alternativk\"orper vorausgesetzt
 wird.
 \medskip\noindent
 {\bf 6.2. Satz.} {\it Ist $M$ eine Moufangebene, so ist die von
 allen Elationen erzeugte Kollineationsgruppe von $M$ auf der Menge
 der Tripel $(P, Q, G)$, wobei $P$ und $Q$ zwei verschiedene Punkte
 und $G$ eine Ge\-ra\-de von $M$ ist, die weder $P$ noch $Q$ enth\"alt,
 transitiv.}
 \smallskip
       Beweis. Es seien $G$ und $G'$ Geraden von $M$ und $P$, $Q$ seien
 zwei verschiedene Punkte von $M_G$ und $P', Q'$ seien zwei
 verschiedene Punkte von $M_{G'}$. Weil $\Rg(M) = 3$ ist, gibt es
 einen Punkt $R \leq G \cap G'$. Ferner gibt es eine von $G$ und
 $G'$ verschiedene Gerade $H$ durch $R$. Es sei $C$ ein von $R$
 verschiedener Punkt auf $H$ und $I$ eine von $H$ verschiedene
 Gerade durch $C$. Definiere $X$ und $Y$ durch $X := G \cap I$ und
 $Y := G' \cap I$. Dann sind $X$ und $Y$ von $C$ verschiedene
 Punkte. Weil $M$ eine Moufangebene ist, gibt es ein $\tau \in \E(R,H)$ mit
 $X^\tau = Y$. Es folgt
 $$ G^\tau = (X + R)^\tau = Y + R = G'. $$
 Wir d\"urfen daher annehmen, dass $G = G'$ ist.
 Weil $\E(G)$ auf der Menge der Punkte von $E_G$ transitiv
 operiert, d\"urfen wir weiterhin annehmen, dass auch $P = P'$ ist.
 Wir haben nun zwei F\"alle zu unterscheiden.
 \par
       1. Fall: $P$, $Q$ und $Q'$ sind nicht kollinear. Dann ist
 insbesondere $Q \neq Q'$. Daher ist $Q + Q'$ eine Gerade, die wegen
 $Q \not\leq G$ von $G$ verschieden ist. Somit ist $C := (Q + Q') \cap G$ ein
 Punkt. Weil $P \not\leq G$ ist, ist auch $P + C$ eine Gerade,
 die von $G$ verschieden ist. Sie ist aber auch von $Q + Q'$
 verschieden, da $P$, $Q$, $Q'$ ja nicht kollinear sind. Es gibt
 folglich ein $\tau \in \E(C,P + C)$ mit $Q^\tau = Q'$. Weil $P^\tau = P$ und
 wegen $C \leq G$ auch $G^\tau = G$ ist, ist
 $(P^\tau,Q^\tau,G^\tau) = (P,Q',G) = (P',Q',G')$.
 \par
       2. Fall: $P$, $Q$ und $Q'$ sind kollinear. In diesem Fall w\"ahlen
 wir einen Punkt $R$ mit $R \not\leq P + Q$. Dann sind weder $P$, $Q$, $R$
 noch $P$, $Q'$, $R$ kollinear. Wie wir im Falle 1 gesehen haben, gibt
 es dann Elationen $\sigma$ und $\tau$ mit $P^\sigma = P$, $Q^\sigma = R$,
 $G^\sigma = G$ und $P^\tau = P$, $R^\tau = Q'$, $G^\tau = G$. Dann ist
 $\sigma \tau$ eine Kollineation aus der von allen Elationen erzeugten Gruppe,
 die $(P,Q,G)$ auf $(P',Q',G')$ abbildet.  Damit ist 6.2 bewiesen.
 \medskip\noindent
 {\bf 6.3. Satz.} {\it Ist $M$ eine Moufangebene, ist $G$ eine Gerade von $M$ und
 sind $P$ und $Q$ zwei verschiedene Punkte von $M_G$, so gibt es eine 
 involutorische Perspektivit\"at\index{involutorische Perspektivit\"at}{} mit der
 Achse $G$, die $P$ und $Q$ vertauscht.}
 \smallskip
       Beweis. Wir zeigen zun\"achst, dass $M$ eine involutorische
 Perspektivit\"at mit der Achse $G$ besitzt. Dazu d\"urfen wir annehmen, dass
 alle Elationen mit der Achse $G$ eine von $2$ verschiedene Ordnung haben. Der
 Kern $K(G)$ von $\E(G)$, der nach II.2.1 ein K\"orper ist, hat dann von 2
 verschiedene Charakteristik. Hieraus folgt, dass die multiplikative Gruppe von
 $K$ ein Element der Ordnung 2 enth\"alt, n\"amlich $-1$. Nach II.2.1 gibt es
 folglich eine involutorische Streckung mit der Achse $G$.
 \par
       Es sei nun $\rho$ eine involutorische Perspektivit\"at mit der
 Achse $G$. Es gibt dann zwei verschiedene Punkte $X$ und $Y$ von
 $M_G$ mit $X^\rho = Y$ und $Y^\rho  = X$. Nach 6.2 gibt es eine
 Kollineation $\gamma$ mit $X^\gamma = P, Y^\gamma = Q$ und
 $G^\gamma = G$. Die Kollineation $\gamma^{-1} \rho \gamma$ ist dann
 eine involutorische Perspektivit\"at mit der Achse $G$, die $P$
 mit $Q$ vertauscht.
 
 \bigskip

       $M$ sei wieder eine Moufangebene. Ferner sei $(U,V,O,E)$ ein Rahmen von
 $M$. Die Menge der Punkte auf $O + V$, die von $V$ verschieden sind, bezeichnen
 wir mit $K$. Au\ss erdem bezeichnen wir die Elemente von $K$ mit kleinen
 lateinischen Buchstaben.  Schlie\ss lich schreiben wir $PQ$ f\"ur die
 Verbindungsgerade $P + Q$ zweier verschiedener Punkte $P$ und $Q$, da wir dem
 Pluszeichen in diesem Abschnitt eine andere Bedeutung unterlegen werden, wie
 gleich erl\"autert wird. Ist $P$ ein Punkt von $M_{UV}$, so setzen wir
 $$ x_P := (PV \cap OE) U \cap OV $$
 und
 $$ y_P := PU \cap OV. $$
 Ist umgekehrt $(x,y) \in K \times K$, so definieren wir $\tau(x,y)$ durch
 $$ \tau(x,y) := (xU \cap OE)V \cap yU. $$
 Offensichtlich gilt $\tau(x_P, y_P) = P$ und $x_{\tau(x,y)} = x$ sowie
 $y_\tau(x,y) = y$.  Folglich ist $\tau$ eine Bijektion von $K \times K$ auf die
 Menge der Punkte von $M_{UV}$. Ferner gilt $\tau^{-1}(P) = (x_P,y_P)$.
 \par 
       Ist $G$ eine Gerade von $M_{UV}$, die nicht durch $V$ geht, so
 setzen wir
 $$ m_G := \bigl((G \cap UV)O \cap VE\bigr)U \cap OV $$
 und
 $$ b_G := G \cap OV. $$
 Ist umgekehrt $[m,b] \in K \times K$, so setzen wir
 $$ \gamma [m,b] := \bigl((mU \cap EV)O \cap UV\bigr)b. $$
 Man sieht wiederum ohne M\"uhe, dass die Abbildung $\gamma$ eine Bijektion von
 $K \times K$ auf die Menge der nicht durch $V$ gehenden Geraden von $M$ ist und
 dass $\gamma^{-1}(G) = [m_G,b_G]$ ist.
 \par
       Schlie\ss lich definieren wir auf $K$ eine Addition durch
 $$ x + y := \bigl((xU \cap OE)V \cap (EO \cap UV)y\bigr)U \cap OV $$
 und eine Multiplikation durch
 $$ xy := \bigl((xU \cap EV)O \cap (yU \cap OE)V\bigr)U \cap OV. $$
 \medskip\noindent
 {\bf 6.4. Satz.} {\it Ist $\tau \in \E(UV)$ und ist $O^\tau = \pi(a,b)$, so ist
 $$ \tau(x,y)^\tau = \tau(x+a, y+b) $$
 f\"ur alle $x,y \in K$. Insbesondere ist $K$ bez.\ der Addition eine zu
 $\E(V,UV)$ isomorphe Gruppe und $O$ ist das bez.\ der Addition
 neutrale Element von $K$. Wir setzen daher im Folgenden $0 := O$.}
 \smallskip
       Beweis. Es sei zun\"achst $a = 0$. Dann ist $\tau (a,b) = \tau (0,b) = b$.
 Ferner ist $y = (yU \cap OE)U \cap OV$ und daher
 $$ y^\tau = (yU \cap OE)^\tau U \cap OV, $$
 da ja $\tau \in \E (V, UV)$ ist. Nun ist
 $$ (yU \cap OE)^\tau \leq (yU \cap OE)V. $$
 Ferner ist
 $$ yU \cap OE \leq OE = O(OE \cap UV) $$
 und somit
 $$ (yU \cap OE)^\tau \leq O^\tau(OE \cap UV) = b(OE \cap UV). $$
 Insgesamt erhalten wir
 $$ y^\tau = \bigl((yU \cap OE)V \cap b(OE \cap UV)\bigr) U \cap OV = y + b. $$
 Nun ist $\pi(x,y) = (xU \cap OE)V \cap yU$ und daher
 $$\eqalign{\pi(x,y)^\tau &= \bigl((xU \cap OE)V\bigr)^\tau \cap y^\tau U^\tau \cr
                          &= (xU \cap OE)V \cap (y + b)U               \cr
                          &= \pi(x,y + b). \cr}$$
 Es sei $\rho$ die nach 6.3 vorhandene involutorische Perspektivit\"at mit der
 Achse $OE$, die $U$ mit $V$ vertauscht. Dann ist
 $$ \eqalign{\pi(x,y)^\rho &= \bigl((xU \cap OE)V \cap yU\bigr)^\rho    \cr
                           &= \bigl((xU \cap OE)V \cap (yU \cap OE)U\bigr)^\rho
									\cr
			   &= (xU \cap OE) U \cap (yU \cap OE)V         \cr
                           &= (yU \cap OE)V \cap xU                     \cr
			   &= \pi(y,x). \cr} $$
 Es sei nun $b = 0$. Dann ist $\tau \in \E(U,UV)$ und daher
 $\rho^{-1} \tau \rho \in \E(V,UV)$. Ferner ist
 $$ O^{\rho^{-1}} = O^{\tau \rho} = \pi(a,0)^\rho = \pi(0,a). $$
 Nach dem bereits erledigten Fall ist also
 $$ \pi(x,y)^{\tau\rho} = \pi(y,x)^{\rho^{-1}\tau \rho} = \pi(y,x + a)
			= \pi(x + a,y)^\rho. $$
 Also ist $\pi(x,y)^\tau = \pi(x + a,y)$.
 \par
       Bezeichnet man die Elation aus $\E(U,V)$, die O auf $\pi(a,b)$
 abbildet, mit $\tau(a,b)$, so folgt mit dem bereits Bewiesenem
 $$ O^{\tau(a,0)\tau(0,b)} = \pi(a,0)^{\tau(0,b)} = \pi(a,0 + b)
                           = \pi(a,b) = O^{\tau(a,b)}, $$
 so dass $\tau(a,b) = \tau(a,0)\tau(0,b)$ ist. Hieraus folgt schlie\ss lich
 $$\eqalign{ \pi(x,y)^{\tau(a,b)} &= \pi(x,y)^{\tau (a,0)\tau(0,b)}   \cr
				  &= \pi(x + a,y)^{\tau(0,b)}         \cr
				  &= \pi(x + a,y + b). \cr} $$
 Schlie\ss lich ist
 $$ O^{\tau(0,a)\tau(0,b)} = \pi(0,a + b) = O^{\tau(0,a + b)}. $$
 Daher ist $\tau(0,a)\tau(0,b) = \tau(0,a + b)$. Weil die Abbildung
 $a \to \tau(0,a)$ offensichtlich eine Bijektion von $K$ auf $\E(V,UV)$
 ist, folgt somit, dass $(K,+)$ eine zu $\E(V,UV)$ isomorphe Gruppe ist.
 \medskip\noindent
 {\bf 6.5. Satz.} {\it Genau dann ist $\pi(x,y) \leq \gamma [m,b]$,
 wenn $y = mx + b$ ist. Ist $G$ eine von $UV$ verschiedene Gerade durch
 $V$, so ist genau dann $\pi(x,y) \leq G$, wenn $x = (G \cap OE)U \cap OV$ ist.}
 \smallskip
       Beweis. $\gamma[m,0]$ ist eine Gerade durch $O$ und aus der
 Definition der Multiplikation in $K$ folgt unmittelbar, dass genau
 dann $\pi(x,y) \leq \gamma[m,0]$ ist, wenn $y = mx$ ist. Es ist
 $\gamma[m,b]^{\tau(0,-b)} = \gamma[m,0]$. Daher ist wegen 6.4
 genau dann $\pi(x,y) \leq \gamma[m,b]$, wenn $\pi(x,y - b) \leq \gamma[m,0]$
 ist. Hieraus folgt, dass $\pi(x,z)$ genau dann auf $\gamma[m,b]$ liegt, wenn
 $y - b = mx$ ist. Weil $K$ bez\"uglich der Addition eine Gruppe ist,
 folgt somit, dass $y = mx + b$ mit $P(x,y) \leq \gamma[m,b]$
 gleichbedeutend ist.
 \par
       Die zweite Behauptung des Satzes ist trivial.
 \medskip
       Wir setzen im Folgenden $1 := UE \cap OV$.
 \medskip\noindent
 {\bf 6.6. Satz.} {\it Es ist $0x = 0 = 0x$ und $1x = x = x1$ f\"ur alle
 $x \in K$.}
 \smallskip
       Beweis. Es ist
 $$ 0x = \bigl(OU \cap (xU \cap OE)V\bigr)U \cap OV = OU \cap OV = 0 $$
 und
 $$ x0 = \bigl((xU \cap EV)O \cap OV\bigr)U \cap OV = 0. $$
 Also gilt $0x = x0 = 0$ f\"ur alle $x \in K$.
 \par
       Es ist
 $$\eqalign{1x &= \bigl(EO \cap (xU \cap OE)V\bigr)U \cap OV    \cr
	      &= (xU \cap OE)U \cap OV                         
              = xU \cap OV = x \cr} $$
 und
 $$\eqalign{x1 &= \bigl((xU \cap EV)O \cap EV\bigr)U \cap OV     \cr
	       &= (xU \cap EV)U \cap OV                       
               =  xU \cap OV = x. \cr} $$
 Damit ist alles bewiesen.
 \medskip\noindent
 {\bf 6.7. Satz.} {\it Sind $a,b \in K$ und ist $a \neq 0$, so gibt es eindeutig
 bestimmte Elemente $x,y \in K$ mit $ax = b$ und $ya = b$.}
 \smallskip
       Beweis. Wir betrachten die Geraden $\gamma[a,0]$ und $\gamma[0,b]$. Weil
 $a \neq 0$ ist, sind diese beiden Geraden nicht parallel. Es gibt also einen
 Punkt $\pi(x,z)$, der sowohl auf $\gamma[a,0]$ als auch auf $\gamma[0,b]$ liegt.
 Nach 6.6 und 6.5 sowie 6.4 ist daher $z = 0x + b = b$ und $z = ax + 0 = ax$.
 Also ist $ax = b$. Ist $ax' = b$, so liegt wegen $b = 0x' + b$ der Punkt
 $\pi(x',b)$ ebenfalls auf den beiden Geraden $\gamma[a,0]$ und
 $\gamma[0,b]$. Daher ist $\pi(x,b) = \pi(x',b)$ und weiter $x = x'$.
 Dies zeigt, dass es im Falle $a \neq 0$ genau ein $x$ gibt mit $ax = b$.
 \par
       Da $a \neq 0$ ist, ist $\pi(a,b) \neq \pi(0,0)$. Es gibt folglich genau
 eine Ge\-ra\-de $G$, die $\pi(0,0)$ mit $\pi(a,b)$ verbindet.  Aus 6.5 und
 $a \neq 0$ folgt, dass $G$ nicht durch $V$ geht. Also ist $G = \gamma[y,z]$.
 Weil $O$ auf $G$ liegt, ist $z = 0$. Daher ist $b = ya$. Somit existiert also
 ein $y$ mit $b = ya$. Die Einzigkeit der L\"osung $y$ folgt aus der
 Eindeutigkeit der Verbindungsgeraden von $O$ mit $\pi(a,b)$.
 \medskip\noindent
 {\bf 6.8. Satz.} {\it Es ist $a(b + c) = ab + ac$ und $(a + b)c = ac + bc$
 f\"ur alle $a$, $b$, $c \in K$.}
 \smallskip
       Beweis. Wir betrachten die durch $\pi(x,y)^\tau := \pi(x + c,y)$
 definierte Abbildung $\tau$. Nach 6.4 ist $\tau \in \E(U,UV)$. Da $\tau$ jede
 Gerade auf eine zu ihr parallele abbildet, ist $\gamma[a,0]^\tau = \gamma[a,d]$.
 Nun ist $\pi(x,ax) \leq \gamma[a,0]$. Daher ist
 $\pi(x + c, ax] \leq \gamma[a,d]$. Folglich ist $ax = a(x + c) + d$ f\"ur alle
 $x \in K$. Mit $x := 0$ folgt, dass $ac + d = 0$ ist. Folglich ist auch
 $d + ac = 0$ und weiter
 $$ ax + ac = a(x + c) + d + ac = a(x + c). $$
 Mit $x := b$ erh\"alt man hieraus die G\"ultigkeit der ersten Behauptung des
 Satzes.
 \par
       Es sei $a \in K$. Weil $\E(V, OV)$ die von $V$ verschiedenen Punkte von
 $UV$ transitiv permutiert, gibt es ein $\tau \in \E(V,OV)$ mit
 $\gamma[0,0]^\tau = \gamma[\alpha, r]$ mit einem geeigneten $r \in K$. Nun ist
 $O \leq \gamma[0,0]$ und $O^\tau = O$. Daher ist $r = 0$, dh., es ist
 $\gamma[0,0]^\tau = \gamma[a,0]$.
 \par
       Die Geraden $\gamma[0,0]$ und $\gamma[0,y]$ sind parallel. Daher
 sind auch die Ge\-ra\-den $\gamma[0,0]^\tau$ und $\gamma[0,y]^\tau$ parallel.
 Hieraus folgt,
 dass $\gamma[0,y]^\tau = \gamma[a,y']$ ist. Nun ist $\pi(0,y) \leq \gamma[0,y]$
 und $\pi(0,y)^\tau = \pi(0,y)$. Folglich ist $y = a0 + y'$ und somit $y = y'$.
 Es ist also $\gamma[0,y]^\tau=\gamma[a,y]$.
 \par
       Ist $\pi(x,y)^\tau = \pi(x',y')$, so ist $x = x'$, da $\tau$ die Geraden
 durch $V$ einzeln festl\"asst. Nun ist $\pi(x,y) \leq \gamma[0,y]$ und daher
 $\pi(x,y') \leq \gamma[a,y]$. Es ist also $\gamma' = ax + y$, dh., es ist
 $\pi(x,y)^\tau = \pi(x,ax + y)$.
 \par
       Es sei $\gamma[b,0]^\tau = \gamma[b', t]$. Wegen
 $O^\tau = O \leq \gamma[b,0]$ ist $t = 0$. Nun ist $\pi(x,bx) \leq \gamma[b,0]$.
 Hieraus folgt, dass $\pi(x,ax + bx) \leq \gamma[b',0]$ ist, so dass
 f\"ur alle $x \in K$ die Gleichung $ax + bx = b'x$ gilt. Mit $x = 1$
 erh\"alt man $a + b = b'$ und mit $x = c$ schlie\ss lich $ac + bc = (a + b)c$.
 Damit ist alles bewiesen.
 \medskip\noindent
 {\bf 6.9. Satz.} {\it Ist $a'a = 1$, so ist $a'(ab) = b$ f\"ur alle $b \in K$.}
 \smallskip
       Beweis. $\gamma[a,0]$ ist die Verbindungsgerade von $O$ und $\pi(1,a)$ und
 $\gamma[a',0]$ ist die Verbindungsgerade von $O$ und $\pi(a,1)$. Ist $\rho$ die
 nach 6.3 existierende Perspektivit\"at mit der Achse $OE$, die $U$ mit $V$
 vertauscht, so ist, wie wir beim Beweise von 6.4 gesehen haben,
 $\pi(x,y)^\rho = \pi(y,x)$. Wegen $O^\rho = O$ und $\pi(1,a)^\rho = \pi(a,1)$
 ist folglich $\gamma[a,0]^\rho = \gamma[a',0]$. Nun ist
 $\pi(b,ab) \leq \gamma[a,0]$ und daher $\pi(ab,b) \leq \gamma[a',0]$, so dass
 also $b = a'(ab)$ gilt, q. e. d.
 \medskip\noindent
 {\bf 6.10. Satz.} {\it Ist $a'a = 1$, so ist auch $aa' = 1$. Ist ferner
 $aa'' = 1$, so ist $a' = a''$. Jedes von Null verschiedene Element von $K$ hat
 also genau ein Inverses.}
 \smallskip
       Beweis. Nach 6.9 ist $a'1 = a' = a'(aa')$. Wegen $a' \neq 0$ ist daher
 $aa' = 1$ Also ist $aa' = aa''$. Weil $a \neq 0$ ist, ist somit $a' = a''$.
 \par
       Die Existenzaussage des Satzes folgt mit 6.7.
 \medskip\noindent
 {\bf 6.11. Satz.} {\it Ist $a \neq 0$, so ist $(ba)a^{-1} = b$ f\"ur alle
 $b \in K$.}
 \smallskip
       Beweis. Es sei $\rho$ die involutorische Perspektivit\"at mit der
 Achse $EV$, die $O$ mit $U$ vertauscht. Da $V$ ein Fixpunkt von
 $\rho$ ist, permutiert $\rho$ die Geraden durch $V$ unter sich.
 Ist $\pi(x,y)$ ein Punkt von $M_{UV}$ mit $x \neq 0$, so ist daher
 $\pi(x,y)^\rho = \pi(\varphi(x),\psi(x,y))$ mit $\varphi(x)$, $\psi(x,y) \in K$.
 Nun ist $\gamma[0,y] = U\pi(1,y)$. Daher ist, da
 ja $\pi(1,y) \leq EV$ ist, $\gamma[0,y]^\rho = O \pi(1,y)$. Somit
 ist $\gamma[0,y]^\rho = \gamma[y,0]$. Ferner ist $\pi(x,y) \leq \gamma[0,y]$ und
 folglich $\pi(\varphi(x),\psi(x,y)) \leq \gamma[y,0]$. Hieraus folgt, dass
 $\psi(x,y) = y \varphi(x)$ ist, falls nur $x \neq 0$ ist. Es ist
 $\gamma[0,1]^\rho = \gamma[1,0]$ und daher, da $\rho$ involutorisch ist,
 $\gamma[1,0]^\rho = \gamma[0,1]$. Wegen $\pi(x,x) \in \gamma[1,0]$ ist
 $\pi(\varphi(x),\psi(x,x)) \leq \gamma[0,1]$, so dass $\psi(x,x) = 1$ ist f\"ur
 alle $x\in K$. Ist $x \neq 0$ so ist $\psi(x,x) = x\varphi(x)$, so dass in
 diesem Falle $\varphi(x) = x^{-1}$ ist. Also ist $\psi(x,y) = yx^{-1}$ f\"ur
 alle von Null verschiedenen $x \in K$.
 \par
       Es sei nun $a \neq 0$. Dann ist $\pi(a,ba) \leq \gamma[b,0]$.
 Folglich ist
 $$ \pi\bigl(a^{-1},(ba)a^{-1}\bigr) \leq \gamma[0,b]. $$
 Hieraus folgt, dass $(ba)a^{-1} = 0a^{-1} + b = b$ ist, q. e. d.
 \medskip\noindent
 {\bf 6.12.} {\it Es ist $a(ab) = a^2b$ und $(ba)a = ba^2$ f\"ur alle $a$, $b \in
 K$, dh., $K$ ist ein {\rm Alternativk\"orper\/}.}\index{Alternativk\"orper}{}
 \smallskip
       Beweis. Ist $a = 0$, so ist sicherlich $(ba)a = ba^2$.
       Es gilt
 $$ 0 = 0 \cdot (-1) = (1 - 1)(-1) = 1 \cdot (-1) + (-1)^2 = -1 + (-1)^2. $$
 Also ist $(-1)^2 = 1$.
 \par
       Weil beide Distributivgesetze gelten, gilt auch $x(-y) = (-x)y = -xy$.
 Es folgt
 $$ \bigl(b(-1)\bigr)(-1) = -(-b) = b = b(-1)^2. $$
 Die Aussage des Satzes gilt also f\"ur $a = 0$ und $a = -1$.
 \par
       Es sei also $a \neq 0$, $-1$. Dann gilt mit $c \in K$
 $$\eqalign{
    \bigl((ca)(a^{-1} - (a + 1)^{-1}\bigr)(a + 1)
     &= \bigl(c - (ca)(a + 1)^{-1}\bigr)(a + 1) \cr
     &= ca + c - ca = c. \cr} $$
 Hierbei haben wir benutzt, dass $(K,+)$ auf Grund von 6.4 und
 II.1.15 abelsch ist. Hieraus folgt
 $$ ca = \bigl(c(a + 1)^{-1}\bigr)\bigl(a^{-1} - (a + 1)^{-1}\bigr)^{-1} $$
 f\"ur alle $c\in K$. Mit $c = a + 1$ folgt
 $$ (a + 1)a = (a^{-1} - (a + 1)^{-1})^{-1}. $$
 Also ist
 $$ ca = \bigl(c(a + 1)^{-1}\bigr)\bigl((a + 1)a\bigr) $$
 f\"ur alle $c \in K$. Setzt man $c := b(a + 1)$, so folgt hieraus
 $$ (ba)a + ba = b(a^2 + a) = ba^2 + ba $$
 und damit $(ba)a = ba^2$.
       Die andere Identit\"at beweist man entsprechend.
 \bigskip
       Sind $a$, $b$, $c \in K$, so definieren wir den
 {\it Assoziator\/}\index{Assoziator}{} $[a,b,c]$
 durch $[a,b,c] := (ab)c - a(bc)$. Man \"uberzeugt sich
 leicht, dass der Assoziator in allen Argumenten additiv ist.
       6.12 besagt nun gerade, dass $[a,a,b] = 0 = [a,b,b]$ ist. Ferner ist
 $$ 0 = [a,b + c,b + c] = [a,b,b] + [a,b,c] + [a,c,b] + [a,c,c]. $$
 Also ist $[a,b,c] = -[a,c,b]$. Insbesondere ist also auch $[a,b,a] = 0$.
 Ebenso leicht rechnet man nach, dass $[a,b,c] = -[b,a,c]$ ist. Dann
 folgt aber auch $[a,b,c] = -[a,c,b] = [c,a,b] = -[c,b,a]$. Der
 Assoziator \"andert also sein Vorzeichen, wenn man zwei seiner
 Argumente vertauscht.
 \par
       Au\ss er dem Assoziator ben\"otigen wir noch den
 {\it Kommutator\/}\index{Kommutator}{}
 $$ [a,b] := ab - ba $$
 zweier Elemente $a$, $b \in K$ und die durch
 $$ k(w,x,y,z) := [wx,y,z] - [x,y,z]w - x[w,y,z] $$
 definierte Funktion $k$. Man \"uberzeugt sich wiederum sehr
 leicht, dass sowohl der Kommutator als auch die Funktion $k$ in
 allen Argumenten additiv ist.
 \medskip\noindent
 {\bf 6.13. Satz.} {\it Ist $\rho$ eine Permutation von $\{w,x,y,z\}$, so ist}
 $$k(w,x,y,z) = \sgn(\rho)k(w^\rho, x^\rho, z^\rho).$$
 \par
       Beweis. Es ist
 $$\eqalign{
   [wx,y,z] - [xy,z,w] + [yz,w,x]                     
                       &= [wx, y,z] - [w,xy,z] + [w,x,yz] \qquad\qquad\quad \cr
                       &= \bigl((wx)y\bigr)z - (wx)(yz) - \bigl(w(xy)\bigr)z  \cr
		       & \hskip .5cm + w\bigl((xy)z\bigr) + (wx)(yz)
			      - \bigl(w(x(yz)\bigr) \cr
		       &= [w,x,y]z + w[x,y,z]. \cr}$$
 Daher ist
 $$\eqalign{
  k(w,x,y,z) &- k(x,y,z,w) + k(y,z,w,x)                            \cr
        &= [wx,y,z] - [x,y,z]w - x[w,y,z] - [xy,z,w] + [y,z,w]x    \cr 
	& \hskip .4cm + y[x,z,w] + [yz,w,x] - [z,w,x]y - z[y,w,x]  \cr
        &= [w,x,y]z + w[x,y,z] - [x,y,z]w - x[w,y,z] + [y,z,w]x    \cr
	& \hskip .4cm + y[x,z,w] - [z,w,x]y - z[y,w,x]             \cr
        &= \bigl[w,[x,y,z]\bigr] + \bigl[[w,x,y],z\bigr]
	      + \bigl[[y,z,w],x\bigr] + \bigl[y,[x,z,w]\bigr]. \cr} $$
 Es ist also
 $$\eqalign{
      k(w,x,y,z) &- k(x,y,z,w) + k(y,z,w,x)   \cr
            &= \bigl[w,[x,y,z]\bigr] + \bigl[[w,x,y],z\bigr]
	      + \bigl[[y,z,w],x\bigr] + \bigl[y,[x,z,w]\bigr]. \cr}$$
 Ersetzt man hierin $w$, $x$, $y$, $z$ durch $x$, $y$, $z$, $w$, so erh\"alt man
 $$\eqalign{
     k(x,y,z,w) &- k(y,z,w,x) + k(z,w,x,y)           \cr
           &= \bigl[x,[y,z,w]\bigr] + \bigl[[x,y,z],w\bigr]
	      + \bigl[[z,w,x],y\bigr] + \bigl[z,[y,w,x]\bigr]. \cr}$$
 Addiert man diese beiden Gleichungen, so ergibt sich
 $$ k(w,x,y,z) + k(z,w,x,y) = 0. $$ Also ist
 $$ k(w,x,y,z) = -k(z,w,x,y). $$
 Ferner ist
 $$\eqalign{
        k(w,x,y,z) &= [wx,y,z] - [x,y,z]w - x[w,y,z]   \cr
	           &= -\bigl([wx,z,y] - [x,z,y]w - x[w,z,y]\bigr) \cr
		   &= -k(w,x,z,y). \cr} $$
 Definiert man $\sigma$ durch $w^\sigma := z$, $x^\sigma := w$, $y^\sigma := x$,
 $z^\sigma := y$ und $\tau$ durch $w^\tau := w$,
 $x^\tau := x$, $y^\tau := z$ und $z^\tau := y$, so ist $\sigma$ ein
 4-Zyklus und $\tau$ eine Transposition, daher ist \sgn$(\sigma) = -1 =
 \sgn (\tau)$. Also ist
 $$ k(w,x,y,z) = \sgn(\sigma)k(w^\sigma,x^\sigma,y^\sigma,z^\sigma)
      = \sgn(\tau)k(w^\tau,x^\tau,y^\tau,z^\tau). $$
 Weil $\sigma$ und $\tau$ die symmetrische Gruppe auf $\{x,y,z,w\}$ erzeugen,
 folgt damit die Behauptung von 6.13.
 \medskip
       Nun ist $k(w,x,y,y) = 0$. Aus 6.13 folgt, 
 dass sicher dann
 $k(w,x,y,z) = 0$ ist, falls zwei der Argumente gleich sind.
 \par
       Der n\"achste Satz ist der Schl\"ussel zum Beweise des Satzes von
 Artin-Zorn, dass n\"amlich endliche Alternativk\"orper K\"orper
 sind. Da\-bei spielt auch eine Rolle, dass die multiplikative Gruppe
 eines end\-li\-chen K\"orpers zyklisch ist, so dass endliche K\"orper
 von einem Element erzeugt werden.
 \par
       Alternativk\"orper\index{Alternativk\"orper}{} beschreiben stets
 Moufangebenen\index{Moufangebene}{}. Das werden wir hier nicht beweisen, da
 es zum Beweise von 6.14 und dem Satz von Artin-Zorn nicht ben\"otigt wird. Der
 interessierte Leser sei f\"ur die Theorie der Alternativk\"orper auf Pickert
 1955, Kapitel 6 verwiesen.
 \medskip\noindent
 {\bf 6.14. Satz.} {\it Sind $a$ und $b$ Elemente des Alternativk\"orpers $K$, so
 gibt es einen assoziativen Unterring von $K$, der $a$ und $b$ enth\"alt.}
 \smallskip
       Beweis. Wir setzen $E:= \{a,b\}$. Sind $e$, $f$, $g \in E$, so ist
 $[e,f,g] = 0$, da wenigstens zwei der Elemente $e$, $f$, $g$ gleich sind.
 Wir definieren die Menge $X$ durch
 $$ X := \{x \mid x \in K,\ [e,f,x]=0\ \hbox{\rm f\"ur alle}\ e, f \in E\}. $$
 Dann ist $E \subseteq X$, wie wir gerade gesehen haben. Sind $e$, $f$, $g \in E$
 und ist $x \in E$, so ist
 $$ [ex,f,g] = k(e,x,f,g) + [x,f,g]e + x[e,f,g] = 0, $$
 auf Grund der Definition von $X$ bzw. der Tatsache, dass zwei der Elemente
 $e$, $f$, $g$ gleich sind. Also ist $EX \subseteq X$. Wir definieren $Y$
 durch
 $$ Y := \bigl\{y \mid y \in K, [e,x,y] = 0\ \hbox{\rm f\"ur alle}\ e \in E\ 
	     \hbox{und alle}\ x \in X\bigr\}. $$
 Wegen $E \subseteq X$ ist $Y \subseteq X$. Es sei $y \in Y$ und $x \in X$ und
 $e$, $f \in E$. Dann ist
 $$\eqalign{
     [e,f,yx] &= [yx,e,f]     \cr
	      &= k(y,x,e,f) + [x,e,f]y + x[e,y,f]  \cr
	      &= k(y,x,e,f)                       
	      = -k(e,x,y,f)                       \cr
	      &= -[ex,y,f] + [x,y,f]e + x[e,y,f]  
	      = 0, \cr} $$
 da auf Grund der Definition von $Y$ und wegen $EX \subseteq X$ die
 Assoziatoren $[ex,y,f]$, $[x,y,f]$ und $[e,y,f]$ alle gleich Null
 sind. Somit ist auch $YX \subseteq X$. Schlie\ss lich sei
 $$ R := \bigl\{r \mid r \in K, [x,y,r] = 0\ \hbox{\rm f\"ur alle}\ x \in X\ 
     \hbox{\rm und alle}\ y \in Y\bigr\}. $$
 Aus der Definition von $X$ und $Y$ folgt $E \subseteq R$. Weil $E \subseteq Y$
 ist, ist $R \subseteq Y$. Also ist $E \subseteq R \subseteq Y \subseteq X$.
 Folglich ist $[r,s,t] = 0$ f\"ur alle $r$, $s$, $t \in R$. Weil $R$ wegen
 $E \subseteq R$ nicht leer ist und offensichtlich $r - s \in R$ gilt f\"ur alle
 $r$, $s \in R$, ist $(R,+)$ eine Gruppe. Es seien $r$, $s \in R$ und $x \in X$
 und $y \in Y$. Dann ist
 $$ k(x,y,r,s) = [yx,r,r] - [x,r,s]y - x[y,r,s] = 0, $$
 da $R \subseteq Y$ und $YX\subseteq X$ ist. Nach 6.13 ist folglich auch
 $k(r,s,y,x)=0$. Somit ist
 $$ 0 = [rs,y,x] - [s,x,y]r - s[r,y,x] = [rs,y,x]. $$
 Also ist auch $rs \in R$, so dass $R$ in der Tat ein assoziativer Ring ist, der
 $a$ und $b$ enth\"alt.\index{Satz von Artin--Zorn}{}
 \medskip\noindent
 {\bf 6.15. Satz von Artin-Zorn.} {\it Ist $K$ ein endlicher
 Alternativk\"orper, so ist $K$ ein K\"orper.}
 \smallskip
       Beweis. Es sei $S$ ein maximaler assoziativer Teilring von $K$.
 Weil $K$ endlich ist und keine Nullteiler hat, folgt, dass $S$ ein
 endlicher, nullteilerfreier Ring ist. Folglich ist $S$ ein
 K\"orper. Als endlicher K\"orper ist $S$ nach dem Satz von
 Wedderburn kommutativ. Insbesondere ist die multiplikative Gruppe
 von $S$ zyklisch. Es sei $a$ ein erzeugendes Element der
 multiplikativen Gruppe von $S$. Es sei $b \in K$. Nach 6.14 gibt
 es einen assoziativen Teilring $R$ von $K$ mit $a,b \in R$. Es
 folgt $S \subseteq R$ und damit $S = R$, da $S$ ja ein maximaler
 assoziativer Teilring von $K$ ist. Also ist $b \in S$ und somit $K
 = S$, so dass $K$ in der Tat ein K\"orper ist.
 \medskip
       Wir sind nun in der Lage, Satz 6.1 zu beweisen. Weil $M$ endlich
 ist, ist $K$ nach 6.15 ein K\"orper. Es sei $k \in K^*$. Wir
 definieren $\sigma(k)$ durch
 $$ \pi(x,y)^{\sigma(k)} := \pi (xk,yk). $$
 Dann ist $\sigma (k)$ eine Bijektion der Punktmenge von $M_{UV}$ auf sich und es
 gilt $\sigma (k)^{-1} = \sigma(k^{-1})$. Ferner folgt aus
 $\pi(c,y)^{\sigma(k)} = \pi(ck,yk)$ und
 $$ \pi(x,mx + b)^{\sigma(k)} = \pi\bigl(xk,(mx + b)k\bigr)
			      = \pi\bigl(xk,m(xk) + b\bigr), $$
 dass $\sigma (k)$ kollineare Punkte auf kollineare Punkte abbildet.
 Weil dann auch $\sigma(k^{-1})$ kollineare Punkte auf kollineare
 Punkte abbildet, ist $\sigma (k)$ eine Kollineation von $M_{UV}$.
 Weil $\sigma (k)$ offensichtlich jede Gerade auf eine zu ihr
 parallele Gerade abbildet und weil $O^{\sigma(k)} = O$ ist, ist
 $\sigma(k) \in \Delta(O,UV)$. Aus
 $$ E^{\sigma(k)} = \pi(1,1)^{\sigma(k)} = \pi(k,k) $$
 folgt schlie\ss lich, dass $M$ eine $(O, UV)$-transitive Ebene ist. Da dies
 f\"ur alle nicht kollinearen Punkte $O$, $U$, $V$ gilt, ist $M$ nach II.1.8
 desarguessch.

\mysection{7. Ebenen mit vielen Perspektivit\"aten}

\noindent
 Ist $K$ eine Kollineationsgruppe
 eines projektiven Verbandes $L$ und ist $H$ eine Hyperebene von $L$, so setzen
 wir $K(H) := K \cap \E(H)$. Ist au\ss erdem $P$ ein Punkt von $L$, so setzen wir
 $K(P,H) := K \cap \Delta(P,H)$. Es ist also $K(H)$ die Menge der in $K$
 enthaltenen Elationen mit der Achse $H$ und $K(P,H)$ die Menge der in $K$
 enthaltenen Perspektivit\"aten mit der Achse $H$ und dem Zentrum $P$. Sowohl
 $K(H)$ als auch $K(P,H)$ sind Untergruppen von $K$. Der Leser beachte, dass
 $K(H)$ immer eine Gruppe von Elationen ist, w\"ahrend $K(P,H)$ auch eine Gruppe
 von Streckungen sein kann.
 \medskip\noindent
 {\bf 7.1. Satz.} {\it Es sei $L$ ein endlicher projektiver Verband mit
 $\Rg(L) \geq 3$ und $K$ sei eine Kollineationsgruppe von $L$. Ist $H$ eine
 Hyperebene von $L$ und ist $B$ die Menge der Punkte von $L_H$ mit
 $K(P,H) \neq \{1\}$, so ist $B$ eine Bahn von $K(H)$.}
 \smallskip
       Beweis. Es sei $\Gamma$ die Menge aller in $K$ enthaltenen
 Perspektivit\"aten mit der Achse $H$. Dann ist $\Gamma$ eine
 Untergruppe von $K$. Setze $a_P := |K(P,H)|$. Dann ist
 $$ |\Gamma| = \big|K(H)\big| + \sum_{P \in B}(a_P - 1), $$
 da ja jedes von $1$ verschiedene Elemente von $\Gamma$ genau ein Zentrum hat.
 Setze $k := |B|$.
 Weil $|K(H)| \geq 1$ und $a_P \geq 2$ ist, ist $|\Gamma| \geq 1 + k > k$. Es
 seien $B_1$, \dots, $B_r$, die Punktbahnen von $\Gamma$ mit
 $|B_i| < |\Gamma|$. Dann ist $B = \bigcup^r_{i:=1} B_i$. Ist $P \in B_i$, so ist
 $|\Gamma| = |B_i||K(P,H)|$. Also ist
 $$\eqalign{
    |\Gamma| &= \big|K(H)\big| - k + \sum_{P \in B} \big|K(P,H)\big| \cr
             &= \big|K(H)\big| - k + \sum_{i:=1}^r \sum_{P \in B_i}
	                     \big|K(P,H)\big|                        \cr
	     &= \big|K(H)\big| - k + \sum_{i:=1}^r \sum_{P \in B_i}
			    |\Gamma||B_i|^{-1}                       \cr
	     &= \big|K(H)\big| - k + \sum_{i:=1}^r |\Gamma|.         \cr
	     &= \big|K(H)\big| - k + r|\Gamma|. \cr} $$
 Folglich ist
 $$ 0 \leq (r - 1)|\Gamma| = k - \big|K(H)\big| < k < |\Gamma|. $$
 Hieraus folgt, dass $0 \leq r - 1 < 1$ ist. Dies impliziert wiederum
 $r = 1$ und dann $k = |K(H)|$.
 \par
       $K(H)$ zerlegt $B$ in Bahnen. Weil $K(H)$ auf der Menge der Punkte
 von $L_H$ regul\"ar operiert, folgt aus $|K(H) = k = |B|$, dass $B$
 eine Bahn von $K(H)$ ist. Damit ist der Satz bewiesen.

 \medskip\noindent
 {\bf 7.2. Korollar.} {\it Es sei $L$ ein endlicher projektiver
 Verband, dessen Rang mindestens $3$ sei. Sind dann $\rho$ und
 $\sigma$ von $1$ verschiedene Stre\-ckun\-gen von $L$ mit den Zentren
 $P$ und $Q$ und haben $\rho$ und $\sigma$ die gleiche Achse $H$,
 so gibt es in der von $\rho$ und $\sigma$ erzeugten Gruppe eine
 Elation $\tau$, die nat\"urlich ebenfalls die Achse $H$ hat, mit $P^\tau = Q$.}
 \smallskip
       Dies folgt unmittelbar aus 7.1.
 \medskip\noindent
 {\bf 7.3. Satz.} {\it Es sei $E$ eine endliche projektive Ebene. Ist $K$ eine
 Kol\-li\-ne\-a\-ti\-ons\-grup\-pe von $E$ und ist $K(P,G) \neq \{1\}$ f\"ur alle
 nicht inzidenten Punkt-Geraden\-paare $(P,G)$ von $E$, so ist $E$ desarguessch und
 $K$ enth\"alt die kleine projektive Gruppe von $E$.}
 \medskip
       Beweis. Ist $G$ eine Gerade von $E$, so ist $K(G)$ nach 7.1 auf
 der Menge der Punkte von $E_G$ transitiv. Dies besagt einmal, dass
 $E$ eine Moufangebene und somit nach 6.1 desarguessch ist, und zum
 anderen, dass $K$ alle Elationen von $E$ und folglich die kleine
 projektive Gruppe von $E$ enth\"alt.
 \medskip
       Der Satz, dass eine Korrelation eines endliche projektiven
 Verbandes stets einen absoluten Punkt hat, hat eine
 \"uberraschende Konsequenz, n\"amlich den folgenden Satz.
 \medskip\noindent
 {\bf 7.4. Satz.} {\it Es sei $L$ ein endlicher projektiver Verband
 mit $\Rg(L) \geq 3$. Ist $K$ eine Kollineationsgruppe von $L$ und
 ist jeder Punkt von $L$ Zentrum einer nicht trivialen Streckung
 aus $K$, so enth\"alt $K$ eine nicht triviale Elation.}
 \smallskip
       Beweis. Ist $L$ ein Gegenbeispiel zu 7.4, so folgt aus dem zu 7.2
 dualen Satz, dass es zu jedem Punkt $P$ von $L$ genau eine
 Hyperebene $P^\pi$ gibt mit $P \not\in P^\pi$ und $K(P,P^\pi) \neq \{1\}$. Aus
 7.2 folgt ferner, dass $\pi$ injektiv ist. I.7.7
 impliziert schlie\ss lich, dass $\pi$ eine Bijektion der Menge der
 Punkte auf die Menge der Hyperebenen von $L$ ist.
 \par
       Es sei $r := \Rg(L)$. Wir zeigen als N\"achstes, dass das
 Ab\-bil\-dungs\-paar $\kappa := (\pi, \pi^{-1})$ eine Korrelation von
 $L_{1,r-1}$ ist. Es sei $P$ ein Punkt und $H$ eine Hyperebene von
 $L$. Ferner sei $Q := H^{\pi^{-1}}$. Ist $P \leq H$, so folgt
 $H^\sigma = H$ f\"ur alle $\sigma \in K(P,P^\pi)$. Hieraus folgt
 wiederum, dass $Q^\sigma = Q$ ist f\"ur alle $\sigma \in K(P,P^\pi)$. Wegen
 $Q \not\in Q^\pi = H$ und $P \leq H$ ist $Q \neq P$. Daher ist $Q \leq P^\pi$
 nach II.1.3. Damit ist gezeigt, dass
 $Q \leq P^\pi$ aus $P \leq H = Q^\pi$ folgt. Wegen $\kappa^2 = 1$ gilt
 daher $P \leq Q^\pi$ genau dann, wenn $Q \leq P^\pi$ ist. Damit
 ist $\kappa$ als Korrelation von $L_{1,r-1}$ erkannt. Da $\kappa$
 eine Korrelation von $L_{1,r-1}$ ist, gibt es nach 4.3a) einen
 Punkt $W$ mit $W \leq W^\kappa$. Andererseits ist $P \not\leq P^\pi$
 f\"ur alle Punkte $P$ von $L$, ein Widerspruch.
 \medskip\noindent
 {\bf 7.5. Satz} {\it Es sei $E$ eine endlich projektive Ebene und
 $K$ eine Kol\-li\-ne\-a\-ti\-ons\-grup\-pe von $E$. Ist $(P,G)$ eine Fahne von
 $E$, so ist entweder $|K(P,G)|^2 \geq |K(G)|$ oder alle Elationen
 aus $K$, deren Zentrum $P$ ist, haben $G$ als Achse.}
 \smallskip
       Beweis. Es sei $H$ eine von $G$ verschiedene Gerade durch $P$ und
 es sei $1 \neq \tau \in K(P,H)$. Ist $\sigma \in K(G)$, so ist
 $\sigma^{-1}\tau\sigma \in K(P,H^\sigma)$. Folglich hat
 $\tau^{-1}\sigma^{-1}\tau\sigma$ das Zentrum $P$. Andererseits
 ist $\tau^{-1}\sigma^{-1}\tau \in K(G^\tau) = K(G)$. Daher hat
 $\tau^{-1}\sigma^{-1}\tau\sigma$ die Achse $G$. Also ist
 $$ \tau^{-1}\sigma^{-1}\tau \sigma \in K(P,G). $$
 Es sei $\rho \in K(G)$. Genau dann ist $\tau^{-1}\sigma^{-1}\tau\sigma
 = \tau^{-1}\rho^{-1}\tau\rho$, wenn
 $\rho\sigma^{-1}\tau = \tau \rho \sigma^{-1}$ ist. Wegen $\tau \neq 1$ ist dies
 genau dann der Fall, wenn $H^{\rho\sigma^{-1}} = H$ ist, dh., wenn
 $\rho\sigma^{-1} \in K(P,G)$ ist. Es gibt somit genau $|K:K(P,G)|$ Elemente der
 Form $\tau^{-1}\sigma^{-1}\tau\sigma$. Da alle diese Elemente in $K(P,G)$
 liegen, ist $|K(P,G)| \geq |K(G) : K(P,G)|$, q. e. d.
 \medskip\noindent
 {\bf 7.6. Satz.} {\it Es sei $E$ eine endliche projektive Ebene
 der Ordnung $n$. Ist $G$ eine Gerade von $E$ und $K$ eine
 Kollineationsgruppe mit $|K(G)| > n$, so hat $K(G)$ eine Geradenbahn
 der L\"ange $n$. Insbesondere ist $n$ Teiler von $|K(G)|$.}
 \smallskip
       Beweis. Es sei $H$ eine Gerade von $E$ mit der Eigenschaft, dass
 die L\"ange $a$ der Bahn, zu der $H$ geh\"ort, maximal sei. Ist
 $t$ die Anzahl der Punktbahnen von $K(G)$ in $E_G$, so folgt, da
 $K$ nach 3.4 in $E$ ebenso viele Punkt- wie Ge\-ra\-den\-bahnen hat,
 dass $K(G)$ in $E_G$ genau $t + n$ Geradenbahnen hat, da $K(G)$ ja
 alle Punkte von $G$ festl\"asst und folglich in $E$ genau
 $t + n + 1$ Punktbahnen hat. Nun ist $n(n + 1)$ die Anzahl der von $G$
 verschiedenen Geraden von $E$. Daher folgt aus der Maximalit\"at von $a$, dass
 $a(t + n) \geq n(n+1)$ ist. Andererseits ist $n^2 = t|K(G)| \geq tn$, so dass
 $t \leq n$ ist. Hieraus folgen die Ungleichungen
 $$ 2an = a(n + n) \geq a(t + n) \geq n(n + 1) $$
 und somit $2a \geq n + 1 > n$. Ist $P := G \cap H$ und ist $L$ eine von $G$
 verschiedene Gerade durch $P$, so ist $K(G)_L = K(P,G) = K(G)_H$.  Folglich ist
 $$ \big|K(G) : K(G)_L\big| = \big|K(G) : K(G)_H\big| = a. $$
 Daher wird die Menge der $n$ von $G$ verschiedenen Geraden durch $P$ in
 lauter Bahnen der L\"ange $a$ zerlegt, so dass $a$ Teiler von $n$
 ist. Weil aber $2a > n$ gilt, kann dies nur so sein, dass $a = n$ ist.
 Damit ist der Satz bewiesen.
 \medskip\noindent
 {\bf 7.7. Satz.} {\it Es sei $E$ eine endliche projektive Ebene der Ordnung $n$
 und $K$ sei eine Kollineationsgruppe von $E$. Ist $G$ eine Gerade von $E$ und
 ist $|K(G)| \geq n + 1$, so ist $K(P,G) \neq \{1\}$ f\"ur alle Punkte $P$ auf
 $G$. \"Uberdies ist $n$ Potenz einer Primzahl.}
 \smallskip
       Beweis. Es sei $P$ ein Punkt auf $G$ und $H$ sei eine von $G$
 verschiedene Gerade durch $P$. Ist $a$ die L\"ange der Bahn von $H$ unter
 $K(G)$, so ist wegen $K(G)_H = K(P,G)$ dann
 $$ n\big|K(P,G)\big| \geq a\big|K(P,G)\big| = \big|K(G)\big| \geq n + 1. $$
 Hieraus folgt $K(P,G) \neq \{1\}$.
 \par
       Weil $K(P,G) \neq \{1\}$ ist f\"ur alle $P \leq G$ ist $K(G)$ nach
 II.1.13 eine elementar\-abelsche $p$-Gruppe. Schlie\ss lich ist $n$
 wegen $|K(G)| \geq n + 1$ nach 7.6 ein Teiler von $|K(G)|$, so dass
 also auch $n$ eine Potenz von $p$ ist. Damit ist alles bewiesen.
       
 \medskip\noindent
 {\bf 7.8. Satz.} {\it Es sei $E$ eine projektive Ebene der
 Ordnung $n$. Ferner sei $(P,G)$ eine Fahne und $K$ eine
 Kollineationsgruppe von $E$. Ist dann $h$ eine nat\"urliche Zahl
 mit $h >1$ und ist $|K(Q,G)| = h$ f\"ur alle von $P$ verschiedenen
 Punkte $Q$ auf $G$, so ist $|K(P,G)| = n$.}
 \smallskip
       Beweis. Es ist
 $$ \big|K(G)\big| = \big|K(P,G)\big| + n(h - 1). $$
 Aus $h > 1$ folgt $|K(G)| \geq n + 1$. Aus 7.6 erschlie\ss en wir daher, dass
 $n$ Teiler von $|K(G)|$ ist. Dann ist $n$ aber auch Teiler von
 $|K(P,G)|$. Hier\-aus folgt zusammen mit den trivialerweise
 geltenden Ungleichungen $1 \leq |K(P,G)| \leq n$ die Behauptung.
 \medskip
       Eine unmittelbare Folgerung aus 7.8 ist
 \medskip\noindent
 {\bf 7.9. Korollar.} {\it Es sei $E$ eine endliche projektive
 Ebene der Ordnung $n$ und $G$ sei eine Gerade von $E$. Ist $K$
 eine Kollineationsgruppe von $E$ und ist $h > 1$ eine nat\"urliche
 Zahl, ist ferner $|K(P,G)| = h$ f\"ur alle Punkte $P$ auf $G$, so
 ist $h = n$ und $E$ ist eine Translationsebene bez. $G$. Au\ss erdem
 enth\"alt $K$ alle Elationen mit der Achse $G$.}
 \medskip\noindent
 {\bf 7.10. Satz.} {\it Es sei $\Omega$ eine endliche Menge und
 $\Delta$ sei eine Teilmenge von $\Omega$. Ferner sei $G$ eine
 Permutationsgruppe auf $\Omega$ und $p$ sei eine Primzahl. Gibt es
 dann zu jedem $\delta \in \Delta$ eine $p$-Untergruppe $G(\delta)$
 von $G_\delta$, die au\ss er $\delta$ kein weiteres Fixelement
 hat, so ist $\Delta$ in einer Bahn von $G$ enthalten.}
 \smallskip
       Beweis. Es sei $\delta \in \Delta$ und $\Sigma$ sei die Bahn von
 $\delta$ unter $G$. Dann wird $\Sigma$ von $G(\delta)$ in Bahnen
 $\Sigma_1 = \{\delta\}$, $\Sigma_2$, \dots, $\Sigma_t$ zerlegt. Auf
 Grund unserer Annahme ist $|\Sigma_i| > 1$ f\"ur alle $i \geq 2$.
 Weil $G(\delta)$ eine $p$-Gruppe ist, ist daher $|\Sigma_i|$ als
 Teiler von $|G(\delta)|$ durch $p$ teilbar. Daher ist
 $$ |\Sigma| = \sum_{i:=1}^t|\Sigma_i| \equiv |\sigma_1| \equiv 1 \mod p. $$
 Es sei $\eta \in \Delta$ und $\eta \not\in \Sigma$. Die Gruppe $G(\eta)$ zerlegt
 $\Sigma$ in Bahnen $\Sigma'_1$, \dots, $\Sigma'_s$. Weil $\eta \not\in \Sigma$
 ist, folgt $|\Sigma'_i| > 1$ und damit $|\Sigma'_i| \equiv 0 \mod p$
 f\"ur alle $i$. Hiermit erhalten wir den Widerspruch
 $$ |\Sigma| = \sum_{i:=1}^t |\Sigma'_i| \equiv  0 \mod p. $$
 Dieser Widerspruch zeigt, dass doch $\Delta \subseteq \Sigma$ ist.
 \medskip\noindent
 {\bf 7.11. Satz.} {\it Es sei $E$ eine endliche projektive Ebene
 und $K$ sei eine Kollineationsgruppe von $E$. Ist dann jeder Punkt
 von $E$ Zentrum und jede Gerade von $E$ Achse einer nicht
 trivialen Elation aus $K$, so ist $E$ desarguessch und $K$
 enth\"alt die kleine projektive Gruppe von $E$.}
 \smallskip
       Beweis. Wir zeigen zun\"achst, dass $K(P,G) \neq \{1\}$ ist f\"ur
 alle Fahnen $(P,G)$ von $E$. Angenommen es sei $(P,G)$ eine Fahne
 mit $K(P,G) = \{1\}$. Auf Grund unserer Annahme ist $K(G) \neq \{1\}$
 und $K(P,H) \neq \{1\}$ f\"ur wenigstens eine von $G$ verschiedene
 Gerade $H$ durch $P$. Mit 7.5 erhalten wir daher den Widerspruch
 $$ 1 = \big|K(P,G)\big|^2 \geq \big|K(G)\big| > 1. $$
 Also ist doch $K(P,G) \neq \{1\}$.
 \par
       Es sei $n$ die Ordnung von $E$. Wegen $K(P,G) \neq \{1\}$ ist
 $|K(G)| > n + 1$, so dass $n$ nach 7.7 Potenz einer Primzahl ist.
 \"Uberdies ist $n$ nach 7.6 Teiler von $|K(G)|$. Nach II.1.15
 ist $K(G)$ eine elementarabelsche $p$-Gruppe. Da $n$ auch ein
 Teiler von $|K(H)|$ ist, wenn $H$ eine Gerade von $E$ ist, so
 folgt, dass $K(H)$ f\"ur alle Geraden $H$ von $E$ eine
 elementarabelsche $p$-Gruppe ist.
 \par
       Es sei weiterhin $G$ eine Gerade von $E$. Ist $P$ ein Punkt auf
 $G$ und ist $H$ eine von $G$ verschiedene Gerade durch $P$, so ist
 $K(P,H)$ eine $p$-Gruppe, die auf $G$ eine Permutationsgruppe
 induziert, die $P$ und nur $P$ zum Fixpunkt hat. Mit 7.10 folgt
 daher, dass der Stabilisator $K_G$ von $G$ in $K$ auf der Menge
 der Punkte von $G$ transitiv operiert. Hieraus folgt wiederum,
 dass die Gruppen $K(P,G)$ mit $P \leq G$ alle in $K$ konjugiert
 sind und folglich alle die Ordnung $h > 1$ haben. 7.9 impliziert
 nun, dass $E$ eine Translationsebene bez.\ $G$ ist und dass $\E(G)$
 in $K$ enthalten ist. Somit ist $E$ eine Moufangebene und daher
 nach 6.1 desarguessch. \"Uberdies enth\"alt $K$ die kleine
 projektive Gruppe von $E$.
 \medskip\noindent
 {\bf 7.12. Korollar.} {\it Es sei $E$ eine endliche projektive
 Ebene und $K$ sei eine auf der Menge der Punkte von $E$ transitive
 Kollineationsgruppe von $E$. Enth\"alt $K$ eine von $1$ verschiedene
 Perspektivit\"at, so ist $E$ desarguessch und $K$ enth\"alt die
 kleine projektive Gruppe.}
 \smallskip
       Beweis. Enth\"alt $K$ eine nicht triviale Streckung, so folgt aus
 der Punkttransitivit\"at von $K$, dass jeder Punkt von $E$ Zentrum
 einer nicht trivialen Streckung aus $K$ ist. Nach 7.4 enth\"alt
 $K$ dann auch eine nicht triviale Elation, so dass wir von
 vornherein annehmen d\"urfen, dass $K$ eine nicht triviale Elation
 enth\"alt. Wir folgern wie\-d\-erum aus der Punktransitivit\"at von
 $K$, dass jeder Punkt von $E$ Zentrum einer nicht trivialen
 Elation aus $K$ ist. Weil $K$ nach 3.4 auch auf der Menge der
 Geraden von $E$ transitiv ist, ist auch jede Gerade von $E$ Achse
 einer nicht trivialen Elation aus $K$. Anwendung von 7.11 liefert
 nun die Behauptung.

\mysection{8. Einiges \"uber Permutationsgruppen}

\noindent
 In 
 diesem Abschnitt beweisen wir einige S\"atze \"uber Permutationsgruppen, die
 wir im n\"achsten Abschnitt ben\"otigen werden.
 \medskip\noindent
 {\bf 8.1. Satz.} {\it Es sei $\Gamma$ eine Permutationsgruppe auf
 der Menge $M$ und $B$ und $C$ seien zwei Bahnen von $\Gamma$ mit
 $\ggT(|B|,|C|) = 1$. Ist $b \in B$, so ist $C$ auch eine Bahn von
 $\Gamma_b$.}
 \smallskip
       Beweis. Es sei $c \in C$ und $C'$ sei die Bahn von $c$ unter
 $\Gamma_b$. Ferner sei $B'$ die Bahn von $b$ unter $\Gamma_c$.
 Dann ist
 $$ |\Gamma| = |B||\Gamma_b| = |B||C'||\Gamma_{b,c}| $$
 und
 $$ |\Gamma| = |C||\Gamma_c| = |C||B'||\Gamma_{c,b}|. $$
 Wegen $\Gamma_{b,c} = \Gamma_{c,b}$ ist $|B||C'| = |C||B'|$. Weil $|B|$ und
 $|C|$ teilerfremd sind, folgt, dass $|C|$ Teiler von $|C'|$
 ist. Weil andererseits $C'$ in $C$ enthalten ist, folgt $C = C'$, q. e. d.
 \medskip\noindent
 {\bf 8.2. Satz.} {\it Es sei $\Gamma$ eine auf der Menge $M$ transitiv
 operierende Permutationsgruppe. Ferner sei $n := |M|$ endlich und $p$ sei eine
 in $|\Gamma|$ aufgehende Primzahl.  Schlie\ss lich sei $p^a$ mit $a \geq 0$ die
 h\"ochste Potenz von $p$, die in $n$ aufgeht. Ist $\Pi$ eine $p$-Sylowgruppe von
 $\Gamma$, so ist die L\"ange einer jeden Bahn von $\Pi$ durch $p^a$ teilbar.
 Insbesondere ist die L\"ange jeder k\"urzesten Bahn von $\Pi$ gleich $p^a$.}
 \smallskip
       Beweis. Es sei $x \in M$ und $p^b$ sei die h\"ochste in $|\Gamma_x|$
 aufgehende Potenz von $p$. Dann ist $|\Pi| = p^{a+b}$. Es sei $B$ eine Bahn von
 $\Pi$. Dann ist
 $$ p^{a+b} = |\Pi| = |B||\Pi_y|, $$
 wobei $y \in B$ ist. Wegen $\Pi_y \subseteq \Gamma_y$ und
 $|\Gamma_y| = |\Gamma_x|$, letzteres weil die beiden Gruppen wegen
 der Transitivit\"at von $\Gamma$ konjugiert sind, ist $|\Pi_y|$
 Teiler von $p^b$. Somit ist $|B|$ durch $p^a$ teilbar.
 \par
       Sind $B_1$, \dots, $B_r$ alle Bahnen von $\Pi$, so ist $n =
 \sum_{i:=1}^r|B_i|$. Weil $n$ nicht durch $p^{a+1}$ teilbar ist,
 k\"onnen nicht alle $|B_i|$ durch $p^{a+1}$ teilbar sein. Folglich
 hat $\Pi$ eine Bahn der L\"ange $p^a$.
 \medskip
       Der n\"achste Satz ist einer der \"altesten S\"atze aus der Theorie der
 Permutationsgruppen.\index{Permutationsgruppe}{} Er stammt von Evariste
 Galois.\index{Galois, E.}{}
 \medskip\noindent
 {\bf 8.3. Satz.} {\it Es sei $\Gamma$ eine auf $M$ primitiv
 operierende endliche Permutationsgruppe. Enth\"alt $\Gamma$ einen
 von $\{1\}$ verschiedenen, aufl\"osbaren Normalteiler, so enth\"alt
 $\Gamma$ einen und nur einen minimalen Normalteiler $N$. Der
 Normalteiler $N$ ist eine elementarabelsche $p$-Gruppe und es ist $|N| = |M|$.}
 \smallskip
       Beweis. $\Gamma$ enth\"alt einen nicht trivialen
 aufl\"osbaren Normalteiler $B$. Es gibt daher eine nat\"urliche
 Zahl $i$ mit $N := B^{(i)} \neq \{1\}$ und $B^{(i+1)} = \{1\}$. Weil
 $N$ in $B$ charakteristisch ist, ist $N$ in $\Gamma$ normal. Weil $\Gamma$
 auf $M$ primitiv operiert, ist $N$ folglich transitiv auf $M$.
 \"Uberdies ist $N$ abelsch, da ja $N' = B^{(i+1)} = \{1\}$ ist.
 Hieraus folgt weiter, dass $N$ auf $M$ scharf transitiv
 operiert, so dass $|M| = |N|$ ist. Dies impliziert, dass $N$ ein
 minimaler Normalteiler ist. Folglich ist $N$ charakteristisch
 einfach und als abelsche Gruppe dann eine elementarabelsche
 $p$-Gruppe.
 \par
       Es sei $N^*$ ein von $N$ verschiedener minimaler Normalteiler.
 Dann ist $N^* \cap N = \{1\}$ und folglich $N^* \subseteq C_\Gamma(N) = N$.
 Dieser Widerspruch beweist die Einzigkeit von $N$. Damit ist alles bewiesen.
 \medskip\noindent
 {\bf 8.4. Satz.} {\it Es sei $n$ eine ungerade nat\"urliche Zahl
 und $\Gamma$ sei eine Permutationsgruppe der Ordnung $2n$, die auf
 der Menge $M$ der L\"ange $n$ transitiv operiere. Ist der
 Stabilisator irgend zweier Elemente von $M$ in $\Gamma$ stets
 trivial, so gilt:
 \item{a)} Sind $a$ und $b$ zwei verschiedene Elemente von $M$, so
 gibt es eine Involution in $\Gamma$, welche $a$ und $b$
 vertauscht.
 \item{b)} Ist $\Delta$ eine Untergruppe von $\Gamma$, die genau $r \geq 1$
 Involutionen enth\"alt, so ist $|\Delta| = 2r$.\par}
 \smallskip
       Beweis. Da $n$ ungerade ist und weil der Stabilisator zweier
 verschiedener Elemente von $M$ in $\Gamma$ stets trivial ist, hat
 jede Involution aus $\Gamma$ genau einen Fixpunkt. Wegen
 $|\Gamma_x| = 2$ f\"ur alle $x \in M$ enth\"alt $\Gamma$ daher
 genau $n$ Involutionen. Sind $\gamma$ und $\delta$ zwei
 Involutionen aus $\Gamma$ und ist $a^\gamma = a^\delta$ f\"ur ein $a \in M$, so
 ist $\gamma = \delta$. Ist n\"amlich $a = a^\gamma = a^\delta$, so ist
 $\gamma$, $\delta \in \Gamma_a$ und
 aus $|\Gamma_a| = 2$ folgt $\gamma = \delta$. Ist $a \neq a^\gamma$,
 so folgt $a^{\gamma\delta} = a^{\delta^2} = a$ und
 $(a^\gamma)^{\gamma \delta} = a^\delta = a^\gamma$. Somit l\"asst
 $\gamma \delta$ sowohl $a$ als auch $a^\gamma$ fest. Nach
 Voraussetzung ist daher $\gamma \delta = 1$ und folglich $\gamma = \delta$.
 Hieraus folgt
 $$ \big|\{a^\gamma \mid \gamma\ \hbox{\rm ist Involution aus}\ \Gamma\}\big|
	 = n. $$
 Daher ist 
 $$ M = \{a^\gamma \mid \gamma\ \hbox{\rm ist Involution aus}\ \Gamma\}, $$
 so dass es in der Tat eine Involution gibt, die $a$ mit $b$ vertauscht.
 \par
       Es sei $\delta$ eine Involution aus $\Delta$ und $a$ sei der
 Fixpunkt von $\Delta$. Mit $T$ bezeichnen wir die Bahn von $a$
 unter $\Delta$. Dann ist $|\Delta| = |T||\Delta_a|$. Wegen
 $1 < |\Delta_a|\leq|\Gamma_a| = 2$ ist daher $|\Delta| = 2|T|$. Nun hat
 $\delta$ keinen Fixpunkt in $T - \{a\}$, so dass $|T|$ ungerade ist.
 Das Paar $(\Delta,T)$ erf\"ullt somit ebenfalls die
 Voraussetzungen von 8.4. Wie der Beweis von a) zeigt, enth\"alt
 $\Delta$ daher genau $|T|$ Involutionen, so dass also $|T| = r$ ist.
 Damit ist der Satz bewiesen.
 \medskip\noindent
 {\bf 8.5. Satz.} {\it Es sei $M$ eine Menge mit $n + 1$ Elementen
 und $n$ sei ungerade. Ferner sei $\Gamma$ eine auf $M$ zweifach
 transitive Per\-mu\-ta\-ti\-ons\-grup\-pe. Ist $a \in M$ und enth\"alt
 $\Gamma_a$ eine Untergruppe $F_a$ gerader Ordnung, die auf
 $M-\{a\}$ transitiv operiert und die ferner die Eigenschaft hat,
 dass $F_{a,b,c} = \{1\}$ ist f\"ur alle $b$, $c \in M - \{a\}$ mit $b \neq c$,
 so operiert $\Gamma_a$ auf $M - \{a\}$ primitiv.}
 \smallskip
       Beweis. $F_a$ operiert als Frobeniusgruppe auf $M - \{a\}$. Daher
 besitzt $F_a$ eine charakteristische Untergruppe $K_a$, den
 Frobeniuskern von $F_a$, die auf $M - \{a\}$ scharf transitiv
 operiert (siehe etwa Huppert 1967, Satz 8.2, S. 496). Es gibt ferner
 eine Involution $\gamma \in F_a$, da die Ordnung von $F_a$
 gerade ist. Wegen $|K_a| = n \equiv 1 \mod 2$ ist $\gamma \in K_a$.
 Folglich ist $K_a\langle \gamma \rangle$ eine Untergruppe der
 Ordnung $2n$, die ebenfalls als Frobeniusgruppe auf $M - \{a\}$
 operiert. Wir d\"urfen daher annehmen, dass $|F_a| = 2n$ ist.
 Aufgrund der Transitivit\"at von $\Gamma$ auf $M$ enth\"alt dann
 $\Gamma_b$ f\"ur alle $b \in M$ eine solche Frobeniusgruppe $F_b$.
 \par
      Wir nehmen nun an, dass $\Gamma_a$ auf $M - \{a\}$ imprimitiv
 operiert. Es sei $I := \{a_1, \dots, a_t\}$ ein Imprimitivit\"atsgebiet von
 $\Gamma_a$. Dann ist $1 < t < n$. Ferner ist $t$ Teiler von $n$,
 so dass $t$ insbesondere ungerade ist. Wir setzen $T := I \cup \{a\}$. Es sei
 $b \in M - T$. Nach 8.4a) gibt es f\"ur $i := 1, \dots, t$ eine Involution
 $\gamma_i \in F_b$ mit
 $a^{\gamma_i} = a_i$. Wir zeigen, dass $T^{\gamma_i} = T$ ist.
 Offensichtlich ist $a^{\gamma_i} \in T$. Es sei $c := a_j^{\gamma_i}
 \not\in T$. Dann ist $j \neq i$ da andernfalls $c = a \in T$ w\"are.
 Wiederum nach 8.4a) gibt es eine Involution $\gamma \in F_c$ mit
 $a^\gamma = a_i$. Hieraus folgt, $a^{\gamma_i \gamma} = a_i^\gamma = a$ und
 somit $\gamma_i \gamma \in \Gamma_a$. Ferner ist
 $a_i^{\gamma_i \gamma} = a^\gamma = a_i$. Folglich ist $I^{\gamma_i
 \gamma} \cap I \neq \emptyset$ und daher $I^{\gamma_i \gamma} = I$. Hieraus
 folgt $c = c^\gamma = a_j^{\gamma_i \gamma} \in I$ und
 damit der Widerspruch $c \in T$. Dieser Widerspruch zeigt, dass
 $T^{\gamma_i} \subseteq T$ und damit, dass $T^{\gamma_i} = T$ ist.
 \par
       Wir betrachten nun die Gruppe $\Delta := \langle \gamma_i \mid i := 1,
 \dots, t\rangle$. Dann ist $\Delta$ eine Untergruppe von $F_b$,
 die auf $T$ transitiv operiert. Weil $|T| = t + 1$ gerade ist, muss
 $\Delta$ als Untergruppe von $F_b$ scharf transitiv operieren.
 Also ist $|\Delta| = t + 1$. Andererseits enth\"alt $T$ mindestens $t$
 Involutionen. Nach 8.4b) ist daher $|\Delta| \geq 2t$. Folglich
 ist $t + 1 \geq 2t$, woraus der Widerspruch $t = 1$ folgt. Damit ist
 der Satz bewiesen.
 \medskip\noindent
 {\bf 8.6. Satz.} {\it Die Voraussetzungen seien die gleichen wie
 in 8.5. Ist der Frobeniuskern $K_a$ von $F_a$ normal in
 $\Gamma_a$, so ist $K_a$ eine elementarabelsche $p$-Gruppe und $n$
 ist eine Potenz von $p$.}
 \smallskip
       Beweis. Auf Grund von 8.5 und 8.3 gen\"ugt es zu zeigen, dass
 $K_a$ abelsch ist.
       Wie wir wissen, gibt es eine Involution $\sigma \in \Gamma_a$, die
 genau einen Fixpunkt $b \in M - \{a\}$ hat. Weil $K_a$ in $\Gamma_a$ normal ist,
 ist $K^\sigma_a = K_a$. Es sei $\kappa \in K_a$ und $\kappa^\sigma = \kappa$.
 Dann ist $\sigma$ mit $\kappa$ vertauschbar und folglich $b^\kappa = b$, woraus
 $\kappa = 1$ folgt. Ist $\kappa^\sigma \kappa^{-1}
 = \lambda^\sigma \lambda^{-1}$ mit $\kappa$, $\lambda \in K_a$, so ist
 $(\lambda^{-1} \kappa)^\sigma = \lambda^{-1} \kappa$ und daher
 $\kappa = \lambda$. Also ist
 $$ \big|\{\kappa^\sigma \kappa^{-1} \mid \kappa \in K_a\}\big| = n, $$
 so dass
 $$ K_a = \{\kappa^\sigma \kappa^{-1} \mid \kappa \in K_a\} $$
 ist.  Es sei nun $\lambda \in K_a$. Dann gibt es also ein $\kappa \in K_a$ mit
 $\lambda = \kappa^\sigma \kappa^{-1}$. Daher ist
 $$ \lambda^\sigma = \kappa^{\sigma^2} \kappa^{-\sigma}
                   = \kappa \kappa^{-\sigma}
		   = (\kappa^\sigma \kappa^{-1})^{-1} = \lambda^{-1}. $$
 Da dies f\"ur alle $\lambda \in K_a$ gilt, folgt, dass $K_a$ abelsch ist.
 
\bigskip

\mysection{9. Geradenhomogene affine Ebenen}

\medskip

\noindent
 Hauptziel
 dieses Abschnitts ist zu zeigen, dass eine endliche affine Ebene, deren
 Kollineationsgruppe die Geraden transitiv permutiert, eine Translationsebene
 ist. Um dies zu beweisen, ben\"otigen wir noch einige Hilfsmittel, die wir zuvor
 bereitstellen.
 \medskip\noindent
 {\bf 9.1. Satz.} {\it Ist $E$ eine projektive Ebene und ist
 $\gamma$ eine involutorische Kollineation von $E$, so gilt eine
 der beiden folgenden Aussagen:
 \item{a)} $\gamma$ ist eine Perspektivit\"at von $E$.
 \item{b)} Die aus den Fixpunkten und Fixgeraden von $\gamma$ bestehende
 Teilstruktur von $E$ ist eine projektive Ebene $F$. Jede Gerade von $E$
 tr\"agt einen Punkt von $F$ und jeder Punkt von $E$ liegt auf einer Geraden
 von $F$.\par}
 \smallskip
       Beweis. Wir f\"uhren den Beweis in mehreren Schritten.\\
      1) Ist $P$ ein Punkt von $E$ mit $P \neq P^\gamma$, so ist $P + P^\gamma$
 eine Fixgerade von $\gamma$.
 \smallskip
       Es ist ja
 $$ (P + P^\gamma)^\gamma = P^\gamma + P^{\gamma^2} = P + P^\gamma. $$\\
       2) $\gamma$ hat mindestens zwei Fixpunkte.
 \smallskip
       Weil $\gamma \neq 1$ ist, gibt es eine Gerade $G$ mit $G^\gamma \neq G$.
 Dann ist aber $P := G \cap G^\gamma$ nach der zu 1) dualen Aussage ein Fixpunkt
 von $\gamma$. W\"are $P$ der einzige Fixpunkt von $\gamma$, so folgte, dass
 $\gamma$ alle Geraden von $E$, die nicht durch $P$ gehen, fest lie\ss e. Mit
 I.8.8, angewandt auf $(E^d)_P$, folgte $\gamma = 1$. Dieser Widerspruch zeigt
 die G\"ultigkeit von 2).
 \smallskip\noindent
       3) Jeder Punkt von $E$ liegt auf einer Fixgeraden von $\gamma$.
 \smallskip
       Ist $P$ ein Punkt von $E$, der kein Fixpunkt ist, so ist $P + P^\gamma$
 nach 1) eine Fixgerade durch $P$. Ist $P$ Fixpunkt, so gibt es nach 2) einen
 weiteren Fixpunkt $Q$. Dann ist $P + Q$ eine Fixgerade durch $P$.

       a) Wir nehmen nun an, dass $\gamma$ keinen Rahmen punktweise
 festl\"asst. Es seien $P$, $Q$ und $R$ drei nicht
 kollineare Fixpunkte von $\gamma$ und $X$ sei ein Punkt, der auf
 keiner der Geraden $P + Q$, $Q + R$, $R + P$ liegt. Nach 3) gibt es eine
 Fixgerade $G$ durch $X$. Diese ist von den Geraden $P + Q$, $Q + R$, $R + P$
 verschieden, da $X$ auf keiner dieser Geraden liegt. Weil die
 Fixpunktmenge von $\gamma$ keinen Rahmen enth\"alt, geht $G$ durch
 einen der Punkte $P$, $Q$, $R$. Wir d\"urfen annehmen, dass $R \leq G$
 ist. Ist $H$ eine Fixgerade von $\gamma$, die von $P + Q$, $Q + R$, $R + P$
 verschieden ist, so geht $H$ ebenfalls durch einen der Punkte, $P$, $Q$ und $R$
 und aus der Tatsache, dass die Fixpunktmenge von $\gamma$ keinen Rahmen
 enth\"alt, sowie aus $R \leq G$ folgt, dass auch $R \leq H$ ist. Es sei
 schlie\ss lich $J$ eine Gerade durch $R$. Ist $J = P + R$ oder $J = Q + R$, so
 ist $J$ eine Fixgerade von $\gamma$. Ist $J \neq P + Q$, $P + R$, so gibt es
 einen Punkt $Y \leq J$ mit $Y \neq R$ und $Y \leq P + Q$. Nun ist $Y$ nach 3) in
 einer Fixgeraden $L$ enthalten. Ferner ist, wie wir bereits wissen, $R \leq L$,
 woraus $L = J$ folgt. Damit ist gezeigt, dass $\gamma$ eine
 Zentralkollineation mit dem Zentrum $R$ ist.
 \par
       Wir d\"urfen des Weiteren annehmen, dass alle Fixpunkte von
 $\gamma$ kollinear sind. Weil $\gamma $ nach 2) mindestens zwei
 Fixpunkte hat, gibt es eine Fixgerade $G$ von $\gamma$, die alle
 Fixpunkte tr\"agt. Die von $\gamma$ in $E_G$ induzierte
 Kollineation hat dann keine Fixpunkte in $E_G$, so dass jeder
 Punkt von $E_G$ wegen 3) auf genau einer Fixgeraden liegt. Dies
 besagt wiederum, dass die zu $E_G$ geh\"orenden Fixgeraden von
 $\gamma$ gerade die Geraden einer Parallelschar von $E_G$ sind, so
 dass $\gamma$ als Elation erkannt ist.

       b) Die Menge der Fixpunkte von $\gamma$ enthalte einen
 Rahmen. Dann ist klar, dass die Menge $F$ der Fixpunkte und
 Fixgeraden von $\gamma$ eine Unterebene von $E$ bilden. Weil
 $\gamma$ auch eine involutorische Kollineation von $E^d$ ist,
 folgt, dass au\ss er 3) auch die zu 3) duale Aussage gilt, so dass
 also jeder Punkt von $E$ auf einer Geraden von $F$ und jede Gerade
 von $E$ durch einen Punkt von $F$ geht. Damit ist alles bewiesen.
 \medskip
       Ist $E$ eine projektive Ebene und ist $F$ eine von $E$ verschiedene
 Unterebene von $F$, so hei\ss t $F$
 {\it Baerunterebene\/}\index{Baerunterebene}{} von $F$, wenn jeder Punkt von $E$
 auf einer Geraden von $F$ liegt und jede Gerade von $E$ einen Punkt von $F$
 tr\"agt. Eine involutorische Kollineation einer projektiven Ebene,
 die eine Baerunterebene elementweise festl\"asst, hei\ss t
 {\it Baerinvolution\/}.\index{Baerinvolution}{}
 \par
       Baerunterebenen und Baerinvolutionen kann es bei endlichen
 projektiven Ebenen nur geben, wenn die Ordnung der Ebene ein
 Quadrat ist. Dies besagt der folgende Satz.
 \medskip\noindent
 {\bf 9.2. Satz.} {\it Es sei $E$ eine endliche projektive Ebene
 der Ordnung $n$. Ist $F$ eine von $E$ verschiedene Unterebene von
 $E$ der Ordnung $r$ und liegt jeder Punkt von $E$ auf einer
 Geraden von $F$, so ist $n = r^2$.}
 \smallskip
       Beweis. $F$ hat $r^2 + r + 1$ Punkte und ebenso viele Geraden. Da
 jeder Punkt von $E$, der nicht zu $F$ geh\"ort, auf einer und dann
 auch auf nur einer Geraden von $F$ liegt, ist die Anzahl der
 Punkte von $E$ gleich
 $$ r^2 + r + 1 + (r^2 + r + 1)(n + 1 - r - 1). $$
 Andererseits ist die Anzahl auch gleich $n^2 + n + 1$. Also ist
 $$\eqalign{
       0 &= n^2 + n + 1 - r^2 - r - 1 - (r^2 + r + 1)(n - r)  \cr
         &= (n - r)(n + r + 1 - r^2 - r - 1)                  \cr
	 &= (n - r)(n - r^2). \cr} $$
 Weil $F \neq E$ ist, ist $r \neq n$ und daher $n = r^2$, q. e. d.
 \medskip\noindent
 {\bf 9.3. Korollar} {\it Es sei $E$ eine endliche projektive Ebene
 der Ordnung $n$ und $\Gamma$ sei eine Kollineationsgruppe von $E$,
 deren Ordnung eine Potenz von $2$ ist. Sind die Fixpunkte von
 $\Gamma$ die Punkte einer Unterebene der Ordnung $m$ von $E$, so
 gibt es eine ganze Zahl $a$ mit $n = m^{2^a}$.}
 \smallskip
       Beweis. Die Fixpunkte von $\Gamma$ seien die Punkte der Unterebene
 $U$. Ist $\Gamma = \{1\}$, so ist $U = E$ und daher $n = m$. Es sei also
 $\Gamma \neq \{1\}$. Es gibt dann eine Involution $\rho \in
 Z(\Gamma)$. Weil $\Gamma$ einen Rahmen punktweise festl\"asst,
 l\"asst $\rho$ nach 9.1 und 9.2 eine Unterebene $F$ der Ordnung
 $r$ mit $r^2 = n$ punktweise fest. Da $\rho$ im Zentrum von $\Gamma$
 liegt, induziert $\Gamma$ eine Kollineationsgruppe $\Gamma^*$ auf
 $F$. Da alle Fixpunkte von $\Gamma$ auch Fixpunkte von $\rho$
 sind, ist $U$ eine Unterebene von $F$. Nun ist $\rho$ im Kern des
 Homomorphismus $*$ von $\Gamma$ auf $\Gamma^*$ enthalten, so dass
 $|\Gamma^*| < |\Gamma|$ ist. Nach Induktionsannahme gibt es folglich eine
 nat\"urliche Zahl $a$ mit $r = m^{2^{a-1}}$. Daher ist $n = m^{2^a}$, q. e. d.
 \medskip\noindent
 {\bf 9.4. Satz.} {\it Es sei $E$ eine projektive Ebene und $P$, $Q$
 und $R$ seien drei nicht kollineare Punkte von $E$. Ist $\rho$
 eine involutorische Perspektivit\"at aus $\Delta(P,Q + R)$ und $\sigma$ eine
 involutorische Perspektivit\"at aus $\Delta(Q,R + P)$, so ist
 $\rho \sigma = \sigma \rho$ und $\rho \sigma$ ist eine
 involutorische Streckung aus $\Delta(R, P + Q)$.}
 \smallskip
       Beweis. Es ist $\rho^{-1} \sigma \rho \in \Delta(Q^\rho,(R + P)^\rho)
 = \Delta (Q,R + P)$ und daher
 $$ \sigma^{-1} \rho^{-1} \sigma \rho \in \Delta(Q,R + P). $$
 Ferner ist $\sigma^{-1} \rho^{-1} \sigma \in \Delta(P^\sigma,(Q + R)^\sigma)
 = \Delta(P,Q + R)$, woraus
 $$ \sigma^{-1}\rho^{-1}\sigma\rho \in \Delta(P,Q + R) $$
 folgt. Also hat
 $\sigma^{-1} \rho^{-1} \sigma \rho$ die beiden verschiedenen
 Zentren $P$ und $Q$ und ist folglich nach II.1.1 gleich 1. Also
 ist $\rho \sigma = \sigma \rho$. Hieraus folgt weiter $(\rho\sigma)^2 = 1$.
 Weil $P \neq Q$ ist, ist \"uberdies $\rho \sigma \neq 1$. Folglich hat
 $\rho \sigma$ auf $Q + R$ nur die Fixpunkte
 $Q$ und $R$. Hieraus folgt, dass $\rho \sigma$ keine
 Baerinvolution ist, da sie sonst mindestens drei Fixpunkte auf $Q + R$
 haben m\"usste. Also ist $\rho \sigma$ nach 9.1 eine
 Zentralkollineation. Da $Q + R$ eine Fixgerade ist, die von der
 Achse von $\rho \sigma$ verschieden ist, geht sie durch das
 Zentrum von $\rho \sigma$. Nun wirkt $\rho \sigma$ auf den Punkten
 von $R + P$ wie $\rho$, da $\sigma$ alle Punkte dieser Geraden
 festl\"asst. Also ist auch $R + P$ eine Gerade durch das Zentrum
 von $\rho \sigma$, so dass $R$ das Zentrum von $\rho \sigma$ ist.
 Weil schlie\ss lich $P + Q$ eine Fixgerade von $\rho \sigma$ ist,
 die nicht durch $R$ geht, ist $P + Q$ die Achse von $\rho \sigma$.
 Damit ist alles bewiesen.
 \medskip
       Der n\"achste Satz ist wieder rein technischer Natur.
 \medskip\noindent
 {\bf 9.5. Satz.} {\it Es sei $E$ eine endliche projektive Ebene
 und $(P,G)$ sei ein nicht inzidentes Punkt-Geradenpaar von $E$.
 Gibt es f\"ur jedes Punktepaar $A$, $B$ mit $A \neq B$ und $A$, $B \leq G$ eine
 involutorische Streckung mit dem Zentrum $A$ und der
 Achse $B + P$, so hat $E$ Prim\-zahl\-po\-tenz\-ord\-nung.}
 \smallskip
       Beweis. Die Ordnung von $n$ ist ungerade, da $E$ involutorische
 Streckungen besitzt.
       Es sei $\Pi$ die von allen involutorischen Streckungen, deren Zentren auf
 $G$ liegen und deren Achsen durch $P$ gehen, erzeugte
 Kol\-li\-ne\-a\-ti\-ons\-grup\-pe von
 $E$. Nach 7.1 ist dann $\E(B,B + P) \subseteq \Pi$ und $|\E(B, B + P)| = n$
 f\"ur alle Punkte $B$ auf $G$. Somit ist $\Pi$ zweifach transitiv auf der Menge
 der Punkte von $G$. Schlie\ss lich ist, falls $A$ und $B$ zwei verschiedene
 Punkte auf $G$ sind, $\Pi(A,B + P)\E(B,B + P)$ eine Frobeniusgruppe gerader
 Ordnung und $\E(B,B + P)$ ist ein Normalteiler von $\Pi_B$. Aus 8.6 folgt daher
 die Behauptung.
 \medskip\noindent
 {\bf 9.6. Satz.} {\it Es sei $E$ eine endliche projektive Ebene, $U$ sei eine
 Gerade von $E$ und $\Delta$ sei eine Kollineationsgruppe von $E_U$. Gilt dann
 \item{a)} Ist $(P,G)$ eine Fahne von $E_U$, so gibt es eine involutorische
 Streckung $\gamma \in \Delta$ mit $P^\gamma = P$ und $G^\gamma = G$,
 \item{b)} Sind $A$ und $B$ Punkte auf $U$, sind $G$ und $H$ Geraden von $E_U$ mit
 $B \leq G, H$ und gibt es eine involutorische Streckung $\gamma \in \Delta$
 mit $A^\gamma = A$, $B^\gamma = B$ und $G^\gamma = G$, so gibt es auch eine
 involutorische Streckung $\delta \in \Delta$ mit $A^\delta = A$, $B^\delta = B$
 und $H^\delta = H$,\par
 \noindent
 so ist die Ordnung von $E$ eine Primzahlpotenz.}
 \smallskip
       Beweis. Wir betrachten zun\"achst den Fall, dass $\Delta$ keine
 involutorische Stre\-ckung enth\"alt, deren Zentrum ein Punkt von $E_U$ ist. Es
 sei $\gamma$ eine nach a) existierende involutorische Streckung aus $\Delta$.
 Ferner sei $A$ das Zentrum und $G$ die Achse von $\gamma$. Dann ist $A \leq U$
 und $G \neq U$. Setze $B := G \cap U$. Weil $\gamma$ die Punkte $A$ und $B$
 sowie die Gerade $G$ festl\"asst, gibt es nach b) eine involutorische
 Streckung $\delta \in \Delta$ mit $A^\delta = A$, $B^\delta = B$ und
 $H^\delta = H$, wobei $H$ eine beliebig vorgegebene, jedoch von $U$ verschiedene
 Gerade durch $B$ ist.  Wegen $\delta \in \Delta$ ist entweder $A$ oder $B$ das
 Zentrum von $\delta$. W\"are $B$ das Zentrum von $\delta$ so l\"age $A$ nach
 II.1.3 auf der Achse $J$ von $\delta$. Nach 9.4 w\"are dann $\gamma \delta$ eine
 involutorische Streckung mit dem Zentrum $G \cap J$ und der Achse $A + B = U$.
 Somit enthielte $\Delta$ doch eine involutorische Streckung mit einem Zentrum in
 $E_U$. Dieser Widerspruch zeigt, dass $A$ das Zentrum von $\delta$ ist. Dann ist
 aber $H$ nach II.1.3 die Achse von $\delta$. Aus dem zu Satz 7.1 dualen Satz
 folgt daher $|\E(A, U)| = n$ Es sei $C$ ein von $A$ verschiedener Punkt auf
 $U$. Gibt es in $\Delta$ eine involutorische Streckung mit dem Zentrum $C$, so
 ist auch $|\E(C,U)| = n$. Wegen $A \neq C$ ist $E_U$ eine Translationsebene, so
 dass $n$ nach II.1.16 Potenz einer Primzahl ist. Wir d\"urfen daher annehmen,
 dass $A$ Zentrum aller involutorischen Streckungen aus $\Delta$ ist. Es sei nun
 $G$ eine Gerade von $E$, die nicht durch $A$ geht. Nach a) gibt es eine
 involutorische Streckung $\rho \in \Delta$, die $G$ invariant l\"asst. Weil
 $A$ das Zentrum von $\rho$ ist, ist $G$ nach II.1.3 die Achse von $\rho$. Nach
 7.1 ist daher $E^d$ eine Translationsebene bez. $A$, so dass $n$ auch in
 diesem Falle Potenz einer Primzahl ist.
 \par
       Wir betrachten nun den Fall, dass es genau einen Punkt in $E_U$
 gibt, der Zentrum einer involutorischen Streckung aus $\Delta$
 ist. Es sei $P$ dieser Punkt. Dann ist $P$ ein Fixpunkt von
 $\Delta$. Es seien $A$ und $B$ zwei verschiedene Punkte auf $U$
 und $X$ sei ein von $P$ und $B$ verschiedener Punkt auf $P + B$.
 Nach a) gibt es eine involutorische Streckung $\rho \in \Delta$
 mit $X^\rho = X$ und $(X + A)^\rho = X + A$. Wegen $P^\rho = P$ und
 $U^\rho = U$ ist dann $(P + X)^\rho = P + X$ sowie
 $$ B^\rho = \bigl(U \cap (P + X)\bigr)^\rho = U^\rho \cap (P + X)^\rho
	   = U \cap (P + X) =  B. $$
 Somit hat $\rho$ drei verschiedene Fixpunkte auf $P + B$, so dass $P + B$
 die Achse von $\rho$ ist. Schlie\ss lich ist
 $$A^\rho = \bigl(U \cap (A + X)\bigr)^\rho = U^\rho \cap (A + X)^\rho = U \cap (A + X)
	  = A$$
 und $A \not\leq P + B$, so dass $A$ nach II.1.3 das Zentrum von $\rho$ ist.
 Anwendung von 9.5 liefert, dass $n$ auch in diesem Falle Potenz
 einer Primzahl ist.
 \par
       Es bleibt der Fall zu betrachten, dass wenigstens zwei Punkte in
 $E_U$ Zentren von involutorischen Streckungen aus $\Delta$ sind.
 \par
       Wir nehmen zun\"achst weiter an, dass jede Gerade von $E_U$ ein in
 $E_U$ liegendes Zentrum einer involutorischen Streckung aus
 $\Delta$ tr\"agt. Mit $k_G$ bezeichnen wir die Anzahl dieser
 Zentren auf der Geraden $G$. Dann ist $k_G \geq 1$ f\"ur alle
 Geraden $G$ von $E_U$ und $k_G \geq 2$ f\"ur wenigstens eine
 Gerade. Da $n(n + 1)$ die Anzahl der Geraden von $E_U$ ist, ist
 $$ \sum_G k_G > n(n + 1). $$
 Andererseits gilt nach 1.2 a) die Gleichung
 $$ v(n + 1) = \sum_G k_G, $$
 wenn $v$ die Anzahl der Punkte von $E_U$ ist, die Zentren von involutorischen
 Streckungen aus $\Delta$ sind. Also ist $v \geq n + 1$. Hieraus folgt mit 7.1
 und 7.7, dass $n$ Potenz einer Primzahl ist.
 \par
       Wir d\"urfen des Weiteren annehmen, dass es eine Gerade $G$ von
 $E_U$ gibt, die kein Zentrum einer involutorischen
 Perspektivit\"at aus $\Delta$ tr\"agt. Setze $A := G \cap U$. Es
 sei $B$ ein von $A$ verschiedener Punkt auf $U$. Ferner sei $P$
 ein Punkt von $E_U$, der Zentrum einer involutorischen Streckung
 aus $\Delta$ ist. Es gibt dann eine involutorische Streckung in $\Delta$ mit
 den Fixelementen $A$, $B$ und $A + P$. Nach b) gibt es
 dann auch eine involutorische Streckung $\rho \in \Delta$ mit den
 Fixelementen $A$, $B$ und $G$. W\"are $U$ die Achse von $\rho$, so
 w\"are wegen $G^\rho = G \neq U$ das Zentrum von $\rho$ ein von
 $A$ verschiedener Punkt auf $G$, im Widerspruch zu der Wahl von
 $G$. Also sind $A$ und $B$ die einzigen Fixpunkte von $\rho$ auf
 $U$. Da $B$ beliebig gew\"ahlt war, folgt mit 7.10, dass
 $\Delta_{A,G}$ die von $A$ verschiedenen Punkte von $U$ transitiv
 permutiert.
 \par
       Enth\"alt $\Delta$ eine von 1 verschiedene Elation mit der
 Achse $U$ und dem Zentrum $B \neq A$, so folgt wiederum
 $|\E(U)|\geq n + 1$, so dass $n$ nach 7.7 Potenz einer Primzahl ist.
 Wir d\"urfen daher annehmen, dass alle von Eins verschiedenen
 Elationen aus $\Delta$, deren Achse $U$ ist, das Zentrum $A$
 haben. Nun enth\"alt $\Delta$ nach 7.2 eine von 1 verschiedene
 Elation mit der Achse $U$. Somit ist $\Delta (A,U) \neq \{1\}$, so
 dass $A$ ein Fixpunkt von $\Delta$ ist. Es sei $P$ ein von $A$
 verschiedener Punkt auf $G$ und $B$ sei ein von $A$ verschiedener
 Punkt auf $U$. Nach a) gibt es eine involutorische Streckung
 $\sigma \in \Delta$ mit $P^\sigma = P$ und $(P +  B)^\sigma  = P + B$.
 Dann ist aber $A^\sigma = A$ und $B^\sigma = B$. Weil $U$ nicht die
 Achse von $\sigma$ sein kann, $P$ w\"are ja sonst Zentrum, sind
 $A$ und $B$ die einzigen Fixpunkte von $\sigma$ auf $U$. Mit Hilfe
 von 7.10 folgt wiederum, dass $\Delta_{A,P}$ auf der Menge der von
 $A$ verschiedenen Punkte von $U$ transitiv operiert. Ist $B$ das
 Zentrum von $\sigma$, so ist $G = A + P$ die Achse von $\sigma$ und
 aus der Transitivit\"at von $\Delta_{A,P}$ auf der Menge der von
 $A$ verschiedenen Punkte von $U$ folgt, dass $|\Delta(A,G)| = n$
 ist. Ist $B$ nicht das Zentrum von $\sigma$, so ist $A$ das
 Zentrum und $P + B$ die Achse von $\sigma$. Die Transitivit\"at von
 $\Delta_{Q,P}$ auf der Menge der von $A$ verschiedenen Punkte von $U$ impliziert
 auch in diesem Falle, dass |$\Delta_{A,G}| = n$ ist. Nun ist aber
 $\Delta(A,U) \neq \{1\}$, wie wir bereits bemerkten, so dass
 $|\Delta(A)|\geq n + 1$ ist. Aus 7.7 folgt somit, dass auch im letzten Falle $n$
 Potenz einer Primzahl ist. Damit ist der Satz bewiesen.
 \medskip\noindent
 {\bf 9.7. Satz.} {\it Ist $E$ eine endliche projektive Ebene, ist\/
 $U$ eine Gerade von $E$ und ist $K$ eine Kollineationsgruppe von
 $E_U$, so sind die folgenden Aussagen \"aquivalent:
 \item{a)} $K$ operiert auf der Menge der Geraden von $E_U$ transitiv.
 \item{b)} $K$ operiert auf der Menge der Punkte von $E_U$ und auch auf der Menge
 der Punkte von $U$ transitiv.
 \item{c)} $K$ ist transitiv auf der Menge der Fahnen von $E_U$.}
 \smallskip
       Beweis. a) impliziert b): Dass die Geradentransitivit\"at von $K$
 auf $E_U$ die Transitivit\"at von $K$ auf der Menge der Punkte
 nach sich zieht, ist unmittelbar einsichtig. Dass $K$ auch auf der
 Menge der Punkte von $E_U$ transitiv ist, folgt aus 5.5.
 \par
       b) impliziert c): Es seien $(P,G)$ und $(Q,H)$ zwei Fahnen von $E_U$. Auf
 Grund der Transitivit\"at von $K$ auf der Menge der Punkte von $E_U$ d\"urfen
 wir annehmen, dass $P = Q$ ist. Es ist dann $G = P + (G \cap U)$ und
 $H = P + (H \cap U)$. Die Anzahl der Punkte von $E_U$ ist $n^2$ und die Anzahl
 der Punkte von $U$ ist $n + 1$.  Weil $n^2$ und $n + 1$ teilerfremd sind, ist
 $K_P$ nach 8.1 transitiv auf der Menge der Punkte von $U$. Es gibt also ein
 $\kappa \in K_P$ mit $(G \cap U)^\kappa = H \cap U$. Es folgt $P^\kappa = P$ und
 $G^\kappa = H$.
 \par
       c) impliziert a): Banal.
 \medskip
       Eine unmittelbare Folgerung ist
 \medskip\noindent
 {\bf 9.8. Korollar} {\it Es sei $E$ eine endliche projektive Ebene
 der Ordnung $n$. Ferner sei $U$ eine Gerade von $E$ und $K$ eine
 Kollineationsgruppe von $E_U$, die die Geraden von $E_U$ transitiv
 permutiert. Ist dann $(P,G)$ eine Fahne von $E_U$ und $A$ ein
 Punkt auf\/ $U$ und setzt man $k := |K_{P,G}|$, so ist $|K| = n^2(n + 1)k$,
 $|K_A| = n^2k$, $|K_G| = nk$, $|K_P| = (n + 1)k$ und $|K_{A,P}| = k$.}
 \medskip
 Ferner gilt
 \medskip\noindent
 {\bf 9.9. Satz.} {\it Es sei $E$ eine endliche projektive Ebene,
 $U$ sei eine Gerade von $E$ und $K$ sei eine Kollineationsgruppe
 von $E_U$. Genau dann ist $K$ auf der Menge der Punkte von $E_U$
 transitiv, wenn $K_A$ f\"ur jeden Punkt $A$ auf $U$ auf der Menge
 der mit $A$ inzidierenden Geraden von $E_U$ transitiv ist.}
 \smallskip
       Beweis. Es sei $s$ die Anzahl der in $U$ enthaltenen Punktbahnen von $K$.
 \par
       Es sei $K$ transitiv auf der Menge der Punkte von $E_U$. Dann hat
 $K$ genau $s + 1$ Punktbahnen in $E$ und daher nach 3.4 auch genau
 $s + 1$ Geradenbahnen. Weil $\{U\}$ eine dieser Geradenbahnen ist,
 zerf\"allt die Menge der Geraden von $E_U$ unter $K$ in genau $s$
 Bahnen. Es seien $\Pi_1$, \dots, $\Pi_s$ die in $U$ enthaltenen
 Punktbahnen von $K$. Wir setzen
 $$ \Gamma_i := \{G \mid G\ \hbox{\rm ist Gerade von}\ E_U\ \hbox{\rm mit}\
       G \cap U \in \Pi_i\} $$
 f\"ur $i := 1$, \dots, $s$. Offensichtlich ist jedes
 $\Gamma_i$ invariant unter $K$ und somit Vereinigung von
 Geradenbahnen von $K$. Weil $K$ in $E_U$ genau $s$ Geradenbahnen
 hat, folgt, dass die $\Gamma_i$ allesamt Geradenbahnen von $K$
 sind. Es sei schlie\ss lich $A$ ein Punkt auf $U$ und $G$ und $H$
 seien zwei mit $A$ inzidierende Geraden von $E_U$. Es gibt ein $i$
 mit $A \in \Pi_i$. Es folgt $G$, $H \in \Gamma_i$. Es gibt daher ein
 $\kappa \in K$ mit $G^\kappa = H$. Es folgt
 $$ A^\kappa = (G \cap U)^\kappa = G^\kappa \cap U^\kappa = H\cap U = A, $$
 so dass sogar $\kappa \in K_A$ gilt.
 \par
       Wir nehmen umgekehrt an, dass $K_A$ f\"ur alle $A \leq U$ die
 Menge der Geraden von $E_U$, die durch $A$ gehen, transitiv
 permutiert. $K$ zerlegt dann offensichtlich die Menge der Geraden
 von $E_U$ in $s$ Bahnen, so dass $K$ genau $s + 1$ Geradenbahnen in
 $E$ hat. Folglich hat $K$ auch genau $s + 1$ Punktbahnen in $E$. Da
 $s$ von diesen Punktbahnen in $U$ enthalten sind, besteht die
 restliche Bahn aus den Punkten von $E_U$. Damit ist der Satz bewiesen.
 \medskip\noindent
 {\bf 9.10. Korollar.} {\it Es sei $E$ endliche projektive Ebene,
 $U$ sei eine Ge\-ra\-de von $E$ und $K$ sei eine Kollineationsgruppe
 von $E_U$, die auf der Menge der Geraden von $E_U$ transitiv ist.
 Sind $A$ und $B$ Punkte von $U$ und sind $G$ und $H$ durch $B$
 gehende Geraden von $E_U$, so ist $|K_{A,B,G}| = |K_{A,B,H}|$.}
 \smallskip
       Beweis. Wegen $\ggT(n^2,n + 1) = 1$ folgt aus 9.7 und 8.1, dass $K_A$
 auf der Punktmenge von $E_U$ transitiv operiert. Aus 9.9 folgt
 daher die Transitivit\"at von $K_{A,B}$ auf der Menge der durch
 $B$ gehenden Geraden von $E_U$. Hieraus folgt die Behauptung.
 \medskip
       Die n\"achsten drei S\"atze sind die H\"ohepunkte dieses Abschnitts.
 9.12 wurde von Ostrom \& Wagner\index{Ostrom, T. G.}{} vor dem von
 Wagner\index{Wagner, A.}{} stammenden n\"achsten Satz
 bewiesen, aus dem er wiederum sehr einfach folgt.
 \medskip\noindent
 {\bf 9.11. Satz.} {\it Es sei $E$ eine endliche projektive Ebene
 und\/ $U$ sei eine Gerade von $E$. Ist $K$ eine Kollineationsgruppe
 von $E_U$, die die Ge\-ra\-den von $E_U$ transitiv permutiert, so ist
 $E$ eine Translationsebene bez.\ $U$ und $\E(U)$ ist in $K$
 enthalten.}
 \smallskip
       Beweis. Es sei $n$ die Ordnung von $E$. Ist $(P,G)$ eine Fahne von
 $E_U$, so setzen wir $k := |K_{P,G}|$.
 \par

 \smallskip

       1.\ Fall: $n$ ist gerade. Es sei $2^a$ die h\"ochste in $n$
 aufgehende Potenz von $2$. Dann ist $a \geq 1$. Es sei $\Sigma$
 eine 2-Sylowgruppe von $K$. Nach 9.8 ist $|\Sigma| = 2^{2a+b}$.
 Schlie\ss lich sei $\sigma$ eine Involution aus dem Zentrum
 $Z(\Sigma)$ von $\Sigma$. Ist $\sigma$ eine Perspektivit\"at, so
 ist $\sigma$, da $n$ gerade ist, eine Elation. Nach 9.1 ist
 $\sigma$ daher entweder eine Elation oder aber die Fixpunkte von
 $\sigma$ sind die Punkte einer Baerunterebene $F$ von $E$. Hat $F$
 die Ordnung $r$, so ist $r^2 = n$ nach 9.3. Sind die Fixpunkte von
 $\sigma$ die Punkte einer Baerunterebene, so hat $\sigma$ genau
 $r^2 + r + 1 = n + r + 1$ Fixpunkte. Wegen $U^\sigma = U$ ist $U$ eine
 Gerade von $F$, so dass auf $U$ genau $r + 1$ Fixpunkte von $\sigma$
 liegen, dh., genau $n$ der Fixpunkte von $\sigma$ liegen in
 $E_U$. Ist $\sigma$ eine Elation und ist $U$ nicht die Achse von
 $\sigma$, so hat $\sigma$ genau $n + 1$ Fixpunkte, von denen
 wiederum genau $n$ in $E_U$ liegen. Es gilt also, dass genau $n$
 der Fixpunkte von $\sigma$ in $E_U$ liegen, wenn wenigstens einer
 der Punkte von $E_U$ ein Fixpunkt von $\sigma$ ist.
 \par
       Es sei $S$ die Menge der in $E_U$ liegenden Fixpunkte von
 $\sigma$. Dann ist $S = \emptyset$ oder $|S| = n$, wie wir gerade
 gesehen haben. Wir nehmen an, dass $|S| = n$ ist. Wegen $\sigma \in Z(\Sigma)$
 ist $S$ unter $\Sigma$ invariant. Folglich ist $S$
 Vereinigung von Bahnen von $\Sigma$. Nach 8.2 ist die L\"ange
 jeder Bahn von $\Sigma$ durch $2^{2a}$ teilbar, da $K$ ja auf der
 Menge der $n^2$ Punkte von $E_U$ transitiv operiert. Also ist
 $2^{2a}$ Teiler von $|S| = n$. Hieraus folgt $2a \leq a$, was wegen
 $a \geq 1$ ein Widerspruch ist. Also ist $S = \emptyset$, was
 wiederum besagt, dass $\sigma$ eine Elation mit der Achse $U$ ist.
 Weil $K$ auf der Menge der Punkte von $U$ transitiv operiert, sind
 die Gruppen $K(P, U)$ mit $P \leq U$ in $K$ konjugiert und haben
 daher alle die Ordnung $h \geq 2$. Nach 7.9 ist folglich $h = n$
 und $K(P,U) = \E(P,U)$ f\"ur alle $P \leq U$. Damit ist der Satz im
 Falle, dass $n$ gerade ist, bewiesen.
 \par

 \smallskip

       2.\ Fall: $n$ ist ungerade. Wir zeigen zun\"achst, dass dann $n$ Potenz
 einer Primzahl ist.

 \smallskip

      Es sei $F$ eine Unterebene von $E$ mit den beiden Eigenschaften:
 \item{j)} Es gibt eine 2-Untergruppe von $K$, so dass die
 Fixpunkte dieser Untergruppe gerade die Punkte von $F$ sind.
 \item{ij)} Es gibt keine echte Unterebene von $F$, so dass die
 Punkte dieser Unterebene gerade die Fixpunkte einer 2-Untergruppe von $K$ sind.

 \noindent
       Weil die Gruppe, die nur aus der $1$ besteht, eine 2-Untergruppe
 von $K$ ist, hat $E$ die Eigenschaft j). Weil $E$ endlich ist,
 folgt daher, dass $E$ auch eine Unterebene $F$ besitzt, die die
 beiden Eigenschaften j) und ij) hat.
 
 Es sei $\Sigma$ eine 2-Untergruppe von $K$, die maximal ist
 bez\"uglich der Eigenschaft, dass die Fixpunkte von $\Sigma$ gerade
 die Punkte von $F$ sind. Ferner sei $2^a$ die h\"ochste in $n + 1$
 und $2^b$ die h\"ochste in $k$ aufgehende Potenz von $2$.
 Schlie\ss lich sei $|\Sigma| = 2^c$. Weil $n$ ungerade ist, ist $n + 1$
 gerade und daher $a \geq 1$.
 
       a) Ist $(P,G)$ eine Fahne von $F_U$, so gibt es eine
 involutorische Streckung $\gamma$ von $F_U$ mit $P^\gamma = P$ und
 $G^\gamma = G$, die von einer Kol\-li\-ne\-a\-ti\-on aus $K$ induziert wird.
 
       $\Sigma$ ist eine Untergruppe von $K_{P,G}$. Wegen $|K_{P,G}| = k$
 ist die Ordnung $2^c$ von $\Sigma$ ein Teiler von $2^b$. Nun ist
 $2^{a+b}$ die Ordnung einer 2-Sylowgruppe von $K$ und es ist $a
 \leq 1$. Daher ist $\Sigma$ keine 2-Sylowgruppe von $K$. Es gibt
 folglich eine 2-Untergruppe $\Sigma^*$ von $K$, die $\Sigma$
 enth\"alt und f\"ur die $|\Sigma^* : \Sigma| = 2$ ist. Hieraus
 folgt, dass $\Sigma$ ein Normalteiler von $\Sigma^*$ ist, was
 wiederum impliziert, dass $F$ von $\Sigma^*$ invariant gelassen
 wird. Aus der Maximalit\"at von $\Sigma$ folgt, dass $\Sigma^*$
 eine Kollineationsgruppe der Ordnung $2$ in $F$ induziert. Aus der
 Minimalit\"at von $F$ folgt weiter, dass das erzeugende Element
 $\sigma$ dieser Gruppe eine Perspektivit\"at ist. Ist $m$ die
 Ordnung von $F$, so ist $n$ nach 9.3 eine Potenz von $m$, so dass
 mit $n$ auch $m$ ungerade ist. Daher ist $\sigma$ eine Stre\-ckung
 von $F_U$. Hieraus folgt wiederum die Existenz einer Fahne $(Q,J)$
 von $F_U$ mit $Q^\sigma = Q$ und $J^\sigma = J$. Also ist $\Sigma^*
 \subseteq K_{Q,J}$, so dass wegen $|K_{Q,J}| = k$ die Ungleichung
 $2^c < 2^b$ gilt. Dies impliziert, dass $\Sigma$ keine
 2-Sylowgruppe von $K_{P,G}$ ist. Es gibt daher bereits in
 $K_{P,G}$ eine Gruppe $\Sigma^*$ mit $|\Sigma^* : \Sigma|=2$. Dann
 gibt es aber, wie wir gesehen haben, eine involutorische Streckung
 von $F_U$, die die Fahne $(P,G)$ invariant l\"asst und die von
 einer Kollineation aus $K$ induziert wird.
 \par
       b) $A$ und $B$ seien zwei auf $U$ liegende Punkte von $F$
 und $G$ und $H$ seien zwei Geraden von $F_U$ mit $B \leq G$, $H$.
 Gibt es dann eine involutorische Streckung von $F_U$ mit den
 Fixelementen $A$, $B$ und $G$, die von einer Kollineation aus $K$
 induziert wird, so gibt es eine ebensolche Streckung mit den
 Fixelementen $A$, $B$ und $H$.
 \par
       Es sei $\Lambda$ der Stabilisator von $F$ in $K$ und
 $\bar{\Lambda}$ sei die Untergruppe von $\Lambda$, die $F$
 elementweise festl\"asst. Dann ist $\Sigma \subseteq \bar\Lambda$
 und $\bar\Lambda$ ist ein Normalteiler von $\Lambda_{A,B,G}$. Aufgrund
 unserer Annahme enth\"alt $\Lambda_{A,B,G}/\bar\Lambda$ ein Element
 der Ordnung $2$. Wegen $\Lambda_{A,B,G} \subseteq K_{A,B,G}$ ist
 somit $\Sigma$ keine 2-Sylowgruppe von $K_{A,B,G}$. Nach 9.10 ist
 $|K_{A,B,G}| = |K_{A,B,H}|$, so dass $\Sigma$ auch keine
 2-Sylowgruppe von $K_{A,B,H}$ ist. Es gibt folglich eine
 2-Untergruppe $\Sigma^*$ von $K_{A,B,H}$ mit
 $|\Sigma^* : \Sigma| = 2$. Wie der Beweis von a) zeigt, induziert
 $\Sigma^*$ eine Kollineationsgruppe der Ordnung $2$ in $F_U$,
 deren erzeugendes Element eine involutorische Streckung mit den
 Fixelementen $A$, $B$ und $H$ ist. Damit ist auch b) bewiesen. Die
 G\"ultigkeit von a) und b) impliziert nun nach 9.6, dass $m$
 Potenz einer Primzahl ist. Weil schlie\ss lich $n$, wie wir
 bereits bemerkten, Potenz von $m$ ist, ist auch $n$ Potenz einer Primzahl.
 \par

 \smallskip

       Wir sind nun in der Lage den Beweis von 9.11 zu beenden. Wie wir
 gesehen haben, ist $n = p^r$ mit einer Primzahl $p$. Es sei $p^s$
 die h\"ochste in $k$ aufgehende Potenz von $p$. Ist dann $\Pi$ eine
 $p$-Sylowgruppe von $K$, so ist $|\Pi| = p^{2r+s}$. Weil $K$ auf der
 Menge der $p^{2r}$ Punkte von $E_U$ transitiv ist, ist auch $\Pi$
 nach 8.2 auf der Menge der Punkte von $E_U$ transitiv.
 Andererseits hat $\Pi$ auf $U$ einen Fixpunkt $A$, da die Anzahl
 $n + 1$ der Punkte auf $U$ nicht durch $p$ teilbar ist. Wegen
 $|\Pi| = p^{2r+s} > p^r$ gibt es einen von $A$ verschiedenen Punkt
 $B$ auf $U$ mit $\Pi_B \neq \{1\}$. Der Punkt $B$ sei so
 gew\"ahlt, dass $|\Pi_B| \geq |\Pi_Q|$ ist f\"ur alle von $A$
 verschiedenen Punkte $Q$ auf $U$. Es sei $C$ ein von $A$
 verschiedener Fixpunkt von $\Pi_B$. Wegen $\Pi_B = \Pi_{B,C} \subseteq \Pi_C$
 und $|\Pi_B| \geq |\Pi_C|$ ist dann $\Pi_B = \Pi_C$.
 Weil $\Pi$ auf der Menge der Punkte von $E_U$ transitiv ist, ist
 $\Pi_B$ auf der Menge der durch $C$ gehenden Geraden von $E_U$
 nach 9.9 transitiv. Es sei $\bar{\Pi}$ eine $p$-Sylowgruppe
 von $K_B$, die $\Pi_B$ enth\"alt. Dann ist $\Pi_B \subseteq \bar\Pi_A$. Weil
 $|K_B| = n^2k$ ist, ist $\bar\Pi$ sogar eine $p$-Sylowgruppe von $K$. Daher sind
 $\Pi$ und $\bar\Pi$ konjugiert. Aus der Maximalit\"at von $\Pi_B$
 folgt daher $\Pi_B = \bar{\Pi}_A$, so dass $\Pi_B$ auch auf der
 Menge der durch $A$ gehenden Geraden von $E_U$ transitiv ist.
 \par
       Es sei $1 \neq \tau \in Z(\Pi_B)$. Weil $\tau$ zwei Fixpunkte hat,
 n\"amlich $A$ und $B$, hat $\tau$ nach 3.2 auch zwei Fixgeraden.
 Es gibt daher eine Gerade $G$ von $E_U$ mit $G^\tau = G$. Gibt es
 ein Element $\mu \in \Pi_B$ mit $(G \cap U)^\pi \neq G \cap U$, so
 ist $G \cap G^\mu$ ein Punkt von $E_U$. Wegen $\tau \in Z(\Pi_B)$
 ist
 $$ (G \cap G^\mu)^\tau = G^\tau \cap G^{\mu \tau} = G \cap G^{\tau\mu}
                        = G \cap G^\mu. $$
 Daher ist $A + (G \cap G^\mu)$ eine Fixgerade von $\tau$, die $U$ in einem 
 Fixpunkt von $\Pi_B$, n\"amlich $A$, schneidet. Wir d\"urfen daher annehmen,
 dass $C := G \cap U$ ein Fixpunkt von $\Pi_B$ ist. Aus der Transitivit\"at von
 $\Pi_B$ auf der Menge der durch $C$ gehenden Geraden von $E_U$,
 folgt, dass $\tau$ alle Geraden durch $C$ einzeln invariant
 l\"asst, so dass $\tau$ zentral mit dem Zentrum $C$ ist. Ist
 $\lambda$ die Anzahl der Fixpunkte von $\Pi_B$ auf $U$, so ist
 $\lambda \equiv 1 \mod p$ und $\lambda \geq 2$. Daher
 ist sogar $\lambda \geq 3$. Somit hat $\tau$ auf $U$ mindestens
 drei Fixpunkte, so dass $U$ die Achse von $\tau$ ist. Also ist
 $K(C,U) \neq \{1\}$. Hieraus folgt wiederum genauso wie im Falle
 eines geraden $n$, dass $E$ eine Translationsebene bez. $U$ ist
 und dass $\E(U)$ in $K$ enthalten ist. Damit ist der Satz bewiesen.
 \medskip
       Wir benutzen nun den gerade bewiesenen Satz, um den folgenden von
 Ostrom \& Wagner stammenden Satz zu bewiesen, der, wie schon
 gesagt, der \"altere der beiden S\"atze ist.\index{Satz von Ostrom--Wagner}{}
 \medskip\noindent
 {\bf 9.12. Satz.} {\it Ist $E$ eine endliche projektive Ebene und
 ist $K$ eine Kollineationsgruppe von $E$, die auf der Menge der
 Punkte von $E$ zweifach transitiv operiert, so ist $E$
 desarguessch und $K$ enth\"alt die kleine projektive Gruppe von $E$.}
 \smallskip
       Beweis. Nach 3.5 ist $K$ auch auf der Menge der Geraden von $E$
 zweifach transitiv. Ist $U$ eine Gerade von $E$, so ist also $K_U$
 auf der Menge der Geraden von $E_U$ noch transitiv. Nach 9.11 ist
 $E_U$ folglich eine Translationsebene und $E(U)$ ist in $K_U$
 enthalten. Somit ist $E$ eine Moufangebene und $K$ enth\"alt die
 kleine projektive Gruppe von $E$. Weil $E$ endlich ist, ist $E$
 nach 6.1 desarguessch.
 \medskip\noindent
 {\bf 9.13. Satz.} {\it Es sei $E$ eine endliche projektive Ebene
 und $K$ sei
 eine Kollineationsgruppe von $E$. Ist $K$ auf der
 Menge der nicht inzidenten Punkt-Geradenpaare von $E$ transitiv,
 so ist $E$ desarguessch und $K$ enth\"alt die kleine projektive Gruppe von $E$.}
 \smallskip
       Beweis. Es sei $(P,G)$ eine Fahne von $E$. Dann ist $K_G$ auf der
 Menge der Punkte von $E_G$ transitiv, so dass $K_{P,G}$ nach 9.9
 auf der Menge der von $G$ verschiedenen Geraden durch $P$
 transitiv ist. Weil durch $P$ mindestens drei Geraden gehen und
 $G$ eine beliebige Gerade durch $P$ ist, folgt, dass $K_P$ auf der
 Menge aller Geraden durch $P$ transitiv operiert. Weil $K_P$ nach
 Voraussetzung auch auf der Menge der nicht durch $P$ gehenden
 Geraden transitiv operiert, hat $G_P$ genau zwei Geradenbahnen.
 Nach 3.4 hat $G_P$ dann auch genau zwei Punkt\-bah\-nen. Eine davon
 ist $\{P\}$. Dies besagt, dass $K$ auf der Menge der Punkte von
 $E$ zweifach transitiv operiert. 9.13 ist daher eine Folge von 9.12.

\mysectionten{10. Endliche projektive R\"aume}

\noindent
 Es sei 
 $T$ ein $2$-$(v,k,\lambda)$-Blockplan. Sind $P$ und $Q$ zwei
 verschiedene Punkte von $T$ und sind $c_1$, \dots, $c_\lambda$ die
 mit $P$ und $Q$ inzidierenden Bl\"ocke, so definieren wir die
 {\it Gerade\/}\index{Gerade}{} $P + Q$ durch
 $$ P + Q := \{X \mid X \I c_1,\ \dots,\ c_\lambda\}. $$
 Da durch zwei verschiedene Punkte von $T$ stets
 genau $\lambda$ Bl\"ocke gehen, ist jede Gerade durch irgend zwei
 auf ihr liegende Punkte eindeutig bestimmt.
 \par
       Sind $P$, $Q$ und $R$ drei nicht kollineare Punkte von $T$, so
 definieren wir die durch $P$, $Q$ und $R$ bestimmte {\it Ebene\/}\index{Ebene}{}
 $P + Q + R$ durch
 $$ P + Q + R := \{X \mid X \I c_1,\ \dots,\ c_i\}, $$
 wobei $c_1$, \dots, $c_i$ die mit $P$, $Q$ und $R$ inzidierenden Bl\"ocke sind.
 Ist $i = 0$, so ist $P + Q + R$ die Menge aller Punkte von $T$. Sind $P'$,
 $Q'$ und $R'$ drei nicht kollineare Punkte in $P + Q + R$, so ist
 $P' + Q' + R' \subseteq P + Q + R$ und es kann vorkommen, dass $P' + Q' + R'$
 echt in $P + Q + R$ enthalten ist.
 \medskip\noindent
 {\bf 10.1. Satz.} {\it Ist $T$ ein $2$-Blockplan mit den Parametern
 $v$, $b$, $k$, $r$ und $\lambda$, so sind die folgenden Bedingungen
 \"aquivalent:
 \item{a)} Es gibt einen endlichen projektiven Verband $L$ mit $\Rg(L) \geq 3$, so
 dass\/ $T$ und $L_{1,\Rg(L) - 1}$ isomorph sind.
 \item{b)} Es ist $v \geq k + 2$ und jede Gerade von $T$ hat mit
 jedem Block von $T$ einen Punkt gemeinsam.
 \item{c)} Es ist $v \geq k + 2$ und auf jeder Geraden von $T$ liegen
 genau ${b-\lambda \over r - \lambda}$ Punkte.}
 \smallskip
       Beweis. In jeder projektiven Geometrie, deren Rang mindestens 3
 ist, trifft jede Gerade jede Hyperebene und au\ss erhalb jeder
 Hyperebene liegt mehr als ein Punkt. Daher ist b) eine Folge von a).
 \par
       Um die \"Aquivalenz von b) und c) zu beweisen, sei $g$ eine Gerade
 von $T$. Durch jeden Punkt $P$ von $g$ gehen $r$ Bl\"ocke, von
 denen $\lambda$ die Gerade $g$ enthalten. Die restlichen $r - \lambda$
 Bl\"ocke treffen $g$ nur in $P$. Die Zahl der $g$ treffenden Bl\"ocke ist daher
 $$ |g|(r - \lambda) + \lambda. $$
 Also trifft $g$ genau dann jeden Block, wenn
 $$ |g|(r - \lambda) + \lambda = b $$
 ist, dh. genau dann, wenn
 $$ |g| = {b - \lambda \over r - \lambda} $$
 ist. Dies zeigt die Gleichwertigkeit von b) und c).
 \par
       Wir setzen jetzt die G\"ultigkeit von b) und c) voraus. Als erstes
 zeigen wir, dass die aus den Punkten und Geraden von $T$
 bestehende Inzidenzstruktur eine irreduzible projektive Geometrie
 ist. Durch zwei verschiedene Punkte von $T$ geht nat\"urlich genau
 eine Gerade. Auf jeder Geraden liegen mindestens drei Punkte.
 Andernfalls w\"are n\"am\-lich ${b - \lambda \over r - \lambda} = 2$, so
 dass jede Gerade genau zwei Punkte tr\"uge. Dann g\"abe es aber wegen
 $v \geq k + 2$ eine Gerade, die einen gegebenen Block nicht trifft.
 \par
       Es seien $P$, $Q$ und $R$ drei nicht kollineare Punkte und $\rho$
 sei die Anzahl der Bl\"ocke, die $P$, $Q$ und $R$ enthalten. Dann
 ist $r - \rho$ die Anzahl der Bl\"ocke durch $P$, die $P + Q + R$ nicht
 enthalten. Wegen b) trifft jeder dieser Bl\"ocke die Gerade $Q + R$
 in genau einem Punkt. Somit ist
 $$ r - \rho = |Q + R|(\lambda - \rho), $$
 so dass $\rho$ wegen $|Q + R| = {b - \lambda \over r - \lambda}$ von der Auswahl
 von $P$, $Q$ und $R$ unabh\"angig ist. Daher ist jede Ebene durch irgend drei
 in ihr liegende, nicht kollineare Punkte eindeutig bestimmt.
 \par
       Um zu zeigen, dass auch das Veblen-Young-Axiom gilt, gen\"ugt es
 zu zeigen, dass zwei Geraden, die in einer Ebene liegen, einen
 Schnittpunkt haben. Es seien also $g$ und $h$ zwei Geraden, die in
 einer Ebene $E$ und $T$ liegen. Wegen $\rho < \lambda$ gibt es
 einen Block $c$ mit $h \subseteq c$ und $E \not\subseteq c$. Weil
 eine Ebene durch irgend drei ihrer Punkte eindeutig bestimmt ist,
 folgt dass $E \cap c = h$ ist. Wegen b) gibt es einen Punkt $P \in g$ mit
 $P \I c$. Wegen $g \subseteq E$ ist daher $P \in E \cap c = h$. Somit haben $g$
 und $h$ den Punkt $P$ gemeinsam. Damit ist gezeigt, dass die Inzidenzstruktur
 aus den Punkten und Geraden von $T$ eine projektive Geometrie ist.
 \par
       Die Bl\"ocke von $T$ sind offensichtlich Unterr\"aume der eben
 konstruierten projektiven Geometrie $L$ und da die Hyperebenen von
 $\sigma$ die einzigen Unterr\"aume von $L$ sind, die von allen
 Geraden getroffen werden, folgt, dass die Bl\"ocke von $T$
 Hyperebenen von $L$ sind. Nun besagt die fisher'sche Ungleichung,
 dass $v \leq b$ ist, weil andererseits nach I.7.7 die Anzahl der
 Hyperebenen von $L$ gleich der Anzahl der Punkte von $L$ ist, ist
 $b \leq v$. Also ist $b = v$, so dass $T$ in der Tat zu $L_{1,r-1}$
 isomorph ist, wenn $r$ der Rang von $L$ ist. Schlie\ss lich folgt
 aus $v \geq k + 2$, dass $r \geq 3$ ist. Damit ist der Satz bewiesen.
 \medskip
       Wenn man in 10.1 nur $v \geq k + 1$ voraussetzt, so erh\"alt man als
 weitere Beispiele nur noch die Blockpl\"ane $(M, P_{|M|-1}(M), \in)$
 mit einer endlichen Menge $M$, wie der Beweis von 10.1 klar zeigt. Dabei ist
 $P_{|M|-1}(M)$ die Menge der $(|M| - 1)$-Teilmengen von $M$.
 \medskip\noindent
 {\bf 10.2. Satz.} {\it Ist $T$ ein projektiver Blockplan mit $v \geq k + 2$, so
 sind die folgenden Bedingungen \"aquivalent:
 \item{a)} Es gibt einen Vektorraum $V$ des Ranges $r \geq 3$, so
 dass $T$ und $L(V)_{1,r-1}$ isomorph sind.
 \item{b)} $T$ besitzt eine auf der Menge der geordneten Tripel
 nicht kol\-li\-ne\-a\-rer Punkte transitive Kollineationsgruppe.
 \item{c)} $T$ besitzt eine auf der Menge der nicht inzidenten
 Punkt-Ge\-ra\-den\-paa\-re transitive Kollineationsgruppe.\par}
 \smallskip
       Beweis. Es gelte a). Da die Geraden von $L(V)$ gleich dem Schnitt
 der sie umfassenden Hyperebenen sind, bedeutet Kollinearit\"at in
 $T$ dasselbe wie Kollinearit\"at in $L(V)_{1,r-1}$. Daher folgt
 aus der Transitivit\"at der Kollineationsgruppe von $L(V)$ auf der
 Menge der Rahmen von $L(V)$ die Transitivit\"at dieser Gruppe auf
 der Menge der geordneten Tripel nicht kollinearer Punkte. Aus a) folgt also b).
 \par
       c) folgt trivialerweise aus b).
 \par
       Um zu zeigen, dass a) von c) impliziert wird, zeigen wir
 zun\"achst, dass jede Gerade von $T$ jeden Block von $T$ trifft.
 Es sei $G$ eine Gerade von $T$ und $\Gamma$ sei der Stabilisator
 von $G$ in der Kollineationsgruppe von $T$. Dann ist $\Gamma$ auf
 der Menge der nicht mit $G$ inzidierenden Punkte transitiv.
 Zerlegt $\Gamma$ die Menge der Punkte auf $G$ in $s$ Bahnen, so
 hat $\Gamma$ insgesamt $s + 1$ Punktbahnen und daher nach 3.4 auch
 genau $s + 1$ Blockbahnen. Sind $c$ und $d$ Bl\"ocke, die $G$
 treffen, jedoch nicht ent\-hal\-ten, so gibt es genau zwei Punkte $P$
 und $Q$ auf $G$ mit $P \I c$ und $Q \I d$. Liegen $P$ und $Q$ in
 verschiedenen Punktbahnen, so liegen $c$ und $d$ in verschiedenen
 Blockbahnen von $\Gamma$. Somit zerf\"allt die Menge der Bl\"ocke,
 die $G$ treffen, aber nicht enthalten, in mindestens $s$ Bahnen
 unter $\Gamma$. Da die Menge der Bl\"ocke, die $G$ enthalten,
 ebenfalls unter $\Gamma$ invariant ist, zerf\"allt die Menge der
 Bl\"ocke, die $G$ treffen, in mindestens $s + 1$ Bahnen unter
 $\Gamma$. Da $\Gamma$ aber genau $s + 1$ Blockbahnen hat, trifft die
 Gerade $G$ jeden Block von $T$. Nach 10.1 gibt es daher einen
 endlichen projektiven Verband $L$ mit $\Rg(L) = r \geq 3$, so dass
 $T$ und $L_{1,r-1}$ isomorph sind. Ist $r = 3$, so ist $T$ eine
 projektive Ebene, die nach 9.13 desarguessch ist. Also ist $L$ in
 jedem Fall eine desarguessche projektive Geometrie. Nach II.6.1
 gibt es somit einen Vektorraum, so dass $L$ und $L(V)$ isomorph
 sind. Damit ist alles bewiesen.

\mysectionten{11. Ein Satz von N. Ito}

\noindent
 Eine
 Kollineation eines Blockplanes $T$ hei\ss t {\it Zentralkollineation\/}
 von\index{Zentralkollineation}{} $T$ mit {\it Zentrum\/}\index{Zentrum}{} $P$,
 wenn sie jeden Block durch $P$ invariant l\"asst.
 \medskip\noindent
 {\bf 11.1. Satz.} {\it Es sei $T$ ein Blockplan mit der
 Eigenschaft, dass zwei verschiedene Bl\"ocke von $T$ h\"ochstens
 $k - 2$ Punkte gemeinsam haben. Ist dann $\sigma$ eine
 Zentralkollineation von $T$ und hat $\sigma$ zwei verschiedene
 Zentren, so ist $\sigma = 1$.}
 \smallskip
       Beweis. $\sigma$ habe die beiden Zentren $P$ und $Q$. Ist $G$ eine
 Gerade durch $P$ oder $Q$, so ist $G^\sigma = G$, da $\sigma$ ja
 alle Bl\"ocke durch $P$ bzw. $Q$ festl\"asst. Ist $R$ ein Punkt
 mit $R \not\leq P + Q$, so ist $R = (P + R) \cap (Q + R)$, woraus $R^\sigma = R$
 folgt. Ist $R \leq P + Q$, so gibt es einen Block $c$ mit
 $R = (P + Q) \cap c$. Nun enth\"alt $c$ mindestens $k - 1$ Fixpunkte, so
 dass $c$ und $c^\sigma$ mindestens $k - 1$ Punkte gemeinsam haben.
 Auf Grund unserer Annahme ist daher $c^\sigma = c$. Also ist
 $$ R^\sigma = \bigl((P + Q) \cap c\bigr)^\sigma
	     = (P + Q)^\sigma \cap c^\sigma = (P + Q)\cap c = R. $$
 Damit ist gezeigt, dass $\sigma$ alle Punkte von $T$ festl\"asst, woraus
 $\sigma = 1$ folgt. 
 \medskip\noindent
 {\bf 11.2. Satz.} {\it Es sei $G$ eine auf $\Omega$ zweifach transitiv
 operierende, end\-li\-che Permutationsgruppe. Ist $U$ eine Untergruppe von $G$
 mit $|G : U| < |\Omega|$, so ist $U$ transitiv auf $\Omega$.}
 \smallskip
       Beweis. $U$ sei intransitiv. Weil $U \neq \{1\}$ ist, gibt es dann
 eine Bahn $c$ von $U$ mit $2 \leq |c| < |\Omega|$. Es sei $B$ die
 Menge aller Bilder von $c$ unter $G$. Wegen der zweifachen
 Transitivit\"at von $G$ auf $\Omega$ ist dann $T := (\Omega,B,\in )$ ein
 2-Blockplan mit $k = |c| < |\Omega| = v$. Nach 2.4 ist
 daher $v \leq |B|$. Nun ist $|B| = |G : H|$, wobei $H$ der
 Stabilisator von $c$ in $G$ ist. Weil $c$ eine Bahn von $U$ ist,
 ist $U \subseteq H$. Somit ist $|G : H| \leq |G : U|$ und daher
 $$ v \leq |G : H| \leq |G : U| < v, $$
 q. e. a.
 \medskip
       Nun haben wir alles beisammen, um den folgenden Satz von
 Ito\index{Satz von Ito}{} zu beweisen.
 \medskip\noindent
 {\bf 11.3. Satz.} {\it Es sei $\Pi$ eine Menge von $v$ Punkten und
 $G$ sei eine auf $\Pi$ zweifach transitive Permutationsgruppe.
 Ferner sei $c$ eine Teilmenge von $\Pi$ mit $2 \leq |c| < v$ und es sei
 $B := \{c^\gamma \mid \gamma \in G\}$. Ist dann $c$ Bahn einer
 Untergruppe $\Delta$ vom Index $v$ in $G$ und operiert $\Delta$
 auf $c$ nicht treu, so besteht $T := (\Pi,B,\in)$ aus den Punkten
 und Hyperebenen  einer desarguesschen projektiven Geometrie $L(V)$
 mit $\Rg(L(V)) \geq 3$ und es ist
 $\PSL(V) \subseteq G \subseteq \PGammaL(V)$.}
 \smallskip
       Beweis. Weil $G$ auf $\Pi$ zweifach transitiv operiert, ist $T$
 ein Blockplan. Ist $\bar{\Delta}$ der Stabilisator von $c$ in
 $G$, so ist $\Delta \subseteq \bar{\Delta}$ und somit
 $$ |B| = |G : \bar\Delta| \leq |G : \Delta| = v. $$
 Hieraus folgt mittels der fisherschen Ungleichung, dass $|B| = v$ ist. Folglich
 ist $\bar\Delta = \Delta$ und der Blockplan $T$ ist projektiv.
 \par
       Es sei $N$ die Untergruppe von $\Delta$, die $c$ punktweise
 festl\"asst. Da $\Delta$ auf $c$ nicht treu operiert, ist $N \neq \{1\}$.
 Hieraus folgt, dass $|\Pi - c| \geq 2$ ist. Setzt man
 $k := |c|$, so ist also $v \geq k + 2$. Es sei $\lambda$ die Anzahl
 der mit zwei verschiedenen Punkten inzidierenden Bl\"ocke. Sind
 dann $d$ und $c$ zwei verschiedene Bl\"ocke von $T$, so ist nach
 2.6 auch $|d \cap e| = \lambda$. Nun ist $k(k - 1) = \lambda(v - 1)$, so
 dass wegen $k < v - 1$ auch $\lambda < k - 1$ gilt. Wir k\"onnen
 daher im Folgenden 11.1 anwenden.
 \par
       $G$ ist eine Gruppe von Kollineationen von $T$, die auf $\Pi$
 zweifach transitiv operiert. Nach 3.5 operiert $G$ folglich auch
 zweifach transitiv auf $B$.
 \par
       Wir zeigen nun, dass $N$ auf $\Pi - c$ transitiv ist. Dazu nehmen
 wir an, dass dies nicht der Fall sei. Wir zeigen dann
 \smallskip
       a) $N$ operiert regul\"ar auf $\Pi -c$.
 \smallskip
       Es ist $\Delta = G_c$, so dass $\Delta$ wegen der zweifachen
 Transitivit\"at von $G$ auf $B$ zwei Blockbahnen hat. Nach 3.4 hat
 $\Delta$ auch genau zwei Punktbahnen. Folglich ist $\Delta$ auf
 $\Pi - c$ transitiv. Weil $N$ in $\Delta$ normal ist, haben alle
 in $\Pi - c$ enthaltenen Bahnen von $N$ die gleiche L\"ange. Die
 all diesen Bahnen gemeinsame L\"ange sei $t$. Es sei $1 \neq \nu \in N$. Ferner
 sei $P \in \Pi -c$ und es gelte $P^\nu = P$. Ist
 $d$ ein Block durch $P$, so ist $|c \cap d| = \lambda$. Daher
 enth\"alt $d$ mindestens $\lambda + 1$ Fixpunkte von $\nu$. Weil
 zwei verschiedene Bl\"ocke nur $\lambda$ Punkte gemeinsam haben,
 folgt $d^\nu = d$. Somit ist $P$ Zentrum von $\nu$. Weil $\nu \neq 1$ ist, folgt
 nach 11.1 weiter, dass $P$ der einzige Fixpunkt von
 $\nu$ in $\Pi - c$ ist. Es sei $\Phi$ die Bahn von $N$, die $P$
 enth\"alt. Dann ist $|\Phi| = t$ und $\Phi^\nu = \Phi$. Weil $\nu^i$
 f\"ur alle nat\"urliche Zahlen $i$ zentral mit dem Zentrum $P$
 ist, folgt mit 11.1, dass die Ordnung von $\nu$ Teiler von $t - 1$
 ist. Es gibt eine von $\Phi$ verschiedene, in $\Pi - c$ enthaltene
 Bahn $\Psi$ von $N$. Wegen $|\Psi| = t$ und $\ggT (t,t - 1) = 1$ folgt,
 dass $\nu$ einen Fixpunkt in $\Psi$ hat. Wegen $\Phi \cap \Psi = \emptyset$ ist
 $P \neq Q$, so dass doch $\nu = 1$ ist. Dieser
 Widerspruch zeigt die G\"ultigkeit von a).
 \par
       Setze $t := |N|$. Nach a) ist dann $(v - k)t^{-1}$ die Anzahl der
 Bahnen von $N$ in $\Pi - c$. Daher ist $k + (b - k)t^{-1}$ die Anzahl
 der Punktbahnen von $N$.
 \par
       Die Gruppe $\Delta$ ist auf der Menge der von $c$ verschiedenen
 Bl\"ocke transitiv. Weil $N$ in $\Delta$ normal ist, zerlegt $N$
 daher diese Menge in lauter Bahnen gleicher L\"ange. Ist $d$ ein
 von $c$ verschiedener Block und ist $u := |N : N_d|$, so ist $u$ die
 L\"ange einer solchen Bahn. Daher ist $(v - 1)u^{-1}$ die Anzahl der
 von $\{c\}$ verschiedenen Blockbahnen von $N$. Insgesamt ist also
 $1 + (v - 1)u^{-1}$ die Anzahl der Blockbahnen von $N$. Wegen 3.4 gilt somit
 \smallskip
       b) Es ist $k + (v - k)t^{-1} = 1 + (v - 1)u^{-1}$.
 \smallskip
       Unser n\"achstes Ziel ist zu zeigen, dass $N$ eine $p$-Gruppe ist.
 Dazu sei $p$ ein Primteiler von $|N|$ und $p^a$ sei die Ordnung
 einer $p$-Sy\-low\-grup\-pe $\Sigma$ von $N$. Schlie\ss lich sei
 $N_\Delta (\Sigma)$ der Normalisator von $\Sigma$ in $\Delta$.
 Weil $N$ ein Normalteiler von $\Delta$ ist, liegen alle
 Konjugierten von $\Sigma$ in $N$. Daher gilt
 $$ \big|\Delta : N_\Delta(\Sigma)\big|
		  = \big|N : N_\Delta(\Sigma) \cap N\big|. $$
 Hieraus folgt,
 $$ |N_\Delta(\Sigma) N|
	= {|N_\Delta(\Sigma)||N|\over |N_\Delta(\Sigma) \cap N|}
	= {|N_\Delta(\Sigma)||\Delta| \over |N_\Delta(\Sigma)|} = |\Delta| $$
 und damit $\Delta = N_\Delta(\Sigma)N$. (Dies ist das sogenannte
 Frattiniargument.) Hieraus folgt, dass
 $N_\Delta(\Sigma)$ auf der Menge der von $\{c\}$ verschiedenen Blockbahnen von
 $N$ transitiv operiert, da $\Delta$, wie wir schon bemerkten, auf der Menge der
 von $c$ verschiedenen Bl\"ocke transitiv ist. Es sei $x$ die
 Anzahl der in einer von $\{c\}$ verschiedenen Blockbahn von $N$
 enthaltenen Bahnen von $\Sigma$, so folgt daher, dass
 $(v - 1)xu^{-1}$ die Anzahl der von $\{c\}$ verschiedenen
 Blockbahnen von $\Sigma$ ist. Also hat $\Pi$ genau
 $1 + (v - 1)xu^{-1}$ Blockbahnen. Weil $\Sigma$ auf $\Pi - c$ regul\"ar
 operiert, ist $k + (v - k)p^{-a}$ die Anzahl der Punktbahnen von
 $\Pi$. Nach 3.4 gilt daher
 \smallskip
       c) $k + (v - k)p^{-a} = 1 + (v - 1)xu^{-1}.$
 \smallskip
       Es sei $t = p^at'$, $u = p^bu'$ mit $\ggT(p,u') = 1$, sowie
 $|N_d| = |N|u^{-1} = y = p^cy'$ mit $\ggT(p,y') = 1$. Dann ist $t = |N| = uy$
 und daher
 \smallskip
       d) $a=b+c$ und $t' = u'y'$.
 \smallskip
 Aus a), b) und c) folgt nun
 $$ k + t'(v - k)t^{-1} = 1 + \bigl(k + (v - k)t^{-1} - 1\bigr)x. $$
 Daher ist
 \smallskip
       e) $(v - k)(t' - x)t^{-1} = (k - 1)(x - 1)$.
 \smallskip
       Ist $P$ ein Punkt von $c$, so gibt es $k - 1$ von $c$ verschiedene
 Bl\"ocke durch $P$. Da die Menge dieser Bl\"ocke unter $N$ invariant bleibt, ist
 \smallskip
       f) $k - 1 =  du$,
 \smallskip\noindent
 wobei $d$ eine nat\"urliche Zahl ist. Aus b) und f) folgt
 $$ (v - k)t^{-1} - (v - k)u^{-1} = -(k - 1) + (k - 1) + (k - 1)u^{-1}
	   = -du + d. $$
 Mit Hilfe von d) erhalten wir daher
 \smallskip
       g) $(v - k)(y - 1)t^{-1} = d(u - 1)$.
 \smallskip\noindent
 Multipliziert man die Gleichung g) mit $t' - x$, so erh\"alt man mit Hilfe von
 e) und f)
 $$\eqalign{
    d(u - 1)(t' - x) &= (v - k)(t' - x)t^{-1}(y - 1)  \cr
                    &= (k - 1)(x - 1)(y - 1)          \cr
		    &= du(x - 1)(y - 1). \cr} $$
 Also ist
 \smallskip
       h) $(u - 1)(t' - x) = u(x - 1)(y - 1),$
 \smallskip\noindent
 so dass $t'-x \equiv 0$ mod $u'$ ist. Wegen d) gilt daher
 \smallskip
       i) $x \equiv 0 \mod u'$.
 \smallskip
       Es sei $\Gamma$ eine von $\{c\}$ verschiedene Blockbahn von $N$
 und $\Gamma_1, \dots, \Gamma_x$ seien die in $\Gamma$ enthaltenen
 Blockbahnen von $\Sigma$. Ferner sei $c_i$ ein Block aus
 $\Gamma_i$ und $N_i$ bzw.\ $\Sigma_i$ sei der Stabilisator von
 $c_i$ in $N$ bzw.\ $\Sigma$. Wie wir wissen, ist $|N : N_i| = u = p^bu'$.
 Weil $\Sigma_i$ in einer $p$-Sylowgruppe von $N_i$ liegt, folgt,
 dass $p_b$ ein Teiler von $|\Sigma : \Sigma_i|$ ist, dh., es ist
 $|\Sigma : \Sigma_i| = p^b u_i$ mit einer nat\"urlichen Zahl $u_i$.
 Schlie\ss lich gilt wegen
 $$ |\Gamma_i| = |\Sigma : \Gamma_i| = p^bu_i $$
 die Gleichung
 $$ \sum_{i:=1}^x p^b u_i = |\Gamma| = u. $$
 Also ist $\sum_{i:=1}^x u_i = u'$, woraus wegen $u_i \geq 1$ folgt, dass
 $x \leq u'$ ist.  Zusammen mit i) ergibt das $u' = x$. Hieraus, aus h) und d)
 folgt nun die Gleichung
 $$ (u - 1)(u'y' - u') = p^bu'(u' - 1)(y - 1).$$
 Daher gilt
 \smallskip
       j) $(u - 1)(y' - 1) = p^b(u' - 1)(y - 1).$
 \smallskip\noindent
 W\"are $u = 1$, so w\"are nach b) auch $t = 1$, dh., es
 w\"are $N = \{1\}$. Dieser Widerspruch zeigt, dass $u > 1$ ist.
 W\"are $y = 1$, so w\"are $t = u$ und aus b) folgte wiederum $t = 1$.
 Also ist auch $y > 1$. Daher folgt aus j), dass genau dann $y' = 1$
 ist, wenn $u' = 1$ ist. Also folgt $t' = 1$ sowohl aus $y' = 1$ als auch
 aus $u' = 1$, da ja $t' = y' u'$ ist. Wir d\"urfen daher annehmen,
 dass $u'$, $y' \neq 1$ ist. Dann ist aber
 \smallskip
       k) $(p^b u' - 1)(a' - 1)^{-1} = p^b(p^cy' - 1)(y' - 1)^{-1}.$
 \smallskip\noindent
       W\"are $b = 0$, so folgt aus k) auch $c = 0$, woraus mit d)
 der Widerspruch $0 = b + c = a \geq 1$ folgte. Also ist $b \neq 0$.
 Wegen $u' > 1$ ist
 $$ p^b - 1 \geq (p^b - 1)(u' - 1)^{-1}. $$
 Daher ist
 $$ 2p^b - 1 \geq b^b + (p^b - 1)(u' - 1)^{-1} = (p^b u' - 1)(u' - 1)^{-1}. $$
 Wegen k) ist also
 $$ 2p^b - 1 \geq p^b (p^c y' - 1)(y' - 1)^{-1}. $$
 Folglich ist
 $$ 2 > (p^cy' - 1)(y' - 1)^{-1} = p^c + (p^c - 1)(y' - 1)^{-1} \geq p^c, $$
 so dass $c = 0$ ist. Aus j) folgt weiter, dass $p^bu' - 1= p^b(u' - 1)$ ist.
 Hieraus folgt schlie\ss lich der Widerspruch $-1 \equiv 0 \mod p$, so dass also
 doch $u' = 1 = y'$ ist. Dies besagt wiederum, dass
 $t = p^a$ ist, womit gezeigt ist, dass $N$ eine $p$-Gruppe ist.
 \par
       Die Formel b) erh\"alt nun die Gestalt
 \smallskip
       l) $k + (v - k)p^{-a} = 1 + (v - 1)p^{-b}.$
 \smallskip\noindent
       Weil $T$ ein projektiver Blockplan ist, ist $k(k - 1) = \lambda(v - 1)$.
 Ferner ist nach l)
 $$ \lambda(v - 1) = \lambda p^b\bigl(k + (v - k)p^{-a} - 1\bigr), $$
 so dass also $k(k - 1) = \lambda p^b((v - k)p^{-a} + k - 1)$ ist. Somit ist
 $\lambda(v - k)p^{-a}$ und wegen $c \leq a$ dann auch $\lambda(v - k)p^{-c}$
 durch $k - 1$ teilbar. Wir definieren $e$ durch
 \smallskip
       m) $\lambda(v - k)p^{-c} = e(k - 1).$
 \smallskip\noindent
       Multipliziert man 1) mit $\lambda p^b$ und beachtet man,
 dass $\lambda(v - 1) = k(k - 1)$ ist, so erh\"alt man
 $$ k(k - 1) = \lambda p^b(k - 1) + \lambda(v - k)p^{-c}. $$
 Hieraus und aus m) erh\"alt man nun
 \smallskip
       n) $k = \lambda p^b + e$.
 \smallskip\noindent
 Aus f) folgt schlie\ss lich, dass
 $$ e - 1 = k - 1 - \lambda p^b = du - \lambda p^b = p^b(d - \lambda) $$
 ist. Folglich ist $f := d - \lambda$  eine nicht negative ganze Zahl. \"Uberdies
 ist
 \smallskip
       o) $e=p^b f+1$.
 \smallskip\noindent
 Aus $k(k - 1) = \lambda(v - 1)$ folgt $\lambda(v - k) = (k - \lambda)(k - 1)$.
 Mit m) erh\"alt man daher $k = ep^c + \lambda$.
 Hieraus und aus o) folgt $k = fp^a + p^c + \lambda$. W\"are $f = 0$, so
 w\"are $d = \lambda$, dh., es w\"are $e = 1$. Aus m) und n) folgte
 daher $\lambda(v - k) = p^c \lambda p^b$ und somit $v - k = p^a = |N|$, so
 dass $N$ auf $\Pi - c$ transitiv w\"are. Dieser Widerspruch zeigt,
 dass $f \neq 0$ ist. Dann ist aber
 \smallskip
       p) $k > p^a$.
 \smallskip\noindent
       Weil $y = p^b > 1$ ist, ist $c$ der einzige Fixblock von
 $N$. Sind nun $P$ und $Q$ zwei verschiedene Punkte auf $c$ und
 sind $c = c_1$, $c_2$, \dots, $c_\lambda$ alle Bl\"ocke durch $P$ und
 $Q$, so folgt, da die Stabilisatoren von $c_1$, \dots, $c_\lambda$
 alle in $G$ konjugiert sind, dass es $p$-Gruppen $N_i$ gibt, die
 $c_i$ und nur $c_i$ zum Fixblock haben. Wegen $N_i \subseteq G_{P,Q}$ f\"ur alle
 $i$ ist $G_{P,Q}$ nach 7.10 auf $\{c_1, \dots, c_\lambda\}$ transitiv. Hieraus
 folgt die Gleichung
 $$ |G_{P,Q} : G_{P,Q} \cap \Delta| = \lambda, 
 $$ da ja $\Delta = G_c$ ist. Somit ist
  $$
    |G_P : G_{P,Q} \cap \Delta|
	 = |G_P : G_{P,Q}||G_{P,Q} : G_{P,Q} \cap \Delta| 
	 = \lambda(v - 1)                                
	 = k(k - 1). 
 $$
 Weil $G_P$ genau zwei Punktbahnen hat, hat $G_P$ nach
 3.4 auch genau zwei Blockbahnen. Folglich ist $G_P$ auf der Menge
 der Bl\"ocke durch $P$ transitiv. Also ist $|G_P : G_P \cap \Delta| = k$.
 Hieraus folgt wiederum
 $$\eqalign{
      k(k - 1) &= |G_P : G_{P,Q} \cap \Delta |      \cr
         &= |G_P :G_P \cap \Delta||G_P \cap \Delta : G_{P,Q} \cap \Delta| \cr
	 &= k|G_P \cap \Delta : G_{P,Q} \cap \Delta|. \cr} $$
 Also ist $|G_P \cap \Delta : G_P,Q \cap \Delta| = k - 1$, so
 dass $G_P \cap \Delta$ auf $c - \{P\}$ transitiv ist. Weil $P$ ein
 beliebiger Punkt von $c$ ist, operiert $\Delta$ auf $c$ zweifach
 transitiv. Es sei nun $1 \neq \nu \in N$ und $Z$ sei der
 Zentralisator von $\nu$ in $G$. Weil die Fixpunkte von $\nu$
 gerade die Punkte von $c$ sind, ist $c^z = c$, dh., es ist $Z \subseteq \Delta$.
 Nun ist $\nu^\Delta \subseteq N$, da $N$ ein
 Normalteiler von $\Delta$ ist. Also ist
 $$ |\Delta : Z| \leq |N| = p^a, $$
 so dass nach p) die Ungleichung $|\Delta : Z| < k$ gilt. Nach 11.2 ist $Z$
 folglich transitiv auf $c$. Nun hat
 $\nu$ genau $k$ Fixpunkte. Daher hat $\nu$ nach 3.2 auch genau
 $k$ Fixbl\"ocke. Es seien $c_1 := c$, $c_2$, \dots, $c_k$ diese
 Fixbl\"ocke. Weil $Z$ auf $c$ transitiv operiert, ist $(c, \{c_1,
 \dots, c_k\}, \in)$ eine linksseitige taktische Konfiguration mit
 $k$ Punkten. Weil nach 2.6 die Gleichung $|c_1 \cap c_i| =
 \lambda$ gilt f\"ur $i := 2$, $3$, \dots, $k$, folgt aus 1.2 a) die
 Kongruenz $k + \lambda (k - 1) \equiv 0 \mod k$. Hieraus folgt
 wiederum $\lambda \equiv 0 \mod k$, so dass wegen $1 \leq \lambda \leq k$ sogar
 $k = \lambda$ gilt. Dann ist aber wegen $k(k - 1) = \lambda (v - 1)$ auch
 $k = v$. Dieser letzte Widerspruch zeigt schlie\ss lich,
 dass $N$ doch auf $\Pi - c$ transitiv operiert.
 \par
       Wir zeigen nun, dass $G$ auf der Menge der nicht inzidenten
 Punkt-Geraden\-paare von $T$ transitiv operiert. Weil $G$ zweifach
 transitiv ist, gen\"ugt es zu zeigen, dass der Stabilisator $G_h$
 einer Geraden $h$ auf $\Pi - h$ noch transitiv ist. Es seien $c_1$,
 \dots, $c_\lambda$ die Bl\"ocke, die $h$ enthalten. Ist dann $N_i$
 die Untergruppe von $G$, die $c_i$ punktweise festl\"asst, so ist
 $N_i$ auf $\Pi - c$ transitiv. Au\ss erdem ist $N_i \subseteq G_h$. Weil
 $$ \Pi - h = \bigcup_{i:=1}^\lambda (\Pi - c_i) $$
 ist, folgt, dass $\Pi - h$ eine Bahn von $G_h$ ist, es sei
 denn, es gibt zwei Indizes $i$ und $j$ mit
 $(\Pi - c_i) \cap (\Pi - c_j)= \emptyset$. Dann ist aber $\Pi = c_i \cup c_j$
 und daher $v = 2k - \lambda$. Aus
 $$ k(k - 1) = \lambda(v - 1) = \lambda(2k - \lambda - 1) $$
 folgt die Gleichung $(k - \lambda)^2 = k - \lambda$. Wegen $k > \lambda$ ist
 somit $k = \lambda + 1$, woraus wiederum $v = k + 1$ folgt. Wie wir jedoch schon
 fr\"uher bemerkten, ist $v \geq k + 2$, so dass $\Pi - h$ doch eine
 Bahn von $G_h$ ist. Nach 10.2 gibt es folglich einen Vektorraum
 $V$ mit $r := \Rg(V) \geq 3$, so dass $T$ und $L(V)_{1,r-1}$
 isomorph sind.
 \par
       Nun ist $N$ auf $\Pi -c$ transitiv. Hieraus folgt mit 7.1, dass
 $N$ alle Elationen mit der Achse $c$ enth\"alt, was wiederum
 impliziert, dass $G$ alle Elationen von $T = L(V)_{1, r-1}$
 enth\"alt. Somit gilt auch $\PSL(V) \subseteq G \subseteq \PGammaL(V)$. Damit
 ist der Satz von Ito bewiesen.


 \newpage
       
 \mychapter{V}{Polarit\"aten}

 \noindent
      Von besonderem Interesse unter den Korrelationen projektiver
 R\"aume sind die involutorischen Korrelationen. Sie werden gemeinhin
 Polarit\"aten genannt. Ihrem Studium und dem Studium ihrer
 Zen\-tra\-li\-sa\-to\-ren
 ist das vorliegende Kapitel gewidmet. Dieses Studium werden wir in diesem
 Kapitel jedoch noch nicht beenden. Den orthogonalen Gruppen, die in diesen
 Kontext geh\"oren, werden wir ein eigenes Kapitel widmen.\footnote{Anmerkung
 der Herausgeber: Dieses Kapitel \"uber orthogonale Gruppen fehlt.}

\mysection{1. Darstellung von Polarit\"aten}

\noindent
       Es sei $V$ ein $K$-Vektorraum mit $\Rg_K(V) \geq 3$. Eine Korrelation
 $\pi$ von $\La(V)$ hei\ss t {\it Polarit\"at\/},\index{Polarit\"at}{} falls
 $\pi^2 = 1$ ist.
 \par
       Wenn wir sagen, dass $\pi$ eine Korrelation von $\La(V)$ sei, so
 beinhaltet das implizit, dass der Rang von $V$ endlich ist, da $\La(V)$
 ja h\"ochstens in diesem Fall eine Korrelation gestattet, wie wir fr\"uher
 sahen.
 \medskip\noindent
 {\bf 1.1. Satz.} {\it Es sei $V$ ein Vektorraum \"uber dem K\"orper $K$ mit
 $\Rg(V) \geq 3$ und $\pi$ sei eine Korrelation von $\La(V)$. Ferner werde
 $\pi$ durch die $\alpha$-Semibilinearform $f$ dargestellt. Genau dann ist
 $\pi$ eine Polarit\"at, wenn aus $u$, $v \in V$ und $f(u,v) = 0$ stets
 $f(v,u) = 0$ folgt.}
 \smallskip
       Beweis. Wir definieren
 $g$ durch $g(u,v) := f(v,u)^{\alpha^{-1}}$. Dann ist
 $g$ nach II.8.7 eine $\alpha^{-1}$-Semibilinearform und diese Form stellt
 $\pi^{-1}$ dar. Nach II.8.8b) gilt daher genau dann $\pi^{-1} = \pi$,
 dh.\ $\pi^2 = 1$, wenn es ein $k \in K^*$ gibt mit $kg(u,v) = f(u,v)$
 f\"ur alle $u$, $v \in V$. Ist $\pi^2 = 1$, so folgt also $f(v,u) = 0$
 aus $f(u,v) = 0$.
 \par
       Folgt andererseits $f(v,u) = 0$ aus $f(u,v) = 0$, so ist offensichtlich
 $U \leq U^{\pi^2}$ f\"ur alle $U \in \La(V)$. Weil $V$ endlichen Rang hat,
 folgt $U = U^{\pi^2}$, so dass $\pi^2 = 1$ ist. Damit ist alles bewiesen.
 \medskip\noindent
 {\bf 1.2. Satz.} {\it Es sei $V$ ein Vektorraum \"uber $K$ mit $\Rg_K(V) \geq 3$.
 Wird die Polarit\"at $\pi$ von $\La(V)$ durch die $\alpha$-Semibilinearform
 $f$ dargestellt und gibt es einen Vektor $w \in V$ mit $f(w,w) = 1$, so
 ist $\alpha^2 = 1$ und $f(u,v) = f(v,u)^\alpha$ f\"ur alle $u$, $v \in V$.}
 \smallskip
       Beweis. Wie wir beim Beweise von 1.1 gesehen haben, gibt es ein
 $k \in K^*$ mit $f(u,v) = kf(v,u)^\alpha$ f\"ur alle $u$, $v \in V$. Es folgt
 $$ 1 = f(w,w) = kf(w,w)^\alpha = k, $$
 so dass also tats\"achlich $f(u,v) = f(v,u)^\alpha$ ist f\"ur alle
 $u$, $v \in V$.
 \par
       Weil $f(u,v) = f(v,u)^\alpha$ ist f\"ur alle $u$ und $v$, ist auch
 $f(v,u) = f(u,v)^\alpha$ f\"ur alle $u$, $v \in V$. Somit ist $f(u,v) =
 f(u,v)^{\alpha^2}$ f\"ur alle $u$, $v \in V$. Hieraus folgt
 $$ x = f(w,w)x = f(w,wx)
      = f(w,wx)^{\alpha^2} = f(w,w)^{\alpha^2}x^{\alpha^2}
      = x^{\alpha^2}, $$
 so dass $\alpha^2 = 1$ ist. Damit ist alles bewiesen.
 \medskip
       Wir nennen die $\alpha$-Semibilinearform {\it f
 symmetrisch\/},\index{symmetrische Semibilinearform}{} falls
 $f(u,v) = f(v,u)^\alpha$ f\"ur alle $u$, $v \in V$ gilt. Wie wir gerade
 gesehen haben, ist dann $f(u,v) = f(u,v)^{\alpha^2}$. Gibt es Vektoren
 $u$, $v$ mit $f(u,v) \neq 0$, so gibt es auch Vektoren $u$, $v$ mit
 $f(u,v) = 1$. Dann ist aber
 $$    x = f(u,v)x = f(u,vx)
         = f(u,vx)^{\alpha^2} = f(u,v)^{\alpha^2}x^{\alpha^2}
         = x^{\alpha^2} $$
 f\"ur alle $x \in K$, so dass $\alpha^2 = 1$ ist. Ist also $f$ nicht die
 Nullform und ist $f$ symmetrisch, so ist $\alpha^2 = 1$.
 \par
       In Abschnitt 4 von Kapitel II hatten wir die Begriffe {\it absoluter
 Punkt\/} und {\it absolute Hyperebene\/} in Bezug auf eine Korrelation
 de\-fi\-niert. Wir wiederholen hier diese Definition f\"ur den Fall, dass
 $\pi$ eine Polarit\"at des projektiven Verbandes $L$ ist. Ist $P$ ein Punkt
 von $L$, so hei\ss t $P$ {\it absoluter Punkt\/}\index{absoluter Punkt}{} von
 $\pi$, falls $P \leq P^\pi$ ist. Analog hei\ss t die Hyperebene $H$ von $L$
 {\it absolut\/},\index{absolute Hyperebene}{}
 falls $H^\pi \leq H$ ist.
 \medskip\noindent
 {\bf 1.3. Satz.} {\it Es sei $V$ ein $K$-Vektorraum mit
 $\Rg_K(V) \geq 3$ und $\pi$ sei eine Polarit\"at von $V$. Sind nicht alle
 Punkte von $\La(V)$ absolut, so l\"asst sich $\pi$ durch eine symmetrische
 $\alpha$-Semibilinearform darstellen.}
 \smallskip
       Beweis. Nach dem dritten Struktursatz (II.8.3) gibt es eine
 $\beta$-Semibilinear\-form $g$, welche $\pi$ darstellt. Da nicht alle Punkte
 absolut sind, gibt es ein $w \in V$ mit $g(w,w) \neq 0$. Setze $k := g(w,w)$
 und definiere $f$ durch $f(u,v) := k^{-1}g(u,v)$. Dann stellt auch $f$
 nach II.8.7a) die Polarit\"at $\pi$ dar. Wegen $f(w,w) = 1$ ist $f$ nach
 1.2 und der dem Beweise von 1.2 nachfolgenden Definition symmetrisch.
 \medskip
       Der n\"achste Satz gibt Auskunft, was geschieht, wenn alle Punkte
 der Polarit\"at $\pi$ absolut sind.
 \medskip\noindent
 {\bf 1.4. Satz.} {\it Es sei $V$ ein $K$-Vektorraum mit
 $\Rg_K(V) \geq 3$. Ferner sei $\pi$ eine Korrelation von $\La(V)$ und $f$
 eine $\pi$ darstellende $\alpha$-Semi\-bi\-li\-ne\-ar\-form.
 Sind alle Punkte von
 $\La(V)$ absolut bez\"uglich $\pi$, so gilt:
\item{a)} Es ist $f(v,v) = 0$ f\"ur alle $v \in V$.
\item{b)} Es ist $f(u,v) = -f(v,u)$ f\"ur alle $u$, $v \in V$.
\item{c)} $\pi$ ist eine Polarit\"at.
\item{d)} Es ist $\alpha = 1$. Insbesondere ist $K$ kommutativ.}
\smallskip
       Beweis. a) Ist $v = 0$, so ist $f(v,v) = 0$. Ist $v \neq 0$, so ist
 $vK$ ein Punkt und daher $vK \leq (vK)^\pi$. Es folgt $f(v,v) = 0$.
 \par
       b) Nach a) ist
 $$
       0 = f(u + v,u + v) 
         = f(u,u) + f(u,v) + f(v,u) + f(v,v) 
         = f(u,v) + f(v,u), 
 $$
 woraus $f(u,v) = -f(v,u)$ folgt.
 \par
       c) Wegen $f(u,v) = -f(v,u)$ ist $f(v,u) = 0$ eine Folge von $f(u,v)
 = 0$. Daher ist $\pi$ nach 1.1 eine Polarit\"at.
 \par
       d) Weil $f$ nicht ausgeartet ist, gibt es Vektoren $u$ und $v$ mit
 $f(u,v) \neq 0$. Setze $k := f(u,v)$. Ist $x \in K$, so folgt
 $$
       kx = f(u,v)x = f(u,vx) = -f(vx,u) 
          = -x^\alpha f(v,u) = x^\alpha f(u,v) 
          = x^\alpha k. 
 $$
 Somit ist $x^\alpha = kxk^{-1}$, so dass $\alpha$ ein Automorphismus von
 $K$ ist. Weil $\alpha$ andererseits auch ein Antiautomorphismus von $K$ ist,
 ist $K$ kommutativ. Hieraus folgt schlie\ss lich $x^\alpha = x$, so dass
 $\alpha$ die Identit\"at ist.
 \medskip
       Ist $\kappa$ eine Korrelation des projektiven Verbandes $L$, so ist
 $\kappa$ nichts anderes als ein Isomorphismus von $L$ auf $L^d$. Daher ist
 $\Rg(X) = \Ko(X^\kappa)$, dh. es ist
 $$ \Rg(L) = \Ko(X^\kappa) + \Rg(X^\kappa) = \Rg(X) + \Rg(X^\kappa) $$
 f\"ur alle $X \in L$. Aus dieser Bemerkung folgt unmittelbar der n\"achste
 Satz.
 \medskip\noindent
 {\bf 1.5. Satz.} {\it Es sei $V$ ein Vektorraum und $\kappa$ sei eine
 Korrelation von $\La(V)$. Ist $X \in \La(V)$, so sind die folgenden Bedingungen
 \"aquivalent:
\item{a)} Es ist $X \cap X^\kappa = 0$.
\item{b)} Es ist $X + X^\kappa = V$.
\item{c)} Es ist $V = X \oplus X^\kappa$.}
 \medskip
       Eine $\alpha$-Semibilinearform mit $\alpha = 1$ nennen wir
 {\it Bilinearform\/}.\index{Bilinearform}{}
 \medskip\noindent
 {\bf 1.6. Satz.} {\it Es sei $V$ ein endlich erzeugter Vektorraum \"uber $K$ mit
 $\Rg_K(V) \geq 2$ und $f$ sei eine nicht ausgeartete Bilinearform auf $V$.
 Ist $f(v,v) = 0$ f\"ur alle $v \in V$, so ist $\Rg_K(V) = 2n$ und es gibt
 eine Basis $b_1$, \dots, $b_{2n}$ von $V$ mit
 $$
 f\bigl(\sum_{i:=1}^{2n} b_ix_i,\sum_{j:=1}^{2n} b_jy_j \bigr) =
       \sum_{i:=1}^{n} (x_{2i-1}y_{2i} - x_{2i}y_{2i-1}) $$
 f\"ur alle $n$-Tupel $x$ und $y$ \"uber $K$.}
 \smallskip
       Beweis. Weil $f$ nicht ausgeartet ist, gibt es Vektoren $u$, $v \in V$
 mit $f(u,v) \neq 0$. Es gibt daher auch Vektoren $b_1$, $b_2 \in V$ mit
 $f(b_1,b_2) = 1$. W\"aren $b_1$ und $b_2$ linear abh\"angig, so w\"are $b_2
 = b_1k$ mit einem $k \in K$ und daher
 $$ 1 = f(b_1,b_2) = f(b_1,b_1)k = 0. $$
 Dieser Widerspruch zeigt, dass $b_1$ und $b_2$ linear unabh\"angig sind.
 Gilt nun $\Rg_K(V) = 2$, so ist $\{b_1,b_2\}$ eine Basis von $V$ und es ist,
 da ja offenbar $f(u,v) = -f(v,u)$ ist,
 $$ f(b_1x_1 + b_2x_2,b_1y_1 + b_2y_2) = x_1y_2 - x_2y_1. $$
 Wir k\"onnen daher annehmen, dass $\Rg_K(V) > 2$ ist. Weil $f$ nicht
 ausgeartet ist, induziert $f$ nach II.8.4 eine Korrelation $\pi$ in $\La(V)$,
 die wegen $f(u,v) = -f(v,u)$ nach 1.1 sogar eine Polarit\"at ist. Es sei
 $U := b_1K + b_2K$. Ist $v \in U \cap U^\pi$, so ist $v = b_1x_1 + b_2x_2$
 und $0 = f(b_1,v) = x_2$ und $0 = f(b_2,v) = x_1$, so dass $v = 0$ ist.
 Somit ist $U \cap U^\pi = \{0\}$. Nach 1.5 ist folglich $V = U \oplus U^\pi$.
 Weil $U^{\pi^2} = U$ ist, ist die Einschr\"ankung von $f$ auf $U^\pi$ nicht
 ausgeartet. Nach Induktionsannahme ist daher $\Rg_K(U^\pi) = 2(n - 1)$ und es
 gibt eine Basis $b_3$, $b_4$, \dots, $b_{2n-1}$, $b_{2n}$ von $U^\pi$ mit
 $$\textstyle f\bigl(\sum_{i:=3}^{2n} b_ix_i,\sum_{j:=3}^{2n} b_jy_j\bigr) =
       \sum_{i:=2}^n (x_{2i-1}y_{2i} - x_{2i}y_{2i-1}). $$
 Daher ist $\Rg_K(V) = 2n$ und
 $$\textstyle f\bigl(\sum_{i:=1}^{2n} b_ix_i,\sum_{j:=1}^{2n} b_jy_j\bigr) =
       \sum_{i:=1}^n (x_{2i-1}y_{2i} - x_{2i}y_{2i-1}). $$
 Damit ist alles bewiesen.
 \medskip
       Ist $\pi$ eine Polarit\"at von $\La(V)$ und ist jeder Punkt von
 $\La(V)$ absolut, so nennen wir $\pi$
 {\it symplektisch\/}.\index{symplektische Polarit\"at}{}
 \medskip\noindent
 {\bf 1.7. Satz.} {\it Es sei $V$ ein Vektorraum \"uber $K$ und es gelte $\Rg_K(V)
 \geq 3$. Genau dann besitzt $\La(V)$ eine symplektische Polarit\"at, wenn
 $K$ kommutativ und der Rang von $V$ endlich und gerade ist. Sind $\pi$
 und $\pi'$ symplektische Polarit\"aten von $\La(V)$, so gibt es ein $\gamma
 \in \PGL(V)$ mit $\gamma^{-1}\pi\gamma = \pi'$.}
 \smallskip
       Beweis. Besitzt $\La(V)$ eine symplektische Polarit\"at, so folgt mit
 I.5.7, dass der Rang von $V$ endlich ist. Aus 1.4 folgt die Kommutativit\"at
 von $K$ und aus 1.6, dass $\Rg_K(V)$ gerade ist.
 \par
       Es sei umgekehrt $\Rg_K(V) = 2n$ und $K$ sei kommutativ. Ist $b_1$,
 \dots, $b_{2n}$ eine Basis von $V$ und definiert man $f$ durch
 $$\textstyle f\bigl(\sum_{i:=1}^{2n} b_ix_i,\sum_{j:=1}^{2n} b_jy_j\bigr) :=
       \sum_{i:=1}^n (x_{2i-1}y_{2i} - x_{2i}y_{2i-1}), $$
 so verifiziert man leicht, dass $f$ eine nicht ausgeartete Bilinearform
 mit $f(v,v) = 0$ f\"ur alle $v \in V$ ist, so dass $\La(V)$ nach II.8.4
 und 1.4 eine symplektische Polarit\"at besitzt.
 \par
       Es seien schlie\ss lich $\pi$ und $\pi'$ symplektische Polarit\"aten
 von $\La(V)$. Ferner seien $f$ und $f'$ Bilinearformen, die $\pi$ bzw. $\pi'$
 darstellen. Nach 1.4 und 1.6 gibt es dann Basen $b_1$, \dots, $b_{2n}$
 bzw. $b_1'$, \dots, $b_{2n}'$ von $V$ mit
 $$\textstyle f\bigl(\sum_{i:=1}^{2n} b_ix_i,\sum_{j:=1}^{2n} b_jy_j\bigr) :=
       \sum_{i:=1}^n (x_{2i-1}y_{2i} - x_{2i}y_{2i-1}) $$
 bzw.  
 $$\textstyle f'\bigl(\sum_{i:=1}^{2n} b'_ix_i,\sum_{j:=1}^{2n} b'_jy_j\bigr) :=
       \sum_{i:=1}^n (x_{2i-1}y_{2i} - x_{2i}y_{2i-1}). $$
 Es gibt ein $\beta \in \GL(V)$ mit $b_i^\beta = b'_i$ f\"ur $i := 1$, \dots,
 $2n$. Dann ist aber $f(u,v) = f'(u^\beta,v^\beta)$ f\"ur alle $u,v \in V$.
 Daher ist genau dann $v \in U^\pi$, wenn $v^\beta \in U^{\beta\pi'}$ ist.
 Andererseits ist genau dann $v \in U^\pi$, wenn $v^\beta \in U^{\pi\beta}$
 ist. Folglich ist $U^{\pi\beta} = U^{\beta\pi'}$ f\"ur alle $U \in \La(V)$.
 Ist $\gamma$ die von $\beta$ in $\La(V)$ induzierte Kollineation, so ist also
 $\pi\gamma = \gamma\pi'$. Es folgt $\pi' = \gamma^{-1}\pi\gamma$, womit auch
 die letzte noch ausstehende Behauptung bewiesen ist.
 \medskip
       Der gerade bewiesene Satz gibt Anlass f\"ur eine Definition. Ist
 $\pi$ eine symplektische Polarit\"at von $\La(V)$ und ist $f$ eine $\pi$
 darstellende Bilinearform, ist ferner $b_1$, \dots, $b_{2n}$ eine Basis von
 $V$ mit
 $$\textstyle f\bigl(\sum_{i:=1}^{2n} b_ix_i,\sum_{j:=1}^{2n} b_jy_j\bigr) =
       \sum_{i:=1}^n (x_{2i-1}y_{2i} - x_{2i}y_{2i-1}), $$
 so nennen wir $b_1$, \dots, $b_{2n}$ eine {\it symplektische
 Basis\/}\index{symplektische Basis}{} von $V$.
 \medskip\noindent
 {\bf 1.8. Satz.} {\it Es sei $V$ ein Vektorraum \"uber $K$ mit
 $\Rg_K(V) \geq 3$.
 Genau dann besitzt $\La(V)$ eine Polarit\"at, wenn $\Rg_K(V)$ endlich ist
 und wenn $K$ einen Antiautomorphismus $\alpha$ mit $\alpha^2 = 1$ besitzt.}
 \smallskip
       Beweis. Es sei $\pi$ eine Polarit\"at von $V$. Nach I.5.7 ist dann
 $\Rg_K(V)$ wegen $\Rg_K(V) = \Rg(\La(V))$ endlich. Ist $\pi$ eine
 symplektische Polarit\"at, so ist $K$ kommutativ, und die Identit\"at ist
 ein Antiautomor\-phismus von $K$. Ist $\pi$ nicht symplektisch, so l\"asst
 sich $\pi$ nach 1.3 mittels einer symmetrischen $\alpha$-Semibilinearform
 darstellen. Nach der Bemerkung vor 1.3 ist dann $\alpha^2 = 1$.
 \par
       Es sei umgekehrt $\Rg_K(V) = n$ endlich und $\alpha$ sei ein
 Antiautomorphismus von $K$ mit $\alpha^2 = 1$. Ist $b_1$, \dots, $b_n$ eine
 Basis von $V$, ist $u = \sum_{i:=1}^n b_ix_i$ und $v = \sum_{i:=1}^n b_ny_n$,
 so setzen wir
 $$\textstyle f(u,v) := \sum_{i:=1}^n x_i^\alpha y_i. $$
 Man verifiziert m\"uhelos, dass $f$ eine symmetrische, nicht ausgeartete
 $\alpha$-Semi\-bi\-line\-ar\-form ist, so dass $\La(V)$ nach II.8.4 und 1.1
 eine Polarit\"at besitzt.
 \medskip
       Aus 1.8 folgt insbesondere, dass alle projektiven Geometrien
 endlichen Rang\-es, deren Koordinatenk\"orper kommutativ ist,
 Po\-la\-ri\-t\"a\-ten
 besitzen, da ja in diesem Falle die Identit\"at die Rolle des
 Antiautomorphismus $\alpha$ mit $\alpha^2 = 1$ spielen kann.
 \par
       Es sei $\pi$ eine Korrelation und $f$ eine $\pi$ darstellende
 $\alpha$-Semi\-bi\-li\-ne\-ar\-form.
 Ist $g$ eine ebenfalls $\pi$ darstellende
 $\beta$-Semibilinearform, so gibt es nach II.8.7 ein Element $k \in K^*$
 mit $x^\beta = k^{-1}x^\alpha k$ f\"ur alle $x \in K$. Ist nun $\alpha = 1$,
 so ist $K$ kommutativ, so dass auch $\beta = 1$ ist. Hieraus folgt, dass
 jede $\pi$ darstellende Semibilinearform eine Bilinearform ist, falls nur
 eine von ihnen es ist. Ist dies der Fall, so nennen wir die Korrelation
 $\pi$ {\it projektiv\/}.\index{projektive Korrelation}{} Symplektische
 Polarit\"aten sind also stets projektiv. Nicht projektive Polarit\"aten
 hei\ss en {\it unit\"ar\/}.\index{unitare@unit\"are Polarit\"at}{}
 \par
       Wie wir gesehen haben, l\"asst sich jede nicht symplektische
 Polarit\"at $\pi$ durch eine symmetrische $\alpha$-Semibilinearform
 darstellen und im Falle, dass $\pi$ projektiv ist, folgt aus II.8.7,
 dass alle $\pi$ darstellenden Bilinearformen symmetrisch sind. F\"ur
 unit\"are Polarit\"aten gibt es noch eine Darstellungen durch eine
 {\it antisymmetrische\/}
 Semibilinearform,\index{antisymmetrische Form}{}
 die gelegentlich von Nutzen ist. Hier die genaue
 Formulierung.
 \medskip\noindent
 {\bf 1.9. Satz.} {\it Es sei $V$ ein Vektorraum \"uber $K$ mit
 $\Rg_K(V) \geq 3$.
 Ferner sei $\pi$ eine unit\"are Polarit\"at von $\La(V)$ und $f$ sei eine
 $\pi$ darstellende, symmetrische $\alpha$-Semibilinearform. Es sei $l \in K$
 und es gelte $l^\alpha \neq l$. Setze $k := l - l^\alpha$ und definiere
 $\beta$ durch $x^\beta := kx^\alpha k^{-1}$. Setze schlie\ss lich
 $$ g(u,v) := kf(u,v) $$
 f\"ur alle $u$, $v \in V$. Dann ist $g$ eine $\pi$ darstellende
 $\beta$-Semibilinearform und es gilt
 $$ g(u,v) = -g(v,u)^\beta $$
 f\"ur alle $u$, $v \in V$. Au\ss erdem gilt $\beta^2 = 1$.}
 \smallskip
       Beweis. Dass $g$ eine $\beta$-Semibilinearform ist, die $\pi$
 darstellt, folgt mit II.8.7a). Dass $\beta^2 = 1$ ist, ist einfach
 nachzurechnen.
 \par
       Es ist
 $$ k^\alpha = l^\alpha - l^{\alpha^2} = l^\alpha - l = -k $$
 und
 $$ k^\beta = kk^\alpha k^{-1} = -kkk^{-1} = -k = k^\alpha, $$
 so dass auch $k^\beta = -k$ gilt.
 Sind nun $u$, $v \in V$, so ist also
 $$
 g(u,v) = kf(u,v) = kf(v,u)^\alpha k^{-1}k       
        = -f(v,u)^\beta k^\beta = -(kf(u,v))^\beta 
        = -g(v,u)^\beta.                         
 $$
 Damit ist alles bewiesen.

\mysection{2. Zentralisatoren von Polarit\"aten}

\noindent
       Von gro\ss em geometrischem und algebraischem Interesse sind die
 Zentralisatoren von Polarit\"aten. Das geometrische Interesse
 an diesen Gruppen r\"uhrt daher, dass man mit ihrer Hilfe viel \"uber die
 Struktur eines projektiven Raumes mit vorgegebener Polarit\"at herausfinden
 kann, w\"ahrend das Interesse der Algebraiker an diesen Gruppen daher
 kommt, dass viele unter den Zentralisatoren von Polarit\"aten einfach sind,
 bzw. gro\ss e einfache Untergruppen enthalten. In diesem Abschnitt beginnen
 wir das Studium dieser Gruppen.
 \medskip\noindent
 {\bf 2.1. Satz.} {\it Es sei $V$ ein Vektorraum \"uber $K$ mit
 $\Rg_K(V) \geq 3$ und $\kappa$ sei eine Korrelation von $\La(V)$, die durch
 die $\alpha$-Semibilinearform $f$ dargestellt werde. Schlie\ss lich sei
 $\sigma$ eine Kollineation von $\La(V)$ und $\beta$ sei der begleitende
 Automorphismus von $\sigma$. {\rm(}Wir interpretieren hier $\sigma$ einmal als
 Kollineation und einmal als die sie induzierende semilineare Abbildung.{\rm)} Genau
 dann ist $\kappa\sigma = \sigma\kappa$, wenn es ein $r \in K^*$ gibt, so
 dass $f(x^\sigma,y^\sigma) = rf(x,y)^\beta$ ist f\"ur alle $x$, $y \in V$.}
 \smallskip
       Beweis. Es sei $\kappa\sigma = \sigma\kappa$. Wir definieren $g$ durch
 $g(x,y) := f(x^\sigma,y^\sigma)^{\beta^{-1}}$ f\"ur alle $x$, $y \in V$.
 Dann ist $g$ eine $\beta\alpha\beta^{-1}$-Semi\-bi\-li\-ne\-ar\-form.
 Weil $f$ nicht
 ausgeartet ist, ist auch $g$ nicht ausgeartet. Daher stellt $g$ nach II.8.4
 eine Korrelation $\kappa'$ dar. Nun ist
 $$\eqalign{
 U^{\kappa'} &= \bigl\{v \mid v \in V, g(u,v) = 0\ \hbox{f\"ur alle\ }
		       u \in U\bigr\} \cr
             &= \bigl\{v \mid v \in V, f(u^\sigma,v^\sigma) = 0\
                              \hbox{f\"ur alle\ } u \in U\bigr\}. \cr} $$
 Ersetzt man in dem letzten Ausdruck $u$ und $v$ durch $u^{\sigma^{-1}}$
 bzw. $v^{\sigma^{-1}}$, so folgt weiter
 $$\eqalign{
 U^{\kappa'} &= \bigl\{v^{\sigma^{-1}} \mid v \in V, f(u,v) = 0\ \hbox{f\"ur
			alle\ } u \in U^\sigma\bigr\} \cr
             &= \bigl\{v \mid v \in V, f(u,v) = 0\ \hbox{f\"ur alle\ } u \in
                     U^\sigma\bigr\}^{\sigma^{-1}} 
             = U^{\sigma\kappa\sigma^{-1}}. \cr} $$
 Hieraus folgt schlie\ss lich, da $\sigma$ mit $\kappa$ vertauschbar ist,
 dass $U^{\kappa'} = U^\kappa$ ist. Da dies f\"ur alle $U$ gilt, ist
 $\kappa' = \kappa$. Weil $\kappa$ also auch durch $g$ dargestellt wird,
 gibt es nach II.8.7b) ein $s \in K^*$ mit $g(x,y) = sf(x,y)$ f\"ur alle
 $x$, $y \in V$. Setzt man schlie\ss lich $r := s^\beta$, so folgt
 $f(x^\sigma,y^\sigma) = rf(x,y)^\beta$.
 \par
       Es sei umgekehrt $r \in K^*$ und es gelte $f(x^\sigma,y^\sigma) =
 rf(x,y)^\beta$. Dann ist genau dann $f(x,y) = 0$, wenn $f(x^\sigma,y^\sigma)
 = 0$ ist. Hieraus folgt
 $$\eqalign{
 U^{\kappa\sigma} &= \bigl\{v^\sigma \mid v \in V, f(u,v) = 0\ 
	       \hbox{f\"ur alle\ } u \in U\bigr\} \cr
       &= \bigl\{v^\sigma \mid v \in V, f(u^\sigma,v^\sigma) = 0\ \hbox{f\"ur
			     alle\ } u \in U\bigr\} \cr
       &= \bigl\{v \mid v \in V, f(u^\sigma,v) = 0\ \hbox{f\"ur alle\ }i
			       u \in U\bigr\} \cr
       &= \bigl\{v \mid v \in V, f(u,v) = 0\ \hbox{f\"ur alle\ } u \in
				U^\sigma\bigr\} 
       = U^{\sigma\kappa}, \cr} $$
 so dass $\kappa\sigma = \sigma\kappa$ ist. Damit ist alles bewiesen.
 \medskip
       Ist $\pi$ eine Polarit\"at von $\La(V)$, so bezeichnen wir mit
 $PC^*(\pi)$ den Zentralisator von $\pi$ in $\PGL(V)$ und mit $C^*(\pi)$
 das volle Urbild von $PC^*(\pi)$ in $\GL(V)$.
 \medskip\noindent
 {\bf 2.2. Satz.} {\it Es sei $V$ ein Vektorraum \"uber $K$ mit
 $\Rg_K(V) \geq 3$ und $\pi$ sei eine Polarit\"at von $\La(V)$, die durch
 die $\alpha$-Semibilinearform $f$ dargestellt werde. Ist $\sigma \in
 C^*(\pi)$, so gibt es ein $r_\sigma \in K^*$ mit
 $$ f(x^\sigma,y^\sigma) = r_\sigma f(x,y) $$
 f\"ur alle $x$, $y \in V$. Die Abbildung $r$ ist ein Homomorphismus von
 $C^*(\pi)$ in $Z(K^*)$, der von der Wahl von $f$ unabh\"angig ist.}
 \smallskip
       Beweis. Es sei $a \in K$. Dann ist
 $$\eqalign{
 a^\alpha r_\sigma f(x,y) &= a^\alpha f(x^\sigma,y^\sigma)
                      = f(x^\sigma a,y^\sigma)\cr
                   &= f\bigl((xa)^\sigma,y^\sigma\bigr) = r_\sigma f(xa,y)
                      = r_\sigma a^\alpha f(x,y). \cr} $$
 Weil $f$ nicht ausgeartet ist, folgt $a^\alpha r_\sigma = r_\sigma a^\alpha$
 f\"ur alle $a \in K$, so dass $r_\sigma \in Z(K^*)$ ist. Hieraus folgt
 weiter, dass
 $$  r_{\sigma\tau}f(x,y) = f(x^{\sigma\tau},y^{\sigma\tau})
                          = r_\tau r_\sigma f(x,y)
                          = r_\sigma r_\tau f(x,y). $$
 Wiederum weil $f$ nicht ausgeartet ist, folgt $r_{\sigma\tau} = r_\sigma
 r_\tau$, so dass $r$ in der Tat ein Homomorphismus von $C^*(\pi)$ in
 $Z(K^*)$ ist.
 \par
       Es bleibt zu zeigen, dass $r$ nicht von $f$ abh\"angt. Es sei
 also $g$ eine zweite, $\pi$ darstellende Semibilinearform. Nach II.8.7b)
 gibt es ein $k \in K^*$ mit $g(x,y) = kf(x,y)$ f\"ur alle $x$, $y \in V$.
 Ist $\sigma \in C^*(\pi)$, so ist also, da ja $r_\sigma \in Z(K)$ gilt,
 $$ g(x^\sigma,y^\sigma) = kf(x^\sigma,y^\sigma) = kr_\sigma f(x,y)
                         = r_\sigma kf(x,y) = r_\sigma g(x,y), $$
 so dass $r$ in der Tat nicht von der Darstellung von $\pi$ abh\"angt.
 \medskip
       Den Kern der Abbildung $r$ bezeichnen wir mit $C(\pi)$. Er besteht
 aus allen Abbildungen $\sigma \in \GL(V)$ mit $f(x^\sigma,y^\sigma) = f(x,y)$
 f\"ur alle $x$, $y \in V$. Die Elemente aus $C(\pi)$ hei\ss en auch
 {\it Isometrien\/}\index{Isometrie}{} des Paares $(V,\pi)$. Die von $C(\pi)$ in
 $\La(V)$ induzierte Kollineationsgruppe bezeichnen wir mit $PC(\pi)$. Sie ist
 normal in $PC^*(\pi)$.
 \medskip\noindent
 {\bf 2.3. Satz.} {\it Es sei $V$ ein Vektorraum \"uber $K$ und $\pi$ sei eine
 Polarit\"at von $\La(V)$. Ferner sei $r$ der in 2.2 beschriebene
 Homomorphismus von $C^*(\pi)$ in $Z(K^*)$. Wird $\pi$ durch die
 $\alpha$-Semibilinearform $f$ dargestellt und ist $\sigma \in C^*(\pi)$,
 so ist $r_\sigma^\alpha = r_\sigma$.}
 \smallskip
       Beweis. Ist $\alpha = 1$, so ist nichts zu beweisen. Es sei also
 $\alpha \neq 1$. Dann ist $\pi$ nach Satz 1.4d) nicht symplektisch. Nach
 1.3 gibt es folglich eine symmetrische $\beta$-Semibilinearform $g$, die
 $\pi$ darstellt. F\"ur $g$ und $\beta$ gilt daher, da $r_\sigma$ im Zentrum
 von $K$ liegt,
 $$\eqalign{
 g(x,y)r_\sigma &= r_\sigma g(x,y) = g(x^\sigma,y^\sigma)
                 = g(y^\sigma,x^\sigma)^\beta \cr
                &= \bigl(r_\sigma g(y,x)\bigr)^\beta
                 = g(y,x)^\beta r_\sigma^\beta = g(x,y)r_\sigma^\beta. \cr}$$
 Weil $g$ nicht ausgeartet ist, folgt $r_\sigma = r_\sigma^\beta$. Nun
 unterscheiden sich $\alpha$ und $\beta$ nach II.8.7a) nur um einen
 inneren Automorphismus von $K$, so dass auch $r_\sigma = r_\sigma^\alpha$
 gilt, da $r_\sigma$ ja im Zentrum von $K$ liegt.
 \medskip\noindent
 {\bf 2.4. Satz.}  {\it Ist $\pi$ eine Polarit\"at des projektiven Verbandes
 $\La(V)$, so ist der Quotient $PC^*(\pi)/PC(\pi)$ eine elementarabelsche $2$-Gruppe.}
 \smallskip
       Beweis. Es sei $\sigma \in C^*(\pi)$ und $f$ sei eine $\pi$ darstellende
 $\alpha$-Semibilinearform. Wir definieren $\rho$ durch $v^\rho :=
 vr_\sigma^{-1}$ f\"ur $v \in V$. Dann ist, da ja $r_\sigma^\alpha =
 r_\sigma \in Z(K)$ gilt,
 $$ f(u^{\sigma^2\rho},v^{\sigma^2\rho}) = r_\sigma^{-\alpha} f(u^{\sigma^2},
                                          v^{\sigma^2})r_\sigma^{-1}
                                       = r_\sigma^{-2}r_\sigma^2 f(u,v)
                                       = f(u,v). $$
 Somit ist $\sigma^2\rho \in C(\pi)$. Weil $\sigma^2$ und $\sigma^2\rho$ in
 $\La(V)$ die gleiche Kol\-li\-ne\-a\-ti\-on induzieren, folgt, dass
 $PC^*(\pi)/PC(\pi)$ eine elementarabelsche 2-Grup\-pe ist.
 \medskip
       In einigen F\"allen k\"onnen wir noch weitergehende Aussagen ma\-chen.
 \medskip\noindent
 {\bf 2.5. Satz.} {\it Ist $V$ ein Vektorraum mit $\Rg_K(V) \geq 3$ \"uber dem
 endlichen K\"orper $K$ und ist $\pi$ eine unit\"are Polarit\"at von
 $\La(V)$, so ist $PC^*(\pi) = PC(\pi)$.}
 \smallskip
       Beweis. Weil $\pi$ nicht symplektisch ist, gibt es eine symmetrische
 $\alpha$-Semi\-bi\-line\-ar\-form $f$, die $\pi$ darstellt. Insbesondere ist dann
 $\alpha^2 = 1$ jedoch $\alpha \neq 1$. Hieraus folgt, dass $K = \GF(q^2)$
 ist und dass f\"ur alle $k \in K$ die Gleichung $k^\alpha = k^q$ gilt.
 \par
       Es sei $\sigma \in C^*(\pi)$. Dann ist $f(u^\sigma,v^\sigma) =
 r_\sigma f(u,v)$ f\"ur alle $u$, $v \in V$. Wegen $r^\alpha_\sigma =
 r_\sigma$ ist $r_\sigma \in \GF(q)$ und aus $K^{q+1} = \GF(q)$ folgt die
 Existenz eines $a \in K^*$ mit $r_\sigma = a^{q+1} = a^{\alpha + 1}$. Wir
 definieren $\rho$ durch $v^\rho := va^{-1}$ f\"ur alle $v \in V$. Dann ist
 $$ f(u^{\sigma\rho},v^{\sigma\rho}) = a^{-\alpha}f(u^\sigma,v^\sigma)a^{-1}
       = a^{-\alpha-1}r_\sigma f(u,v)
       = f(u,v),$$
 so dass $\sigma\rho \in C(\pi)$ ist. Weil $\sigma$ und $\sigma\rho$ in
 $\La(V)$ die gleiche Kollineation induzieren, folgt $PC^*(\pi) = PC(\pi)$,
 q. o. o.
 \medskip\noindent
 {\bf 2.6. Satz.} {\it Es sei $V$ ein Vektorraum \"uber $K$ mit $\Rg_K(V)
 \geq 3$.
 Ist $\pi$ eine symplektische Polarit\"at von $\La(V)$, so ist der in 2.2
 definierte Homomorphismus $r$ surjektiv. Insbesondere ist also
 $C^*(\pi)/C(\pi) \cong K^*$ und $PC^*(\pi)/PC(\pi) \cong K^*/K^{*2}$.}
 \smallskip
       Beweis. Es sei $f$ eine $\pi$ darstellende Bilinearform. Nach 1.7 ist
 $\Rg_K(V) = 2n$ und nach 1.6 besitzt $V$ eine symplektische Basis $b_1$,
 \dots, $b_{2n}$. F\"ur $a \in K^*$ definieren wir die Abbildung $\sigma \in
 \GL(V)$ verm\"oge $b^\sigma_{2i-1} := -b_{2i}$ und $b_{2i}^\sigma :=
 b_{2i-1}a$ f\"ur $i := 1$, \dots, $n$. Eine leichte, wenn auch l\"angliche
 Rechnung --- man nehme ein Blatt in Querformat --- zeigt, dass
 $f(u^\sigma,v^\sigma) = af(u,v)$ ist. Folglich ist $\sigma \in C^*(\pi)$ und
 $r_\sigma = a$. Dies zeigt die Surjektivit\"at von $r$. Hieraus folgt
 unmittelbar die vorletzte Aussage des Satzes.
 \par
       Um die letzte Aussage zu beweisen, sei $\sigma \in C^*(\pi)$ und
 $N$ be\-zeich\-ne den Kern des Homomorphismus $P$ von $C^*(\pi)$ auf
 $PC^*(\pi)$. Wir setzen
 $$ \varphi(\sigma N) := r_\sigma K^{*2}. $$
 Ist $\tau \in \sigma N$, so ist $v^{\sigma\tau^{-1}} = va$ f\"ur alle
 $v \in V$ mit einem von $v$ unabh\"angigem $a \in K^*$. Es folgt
 $$
 r_\sigma r_\tau^{-1}f(u,v) = r_{\sigma\tau^{-1}}f(u,v) 
       = f(u^{\sigma\tau^{-1}},v^{\sigma\tau^{-1}}) = f(ua,va) 
       = a^2f(u,v) 
 $$
 Also ist $r_\sigma = r_\tau a^2$, so dass $\tau$ wohldefiniert ist.
 Nach dem bereits Bewiesenen ist $\varphi$ surjektiv.
 \par
       Genau dann ist $\varphi(\sigma N) = K^{*2}$, wenn $r_\sigma = k^2$ ist
 mit einem $k \in K^*$. Insbesondere ist also $PC(\pi) \subseteq \Kern(
 \varphi)$. Es sei $\sigma N \in \Kern(\varphi)$. Dann ist $r_\sigma = k^2$.
 Definiere $\tau$ durch $v^\tau := v^\sigma k^{-1}$. Dann ist
 $$ f(u^\tau,v^\tau) = k^{-2}f(u^\sigma,v^\sigma) = k^{-2}r_\sigma f(u,v)
       = f(u,v) $$
 f\"ur alle $u$, $v \in V$. Es folgt $\sigma N = \tau N \in PC(\pi)$. Damit
 ist alles gezeigt.
 \medskip\noindent
 {\bf 2.7. Satz.} {\it Es sei $V$ ein Vektorraum des Ranges $n \geq 3$ \"uber
 dem kommutativen K\"orper $K$. Ferner sei $\pi$ eine Polarit\"at von
 $\La(V)$, die durch die $\alpha$-Semibilinearform $f$ dargestellt werde.
 Ist $\sigma \in C^*(\pi)$, so gibt es ein $a \in K^*$ mit $r_\sigma^n
 = a^{\alpha + 1}$.}
 \smallskip
       Beweis. Es sei $b_1$, \dots, $b_n$ eine Basis von $V$ und $F$ sei
 die durch
 $ F_{ij} := f(b_i,b_j) $
 definierte Matrix. Ferner sei $c$ die durch
 $ b_i^\sigma := \sum_{r:=1}^n b_rc_{ri} $
 definierte Matrix. Schlie\ss lich sei $a := \det(F)$. Dann ist
 $$
 r_\sigma f(b_i,b_j) = f(b_i^\sigma,b_j^\sigma)   
                     = f\biggl(\sum_{r:=1}^n b_rc_{ri},\sum_{s:=1}^n b_sc_{sj}
                                                        \biggr) 
                     = \sum_{r:=1}^n \sum_{s:=1}^n c^\alpha_{ri}
                            f(b_r,b_s)c_{sj}. 
 $$
 Also ist $r_\sigma F = c^{\alpha t}Fc$ --- $t$ bedeutet transponieren. Daher ist
 $$ r_\sigma^n\det(F) = \det(c^\alpha)\det(c)\det(F). $$
 Nun ist $f$ nicht ausgeartet, woraus folgt, dass $\det(F) \neq 0$ ist.
 Somit ist $r_\sigma^n = \det(c)^\alpha\det(c) = a^{\alpha + 1}$, q. e. d.
 \medskip\noindent
 {\bf 2.8. Satz.} {\it Es sei $V$ ein Vektorraum \"uber dem kommutativen K\"orper
 $K$ mit $\Rg_K(V) \geq 3$, und $\pi$ sei eine Polarit\"at von $\La(V)$, die
 von der $\alpha$-Semibilinearform $f$ dargestellt werde. Ist der Rang von
 $V$ ungerade, so ist $C^*(\pi)/C(\pi) \cong (K^*)^{1+\alpha}$ und daher
 $PC^*(\pi) = PC(\pi)$.}
 \smallskip
       Beweis. Setze $n := \Rg_K(V)$. Dann ist $n = 2k + 1$. Nach 2.7 gilt
 $r_\sigma^n = a^{\alpha+1}$ mit einem $a \in K^*$. Wir setzen $b :=
 r_\sigma^{-k}$. Nach 2.3 ist dann $b^\alpha = b$. Folglich ist
 $$r_\sigma = r_\sigma^{2k+1}b^2 = a^{\alpha+1}b^\alpha b = (ab)^{\alpha+1}. $$
 Dies besagt, dass $r$ ein Homomorphismus von $C^*(\pi)$ in $(K^*)^{1+
 \alpha}$ ist. Weil andererseits $(K^*)^{\alpha+1}$ im Bild von $r$ enthalten
 ist, folgt, dass $r$ ein Epimorphismus von $C^*(\pi)$ auf $(K^*)^{\alpha+1}$
 ist. Ist nun $\sigma \in C^*(\pi)$, so gibt es also ein $c \in K^*$ mit
 $r_\sigma = c^{\alpha+1}$. Definiert man $\rho$ verm\"oge $v^\rho :=
 vc^{-1}$ f\"ur alle $v \in V$, so folgt $\sigma\rho \in C(\pi)$. Weil $\sigma$
 und $\sigma\rho$ die gleiche Kollineation in $\La(V)$ hervorrufen, ist daher
 $PC^*(\pi) = PC(\pi)$. Damit ist alles bewiesen.
 \medskip\noindent
 {\bf 2.9. Satz.} {\it Es sei $K$ ein endlicher K\"orper und $Q$ sei die Menge
 der Quadrate von $K$. Sind $a$ und $b$ zwei von $0$ verschiedene Elemente von
 $K$, so ist $K = aQ + bQ$.}
 \smallskip
       Beweis. Ist die Charakteristik von $K$ gleich 2, so ist $K = Q$ und
 weiter nichts zu beweisen. Es sei also die Charakteristik von $K$ von 2
 verschieden. Dann enth\"alt $K^*$ genau ${1 \over 2}(q - 1)$ Quadrate, so
 dass $Q$ wegen $0^2 = 0$ genau ${1 \over 2}(q + 1)$ Elemente enth\"alt. Ist
 nun $k \in K$, so gelten wegen $a$, $b \neq 0$ die Gleichungen
 $$ |aQ| =\ {1 \over 2}(q + 1) = |-{}bQ + k|. $$
 Somit ist
 $$ \eqalign{
    q \geq \big|(aQ) \cup (-{}bQ + k)\big|                 
      & = |aQ| + |-{}bQ + k| - \big|(aQ) \cap (-bQ + k)\big| \cr
      &= q + 1 - \big|(aQ) \cap (-bQ + k)\big|.            \cr} $$
 Hieraus folgt, dass $(aQ) \cap (-bQ + k) \neq \emptyset$ ist. Es gibt also
 Elemente $g$, $h \in Q$ mit $ag = -bh + k$, so dass $k = ag + bh$ ist.
 \medskip\noindent
 {\bf 2.10. Satz.} {\it Es sei $V$ ein Vektorraum geraden Ranges \"uber dem
 endlichen K\"orper $K$. Ist $\pi$ eine projektive Polarit\"at von $\La(V)$,
 so ist $C^*(\pi)/C(\pi) \cong K^*$ und daher $PC^*(\pi)/PC(\pi) \cong
 K^*/K^{*2}$.}
 \smallskip
       Beweis. Dies ist nach 2.6 richtig, falls $\pi$ symplektisch ist.
 Ist die Charakteristik von $K$ gleich 2, so ist jedes Element von $K$ ein
 Quadrat. Weil das Bild von $r$ die Gruppe der Quadrate umfasst, folgt
 auch in diesem Falle, dass $C^*(\pi)/C(\pi) \cong K^*$ ist, was wiederum
 $PC^*(\pi)/PC(\pi) \cong K^*/K^{*2}$ nach sich zieht. Wir d\"urfen daher
 annehmen, dass $\pi$ nicht symplektisch ist und dass die Charakteristik
 von $K$ von 2 verschieden ist.
 \par
       Weil $\pi$ nicht symplektisch ist, gibt es einen Punkt $P$ mit $P
 \not\leq P^\pi$. Nach 1.5 ist daher $V = P \oplus P^\pi$. Es sei $f$ eine
 $\pi$ darstellende symmetrische Bilinearform. W\"are $f(u,u) = 0$ f\"ur
 alle $u \in P^\pi$, so w\"are
 $$ 0 = f(u + u',u + u') = 2f(u,u') $$
 f\"ur alle $u$, $u' \in P^\pi$. Weil die Charakteristik von $K$ von 2
 verschieden ist, ist daher $f(u,u') = 0$ f\"ur alle $u$, $u' \in P^\pi$.
 Ist $v \in V$, so ist $v = p + u$ mit $p \in P$ und $u \in P^\pi$. Ist $x \in
 P^\pi$, so folgt
 $$ f(x,v) = f(x,p + u) = f(x,p) + f(x,u) = 0. $$
 Da dies f\"ur alle $v \in V$ gilt, folgt $x = 0$, da $f$ ja nicht ausgeartet
 ist. Also ist $P^\pi = \{0\}$ und somit $V = P$. Dieser Widerspruch zeigt,
 dass es einen Punkt $Q \leq P^\pi$ gibt, der nicht absolut ist. Es sei
 $G := P + Q$ und $P = pK$ sowie $Q = qK$. Indem man $f$ notfalls mit einem
 Skalarfaktor ab\"andert --- Hier benutzen wir, dass $r$ nicht von der
 Auswahl von $f$ abh\"angt ---, kann man erreichen, dass $f(p,p) = 1$ ist.
 Ist $x := pk + ql \in G \cap G^\pi$, so ist $0 = f(p,x) = k$ und $0 = f(q,x)
 = f(q,q)l$, so dass $k = l = 0$ ist. Also ist $G \cap G^\pi = \{0\}$
 und daher $V = G \oplus G^\pi$.
 \par
       Es sei $s \in K^*$ und $s$ sei kein Quadrat.
 \par
       1. Fall: Es sei $f(q,q)$ kein Quadrat. Dann ist $sf(q,q)$ ein Quadrat.
 Es gibt also Elemente $b$, $c \in K^*$ mit $s = c^2f(q,q)$ und $sf(q,q) =
 b^2$. Setze ferner $a := 0$ und $d := 0$.
 \par
       2. Fall: Es sei $f(q,q)$ ein Quadrat. In diesem Falle d\"urfen wir
 annehmen, dass $f(q,q) = 1$ ist. Nach 2.9 gibt es Elemente $a$, $c \in K$
 mit $s = a^2 + c^2$. In diesem Falle setzen wir $b := c$ und $d := -a$.
 \par
       Wir definieren nun eine Abbildung $\sigma'$ auf $G$ verm\"oge
 $p^{\sigma'} := pa + qc$ und $q^{\sigma'} := pb + qd$. Dann ist, wie leicht
 nachzurechnen, $f(g^{\sigma'},h^{\sigma'}) = sf(g,h)$ f\"ur alle $g$, $h
 \in G$. Wegen $\Rg_K(V) = 2 + \Rg_K(G^\pi)$ ist auch $\Rg_K(V^\pi)$
 gerade. Es gibt daher nach In\-duk\-ti\-ons\-an\-nah\-me eine Abbildung
 $\sigma'' \in \GL(G^\pi)$ mit $f(u^{\sigma''}, v^{\sigma''}) = sf(u,v)$ f\"ur
 alle $u$, $v \in G^\pi$. Ist schlie\ss lich $v \in V$ und $v = x + y$ mit
 $x \in G$ und $y \in G^\pi$ und definiert man $\sigma$ durch $v^\sigma :=
 x^{\sigma'} + x^{\sigma''}$, so zeigt eine triviale Rechnung, dass
 $f(u^\sigma,v^\sigma) = sf(u,v)$ f\"ur alle $u$, $v \in V$ ist. Da die
 Quadrate stets im Bilde von $r$ liegen, ist damit gezeigt, dass $r$
 surjektiv ist.

\mysection{3. Symplektische Polarit\"aten und ihre Zentralisatoren}

\noindent
       Ist $\pi$ eine symplektische Polarit\"at von $\La_K(V)$, so setzen
 wir $\Sp(V) := C(\pi)$ und $\PSp(V) := PC(\pi)$. Soll der Rang aus irgendeinem
 Grunde erw\"ahnt werden, so schreiben wir statt $\Sp(V)$ und $\PSp(V)$ auch
 $\Sp(n,K)$ bzw. $\PSp(n,K)$. Diese Bezeichnungsweisen sind gerechtfertigt,
 da die Struktur von $\Sp(V)$ und $\PSp(V)$ nach 1.7 nur von $V$ abh\"angt.
 Ist $f$ eine Form, die $\pi$ darstellt, so ist $f$ also eine nicht
 ausgeartete, alternierende Bilinearform und $\Sp(V)$ ist die Gruppe aller
 $\sigma \in \GL(V)$ mit $f(u^\sigma,v^\sigma) = f(u,v)$ f\"ur alle $u$, $v \in
 V$. Wir nennen beide Gruppen, $\Sp(V)$ und $\PSp(V)$
 {\it symplektisch\/}.\index{symplektische Gruppe}{}
 \medskip\noindent
 {\bf 3.1. Satz.} {\it Es sei $V$ ein Vektorraum geraden Ranges \"uber dem
 kommutativen K\"orper $K$. Ist $\sigma \in \GL(V)$, so sind die folgenden
 Aussagen \"aquivalent:
 \item{a)} Es ist $\sigma \in \Sp(V)$.
 \item{b)} Symplektische Basen werden von $\sigma$ auf symplektische Basen
 abgebildet.
 \item{c)} Es gibt eine symplektische Basis, die von $\sigma$ auf eine
 symplektische Basis abgebildet wird.

 \noindent
       Ferner ist die Gruppe $\Sp(V)$ auf der Menge der symplektischen
 Basen scharf transitiv.}
 \smallskip
       Beweis. Ist $b_1, \dots, b_{2n}$ eine Basis von $V$, so ist
 $$ f(b_{2i-1},b_{2i}) = 1 = -f(b_{2i},b_{2i-1}) $$
 f\"ur $i := 1, \dots, n$ und $f(b_i,b_j) = 0$ f\"ur alle \"ubrigen
 Indexpaare kennzeichnend daf\"ur, dass $b_1, \dots, b_{2n}$ eine
 symplektische Basis von $V$ ist. Dabei ist $f$ eine Bilinearform, die
 eine symplektische Polarit\"at darstellt. Ist $\sigma \in \Sp(V)$, so ist
 also mit $b_1, \dots, b_{2n}$ auch $b_1^\sigma, \dots, b_{2n}^\sigma$
 eine symplektische Basis. Es ist also b) eine Folge von a).
 \par
       Da es nach 1.6 stets symplektische Basen gibt, ist c) eine Folge von b).
 \par
       Es sei $\sigma \in \GL(V)$ und $b_1, \dots, b_{2n}$ und
 $b_1^\sigma, \dots, b_{2n}^\sigma$ seien symplektische Basen von $V$.
 Ist $u = \sum_{i:=1}^{2n} b_ix_i$ und $v = \sum_{i:=1}^{2n} b_iy_i$, so ist
 folglich
 $$ f(u,v) = \sum_{i:=1}^n (x_{2i-1}y_{2i} - x_{2i}y_{2i-1}) $$
 und
 $$ f(u^\sigma,v^\sigma) = \sum_{i:=1}^n (x_{2i-1}y_{2i} - x_{2i}y_{2i-1}), $$
 so dass $f(u^\sigma,v^\sigma) = f(u,v)$ ist. Dies besagt aber gerade,
 dass $\sigma$ ein Element der symplektischen Gruppe ist. Also ist a) eine Folge
 von c).
 \par
       Sind schlie\ss lich $b_1, \dots, b_{2n}$ und $c_1, \dots, c_{2n}$
 zwei symplektische Basen von $V$, so gibt es genau ein $\sigma \in \GL(V)$
 mit $b_i^\sigma = c_i$ f\"ur alle $i$, so dass $\Sp(V)$ nach dem bereits
 Bewiesenen auf der Menge der symplektischen Basen scharf transitiv
 operiert. Damit ist alles bewiesen.
 \medskip\noindent
 {\bf 3.2. Satz.} {\it Es sei $V$ ein Vektorraum \"uber dem kommutativen K\"orper
 $K$ und $f$ sei eine nicht ausgeartete, alternierende Bilinearform auf $V$.
 Sind dann $\{b_1, b_3, \dots,$ $ b_{2r+1}\}$ und $\{b_2, b_4, \dots,
 b_{2s}\}$ Mengen von linear unabh\"angigen Vektoren mit
 $$ f(b_{2i-1},b_{2i}) = 1 = -f(b_{2i},b_{2i-1}) $$
 f\"ur $i := 1, \dots, \min\{r,s\}$ sowie $f(b_i,b_j) = 0$ f\"ur alle
 \"ubrigen Indexpaare, so gibt es eine symplektische Basis $c_1, \dots,
 c_{2n}$ mit $c_{2i-1} = b_{2i-1}$ f\"ur $i := 1, \dots, r$ und
 $c_{2i} = b_{2i}$ f\"ur $i := 1, \dots, s$.}
 \smallskip
       Beweis. Wir d\"urfen annehmen, dass weder $\{b_1, \dots, b_{2r+1}\}$
 noch $\{b_2, \dots, b_{2s}\}$ leer ist. Sind n\"amlich beide Mengen leer,
 so besagt der Satz nichts anderes als die Existenz von symplektischen
 Basen. Ist eine der beiden Mengen, etwa die zweite, leer, so ist
 $$ b_3K + b_5K + \dots b_{2r+1}K < b_1K + b_3K + \dots b_{2r+1}K, $$
 so dass also
 $$ (b_1K + b_3K + \dots b_{2r+1}K)^\pi <
                                 (b_3K + b_5K + \dots b_{2r+1}K)^\pi $$
 ist. Dabei bezeichne $\pi$ die von $f$ induzierte symplektische
 Polarit\"at.
 Es gibt also einen Vektor $y \in (b_3K + \dots + b_{2r+1}K)^\pi$,
 der nicht in $(b_1K + \dots b_{2r+1}K)^\pi$ liegt. Setze
 $b_2 := f(b_1,y)^{-1}y$. Dann ist $f(b_1,b_2) = 1$ und $f(b_{2i+1},b_2) = 0$
 f\"ur $i := 1, \dots, r$. Damit ist die Zwischenbehauptung bewiesen.
 \par
       Ist nun $\Rg_K(V) = 2$, so ist $b_1, b_2$ eine symplektische Basis
 von $V$, so dass der Satz in diesem Falle bewiesen ist. Es sei also
 $\Rg_K(V) > 2$. Nun ist $G := b_1K + b_2K$ eine Gerade von $\La(V)$ mit
 $V = G \oplus G^\pi$, wie wir schon einmal bemerkten. Auf Grund unserer
 Annahme gilt
 $$ \{b_3, \dots, b_{2r+1}\} \cup \{b_4, \dots, b_{2s}\} \subseteq G^\pi. $$
 Wegen $\Rg_K(G) = \Rg_K(V) - 2$ f\"uhrt Induktion zum Ziele.
 \medskip
       Statt $\Sp(2n,\GF(q))$ bzw. $\PSp(2n,\GF(q))$ schreiben wir im Folgenden auch
 $\Sp(2n,q)$ bzw. $\PSp(2n,q)$. Man beachte, dass der Rang eines projektiven
 Geometrie mit einer symplektischen Polarit\"at stets gerade ist.
 \medskip\noindent
 {\bf 3.3. Satz.} {\it Es ist
 $$ \big|\Sp(2n,q)\big| = q^{n^2} \prod_{i:=1}^n (q^{2i} - 1) $$
 und
 $$ \big|\PSp(2n,q)\big| = {1 \over \ggT(2,q - 1)}\big|\Sp(2n,q)\big|. $$}
 \smallskip
       Beweis. Weil $\Sp(2n,q)$ auf den symplektischen Basen des
 zu\-ge\-h\"o\-ri\-gen
 Vektorraumes scharf transitiv operiert, ist $|\Sp(2n,q)|$ gleich der Anzahl
 dieser Basen. Ist $b_1, b_2, \dots, b_{2n}$ eine symplektische Basis
 von $V$, so ist $b_3, \dots, b_{2n}$ eine symplektische Basis von
 $(b_1K + b_2K)^\pi$, so dass die Anzahl der symplektischen Basen von $V$
 gleich der Anzahl der Paare $b_1, b_2$ mit $f(b_1,b_2) = 1$ mal der
 Anzahl der symplektischen Basen eines Vektorraumes vom Range $2n - 2$ ist.
 Ist nun $b_1$ ein vom Nullvektor verschiedener Vektor und ist $xK$ ein Punkt,
 der nicht in $(b_1K)^\pi$ liegt, so ist $f(b_1,x) \neq 0$. Es gibt also
 genau ein $b_2 \in xK$ mit $f(b_1,b_2) = 1$. Die Anzahl der Vektoren $b_2$
 mit $f(b_1,b_2) = 1$ ist somit gleich der Anzahl der Punkte von $\La(V)$,
 die nicht in der Hyperebene $(b_1K)^\pi$ liegen, dh. gleich $q^{2n-1}$.
 Folglich ist die Anzahl der Paare $b_1$, $b_2$ gleich $(q^{2n} - 1)q^{2n-1}$.
 Daher ist die Anzahl der symplektischen Basen von $V$ gleich
 $$ (q^{2n} - 1)q^{2n-1}q^{(n-1)^2}\prod_{i:=1}^{n-1} (q^{2i} - 1). $$
 Also ist
 $$ \big|\Sp(2n,q)\big| = q^{n^2}\prod_{i:=1}^n (q^{2i} - 1). $$
 Ist $\sigma$ im Kern des Homomorphismus von $\Sp(2n,q)$ auf $\PSp(2n,q)$, so ist
 $v^\sigma = vk$ f\"ur alle $v \in V$ und einem geeigneten $k \in K^*$. Es
 folgt $f(u,v) = f(u^\sigma,v^\sigma) = f(u,v)k^2$ und weiter $k^2 = 1$.
 Also ist $k = 1$ oder $k = -1$. Hieraus folgt auch noch die letzte
 Behauptung des Satzes.
 \medskip
       Es sei $\pi$ eine Polarit\"at von $\La(V)$. Ist $U \in \La(V)$ und gilt
 $U \cap U^\pi \neq \{0\}$, so hei\ss t $U$ {\it isotrop\/}.\index{isotrop}{}
 Ist sogar $U \leq U^\pi$, so hei\ss t $U$ {\it vollst\"andig
 isotrop\/}.\index{vollst\"andig isotrop}{} Ist $U$ vollst\"andig isotrop, so ist
 $$ 2\Rg_K(U) \leq \Rg_K(U) + \Rg_K(U^\pi) = \Rg_K(V) $$
 und daher $\Rg_K(U) \leq {1 \over 2}\Rg_K(V)$.
 \medskip\noindent
 {\bf 3.4. Satz.} {\it Es sei $V$ ein Vektorraum des Ranges $2n$ \"uber dem
 kommutativen K\"orper $K$ und $\pi$ sei eine symplektische Polarit\"at
 von $\La(V)$. Ist $U \in \La(V)$ vollst\"andig isotrop, so gibt es einen
 vollst\"andig isotropen Unterraum $W$ von $V$ mit $U \leq W$ und $W = W^\pi$,
 dh. mit $\Rg_K(W) = n$.}
 \smallskip
       Beweis. Es sei $\Rg_K(U) = r + 1$ und $b_1$, $b_3$, \dots, $b_{2r+1}$
 sei eine Basis von $U$. Weil $U$ vollst\"andig isotrop ist, erf\"ullen dann
 $\{b_1, b_3, \dots, b_{2r+1}\}$ und $\emptyset$ die Voraussetzungen von
 3.2. Es gibt daher eine symplektische Basis $c_1$, $c_2$, \dots, $c_{2n}$ von
 $V$ mit $b_{2i+1} = c_{2i+1}$ f\"ur $i := 0$, \dots, $r$. Setze $W :=
 \sum_{i:=1}^n c_{2i-1}K$. Dann ist $U \leq W$ und $W = W^\pi$, q. e. d.
 \medskip\noindent
 {\bf 3.5. Satz.} {\it Es sei $V$ ein Vektorraum und $\pi$ sei eine symplektische
 Polarit\"at von $\La(V)$. Ferner seien $U$, $W \in \La(V)$. Genau dann
 gibt es ein $\sigma \in \Sp(V)$ mit $U^\sigma = W$, wenn $\Rg_K(U) = \Rg_K(W)$
 und $\Rg_K(U \cap U^\pi) = \Rg_K(W \cap W^\pi)$ ist.}
 \smallskip
       Beweis. Ist $\sigma \in \Sp(V)$ und ist $U^\sigma = W$, so folgt aus
 $\sigma\pi = \pi\sigma$, dass $U^{\pi\sigma} = W^\pi$ ist. Folglich ist
 auch $(U \cap U^\pi)^\sigma = W \cap W^\pi$,  so dass in der Tat
 $\Rg_K(U) = \Rg_K(W)$ und $\Rg_K(U \cap U^\pi) = \Rg_K(W \cap W^\pi)$
 ist.
 \par
       Es sei umgekehrt $\Rg_K(U) = \Rg_K(W)$ und $\Rg_K(U \cap U^\pi) =
 \Rg_K(W \cap W^\pi)$. Ist nun $U = U_0 \oplus (U \cap U^\pi)$ und
 $W = W_0 \oplus (W \cap W^\pi)$, so ist also $\Rg_K(U_0) = \Rg_K(W_0)$ und
 $U_0$ und $W_0$ sind beide nicht isotrop. Es gibt also symplektische Basen
 $b_1, \dots, b_{2s}$ von $U_0$ bzw. $b'_1, \dots, b'_{2s}$ von $W_0$.
 Es sei ferner $b_{2s+1}, b_{2s+3}, \dots, b_{2r+1}$ eine Basis von
 $U \cap U^\pi$ und $b'_{2s+1}, b'_{2s+3}, \dots, b'_{2r+1}$ eine solche
 von $W \cap W^\pi$. Dann erf\"ullen $\{b_1, b_3, \dots, b_{2r+1}\}$ und
 $\{b_2, b_4, \dots, b_{2s}\}$ bzw.
 $\{b'_1, b'_3, \dots, b_{2r+1}\}$ und
 $\{b'_2, b'_4, \dots, b'_{2s}\}$ die Voraussetzungen von 3.2, so dass
 aus diesem Satz und aus 3.1 die Existenz eines $\sigma \in \Sp(V)$ mit
 $U^\sigma = W$ folgt.
 \medskip\noindent
 {\bf 3.6. Satz.} {\it Ist $V$ ein Vektorraum des Ranges $2$ \"uber dem
 kommutativen
 K\"orper $K$, so ist $\Sp(V) = \SL(V)$ und somit auch $\PSp(V) = \PSL(V)$.}
 \smallskip
       Beweis. Es sei $b_1$, $b_2$ eine symplektische Basis von $V$ und
 $\sigma \in \GL(V)$. Ist $b_i^\sigma = b_1a_{1i} + b_2a_{2i}$, so ist
 $$ f(b_1^\sigma,b_2^\sigma) = a_{11}a_{22} - a_{12}a_{21} = \det(\sigma). $$
 Folglich ist $b_1^\sigma$, $b_2^\sigma$ genau dann eine symplektische
 Basis von $V$, wenn $\det(\sigma) = 1$, dh. wenn $\sigma \in \SL(V)$ ist
 (Satz III.1.12b)). Aus 3.1 folgt daher die Behauptung.
 \medskip\noindent
 {\bf 3.7. Satz.} {\it Es sei $V$ ein Vektorraum \"uber dem kommutativen
 K\"orper $K$ und $f$ sei eine nicht ausgeartete, alternierende Bilinearform
 auf $V$.
 \item{a)} Ist $a \in V$ und $k \in K$, so liegt die durch $x^\tau := x + 
 af(a,x)k$ definierte Transvektion $\tau$ in $\Sp(V)$.
 \item{b)} Ist $\varphi$ eine lineare Abbildung von $V$ in $K$, ist ferner
 $a \in \Kern(\varphi)$ und liegt die durch $x^\tau := x + a\varphi(x)$
 definierte Transvektion $\tau$ in $\Sp(V)$, so gibt es ein $k \in K$ mit
 $\varphi(x) = f(a,x)k$ f\"ur alle $x \in V$.\par}
 \smallskip
       Beweis. a) Es ist klar, dass die Abbildung $x \to f(a,x)k$ linear
 ist, da $K$ ja kommutativ ist. Wegen $f(a,a) = 0$ ist $\tau$ folglich eine
 Transvektion. Schlie\ss lich ist
 $$
 f(x^\tau,y^\tau) = f(x,y) + f(x,a)f(a,y)k + f(a,y)f(a,x)k 
                           + f(a,a)f(a,x)f(a,y)k^2, 
 $$
 so dass wegen $f(a,a) = 0$ und $f(a,x) = -f(x,a)$ f\"ur alle $x$, $y \in V$
 die Gleichung $f(x^\tau,y^\tau) = f(x,y)$ gilt. Folglich ist $\tau \in
 \Sp(V)$.
 \par
       b) Wir d\"urfen annehmen, dass $\tau \neq 1$ ist. Dann ist $aK$ das
 Zentrum der von $\tau$ induzierten Elation. Ist $\pi$ die von $f$ induzierte
 symplektische Polarit\"at, so folgt aus $\pi\tau = \tau\pi$, dass
 $(aK)^\pi$ die Achse von $\tau$ ist. Daher ist
 $$ \Kern(\varphi) = (aK)^\pi = \big\{x \mid x \in V, f(a,x) = 0\big\}. $$
 Hieraus folgt, dass es ein $k \in K$ gibt mit $\varphi(x) = f(a,x)k$
 f\"ur alle $x \in V$. Damit ist der Satz bewiesen.
 \medskip\noindent
 {\bf 3.8. Satz.} {\it Die Gruppe $\Sp(2n,K)$ wird von Transvektionen erzeugt.
 Insbesondere ist $\det(\sigma) = 1$ f\"ur alle $\sigma \in \Sp(2n,K)$.}
 \smallskip
       Beweis. Es sei $T$ die von allen in $\Sp(2n,K)$ liegenden Transvektionen
 er\-zeug\-te Untergruppe von $\Sp(2n,K)$. Wir zeigen zun\"achst, dass $T$ auf
 der Menge der von $0$ verschiedenen Vektoren von $V$ transitiv operiert.
 Es seien also $a$, $b$ zwei von 0 verschiedene Vektoren aus $V$. Ist $f(a,b)
 \neq 0$, so gibt es ein $k \in K$ mit $f(a,b)k = -1$. Die durch
 $$ v^\tau := v + (b - a)f(b - a,v)k $$
 definierte Transvektion liegt nach 3.7a) in $T$. Ferner ist
 $$ a^\tau = a + (b - a)f(b - a,a)k = a + (b - a)f(b,a)k = b $$
 In diesem Falle gibt es also sogar eine Transvektion, die $a$ auf $b$
 abbildet. Ist $f(a,b) = 0$, so gibt es einen Vektor $c$ mit $f(a,c) \neq 0
 \neq f(b,c)$. Nach dem eben Bewiesenen gibt es dann zwei Transvektionen
 $\sigma$, $\tau \in T$ mit $a^\sigma = c$ und $c^\tau = b$, so dass es
 auch in diesem Falle ein Element in $T$ gibt, n\"amlich $\sigma\tau$, welches
 $a$ auf $b$ abbildet.
 \par
       Es seien nun $a$, $b$ und $a'$, $b'$ zwei Paare von Vektoren mit
 $f(a,b) = 1$ und $f(a',b') = 1$. Wir zeigen, dass es ein $\rho \in T$ gibt,
 mit $a^\rho = a'$ und $b^\rho = b'$. Weil $T$ auf $V - \{0\}$ transitiv
 operiert, d\"urfen wir $a = a'$ annehmen. Ist $f(b',b) \neq 0$, so gibt es
 ein $k \in K$ mit $f(b',b)k = 1$. In diesem Falle definieren wir $\tau$
 durch
 $$ x^\tau := x + (b' - b)f(b' - b,x)k. $$
 Dann ist
 $$\eqalign{
 a^\tau &= a + (b' - b)\bigl(f(b',a) - f(b,a)\bigr)k  \cr
        &= a' + (b' - b)\bigl(f(b',a') - f(b,a)\bigr)k \cr
        &= a' + (b' - b)\bigl(-1 + 1\bigr)k            
        = a' \cr} $$
 und
 $$ b^\tau = b + (b' - b)\bigl(f(b',b) - f(b,b)\bigr)k = b'. $$
 Ist $f(b',b) = 0$, so betrachten wir die beiden Paare $a$, $b$ und $a$,
 $a + b$ sowie $a$, $a + b$ und $a$, $b'$. Dann ist $f(b,a + b) = -1 \neq 0$
 und $f(a + b,b') = f(a,b') = f(a',b') = 1 \neq 0$, so dass es zwei
 Transvektionen $\sigma$, $\tau \in T$ gibt mit $a^\sigma = a = a^\tau$ und
 $b^\sigma = a + b$ und $(a + b)^\tau = b'$. Es folgt $a^{\sigma\tau} = a'$
 und $b^{\sigma\tau} = b'$.
 \par
       Ist nun $n = 1$, so ist 3.8 wegen 3.6 richtig. Es sei also $n > 1$ und
 $\gamma \in \Sp(2n,K)$. Es gibt $a$, $b \in V$ mit $f(a,b) = 1$. Dann ist
 auch $f(a^\gamma,b^\gamma) = 1$, so dass es ein $\rho \in T$ gibt mit
 $a^{\gamma\rho} = a$ und $b^{\gamma\rho} = b$. Wegen $f(a,b) = 1$ ist
 $aK + bK$ nicht isotrop. Ist $\pi$ die durch $f$ induzierte symplektische
 Polarit\"at und $U := aK + bK$, so ist also $V = U \oplus U^\pi$ und
 $U^{\pi\gamma\rho} = U^\pi$, da ja $\gamma\rho \in \Sp(2n,K)$. Wegen
 $\Rg_K(U^\pi) = 2(n - 1)$ gibt es dann Transvektionen $\lambda_1, \dots,
 \lambda_r \in \Sp(2n - 2,K)$ mit $x^{\gamma\rho} = x^\sigma$ f\"ur alle
 $x \in U^\pi$, wobei $\sigma := \lambda_1 \cdots \lambda_r$ gesetzt wurde.
 Definiert man nun $\tau_i$ durch
 $$ (u + x)^{\tau_i} := u + x^{\lambda_i}, $$
 so ist $\tau_i$ eine Transvektion aus $\Sp(2n,K)$ und es ist
 $\gamma\rho = \tau_1 \cdots \tau_r$, so dass also $\gamma \in T$ gilt.
 Damit ist der Satz bewiesen.
 \medskip
       Zur Erinnerung: Ist $G$ eine Gruppe, so bezeichnen wir mit $G'$ die
 Kommutatorgruppe von $G$.
 \medskip\noindent
 {\bf 3.9. Satz.} {\it Ist $K$ ein kommutativer K\"orper und $n$ eine
 nat\"urliche Zahl, so ist $\Sp(2n,K) = \Sp(2n,K)'$, es sei denn, es ist $n = 1$
 und $|K| \leq 3$ oder $n = 2$ und $|K| = 2$.}
 \smallskip
       Beweis. Wir betrachten zun\"achst den Fall $|K| > 3$. Es gibt dann
 ein $k \in K$ mit $k^2 \neq 0$, 1. Es sei $0 \neq a \in V$. Es gibt dann ein
 $\sigma \in \Sp(2n,K)$ mit $a^\sigma = ak$. Es sei $\tau$ die durch
 $x^\tau := x + af(a,x)l$ mit $l \in K$ definierte Transvektion. Dann ist
 $$\eqalign{
 x^{\sigma^{-1}\tau^{-1}\sigma\tau}
       &= \bigl(x^{\sigma^{-1}} - af(a,x^{\sigma^{-1}})l\bigr)^{\sigma\tau} \cr
       &= \bigl(x - a^\sigma f(a,x^{\sigma^{-1}})l\bigr)^\tau               \cr
       &= \bigl(x - akf(a,x^{\sigma^{-1}})l\bigr)^\tau                      \cr
       &= x^\tau - a^\tau kf(a,x^{\sigma^{-1}})l                 \cr
       &= x + af(a,x)l - ak^2f(a^{\sigma^{-1}},x^{\sigma^{-1}})l.\cr} $$
 Weil $\sigma \in \Sp(2n,K)$ ist, ist $f(a^{\sigma^{-1}},x^{\sigma^{-1}}) =
 f(a,x)$. Also ist
 $$x^{\sigma^{-1}\tau^{-1}\sigma\tau} = x + af(a,x)(1 - k^2)l. $$
 Da $a$ irgendein Vektor ist und weil $1 - k^2 \neq 0$ ist, folgt mit 3.7,
 dass jede Transvektion aus $\Sp(2n,K)$ ein Kommutator ist. Weil $\Sp(2n,K)$
 nach 3.8 von ihren Transvektionen erzeugt wird, ist $
 \Sp(2n,K)$ also in
 diesem Falle gleich ihrer Kommutatorgruppe.
 \par
       Als N\"achstes betrachten wir den Fall $|K| = 3$ und $n \geq 2$. Dann
 ist $\Rg_K(V) \geq 4$. Jede symplektische Basis hat also mindestens vier
 Vektoren. Es seien $v_1$, $v_2, v_3$, $v_4$ die ersten vier Vektoren
 einer solchen Basis. Setze $U := (\sum_{i:=1}^4 v_iK)^\pi$. Dann ist
 $V = (\sum_{i:=1}^4 v_iK) \oplus U$. Ferner ist auch $v_1$, $v_2 + v_3$,
 $v_3$, $v_1 + v_4$ eine symplektische Basis von $\sum_{i:=1}^4 v_iK$.
 Es gibt zwei Abbildungen $\sigma$, $\tau \in \GL(V)$ mit $v_1^\sigma = v_2$,
 $v_2^\sigma = -v_1$, $v_3^\sigma = v_3$, $v_4^\sigma = v_4$ und $u^\sigma = u$
 f\"ur alle $u \in U$, bzw. $v_1^\tau = v_2 + v_3$, $(v_2 + v_3)^\tau =
 -v_1$, $v_3^\tau = v_3$, $(v_1 + v_4)^\tau = v_1 + v_4$ und wiederum $u^\tau
 = u$ f\"ur alle $u \in U$. Mit 3.1 erschlie\ss en wir, dass $\sigma$ und
 $\tau$ in $\Sp(2n,K)$ liegen. Man rechnet ferner leicht nach, dass
 $\sigma^4 = \tau^4 = 1$ ist. Nun wird $\Sp(2n,K)$ von Transvektionen erzeugt,
 die wegen $|K| = 3$ alle die Ordnung 3 haben. Daher ist $\Sp(2n,K)/\Sp(2n,K)'$
 eine elementarabelsche 3-Gruppe, so dass wegen $\sigma^4 = \tau^4 = 1$
 die Abbildungen $\sigma$ und $\tau$ in $\Sp(2n,K)'$ liegen. Wir betrachten
 die Abbildung $(\tau\sigma)^2$. Es ist, wie eine einfache Rechnung zeigt,
 $v_i^{(\tau\sigma)^2} = v_i$ f\"ur $i := 1$, 2, 3. Ferner ist --- hier
 benutzen wir, dass wegen $|K| = 3$ die Gleichung $-v_3 = 2v_3$ gilt ---,
 $$\eqalign{
 v_4^{(\tau\sigma)^2} &= (v_1 + v_4 - v_1)^{(\tau\sigma)^2}
                       = (v_1 + v_4 - v_2 - v_3)^{\sigma\tau\sigma} \cr
                      &= (v_2 + v_4 + v_1 - v_3)^{\tau\sigma}
                       = (v_2 + v_3 + v_1 + v_4 + v_3)^{\tau\sigma} \cr
                      &= (-v_1 + v_1 + v_4 + v_3)^\sigma
                       = v_3 + v_4.                                 \cr} $$
 Schlie\ss lich ist $u^{(\tau\sigma)^2} = u$ f\"ur alle $u \in U$. Definiert
 man $\rho$ durch $x^\rho := x + v_3f(v_3,x)$, so ist $\rho$ eine
 Transvektion. Ferner sieht man, dass $\rho$ mit $(\tau\sigma)^2$ auf der
 Menge $\{v_1,v_2,v_3,v_4\} \cup U$ \"ubereinstimmt, so dass $\rho =
 (\tau\sigma)^2$ ist. Also ist $\rho \in \Sp(2n,K)'$.
 \par
       Es sei nun $\tau$ eine Transvektion ungleich 1 aus $\Sp(2n,K)$.
 Nach 3.7 gibt es ein $a \in V$ und ein $k \in K$ mit $x^\tau = x + af(a,x)k$
 f\"ur alle $x \in V$. Wegen $\tau \neq 1$ ist $a \neq 0$ und $k \neq 0$.
 Es gibt also ein $\lambda \in \Sp(2n,K)$ mit $a^\lambda = v_3$. Es folgt
 $x^{\lambda^{-1}\tau\lambda} = x + v_3f(a,x)k$. Wegen $k \neq 0$ und $|K| = 3$
 ist $k = 1$ oder $-1$. Dies besagt, dass $\tau$ zu $\rho$ oder $\rho^{-1}$
 konjugiert ist. Weil $\Sp(2n,K)'$ ein Normalteiler von $\Sp(2n,K)$ ist, folgt
 $\rho \in \Sp(2n,K)'$. Hieraus folgt wieder, dass $\Sp(2n,K) = \Sp(2n,K)'$ ist.
 \par
       Schlie\ss lich sei $|K| = 2$ und $n \geq 3$. In diesem Falle ist
 $\Rg_K(V) \geq 6$. Es gibt daher sechs Vektoren $v_1$, \dots, $v_6$, die
 den Anfang einer symplektischen Basis bilden. Wir setzen wieder
 $U := (\sum_{i:=1}^6 v_iK)^\pi$. Es folgt $V = (\sum_{i:=1}^6 v_iK) \oplus
 U$.   
 \par
       Wir definieren eine Abbildung $\rho \in \Sp(2n,K)$ verm\"oge $v_i^\rho
 := v_i$ f\"ur $i:= 1$, 2, 5, 6 und $v_3^\rho := v_3 + v_4$ und $v_4^\rho :=
 v_3$ sowie $u^\rho := u$ f\"ur alle $u \in U$. (Dass $\rho$ tats\"achlich
 in $\Sp(2n,K)$ liegt, erschlie\ss t man wieder mit Satz 3.1.)
 Dann ist $\rho^3 = 1$.
 Weil $\Sp(2n,K)$ von Transvektionen erzeugt wird, die wegen $|K| = 2$ die
 Ordnung 2 haben, ist $\Sp(2n,K)/\Sp(2n,K)'$ eine elementarabelsche 2-Gruppe.
 Folglich ist $\rho \in \Sp(2n,K)'$
 \par
       Wir definieren eine weitere Abbildung $\sigma \in \Sp(2n,K)$ verm\"oge
 $v_i^\sigma := v_i$ f\"ur $i := 1$, 5, 6 sowie $(v_2 + v_4)^\sigma = v_2 +
 v_4$, $(v_1 + v_3)^\sigma = v_4$ und $v_4^\sigma = v_1 + v_3 + v_4$ sowie
 $u^\sigma = u$ f\"ur alle $u \in U$. Dann ist $\sigma^3 = 1$, so dass
 auch $\sigma \in \Sp(2n,K)'$ gilt. Folglich ist $\mu_1 := \rho\sigma \in
 \Sp(2n,K)'$. Man rechnet leicht nach, dass
 $$ v_1^{\mu_1} = v_1,\ v_3^{\mu_1} = v_3,\ v_5^{\mu_1} = v_5 $$
 und
 $$v_2^{\mu_1} = v_1 + v_2 + v_3,\ v_4^{\mu_1} = v_1 + v_4,\ v_6^{\mu_1} = v_6
       $$
 sowie $u^{\mu_1} = u$ f\"ur alle $u \in U$ ist. Ersetzt man hierin der Reihe
 nach $v_1$, $v_2$, $v_3$, $v_4$, $v_5$, $v_6$ durch $v_1$, $v_2$, $v_5$,
 $v_6$, $v_3$, $v_4$, so sieht man, dass es auch eine Abbildung $\mu_2 \in
 \Sp(2n,K)'$ gibt mit
 $$ v_1^{\mu_2} = v_1,\ v_5^{\mu_2} = v_5,\ v_3^{\mu_2} = v_3 $$
 und
 $$v_2^{\mu_2} = v_1 + v_2 + v_5,\ v_6^{\mu_2} = v_1 + v_6,\ v_4^{\mu_2} = v_4
       $$
 sowie $u^{\mu_2} = u$ f\"ur alle $u \in U$ ist. Ersetzt man wiederum in der
 ersten Gruppe von Gleichungen $v_1$, $v_2$, $v_3$, $v_4$, $v_5$, $v_6$ der
 Reihe nach durch $v_1$, $v_2$, $v_3 + v_5$, $v_4$, $v_5$, $v_4 + v_6$, so
 folgt, dass es ein $\mu_3 \in \Sp(2n,K)'$ gibt mit
 $$ v_1^{\mu_3} = v_1,\ v_3^{\mu_3} = v_3,\ v_5^{\mu_3} = v_5 $$
 und
 $$v_2^{\mu_3} = v_1 + v_2 + v_3 + v_5,\ v_4^{\mu_3} = v_1 + v_4,\ 
     v_6^{\mu_3} = v_1 + v_6 $$
 sowie $u^{\mu_3} = u$ f\"ur alle $u \in U$ ist. Setzt man nun $\mu :=
 \mu_1\mu_2\mu_3$, so zeigen einfache Rechnungen, dass f\"ur $i \neq 2$ die
 Gleichung $v_i^\mu = v_i$ gilt, dass $v_2^\mu = v_1 + v_2$  und dass
 $u^\mu = u$ f\"ur alle $u \in U$ ist. Folglich ist $\mu$ eine in $\Sp(2n,K)'$
 liegende, von 1 verschiedene Transvektion. Weil die Charakteristik von $K$
 gleich 2 ist, sind alle Transvektionen aus $\Sp(2n,K)$ konjugiert, so dass
 auch in diesem Falle $\Sp(2n,K) = \Sp(2n,K)'$ gilt. Damit ist der Satz bewiesen.
 \medskip
       Die Gruppen $\SL(2,2)$ und $\SL(2,3)$ sind aufl\"osbar, so dass im
 Falle $|K| \leq 3$ und $n = 1$ tats\"achlich $\Sp(2,K)' \neq \Sp(2,K)$ ist.
 Dass auch $\Sp(4,2)' \neq \Sp(4,2)$ ist, werden wir am Ende dieses Abschnitts
 sehen, wo wir zeigen werden, dass $\Sp(4,2)$ und $S_6$ isomorph sind.
 \medskip\noindent
 {\bf 3.10. Satz.} {\it Es sei $V$ ein Vektorraum \"uber $K$ und $\pi$ sei eine
 symplektische Polarit\"at von $\La(V)$. Ist $G := \PSp(V)$ die zugeh\"orige
 projektive symplektische Gruppe, so gilt:
 \item{a)} $G$ ist auf der Menge der Punkte von $\La(V)$ transitiv.
 \item{b)} Ist $P$ ein Punkt von $\La(V)$, so hat $G_P$ genau drei Punktbahnen,
 n\"amlich: $\{P\}$, die Menge der von $P$ verschiedenen Punkte von $P^\pi$
 (diese Menge ist im Falle $\Rg_K(V) = 2$ leer)
 und die Menge der Punkte von $\La(V)_{P^\pi}$.
 \item{c)} $G$ operiert auf der Menge der Punkte von $\La(V)$ primitiv.}
 \smallskip
       Beweis. a) folgt aus der Transitivit\"at von $\Sp(V)$ auf der Menge der
 von Null verschiedenen Vektoren von $V$.
 \par
       b) Es sei $P = pK$. Ferner seien $uK$ und $vK$ zwei Punkte von $P^\pi$.
 Dann ist $f(p,u) = 0 = f(p,v)$. Aus 3.2 und 3.1 folgt die Existenz einer
 Abbildung $\rho \in \Sp(V)$ mit $p^\rho = p$ und $u^\rho = v$, so dass die
 von $P$ verschiedenen Punkte von $P^\pi$ eine Bahn von $G_P$ bilden.
 \par
       Es seien $uK$ und $vK$ zwei Punkte, die nicht in $P^\pi$ liegen. Dann
 ist $f(p,u)$, $f(p,v) \neq 0$. Wir d\"urfen daher, indem wir ggf. $u$ und
 $v$ durch skalare Vielfache ersetzen, annehmen, dass $f(p,u) = 1 = f(p,v)$
 ist. Wiederum mit 3.2 und 3.1 erschlie\ss en wir die Existenz eines $\sigma
 \in \Sp(V)$ mit $p^\sigma = p$ und $u^\sigma = v$, so dass auch die Punkte
 au\ss erhalb von $P^\pi$ eine Bahn von $G_P$ bilden.
 \par
       Angenommen $G$ sei imprimitiv. Ist $\Delta$ ein Imprimitivit\"atsgebiet
 von $G$, so ist also $|\Delta| \geq 2$, jedoch enth\"alt $\Delta$ nicht alle
 Punkte von $\La(V)$. Es sei $P \in \Delta$. Ist $H := G_P$, so ist
 $\Delta^H = \Delta$. Folglich zerf\"allt $\Delta$ in genau zwei Bahnen von
 $H$. Aus b) folgt somit, dass $\Delta$ entweder aus den Punkten von
 $P^\pi$ oder aus $P$ und den Punkten des affinen Raumes $\La(V)_{P^\pi}$
 besteht. Beides f\"uhrt aber auf einen Widerspruch. Nach a) ist $G$ ja
 auf den Punkten von $\La(V)$ transitiv. Weil nun jede Hyperebene von $V$
 von der Form $P^\pi$ ist, wobei $P$ ein Punkt ist, ist $G$ auch auf der
 Menge der Hyperebenen transitiv. Folglich w\"aren alle Hyperebenen
 Imprimitivit\"atsgebiete. Dies kann aber nicht sein, da Hyperebenen stets
 Punkte gemeinsam haben, Imprimitivit\"atsgebiete aber disjunkt liegen. Da es
 andererseits zu zwei verschiedene Hyperebenen stets einen Punkt gibt, der
 auf keiner der beiden Hyperebenen liegt, kann auch der zweite Fall nicht
 eintreten. Damit ist auch c) bewiesen.
 \medskip\noindent
 {\bf 3.11. Satz.} {\it Die Gruppe $\PSp(2n,K)$ ist einfach, es sei denn es ist
 $n = 1$ und $|K| \leq 3$ oder $n = 2$ und $|K| = 2$.}
 \smallskip
       Beweis. Setze $G := \PSp(2n,K)$. Es sei $P$ ein Punkt der zugrunde
 liegenden Geometrie. Die in $G$ liegenden Elationen mit dem Zentrum $P$
 bilden einen abelschen Normalteiler $A_P$ von $G_P$. Ist $Q$ ein weiterer
 Punkt, so folgt mit 3.10a), dass $A_P$ und $A_Q$ in $G$ konjugiert sind.
 Nach 3.8 wird $G$ von seinen Elationen erzeugt. Nach 3.10c) operiert $G$
 auf der Punktmenge von $\La(V)$ primitiv und nach 3.8 ist $G = G'$,
 falls nicht einer der drei im Satz erw\"ahnten Ausnahmef\"alle auftritt.
 Folglich ist $G$, von den Ausnahmef\"allen abgesehen, nach Satz III.2.1,
 dem Satz von Iwasawa, einfach.
 \medskip
 Die Gruppen $\PSp(2,2)$ und $\PSp(2,3)$ sind, wie wir schon bemerkten,
 auf\-l\"os\-bar. Dass auch $\PSp(4,2)$ nicht einfach ist, zeigt der folgende
 Satz.\index{Ausnahmeisomorphien}{}
 \medskip\noindent
 {\bf 3.12. Satz.} {\it Die Gruppen $\PSp(4,2)$ und $S_6$ sind isomorph.}
 \smallskip
       Beweis. Wir betrachten die Menge $M := \{1,2,3,4,5,6\}$ und bezeichnen
 mit $\Pi$ die Menge ihrer 2-Teilmengen. Ist $P \in \Pi$, so sei $b_P$ die
 Menge aus $P$ und den 2-Teilmengen von $M - P$. Schlie\ss lich sei
 $\Gamma := \{b_P \mid P \in \Pi\}$. Wir betrachten die Inzidenzstruktur
 $$ \Lambda := (\Pi,\Gamma,\in). $$
 Die n\"achste Aussage folgt unmittelbar aus der Definition.
 \smallskip
 \item{a)} Sind $P$, $Q \in \Pi$, so ist genau dann $P \in b_Q$, wenn $P = Q$
 oder $P \cap Q = \emptyset$ ist.
 \smallskip
       Aus a) folgt, dass die Anzahl $r_P$ der Bl\"ocke durch den Punkt $P$
 gleich $1 + {4 \choose 2} = 7$ und dass die Anzahl $k_b$ von Punkten auf
 einem Block ebenfalls gleich $1 + {4 \choose 2} = 7$ ist. Also ist
 $\Lambda$ eine taktische Konfiguration mit den Parametern $v = b = {6 \choose
 2} = 15$ und $r = k = 7$.
 \smallskip
 \item{b)} Es sei $P := \{a,b\}$ und $Q := \{a,c\}$ und $b \neq c$. Es sei
 ferner $X \in \Pi$. Genau dann ist $P$, $Q \in B_X$, wenn $P \cap X =
 \emptyset = Q \cap X$ ist.
 \smallskip
       Ist $P \cap X = Q \cap X = \emptyset$, so ist $P$, $Q \in b_X$ nach a).
 Es seien $P$, $Q \in b_X$. W\"are $P = X$, so folgte $a \in Q \cap X$. Mit
 a) folgt weiter $Q = X$ und damit $Q = P$ im Widerspruch zu $c \neq b$.
 \smallskip
 \item{c)} Sind $P$, $Q \in \Pi$ und ist $P \cap Q = \emptyset$, so sind
 $b_P$, $b_Q$ und $b_{P-(P \cup Q)}$ die einzigen Bl\"ocke, die gleichzeitig
 mit $P$ und $Q$ inzidieren.
 \smallskip
       Dies folgt unmittelbar aus a).
 \smallskip
       Aus b) und c) folgt, dass zwei verschiedene Punkte von $\Lambda$ mit
 genau drei Bl\"ocken inzidieren. $\Lambda$ ist also ein
 2-(15,7,3)-Blockplan.
 \par
       Es seien $P$ und $Q$ zwei verschiedene Punkte von $\Lambda$. Mit
 $P + Q$ bezeichnen wir ihre Verbindungsgerade, das ist der Schnitt \"uber
 die drei Bl\"ocke, die $P$ und $Q$ enthalten. Ist $P \cap Q \neq \emptyset$,
 so gilt nach b) genau dann $Y \in P + Q$, wenn $Y \cap X = \emptyset$ ist
 f\"ur alle 2-Teilmengen $X$ von $M - (P \cap Q)$, dh. genau dann, wenn $Y$
 eine 2-Teilmenge von $P \cup Q$ ist. In diesem Falle ist also $|P + Q| = 3$.
 Ist $P \cap Q = \emptyset$, so folgt aus c), dass $P$, $Q$ und
 $M - (P \cup Q)$ die einzigen Punkte auf $P + Q$ sind. Also gilt auch in
 diesem Falle $|P + Q| = 3$. Nun ist
 $$ {b - \lambda \over r - \lambda} = {15 - 3 \over 7 - 3} = 3, $$
 so dass $(\Pi,\Gamma,\in)$ nach IIII.10.1 die Geometrie aus den
 Punkten und Ebenen der projektiven Geometrie $L$ des Ranges 4 \"uber $\GF(2)$
 ist. Ferner folgt, dass $S_6$ eine Kollineationsgruppe von $L$ ist. Wir
 de\-fi\-nie\-ren nun $\pi$ durch $P^\pi := b_P$ und $b_P^\pi := P$. Weil genau
 dann $P \in b_Q$ gilt, wenn $P = Q$ oder $P \cap Q = \emptyset$ ist, und
 da dies mit $Q = P$ und $Q \cap P = \emptyset$ gleichbedeutend ist, gilt
 genau dann $P \in b_Q$, wenn $Q \in b_P$ gilt. Daher ist $\pi$ eine
 Polarit\"at, die wegen $P \in b_P = P^\pi$ symplektisch ist. Es folgt
 $S_6 \subseteq PC^*(\pi)$. Weil $L$
 eine projektiver Raum \"uber $\GF(2)$ ist, gilt $PC^*(\pi) = PC(\pi) =
 \PSp(4,2)$ nach 2.6. Also ist sogar $S_6 \subseteq \PSp(4,2)$. Aus 3.3 folgt
 schlie\ss lich, dass $|\PSp(4,2)| = 6! = |S_6|$ ist, so dass in der
 Tat $S_6 = \PSp(4,2)$ gilt.

\mysection{4. Polarit\"aten bei Charakteristik 2}

\noindent
       Die Zentralisatoren von Polarit\"aten operieren in der Regel
 irreduzibel auf den zugrunde liegenden Vektorr\"aumen. Ausnahmen treten
 nur bei Charakteristik 2 auf. Einen der zwei m\"oglichen Ausnahmef\"alle
 wollen wir in diesem Abschnitt betrachten. Der hier nicht behandelte Fall
 ist der einer unit\"aren Polarit\"at $\pi$ von $\La(V)$.
 \par
       Im Folgenden sei $V$ ein Vektorraum \"uber dem K\"orper $K$ der
 Charakteristik 2 und es sei $\Rg_K(V) \geq 3$. Ferner sei $\pi$ eine
 projektive Polarit\"at von $\La(V)$ und $f$ sei eine $\pi$ darstellende
 Bilinearform. Wie wir wissen, ist $f$ symmetrisch. Sind $u$, $v \in V$, so
 ist
 $$
       f(u + v,u + v) = f(u,u) + 2f(u,v) + f(v,v) 
                      = f(u,u) + f(v,v).         
 $$
 Hieraus folgt, dass die Menge
 $$ H := \bigl\{v \mid v \in V, f(v,v) = 0\bigr\} $$
 ein Unterraum von $V$ ist. Wir nehmen weiterhin an, dass $\pi$ nicht
 symplektisch ist. Dann ist also $H \neq V$.
 \medskip\noindent
 {\bf 4.1. Satz.} {\it Es sei $V$ ein Vektorraum des Ranges mindestens $3$ \"uber
 dem kommutativen K\"orper $K$ der Charakteristik $2$. Ferner sei $\pi$ eine
 projektive Polarit\"at von $\La(V)$, die nicht symplektisch sei. Ist dann
 $f$ eine $\pi$ darstellende Bilinearform und setzt man
 $$ H := \bigl\{v \mid v \in V, f(v,v) = 0\bigr\}, $$
 so ist $1 \leq \Rg_K(V/H) \leq [K:K^2]$. Insbesondere ist $H$ eine
 Hyperebene, falls $K$ perfekt ist.}
 \smallskip
       Beweis. Weil $H$ von $V$ verschieden ist, ist der Rang von $V/H$
 mindestens 1. Die Abbildung $v \to f(v,v)$ ist eine semilineare Abbildung
 des $K$-Vektorraumes $V$ in den $K^2$-Vektorraum $K$. Daher ist
 $\Rg_K(V/H) \leq [K:K^2]$. Ist $K$ perfekt, so ist $K = K^2$, so dass in
 diesem Falle $\Rg_K(V/H) = 1$ ist. Damit ist alles bewiesen.
 \medskip
       Wir setzen $H_0 := H \cap H^\pi$. Ferner seien $H_1$, $H_2$, $H_3 \in
 \La(V)$ mit $H = H_0 \oplus H_1$, $H^\pi = H_0 \oplus H_2$ und
 $H_1^\pi \cap H_2^\pi = H_0 \oplus H_3$. Dann ist
 $$ H_0 \oplus H_1 \oplus H_2 = H + H^\pi $$
 und somit
 $$ H_0^\pi \cap (H_1 + H_2)^\pi = H_0. $$
 Wegen $H_0^\pi = H + H^\pi$ ist daher nach dem Modulargesetz
 $$\eqalign{
 H_0 &= \bigl(H_0 + (H_1 + H_2)\bigr) \cap (H_1 + H_2)^\pi    \cr
     &= H_0 + \bigl((H_1 + H_2) \cap (H_1 + H_2)^\pi\bigr)    \cr} $$
 Also ist
 $$ (H_1 + H_2) \cap (H_1 + H_2)^\pi \leq H_0. $$
 Hieraus folgt
 $$ (H_1 + H_2) \cap (H_1 + H_2)^\pi \leq H_0 \cap (H_1 + H_2) = \{0\}, $$ 
 so dass also
 $$ (H_1 + H_2) \cap (H_1 + H_2)^\pi = \{0\} $$ 
 ist. Nach 1.5 ist daher
 $$\eqalign{
  V &= H_1 \oplus H_2 \oplus (H_1 + H_2)^\pi    \cr
    &= H_1 \oplus H_2 \oplus H_0 \oplus H_3     \cr
    &= (H + H^\pi) \oplus H_3.                  \cr} $$
 Es ist $H_0 + H_2 + H_3 \leq H_1^\pi$, so dass $V = H_1 + H_1^\pi$ ist.
 Mit 1.5 erhalten wir daher $H_1 \cap H_1^\pi = \{0\}$. Folglich induziert
 $f$ in $H_1$ eine nicht ausgeartete, alternierende Bilinearform. Mit
 $\Sp(H_1)$ bezeichnen wir die zu der von $f$ in $H_1$ induzierten Form
 geh\"orige symplektische Gruppe.
 \medskip\noindent
 {\bf 4.2. Satz.} {\it Es sei $V$ ein Vektorraum des Ranges mindestens $3$ \"uber
 dem kommutativen K\"orper $K$ der Charakteristik $2$. Ferner sei $\pi$ eine
 projektive Polarit\"at von $\La(V)$, die nicht symplektisch sei. Ferner sei
 $f$ eine $\pi$ darstellende Bilinearform. Schlie\ss lich haben $H$, $H_0$,
 $H_1$, $H_2$ und $H_3$ die oben eingef\"uhrten Bedeutungen. Ist $\sigma \in
 C(\pi)$, so gibt es Abbildungen $\alpha \in \Hom_K(H_3,H_0)$, $\beta \in
 \Hom_K(H_3,H_1)$, $\gamma \in \Sp(H_1)$ und $\delta \in \Hom_K(H_1,H_0)$, so
 dass
 \item{(j)} $f(y,z^\alpha) + f(y^\alpha,z) + f(y^\beta,z^\beta) = 0$ f\"ur
 alle $y$, $z \in H_3$
 und
 \item{(ij)} $f(x^\gamma,h^\beta) + f(x^\delta,h) = 0$ f\"ur alle $x \in H_1$
 und alle $h \in H_3$
 sowie 
 \item{(iij)} $(h_0 + h_1 + h_2 + h_3)^\sigma = h_0 + h_1^\gamma + h_1^\delta
       + h_2 + h_3 + h_3^\alpha + h_3^\beta$
 f\"ur alle $h_i \in H_i$
 \par\noindent
 gilt.
 \par
       Ist umgekehrt
 $\alpha \in \Hom_K(H_3,H_0)$, $\beta \in \Hom_K(H_3,H_1)$, $\gamma \in \Sp(H_1)$
 und $\delta \in \Hom_K(H_1,H_0)$ und gelten {\rm (j)} und {\rm (ij)} f\"ur diese
 Abbildungen, so liegt die
 durch {\rm (iij)} erkl\"arte Abbildung $\sigma$ in $C(\pi)$. Die Zuordnung
 $\sigma \to (\alpha,\beta,\gamma,\delta)$ ist bijektiv.}
 \smallskip
       Beweis. Es sei $\sigma \in C(\pi)$. Dann ist offenbar $H^\sigma = H$
 und damit auch $H^{\pi\sigma} = H^\pi$ und $H_0^\sigma = H_0$. Es sei
 $y \in H_2 + H_3$. Es gibt dann Homomorphismen $\eta \in \Hom_K(H_2 + H_3,H_2
 + H_3)$ und $\zeta \in \Hom_K(H_2 + H_3,H_0 + H_1)$ mit $y^\sigma =
 y^\eta + y^\zeta$. Daher ist
 $$ f(y,y) = f(y^\sigma,y^\sigma) = f(y^\eta,y^\eta) + (f(y^\zeta,y^\zeta). $$
 Weil $H_0 + H_1 = H$ ist, ist $f(y^\zeta,y^\zeta) = 0$. Somit ist
 $f(y,y) = f(y^\eta,y^\eta)$, was wiederum $f(y^\eta + y,y^\eta + y) = 0$
 nach sich zieht. Daher ist
 $$ y^\eta + y \in H \cap (H_2 + H_3) = \{0\}, $$
 so dass f\"ur alle $y \in H_2 + H_3$ die Gleichung $y^\eta = y$ gilt.
 \par
       Es sei $y \in H_2$. Nach dem soeben Bewiesenen ist dann $y^\sigma =
 y + y^\zeta$ mit $y^\zeta \in H$. Wegen $H_2 \leq H^\pi = H^{\pi\sigma}$
 ist daher
 $$ y^\zeta = y^\sigma + y \in (H_2^\sigma + H_2) \cap H \leq H^\pi \cap H
              = H_0. $$
 Mit $u \in H_3$ folgt somit
 $$\eqalign{
 f(u,y) &= f(u^\sigma,y^\sigma) = f(u + u^\zeta,y + y^\zeta)          \cr
        &= f(u,y) + f(u,y^\zeta) + f(u^\zeta,y) + f(u^\zeta,y^\zeta). \cr} $$
 Nun ist $u^\zeta \in H$ und $y \in H_2 \leq H^\pi$, so dass $f(u^\zeta,
 y) = 0$ ist. Ebenso folgt $f(u^\zeta,y^\zeta) = 0$. Folglich ist
 $$ f(u,y) = f(u,y) + f(u,y^\zeta), $$
 so dass auch $f(u,y^\zeta) = 0$ ist f\"ur alle $u \in H_3$ und alle
 $y \in H_2$. Hieraus folgt
 $$ y^\zeta \in H_3^\pi \cap H_0 = H_3^\pi \cap H \cap H^\pi =
 (H_3 + H + H^\pi)^\pi = V^\pi = \{0\} $$
 f\"ur alle $y \in H_3$. Somit ist $y^\sigma = y$ f\"ur alle $y \in H_2$.
 \par
       Es sei $y \in H_0$ und $x \in H_3$. Dann ist
 $$ f(y,x) = f(y^\sigma,x^\sigma) = f(y^\sigma,x + x^\zeta) = f(y^\sigma,x)
 + f(y^\sigma,x^\zeta). $$
 Nun ist $y^\sigma \in H_0^\sigma = H_0$ und $x^\zeta \in H$, so dass
 $f(y^\sigma,x^\zeta) = 0$ ist. Hieraus folgt
 $$ f(y - y^\sigma,x) = 0 $$
 f\"ur alle $y \in H_0$ und alle $x \in H_3$. Folglich ist
 $$ y - y^\sigma \in H_3^\pi \cap H_0 = \{0\}. $$
 Also induziert $\sigma$ auch in $H_0$ die Identit\"at. Wegen $H^\pi =
 H_0 + H_2$ ist $y^\sigma = y$ f\"ur alle $y \in H^\pi$.
 \smallskip
       Bevor wir den Beweis zu Ende f\"uhren, machen wir noch die folgende
 Bemerkung: ist $H = \{0\}$, so ist $H^\pi = V$ und daher
 $C(\pi) = \{1\}$. Dieser Fall kann durchaus eintreten. Ist n\"amlich $K$ nicht
 perfekt und ist $[K:K^2] = n \geq 3$, ist ferner $k_1, \dots, k_n$ eine Basis
 von $K$ bez. $K^2$ und $V$ ein Vektorraum \"uber $K$ mit $3 \leq \Rg_K(V) = m
 < n$, ist schlie\ss lich $b_1, \dots, b_m$ eine Basis von $V$ \"uber
 $K$ und definiert man $f$ verm\"oge
 $$ f\biggl(\sum_{i:=1}^m b_ix_i,\sum_{i:=1}^m b_iy_i \biggr)
                        := \sum_{i:=1}^m k_ix_iy_i, $$
 so ist $f$ eine nicht entartete Bilinearform auf $V$, f\"ur die genau dann
 $f(x,x) = 0$ gilt, wenn $x = 0$ ist.
 \smallskip
       Es sei $x \in H_3$. Dann ist $x^\zeta \in H = H_0 \oplus H_1$. Es gibt
 daher ein $\alpha \in \Hom_K(H_3,H_0)$ und ein $\beta \in \Hom_K(H_3,H_1)$ mit
 $x^\zeta = x^\alpha + x^\beta$.
 \par
       Es sei $x \in H_1$. Dann ist $x^\sigma = x^\gamma + x^\delta$ mit
 $\gamma \in \Hom_K(H_1,H_1)$ und $\delta \in \Hom_K(H_1,H_0)$, da ja
 $$ x^\sigma \in H_1^\sigma \leq H^\sigma = H = H_0 \oplus H_1 $$
 ist. Sind $x$, $y \in H_1$, so ist also
 $$ f(x,y) = f(x^\sigma,y^\sigma) = f(x^\gamma + x^\delta,y^\gamma + y^\delta)
           = f(x^\gamma,y^\gamma). $$
 Ist $x^\gamma = 0$, so ist $f(x,y) = 0$ f\"ur alle $y \in H_1$. Weil $f$,
 wie wir gesehen haben, in $H_1$ eine nicht ausgeartete, alternierende
 Bilinearform induziert, folgt weiter, dass $x = 0$ ist. Folglich ist
 $\gamma$ injektiv und wegen der Endlichkeit von $\Rg_K(H_1)$ dann auch
 bijektiv. Hieraus und aus $f(x^\gamma,y^\gamma) = f(x,y)$ f\"ur alle $x$,
 $y \in H_1$ folgt schlie\ss lich $\gamma \in \Sp(H_1)$.
 \par
       Ist nun $h_i \in H_i$ f\"ur $i := 0$, 1, 2, 3, so ist also
 $$\eqalign{
 (h_0 + h_1 + h_2 + h_3)^\sigma &= h_0^\sigma + h_1^\sigma + h_2^\sigma +
                                                        h_3^\sigma \cr
       &= h_0 + h_1^\gamma + h_1^\delta + h_2 + h_3^\alpha + h_3^\beta. \cr} $$
 Damit ist (iij) nachgewiesen. Um (j) zu beweisen, seien $y$, $z \in H_3$.
 Dann ist
 $$\eqalign{
 f(y,z) &= f(y^\sigma,z^\sigma) = f(y + y^\alpha + y^\beta,z + z^\alpha +
                                   z^\beta)                          \cr
        &= f(y,z) + f(y,z^\alpha) + f(y,z^\beta) + f(y^\alpha,z)
                                          + f(y^\alpha,z^\alpha)     \cr
        &\ \ \ +\ f(y^\alpha,z^\beta) + f(y^\beta,z) + f(y^\beta,z^\alpha)
                          + f(y^\beta,z^\beta).                      \cr} $$
 Wegen $H_3 \leq H_1^\pi$ ist
 $$ f(y,z^\beta) = 0 = f(y^\beta,z). $$
 Ferner ist $H_0 \leq H_0^\pi$, $H_1^\pi$, so dass auch
 $$ f(y^\alpha,z^\alpha) = 0 = f(y^\alpha,z^\beta) = f(y^\beta,z^\alpha) $$
 ist. Aus obiger Gleichung f\"ur $f(y,z)$ folgt somit
 $$ f(y,z^\alpha) + f(y^\alpha,z) + f(y^\beta,z^\beta) = 0, $$
 dh. (j).
 \par
       Schlie\ss lich sei $x \in H_1$ und $h \in H_3$. Wegen $H_3 \leq H_1^\pi$
 ist dann $f(x,h) = 0 = f(x^\gamma,h)$. Daher ist
 $$\eqalign{
   0 &= f(x,h) = f(x^\sigma,h^\sigma) = f(x^\gamma + x^\delta,h + h^\alpha
                                                 + h^\beta) \cr
     &= f(x^\gamma,h) + f(x^\gamma,h^\alpha) + f(x^\gamma,h^\beta) 
    + f(x^\delta,h) + f(x^\delta,h^\alpha) + f(x^\delta,h^\beta). \cr}
 $$
 Wegen $h^\alpha \in H_0$ ist
 $f(x^\gamma,h^\alpha) = 0 = f(x^\delta,h^\alpha)$.
 Ferner ist $x^\delta \in H_0$ und $h^\beta \in H_1$, so dass auch
 $f(x^\delta,h^\beta) = 0$ ist. Daher gilt
 $$ 0 = f(x^\gamma,h^\beta) + f(x^\delta,h), $$
 dh. (ij).
 \par
       Es bleibe dem Leser \"uberlassen zu verifizieren, dass umgekehrt
 jedes Quadrupel $(\alpha,\beta,\gamma,\delta)$, welches (j) und (ij)
 erf\"ullt, verm\"oge (iij) ein $\sigma \in C(\pi)$ liefert und dass die
 Zuordnung $\sigma \to (\alpha,\beta,\gamma,\delta)$ injektiv ist.
 \medskip\noindent
 {\bf 4.3. Satz.} {\it Es sei $V$ ein Vektorraum des Ranges mindestens $3$ \"uber
 dem kommutativen K\"orper $K$ der Charakteristik $2$. Ferner sei $\pi$ eine
 projektive Polarit\"at von $\La(V)$, die nicht symplektisch sei. Ferner seien
 $\rho$, $\sigma \in C(\pi)$. Geh\"ort zu $\rho$ gem\"a\ss\ 4.2 das
 Quadrupel $(\alpha,\beta,\gamma,\delta)$  und zu $\sigma$ das Quadrupel
 $(\alpha',\beta',\gamma',\delta')$, so geh\"ort zu $\rho\sigma$ das
 Quadrupel}
 $$ (\alpha + \alpha' + \beta\delta',\beta' + \beta\gamma',\gamma\gamma',
       \gamma\delta' + \delta). $$
 \par
       Beweis. Es sei $(\alpha'',\beta'',\gamma'',\delta'')$ das zu
 $\rho\sigma$ geh\"orende Quadrupel. Ist $h_1 \in H_1$, so ist $h_1^\rho =
 h_1^\gamma + h_1^\delta$. Daher ist
 $$ h_1^{\rho\sigma} = h_1^{\gamma\sigma} + h_1^{\delta\sigma}
       = h_1^{\gamma\gamma'} + h_1^{\gamma\delta'} + h_1^\delta. $$
 Somit ist $\gamma'' = \gamma\gamma'$ und $\delta'' = \gamma\delta' + \delta$.
 Ist $h_3 \in H_3$, so ist
 $$
 h_3^{\rho\sigma} = h_3^\sigma + h_3^{\alpha\sigma} + h_3^{\beta\sigma}
        = h_3 + h_3^{\alpha'} + h_3^{\beta'} + h_3^\alpha
              + h_3^{\beta\gamma'} + h_3^{\beta\delta'}. $$
 Also ist $\alpha'' = \alpha + \alpha' + \beta\delta'$ und $\beta'' =
 \beta' + \beta\gamma'$. Damit ist alles bewiesen.
 \medskip\noindent
 {\bf 4.4. Satz.} {\it Es sei $V$ ein Vektorraum des Ranges mindestens $3$ \"uber
 dem kommutativen K\"orper $K$ der Charakteristik $2$. Ferner sei $\pi$ eine
 projektive Polarit\"at von $\La(V)$, die nicht symplektisch sei. Ferner sei
 $\rho \in C(\pi)$ und zu $\rho$ geh\"ore gem\"a\ss\ 4.2 das
 Quadrupel $(\alpha,\beta,\gamma,\delta)$. Genau dann liegt $\rho$ im Zentrum
 von $C(\pi)$, wenn $\beta = 0$, $\gamma = 1$ und $\delta = 0$ ist.}
 \smallskip
       Beweis. Es sei $\sigma \in C(\pi)$ und
 $(\alpha',\beta',\gamma',\delta')$ sei das zugeh\"orige Quadrupel. Nach
 4.3 ist genau dann $\rho\sigma = \sigma\rho$, wenn
 $$\eqalign{
  \alpha + \alpha'  + \beta\delta' &= \alpha' + \alpha + \beta'\delta, \cr
             \beta' + \beta\gamma' &= \beta + \beta'\gamma,            \cr
                     \gamma\gamma' &= \gamma'\gamma,                   \cr
            \gamma\delta' + \delta &= \gamma'\delta + \delta'          \cr} $$
 ist. Nach 4.3 entspricht jedem Quadrupel der Form $(0,0,\gamma',0)$ ein
 $\sigma \in C(\pi)$. Ist nun $\rho \in Z(C(\pi))$, so ist also
 $\beta\gamma' = \beta$, $\gamma\gamma' = \gamma'\gamma$ und $\delta = \gamma'
 \delta$ f\"ur alle $\gamma' \in \Sp(H_1)$. Hieraus folgt $\beta = 0$ und
 $\delta = 0$ sowie $\gamma \in Z(\Sp(\pi))$. Weil die Charakteristik von
 $K$ gleich 2 ist, ist $Z(\Sp(H_1)) = \{1\}$, so dass $\gamma = 1$ ist.
 \par
       Ist umgekehrt $(\alpha,\beta,\gamma,\delta) = (\alpha,0,1,0)$, so ist
 $\rho \in Z(\Sp(\pi))$.
 \medskip\noindent
 {\bf 4.5. Satz.} {\it Es sei $V$ ein Vektorraum des Ranges mindestens $3$ \"uber
 dem kommutativen K\"orper $K$ der Charakteristik $2$. Ferner sei $\pi$ eine
 projektive Polarit\"at von $\La(V)$, die nicht symplektisch sei. Dann
 ent\-h\"alt $C(\pi)$ einen Normalteiler $M$ mit $Z(C(\pi)) \subseteq M
 \subseteq C(\pi)$ mit $\Sp(H_1) \cong C(\pi)/M$ und $M' \subseteq Z(C(\pi))$.}
 \smallskip
       Beweis. Es sei $M$ die Menge aller $\rho \in C(\pi)$, f\"ur die das
 zugeh\"orige Quadrupel die Form $(\alpha,\beta,1,\delta)$ hat. Aus 4.3 folgt,
 dass $M$ ein Normalteiler von $C(\pi)$ ist. Ist ferner $G$ die Menge
 aller $\rho \in C(\pi)$, deren zugeh\"origes Quadrupel die Form $(0,0,\gamma,
 0)$ hat, so ist $G$ eine zu $\Sp(H_1)$ isomorphe Untergruppe von $C(\pi)$ mit
 $C(\pi) = GM$ und $G \cap M = \{1\}$. Also ist $C(\pi)/M$ zu $\Sp(H_1)$
 isomorph. Aus 4.4 folgt $Z(C(\pi)) \subseteq M$ und mit 4.3 folgt
 schlie\ss lich $M' \subseteq Z(C(\pi))$.
 \medskip\noindent
 {\bf 4.6. Satz.} {\it Die Voraussetzungen seien wie in 4.5. Ist $K$ perfekt, so
 gilt: 
 \item{a)} Ist $\Rg_K(V) = n$ ungerade, so ist $M = \{1\}$ und $C(\pi) \cong
 \Sp(n - 1,K)$.
 \item{b)} Ist $\Rg_K(V) = n$ gerade, so ist $Z(C(\pi))$ zur additiven Gruppe
 von $K$ und $M/Z(C(\pi))$ zu $\Hom_K(H_3,H_1)$ isomorph.\par}
 \smallskip
       Beweis. Weil $K$ perfekt ist, ist $H$ nach 4.1 eine Hyperebene.
 Folg\-lich ist $H^\pi$ ein Punkt. Wegen $H_0 \leq H^\pi$ folgt daher, dass
 $\Rg_K(H_0) \leq 1$ ist. Ferner ist $\Rg_K(H_1)$ gerade.
 \par
       a) Ist $\Rg_K(V)$ ungerade, so ist $\Rg_K(H)$ gerade. Wegen
 $$ \Rg_K(H) = \Rg_K(H_0) + \Rg_K(H_1) $$
 ist daher auch $\Rg_K(H_0)$ gerade, so dass
 $\Rg_K(H_0) = 0$ und damit $H_0 = \{0\}$ ist. Also ist $H = H_1$ und
 $H^\pi = H_2$. Hieraus folgt mit
 $$ H_3 \oplus H_0 = H_1^\pi \cap H_2^\pi = H^\pi \cap H = H_0 = \{0\}, $$
 dass auch $H_3 = \{0\}$ ist. Daher ist $\Hom_K(H_3,H_0) = \{0\}$ und 
 $\Hom_K(H_3,H_1) = \{0\}$. Wegen $H_0 = \{0\}$ ist auch
 $\Hom_K(H_1,H_0) = \{0\}$. Dies impliziert
 wiederum $M = \{1\}$. Wegen $H_1 = H$ ist daher
 $$ C(\pi) \cong \Sp(H_1) \cong \Sp(n - 1,K). $$
 \par
       b) Ist $\Rg_K(V)$ gerade, so ist $\Rg_K(H)$ ungerade. Wegen
 $$ \Rg_K(H_0) = \Rg_K(H) - \Rg_K(H_1) $$
 ist dann auch $\Rg_K(H_0)$ ungerade.
 Weil der Rang von $H_0$ h\"ochstens Eins ist, ist also
 $\Rg_K(H_0) = 1$. Insbesondere ist dann $\Rg_K(H_1) = n - 2$,
 was seinerseits $\Rg(H_1^\pi) = 2$ impliziert. Weil $\Rg_K(H_0) = 1 =
 \Rg_K(H^\pi)$ ist, ist $H^\pi = H_0$, so dass $H_2 = \{0\}$ ist. Dies
 besagt wiederum $H_2^\pi = V$. Also ist
 $$ H_1^\pi = H_1^\pi \cap H_2^\pi = H_0 \oplus H_3, $$
 so dass
 $$ \Rg_K(H_3) = \Rg_K(H_1^\pi) - \Rg_K(H_0) = 2 - 1 = 1 $$
 ist. Wegen
 $$ V = H_0 \oplus H_1 \oplus H_3 = H_1 \oplus H_1^\pi $$
 ist $H_1^\pi$ nicht isotrop. Es gibt daher Vektoren $p$ und $q$ mit $H_0 =
 pK$ und $H_3 = qK$ und $f(p,q) = 1$. Es sei nun $\alpha \in \Hom_K(H_3,H_0)$
 und $\beta \in \Hom_K(H_3,H_1)$. Sind $y$, $z \in H_3$, so sind $y$ und $z$
 linear abh\"angig. Dann sind aber auch $y^\beta$ und $z^\beta$ linear
 abh\"angig, so dass $f(y^\beta,z^\beta) = 0$ ist. Es sei $q^\alpha = pr$
 und $y = qs$ sowie $z = qt$. Dann ist
 $$ f(y,z^\alpha) + f(y^\alpha,z) + f(y^\beta,z^\beta) =
       f(q,p)srt + f(p,q)srt = 0, $$
 so dass die Bedingung (j) von 4.2 stets erf\"ullt ist.
 \par
       Es sei $h = qr \in H_3$. Genau dann ist
 $$0 = f(x\gamma,q^\beta)r + f(x\delta,q)r $$
 f\"ur alle $r \in K$, wenn
 $$ 0 = f(x^\gamma,q^\beta) + f(x^\delta,q) $$
 ist. Definiert man nun $\varphi$ durch $x^\delta = p\varphi(x)$, so ist
 $\varphi \in \Hom_K(H_1,K)$ und es gilt
 $$ 0 = f(x^\gamma,q^\beta) + f(p,q)\varphi(x) = f(x^\gamma,q^\beta)
                                +  \varphi(x), $$
 so dass $\varphi(x) = f(x^\gamma,q^\beta)$ ist. Dies zeigt, dass
 $$ x^\delta = pf(x^\gamma,q^\beta) $$
 ist.  
 \par
       Sind umgekehrt $\alpha$, $\beta$, $\gamma$ gegeben und definiert man
 $\delta$ durch
 $$ x^\delta := pf(x^\gamma,q^\beta), $$
 so erf\"ullt das Quadrupel $(\alpha,\beta,\gamma,\delta)$ die Bedingung
 (j) sowieso und auch die Bedingung (ij), wie man leicht nachrechnet. Hieraus
 folgt alles Weitere.

\mysection{5. Quadratische Formen}

\noindent
       Es sei $V$ ein Vektorraum endlichen Ranges \"uber dem kommutativen
 K\"orper $K$. Eine Abbildung $Q$ von $V$ in $K$ nennen wir, wie schon in
 Kapitel IIII, Abschnitt~4, {\it quadratische Form\/} auf $V$, falls $Q$ die
 beiden folgenden Bedingungen erf\"ullt:\index{quadratische Form}{}
 \item{1)} Es ist $Q(xk) = Q(x)k^2$ f\"ur alle $x \in V$ und alle $k\in K$.
 \item{2)} Die durch $f(x,y) := Q(x + y) - Q(x) - Q(y)$ erkl\"arte
 Abbildung $f$ von $V \times V$ in $K$ ist eine Bilinearform.
 \par\noindent
       Es ist
 $$ f(x,x) = Q(2x) - 2Q(x) = 4Q(x) - 2Q(x) = 2Q(x), $$
 so dass $f$ alternierend ist, falls $\Char(K) = 2$ ist. Ferner
 ist klar, dass $f(x,y) = f(y,x)$ ist, so dass $f$ also in jedem Fall
 symmetrisch ist. Man nennt $f$ die zu $Q$ {\it assoziierte
 Bilinearform\/}.\index{assoziierte Bilinearform}{} Ist $\Char(K) \neq 2$, so ist
 $Q(x) = {1 \over 2}f(x,x)$. Die quadratische Form $Q$ hei\ss t {\it nicht
 ausgeartet\/},\index{nicht ausgeartet}{} falls die zu $Q$
 assoziierte Bilinearform $f$ nicht
 ausgeartet ist. Ist $\Char(K) = 2$, so ist $f$ alternierend, so dass $Q$
 h\"ochstens dann nicht ausgeartet ist, wenn $\Rg_K(V)$ gerade ist.
 \par
       Ist $Q$ eine quadratische Form auf $V$ und ist $f$ die zu $Q$
 assoziierte Bilinearform, so hei\ss t die Menge
 $$ \Kern(Q) := \{x \mid x \in V, Q(x) = 0\ \hbox{\rm und}\ f(x,y) = 0\ 
                     \hbox{f\"ur alle}\ y \in V\} $$
 {\it Kern\/} von $Q$. Offenbar ist $\Kern(Q)$ ein Unterraum von $V$. Ist
 n\"amlich $x \in \Kern(Q)$ und $k \in K$, so ist
 $$ Q(xk) = Q(x)k^2 = 0 $$
 und
 $$ f(xk,y) = f(x,y)k = 0 $$
 f\"ur alle $y \in V$. Also ist $xk \in \Kern(Q)$. Sind $x$, $x' \in \Kern(Q)$,
 so ist
 $$ Q(x + x') = f(x,x') + Q(x) + Q(x') = 0 $$
 und
 $$ f(x + x',y) = f(x,y) + f(x',y) = 0 $$
 f\"ur alle $y \in V$. Daher ist auch $x + x' \in \Kern(Q)$.
 \par
       Die quadratische Form $Q$ hei\ss t
 {\it singul\"ar\/},\index{singul\"are quadratische Form}{} falls $\Kern(Q) \neq \{0\}$ ist.
 \medskip\noindent
 {\bf 5.1. Satz.} {\it Es sei $V$ ein Vektorraum endlichen Ranges \"uber dem
 kommutativen K\"orper $K$ und $Q$ sei eine quadratische Form auf $V$.
 \item{a)} Ist $Q$ singul\"ar, so ist $Q$ ausgeartet.
 \item{b)} Ist $\Char(K) \neq 2$ und ist $Q$ ausgeartet, so ist $Q$ singul\"ar.}
 \smallskip
       Beweis. a) Es sei $Q$ singul\"ar. Es gibt dann einen von 0
 verschiedenen Vektor $x \in \Kern(Q)$. F\"ur diesen Vektor gilt $f(x,y) = 0$
 f\"ur alle $y \in V$, so dass $Q$ ausgeartet ist.
 \par
       b) $Q$ sei ausgeartet. Es gibt dann einen Vektor $x \in V$ mit
 $x \neq 0$ und $f(x,y) = 0$ f\"ur alle $y \in V$. Insbesondere ist dann auch
 $$ 0 = f(x,x) = 2Q(x). $$
 Weil $\Char(K) \neq 2$ ist, folgt $Q(x) = 0$, so dass $x \in \Kern(Q)$
 ist. Folglich ist $Q$ singul\"ar.
 \medskip
       Dass es im Falle der Charakteristik 2 ausgeartete quadratische Formen
 gibt, die nicht singul\"ar sind, zeigt folgendes Beispiel. Es
 sei $K$ ein K\"orper der Charakteristik 2 und $V$ sei ein Vektorraum des
 Ranges 3 \"uber $K$ und $b_1$, $b_2$, $b_3$ sei eine Basis von $V$.
 Ist $v = b_1k_1 + b_2k_2 + b_3k_3$, so setzen wir
 $$ Q(v) := k_1^2 + k_2k_3 + k_3^2. $$
 Ist dann auch noch $w = b_1l_1 + b_2l_2 + b_3l_3$, so ist --- hier wird
 benutzt, dass $\Char(K) = 2$ ist ---,
 $$ f(v,w) = l_2k_3 + k_2l_3. $$
 Hieraus folgt, dass $f(v,w) = 0$ f\"ur alle $w \in V$
 genau dann gilt, wenn $v \in b_1K$ ist.
 Also ist $\Kern(Q) \leq b_1K$. Nun ist $Q(b_1) = 1$, so dass
 $\Kern(Q) = \{0\}$ ist.
 \medskip
       Ist $Q$ eine quadratische Form auf dem Vektorraum $V$ und ist $f$ die
 zu $Q$ assoziierte Bilinearform, so definieren wir das
 {\it Radikal\/}\index{Radikal}{}
 $\rad(Q)$ von $Q$ durch
 $$ \rad(Q) := \bigl\{x \mid x \in V, f(x,y) = 0\ \hbox{f\"ur alle}\
	y \in V\bigr\}. $$
 Offenbar ist $\rad(Q)$ ein Unterraum.
 \medskip\noindent
 {\bf 5.2. Satz.} {\it Es sei $K$ ein perfekter K\"orper der Charakteristik $2$.
 Ferner sei $V$ ein Vektorraum \"uber $K$ und $Q$ sei eine quadratische Form
 auf $V$. Dann ist $Q$ sicher dann singul\"ar, wenn $\Rg_K(\rad(Q)) \geq 2$
 ist.}
 \smallskip
       Beweis. Es seien $a$ und $b$ zwei linear unabh\"angige Vektoren aus
 $\rad(Q)$. Ist $Q(a) = 0$ oder $Q(b) = 0$, so ist $Q$ singul\"ar. Wir
 k\"onnen daher annehmen, dass $Q(a)$, $Q(b) \neq 0$ ist. Es gibt dann, da
 $K$ perfekt ist, ein $k \in K$ mit
 $$ k^2 = {Q(a) \over Q(b)}. $$
 Weil $a$ und $b$ linear unabh\"angig sind, ist $a + bk \neq 0$. Ferner
 ist $a + bk \in \rad(Q)$, da $\rad(Q)$ ja ein Unterraum von $V$ ist.
 Schlie\ss lich gilt
 $$ Q(a + bk) = f(a,b)k + Q(a) + Q(b)k^2 = 0, $$
 so dass $0 \neq a + bk \in \Kern(Q)$ ist, q. e. d.
 \medskip\noindent
 {\bf 5.3. Satz.} {\it Es sei $K$ ein perfekter K\"orper der Charakteristik $2$
 und $V$ sei ein Vektorraum geraden Ranges \"uber $K$. Ist $Q$ eine
 qua\-dra\-ti\-sche
 Form auf $V$, so ist $Q$ genau dann nicht singul\"ar, wenn $Q$ nicht
 entartet ist.}
 \smallskip
       Beweis. Ist $Q$ nicht entartet, so ist $Q$ nach 5.1 a) nicht
 singul\"ar.
 \par
       Ist $Q$ entartet, so ist $\rad(Q) \neq \{0\}$. Weil die zu $Q$
 assoziierte Bi\-li\-ne\-ar\-form $f$ alternierend ist, induziert $f$ in
 $V/\rad(Q)$
 eine alternierende Bilinearform, die \"uberdies nicht ausgeartet ist. Nach
 1.6 ist somit $\Rg_K(V/\rad(Q))$ gerade. Weil $\Rg_K(V)$ gerade ist, ist
 auch $\Rg_K(\rad(Q))$ gerade. Wegen $\rad(Q) \neq \{0\}$ ist somit
 $\Rg_K(\rad(Q)) \geq 2$, so dass $Q$ nach 5.2 singul\"ar ist. Damit ist
 alles bewiesen.
 \medskip
       Es sei $Q$ eine quadratische Form auf $V$ und $U$ sei der Kern von
 $Q$. Ist $v \in V$ und $u \in U$, so ist
 $$ Q(v + u) = f(v,u) + Q(v) + Q(u) = Q(v). $$
 Definiert man $Q'$ auf $V/U$ durch $Q'(v + U) := Q(v)$, so ist also $Q'$ eine
 Abbildung von $V/U$ in $K$. Ferner ist
 $$ Q'(vk + U) = Q(vk) = Q(v)k^2 = Q'(v + U)k^2. $$
 Definiert man $f'$ durch
 $$ f'(v + U,w + U) := Q'(v + w + U) - Q'(v + U) - Q'(w + U), $$
 so folgt
 $$ f'(v + U,w + U) = f(v,w). $$
 Dies impliziert, dass auch $f'$ eine Bilinearform ist. Folglich ist $Q'$
 eine quadratische Form auf $V/U$. Ist $f'(v + U,w + U) = 0$ f\"ur alle
 $w + U \in V/U$ und ist $Q'(v + U) = 0$, so ist $f(v,w) = 0$ f\"ur alle
 $w \in V$ und $Q(v) = 0$. Also ist $v \in U$. Somit ist $Q'$ nicht singul\"ar.
 Hieraus folgt, dass man eine \"Ubersicht \"uber alle quadratischen
 Formen hat, wenn man alle nicht singul\"aren quadratischen Formen kennt.
 \par
       Ein Vektor $v \in V$ hei\ss t
 {\it singul\"ar\/},\index{singul\"arer Vektor}{} falls $v \neq 0$ und
 $Q(v) = 0$ ist. Ist $v$ singul\"ar, so ist $vk$ singul\"ar f\"ur alle
 $k \in K^*$. Ferner ist $f(w,w') = 0$ f\"ur alle $w$, $w' \in vK$. Ist $U \in
 \La(V)$ und ist $Q(u) = 0$ und $f(u,w) = 0$ f\"ur alle $u$, $w \in U$, so
 hei\ss t $U$
 {\it vollst\"andig singul\"ar\/}.\index{vollst\"andig singul\"arer Raum}{}
 Die von singul\"aren Vektoren aufgespannten Punkte sind also vollst\"andig singul\"ar.
 \par
       Die Menge der zu $Q$ geh\"orenden vollst\"andig singul\"aren Punkte
 nennen wir die zu $Q$ {\it geh\"orende\/}, bzw. die {\it durch Q definierte
 Quadrik\/}.\index{Quadrik}{} Ist $Q$ nicht singul\"ar, so nennen wir
 die zugeh\"orige Quadrik, auch wenn sie leer sein sollte, {\it nicht
 ausgeartet\/}.\index{nicht ausgeartet}{} Ist
 $Q$ singul\"ar, so nennen wir die zugeh\"orige Quadrik, die dann niemals
 leer ist, auch {\it Kegel\/}.\index{Kegel}{} Ist $U$ der Kern von $Q$, so ist
 der durch $Q$
 definierte Kegel offensichtlich die Projektion der zu $Q'$ geh\"orenden
 Quadrik in $\La(V/U)$ aus $U$. Dies rechtfertigt den Namen Kegel.
 \par
       Sind $\Sigma$ und $\Sigma'$ Quadriken in $\La(V)$, so hei\ss en
 $\Sigma$ und $\Sigma'$ {\it projektiv
 \"aquivalent\/},\index{projektiv \"aquivalent}{} falls es eine
 projektive Kollineation von $\La(V)$ gibt, welche $\Sigma$ auf $\Sigma'$
 abbildet. Sind $Q$ und $Q'$ quadratische Formen, so hei\ss en sie
 {\it projektiv
 \"aquivalent\/},\index{projektiv \"aquivalent}{}
 falls es ein $k \in K^*$ und ein $\sigma
 \in \GL(V)$ gibt mit $Q(x) = kQ'(x^\sigma)$ f\"ur alle $x \in V$. Sind
 $\Sigma$ und $\Sigma'$ projektiv \"aquivalente Quadriken und wird $\Sigma$
 durch die quadratische Form $Q$ dargestellt, so definieren wir die
 quadratische Form $Q'$ durch $Q'(x) := Q(x^{\gamma^{-1}})$ f\"ur alle
 $x \in V$, wobei $\gamma \in \GL(V)$ so gew\"ahlt sei, dass $\Sigma^\gamma
 = \Sigma'$ ist. Nun ist genau dann $Q'(x) = 0$, wenn $Q(x^{\sigma^{-1}}) = 0$
 ist. Hieraus folgt, dass $\Sigma'$ durch $Q'$ dargestellt wird. Damit
 ist gezeigt, dass projektiv \"aquivalente Quadriken durch projektiv
 \"aquivalente quadratische Formen dargestellt werden. Die Umkehrung, dass
 projektiv \"aquivalente quadratische Formen projektive \"aquivalente
 Quadriken darstellen ist ebenso einfach zu zeigen.
 \medskip\noindent
 {\bf 5.4. Satz.} {\it Es sei $Q$ eine quadratische Form auf dem Vektorraum $V$.
 Ist $U$ ein Teilraum von $V$, dessen Punkte alle auf der durch $Q$
 definierten Quadrik liegt, so ist\/ $U$ vollst\"andig singul\"ar.}
 \smallskip
       Beweis. Es seien $x$, $y \in U$. Dann ist $Q(x) = Q(y) = Q(x + y) = 0$.
 Daher ist
 $$ f(x,y) = Q(x + y) - Q(x) - Q(y) = 0. $$
 Damit ist bereits alles bewiesen.

\mysection{6. Die wittsche Zerlegung}

\noindent
       Es sei $f$ eine symmetrische $\alpha$-Semibilinearform auf $V$. Die
 Form $f$ hei\ss t {\it spurwertig\/},\index{spurwertig}{} falls es zu jedem
 $x \in V$ ein $k \in K$ gibt mit $f(x,x) = k + k^\alpha$.
 \par
       Ist die Charakteristik von
 $K$ ungerade, so gilt wegen $f(x,x)^\alpha = f(x,x)$ die Gleichung
 $$ f(x,x) = {1 \over 2}f(x,x) + \left({1 \over 2}f(x,x)\right)^\alpha, $$
 so dass $f$ spurwertig ist.
 \par
       Ist $f$ alternierend, so ist
 $ f(x,x) = 0 + 0^\alpha$,
 so dass $f$ auch in diesem Falle spurwertig ist.
 \par
       Ist $\Char(K) = 2$ und induziert $\alpha$ auf $Z(K)$ nicht die
 Identit\"at, so gibt es nach III.6.3b) ein $z \in Z(K)$ mit $z + z^\alpha = 1$.
 Dann ist
 $$ f(x,x) = f(x,x)(z + z^\alpha) = f(x,x)z + (f(x,x)z)^\alpha, $$
 so dass $f$ auch in diesem Falle spurwertig ist. Dieser Fall liegt
 insbesondere dann vor, wenn $K$ kommutativ und $\alpha$ ein involutorischer
 Automorphismus von $K$ ist.
 \par
       Ist $f$ die zu einer quadratischen Form assoziierte Bilinearform, so ist
 $f$ ebenfalls spurwertig. Dies ist f\"ur $\Char(K) \neq 2$ nach obigem
 selbstverst\"andlich und folgt f\"ur $\Char(K) = 2$ daraus, dass $f$ in diesem
 Falle alternierend ist.
 \par
       Dass es auch nicht spurwertige symmetrische Bilinearformen gibt, zeigen
 die Beispiele von projektiven Polarit\"aten, die wir in
 Abschnitt 4 betrachtet haben.
 \par
       Wir setzen im Folgenden stets voraus, dass $f$ eine spurwertige,
 symmetrische $\alpha$-Semibilinearform auf $V$ ist, wobei wir statt
 $\alpha$-Se\-mi\-bi\-li\-near\-form auch kurz $\alpha$-Form sagen werden. Nach
 obiger Bemerkung kann $f$ also auch die zu einer quadratischen Form assoziierte
 Bilinearform sein. Sind $x$ und $y$ zwei Vektoren aus $V$ und ist $f(x,y) = 0$,
 so nennen wir $x$ und $y$ {\it orthogonal\/}.\index{orthogonale Vektoren}{}
 Ist $U \in \La(V)$, so bezeichnen
 wir mit $U^\bot$ die Menge aller zu allen Vektoren aus $U$ orthogonalen
 Vektoren. Zwei Unterr\"aume $U$ und $W$ hei\ss en {\it orthogonal\/}, falls
 $W \leq U^\bot$ ist.  Weil $f$ symmetrisch ist, ist auch die
 Orthogonalit\"atsrelation symmetrisch.
 \par
       Ist $U \cap U^\bot \neq \{0\}$, so hei\ss t $U$, wie schon zuvor,
 {\it isotrop\/}.\index{isotroper Raum}{} R\"uhrt $f$ von einer quadratischen
 Form $Q$ her, so hei\ss t $U$ {\it singul\"ar\/},\index{singul\"arer Raum}{}
 falls $U$ isotrop ist und au\ss erdem $Q(u) = 0$ ist
 f\"ur alle $u \in U \cap U^\bot$. Wir nennen $U$ {\it vollst\"andig
 singul\"ar\/},\index{vollst\"andig singul\"arer Raum}{}
 falls $U$ vollst\"andig isotrop und singul\"ar ist. Diesen Begriff hatten wir
 schon fr\"uher definiert. Man beachte, dass nicht isotrope R\"aume
 nicht singul\"ar sind.
 \medskip\noindent
 {\bf 6.1. Satz.} {\it Es sei $V$ ein Vektorraum \"uber $K$ und $f$ sei entweder
 eine spurwertige $\alpha$-Form oder die zu der quadratischen Form $Q$ auf
 $V$ assoziierte Bilinearform. Ferner sei $U$ ein von $\{0\}$ verschiedener,
 vollst\"andig isotroper, bzw. vollst\"andig singul\"arer Unterraum von $V$.
 Ist $x \in V$ und $x \not\in U^\bot$, ist ferner $k \in K$, so gibt es
 ein $y \in U$ mit $f(x + y,x + y) = k + k^\alpha$, bzw. mit $Q(x + y) = k$.}
 \smallskip
       Beweis. Wir betrachten zuerst den Fall, dass $f$ eine $\alpha$-Form
 ist. Da $f$ spurwertig ist, gibt es ein $l \in K$ mit $f(x,x) = l + l^\alpha$.
 Ist $y \in U$, so ist $f(y,y) = 0$, da ja $U \leq U^\bot$ ist. Daher ist
 $$\eqalign{
 f(x + y,x + y) &= f(x,x) + f(x,y) + f(y,x) + f(y,y)    \cr
                &= l + l^\alpha + f(x,y) + f(x,y)^\alpha  \cr
                &= l + f(x,y) + (l + f(x,y))^\alpha .    \cr} $$
 Wegen $x \not\in U^\bot$ ist die Abbildung $y \to f(x,y)$ eine von Null
 verschiedene lineare Abbildung von $U$ in und damit auf $K$. Es gibt also
 ein $y \in U$ mit $f(x,y) = k - l$, womit der Satz im vorliegenden Falle
 bewiesen ist.
 \par
       Es sei nun $f$ die zu der quadratischen Form $Q$ assoziierte
 Bi\-li\-ne\-ar\-form. Wir setzen $l := Q(x)$. Ist $y \in U$, so ist $Q(y) = 0$,
 da $U$ ja vollst\"andig singul\"ar ist. Es folgt also, dass
 $$ Q(x + y) = Q(x) + Q(y) + f(x,y) = l + f(x,y) $$
 ist. Es ist wieder $f(x,y) \neq 0$, so dass wir wie im ersten Fall die
 Existenz eines $y \in U$ erschlie\ss en, welches $f(x,y) = k - l$ erf\"ullt.
 Dann ist $Q(x + y) = k$, so dass die Behauptung auch im zweiten Falle
 bewiesen ist.
 \medskip\noindent
 {\bf 6.2. Satz.} {\it Es sei $f$ eine spurwertige, nicht entartete, symmetrische
 $\alpha$-Form auf dem $K$-Vektorraum $V$, bzw. die zu einer nicht entarteten
 quadratischen Form $Q$ assoziierte Bilinearform. Ferner sei $X$ ein
 voll\-st\"an\-dig isotroper, bzw. vollst\"andig singul\"arer Unterraum von $V$.
 \item{a)} Ist $Y \in \La(V)$ ein vollst\"andig isotroper, bzw. vollst\"andig
 singul\"arer Unterraum mit $\Rg_K(Y) = \Rg_K(X) = r$ und $Y \cap X^\bot =
 \{0\}$, so ist $X + Y$ nicht isotrop und zu jeder
 Basis $b_1$, \dots, $b_r$ von $X$ gibt es eine Basis $c_1$, \dots, $c_r$
 von $Y$ mit $f(b_i,c_j) = \delta_{ij}$ f\"ur $i$, $j := 1$, \dots, $r$.
 \item{b)} Ist $U$ ein vollst\"andig isotroper bzw. vollst\"andig singul\"arer
 Unterraum mit $\Rg_K(U) \leq \Rg_K(X)$ und $U \cap X^\bot = \{0\}$, so gibt es
 einen vollst\"andig
 isotropen bzw. einen vollst\"andig singul\"aren Unterraum $Y$ mit
 $U \leq Y$, $\Rg_K(Y) = \Rg_K(X)$ und $Y \cap X^\bot = \{0\}$.\par}
 \smallskip
       Beweis. a) Weil $f$ nicht entartet ist, ist die Abbildung $\bot$
 eine Polarit\"at von $\La(V)$. Daher ist
 $$(X + Y) \cap (X + Y)^\bot = (X + Y) \cap X^\bot \cap Y^\bot. $$
 Wegen $X \leq X^\bot$ gilt
 daher auf Grund des Modulargesetzes, und weil $X \cap Y^\bot = \{0\}$
 vorausgesetzt ist,
 $$ (X + Y) \cap X^\bot \cap Y^\bot = \bigl(X + (Y \cap X^\bot)\bigr) \cap Y^\bot
                                    = X \cap Y^\bot. $$
 Nun ist $X^\bot \cap Y = \{0\}$ und daher $X + Y^\bot = V$. Weil $f$ nicht
 ent\-ar\-tet ist, ist $\Rg_K(Y) + \Rg_K(Y^\bot) =\Rg_K(V)$ und wegen
 $\Rg_K(X) = \Rg_K(Y)$ ist $\Rg_K(X) + \Rg_K(Y^\bot) = \Rg_K(V)$. Hieraus
 folgt, dass $V = X \oplus Y^\bot$ ist. Dies besagt wiederum $X \cap
 Y^\bot = \{0\}$. Also ist $(X + Y) \cap (X + Y)^\bot = \{0\}$, so dass
 $X + Y$ nicht isotrop, bzw. nicht singul\"ar ist.
 \par
       Wir definieren f\"ur $i := 1$, \dots, $r$ die Abbildung $\varphi_i
 \in \Hom_K(Y,K)$ durch
 $$ \varphi_i(y) := f(b_i,y) $$
 f\"ur alle $y \in Y$. Ferner setzen wir $H_i := \Kern(\varphi_i)$. Dann
 ist
 $$ \bigcap_{i:=1}^r H_i \leq X^\bot \cap Y = \{0\}. $$
 Daher ist
 $$ P_j := \bigcap_{i:=1,\>i\neq j}^r H_i $$
 ein Punkt, da ja $\Rg_K(Y) = r$ ist. Ist $y \in P_j$, so ist $f(b_i,y) = 0$
 f\"ur $i \neq j$ und $f(b_j,y) \neq 0$. Es gibt daher auch ein $c_j \in P_j$
 mit $f(b_i,c_j) = \delta_{ij}$. Offenbar ist $c_1$, \dots, $c_r$ eine
 Basis von $Y$. Damit ist a) bewiesen.
 \par
       b) Ist $\Rg_K(U) = r$, so ist nichts zu beweisen. Es sei also
 $$ \Rg_K(U) < r. $$
 Wegen $U \cap X^\bot = \{0\}$ ist $V = U^\bot + X$. Daher ist
 $$\eqalign{
   \Rg_K(V) &= \Rg_K(X) + \Rg_K(U^\bot) - \Rg_K(X \cap U^\bot) \cr
            &= \Rg_K(X) + \Rg_K(V) - \Rg_K(U) - \Rg_K(X \cap U^\bot) \cr} $$
 Hieraus folgt
 $$ \Rg_K(X) = \Rg_K(U) + \Rg_K(X \cap U^\bot), $$
 so dass wegen $\Rg_K(U) < \Rg_K(X)$ die Ungleichung $X \cap U^\bot \neq
 \{0\}$ gilt. W\"are $U^\bot \leq U + X^\bot$, so w\"are $U \geq U^\bot
 \cap X$ und folglich w\"are
 $$\{0\} = U \cap X^\bot \geq U^\bot \cap X \cap X^\bot = U^\bot \cap X \neq
       \{0\}. $$
 Dieser Widerspruch zeigt, dass es ein $y \in U^\bot$ gibt mit $y \not\in
 U + X^\bot$. Ist $x \in U^\bot \cap X$, so ist $y + x \in U^\bot$ und wegen
 $X \leq X^\bot$ auch $y + x \not\in U + X^\bot$. Weil $U^\bot \cap X \neq
 \{0\}$ und au\ss erdem vollst\"andig isotrop bzw. vollst\"andig singul\"ar
 ist, gibt es nach 6.1 ein $x \in U^\bot \cap X$, so dass $y + x$ isotrop,
 bzw. singul\"ar ist. Wegen $y + x \not\in U + X^\bot$ ist $y + x \neq 0$.
 Setze
 $$ U' := U \oplus (y + x)K. $$
 Dann ist $\Rg_K(U') = \Rg_K(U) + 1$. Ferner ist $U'$ vollst\"andig
 isotrop, bzw. vollst\"andig singul\"ar. Wegen $y + x \not\in U + X^\bot$ ist
 $$\eqalign{
 \Rg_K(U' + X^\bot) &= \Rg_K(U + X^\bot) + 1     \cr
       &= \Rg_K(U) + \Rg_K(X^\bot) + 1           \cr
       &= \Rg_K(U') + \Rg_K(X^\bot),             \cr} $$
 so dass $\Rg_K(U' \cap X^\bot) = 0$ ist. Folglich ist $U' \cap X^\bot =
 \{0\}$. Vollst\"andige Induktion liefert nun die Behauptung unter b).
 \medskip\noindent
 {\bf 6.3. Korollar.} {\it Ist $X$ ein vollst\"andig isotroper bzw. vollst\"andig
 singul\"arer Unterraum von $V$, so gibt es einen vollst\"andig isotropen
 bzw. vollst\"andig singul\"aren Unterraum $Y$ von $V$ mit $\Rg_K(Y) =
 \Rg_K(X)$ und $X \cap Y = \{0\}$, so dass $X + Y$ nicht isotrop bzw. nicht
 singul\"ar ist.}
 \smallskip
       Dies folgt mit $U := \{0\}$ aus 6.2.
 \medskip
 Aus 6.3 folgt wiederum, dass $2\Rg_K(X) \leq \Rg_K(V)$ ist, falls $X$
 ein vollst\"andig isotroper bzw. vollst\"andig singul\"arer Unterraum von
 $V$ ist.
 \medskip\noindent
 {\bf 6.4. Satz.} {\it Es sei $f$ eine spurwertige, nicht entartete, symmetrische
 $\alpha$-Form auf dem $K$-Vektorraum $V$, bzw. die der nicht entarteten
 quadratischen Form $Q$ auf $V$ assoziierte Bilinearform. Ferner seien $X_1$
 und $X_2$ zwei maximale vollst\"andig isotrope bzw. vollst\"andig singul\"are
 Unterr\"aume von $V$. Setze $X := X_1 \cap X_2$. Es seien $S_1$ und $S_2$
 Unterr\"aume mit $X_i = X \oplus S_i$. Schlie\ss lich setzen wir $S :=
 S_1 + S_2$. Es gibt dann zwei Unterr\"aume $U$ und $W$ von $V$ mit:
 \item{a)} Die Unterr\"aume $U + X$, $S$ und $W$ sind nicht isotrop und
 paarweise orthogonal.
 \item{b)} Es ist $V = X \oplus S \oplus U \oplus W$.
 \item{c)} Es gibt keinen von 0 verschiedenen isotropen bzw. singul\"aren
 Vektor in $W$.
 \item{d)} $U$ ist vollst\"andig isotrop bzw. vollst\"andig singul\"ar.
 \par\noindent
       Dar\"uber hinaus ist $\Rg_K(X_1) = \Rg_K(X_2)$, $\Rg_K(U) = \Rg_K(X)$,
 $\Rg_K(S_1) = \Rg_K(S_2)$ und $\Rg_K(W) = \Rg_K(V) - 2\Rg_K(X_1)$.}
 \smallskip
       Beweis. Ist $M$ ein maximaler vollst\"andig isotroper bzw.\
 vollst\"an\-dig
 singul\"arer Unterraum, ist $x \in M^\bot$ und ist $x$ isotrop bzw.\
 singul\"ar, so ist $x \in M$, da andernfalls $M + xK$ ein vollst\"andig
 isotroper bzw.\ vollst\"an\-dig singul\"arer Unterraum von $V$ w\"are, der $M$
 echt umfasste. Ist also $x \in X_i^\bot$ und ist $x$ isotrop bzw.\
 singul\"ar, so ist $x \in X_i$.
 \par
       Ist $y \in S_1 \cap S_2^\bot$, so ist $y \in X_1^\bot$, da ja
 $S_1 \leq X_1 \leq X_1^\bot$ ist. Weil $X \leq X_1$ gilt, ist also
 $y \in X^\bot$. Somit ist
 $$ y \in X^\bot \cap S_2^\bot = (X + S_2)^\bot = X_2^\bot. $$
 Da $y$ isotrop bzw. singul\"ar ist, ist also $y \in X_2$. Also ist
 $$ y \in S_1 \cap X_2 = S_1 \cap X_1 \cap X_2 = S_1 \cap X = \{0\}. $$
 Somit ist $S_1 \cap S_2^\bot = \{0\}$. Ebenso folgt $S_2 \cap S_1^\bot =
 \{0\}$.
 \par
       Weil $f$ nicht entartet ist, ist also, da ja $S_1 \cap S_2^\bot = 
 \{0\}$ ist, $V = S_1 \oplus S_2^\bot$. Also ist
 $$ \Rg_K(S_1) + \Rg_K(S_2^\bot) = \Rg_K(V). $$
 Andererseits ist
 $$ \Rg_K(S_2) + \Rg_K(S_2^\bot) = \Rg_K(V), $$
 so dass
 $$ \Rg_K(S_1) = \Rg_K(S_2) $$
 ist. Hieraus und aus $S_1 \cap S_2^\bot = \{0\}$ folgt nach 6.2 a), dass
 $S = S_1 + S_2$ nicht isotrop ist. Aus 1.5 folgt weiter, dass $S^\bot$
 nicht isotrop ist. Es ist $X \leq S^\bot$. Nach 6.3, angewandt auf $X$ und
 $S^\bot$, gibt es daher einen vollst\"andig isotropen bzw. vollst\"andig
 singul\"aren Unterraum $U \leq S^\bot$ mit $\Rg_K(X) = \Rg_K(U)$ und
 $X \cap U = \{0\}$, so dass $X + U$ nicht singul\"ar ist. Wir setzen
 $$ W := (X + U)^\bot \cap S^\bot. $$
 Weil $S$ nicht isotrop ist, ist $V = S \oplus S^\bot$. Daher ist, da ja
 $S \leq (X + U)^\bot $
 gilt,
 $$\eqalign{
 (X + U)^\bot &= (X + U)^\bot \cap (S \oplus S^\bot) \cr
       &= S \oplus \bigl((X + U)^\bot \cap S^\bot\bigr)        \cr
       &= S \oplus W.                                \cr} $$
 Weil $X + U$ und damit $(X + U)^\bot$ nicht isotrop sind und dies auch von
 $S$ gilt, ist auch $W$ nicht isotrop. Weil ferner $S \leq (X + U)^\bot$,
 $W \leq S^\bot$ und $X + U \leq W^\bot$ gilt, ist a) erf\"ullt. Ferner ist
 $$ V = X \oplus U \oplus (X \oplus U)^\bot = X \oplus U \oplus S \oplus W, $$
 so dass auch b) gilt
 \par
       Es sei $w \in W$ und $w$ sei isotrop bzw. singul\"ar. Dann ist
 $$ w \in (X + S)^\bot \leq X_1^\bot, $$
 so dass nach unserer Bemerkung zu Anfang des Beweises $w \in X_1$ gilt.
 Daher ist
 $$ w \in X_1 \cap W \leq (X + S) \cap W = \{0\}. $$
 Damit ist auch c) bewiesen.
 \par
       Wir haben bereits gezeigt, dass
 $$ \Rg_K(S_1) = \Rg_K(S_2) $$
 und
 $$ \Rg_K(X) = \Rg_K(U) $$
 ist. Alle \"ubrigen Rangaussagen folgen aus diesen.
 \medskip
       Es sei $V$ ein Vektorraum und $f$ sei eine spurwertige $\alpha$-Form
 auf $V$ oder die einer quadratischen Form assoziierte Bilinearform. Sind
 $X$, $Y$ und $Z$ Unterr\"aume von $V$ und ist
 $V = X \oplus Y \oplus Z$, so nennt man diese Zerlegung von $V$ in eine
 direkte Summe {\it wittsche Zerlegung von V\/},\index{wittsche Zerlegung}{}
 falls $X$ und $Y$
 vollst\"andig isotrop bzw. vollst\"andig singul\"ar sind, falls $Z$ nicht
 isotrop ist und falls $Z \leq (X + Y)^\bot$ gilt.
 \medskip\noindent
 {\bf 6.5. Satz.} {\it Es sei $f$ eine spurwertige, nicht ausgeartete
 symmetrische $\alpha$-Form auf dem $K$-Vektorraum $V$ bzw. die der nicht
 ausgearteten quadratischen Form $Q$ auf $V$ assoziierte Bilinearform. Dann
 gilt: 
 \item{(a)} Sind $X$ und $Y$ zwei maximale vollst\"andig isotrope bzw.
 voll\-st\"an\-dig singul\"are Teilr\"aume von $V$, so ist
 $ \Rg_K(X) = \Rg_K(Y)$. 
 \item{(b)} Ist $X$ ein maximaler vollst\"andig isotroper bzw. vollst\"andig
 singul\"arer Teilraum von $V$, so gibt es einen maximalen vollst\"andig
 isotropen bzw. voll\-st\"an\-dig singul\"aren Teilraum $Y$ von $V$ mit
 $X \cap Y = \{0\}$.
 \item{(c)} Sind $X$ und $Y$ maximale vollst\"andig isotrope bzw. vollst\"andig
 singul\"are Teil\-r\"aume mit $X \cap Y = \{0\}$, so ist $X + Y$ nicht
 isotrop und $(X + Y)^\bot$ enth\"alt keinen von $0$ verschiedenen isotropen
 bzw. singul\"aren Vektor. Insbesondere ist $X$, $Y$, $(X + Y)^\bot$ eine
 wittsche Zerlegung von $V$.
 \item{(d)} $V$ besitzt eine wittsche Zerlegung.}
 \smallskip
       Dies folgt unmittelbar aus 6.4 und 6.3.
 \medskip
       Ist $\nu$ der Rang eines maximalen vollst\"andig isotropen bzw.
 eines maximalen vollst\"andig
 singul\"aren Teilraumes von $V$, so ist $\nu$ nach 6.5a) der Rang
 eines jeden maximalen vollst\"andig isotropen bzw. singul\"aren Teilraumes.
 Die Zahl $\nu$ hei\ss t {\it Index\/}\index{Index}{} von $f$ bzw. von $Q$. Nach
 6.5b) ist $2\nu \leq \Rg_K(V)$. Ist $2\nu = \Rg_K(V)$, so hei\ss t $f$ bzw. $Q$
 von {\it maximalem Index\/}.\index{maximaler Index}{} Formen von maximalem Index
 gibt es also h\"ochstens dann, wenn der Rang von $V$ gerade ist. In diesem Falle
 gibt es aber auch immer eine. Genauer gilt:
 \medskip\noindent
 {\bf 6.6. Satz.} {\it Es sei $V$ ein Vektorraum des Ranges $2\nu$. Dann gibt es
 bis auf Isometrie genau eine quadratische Form des Index $\nu$ auf $V$.}
 \smallskip
       Beweis. Es sei $Q$ eine quadratische Form vom Index $\nu$ auf $V$ und
 $f$ sei die zugeh\"orige Bilinearform. Nach 6.3 gibt es dann zwei
 vollst\"andig singul\"are Unterr\"aume $X$ und $Y$ des Ranges $\nu$ mit
 $X \cap Y = \{0\}$, so dass $X + Y$ nicht isotrop ist. Dann ist
 $$ V = X \oplus Y, $$
 da ja $\Rg_K(V) = \Rg_K(X) + \Rg_K(Y)$ ist. Es sei
 $v \in V$. Es gibt dann $x \in X$ und $y \in Y$ mit $v = x + y$. Es folgt
 $$ Q(v) = Q(x) + Q(y) + f(x,y) = f(x,y). $$
 Wegen $X = X^\bot$ ist $X^\bot \cap Y = \{0\}$, so dass es nach 6.2
 eine Basis $b_1, \dots, b_\nu$ von $X$ und eine Basis $c_1, \dots, c_\nu$
 von $Y$ gibt mit $f(b_i,c_j) = \delta_{ij}$ f\"ur alle $i$ und $j$.
 Hieraus und aus $Q(v) = f(x,y)$ folgt, dass $Q$ bis auf Isometrie eindeutig
 festliegt.
 \par
       Es seien umgekehrt $X$ und $Y$ Unterr\"aume des Ranges $\nu$ von $V$
 und es gelte $V = X \oplus Y$. Ist $b_1$, \dots, $b_\nu$ eine Basis von
 $X$ und $c_1$, \dots, $c_\nu$ eine solche von $Y$, so definieren wir $f$
 durch
 $$ f(x,y) := \sum_{i:=1}^\nu x_iy_i, $$
 falls $x = \sum_{i:=1}^\nu b_ix_i$ und $y = \sum_{i:=1}^\nu c_iy_i$ ist,
 und weiter durch
 $$ f(x + y,x' + y') = f(x,y') + f(x',y) $$
 f\"ur alle $x$, $x' \in X$ und alle $y$, $y' \in Y$. Dann ist
 $$ f(y,x) = f(0 + y,x + 0) = f(0,0) + f(x,y) = f(x,y). $$
 Ferner ist
 $$ Q((x + y)k) = f(xk,yk) = f(x,y)k^2 = Q(x + y)k^2 $$
 und
 $$\eqalign{
 Q(x+ y + x' + y') - Q(x + y) &- Q(x' + y')     \cr
       &= f(x + x',y + y') - f(x,y) - f(x',y')  \cr
       &= f(x,y') + f(x',y)                    
       = f(x + y,x' + y').                     \cr} $$
 Damit ist gezeigt, dass $Q$ eine quadratische Form ist, die \"uberdies
 nicht ausgeartet ist. Weil $X$ und $Y$ vollst\"andig singul\"ar sind, ist
 $\nu$ der Index von $Q$. Damit ist alles bewiesen.

\mysection{7. Der Satz von Witt}

\noindent
       Der Satz von Witt, den wir in diesem Abschnitt beweisen werden, ist
 von gro\ss er Wichtigkeit beim Studium der Zentralisatoren von Polarit\"aten.
 Bevor wir ihn formulieren, jedoch noch einige Bemerkungen. Zun\"achst
 erweitern wir den Begriff der Isometrie etwas.
 \par
       Ist $f$ eine $\alpha$-Form auf dem $K$-Vektorraum $V$ und $f'$ eine
 $\alpha$-Form auf dem $K$-Vektorraum $V'$, ist $\sigma$ eine bijektive
 lineare Abbildung von $V$ auf $V'$ und gilt $f'(x^\sigma,y^\sigma) = f(x,y)$
 f\"ur alle $x$, $y \in V$, so hei\ss t $\sigma$
 {\it Isometrie\/}\index{Isometrie}{} von $V$
 auf $V'$. Ferner nennen wir $\sigma$ auch dann eine Isometrie von $V$ auf
 $V'$, wenn $Q$ eine quadratische Form auf $V$ und $Q'$ eine quadratische
 Form auf $V'$ ist und wenn f\"ur alle $x \in V$ die Gleichung $Q'(x^\sigma)
 = Q(x)$ gilt.
 \par
       Es sei $f$ eine nicht entartete $\alpha$-Form auf $V$ und $f'$ sei
 eine $\alpha$-Form auf $V'$. Ferner sei $\sigma$ eine lineare Abbildung
 von $V$ in $V'$ mit $f'(x^\sigma,y^\sigma) = f(x,y)$ f\"ur alle $x$, $y \in
 V$. Ist $y^\sigma = 0$, so ist
 $$ f(x,y) = f'(x^\sigma,y^\sigma) = f'(x^\sigma,0) = 0 $$
 f\"ur alle $x \in V$, so dass $y = 0$ ist, da $f$ als nicht entartet
 vorausgesetzt wurde. Also ist $\sigma$ injektiv. Ist $\Rg_K(V) = \Rg_K(V')$,
 so ist $\sigma$ also eine Bijektion und damit eine Isometrie.
 \medskip\noindent
 {\bf 7.1. Satz.} {\it Ist $f$ eine symmetrische $\alpha$-Form auf $V$ bzw. die
 zu einer quadratischen Form $Q$ assoziierte Bilinearform, sind ferner
 $X_1$ und $X_2$ Teilr\"aume von $V$ sowie $\sigma_1$ und $\sigma_2$
 Isometrien von $X_1$ bzw. $X_2$ in $V$, gilt
 $$ X_1 \cap X_2 = \{0\} = X_1^{\sigma_1} \cap X_2^{\sigma_2} $$
 und
 $$ f(x_1^{\sigma_1},x_2^{\sigma_2}) = f(x_1,x_2) $$
 f\"ur alle $x_1 \in X_1$ und alle $x_2 \in X_2$ {\rm (}dies ist eine zus\"atzliche
 Eigenschaft zu den Isometriebedingungen{\rm )}, so ist die durch
 $$ (x_1 + x_2)^\nu := x_1^{\sigma_1} + x_2^{\sigma_2} $$
 definierte lineare Abbildung $\nu$ von $X_1 \oplus X_2$ auf $X_1^{\sigma_1}
 \oplus X_2^{\sigma_2}$ eine Isometrie von $X_1 \oplus X_2$ in $V$.}
 \smallskip
       Beweis. Es ist banal, dass $\nu$ eine bijektive lineare Abbildung von
 $X_1 \oplus X_2$ auf $X_1^{\sigma_1} \oplus X_2^{\sigma_2}$ ist.
       Ferner ist
 $$\eqalign{
 f\bigl((x_1 + x_2)^\nu,(y_1 &+ y_2)^\nu\bigr)
			= f(x_1^{\sigma_1} + x_2^{\sigma_2},
                                 y_1^{\sigma_1} + \sigma_2^{\sigma_2})     \cr
   &= f(x_1^{\sigma_1},y_1^{\sigma_1}) + f(x_1^{\sigma_1},y_2^{\sigma_2})
       + f(x_2^{\sigma_2},y_1^{\sigma_1})                                 
           + f(x_2^{\sigma_2},y^{\sigma_2})                             \cr
   &= f(x_1,y_1) + f(x_1,y_2) + f(x_2,y_1)                              
      + f(x_2,y_2)                                                 \cr
   &= f(x_1 + x_2,y_1 + y_2),                                         \cr} $$
 bzw.
 $$\eqalign{
 Q\bigl((x_1 + x_2)^\nu\bigr) &= Q(x_1^{\sigma_1} + x_2^{\sigma_2})      \cr
   &= Q(x_1^{\sigma_1}) + Q(x_2^{\sigma_2}) + f(x_1^{\sigma_1},x_2^{\sigma_2})
                                                               \cr
   &= Q(x_1) + Q(x_2) + f(x_1,x_2)                             \cr
   &= Q(x_1 + x_2).                                            \cr} $$
 Damit ist gezeigt, dass $\nu$ eine Isometrie ist.
 \medskip\noindent
 {\bf 7.2. Satz von Witt.} {\it Es sei $f$ eine spurwertige, symmetrische,
 nicht entartete $\alpha$-Form auf dem $K$-Vektorraum $V$, bzw. die der
 nicht entarteten quadratischen Form $Q$ auf $V$ assoziierte Bilinearform.
 Ist $X$ ein Teilraum von $V$ und ist $\sigma$ eine injektive lineare
 Abbildung von $X$ in $V$ mit $f(x^\sigma,y^\sigma) = f(x,y)$ f\"ur alle
 $x$, $y \in X$, bzw. $Q(x^\sigma) = Q(x)$ f\"ur alle $x \in X$, so gibt es
 eine Isometrie $\tau$ von $V$ auf sich mit $x^\tau = x^\sigma$ f\"ur alle
 $x \in X$.}
 \smallskip
       Beweis. Setze
 $$ Y := \{y \mid y \in X, y^\sigma = y\} \quad
 \mathrm{und}\quad
  P := \{x^\sigma - x \mid x \in X\}. $$
 Dann ist $P \cong X/Y$.
 \par
       Wir betrachten zun\"achst den Fall, wo $Y$ eine Hyperebene von $X$ ist.
 Wegen $P \cong X/Y$ ist $P$ dann ein Punkt.
 \par
       Ist $X' \leq P^\bot$ und $X \cap X' = X^\sigma \cap X'
 = \{0\}$, so gilt f\"ur $x \in X$ und $y \in X'$ wegen $X' \leq P^\bot$
 die Gleichung
 $$ 0 = f(x^\sigma - x,y) = f(x^\sigma,y) - f(x,y), $$
 so dass $f(x^\sigma,y) = f(x,y)$ ist. Definiert man $\nu$ durch
 $$ (x + y)^\nu := x^\sigma + y, $$
 so ist $\nu$ nach 7.1 eine Isometrie von
 $X \oplus X'$ in $V$, die \"uberdies die Vektoren aus $X'$ festl\"asst.
 Hieraus folgt weiter
 $$ P = \{x^\nu - x \mid x \in X + X'\}. $$
 \par
       Sind $x$, $y \in X$, so ist
 $$ f(x^\sigma,y^\sigma - y) = f(x^\sigma,y^\sigma) - f(x^\sigma,y)     = f(x,y) - f(x^\sigma,y)
   = f(x - x^\sigma,y). \leqno(*) $$
 Ist \"uberdies $x \in Y$, so ist nach $(*)$ also
 $$ f(x,y^\sigma - y) = f(x^\sigma,y^\sigma - y) = f(x - x^\sigma,y) = 0, $$
 da in diesem Falle ja $x = x^\sigma$ ist. Also ist $Y \leq P^\bot$.
 \par
       1. Fall: Es ist $X \not\leq P^\bot$. Aus $(*)$ folgt, dass dann auch
 $X^\sigma \not\leq P^\bot$ ist. Weil $Y$ eine Hyperebene von $X$ und damit
 auch von $X^\sigma$ ist und weil
 $$ Y \leq X \cap P^\bot,\ X^\sigma \cap P^\bot $$
 gilt, ist also
 $$ X \cap P^\bot = Y = X^\sigma \cap P^\bot. $$
 W\"ahle ein Komplement $X'$ von $Y$ in $P^\bot$. Dann ist $X \cap X'
 = \{0\}$ und $X^\sigma \cap X' = \{0\}$, so dass sich $\sigma$, wie
 eingangs bemerkt, zu einer Isometrie $\nu$ von $X + X'$ in $V$ fortsetzen
 l\"asst. Nun ist
 $$ P^\bot = Y \oplus X' < X + X' \leq V, $$
 so dass $X + X' = V$ ist, da $P^\bot$ ja eine Hyperebene von $V$ ist.
 In diesem Falle setze man $\tau := \nu$. Dann ist $\tau$ eine m\"ogliche
 Fortsetzung von $\sigma$.
 \par
       2. Fall: Es ist $X \leq P^\bot$. Aus $(*)$ folgt, dass dann auch
 $X^\sigma \leq P^\bot$ ist. Hieraus folgt weiter
 $$ P \leq X + X^\sigma \leq P^\bot, $$
 so dass $P$ isotrop bzw. singul\"ar ist. F\"ur Letzteres muss man noch
 beachten, dass f\"ur $p \in P$ die Gleichung $Q(p) = 0$ gilt. Ist $p \in P$,
 so gibt es ein $x \in X$ mit $p = x^\sigma - x$ und es folgt
 $$
 Q(x^\sigma - x) = Q(x^\sigma) + Q(x) - f(x^\sigma,x)     
       = 2Q(x) - f(x^\sigma - x,x) - f(x,x) = 0.        
 $$
 \par
       Nach I.7.2 gibt es einen Teilraum $X'$ von $P^\bot$ mit
 $$ X \oplus X' = P^\bot = X^\sigma \oplus X'. $$
 Nach 7.1 l\"asst sich $\sigma$ zu einer
 Isometrie von $X \oplus X'$ auf $X^\sigma \oplus X'$, dh. zu einer
 Isometrie von $P^\bot$ fortsetzen, die \"uberdies die Hyperebene $Y + X'$ von
 $P^\bot$ elementweise festl\"asst. Wir d\"urfen daher im weiteren Verlauf
 des Beweises annehmen, dass $X = P^\bot$ ist und dass $\sigma$ eine
 Isometrie von $X$ auf sich ist. Dies wird im Folgenden wesentlich benutzt.
 \par
       Ist $z \in V$, so ist die Abbildung $x \to f(x^{\sigma^{-1}},z)$
 f\"ur $x \in X$ eine lineare Abbildung von $X$ in $K$. Es gibt ferner eine
 lineare Abbildung $\varphi$ von $V$ in $K$ mit $\varphi(x) =
 f(x^{\sigma^{-1}},z)$ f\"ur alle $x \in X$. Weil $f$ nicht ausgeartet ist,
 gibt es ein $z' \in V$ mit $f(v,z') = \varphi(v)$ f\"ur alle $v \in V$.
 Ersetzt man $x$ durch $x^\sigma$, so sieht man,
 dass es zu jedem $z \in V$ ein $z' \in V$ gibt mit
 $$ f(x,z) = \varphi(x^\sigma) = f(x^\sigma,z') $$
 f\"ur alle $x \in X$. Ist $y \in X$, so gilt wegen $X = X^\sigma = P^\bot$ die
 Gleichung
 $$ f(x^\sigma,z') = f(x^\sigma,z' + y^\sigma - y). \leqno(**) $$
 Ist $z \not\in X$, so ist auch $z' \not\in X$. W\"are n\"amlich $z' \in X$,
 so w\"are
 $$ f(x^\sigma - x,z') = 0 $$
 f\"ur alle $x \in X$, woraus
 $f(x^\sigma,z') = f(x,z)$ f\"ur alle $x \in X$ folgte. Daher g\"alte nach
 $(**)$ die Gleichung
 $f(x,z') = f(x,z)$ f\"ur alle $x \in X$, was wiederum $f(x,z - z') = 0$
 nach sich z\"oge, so dass
 $$ z - z' \in X^\bot = P^{\bot\bot} = P \leq X $$
 w\"are im Widerspruch zu $z \not\in X$. Da $P$ vollst\"andig isotrop bzw.
 vollst\"andig singul\"ar ist und weil $f$ als spurwertig vorausgesetzt wurde,
 gibt es nach 6.1 ein $y \in X$ mit
 $$ f(z' + y^\sigma - y,z' + y^\sigma - y) = f(z,z) $$
 bzw.
 $$ Q(z' + y^\sigma - y) = Q(z). $$
 Weil schlie\ss lich
 $$ f(x^\sigma,z') = f(x^\sigma,z' + y^\sigma - y) $$
 f\"ur alle $x \in X$ gilt, ersetzen wir $z'$ durch $z' + y^\sigma - y$ um
 zu sehen, dass es zu $z \in V$ mit $z \not\in X$ ein $z' \in V$ mit
 $z' \not\in X$ gibt, so dass
 $$ f(x^\sigma,z') = f(x,z') $$
 f\"ur alle $x \in X$ gilt, sowie
 $$ f(z',z') = f(z,z) $$
 bzw.
 $$ Q(z') = Q(z). $$
 Hieraus folgt, dass sich $\sigma$ zu einer Isometrie von $X \oplus zK = V$
 auf $X \oplus z'K = V$ fortsetzen l\"asst.
 \par
       Der Rest des Beweises folgt nun mit Induktion nach $\Rg_K(X)$.
 Ist $\Rg_K(X) = 0$, so ist der Satz trivial. Es sei also $n := \Rg_K(X) > 0$
 und $Z \leq X$ sowie $\Rg_K(Z) = n - 1$. Schr\"ankt man $\sigma$ auf $Z$
 ein, so gibt es nach Induktionsannahme eine Isometrie $\rho$ von $V$ mit
 $y^\sigma = y^\rho$ f\"ur alle $y \in Z$.
 Ist sogar $x^\sigma = x^\rho$ f\"ur alle $x \in X$, so ist $\tau := \rho$
 eine Fortsetzung von $\sigma$ auf $V$. Es sei also $x^\rho \neq x^\sigma$
 f\"ur wenigstens ein $x \in X$. Dann ist $\sigma\rho^{-1}$ eine Isometrie
 von $X$ in $V$ mit
 $$ Z = \{x \mid x \in X, x^{\sigma\rho^{-1}} = x\}. $$
 Nach dem bereits Bewiesenen --- $Z$ spielt jetzt die Rolle von $Y$ ---
 gibt es dann eine Isometrie $\nu$ von $V$ mit
 $x^\nu = x^{\sigma\rho^{-1}}$ f\"ur alle $x \in X$. Die Abbildung $\tau :=
 \nu\rho$ ist dann eine Isometrie von $V$ mit $x^\sigma = x^\tau$ f\"ur alle
 $x \in X$. Damit ist 7.2 bewiesen.
 \medskip\noindent
 {\bf 7.3. Korollar.} {\it Es sei $f$ eine spurwertige, nicht entartete
 $\alpha$-Form auf dem Vektorraum $V$, bzw. die der nicht entarteten
 quadratischen Form $Q$ auf $V$ assoziierte Bilinearform. Ist dann $G$ die
 Gruppe der zugeh\"origen Isometrien, so gilt:
 \item{a)\,} Ist $M$ die Menge der vollst\"andig isotropen, bzw. vollst\"andig
 singul\"aren Un\-ter\-r\"au\-me vom Rang $i$ von $V$, so ist $G$ auf $M$ transitiv.
 \item{b)\,} Ist $U$ ein vollst\"andig isotroper bzw. vollst\"andig singul\"arer
 Unterraum von $V$, so wird jede Isometrie von $U$ auf sich von einem
 $\sigma \in G$ induziert.\par}
 \smallskip
       Dies folgt unmittelbar aus 7.2.

\mysection{8. Unit\"are Gruppen}

\noindent
       Als N\"achstes studieren wir die Zentralisatoren von unit\"aren
 Po\-la\-ri\-t\"a\-ten, die durch spurwertige $\alpha$-Formen dargestellt werden.
 Das uns hier vor allem interessierende Resultat ist, dass die von
 Elationen erzeugte Untergruppe einer solchen Gruppe bis auf drei Ausnahmen
 einfach ist. Wir werden diesen Satz nicht in voller Allgemeinheit
 beweisen, aber nur solche Geometrien ausschlie\ss en, deren
 Koordinatenk\"orper nicht kommutativ sind und ein Zentrum haben, welches
 zu $\GF(2)$, $\GF(3)$, $\GF(4)$ oder $\GF(9)$ isomorph ist. Dass wir diese
 sehr fremdartigen K\"orper aus der Betrachtung ausschlie\ss en, liegt an
 unserer Beweismethode. Es gibt zwar Schiefk\"orper, deren Zentrum zu einem
 der genannten Galoisfelder isomorph ist (P. M. Cohn 1977, S. 116/117),
 ich wei\ss\ aber nicht, ob es unter ihnen auch solche gibt, die einen
 Antiautomorphismus gestatten. Das Studium der unit\"aren Gruppen wird uns
 drei Abschnitte lang besch\"aftigen.
 \par
       Bei den Studien dieses Abschnitts werden wir uns des Satzes 1.9
 bedienen, der besagt, dass sich unit\"are Polarit\"aten
 nicht nur durch symmetrische $\alpha$-Formen, sondern auch durch
 antisymmetrische $\alpha$-Formen\index{antisymmetrische Form}{} darstellen
 lassen. Dabei hei\ss t eine $\alpha$-Form {\it antisymmetrisch\/}, falls
 $$ f(u,v)^\alpha = -f(v,u) $$
 ist f\"ur alle $u$, $v \in V$ und falls im Falle $\alpha = 1$ au\ss erdem
 $f(x,x) = 0$ f\"ur alle $x \in V$ gilt. Dies hat den Vorteil, dass die
 hier vorzutragenden S\"atze auch f\"ur die symplektischen Gruppen gelten.
 Dies werden wir nicht weiter ausnutzen, da wir die Einfachheit der
 symplektischen Gruppen\index{symplektische Gruppe}{} ja bereits bewiesen haben.
 \medskip\noindent
 {\bf 8.1. Satz.} {\it Es sei $V$ ein Vektorraum \"uber dem K\"orper $K$ und $f$
 sei eine nicht ausgeartete, antisymmetrische $\alpha$-Semibilinearform auf
 $V$. Ferner sei $\varphi$ eine lineare Abbildung von $V$ in $K$ und $h$
 sei ein Element des Kernes von $\varphi$. Ist $\tau$ die durch
 $$ x^\tau := x + h\varphi(x) $$
 definierte Transvektion von $V$, so ist $\tau$ genau dann eine Isometrie
 von $V$, wenn es ein $\lambda \in K$ gibt mit $\lambda^\alpha = \lambda$
 und
 $$ \varphi(x) = \lambda f(h,x) $$
 f\"ur alle $x \in V$.}
 \smallskip
 Transvektionen dieser Art nennen wir in Zukunft {\it unit\"ar\/}.
 \smallskip
       Beweis. Wir d\"urfen annehmen, dass $\tau \neq 1$ ist.
       Es sei $\tau$ eine Isometrie. Dann ist
 $$\eqalign{
  f(u,v) &= f(u^\tau,v^\tau)                                 
       = f(u + h\varphi(u),v + h\varphi(v))                   \cr
       &= f(u,v) + f(u,h)\varphi(u) + \varphi(u)^\alpha f(h,v)
              + \varphi(u)^\alpha f(h,h)\varphi(v)             \cr} $$
 f\"ur alle $u$, $v \in V$. Es folgt
 $$ 0 = f(u,h)\varphi(v) + \varphi(u)^\alpha f(h,v)
                     + \varphi(u)^\alpha f(h,h)\varphi(v) $$
 f\"ur alle $u$, $v \in V$. Setze $u := h$. Dann ist
 $$
 0 = f(h,h)\varphi(v) + 0^\alpha f(h,v) + 0^\alpha f(h,h)\varphi(v) 
   = f(h,h)\varphi(v) 
 $$
 f\"ur alle $v \in V$. Wegen $\tau \neq 1$ gibt es ein $v \in V$ mit
 $\varphi(v) \neq 0$, so dass $f(h,h) = 0$ ist.
 \par
       Wegen $f(h,h) = 0$ ist
 $$ 0 = f(u,h)\varphi(v) + \varphi(u)^\alpha f(h,v) $$
 f\"ur alle $u$, $v \in V$. Da $f$ nicht entartet ist und wegen $\tau \neq 1$
 auch $h \neq 0$ gilt, gibt es ein $u \in V$ mit $f(u,h) = 1$. Es folgt
 $$ \varphi(v) = -\varphi(u)^\alpha f(h,v) $$
 f\"ur alle $v \in V$. Setze $\lambda := -\varphi(u)^\alpha$. Dann ist
 also $\varphi(v) = \lambda f(h,v)$
 f\"ur alle $v \in V$. Es folgt
 $$ \varphi(u) = \lambda f(h,u) = \lambda \bigl(-f(u,h)^\alpha\bigr) = -\lambda
       = \varphi(u)^\alpha. $$
 Hieraus folgt wiederum $\lambda^\alpha = \lambda$.
 \par
       Ist umgekehrt $\lambda^\alpha = \lambda$ und $f(h,h) = 0$, so zeigt
 eine banale Rechnung, dass die durch
 $ x^\tau := x + h\lambda f(h,x) $
 definierte Transvektion $\tau$ eine Isometrie ist.
 \medskip\noindent
 {\bf 8.2. Satz.} {\it Es sei $V$ ein Vektorraum \"uber dem K\"orper $K$ und $f$
 sei eine nicht ausgeartete, antisymmetrische $\alpha$-Form auf $V$. Ist $P$
 ein isotroper Punkt, ist $G$ eine nicht isotrope Gerade durch $P$
 und sind $Q$ und $R$ von $P$ verschiedene isotrope Punkte auf $G$, so
 gibt es eine unit\"are Transvektion $\tau$ mit dem Zentrum $P$, f\"ur die
 $Q^\tau = R$ gilt.}
 \smallskip
       Beweis. Es sei $P = pK$, $Q = qK$ und $R = rK$. Dann sind die Vektoren
 $p$, $q$ und $r$ isotrop. Da $G$ nicht isotrop ist, ist $f(p,q)$, $f(p,r)
 \neq 0$. Indem man $q$ und $r$ ggf. mit einem Skalar multipliziert, kann
 man erreichen, dass $f(p,q) = f(p,r) = 1$ ist. Weil $p$ und $q$ linear
 unabh\"angig sind, gibt es $\mu$ und $\lambda$ mit $r = p\lambda + q\mu$. Es
 folgt
 $$ 1 = f(p,r) = f(p,p\lambda + q\mu) = \mu $$
 und
 $$ 0 = f(r,r) = f(p\lambda + q,p\lambda + q) = \lambda^\alpha f(p,q) +
       f(q,p)\lambda = \lambda^\alpha - \lambda. $$
 Also ist $\lambda^\alpha = \lambda$. Definiert man nun $\tau$ durch
 $$ x^\tau := x + p\lambda f(p,x), $$
 so ist $\tau$ nach 8.1 eine unit\"are Transvektion mit dem Zentrum $P$ und
 es gilt
 $$ q^\tau = q + p\lambda = r, $$
 so dass $Q^\tau = R$ ist. Damit ist alles bewiesen.
 \medskip
       Ist $g$ eine spurwertige, nicht ausgeartete, symmetrische
 $\alpha$-Form mit $\alpha \neq 1$, so gibt es nach Satz 1.9 eine nicht
 ausgeartete, antisymmetrische $\beta$-Form $f$
 und ein $k \in K^*$ mit $f(u,v) = kg(u,v)$. Die
 Isometrien von $f$ sind offenbar die gleichen wie die von $g$. Ist $g$
 spurwertig, so k\"onnen wir die S\"atze \"uber spurwertige $\alpha$-Formen
 auch auf $f$ anwenden. Um dies auszudr\"ucken werden wir auch $f$ {\it
 spurwertig\/} nennen, wenn der gerade geschilderte Zusammenhang besteht.
 \medskip\noindent
 {\bf 8.3. Satz.} {\it Es sei $f$ eine nicht ausgeartete, antisymmetrische,
 spur\-wer\-ti\-ge $\alpha$-Form auf dem $K$-Vektorraum $V$. Mit $\Sigma$
 bezeichnen wir die Elemente $s \in K$, f\"ur die $s^\alpha = s$ gilt.
 Es sei $G$ eine nicht isotrope Gerade von $\La(V)$ und $P$ sei ein isotroper
 Punkt auf $G$. Es gibt dann einen weiteren isotropen Punkt $Q$ auf $G$.
 Ist dann $P = pK$ und $Q = qK$, so ist
 $$ \bigl\{(q + ps)K \mid s \in \Sigma\bigr\} $$
 die Menge der von $P$ verschiedenen isotropen Punkte auf $G$. \"Uberdies
 gilt: Ist $(q + ps)K = (q + pt)K$, so ist $s = t$. Insbesondere gilt,
 dass $G$ mindestens drei isotrope Punkte tr\"agt.}
 \smallskip
       Beweis. Weil $G$ nicht isotrop ist,
 ist die Einschr\"ankung von $f$ auf $G$ nicht ausgeartet. Nach 6.3 gibt
 es daher einen isotropen Punkt $Q$ von $G$ mit $G = P + Q$. Wiederum weil
 $G$ nicht isotrop ist, gilt $f(p,q) \neq 0$, falls $0 \neq p \in P$ und
 $0 \neq q \in Q$ ist. Daher gibt es ein $p \in P$ und ein $q \in Q$ mit
 $f(p,q) = 1$. Ist nun $R$ ein von $P$ verschiedener isotroper Punkt auf
 $G$, so gibt es, wie wir beim Beweise von 8.2 gesehen haben, ein $r \in
 \Sigma$ mit $R = (q + pr)K$. Daher ist $R$ in der Menge
 $$ \bigl\{(q + ps)K \mid s \in \Sigma\bigr\} $$
 enthalten. Ist andererseits $s \in \Sigma$, so folgt
 $$\eqalign{
 f(q + ps,q + ps) &= f(q,ps) + f(ps,q) 
       = f(q,ps) - f(q,ps)^\alpha     \cr
       &= f(p,q)s - \bigl(f(p,q)s\bigr)^\alpha   
       = s - s^\alpha                
       = 0. \cr} $$
 Also ist $\{(q + ps)K \mid s \in \Sigma\}$ in der Tat die Menge der
 isotropen Punkte auf $G$.
 \par
       Sind $s$, $t \in \Sigma$ und gilt $(q + ps)K = (q + pt)K$, so folgt
 $s = t$, da $p$ und $q$ ja linear unabh\"angig sind.
 \par
       Die letzte Aussage des Satzes folgt schlie\ss lich daraus, dass
 $\Sigma$ zumindest die beiden Elemente 0 und 1 enth\"alt.
 \medskip
       Ist $f$ eine nicht ausgeartete, antisymmetrische oder auch eine nicht
 ausgeartete, symmetrische $\alpha$-Form auf dem Vektorraum $V$, so bezeichnen
 wir mit $\Un(V,f)$ die Gruppe aller Isometrien von $f$ und
 mit $\TU(V,f)$ die von allen bez.\ $f$ unit\"aren Transvektionen
 erzeugte Untergruppe von $\Un(V,f)$.
 \medskip\noindent
 {\bf 8.4. Satz.} {\it Es sei $f$ eine nicht ausgeartete, antisymmetrische,
 spur\-wer\-ti\-ge $\alpha$-Form auf dem $K$-Vektorraum $V$, deren Index
 mindestens $1$ sei. Dann ist die Gruppe $\TU(V,f)$ auf der Menge
 der isotropen Punkte von $\La(V)$ transitiv.}
 \par
       Beweis. Es seien $P = pK$ und $Q = qK$ zwei verschiedene
 isotrope Punkte von $\La(V)$. Ist $f(p,q) \neq 0$, so ist
 $P + Q$ eine nicht isotrope Gerade. Es  gibt daher nach 8.3 einen dritten
 isotropen Punkt $R$ auf $G$. Nach 8.2 gibt es dann eine Transvektion mit
 dem Zentrum $R$, die $P$ auf $Q$ abbildet.
 \par
       Ist $f(p,q) = 0$, so ist $P + Q$ vollst\"andig isotrop. Nach dem
 Korollar 6.3 gibt es einen vollst\"andig isotropen Unterraum $U$ von $V$
 mit $U \cap (P + Q)^\bot = \{0\}$. Nach 8.2 a) gibt es ferner eine Basis
 $p'$, $q'$ von $U$ mit $f(p,p') = 1 = f(q,q')$ und $f(p,q') = 0 = f(q,p')$.
 Dann ist $(p' + q')K$ ein isotroper Punkt und es gilt
 $$ f(p,p' + q') = 1 = f(q,p' + q'), $$
 so dass es, wie soeben gezeigt, eine unit\"are Transvektion gibt, die
 $P$ auf $(p' + q')K$ abbildet, und eine unit\"are Transvektion, die
 $(p' + q')K$ auf $Q$ abbildet. Damit ist alles gezeigt.
 \medskip
       Im Folgenden benutzen wir Techniken, die wir schon beim Beweise der
 Einfachheit der kleinen projektiven Gruppe in Abschnitt 2 des Kapitels II
 benutzt haben. Hier sind es insbesondere die S\"atze II.2.4 und II.2.5,
 die sich ohne M\"uhe auf den vorliegenden Fall \"uber\-tra\-gen lassen.
 \medskip\noindent
 {\bf 8.5. Satz.} {\it Es sei $f$ eine nicht ausgeartete,
 antisymmetrische $\alpha$-Form auf dem
 $K$-Vektorraum $V$. Es seien $u$ und $v$ zwei isotrope Vektoren mit
 $f(u,v) = 1$. Ferner sei $a \in K^*$ und es gelte $a^\alpha = a$.
 Nach 8.1 werden dann durch
 $$\eqalign{
 x^{\tau_1} &:= x - uf(u,x),              \cr
 x^{\tau_2} &:= x - v(1 - a)a^{-1}f(v,x), \cr
 x^{\tau_3} &:= x - uaf(u,x),             \cr
 x^{\tau_4} &:= x - v(1 - a)a^{-2}f(v,x)    \cr} $$
 vier unit\"are Transvektionen $\tau_i$ definiert. Setze
 $$ \sigma := \tau_1\tau_2\tau_3\tau_4. $$
 Dann ist $u^\sigma = ua$, $v^\sigma = va^{-1}$ und $x^\sigma = x$
 f\"ur alle $x \in (uK + vK)^\bot$.}
 \smallskip
       Beweis. Wie schon bei Satz III.2.4 sage ich auch hier: Rechnen!
 \medskip\noindent
 {\bf 8.6. Satz.} {\it Es sei $f$ eine nicht ausgeartete, antisymmetrische
 $\alpha$-Form auf dem $K$-Vektorraum $V$. Es seien $u$ und $v$ zwei
 isotrope Vektoren aus $V$ mit $f(u,v) = 1$.
 Ferner seien $a$, $b \in K^*$ und es gelte
 $a^\alpha = a$ und $b^\alpha = b$. Wir definieren $\sigma$, $\tau \in
 \TU(V,f)$ durch
 $$ (uk + vl + x)^\sigma := uak + va^{-1}l + x $$
 f\"ur alle $k$, $l \in K$ und alle $x \in (uK + vK)^\bot$ sowie
 $$ x^\tau = x + vbf(v,x) $$
 f\"ur alle $x \in V$. Dann ist
 $$ x^{\sigma^{-1}\tau^{-1}\sigma\tau} = 
       x + v(b - a^{-1}ba^{-1})f(v,x). $$}
 \par
       Beweis. Dass $\tau$ eine unit\"are Transvektion ist, folgt aus
 Satz 8.1, und dass $\sigma \in \TU(V,f)$ ist, aus Satz 8.5.
 \par
       Es ist $(a^{-1})^\alpha = a^{-1}$ und $\sigma$ ist eine Isometrie.
 Daher ist
 $$\eqalign{
 x^{\sigma^{-1}\tau^{-1}\sigma\tau}
       &= (x^{\sigma^{-1}} - vbf(v,x^{\sigma^{-1}}))^{\sigma\tau}      \cr
       &= (x - va^{-1}bf(v,x^{\sigma^{-1}}))^\tau                      \cr
       &= \bigl(x - va^{-1}ba^{-1}f(v^{\sigma^{-1}},x^{\sigma^{-1}})\bigr)^\tau
		    \cr
       &= x + vbf(v,x) - va^{-1}ba^{-1}f(v,x), \cr} $$
 q. o. o.
 \medskip
       Um diesen Satz erfolgreich anwenden zu k\"onnen, m\"ussen wir wissen,
 wann es in einem K\"orper $K$ von Null verschiedene Elemente $a$ gibt mit
 $a^\alpha = a$ und $a^2 \neq 1$. Dar\"uber gibt der n\"achste Satz Auskunft.
 \medskip\noindent
 {\bf 8.7. Satz.} {\it Es sei $K$ ein K\"orper und $\alpha$ sei ein
 Antiautomorphismus mit $\alpha^2 = 1$. Setze
 $$ L := \{ a \mid a \in K, a^\alpha = a\}. $$
 Gilt $a^2 = 1$ f\"ur alle $a \in L - \{0\}$, so ist $L = \GF(2)$ oder $\GF(3)$
 und $[K : L] \leq 2$.}
 \smallskip
       Beweis. Es sei $0 \neq a \in L$.
 Dann ist $(a - 1)(a +1) = a^2 - 1 = 0$ und daher $a = 1$ oder $-1$.
 Also ist $L = \{0, 1, -1\}$. Folglich liegt $L$ im Zentrum $Z(K)$ von
 $K$. Weil $\alpha$ in $Z(K)$ einen Automorphismus induziert, ist $L$
 ein Teilk\"orper von $Z(K)$. Ist nun $k \in K$, so ist
 $$ (k + k^\alpha)^\alpha = k^\alpha + k^{\alpha^2} = k + k^\alpha, $$
 so dass $k + k^\alpha \in L$ gilt f\"ur alle $k \in K$. Weil $\alpha$ ein
 Antiautomorphismus ist, folgt
 $$ (k^\alpha k)^\alpha = k^\alpha k^{\alpha^2} = k^\alpha k, $$
 so dass auch $k^\alpha k \in L$ gilt.
 Nun ist\footnote{Anmerkung der Herausgeber: F\"ur $a,b\in K$ folgt aus dieser Beziehung
 $ab + ba = (a+b)^2 - a^2 -b^2  \in L + La + Lb$, und $L+La+Lb+Lab$
 enth\"alt $ba$ und ist ein Teilring, also ein endlicher
 Teilk\"orper von $K$ (mit $2^4$ oder $3^4$ Elementen) und daher kommutativ (nach
 Wedderburn, oder man betrachte ein Element der multiplikativen Ordnung $5$). Also
 ist $K$ kommutativ.}
 $$ k^2 - (k + k^\alpha)k + k^\alpha k = 0. $$
 Ist also $k \not\in L$, so hat das Minimalpolynom von $k$ \"uber $L$ den
 Grad 2. Nach Satz II.10.6 ist $K$ dann kommutativ oder es ist $Z(K) = L$
 und $K$ ist ein Quaternionenschiefk\"orper \"uber $L$. Weil $L$ endlich ist,
 w\"are im letzteren Fall $K$ ein endlicher, nicht kommutativer K\"orper im
 Widerspruch zum Satz von Wedderburn, dass alle endlichen K\"orper
 kommutativ sind. Also ist $K$ kommutativ, und es folgt $[K:L] \leq 2$,
 da $L$ ja der Fixk\"orper von $\alpha$ ist.
 Nun ist aber $L = \{0, 1, -1\}$
 und folglich $L = \GF(2)$ oder $\GF(3)$. Damit ist alles bewiesen.
 \medskip\noindent
 {\bf 8.8. Satz.} {\it Es sei $K$ ein K\"orper und $f$ sei eine nicht
 ausgeartete, spurwertige, antisymmetrische $\alpha$-Form auf dem
 $K$-Vektorraum $V$. Der Fixk\"orper des von $\alpha$ auf $Z(K)$ induzierten
 Automorphismus sei von $\GF(2)$ und $\GF(3)$ verschieden. Ist dann der
 Index von $f$ mindestens $1$, so ist
 $\TU(V,f)' = \TU(V,f)$, dh., $\TU(V,f)$ ist perfekt.}
 \smallskip
       Die Voraussetzungen an $K$ sind sicher dann erf\"ullt, wenn $Z(K)$ von
 $\GF(2)$, $\GF(3)$, $\GF(4)$ und $\GF(9)$ verschieden ist.
 \smallskip
       Beweis. Es sei $\tau$ eine von 1 verschiedene unit\"are Transvektion.
 Es gibt dann einen isotropen Vektor $v$ und ein $c \in K^*$ mit
 $c^\alpha = c$, so dass
 $$ x^\tau = x + vcf(v,x) $$
 ist f\"ur alle $x \in V$. Weil $f$ nicht ausgeartet ist, gibt es eine nicht
 isotrope Gerade $G$ durch $vK$, und weil $f$ spurwertig ist, gibt
 es einen von $vK$ verschiedenen isotropen Punkt $uK$ auf $G$. Dann ist
 $f(u,v) \neq 0$, so dass wir annehmen d\"urfen, dass $f(u,v) = 1$ ist.
 \par
       Weil $\alpha$ das Zentrum von $K$ invariant l\"asst, folgt
 mit 8.7, dass es ein $a \in Z(K)^*$ gibt mit $a^\alpha = a$ und
 $a^2 \neq 1$. Setze
 $$ b := c(1 - a^{-2})^{-1}. $$
 Dann ist auch $b^\alpha = b$. Nun ist
 $$
 x^\tau = x + vcf(v,x) = x + vb(1 - a^{-2})f(v,x) 
        = x + v(b - a^{-1}ba^{-1})f(v,x), 
 $$
 so dass $\tau$ nach 8.6
 ein Kommutator ist. Weil also alle Transvektionen Kommutatoren sind,
 ist $\TU(V,f)$ perfekt.
 \medskip
       Dass $\TU(V,f)$ nicht immer perfekt ist, zeigt der folgende Satz.
 \medskip\noindent
 {\bf 8.9. Satz.} {\it Es sei $K$ ein kommutativer K\"orper und $V$ sei ein
 Vektorraum des Ranges $2$ \"uber $K$. Ist dann $f$ eine nicht ausgeartete,
 spurwertige, antisymmetrische $\alpha$-Form auf $V$, deren Index mindestens
 $1$ ist, so ist $\TU(V,f)$ zu $\SL(2,L)$ isomorph. Dabei ist $L$ der durch
 $$ L := \{a \mid a \in K, a^\alpha = a\} $$
 definierte Teilk\"orper von $K$.}
 \smallskip
       Beweis. Weil $K$ kommutativ ist, ist $\alpha$ ein Automorphismus
 von $K$, so dass $L$ in der Tat ein Teilk\"orper ist.
 \par
       Weil $f$ spurwertig ist und weil der Index von $f$ mindestens 1 ist,
 gibt es zwei isotrope Vektoren $u$ und $v$ mit $f(u,v) = 1$. Die von
 den Transvektionen $\tau$ der Form
 $$ x^\tau = x + uaf(u,x) $$
 bzw.
 $$ x^\tau = x + vaf(u,x) $$
 mit $a \in L$
 erzeugte Untergruppe $U$ von $\TU(V,f)$ ist offensichtlich isomorph zu
 $\SL(2,L)$. Andererseits operiert sie auf den Punkten von $\La(V)$ transitiv
 und enth\"alt daher alle Elationen von $\TU(V,f)$. Also ist $U = TU(V,f)$.
 Damit ist alles bewiesen.
 \medskip
       Weil die Gruppen $\SL(2,2)$ und $\SL(2,3)$ aufl\"osbar sind, zeigt 8.9,
 dass 8.8. nicht ohne Einschr\"ankung an den K\"orper g\"ultig
 ist. Wir werden noch sehen, dass $\TU(V,f)$ auch dann nicht perfekt ist,
 wenn $K = \GF(4)$ und $\Rg_{\GF(4)}(V) = 3$ ist. Dies sind im \"ubrigen
 alle Ausnahmen.
 \medskip\noindent
 {\bf 8.10. Satz.} {\it Es sei $f$ eine nicht ausgeartete, antisymmetrische,
 spurwertige $\alpha$-Form auf dem $K$-Vektorraum $V$. Ferner sei
 $\Rg_K(V) \geq 3$ und der Index von $f$ sei mindestens $1$.
 Ist $P$ ein nicht isotroper Punkt von $\La(V)$, so gibt es zwei nicht
 isotrope Geraden durch $P$, die je zwei isotrope Punkte tragen.}
 \smallskip
       Beweis. Es gibt nach Voraussetzung einen isotropen Punkt $Q$. Es ist
 $Q \neq P$, so dass $P + Q$ eine Gerade ist. Es gibt einen Punkt $R$ mit
 $R \not\leq P + Q$, $P^\bot$, $Q^\bot$. W\"are dies nicht der Fall, so
 enthielten $P + Q$, $P^\bot$ und $Q^\bot$ alle Punkte von $\La(V)$. Dies
 h\"atte $\Rg_K(V) = 3$ und $|K| = 2$ zur Folge. Wegen $|K| = 2$ w\"are
 $\alpha = 1$, so dass $f$ eine symplektische Polarit\"at induzierte.
 Mit Satz 1.7 erhielten wir den Widerspruch, dass $3 = \Rg_K(V)$ gerade
 w\"are. Wir betrachten die Gerade $Q + R$. Wegen $Q \leq Q^\bot$ und
 $R \cap Q^\bot = \{0\}$ ist dann
 $$ (Q + R)^\bot \cap (Q + R) = R^\bot \cap Q^\bot \cap (Q + R) =
     R^\bot \cap Q. $$
 Wegen $R \cap Q^\bot = \{0\}$ ist $Q + R^\bot = V$. Da 
 $Q$ ein Punkt und
 $R^\bot$ eine Hyperebene ist, folgt $Q \cap R^\bot = \{0\}$. Damit ist
 gezeigt, dass $Q + R$ nicht isotrop ist. Weil $Q$ ein isotroper Punkt
 auf $Q + R$ ist, folgt mit 8.3, dass $Q + R$ noch zwei von $Q$ verschiedene
 isotrope Punkte $Q'$ und $Q''$ tr\"agt. Von den drei Punkten $Q$, $Q'$
 und $Q''$ liegen wenigstens zwei nicht in $P^\bot$. Wir d\"urfen annehmen,
 dass dies die Punkte $Q$ und $Q'$ sind. Dann sind aber die Geraden
 $P + Q$ und $P + Q'$ nicht isotrop. Es ist ja
 $$\eqalign{
 (P + Q)^\bot \cap (P + Q) &= P^\bot \cap Q^\bot \cap (P + Q) \cr
      &= P^\bot \cap (Q + (P \cap Q^\bot))                  
      = P^\bot \cap Q = \{0\}. \cr} $$
 Ebenso folgt, dass auch
 $$ (P + Q')^\bot \cap (P + Q') = \{0\} $$
 ist. Da die beiden Geraden $P + Q$ und $P + Q'$ nicht isotrop sind, enthalten
 sie nach 8.3 noch je einen weiteren isotropen Punkt. Damit ist der Satz
 bewiesen.
 \medskip\noindent
 {\bf 8.11. Satz.} {\it Es sei $f$ eine nicht ausgeartete, antisymmetrische,
 spurwertige $\alpha$-Form auf dem $K$-Vektorraum $V$. Ferner sei
 $\Rg_K(V) \geq 2$ und der Index von $f$ sei mindestens $1$. Es sei $\Pi$ die
 Menge der isotropen Punkte von $f$. Ist dann $N$ der Normalteiler von
 $\TU(V,f)$, der aus allen Abbildungen von $\TU(V,f)$ besteht, die alle Punkte
 von $\Pi$ in sich abbilden, so ist}
 $$ Z(\TU(V,f)) = N = \TU(V,f) \cap Z(\GL(V)). $$
 \par
       Beweis. Nach III.1.2 gilt $\TU(V,f) \cap Z(\GL(V)) \subseteq N$.
 Um die umgekehrte Inklusion zu etablieren,
 sei zun\"achst $\Rg_K(V) \geq 3$. Es sei $\nu \in N$ und $P$ sei ein Punkt
 von $\La(V)$. Ist $P \in \Pi$, so ist
 $P^\nu = P$ nach Voraussetzung. Es sei also $P \not\in \Pi$. Nach 8.10
 gibt es dann zwei Geraden $G$ und $H$ durch $P$, die beide zwei isotrope
 Punkte tragen. Es folgt $G^\nu = G$ und $H^\nu = H$ und damit $P^\nu = P$.
 Damit ist gezeigt, dass $\nu$ alle Punkte von $\La(V)$ invariant
 l\"asst. Mit Satz III.1.2 folgt nun die Behauptung in diesem Falle.
 \par
       Es sei $\Rg_K(V) = 2$. Es gibt dann eine Basis $b_1$ und $b_2$ aus
 isotropen Vektoren mit $f(b_1,b_2) = 1$. Ist $\lambda \in K$ und
 $\lambda^\alpha = \lambda$, so wird durch
 $$ v^\tau := v - b_1\lambda f(b_1,v) $$
 eine unit\"are Transvektion definiert, wie wir wissen. Es folgt, dass
 $$ b_2^\tau = b_2 - b_1\lambda f(b_1,b_2) = b_2 - b_1\lambda $$
 ein isotroper Vektor ist.
 \par
       Es sei $\nu \in N$. Es gibt dann $r_1$, $r_2$, $r \in K^*$ mit
 $b_1^\nu = b_1r_1$, $b_2^\nu = b_2r_2$ und $b_2^{\tau\nu} = b_2^\tau r
 = (b_2 - b_1\lambda)r$. Es folgt
 $$ b_2r_2 - b_1r_1\lambda = b_2^\nu - b_1^\nu\lambda = (b_2 - b_1\lambda)^\nu
       = b_2^{\tau\nu} = b_2r - b_1\lambda r $$
 und damit $r_2 = r$ und $r_1\lambda = \lambda r$. Daher ist
 $$ r_1\lambda = \lambda r_2. $$
 Weil $r_1$ und $r_2$ von $\lambda$ unabh\"angig sind, gilt diese Gleichung
 f\"ur alle $\lambda \in K$, f\"ur die $\lambda^\alpha = \lambda$ gilt.
 Mit $\lambda = 1$ folgt $r_2 = r_1 = r$, womit gleichzeitig gezeigt ist,
 dass $r$ nicht von $\lambda$ abh\"angig ist. Wegen $r\lambda = \lambda r$
 f\"ur alle $\lambda$ mit $\lambda^\alpha = \lambda$ ist $r$ nach III.10.11
 ein Element des Zentrums von $K$, es sei denn, es ist $K$ ein
 Quaternionenschiefk\"orper mit von 2 verschiedener Charakteristik und es
 ist
 $$ Z(K) = \{\lambda \mid \lambda \in K, \lambda^\alpha = \lambda\}. $$
 Ist nun $T_1$ die Menge aller unit\"aren Transvektionen der Form
 $v^\rho = v - b_1\lambda f(b_1,v)$ und $T_2$ die Menge aller Transvektionen
 der Form $v^\tau = v - b_2\mu f(b_2,v)$ mit $\lambda$, $\mu \in Z(K)$,
 so folgt mit 8.2, dass die von
 $T_1$ und $T_2$ erzeugte Untergruppe von $\TU(V,f)$ auf $\Pi$ zweifach
 transitiv operiert und folglich alle unit\"aren Transvektionen enth\"alt.
 Also ist sie gleich $\TU(V,f)$. Hieraus folgt, dass die
 $(2 \times 2)$-Matrizen, die Abbildungen aus $T\Un(V,f)$ bez\"uglich der
 Basis $b_1$, $b_2$ darstellen, nur Koeffizienten aus $Z(K)$ haben. Daher
 ist auch $r \in Z(K)$. Damit ist gezeigt, dass
 $$ N = \TU(V,f) \cap Z\bigl(\GL(V)\bigr) \subseteq Z\bigl(\TU(V,f)\bigr) $$
 ist.
 \par
       Es sei $\zeta \in Z(\TU(V,f))$ und $aK$ sei ein isotroper Punkt.
 Dann wird durch $v^\tau := v - af(a,v)$ eine unit\"are Transvektion
 definiert. Es folgt
 $$ v - af(a,v) = v^\tau = v^{\zeta^{-1}\tau\zeta} =
       v - a^\zeta f(a,v^{\zeta^{-1}}). $$
 Hieraus folgt die Existenz eines $k \in K^*$ mit $a^\zeta = ak$. Daher ist
 $$ \zeta \in T\Un(V,f) \cap Z\bigl(\GL(V)\bigr), $$
 so dass $Z(T\Un(V,f)) = N$ ist. Damit ist der Satz in allen seinen Teilen
 bewiesen.
 \medskip\noindent
 {\bf 8.12. Satz.} {\it Es sei $K$ ein K\"orper, dessen Zentrum von $\GF(2)$,
 $\GF(3)$, $\GF(4)$ und $\GF(9)$ verschieden sei. Es sei $f$ eine nicht
 ausgeartete, antisymmetrische, spurwertige $\alpha$-Form auf dem
 $K$-Vektorraum $V$. Ferner sei $\Rg_K(V) \geq 2$ und der Index von $f$
 sei mindestens $1$ und $\Pi$ sei die Menge der isotropen Punkte von $V$. Ist
 dann $N$ der Normalteiler von $\TU(V,f)$, der aus allen Abbildungen von
 $\TU(V,f)$ besteht, die alle Punkte von $\Pi$ in
 sich abbilden, ist ferner $M$ ein Normalteiler von $\TU(V,f)$, der nicht
 in $N$ enthalten ist, so ist $M = \TU(V,f)$. Insbesondere ist die Gruppe
 $$ \PTU(V,f) := \TU(V,f)/N $$
 einfach.}
 \smallskip
       Beweis. Weil der Index von $f$ mindestens gleich 1 ist, gibt es eine
 nicht isotrope Gerade $H$, die mindestens zwei isotrope Punkte tr\"agt.
 Es sei $U$ die Untergruppe von $\TU(V,f)$ die von den Transvektionen von
 $\TU(V,f)$ erzeugt wird, deren Zentren auf $H$ liegen. Jedes Element von $U$
 l\"asst $H^\bot$ vektorweise fest.
 Wie der Beweis von
 Satz 8.8 zeigt --- hier benutzen wir die Voraussetzungen \"uber das Zentrum
 von $K$ ---, ist jede Transvektion aus $U$ ein Kommutator in $U$, so dass
 $U = U'$ ist. Dies werden wir zweimal benutzen.
 \par
       Es sei $U^*$ die Gruppe, die von $U$ auf der Menge der auf $H$
 liegenden Punkte induziert wird.
 Ist dann $f^*$ die Einschr\"ankung von $f$ auf $H$, so folgt mit 8.11,
 dass $U^*$ zu $\PTU(G,f^*)$ isomorph ist. Diese Gruppe operiert auf der
 Menge der
 isotropen Punkte von $H$ zweifach transitiv, also insbesondere primitiv,
 wie wir wissen. Au\ss erdem ist $U^*$ perfekt, da $U$ perfekt ist.
 Da der Stabilisator eines isotropen Punktes von $H$ in $U^*$ den
 abelschen
 Normalteiler enth\"alt, der von den Transvektionen mit diesem Punkt als
 Zentrum induziert wird, und da diese Normalteiler die Gruppen $U^*$
 erzeugen, ist $U^*$ nach dem Satz III.2.1 von Iwasawa einfach.
 \par
       Es sei $M$ ein Normalteiler von $U$, der nicht im Zentrum von $U$
 enthalten ist. Ferner sei $Z$ das Zentrum von $U$. Dann operieren die
 Elemente von $Z$ als Skalarmultiplikationen auf $H$ und als Identit\"at
 auf $H^\bot$. Daher operiert die Einschr\"ankung von $M$ auf die Menge der
 Punkte von $H$ nicht trivial. Weil $U^*$ einfach ist, ist also $M^* = U^*$.
 Es folgt $MZ = U$. Weiter folgt
 $$ U/M = (MZ)/M \cong Z/(M \cap Z), $$
 so dass $U/M$ abelsch ist. Hieraus folgt
 $ U = U' \subseteq M $
 und damit $U = M$. --- Diesem Argument sind wir in Kapitel III schon einmal
 begegnet.
 \par
       Es sei nun $M$ ein Normalteiler von $\TU(V,f)$, der nicht im Zentrum
 von $\TU(V,f)$ enthalten ist. Nach Satz 8.8 gibt es dann einen isotropen
 Punkt $P$ und ein $\sigma \in N$ mit $P \neq P^\sigma$.  Wir zeigen, dass
 es sogar einen isotropen Punkt $Q$ gibt mit $Q^\sigma \not\leq Q^\bot$.
 Dazu d\"urfen wir annehmen, dass $P \leq P^\sigma$ ist. Es sei $P = vK$.
 Dann ist also $f(v,v^\sigma) = 0$. Wegen $P \neq P^\sigma$ ist
 $$ P^\bot \neq (P^\sigma)^\bot. $$
 Es gibt also ein $z \in (P^\sigma)^\bot - P^\bot$. Setze $H := zK + P$.
 Wegen $P \leq P^\bot$ und $zK \cap P^\bot = \{0\}$ folgt
 $$\eqalign{
 H \cap H^\bot &= (zK + P) \cap P^\bot \cap (zK)^\bot  \cr
       &= \bigl((zK \cap P^\bot) + P\bigr)\cap (zK)^\bot         \cr
       &= P \cap (zK)^\bot.                             \cr} $$
 Aus $zK \cap P^\bot = \{0\}$ folgt auch $(zK)^\bot + P = V$. Da $(zK)^\bot$
 eine Hyperebene ist, folgt schlie\ss lich $(zK)^\bot \cap P = \{0\}$.
 Also ist $H \cap H^\bot = \{0\}$, so dass $H$ nicht isotrop ist. Daher gibt
 es neben $P$ noch zwei weitere isotrope Punkte $yK$ und $y'K$ auf $H$.
 Es gibt ferner eine Transvektion $\tau$ mit Zentrum $y'K$, die $vK$ auf
 $yK$ abbildet. Es gibt also ein $l \in K$ mit $v^\sigma = yl$.
 \par
       Es gibt $a$, $b \in K$ mit $y' = va + zb$. Es folgt
 $$ f(y',v^\sigma) = f(va + zb,v^\sigma)
       = a^\alpha f(v,v^\sigma) + b^\alpha f(z,v^\sigma)
       = 0, $$
 da ja $f(v,v^\sigma) = 0 = f(z,v^\sigma)$ ist. Dies zeigt, dass $v^\sigma$
 in der Achse $(y'K)^\bot$ von $\tau$ liegt. Daher ist
 $v^{\sigma\tau} = v^\sigma$. Es folgt
 $$ v^{\sigma\tau^{-1}\sigma^{-1}\tau} = v^{\sigma\sigma^{-1}\tau}
       = v^\tau = yl. $$
 Nun ist aber $\sigma(\tau^{-1}\sigma^{-1}\tau) \in N$ und
 $f(v,yl) \neq 0$. Letzteres, weil $H = vK + yK$ ist und weil $H$ nicht
 isotrop ist. Damit ist gezeigt, dass es eine Abbildung in $G$ gibt, die
 wir wieder $\sigma$ nennen, so dass $P^\sigma \not\leq P^\bot$ gilt.
 \par
       Es sei $\tau$ eine von $1_V$ verschiedene Transvektion mit dem
 Zentrum $P$. Dann ist $\sigma^{-1}\tau\sigma$ eine Transvektion mit dem
 Zentrum $P^\sigma$. Wegen $P \neq P^\sigma$ und $P \not\leq (P^\sigma)^\bot$
 sind $\tau$ und $\sigma^{-1}\tau\sigma$ nicht vertauschbar. Dann sind
 auch $\tau^{-1}$ und $\sigma^{-1}\tau$ nicht vertauschbar. Daher ist
 $$ \tau^{-1}\sigma^{-1}\tau\sigma $$
 ein Element der Gruppe $U$, die von den Elationen mit einem Zentrum auf
 $H$ erzeugt wird, welches nicht im Zentrum von $U$ liegt. Also ist
 $M \cap U$ ein Normalteiler von $U$, der nicht im Zentrum von $U$ liegt.
 Nach der Eingangs gemachten Bemerkung ist folglich
 $$ U = M \cap U \subseteq M. $$
 Daher enth\"alt $M$ alle Transvektionen mit dem Zentrum $P$ und wegen der
 Transitivit\"at von $\TU(V,f)$ auf der Menge aller isotropen Punkte
 \"uberhaupt alle unit\"aren Transvektionen. Folglich ist $M = \TU(V,f)$.
 Hieraus folgt wiederum, dass $\PTU(V,f)$ einfach ist.
 \medskip
       Wie schon gesagt, gibt es nicht kommutative K\"orper, deren Zentren
 zu $\GF(2)$, $\GF(3)$, $\GF(4)$ oder $\GF(9)$ isomorph sind. Ob es auch
 solche gibt, die einen involutorischen Antiautomorphismus besitzen, wei\ss\
 ich nicht. Wenn ja, so ist die $\PTU(V,f)$ auch in diesen F\"allen einfach,
 wie Dieudonn\'e zeigte. Wir werden dies hier aber nicht weiter verfolgen.
 Ist $K$ kommutativ, so kommen die F\"alle $\GF(2)$ und $\GF(3)$ nicht vor,
 da diese beiden K\"orper keinen involutorischen Antiautomorphismus
 besitzen. Es bleiben also nur noch die unit\"aren Gruppen \"uber $\GF(4)$ und
 $\GF(9)$ zu untersuchen. Dies werden wir in einem etwas allgemeineren Rahmen
 im \"ubern\"achsten Abschnitt tun.

 \mysection{9. Endliche unit\"are Gruppen}
 
\noindent
       Unser erstes Ziel in diesem Abschnitt ist, die Ordnung der endlichen
 unit\"aren Gruppen zu bestimmen. Anschlie\ss end werden wir die Geometrie der
 vollst\"andig isotropen Unterr\"aume einer unit\"aren Polarit\"at eines endlichen
 projektiven Raumes studieren, um im letzten Abschnitt dann die Einfachheit
 der $\PTU(V,f)$ auch in den F\"allen, dass der zugrunde liegende K\"orper
 $\GF(4)$ oder $\GF(9)$ ist, zu beweisen, wobei es, wie schon erw\"ahnt, drei
 Ausnahmen gibt, von denen wir zwei bereits kennen und die dritte in diesem
 Abschnitt kennen lernen werden.
 \par
       Ist $K$ ein kommutativer K\"orper und ist $\alpha$ ein involutorischer
 Automorphismus von $K$, so ist $K$ eine quadratische Erweiterung des
 Fixk\"orpers von $\alpha$. Ist $K$ endlich, so ist also $K = \GF(q^2)$. Da
 die Automorphismengruppe eines endlichen K\"orpers zyklisch ist, hat
 $\GF(q^2)$ nur einen involutorischen Automorphismus $\alpha$. F\"ur diesen
 gilt $k^\alpha = k^q$ f\"ur alle $k \in K$.
 \medskip\noindent
 {\bf 9.1. Satz.} {\it Es sei $V$ ein Vektorraum des Ranges $n \geq 1$ \"uber
 $\GF(q^2)$ und $\alpha$ sei der involutorische Automorphismus von $\GF(q^2)$.
 Ist dann $f$ eine nicht ausgeartete symmetrische $\alpha$-Form auf $V$,
 so gibt es eine Basis $b_1$, \dots, $b_n$ von $V$ mit
 $$\textstyle f\bigl(\sum_{i:=1}^n b_ik_i,\sum_{j:=1}^n b_jl_j\bigr)
               = \sum_{i:=1}^n k_i^\alpha l_i $$
 f\"ur alle $k_i$, $l_j \in \GF(q^2)$.}
 \smallskip
 Jede solche Basis hei\ss t {\it Orthonormalbasis\/}.\index{Orthonormalbasis}{}
 \smallskip
       Beweis. Weil $f$ nicht ausgeartet ist, gibt es ein $v \in V$ mit
 $f(v,v) \neq 0$. Wegen $f(v,v) = f(v,v)^\alpha$ ist $f(v,v) \in \GF(q)$.
 Weil $f(v,v) \neq 0$ und weil die multiplikative Gruppe von $\GF(q^2)$
 zyklisch ist, gibt es ein $k \in \GF(q^2)^*$ mit
 $$ k^{\alpha+1} = k^{q+1} = f(v,v). $$
 Setze $b_1 := vk^{-1}$. Dann ist
 $$ f(b_1,b_1) = k^{-\alpha} f(v,v)k^{-1}
       = k^{-\alpha}k^{\alpha+1}k^{-1} = 1. $$
 Es folgt $V = b_1K \oplus (b_1K)^\bot$. Dies impliziert wiederum, dass die
 Einschr\"ankung von $f$ auf $(b_1K)^\bot$ nicht entartet ist. Daher gibt es
 eine Basis $b_2$, \dots, $b_n$ von $(b_1K)^\bot$ mit $f(b_i,b_i) = 1$
 f\"ur $i := 2$, \dots, $n$ und $f(b_i,b_j) = 0$ f\"ur $i$, $j := 1$, \dots,
 $n$ und $i \neq j$. Damit ist $9.1$ bewiesen.
 \medskip
       Dieser Satz ist grundlegend f\"ur alles Folgende. Er besagt unter
 anderem, dass ein endlicher projektiver Raum \"uber $\GF(q^2)$ im
 We\-sent\-li\-chen nur eine unit\"are Polarit\"at besitzt. Um 9.1 ausnutzen zu
 k\"onnen, beweisen wir zun\"achst
 \medskip\noindent
 {\bf 9.2. Satz.} {\it Es sei $0 \neq k \in \GF(q)$.
 Ist $A_{n,q}$ die Anzahl der $n$-Tupel $y_1$, \dots, $y_n$
 mit $y_i \in \GF(q^2)$ f\"ur alle $i$, die die Gleichung
 $$ \sum_{i:=1}^n y_i^{q + 1} = k $$
 erf\"ullen, so ist
 $$ A_{n,q} = q^{2n-1} - (-1)^nq^{n-1}. $$
 Insbesondere ist $A_{n,q}$ von $k$ unabh\"angig.}
 \smallskip
       Beweis. Die multiplikative Gruppe von $\GF(q^2)$ ist zyklisch und hat die
 Ordnung $q^2 - 1$. Daher ist das Potenzieren mit $q + 1$ ein Epimorphismus
 dieser Gruppe auf die multiplikative Gruppe von $\GF(q)$. Der Kern dieses
 Epimorphismus hat die Ordnung $q + 1$. Folglich hat jedes von $0$
 verschiedene Element von $\GF(q)$ genau $q + 1$ Urbilder unter diesem
 Epimorphismus, w\"ahrend die Gleichung $k^{q+1} = 0$ genau eine L\"osung
 hat, n\"amlich $k = 0$.
 \par
       Es sei $n > 1$. Ist $\sum_{i:=1}^{n-1} y_i^{q+1} = k$,
 so ist $y_n = 0$. Daher ist die Anzahl der $n$-Tupel dieser Art gleich
 $$ A_{n-1,q}. $$
 Ist $\sum_{i:=1}^{n-1} y_i^{q+1} \neq k$,
 so gibt es nach der Eingangs gemachten Bemerkung genau $q + 1$ Werte $y_n$
 mit $\sum_{i:=1}^n y_i^{q+1} = k$, da ja
 $k - \sum_{i:=1}^{n-1} y_i^{q+1} \in \GF(q)$
 gilt. Da es von dieser Sorte $(n - 1)$-Tupel genau
 $$ q^{2(n-1)} - A_{n-1,q} $$
 St\"uck gibt, ist
 $$ 
 A_{n,q} = A_{n-1,q} + (q + 1)(q^{2(n-1)} - A_{n-1,q}) 
       = q^{2n-1} + q^{2(n-1)} - qA_{n-1,q}.  $$
 \par
       Es sei nun $n = 1$. Nach der eingangs gemachten Bemerkung hat die
 Gleichung $y_1^{q+1} = k$ genau $q + 1$ L\"osungen, so dass der Satz in
 diesem Falle korrekt ist. Es sei $n > 1$ und der Satz gelte f\"ur $n - 1$.
 Dann ist
 $$
 A_{n,q} = q^{2n-1} + q^{2(n-1)} - q\bigl(q^{2n-3} - (-1)^{n-1}q^{n-2}\bigr) 
         = q^{2n-1} - (-1)^nq^{n-1},  $$
 q. e. d.
 \medskip
       Es sei $K$ ein K\"orper und $\alpha$ sei ein involutorischer
 Antiautomorphismus von $K$. Ferner sei $V$ ein $K$-Vektorraum und $f$ sei
 eine nicht ausgeartete, spurwertige, symmetrische $\alpha$-Form auf $V$.
 Wir setzen $\SU(V,f) := U \cap \SL(V)$ und bezeichnen mit $\PGU(V,f)$ und
 $\PSU(V,f)$ die von $\Un(V,f)$ bzw. $\SU(V,f)$ auf $\La(V)$ induzierten
 Kol\-li\-ne\-a\-ti\-ons\-grup\-pen. Weil $\SL(V)$ alle Transvektionen
 enth\"alt, gilt $$ \TU(V,f) \subseteq \SU(V,f) $$
 und dann auch $\PTU(V,f) \subseteq \PSU(V,f)$.
 Gleichheit gilt sicher dann nicht, wenn der Index von $f$ Null ist. Sie
 gilt aber auch im Falle, dass der Index mindestens 1 ist nicht immer. Ist
 $K$ kommutativ, so gibt es aber nur eine Ausnahme, wie wir noch sehen werden.
 \par
       Ist $K = \GF(q^2)$, so gibt es nach Satz 9.1 bis auf Isometrie nur
 eine unit\"are Polarit\"at. Daher bezeichnen wir in diesem Falle die Gruppen
 $\Un(V,f)$ etc. auch mit $\Un(n,q^2)$ etc., wenn $n$ der Rang von $V$ ist.
 \medskip\noindent
 {\bf 9.3. Satz.} {\it Es sei $n$ eine nat\"urliche Zahl und $q$ sei Potenz einer
 Primzahl $p$. Dann gilt:
 \smallskip
 \item{a)} Es ist
 $$ \big|\Un(n,q^2)\big| = \prod_{i:=1}^n \bigl(q^i - (-1)^i\bigr)q^{i-1}. $$
 \smallskip
 \item{b)} Es ist
 $$ \big|\SU(n,q^2)\big| = \prod_{i:=2}^n \bigl(q^i - (-1)^i\bigr)q^{i-1}. $$
 \smallskip
 \item{c)} Es ist
 $$ \big|\PGU(n,q^2)\big| = \prod_{i:=2}^n \bigl(q^i - (-1)^i\bigr)q^{i-1}. $$
 \smallskip
 \item{d)} Es sei $n \geq 2$. Dann ist
 $$\big|\PSU(n,q^2)\big| = {1 \over \ggT(n,q+1)}\prod_{i:=2}^n \bigl(q^i -
	      (-1)^i\bigl)q^{i-1}. $$
 \par\noindent
       Die Ordnung einer $p$-Sylowgruppe irgendeiner dieser Gruppen ist
 $q^{{1 \over 2}n(n-1)}$.}
 \smallskip
       Beweis. a) Es sei $V$ ein Vektorraum des Ranges $n$ \"uber $\GF(q^2)$,
 es sei $\alpha$ der involutorische Automorphismus von $\GF(q^2)$ und $f$
 sei eine nicht ausgeartete, symmetrische $\alpha$-Form auf $V$. Ferner sei
 $b_1$, \dots, $b_n$ eine nach 9.1 existierende Orthonormalbasis. Ist dann
 $v = \sum_{i:=1}^n b_ik_i$, so ist genau dann $f(v,v) = 1$, wenn
 $\sum_{i:=1}^n k^{q+1} = 1$ ist. Mit 9.2 folgt, dass
 $$ A_{n,q} = \bigl(q^n - (-1)^n\bigr)q^{n-1} $$
 die Anzahl dieser Vektoren ist. Aufgrund des Satzes von Witt (Satz 7.2) ist
 $\Un(n,q^2)$ auf der Menge der Vektoren der L\"ange 1 transitiv. Daher ist
 $$ \big|\Un(n,q^2)\big| = A_{n,q}\big|\Un(n,q^2)_v\big|, $$
 wenn $v$ ein Vektor der L\"ange 1 ist. Die Gruppe $\Un(n,q^2)_v$ operiert
 treu auf $(v\GF(q^2))^\bot$, weil sie den Punkt $v\GF(q^2)$ vektorweise
 festl\"asst und weil $V = v\GF(q^2) \oplus (v\GF(q^2))^\bot$.
 Da die Einschr\"ankung von $f$ auf $(v\GF(q^2))^\bot$ nicht
 ausgeartet ist, induziert $\Un(n,q^2)_v$ auf $(v\GF(q^2))^\bot$ eine
 Untergruppe von $\Un(n-1,q^2)$, die wegen des Satzes von Witt sogar gleich
 $\Un(n-1,q^2)$ ist. Daher ist
 $$ \big|\Un(n,q^2)\big| = A_{n,q}\big|\Un(n-1,q^2)\big|. $$
 Hieraus folgt mit Induktion die Aussage a).
 \par
       b) Es sei $\sigma \in \Un(n,q^2)$. Ferner sei $a$ die Matrix, die
 $\sigma$ bez\"uglich der Orthonormalbasis $b_1$, \dots, $b_n$ darstellt.
 Aus $f(b^\sigma_i,b^\sigma_j) = f(b_i,b_j)$ folgt dann, dass
 $$ (a^{-1})_{ij} = a_{ji}^q $$
 ist. Hieraus folgt weiter, dass
 $$ 1 = \det(a)^{q+1} $$
 ist. Ist andererseits $k \in \GF(q^2)$, gilt $k^{q+1} = 1$ und definiert
 man $\sigma$ durch $b_1^\sigma := b_1k$ und $b_i^\sigma := b_i$ f\"ur
 $i \geq 2$, so ist $\sigma \in \Un(n,q^2)$ und es gilt $\det(\sigma) = k$.
 Also ist $\det$ ein Epimorphismus von $\Un(n,q^2)$ auf die Gruppe
 $$ \bigl\{k \mid k \in \GF(q^2), k^{q+1} = 1\bigr\}. $$
 Da diese Gruppe die Ordnung $q + 1$ hat und $\SU(n,q^2)$ der Kern von $\det$
 ist, gilt auch b).
 \par
       c) F\"ur $n = 1$ ist nichts zu beweisen. Es sei also $n \geq 2$.
 Induziert dann $\delta \in \Un(n,q^2)$ auf $\La(V)$ die Identit\"at, so
 gibt es ein $k \in \GF(q^2)$ mit $v^\delta = vk$. Es folgt $k^{q+1} = 1$.
 Umgekehrt definiert jedes solche $k$ ein $\delta \in \Un(n,q^2)$, welches
 auf $\La(V)$ die Identit\"at induziert. Weil es genau $q + 1$ solcher
 $k$ gibt, gilt auch c).
 \par
       d) Es sei $\delta \in \SU(n,q^2)$ und $\delta$ induziere auf $\La(V)$
 die Identit\"at. Dann ist einmal $v^\delta = vk$ mit einem $k \in \GF(q^2)$,
 f\"ur das $k^{q+1} = 1$ gilt. Andererseits ist $1 = \det(\delta) = k^n$.
 Folglich ist $k^{\ggT(n,q+1)} = 1$. Gilt umgekehrt diese Gleichung und
 definiert man $\delta$ durch $v^\delta := vk$, so ist $\delta \in \SU(n,q^2)$
 und $\delta$ induziert die Identit\"at in $\La(V)$. Also gilt auch d).
 \par
       Die letzte Aussage ist banal. Damit ist alles bewiesen.
 \medskip\noindent
 {\bf 9.4. Satz.} {\it Es sei $V$ ein Vektorraum des Ranges $n \geq 2$ \"uber
 $\GF(q^2)$ und $\alpha$ sei der involutorische Automorphismus von $\GF(q^2)$.
 Ist $f$ eine nicht entartete, symmetrische $\alpha$-Form auf $V$ und
 bezeichnet $T^1_{n,q}$ die Anzahl der bez\"uglich $f$ isotropen Punkte von
 $\La(V)$, so gilt:
 \item{a)} Es ist $T^1_{n,q} = T^1_{n-1,q} + A_{n-1,q}$, wobei $A$ die in
 Satz 9.2 definierte Matrix ist.
 \item{b)} Es ist}
 $$ T^1_{n,q} = \frac 1 {q^2-1} (q^n - (-1)^n)(q^{n-1} - (-1)^{n-1}) 
 . $$
 \par
       Beweis. a) Nach 9.1 gibt es eine Basis $b_1$, \dots, $b_n$ mit
 $$\textstyle f\bigl(\sum_{i:=1}^n b_ix_i,\sum_{j:=1}^n b_jy_j\bigr) =
       \sum_{i:=1}^n x^q_iy_i. $$
 Es sei
 $$\textstyle P = \bigl(\sum_{i:=1}^n b_ix_i\bigr)\GF(q^2) $$
 ein Punkt von $V$. Genau dann ist $P$ isotrop, wenn
 $$\textstyle \sum_{i:=1}^n x_i^{q+1} = 0 $$
 ist. Die Punkte mit $x_n = 0$ f\"ullen eine nicht isotrope Hyperebene $H$.
 Weil $H$ nicht isotrop ist, ist $T^1_{n-1,q}$ die Anzahl der auf $H$
 liegenden isotropen Punkte bez\"uglich $f$.
 \par
       Es sei $P$ ein isotroper Punkt, der nicht auf $H$ liegt. Dann ist
 $x_n \neq 0$. Wir setzen $y_i := x_ix_n^{-1}$. Dann ist
 $$ P = (b_1y_1 + \dots + b_{n-1}y_{n-1} + b_n)\GF(q^2). $$
 Es folgt
 $$\textstyle \sum_{i:=1}^{n-1} y_i^{q+1} = -1. $$
 Umgekehrt definiert jedes $(n - 1)$-Tupel $y$, welches diese Gleichung
 erf\"ullt, einen isotropen Punkt und verschiedene solche $(n - 1)$-Tupel
 definieren verschiedene Punkte. Mit 9.2 folgt daher, dass $A_{n-1,q}$
 die Anzahl der bez\"uglich $f$ isotropen Punkte ist, die nicht auf $H$
 liegen. Damit ist a) bewiesen.
 \par
       b) Weil $f$ eine unit\"are Polarit\"at darstellt und unit\"are
 Polarit\"aten nach Satz 1.9 sich auch durch schiefsymmetrische
 $\alpha$-Semibilinearfor\-men darstellen lassen, folgt aus Satz 8.3,
 dass $T^1_{2,q} = q + 1$ ist. Daher gilt b) f\"ur $n = 2$. Der Rest
 folgt mit einer einfachen Induktion unter Zuhilfenahme von a).
 \medskip
       Ist $V$ ein Vektorraum des Ranges $n$ \"uber $\GF(q^2)$, ist
 $\alpha$ der involutorische Automorphismus von $\GF(q^2)$ und ist $f$ eine
 nicht ausgeartete, symmetrische $\alpha$-Semibilinearform auf $V$,
 so bezeichnen wir mit $T^r_{n,q}$ die Anzahl der bez\"uglich $f$
 vollst\"andig isotropen Teilr\"aume des Ranges $r$ von $V$.
 Auf Grund von Satz 9.1 ist $T^r_{n,q}$
 nicht von der Wahl von $f$ abh\"angig. Mittels der Bemerkung nach 6.3 folgt
 $T^r_{n,q} = 0$, falls $r > {n \over 2}$ ist.
 \medskip\noindent
 {\bf 9.5. Satz.} {\it Es sei $V$ ein Vektorraum des Ranges $n \geq 2$ \"uber
 $\GF(q^2)$ und $\alpha$ sei der involutorische Automorphismus dieses
 K\"orpers.
 Ist dann $f$ eine nicht ausgeartete symmetrische $\alpha$-Form auf $V$,
 so gilt f\"ur die Anzahlen der bez\"uglich $f$ vollst\"andig isotropen Teilr\"aume:
 \item{a)} Ist $1 \leq s \leq r \leq {n \over 2}$, so ist
 $$ T^s_{n,q}T^{r-s}_{n-2s,q} = T^r_{n,q}{r \choose s,q^2}. $$
 \item{b)} Ist $1 \leq r \leq {n \over 2}$, so ist
 $$ T^r_{n,q} = {\prod_{i:=n+1-2r}^n (q^i - (-1)^i) \over
                     \prod_{j:=1}^r (q^{2j} - 1)}. $$}
 \par
       Beweis. a) Wir betrachten die folgende Inzidenzstruktur: die Punkte sind die
 vollst\"andig isotropen Unterr\"aume des Ranges $s$ und die Bl\"ocke sind die vollst\"andig
 isotropen Unterr\"aume des Ranges $r$, wobei die Inzidenz mit der
 in $\La(V)$ gegebenen Relation $\leq$ identisch sei. Dann ist
 $T^s_{n,q}$ die Anzahl der Punkte und $T^r_{n,q}$ die Anzahl der
 Bl\"ocke dieser Inzidenzstruktur.
 Ferner ist $r \choose s,q^2$ die Anzahl der Punkte pro Block, da ja
 jeder Teilraum eines vollst\"andig isotropen Teilraums vollst\"andig isotrop ist. Ist
 $U$ ein vollst\"andig isotroper Teilraum des Ranges $s$, so induziert $f$ wegen
 $\bot^2 = 1_{\La(V)}$ auf $U^\bot/U$ eine nicht ausgeartete symmetrische
 $\alpha$-Semibilinearform. Ferner gilt, dass jeder vollst\"andig isotrope
 Teilraum,
 der $U$ umfasst, in $U^\bot$ enthalten ist. Hieraus folgt, dass die
 Anzahl der vollst\"andig isotropen Teilr\"aume des Ranges $r$ von $V$, die $U$
 umfassen, gleich der Anzahl der vollst\"andig isotropen Teilr\"aume des Ranges
 $r - s$ von $U^\bot/U$, dh. gleich $T^{r-s}_{n-2s,q}$ ist.
 Mit IIII.1.2c) folgt nun die Behauptung.
 \par
       b) Mit $s = 1$ folgt aus a) mit I.7.6 und 9.3b)
 $$\eqalign{
 {q^{2r}- 1 \over q^2 - 1}T^r_{n,q} &= {r \choose 1,q^2}T^r_{n,q}
       = T^1_{n,q}T^{r-1}_{n-2,q}                                 \cr
       &= \frac 1 {q^2-1} (q^n - (-1)^n)(q^{n-1} - (-1)^{n-1}) 
                     T^{r-1}_{n-2,q}. \cr} $$
 Hieraus folgt per Induktion die Behauptung.
 \medskip
       Ist $n = 3$, so gibt es auf Grund der Bemerkung nach 6.3 keine
 vollst\"andig isotropen Geraden. Dies macht die Menge der isotropen Punkte
 in diesem Falle besonders homogen und damit besonders interessant. Diese
 Sonderrolle der projektiven Ebenen werden wir nun --- vom Hauptweg unserer
 Untersuchungen abweichend --- etwas eingehender betrachten. Dabei spielt die
 Endlichkeit nur bei den Anzahlbestimmungen eine Rolle.
 \par
       Es sei $K$ ein K\"orper und $\alpha$ sei ein
 involutorischer Antiautomorphismus von $K$. Ferner sei $V$ ein Vektorraum
 des Ranges 3 \"uber $K$ und $f$ sei eine nicht ausgeartete, symmetrische
 $\alpha$-Semibilinearform auf $V$, deren Index mindestens 1 sei. Wie wir
 wissen, ist die
 letzte Bedingung \"uberfl\"ussig, da automatisch erf\"ullt, wenn $K$ endlich
 ist. Mit $\U(V,f)$ bezeichnen wir die Geometrie aus den bez\"uglich $f$
 isotropen Punkten und den nicht isotropen Geraden, die wenigstens einen
 isotropen Punkt tragen. Man nennt $\U(V,f)$ {\it Unital\/}.\index{Unital}{}
 \medskip\noindent
 {\bf 9.6. Satz.} {\it Es sei $K$ ein K\"orper und $\alpha$ sei ein
 involutorischer
 Antiautomorphismus von $K$. Ferner sei $V$ ein Vektorraum des Ranges $3$
 \"uber $K$ und $f$ sei eine nicht ausgeartete, symmetrische, spurwertige
 $\alpha$-Form auf $V$, deren Index $1$ sei. Dann gilt: 
 \item{a)} Durch zwei verschiedene Punkte von $\U(V,f)$ geht genau eine
 Gerade von $\U(V,f)$.
 \item{b)} Ist $P$ ein Punkt von $\U(V,f)$ und ist $G$ eine von $G^\bot$
 verschiedene Gerade durch $P$, so ist $G$ eine Gerade von $\U(V,f)$.
 \item{c)} Ist $K$ endlich, also $K = \GF(q^2)$, so ist $\U(V,f)$ ein
 $2$-$(q^3 + 1,q + 1,1)$ Blockplan.
 \item{d)} Die Gruppe $\PGU(V,f)$ operiert zweifach transitiv auf der Menge der
 Punkte von $\U(V,f)$.\par}
 \smallskip
       Beweis. a) Es seien $P$ und $Q$ zwei verschiedene Punkte von $\U(V,f)$.
 Weil $\bot$ eine Polarit\"at ist, folgt mit Satz II.4.1, dass $P + Q$
 keine absolute Gerade ist. Da es keine vollst\"andig isotropen Ge\-ra\-den gibt, ist
 $P + Q$ also nicht isotrop und folglich eine Gerade von $\U(V,f)$. Durch zwei
 verschiedene Punkte von $\U(V,f)$ geht also mindestens und damit genau eine
 Gerade von $\U(V,f)$.
 \par
       b) Nach dem zu II.4.1 dualen Satz ist $P^\bot$ die einzige absolute
 Gerade durch $P$. Folglich ist $G$ nicht isotrop und damit eine Gerade von
 $\U(V,f)$.
 \par
       c) Dies folgt mit a) und Satz 9.3 b).
 \par
       d) Ist $G$ eine Gerade von $\U(V,f)$, so gibt es einen und nach Satz
 8.3 dann mindestens drei isotrope Punkte auf $G$. Dann gibt es aber eine
 Basis $a$, $b$ von $G$ mit isotropen Vektoren $a$ und $b$, so dass
 $f(a,b) = 1$ ist. Dies zeigt, dass alle Geraden von $\U(V,f)$ isometrisch
 sind. Mit dem Satz von Witt (Satz 7.2) folgt daher, dass $\PGU(V,f)$
 auf der Menge der Geraden von $\U(V,f)$ transitiv operiert. Ist $G$ wieder
 eine Gerade von $\U(V,f)$, so haben wir fr\"uher schon gesehen, dass die
 Gruppe, die von allen Transvektionen mit Zentren auf $G$ erzeugt wird, auf
 der Menge der Punkte von $G$ zweifach transitiv operiert. Mit a) folgt
 daher auch die letzte Behauptung.
 \medskip
       Ist $K = \GF(q^2)$, so ist $\U(V,f)$ also ein 2-$(q^3 + 1,q + 1,1)$
 Blockplan. Mit b) etwa folgt, dass die Anzahl der Geraden von $\U(V,f)$
 durch eine Punkt dieser Geometrie gleich $q^2$ ist. Es folgt, dass die
 Anzahl der Geraden von $\U(V,f)$ gleich $q^2(q^2 - q + 1)$ ist. Also tr\"agt
 jede Gerade von $\La(V)$ einen Punkt von $\U(V,f)$.
 \par
       Ist $V$ ein Vektorraum des Ranges 3 \"uber $\GF(4)$, ist $\alpha$ der
 involutorische Automorphismus von $\GF(4)$ und ist $f$ eine nicht entartete,
 symmetrische $\alpha$-Form auf $V$, so ist $\U(V,f)$ nach Satz 9.5c) ein
 2-$(9,3,1)$ Blockplan, so dass $\U(V,f)$ nichts anderes ist als eine affine
 Ebene der Ordnung 3. Diese Situation, n\"amlich die Einbettung der affinen
 Ebene der Ordnung 3 in eine projektive Ebene, haben wir in Abschnitt III.7
 ausf\"uhrlich studiert. Es steht zu erwarten, dass die hessesche Gruppe
 im vorliegenden Falle als die Gruppe $\PGU(V,f)$ zu interpretieren ist. Dies
 ist tats\"achlich der Fall, wie wir uns nun \"uberlegen werden.
 \par
       Dazu sei $P$ ein Punkt von $\La(V)$. Ist $P$ ein Punkt von $\U(V,f)$,
 so ist $P^\bot$ die einzige Gerade durch $P$, die mit $\U(V,f)$ nur den
 Punkt $P$ gemeinsam hat, wie 9.5b) sagt. Ist $P$ kein Punkt von $\U(V,f)$,
 so ist $P \not\leq P^\bot$. Daher ist  $P^\bot$ eine Gerade von
 $\U(V,f)$. Ist $Q$ ein Punkt von $\U(V,f)$, der auf $P^\bot$ liegt, so
 ist $Q^\bot = P + Q$. Es gibt also drei Geraden durch $P$, die $\U(V,f)$
 in genau einem Punkt treffen und diese drei Punkte sind kollinear. Die
 restlichen zwei Geraden durch $P$ sind nicht isotrop und sind daher Geraden
 von $\U(V,f)$. Diese beiden Geraden und die Gerade $P^\bot$ sind die
 drei Geraden einer Parallelenschar der affinen Ebene $\U(V,f)$. Dies zeigt,
 dass $\bot$ vollst\"andig durch die Einbettung von $\U(V,f)$ in $\La(V)$
 beschrieben wird. Daher ist die hessesche Gruppe im Zentralisator von
 $\bot$ enthalten. Da die hessesche Gruppe der Stabilisator von $\La(V)$
 in $\PGL(V)$ ist, ist sie sogar gleich $\PGU(V,f)$. Nach III.7.8 enth\"alt
 die hessesche Gruppe den Normalteiler $T$ der Translationen der affinen
 Ebene $\U(V,f)$. Ist $\sigma$ eine Elation ungleich 1 aus $\PGU(V,f)$, so
 enth\"alt die Gruppe $T\{1,\sigma\}$ alle Elationen aus $\PGU(V,f)$, da $T$
 auf der Menge der Punkte von $\U(V,f)$ transitiv ist.
 Also ist
 $$ \PTU(V,f) = T\{1,\sigma\}. $$
 Somit ist $\PTU(V,f)$ auch in diesem Falle nicht einfach. \"Uber\-dies gilt hier
 $$ \big|\PTU(V,f)\big| = 9 \cdot 2 \neq 9 \cdot 8 = \big|\PSU(V,f)\big|, $$
 so dass $\PTU(V,f) \neq \PSU(V,f)$ ist.
 \medskip\noindent
 {\bf 9.7. Satz.} {\it Es sei $V$ ein Vektorraum des Ranges $n \geq 2$ \"uber
 $\GF(q^2)$ und $\alpha$ sei der involutorische Automorphismus von $\GF(q^2)$.
 Ferner sei $f$ eine nicht entartete, symmetrische $\alpha$-Form auf $V$.
 Ist $P$ ein isotroper Punkt, so bezeichnen wir mit $O(P)$ die Menge der von
 $P$ verschiedenen, zu
 $P$ orthogonalen und mit $N(P)$ die Menge der zu $P$ nicht orthogonalen
 isotropen Punkte von $V$. Dann ist}
 $$ \big|O(P)\big| = q^2T^1_{n-2,q} \qquad
 \mathit{und}\qquad
  \big|N(P)\big| = q^{2n-3}. $$
 \par
       Beweis. Es ist $O(P)$ die Menge der von $P$ verschiedenen isotropen
 Punkte auf $P^\bot$. Ist $Q \in O(P)$, so ist $P$, $Q \leq Q^\bot$ und
 $P$, $Q \leq P^\bot$. Also ist
 $$ P + Q \leq P^\bot \cap Q^\bot = (P + Q)^\bot. $$
 Folglich ist $P + Q$ vollst\"andig isotrop. Ist andererseits $G$ eine
 voll\-st\"an\-dig
 isotrope Gerade durch $P$, so ist $G \leq P^\bot$ und alle von $P$
 verschiedenen Punkte auf $G$ geh\"oren zu $O(P)$. Da die vollst\"andig
 isotropen
 Geraden durch $P$ den isotropen Punkten von $P^\bot/P$ ent\-spre\-chen, gilt
 also
 $$ \big|O(P)\big| = q^2T^1_{n-2,q}. $$
 Hieraus folgt mit 9.5 weiter, dass die Anzahl der Punkte in $N(P)$ gleich
 $$ T^1_{n,q} - q^2T^1_{n-2,q} - 1 = q^{2n - 3} $$
 ist.
 \medskip\noindent
 {\bf 9.8. Satz.} {\it Es sei $V$ ein Vektorraum des Ranges $n \geq 2$ \"uber
 dem kommutativen K\"orper $K$ und $\alpha$ sei ein involutorischer
 Automorphismus von $K$.
 Ferner sei $f$ eine nicht entartete, symmetrische, spurwertige
 $\alpha$-Form auf $V$, deren Index mindestens $1$ sei. Ist $P$ ein isotroper
 Punkt, so bezeichnen wir mit $O(P)$ die Menge der von $P$ verschiedenen,
 zu $P$ orthogonalen und mit $N(P)$ die Menge der zu $P$ nicht
 orthogonalen isotropen Punkte von $V$. Dann sind die Mengen $O(P)$ und
 $N(P)$ Bahnen von $\PSU(V,f)_P$.}
 \smallskip
       Beweis. Ist $n = 2$, so ist $O(P) = \emptyset$ und die Gruppe
 $\PTU(V,f)$ ist auf der Menge der isotropen Punkte von $\La(V)$ zweifach
 transitiv, wie wir schon verschiedentlich bemerkten.
 Daher ist $N(P)$ eine Bahn von $\PTU(V,f)_P$ und damit auch von
 $\PSU(V,f)_P$. Wir d\"urfen daher annehmen, dass $n \geq 3$ ist.
 \par
       Zun\"achst beachten wir jedoch, unabh\"angig vom Rang und vom
 Index, dass $V$
 eine Orthogonalbasis hat. Es gibt ja sicherlich ein $b_1 \in V$ mit
 $f(b_1,b_1) \neq 0$. Dann ist $(b_1K)^\bot$ nicht isotrop.
 Nach Induktionsannahme
 gibt es daher eine Orthogonalbasis $b_2$, \dots, $b_n$ von $(b_1K)^\bot$,
 so dass $b_1$, $b_2$, \dots, $b_n$ eine Orthogonalbasis von $V$ ist.
 \par
       Es sei $\sigma \in U(V,f)$ und es gelte $b_j = \sum_{i:=1}^n b_ia_{ij}$.
 Definiere ferner die Matrix $c$ durch $c_{ij} := f(b_i,b_j)$. Dann ist
 $$ c = (a^\alpha)^tca. $$
 Wegen $\det(c) = \prod_{i:=1}^n f(b_i,b_i) \neq 0$ ist daher
 $$ \det(\sigma)^{\alpha + 1} = \det(a)^{\alpha + 1} = 1. $$
 Entsprechend interpretiert gilt dieser Sachverhalt auch bei nicht kommutativem
 $K$.
 \par
       Zur\"uck zum Beweise des Satzes. Wir zeigen zun\"achst, dass $N(P)$
 eine Bahn von $\PSU(V,f)_P$ ist. Dazu seien $Q$, $R \in N(P)$.
 Es gibt dann ein $q \in Q$ und ein $r \in R$ mit
 $$ f(p,q) = f(p,r) = 1. $$
 Weil auch
 $$ f(q,q) = f(r,r) = 0 $$
 gilt, gibt es nach dem Satz von Witt ein $\tau \in \Un(V,f)$ mit $p^\tau = p$
 und $q^\tau = r$. Wir setzen $G := P + Q$. Dann ist $G$ eine nicht isotrope
 Gerade, so dass $V = G \oplus G^\bot$ ist. Es gibt einen nicht isotropen
 Vektor $s \in G^\bot$. Setze $k := \det(\tau)$. Nach unserer zuvor
 gemachten Bemerkung ist $k^{\alpha + 1} = 1$ und daher auch $k^{-\alpha - 1}
 = 1$. Definiere $\sigma$ durch $s^\sigma := sk^{-1}$ und $u^\sigma = u$
 f\"ur alle $u \in (sK)^\bot$. Es folgt $\sigma \in \Un(V,f)$, da $K$ kommutativ
 ist. Wegen
 $p$, $r \in (sK)^\bot$ folgt $p^{\tau\sigma} = p$ und $q^{\tau\sigma} = r$,
 so dass $P^{\tau\sigma} = P$ und $Q^{\tau\sigma} = R$ ist. Schlie\ss lich
 ist $\det(\tau\sigma) = 1$, so dass $\tau\sigma \in \SU(V,f)$ gilt. Dies
 zeigt, dass $N(P)$ eine Bahn von $\PSU(V,f)_P$ ist, da $N(P)$ ja unter
 $\PSU(V,f)_P$ invariant ist.
 \par
       Nun zeigen wir, dass auch $O(P)$ eine Bahn von $\PSU(V,f)_P$ ist.
 Da es gleichg\"ultig ist, mit welcher Form $\bot$ dargestellt wird,
 nehmen wir f\"ur den Augenblick an, dass $f$ schiefsymmetrisch sei.
 Sind $x$, $x'$, $y$, $y' \in P^\bot$ und gilt $x - x'$, $y - y' \in P$,
 so ist $f(x,y) = f(x',y')$, wie man leicht nachrechnet. Daher wird durch
 $f(x + P,y + P) := f(x,y)$ eine schiefsymmetrische, ebenfalls nicht
 ausgeartete $\alpha$-Form auf $P^\bot/P$ definiert. Ist $q + P$ ein
 isotroper Vektor von $P^\bot/P$, ist $l \in K$ und gilt $l^\alpha = l$, so
 wird durch
 $$ (v + P)^\tau := v + P + (q + P)lf(q + P,v + P) $$
 eine unit\"are Transvektion auf $P^\bot/P$ definiert. Wegen
 $$ (v + P)^\tau = v + qlf(q,v) + P $$
 wird $\tau$ von einer Transvektion aus $\PSU(V,f)_P$ induziert. Weil die
 isotropen Punkte von $P^\bot/P$ den vollst\"andig isotropen Geraden durch $P$
 entsprechen, folgt mittels 8.4, dass $\PSU(V,f)_P$ auf der Menge der
 vollst\"andig isotropen Geraden durch $P$ transitiv operiert.
 \par
       Es sei $f$ wieder die urspr\"ungliche Form. Sind nun $Q$, $R \in O(P)$,
 so gibt es also ein $\gamma \in \PSU(V,f)_P$ mit $Q^\gamma \leq P + R$.
 Wir d\"urfen also des weiteren annehmen, dass $P + Q = P + R$ ist.
 Die Gerade $P + Q$ ist vollst\"andig isotrop. Nach 6.3 gibt es eine vollst\"andig
 isotrope Gerade $G$, so dass $(P + Q) \cap G = \{0\}$ ist und
 $P + Q + G$ nicht isotrop ist. Setze $U := P + Q + G$. Dann ist also $U$
 ein nicht isotroper Raum des Ranges 4. Die Einschr\"ankung von $f$ auf $U$
 ist nicht ausgeartet und definiert eine Polarit\"at, die wir mit $\bot(U)$
 bezeichnen.
 Wegen $\Rg_K(U) = 4$ ist $P + Q = (P + Q)^{\bot(U)}$ und $G = G^{\bot(U)}$.
 \par
       Setze $T := Q^{\bot(U)} \cap G$. Dann ist $T$ ein Punkt auf $G$.
 Es folgt
 $$
 P^{\bot(U)} \cap T = P^{\bot(U)} \cap Q^{\bot(U)} \cap G =
                            (P + Q)^{\bot(U)} \cap G 
                    = (P + Q) \cap G = \{0\}.  $$
 Somit ist $P^{\bot(U)} \cap G$ ein von $T$ verschiedener Punkt auf $G$.
 Es gibt noch
 einen Punkt $S$ auf $G$, der von $T$ und $P^{\bot(U)} \cap G$ verschieden ist.
 Es sei $T = tK$. Dann ist $f(q,t) = 0$. Wegen $T \not\leq P^{\bot(U)}$ gibt
 es ein $p \in P$ mit $f(p,t) = 1$. Ebenso gibt es ein $s \in S$ mit
 $f(p,s) = 1$. Wegen $S \neq T$ gibt es schlie\ss lich ein $q \in Q$ mit
 $f(q,s) = 1$. Es sei $R = rK$. Es gibt dann ein $\lambda \in K$ mit
 $r = q + p\lambda$. Definiert man $\sigma$ durch
 $$\eqalign{
       p^\sigma &:= p                                         \cr
       q^\sigma &:= p\lambda + q                              \cr
       s^\sigma &:= s(1 - \lambda^\alpha) + t\lambda^\alpha   \cr
       t^\sigma &:= -s\lambda^\alpha + t(1 + \lambda^\alpha)  \cr} $$
 und $u^\sigma := u$ f\"ur alle $u \in (P + Q + G)^\bot$, so zeigt eine
 einfache Rechnung, dass $\sigma \in U(V,f)$ und $\det(\sigma) = 1$ ist.
 Also ist $\sigma \in \SU(V,f)_P$ und es gilt $Q^\sigma = R$. Damit ist
 alles bewiesen, da auch $O(P)$ unter $\PSU(V,f)_P$ invariant ist.
 \medskip\noindent
 {\bf 9.9. Satz.} {\it Es sei $K$ ein kommutativer K\"orper und $\alpha$ sei ein
 involutorischer Automorphismus von $K$. Ist $V$ ein $K$-Vektorraum des
 Ranges $n \geq 2$ und ist $f$ eine nicht entartete, symmetrische und
 spurwertige $\alpha$-Form auf $V$, deren Index mindestens $1$ ist, so
 operiert die Gruppe $\PSU(V,f)$ auf der Menge der isotropen Punkte primitiv.}
 \smallskip
       Beweis. Ist der Index von $f$ gleich 1, so ist $N(P) = \emptyset$
 f\"ur alle isotropen Punkte $P$ von $\La(V)$. In diesem Falle ist $\PSU(V,f)$
 nach 9.8 auf der Menge der isotropen Punkte zweifach transitiv, also erst recht
 primitiv.
 \par
       Wir nehmen nun an, der Satz sei falsch. Dann ist der Index von $f$
 mindestens 2 und folglich $n \geq 4$. Es sei $\Delta$ ein
 Imprimitivit\"atsgebiet von $\PSU(V,f)$. Ist $P \in \Delta$, so gibt es
 noch einen weiteren Punkt in $\Delta$, so dass mit 9.8 folgt, dass
 $\Delta = \{P\} \cup O(P)$ oder $\Delta = \{P\} \cup N(P)$ ist.
 \par
       1. Fall: Es sei $\Delta = \{P\} \cup O(P)$. Wegen $n \geq 4$ ist
 der Rang von $P^\bot/P$ mindestens 2. Weil $f$ auf $P^\bot/P$ eine nicht
 entartete Form induziert, gibt es eine Ebene $E$ mit $P \leq E \leq P^\bot$,
 so dass $E/P$ eine Gerade von $P^\bot/P$ ist, die mindestens zwei
 isotrope Punkte enth\"alt. Ist dann $G$ ein Komplement von $P$ in $E$,
 so folgt aus $E \leq P^\bot$, dass $E/P$ und $G$ isometrisch sind.
 Es gibt daher zwei isotrope Vektoren $v$, $w \in G$ mit $f(v,w) = 1$.
 Es gilt
 $$ vK,\ wK \in O(P) \subseteq \Delta, $$
 und daher
 $$ \{vK\} \cup O(vK) = \Delta = \{wK\} \cup O(wK). $$
 Wegen $vK \neq wK$ ist daher $wK \in O(vK)$, so dass wir den Widerspruch
 $0 = f(v,w) = 1$ erhalten.
 \par
       2. Fall: Es sei $\Delta = \{P\} \cup N(P)$. Es sei $P = vK$.
 Weil der Index von $f$
 mindestens 2 ist, gibt es eine vollst\"andig isotrope Gerade $G$ und in jedem
 Falle auch eine nicht isotrope Gerade $H$ durch $P$ (Satz 8.10).
 Es gibt dann einen isotropen Vektor $u \in H$ mit $f(v,u) = 1$. Es sei
 $w \in G - P$. Es ist $uK \in N(P) \subseteq \Delta$. Also ist
 $$ \Delta = \{uK\} \cup N(uK). $$
 Setze $s := v + w$. Dann ist $P \neq sK$. Au\ss erdem ist
 $f(v,s) = f(v,v) + f(v,w) = 0$. Also ist $sK \in O(P)$ und damit
 $$ sK \not\in \Delta. $$
 Andererseits ist $f(u,s) = f(u,v) + f(u,u) = 1$ und daher
 $$ sK \in N(uK) \subseteq \Delta, $$
 so dass wir auch in diesem Falle einen Widerspruch erhalten. Damit
 ist der Satz bewiesen.
 \medskip
       Wir beschlie\ss en diesen Abschnitt, indem wir einen weiteren
 Ausnahmeisomorphismus\index{Ausnahmeisomorphien}{} etablieren.
 \medskip\noindent
 {\bf 9.10. Satz.} {\it Die Gruppen $\PSU(4,4)$ und $\PSp(4,3)$ sind isomorph.}
 \smallskip
       Beweis. Es sei $V$ ein Vektorraum des Ranges 4 \"uber $\GF(4)$ und
 $f$ sei eine nicht ausgeartete $\alpha$-Form auf $V$, wobei $\alpha$ der
 involutorische Automorphismus von $\GF(4)$ sei. Wir definieren eine
 Inzidenzstruktur $\Delta = (\Pi,\Beta,\I)$ wie folgt: ist $P$ ein Punkt von
 $\La(V)$, so geh\"ore $P$ genau dann zu $\Pi$, wenn $P$ nicht isotrop ist.
 Ferner sei $\Beta := \{P^\bot \mid P \in \Pi\}$. Es gelte
 $Q\>\I\>P^\bot$ genau dann, wenn $Q = P$ oder wenn $Q \leq P^\bot$ ist.
 \par
       Die Anzahl der Punkte von $\La(V)$ ist $4^3 + 4^2 + 4 + 1 = 85$. Nach
 9.4b) ist die Anzahl der isotropen Punkte gleich $(4 + 1)i \cdot 9 = 45$. Daher
 ist die Anzahl $v$ der Punkte von $\Delta$ gleich 40. Dann ist aber auch
 $|\Beta| = 40$, da $\bot$ ja bijektiv ist. Ist $P \in \Pi$, so ist
 die Anzahl aller Punkte von $P^\bot$ gleich $4^2 + 4 + 1 = 21$.
 Daher inzidieren
 in $\Delta$ mit $P^\bot$ genau $21 - 9 + 1 = 13$ Punkte. Ebenso sieht man,
 dass jedes $P \in \Pi$ in $\Delta$ mit genau 13 Bl\"ocken inzidiert.
 Also ist $\Delta$ eine taktische Konfiguration mit den Parametern
 $v = b = 40$ und $k = r = 13$.
 \par
       Es seien $P$ und $Q$ zwei verschiedene Punkte von $\Delta$.
 \par
       1. Fall: Es ist $P \leq Q^\bot$. Dann ist auch $Q \leq P^\bot$.
 Mittels des Modulargesetzes folgt
 $$\eqalign{
 (P + Q) \cap (P + Q)^\bot &= (P + Q) \cap P^\bot \cap Q^\bot \cr
       &= \bigl((P \cap P^\bot) + Q\bigr) \cap Q^\bot                 
       = Q \cap Q^\bot                                       
       = \{0\}.                                              \cr} $$
 Somit ist $P + Q$ nicht isotrop, tr\"agt also genau drei isotrope und dann
 genau zwei nicht isotrope Punkte, n\"amlich $P$ und $Q$. Dann tr\"agt aber
 auch die Gerade $(P + Q)^\bot$ genau zwei nicht isotrope Punkte, die wir
 $P_1$ und $Q_1$ nennen. Es folgt, dass $P$ und $Q$ simultan mit genau den
 Bl\"ocken $P^\bot$, $Q^\bot$, $P_1^\bot$, $Q_1^\bot$ inzidieren. Auf all
 diesen Bl\"ocken liegen aber auch die Punkte $P_1$ und $Q_1$, so dass
 die $\Delta$-Gerade durch $P$ und $Q$ genau vier Punkte enth\"alt.
 \par
       2. Fall: Es ist $P \not\leq Q^\bot$. Dann ist auch $Q \not\leq P^\bot$.
 Weil alle nicht isotropen Geraden isometrisch sind, kann $P + Q$ nicht
 nicht isotrop sein. Daher gibt es auf $P + Q$ genau einen isotropen Punkt,
 n\"amlich $\rad(P + Q)$ und neben $P$ und $Q$ noch zwei weitere nicht
 isotrope Punkte $R$ und $S$. Weil $P + Q$ und $(P + Q)^\bot$ symmetrische
 Rollen spielen, gibt es auf $(P + Q)^\bot$ ebenfalls vier nicht isotrope
 Punkte $P_1$, $Q_1$, $R_1$ und $S_1$. Die Punkte $P$ und $Q$ inzidieren in
 $\Delta$ simultan also mit genau den Bl\"ocken $P_1^\bot$, $Q_1^\bot$,
 $R_1^\bot$, $S_1^\bot$. Dar\"uber hinaus liegen auch $R$ und $S$ auf der
 $\Delta$-Geraden durch $P$ und $Q$.
 \par
       Damit ist gezeigt, dass $\Delta$ ein 2-$(40,13,4)$ Blockplan ist,
 dessen Geraden alle
 $$ 4 = {40 - 4 \over 13 - 4} = {b - \lambda \over r - \lambda} $$
 Punkte tragen. Nach IIII.10.1 besteht $\Delta$ also aus den Punkten und
 Hyperebenen einer projektiven Geometrie, da ja auch $v - k = 40 - 13 \geq 2$
 ist. Diese Geometrie hat die Ordnung $q = 3$ und dann wegen
 $$ 40 = 1 + 3 + 3^2 + 3^3 $$
 --- einer Zerlegung der 40, die schon bei Fibonacci eine Rolle spielte ---
 den Rang 4. Also ist $\Delta$ zur Geometrie der Punkte und Ebenen eines
 Vektorraums vom Rang 4 \"uber $\GF(3)$ isomorph. Die Abbildung $\pi$, die
 durch $P^\pi := P^\bot$ und $(P^\bot)^\pi := P$ definiert wird, ist
 offenbar eine Polarit\"at von $\Delta$, die von einer symplektischen
 Polarit\"at der zugeh\"origen projektiven Geometrie herr\"uhrt. Es folgt
 $\PGU(4,4) \subseteq \PSp(4,3)$. Nach 9.3 ist
 $$ \big|\PSU(4,4)\big| = \big|\PGU(4,4)\big| = 2^6 \cdot 3 \cdot 9 \cdot 15 $$
 und nach 3.3 ist
 $$ \big|\PSp(4,3)\big| = {1 \over 2} \cdot 3^4 \cdot 8 \cdot 80 =
	  2^6 \cdot 3 \cdot 9 \cdot 15. $$
 Also ist $\PSU(4,4) = \PSp(4,3)$, q. e. d.
 \medskip\noindent
 {\bf 9.11. Korollar.} {\it Die Gruppe $\PSU(4,4)$ ist einfach.}
 \smallskip
       Beweis. Dies folgt mit 9.10 aus 3.11.

\mysectionten{10. Die speziellen unit\"aren Gruppen}

\noindent
       Unser n\"achstes Ziel ist zu zeigen, dass bis auf eine Ausnahme,
 die wir schon kennen, $\TU(V,f) = \SU(V,f)$ gilt, falls der $V$ zugrunde
 liegende K\"orper kommutativ ist. Dies werden wir dann benutzen, um die
 Einfachheit der $\PSU(V,f)$ in den noch fehlenden F\"allen nachzuweisen.
 \medskip\noindent
 {\bf 10.1. Satz.} {\it Es sei $K$ ein K\"orper und $\alpha$ sei ein
 involutorischer Antiautomorphismus von $K$. Ist $V$ ein Vektorraum des
 Ranges $n \geq 2$ \"uber $K$ und ist $f$ eine nicht entartete, symmetrische,
 spurwertige $\alpha$-Form auf $V$, deren Index mindestens gleich $1$ ist,
 so gibt es zu jedem $v \in V$ eine nicht isotrope, isotrope Punkte tragende
 Gerade $G$ mit $v \in G$.}
 \smallskip
       Beweis. Ist $f(v,v) = 0$, so folgt dies mit 8.10. Es sei also
 $f(v,v) \neq 0$. Weil $f$ spurwertig ist, gibt es ein $b \in K$ mit
 $f(v,v) = b + b^\alpha$. Weil der Index von $f$ mindestens 1 ist, gibt es
 eine Gerade $H = v_1K + v_2K$ mit isotropen $v_i$, so dass $f(v_1,v_2) = 1$
 ist. Es folgt
 $$ f(v_1 + v_2b,v_1 + v_2b) = b^\alpha + b = f(v,v). $$
 Nach dem Satz von Witt gibt es ein $\sigma \in \Un(V,f)$ mit 
 $(v_1 + b_2b)^\sigma = v$. Setze $G := H^\sigma$. Dann ist $G$ eine Gerade
 der verlangten Art.
 \medskip\noindent
 {\bf 10.2. Satz.} {\it Es sei $K$ ein kommutativer K\"orper, der mehr als $4$
 Elemente enthalte, und $\alpha$ sei
 ein involutorischer Automorphismus von $K$. Ferner sei $V$ ein Vektorraum
 des Ranges $n \geq 3$ \"uber $K$ und $f$ sei eine nicht ausgeartete,
 symmetrische, spurwertige $\alpha$-Form auf $V$, deren Index mindestens $1$
 sei. Sind $v_1$, $v_2 \in V$, gilt $f(v_1,v_1) = f(v_2,v_2) \neq 0$, ist
 $v_1K \neq v_2K$ und ist die Gerade $v_1K + v_2K$ isotrop, so gibt es ein
 $v_3 \in V$ mit $f(v_3,v_3) = f(v_1,v_1)$, so dass $v_1K + v_3K$ und
 $v_2K + v_3K$ nicht isotrope Geraden sind.}
 \smallskip
       Beweis. Setze $G := v_1K + v_2K$. Weil $f(v_1,v_1) \neq 0$ ist, ist
 $$ G = v_1K \oplus \rad(G), $$
 wobei $\rad(U) := U \cap U^\bot$ f\"ur alle $U \in \La(V)$ gesetzt sei.
 Es gibt einen Vektor $r \in \rad(G)$ und ein $a \in K$ mit
 $$ v_2 = r + v_1a. $$
 Es folgt
 $$ f(v_1,v_1) = f(v_2,v_2) = a^\alpha f(v_1,v_1)a. $$
 Weil $f(v_1,v_1) \neq 0$ und weil $K$ kommutativ ist, folgt
 $$ 1 = a^{\alpha+1}. $$
 Wegen $v_2 \not\in v_1K$ ist $r \neq 0$. Ferner ist $r \in (v_1K)^\bot$.
 Weil $(v_1K)^\bot$ nicht isotrop ist, gibt es einen weiteren isotropen Vektor
 $r' \in (v_1K)^\bot$ mit $f(r,r') = 1$.
 \par
       Setze $L := \{k \mid k \in K, k^\alpha = k\}$. Dann ist $L$ ein
 Teilk\"orper von $K$ mit $[K:L] = 2$. Es folgt, dass $L$ mindestens drei
 Elemente enth\"alt, da $K$ mehr als vier Elemente enth\"alt. Die Norm der
 Erweiterung $K:L$ ist bekanntlich die durch $N(k) := k^{1+\alpha}$
 definierte Abbildung $N$. Die Menge $N(K^*)$ enth\"alt alle Quadrate
 von $L$ und und damit mindestens ein von 1 verschiedenes Element, falls
 $L$ mehr als drei Elemente enth\"alt. Dies ist aber auch richtig,
 falls $|L| = 3$ ist, da in diesem Falle $N(K^*) = L$ ist. Es gibt also ein
 $z \in K^*$ mit $N(z) \neq 1$.
 \par
       Ist $k \in K - L$, so ist $k^{1-\alpha} \neq 1$. Ferner gilt
 $N(k^{1-\alpha}) = 1$. Folglich ist $\Kern(N)$ nicht trivial. Weil
 andererseits auch $N(K^*)$ mindestens zwei Elemente enth\"alt, enth\"alt
 $N(K^*)$ zun\"achst mindestens ein, dann aber mindestens zwei Elemente,
 die von $f(v_1,v_1)$ verschieden sind. Es gibt daher ein $t \in N(K^*)$ mit
 $t \neq f(v_1,v_1)$, $a^\alpha zf(v_1,v_1)$. Wegen $t \in N(K^*)$ ist
 $t^\alpha = t$, was wir noch benutzen werden.
 \par
       Setze $y := t - a^\alpha zf(v_1,v_1)$. Dann ist $y \neq 0$.
 \par
       Weil $f$ spurwertig ist, gibt es ein $k \in K$ mit
 $$ k + k^\alpha = -zz^\alpha f(v_1,v_1) + f(v_1,v_1). $$
 Setze $x := y^{-\alpha}k$. Dann ist
 $$ xy^\alpha + x^\alpha y = -zz^\alpha f(v_1,v_1) + f(v_1,v_1). $$
 \par
       Setze nun $v_3 := rx + r'y + v_1z$. Dann ist
 $$ f(v_3,v_3) = x^\alpha y + y^\alpha x + zz^\alpha f(v_1,v_1)
                      = f(v_1,v_1). $$
 \par
       Um nachzuweisen, dass $v_1K + v_3K$ nicht isotrop ist, gen\"ugt es
 nachzuweisen, dass die Determinante der Matrix
 $$ \pmatrix{ f(v_1,v_1) & f(v_1,v_3) \cr
              f(v_3,v_1) & f(v_3,v_3) \cr} $$
 nicht Null ist. Wegen $f(v_1,v_1) = f(v_3,v_3)$ und $f(v_1,v_3) =
 zf(v_1,v_1)$ ist diese Determinante gleich
 $$ f(v_1,v_1)^2(1 - zz^\alpha), $$
 also in der Tat von Null verschieden, da ja $zz^\alpha \neq 1$ ist.
 \par
       Um nachzuweisen, dass $v_2K + v_3K$ nicht isotrop ist, muss man
 entsprechend nachweisen, dass die Determinante der Matrix
 $$ \pmatrix{ f(v_2,v_2) & f(v_2,v_3) \cr
              f(v_3,v_2) & f(v_3,v_3) \cr} $$
 ungleich Null ist. Wegen
 $$ f(v_3,v_3) = f(v_1,v_1) = f(v_2,v_2) $$
 und
 $$ f(v_2,v_3) = y + a^\alpha zf(v_1,v_1) = t $$
 ist diese Determinante gleich
 $$ f(v_1,v_1)^2 - tt^\alpha = f(v_1,v_1)^2 - t^2 \neq 0. $$
 Damit ist der Satz bewiesen.
 \medskip\noindent
 {\bf 10.3. Satz.} {\it Es sei $K$ ein kommutativer K\"orper und $\alpha$ sei
 ein involutorischer Automorphismus von $K$. Ferner sei $V$ ein
 $K$-Vektorraum des Ranges $n \geq 3$. Ist $|K| = 4$, so sei $n \geq 4$.
 Ist $f$ eine nicht ausgeartete, symmetrische, spurwertige $\alpha$-Form
 auf $V$, deren Index mindestens $1$ ist, sind $v_1$, $v_2 \in V$ und gilt
 $f(v_1,v_1) = f(v_2,v_2) \neq 0$, so gibt es ein $\gamma \in \TU(V,f)$ mit
 $v_1^\gamma = v_2$.}
 \smallskip
       Beweis. Es sei zun\"achst $v_1K = v_2K$. Dann ist $v_2 = v_1a$ mit einem
 $a \in K$. Es folgt $f(v_1,v_1) = a^\alpha f(v_1,v_1)a$ und daher
 $a^{1+\alpha} = 1$,
 da ja $f(v_1,v_1) \neq 0$ ist. Nach 10.1 gibt es eine nicht isotrope Gerade
 durch $v_1K$, die zwei isotrope Punkte tr\"agt. Es sei
 $0 \neq u \in G \cap (v_1K)^\bot$. Weil $v_1$ nicht isotrop ist, ist $u$,
 $v_1$ eine Basis von $G$. Definiere nun $\sigma$ durch
 $$x^\sigma := \cases{v_2, &falls $x = v_1$        \cr
                      ua^{-\alpha}, &falls $x = u$ \cr
                      x, &falls $x \in G^\bot$.    \cr} $$
 Dann ist $\sigma$ ein Element von $\SU(V,f)$. Die unit\"aren Transvektionen,
 deren Zentren auf $G$ liegen, erzeugen eine Untergruppe $S$ von $\SU(V,f)$,
 die auf $G$ die $\SU(G,f)$ induziert, wie aus Satz 8.9 folgt. Weil $\sigma$
 auf $G$ ein Element der $\SU(G,f)$ induziert und auf $G^\bot$ gleich der
 Identit\"at ist, folgt $\sigma \in S$, so dass $\sigma \in \TU(V,f)$ ist.
 \par
       Es sei also $v_1K \neq v_2K$. Setze $G := v_1K + v_2K$.
 \par
       1. Fall: $G$ ist nicht isotrop, enth\"alt aber isotrope Punkte.
 Hier definieren wir $\sigma$ durch
 $$ x^\sigma := \cases{v_2,  &falls $x = v_1$    \cr
                       -v_1, &falls $x = v_2$    \cr
                       x,  &falls $x \in G^\bot$ \cr} $$
 und schlie\ss en wie eben weiter.
 \par
       2. Fall: $G$ enth\"alt keine isotropen Punkte. Dies hat zur Konsequenz,
 dass $K$ nicht endlich ist. Dies werden wir noch auszunutzen haben.
 \par
       Nach 10.1 gibt es eine
 nicht isotrope Gerade $H$ durch den Punkt $(v_1 - v_2)K$, die isotrope Punkte
 tr\"agt. Es ist also $V = H \oplus H^\bot$. Es gibt daher $w_i \in H$ und
 $s_i \in H^\bot$ mit $v_i = w_i + s_i$. Es folgt
 $$ s_1 - s_2 = v_1 - v_2 + w_2 - w_1 \in H \cap H^\bot = \{0\}. $$
 Setzt man $s := s_1$, so ist also $v_i = w_i + s$ f\"ur $i := 1$, 2. Ferner
 gilt
 $$ f(v_1,v_1) = f(w_1,w_1) + f(s,s) $$
 und
 $$ f(v_2,v_2) = f(w_2,w_2) + f(s,s) $$
 und folglich
 $$ f(w_1,w_1) = f(w_2,w_2). $$
 \par
       Ist $f(w_1,w_1) \neq 0$, so gibt es, wie im ersten Fall gesehen,
 ein $\sigma \in \TU(V,f)$ mit $w_1^\sigma = w_2$ und $s^\sigma = s$. Es folgt
 $v_1^\sigma = v_2$.
 \par
       Es seien $w_1$ und $w_2$ isotrop. Weil $\alpha$ nicht die Identit\"at
 ist, gibt es ein $a \in K$ mit $a^\alpha \neq a$. Setze $k := a - a^\alpha$.
 Dann ist $k \neq 0$ und $k^\alpha = -k$. Setze ferner $g := kf$. Dann ist
 $g$ schiefsymmetrisch und es gilt $\SU(V,f) = \SU(V,g)$ und $\TU(V,f) =
 \TU(V,g)$.
 Ferner gilt $f(v,v) = f(w,w)$ genau dann, wenn $g(v,v) = g(w,w)$ ist. Also
 d\"urfen wir bei den nun folgenden Betrachtungen $f$ durch $g$ ersetzen.
 Wir setzen
 $$ L := \{k \mid k \in K, k^\alpha = k\}. $$
 Ferner setzen wir $a := g(w_1,w_2)$. Weil $G$ keine isotropen Punkte tr\"agt,
 ist $g(v_1 - v_2,v_1 - v_2) \neq 0$. Andererseits ist
 $$\eqalign{
 g(v_1 - v_2,v_1 - v_2) &= g(w_1 - w_2,w_1 - w_2)  \cr
       &= -g(w_1,w_2) - g(w_2,w_1)               
       = -a + a^\alpha,                         \cr} $$
 so dass $a^\alpha \neq a$ ist. Somit ist $a \not\in L$. Setze
 $w := w_1a + w_2$. Dann ist
 $$ g(w,w) = a^\alpha g(w_1,w_2) + g(w_2,w_1)a = a^\alpha a - a^\alpha a = 0 $$
 und
 $$ g(v_1,w) = g(w_1 + s,w_1a + w_2) = g(w_1,w_2) = a \neq 0 $$
 sowie
 $$ g(v_2,w) = g(w_2 + s,w_1a + w_2) = g(w_2,w_1)a = -a^\alpha a \neq 0. $$
 \par
       F\"ur $x \in L$ setze
 $$\eqalign{
 p(x) := x^2(a^{1+\alpha})^3 &+ xa^{1+\alpha}\bigl(a^\alpha g(v_1,v_2) +
                                          ag(v_1,v_2)^\alpha\bigr)       \cr
                     &+ g(v_1,v_1)^2 + g(v_1,v_2)g(v_1,v_2)^\alpha. \cr} $$
 Wegen
 $$ \bigl(g(v_1,v_1)^2\bigr)^\alpha = \bigl(g(v_1,v_1)^\alpha\bigr)^2 =
	      \bigl(-g(v_1,v_1)\bigr)^2 = g(v_1,v_1)^2 $$
 ist $p$ eine Abbildung von $L$ in sich.
 \par
       F\"ur $y \in L$ setze
 $$ q(y) := y^2g(v_1,v_1) + y\bigl(g(v_1,v_2) + g(v_2,v_1)\bigr) + g(v_1,v_1). $$
 Hat $p$ keine Nullstelle in $L$, so sei $\eta$ irgendein Element von $L^*$.
 Hat $p$ eine Nullstelle in $L$, so hat $p$ zwei Nullstellen $\nu_1$ und
 $\nu_2$, falls man Vielfachheiten ggf. mitz\"ahlt. Weil $K$ unendlich ist, ist
 auch $L$ unendlich. Es gibt daher ein $\eta \in L^*$ mit
 $$ q(\eta) + \nu_i\eta aa^\alpha(a - a^\alpha) \neq 0. $$
 Setze
 $$ \xi := {-q(\eta) \over \eta aa^\alpha(a - a^\alpha)}. $$
 Dann ist $\xi \in L$, da $\xi$ unter $\alpha$ invariant bleibt, wie man
 sich rasch \"uberzeugt. Aus
 $$q(\eta) + \xi\eta aa^\alpha(a - a^\alpha) = 0 $$
 folgt,
 dass $\xi$ keine Nullstelle von $p$ ist. Also ist $p(\xi) \neq 0$.
 Wir definieren nun $\tau$ durch
 $$ v^\tau := v - w\xi g(w,v).$$
 Dann ist
 $$ v_1^\tau = v_1 - w\xi f(w,v_1) = v_1 + w\xi a^\alpha. $$
 Wir zeigen, dass $v_2K + v_1^\tau K$ nicht isotrop ist, aber isotrope Punkte
 tr\"agt. Eine einfache Rechnung zeigt, dass
 $$ \det\pmatrix{g(v_2,v_2)      & g(v_2,v_1^\tau)      \cr
                 g(v_1^\tau,v_2) & g(v_1^\tau,v_1^\tau) \cr}
                     = p(\xi) \neq 0 $$
 ist. Dies zeigt, dass $v_2K + v_1^\tau K$ nicht isotrop ist. Ferner gilt
 $v_1^\tau + v_2\eta \in v_2K + v_1^\tau K$ und
 $$ g(v_1^\tau + v_2\eta,v_1^\tau + v_2\eta) =
       q(\eta) + \eta\xi aa^\alpha(a - a^\alpha) = 0. $$
 Daher ist $v_2K + v_1^\tau K$ eine isotrope Punkte tragende nicht isotrope
 Gerade. Nach Fall 1 gibt es daher ein $\sigma \in \TU(V,g)$
 mit $v_1^{\tau\sigma} = v_2$. Damit ist der fragliche Sachverhalt auch in
 diesem Falle bewiesen.
 \par
       3. Fall: $G$ ist isotrop und $K$ enth\"alt mehr als vier Elemente.
 Nach 10.2 gibt es dann ein $v_3 \in V$ mit $g(v_1,v_1) = g(v_3,v_3)$, so
 dass die Geraden $v_1K + v_3K$ und $v_3K + v_2K$ nicht singul\"ar sind.
 Mittels der bereits erledigten F\"alle 1 und 2 gelangt man daher ans Ziel.
 \par
       4. Fall: $G$ ist isotrop und $K$ enth\"alt genau vier Elemente. Nach
 Voraussetzung ist dann $n \geq 4$. Es folgt $\Rg_K(G^\bot) = n - 2 \geq 2$.
 W\"are $G^\bot$ vollst\"andig isotrop, so folgte
 $$ G^\bot \leq G^{\bot\bot} = G $$
 und damit $n - 2 \leq 2$, so dass $G = G^\bot$ w\"are. Somit w\"are 
 $G$ vollst\"andig isotrop, was nicht der Fall ist, da $g(v_1,v_1) \neq 0$ ist.
 Weil $G^\bot$ also nicht isotrop ist, gibt es ein $v_3 \in G^\bot$ mit
 $g(v_3,v_3) \neq 0$. Es folgt, dass die Geraden $v_1K + v_3K$ und
 $v_2K + v_3K$ nicht singul\"ar sind, so dass auch in diesem Falle mit
 Fall 1 die Behauptung folgt.
 \medskip\noindent
 {\bf 10.4. Satz.} {\it Es sei $K$ ein kommutativer K\"orper und $\alpha$ sei ein
 involutorischer Automorphismus von $K$. Ist $V$ ein Vektorraum des Ranges
 $n \geq 2$ \"uber $K$ und ist $f$ eine nicht ausgeartete, symmetrische,
 spurwertige $\alpha$-Form auf $V$, deren Index mindestens $1$ sei, so ist
 $SU(V,f) = TU(V,f)$, es sei denn, es ist $n = 3$ und $|K| = 4$.}
 \smallskip
       Beweis. Ist $n = 2$, so folgt dies aus Satz 8.9. Dies nehmen wir
 als Induktionsverankerung, falls $K$ mehr als vier Elemente hat. Hat
 $K$ nur vier Elemente, so ist die Gruppe $\PSU(V,f)$ nach 9.11 einfach.
 Mit 9.3b) und d) folgt, dass die Gruppen $\PSU(V,f)$ und $\SU(V,f)$
 wegen $\ggT(4,2 + 1) = 1$ isomorph sind. Daher ist auch $\SU(V,f)$ in diesem
 Falle einfach. Weil $\TU(V,f)$ ein nicht trivialer Normalteiler von $\SU(V,f)$
 ist, ist also $\TU(V,f) = \SU(V,f)$
 \par
       Es sei nun $n \geq 3$ und $K$ enthalte mehr als vier Elemente oder
 $n \geq 5$ und $|K| = 4$. Ferner $\sigma \in \SU(V)$ und $v \in V$ mit
 $f(v,v) \neq 0$. Nach 10.3 gibt es dann ein $\tau \in TU(V,f)$ mit
 $v^{\sigma\tau} = v$. Die Einschr\"ankung von $\sigma\tau$ auf $(vK)^\bot$ ist
 nach Induktionsannahme ein Produkt von Transvektionen mit Zentrum in
 $(vK)^\bot$. Da diese den Punkt $vK$ vektorweise festlassen, ist dieses
 Produkt auf ganz $V$ gleich $\sigma\tau$. Daher ist $\sigma\tau \in \TU(V,f)$
 und folglich $\sigma \in \TU(V,f)$. Damit ist der Satz bewiesen.
 \medskip\noindent
 {\bf 10.5. Satz.} {\it Es sei $K$ ein kommutativer K\"orper und $\alpha$ sei ein
 involutorischer Automorphismus von $K$. Ist $V$ ein Vektorraum des Ranges
 $n \geq 2$ \"uber $K$ und ist $f$ eine nicht entartete, symmetrische,
 spur\-wer\-ti\-ge $\alpha$-Form auf $V$, deren Index mindestens $1$ sei, so ist
 $\SU(V,f) = \SU(V,f)'$, es sei denn es ist $(n,|K|) = (2,4)$, $(2,9)$, $(3,4)$.}
 \smallskip
       Beweis. Ist $n = 2$, so folgt dies mit den S\"atzen 8.9 und III.2.7.
 Es sei also $n \geq 3$. Es sei $\tau$ eine unit\"are Transvektion.
 Es gibt dann einen isotropen Vektor $v_1$ und ein $c \in K$ mit $c^\alpha =
 -c$, so dass
 $$ v^\tau = v - v_1cf(v_1,v) $$
 ist. Nach 8.10 gibt es einen isotropen Vektor $v_3$, so dass die Gerade
 $G := v_1K + v_3K$ nicht isotrop ist. Es ist
 $V = G \oplus G^\bot$ und $G^\bot \neq \{0\}$. Daher gibt es einen nicht
 isotropen Vektor $v_2 \in G^\bot$. Indem man $v_3$ gegebenenfalls durch einen
 Skalar ab\"andert, kann man erreichen, dass $f(v_1,v_3) = f(v_2,v_2)$ ist.
 Weil $T$ auf $(G + v_2K)^\bot$ die Identit\"at induziert und
 $$ V = G \oplus v_2K \oplus (G + v_2K)^\bot $$
 gilt, d\"urfen wir, um zu zeigen, dass $\tau$ ein Kommutator ist,
 annehmen, dass $n = 3$ ist. Dann ist $v_1$, $v_2$, $v_3$ eine Basis von
 $V$ und es gilt
 $$\eqalign{
       v_1^\tau &= v_1,        \cr
       v_2^\tau &= v_2,        \cr
       v_3^\tau &= v_3 - v_1c. \cr} $$
 \par
       F\"ur $x$, $y \in K$ definieren wir eine Abbildung $\rho(x,y) \in
 \SL(V)$ durch
 $$\eqalign{
       v_1^{\rho(x,y)} &:= v_1,                     \cr
       v_2^{\rho(x,y)} &:= v_1x + v_2,              \cr
       v_3^{\rho(x,y)} &:= v_1y - v_2x^\alpha + v_3.\cr} $$
 \par
       Wir beachten zun\"achst, dass
 $ \tau = \rho(0,-c) $
 ist. Ferner gilt, wie einfache Rechnungen zeigen,
 $$ \rho(x',y')\rho(x,y) = \rho(x' + x,y' + y - x'^\alpha x). $$
 Hieraus folgt
 $$ \rho(x,y)^{-1} = \rho(-x,-y - xx^\alpha) $$
 und weiter
 $$ \rho(x',y')^{-1}\rho(x,y)^{-1}\rho(x',y')\rho(x,y) =
       \rho(0,x'x^\alpha - xx'^\alpha). $$
 \par
       Als N\"achstes suchen wir Bedingungen daf\"ur, dass die Abbildung
 $\rho(x,y)$ zu $\SU(V,f)$ geh\"ort.
 Sofort zu sehen ist, dass $\det(\rho(x,y)) = 1$ ist. Ferner ist
 $$\eqalign{
 f(v_1^{\rho(x,y)},v_1^{\rho(x,y)}) &= f(v_1,v_1) \cr
 f(v_1^{\rho(x,y)},v_2^{\rho(x,y)}) &= f(v_1,v_1x + v_2) = f(v_1,v_2) \cr
 f(v_1^{\rho(x,y)},v_3^{\rho(x,y)}) &= f(v_1,v_1y - v_2c^\alpha + v_3)
                      = f(v_1,v_3) \cr
 f(v_2^{\rho(x,y)},v_2^{\rho(x,y)}) &= f(v_1x + v_2,v_1x + v_2) = f(v_2,v_2).
       \cr} $$
 Weiter ist, da ja $f(v_1,v_3) = f(v_2,v_2)$ ist,
 $$\eqalign{
 f(v_2^{\rho(x,y)},v_3^{\rho(x,y)}) &= f(v_1x + v_2,v_1y - v_2x^\alpha + v_3)
                                          \cr
       &= x^\alpha f(v_1,v_3) - f(v_2,v_2)x^\alpha 
       = 0 = f(v_2,v_3).                          \cr} $$
 Bis hierher ergeben sich also noch keine Einschr\"ankungen daf\"ur, dass
 $\rho(x,y)$ ein Element von $\SU(V,f)$ ist. Dies \"andert sich jetzt. Es
 gilt n\"amlich
 $$\eqalign{
 f(v_3^{\rho(x,y)},v_3^{\rho(x,y)})
       &= f(v_1y - v_2x^\alpha + v_3,v_1y - v_2x^\alpha + v_3)     \cr
       &= y^\alpha f(v_1,v_3) + xx^\alpha f(v_2,v_2) + yf(v_3,v_1) \cr
       &= (y + y^\alpha + xx^\alpha)f(v_2,v_2).                     \cr} $$
 Bei dieser Rechnung ist zu beachten, dass
 $$ f(v_3,v_1) = f(v_1,v_3)^\alpha = f(v_2,v_2)^\alpha = f(v_2,v_2) $$
 ist. Es gilt also $\rho(x,y) \in i\SU(V,f)$ genau dann, wenn
 $$ y + y^\alpha + xx^\alpha = 0 $$
 ist.
 \par
       Wir beachten weiter: Es gibt ein $k \in K$ mit $k^\alpha + k \neq 0$
 (Satz III.6.3).
 Ist $m \in L$, so setzen wir $l := (k^\alpha + k)^{-1}m$. Dann folgt
 $$ (lk)^\alpha + lk = m. $$
 Somit ist die Abbildung $x \to x^\alpha + x$ surjektiv.
 \par
       Es gibt ein $t \in K$ mit $t + t^\alpha \neq 0$. Setze
 $$ x := {ct^\alpha \over t + t^\alpha}. $$
 Wegen $c + c^\alpha = 0$ ist dann
 $$ x^\alpha - x = {c^\alpha t - ct^\alpha \over t + t^\alpha}
       = {c^\alpha t + ct - ct - ct^\alpha \over t + t^\alpha}
       = -c. $$
 (Dies ist ein Spezialfall der additiven Form von Hilberts Satz
 90.)\index{Hilberts Satz 90}{}
 Wie wir gerade bemerkten, gibt es ein $y \in K$ mit $y + y^\alpha +
 xx^\alpha = 0$. Ebenso gibt es ein $z \in K$ mit $z + z^\alpha + 1 = 0$.
 Dann sind $\rho(1,z)$, $\rho(x,y) \in \SU(V,f)$ und es gilt
 $$ \rho(1,z)^{-1}\rho(x,y)^{-1}\rho(1,z)\rho(x,y) = \rho(0,x^\alpha - x)
       = \rho(0,-c) = \tau. $$
 Damit ist gezeigt, dass $\tau$ ein Kommutator ist. (Bis hierhin gilt der
 Beweis auch im Falle, dass $n = 3$ und $|K| = 4$ ist.) Weil $\SU(V,f)$ nach
 Satz 10.4 von seinen Transvektionen erzeugt wird, ist also $\SU(V,f) =
 \SU(V,f)'$, q. e. d.
 \medskip\noindent
 {\bf 10.6. Satz.} {\it Es sei $K$ ein kommutativer K\"orper und $\alpha$ sei ein
 involutorischer Automorphismus von $K$. Ferner sei $V$ ein Vektorraum
 des Ranges $n \geq 2$ \"uber $K$. Ist dann $f$ eine nicht ausgeartete,
 symmetrische, spurwertige $\alpha$-Form auf $V$, deren Index mindestens $1$
 ist, so ist $\PSU(V,f)$ einfach,
 es sei denn, es ist $n = 2$ und $K = \GF(4)$ oder $\GF(9)$ oder es ist
 $n = 3$ und $K = \GF(4)$.}
 \smallskip
       Beweis. Es sei $\Omega$ die Menge der isotropen Punkte. Dann operiert
 $\PSU(V,f)$ nach Satz 9.9 auf $\Omega$ primitiv. Mit 10.5 folgt $\PSU(V,f)' =
 \PSU(V,f)$, wenn man von den Ausnahmef\"allen absieht. Ist $P \in \Omega$,
 so bilden die unit\"aren Elationen mit Zentrum $P$ einen abelschen
 Normalteiler von $\PSU(V,f)_P$. Dieser Normalteiler erzeugt zusammen mit
 seinen Konjugierten nach 10.4 die Gruppe $\PSU(V,f)$. Nach dem Satz III.2.1 von
 Iwasawa ist $\PSU(V,f)$ also einfach.
 \medskip
       Damit haben wir in allen uns interessierenden F\"allen die
 Einfachheit von $\PTU(V,f)$ nachgewiesen.
 \medskip
       Wie schon gesagt, werden wir uns zu einem sp\"ateren Zeitpunkt um
 die orthogonalen Gruppen k\"ummern.\footnote{Anmerkung der Herausgeber:
 Das Kapitel \"uber orthogonale Gruppen fehlt.}
 Dabei wird f\"ur sie das Glei\-che
 gelten wie f\"ur die unit\"aren Gruppen, dass wir n\"amlich die Gruppen
 vom Index 0 nicht untersuchen werden. Ihre Struktur ist in beiden F\"allen
 so eng mit der
 arithmetischen Struktur des Koordinatenk\"orpers verwoben, dass man
 tief in die Zahlentheorie eintauchen muss, will man Aussagen \"uber ihre
 Struktur gewinnen.


 \newpage
       
 \mychapter{VI}{Segresche Mannigfaltigkeiten}

 \noindent
 Seit meiner Komputeralgebrazeit liebe ich freie
 Konstruktionen. Ich definiere daher das Tensorprodukt zweier
 Moduln hier in der Allgemeinheit, in der es gemeinhin in der
 Algebra ben\"otigt wird, dh., allgemeiner als wir es brauchen
 werden, da man in dieser allgemeineren Situation nicht umhin
 kommt, sich freier Konstruktionen zu bedienen. Hat man dann
 Tensorprodukte so allgemein definiert, wie wir es tun werden, so
 stellt man fest, dass diese Allgemeinheit der Geometrie wiederum
 zugute kommt. Man kann n\"amlich Tensorprodukte gut dazu benutzen,
 Homomorphismen projektiver Verb\"ande zu be\-schrei\-ben. Darauf werden
 wir im zweiten Abschnitt dieses Kapitels eingehen. Dies ist das
 Extra, das dieses Kapitel birgt. Sein eigentliches Thema sind die
 Segreschen Mannigfaltigkeiten, von denen man in B\"uchern \"uber
 projektive Geometrie auch nur selten etwas erf\"ahrt.
 \par
       Investieren wir zun\"achst in weiteres Werkzeug.

 \mysection{1. Tensorprodukte}
 
\noindent
 Es sei $R$ ein Ring mit Eins. Ferner sei $M$ ein unit\"arer
 $R$-Rechtsmodul und $N$ ein unit\"arer $R$-Linksmodul. (Das Wort
 \anff unit\"ar`` werden wir uns im folgenden schenken). Es
 sei weiterhin $A$ eine abelsche Gruppe, deren Verkn\"upfung wir
 als Addition notieren. Die Abbildung $f$ von $M \times N$ in $A$
 hei\ss e genau dann {\it tensoriell\/}, wenn gilt:
 \item{a)} Es ist $f(m + m',n) = f(m,n) + f(m',n)$ f\"ur alle $m$, $m' \in M, n \in N$.
 \item{b)} Es ist $f(m,n + n') = f(m,n) + f(m,n')$ f\"ur alle $m \in M$ und alle
 $n, n' \in N$.
 \item{c)} Es ist $f(mr,n) = f(m,rn)$ f\"ur alle $m \in M$, f\"ur alle $n \in N$
 und alle $r \in R$.
 \par\noindent
       Wie es so h\"aufig geschieht, erscheint von $A$ nur die Verkn\"upfung
 explizit in dieser Definition.
 \par
       Zu den bisherigen Daten nehmen wir noch eine abelsche Gruppe $T$ hinzu
 sowie eine tensorielle Abbildung $\tau$ von $M \times N$ in $T$. Das Paar
 $(T,\tau)$ hei\ss e ein {\it Tensorprodukt\/} von $M$ mit $N$, falls gilt:
 \item{1)} Ist $A$ eine abelsche Gruppe und ist $f$ eine tensorielle Abbildung von
 $M \times N$ in $A$, so gibt es einen Homomorphismus $\varphi$ von $T$ in $A$
 mit $f = \varphi \tau$.
 \item{2)} Die Gruppe $T$ wird von der Menge $\{\tau (m,n) \mid m \in M,\ 
 n \in N\}$ erzeugt.
 \par
       Wer sich mit universellen Objekten auskennt, wird sich \"uber die
 n\"achsten beiden S\"atze nicht wundern.
 \medskip\noindent
 {\bf 1.1. Satz.} {\it Es sei $R$ ein Ring mit Eins. Ferner sei $M$ ein
 $R$-Rechts\-mo\-dul und $N$ ein $R$-Linksmodul. Schlie\ss lich sei $(T,\tau)$
 ein Tensorprodukt von $M$ mit $N$. Ist $A$ eine abelsche Gruppe und ist $f$ eine
 tensorielle Abbildung von $M \times N$ in $A$, so gibt es genau einen
 Homomorphismus $\varphi$ von $T$ in $A$ mit $f=\varphi \tau$.}
 \smallskip
       Beweis. Dass es ein solches $\varphi$ gibt, besagt die Definition
 des Tensorproduktes. Es sei $\psi$ ein zweiter Homomorphismus mit
 $f = \psi \tau$. Ist dann $(m,n) \in M \times N$, so ist
 $$ \varphi(\tau(m,n)) = f(m,n) = \psi(\tau(m,n)). $$
 Dies zeigt, dass $\varphi$ und $\psi$ auf einem Erzeugendensystem von $A$
 \"u\-ber\-ein\-stim\-men, so dass, wie behauptet, $\psi = \varphi$ ist.
 \medskip\noindent
 {\bf 1.2. Satz.} {\it Es sei $R$ ein Ring mit Eins. Ferner sei $M$
 ein $R$-Links\-mo\-dul und $N$ ein $R$-Rechtsmodul. Sind $(T, \tau)$
 und $(T', \tau')$ Tensorprodukte von $M$ mit $N$, so sind $(T,
 \tau)$ und $(T', \tau')$ isomorph.}
 \smallskip
       Beweis. Nach Voraussetzung gibt es einen Homomorphismus $\varphi$
 von $T$ in $T'$ mit $\tau' = \varphi\tau$ und einen Homomorphismus
 $\varphi'$ von $T'$ in $T$ mit $\tau = \varphi'\tau'$. Hieraus
 folgt
 $$ \tau' = \varphi\varphi'\tau'. $$
 Nach 1.1 ist daher $\varphi \varphi' = 1_{T'}$. Ebenso folgt, dass
 $\tau'\tau = 1_T$ ist. Daher ist $\varphi$ ein Isomorphismus von $T$ auf $T'$
 mit $\tau' = \varphi\tau$, so dass alles bewiesen ist.
 \medskip
       Alle bisherigen Anstrengungen w\"aren umsonst, g\"abe es keine
 Tensorprodukte. Doch der mit freien Konstruktionen Vertraute
 wei\ss\ na\-t\"ur\-lich, wie er vorzugehen hat, um die Existenz des
 Tensorproduktes si\-cher\-zu\-stel\-len. Wir nehmen an, dass der Leser
 zumindest wei\ss, wie man sich freie abelsche Gruppen verschafft.
 Wer keine Kenntnis der grundlegenden Konstruktionen freier Objekte
 hat, sei f\"ur diese auf mein Buch \anff Tools and Fundamental
 Constructions of Combinatorial Mathematics`` 
 verwiesen.

\medskip

\noindent
{\bf 1.3. Satz.}{\it Es sei $R$ ein Ring mit Eins. Ferner sei $M$
ein Rechts- und $N$ ein Linksmodul \"uber $R$. Schlie\ss lich sei
$F$ die freie abelsche Gruppe \"uber dem freien Erzeugendensystem
$M \times N$ und $T$ sei die von allen Elementen der Form 
$$(m+m', n) - (m;n)-(m',n)$$ bzw. 
$$(m,n+n')-(m,n)-(m,n')$$ bzw. $$(mr,n)-(m,rn)$$ erzeugte
Untergruppe von $F$. Setzt man  $$M \otimes_R N := F/T$$
und definiert man $\tau$ durch  $$\tau(m,n) := (m,n)+T,$$
so ist $(M \otimes_R N,\tau)$ ein Tensorprodukt von $M$ mit $N$.}

\smallskip

Beweis. Aus der Definition von $T$ folgt, dass $\tau$ eine
tensorielle Abbildung von $M \times N$ in $M \otimes_R N$ ist.
Weil $M \times N$ ein Erzeugendensystem von $F$ ist, wird $M
\otimes_R N$ von der Menge der $\tau (m,n)$ erzeugt. Es bleibt zu
zeigen, dass jede tensorielle Abbildung von $M \times N$ in eine
abelsche Gruppe sich durch $\tau$ faktorisieren l\"asst. Dies
ist aber auch banal. Jede faktorielle Abbildung $f$ von $M \times
N$ in eine abelsche Gruppe $A$ ist insbesondere eine Abbildung
jener Menge in $A$, so dass $f$ sich zu einem Homomorphismus $g$
von $F$ in $A$ fortsetzen l\"asst. Weil $f$ tensoriell ist, liegt
$T$ im Kern von $g$, so dass mit bekannten S\"atzen folgt, dass es
einen Homomorphismus $\varphi$ von $F/T = M \otimes_R N$ in $A$
gibt mit $f(m,n)=g(m,n)=(\varphi \tau)(m,n)$ f\"ur alle $(m,n) \in
M \times N$. Also ist $f=\varphi \tau$, was noch zu beweisen war.

\smallskip

Es ist \"ublich, das Bild von $(m,n)$ unter $\tau$ mit $m \otimes
n$ zu bezeichnen. Dieser Konvention werden wir uns im folgenden
anschlie\ss en.

Der n\"achste Satz ist ebenfalls von gro\ss er Bedeutung.

\medskip

\noindent
{\bf 1.4. Satz.} {\it Es sei $R$ ein Ring mit Eins. Ferner seien
$M$ und $M'$ zwei Rechts- und $N$ und $N'$ zwei Linksmoduln \"uber
$R$. Ist dann $f$ ein Homomorphismus von $M$ in $M'$ und $g$ ein
Homomorphismus von $N$ in $N'$, so gibt es genau einen
Homomorphismus, den wir mit $f \otimes g$ bezeichnen, von $M
\otimes_R N$ in $M' \otimes_R N'$ mit
$$ (f \otimes g)(m
\otimes n) = f(m) \otimes g(n) $$
f\"ur alle $(m,n) \in M \times N$.}
\smallskip
Beweis. Die durch $\sigma(m,n):=f(m) \otimes g(n)$ definierte
Abbildung $\sigma$ ist tensoriell, so dass die Existenz von $f
\otimes g$ aus der Definition des Tensorproduktes und die
Einzigkeit aus Satz 1.1 folgt.
\smallskip
Hat man drei $R$-Rechtsmoduln $M, M'$ und $M''$ sowie drei
$R$-Linksmoduln $N, N'$ und $N''$, ist $f$ ein Homomorphismus von
$M$ in $M'$ und $f'$ ein solcher von $M'$ in $M''$, ist ferner $g$
ein Homomorphismus und $f'$ ein solcher von $M'$ in $M''$, ist
ferner $g$ ein Homomorphismus von $N$ in $N'$ und $g'$ ein
Homomorphismus von $N'$ in $N''$, so ist  $$(f' \otimes
g')(f \otimes g) = (f' f) \otimes (g' g),$$ wie unschwer zu sehen
ist.
\par
Wie n\"utzlich Satz 1.4 ist, sieht man schon beim Beweise des
n\"achsten Satzes.

\medskip
\noindent
{\bf 1.5. Satz.} {\it Es sei $R$ ein Ring mit Eins und $M$ sei ein
$R$-Rechts- und $N$ ein $R$-Linksmodul. Ferner sei $(D_i\mid i \in I)$
eine Familie von Teilmoduln von $M$ und es gelte  $$M=
\bigoplus_{i \in I} D_i.$$ F\"ur $i \in I$ sei $\pi_i$ die
Projektion von $M$ auf $D_i$ mit $D_j \subseteq \Kern (\pi_i)$
f\"ur alle $j \in \I - \{i\}$. Setzt man $\Delta_i := (\pi_i
\otimes 1)(M \otimes_R N)$ und bezeichnet man mit $\otimes_i$ die
Einschr\"ankung von $D_i \times N$, so gilt  $$M \otimes_R
N = \bigoplus_{i \in I} \Delta_i$$ und  $(\Delta_i, \otimes_i)$
ist f\"ur alle $i \in I$ ein Tensorprodukt von $D_i$ mit $N$.}

\smallskip

Beweis. Es sei $(m,n) \in M \times N$. Es gibt dann eine endliche
Teilmenge $J$ von $I$ mit $m \in \sum_{j \in J} D_j$, dh., es ist
 $$m = \sum_{j \in J} m_j$$ mit $m_j \in D_j$. Weil nach
Konstruktion $\pi_i^2 = \pi_i$ f\"ur alle $i$ und $\pi_i \pi_j =
0$ f\"ur alle $i$ und $j$ mit $i \neq j$ gilt, folgt $m_j =
\pi_j(m)$ f\"ur alle $j \in J$ und damit  $$m = \sum_{j \in
J} \pi_j (m).$$ Hiermit folgt  $$m \otimes n = \sum_{j\in
J} (\pi_j (m) \otimes n) = \sum_{j \in J} (\pi_j \otimes 1)(m
\otimes n),$$ so dass $m \otimes n \in \sum_{i \in I} \Delta_i$
ist. Weil $M \otimes_R N$ von der Menge der $m \otimes n$ erzeugt
wird, folgt  $$M \otimes_R N= \sum_{i \in I} \Delta_i.$$ Es
ist  $$(\pi_i \otimes 1) (\pi_j  \otimes 1) = (\pi_i \pi_j)
\otimes 1,$$ so dass die $\pi_i \otimes 1$ paarweise orthogonale
Idempotente sind. Hieraus folgt, dass sogar  $$M \otimes_R
N= \bigoplus_{i \in I} \Delta_i$$ gilt.

Es bleibt zu zeigen, dass $(\Delta_i, \otimes_i)$ ein
Tensorprodukt von $D_i$ mit $N$ ist. Dazu sei $f$ eine tensorielle
Abbildung von $D_i \times N$ in eine abelsche Gruppe $A$. Wir
definieren dann eine Abbildung $F$ von $M \times N$ in $A$ durch
$$F(m,n) := f(\pi_i(m), n).$$ Offenbar ist auch $F$ tensoriell.
Es gibt also einen Homomorphismus $\varphi$ von $M \otimes_R N$ in
$A$ mit $F(m,n)=\varphi (m \otimes n)$. Setze $\psi :=\varphi
(\pi_i \otimes 1)$. Ist dann $(y,n) \in D_i \times N$, so gilt
$$\eqalign{
f(y,n) &= F(y,n)=\varphi(y \otimes n)\cr &=\varphi(\pi_i(y)
\otimes n) = \varphi ((\pi_i \otimes 1)(y \otimes n))\cr &=\psi (y
\otimes_i n).}$$
Schlie\ss lich ist klar, dass $\Delta_i$ von der Menge der $y
\otimes_i n$ erzeugt wird. Damit ist alles bewiesen.

\medskip

Dieser Satz findet sich in der Literatur meist wie folgt
formuliert: Unter den gemachten Voraussetzungen gilt  $$M
\otimes_R N \cong \bigoplus_{i \in I} (D_i \otimes N).$$ Dies ist
zwar richtig --- und wir werden ihn auch wohl immer in dieser
Formulierung benutzen ---, doch ist dies viel weniger
aussagekr\"aftig als obige Formulierung.

Wir m\"ochten Tensorprodukte von Vektorr\"aumen studieren und
dabei erreichen, dass die betrachteten Tensorprodukte ebenfalls
Vektorr\"aume sind. Um dorthin zu gelangen, definieren wir
zun\"achst den Begriff des $(R,S)$-Bimoduls. Es seien $R$ und $S$
Ringe mit Eins. Die abelsche Gruppe hei\ss e genau dann
$(R,S)$-\emph{Bimodul}, falls $M$ ein $R$-Links- und ein
$S$-Rechtsmodul ist und dar\"uber hinaus $r(ms)=(rm)s$ f\"ur alle
$r \in R$, alle $m \in M$ und alle $s \in S$ gilt.

\medskip
\noindent
{\bf 1.6. Satz.} {\it Es seien $R$ und $S$ zwei Ringe mit Eins.
Ferner sei $M$ ein $R$-Rechtsmodul und $N$ ein $(R,S)$ Bimodul.
Dann tr\"agt $M \otimes_R N$ genau eine $S$-Rechtsmodulstruktur
mit $(m \times n) s=m \times ns$ f\"ur alle $m \in M$, alle $n \in
N$ und alle $s \in S$.}

\smallskip

Beweis. Die Einzigkeit folgt wieder daraus, dass die Gesamtheit
der $m \otimes n$ das Tensorprodukt $M \otimes_R N$ erzeugen.

Es sei $s \in S$. Wir definieren die Abbildung $f$, von $M \otimes
N$ in $M \otimes_R N$ durch  $$f_s(m,n) := m \otimes ns.$$
Banale Rechnungen zeigen, dass $f_s$ tensoriell ist. Es gibt also
einen Endomorphismus $\varphi_s$ von $M \otimes_R N$ mit 
$$\varphi_s(m \otimes n) = m \otimes ns$$ f\"ur alle $(m,n) \in M
\times N$. Definiert man nun $ys$ f\"ur $y \in M \otimes_R N$
durch $ys := \varphi_s(y)$, so zeigen Routinerechnungen, dass $M
\otimes_R N$ auf diese Weise zu einem $S$-Rechtsmodul wird, f\"ur
den \"uberdies $(m \otimes n)s = m \otimes ns$ f\"ur alle $m \in
M$, alle $n \in N$ und alle $s \in S$ gilt.

Jeder Ring $R$ ist nat\"urlich ein $(R,R)$-Bimodul. Also tr\"agt
$M \otimes_R R$ gem\"a\ss\ dem gerade bewiesenen Satz eine
Struktur als $R$-Rechtsmodul, falls $M$ ein $R$-Rechtsmodul ist.
\"Uber diesen Modul gilt der folgende Satz.

\medskip
\noindent
{\bf 1.7. Satz.} {\it  Es sei $R$ ein Ring mit Eins und $M$ sei
ein $R$-Rechtsmodul. Es gibt dann einen Isomorphismus $\sigma$ des
$R$-Rechtsmoduls $M \otimes_R R$ auf $M$ mit $\sigma(m \otimes r)
= mr$ f\"ur alle $m \in M$ und alle $r \in R$.}

\smallskip

Beweis. Die Abbildung, die $(m,r)$ auf $mr$ abbildet, ist
tensoriell. Es gibt daher einen Homomorphismus $\sigma$ der
abelschen Gruppe $M \otimes_R R$ in $R$ mit $\sigma (m \otimes
r)=mr$. Hieraus folgt unmittelbar, dass $\sigma$ sogar ein
Modulhomomorphismus ist. Es bleibt zu zeigen, dass $\sigma$
bijektiv ist.

Ist $m \in M$, so sei $\tau(m)$ durch $\tau(m) := m \otimes 1$
gegeben. Die Abbildung $\tau$ ist sicherlich additiv. Wegen
$$\tau(mr) = mr \otimes 1=m \otimes r=(m \otimes 1)r$$ ist
$\tau$ ein Modulhomomorphismus. Nun ist aber  $$\sigma \tau
(m) = m1 = m$$ und  $$\tau \sigma (m,r) = mr \otimes 1 =m
\otimes r.$$ Hieraus folgt, dass $\sigma$ bijektiv und dass $\tau$
die zu $\sigma$ inverse Abbildung ist.

\smallskip

F\"ur diesen Satz braucht man zum ersten Male, dass $R$ eine Eins
hat. Zuvor war diese Annahme immer \"uberfl\"ussig.

Es sei $R$ ein Ring mit Eins und $M$ sei ein Rechts- und $N$ ein
Linksmodul \"uber $R$. Ist $D$ ein direkter Summand von $M$, so
zeigt Satz 1.5, dass $M \otimes_R N$ einen zu $D \otimes_R N$
einen zu $D \otimes_R N$ isomorphen direkten Summanden besitzt.
Ist $D$ nur ein Teilmodul von $M$, so kann man nicht schlie\ss en,
dass $D \otimes_R N$ zu einem Teilmodul von $M \otimes_R N$
isomorph ist. Um dies zu belegen, sei $Z$ der Ring der ganzen
Zahlen, $Z_2$ die zyklische Gruppe der Ordnung $2$ und $Q$ der
K\"orper der rationalen Zahlen. Nach Satz 1.7 ist dann
\emph{mutatis mutandis}  $$Z \otimes_Z Z_2 \cong Z_2.$$
Andererseits ist  $$Q \otimes_Z Z_2 = \{0\},$$ wie wir
jetzt zeigen werden. Ist n\"amlich $z \in Z_2$ und $r \in Q$, so
folgt  $$r \otimes z = \frac{r}{2} \otimes 2 z = 0$$ und
damit die Behauptung.

In Satz 1.6 haben wir gesehen, wie man aus dem Tensorprodukt $M
\otimes_R N$ einen $S$-Rechtsmodul macht, wenn $N$ ein
$(R,S)$-Bimodul ist. Ganz analog sieht man, dass $M \otimes_R N$
genau eine $S$-Linksmodulstruktur mit $$s(m \otimes n) = sm
\otimes n$$ f\"ur alle fraglichen $s, m$ und $n$ tr\"agt, wenn nur
$M$ ein $(S,R)$-Bimodul und $N$ ein $R$-Linksmodul ist. Diese
Rechts- bzw. Linksmodulstrukturen sind im folgenden gemeint, wenn
davon die Rede ist, dass $M \otimes_R N$ eine Rechts- bzw.
Linksmodulstruktur tr\"agt.

Dies bemerkt, formuliert und beweist man den folgenden Satz.

\medskip
\noindent
{\bf 1.8. Satz.} {\it Es seien $R$ und $S$ zwei Ringe mit Eins.
Ferner sei $M$ ein $R$-Rechtsmodul, $N$ ein $(R,S)$-Bimodul und
$O$ ein $S$-Linksmodul. Dann ist $M \otimes_R N$ ein
$S$-Rechtsmodul und $N \otimes O$ ein $R$-Linksmodul, so dass die
Tensorprodukte $M \otimes_R(N \otimes_S O)$ und $(M \otimes_R N)
\otimes_S O$ definiert sind. Dar\"uber hinaus gibt es einen
Isomorphismus $\sigma$ von  $$M \otimes_R(N\otimes_S O)$$
auf  $$(M \otimes_R N) \otimes_S O$$ mit  $$\sigma(m
\otimes(n \otimes o))= (m \otimes n) \otimes o$$ f\"ur alle
$(m,n,o) \in M \times N \times O.$}

\smallskip

Beweis. Es sei $m \in M$. Wir definieren die Abbildung $a_m$ von
$N$ in $M \otimes_R N$ durch $a_m(x) := m \otimes x$ f\"ur alle $x
\in N$. Dann ist $a_m$ ein Homomorphismus des $S$-Rechtsmoduls $N$
in den $S$-Rechtsmodul $M \otimes_R N$. Wir setzen $b_m := a_m
\otimes 1_O$. Dann ist $b_m$ ein Homomorphismus von $N \otimes_S
O$ in $(M \otimes_R N) \otimes_S O$. Wir definieren nun die
Abbildung $f$ von $M \times (N \otimes_S O)$ in $(M \otimes_R N)
\otimes_S O$ durch  $$f(m,z) := b_m(z)$$ f\"ur alle $m \in
M$ und alle $z \in N \otimes_S O$. Die Abbildung $f$ ist offenbar
in beiden Argumenten additiv. Es gilt aber auch $f(mr,z) =
f(m,rz)$ f\"ur alle in Frage kommenden $m,r$ und $z$. Um dies zu
beweisen, d\"urfen wir annehmen, dass $z = n \otimes o$ ist. Dann
ist
$$\eqalign{
f(mr,z) &= (a_{mr} \otimes 1_O)(n \otimes o) = (mr \otimes n)
\otimes o = (m \otimes rn) \otimes o\cr &= (a_m \otimes 1_O)(rn
\otimes o) = f(m,rn \otimes o) = f(m,rz).}
$$
Damit ist gezeigt, dass $f$ tensoriell ist. Es gibt also einen
Homomorphismus $\sigma$ von $M \otimes_R(N \otimes_S O)$ in $(M
\otimes_R N) \otimes_S O$ mit $f(m,z)=\sigma (m \otimes z)$ f\"ur
alle $m \in M$ und alle $z \in N \otimes_S O$. Hieraus folgt
insbesondere, dass  $$\sigma (m \otimes (n \otimes o)) = (m
\otimes n) \otimes o$$ f\"ur alle in Frage kommenden $m,n$ und $o$
gilt. Nach dem bislang Bewiesenen wird es dem Leser nicht schwer
fallen zu zeigen, dass es einen Homomorphismus $\tau$ von $(M
\otimes_R N) \otimes_S O$ in $M \otimes_R (N \otimes_S O)$ gibt
mit  $$\tau ((m \otimes n) \otimes o) = m \otimes (n
\otimes o)$$ f\"ur alle $m,n$ und $o$. Hieraus folgt schlie\ss
lich, dass $\sigma$ und $\tau$ invers zueinander sind, so dass
$\sigma$ in der Tat ein Isomorphismus ist.

\medskip

Falls es dem Leser wider Erwarten nicht gelungen sein sollte, die
Existenz von $\tau$ nachzuweisen, so findet er einen solchen
Nachweis bei Bourbaki, der es allerdings dem Leser \"uberl\"a\ss
t, die Existenz von $\sigma$ zu beweisen.

Es sei $K$ ein kommutativer K\"orper und $V$ sei ein
Rechtsvektorraum \"uber $K$. Ist dann $k \in K$ und $v \in V$, so
setzen wir $kv := vk$. Weil $K$ kommutativ ist, wird $V$ auf diese
Weise zu einem $K$-Linksvektorraum. Es gilt sogar, dass $V$ ein
$(K,K)$-Bivektorraum ist. Sind n\"amlich $k,l \in K$ und $v \in
V$, so ist  $$k(vl)=(vl)k=v(lk)=v(kl)=(vk)l=(kv)l.$$
Sprechen wir in Zukunft von einem Vektorraum $V$ \"uber einem
kommutativen K\"orper $K$, so verstehen wir darunter immer einen
$(K,K)$-Bivektorraum mit $kv = vk$ f\"ur alle $v \in V$ und alle
$k \in K$. Mit dieser Verabredung wird dann auch das Tensorprodukt
zweier Vektorr\"aume \"uber einem kommutativen K\"orper zu einem
Vektorraum \"uber eben diesem K\"orper und die Verabredung, die
wir \"uber gegebene Vektorr\"aume \"uber kommutativen K\"orpern
getroffen haben, dass n\"amlich $kv=vk$ f\"ur alle in Frage
kommenden $k$ und $v$ gilt, gilt auch f\"ur dass Tensorprodukt der
beiden Vektorr\"aume, wie Satz 1.6 zeigt.

Macht man solche Generalvoraussetzungen, wie wir es gerade taten,
so muss man nat\"urlich vorsichtig sein, da man m\"oglicherweise
ein Objekt von der Art konstruiert, \"uber die man die
Generalvoraussetzung gemacht hat. Dann muss man zeigen, dass das
Objekt auch die gew\"unschten Eigenschaften hat, was sicher nicht
immer zutrifft. Das erinnert mich an eine Episode aus der Zeit, da
meine Kinder zur Schule gingen. Meine Tochter Suzanne kam eines
Tages zu mir, der ich Zeitung lesend auf dem Sofa lag, und sagte:
\glqq Papa, wir schreiben morgen eine Mathe-Arbeit. Kannst du mich
abfragen?\grqq{} Ich konnte: Sie gab mir ihr Heft: Ich schaute in
ihr Heft und tat das, was ich immer in zweifelhaften Situation
tue, ich sagte: \glqq Ach du meine G\"ute!\grqq, und fuhr fort,
\glqq Gib mir, bitte, dein Buch.\grqq{} Im Buch stand der gleiche
Unfug. Es brachte eine zweifelhafte Definition von Kardinalzahl
und fuhr fort, dass jede Menge eine Kardinalzahl habe und dass die
Menge der Kardinalzahlen die Menge der nat\"urlichen Zahlen sei.
Das war f\"ur mich ein gefundenes Fressen. Ich fragte sie also:
\glqq Jede Menge hat eine Kardinalzahl?\grqq \glqq Ja.\grqq \glqq
Also auch die Menge der nat\"urlichen Zahlen?\grqq \glqq Ja.\grqq
\glqq Jede Kardinalzahl ist eine nat\"urliche Zahl?\grqq \glqq
Ja.\grqq \glqq Ist 10395 die Kardinalzahl der Menge der
nat\"urlichen Zahlen?\grqq \glqq Nein.\grqq \glqq Ist 23576 die
Kardinalzahl der Menge der nat\"urlichen Zahlen?\grqq \glqq
Nein.\grqq \glqq Ja, was ist denn die Kardinalzahl der Menge der
nat\"urlichen Zahlen?\grqq \glqq Unendlich.\grqq \glqq Ist
Unendlich eine nat\"urliche Zahl?\grqq \glqq Nein.\grqq \glqq Was
ist denn Unendlich?\grqq \glqq So'n Symbol.\grqq{} Ich rief also
den Lehrer an und \"au\ss erte ihm meine Bedenken. Unter anderem
erw\"ahnte ich, dass es neben den nat\"urlichen Zahlen ja auch
noch andere Kardinalzahlen g\"abe. Darauf der Lehrer: \glqq Das
ist nat\"urlich richtig, aber die Mengen, die wir betrachten, sind
alle endlich.\grqq \glqq So!\grqq, sagte ich, \glqq Und die Menge
der nat\"urlichen Zahlen?\grqq \glqq Dar\"uber habe ich noch nicht
nachgedacht.\grqq{} Die Arbeit wurde erst am \"ubern\"achsten Tag
geschrieben. Bei dem Buch handelte es sich um den
Lambacher--Schweitzer, den Namen des Lehrers zu nennen verbietet
des S\"angers H\"oflichkeit.

\smallskip

Doch zur\"uck zu unserem Thema.

\medskip
\noindent
{\bf 1.9. Satz.} {\it Es sei $K$ ein kommutativer K\"orper und $V$
und $W$ seien zwei Vektorr\"aume \"uber $K$. Ist $B$ eine Basis
von $V$ und $C$ eine Basis von $W$, so ist  $$\{b \otimes
c\mid b \in B, c \in C\}$$ eine Basis von $V \otimes_K W$.}

\smallskip

Beweis. Um diesen Satz zu beweisen, kann man sich nat\"urlich des
Satzes 1.5 bedienen. Da Wiederholung f\"ur den Lernenden aber
n\"utzlich ist, f\"uhren wir den Beweis des vorliegenden Satzes
unabh\"angig von Satz 1.5 aus, auch wenn wir daf\"ur einiges
doppelt machen.

Es sei $v \in V$ und $w \in W$. Es gibt dann $b_1, \dots, b_m \in
B$ und $k_1, \dots, k_m \in K$ sowie $c_1, \dots, c_n \in C$ und
$l_1, \dots, l_n \in K$ mit $v = \sum^{m}_{i:=1} k_i b_i$ und $x =
\sum^n_{j:=1} l_jc_j$. Es folgt  $$v \otimes w =
\sum^m_{i:=1} \sum^n_{j:=1} k_i l_j (b_i \otimes c_j).$$ Weil $V
\otimes_K W$ als abelsche Gruppe von der Menge der $v \otimes w$
erzeugt wird, wird $V \otimes_K W$ auch als Vektorraum von dieser
Menge erzeugt. Daher ist die Menge der $b \otimes c$ ein
Erzeugendensystem des $K$-Vektorraumes $V \otimes_K W$.

Es sei $b \in B$ und $c \in C$. Wir definieren die linearen
Abbildungen $\pi$ und $\rho$ von $V$ bzw. $W$ in $K$ durch $\pi(b)
:= 1$ und $\pi (x) := 0$ f\"ur $x \in B - \{b\}$, bzw., $\rho (c)
:= 1$ und $\rho (y) := 0$ $y \in C - \{c\}$. Ferner definieren wir
die Abbildung $f$ von $V \times W$ in $K$ durch $f(v,w) := \pi (v)
\rho (w)$. Dann ist $f$ tensoriell, so dass es eine lineare
Abbildung $\varphi$ von $V \otimes_K W$ in $K$ gibt mit $$\varphi
(v \otimes w) = \pi (v) \rho (w)$$ f\"ur alle $v \in V$ und alle
$w \in W$. Es folgt $\varphi (b \otimes c) = 1$ und $\varphi (x
\otimes y) = 0$ f\"ur alle $(x,y) \in B \times C$ --- $\{(b,c)\}$,
so dass $b \otimes c$ von $B \times C$ --- $\{(b,c)\}$ linear
unabh\"angig ist. Damit ist der Satz bewiesen.

\medskip
\noindent
{\bf 1.10. Korollar.} {\it Es sei $K$ ein kommutativer K\"orper.
Sind $V$ und $W$ Vektorr\"aume endlichen Ranges \"uber $K$, so hat
auch $V \otimes_K W$ endlichen Rang und es gilt}  $$\Rg_K (V
\otimes_K W) = \Rg_K (V) \Rg_K (W).$$

Das Argument, welches wir
zum Beweise von 1.9 verwandten, l\"asst sich noch einmal
verwenden, was daraufhin deutet, dass eigentlich ein Satz zu
formulieren sei, der dieses Argument institutionalisiert.

\medskip
\noindent
{\bf 1.11. Satz.} {\it Es sei $K$ ein kommutativer K\"orper und
$V$ und $W$ seien Vektorr\"aume \"uber $K$. Sind $b_1, \dots, b_n$
linear unabh\"angige Vektoren aus $V$, sind $y_1, \dots, y_n \in
W$ und gilt  $$\sum^n_{i:=1} (b_i \otimes y_i) = 0,$$ so
ist $y_i = 0$ f\"ur $i:= 1, \dots, n$.}

\smallskip

Beweis. Es sei $i \in \{1, \dots, n\}$. Es sei $i \in \{1, \dots,
n\}$. Es gibt dann eine lineare Abbildung $f$ von $V$ in $K$ mit
$f(b_i)=1$ und $f(b_j)=0$ f\"ur alle von $i$ verschiedenen $j$.
Ist $g$ irgendeine lineare Abbildung von $W$ in $K$, so ist die
durch $\varphi (v,w):= f(v) g(w)$ definierte Abbildung $\varphi$
tensoriell. Es gibt daher eine lineare Abbildung $\pi$ von $V
\otimes_K W$ in $K$ mit $\pi(v \otimes w) = f(v) g(w)$. Hieraus
folgt  $$0=\pi \Bigg( \sum^n_{j:=1} b_j \otimes y_j \Bigg)
= \sum^n_{i:=1} f(b_j)g(y_j)=g(y_i).$$ Also ist $g(y_i)=0$ f\"ur
alle linearen Abbildungen $g$ von $W$ in $K$, so dass, wie
behauptet, $y_i=0$ ist.

\mysection{2. Homomorphismen projektiver R\"aume}

\noindent
Da wir Tensorprodukte allgemeiner definiert haben, als wir sie
eigentlich brau\-chen, werden wir von unserem Thema abweichen und
zeigen, dass die allgemeinere Version des Tensorproduktes auch
f\"ur die projektive Geometrie von Nutzen ist.

Ist $K$ ein K\"orper und ist $R$ ein Teilring von $K$ mit $1 \in
R$, so ist $K$ nat\"urlich ein $(R,K)$-Bimodul. Ist $M$ ein
$R$-Rechtsmodul, so tr\"agt $M \otimes_R K$ nach 1.6 daher eine
Struktur als $K$-Rechtsvektorraum, die mit der Struktur des
Tensorproduktes vertr\"aglich ist. Dies machen wir uns jetzt zu
Nutze.

\medskip
\noindent
{\bf 2.1. Satz.} {\it Es sei $K$ ein K\"orper und $R$ sei ein
Teilring von $K$ mit $1 \in R$. Ist $M$ ein $R$-Rechtsmodul, so
tr\"agt $M \otimes_R K$ genau eine $K$-Rechtsvektorraumstruktur
mit  $$(m \otimes k)l = m \otimes (kl)$$ f\"ur alle $m \in
M$ und alle $k,l \in K$. Ist $M$ ein freier $R$-Modul und ist $B$
eine $R$-Basis von $M$, so ist  $$\{b \otimes 1\mid b \in B\}$$
eine $K$-Basis von $M \otimes_R K$.}

\smallskip

Beweis. Die erste Aussage folgt, wie schon bemerkt, aus Satz 1.6.

Weil $B$ eine Basis von $M$ ist, ist $M = \oplus_{b \in B} b R$.
F\"ur $b \in B$ sei $\pi_b$ die durch $\pi_b (b) := b$ und
$\pi_b(c) := 0$ f\"ur $c \in B - \{b\}$ definierte Projektion von
$M$ auf $bR$. Nach 1.5 ist dann  $$M \otimes_R K=
\bigoplus_{b \in B} (\pi_b \otimes 1_K)(M \otimes_R K).$$ Weil $B$
eine Basis von $M$ ist, wird $M \otimes_R K$ von der Menge der
Elemente $c \otimes k$ mit $c \in B$ und $k \in K$ erzeugt. Daher
wird  $$(\pi_b \otimes 1_K)(M \otimes_R K)$$ von den
Bildern dieser Elemente unter $\pi_b \otimes 1_K$ erzeugt. Ist nun
$c \in B - \{b\}$, so ist $(\pi_b \otimes 1_K)(c \otimes k) = 0$.
Ferner ist  $$(\pi_b \otimes 1_K)(b \otimes k) = (b \otimes
1)k.$$ Damit ist gezeigt, dass  $$(\pi_b \otimes 1_K)(M
\otimes_R K) = (b \otimes 1)K$$ ist. Um zu zeigen, dass $\{b
\otimes 1\mid b \in B\}$ eine Basis ist, ist also nur noch zu zeigen,
dass die direkten Summanden allesamt von $\{0\}$ verschieden sind.

Zu diesem Zweck zitieren wir zun\"achst noch einmal 1.5. Nach
diesem Satz ist $(\pi_b \otimes 1_K)(M \otimes_R K)$ zu $b R
\otimes_R K$ isomorph. Da die durch $\beta(br, k) := rk$
definierte Abbildung $\beta$ tensoriell ist, gibt es einen
Homomorphismus $\eta$ von $b R \otimes_R K$ in $K$ mit $\eta(br
\otimes k) = rk$ f\"ur alle $r \in R$ und alle $k \in K$. Es folgt
$\eta (b \otimes 1) = 1$, so dass $(bR \otimes_R) K)$ nicht nur
aus der Null besteht. Damit ist alles bewiesen.

\medskip

Der Leser beachte, dass $b \otimes 1$ im letzten Absatz ein
Element aus $bR \otimes_R K$ bezeichnet, w\"ahrend zuvor ein
Element aus $M \otimes_R K$ mit diesem Ausdruck gemeint war.

\medskip
\noindent
{\bf 2.2. Korollar.} {\it Es sei $K$ ein K\"orper und $R$ sei ein
Teilring von $K$ mit $1 \in R$. Ist dann $F$ ein freier Modul
\"uber $R$, ist $f \in F$ und gilt $f \otimes 1 = 0$, so ist $f =
0$.}

\smallskip

Beweis. Weil $F$ ein freier $R$-Modul ist, gibt es eine Basis $B$
von $F$. Es gibt dann weiter $b_1, \dots, b_m \in B$ und $r_1,
\dots, r_m \in R$ mit $f = \sum^m_{i:=1} b_i r_i$. Es folgt
$$0 = f \otimes 1= \sum^m_{i:=1} (b_i r_i \otimes 1)=
\sum^m_{i:=1} (b_i \otimes 1) r_i,$$ so dass nach 2.1 gilt, dass
$r_i = 0$ ist f\"ur alle $i$. Somit ist $f = 0$, wie behauptet.

\medskip
\noindent
{\bf 2.3. Satz.} {\it Es sei $K$ ein K\"orper und $R$ sei ein
Teilring von $K$ mit $1 \in R$. Ferner sei $F$ ein freier
Rechtsmodul \"uber $R$. Ist $U \in \La_K(F \otimes_R K)$ und ist
 $$D := \{f\mid f \in F, f \otimes 1 \in U\},$$ so gilt: Ist $f
\in F$ und $r \in R^*$ und gilt $fr \in D$, so ist $f \in D$,
m.a.W., der Faktormodul $F/D$ ist \emph{torsionsfrei}.}

\smallskip

Beweis. Wegen $fr \in D$ ist $fr \otimes 1 \in U$. Es folgt
$$f \otimes 1 = (f \otimes r)r^{-1} = (fr \otimes 1)^{-1}
\in U,$$ so dass in der Tat $f \in D$ gilt.

\smallskip

Der Modul $D$ kann nat\"urlich, gemessen an $U$, sehr klein sein.
Wir k\"ummern uns nun zun\"achst um Bedingungen, die erzwingen,
dass $D$ so gro\ss\ wie m\"oglich wird.

Die Bemerkung, dass im folgenden Satz a) eine Konsequenz von d)
ist, verdanke ich Herrn U.~Dempwolff. Wir ben\"otigen des weiteren
jedoch nur, dass b) aus a) folgt.

\medskip
\noindent
{\bf 2.4. Satz.} {\it Ist $K$ ein K\"orper und ist $R$ ein
Teilring von $K$ mit $1 \in R$, so sind die folgenden Aussagen
\"aquivalent:
\item{a)} Zu jeder endlichen Teilmenge $E$ von $K$ gibt es ein $s
\in R^*$ mit $es \in R$ f\"ur alle $e \in E$.
\item{b)} Ist $M$ ein Rechtsmodul \"uber $R$, sind $f_1, \dots,
f_n \in M$ und $k_1, \dots, k_n \in K$, so gibt es $r_1, \dots,
r_n \in R$ und ein $y \in K^*$ mit  $$\sum^n_{i:=1} (f_i
\otimes k_i) = \Bigg( \sum^n_{i:=1} f_i r_i \Bigg) \otimes y$$ und
der weiteren Eigenschaft, dass $r_i = 0$ genau dann gilt, wenn
$k_i = 0$ ist.
\item{c)} Ist $M$ ein Rechtsmodul \"uber $R$ und ist $v \in M
\otimes_R K$, so gibt es ein $m \in M$ und ein $k \in K^*$ mit $v
= m \otimes k$.
\item{d)} Es gibt einen nicht endlich erzeugten freien Rechtsmodul
$F$ \"uber $R$, so dass es zu jedem $v \in F \otimes_R K$ ein $f
\in F$ und ein $k \in K^*$ gibt mit $v = f \otimes k$.\par}

\smallskip

Beweis. a) impliziert b): Nach Voraussetzung gibt es ein $s \in
R^*$ mit $k_i s \in R$ f\"ur $i := 1, \dots, n$. Setze $r_i :=
k_is$ und $y:=s^{-1}$. Dann gilt in der Tat 
$$\sum^n_{i:=1} (f_i \otimes k_i) = \Bigg( \sum^n_{i:=1} f_i r_i
\Bigg) \otimes y.$$ Wegen $r_i = k_i s$ und $s \neq 0$ gilt auch
die Aussage \"uber das Verschwinden der $r_i$.

Es ist trivial, dass c) eine Folge von b) und dass d) eine Folge
von c) ist.

d) impliziert a): Es sei $B$ eine Basis von $F$ und $E$ sei eine
endliche Teilmenge von $K$. Wir d\"urfen annehmen, dass $1 \in E$
ist. Weil $F$ nicht endlich erzeugt ist, ist $B$ nicht endlich, so
dass es eine injektive Abbildung $b$ von $E$ in $B$ gibt. Setze
 $$v := \sum_{e \in E} (b_e \otimes e).$$ Es gibt dann ein
$f \in F$ und ein $k \in K^*$ mit $v = f \otimes k$. Es sei $f =
\sum_{c \in B} cr(c)$. Dann ist  $$v = f \otimes k =
\sum_{c\in B}(c \otimes 1) r (c) k.$$ Andererseits ist 
$$v= \sum_{e \in E} (b_e \otimes 1)e.$$ Nach 2.1 gilt folglich $e
= r(b_e)k$ f\"ur alle $e \in E$. Dann ist aber $ek^{-1} \in R$
f\"ur alle $e \in E$. Wegen $1 \in E$ folgt $k^{-1} = 1k^{-1} \in
R$. Damit ist alles bewiesen.

\smallskip

Wir treffen des weiteren die Verabredung, mit $\Delta_R (M)$ die
Menge aller Teilmoduln $N$ des $R$-Rechtsmoduls $M$ zu bezeichnen,
f\"ur die $M/N$ torsionsfrei ist. Ist $\Phi$ eine Teilmenge von
$\Delta_R(M)$, so ist auch  $$\bigcap_{X \in\Phi} X \in
\Delta_R(M).$$ Nach I.2.1 ist $(\Delta_R(M), \subseteq)$ also ein
vollst\"andiger Verband. Zu bemerken ist, dass die obere Grenze
zweier Elemente aus $\Delta_R(M)$ auch in gut sich stellenden
F\"allen meist nicht die Summe dieser beiden Elemente ist. (Diese
geschraubt klingende Formulierung habe ich gew\"ahlt, um darauf
hinweisen zu k\"onnen, dass \glqq ill posed problems\grqq{} keine
\glqq schlecht gestellten Probleme\grqq --- so etwas gibt es nicht
---, sondern \glqq schlecht sich stellende Probleme\grqq{} sind.
Wer von schlecht gestellten Problemen redet, darf sich nicht
\"uber das Deutsch unserer Studenten beklagen.)

\medskip
\noindent
{\bf 2.5. Satz.} {\it Es sei $K$ ein K\"orper und $R$ sei ein
Teilring von $K$ mit $1 \in R$. Ferner gebe es zu jeder endlichen
Teilmenge $E$ von $K$ ein $r \in R^*$ mit $er \in R$ f\"ur alle $e
\in E$. Schlie\ss lich sei $F$ ein freier $R$-Rechtsmodul.
Definiert man die Abbildung $\varphi$ von $\Delta_R(F)$ in
$\La_K(F \otimes_R K)$ durch  $$\varphi (D) := \sum_{f \in
D} (f \otimes 1) K$$ f\"ur alle $D \in \Delta_R (F)$, so ist
$\varphi$ ein Isomorphismus von $(\Delta_R(F), \subseteq)$ auf
$(\La_K (F \otimes_R K), \leq)$.}

\smallskip

Beweis. Setze $V := F \otimes_R K$. F\"ur $U \in \La_K (V)$ werde
$\psi (U)$ definiert durch  $$\psi (U) := \{f\mid f \in F, f
\otimes 1 \in U\}.$$ Zun\"achst folgt  $$D \subseteq \psi
\varphi (D)$$ f\"ur alle $D \in \Delta_R (F)$.

Es sei $D \in \Delta_R(F)$. Ferner sei $f \in \psi \varphi (D)$.
Es gibt dann $d_1, \dots, d_n \in D$ und $k_1, \dots, k_n \in K$
mit   $$f \otimes 1 = \sum^n_{i := 1} (d_i \otimes k_i).$$
Nach 2.4 gibt es $r_1, \dots, r_n \in R$ und ein $y \in K^*$ mit
 $$f \otimes 1 = \Bigg( \sum^n_{i := 1} d_i r_i \Bigg)
\otimes y.$$ Nach Voraussetzung gibt es ein $a \in R^*$ mit $ya
\in R$. Es folgt  $$fa \otimes 1 = \Bigg( \sum^n_{i:=1} d_i
r_i ya \Bigg) \otimes 1.$$ Mit 2.2 erhalten wir folglich 
$$fa = \sum^n_{i:=1} d_ir_i ya \in D.$$ Wegen $a \neq 0$ und $D
\in \Delta_R (F)$ folgt daher $f \in D$. Somit ist $D =
\psi\varphi (D)$.

Es sei $U \in \La_K (V)$. Dann ist  $$\varphi \psi (U) =
\sum_{f \in F, f \otimes 1 \in U} (f \otimes 1) K.$$ Banal ist,
dass $\varphi \psi (U) \leq U$ gilt. Ist andererseits $u \in U$,
so folgt mit 2.4 die Existenz von $f \in F$ und $k \in K^*$ mit $u
= f \otimes k$. Hieraus folgt zun\"achst $f \otimes 1 =uk^{-1}
\in U$ und dann $u = (f \otimes 1)k \in \varphi \psi (U)$, so dass
$U = \varphi\psi (U)$ ist. Damit ist gezeigt, dass $\varphi$ eine
Bijektion von $\Delta_R (F)$ auf $\La_K (V)$ ist.

Sind $D, D' \in \Delta_R (F)$ und ist $D \subseteq D'$, so ist
nat\"urlich $\varphi (D) \leq \varphi (D')$. Sind $U, U' \in \La_K
(V)$ und ist $U \leq U'$, so folgt genauso einfach die
G\"ultigkeit der Inklusion $\psi(U) \subseteq \psi(U')$. Damit ist
der Satz bewiesen.

\medskip

Dieser Satz ist bestm\"oglich, wie wir noch sehen werden.
Andererseits kann man mehr beweisen, wenn der Rang von $F$ endlich
und $R$ etwa ein Hauptidealbereich ist. Davon sp\"ater mehr.
Zun\"achst wollen wir nach Beispielen von Ringen suchen, die die
Voraussetzungen von 2.5 erf\"ullen.

Ist $K$  ein kommutativer K\"orper und ist $R$ ein Teilring von
$K$ mit $1 \in R$, so erf\"ullt dieses Paar genau dann die
Voraussetzungen von 2.5, wenn $K$ der Quotientenk\"orper von $R$
ist, so dass dieser Satz also von realen Verh\"altnissen handelt.

Eine weitere Situation, die im Kommutativen ein Spezialfall der
gerade beschriebenen Situation ist, ist die folgende: Es sei $K$
ein K\"orper und $R$ sei ein Teilring von $K$ mit $1 \in R$. Wir
nennen $R$ \emph{Bewertungsring} von $K$, falls gilt: Ist $k \in K
- R$, so ist $k^{-1} \in R$. Zeigen wir zun\"achst, dass f\"ur
K\"orper mit Bewertungsringen der Satz 2.5 gilt.

\medskip
\noindent
{\bf 2.6. Satz.} {\it Es sei $K$ ein K\"orper und $R$ sei ein
Bewertungsring von $K$. Ist $E$ eine endliche Teilmenge von $K$,
so gibt es ein $r \in R^*$ mit $er \in R$ f\"ur alle $e \in E$.}

\smallskip

Beweis. Ist $E$ eine Teilmenge  von $R$, so setzen wir $r : = 1$.
Wir d\"urfen daher annehmen, dass es ein $e \in E$ gibt, welches
nicht in $R$ liegt. Dann liegt aber $e^{-1} \in R$. Ist $e$ das
einzige Element von $E$, so setzen wir $r := e^{-1}$. Ist $e$
nicht das einzige Element von $E$, so gibt es nach
Induktionsannahme ein $s \in R^*$ mit $fe^{-1} s \in R$ f\"ur alle
$f \in E - \{e\}$. Setzt man $r := e^{-1} s$, so gilt $fr \in R$
f\"ur alle $f \in E$.

Weitere Eigenschaften von Bewertungsringen sind in den n\"achsten
beiden S\"atzen notiert.

\medskip
\noindent
{\bf 2.7. Satz.} {\it Es sei $K$ ein K\"orper und $R$ sei ein
Bewertungsring von $K$. Ist dann $E$ eine endliche Teilmenge von
$R^* := R - \{0\}$, so gibt es $s, t \in E$ mit $fs^{-1}, t^{-1} f
\in R$ f\"ur alle $f \in E$.}

\smallskip

Beweis. Da $1 \in R$ gilt, ist der Satz richtig, falls $E$ nur ein
Element enth\"alt. Es sei $|E| \geq 1$ und $e \in E$. Es gibt dann
ein $\lambda \in E - \{e\}$ mit $f \lambda^{-1} \in R$ f\"ur alle
$f \in E - \{e\}$. Ist $e \lambda^{-1} \in R$, so tut's $s :=
\lambda$. Ist $e \lambda^{-1} \notin R$, so ist $\lambda e^{-1}
\in R$. In diesem Falle tut's $s := e$.

Die Existenz von $t$ beweist sich analog.

\medskip
\noindent
{\bf 2.8. Satz.} {\it Es sei $K$ ein K\"orper und $R$ sei ein
Bewertungsring von $K$. Ist dann $M:= \{f\mid f \in R, f^{-1} \notin
R\}$, so ist $M$ ein zweiseitiges Ideal von $R$ und alle von $R$
verschiedenen Rechts- wie Linksideale sind in $M$ enthalten.}

\smallskip

Beweis. $M$ ist die Menge der Elemente von $R$, die keine
Einheiten von $R$ sind. Um zu zeigen, dass $M$ ein Ideal ist, ist
nur kritisch nachzuweisen, dass $M$ additiv abgeschlossen ist.
Dazu seien $a, b \in M$. Wegen Satz 2.5 d\"urfen wir annehmen,
dass $b = ac$ ist mit $c \in R$. Dann ist aber, falls $a + b \neq
0$ ist,  $$1 = (a+b)(a+b)^{-1} = a(1+c)(a+b)^{-1}.$$ Also
ist  $$a^{-1} = (1+c)(a+b)^{-1}.$$ Hieraus folgt
$(a+b)^{-1} \notin R$, da andernfalls $a^{-1} \in R$ w\"are. Also
ist $a+b \in M$, was nat\"urlich auch richtig ist, wenn $a +b=0$
ist.

Die restlichen Aussagen folgen daraus, dass $R-M$ nur aus
Einheiten von $R$ besteht.

\smallskip

F\"ur den Kenner sei hier folgendes erw\"ahnt. Ist $R$ ein
Integrit\"atsbereich und ist $F$ ein freier $R$-Modul, so ist
$\Delta_R (F)$ gerade die Menge der reinen Teilmoduln von $F$.
Dies gesagt, ist klar, dass man im Falle von Hauptidealbereichen
einige der folgenden S\"atze allgemeiner formulieren kann, als wir
es tun werden.

\medskip
\noindent
{\bf 2.9. Satz.} {\it Es sei $R$ ein Integrit\"atsbereich oder ein
Bewertungsring. Ferner sei $M$ ein Modul \"uber $R$. Ist $T(M)$
die Menge der Torsionselemente von $M$, so ist $M$ ein Teilmodul
von $M$.}

\smallskip

Beweis. F\"ur Integrit\"atsbereiche ist die Aussage v\"ollig
banal. Es sei $R$ also ein Bewertungsring. Ferner seien $m$ und
$n$ Torsionselemente von $M$. Es gibt dann $r, s \in R^*$ mit $mr
= 0 = ns$. Nach 2.7 d\"urfen wir annehmen, dass $s=rt$ ist mit
einem $t \in R$. Es folgt  $$(m+n)s=mrt+ns=0,$$ so dass $m
+ n \in T(M)$ gilt.

Es sei weiterhin $m \in T(M)$ und es sei $0 \neq s \in R$. Es gibt
ein $r \in R^*$ mit $mr = 0$. Ist $r = st$ mit $t \in R$, so ist
einmal $t \neq 0$ und zum anderen $(ms)t = mr = 0$, so dass $ms
\in T(M)$ gilt. Ist $s$ kein Linksteiler von $r$, so gibt es nach
2.7 ein $t \in R$ mit $s = rt$. In diesem Falle ist aber $ms=0$.
Damit ist alles bewiesen.

\medskip
\noindent
{\bf 2.10. Satz.} {\it Es sei $R$ ein Hauptidealbereich oder ein
Bewertungsring. Es sei ferner $M$ ein torsionsfreier Rechtsmodul
\"uber $R$ und $0 \neq m \in M$. Der Teilmodul $U$ von $M$ werde
definiert durch $U/mR=T(M/mr)$. Ist $U$ endlich erzeugt, so gibt
es ein $u \in U$ mit $U = uR$.}

\smallskip

Beweis. Es sei $U=\sum^n_{i:=1} f_iR$. Es gibt dann $r_i \in  R*$
und $s_i \in R$ mit $f_i r_i = ms_i$ f\"ur $i := 1, \dots, n$. Wir
d\"urfen nat\"urlich annehmen, dass $f_i \neq 0$ ist. Weil $M$
torsionsfrei ist, ist dann $f_ir_i \neq 0$ und folglich $s_i \neq
0$ f\"ur alle $i$.

Wir betrachten zun\"achst den Fall, dass $R$ ein Hauptidealbereich
ist. Dann ist $R$ insbesondere kommutativ. Setze $a :=
\prod^n_{i:=1} r_i$. Ist $x \in U$, so gibt es $\lambda_i \in R$
mit $x=\sum^n_{i:=1} f_i \lambda_i$. Weil $r_i$ auf Grund der
Kommutativit\"at von $R$ ein Teiler von $a$ ist, gibt es $g_i \in
R$ mit $\lambda_i a =r_ig_i$. Es folgt  $$xa=\sum^n_{i:=1}
f_ir_ig_i=m \sum^n_{i:=1} s_ig_i.$$ Somit ist die durch $\sigma(x)
:= xa$ definierte Abbildung $\sigma$ zun\"achst eine Abbildung von
$U$ in $mR$. Da $R$ kommutativ ist, ist $\sigma$ sogar ein
Homomorphismus von $U$ in $mR$. Weil $a$ nicht Null und $M$
torsionsfrei ist, ist $\sigma$ sogar ein Monomorphismus. Nun ist
$mR$ aber ein epimorphes Bild von $R$. Weil $R$ ein
Hauptidealbereich ist, wird folglich jeder Teilmodul von $mR$ von
einem Element erzeugt. Daher wird $\sigma(U)$ und damit $U$ von
einem Element erzeugt.

Es sei nun $R$ ein Bewertungsring. Ist $n = 1$, so ist nichts zu
beweisen. Es sei also $n > 1$. Auf Grund von 2.7 d\"urfen wir
annehmen, dass $s_{n-1} = s_nt$ mit $t \in R$ ist. Dann ist aber
 $$f_{n-1} r_{n-1} = ms_nt= f_nr_nt.$$ nach 2.7 ist
$r_{n-1}$ ein Rechtsteiler von $r_n t$ oder $r_n t$ ein
Rechtsteiler von $r_{n-1}$. Wegen der Torsionsfreiheit von $M$
gilt im ersten Falle $f_{n-1} \in f_n R$ und im zweiten Falle $f_n
\in f_{n-1} R$. Induktion f\"uhrt nun zum Ziele.

\medskip

Der Leser, der mit den Feinheiten von $\ggT$-Bereichen vertraut
ist, wird sehen, dass man bei der folgenden Definition und einigen
der weiteren S\"atze statt Hauptidealbereich auch $\ggT$-Bereich
sagen k\"onnte. Da die wesentlichen S\"atze f\"ur $\ggT$-Bereiche
jedoch ihre G\"ultigkeit verlieren, verzichten wir bei den
fraglichen S\"atzen auf die gr\"o\ss ere Allgemeinheit.

Es sei $R$ ein Hauptidealbereich oder ein Bewertungsring. Ferner
sei $F$ ein freier $R$-Rechtsmodul und $B$ sei eine Basis von $F$.
Ist $f \in F$, so gibt es eine Abbildung $r$ von $B$ in $R$ mit
endlichem Tr\"ager, so dass $f = \sum_{b \in B} br_b$ gilt. Ist
$R$ ein Hauptidealbereich, so setzen wir 
$$\cont(f) := \ggT(r_b\mid b \in B).$$ Ist $R$ ein Bewertungsring, so
gibt es nach 2.7 ein $b \in B$ mit $r r_b^{-1} \in R$ f\"ur alle
$c \in B$, da der Tr\"ager von $r$ ja endlich ist. In diesem Falle
setzen wir  $$\cont (f) := r_b.$$ Wir nennen $\cont (f)$
den \emph{Inhalt} von $f$. Der Inhalt $\cont(f)$ von $f$ ist bis
auf Einheiten eindeutig bestimmt. Ist $\cont(f) = 1$, so nennen
wir $f$ \emph{primitiv}.

Es sei $C$ eine weitere Basis von $F$. Ist $c = \sum_{b \in B}
bA_{bc}$, und $f=\sum_{c \in C} c \lambda_c$, so ist  $$f=
\sum_{b \in B} b \sum_{c \in C} A_{bc} \lambda_c.$$ Hieraus folgt,
dass der mit Hilfe von $C$ definierte Inhalt von $f$ ein Teiler
des mit Hilfe von $B$ definierten Inhalts von $f$ ist. Vertauscht
man in diesem Argument die Rollen von $B$ und $C$, so sieht man,
dass auch der mittels $B$ definierte Inhalt ein Teiler des mittels
$C$ definierten Inhalts ist. Somit h\"angt die Funktion $\cont$
nicht von der Wahl der Basis ab.

\medskip
\noindent
{\bf 2.11. Satz.} {\it Es sei $R$ ein Hauptidealbereich oder ein
Bewertungsring. Ist $F$ ein freier Rechtsmodul \"uber $R$ und ist
$f \in F$ primitiv, so ist $f R \in \Delta_R (F)$. Ist $0 \neq g
\in F$, so gibt es ein primitives $f \in F$ mit $g = f \cont
(g)$.}

\smallskip

Beweis. Es sei $B$ eine Basis von $F$. Dann ist $f = \sum_{b \in
B} b \alpha_b$, wobei $\alpha$ eine Abbildung von $B$ in $R$ ist,
deren Tr\"ager endlich ist. Es sei $0 \neq h \in F$ und $s \in
R^*$ und es gelte $hs \in fR$. Es gibt dann wegen der
Torsionsfreiheit von $F$ (Satz 2.2) ein $t \in R^*$ mit $hs=ft$.
Ferner ist $h = \sum_{b \in B} b \beta_b$ und es folgt 
$$\sum_{b\in B} b \beta_b s = hs = ft = \sum_{b \in B} b \alpha_b
t$$ und damit $\beta_b s = \alpha_b t$ f\"ur alle $b \in B$.

Ist $R$ ein Hauptidealbereich, so ist  $$\cont(h)s= \ggT
(\beta_b s\mid b \in B) = \ggT(\alpha_b t\mid b \in B) = \cont(f)t=t.$$

Ist $R$ ein Bewertungsring, so gibt es ein $b$, so dass $\alpha_b$
eine Einheit ist, da $f$ ja primitiv ist. Es folgt
$t=\alpha_b^{-1} \beta_b s$ und weiter $h=f\alpha_b^{-1} \beta_b
\in fR$. Damit ist gezeigt, dass $fR \in \Delta_R (F)$ gilt.

Es sei $g=\sum_{b \in B} br_b$. Setze $f := \sum_{b \in B} br_b
\cont (g)^{-1}$. Dann ist $f$ primitiv und $g=f \cont (g)$.

\medskip
\noindent
{\bf 2.12. Satz.}{\it Es sei $R$ ein Ring. Ferner sei $M$ ein
Rechtsmodul \"uber $R$ und $U$ sei ein Teilmodul von $M$. Ist
$M/U$ ein freier $R$-Modul, so ist $U$ ein direkter Summand von
$M$.}

\smallskip

Beweis. Mit Hilfe des Auswahlaxioms erhalten wir eine Familie $B$
von Elementen von $M$, so dass  $\{b+U\mid b\in B\}$ eine
Basis von $M/U$ ist. Wir setzen  $$V:=\sum_{b \in B} b R.$$
Nat\"urlich gilt $M = U+V$. Es sei $u \in U \cap V$. Es gibt dann
$b_1, \dots, b_n \in B$ und $r_1, \dots, r_n \in R$ mit $u =
\sum^n_{i:=1} b_ir_i$. Es folgt  $$U=u+U=\sum^n_{i:=1}
(b_i+u)r_i.$$ Also ist $r_i=0$ f\"ur alle $i$, so dass $u=0$ und
folglich $M=U \oplus V$ gilt.

\medskip
\noindent
{\bf 2.13. Satz.} {\it Es sei $R$ ein Hauptidealbereich oder ein
Bewertungsring. Ist $M$ ein endlich erzeugter, torsionsfreier
Rechtsmodul \"uber $R$, so ist $M$ frei in endlich vielen
Erzeugenden.}

\smallskip

Beweis. Es sei $M = \sum^n_{i:=1} f_i R$. Ist $n=1$, so ist $M$
frei in einer Erzeugenden. Es sei also $n > 1$. Der Teilmodul $U$
von $M$ werde definiert durch $U/f_n R= T(M/f_nR)$. Dann ist $M/U$
ein torsionsfreier $R$-Rechtsmodul, der von $f_1+U, \dots, f_{n-1}
+ U$ erzeugt wird. Nach Induktionsannahme ist $M/U$ frei in
h\"ochstens $n-1$ Erzeugenden. Weil $M/U$ frei ist, ist $U$ nach
2.12 ein direkter Summand von $M$. Es sei $C$ ein Komplement von
$U$. Dann ist also $M=U \oplus C$ und somit $M/C \cong U$. Dies
besagt, dass $U$ endlich erzeugt ist. Mit 2.10 folgt, dass es ein
$u \in U$ gibt mit $U = uR$. Da $C$ wegen $C \cong M/U$ frei ist,
ist auch $M$ frei und zwar in h\"ochstens $n$ Erzeugenden.

\medskip
\noindent
{\bf 2.14. Satz.} {\it Es sei $R$ ein Hauptidealbereich oder ein
Bewertungsring. Ferner sei $F$ ein freier Rechtsmodul \"uber $R$.
Ist $D$ ein direkter Summand von $F$, so ist $D \in \Delta_R (F)$.

Genau dann besteht $\Delta_R(F)$ genau aus den direkten Summanden
von $F$, wenn $F$ endlich erzeugt ist.}

\smallskip

Beweis. Es sei $D$ ein direkter Summand von $F$. Weil $F$ frei und
somit torsionsfrei ist, ist jedes Komplement von $D$ torsionsfrei.
Dies impliziert, dass auch $F/D$ torsionsfrei ist. Folglich gilt
$D \in \Delta_R(F)$.

Der Modul $F$ sei endlich erzeugt. Ferner sei $D \in \Delta_R(F)$.
Dann ist $F/D$ endlich erzeugt, da $F$ es ist, und nach Definition
von $\Delta_R(F)$ ist $F/D$ torsionsfrei. Nach 2.13 ist $F/D$
daher frei, was nach 2.12 impliziert, dass $D$ ein direkter
Summand von $F$ ist. Ist $F$ endlich erzeugt, so besteht
$\Delta_R(F)$ also genau aus den direkten Summanden von $F$.

Ist $R$ Bewertungsring, so sei $R$ Bewertungsring des K\"orpers
$K$. Ist $R$ Hauptidealbereich, so sei $K$ der Quotientenk\"orper
von $R$. In beiden F\"allen sei $p$ ein Element von $R$, welches
keine Einheit sei. Wir setzen  $$M := \sum^\infty_{i:=0}
p^{-i} R.$$ Der Modul $F$ sei nicht endlich erzeugt. Ist $B$ eine
Basis von $F$, so ist $B$ nicht endlich, enth\"alt also eine
abz\"ahlbare Teilmenge $C$. Es sei $a$ eine mit $0$ beginnende
Abz\"ahlung von $C$. Es gibt dann einen Epimorphismus $\epsilon$
von $F$ auf $M$ mit $\epsilon (a_i):=p^{-i}$ f\"ur alle nicht
negativen ganzen Zahlen $i$ und $\epsilon(b) =0$ f\"ur alle $b\in
B-C$. Es sei $D$ der Kern von $\epsilon$. Dann ist $F/D$ zu $M$
isomorph, also nicht torsionsfrei. Daher ist $D \in \Delta_R(F)$.

Wir nehmen nun an, $D$ h\"atte ein Komplement $G$, und fahren im
Indikativ fort. Dann ist $G$ zu $M$ isomorph. Es gibt daher eine
Folge $z$ auf $G$ mit $G=\sum^{\infty}_{i:=0} z_iR$ und $z_{i+1}
p=z_i$ f\"ur alle $i$. Nach 2.11 gibt es ein primitives $f \in F$
mit $z_0=f$ cont $(z_0)$. Weil $G$ ein direkter Summand ist, ist
$G \in \Delta_R(F)$, so dass $f \in G$ ist. Weil $f$ primitiv ist,
ist $fR$ nach 2.11 ein Element von $\Delta_R(F)$. Nun ist
$z_ip^i=z_0 \in fR$ und folglich $z_i \in f R$ f\"ur alle $i$.
Also ist $G = fR$. Weil $G$ und $M$ isomorph sind, gibt es ein $v
\in M$ mit $M=vR$. Es gibt weiter $r_0, \dots, r_n \in R$ mit
 $$v=\sum^n_{i:=0} p^{-i} r_i.$$ Es gibt au\ss erdem ein $s
\in R$ mit  $$p^{-n-1} = vs.$$ Hieraus folgt 
$$1=p^{n+1} p^{-n-1} = p \sum^n_{i:=0} p^{n-i} r_is,$$ so dass $p$
eine Einheit ist. Dies widerspricht aber der Wahl von $p$, so dass
$D$ doch kein direkter Summand von $F$ ist. Damit ist alles
bewiesen.

\medskip

Der Leser hat hoffentlich bemerkt, dass der Satz falsch ist.
Unterr\"aume von Vektorr\"aumen sind stets direkte Summanden, so
dass man also noch voraussetzen muss, dass $R$ kein K\"orper ist.
Beim Beweis des Satzes haben wir ja ein von Null verschiedenes
Element von $R$ ben\"otigt, welches in $R$ keine Einheit ist.

Der Kern $D$ von $\epsilon$ ist eine Hyperebene von $\Delta_R(F)$.
Um dies zu beweisen, \"uberlege sich der Leser, dass je zwei
Elemente des Moduls $M$ und damit je zwei Elemente des Moduls
$R/D$ \"uber $R$ linear abh\"angig sind. Hieraus folgt, dass $D$
unter der in 2.5 definierten Abbildung $\varphi$ auf eine
Hyperebene von $F \otimes_R K$ abgebildet wird, wobei $K$ der
Quotientenk\"orper von $R$ ist, falls $R$ ein Hauptidealbereich
ist, bzw. ein K\"orper, von dem $R$ ein Bewertungsring ist. Es
gibt also Hyperebenen und damit Unterr\"aume endlichen Ko-Ranges
in $\Delta_R(F)$, die keine direkten Summanden sind, falls $F$
nicht endlich erzeugt ist. Am unteren Ende des Verbandes
$\Delta_R(F)$ herrschen andere Verh\"altnisse, wie der n\"achste
Satz lehrt.

\medskip
\noindent
{\bf 2.15. Korollar.} {\it Es sei $R$ ein Hauptidealbereich oder
ein Bewertungsring. Ferner sei $F$ ein freier Rechtsmodul \"uber
$R$. Ist $D \in \Delta_R(F)$ endlich erzeugt, so ist $D$ ein
direkter Summand von $F$.}

\smallskip

Beweis. Es sei $B$ eine Basis von $F$. Weil $D$ endlich erzeugt
ist, gibt es eine endliche Teilmenge $E$ von $B$ mit  $$D
\subseteq \sum_{b \in E} bR.$$ Setzt man die Summe rechter Hand
gleich $U$, so ist $U$ ein direkter Summand von $F$, da das
Komplement von $E$ in $B$ ein Komplement von $U$ in $F$ erzeugt.
Ferner ist $D \in \Delta_R (U)$, so dass $D$ nach 2.14 ein
direkter Summand von $U$ ist. Dann ist aber auch ein direkter
Summand von $F$.

\medskip
\noindent
{\bf 2.16. Korollar.} {\it Es sei $R$ ein Hauptidealbereich oder
ein Bewertungsring. Ferner sei $F$ ein freier Rechtsmodul \"uber
$R$. Ist $D$ ein endlich erzeugter direkter Summand von $F$, so
ist auch jedes Komplement von $D$ in $F$ frei.}

\smallskip

Beweis. Wie beim Beweise von 2.15 sehen wir, dass $D$ in einem
endlich erzeugten direkten Summanden $U$ von $F$ liegt, der ein
Komplement besitzt, welches frei ist. Wegen $D \in \Delta_R (F)$
folgt mit 2.14, dass $D$ auch ein direkter Summand von $U$ ist.
Nach 2.13 ist jedes Komplement von $D$ in $U$ frei. St\"uckelt man
ein solches mit dem freien Komplement von $U$ zusammen, so
erh\"alt man ein Komplement von $D$, welches ein freier Modul ist.
Dann ist aber auch $M/D$ frei und folglich jedes Komplement von
$D$.

\medskip
\noindent
{\bf 2.17. Satz.} {\it Es sei $R$ ein Ring mit Eins und $M$ sei
ein zweiseitiges Ideal von $R$ und $R/M$ sei ein K\"orper. Ist $F$
ein Rechtsmodul \"uber $R$, so bezeichnen wir mit $FM$ den von
allen $fm$ mit $f \in F$ und $m \in M$ erzeugten Teilmodul von
$F$. Dann ist $F/FM$ ein Rechtsvektorraum \"uber $R/M$. Ist $F$
frei und ist $B$ eine Basis von $F$, so ist  $$\{b+FM\mid b \in
B\}$$ eine Basis von $F/FM$.}

\smallskip

Beweis. Es ist nat\"urlich klar, dass $\{b+FM\mid b \in B\}$ ein
Erzeugendensystem von $F/FM$ ist. Es ist also nur zu zeigen, dass
dieses Erzeugendensystem auch linear unabh\"angig ist. Dazu sei
$C$ eine endliche Teilmenge von $B$ und es gelte  $$\sum_{c
\in C} (c+FM)(r_c+M)=0,$$ wobei die $r_c$ Elemente von $R$ sind.
Dann ist  $$\sum_{c \in C} cr_c \in FM.$$ Es gibt also
$f_1, \dots, f_n \in F$ und $m_1, \dots, m_n \in M$ mit 
$$\sum_{c \in C} cr_c = \sum^n_{i:=1} f_i m_i.$$ Stellt man nun
die $f_i$ mittels der Basis $B$ dar, so erh\"alt man  $$f_i
= \sum_{b \in B} b \lambda_{ib}$$ mit $\lambda_{ib} \in R$ f\"ur
alle $i$ und alle $b$. Es folgt  $$\sum_{c \in C}
cr_c=\sum^n_{i:=1} \sum_{b\in B} b\lambda_{ib} m_i=\sum_{b \in B}
b \sum^n_{i:=1} \lambda_ib m_i.$$ Koeffizientenvergleich --- hier
benutzen wir, dass $B$ eine Basis ist --- liefert $r_c \in M$
f\"ur alle $c \in C$. Damit ist alles bewiesen.

\medskip

Hat ein Ring $R$ ein zweiseitiges Ideal $M$, so dass $R/M$ ein
K\"orper ist, so ist der Rang eines freien Moduls \"uber $R$ eine
Invariante, wie der gerade bewiesene Satz zeigt.

Bei einem Vektorraum ist der geometrische Rang eines Unterraumes
gleich seinem Rang als Vektorraum. Bei freien Moduln \"uber
Ringen, wie wir sie betrachten, ist der entsprechende Sachverhalt
ganz und gar nicht evident. Daher liest sich der n\"achste Satz
etwas umst\"andlich.

\medskip
\noindent
{\bf 2.18. Satz.} {\it Es sei $R$ ein Hauptidealbereich oder ein
Bewertungsring und $M$ sei ein maximales Ideal von $R$. Dann ist
$R/M$ ein K\"orper. Es sei weiterhin $F$ ein freier Rechtsmodul
\"uber $R$ und $FM$ bezeichne wieder den von allen $fm$ mit $f \in
F$ und $m \in M$ erzeugten Teilmodul von $F$. F\"ur $D \in
\Delta_R (F)$ setzen wir  $$\epsilon (D) := (D+FM)/FM.$$
Dann ist $\epsilon$ ein Epimorphismus des projektiven Verbandes
 $$(\Delta_R(F), \subseteq)$$ auf den projektiven Verband
 $$(\La_{R/M}(F/FM), \leq)$$ mit den folgenden
Eigenschaften:
\item{a)} Ist $k$ eine nat\"urliche Zahl, ist $D \in \Delta_R(F)$
als Modul endlich erzeugt und hat $D$ als Modul den Rang $k$, so
hat auch $\epsilon (D)$ den Rang $k$.
\item{b)} Ist $k$ eine nat\"urliche Zahl und hat $U \in \La_{R/M}
(F/FM)$ den Rang $k$, so gibt es ein $D \in \Delta_R(F)$ mit
$\epsilon(D)=U$. Ist $D \in \Delta_R(F)$ und gilt $\epsilon(D)=U$,
so ist $D$ als Modul endlich erzeugt und $k$ ist auch der Rang von
$D$.
\item{c)} Ist $G$ eine Gerade von $\Delta_R(F)$, so gibt es drei
Punkte auf $G$, deren Bilder unter $\epsilon$ paarweise
verschieden sind.\par}

\smallskip

Beweis. Nat\"urlich ist $\epsilon$ ein Homomorphismus. Dass
$\epsilon$ surjektiv ist, beweisen wir erst zum Schlu\ss.

a) Da $D$ erzeugt ist, ist $D$ frei, so dass klar ist, was
es bedeutet, dass der Rang von $D$ gleich $k$ ist, n\"amlich, dass
$D$ eine Basis $C$ der Kardinalit\"at $k$ hat. (Man kann Rang auch
f\"ur torsionsfreie Moduln ohne M\"uhe definieren, was wir hier
jedoch nicht tun.)

 Nach 2.15 ist $D$ ein direkter Summand
und nach 2.16 ist jedes Komplement von $D$ frei. Es gibt daher
eine Basis $B$ von $F$ mit $C \subseteq B$. Nach 2.17 ist
$\{b+FM\mid b\in B\}$ eine Basis von $F/FM$, so dass $\{\gamma+FM\mid \gamma
\in C\}$ linear unabh\"angig ist. Da $\epsilon(D)$ von dieser
Menge erzeugt wird, gilt die Behauptung~a).

b) Es sei zun\"achst $k=1$. Dann gibt es ein $p \in F$, so
dass $U$ von $p+FM$ erzeugt wird. Nach 2.11 gibt es ein primitives
$f$ mit $p=f \cont(p)$. Es folgt $\cont (p) \notin M$, so dass $U$
auch von $f+FM$ erzeugt wird. Nach 2.11 ist $fR \in \Delta_R(F)$.
Also ist $\epsilon$ auf $f R$ anwendbar und wir erhalten $\epsilon
(fR)=U$.

Nun sei $k$ beliebig. Es gibt dann primitive $f_1, \dots, f_k \in
F$, so dass $\epsilon (f_1R), \dots$, $\epsilon (f_kR)$ eine Basis
von $U$ ist. (Projektive Interpretation!) Es sei $W$ das Supremum
der $f_1R_1, \dots, f_kR$ in $\Delta_R(F)$. Weil $W$ von $k$
Punkten erzeugt wird, ist der Rang von $W$ als projektiver Raum
h\"ochstens gleich $k$. Dann ist aber auch der Rang von
$\epsilon(W)$ h\"ochstens gleich $k$. Nun ist $U \leq \epsilon
(W)$, da ja $\epsilon (f_iR) \leq \epsilon (W)$ f\"ur alle $i$
gilt. Daher ist $U = \epsilon (W)$ und der Rang von $W$ als
projektiver Raum ist gleich~$k$.

Es sei $B$ eine Basis von $F$. Dann h\"angt die Menge der $f_i$
von einer endlichen Teilmenge von $B$ ab. Daher liegen die $f_iR$
in einem endlich erzeugten direkten Summanden von $F$. Weil $W$
das Supremum der $f_iR$ ist, liegt auch $W$ in diesem direkten
Summanden. Dann ist aber $W$ auch als Modul endlich erzeugt und
$k$ ist der Rang von $W$ als Modul.

c) Es sei $G$ eine Gerade von $\Delta_R(F)$. Wie der Beweis
von b) zeigt, ist $G$ ein freier Modul des Ranges $2$. Es gibt
also eine Basis $b_1, b_2$ von $G$. Dann sind aber auch $b_1+b_2,
b_1$ und $b_1+b_2, b_2$ Basen von $G$. Weil $G$ ein direkter
Summand ist, gibt es Basen $B_0, B_1$ und $B_2$ von $F$, die der
Reihe nach die drei zuvor genannten Basen von $G$ enthalten.
Hieraus folgt dann mit 2.17, dass $\epsilon(b_1R),
\epsilon((b_1+b_2)R)$ und $\epsilon(b_2R)$ drei verschiedene
Punkte sind, deren Urbilder auf $G$ liegen.

Es bleibt die Surjektivit\"at von $\epsilon$ nachzuweisen. Dazu
sei $U$ ein Teilraum von $F/FM$. Ferner sei $B$ eine Basis des
projektiven Raumes $\La_{R/M} (U)$. Es gibt dann eine Punktmenge
$C$ von $\Delta_R (F)$, so dass die Einschr\"ankung von $\epsilon$
auf $C$ eine Bijektion von $C$ auf $B$ ist. Es sei $D$ das
Erzeugnis von $C$ in $\Delta_R (F)$.  Ist $P$ ein Punkt auf $D$,
so gibt es eine endliche Teilmenge $E$ von $C$, so dass $P$ von
$E$ abh\"angt. Ist $X$ der von $E$ aufgespannte Teilraum von
$\Delta_R(F)$, so hat $X$ h\"ochstens den Rang $|E|$. Andererseits
hat $Y:= \sum_{e \in E} \epsilon (eR)$ genau den Rang $|E|.$ Wegen
$Y \subseteq \epsilon (X)$ hat $X$ daher ebenfalls den Rang $|E|$
und es gilt $Y = \epsilon (X)$. Hieraus folgt  $$\epsilon
(P) \subseteq Y \subseteq U.$$ Weil $\epsilon (D)$ von seinen
Punkten erzeugt wird, folgt mit a) und b), dass $\epsilon(D)
\subseteq U$ gilt. Andererseits ist banalerweise $U \subseteq
\epsilon(D)$, so dass in der Tat $U=\epsilon(D)$ ist. Damit ist
$\epsilon$ auch als surjektiv erkannt.

\medskip

Es gilt auch die Umkehrung dieses Satzes: Ist $V$ ein Vektorraum
\"uber dem K\"orper $K$ und ist der Rang von $V$ mindestens $3$,
so l\"asst sich jeder Epimorphismus von $\La_K(V)$ auf einen
projektiven Verband $L$, der die Menge der Punkte von $\La_K(V)$
auf die Menge der Punkte von $\La_K(V)$ auf die Menge der Punkte
von $L$ und die Menge der Geraden von $\La_K(V)$ auf die Menge der
Geraden von $L$ abbildet, mittels eines Bewertungsringes von $K$
auf die obige Weise be\-schrei\-ben. Einen Beweis f\"ur diesen
Sachverhalt findet der Leser in Machala (1975). Machalas Beweis
ist mehrere Seiten lang. Ich kann ihn zwar nachvollziehen, da er
mir sein Geheimnis, weshalb er korrekt ist, aber nicht preisgibt,
ist er hier nicht abgedruckt.

\mysection{3. Segresche Mannigfaltigkeiten}

\noindent
Nach diesem Intermezzo \"uber Epimorphismen von projektiven
Verb\"anden  wenden wir uns wieder unserem eigentlichen Thema zu.

Ist $R$ ein kommutativer Ring mit Eins und ist $M$ ein $R$-Modul,
so d\"urfen und werden wir annehmen, dass $M$ sowohl ein Rechts-
als auch ein Linksmodul \"uber $R$ ist und dass dar\"uber hinaus
$rm = mr$ f\"ur alle $r \in R$ und alle $m \in M$ gilt. Dann ist
$M$ sogar ein $R$-Bimodul.

Es sei $R$ ein kommutativer Ring mit Eins und $M_1, \dots, M_t$
seien $R$-Moduln. Eine Abbildung $f$ von $M_1 \times \cdots \times
M_t$ in den $R$-Modul $N$ hei\ss t \emph{t-fach linear} oder kurz
\emph{multilinear}, falls $f$ in jedem Argument linear ist.

War das Tensorprodukt bislang eine bin\"are Operation, so f\"uhren
wir nun das Tensorprodukt als $t$-\"are Operation ein, um dann zu
sehen, dass sich die neue Operation auf die alte zur\"uckf\"uhren
l\"asst. Dazu sei $R$ ein kommutativer Ring mit Eins und $M_1,
\dots, M_t$ und $T$ seien Moduln \"uber $R$. Es sei ferner $\tau$
$t$-fach lineare Abbildung von $M_1 \times \cdots M_t$ in $T$. Wir
nennen $(T,\tau)$ ein \emph{Tensorprodukt} von $M_1, \dots, M_t$,
falls gilt:

\item{1)} Der Modul $T$ wird von der Menge der $\tau (m_1, \dots, m_t)$
mit $m_i \in M_i$ f\"ur $i:=1, \dots, t$ erzeugt.

\item{2)} Ist $N$ ein $R$-Modul und ist $f$ eine multilineare
Abbildung von $M_1 \times \cdots M_t$ in $N$, so gibt es eine
lineare Abbildung $g$ von $T$ in $N$ mit $f = g\tau$.
\par\noindent
Die Frage nach Existenz und Eindeutigkeit von Tensorprodukten
beantwortet der folgende Satz.

\medskip
\noindent
{\bf 3.1. Satz.} {\it Es sei $R$ ein kommutativer Ring mit Eins.
Sind $M_1, \dots, M_t$ Moduln \"uber $R$, so gibt es bis auf
Isomorphie genau ein Tensorprodukt $(T, \tau)$ von $M_1, \dots,
M_t$.}

\smallskip

Beweis. Dass es bis auf Isomorphie h\"ochstens ein Tensorprodukt
gibt, beweist man wie bei solch universellen Objekten \"ublich.

Nat\"urlich ist $(M_1, \mathrm{id})$ ein Tensorprodukt von $M_1$, so dass
der Satz f\"ur $t=1$ richtig ist. Es sei $t > 1$ und $1 \leq i <
t$. Nach Induktionsannahme gibt es ein Tensorprodukt $(A, \alpha)$
von $M_1, \dots, M_i$ und ein Tensorprodukt $(B, \beta)$ von
$M_{i+1}, \dots, M_t$. Wir zeigen, dass $(A \otimes_R B, \tau)$,
wobei $\tau$ die durch  $$\tau(x_1, \dots, x_t):=
\alpha(x_1, \dots, x_i) \otimes \beta (x_{i+1}, \dots, x_t)$$
definierte Abbildung ist, ein Tensorprodukt von $M_1, \dots, M_t$
ist.

Um dies zu zeigen, sei $f$ eine multilineare Abbildung von $M_1
\times \cdots \times M_t$ in $W$. Ferner sei $x \in M_1 \times
\cdots \times M_i$ und $y \in M_{i+1} \times \cdots \times M_t$.
Wir definieren $g_x$ durch  $$g_x (y) := f(x,y).$$ Dann ist
$g_x$ multilinear, so dass es eine lineare Abbildung $\psi_x$ von
$B$ in $W$ gibt mit $g_x=\psi_x \beta$.

Es sei $b \in B$. Wir definieren $h_b$ durch  $$h_b(x) :=
\psi_x(b)$$ f\"ur alle $x \in M_1 \times \cdots \times M_i$. Weil
$\psi_x$ ein Homomorphismus ist und weil $B$ von der Menge der
$\beta(y)$ erzeugt wird, folgt, dass auch $h_b$ multilinear ist.
Es gibt daher einen Homomorphismus $\varphi_b$ von $A$ in $W$ mit
$h_b= \varphi_b \alpha$. Wir definieren nun die Abbildung $F$ von
$A \times B$ in $W$ durch  $$F(a,b) := \varphi_b(a).$$ Dann
ist $F$ gewiss linear in $a$. Ist $x \in M_1 \times \cdots
\times M_i$, so folgt  $$F(\alpha(x), b) = \varphi_b
\alpha(x) = h_b (x) =\psi_x(b),$$ so dass $F$ auch in $b$ linear
ist, da $F$ in $a$ linear ist und die $\alpha (x)$ den Modul $A$
erzeugen. Schlie\ss lich ist   $$F(a,
rb)=F(a,br)=F(a,b)r=F(ar,b),$$ so dass $F$ tensoriell ist. Es gibt
daher eine lineare Abbildung $G$ von $A \otimes_R B$ in $W$ mit
$F(a,b)=G(a \otimes b)$ f\"ur alle $a \in A$ und alle $b \in B$.
Mit diesem $G$ folgt schlie\ss lich
$$\eqalign{
(G\tau)(x,y) &= G(\alpha(x) \otimes \beta (y))
= F(\alpha(x),\beta(y))\cr &= \psi_x(\beta(y)) = g_x(y) = f(x,y).}
$$
Um zu erkennen, dass $(A \otimes_R B,\tau)$ ein Tensorprodukt von
$M_1, \dots, M_t$ ist, ist nur noch zu bemerken, dass die Menge
der $\tau(x)$ mit $x \in M_1 \times \cdots \times M_t$ ein
Erzeugendensystem von $A \bigotimes_R B$ ist. Damit ist Satz 3.1
bewiesen.

\medskip

Wir haben 
 viel mehr bewiesen als im Satz formuliert. Um dies
deutlich zu machen, vereinbaren wir zun\"achst, statt $\tau(x_1,
\dots, x_t)$ wie 
 \"ublich  $x_1 \otimes \cdots \otimes
x_t$ zu schreiben und das Tensorprodukt selbst mit
$\bigotimes^t_{i:=1} M_i$ zu bezeichnen.

\medskip
\noindent
{\bf 3.2. Korollar.} {\it Es sei $R$ ein kommutativer Ring mit
Eins und $M_1, \dots, M_t$ seien Moduln \"uber $R$. Ferner sei $1
\leq i \leq t$. Es gibt dann genau einen Isomorphismus $\sigma$
des Tensorproduktes $\bigotimes^t_{k:=1} M_k$ auf das
Tensorprodukt $(\bigotimes^i_{k:=1} M_k) \otimes
(\bigotimes^t_{k:=i+1} M_k)$ mit  $$\sigma (x_1 \otimes
\cdots \otimes x_t)=(x_1 \otimes \cdots \otimes x_i) \otimes
(x_{i+1} \otimes \cdots \otimes x_t).$$ Dabei ist das
Tensorprodukt \"uber eine leere Indexmenge (hier der Fall $i=t$)
als der Ring $R$ zu interpretieren.}

\smallskip

Beweis. Nur der Fall $i=t$ bedarf eines Hinweises. In diesem Falle
wende man Satz 1.7 an.

\medskip

Dieses Korollar macht explizit, dass sich das Tensorprodukt von
mehreren Moduln \"uber einem kommutativen Ring mit Hilfe des
Tensorproduktes von zwei Moduln ausdr\"ucken l\"asst.

Wir betrachten nun Vektorr\"aume $V_1, \dots, V_t$ \"uber einem
kommutativen K\"or\-per $K$. Die Vektoren der Form $v_1 \otimes
\cdots \otimes v_t$ des Tensorproduktes $\bigotimes^t_{i:=1} V_i$
nennen wir \emph{reine} oder auch \emph{zerlegbare Tensoren}. Die
Menge der von reinen Tensoren aufgespannten Punkte von
$\La_K(\bigotimes^t_{i:=1} V_i)$ bezeichnen wir mit $S(V_1, \dots,
V_t)$. Wir nennen diese Menge \emph{Segresche Mannigfaltigkeit mit
den Parameterr\"aumen} $V_1, \dots, V_t$. Es gilt nun

\medskip
\noindent
{\bf 3.3. Satz.} {\it Es sei $K$ ein kommutativer K\"orper und
$V_1, \dots, V_t$ seien Vektorr\"aume \"uber $K$. Ist $\pi \in
S_t$, so gibt es eine projektive Kollineation $\kappa$ von $\La_K
(\bigotimes^t_{i:=1} V_i)$ auf $\La_K(\bigotimes^t_{i:=1}
V_{\pi(i)})$ mit $S(V_1, \dots, V_t)^\kappa = S(V_{\pi(1)}, \dots,
V_{\pi(t)})$.}

\smallskip

Beweis. Die Abbildung, die dem Element $(v_1, \dots, v_t)$ das
Element $v_{\pi(1)} \otimes \cdots \otimes v_{\pi(t)}$ zuordnet,
ist multilinear. Es gibt daher eine lineare Abbildung $\lambda$
von $\bigotimes^t_{i:=1} V_i$ auf $\bigotimes^t_{i:=1} V_{\pi(i)}$
mit  $$(v_1 \otimes \cdots \otimes v_t)^\lambda =
v_{\pi(1)} \otimes \cdots \otimes v_{\pi(t)}.$$ Vertauschen der
Rollen der beiden Tensorprodukte liefert die Existenz der zu
$\lambda$ inversen Abbildung. Folglich ist $\lambda$ ein
Isomorphismus und der von $\lambda$ induzierte Isomorphismus
$\kappa$ von $\La_k(\bigotimes^t_{i:=1} V_i)$ auf
$\La_K(\bigotimes^t_{i:=1} V_{\pi(i)})$ bildet $S(V_1, \dots,
V_t)$ auf $S(V_{\pi(1)}, \dots, V_{\pi(t)})$ ab.

\medskip

Es sei $P$ ein Punkt der Segremannigfaltigkeit $S(V_1, \dots,
V_t)$ mit den Parameterr\"aumen $V_1, \dots, V_t$. Es gibt dann
Vektoren $p_i \in V_i$ mit $P=(p_1 \otimes \cdots \otimes p_t)K$.
F\"ur $i := 1, \dots, t$ setzen wir  $$U_i(P) :=\{p_1
\otimes \cdots \otimes p_{i-1} \otimes x \otimes p_{i+1} \otimes
\cdots \otimes p_t \mid x \in V_i\}.$$ Dann ist $U_i(P)$ ein zu $V_i$
isomorpher Unterraum von $\La_K (\bigotimes ^t_{i:=1} V_i)$, der
$P$ enth\"alt und dessen Punkte alle zu $S(V_1, \dots, V_t)$
geh\"oren. Um dies anzudeuten, schreiben wir auch $U_i (P) \leq
S(V_1, \dots, V_t)$, obgleich das formal nicht korrekt ist. Wir
setzen ferner  $$T(P) := \sum^t_{i:=1} U_i (P)$$ und nennen
$T(P)$ den \emph{Tangentialraum} von $S(V_1, \dots, V_t)$ in $P$.

\medskip
\noindent
{\bf 3.4. Satz.} {\it Es sei $K$ ein kommutativer K\"orper und
$V_1, \dots, V_t$ seien Vektorr\"aume \"uber $K$. Ist $X$ eine
nicht leere Teilmenge von $\{1, \dots, t\}$ und ist $j \in \{1,
\dots, t\} - X$, so ist} $$U_j(P) \cap \sum_{i \in X}
U_i(P)=P.$$


Beweis. Auf Grund von 3.3 d\"urfen wir annehmen, dass $j = 1$ und
$X=\{2, \dots, s\}$ ist. Es sei $P=(p_1 \otimes \cdots \otimes
p_t)K$. Ferner sei $u \in U_1(P) \cap \sum^s_{i:=2} U_i (P)$. Es
gibt dann $x_i \in V_i$ mit  $$u=-x_1 \otimes p_2 \otimes
\cdots \otimes p_t$$ einerseits und  $$u=p_1 \otimes x_2
\otimes \cdots \otimes p_t + \cdots + p_1 \otimes \cdots \otimes
x_s \otimes \cdots \otimes p_t$$ andererseits. Hieraus folgt
weiter unter Zuhilfenahme von 3.2, dass  $$0=x_1 \otimes
(p_2 \otimes \cdots \otimes p_t)+p_1 \otimes y$$ ist, wobei $y$
eine naheliegende Abk\"urzung ist. Weil $p_2 \otimes \cdots
\otimes p_t$ nicht Null ist, sind $x_1$ und $p_1$ nach 1.11 linear
abh\"angig. Daher ist $x_1 \in p_1 K$, so dass der Satz bewiesen
ist.

\medskip
\noindent
{\bf 3.5. Korollar.} {\it Die Voraussetzungen seien die gleichen
wie in Satz 3.4. Sind \"uberdies die R\"ange der $V_i$ endlich und
setzt man $r_i := \Rg_K (V_i)$, so ist} $$\Rg_K (T(P)) =
1-t + \sum^t_{i:=1} r_i.$$ 

Dies folgt unter Benutzung der
Rangformel mittels Induktion aus 3.4.

\medskip

Der n\"achste Satz wird seine Bedeutung verlieren, sobald Satz 3.9
bewiesen ist. F\"ur den Beweis dieses Satzes wird er jedoch
ben\"otigt.

\medskip
\noindent
{\bf 3.6. Satz.} {\it Es sei $K$ ein kommutativer K\"orper, 
$V_1, \dots, V_t$ seien Vektorr\"aume \"uber $K$, und $P$ sei
ein Punkt von $S(V_1, \dots, V_t)$. Ist $U \in
\La_K(\otimes^t_{i:=1} V_i)$, ist $U \leq T(P)$ und liegt jeder
Punkt von $U$ auf $S(V_1, \dots, V_t)$, so gibt  es ein $i$ mit $U
\leq U_i(P)$.}

\smallskip

Beweis. Es sei $P=(p_1 \otimes \cdots \otimes p_t)K$ und $Q=(q_1
\otimes \cdots \otimes q_t)K$ sei ein Punkt von $U$. Es gibt dann
$u_i \in V_i$ mit  $$-q_1 \otimes \cdots \otimes q_t =
\sum^t_{i:=1} p_1 \otimes \cdots \otimes u_i \otimes \cdots
\otimes p_t.$$ Wir d\"urfen $P \neq Q$ annehmen. Es gibt dann ein
$i$ mit $q_i \notin p_i K$. Wegen 3.3 d\"urfen wir annehmen, dass
$i=1$ ist. Dann folgt  $$0=q_1 \otimes (q_2 \otimes \cdots
\otimes q_t) + u_1 \otimes (p_2 \otimes \cdots \otimes p_t)+p_1
\otimes y,$$ wobei $y$ wieder f\"ur einen l\"angeren Ausdruck
steht, dessen Einzelheiten uns nicht interessieren. Weil $q_2
\otimes \cdots \otimes q_t$ nicht Null ist, sind die Vektoren
$q_1, p_1$ und $u_1$ nach 1.11 linear abh\"angig. W\"are $u_1 =
p_1 k$ mit einem $k \in K$, so folgte  $$0=q_1 \otimes (q_2
\otimes (q_2 \otimes \cdots \otimes q_t)+ p_1 \otimes (y + (p_2
\otimes \cdots \otimes p_t)k)$$ im Widerspruch zur linearen
Unabh\"angigkeit von $q_1$ und $p_1$. Somit sind $u_1$ und $p_1$
linear unabh\"angig. Es gibt daher $a,b \in K$ mit $q_1=u_1a+p_1
b$. Es folgt $a \neq 0$ und  $$0=u_1 \otimes ((q_2 \otimes
\cdots \otimes q_t)a+p_2 \otimes \cdots \otimes p_t)+p_1 \otimes
z$$ mit einem geeigneten $z$. Weil $u_1$ und $p_1$ linear
unabh\"angig sind, folgt  $$0=(q_2 \otimes \cdots \otimes
q_t)a + p_2 \otimes \cdots \otimes p_t.$$ Hieraus folgt nun
mittels Induktion, dass $q_i \in p_i K$ ist f\"ur $i:= 2, \dots,
t$. Dies zeigt, dass $Q \leq U_1 (P)$ ist.

Es seien $Q$ und $R$ zwei Punkte von $U$, die beide von $P$
verschieden seien. Es gibt dann $i$ und $j$ mit $Q \leq U_i (P)$
und $R \leq U_j(P)$, wie wir gerade gesehen hatten. Wir zeigen,
dass $i=j$ ist. Dazu nehmen wir an, dass dies nicht der Fall sei.
Wir d\"urfen dann annehmen, dass $i=1$ und $j=2$ ist. Es folgt
$Q=(u_1 \otimes p_2 \otimes \cdots \otimes p_t)K$ und $R= (p_1
\otimes u_2 \otimes \cdots \otimes p_t)K$ mit $u_i \in V_i$ f\"ur
$i:= 1,2$. Es folgt, dass auch der durch  $$T:= (u_1
\otimes p_2 \otimes \cdots \otimes p_t + p_1 \otimes u_2 \otimes
\cdots \otimes p_t)K$$ definierte Punkt $T$ ein Punkt auf $S(V_1,
\dots, V_t)$ ist. Es gibt daher ein $k$ mit $T \leq U_k(P)$. Wegen
$Q+R=R+T=T+Q$ und $U_1 (P) \cap U_2 (P) = P$ (Satz 3.4) ist $k
\neq 1,2$. Daher d\"urfen wir annehmen, dass $k = 3$ ist. Es folgt
$T=(p_1 \otimes p_2 \otimes u_3 \otimes \cdots p_t)K$ mit $u_3 \in
V_3$, so dass es ein $k \in K$ gibt mit  $$(u_1 \otimes p_2
\otimes p_3 +p_1 \otimes u_2\otimes p_3 + p_1 \otimes p_2 \otimes
u_3k) \otimes \cdots \otimes p_t = 0.$$ Hieraus folgt 
$$u_1 \otimes p_2 \otimes p_3 + p_1 \otimes u_2 \otimes p_3 + p_1
\otimes p_2 \otimes u_3 k = 0$$ und weiter  $$u_1 \otimes
(p_2 \otimes p_3) + p_1 \otimes (u_2 \otimes p_3 + p_2 \otimes u_3
k) = 0,$$ so dass $u_1$ und $p_1$ linear abh\"angig sind. Dies hat
aber $Q = P$ zur Folge. Damit ist gezeigt, dass es ein $i$ gibt,
so dass alle von $P$ verschiedenen Punkte von $U$ in $U_i(P)$
liegen. Also ist $U \leq U_i (P)$.

\medskip
\noindent
{\bf 3.7. Satz.} {\it Es sei $K$ ein kommutativer K\"orper und
$V_1, \dots, V_t$ seien Vektorr\"aume \"uber $K$. Sind $P$ und $Q$
Punkte auf $S(V_1, \dots, V_t)$ und ist $U_i(P) \cap U_i(Q) \neq
\{0\}$, so ist $U_i(P) = U_i(Q)$.}

\smallskip

Beweis. Es sei $P=(p_1 \otimes \cdots \otimes p_t)K$ und $Q = (q_1
\otimes \cdots \otimes q_t)K$. Auf Grund unserer Voraussetzung
gibt es $x_i, y_i \in V_i$ mit  $$0 \neq p_1 \otimes \cdots
\otimes x_i \otimes  \cdots \otimes p_t = q_1 \otimes \cdots
\otimes y_i \otimes \cdots \otimes q_t.$$ Hieraus folgt $p_j
K=q_iK$ f\"ur alle von $i$ verschiedenen $j$. Dies zeigt, dass
$U_i(P) = U_i(Q)$ ist.

\medskip

Besonders einfach sind Segresche Mannigfaltigkeiten mit nur zwei
Parameterr\"aumen. F\"ur diese gilt der folgende Satz, der uns
sogleich von Nutzen sein wird.

\medskip
\noindent
{\bf 3.8. Satz.} {\it Es sei $K$ ein kommutativer K\"orper und
$V_1$ und $V_2$ seien Vektorr\"aume \"uber $K$. Sind $P$ und $Q$
zwei Punkte von $S(V_1, V_2)$, so ist $U_1(P) \cap U_2(Q)$ ein
Punkt.}

\smallskip

Beweis. Es sei $U_1(P)=\{x_1 \otimes p_2\mid x_1 \in V_1\}$ und
$U_2(Q)=\{q_1 \otimes x_2\mid x_2 \in V_2\}$. Setze $R:= (p_1 \otimes
q_2)K$. Dann ist $R$ ein Punkt mit $R \leq U_1(P) \cap U_2(Q)$.
Mit 3.7 und 3.4 folgt $$U_1(P) \cap U_2(Q)=U_1(R)\cap
U_2(R)=R.$$ Damit ist der Satz bewiesen.

\medskip
\noindent
{\bf 3.9. Satz.} {\it Es sei $K$ ein kommutativer K\"orper und
$V_1, \dots, V_t$ seien Vektorr\"aume \"uber $K$. Ist $U$ ein
Teilraum von $\bigotimes^t_{i:=1} V_i$ und liegen alle Punkte von
$U$ auf $S(V_1, \dots, V_t)$, so gibt es einen Punkt $P$ auf
$S(V_1, \dots, V_t)$ und ein $i$ mit $U \leq U_i(P)$.}

\smallskip

Beweis. Wir machen Induktion nach $t$. Sei also zun\"achst $t=2$.
Ferner seien $Q$ und $R$ zwei verschiedene Punkte auf $U$. Setze
$P := U_1(Q)\cap U_2(R)$. Nach 3.8 ist $P$ ein Punkt. Mit 3.7
folgt weiter $U_1   (Q)=U_1(P)$ und $U_2(R) = U_2(P)$. Daher ist
 $$Q+R\leq T(P).$$ Weil alle Punkte von $Q+R$ zu $S(V_1,
V_2)$ geh\"oren, folgt nach 3.6 dass es ein $i \in \{1,2\}$ gibt
mit  $$Q+R \leq U_i(P).$$ Hieraus folgt $U_i(Q)=U_i (R)$.
Es sei $S$ ein von $Q$ und $R$ verschiedener Punkt von $U$. Dann
folgt mit dem gleichen Argument, dass es $k, l \in \{1,2\}$ gibt
mit $U_k (Q)=U_k(S)$ und $U_k(S)$ und $U_l(R)=U_l(S)$ ist. Ist
$k=1$, so folgt  $$Q+R \leq U_k(S) \cap U_i(P)$$ und damit
$k=i$, da ja $Q+R$ eine Gerade ist. In diesem Falle ist $S \leq
U_i(P)$. Ist $k \neq l$, so ist, da wir nur zwei Indizes zur
Auswahl haben, $k=i$ oder $l=i$, so dass auch hier $S \leq U_i(P)$
gilt. Weil $U$ die obere Grenze der in $U$ enthaltenen Punkte ist,
ist daher $U \leq U_i (P)$.

Es sei nun $t > 2$. In diesem Falle identifizieren wir
$\otimes^t_{i:=1} V_i$ mit $(\otimes^{t-1}_{i:=1} V_i) \otimes
V_t$. Dann liegen die Punkte von $U$ alle auf 
$$S\big( \bigotimes^{t-1}_{i:=1} V_i, V_t \big).$$ Nach dem bereits
Bewiesenen liegen die Punkte von $U$ entweder in einem Teilraum
der Form  $$\{v \bigotimes x_t\mid x_t \in V_t\}$$ mit einem $v
\in \bigotimes^{t-1}_{i:=1} V_i$ oder in einem Teilraum der Form
$$\{y \otimes v_t\mid y \in \bigotimes^{t-1}_{i:=1} V_i\}$$ mit einem
$v_t \in V_t$. Im ersten Falle folgt, weil die Punkte von $U$ ja
in der Segreschen Mannigfaltigkeit liegen, dass $v$ zerlegbar ist.
Das bedeutet aber, dass es einen Punkt $P$ auf der Segreschen
Mannigfaltigkeit gibt mit $U \leq U_t(P)$. Im zweiten Falle folgt
aus der Induktionsannahme, dass es einen Punkt $P$ sowie ein $i$
gibt mit $1 \leq i \leq t-1$ und $U \leq U_i(P)$. Damit ist auch
dieser Satz bewiesen.

\medskip

Die R\"aume $U_i(P)$ spielen eine hervorragende Rolle, wie wir
jetzt schon sehen. Daher setzen wir  $$E_i:= \{U_i(P) \mid P \in
S(V_1, \dots, V_t)\}$$ und nennen $E_i$ die \emph{i}te \emph{Schar
von Erzeugenden} von $S(V_1, \dots, V_t)$. Ist $t=2$ und
$\Rg_K(V_i)=2$ f\"ur beide $i$, so hei\ss en $E_1$ und $E_2$ auch
\emph{Regelscharen} bzw. \emph{reguli} (Einzahl: \emph{regulus}.
Blaschke bildet den Plural \glqq Regulusse\grqq.) V\"ollig unklar
ist mir die Herkunft von \emph{regulus}. Im Lateinischen bedeutet
es \glqq K\"onig eines kleinen Landes\grqq, weiter einen kleinen
Vogel, hinter dem man den Zaunk\"onig vermutet und, als
\"Ubersetzung von $\beta \alpha \sigma \iota \lambda \imath \sigma
\kappa o \varsigma$, \glqq Basilisk\grqq, eine Eidechsenart. Ich
wei\ss{} auch nicht, welches der beiden W\"orter \glqq
Regelschar\grqq{} und \emph{regulus} im Sinne von Regelschar
\"alter ist.

\medskip
\noindent
{\bf 3.10. Satz.} {\it Es seien $V_1$ und $V_2$ Vektorr\"aume
\"uber dem kommutativen K\"orper $K$. Sind $x, y, z \in V_1$ und
$u, v, w \in V_2$ und gilt $x \otimes u+y \otimes v = z \otimes
w$, so sind $x$ und $y$ oder $u$ und $v$ linear abh\"angig.}

\smallskip

Beweis. Die Vektoren $x$ und $y$ seien linear unabh\"angig. Ist
$u=0$ oder $v=0$, so ist nichts zu beweisen. Es sei also $u \neq
0$ und $v \neq 0$. Es sei ferner $B$ eine Basis von $V_1$ mit $x,y
\in B$ und $C$ sei eine Basis von $V_2$. Schlie\ss lich sei
$z=\sum_{b \in B} b \zeta_b, u=\sum_{c \in C} c \alpha_c,
v=\sum_{c \in C} c \beta_c$ und $w=\sum_{c \in C} c \gamma_c$.
Dann ist  $$\sum_{c \in C} (x \otimes c) \alpha_c + \sum_{c
\in C} (y \otimes c) \beta_c = \sum_{b \in B} \sum_{c \in C} (b
\otimes c) \zeta_b \gamma_c.$$ Da die $b \otimes c$ eine Basis von
$V_1 \otimes_K V_2$ bilden, folgt $\alpha_c = \zeta_x \gamma_c$
und $\beta_c = \zeta_y \gamma_c$ f\"ur alle $c \in C$. Hieraus
folgt $u=w \zeta_x$ und $v = w \zeta_y$. Weil $u$ und $v$ nicht
Null sind, sind $\zeta_x$ und $\zeta_y$ nicht Null, so dass $u$
und $v$ in der Tat linear abh\"angig sind.

\medskip

Ich wei\ss\ nicht, ob man beim n\"achsten Satz auf die
Voraussetzung der Endlichkeit der R\"ange verzichten kann.

\medskip
\noindent
{\bf 3.11. Satz.} {\it Es sei $K$ ein kommutativer K\"orper und
$V_1, \dots, V_t$ seien Vektorr\"aume \"uber $K$ mit $2 \leq
\Rg_K(V_i) < \infty$. Ist $\sigma \in GL (\bigoplus^t_{i:=1} V_i)$
und gilt  $$S(V_1, \dots, V_t)^\sigma = S(V_1, \dots,
V_t),$$ so gibt es ein $\pi \in S_t$ mit $E^\sigma_i = E_{\pi(i)}$
f\"ur $i:=1, \dots, t$.}

\smallskip

Beweis. Wegen 3.3 d\"urfen wir annehmen, dass  $$2 \leq
\Rg_K(V_1) \leq \dots \leq \Rg_K(V_t)$$ ist.

Es sei $U$ ein Element von $E_t$. Da die Segresche
Mannigfaltigkeit $S(V_1, \dots, V_t)$ invariant unter $\sigma$ ist, liegen
alle Punkte von $U^\sigma$ auf $S(V_1, \dots, V_t)$. Nach 3.9 gibt
es daher ein $j$ und ein $W \in E_j$ mit $U^\sigma \leq W$. Auf
Grund unserer Annahme \"uber die R\"ange der $V_i$ folgt $U^\sigma
= W$. Wir wollen zeigen, dass hieraus folgt, dass $E^\sigma_t =
E_j$ ist. Dazu d\"urfen wir wegen 3.3 annehmen, dass $j = t$ ist.

Es ist  $$U=\{p_1 \otimes \cdots \otimes p_{t-1} \otimes
x\mid x \in V_t\}$$ und  $$U^\sigma = \{q_1 \otimes \cdots
\otimes q_{t-1} \otimes y\mid y \in V_t\}.$$ Es sei $W := \{w_1
\otimes p_2 \cdots \otimes p_{t-1} \otimes x\mid x \in V_t\}$, aber
$W^\sigma \notin E_t$. Wir k\"onnen annehmen, dass $W^\sigma \in
E_{t-1}$ ist. Es gibt dann Vektoren $r_i$ mit  $$W^\sigma =
\{r_1 \otimes \cdots \otimes r_{t-2} \otimes x \otimes r_t\mid x \in
V_{t-1}\}.$$ Weil der Rang von $V_t$ mindestens $2$ ist, gibt es
ein $y \in V_t$, so dass $y$ und $r_t$ linear unabh\"angig sind.
Es gibt ferner ein $x \in V_t$ mit  $$(p_1 \otimes \cdots
\otimes p_{t-1} \otimes x)^\sigma = q_1 \otimes \cdots \otimes
q_{t-1} \otimes y$$ und ein $z \in V_{t-1}$ mit  $$(w_1
\otimes p_2 \otimes \cdots \otimes p_{t-1} \otimes x)^\sigma = r_1
\otimes \cdots \otimes r_{t-2} \otimes z \otimes r_t.$$ 
Nun ist
 $$q_1 \otimes \cdots \otimes q_{t-1} \otimes y+r_1 \otimes
\cdots \otimes r_{t-2} \otimes z \otimes r_t = ((p_1+w_1) \otimes
p_2 \otimes \cdots p_{t-1} \otimes x)^\sigma$$ Weil $\sigma$
zerlegbare Tensoren auf zerlegbare Tensoren abbildet und weil $y$
und $r_n$ linear unabh\"angig sind, folgt mit 3.10, dass es ein $k
\in K$ gibt mit  $$q_1 \otimes \cdots \otimes q_{t-1} =
(r_1 \otimes \dots \otimes r_{t-2} \otimes z)k.$$ Hieraus folgt
$q_iK=r_iK$ f\"ur $i:=1, \dots, t-2$ und $q_{t-1} K=zK$. Dies
impliziert die Existenz eines $l \in K$ mit  $$q_1 \otimes
\cdots q_{t-1} \otimes y+r_1 \otimes \cdots \otimes r_{t-2}
\otimes z \otimes r_t = q_1 \otimes \cdots \otimes q_{t-1} \otimes
(q+r_tl).$$ Also ist  $$((p_1+w_1) \otimes p_2 \otimes
\cdots \otimes p_{t-1} \otimes x)^\sigma = q_1 \otimes \cdots
\otimes q_{t-1} \otimes (y+r_tl).$$ Der letzte Vektor liegt aber
im Bild von $U$ unter $\sigma$. Es gibt also ein $u \in V_t$ mit
 $$p_1 \otimes \cdots \otimes p_{t-1} \otimes u = (p_1+w_1)
\otimes p_2 \otimes \cdots \otimes p_{t-1} \otimes x.$$ Hieraus
folgt schlie\ss lich, dass $p_1$ und $w_1$ linear abh\"angig sind,
so dass $U = W$ ist. Dieser Widerspruch zeigt, dass doch $W^\sigma
\in E_t$ ist.

Ist $i \geq 1$ und ist bereits gezeigt, dass  $$\{w_1
\otimes \cdots \otimes w_i \otimes p_{i+1} \otimes \cdots \otimes
p_{t-1} \otimes x\mid x \in V_t\}^\sigma$$ in $E_t$ liegt, so folgt
mittels 3.3 und dem bereits Bewiesenen, dass auch  $$\{w_1
\otimes \cdots \otimes w_{i+1} \otimes p_{i+2} \otimes \cdots
\otimes p_{t-1} \otimes x\mid x \in V_t\}^\sigma$$ in $E_t$ liegt.
Damit ist gezeigt, dass $E^\sigma_t=E_t$ ist.

Die Behauptung des Satzes folgt nun mittels Induktion nach $t$.

\medskip

Wir sind nun in der Lage, den Stabilisator $\hat{G}(V_1, \dots,
V_t)$ von $S(V_1, \dots, V_t)$ in $PGL(\bigotimes^t_{i:=1} V_i)$
zu bestimmen.

\medskip
\noindent
{\bf 3.12. Satz.} {\it Es sei $K$ ein kommutativer K\"orper und
$V_1, \dots, V_t$ seien Vektorr\"aume endlichen Ranges \"uber $K$
mit $2 \leq \Rg_K(V_i)$ f\"ur alle $i$. Ist $\sigma$ eine
Kollineation von $\La_K(\bigotimes^t_{i:=1} V_i)$, die $S(V_1,
\dots, V_t)$ invariant l\"asst, und gibt es einen Punkt $P \in
S(V_1, \dots, V_t)$, so dass $\sigma$ den Tangentialraum $T(P)$
punktweise festl\"asst, so ist $\sigma = 1$.}

\smallskip

Beweis. Weil der Rang von $T(P)$ mindestens $2$ ist, wird $\sigma$
durch eine lineare Abbildung induziert. Dies werden wir sp\"ater
verwenden.

Eine zweite Bemerkung, die uns gleich n\"utzlich sein wird, ist
die, dass $S(V_1, \dots, V_t)$ einen Rahmen von
$\bigotimes^t_{i:=1} V_i$ enth\"alt, woraus nach III.1.10 folgt,
dass einzig die Identit\"at unter den projektiven Kollineationen
$S(V_1, \dots, V_t)$ punktweise festl\"asst. Um die Existenz des
Rahmens zu zeigen, sei $B_i$ eine Basis von $V_i$. Setzt man $u_i
:= \sum_{b \in B_i} b$, und definiert man $b_\alpha$ durch
$b_\alpha := \alpha (1) \otimes \cdots \otimes \alpha (t)$ f\"ur
alle $\alpha  \in \cart^t_{i:=1} B_i$, so bilden die Punkte $P_0
:= (u_1 \otimes \cdots \otimes u_t)K, P_\alpha := b_\alpha K$
einen Rahmen von $\bigotimes^t_{i:=1} V_i$.

Es sei $t=2$. Ferner sei $Q$ ein Punkt von $S(V_1, V_2)$. Dann ist
$R := U_1(Q) \cap U_2(P)$ nach 3.8 ein Punkt. Es folgt $U_1 (Q) =
U_1 (R)$ und $U_2(P)=U_2(R)$. Es folgt, dass $R$ und $U_2(R)$
unter $\sigma$ festbleiben. Dann bleibt aber auch $U_1(R)$ unter
$\sigma$ fest. Wegen $U_1(R)=U_1(Q)$ ist also $U_1(Q)$ unter
$\sigma$ fest. Genauso zeigt man, dass auch $U_2(Q)$ unter
$\sigma$ festbleibt. Wegen $Q=U_1(Q) \cap U_2(Q)$ ist folglich $Q$
bei $\sigma$ fest. Dies zeigt, dass $\sigma$ alle Punkte von
$S(V_1, V_2)$ festl\"asst. Weil $\sigma$ von einer linearen
Abbildung induziert wird, ist $\sigma$ nach obiger Bemerkung die
Identit\"at.

Es sei nun $t >2$ und $P = (p_1 \otimes \cdots \otimes p_t) K$.
Ist $0 \neq x_t \in V_t$, so ist $X := (p_1 \otimes \cdots \otimes
p_{t-1} \otimes x_t)K$ ein Punkt von $U_t(P)$. Hieraus folgt
$X^\sigma = X$. Dies impliziert, dass $U_{t-1} (X)^\sigma
=U_{t-1}(X)$ ist. Da dies f\"ur jedes $x_t$ gilt, folgt, dass der
Raum $W$, der von  $$\{p_1 \otimes \cdots \otimes p_{t-2}
\otimes x_{t-1} \otimes x_t\mid x_{t-1} \in V_{t-1}, x_t \in
V_t\}$$erzeugt wird, unter $\sigma$ invariant ist. Die Punkte von
$S(V_1, \dots, V_t)$, die in $W$ liegen, bilden eine zu
$S(V_{t-1}, V_t)$ isomorphe Segremannigfaltigkeit. Da $P$ zu
dieser Mannigfaltigkeit geh\"ort und da der Tangentialraum
$U_{t-1} (P) + U_t(P)$ an diese Mannigfaltigkeit im Punkte $P$ als
Unterraum von $T(P)$ ebenfalls punktweise festbleibt, folgt, dass
$W$ von $\sigma$ punktweise festgelassen wird. Nun ist
$\bigotimes^t_{i:=1} V_i$ zu  $$\bigg(
\bigotimes^{t-2}_{i:=1} V_i \bigg) \otimes (V_{t-1} \otimes V_t)$$
in kanonischer Weise isomorph, so dass wir diese beiden
Vektorr\"aume identifizieren d\"urfen. Dann l\"asst aber
$\sigma$ die R\"aume $\sum^{t-2}_{i:=1} U_i (P)$ und $W$
punktweise fest. Die lineare Abbildung, von der $\sigma$ induziert
wird, wirkt auf diesen beiden R\"aumen als Skalarmultiplikation.
Da $P$ im Durchschnitt dieser beiden R\"aume liegt, wird sie auf
beiden R\"aumen durch die gleiche Skalarmultiplikation
dargestellt. Somit l\"asst $\sigma$ den Tangentialraum $T' (P)$
von $S(V_1, \dots, V_{t-2}, V_{t-1} \otimes V_t)$ punktweise fest.
Nach Induktionsannahme ist daher $\sigma=1$.

\medskip
\noindent
{\bf 3.13. Satz.} {\it Es sei $K$ ein kommutativer K\"orper und
$V_1, \dots, V_t$ seien Vektorr\"aume endlichen Ranges \"uber $K$.
Es sei  $$\{r_1, \dots, r_s\} := \{\Rg_K(V_i)\mid i := 1,
\dots, t\}$$ und  $$2 \leq r_1 < r_2 < \dots < r_s.$$
Ferner sei $\alpha_j$ die Anzahl der $i$ mit $\Rg_K(V_i)=r_j$. Ist
$G(V_1, \dots, V_t)$ diejenige Untergruppe von $\hat{G}(V_1,
\dots, V_t)$, die die Erzeugendenscharen $E_i$ von $S(V_1, \dots,
V_t)$ jede f\"ur sich invariant l\"asst, so ist 
$$\hat{G}(V_1, \dots, V_t)/G(V_1, \dots, V_t) \cong S_{\alpha_1}
\times \cdots \times S_{\alpha}$$ und} $$G(V_1, \dots, V_t)
\cong PGL (V_1) \times \cdots \times PGL (V_t).$$


Beweis. Die erste Aussage folgt unmittelbar aus 3.11.

Es sei $\rho_i \in GL (V_i)$ f\"ur $i:=1, \dots, t$. Mit Hilfe von
3.2 und 1.3 erhalten wir die Existenz genau einer Abbildung
$\rho_1 \otimes \cdots \otimes \rho_t$ mit  $$(v_1 \otimes
\cdots \otimes v_t)^{\rho_1 \otimes \cdots \otimes \rho_t} =
v_1^{\rho_1} \otimes \cdots \otimes v^{\rho_t}_t.$$ Es sei $\eta$
die Abbildung, die dem Element $(\rho_1, \dots, \rho_t)$ die von
$\rho_1 \otimes \cdots \rho_t$ in $\La_K (\bigotimes^t_{i:=1}
V_i)$ induzierte Kollineation zuordnet. Dann ist $\eta$ ein
Homomorphismus von $GL(V_1) \times \cdots \times GL (V_t)$ in
$G(V_1, \dots, V_t)$. Es sei $(\rho_1, \dots, \rho_t)$ im Kern von
$\eta$. Es gibt dann ein $k \in K^*$ mit  $$v^{\rho_1}_1
\otimes \cdots \otimes v^{\rho_t}_t = (v_1 \otimes \cdots \otimes
v_t)k.$$ Hieraus folgt die Existenz von $k_i \in K^*$ mit
$v^{\rho_i}_i = v_ik_i$ f\"ur $i:=1, \dots, t$. Dies zeigt, dass
$\eta$ einen Monomorphismus von $PGL(V_1) \times \cdots \times PGL
(V_t)$ in $G(V_1, \dots, V_t)$ induziert.

Es bleibt zu zeigen, dass $\eta$ surjektiv ist. Dazu sei $\sigma$
eine lineare Abbildung, die eine Kollineation aus $G(V_1, \dots,
V_t)$ induziert. Es sei $P=(p_1 \otimes \cdots \otimes p_t)K$ ein
Punkt von $S(V_1, \dots, V_t)$. Dann ist auch $P^\sigma$ ein Punkt
von $S(V_1, \dots, V_t)$. Es gibt folglich $Q_i \in V_i$ mit
 $$(p_1 \otimes \dots, \otimes p_t)^\sigma = q_1 \otimes
\dots \otimes q_t.$$ Es gilt dann $U_i(P)^\sigma = U_i(P^\sigma)$
f\"ur alle $i$. Es gibt daher eine Abbildung $\sigma_i \in
GL(V_i)$ mit  $$(p_1 \otimes \cdots \otimes x_i \otimes
\cdots \otimes p_t)^\sigma = q_1 \otimes \cdots \otimes
x^{\sigma_i}_i \otimes \cdots \otimes q_t$$ und es folgt
$p^{\sigma_i}_i = q_i$. Nach dem, was wir bereits gezeigt haben,
induziert $\sigma(\sigma_1 \otimes \cdots \otimes \sigma _t)^{-1}$
eine Kollineation, die zu $G(V_1, \dots, V_t)$ geh\"ort. Diese
Kollineation l\"asst aber $T(P)$ punktweise fest, da die lineare
Abbildung $\sigma(\sigma_1 \otimes \cdots \otimes \sigma_t)^{-1}$
die Unterr\"aume $U_i(P)$ vektorweise festl\"asst. Mit 3.12
folgt daher, dass die fragliche Kollineation die Identit\"at ist.
Also ist $\sigma = \sigma_1 \otimes \cdots \otimes \sigma_t$, so
dass $\eta$ auch surjektiv ist. Damit ist alles bewiesen.

\mysection{4. Geometrische Erzeugung Segrescher Mannigfaltigkeiten}

\noindent
Die Definition der Segreschen Mannigfaltigkeiten, die wir im
letzten Abschnitt gegeben haben, l\"asst nat\"urlich viel zu
w\"unschen \"ubrig. Sie gibt eine rein algebraische Beschreibung
dieser Objekte, so dass sich die Frage erhebt, ob es geometrische
Kennzeichnungen dieser Mannigfaltigkeiten gibt. Die Antwort lautet
nat\"urlich \glqq ja\grqq, zumindest, wenn wir die Endlichkeit der
R\"ange der Parameterr\"aume voraussetzen.

\medskip
\noindent
{\bf 4.1. Satz.} {\it Es sei $L$ ein projektiver Verband und $U_1,
\dots, U_n$ seien unabh\"angige Elemente dieses Verbandes, dh., es
gelte  $$U_k \cap \sum_{i \in X} U_i = 0$$ f\"ur alle
echten Teilmengen $X$ von $\{1, \dots, n\}$ und alle Indizes $k$,
die nicht in $X$ liegend. Ist $P$ ein Punkt von $\sum^n_{i:=1}
U_i$ und gilt  $$P \nleq \sum_{i \in X} U_i$$ f\"ur alle
echten Teilmengen $X$ von $\{1, \dots, n\}$, so gibt es genau
einen Teilraum $T$ des Ranges $n$ von $L$, der $P$ enth\"alt und
der jeden der R\"aume $U_i$ nicht trivial schneidet.

Es ist $T \cap U_i$ f\"ur alle $i$ ein Punkt und  $\{T
\cap U_i\mid i := 1, \dots, n\}$ ist eine Basis von $T$.}

\smallskip

Beweis. Wir machen Induktion nach $n$. Es sei zun\"achst $n=2$. In
diesem Falle setzen wir  $$T:=(U_1+P)\cap(U_2+P).$$ Dann
ist $P \leq T$. Auf Grund des Modulargesetzes gilt  $$U_1 +
((U_1+P)\cap U_2)=(U_1+P)\cap(U_1+U_2)=U_1+P,$$ da ja $P \leq
U_1+U_2$ ist. Hieraus folgt, dass $(U_1+P)\cap U_2$ ein Punkt ist,
da $P$ ja nicht in $U_1$ liegt. Wegen  $$T \cap
U_2=(U_1+P)\cap (U_2+))\cap U_2=(U_1+P)\cap U_2$$ ist $T \cap U_2$
also ein Punkt. Ebenso sieht man, dass $T \cap U_1$ ein Punkt ist.
Nochmalige Benutzung des Modulargesetzes liefert  $$ T =
(U_1+P) \cap(P+U_2)= P+ ((U_1+P) \cap U_2),$$ so dass $T$ die
Summe zweier Punkte ist. Folglich ist $T$ eine Gerade. Damit ist
die Existenz von $T$ gezeigt. Ist andererseits $T'$ eine Gerade
durch $P$, die $U_1$ und $U_2$ trifft, so ist  $$T' \leq
(U_1+P)\cap(U_2+P)=T,$$ womit auch die Einzigkeit von $T$
nachgewiesen ist.

Es sei nun $n \neq 3$. Wir setzen $V:= \sum^{n-1}_{i:=1} U_i$.
Dann gilt $P \leq V + U_n$ und $P \nleq V, U_n$. Nach dem bereits
Bewiesenen gibt es genau eine Gerade $G$ durch $P$, die $V$ und
$U_n$ jeweils in einem Punkte trifft. Setze $Q:=V \cap G$. Ist
dann $X$ eine echte Teilmenge von $\{1, \dots, n-1\}$, so folgt $Q
\nleq \sum_{i \in U_i}$, da andernfalls $P \leq \sum_{i \in X \cup
\{n\}} U_i$ w\"are. Es gibt also einen Teilraum $H$ des Ranges
$n-1$, der durch $Q$ geht und jeden der Teilr\"aume $U_i$ mit $i <
n$ nicht trivial trifft. Setzt man $T := G+H$, so liegt $P$ auf
$T$ und $T$ schneidet jeden der R\"aume $U_i$ nicht trivial. Nun
ist offenbar $G \cap H = Q$ und daher  $$\Rg_L(T) =
\Rg_L(G) + \Rg_L(H) - 1 = 2+n-1-1=n,$$ so dass die Existenz von
$T$ bewiesen ist.

Die Bedeutungen von $G, H, Q$ und $V$ werden beibehalten. Es sei
$T'$ ein Raum des Ranges $n$, der $P$ enth\"alt und alle $U_i$
nicht trivial schneidet. Weil die $U_i$ unabh\"angig sind, folgt
 $$n \geq \Rg_L\bigg(\sum^n_{i:=1} (U_i \cap T') \bigg) =
\sum^n_{i:=1} \Rg_L(U_i \cap T') \geq n$$ und weiter $\Rg_L(U_i
\cap T') = 1$. Setze $R_i := U_i \cap T'$. Dann ist $R_n \leq P +
U_n$. Wegen $P \nleq \sum^{n-1}_{i:=1} R_i$ ist
$T'=P+\sum^{n-1}_{i:=1} R_i$. Es folgt $R_n \leq P+V$. Dies
besagt, dass $R_n \leq (P+U_n) \cap (P+V) =G$ ist. Also ist $G
\leq T'$. Es folgt, dass auch $Q \leq T'$ ist. Also ist $T' \cap
V$ ein Raum des Ranges $n-1$ durch $Q$, der die $U_i$ mit $i \leq
n-1$ nicht trivial schneidet. Nach Induktionsannahme ist daher $T'
\cap V = H$ und weiter $T = T'$. Damit ist die Einzigkeit von $T$
bewiesen. Gleichzeitig sehen wir, dass $T \cap U_i$ ein Punkt ist
und dass die Menge dieser Punkte eine Basis von $T$ ist. Damit ist
alles bewiesen.

\medskip

Da $T$ die R\"aume $U_i$ \emph{transversal} schneidet, nennen wir
$T$ \emph{Transversale} der $U_i$ durch den Punkt $P$.

\medskip
\noindent
{\bf 4.2. Satz.} {\it Es sei $L$ ein projektiver Verband endlichen
Ranges mit dem gr\"o\ss ten Element $\Pi$. Ferner seien $U_1,
\dots, U_n, D \in L$ und $\Pi$ sei die direkte Summe von je $n$
von ihnen. Ist $P$ ein Punkt auf $D$, so bezeichne $T_P$ die
Transversale der $U_1, \dots, U_n$ durch den Punkt $P$. Es gilt
nun:

\item{a)} Ist $B$ eine Basis von $D$, so ist   $$B_i :=
\{T_P \cap U_i\mid P \in B\}$$ eine Basis von $U_i$.

\item{b)} Ist $C$ ein Rahmen von $D$, so ist  $$C_i :=
\{T_P \cap U_i\mid P \in C\}$$ eine Rahmen von $U_i$.

\item{c)} Ist $B$ eine Basis von $D$ und ist $P \in B$, so ist
 $$\{P\} \cup \{T_P \cap U_i\mid i := 1, \dots, n\}$$ ein
Rahmen von $T_P$.}

\smallskip

Beweis. a) Wir setzen $V_i=\sum_{P\in B_i} P$. Dann ist $V_i \leq
U_i$. Es folgt  $$T_P \leq U_1+ \dots + U_{i-1}+V_i+U_{i+1}
+ \dots + U_n.$$ Wegen $P \leq T_P$ ist also  $$D \leq U_1
+ \dots + U_{i-1} + V_i + U_{i+1} + \dots + U_n.$$ Dann ist aber
$$\eqalign{
\Pi &= U_1 + \dots + U_{i-1}+D+U_{i+1} + \dots + U_n\cr
&\leq U_1 +\dots + U_{i-1} + V_i U_{i+1} + \dots + U_n,}
$$
so dass wegen der Unabh\"angigkeit der $U_j$ folgt, dass $U_i =
V_i$ ist. Somit ist $B_i$ ein Erzeugendensystem von $U_i$. Weil
$\Pi$ von je $n$ der R\"aume $U_1, \dots, U_n, D$ erzeugt wird,
sind diese R\"aume alle isomorph. Daher ist $B_i$ ein minimales
Erzeugendensystem, dh., eine Basis von $U_i$.

b) folgt unmittelbar aus a).

c) folgt mit 4.1.

\medskip
\noindent
{\bf 4.3. Korollar.} {\it Es sei $K$ ein K\"orper und $V$ sei ein
Vektorraum endlichen Ranges \"uber $K$. Ferner seien $U_1, \dots,
U_n, D \in \La_K(V)$ und $V$ sei die direkte Summe von je $n$
dieser R\"aume. Ebenso sei $V$ die direkte Summe von je $n$ der
$U'_1, \dots, U'_n, \dots, U'_n, D' \in \La_K(V)$. Ist dann $P_0,
\dots, P_r$ ein Rahmen von $D$ und $P'_0, \dots, P'_r$ ein Rahmen
von $D'$, so gibt es ein $\sigma \in PGL (V)$ mit 
$$U^\sigma_i = U'_i$$ f\"ur $i:= 1, \dots, n$ und 
$$P^\sigma_j = P'_j$$ f\"ur $j := 0, \dots, r$.}

\smallskip

Dieses Korollar zu beweisen, sei dem Leser als \"Ubungsaufgabe
\"uberlassen.

\medskip

Unsere bisherige Definition der Segreschen Mannigfaltigkeiten
machte explizit Gebrauch von den Parameterr\"aumen $V_1, \dots,
V_t$. Davon wollen wir nun loskommen. Dazu definieren wir jetzt
allgemeiner: Es sei $W$ ein Vektorraum endlichen Ranges \"uber dem
kommutativen K\"orper $K$ und $S$ sei eine Menge von Punkten von
$\La_K(W)$. Wir nennen $S$ eine \emph{Segremannigfaltigkeit},
falls es $K$-Vektorr\"aume $V_1, \dots, V_t$ gibt und einen
Isomorphismus $\sigma$ von $\La_K(\bigotimes^t_{i:=1} V_i)$ auf
$\La_K(W)$ mit $S(V_1, \dots, V_t)^\sigma = S$. Geometrisch
relevant ist nur der Fall, dass $t\geq 2$ ist und die R\"ange der
$V_i$ allesamt ebenfalls mindestens gleich $2$ sind. Dann ist der
Rang von $W$ mindestens gleich $4$, so dass es nach dem zweiten
Struktursatz eine semilineare Abbildung von $\bigotimes^t_{i:=1}
V_t$ gibt, durch die $\sigma$ induziert wird. Man sieht leicht,
dass es dann auch eine lineare Abbildung von $\bigotimes^t_{i:=1}
V_i$ auf $W$ gibt, die einen Isomorphismus der Verb\"ande
induziert, welcher $S(V_1, \dots, V_t)$ auf $S$ abbildet. Da die
R\"ange der $V_i$ den Isomorphietyp von $S$ v\"ollig festlegen,
sagen wir auch, falls es die Deutlichkeit erfordert, dass $S$ eine
Segresche Mannigfaltigkeit mit den \emph{Invarianten} $n_1, \dots,
n_t$ sei. Diese Definition ist sorgf\"altig zu lesen. Sie
beinhaltet n\"amlich auch, dass $\Rg_K(W) = \Pi^t_{i:=1} n_i$ ist,
wenn es in $\La_K(W)$ eine Segresche Mannigfaltigkeit mit den
Invarianten $n_1, \dots, n_t$ gibt.

Der n\"achste Satz zeigt uns, wie man Segresche Mannigfaltigkeiten
geo\-me\-trisch erzeugen kann.

\medskip
\noindent
{\bf 4.4. Satz.} {\it Es seien $n_1, \dots, n_t$ nat\"urliche
Zahlen mit $n_i \geq 2$ f\"ur alle $i$. Es sei $K$ ein
kommutativer K\"orper und $W$ sei ein $K$-Vektorraum. Schlie\ss
lich seien $W_1, \dots, W_{n_t}, D$ Unterr\"aume von $W$, so dass
$W$ die direkte Summe von je $n_t$ von ihnen ist. Ist dann $S$
eine Segresche Mannigfaltigkeit mit den Invarianten $n_1, \dots,
n_{t-1}$ von $\La(D)$ und ist $S$ die Menge der Punkte $P$ von
$\La_K(W)$, zu denen es einen Punkt $Q \in S$ gibt, so dass $P$
auf der Transversalen $T_Q$ von $W_1, \dots, W_{n_t}$ durch $Q$
liegt, so ist $\hat{S}$ eine Segresche Mannigfaltigkeit mit den
Invarianten $n_1, \dots, n_t$.}

\smallskip

Beweis. Es seien $V_1, \dots, V_t$ Vektorr\"aume \"uber $K$ mit
$\Rg_K(V_i)=n_i$. Diese werden wir gleich benutzen.

Aus der Annahme, dass $W$ die direkte Summe von je $n_t$ der
R\"aume $W_1, \dots$, $W_{n_t}, D$ ist, folgt, dass die $W_i$
allesamt zu $D$ isomorph sind. Hieraus folgt  $$\textstyle \Rg_K
(W)=n_t\Rg_K (D) = \prod^t_{i:=1} n_i,$$ so dass $W$ als
Vektorraum zu $\bigotimes^t_{i:=1} V_i$ isomorph ist. Wir d\"urfen
daher $W = \bigotimes^t_{i:=1} V_i$ annehmen.

Es sei $B_i$ eine Basis von $V_i$ und $A:=\cart^t_{i:=1} B_i$.
Ferner sei $w_i := \sum_{b \in B_i} b$ und $b_\alpha :=\alpha (1)
\otimes \cdots \otimes \alpha (t)$ f\"ur $\alpha \in A$. Schlie\ss
lich sei $P_0 := (w_1 \otimes \cdots \otimes w_t)K$ und $P_\alpha
:= b_\alpha K$. Dann ist  $$\{P_0\} \cup \{P_\alpha\mid \alpha
\in A\}$$ ein Rahmen von $W$, der in $S(V_1, \dots, V_t)$
enthalten ist, wie wir schon einmal feststellten. F\"ur $b \in
B_t$ setzen wir  $$M_b := \sum_{\alpha \in A, \alpha (t) =
b} b_\alpha K$$ und  $$N:= \sum_{i:=1, \dots, t-1, x_i \in
V_i} (x_1 \otimes \cdots \otimes x_{t-1} \otimes w_t)K.$$ Alle
diese R\"aume haben den Rang $\prod^{t-1}_{i:=1} n_i$, sind also
isomorph zueinander. Ferner folgt, weil $\{b_\alpha\mid \alpha \in
A\}$ eine Basis von $W$ ist, dass $W$ die direkte Summe der $M_b$
ist. Es sei $c \in B_t$ und $y \in N \cap \sum_{b \in B_{t-\{c\}}}
M_b$. Wegen $y \in N$ ist $y=v \otimes w_t$ mit $v \in
\bigotimes^{t-1}_{i:=1} V_i$. Andererseits ist $y = \sum_{\alpha
\in A, \alpha (t) \neq c} b_\alpha k_\alpha$. Hieraus folgt mit
Umsortieren der Summanden $y = \sum_{b \in B_{t}-\{c\}} v_b
\otimes b$ mit $v_b \in \bigotimes^{t-1}_{i:=1} V_i$. Also ist
 $$y=v \otimes w_t = \sum_{b \in B_t-\{c\}} v_b \otimes
b.$$ Weil $w_t$ die Summe aller $b \in B_t$ ist, folgt, dass auch
$(B_t -\{c\}) \cup \{w_t\}$ eine Basis von $V_t$ ist. Mit 1.11
folgt daher, dass unter anderem $v=0$ ist. Also ist $y=0$. Hieraus
folgt schlie\ss lich, dass $W$ auch die direkte Summe von $D$ und
$\sum_{b \in B_{t-\{c\}}} M_b$ ist. Wegen 4.3 d\"urfen wir also
annehmen, dass $D = N$ und dass  $$\{W_i\mid i := 1, \dots,
n_t\} = \{M_b\mid b \in B_t\}$$ ist. Die Untergruppe der projektiven
Gruppe, die die R\"aume $N$ und $M_\alpha$ jeden f\"ur sich
invariant l\"asst, induziert nach 4.3 und III.1.11 die $PGL (D)$
in $\La_K(D)$. Daher d\"urfen wir auch noch annehmen, dass $S$ aus
den Punkten der Form $(x_1 \otimes \cdots \otimes x_{t-1} \otimes
w_t)K$ besteht. Wir m\"ussen schlie\ss lich noch zeigen, dass
$\hat{S} = S(V_1, \dots, V_t)$ ist. Dazu sei $P := (x_1 \otimes
\cdots \otimes x_t)K$ ein Punkt von $S$ und es gelte  $$P
\leq U_t (Q).$$ Ist $b \in B_t$, so ist $(x_1 \otimes \cdots
\otimes x_{t-1} \otimes b)K$ ein Punkt auf $U_t (Q)$, aber auch
auf $M_b$. Dies zeigt, dass $U_t (Q)$ eine Transversale der $M_b$
ist. Also ist $P \in \hat{S}$ und folglich $S(V_1, \dots, V_t)
\subseteq \hat{S}$. Ist andererseits $R \in \hat{S}$, so gibt es
einen Punkt $Q \in S$, so dass $R$ auf der Transversalen $T_Q$ der
R\"aume $M_b$ durch den Punkt $Q$ liegt. Wie wir bereits
feststellten, ist auch $U_t(Q)$ eine Transversale der $M_b$ durch
$Q$. Nach 4.3 ist daher $T_Q = U_t (Q)$, so dass $P \in S(V_1,
\dots, V_t)$ ist. Damit ist alles bewiesen.

\medskip

Da $S(V_1)$ nichts anderes als die Menge der Punkte von
$\La_K(V_1)$ ist, ist nach diesem Satz klar, wie man sich rekursiv
eine Segremannigfaltigkeit mit vorgegebenen Invarianten rein
geometrisch konstruieren kann. Besonders zu\-frie\-den\-stel\-lend ist das
Verfahren, wenn man eine Segremannigfaltigkeit mit nur zwei
Invarianten konstruieren will. Sind beide Invarianten $2$, so
erh\"alt man eine Segremannigfaltigkeit, die man in der sch\"onen
alten Zeit bereits in der Vorlesung \glqq Analytische
Geometrie\grqq{} unter dem Namen \emph{Hyperboloid} kennenlernte.

\medskip
\noindent
{\bf 4.5. Satz.} {\it Es sei $K$ ein kommutativer K\"orper und $V$
sei ein Vektorraum des Ranges $4$ \"uber $K$. Ist $S$ eine
Segremannigfaltigkeit mit den Invarianten $2,2$ in $\La_K(V)$, so
ist $S$ eine Quadrik, die durch eine quadratische Form von
maximalen Index dargestellt wird.}

\smallskip

Beweis. Wir k\"onnen annehmen, dass $V=V_1 \otimes_K V_2$ ist mit
zwei Vektorr\"aumen $V_i$ des Ranges $2$ \"uber $K$ und dass die
fragliche Segremannigfaltigkeit gleich $S(V_1, V_2)$ ist. Es sei
$b_1, b_2$ eine Basis von $V_1$ und $b_3, b_4$ eine Basis von
$V_2$ Ist dann $u=b_1 k_1 + b_2 k_2$ und $v=b_3 k_3 + b_4 k_4$, so
ist $$u \otimes v=(b_1 \otimes b_3)k_1 k_3 + (b_1 \otimes
b_4)k_1 k_3 + (b_2 \otimes b_3)k_2 k_3+(b_2 \otimes b_4)k_2 k_4.$$
Bezeichnet man die Koordinaten eines Vektors bez\"uglich der Basis
aus den $b_i \otimes b_j$ mit $x_{ij}$, so folgt, dass die
Koordinaten des Vektors $u \otimes v$ der quadratischen Gleichung
 $$x_{13} x_{24} - x_{14} x_{23} = 0$$ gen\"ugen, so dass
$S$ also in der durch diese Gleichung definierten Quadrik
enthalten ist. Es sei umgekehrt $(x_{13}, x_{14}, x_{23},
x_{24})K$ ein Punkt dieser Quadrik. Da $b_1$ und $b_2$ v\"ollig
gleichberechtigt sind wie auch $b_3$ und $b_4$ und da auch $V_1$
und $V_2$ gleiche Rollen spielen, d\"urfen wir annehmen, dass
$x_{13} \neq 0$ ist. Wir setzen  $$k_1 := 1, k_3 := x_{13},
k_2 := \frac{x_{23}}{k_3}, k_4 := x_{14}.$$ Dann ist 
$$0=x_{13} x_{24} - x_{14} x_{23} = k_3 x_{24} - k_4 k_2 k_3.$$
Hieraus folgt $x_{24} = k_2 k_4$, da $k_3$ ja von Null verschieden
ist. Dies zeigt, dass alle Punkte der Quadrik zu $S$ geh\"oren. Da
durch jeden Punkt von $S$ zwei Geraden gehen, deren Punkte alle zu
$S$ geh\"oren, und da diese Geraden nach V.5.4 vollst\"andig
isotrop sind, ist der Index der die Quadrik darstellenden
quadratischen Form mindestens 2. Weil der Rang von $V$ gleich $4$
ist, ist er sogar gleich $2$. Damit ist alles gezeigt.

\medskip
\noindent
{\bf 4.6. Satz.} {\it Es sei $V$ ein Vektorraum des Ranges $4$
\"uber einem kommutativen K\"orper. Sind $G_1, G_2$ und $G_3$
paarweise windschiefe Geraden von $V$ und ist $Q$ die Menge der
Punkte, die auf Transversalen von $G_1, G_2$ und $G_3$ liegen, so
ist $Q$ eine Quadrik von maximalem Index.}

\smallskip

Beweis. Mittels 4.4 folgt, dass $Q$ eine Segresche
Mannigfaltigkeit auf den Invarianten $2, 2$ ist. Nach 4.5 ist $Q$
dann eine Quadrik, die durch eine quadratische Form von maximalen
Index dargestellt wird.

\medskip

Wir sind nun in der Lage, Satz I.10.2 zu beweisen, den wir hier
noch einmal notieren werden, damit der Leser ihn vor Augen hat.

\medskip
\noindent
{\bf I.10.2. Satz.} {\it Es sei $L$ ein irreduzibler projektiver
Verband, dessen Rang mindestens $4$ sei. Genau dann gilt in $L$
der Satz von Pappos, wenn f\"ur jedes Hexagramme mystique $G_1,
G_2, G_3, H_1, H_2, H_3$ gilt: ist $H$ eine Transversale der $G_i$
und $G$ eine Transversale der $H_j$, so ist $G \cap H \neq 0$.}

\smallskip

Beweis. Schlie\ss t sich jedes Hexagramme mystique, so gilt in $L$
der Satz von Pappos. Dies ist gerade die Aussage des Satzes von
Dandelin (Satz I.10.1), den wir nicht noch einmal beweisen werden.

Um die Umkehrung zu beweisen, d\"urfen wir auf Grund des ersten
Struktursatzes annehmen, dass es einen Vektorraum $V$ gibt, so
dass $L=\La(V)$ ist. Weil $L$ pappossch ist, ist der
Koordinatenk\"orper von $L$ nach Satz II.6.3 kommutativ. Es sei
nun $G_1, G_2, G_3, H_1, H_2, H_3$ ein Hexagramme mystique. Dann
ist $U:= G_1+G_2$ ein Unterraum des Ranges $4$ und das Hexagramme
mystique ist in $H$ enthalten. Wir d\"urfen daher annehmen, dass
$U=V$ ist. Dann ist $\Rg_K(V) = 4$. Nach Satz 4.4 ist die Menge
$S$ der Punkte, die auf Transversalen der $G_i$ liegen, eine
Segremannigfaltigkeit mit den Invarianten $2,2$. Dann ist $S$ aber
auch die Menge der Punkte, die auf Transversalen der $H_i$ liegen,
da durch einen Punkt von $H_1$ nach Satz 4.1 genau eine
Transversale der $H_i$ geht und die Punkte von $S$ alle auf
Transversalen der $H_i$ liegen. Ist nun $H$ eine Transversale der
$G_i$ und $G$ eine solche der $H_i$, so folgt mit 3.8, dass $G
\cap H$ ein Punkt ist. Damit ist der Beweis von I.10.2
nachgetragen.

\medskip

Die Segreschen Mannigfaltigkeiten mit den Invarianten $2, 2$
gestatten uns auch eine Frage zu beantworten, die A. F. M\"obius
im Jahre 1828 im crelleschen Journal stellte und mit Hilfe seines
baryzentrischen Kalk\"uls auch sogleich beantwortete, die Frage
n\"amlich: \glqq Kann von zwei dreiseitigen Pyramiden eine jede in
Bezug auf die andere um- und einbeschrieben zugleich
heissen?\grqq{} Solche \emph{M\"o\-bi\-us\-paare}, wie sie heute hei\ss{}en, 
gibt es in der Tat in jeder dreidimensionalen projektiven
Geometrie \"uber einem kommutativen K\"orper, der wenigstens drei
Elemente besitzt.

\medskip
\noindent
{\bf 4.7. Satz.} {\it Es sei $K$ ein kommutativer K\"orper, der
von $GF(2)$ verschieden sei. Es sei ferner $V$ ein Vektorraum des
Ranges $4$ \"uber $K$ und $S$ sei eine Segresche Mannigfaltigkeit
mit den Invarianten $2,2$ in $\La(V)$ und $E_1$ und $E_2$ seien
die beiden Geraden aus $E_1$ und $H_1, H_2, H_3, H_4$ vier
verschiedene Geraden aus $E_2$. {\rm (}Damit es solche vier Geraden gibt,
braucht man, dass $K$ mindestens drei Elemente enth\"alt.{\rm )} Setze
$$\eqalign{
P_0 &:= G_0 \cap H_1\qquad Q_0 := G_2\cap H_0\cr
P_1 &:= G_0 \cap H_0\qquad Q_1 := G_2 \cap H_1\cr
P_2 &:= G_1 \cap H_3\qquad Q_2 := G_3 \cap H_2\cr
P_3 &:= G_1 \cap H_2\qquad Q_3 := G_3 \cap H_3.\cr}
$$
Dann sind die Tetraeder $P_0, P_1, P_2, P_3$ und $Q_0, Q_1, Q_2,
Q_3$ ein M\"obiuspaar.}

\smallskip

Beweis. Es ist, da die $G_i \cap H_j$ ja Punkte sind, die f\"ur
verschiedene Indexpaare auch verschieden sind,
$$\eqalign{
P_0+P_1+P_2 &= G_0 \cap H_1 + G_0 \cap H_0 + G_1 \cap H_3\cr
&= G_0 + G_1 \cap H_3\cr &=G_0 + G_0 \cap H_3 + G_1 \cap H_3
= G_0 + H_3,}$$
so dass $Q_3 \leq P_0 + P_1+P_2$ gilt. Ferner ist
$$\eqalign{
P_1+P_2+P_3 &=G_0 \cap H_0 + G_1 \cap H_3 + G_1 \cap H_2\cr
&= G_0 \cap H_0 + G_1\cr &=G_0 \cap H_0 + G_1 \cap H_0 +  G_1
= H_0 + G_1,}$$
so dass $Q_0 \leq P_1 + P_2 + P_3$ ist. Ganz analog beweist man,
dass $Q \leq P_2 + P_3 +P_0$ und $Q_3 \leq P_2 + P_0 + P_1$ gilt.

Ferner gilt
$$\eqalign{
Q_0 + Q_1 + Q_2 &= H_0 \cap G_2 + H_1 \cap G_2 + H_2 \cap G_3\cr &=
G_2 + H_2 \cap G_2 + H_2 \cap G_3
= G_2 + H_2,}$$
so dass $P_3 \leq Q_0 + Q_1+Q_2$ ist. Entsprechend beweist man,
dass $P_0 \leq Q_1 + Q_2 + Q_3, P_1 \leq Q_2 + Q_3 + Q_0$ ist.
Entsprechend beweist man, dass $P_0 \leq Q_1+Q_2+Q_3, P_1 \leq
Q_2, Q_3+Q_0$ und $P_2 \leq  Q_3 + Q_0 + Q_1$ ist.

\medskip

Zum Schlu\ss\ ziehen wir noch eine Folgerung \"uber orthogonale
Gruppen aus dem Vorstehenden.

\medskip
\noindent
{\bf 4.8. Satz.} {\it Es sei $K$ ein kommutativer K\"orper und $V$
sei ein Vektorraum des Ranges $4$ \"uber $K$. Ist $Q$ eine nicht
ausgeartete, spurwertige quadratische Form vom Index $2$ auf $V$,
so enth\"alt die von $O(V,Q)$ auf $\La(V)$ induzierte
Kol\-li\-ne\-a\-ti\-ons\-grup\-pe $PO(V,Q)$ einen Normalteiler vom Index $2$,
welcher zu $PGL(2,K) \times PGL(2,K)$ isomorph ist.}
\smallskip
      Beweis. Dies folgt aus der Isometrie der quadratischen Formen von
maximalem Index (Satz V.6.6) zusammen mit Satz 4.5 und Satz 3.13.


 \newpage
       
 \mychapter{VII}{Gra\ss{}mannsche Mannigfaltigkeiten}

\noindent
       In diesem Kapitel definieren und studieren wir die gra\ss mannschen
 Mannigfaltigkeiten eines Vektorraumes. Das sind die Gebilde aus den
 Unterr\"aumen des Ranges $k$. Einen sehr wesentlichen Satz \"uber diese
 Gebilde haben wir schon bewiesen, n\"amlich den Satz von Chow (Satz I.8.4).
 Im vorliegenden Kapitel setzen wir nun voraus, dass die
 Koordinatenk\"orper kommutativ sind. Das versetzt uns in die Lage, die
 Unterr\"aume des Ranges $k$ einer projektiven Geometrie des Ranges $n$
 durch Punktmannigfaltigkeiten in einem projektiven Raum des Ranges
 $n \choose k$ darzustellen. Um dieses Programm durchzuf\"uhren, ben\"otigen
 wir die Gra\ss mannalgebra eines Vektorraumes, die wir mit Hilfe der
 Tensoralgebra des Vektorraumes definieren werden, so dass wir uns zun\"achst
 mit der Tensoralgebra eines Moduls \"uber einem kommutativen Ring
 besch\"aftigen werden.
 \par
       Bei all diesen Untersuchungen ist es ganz wesentlich vorauszusetzen,
 dass die Moduln Moduln \"uber kommutativen Ringen sind. Ist n\"amlich $M$
 ein Rechts\-mo\-dul \"uber einem kommutativen Ring $R$, so wird $M$ durch die
 Vorschrift $rm := mr$ f\"ur $r \in R$ und $m \in M$ zu einem zweiseitigen
 $R$-Modul, f\"ur den \"uberdies $rm = mr$ f\"ur alle $r \in R$ und alle
 $m \in M$ gilt. Diese Eigenschaften werden wir st\"andig auszunutzen haben.
 Setzt man andererseits voraus, dass $M$ ein zweiseitiger $R$-Modul ist, f\"ur
 den $rm = mr$ f\"ur alle $m$ und $r$ gilt, so folgt mit
 $$ (rs)m = m(rs) = (mr)s = s(rm) = (sr)m, $$
 dass $R/\ann(M)$ kommutativ ist.

 \mysection{1. Die Gra\ss mannalgebra eines Moduls}

\noindent
       Im Folgenden sei $R$ ein kommutativer Ring mit Eins. Ist $M$ ein
 $R$-Modul, so nehmen wir, wie schon verabredet, an, dass $M$ ein unit\"arer
 $R$-Bimodul sei mit $rm = mr$ f\"ur alle $r \in R$ und alle $m \in M$.  Wir
 definieren eine Rekursionsregel $\rho$ durch $\rho(A,M) := A \otimes_R M$
 f\"ur alle $R$-Moduln $A$. Dann ist, wie vor Satz VI.1.9 beschrieben, auch
 $\rho(A,M)$ ein $R$-Bimodul. Nach dem Dedekindschen Rekursionssatz gibt
 es daher eine Abbildung $T_R$ mit $T_R^0(M) = R$ und
 $$ T_R^{n+1}(M) = T_R^n(M) \otimes_R M. $$
 Mittels VI.1.8 und VI.1.7 folgt die Existenz eines Isomorphismus
 $\sigma_{m,n}$ von $T_R^m(M) \otimes T_R^n(M)$ auf $T_R^{m+n}(M)$ mit
 $$ \sigma_{m,n}((x_1 \otimes\cdots\otimes x_m) \otimes (y_1 \otimes\cdots
              \otimes y_n)) = x_1 \otimes\cdots\otimes x_m \otimes y_1
              \otimes\cdots\otimes y_n $$
 f\"ur alle $x_1$, \dots, $x_m$, $y_1$, \dots, $y_n \in M$. Dabei ist etwa
 $x_1 \otimes x_2 \otimes x_3 \otimes x_4$ als
 $((x_1 \otimes x_2) \otimes x_3) \otimes x_4$ zu lesen. Wir setzen nun
 $$ T_R(M) := \bigoplus_{i:=0}^\infty T_R^i(M). $$
 Sind $f$, $g \in T_R(M)$, so definieren wir $fg$ durch
 $$ fg := \sum_{n:=0}^\infty \sum_{i:=0}^n \sigma_{i,n-i}(f_i \otimes g_{n-i}).
       $$
 Dann ist $T_R(M)$ eine $R$-Algebra mit Eins, die sog. {\it Tensoralgebra\/}
 von $M$. Die Multiplikation in $T_R(M)$ ist assoziativ.
 \par
       Wir haben stillschweigend das Element $y \in T_R^i(M)$ mit der Folge
 aus $T_R(M)$ identifiziert, die an der $i$ten Stelle den Wert $y$ hat und
 an allen \"ubrigen Stellen Null ist. Diese Identifizierung machen wir nun
 der Deutlichkeit halber zum Teil wieder r\"uckg\"angig, indem wir mit
 $\epsilon(y)$ die Folge $f$ bezeichnen, f\"ur die $f_1 = y$ und $f_i = 0$
 gilt f\"ur alle von 1 verschiedenen $i$. Dabei ist mit $y$ nat\"urlich
 ein Element aus $T_R^1(M)$ gemeint.
 \medskip\noindent
 {\bf 1.1. Satz.} {\it Es sei $R$ ein kommutativer Ring mit Eins und $M$ sei ein
 $R$-Modul und $A$ eine assoziative $R$-Algebra mit Eins. Ist $\mu$ ein
 Modulhomomorphismus von $M$ in $A$, so gibt es genau einen
 Algebrenhomomorphismus $\alpha$ von $T_R(M)$ in $A$ mit
 $\mu = \alpha\epsilon$.}
 \smallskip
       Beweis. Auf Grund der Definition der Tensoralgebra ist
 $$ x_1 \otimes\cdots\otimes x_n = \epsilon(x_1) \dots \epsilon(x_n) $$
 f\"ur alle $n$ und alle $n$-Tupel $(x_1, \dots, x_n)$ mit $x_i \in M$. Sind
 $\alpha$ und $\beta$ nun Algebrenhomomorphismen von $T_R(M)$ in $A$ mit
 $\alpha\epsilon = \beta\epsilon$, so folgt
 $$\eqalign{
 \alpha(x_1 \otimes\cdots\otimes x_n)
       &= \alpha(\epsilon(x_1) \dots \epsilon(x_n)) \cr
       &= \alpha\epsilon(x_1) \dots \alpha\epsilon(x_n) \cr
       &= \beta\epsilon(x_1) \dots \beta\epsilon(x_n) \cr
       &= \beta(\epsilon(x_1) \dots \epsilon(x_n) \cr
       &= \beta(x_1 \otimes\cdots\otimes x_n). \cr} $$
 Da $T_R(M)$ von der Menge der $x_1 \otimes\cdots\otimes x_n$ als Modul
 erzeugt wird, folgt $\alpha = \beta$, so dass die Einzigkeit von $\alpha$
 bewiesen ist.
 \par
       Es sei $n > 0$. Die Abbildung $\lambda_n$, die durch
 $$ \lambda_n(x_1, \dots, x_n) := \mu(x_1) \dots \mu(x_n) $$
 definiert wird, ist multilinear. Es gibt daher eine $R$-lineare Abbildung
 $\gamma_n$ von $T_R^n(M)$ in $A$ mit
 $$ \gamma_n(x_1 \otimes\cdots\otimes x_n) = \mu(x_1) \dots \mu(x_n). $$
 Ist $e$ die Eins von $A$, so setzen wir noch
 $$ \gamma_0(r) := re $$
 f\"ur alle $r \in R$. Mittels dieser Abbildungen definieren wir $\alpha$
 durch
 $$ \alpha(f) := \sum_{i:=0}^\infty \gamma_i(f_i) $$
 f\"ur alle $f \in T_R(M)$. Es ist klar, dass $\alpha$ eine lineare
 Abbildung von $T_R(M)$ in $A$ ist. Nach Konstruktion von $\alpha$ gilt
 ferner $\alpha\epsilon = \mu$ und $\alpha(1) = e$. Es bleibt zu zeigen, dass
 $\alpha$ multiplikativ ist. Sind $x_1$, \dots, $x_m$, $y_1$, \dots,
 $y_n \in M$, so folgt (hier erst benutzen wir die Assoziativit\"at der
 Multiplikation in A)
 $$\eqalign{
 \alpha((x_1 \otimes\cdots\otimes x_m)(y_1 \otimes\cdots\otimes y_n))  
       &= \alpha(x_1 \otimes\cdots\otimes x_m \otimes y_1 \otimes\cdots
              \otimes y_n)  \cr
       &= \mu(x_1) \dots \mu(x_m)\mu(y_1) \dots \mu(y_n)  \cr
       &= \alpha(x_1 \otimes\cdots\otimes x_m)
            \alpha(y_1 \otimes\cdots\otimes y_n). \cr} $$
 Hieraus folgt zusammen mit der Linearit\"at von $\alpha$ die
 Multiplikativit\"at von $\alpha$.
 \medskip
       Der n\"achste Satz zeigt, dass die Tensoralgebra eines freien Moduls
 ein wohl\-be\-kann\-tes Objekt ist. Er wird uns im n\"achsten Kapitel an einer
 entscheidenden Stelle weiterhelfen.
 \medskip\noindent
 {\bf 1.2. Satz.} {\it Es sei $R$ ein kommutativer Ring mit Eins und $F$ sei ein
 freier $R$-Modul. Ferner sei $\epsilon$ die vor Satz 1.1 definierte Abbildung
 von $F$ in $T_R(F)$. Ist $B$ eine Basis von $F$ und bezeichnet $\pol_R(B)$ die
 freie assoziative Algebra \"uber $R$ in der Menge der Unbestimmten $B$, so gibt
 es einen Isomorphismus $\varphi$ von $\pol_R[B]$ auf $T_R(F)$ mit
 $\varphi(b) = \epsilon(b)$ f\"ur alle $b \in B$.}
 \smallskip
       Beweis. Die Definition der freien assoziativen Algebra besagt,
 dass es einen Homomorphismus  $\varphi$ von $\pol_R[B]$ in $T_R(F)$ gibt mit
 $\varphi(b) = \epsilon(b)$ f\"ur alle $b \in B$. Nach 1.1 gibt es einen
 Homomorphismus $\psi$ von $T_R(F)$ in $\pol_R[B]$ mit $\psi\epsilon(b) = b$
 f\"ur alle $b \in B$. Es folgt $\psi\varphi(b) = b$ und damit
 $$ \psi\varphi = 1_{\pol_R[B]}. $$
 Andrerseits ist $\varphi\psi\epsilon(b) = \epsilon(b)$ f\"ur alle $b \in B$
 und damit
 $$ \varphi\psi = 1_{T_R(F)}. $$
 Damit ist alles bewiesen.
 \medskip
       Es sei $R$ ein kommutativer Ring mit Eins und $M$ sei ein Modul \"uber
 $R$. Ferner sei $I$ das zweiseitige Ideal von $T_R(M)$, welches von allen
 Elementen der Form $x \otimes x$ mit $x \in M$ erzeugt wird. Wir setzen
 $$ {\textstyle\bigwedge}_R(M) := T_R(M)/I $$
 und nennen $\bigwedge_R(M)$ die {\it Gra\ss mannalgebra\/} von $M$. Wir
 setzen ferner
 $$ {\textstyle\bigwedge}_R^n(M) := (T_R^n(M) + I)/I. $$
 Die Multiplikation in $\bigwedge_R(M)$ bezeichnen wir mit $\wedge$.
 \par
       Offenbar gilt $I \cap T_R^0(M) = \{0\} = I \cap T_R^1(M)$. Wir
 d\"urfen daher $r \in R$ mit $r + I$ und $m \in M$ mit $m + I$ identifizieren.
 \medskip\noindent
 {\bf 1.3. Satz.} {\it Es sei $R$ ein kommutativer Ring mit Eins und $M$ sei ein
 Modul \"uber $R$. Ferner sei $A$ eine assoziative $R$-Algebra mit Eins. Ist
 $\varphi$ eine $R$-lineare Abbildung von $M$ in $A$ und gilt $\varphi(x)^2
 = 0$ f\"ur alle $x \in M$, so gibt es genau einen Homomorphismus $\psi$ von
 $\bigwedge_R(M)$ in $A$ mit $\psi(x) = \varphi(x)$ f\"ur alle $x \in M$.}
 \smallskip
       Beweis. Dies folgt mit Standardschl\"ussen aus 1.1.
 \medskip\noindent
 {\bf 1.4. Satz.} {\it Es sei $R$ ein kommutativer Ring mit Eins und $M$ und $N$
 seien zwei $R$-Moduln. Ist $\varphi$ ein Homomorphismus von $M$ in $N$, so
 gibt es genau einen Homomorphismus $\psi$ von $\bigwedge_R(M)$ in
 $\bigwedge_R(N)$ mit $\psi(x) = \varphi(x)$ f\"ur alle $x \in M$. Ist
 $\varphi$ surjektiv, so ist auch $\psi$ surjektiv. Ferner gilt
 $\psi(\bigwedge_R^k(M)) \subseteq \bigwedge_R^k(N)$ f\"ur alle nicht negativen
 ganzen Zahlen $k$.}
 \smallskip
       Beweis. Existenz und Eindeutigkeit von $\psi$ folgen mit 1.3. Die
 Aussage, dass $\psi(\bigwedge_R^k(M)) \subseteq \bigwedge_R^k(N)$ ist, gilt
 f\"ur $k = 1$ und folgt dann f\"ur beliebiges $k$ durch Induktion. Weil
 $\bigwedge_R^1(N)$ die Algebra $\bigwedge_R(N)$ erzeugt, ist $\psi$ surjektiv,
 falls $\varphi$ surjektiv ist.
 \medskip
       Die Abbildung $\psi$ ist nicht notwendig injektiv, wenn $\varphi$ es
 ist.
 \medskip\noindent
 {\bf 1.5. Satz.} {\it Es sei $R$ ein kommutativer Ring mit Eins und $M$ sei ein
 $R$-Modul. Sind $x_1$, \dots, $x_n \in M$ und ist $\sigma$ eine Permutation
 aus der symmetrischen Gruppe $S_n$ vom Grade $n$, so gilt
 $$ x_{\sigma(1)} \wedge \dots \wedge x_{\sigma(n)} =
 {\rm sgn}(\sigma) x_1 \wedge \dots \wedge x_n. $$
 Sind $i$ und $j$ zwei verschiedene der Indizes 1, \dots, $n$ und ist
 $x_i = x_j$, so ist}
 $$ x_1 \wedge \dots \wedge x_n = 0. $$
 \par
       Beweis. Es ist
 $$x_1 \wedge x_2 + x_2 \wedge x_1 = (x_1 + x_2) \wedge (x_2 + x_1)
       - x_1 \wedge x_1 - x_2 \wedge x_2 = 0, $$
 so dass die erste Aussage f\"ur $n = 2$ richtig ist. Induktion zeigt ihre
 G\"ultigkeit f\"ur beliebiges $n$, wenn man nur noch beachtet, dass jede
 Permutation Produkt von Transpositionen ist.
 \par
       Nach dem bereits Bewiesenen ist
 $$ x_1 \wedge \dots \wedge x_n = \pm x_i \wedge x_j \wedge y $$
 mit einem geeigneten $y$. Wegen $x_i = x_j$ ist aber $x_i \wedge x_j = 0$.
 Damit ist alles bewiesen.
 \medskip\noindent
 {\bf 1.6. Satz.} {\it Es sei $R$ ein kommutativer Ring mit Eins und $M$ sei ein
 $R$-Modul. Ist $\{b_1, \dots, b_n\}$ ein Erzeugendensystem von $M$, so ist
 $$ {\textstyle\bigwedge}_R^i(M) = \{0\} 
 \quad\mathit{f\ddot ur\ alle\ }i > n.$$}
 \smallskip
       Beweis. Es sei $i > n$ und sei $x_1$, \dots, $x_i \in M$. Es gibt dann
 $a_{jk} \in R$ mit
 $$ x_k = \sum_{j:=1}^n b_ja_{jk}. $$
 Es folgt
 $$
 x_1 \wedge \dots \wedge x_i
       = \sum_{j:=1}^n b_ja_{j1} \wedge \dots \wedge \sum_{j:=1}^n b_j a_{ji}
       = \sum_{\alpha} (b_{\alpha(1)1} \wedge \dots \wedge b_{\alpha(i)i})
                                                        A_\alpha,    
 $$
 wobei $\alpha$ die Menge der Abbildungen von $\{1, \dots, i\}$ in $\{1,
 \dots n\}$ durch\-l\"auft und $A_\alpha \in R$ gilt. Wegen $n < i$ ist
 keines der $\alpha$ injektiv, so dass nach 1.5 f\"ur alle $\alpha$ die
 Gleichung
 $$ b_{\alpha(1)1} \wedge \dots \wedge b_{\alpha(i)i} = 0 $$
 gilt. Hieraus folgt die Behauptung.
 \medskip\noindent
 {\bf 1.7. Satz.} {\it Es sei $R$ ein kommutativer Ring mit 1 und $F$ sei ein
 freier Modul \"uber $R$. Es sei weiter $B$ eine Basis von $F$, die linear
 geordnet sei. Ist $J \in \Fin(B)$, sind $b_1$, \dots, $b_n$ die Elemente
 von $J$ und gilt $b_1 < \dots < b_n$, so setzen wir
 $$ b_J := b_1 \wedge \dots \wedge b_n. $$
 Ferner setzen wir $b_\emptyset := 1$. Dann ist
 $ \{b_J \mid J \in \Fin(B)\} $
 eine Basis von $\bigwedge_R(F)$.}
 \smallskip
       Beweis. Mit 1.5 folgt, dass $\{b_J \mid J \in \Fin(B)\}$ ein
 Erzeugendensystem von $\bigwedge_R(F)$ ist. Wir m\"ussen also nur noch
 zeigen, dass $\{b_J \mid J \in \Fin(B)\}$ auch linear unabh\"angig ist.
 \par
       Es sei $\{w_J \mid J \in \Fin(B)\}$ eine mit $\Fin(B)$ indizierte
 Menge und $A$ sei ein freier $R$-Modul mit der Basis $\{w_J \mid J \in
 \Fin(B)\}$. Sind $i$, $j \in B$, so setzen wir
 $$ \langle i,j \rangle := \cases{0, &falls $i = j$, \cr
                                  1, &falls $i < j$, \cr
                                 -1, &falls $i > j$. \cr} $$
 Ferner setzen wir
 $$ w_Iw_J := \biggl(\prod_{i \in I, j \in J} \langle i,j \rangle\biggr)
                            w_{I \cup J}. $$
 Setzt man die f\"ur je zwei Basiselemente definierte Multiplikation linear
 fort, so wird $A$ eine $R$-Algebra, wie wir nun sehen werden. Dazu gen\"ugt
 es zu zeigen, dass $A$ eine Eins besitzt und dass
 $$ (w_Iw_J)w_L = w_I(w_Jw_L) $$
 ist f\"ur alle $I$, $J$, $L \in \Fin(B)$.
 \par
       Offenbar ist $w_\emptyset$ das Einselement von $A$.
 \par
       Um die Assoziativit\"at zu beweisen, seien $I$, $J$, $L \in \Fin(B)$.
 Ist $I \cap J \neq \emptyset$,
 so ist $w_Iw_J = 0$ und daher $(w_Iw_J)w_L = 0$. Andererseits folgt aus
 $I \cap J \neq \emptyset$, dass auch $I \cap (J \cup L) \neq \emptyset$
 ist. Daher ist
 $$ w_I(w_Jw_L) = \pm w_I w_{J \cup L} = 0. $$
 In diesem Falle ist also $(w_Iw_J)w_L = w_I(w_Jw_L)$. Der Fall $I \cap L \neq
 \emptyset$ erledigt sich analog. Es sei also $I \cap J = \emptyset$ und
 $J \cap L = \emptyset$. Dann ist
 $$ \prod_{k \in I \cup J, l \in L} \langle k,l \rangle =
 \prod_{k \in I,l \in L} \langle k,l \rangle \prod_{k \in J,l \in L}
                            \langle k,l \rangle $$
 und daher
 $$\eqalign{
 (w_Iw_J)w_L &= \prod_{i \in I,j \in J} w_{I\cup J}w_L             \cr
       &= \prod_{i\in I,j \in J} \langle i,j \rangle
              \prod_{k \in I \cup J,l \in L} \langle k,l \rangle 
                            w_{I \cup J \cup L}                    \cr
       &= \prod_{i \in I,j \in J} \langle i,j \rangle
              \prod_{k \in I , l \in L} \langle k,l \rangle
              \prod_{k \in J , l \in L} \langle k,l \rangle
                                   w_{I\cup J \cup L}.             \cr} $$
 Andererseits ist, da ja $J \cap L = \emptyset$ ist,
 $$\eqalign{
 w_I(w_Jw_L) &= \prod_{k \in J,l \in L} \langle k,l \rangle
                            w_Iw_{J \cup L}                        \cr
       &= \prod_{k \in J,l \in L} \langle k,l \rangle \prod_{i \in I,j \in
              J \cup L} \langle i,j \rangle w_{I \cup J \cup L}    \cr
       &= \prod_{k \in J,l \in L} \langle k,l \rangle
              \prod_{i \in I , j \in J} \langle i,j \rangle
              \prod_{k \in J , l \in L} \langle k,l \rangle
                                   w_{I\cup J \cup L}.             \cr} $$
 Also ist $(w_Iw_J)w_L = w_I(w_Jw_L)$, so dass $A$ in der Tat eine
 $R$-Algebra ist.
 \par
       Es gibt einen Homomorphismus $\alpha$ von $F$ in $A$ mit
 $\alpha(b) = w_{\{b\}}$. Ist $f = \sum_{b\in B} bf_b$ so folgt
 $$\eqalign{
 \alpha(f)^2 &= \sum_{b \in B}\sum_{c \in B} w_{\{b\}}w_{\{c\}}  \cr
       &= \sum_{b \in B} w^2_{\{b\}} f_b^2 +
                \sum_{b,c \in B, b<c} (w_{\{b\}}w_{\{c\}} +
                                       w_{\{c\}}w_{\{b\}}) f_bf_c 
       = 0.  \cr} $$
 Es gibt daher einen Homomorphismus $\beta$ von $\bigwedge_R(F)$ in $A$ mit
 $$ \alpha(b) = \beta(b) $$
 f\"ur alle $b \in B$. Hieraus folgt $\beta(b_I) = w_I$ f\"ur alle $I \in
 \Fin(B)$, so dass die $b_I$ in der Tat eine Basis von $\bigwedge_R(F)$
 sind. Es folgt weiter, dass $\beta$ ein Isomorphismus ist.
       Damit ist alles bewiesen.
 \medskip
       Der Beweis des Satzes liefert auch noch die beiden folgenden
 Korollare.
 \medskip\noindent
 {\bf 1.8. Korollar.} {\it Es sei $R$ ein kommutativer Ring mit Eins und $F$ sei
 ein freier $R$-Modul. Dann ist
 $$ {\textstyle\bigwedge}_R(F) = \bigoplus_{k:=0}^\infty
                            {\textstyle\bigwedge}_R^k(F). $$}
 \par\noindent
 {\bf 1.9. Korollar.} {\it Es sei $R$ ein kommutativer Ring mit Eins und $F$ sei
 ein freier $R$-Modul des Ranges $n$. Dann ist
 $$ \Rg_R({\textstyle\bigwedge}_R(F)) = 2^n $$
 und
 $$ \Rg_R({\textstyle\bigwedge}_R^k(F)) = {n \choose k} $$
 f\"ur alle ganzen Zahlen $k$ mit $0 \leq k \leq n$.}
 \medskip
       Ist $B$ unendlich, so hat $\Fin(B)$ und auch die Menge $\Fin_k(B)$ der
 $k$-Teil\-men\-gen von $B$ die gleiche M\"achtigkeit wie $B$, so dass
 $\bigwedge_R(F)$ und $\bigwedge_R^k(F)$ in diesem Falle den gleichen Rang
 wie $F$ haben.
 \medskip\noindent
 {\bf 1.10. Satz.} {\it Es sei $K$ ein kommutativer K\"orper und $V$ sei ein
 Vektorraum \"uber $K$. Sind $v_1$, \dots, $v_r \in V$, so sind
 $v_1$, \dots, $v_r$ genau dann linear unabh\"angig, wenn
 $v_1 \wedge \dots \wedge v_r \neq 0$ ist.}
 \smallskip
       Beweis. Sind $v_1$, \dots, $v_r$ linear unabh\"angig, so sind sie
 Teil einer Basis von $V$. Dann geh\"ort aber $v_1 \wedge \dots \wedge v_r$
 nach 1.7 zu einer Basis von $\bigwedge_R(V)$, ist also von 0 verschieden.
 \par
       Sind $v_1$, \dots, $v_r$ linear abh\"angig, so ist oBdA
 $$ v_r = \sum_{i:=1}^{r-1} v_i\lambda_i $$
 und nach 1.5 folglich
 $$ v_1 \wedge \dots \wedge v_r =
 \sum_{i:=1}^{r-1} ((v_1 \wedge \dots \wedge v_{r-1}) \wedge v_i)\lambda_i
       = 0. $$
 Damit ist alles bewiesen.
 \medskip
       Wir schlie\ss en diesen Abschnitt mit einem Satz, den wir nur so
 formulieren werden, wie wir ihn sp\"ater ben\"otigen. Zuvor jedoch noch eine
 Definition. Ist $u \in \bigwedge_K(V)$, so hei\ss t $v$ {\it zerlegbar\/},
 falls es $v_1$, \dots, $v_r \in V$ gibt mit $u = v_1 \wedge \dots \wedge
 v_r$.
 \medskip\noindent
 {\bf 1.11. Satz.} {\it Es sei $V$ ein Vektorraum \"uber dem
 kommutativen K\"orper $K$. Ist $\sigma$ ein Endomorphismus von $V$, so gilt:
 \item{a)} Es gibt genau einen Endomorphismus $\sigma_{\#r}$ von
 $\bigwedge_K^r(V)$ mit
 $$ \sigma_{\#r}(v_1 \wedge \dots \wedge v_r) = \sigma(v_1) \wedge \dots
                                                    \wedge \sigma(v_t) $$
 f\"ur alle $v_1$, \dots, $v_r \in V$. Insbesondere bildet $\sigma_{\#r}$
 zerlegbare Vektoren auf zerlegbare Vektoren ab.
 \item{b)} Es ist $(\sigma\tau)_{\#r} = \sigma_{\#r}\tau_{\#r}$, falls $\tau$
 ein weiterer Endomorphismus von $V$ ist.
 \item{c)} F\"ur jeden Automorphismus $\sigma$ von $V$\/ ist $\sigma_{\#r}$ ein
 Automorphismus von $\bigwedge_K^r(V)$.
 \item{d)} Ist $u \in \bigwedge_K^r(V)$ und $v \in \bigwedge_K^s(V)$, so ist}
 $$ \sigma_{\#r}(u) \wedge \sigma_{\#s}(v) = \sigma_{\#(r+s)}(u \wedge v). $$
 \par
       Beweis. a) Es seien $v_1$, \dots, $v_r \in V$. Nach Satz 1.4 gibt es
 einen Endomorphismus $\rho$ von $\bigwedge_K(V)$ mit $\sigma(v) = \rho(v)$
 f\"ur alle $v \in V$. Es sei $\sigma_{\#r}$ die Einschr\"ankung von $\rho$ auf
 $\bigwedge_K^r(V)$. Dann ist $\sigma_{\#r}$ nach 1.4 ein Endomorphismus von
 $\bigwedge_K^r(V)$ und es gilt
 $$\eqalign{
 \sigma_{\#r}(v_1 \wedge \dots \wedge v_r)
       &= \rho(v_1 \wedge \dots \wedge v_r)            \cr
       &= \rho(v_1) \wedge \dots \wedge \rho(v_r)      \cr
       &= \sigma(v_1) \wedge \dots \wedge \sigma(v_r). \cr} $$
 \par
       b) Es ist
 $$\eqalign{
 (\sigma\tau)_{\#r}(v_1 \wedge \dots \wedge v_r)
       &= \sigma\tau(v_1) \wedge \dots \wedge \sigma\tau(v_r)      \cr
       &= \sigma_{\#r}(\tau(v_1) \wedge \dots \wedge \tau(v_r))     \cr
       &= \sigma_{\#r}\tau_{\#r}(v_1 \wedge \dots \wedge v_r).       \cr} $$
 Da $\bigwedge_K^r(V)$ von seinen zerlegbaren Vektoren erzeugt wird, ist somit
 $$ (\sigma\tau)_{\#r} = \sigma_{\#r}\tau_{\#r}. $$
 \par
       c) Es ist $(1_V)_{\#r} = 1_{\bigwedge_K^r(V)}$. Ist nun $\sigma$ ein
 Automorphismus von $V$, so folgt mit b)
 $$ 1_{\bigwedge_K^r(V)} = (\sigma\sigma^{-1})_{\#r}
       = \sigma_{\#r}(\sigma^{-1})_{\#r} $$
 und
 $$ 1_{\bigwedge_K^r(V)} = (\sigma^{-1}\sigma)_{\#r}
       = (\sigma^{-1})_{\#r}\sigma_{\#r}. $$
 Also ist auch $\sigma_{\#r}$ ein Automorphismus.
 \par
       d) Weil $u$ Linearkombination von zerlegbaren Vektoren der L\"ange
 $r$ und $v$ Linearkombination von zerlegbaren Vektoren der L\"ange $s$ ist,
 ist $u \wedge v$ Linearkombination von zerlegbaren Vektoren der L\"ange
 $r + s$, so dass $u \wedge v \in \bigwedge_K^{r+s}(V)$ ist. Hat nun $\rho$
 wieder die Bedeutung wie beim Beweise von a), so ist
 $$ \sigma_{\#r}(u) \wedge \tau_{\#s}(v) = \rho(u) \wedge \rho(v) = \rho(u
       \wedge v) = \sigma_{\#(r+s)}(u \wedge v). $$
 Hieraus folgt die Behauptung unter d).

 \mysection{2. Dualit\"at in der Gra\ss mannalgebra}
 
\noindent
       Es sei $K$ ein kommutativer K\"orper und $V$ und $W$ seien zwei
 Vektorr\"aume \"uber $K$. Die Abbildung $f$ des $r$-fachen cartesischen
 Produktes $V^r$ von $V$ mit sich selbst in $W$ hei\ss t {\it r-fach
 alternierende Abbildung\/} von $V$ in $W$, falls $f$ in jedem Argument
 linear ist und \"uberdies $f(v_1, \dots, v_r) = 0$ ist, falls zwei der
 Argumente gleich sind. Die Menge aller $r$-fach linearen Abbildungen von
 $V$ in $W$ bezeichnen wir mit $\Alt_K^r(V,W)$.
 \medskip\noindent
 {\bf 2.1. Satz.} {\it Es seien $V$ und $W$ Vektorr\"aume \"uber dem kommutativen
 K\"orper $K$. Ist $f \in \Alt_K^r(V)$ und sind $v_1$, \dots, $v_r$ linear
 abh\"angige Vektoren aus $V$, so ist $f(v_1, \dots, v_r) = 0$.}
 \smallskip
       Beweis. Es sei oBdA $v_r = \sum_{i:=1}^{r-1} v_i\alpha_i$. Dann ist
 $$ f(v_1, \dots, v_r) = \sum_{i:=1}^{r-1} f(v_1, \dots, v_{r-1},v_i)
                            \alpha_i = 0. $$
 \par\noindent
 {\bf 2.2. Korollar.} {\it Es seien $V$ und $W$ Vektorr\"aume \"uber dem
 kommutativen K\"orper $K$. Ist $V$ endlich erzeugt und ist $r > \Rg_K(V)$,
 so ist $\Alt_K^r(V) = \{0\}$.}
 \medskip
       Der Beweis des n\"achste Satzes ist eine einfache \"Ubungsaufgabe, wenn
 man nur beachtet, dass die symmetrische Gruppe
 $S_r$ von ihren Transpositionen erzeugt wird.
 \medskip\noindent
 {\bf 2.3. Satz.} {\it Es sei $K$ ein kommutativer K\"orper und $V$ und $W$ seien
 Vektorr\"aume \"uber $K$. Sind $v_1$, \dots, $v_r \in V$, ist $\pi \in S_r$
 und ist $f \in \Alt_K^r(V,W)$, so ist}
 $$ f(v_{\pi(1)}, \dots, v_{\pi(r)}) = \sgn(\pi) f(v_1, \dots, v_r). $$

       Mit $\Lin_K^r(V,W)$ bezeichnen wir die Menge der $r$-fach linearen
 Abbildungen von $V$ in $W$. Ist $f \in \Lin_K^r(V,W)$, so setzen wir
 $$ fa(v_1, \dots, v_r) := \sum_{\pi \in S_r} \sgn(\pi) f(v_{\pi(1)}, \dots,
 v_{\pi(r)}). $$
 Die Abbildung $fa$ hei\ss t die {\it Antisymmetrisierte\/} von $f$. Die
 Menge der antisymmetrisierten, $r$-fach linearen Abbildungen von $V$ in
 $W$ bezeichnen wir mit $\Ant_K^r(V,W)$.
 \medskip\noindent
 {\bf 2.4. Satz.} {\it Es seien $V$ und $W$ Vektorr\"aume \"uber dem kommutativen
 K\"orper $K$. Ist $r$ eine nat\"urliche Zahl, so ist}
 $$ \Alt_K^r(V,W) = \Ant_K^r(V,W). $$
 \par
       Beweis. Es sei $f \in \Lin_K^r(V,W)$. Ferner seien $v_1$, \dots,
 $v_r \in V$ und es gebe $i$, $j$ mit $i \neq j$ und $v_i = v_j$. Es sei
 $\tau$ die Transposition, die $i$ mit $j$ vertauscht. Dann ist
 $S_r = A_r \cup \tau A_r$ und $A_r \cap \tau A_r = \emptyset$. Es folgt
 $$ fa(v_1, \dots, v_r)
       = \sum_{\pi \in A_r} (f(v_{\pi(1)}, \dots, v_{\pi(r)}) -
              f(v_{\tau\pi(1)}, \dots, v_{\tau\pi(r)}). $$
 Ist $\pi(k) \neq i$, $j$, so ist $v_{\tau\pi(k)} = v_{\pi(k)}$. Ist
 $\pi(k) = i$, so folgt
 $$ v_{\tau\pi(k)} = v_{\tau(i)} = v_j = v_i = v_{\pi(k)}, $$
 und ist $\pi(k) = j$, so folgt
 $$ v_{\tau\pi(k)} = v_{\tau(j)} = v_i = v_j = v_{\pi(k)}. $$
 Also ist $v_{\tau\pi(k)} = v_{\pi(k)}$ f\"ur alle $k$ und damit
 $$ fa(v_1, \dots, v_r) = 0. $$
 Dies zeigt, dass $\Ant_K^r(V,W) \subseteq \Alt_K^r(V,W)$ ist.
 \par
       Es sei $g \in \Alt_K^r(V,W)$. Es sei ferner
 $B$ eine Basis von $V$, die linear geordnet sei. Schlie\ss lich bezeichne
 $\Gamma$ die Menge aller Abbildungen von $\{1, \dots, r\}$ in $B$. Die
 Menge aller $r$-Tupel $(\gamma_1, \dots, \gamma_r)$ bildet eine Basis von
 $V^r$. Es gibt daher genau ein $f \in \Lin_K^r(V,W)$ mit
 $$ f(\gamma_1, \dots, \gamma_r) := \cases{g(\gamma_1, \dots, \gamma_r), &
                                 falls $\gamma_1 < \dots < \gamma_r$, \cr
                                          0, &sonst  \cr} $$
 f\"ur alle $\gamma \in \Gamma$. Es ist
 $$ fa(\gamma_1, \dots, \gamma_r) =
  \sum_{\pi \in S_r} \sgn(\pi) f(\gamma_{\pi(1)}, \dots, \gamma_{\pi(r)}). $$
 Ist $\gamma$ nicht injektiv, so ist auch $\gamma\pi$ nicht injektiv f\"ur
 alle $\pi \in S_r$. In diesem Falle ist also
 $$ fa(\gamma_1, \dots, \gamma_r) = 0 = g(\gamma_1, \dots, \gamma_r). $$
 Ist $\gamma$ injektiv, so gibt es genau ein $\rho \in S_r$ mit
 $$ \gamma_{\rho(1)} < \dots < \gamma_{\rho(r)}. $$
 Daher ist
 $$\eqalign{
 fa(\gamma_1, \dots, \gamma_r)
       &= \sgn(\rho)f(\gamma_{\rho(1)}, \dots, \gamma_{\rho(r)}) \cr
       &= \sgn(\rho)g(\gamma_{\rho(1)}, \dots, \gamma_{\rho(r)}) \cr
       &= g(\gamma_1, \dots, \gamma_r). \cr} $$
 Somit stimmen $fa$ und $g$ auf einer Basis von $V^r$ \"uberein, so dass
 $fa = g$ ist.
 \medskip\noindent
 {\bf 2.5. Satz.} {\it Es seien $V$ und $W$ Vektorr\"aume \"uber dem kommutativen
 K\"orper $K$. Es sei $B$ eine Basis von $V$, die linear geordnet sei. Ferner
 sei\/ $\Gamma$ die Menge aller streng monoton steigenden Abbildungen von
 $\{1, \dots, r\}$ in $B$. Ist $\eta$ eine Abbildung von $\Gamma$ in $W$,
 so gibt es genau ein $f \in \Alt_K^r(V,W)$ mit
 $$ f(\beta_1, \dots, \beta_r) = \eta(\beta) $$
 f\"ur alle $\beta \in \Gamma$.}
 \smallskip
       Beweis. Sind $v_1$, \dots, $v_r \in V$, so ist
 $$ v_j = \sum_{b\in B} ba_{bj} $$
 mit einer Matrix $a$ \"uber $K$. Diese Bezeichnung wollen wir im Folgenden
 festhalten.
 \par
       Es sei $\Delta$ die Menge aller Abbildungen von $\{1, \dots, r\}$ in
 $B$. Um die Ein\-zig\-keits\-aus\-sage des Satzes zu beweisen, sei
 $f \in \Alt_K^r(V,W)$ und 
 $f(b_1, \dots, b_r) = \eta(\beta)$ f\"ur
 alle $\beta \in \Gamma$. Dann ist
 $$\eqalign{
 f(v_1, \dots, v_r) &= \sum_{\beta \in \Delta} f(\beta_1, \dots \beta_r)
                            \prod_{i:=1}^r a_{\beta_ii} \cr
       &= \sum_{\beta \in \Gamma} \sum_{\pi \in S_r}
              f(\beta_{\pi(1)}, \dots, \beta_{\pi(r)}) \prod_{i:=1}^r
                                   a_{\beta_ii}         \cr
       &= \sum_{\beta \in \Gamma} \eta(\beta) \sum_{\pi \in S_r} \sgn(\pi)
                     \prod_{i:=1}^r a_{\beta_ii}.       \cr} $$
 Da die rechte Seite nur von $\beta$ und $a$ abh\"angt, ist die Einzigkeit
 von $f$ gezeigt.
 \par
       Definiert man umgekehrt $f$ durch
 $$ f(v_1, \dots, v_r)
       := \sum_{\beta \in \Gamma} \eta(\beta) \sum_{\pi \in S_r} \sgn(\pi)
                     \prod_{i:=1}^r a_{\beta_ii}, $$
 so verifiziert man leicht, dass $f$ alle gew\"unschten Eigenschaften
 hat.
 \medskip\noindent
 {\bf 2.6. Satz.} {\it Es seien $V$ und $W$ endlich erzeugte Vektorr\"aume \"uber
 dem kommutativen K\"orper $K$. Setze $n := \Rg_K(V)$ und $m := \Rg_K(W)$.
 Dann ist}
 $$ \Rg_K(\Alt_K^r(V,W)) = {n \choose r}m. $$
 \par
       Beweis. Nach 2.2 gilt diese Formel jedenfalls dann, wenn $r > n$ ist, da
 dann ja ${n \choose r} = 0$ ist.
 \par
       Es sei $r \leq n$. Es sei $B$ eine linear geordnete Basis von $V$ und
 $C$ sei eine Basis von $W$. Es sei ferner $\Gamma$ die Menge der streng
 monoton wachsenden Abbildungen von $\{1, \dots, r\}$ in $B$. Ist $\beta \in
 \Gamma$ und $c \in C$, so definieren wir die Abbildung $\gamma_{\beta,c}$
 von $\Gamma$ in $W$ durch
 $$ \gamma_{\beta,c}(\delta) := \cases{c, &falls $\beta = \delta$,   \cr
                                       0, &falls $\beta \neq \delta$. \cr} $$
 Nach 2.5 gibt es genau ein $f_{\beta,c} \in \Alt_K^r(V,W)$ mit
 $$ f_{\beta,c}(\delta_1, \dots, \delta_r) := \cases{c, &falls
                                                  $\beta = \delta$,   \cr
                                       0, &falls $\beta \neq \delta$. \cr} $$
 Wir zeigen, dass die $f_{\beta,c}$ eine Basis von $\Alt_K^r(V,W)$ bilden.
 \par
       Es sei
 $$ 0 = \sum_{\beta \in \Gamma} \sum_{c \in C} k_{\beta,c}f_{\beta,c}. $$
 Mit $\gamma \in \Gamma$ folgt
 $$
 0 = \sum_{\beta \in \Gamma} \sum_{c \in C}
                f_{\beta,c}(\gamma_1, \dots, \gamma_r)k_{\beta,c} 
   = \sum_{c \in C}\sum_{\beta \in \Gamma}
                f_{\beta,c}(\gamma_1, \dots, \gamma_r)k_{\beta,c} 
   = \sum_{c \in C} ck_{\gamma,c}.                              
 $$
 Daher ist $k_{\gamma,c} = 0$ f\"ur alle $\gamma \in \Gamma$ und alle $c \in
 C$. Somit sind die $f_{\beta,c}$ linear unabh\"angig.
 \par
       Es sei $f \in \Alt_K^r(V,W)$ und $\gamma \in \Gamma$. Dann ist
 $$ f(\gamma_1, \dots, \gamma_r) = \sum_{c \in C} ck_{c,\gamma}. $$
 Setze
 $$ g := \sum_{\gamma \in \Gamma} \sum_{c \in C} f_{\gamma,c}k_{c,\gamma}. $$
 Ist dann $\beta \in \Gamma$, so folgt
 $$ \eqalign{
 g(\beta_1, \dots, \beta_r)
  &= \sum_{\gamma \in \Gamma}\sum_{c \in C} f_{\gamma,c}(\beta_1, \dots,
                                          \beta_r)k_{c,\gamma} 
  = \sum_{c \in C} f_{\beta,c}(\beta_1, \dots,
                                          \beta_r)k_{c,\beta} \cr
  &= \sum_{c \in C} ck_{c,\beta}                              
  = f(\beta_1, \dots, \beta_r). \cr} $$
 Mit 2.5 folgt $f = g$, so dass die $f_{\beta,c}$ den Raum $\Alt_K^r(V,W)$
 auch erzeugen.
 \par
       Weil die $f_{\beta,c}$, wie wir gerade gesehen haben, eine Basis von
 $\Alt_K^r(V,W)$ bilden, ist
 $$\Rg_K(\Alt_K^r(V,W)) = |\Gamma|m = {n \choose r}m, $$
 q. o. o.
 \medskip
       Es sei $K$ ein kommutativer K\"orper und $V$ und $W$ seien zwei
 $K$-Vektorr\"aume. Ist $\sigma$ eine lineare Abbildung von $T_K^r(V)$ in $W$,
 so hei\ss t $\sigma$ eine {\it r-fach alternierende\/} Abbildung, falls
 $(v_1, \dots, v_r) \to (v_1 \otimes \dots \otimes v_r)^\sigma$ eine
 $r$-fach alternierende Abbildung ist. Mit $A(T_K^r(V),W)$ bezeichnen
 wir die Menge der alternierenden Abbildungen von $T_K^r(V)$ in $W$. Da
 jede Abbildung aus $\Alt_K^r(V,W)$ sich auf genau eine Weise durch das
 Tensorprodukt faktorisieren l\"asst, folgt, dass $A(T_K^r(V),W)$ und
 $\Alt_K^r(V,W)$ verm\"oge dieser Faktorisierung kanonisch isomorph sind.
 Ist $W = K$, so schreiben wir $AT_K^r(V)$ an Stelle des
 schwer\-f\"al\-li\-ge\-ren Ausdrucks $A(T_K^r(V),K)$.
 Ist $\Rg_K(V) = n$, so ist  $\Rg(AT_K^r(V)) = {n \choose k}$.
 \par
       Mit $NT_K^r(V)$ bezeichnen wir den Teilraum von $T_K^r(V)$, der
 von allen Vektoren $v_1 \otimes \dots \otimes v_r$ erzeugt wird, bei denen
 mindestens zwei der $v_i$ gleich sind.
 \medskip\noindent
 {\bf 2.7. Satz.} {\it Es sei $K$ ein kommutativer K\"orper und $V$ sei ein
 Vektorraum \"uber $K$. Ist $I$ das Ideal der Tensoralgebra $T_K(V)$, welches
 von allen $v \otimes v$ mit $v \in V$ erzeugt wird, so ist
 $$ NT_K^r(V) = T_K^r(V) \cap I. $$}
 \par
       Beweis. Mit Satz 1.5 folgt $NT_K^r(V) \subseteq I \cap T_K^r(V)$.
 \par
       Um die umgekehrte
 Inklusion zu beweisen, sei $w \in I \cap T_K^r(V)$. Es gibt dann zerlegbare
 Vektoren $z_1$, \dots, $z_t$, in deren Zerlegung Produkte der Form
 $v \otimes v$ vorkommen, sowie $k_1$, \dots, $k_t \in K$ mit
 $$ w = \sum_{i:=1}^t z_ik_i. $$
 Ohne die Allgemeinheit einzuschr\"anken, d\"urfen wir annehmen, dass
 $z_1$, \dots,$z_s \in T_K^r(V)$ und $z_{s+1}$, \dots, $z_t \not\in T_K^r(V)$
 ist. Mit 1.8 folgt dann, setzt man noch
 $$ C := \bigoplus_{i:=0, i \neq r}^\infty T_K^i(V), $$
 dass
 $$w - \sum_{i:=1}^s z_ik_i = \sum_{j:=s+1}^t z_jk_j 
                             \in T_K^r(V) \cap C = \{0\} $$
 ist. Also ist
 $$ w = \sum_{i:=1}^s z_ik_i \in NT_K^r(V). $$
 Damit ist der Satz bewiesen.
 \medskip\noindent
 {\bf 2.8. Satz.} {\it Es sei $V$ ein Vektorraum \"uber dem kommutativen K\"orper
 $K$. Es gibt genau einen Isomorphismus $\varphi$ von $T_K^r(V)/NT_K^R(V)$ auf
 $\bigwedge_K^r(V)$ mit}
 $$ \varphi(v_1 \otimes \dots \otimes v_r + NT_K^r(V))
                     = v_1 \wedge \dots \wedge v_r.$$
 \par
       Beweis. Es ist ja $\bigwedge_K^r(V) = (T_K^r(V) + I)/I$, wobei $I$
 wieder das Ideal der Tensoralgebra von $V$ ist, welches von den Elementen
 $v \otimes v$ mit $v \in V$ erzeugt wird. Bezeichnet man mit $\kappa$ die
 Einschr\"ankung des kanonischen Epimorphismus von $T_K^r(V) + I$ auf
 $T_K^r(V)$, so gilt
 $$ \kappa(v_1 \otimes \dots \otimes v_r) = v_1 \wedge \dots \wedge v_r, $$
 woraus weiter folgt, dass $\kappa$ surjektiv ist. Mit 2.7 folgt
 $$ \Kern(\kappa) = I \cap T_K^r(V) = NT_K^r(V), $$
 woraus mittels des ersten Isomorphiesatzes die Behauptung folgt.
 \medskip
       Es sei $V$ ein Vektorraum \"uber dem kommutativen K\"orper $K$ und
 $V^*$ sei sein Dualraum. Um Klammern zu sparen, schreiben wir wieder
 $fv$ f\"ur das Bild von $v$ unter $f$, falls $f \in V^*$ und $v \in V$ ist.
 Sind nun $f_1$, \dots, $f_r \in V^*$, so ist
 $$ (f_1 \otimes \dots \otimes f_r)(v_1 \otimes \dots \otimes v_r)
       = f_1v_1 \otimes \dots \otimes f_rv_r
       = \prod_{i:=1}^r fv_i. $$
 Dies zeigt, dass $T_K^r(V^*) \subseteq T_K^r(V)^*$ ist. Ist $V$ endlichen
 Ranges, so ist $T_K^r(V^*) = T_K^r(V)^*$, da beide Vektorr\"aume den
 Rang $\Rg_K(V)^r$ haben.
 \medskip\noindent
 {\bf 2.9. Satz.} {\it Es sei $K$ ein kommutativer K\"orper. Ist $V$ ein
 Vektorraum \"uber $K$, so ist}
 $$ NT_K^r(V)^\bot = AT_K^r(V). $$
 \par
       Beweis. Es ist $NT_K^r(V)$ der von allen $v_1 \otimes \dots \otimes
 v_r$, bei denen wenigstens zwei der $v_i$ gleich sind, erzeugte Unterraum
 von $T_K^r(V)$. Nun ist $\sigma \in T_K^r(V)^*$ genau dann alternierend,
 wenn $NT_K^r(V) \subseteq \Kern(\sigma)$ gilt. Daher ist
 $AT_K^r(V) = NT_K^r(V)^\bot$, q. e. d.

 \medskip

       Es sei $f \in T_K^r(V)$. Wir definieren die Abbildung $g$ von
 $V^r$ in $K$ durch
 $$ g(v_1, \dots, v_r) := f(v_1 \otimes \dots \otimes v_r). $$
 Dann ist $g$ nat\"urlich eine multilineare Abbildung von $V$ in $K$. Daher
 gibt es f\"ur die Antisymmetrisierte $ga$ von $g$ ein $fa \in T_K^r(V)^*$
 mit
 $$ ga(v_1, \dots, v_r) = fa(v_1 \otimes \dots \otimes v_r) $$
 f\"ur alle $v_1$, \dots, $v_r \in V$. Dieses $fa$ nennen wir sinngem\"a\ss\
 {\it Antisymmetrisierte\/} von $f$. Es gilt
 $$\eqalign{
 fa(v_1 \otimes \dots \otimes v_r) = ga(v_1, \dots, v_r) 
       &= \sum_{\pi \in S_r} \sgn(\pi) g(v_{\pi(1)}, \dots g_{\pi(r)}) \cr
       &= \sum_{\pi \in S_r} \sgn(\pi) f(v_{\pi(1)} \otimes \dots
                                          \otimes v_{\pi(r)}). \cr} $$
 \par\noindent
 {\bf 2.10. Satz.} {\it Es sei $K$ ein kommutativer K\"orper und $V$ sei ein
 Vektorraum \"uber $K$. Ferner seien $y_1$, \dots, $y_r \in V^*$. Setze
 $$ f := y_1 \otimes \dots \otimes y_r. $$
 Dann ist}
 $$ fa = \sum_{\pi \in S_r} \sgn(\pi) (y_{\pi(1)} \otimes \dots \otimes
                            y_{\pi(r)}). $$
 \par
       Beweis. Nach der zuvor gemachten Bemerkung folgt mittels
 $f = y_1 \otimes \dots \otimes y_r$, dass
 $$\eqalign{
 fa(v_1 \otimes \dots \otimes v_r)
       &= \sum_{\pi \in S_r} \sgn(\pi) f(v_{\pi(1)} \otimes \dots
                                          \otimes v_{\pi(r)})       \cr
       &= \sum_{\pi \in S_r} \sgn(\pi) \prod_{i:=1}^r y_iv_{\pi(i)} 
        = \sum_{\pi \in S_r} \sgn(\pi) \prod_{i:=1}^r y_{\pi(i)}v_i \cr
       &= \sum_{\pi \in S_r} \sgn(\pi) (y_{\pi(1)} \otimes \dots \otimes
              y_{\pi(r)})(v_1 \otimes \dots \otimes v_r) \cr} $$
 ist. Hieraus folgt die Behauptung.
 \medskip\noindent
 {\bf 2.11. Satz.} {\it Es sei $V$ ein Vektorraum endlichen Ranges \"uber dem
 kommutativen K\"orper $K$. Es gibt dann genau ein Skalarprodukt $f$, sodass
 $(\bigwedge_K^r(V^*),\bigwedge_K^r(V),f)$ ein duales Raumpaar ist
 und \"uberdies
 $$ f(y_1 \wedge \dots \wedge y_r,x_1 \wedge \dots \wedge x_r)
       = \det(y_ix_j \mid i, j := 1, \dots r) $$
 gilt f\"ur alle $y_1$, \dots, $y_r \in V^*$ und alle $x_1$, \dots, $x_r \in
 V$.}
 \smallskip
       Beweis. Dass es h\"ochstens ein solches Skalarprodukt gibt, folgt
 daraus, dass $\bigwedge_K^r(V^*)$ von den $y_1 \wedge \dots \wedge y_r$
 und $\bigwedge_K^r(V)$ von den $x_1 \wedge \dots \wedge x_r$ erzeugt werden.
 \par
       Mit 2.8 folgt die Existenz eines Epimorphismus $\psi$ von
 $T_K^r(V^*)$ auf $\bigwedge_K^r(V^*)$ mit
 $$ \psi(y_1 \otimes \dots \otimes y_r) = y_1 \wedge \dots \wedge y_r. $$
 Es ist
 $$ \Kern(\psi) = NT_K^r(V^*). $$
 Da die Abbildung $(y_1, \dots, y_r) \to (y_1 \otimes \dots \otimes y_r)a$
 alternierend ist, ist
 $$ NT_K^r(V^*) \subseteq \Kern(a). $$
 Somit gibt es einen eindeutig bestimmten Homomorphismus $\eta$ von
 $\bigwedge_K^r(V^*)$ in $T_K^r(V^*)$ mit $a = \eta\psi$. Insbesondere ist
 $$ (y_1 \otimes \dots \otimes y_r)a = \eta(y_1 \wedge \dots \wedge y_r). $$
 Weil die $y_1 \otimes \dots \otimes y_r$ den Raum $T_K^r(V^*)$ erzeugen,
 erzeugen die $(y_1 \otimes \dots \otimes y_r)a$ den Raum $T_K^r(V^*)a$.
 Folglich ist $\eta$ ein Epimorphismus von $\bigwedge_K^r(V^*)$ auf
 $T_K^r(V^*)a$. Mittels 2.4 folgt $AT_K^r(V) = T_K^r(V^*)a$ und daher
 $$ \Rg({\textstyle\bigwedge}_K^r(V^*)) = {n \choose r} = \Rg(AT_K^r(V))
       = \Rg(T_K^r(V^*)a). $$
 Hieraus und aus der Surjektivit\"at von $\eta$ folgt, dass $\eta$ auch
 injektiv ist. Somit ist $\eta$ also ein Isomorphismus von $\bigwedge_K^r(V^*)$
 auf $T_K^r(V^*)a$.
 \par
       $T_K^r(V^*)a$ ist gerade der zu $NT_K^r(V)$ orthogonale Unterraum von
 $T_K^r(V^*)$. Ist $y \in T_K^r(V^*)a$ und $x \in T_K^r(V)$, so wird also durch
 die Vorschrift
 $$ g(y,x + NT_K^r(V)) := yx $$
 ein Skalarprodukt $g$ auf  $T_K^r(V^*)a \times T_K^r(V)/NT_K^r(V)$
 erkl\"art. Benutzt man den Isomorphismus $\varphi^{-1}$ von
 $\bigwedge_K^r(V)$ auf $T_K^R(V)/NT_K^r(V)$ und den eben bestimmten
 Isomorphismus $\eta$ von $\bigwedge_K^r(V^*)$ auf $T_K^r(V^*)a$, so wird
 $$ ({\textstyle\bigwedge}_K^r(V^*),{\textstyle\bigwedge}_K^r(V)) $$
 zu einem dualen Raumpaar bez\"uglich $f$, wenn man $f$ durch die Vor\-schrift
 $$ f(y,x) := g(\eta(y),\varphi^{-1}(x)) $$
 definiert. Ist $x = x_1 \wedge \dots \wedge x_r$ und $y = y_1 \wedge \dots
 \wedge y_r$, so folgt mit Satz 2.10
 $$\eqalign{
 f(y,x) &= g((y_1 \otimes \dots \otimes y_r)a,x_1 \otimes \dots \otimes x_r
                     + NT_K^r(V)) \cr
       &= (y_1 \otimes \dots \otimes y_r)a(x_1 \otimes \dots \otimes x_r) \cr
       &= \biggl(\sum_{\pi\in S_r} \sgn(\pi)(y_{\pi(1)} \otimes \dots
              \otimes y_{\pi(r)}\biggr)(x_1 \otimes \dots \otimes x_r)    \cr
       &= \sum_{\pi\in S_r} \sgn(\pi) \prod_{i:=1}^r y_{\pi(i)}x_i        \cr
       &= \det(y_jx_i\mid j,i := 1, \dots, r). \cr} $$
 Damit ist alles bewiesen.
 \medskip\noindent
 {\bf 2.12. Satz.} {\it Es sei $K$ ein kommutativer K\"orper und $V$ sei ein
 Vektorraum endlichen Ranges \"uber $K$. Ferner sei $f_r$ das in 2.11
 beschriebene Skalarprodukt auf $\bigwedge_K^r(V^*) \times \bigwedge_K^r(V)$.
 Ist dann $x = \sum_{i:=0}^n x_i \in \bigwedge_K(V)$ mit
 $x_i \in \bigwedge_K^i(V)$ und
 $y = \sum_{i:=0}^n y_i \in \bigwedge_K(V^*)$ mit $y_i \in \bigwedge_K^i(V^*)$
 und setzt man
 $$ f(y,x) := \sum_{i:=0}^n f_i(y_i,x_i), $$
 so ist $f$ ein Skalarprodukt auf $\bigwedge_K(V^*) \times \bigwedge_K(V)$.
 \par
       Ist $b_1$, \dots, $b_n$ eine Basis von $V$ und $b_1^*$, \dots, $b_n^*$
 die zu dieser Basis duale Basis von $V^*$, so ist
 $$ \{b_I^* \mid I \in \Fin(\{1, \dots, n\})\} $$
 die zu
 $$ \{b_I^* \mid I \in \Fin(\{1, \dots, n\})\} $$
 duale Basis von $\bigwedge_K(V^*)$.}
 \smallskip
       Beweis. Es ist banal, dass $f$ ein Skalarprodukt ist.
 \par
       Es seien $I$, $J \in \Fin(\{1, \dots, n\})$. Ist $|I| \neq |J|$, so
 ist $f(b_I^*,b_J) = 0$. Es sei also $|I| = |J|$. Dann ist
 $$ f(b_I^*,b_J) = \det(b_i^*b_j \mid i \in I, j \in J). $$
 Ist nun $I \neq J$, so gibt es wegen $|I| = |J|$ ein $k \in I - J$. Es
 folgt $b_k^*b_j = 0$ f\"ur alle $j \in J$ und damit $f(b^*_kb_j) = 0$.
 Ist $I = J$, so ist
 $$ b_i^*b_j = \cases{1, &falls $i = j$     \cr
                      0, &falls $i \neq j$  \cr} $$
 und daher
 $$ \det(b_i^*b_j \mid i, j \in I) = 1. $$
 Damit ist der Satz bewiesen.
 \medskip
       Zum Schluss noch eine \"Ubungsaufgabe f\"ur den Leser.
 \medskip\noindent
 {\bf 2.13. Satz.} {\it Es sei $K$ ein kommutativer K\"orper und $V$ sei ein
 Vektorraum endlichen Ranges \"uber $K$. Ist dann $(\bigwedge_K^r(V^*),
 \bigwedge_K^r(V),f)$ das in 2.11 be\-schrie\-be\-ne duale Raumpaar, so gilt
 $(\sigma^*)_{\#r} = (\sigma_{\#r})^*$ f\"ur alle $\sigma \in End_K(V)$.}

 \mysection{3. Innere Produkte}

\noindent
       Es sei $V$ ein Vektorraum endlichen Ranges \"uber dem kommutativen
 K\"orper $K$ und $V^*$ bezeichne wie \"ublich seinen Dualraum. Wir betrachten
 auf $\bigwedge_K(V^*) \times \bigwedge_K(V)$ das in Satz 2.12
 definierte Skalarprodukt $f$. Es sei $y \in \bigwedge_K(V^*)$. Die durch
 $\varphi_y^*(z) := z \wedge y$ definierte Abbildung $\varphi^*$ ist eine
 Endomorphismus von $\bigwedge_K(V^*)$. Die zu $\varphi_y^*$ duale Abbildung
 von $\bigwedge_K(V)$ in sich bezeichnen wir mit $\varphi_y$. Wir definieren
 $\rhak$ durch
 $$ y \rhak x := \varphi_y(x) $$
 f\"ur alle $y \in \bigwedge_K(V^*)$ und alle $x \in \bigwedge_K(V)$.
 Es gilt dann
 $$ f(z,y \rhak x) = f(z,\varphi_y(x)) = f(\varphi_y^*(z),x)
       = f(z \wedge y,x) $$
 f\"ur alle $y$, $z \in \bigwedge_K(V^*)$ und alle $x \in \bigwedge_K(V)$.
 Wir nennen $y \rhak x$ {\it rechtsseitiges inneres Produkt\/} von $y$ mit
 $x$. Entsprechend ist die Abbildung  $y \to y \wedge x$ ein Endomorphismus
 von $\bigwedge_K(V)$. Es gibt also einen Endomorphismus $\lhak$ von
 $\bigwedge_K(V^*)$ in sich mit
 $$ f(z \lhak y,x) = f(z,y \wedge x) $$
 f\"ur alle $z \in \bigwedge_K(V^*)$ und alle $y$, $x \in \bigwedge_K(V)$.
 Wir nennen $z \lhak y$ {\it linksseitiges inneres Produkt\/} von $z$ mit $y$.
 \par
       Wir notieren einige Rechenregeln.
 \medskip\noindent
 {\bf 3.1. Satz.} {\it Es sei $V$ ein Vektorraum endlichen Ranges \"uber dem
 kommutativen K\"orper $K$. Dann gilt:
 \item{a)} Es ist $z \rhak (y + x) = (z \rhak y) + (z \rhak x)$ f\"ur alle
 $z \in \bigwedge_K(V^*)$ und alle $y$, $x \in \bigwedge_K(V)$.
 \item{b)} Es ist $(z + y) \rhak x = (z \rhak x) + (y \rhak x)$ f\"ur alle
 $z$, $y \in \bigwedge_K(V^*)$ und alle $x \in \bigwedge_K(V)$.
 \item{c)} Es ist $(yk) \rhak x = y \rhak (xk) = (y \rhak x)k$ f\"ur alle
 $y \in \bigwedge_K(V^*)$ und alle $x \in \bigwedge_K(V)$.
 \item{d)} Es ist $(z \wedge y) \rhak x = z \rhak (y \rhak x)$ f\"ur alle
 $z$, $y \in \bigwedge_K(V^*)$ und alle $x \in \bigwedge_K(V)$.
 \par\noindent
       Entsprechend gilt:
 \item{a$'$)} Es ist $z \lhak (y + x) = (z \lhak y) + (z \lhak x)$ f\"ur alle
 $z \in \bigwedge_K(V^*)$ und alle $y$, $x \in \bigwedge_K(V)$.
 \item{b$'$)} Es ist $(z + y) \lhak x = (z \lhak x) + (y \lhak x)$ f\"ur alle
 $z$, $y \in \bigwedge_K(V^*)$ und alle $x \in \bigwedge_K(V)$.
 \item{c$'$)} Es ist $(yk) \lhak x = y \lhak (xk) = (y \lhak x)k$ f\"ur alle
 $y \in \bigwedge_K(V^*)$ und alle $x \in \bigwedge_K(V)$.
 \item{d$'$)} Es ist $z \lhak (y \wedge x) = z \lhak (y \lhak x)$ f\"ur alle
 $z$, $y \in \bigwedge_K(V^*)$ und alle $x \in \bigwedge_K(V)$.\par}
 \smallskip
       Beweis. a) Es sei $u \in \bigwedge_K(V^*)$. Dann ist
 $$\eqalign{
 f(u,z \rhak (x + y)) &= f(u \wedge z,x + y) \cr
       &= f(u \wedge z,x) + f(u \wedge z,y)  \cr
       &= f(u,z \rhak x) + f(u,z \rhak y)    \cr
       &= f(u, (z \rhak x) + (z \rhak y)).   \cr} $$
 Da $u$ beliebig war, gilt
 $$ z \rhak (x + y) = (z \rhak x) + (z \rhak y). $$
 \par
       d$'$) Es sei $u \in \bigwedge_K(V)$. Dann ist
 $$\eqalign{
 f( z \lhak (y \wedge x),u) &= f(z,(y \wedge x) \wedge u)   \cr
       &= f(z,y \wedge (x \wedge u))                        \cr
       &= f(z \lhak y, x \wedge u)                          \cr
       &= f((z \lhak y) \lhak x,u).                         \cr} $$
 Die restlichen Aussagen beweisen sich analog.
 \medskip
       Die Elemente aus $\bigwedge_K^r(V)$ nennen wir auch $r$-{\it Vektoren\/}
 und die Elemente aus $\bigwedge_K^r(V^*)$ nennen wir $r$-{\it Formen\/}.
 Ist $x$ ein $r$-Vektor und $y$ ein $s$-Vektor, so ist $x \wedge y = 0$,
 falls $r + s > \Rg_K(V)$ ist. Ist $r + s \leq \Rg_K(V)$, so folgt mit
 1.8, dass $x \wedge y$ ein $(r + s)$-Vektor ist. Entsprechendes gilt
 f\"ur $r$- und $s$-Formen.
 \medskip\noindent
 {\bf 3.2. Satz.} {\it Es sei $V$ ein Vektorraum endlichen Ranges \"uber dem
 kommutativen K\"orper $K$. Ist $x$ ein $r$-Vektor und $y$ eine $s$-Form,
 so gilt:
 \item{a)} F\"ur $r < s$ ist $y \rhak x = 0$. Ist $r \geq s$, so ist
 $y \rhak x$ ein $(r - s)$-Vektor.
 \item{b)} F\"ur $r > s$ ist $y \lhak x = 0$. Ist $r \leq s$, so ist
 $y \lhak x$ eine $(s - r)$-Form.
 \par\noindent
       Ist $r = s$, so ist $y \rhak x = f(y,x) = y \lhak x$. Dabei ist $f$
 wieder das oben definierte Skalarprodukt auf $\bigwedge_K(V^*) \times
 \bigwedge_K(V)$.}
 \smallskip
       Beweis. a) Es sei $z$ eine $t$-Form. Dann ist
 $$ f(z,y \rhak x) = f(z \wedge y,x). $$
 Ist $r \neq s + t$, was insbesondere f\"ur alle $t$ der Fall ist, falls
 $r < s$ gilt, so ist $f(z \wedge y,x) = 0$ aufgrund der Definition von $f$,
 da $z \wedge y$ ja eine $s + t$-Form ist. Ist $r < s$ so ist also
 $y \rhak x = 0$, da ja dann $f(z,y \rhak x) = 0$ f\"ur alle Formen $z$ gilt.
 Ist $r \geq s$, so ist $f(z,y \rhak x)$ h\"ochstens dann von Null
 verschieden, wenn $r = s + t$ ist. Dies impliziert, dass 
 $$ y \rhak x \in \biggl(\bigoplus_{t:=0, t\neq r - s}^n {\textstyle
   \bigwedge}_K^t(V^*)\biggr)^\top = {\textstyle\bigwedge}_K^{s-r}(V) $$
 ist.
 \par
       b) Es sei $u$ ein $t$-Vektor. Dann ist
 $$ f(y \lhak x,u) = f(y,x \wedge u). $$
 Hieraus erschlie\ss t man b) analog zu a).
 \par
       Ist $r = s$, so ist $y \rhak x \in \bigwedge_K^0(V) = K$. Also ist,
 da $f(1,1) = 1$ ist,
 $$ y \rhak x = (y \rhak x)f(1,1) = f(y \rhak x,1) = f(y,x \wedge 1)
        = f(y,x) $$
 und
 $$ y \lhak x = f(1,1)(y \lhak x) = f(1,y \lhak x) = f(1 \wedge y,x)
       = f(y,x). $$
 Damit ist alles bewiesen.
 \medskip\noindent
 {\bf 3.3. Satz.} {\it Es sei $V$ ein Vektorraum des Ranges $n$ \"uber dem
 kommutativen K\"orper $K$. Ferner sei $z$ eine $n$-Form ungleich der Nullform.
 Definiert man $\varphi$ durch
 $$ \varphi(x) := z \lhak x, $$
 so ist $\varphi$ ein Isomorphismus von $\bigwedge_K(V)$ auf
 $\bigwedge_K(V^*)$. Ist $\varphi_r$ die Einschr\"ankung von $\varphi$ auf
 $\bigwedge_K^r(V)$, so ist $\varphi_r$ ein Isomorphismus von
 $\bigwedge_K^r(V)$ auf $\bigwedge_K^{n-r}(V^*)$.
 \par
       Ist $y$ ein $n$-Vektor mit $f(z,y) = 1$, und ein solcher existiert
 stets, so ist 
 $$ \varphi^{-1}(u) = u \rhak y $$
 f\"ur alle $u \in \bigwedge_K(V^*)$.}
 \smallskip
       Beweis. Aus 3.1 a$'$) und c$'$) folgt, dass $\varphi$ eine
 lineare Abbildung von $\bigwedge_K(V)$ in $\bigwedge_K(V^*)$ ist.
       Es sei $b_1$, \dots, $b_n$ eine Basis von $V$. Dann ist
 $$ 0 \neq b_1 \wedge \dots \wedge b_n \in {\textstyle \bigwedge}_K^n(V). $$
 Wegen $\Rg(\bigwedge_K^n(V)) = 1$ ist daher
 $$ {\textstyle\bigwedge}_K^n(V) = (b_1 \wedge \dots \wedge b_n)K, $$
 so dass insbesondere alle $n$-Vektoren zerlegbar sind. Da dies auch f\"ur
 alle $n$-Formen gilt, ist
 $$ z = (b_1^*k) \wedge b_2^* \wedge \dots \wedge b_n^*, $$
 wobei $b_1^*$, \dots, $b_n^*$ die zu $b_1$, \dots, $b_n$ duale Basis ist.
 Indem wir $b_1$ durch $b_1k^{-1}$ ersetzen, erreichen wir, dass
 $$ z = b_1^* \wedge b_2^* \wedge \dots \wedge b_n^* $$
 ist. Setze
 $$ y := b_1 \wedge \dots \wedge b_n. $$
 Dann ist $f(z,y) = 1$.
 \par
       Es seien $L$, $J \subseteq \{1, \dots, n\}$. Dann ist
 $$
 f(\varphi(b_J),b_L) = f(z \lhak b_J,b_L)
        = f(z,b_J \wedge b_L)            
       = \prod_{i \in L, j \in J} \langle j,i \rangle
                        f(b_{\{1, \dots, b_n\}}^*,b_{J\cup L}). 
 $$
 Ist $I$ das Komplement von $J$ in $\{1, \dots, n\}$, so ist
 $$ f(b_{\{1, \dots, b_n\}}^*, b_{J \cup L}) = f(b_I^*,b_L). $$
 Beachtet man noch, dass $f(b_I^*,b_L) = 0$ ist f\"ur $I \neq L$, so sieht
 man die G\"ultigkeit von
 $$ f(\varphi(b_J),b_L) = \prod_{i \in I,j \in J} \langle j,i \rangle
                            f(b_I^*,b_L). $$
 Definiere die Abbildung $\xi$ von $\bigwedge_K(V)$ in $\bigwedge_K(V^*)$
 durch 
 $$ \xi(b_J) := \prod_{i \in I,j \in J} \langle j,i \rangle b_I^*, $$
 wobei $I$ das Komplement von $J$ in $\{1, \dots, n\}$ ist. Dann ist
 $$ f(\xi(b_J),b_L) = \prod_{i \in I,j \in J} \langle j,i \rangle
                            f(b_I^*,b_L). $$
 Also ist $f(\varphi(b_J),b_L) = f(\xi(b_J), b_L)$ f\"ur alle Teilmengen $L$
 und $J$ von $\{1, \dots, n\}$. Weil die $b_L$ eine Basis von $\bigwedge_K^r
 (V)$ bilden, folgt $\varphi(b_J) = \xi(b_J)$ und damit schlie\ss lich
 $\varphi = \xi$. Also ist
 $$ \varphi(b_J) = \prod_{i \in I,j \in J} \langle j,i \rangle b_I^*, $$
 wenn nur $I$ das Komplement von $J$ in $\{1, \dots, n\}$ ist.
 \par
       Definiere die Abbildung $\psi$ von $\bigwedge_K(V^*)$ in
 $\bigwedge_K(V)$ durch
 $$ \psi(u) := u \rhak y $$
 f\"ur alle $u \in \bigwedge_K(V^*)$. Vertauscht man in der vorstehenden
 Argumentation die Rollen von $V$ und $V^*$, was wegen der Endlichkeit des
 Ranges m\"oglich ist, so sieht man, dass $\psi$ eine lineare Abbildung ist,
 f\"ur die
 $$ \psi(b_I^*) = \prod_{i \in I,j \in J} \langle j,i \rangle b_J $$
 gilt, falls nur $J$ das Komplement von $I$ in $\{1, \dots, n\}$ ist.
 (Wer glaubt, es m\"usse in der Formel $\langle i,j \rangle$ hei\ss en ---
 ich war mir einen Augenblick auch unsicher ---, der rechne!)
 Aus all dem folgt schlie\ss lich
 $$ \psi\varphi(b_J) = \biggl(\prod_{i \in I, j\in J} \langle i,j \rangle
       \biggr)^2 b_J = b_J $$
 und
 $$ \varphi\psi(b_I^*) = \biggl(\prod_{i \in I, j\in J} \langle i,j \rangle
       \biggr)^2 b_I^* = b_I^*. $$
 Somit ist $\varphi$ bijektiv und $\psi = \varphi^{-1}$.
 \par
       Weil $\varphi$ ein Isomorphismus ist, folgt mit 3.2 weiter, dass
 $\varphi_r$ ein Isomorphismus von $\bigwedge_K^r(V)$ auf
 $\bigwedge_K^{n-r}(V^*)$ ist.
 \medskip\noindent
 {\bf 3.4. Korollar.} {\it Die Voraussetzungen seien wie bei Satz 3.3. Ferner
 sei $b_1, \dots, b_n$ eine Basis von $V$ und $b_I$ mit $I \in \Fin(\{1, \dots,
 n\})$ bzw. $b_J^*$ mit $J \in \Fin(\{1, \dots, n\})$ die von der gegebenen
 Basis abgeleiteten Basen von
 $\bigwedge_K(V)$ bzw. $\bigwedge_K(V^*)$. Ist dann $I$ eine Teilmenge
 von $\{1, \dots, n\}$ und ist $J$ ihr Komplement, so ist
 $$ \varphi(b_J) = \prod_{i \in I, j \in J} \langle j,i \rangle b_I^* $$
 und
 $$ \varphi^{-1}(b_I^*) =
       \prod_{i\in I, j \in J} \langle j,i \rangle b_J^*. $$}
 \par
       Ist $X$ ein Vektorraum \"uber $K$ und ist $k \in K$, so bezeichnen wir
 mit $\mu_k$ die durch $\mu_k(x) := xk$ definierte lineare Abbildung von
 $X$ in sich.
 \medskip\noindent
 {\bf 3.5. Satz.} {\it Es sei $V$ ein Vektorraum des Ranges $n$ \"uber $K$.
 Ist $\sigma \in GL(V)$, so gibt es $k \in K^*$ mit
 $$ \varphi\sigma_{\#r}\mu_k = \sigma_{\#(n-r)}^{*-1}\varphi. $$
 f\"ur alle $r$ mit $0 \leq r \leq n$. Dabei ist $\varphi$ wieder der in
 Satz 3.3 definierte Isomorphismus von $\bigwedge_K(V)$ auf $\bigwedge_K(V^*)$.}
 \smallskip
       Beweis. Es sei $y$ der bei der Definition von $\varphi$ benutzte
 $n$-Vektor. Ferner sei $v_1$, \dots, $v_n$ eine Basis von $V$ mit
 $y = v_1 \wedge \dots \wedge v_n$ und $v_1^*$, \dots, $v_n^*$ sei ihre duale
 Basis. Weil $\sigma$ ein Automorphismus ist, ist auch $\sigma(v_1)$, \dots,
 $\sigma(v_n)$ eine Basis von $V$. Es folgt
 $$ 0 \neq \sigma(v_1) \wedge \dots \wedge \sigma(v_n) \in \bigwedge_K^n(V). $$
 Es gibt folglich ein $k \in K^*$ mit
 $$ y = \sigma(v_1k) \wedge \sigma(v_2) \wedge \dots \wedge \sigma(v_n). $$
 \par
       Die zu $v_1k$, $v_2$, \dots, $v_n$ duale Basis ist
 $v_1^*k^{-1}$, $v_2^*$, \dots, $v_n^*$. Setze
 $b_1 := v_1k$ und $b_i := v_i$ f\"ur $i > 1$ sowie
 $b_1^* := v_1^*k^{-1}$ und $b_i^* := v_i^*$ f\"ur $i > 1$.
 Dann ist also $b_1^*$, \dots, $b_n^*$ die zu $b_1$, \dots, $b_n$ duale Basis.
 Ferner ist
 $$ f(\sigma^{*-1}(b_J^*),\sigma(b_i)) = f(b_j^*,\sigma^{-1}\sigma(b_i))
       = f(b_j^*,b_i). $$
 Also ist $\sigma^{*-1}(b_1^*)$, \dots, $\sigma^{*-1}(b_n^*)$ die zu
 $\sigma(b_1)$, \dots, $\sigma(b_n)$ duale Basis.
 Ist nun $i$ eine streng monoton steigende Abbildung von $\{1, \dots, r\}$
 in $\{1, \dots, n\}$ und $j$ eine streng monoton steigende Abbildung von
 $\{1, \dots, n - r\}$ auf $\{1, \dots, n\} - \{i_{\{1, \dots, r\}}\}$,
 so folgt
 mit 3.4
 $$\eqalign{
 \varphi\sigma_{\#r}(b_{i_1} \wedge \dots \wedge b_{i_r})
       &= \varphi(\sigma(b_{i_1}) \wedge \dots \wedge \sigma(b_{i_r})) \cr
       &= \prod_{l:=1}^{n-r}\prod_{m:=1}^r \langle j_l,i_m \rangle
  (\sigma^{*-1}(b_{j_1}^*) \wedge \dots \wedge \sigma^{*-1}(b_{j_{n-r}}^*)) \cr
       &= \prod_{l:=1}^{n-r}\prod_{m:=1}^r \langle j_l,i_m \rangle
  \sigma_{\#(n-r)}^{*-1}(b_{j_1}^* \wedge \dots \wedge b_{j_{n-r}}^*). \cr} $$
 \hfuzz=2pt%
 Wegen der Monotonie von $i$ und $j$ ist entweder $i_1 = 1$ oder $j_1 = 1$.
 Diese beiden F\"alle sind nun zu unterschieden.
 \par
       Es sei $i_1 = 1$. Dann ist also
 $$\eqalign{
 \varphi\sigma_{\#r}\mu_k(v_1 \wedge v_{i_2} \wedge \dots \wedge v_{i_r}) 
    &= \prod_{l:=1}^{n-r}\prod_{m:=1}^r\langle j_l,i_m \rangle
       \sigma^{*-1}_{\#(n-r)}(v_{j_1}^* \wedge \dots \wedge v_{j_{n-r}}^*) \cr
    &= \sigma_{\#(n-r)}^{*-1}(v_1 \wedge v_{i_2} \wedge \dots \wedge v_{i_r}).
                                     \cr} $$
 Ist $j_1 = 1$, so ist
 $$\eqalign{
 \varphi\sigma_{\#r}(v_{i_1} \wedge \dots \wedge v_{i_r}) 
       &= \prod_{l:=1}^{n-r}\prod_{m:=1}^r\langle j_l,i_m \rangle
          \sigma_{\#(n-r)}^{*-1}\mu_k^{-1}
          (v_1^* \wedge v_{j_2}^* \wedge \dots \wedge v_{j_{n-r}}) \cr
       &= \sigma_{\#(n-r)}^{*-1}\mu_k^{-1}\varphi
              (v_{i_1} \wedge \dots \wedge v_{i_r}) \cr} $$
 Da die $v_{i_1} \wedge \dots \wedge v_{i_r}$ eine Basis von
 $\bigwedge_K(V)$ bilden und da $\mu_k$ mit $\varphi$ vertauschbar ist,
 gilt also
 $$ \varphi\sigma_{\#r}\mu_k = \sigma_{\#(n-r)}^{*-1}\varphi. $$

 \mysection{4. Zerlegbare Vektoren}

\noindent
       Zerlegbare Vektoren sind f\"ur die Geometrie von besonderer Bedeutung,
 da sich mit ihrer Hilfe die Unterr\"aume von Vektorr\"aumen dar\-stel\-len
 lassen. Daher werden wir nun einige Aussagen \"uber sie beweisen. Einige der
 S\"atze gelten auch in allgemeineren Situationen. Wir be\-schr\"an\-ken uns aber
 im Folgenden auf Vektorr\"aume, da man die Beweise der fraglichen S\"atze
 abk\"urzen kann, indem man Eigenschaften benutzt, die nur Vek\-tor\-r\"au\-me
 besitzen.
 \medskip\noindent
 {\bf 4.1. Satz.} {\it Es sei $K$ ein kommutativer K\"orper und $V$ sei ein
 Vektorraum \"uber $K$. Ist $0 \neq z \in \bigwedge_K^r(V)$, so setzen wir
 $$ V_z := \{x \mid x \in V, z \wedge x = 0\}. $$
 Dann ist $V_z \in \La(V)$. Sind $x_1$, \dots, $x_s$ linear unabh\"angige
 Vektoren aus $V_z$, so ist $s \leq r$ und es gibt ein $w \in
 \bigwedge_K^{r-s}(V)$ mit
 $$ z = w \wedge x_1 \wedge \dots \wedge x_s. $$
 Insbesondere ist $\Rg_K(V_z) \leq r$.}
 \smallskip
       Beweis. Es ist klar, dass $V_z$ ein Unterraum von $V$ ist.
       Es gibt eine Basis $B$ von $V$ mit $x_1$, \dots, $x_s \in B$. Es
 gibt ferner eine lineare Ordnung von $B$ mit
 $$ x_1 < x_2 < \dots < x_s. $$
 Nach 1.7 ist
 $$ z = \sum_{I \in \Fin_r(B)} b_Ik_I. $$
 Wegen $x_i \in V_z$ ist
 $$ 0 = z \wedge x_i = \sum_{I \in \Fin_r(B)} (b_I \wedge x_i)k_I. $$
 Mit 1.7 folgt hieraus
 $$ (b_I \wedge x_i)k_I = 0 $$
 f\"ur alle $I \in \Fin_r(B)$. Ist
 $x_i \not\in I$, so ist $b_I \wedge x_i \neq 0$ und daher $K_I = 0$. Ist
 $k_I \neq 0$, so ist also $x_i \in I$ f\"ur $i := 1$, \dots, $r$. Weil $z$
 nicht Null ist, gibt es ein $I$ mit
 $|I| = r$ und $k_I \neq 0$. F\"ur dieses $I$ gilt dann also $x_i \in I$ f\"ur
 $i := 1$, \dots, $s$. Also ist $s \leq |I| = r$.
 \par
       Ist $k_I \neq 0$, so ist $x_1$, \dots, $x_s \in I$, wie wir gerade
 gesehen haben. Es gibt dann also ein $w_I \in \bigwedge_K^{r-s}(V)$ mit
 $$ b_I = w_I \wedge x_1 \wedge \dots \wedge x_s. $$
 Ist $K_I = 0$, so ersetzen wir $b_I$ und $w_I$ jeweils durch 0 und nennen
 diese Elemente wiederum $b_I$ und $w_I$. Dann ist also in jedem Falle
 $$ b_I = w_I \wedge x_1 \wedge \dots \wedge x_s. $$
 Hiermit folgt
 $$\eqalign{
 z = \sum_{I \in \Fin_r(B)} b_Ik_I 
   &= \sum_{I \in \Fin_r(B)} (w_I \wedge x_1 \wedge \dots \wedge x_s)k_I \cr
   &= \biggl(\sum_{I \in \Fin_r(B)}
               w_Ik_I\biggr) \wedge x_1 \wedge \dots \wedge x_s.\cr} $$
 Damit ist alles bewiesen.
 \medskip
       Ist $z \in \bigwedge_K^r(V)$, so hei\ss t $z$ genau dann {\it
 zerlegbar\/}, wenn es $x_1$, \dots, $x_r \in V$ gibt mit
 $$ z = x_1 \wedge \dots \wedge x_r. $$
 \medskip\noindent
 {\bf 4.2. Korollar.} {\it Es sei $K$ ein kommutativer K\"orper und $V$ 
 ein $K$-Vektorraum. Ist $0 \neq z \in \bigwedge_K^r(V)$, so ist $z$ genau
 dann zerlegbar, wenn $\Rg_K(V_z) = r$ ist. Ist $z = x_1 \wedge \dots \wedge
 x_r$, so ist
 $$ V_z = \sum_{i:=1}^r x_iK. $$}
 \par
       Beweis. Es sei $z = x_1 \wedge \dots \wedge x_r$ mit $x_i \in V$.
 Dann ist $z \wedge x_i = 0$ und somit $x_i \in V_z$ f\"ur alle $i$. Weil
 $x_1$, \dots, $x_r$ nach 1.10 linear unabh\"angig sind, ist
 $r \leq \Rg_K(V_z)$. Nach 1.10 ist andererseits $\Rg_K(V_z) \leq r$,
 so dass $\Rg_K(V_z) = r$ ist. Hieraus folgt wiederum
 $ V_z = \sum_{i:=1}^r x_iK$.
 \par
       Ist umgekehrt $\Rg_K(V_z) = r$ und ist $x_1$, \dots, $x_r$ eine
 Basis von $V_z$, so gibt es nach 4.1 ein $k \in \bigwedge_K^{r-r}(V) = K$
 mit
 $ z = k \wedge x_1 \wedge \dots \wedge x_r = (kx_1) \wedge \dots \wedge x_r
       $.
 Damit ist das Korollar bewiesen.
 \medskip\noindent
 {\bf 4.3. Satz.} {\it Es sei $K$ ein kommutativer K\"orper und $V$ sei ein
 $K$-Vektorraum. Sind $x_1$, \dots, $x_r$ und $y_1$, \dots, $y_r$ je $r$
 linear unabh\"angige Vektoren aus $V$, so ist genau dann
 $$ \sum_{i:=1}^r x_iK = \sum_{i:=1}^r y_iK, $$
 wenn
 $ (x_1 \wedge \dots \wedge x_r)K = (y_1 \wedge \dots \wedge y_r)K $
 ist.}
 \smallskip
       Beweis. Es sei $(x_1 \wedge \dots \wedge x_r)K =
 (y_1 \wedge \dots \wedge y_r)K$. Es gibt dann ein $k \in K$ mit
 $$ y_1 \wedge \dots \wedge y_r = (x_1 \wedge \dots \wedge x_r)k. $$
 Es folgt
 $$ y_1 \wedge \dots \wedge y_r \wedge x_i = 0 $$
 und damit, wenn man noch 4.2 ber\"ucksichtigt,
 $$ x_i \in V_{y_1 \wedge \dots \wedge y_r} = \sum_{j:=1}^r y_jK $$
 f\"ur alle $i$. Also gilt
 $$ \sum_{i:=1}^r x_iK \leq \sum_{i:=1}^r y_iK. $$
 Aus Symmetriegr\"unden gilt auch die umgekehrte Inklusion, so dass in der
 Tat
 $$ \sum_{i:=1}^r x_iK = \sum_{i:=1}^r y_iK $$
 gilt.
 \par
       Es sei nun $\sum_{i:=1}^r x_iK = \sum_{i:=1}^r y_iK$. Nach 4.2
 ist
 $$ \sum_{i:=1}^r y_iK = V_{y_1 \wedge \dots \wedge y_r}. $$
 Nach 1.10 gibt es daher ein $k \in \bigwedge_K^{r-r}(V)$ mit
 $$ y_1 \wedge \dots \wedge y_r = k \wedge x_1 \wedge \dots \wedge x_r
       = (x_1 \wedge \dots \wedge x_r)k. $$
 Es folgt
 $$ (y_1 \wedge \dots \wedge y_r)K \leq (x_1 \wedge \dots \wedge x_r)K. $$
 Aus Symmetriegr\"unden ist dann
 $$ (y_1 \wedge \dots \wedge y_r)K = (x_1 \wedge \dots \wedge x_r)K. $$
 Damit ist alles bewiesen.
 \medskip\noindent
 {\bf 4.4. Satz.} {\it Es sei $K$ ein kommutativer K\"orper und $V$ sei ein
 Vektorraum \"uber $K$. Ferner sei $0 \neq x \in \bigwedge_K^r(V)$ und $0 \neq
 y \in \bigwedge_K^s(V)$ und $x$ und $y$ seien beide zerlegbar. Genau dann ist
 $V_x \leq V_y$, wenn $r \leq s$ ist und es ein $z \in \bigwedge_K^{s-r}(V)$
 gibt mit $y = z\wedge x$.}
 \smallskip
       Beweis. Es sei $x = x_1 \wedge \dots \wedge x_r$. Dann ist
 $ V_x = \sum_{i:=1}^r x_iK$.
 \"Uberdies sind $x_1$, \dots, $x_r$ nach 1.10 linear unabh\"angig.
 \par
       Es sei nun $V_x \leq V_y$. Dann sind also $x_1$, \dots, $x_r$ linear
 unabh\"angige Vektoren aus $V_y$. Nach 3.10 ist daher $r \leq s$ und es
 gibt ein $z \in \bigwedge_K^{s-r}(V)$ mit $y = z \wedge x$.
 \par
       Es sei umgekehrt $y = w \wedge x$ mit einem
 $w \in \bigwedge_K^{s-r}(V)$. Dann ist
 $$ y \wedge x_i = w \wedge x \wedge x_i = 0 $$
 und folglich
 $$ V_x = \sum_{i:=1}^r x_iK \leq V_y. $$
 Damit ist 4.4 bewiesen.
 \medskip\noindent
 {\bf 4.5. Satz.} {\it Es sei $K$ ein kommutativer K\"orper und $V$ sei ein
 Vektorraum \"uber $K$. Ferner seien $x$ und $y$ zerlegbare Vektoren ungleich
 Null aus $\bigwedge_K(V)$. Genau dann ist $V_x \cap V_y \neq \{0\}$, wenn
 $x \wedge y = 0$ ist.}
 \smallskip
       Beweis. Es sei $V_x \cap V_y \neq \{0\}$. Ist $u_1$, \dots, $u_r$ eine
 Basis von $V_x \cap V_y$ und setzt man $z := u_1 \wedge \dots \wedge u_r$,
 so folgt mit 4.2, dass $V_x \cap V_y = V_z$ ist. Nach 4.1 ist
 $y = z \wedge u$ und $x = w \wedge z$ mit gewissen $u$ und $w$. Es folgt
 $$ x \wedge y = w \wedge z \wedge z \wedge u = 0, $$
 da wegen der Zerlegbarkeit von $z$ ja $z \wedge z = 0$ ist.
 \par
       Es sei $V_x \cap V_y = \{0\}$ und $x = x_1 \wedge \dots \wedge x_r$
 sowie $y = y_1 \wedge \dots \wedge y_s$. Dann ist $x_1$, \dots, $x_r$ eine
 Basis von $V_x$ und $y_1$, \dots, $y_s$ eine Basis von $V_y$. Wegen
 $V_x \cap V_y = \{0\}$ ist $x_1$, \dots, $x_r$, $y_1$, \dots, $y_s$ eine
 Basis von $V_x + V_y$. Also ist
 $$ x \wedge y = x_1 \wedge \dots \wedge x_r \wedge y_1 \wedge \dots \wedge y_s
        \neq 0. $$
 Damit ist alles bewiesen.
 \medskip
       Es sei $z$ eine $n$-Form ungleich der Nullform. Wir definieren wie in
 Satz 3.3 den Isomorphismus $\varphi$ von $\bigwedge_K(V)$ auf
 $\bigwedge_K(V^*)$ durch $\varphi(x) := z \lhak x$.
 \medskip\noindent
 {\bf 4.6. Satz.} {\it Es sei $V$ ein Vektorraum endlichen Ranges \"uber dem
 kommutativen K\"orper $K$. Ist $x \in \bigwedge_K(V)$, so ist $x$ genau
 dann zerlegbar, wenn $\varphi(x)$ zerlegbar ist.}
 \smallskip
       Beweis. Ist $x = 0$, so ist $\varphi(x) = 0$ und $x$ und $\varphi(x)$
 sind zerlegbar. Es sei also $x \neq 0$. Ferner sei $x = x_1 \wedge \dots
 \wedge x_r$ mit $x_i \in V$ f\"ur alle $i$. Nach 1.10 sind die $x_i$ linear
 unabh\"angig, da ja $x \neq 0$ ist. Sie k\"onnen also zu einer Basis
 $x_1$, \dots, $x_r$, \dots, $x_n$ von $V$ erg\"anzt werden. Es sei $x_1^*$,
 \dots, $x_n^*$ die zu $x_1$, \dots, $x_n$ duale Basis. Ist $r = n$, so
 ist $\varphi(x) \in \bigwedge_K^0(V)$, so dass auch $\varphi(x)$ zerlegbar
 ist. Es sei also $r < n$. Indem wir $x_n$ gegebenenfalls durch einen
 Skalar ab\"andern, k\"onnen wir erreichen, dass
 $$ z = x_1^* \wedge \dots \wedge x_n^* $$
 ist, wobei $z$ die zur Konstruktion von $\varphi$ benutzte $n$-Form ist.
 Nach 3.4 ist daher
 $$ \varphi(x) = (-1)^{(n-r)r}x_{r+1}^* \wedge \dots \wedge x_n^*, $$
 so dass $\varphi(x)$ zerlegbar ist.
 \par
       Die Umkehrung ist ebenso trivial zu beweisen.
 \medskip\noindent
 {\bf 4.7. Korollar.} {\it Jeder $(n - 1)$-Vektor ist zerlegbar.}
 \smallskip
       Beweis. Dies folgt aus 4.6 und 3.3, wenn man nur noch bemerkt, dass
 $1$-Formen trivialerweise zerlegbar sind.
 \medskip
       Sind $x$, $y \in \bigwedge_K(V)$, so setzen wir
 $$ x \vee y := \varphi^{-1}(\varphi(x) \wedge \varphi(y)). $$
 Mit dieser Bezeichnung gilt:
 \medskip\noindent
 {\bf 4.8. Satz.} {\it Es sei $V$ ein Vektorraum des Ranges $n$ \"uber dem
 kommutativen K\"orper $K$. Ferner sei $r < n$. Ist $u \in \bigwedge_K^r(V)$,
 so ist $u$ genau dann zerlegbar, wenn f\"ur alle zerlegbaren
 $(n - r - 1)$-Vektoren $x$ gilt, dass $u \vee(u \wedge x) = 0$ ist.}
 \smallskip
       Beweis. Es sei $u$ zerlegbar. Ferner sei $x$ ein zerlegbarer
 $(n - r - 1)$-Vektor. Ist $u \wedge x = 0$, so ist nichts zu beweisen. Es sei
 also $u \wedge x \neq 0$ und $u = z_1 \wedge \dots \wedge u_r$ und
 $x = u_{r+1} \wedge \dots \wedge u_{n - 1}$. Nach 1.10 sind die Vektoren
 $u_1$, \dots, $u_{n-1}$ linear unabh\"angig. Es gibt also einen Vektor
 $u_n$, so dass $u_1$, \dots, $u_n$ eine Basis von $V$ ist. Es sei $u_1^*$,
 \dots, $u_n^*$ die zu $u_1$, \dots, $u_n$ duale Basis. Indem wir
 gegebenenfalls $u_n$ um einen Skalarfaktor ab\"andern, d\"urfen wir
 annehmen, dass 
 $$ y = u_1 \wedge \dots \wedge u_n $$
 und
 $$ z = u_1^* \wedge \dots \wedge u_n^* $$
 ist, wobei $y$ der $n$-Vektor und $z$ die $n$-Form aus Satz 3.3 sind. Nach
 3.4 ist
 $$ \varphi(u \wedge x) = \varphi(u_1 \wedge \dots \wedge u_{n-1}) =
         (-1)^{n-1}u_n^* $$
 und
 $$ \varphi(u) = \varphi(u_1 \wedge \dots \wedge u_r) =
       (-1)^{n(n-r)}(u_{r+1}^* \wedge \dots \wedge u_n^*). $$
 Somit ist $\varphi(u) \wedge \varphi(u \wedge x) = 0$ und daher auch
 $u \vee (u \wedge x) = 0$.
 \par
       Es sei nun umgekehrt $u \vee (u \wedge x) = 0$ f\"ur alle
 $(n - r - 1)$-Vektoren $x$. Ferner sei $x_1$, \dots, $x_n$ eine Basis von
 $V$ und $x_I$ mit $I \in \Fin(\{1, \dots n\})$ sei die dieser Basis
 entsprechende Basis von $\bigwedge_K(V)$. Dann ist
 $$ u = \sum_{I \in \Fin_r(V)} x_Ik_I $$
 mit $k_I \in K$ f\"ur alle $I$. Ist $u = 0$, so ist $u$ zerlegbar. Es sei
 also $u \neq 0$. Es gibt dann ein $L \in \Fin_r(\{1, \dots, n\})$ mit
 $k_L \neq 0$. Die Menge
 $$ \{1, \dots, n\} - L $$
 enth\"alt $n - r$ Elemente und daher auch $n - r$ Teilmengen der L\"ange
 der L\"ange $n - r - 1$. Diese seien $J_1$, \dots, $J_{n-r}$. Dann ist
 $$ u \wedge x_{J_i} = \sum_{|I| = r} (x_I \wedge x_{J_i})k_I. $$
 Nun ist $x_I \wedge x_{J_i} = 0$, falls $I \cap J_i \neq \emptyset$.
 Es sei $I \cap J_i = \emptyset$. Ferner sei
 $$ \{1, \dots, n\} = L \cup J_i \cup \{a_i\}. $$
 Ist nun $I \neq L$, so folgt aus
 $$ I = I \cap (L \cup J_i \cup \{a_i\})
      = (I \cap L) \cup (I \cap J_i) \cup (L \cap \{a_i\})
      = (I \cap L) \cup (L \cap \{a_i\}), $$
 dass $a_i \in L$ und $|I \cap L| = r - 1$ ist. Also ist $I = A \cup \{a_i\}$
 mit einer $(r - 1)$-Teilmenge $A$ von $L$. Hieraus folgt, dass
 $$ u \wedge x_{J_i} = (x_L \cap x_{J_i})k_L +
 \sum_{A \in \Fin_r(L)} (x_{A \cup \{a_i\}} \wedge x_{J_i})
                                          k_{A \cup \{a_i\}} $$
 ist. Es sei $0 = \sum_{i:=1}^{n-r} (u \wedge x_{J_i})l_i$. Dann ist
 $$\eqalign{
 0 &= \sum_{i:=1}^{n-r} (x_L \wedge x_{J_i})k_Ll_i + \sum_{i:=1}^{n-r}
       \sum_{A \in \Fin_r(L)} (x_{A \cup \{a_i\}} \wedge x_{J_i})
                                          k_{A \cup \{a_i\}}l_i \cr
       &= \sum_{i:=1}^{n-r} \pm (x_{L \cup J_i})k_Ll_i +
              \sum_B x_B \kappa_B. \cr} $$
 Dabei sind die $B$ von der Form $A \cup \{a_i\} \cup J_i$. Nun folgt aus
 $L \cap J_i = \emptyset$ f\"ur alle $i$, dass $L \cup J_i = L \cup J_j$
 impliziert, dass $i = j$ ist. Ferner ist $L \cup J_i \neq B =
 A \cup \{a_j\} \cup J_j$, da sonst $L \cap (\{a_j\} \cup J_j) \neq \emptyset$
 w\"are. Daher sind die Vektoren $x_{L \cup J_i}$ und $x_M$ linear
 unabh\"angig. Hieraus folgt, dass $\pm k_Ll_i = 0$ ist f\"ur $i := 1$,
 \dots, $n - r$. Wegen $k_L \neq 0$ sind daher alle $l_i = 0$. Dies besagt,
 dass die $n - r$ Vektoren $u \wedge x_{J_i}$ linear unabh\"angig sind. Weil
 $u \wedge x_{J_i}$ ein $(n - 1)$-Vektor ist, ist $\varphi(u \wedge x_{J_i})$
 eine 1-Form, dh. es ist $\varphi(u \wedge x_{J_i}) \in V^*$. Nun ist
 $u \vee (u \wedge x_{J_i}) = 0$ f\"ur alle $i$. Daher ist auch
 $\varphi(u) \wedge \varphi(u \wedge x_{J_i}) = 0$ f\"ur alle $i$ und damit
 $$ \varphi(u \wedge x_{J_i}) \in V^*_{\varphi(u)}. $$
 Weil $\varphi$ ein Isomorphismus ist, folgt
 $$ \Rg_K(V^*_{\varphi(u)} \geq n - r. $$
 Weil $u$ ein $r$-Vektor ist, ist $\varphi(u)$ nach 3.3 eine $(n - r)$-Form.
 Nach 4.1 ist andererseits
 $$ \Rg_K(V^*_{\varphi(u)}) \leq n - r. $$
 Daher ist $\Rg_K(V^*_{\varphi(u)}) = n - r$, so dass $\varphi(u)$ nach 4.2
 eine zerlegbare $(n - r)$-Form ist. Mit 4.6 folgt, dass auch $u$
 zerlegbar ist. Damit ist der Satz bewiesen.

 \mysection{5. Doppelverh\"altnisse}

\noindent
       Doppelverh\"altnisse spielten in der projektiven Geometrie lange Zeit
 eine herausragende Rolle. Wir ben\"otigen sie nur, um projektive
 Kol\-li\-nea\-tio\-nen papposscher R\"aume zu kennzeichnen, da uns diese
 Kennzeichnung
 im n\"achsten Abschnitt bequem sein wird. Wir werden uns hier daher kurz
 fassen. Wer mehr \"uber Doppelverh\"altnisse wissen m\"ochte, der sei an
 Baer 1952, S. 71--94 verwiesen.
 \par
       Es sei $V$ ein Vektorraum des Ranges 2 \"uber dem kommutativen
 K\"orper $K$. Ferner seien $P$, $Q$, $R$ und $S$ vier verschiedene Punkte
 von $\La(V)$. Es gibt dann Vektoren $p \in P$ und $q \in Q$ mit $P = pK$,
 $Q = qK$ und $R = (p + q)K$. Es gibt ferner ein $d \in K$ mit
 $S = (p + qd)K$. Wir nennen $d$ {\it Doppelverh\"altnis\/} der Punkte
 $P$, $Q$, $R$, $S$. Es gilt nun der folgende Satz.
 \medskip\noindent
 {\bf 5.1. Satz.} {\it Es sei $V$ ein Vektorraum des Ranges 2 \"uber dem
 kommutativen K\"orper $K$. Sind $P$, $Q$, $R$ und $S$ vier verschiedene
 Punkte von $\La(V)$, so haben $P$, $Q$, $R$, $S$ genau ein Doppelverh\"altnis.
 Dieses bezeichnen wir mit $DV(P,Q;R,S)$.}
 \smallskip
       Beweis. Die Existenzaussage wurde schon bewiesen. Um die Eindeutigkeit
 zu beweisen, seien $d$ und $d'$ zwei Doppelverh\"altnisse der Punkte $P$,
 $Q$, $R$ und $S$.
 Es gibt dann von Null verschiedene Vektoren $p$, $p' \in P$ und
 $q$, $q' \in Q$ mit
 $$ R = (p + q)K = (p' + q')K $$
 und
 $$ S = (p + qd)K = (p' + q'd')K. $$
 Es gibt also $a$, $b \in K$ mit
 $$ p + q = (p' + q')a $$
 und
 $$ p + qd = (p' + q'd')b. $$
 Ferner gibt es $u$, $v \in K$ mit $p' = pu$ und $q' = qv$. Also ist
 $$ p + q = pua + qva $$
 und
 $$ p + qd = pub + qvd'b. $$
 Es folgt
 $$ 1 = va = ua = ub $$
 und damit $u = v = a^{-1}$ und $b = a$. Daher ist $d = a^{-1}d'a = d'$.
 \medskip
       Ist $K$ nicht kommutativ, so zeigt eine Analyse des gerade ge\-f\"uhr\-ten
 Beweises, dass die Doppelverh\"altnisse von vier Punkten, die genauso
 definiert werden, wie hier geschehen, eine Konjugiertenklasse von Elementen
 aus $K^*$ ausmachen. Diese Analyse ist bei Baer {\it loc.\ cit.\/}
 durchgef\"uhrt.
 \par
       Sind $X$ und $Y$ zwei Unterr\"aume des Ranges $r - 1$ und $r + 1$
 eines Vektorraumes, so hat der Faktorraum $Y/X$ den Rang 2, daher ist auch
 das Doppelverh\"altnis von vier Unterr\"aumen des Ranges $r$ definiert, solange
 diese Unter\-r\"au\-me $X$ enthalten und in $Y$ enthalten sind. Diese
 Bemerkung machen wir uns beim n\"achsten Satz zunutze.
 \medskip\noindent
 {\bf 5.2. Satz.} {\it Es sei $V$ ein Vektorraum endlichen Ranges \"uber dem
 kommutativen K\"orper $K$. Ferner sei $\Rg_K(V) \geq 3$ und $\kappa$ sei eine
 Kollineation oder Korrelation von $\La(V)$. Dann sind \"aquivalent:
 \item{a)} $\kappa$ ist projektiv.
 \item{b)} Sind $X$ und $Y$ Unterr\"aume von $V$ des Ranges $r - 1$ und $r + 1$
 und sind $P$, $Q$, $R$, $S$ verschiedene Unterr\"aume des Ranges $r$, die
 allesamt zwischen $X$ und $Y$ liegen, so ist
 $$ DV(P,Q;R,S) = DV(\kappa(P),\kappa(Q);\kappa(R),\kappa(S)). $$
 \item{c)} Es gibt zwei Unterr\"aume $X$ und $Y$ von $V$ des Ranges $r - 1$ und
 des Ranges $r + 1$, so dass
 $$ DV(P,Q;R,S) = DV(\kappa(P),\kappa(Q);\kappa(R),\kappa(S)) $$
 gilt f\"ur alle Quadrupel verschiedener Unterr\"aume $P$, $Q$, $R$, $S$ des
 Ranges $r$, die zwischen $X$ und $Y$ liegen.\par}
 \smallskip
       Beweis. Es seien $X$ und $Y$ zwei Unterr\"aume des Ranges $r - 1$ bzw.
 $r + 1$ von $V$ und $P$, $Q$, $R$ und $S$ seien vier verschiedene
 Unterr\"aume des Ranges $r$, die alle zwischen $X$ und $Y$ liegen.
 Es gibt dann Vektoren $p$ und $q$ mit $P = X + pK$, $Q = X + qK$, $R =
 X + (p + q)K$ und $S = X + (p + qd)K$. Dann ist $d$ das Doppelverh\"altnis
 der Punkte $P$, $Q$, $R$ und $S$.
 \par
       Es sei zun\"achst $\kappa$ eine Kollineation von $\La(V)$. Nach
 dem zweiten Struktursatz gibt es dann eine semilineare Abbildung $\sigma$
 von $V$, die $\kappa$ induziert. Es sei $\alpha$ der Begleitautomorphismus
 von $\sigma$. Es folgt
 $$\eqalign{
 \kappa(P) &= \kappa(X) + \sigma(p)K                       \cr
 \kappa(Q) &= \kappa(X) + \sigma(q)K                       \cr
 \kappa(R) &= \kappa(X) + (\sigma(p) + \sigma(q))K         \cr
 \kappa(S) &= \kappa(X) + (\sigma(p) + \sigma(q)d^\alpha)K \cr} $$
 und damit
 $$ DV(\kappa(P),\kappa(Q);\kappa(R),\kappa(S)) = d^\alpha. $$
 \par
       Es sei $\kappa$ eine Korrelation und $f$ sei eine $\kappa$
 darstellende $\alpha$-Se\-mi\-bi\-li\-ne\-ar\-form.
 Es gibt dann Vektoren $p'$ und $q'$
 und ein $d' \in K$ mit
 $$\eqalign{
 P^\kappa &= X^\kappa \cap (pK)^\kappa = Y^\kappa + p'K     \cr
 Q^\kappa &= X^\kappa \cap (qK)^\kappa = Y^\kappa + q'K     \cr
 R^\kappa &= X^\kappa \cap ((p + q)K)^\kappa = Y^\kappa + (p' + q')K \cr
 S^\kappa &= x^\kappa \cap ((p + qd)K)^\kappa = Y^\kappa + (p' + q'd')K.\cr} $$
 Es folgt $f(p,p') = 0 = f(q,q')$ und $f(p,q') \neq 0 \neq f(q,p')$. Weiter
 folgt
 $$ 0 = f(p + q,p' + q') = f(p,q') + f(q,p') $$
 und damit $f(p,q') = -f(q,p')$. Hieraus folgt schlie\ss lich
 $$ 0 = f(p + qd,p' + q'd') = f(p,q')(d' - d^\alpha), $$
 so dass wegen $f(p,q') \neq 0$ auch hier
 $$ DV(\kappa(P),\kappa(Q);\kappa(R);\kappa(S)) = d^\alpha $$
 ist.
 \par
       a) impliziert b): Ist $\kappa$ projektiv, so ist in den vorstehenden
 Aus\-f\"uh\-run\-gen
 $\alpha = 1$, so dass b) in der Tat eine Folge von a) ist.
 \par
       b) impliziert c): Banal.
 \par
       c) impliziert a): $\kappa$ werde durch die semilineare Abbildung
 $\sigma$ mit dem Begleitautomorphismus $\alpha$ induziert bzw. durch
 eine $\alpha$-Form dargestellt. Ferner seien $P$ und $Q$ zwei verschiedene
 Unterr\"aume des Ranges $r$, die zwischen $X$ und $Y$ liegen. Es gibt dann
 Vektoren $p$ und $q$ mit $P = X + pK$ und $Q = X + qK$. Setze
 $R := X + (p + q)K$. Es sei $d \in K$ und $d \neq 0$, 1. Setze
 $S := X + (p + qd)K$. Dann sind
 $P$, $Q$, $R$ und $S$ vier verschiedene Unterr\"aume des Range $r$, die
 alle zwischen $X$ und $Y$ liegen. Es folgt
 $$ d = DV(P,Q,R,S) = DV(\kappa(P),\kappa(Q),\kappa(R),\kappa(S)) = d^\alpha.$$
 Es ist also $d^\alpha = d$ f\"ur alle $d \in K$, die von 0 und 1 verschieden
 sind. Es ist aber auch $0^\alpha = 0$ und $1^\alpha = 1$. Daher ist
 $\alpha = 1_K$, so dass $\kappa$ projektiv ist. Damit ist alles bewiesen.
 \medskip
       Das Semikolon in $DV(P,Q;R,S)$ soll andeuten, dass es Symmetrien gibt.
 So ist beispielsweise $DV(P,Q;R,S) = DV(Q,P;S,R)$.
 Hierauf werden wir aber nicht weiter eingehen. Gesagt sei nur, dass
 $DV$ bei den 24 Permutationen der $P$, $Q$, $R$, $S$ h\"ochstens sechs
 verschiedene Werte annimmt.

 \mysection{6. Gra\ss mannsche Mannigfaltigkeiten I}
 
\noindent
       Es sei $V$ ein Vektorraum des Ranges $n$ \"uber dem kommutativen
 K\"orper $K$. Die Menge $G_r(V)$ der Punkte von $\La(\bigwedge_K^r(V))$,
 die durch zerlegbare Vektoren auf\-ge\-spannt werden, nennen wir
 {\it gra\ss mannsche Mannigfaltigkeit\/} mit den Parametern $n$ und $r$.
 Die Bilder von $G_r(V)$ unter der Kollineationsgruppe von
 $\La(\bigwedge_K^r(V))$ nennen wir ebenfalls {\it gra\ss mannsche
 Mannigfaltigkeiten\/}. Zwei gra\ss mannsche Man\-nig\-fal\-tig\-kei\-ten $G$ in
 $\La(X)$ und $H$ in $\La(Y)$ hei\ss en {\it projektiv \"aquivalent\/}, falls
 es einen Isomorphismus $\sigma$ von $\La(X)$ auf $\La(Y)$ gibt mit
 $\sigma(G) = H$.
 \par
       Wie zuvor bezeichnen wir mit $\UR_r(V)$ die Menge der Unterr\"aume
 des Ranges $r$ von $V$. Ist $U \in \UR_r(V)$ und sind $u_1$, \dots, $u_r$
 und $v_1$, \dots, $v_r$ Basen von $U$, so ist
 $$ (u_1 \wedge \dots \wedge u_r)K = (v_1 \wedge \dots \wedge v_r)K $$
 nach 4.3 und nach 1.10 ist $u_1 \wedge \dots \wedge u_r \neq 0$. Setzt man
 nun
 $$ \gamma(U) := (u_1 \wedge \dots \wedge u_r)K, $$
 so ist $\gamma$ also wohldefiniert und \"uberdies eine Abbildung von
 $\UR_r(V)$ in $G_r(V)$. Ist andererseits $(x_1 \wedge \dots \wedge x_r)K$
 ein Punkt von $G_r(V)$, so sind die $x_i$ linear unabh\"angig und es ist
 $$\textstyle \gamma\bigl(\sum_{i:=1}^r x_iK\bigr) = (x_1 \wedge \dots \wedge x_r)K, $$
 so dass $\gamma$ auch surjektiv ist. Wir nennen $\gamma$
 {\it gra\ss mannsche Abbildung\/} von $\UR_r(V)$ auf $G_r(V)$.
 \par
       Sind $X$ und $Y$ Unterr\"aume von $V$ mit $X \leq Y$ und ist
 $\Rg_K(X) \leq r \leq \Rg_K(Y)$, so setzen wir
 $$ \UR_r(Y,X) := \{U \mid U \in \UR_r(V), X \leq U \leq Y\}. $$
 \par
       Die wichtigste Eigenschaft von $\gamma$ wird durch den folgenden
 Satz ausgedr\"uckt.
 \medskip\noindent
 {\bf 6.1. Satz.} {\it Es sei $V$ ein Vektorraum des Ranges $n$ \"uber dem
 kommutativen K\"orper $K$. Ferner seien $X$ und $Y$ zwei Unterr\"aume von
 $V$ und es gelte $X \leq Y$. Setze $t := \Rg_K(X)$ und $s := \Rg_K(Y)$. Es sei
 $x_1$, \dots, $x_t$ eine Basis von $X$. Ist
 dann $t \leq r \leq s$, so gibt es einen Monomorphismus $\eta$ von
 $\bigwedge_K^{r-t}(Y/X)$ in $\bigwedge_K^r(V)$ mit
 $$ \eta((y_1 + X) \wedge \dots \wedge (y_{r-t} + X)) =
       x_1 \wedge \dots \wedge x_t \wedge y_1 \wedge \dots \wedge y_{r-t} $$
 f\"ur alle $y_1$, \dots, $y_{r-t} \in Y$. Setzt man
 $$ G := \{\gamma(U) \mid U \in \UR_r(Y,X)\}, $$
 so ist weiter $G = \eta(G_{s-t,r-t}(Y/X))$. Insbesondere ist $G$ 
 eine zu $G_{s-t,r-t}(Y/X)$ isomorphe Teilmannigfaltigkeit von $G_r(V)$.
 Die Punkte von $G$ spannen einen Teilraum des Ranges $s - t \choose r - t$
 auf.}
 \smallskip
       Beweis. Sind $y_1$ \dots, $y_{r-t} \in Y$, so setzen wir
 $$ \psi(y_1 + X, \dots, y_{r-t} + X) := x_1 \wedge \dots \wedge x_r \wedge
       y_1 \wedge \dots \wedge y_{r-t}. $$
 Ist $u \in X$, so ist
 $$ x_1 \wedge \dots \wedge x_t \wedge u = 0. $$
 Daher ist $\psi$ wohldefiniert. Weil $\psi$ banalerweise $(r - t)$-fach
 alternierend ist, gibt es eine lineare Abbildung $\eta$ von 
 $\bigwedge_K^{r-t}(Y/X)$ in $\bigwedge_K^r(V)$ mit
 $$ \eta((y_1 + X) \wedge \dots \wedge (y_{r-t} + X)) =
       x_1 \wedge \dots \wedge x_t \wedge y_1 \wedge \dots \wedge y_{r-t}. $$
 Mittels einer Basis $z_1 + X$, \dots, $z_{s-t} + X$ von $Y/X$ sieht man,
 dass $\eta$ mindestens den Rang $s - t \choose r - t$ hat. Weil das Urbild
 aber diesen Rang hat, folgt, dass $\eta$ ein Monomorphismus ist. Hieraus
 folgt wiederum, dass
 $$ H := \eta(G_{s-t,r-t}(Y/X)) $$
 eine zu $G_{s-t,r-t}(Y/X)$ isomorphe
 Teilmannigfaltigkeit von $G_r(V)$ ist. Ferner ist $H \subseteq G$.
 Es sei nun $P \in G$ und $U := \gamma^{-1}(P)$. Dann ist $X \leq U \leq Y$.
 Es gibt daher Elemente $y_1$, \dots, $y_{r-t} \in Y$, so dass $x_1$, \dots,
 $x_r$, $y_1$, \dots, $y_{r-t}$ eine Basis von $U$ ist. Wegen
 $$\eqalign{
 P &= \gamma(U) = (x_1 \wedge \dots \wedge x_t \wedge y_1 \wedge \dots
                     \wedge y_{r-t})K \cr
   &= \eta((y_1 + X) \wedge \dots \wedge (y_{r-t} + X))K \cr} $$
 ist sogar $H = G$. Damit ist alles bewiesen.
 \medskip
       Die Vektoren aus $\bigwedge_K^1(V)$ sind alle zerlegbar und nach 4.7
 sind auch alle Vektoren aus $\bigwedge_K^{n-1}(V)$ zerlegbar, falls $n$ der
 Rang von $V$ ist. Hieraus folgt, dass $G_{n,1}(V)$ aus den Punkten von
 $\La(\bigwedge_K^1(V))$ und $G_{n,n-1}(V)$ aus den Punkten von
 $\La(\bigwedge_K^{n-1}(V))$ besteht. Aus diesen Bemerkungen folgen nun sofort
 die Korollare 6.2 und 6.3.
 \medskip\noindent
 {\bf 6.2. Korollar.} {\it Es sei $V$ ein Vektorraum des Ranges $n$ \"uber dem
 kommutativen K\"orper $K$. Ferner seien $X$ und $Y$ zwei Teilr\"aume des
 Ranges $r - 1$ und $s$ von $V$ und es gelte $X \leq Y$. Ist wieder
 $$ G := \{\gamma(U) \mid U \in \La_r(Y,X) \}, $$
 so besteht $G$ aus den Punkten eines Unterraumes $W$ des Ranges $s - r + 1$.
 \"Uberdies wird die Einschr\"ankung $\gamma_0$ von $\gamma$ auf
 $\La_r(Y,X)$ durch eine lineare Abbildung von $Y/X$ auf $W$ induziert.}
 \medskip\noindent
 {\bf 6.3. Korollar.} {\it Es sei $V$ ein Vektorraum des Ranges $n$ \"uber dem
 kommutativen K\"orper $K$. Ferner seien $X$ und $Y$ zwei Teilr\"aume des
 Ranges $t$ und $r + 1$ von $V$ und es gelte $X \leq Y$. Dann besteht
 $$ G := \{\gamma(U) \mid U \in \UR_r(Y,X))\} $$
 aus den Punkten eines Unterraumes $W$ des Ranges $r + 1 - t$. \"Uberdies wird
 die Einschr\"ankung $\gamma_1$ von $\gamma$ auf $\La_r(Y,X)$ durch eine
 lineare Abbildung von $Y/X$ auf $W$ induziert.}
 \medskip
       Setzt man $s := n$ in 6.2, so folgt, dass $G_r(V)$ eine Schar
 $S^I$ von Teilr\"aumen des Ranges $n - r + 1$ enth\"alt. Setzt man $t := 0$
 in 6.3, so sieht man, dass $G_r(V)$ auch eine Schar $S^{II}$
 von Teilr\"aumen des Ranges $r + 1$ enth\"alt. Wir werden sp\"ater sehen,
 dass alle Teilr\"aume von $\La(\bigwedge_K^r(V))$, deren Punkte zu
 $G_r(V)$ geh\"oren, in einem Raum aus $S^I \cup S^{II}$ enthalten sind.
 \medskip\noindent
 {\bf 6.4. Satz.} {\it Es sei $V$ ein Vektorraum endlichen Ranges \"uber dem
 kommutativen K\"orper $K$. Ferner sei $U \in \UR_r(V)$ und $V = U \oplus C$.
 Ist dann $S$ die Menge der Punkte der Form $(c \wedge u_1 \wedge \dots \wedge
 u_{r-1})K$ mit $c \in C$ und $u_i \in U$, so ist $S$ zur Segreschen
 Mannigfaltigkeit $S(C,U^*)$ projektiv \"aquivalent. Ist $W$ ein erzeugender
 Raum der ersten Schar von $S$, so ist $\gamma(U) + W \in S^I$. Ist $W$ ein
 erzeugender Raum der zweiten Schar von $S$, so ist $\gamma(U) + W \in S^{II}$.}
 \smallskip
       Beweis. Es sei $\epsilon$ die Inklusionsabbildung von $U$ in
 $\bigwedge_K(U)$ und $\iota$ sei die Inklusionsabbildung von $U$ in
 $V \leq \bigwedge_K(V)$. Dann ist $\iota(u) \wedge \iota(u) = u \wedge u = 0$
 f\"ur alle $u \in U$. Es gibt also einen Algebrenhomomorphismus $\tau$ von
 $\bigwedge_K(U)$ in $\bigwedge_K(V)$ mit $\tau\epsilon(u) = \iota(u) = u$
 f\"ur alle $u \in U$. Insbesondere ist
 $$ \tau((\epsilon(u_1) \wedge \dots \wedge \epsilon(u_{r-1})) =
     \tau\epsilon(u_1) \wedge \dots \wedge \tau\epsilon(u_{r-1}) =
      u_1 \wedge \dots \wedge u_{r-1}. $$
 Also ist $\tau(\bigwedge_K^{r-1}(U)) \leq \bigwedge_K^{r-1}(V)$. Da
 $\Rg_K(U) = r$ ist, sind die Elemente aus $\bigwedge_K^{r-1}(U)$ alle
 zerlegbar. Es sei $x \in \bigwedge_K^{r-1}(U)$ und es gelte $\tau(x) = 0$.
 Es gibt dann $u_1$, \dots, $u_{r-1} \in U$ mit
 $$ x = \epsilon(u_1) \wedge \dots \wedge \epsilon(u_{r-1}). $$
 Es folgt
 $$ 0 = \tau(\epsilon(u_1) \wedge \dots \wedge \epsilon(u_{r-1})
       = u_1 \wedge \dots \wedge u_{r-1}, $$
 so dass die $u_i$ linear abh\"angig sind. Dann sind aber auch die
 $\epsilon(u_i)$ linear abh\"angig, so dass
 $$ x = \epsilon(u_1) \wedge \dots \wedge \epsilon(u_{r-1}) = 0 $$
 ist. Somit ist die Einschr\"ankung von $\tau$ auf $\bigwedge_K^{r-1}(U)$ ein
 Monomorphismus von $\bigwedge_K^{r-1}(U)$ in $\bigwedge_K^{r-1}(V)$. Setze
 $$ \tilde U := \{u_1 \wedge \dots \wedge u_{r-1} \mid u_i \in U\}. $$
 Dann ist $\tilde U$ ein zu $\bigwedge_K^{r-1}(U)$ isomorpher Unterraum von
 $\bigwedge_K^{r-1}(V)$, wie gerade gezeigt.
 \par
       Die durch
 $$ (u_1 \wedge \dots \wedge u_{r-1})^*u := u_1 \wedge \dots \wedge u_{r-1}
                                                 \wedge u $$
 definierte Abbildung $*$ ist eine lineare Abbildung von $\tilde U$ in $U^*$.
 Aus $(u_1 \wedge \dots \wedge u_{r-1})^* = 0$ folgt, da $\Rg_K(U) = r$ ist,
 dass $u_1 \wedge \dots \wedge u_{r-1} = 0$ ist. Daher ist $*$ injektiv.
 Schlie\ss lich folgt aus
 $$ \Rg_K(\tilde U) = \Rg_K({\scriptstyle\bigwedge}_K^{r-1}(U))
                      = r = \Rg_K(U^*), $$
 dass $*$ sogar bijektiv ist.
 \par
       Die Abbildung
 $$ (c,(u_1 \wedge \dots \wedge u_{r-1})^*) \to c \wedge u_1 \wedge \dots
                                   \wedge u_{r-1} $$
 ist eine bilineare Abbildung von $C \times U^*$ in $\bigwedge_K^r(V)$. Es
 gibt folglich eine lineare Abbildung $\sigma$ von $C \otimes U^*$ in
 $\bigwedge_K^r(V)$ mit
 $$ \sigma(c \otimes (u_1 \wedge \dots \wedge u_{r-1})^*) =
                     c \wedge u_1 \wedge \dots u_{r-1}. $$
 Der von den $c \wedge u_1 \wedge \dots \wedge u_{r-1}$ aufgespannte
 Unterraum $X$ hat mindestens den Rang $(n - r)r$. Ist n\"amlich $c_1$, \dots,
 $c_{n-r}$ eine Basis von $C$ und $b_1$, \dots, $b_r$ eine Basis von $U$,
 so ist wegen $U \cap C = \{0\}$ die Menge der $r$-Vektoren der Form
 $$ c_j \wedge b_1 \wedge \dots \wedge b_{i-1} \wedge b_{i+1} \wedge \dots
              \wedge b_r $$
 nach Satz 1.7 linear unabh\"angig. Die Anzahl dieser Vektoren ist aber
 gleich $(n - r)r$. Weil $X$ ein epimorphes Bild von $C \otimes U^*$ ist
 und weil der Rang dieses Raumes gleich $(n - r)r$ ist, folgt, dass der
 Rang von $X$ gleich $(n - r)r$ ist. Dies hat wiederum zur Folge, dass
 $\sigma$ ein Isomorphismus von $C \otimes U^*$ auf $X$ ist. Offensichtlich
 gilt auch
 $$ \sigma(S(C,U^*)) = S. $$
 Damit ist die erste Aussage des Satzes bewiesen.
 \par
       Als n\"achstes zeigen wir, dass $\gamma(U) \cap X = \{0\}$ ist.
 Dazu nehmen wir an, dies sei nicht der Fall. Da $\gamma(U)$ ein Punkt ist,
 ist dann $\gamma(U) \leq X$, so dass es ein $c \in C$ und $u_i$ in $U$
 gibt mit
 $$ \gamma(U) = (c \wedge u_1 \wedge \dots \wedge u_{r-1})K. $$
 Weil $c \wedge u_1 \wedge \dots \wedge u_{r-1} \neq 0$ ist, sind die $u_i$
 linear unabh\"angig. Es gibt also ein $u_0 \in U$, so dass $u_0$, \dots,
 $u_{r-1}$ eine Basis von $U$ ist. Dann ist aber $\gamma(U) =
 (u_0 \wedge \dots \wedge u_{r-1})K$. Es gibt also ein $k \in K^*$ mit
 $$ c \wedge u_1 \wedge \dots \wedge u_{r-1} =
       (u_0 \wedge \dots u_{r-1})k. $$
 Es folgt
 $$ (c - u_0k) \wedge u_1 \wedge \dots u_{r-1} = 0, $$
 so dass $c - u_0k$, $u_1$, \dots, $u_{r-1}$ linear abh\"angig sind.
 Weil die $u_1$, \dots, $u_{r-1}$ linear unabh\"angig sind, folgt, dass
 $c - u_0k$ von $u_1$, \dots, $u_{r-1}$ linear abh\"angt. Dies hat
 schlie\ss lich $c \in U$ zur Folge. Damit erhalten wir den Widerspruch
 $$ 0 \neq c \in C \cap U = \{0\}. $$
 Also ist doch $\gamma(U) \cap X = \{0\}$.
 \par
       Es sei nun $W$ ein Raum der ersten Schar von $S$. Es gibt dann
 $u_1$, \dots, $u_{r-1} \in U$ mit
 $$ W = \{c \wedge u_1 \wedge \dots \wedge u_{r-1} \mid c \in C\}. $$
 Es gibt ferner ein $u_0$, so dass $u_0$, $u_1$, \dots, $u_{r-1}$ eine
 Basis von $U$ ist. Es folgt
 $$ \gamma(U) = \{(u_0k) \wedge u_1 \wedge \dots \wedge u_{r-1} \mid k \in K\}. $$
 Somit ist
 $$ \gamma(U) + W = \{(u_0k + c) \wedge u_1 \wedge \dots \wedge u_{r-1} \mid
                            k \in K, c \in C\}. $$
 Hieraus folgt, dass $\Rg_K(\gamma(U) + W) = n - r + 1$ ist. Ferner ist klar,
 dass die Punkte von $\gamma(U) + W$ gerade den Unterr\"aumen vom Rang $r$
 entsprechen, die $u_1$, \dots, $u_r$ enthalten. Also ist $\gamma(U) + W \in
 S^I$. 
 \par
       Es sei nun $W$ ein Raum der zweiten Schar von $S$. Es gibt dann ein
 $c \in C$ mit $c \neq 0$ und
 $$ W = \{c \wedge u_1 \wedge \dots \wedge u_{r-1} \mid u_i \in U\}. $$
 Es sei $v_0$, \dots, $v_{r-1}$ eine Basis von $U$. Dann ist $\gamma(U) =
 (v_0 \wedge \dots \wedge v_{r-1})K$. Es sei $P$ ein Punkt von $\gamma(U) + W$
 und es sei $P = xK$. Es gibt dann ein $k \in K$ und $u_1$, \dots, $u_{r-1}
 \in U$ mit
 $$ x = (v_0k) \wedge v_1 \wedge \dots \wedge v_{r-1} +
              c \wedge u_1 \wedge \dots \wedge u_{r-1}. $$
 Wir d\"urfen annehmen, dass $u_1 \wedge \dots \wedge u_{r-1} \neq 0$ ist.
 Dann gibt es ein $u_0$, so dass $u_0$, $u_1$, \dots, $u_{r-1}$ eine
 Basis von $U$ ist. Es folgt $\gamma(U) = (u_0 \wedge \dots \wedge u_{r-1})K$,
 so dass es ein $l \in K$ gibt mit
 $$\eqalign{
   x &= (u_0l) \wedge u_1 \wedge \dots \wedge u_{r-1} +
              c \wedge u_1 \wedge \dots \wedge u_{r-1}      \cr
       &= (u_0l + c) \wedge u_1 \wedge \dots \wedge u_{r-1}. \cr} $$
 Hieraus folgt
 $$ \gamma^{-1}(P) = V_x \leq U \oplus cK. $$
 \par
       Es sei umgekehrt $X$ ein Unterraum des Ranges $r$ von $U \oplus cK$.
 Ist $X = U$, so ist $\gamma(X) \leq \gamma(U) + W$. Es sei also $X \neq U$.
 Dann ist $\Rg_K(X \cap U) = r - 1$. Es sei $u_1$, \dots, $u_{r-1}$ eine
 Basis von $U \cap X$ und $v_0$, $u_1$, \dots, $u_{r-1}$ eine Basis von $X$.
 Es gibt dann ein $u \in U$ und ein $k \in K$ mit $v = u + ck$. Es folgt
 $$ v \wedge u_1 \wedge \dots \wedge u_{r-1} =
 u \wedge u_1 \wedge \dots \wedge u_{r-1} + (c \wedge u_1 \wedge \dots \wedge
       u_{r-1})k .$$
 Folglich ist $\gamma(X) \leq \gamma(U) + W$. Damit ist gezeigt, dass
 die Punkte von $\gamma(U) + W$ gerade die Bilder der Unterr\"aume des
 Ranges $r$ von $U \oplus cK$ sind. Folglich ist $\gamma(U) + W \in S^{II}$.
 Damit ist alles bewiesen.
 \medskip
       Mit $S_{\gamma(U)}$ bezeichnen wir die Menge der R\"aume
 $\gamma(U) + W$ mit $W \in E_1 \cup E_2$, wobei $E_i$ die $i$te Schar von
 Erzeugenden von $S$ ist. Wir nennen $S_\gamma(U)$ {\it Segreschen Kegel mit
 der Spitze\/} $\gamma(U)$ von $G_r(V)$. Den von den R\"aumen aus
 $S_{\gamma(U)}$ aufgespannte Unterraum bezeichnen wir mit $T_{\gamma(U)}$
 und nennen ihn {\it Tangentialraum\/} von $G_r(V)$ in $\gamma(U)$.
 Wie der Beweis von 6.4 zeigt, gilt
 \medskip\noindent
 {\bf 6.5. Korollar.} {\it Es ist $\Rg_K(T_{\gamma(U)}) = (n - r)r + 1$.}
 \medskip\noindent
 {\bf 6.6. Satz.} {\it Es sei $V$ ein Vektorraum des Ranges $n$ \"uber dem
 kommutativen K\"orper $K$. Ist $U$ ein Unterraum von $\bigwedge_K^r(V)$,
 dessen Punkte alle in $G_r(V)$ liegen, so ist $U$ Unterraum eines
 Raumes aus $S^I \cup S^{II}$.}
 \smallskip
       Beweis. Es seien $P$ und $Q$ zwei verschiedene Punkte von $U$.
 Ferner sei $P = (u_1 \wedge \dots \wedge u_r)K$ und $Q = (v_1 \wedge \dots
 \wedge v_r)K$. Weil $P$ und $Q$ verschieden sind, sind $u_1 \wedge \dots
 \wedge u_r$ und $v_1 \wedge \dots \wedge v_r$ linear unabh\"angig. Wegen
 $P + Q \leq U$ und weil die Punkte von $U$ alle zu $G_r(V)$ geh\"oren
 gibt es $w_1$, \dots, $w_r$ mit
 $$ u_1 \wedge \dots \wedge u_r + v_1 \wedge \dots \wedge v_r =
       w_1 \wedge \dots \wedge w_r. $$
 \"Uberdies ist $w_1 \wedge \dots \wedge w_r \neq 0$, da $u_1 \wedge \dots
 \wedge u_r$ und $v_1 \wedge \dots \wedge v_r$ linear unabh\"angig sind.
 Schlie\ss lich folgt, dass $V_{w_1 \wedge \dots \wedge w_r} \neq
 V_{v_1 \wedge \dots \wedge v_r}$ ist. Wir d\"urfen daher oBdA annehmen,
 dass $w_1 \not\in V_{v_1 \wedge \dots \wedge v_r}$ ist. Nun ist
 $$ 0 = w_1 \wedge \dots \wedge w_r \wedge w_1 = u_1 \wedge \dots \wedge u_r
       \wedge w_1 + v_1 \wedge \dots \wedge v_r \wedge w_1 $$
 und daher
 $$ u_1 \wedge \dots \wedge u_r \wedge w_1 =
      -v_1 \wedge \dots \wedge v_r \wedge w_1. $$
 Setze $y := u_1 \wedge \dots \wedge u_r \wedge w_1$. Aus der gerade bewiesenen
 Gleichung folgt dann, da $w_1 \not\in V_{v_1 \wedge \dots \wedge v_r}$ gilt,
 dass $y \neq 0$ ist. Nun ist
 $$ y \wedge u_i = 0 = y \wedge v_i $$
 f\"ur alle $i$. Daher ist
 $$ V_{u_1 \wedge \dots \wedge u_r} \leq V_y $$
 und
 $$ V_{v_1 \wedge \dots \wedge v_r} \leq V_y. $$
 Weil $\Rg_K(V_y) = r + 1$ ist,
 folgt, dass $P + Q$ in einem Raum aus $S^{II}$ liegt, da ja
 $\gamma(V_{u_1 \wedge \dots \wedge u_r}) = P$ und $\gamma(V_{v_1 \wedge
 \dots \wedge v_r}) = Q$ ist. Weil $\Rg_K(V_y) = r + 1$ ist folgt, dass
 $\Rg_K(V_{u_1 \wedge \dots \wedge u_r} \cap V_{v_1 \wedge \dots \wedge v_r})
 = r - 1$ ist. Daher liegt $P + Q$ auch in einem Teilraum aus $S^I$.
 Die beiden R\"aume aus $S^I \cup S^{II}$ geh\"oren zum Segreschen Kegel
 $S_P$ mit der Spitze $P$. Hieraus folgt, falls $X$ der von der zu $P$
 geh\"orenden Segreschen Mannigfaltigkeit $S$ aufgespannte Raum ist, dass
 $U \cap X$ ganz in $S$ liegt. Aus VI.3.9 und 6.4 folgt nun die Behauptung.
 \medskip\noindent
 {\bf 6.7. Satz.} {\it Es sei $V$ ein Vektorraum \"uber dem kommutativen K\"orper
 $K$ und der Rang von $V$ sei endlich. Sind $U$, $U' \in \UR_r(V)$, so sind
 \"aquivalent:
 \item{a)} Es ist $\Rg_K(U + U') = r + 1$.
 \item{b)} Es ist $\Rg_K(U \cap U') = r - 1$.
 \item{c)} Es ist $\gamma(U) \neq \gamma(U')$ und alle Punkte von
 $\gamma(U) + \gamma(U')$ liegen in $G_r(V)$.}
 \smallskip
       Beweis. Es ist
 $$ \Rg_K(U + U') + \Rg_K(U \cap U') = \Rg_K(U) + \Rg_K(U') = 2r. $$
 Hieraus folgt, dass a) und b) \"aquivalent sind.
 \par
       Aus a) und b) zusammen folgt einmal $\gamma(U) \neq \gamma(U')$ und
 mit 6.1 weiterhin, dass
 $$ \gamma(U), \gamma(U') \in G \cong G_1((U + U')/(U \cap U')) $$
 ist. Dies besagt, dass $\gamma(U)$ und $\gamma(U')$ auf einer Geraden
 liegen, deren Punkte alle zu $G_r(V)$ geh\"oren.
 \par
       Aus c) folgen, wie wir beim Beweise von 6.6 gesehen haben, a) und b).
 Damit ist alles gezeigt.
 \medskip\noindent
 {\bf 6.8. Satz.} {\it Es sei $V$ ein Vektorraum des Ranges $n$ \"uber dem
 kommutativen K\"orper $K$. Ferner seien $X$ und $Y$ Unterr\"aume des Ranges
 $r - 1$ bzw. $r + 1$ von $V$. Sind dann $P$, $Q$, $R$, $S$ vier verschiedene
 Unterr\"aume des Ranges $r$ von $Y$, die $X$ enthalten, und ist $\gamma$
 die gra\ss mannsche Abbildung von $\UR_r(V)$ auf $G_r(V)$, so
 gilt
 $$ DV(P,Q;R,S) = DV(\gamma(P),\gamma(Q);\gamma(R),\gamma(S)). $$
 Sind umgekehrt $P'$, $Q'$, $R'$ und $S'$ vier verschiedene, kollineare
 Punkte von $G_r(V)$, so ist}
 $$ DV(\gamma^{-1}(P'),\gamma^{-1}(Q');\gamma^{-1}(R'),\gamma^{-1}(S'))
       = DV(P',Q';R',S'). $$
 \par
       Beweis. Es sei $u_1$, \dots, $u_{r-1}$ eine Basis von $X$. Ferner sei
 $d := DV(P,Q;R,S)$. Es gibt Vektoren $p$ und $q$ mit $P = X + pK$, $Q =
 X + qK$, $R = X + (p + q)K$ und $S = X + (p + qd)K$. Hieraus folgt
 $$\eqalign{
 \gamma(P) &= (u_1 \wedge \dots \wedge u_{r-1} \wedge p)K             \cr
 \gamma(Q) &= (u_1 \wedge \dots \wedge u_{r-1} \wedge q)K             \cr
 \gamma(R) &= (u_1 \wedge \dots \wedge u_{r-1} \wedge (p + q))K       \cr
 \gamma(S) &= (u_1 \wedge \dots \wedge u_{r-1} \wedge (p + qd))K.     \cr} $$
 Setzt man
 $$\eqalign{
 p' &:= u_1 \wedge \dots \wedge u_{r-1} \wedge p \cr
 q' &:= u_1 \wedge \dots \wedge u_{r-1} \wedge q,\cr} $$
 so ist also $\gamma(P) = p'K$, $\gamma(Q) = q'K$, $\gamma(R) = (p' + q')K$
 und $\gamma(S) = (p' + q'd)K$. Also ist
 $$ DV(\gamma(P),\gamma(Q);\gamma(R),\gamma(S)) = d. $$
 \par
       Dass $\gamma^{-1}$ ebenfalls das Doppelverh\"altnis invariant
 l\"asst, beweist sich ebenso einfach.
 \medskip\noindent
 {\bf 6.9. Satz.} {\it Es sei $V$ ein Vektorraum des Ranges $n$ \"uber
 dem kommutativen K\"orper $K$. Ferner sei $\gamma$ die gra\ss mannsche
 Abbildung von $\UR_r(V)$ auf $G_r(V)$ und $\delta$ sei die
 gra\ss mannsche Abbildung von $\UR_s(V)$ auf $G_s(V)$. Ist $\kappa$
 ein Isomorphismus von $\La(\bigwedge_K^r(V))$ auf $\La(\bigwedge_K^s(V))$,
 der $G_r(V)$ auf $G_s(V)$ abbildet, so gibt es
 eine Kollineation oder eine Korrelation $\lambda$ von $\La(V)$
 mit $\lambda(U) = \delta^{-1}\kappa\gamma(U)$ f\"ur alle $U \in \UR_r(V)$.
 Ist $\lambda$ eine Kollineation, so ist $r = s$, und ist $\lambda$ eine
 Korrelation, so ist $r + s = n$. Genau dann ist $\kappa$ projektiv, wenn
 $\lambda$ projektiv ist.}
 \smallskip
       Beweis. Es seien $U$ und $U'$ benachbarte Unterr\"aume des Ranges $r$.
 Dann ist $\Rg_K(U + U') = r + 1$, so dass nach 6.7 alle Punkte von
 $\gamma(U) + \gamma(U')$ in $G_r(V)$ liegen. Daher liegen alle Punkte
 von $\kappa\gamma(U) + \kappa\gamma(U')$ in $G_s(V)$. Nach 6.7 sind
 folglich auch $\delta^{-1}\kappa\gamma(U)$ und $\delta^{-1}\kappa\gamma(U')$
 benachbart, so dass die Existenz von $\lambda$ aus dem Satz von Chow
 (Satz I.8.4) folgt, wie auch die Aussage \"uber $r$ und $s$.
 \par
       Weil es Geraden gibt, deren Punkte allesamt in $G_r(V)$ liegen,
 folgt die letzte Behauptung aus den S\"atzen 6.8 und 5.2.
 \medskip\noindent
 {\bf 6.10. Korollar.} {\it Es sei $V$ ein Vektorraum des Ranges $n$ \"uber
 dem kommutativen K\"orper $K$. Genau dann sind $G_r(V)$ und $G_s(V)$
 projektiv \"aquivalent, wenn $r = s$ oder wenn $r + s = n$ ist.}
 \smallskip
       Beweis. Es ist aufgrund von 6.9 nur noch zu zeigen, dass $G_r(V)$
 und $G_s(V)$ projektiv \"aquivalent sind, wenn $r + s = n$ ist. Dies folgt
 aber aus 3.3 und 4.6, wenn man nur noch $V$ mit $V^*$ \"uber eine Basis
 von $V$ und deren Dualbasis von $V^*$ identifiziert, was ja wegen der
 Kommutativit\"at von $K$ m\"oglich ist.
 \medskip\noindent
 {\bf 6.11. Satz.} {\it Es sei $V$ ein Vektorraum des Ranges $n$ \"uber dem
 kommutativen K\"orper $K$. Der Monomorphismus $\sigma \to \sigma_{\#r}$
 induziert einen Monomorphismus von $PGL(V)$
 auf eine Untergruppe $\Gamma_r^0(V)$ des
 Stabilisators $\Gamma_r(V)$ von $G_r(V)$ in $PGL(\bigwedge_K^r(V))$. Es
 gilt: 
 \item{a)} Ist $n \neq 2r$, so ist $\Gamma_r(V) = \Gamma_r^0(V)$.
 \item{b)} Ist $n = 2r$, so ist $|\Gamma_r(V):\Gamma_r^0(V)| = 2$.}
 \smallskip
       Beweis. Aus Satz 1.11 folgt, dass die Abbildung
 $$ \sigma \to \sigma_{\#r} $$
 einen Homomorphismus von $PGL(V)$ auf eine Untergruppe $\Gamma_r^0(V)$ von
 $\Gamma_r(V)$ induziert. L\"asst $\sigma_{\#r}$ alle Punkte von
 $G_r(V)$ fest, so ist also
 $$ \sigma(v_1) \wedge \dots \wedge \sigma(v_r) =
       \sigma_{\#r}(v_1 \wedge \dots \wedge v_r) = (v_1 \wedge \dots \wedge
              v_r)k $$
 f\"ur alle $(v_1, \dots, v_r) \in V^r$, wobei $k$ von den $v_i$ abh\"angt.
 Es folgt, dass $\sigma$ alle Unterr\"aume des Ranges $r$ von $V$
 festl\"asst, so dass die Abbildung $\sigma$ in $\La(V)$ die Identit\"at
 induziert. Dies zeigt, dass die Abbildung $\sigma \to \sigma_{\#r}$
 einen Isomorphismus von $PGL(V)$ auf $\Gamma_r^0(V)$ induziert.
 \par
       Nach 6.9 gibt es zu jedem $\kappa \in \Gamma_r(V)$ eine Kollineation
 oder Korrelation $\lambda$ von $\La(V)$ mit
 $$ \lambda(U) = \gamma^{-1}\kappa\gamma(U) $$
 f\"ur alle $U \in \UR_r(V)$. Es sei $\Gamma_r^1(V)$ die Untergruppe aller
 $\kappa \in \Gamma_r(V)$, f\"ur die $\lambda$ eine Kollineation ist. Es
 sei $U = \bigoplus_{i:=1}^r u_iK$ und $\kappa \in \Gamma_r^1(V)$. Schlie\ss
 lich sei $\lambda \in GL(V)$ und $\lambda$ induziere die zu $\kappa$
 geh\"orende Kollineation in $\La(V)$, die wir ebenfalls mit $\lambda$
 bezeichnen. Dann ist
 $$ \gamma\lambda(U) = (\lambda(u_1) \wedge \dots \wedge \lambda(u_r))K
    = \lambda_{\#r}(u_1 \wedge \dots \wedge u_r)K = \lambda_{\#r}\gamma(U).$$
 Somit ist $\gamma\lambda = \lambda_{\#r}\gamma$, dh.,
 $$ \gamma^{-1}\kappa\gamma = \lambda = \gamma^{-1}\lambda_{\#r}\gamma, $$
 so dass $\kappa = \lambda_{\#r}$ ist. Folglich ist $\Gamma_r^1(V) =
 \Gamma_r^0(V)$.
 \par
       a) In diesem Falle ist $\Gamma_r^1(V) = \Gamma_r(V)$ und damit
 $\Gamma_r^0(V) = \Gamma_r(V)$.
 \par
       b) Weil das Produkt zweier Korrelationen eine Kollineation ist, ist
 $$ |\Gamma_r(V):\Gamma_r^0| = |\Gamma_r(V):\Gamma_r^1(V)| \leq 2. $$
 Die Kollineationen aus $\Gamma_r^0(V)$ lassen die Scharen $S^I$ und
 $S^{II}$ je f\"ur sich invariant. L\"asst man $r$ in 6.10 die Rollen von
 $r$ und von $s$ spielen, was wegen $n = 2r$ ja m\"oglich ist, so sieht man,
 dass es auch Kollineationen in $\Gamma_r(V)$ gibt, die die Scharen $S^I$
 und $S^{II}$ vertauschen. Damit ist alles bewiesen.

 \mysection{7. Gra\ss mannsche Mannigfaltigkeiten II}

\noindent
       Wir haben in Abschnitt 4 Kriterien f\"ur die Zerlegbarkeit von
 $r$-Vektoren gegeben. In diesem Abschnitt werden wir nun eine geo\-me\-tri\-sche
 Deutung dieser Kriterien geben. Sie lassen sich n\"amlich dahingehend
 interpretieren, dass sich jede gra\ss mannsche Mannigfaltigkeit als Schnitt
 von endlich vielen Quadriken darstellen l\"asst.
 \par
       Es sei $V$ ein Vektorraum \"uber dem kommutativen K\"orper $K$ und $f$ sei
 das in Satz 2.12 beschriebene Skalarprodukt auf
 $\bigwedge_K(V^*) \times \bigwedge_K(V)$. Ist $x$ ein zerlegbarer
 $(n - r - 1)$-Vektor und $w$ eine zerlegbare $(r - 1)$-Form, so definieren
 wir die Abbildung $Q_{w,x}$ von $\bigwedge_K^r(V)$ in $K$ durch
 $$ Q_{w,x}(z) := f(w,z \vee (z \wedge x)). $$
 Banal ist, dass
 $Q_{w,x}(zk) = Q_{w,x}(z)k^2$ ist. Langweilige Routinerechnungen zeigen, wobei
 man auf die Definition von $\vee$ vor Satz 4.8 zur\"uckgreifen muss, dass
 $$ Q_{w,x}(z + z') = Q_{w,x}(z) + Q_{w,x}(z') + f(w,z \vee (z' \wedge x)) +
                     f(w,z' \vee (z \wedge x)) $$
 ist. Ebensolche Routinerechnungen zeigen, dass die Abbildung
 $$ (z,z') \to f(w,z \vee (z' \wedge x)) + f(w,z' \vee (z \wedge x)) $$
 bilinear ist. Daher ist $Q_{w,x}$ eine quadratische Form.
 \medskip\noindent
 {\bf 7.1. Satz.} {\it Ist $\sigma \in GL(V)$, ist $x$ ein zerlegbarer
 $(n - r - 1)$-Vektor und $w$ eine zerlegbare $(r - 1)$-Form, ist
 schlie\ss lich $\bar x = \sigma_{\#(n-r-1)}^{-1}(x)$ und $\bar w =
 \sigma_{\#(r-1)}^*(w)$, so sind $Q_{w,x}$ und $Q_{\bar w,\bar x}$ projektiv
 \"aquivalent.}
 \smallskip
       Beweis. Nach 3.5 gibt es ein $k \in K^*$ mit $\varphi\sigma_{\#r}\mu_k
 = \sigma_{\#(n-r)}^{*-1}\varphi$. Benutzt man dies und 1.11 d), so folgt:
 $$\eqalign{
 \varphi(\sigma_{\#r}(z)) \vee (\sigma_{\#r}(z) \wedge x)
  &= \varphi(\sigma_{\#r}(z)) \wedge \varphi(\sigma_{\#r}(z) \wedge
              \sigma_{\#(n-r-1)}(\bar x)) \cr
  &= \sigma_{\#(n-r)}^{*-1}\varphi(z) \wedge \sigma^{*-1}\varphi(\bar x)k^{-2}
                                         \cr
 &= \sigma_{\#(n-r+1)}^{*-1}(\varphi(z) \wedge \varphi(z \wedge \bar x))k^{-2}.
                                         \cr} $$
 Also ist
 $$ \sigma_{\#r}(z) \vee (\sigma_{\#r} \wedge x) = 
       \sigma^{*-1}(\varphi(z) \wedge \varphi(z \wedge \bar x))k^{-2}. $$
 Nach 3.5 gilt auch $\varphi\sigma_{\#(r-1)}\mu_k = \sigma_{\#(n-r+1)}\varphi$,
 da $k$ von $r$ unabh\"angig ist. Also ist
 $$ \sigma_{\#r}(z) \vee (\sigma_{\#r}(z) \wedge x) =
       \sigma_{\#(r-1)}(z \vee (z \wedge \bar x))k^{-1}. $$
 Hieraus folgt schlie\ss lich, dass
 $$ Q_{w,x}(\sigma_{\#r}(z)) = f(w,z \vee (z \wedge \bar x))k^{-1}
       = k^{-1}Q_{\bar w, \bar x}(z) $$
 ist. Folglich sind $Q_{w,x}$ und $Q_{\bar w,\bar x}$ projektiv \"aquivalent.
 \medskip
       F\"ur den Rest des Abschnitts vereinbaren wir das folgende. Es sei $V$
 ein Vektorraum des Ranges $n$ \"uber dem kommutativen K\"orper $K$. Es sei
  $x_1$, \dots, $x_n$ eine Basis von $V$ und $y = x_1 \wedge \dots
 \wedge x_n$. Ferner sei $x_1^*$, \dots, $x_n^*$ die zu $x_1$, \dots, $x_n$
 duale Basis von $V^*$. Mit $(x_I \mid I \subseteq \{1, \dots, n\})$ bzw.
 $(x_J^* \mid J \subseteq \{1, \dots, n\})$ seien die entsprechenden Basen von
 $\bigwedge_K(V)$ bzw. $\bigwedge_K(V^*)$ bezeichnet. Statt
 $Q_{x_J^*,x_I}$ schreiben wir $Q_{J,I}$.
 \medskip\noindent
 {\bf 7.2. Satz.} {\it Es seien $I$, $I'$, $J$, $J'$ Teilmengen von $\{1,\dots, n
 \}$. Gilt $|I| = n - r - 1 = |I'|$ und $|J| = r - 1 = |J'|$ sowie
 $|I \cap J| = |I' \cap J'|$, so sind $Q_{J,I}$ und $Q_{J',I'}$ projektiv
 \"aquivalent.}
 \smallskip
       Beweis. Aufgrund der Annahme \"uber $I$, $J$, $I'$ und $J'$ gibt es
 ein $\pi \in S_n$ mit $\pi(I) = I'$ und $\pi(J) = J'$. Ferner gibt es genau
 ein $\sigma \in GL(V)$ mit $\sigma(x_i) = x_{\pi^{-1}(i)}$ f\"ur alle $i$.
 Dann ist $\sigma^{-1}(x_i) = x_{\pi(i)}$ und $\sigma^*(x^*_j) = x^*_{\pi(j)}$,
 wie man leicht nachrechnet. Es folgt $\sigma_{\#n}(y) = \sgn(\pi)y$ und
 $\sigma_{\#(n-r-1)}(x_I) = \epsilon_1x_{I'}$ und $\sigma^*_{\#(r-1)}(x^*_J)
 = \epsilon_2x^*_{J'}$, wobei $\epsilon_i = \pm1$ ist. Hieraus folgt, wie
 der Beweis von 7.1 zeigt, dass   
 $$ Q_{J,I}(\sigma_{\#r}(z)) = \epsilon_1\epsilon_2\sigma(\pi)Q_{J',I'}(z) $$
 ist. Damit ist der Satz bewiesen.
 \medskip
       Es sei $z \in \bigwedge_K^r(V)$. Dann ist $z = \sum_{|L| = r} x_La_L$.
 Nach 3.4 ist daher
 $$ \varphi(z) = \sum_{|L| = r} x^*_{L^c}a_L \prod_{\lambda \in L,
       \mu \in L^c} \langle \mu,\lambda \rangle. $$
 Dabei ist $L^c$ das Komplement von $L$ in $\{1, \dots, n\}$. Ferner ist
 $$ z \wedge x_I = \sum_{|L|=r} x_{L\cup I} a_L \prod_{\alpha \in L,\beta \in
                            I} \langle \alpha,\beta \rangle, $$
 und daher
 $$ \varphi(z \wedge x_I) = \sum_{|L|=r} x^*_{L^c \cap I^c} a_L
       \prod_{\alpha \in L,\beta \in I} \langle \alpha,\beta \rangle
       \prod_{\gamma \in L \cup I, \delta \in L^c \cap I^c} \langle \delta,
                            \gamma \rangle. $$
 Somit ist
 $$ \varphi(z) \wedge \varphi(z \wedge x_I) = \sum_{|L|=r} \sum_{|M|=r}
 x^*_{L^c\cup(M^c\cap I)} a_La_MB_{L,M},  $$
 wobei
 $$\eqalign{
    &B_{L,M} = \cr
    &   \prod_{\lambda \in L,\mu \in L^c} \langle \mu,\lambda \rangle
         \prod_{\alpha \in M,\beta \in I} \langle \alpha,\beta \rangle
         \prod_{\gamma \in M \cup I,\delta \in M^c \cap I^c} \langle \gamma,
                            \delta \rangle
         \prod_{\epsilon \in L^c,\zeta \in M^c \cap I^c} \langle \epsilon,
                            \zeta \rangle \cr } $$
 ist. Schlie\ss lich ist
 $$
   z\vee (z \wedge x_I) = 
       \sum_{|L|=r}\sum_{|M|=r} x_{L \cap(M \cup I)}
        a_La_MB_{L,M}\prod_{{\eta \in L \cap (M \cup I)} \atop {\theta \in L^c \cup
                                  (M^c \cap I^c)}} \langle \theta,\eta \rangle.
 $$
 Nun sind $(x_R \mid R)$ und $(x^*_S \mid S)$ duale Basen. Folglich ist
 $$ f(x_{L \cap (M \cup I)},x^*_J) = 0, $$
 falls $L \cap (M \cup I) \neq J$ ist. Also ist
 $$ Q_{J,I} = \epsilon_J \sum_{|L| = |M| = r, L \cap(M \cup I) = J}
                     a_La_MB_{L,M} $$
 mit
 $$ \epsilon_J = \prod_{\eta \in J,\theta \in J^c} \langle \theta,\eta
                            \rangle. $$
 \par
       Wir wollen sehen, wie wir den Ausdruck f\"ur $Q_{J,I}$ noch
 vereinfachen k\"onnen. Dazu fragen wir zun\"achst, wann unter der
 Nebenbedingung $L \cap (M \cup I) = J$ der Koeffizient $B_{L,M} \neq 0$ ist.
 Dazu ist, wie der eben etablierte Ausdruck f\"ur $B_{L,M}$
 zeigt, hinreichend und notwendig, dass $M \cap I = \emptyset$
 und $L \cup M \cup I = \{1, \dots, n\}$ ist. Aus $M \cap I = \emptyset$
 folgt, dass
 $$ |M \cup I| = r + n - r - 1 = n - 1 $$
 ist. Es gibt also ein $i_M$ mit
 $$ M \cup I \cup \{i_M\} = \{1, \dots, n\} $$
 und
 $$ i_M \not\in M \cup I. $$
 Aus $L \cup M \cup I = \{1, \dots, n\}$ folgt, dass $i_M \in L$ ist, und
 $L \cap (M \cup I) = J$ und $i_M \not\in M \cup I$ implizieren, dass
 $i_M \not\in J$ ist. Also ist $|J \cup \{i_M\}| = r$. Nun ist $L
 \subseteq J \cup \{i_M\}$ und $|L| = r$, so dass $L = J \cup \{i_M\}$
 ist. Wir setzen nun
 $$ C_{L,M} := \prod_{\alpha \in M,\beta \in I} \langle \alpha,\beta \rangle
       \prod_{\gamma \in L^c} \langle \gamma,i_M \rangle. \leqno (i) $$
 Dann gilt
 \medskip\noindent
 {\bf 7.3. Satz.} {\it Es sei $|I| = n - r - 1$ und $|J| = r - 1$. Dann ist
 $$ Q_{J,I}(z) = \epsilon_J\epsilon'_J\sum_{|M|=r,M\cap I = \emptyset, 
                                          J \subseteq M \cup I}
                             a_Ma_{J \cup \{i_M\}}C_{J \cup \{i_M\},M}. $$
 Dabei ist $i_M$ das Element aus der Gleichung
 $ M \cup I \cup \{i_M\} = \{1, \dots, n\}$.
 Ferner ist
 $$ \epsilon_J = \prod_{\eta \in J, \delta \in J^c} \langle \theta,\eta
                                                  \rangle\qquad
 \mathit{und}\qquad
  \epsilon'_J = (-1)^{nr + {1 \over 2}r(r - 1) + \Sigma J - 1}, $$
 wobei $\Sigma X$ abk\"urzend f\"ur $\sum_{y\in X} y$ steht.}
 \smallskip
       Beweis. Setze $L := J \cup \{i_M\}$. Es sei $L = \{\lambda_1, \dots,
 \lambda_r\}$ mit $\lambda_1 < \dots < \lambda_r$. Dann ist
 $$\eqalign{
 \prod_{\lambda \in L,\,\mu \in L^c} \langle \mu,\lambda \rangle
       = \prod_{i:=1}^r\prod_{\mu \in L^c} \langle \mu,\lambda_i \rangle 
       &= \prod_{i:=1}^r (-1)^{n-\lambda_i - (r - i)}                     \cr
       &= (-1)^{nr+{1\over2}r(r-1)+i_M+\Sigma J}. \cr} $$
 Ferner ist $M^c \cap I^c = \{i_M\}$, so dass
 $$ \prod_{\gamma \in M \cup I,\delta \in M^c \cap I^c} \langle \delta,
       \gamma \rangle = \sum_{\gamma \in M \cup I} \langle i_M,\gamma \rangle
              = (-1)^{i_M-1} $$
 ist. Schlie\ss lich ist
 $$ \prod_{\epsilon \in L^c,\zeta \in M^c \cap I^c} \langle \epsilon,\zeta
       \rangle = \prod_{\epsilon \in L^c} \langle \epsilon,i_M \rangle. $$
 Daher ist
 $$
 B_{L,M} = (-1)^{mr+{1\over2}r(r-1)+i_M+\Sigma J+i_M-1} C_{L,M} 
         = \epsilon'_JC_{L,M}.  $$
 Hieraus folgt schlie\ss lich die Behauptung.
 \medskip\noindent
 {\bf 7.4. Satz.} {\it Ist $I \cap J = \emptyset$, so ist $Q_{J,I} = 0$.}
 \smallskip
       Beweis. Aufgrund von 7.2 d\"urfen wir annehmen, dass
 $$ I = \{1, \dots, n - r - 1\} $$
 und
 $$ J := \{n - r, \dots, n - 2\} $$
 ist. Aus $J \subseteq M \cup I$ und $I \cap J = \emptyset$ folgt, dass
 $J \subseteq M$ ist. F\"ur $M$ gibt es daher nur die beiden M\"oglichkeiten
 $M_1 = J \cup \{n - 1\}$ und $M_2 = J \cup \{n\}$. Dann ist aber
 $$ L_1 = J \cup \{n\} = M_2 $$
 und
 $$ L_2 = J \cup \{n - 1\} = M_1. $$
 Nun ist $\beta < \alpha$, falls $\beta \in I$ und $\alpha \in M_i$ ist.
 Also ist
 $$ \prod_{\alpha\in M_1,\beta \in I} \langle \alpha,\beta \rangle =
    \prod_{\alpha\in M_2,\beta \in I} \langle \alpha,\beta \rangle. $$
 Diesen gemeinsamen Wert nennen wir $\epsilon$.
 \par
       Es ist $i_{M_1} = n \in L_1$ und daher
 $$ \prod_{\gamma\in L_1^c} \langle \gamma,i_{M_1} \rangle = 1. $$
 Ferner ist $i_{M_2} = n - 1 \in L_2$ und $n \in L_2^c$. Daher ist
 $$ \prod_{\gamma\in L_2^c} \langle \gamma,i_{M_2} \rangle = -1. $$
 Somit ist
 $$ B_{L_1,M_1} = \epsilon = -B_{L_2,M_2}. $$
 Hieraus folgt schlie\ss lich, dass
 $$ Q_{J,I}(z) = \epsilon_J\epsilon'_J\epsilon(a_{L_1}a_{M_1} - a_{L_2}
                     a_{M_2})
       = \epsilon_J\epsilon'_J\epsilon(a_{M_2}a_{M_1} - a_{M_1}a_{M_2})
       = 0 $$
 ist. Damit ist 7.4 bewiesen.
 \medskip
       Wir sind nun in der Lage, den folgenden Satz zu beweisen.
 \medskip\noindent
 {\bf 7.5. Satz.} {\it Ist $V$ ein Vektorraum des Ranges $n$ \"uber dem
 kommutativen K\"orper $K$, so ist $G_r(V)$ Schnitt von
 $$ {n \choose r + 1}\biggl[{n \choose r - 1} - {r + 1 \choose r - 1}\biggr] $$
 Quadriken in $\La(\bigwedge_K^r(V))$.}
 \smallskip
       Beweis. Nach 4.8 ist $z \in \bigwedge_K^r(V)$ genau dann zerlegbar,
 wenn f\"ur alle zerlegbaren $(n - r - 1)$-Vektoren $x$ die Gleichung
 $z \vee (z \wedge x) = 0$ gilt. Dies ist nat\"urlich genau dann der Fall,
 wenn $z \vee (z \wedge x_I) = 0$ f\"ur alle $(n - r - 1)$-Teilmengen $I$
 gilt. Nun ist $z \vee (z \wedge x_I)$ ein $(r - 1)$-Vektor. Also ist
 $z \vee (z \wedge x_J)$ genau dann gleich Null, wenn
 $$ Q_{J,I}(z) = 0 $$
 f\"ur alle $I$ mit $|I| = r - 1$ ist. Es folgt, dass $G_r(V)$ Schnitt der
 Quadriken zu den
 $$ {n \choose n - r - 1}{n \choose r - 1} = {n \choose r + 1}
                            {n \choose r - 1} $$
 quadratischen Formen $Q_{J,I}$ ist. Von diesen tragen aber nach 7.4 diejenigen
 zu dem Schnitt nichts bei, f\"ur die $I \cap J = \emptyset$ ist. Deren Anzahl
 ist
 $$ {n \choose n - r - 1}{r + 1 \choose r - 1}
       = {n \choose r + 1}{r + 1 \choose r - 1}. $$
 Hieraus folgt die Behauptung.\footnote{Anmerkung der Herausgeber: der Autor hat hier
 einen Hinweis auf Fehler bei Bertini, Krull und Lense notiert.}

 \medskip
       Wir betrachten den Spezialfall $n = 4$, $r = 2$, der schon von
 Pl\"ucker untersucht wurde. In diesem Falle ist $n - r - 1 = 1 = r - 1$.
 Hier sind also die quadratischen Formen $Q_{\{k\},\{k\}}$ f\"ur $k := 1$,
 \dots, 4 zu betrachten. Berechnet man sie mit Hilfe von 7.3, so erh\"alt
 man, wobei der unwesentliche Faktor $\epsilon_J\epsilon'_J$ nicht
 ber\"ucksichtigt ist,
 $$\eqalign{
 Q_{1,1}(z) &= a_{12}a_{34} - a_{13}a_{24} + a_{14}a_{23} \cr
 Q_{2,2}(z) &= a_{12}a_{34} - a_{24}a_{13} + a_{23}a_{14} \cr
 Q_{3,3}(z) &= a_{34}a_{12} - a_{13}a_{24} + a_{23}a_{14} \cr
 Q_{4,4}(z) &= a_{34}a_{12} - a_{24}a_{13} + a_{14}a_{23}. \cr} $$
 Die vier Quadriken, die uns Satz 7.5 liefert, sind also nicht voneinander
 verschieden. Dies notieren wir als
 \medskip\noindent
 {\bf 7.6. Satz.} {\it Ist $V$ ein Vektorraum des Ranges $4$ \"uber dem
 kommutativen K\"orper $K$, so ist $G_2(V)$ eine Quadrik, die durch die Form
 $$ Q(z) := a_{12}a_{34} - a_{13}a_{24} + a_{14}a_{23} $$
 dargestellt wird. Sie ist von maximalem Index.}
 \medskip
       Man nennt $G_2(V)$ {\it Pl\"uckerquadrik\/}, falls $\Rg_K(V) = 4$ ist.
 \par
       Wie der Fall der Pl\"uckerquadrik zeigt, kann es sein, dass
 $Q_{J,I}$ und $Q_{J',I'}$ die gleiche Quadrik darstellen, ohne dass
 $I = I'$ und $J = J'$ ist. \"Uber diese Situation gibt der n\"achste Satz
 Auskunft. Um ihn zu formulieren, ben\"otigen wir noch die folgende
 Bezeichnung. Sind $X$ und $Y$ Mengen, so bezeichne $X \bt Y$ ihre
 symmetrische Differenz, das hei\ss t die Menge
 $$ (X \cup Y) - (X \cap Y). $$
 Die Potenzmenge einer Menge versehen mit der symmetrischen Differenz als
 Verkn\"upfung ist bekanntlich eine elementarabelsche 2-Grup\-pe.
 \medskip\noindent
 {\bf 7.7. Satz.} {\it Es sei $V$ ein Vektorraum des Ranges $n$ \"uber dem
 kommutativen K\"orper $K$. Ferner seien $I$ und $I'$ zwei
 $(n - r - 1)$- und $J$ und $J'$ zwei $(r - 1)$-Teilmengen von $\{1, \dots, n
 \}$ mit $I \cap J \neq \emptyset \neq J' \cap I'$. Genau dann gibt es ein
 $k \in K^*$ mit $Q_{J,I} = kQ_{J',I'}$, wenn
 $I = I'$ und $J = J'$ oder wenn $I \bt I' = J \bt J'$ und $|I \bt I'| = 2$
 ist.}
 \smallskip
       Beweis. Es seien $I$ und $J$ Teilmengen von $\{1, \dots, n\}$ mit
 $|J| = r - 1$, $|I| = n - r - 1$ und $I \cap J \neq \emptyset$.
 Nach 7.3 ist bis auf ein Vorzeichen, das f\"ur unsere Untersuchungen
 irrelevant ist,
 $$ Q_{J,I}(z) = \sum_{|M|=r,\,M\cap I = \emptyset,\, 
                                          J \subseteq M \cup I}
                             a_Ma_{J \cup \{i_M\}}C_{J \cup \{i_M\},M}. $$
 Dabei ist $i_M$ das Element aus der Gleichung
 $ M \cup I \cup \{i_M\} = \{1, \dots, n\}$.
 \par
       Es seien $I'$ und $J'$ zwei weitere Mengen der gleichen Art und es
 gebe ein $k \in K^*$ mit $Q_{J,I} = kQ_{J',I'}$. W\"ahle ein $M$ mit
 $|M| = r$, $M \cap I = \emptyset$ und $J \subseteq M \cup I$. Wir definieren
 $z$ durch $a_M := 1$, $a_{J \cup \{i_M\}} := 1$ und $a_L := 0$ f\"ur alle
 von $M$ und $J \cup \{i_M\}$ verschiedenen $L$. Dann ist
 $$ Q_{J,I}(z) = a_Ma_{J \cup \{i_M\}}C_{J \cup \{i_M\},M} \neq 0. $$
 Es folgt, dass auch $Q_{J',I'}(z) \neq \emptyset$ ist. Es gibt also
 ein $P$ mit $|P| = r$, $P \cap I' = \emptyset$ und $J' \subseteq P \cup I'$,
 so dass
 $$ a_Pa_{J' \cup \{i'_P\}} \neq 0 $$
 ist. Wegen $P \neq J' \cup \{i'_P\}$ gibt es nun zwei F\"alle, n\"amlich den
 \item{Fall 1:} Es ist $M = P$ und $J \cup \{i_M\} = J' \cup \{i'_M\}$
 \par\noindent
 und den
 \item{Fall 2:} Es ist $P = J \cup \{i_M\}$ und $M = J' \cup \{i'_P\}$.
 \smallskip
       Wir zeigen, dass genau dann $I = I'$ gilt, wenn $J = J'$ ist. Dazu
 sei zun\"achst $I = I'$. In Fall 1 ergibt sich $i_M = i'_M$ und wegen
 $i_M \not\in J$ und $i_M = i'_M \not\in J'$ und $J' \cup \{i'_M\} =
 J \cup \{i_M\}$ dann $J = J'$. In Fall 2 ergibt sich der Widerspruch
 $$ \emptyset = I \cap M = I' \cap (J' \cup \{i'_P\}) \supseteq I' \cap J'
       \neq \emptyset. $$
 \par
       Ist umgekehrt $J = J'$, so folgt in Fall 1 aus $J \cup \{i_M\} =
 J' \cup \{i'_M\}$, dass $i_M = i'_M$ ist. Dies impliziert wiederum
 $I = I'$. In Fall 2 folgt zun\"achst aus $I' \cap P = \emptyset$,
 dass $I' \cap J = \emptyset$ ist. Damit erhalten wir den Widerspruch
 $$ \emptyset \neq I' \cap J' = I' \cap J = \emptyset. $$
 \par
       Es sei also $I \neq I'$ und $J \neq J'$. Wegen $|I| = |I'|$ ist dann
 $$ |I - I'| = |I' - I| \neq 0. $$
 Es ist
 $$ |I^c \cap J| = |J| - |I \cap J| = r - 1 - |I \cap J| \leq r - 2. $$
 Es sei $x \in I' - I$. Dann ist folglich
 $$ \big|(J - I) \cup \{x\}\big| \leq r - 1. $$
 Es gibt also ein $M$ mit $|M| = r$, $I \cap M = \emptyset$ und
 $(J - I) \cup \{x\} \subseteq M$. Mit diesem $M$ liegt wegen $M \cap I' \neq
 \emptyset$ Fall 2 vor. Es ist also
 $$ M = J' \cup \{i'_P\} $$
 und
 $$ P = J \cup \{i_M\}. $$
 Wegen $x \in I'$ ist $x \neq i'_P$, so dass wegen $x \in M$ dann $x \in J'$
 gilt. Weil $x$ irgendein Element aus $I' - I$ ist, folgt $I' - I \subseteq
 J' \cap I'$. Wegen $\emptyset = I \cap M$ und $J' \subseteq M$ folgt
 $I \cap J' = \emptyset$. Daher ist
 $$ I' - I = J' \cap I'. $$
 entsprechend folgt
 $$ I' \cap J = \emptyset $$
 und  
 $$ I - I' = J \cap I. $$
 Setze $s := |I' - I|$. Dann ist auch $s = |I - I'|$. Es folgt
 $$ |I \cup (J \cap I^c) \cup (I' - I)|
       = n - r - 1 + r - 1 - s + s = n - 2. $$
 W\"are nun $s > 1$, so g\"abe es also ein $M$ mit $|M| = r$, $M \cap I =
 \emptyset$, $J \cap I^c \subseteq M$ und $|M \cap I'| = s - 1$, so dass
 auch $M \cap I' \neq \emptyset$ w\"are. Dann w\"are aber $J' \subseteq M$
 im Widerspruch zur Wahl von $M$. Dieser Widerspruch zeigt, dass doch
 $s = 1$ ist.
 \par
       Es ist
 $$ |I \cup I'| = |I \cup (I' - I)| = n - r - 1 + 1 = n - r. $$
 W\"ahlt man nun $M := (I \cup I')^c$, so ist man in Fall 1 und es gilt
 $$ J \cup \{i_M\} = J' \cup \{i'_M\}. $$
 Es folgt
 $$ J = (J \cap J') \cup \{i'_M\} $$
 und
 $$ J' = (J \cap J') \cup \{i_M\}. $$
 Es ist
 $$ i'_M \in (I' \cup M)^c = (I' \cup (I \cup I')^c)^c = I'^c \cap (I \cup I')
       = I'^c \cap I. $$
 Also ist
 $$ I = (I \cap I') \cup \{i'_M\}. $$
 Ebenso folgt
 $$ I' = (I \cap I') \cup \{i_M\}. $$
 Damit ist gezeigt, dass
 $$ I \bt I' = \{i'_M,i_M\} = J \bt J' $$
 ist. Damit ist die Notwendigkeit der Bedingungen des Satzes gezeigt.
 \par
       Die Bedingungen des Satzes seien erf\"ullt. Ist $I = I'$ und $J = J'$,
 so ist $Q_{I,J} = Q_{J',I'}$. Es sei also $I \bt I' = J \bt J'$ und
 $|I \bt I'| = 2$. Auf Grund von Satz 7.2 d\"urfen wir annehmen, dass
 folgendes gilt:
 $$\eqalign{
  I \cap I' &= \{1, \dots, n - r - 2\}             \cr
          I &= (I \cap I') \cup \{n - r - 1\}      \cr
         I' &= (I \cap I') \cup \{n - r\}          \cr
  J \cap J' &= \{n - r + 1, \dots, n - 2\}         \cr
          J &= \{n - r - 1\} \cup (J \cap J')      \cr
         J' &= \{n - r\} \cup (J \cap J').         \cr} $$
 F\"ur $M$, $i_M$ und $P$ gibt es wegen $I \cap M = \emptyset$ und
 $J \cap J' \subseteq M$ nur die folgenden drei M\"oglichkeiten:
 $$\eqalign{
 M_1 &= (J \cap J') \cup \{n - 1,n\} \cr
 M_2 &= (J \cap J') \cup \{n - r,n\} \cr
 M_3 &= (J \cap J') \cup \{n - r,n - 1\} \cr} $$
 Es folgt $i_{M_1} = n - r$, $i_{M_2} = n - 1$ und $i_{M_3} = n$.
 F\"ur die zugeh\"origen $P$ erhalten wir
 $$\eqalign{
 P_1 &= (J \cap J') \cup \{n - 1,n\} = M_1 \cr
 P_2 &= J \cup \{n - 1\}                   \cr
 P_3 &= J \cup \{n\}.                       \cr} $$
 Es folgt $i'_{P_1} = n - r - 1$, $i'_{P_2} = n$ und $i'_{P_3} = n - 1$.
 Es sind nun die Koeffizienten $C_{J \cup \{i_M\},M}$ und
 $C'_{J' \cup \{i_P\},P}$ auszurechnen, wobei die Koeffizienten von
 $Q_{J',I'}$ mittels $I'$ zu berechnen sind. Zur Erinnerung, es ist
 $$C_{J \cup \{i_M\},M} = \prod_{\alpha \in M, \beta \in I} \langle \alpha,
       \beta \rangle \prod_{\gamma \in (J \cup \{i_M\})^c} \langle \gamma,
       i_M \rangle. $$
 Mit $M := M_1$ und $P := P_1$ ergibt sich
 $$ C_{J \cup \{i_M\},M} = (-1)^{r(n-r-1)}(-1)^2 = (-1)^{rn} $$
 und
 $$ C'_{J' \cup \{i'_P\},P} = (-1)^{r(n-r-1)}(-1)^2 = (-1)^{rn}. $$
 Mit $M := M_2$ und $P := P_2$ ergibt sich
 $$ C_{J \cup \{i_M\},M} = (-1)^{r(n-r-1)}(-1) = (-1)^{rn-1} $$
 und
 $$ C'_{J' \cup \{i'_P\},P} = (-1)^{r(n-r-1)-1} = (-1)^{rn-1}. $$
 Mit $M := M_3$ und $P := P_3$ ergibt sich
 $$ C_{J \cup \{i_M\},M} = (-1)^{r(n-r-1)} = (-1)^{rn} $$
 und
 $$ C'_{J' \cup \{i'_P\},P} = (-1)^{r(n-r-1)-1}(-1) = (-1)^{rn}. $$
 Es ist also in der Tat
 $$ Q_{J,I} = kQ_{J',I'}, $$
 wobei
 $$ k = \epsilon_J\epsilon'_J\epsilon'_{J'}\epsilon''_{J'} $$
 ist.
       Damit ist der Satz bewiesen.
 \medskip
       Beim Beweise dieses Satzes wurde wieder wesentlich von der
 Kommutativit\"at von $K$ Gebrauch gemacht.
 \par
       Gleiche Quadriken zu definieren, ist nat\"urlich eine
 \"Aqui\-va\-lenz\-re\-la\-tion auf der Menge der $Q_{J,I}$. Ist $r = 2$, so ist
 es relativ einfach, die Anzahl der \"Aqui\-va\-lenz\-klas\-sen abzuz\"ahlen. Dies liegt
 daran, dass in diesem Falle $J \subseteq I$ und $|J| = 1$ gilt.
 \medskip\noindent
 {\bf 7.8. Satz.} {\it Es sei $V$ ein Vektorraum des Ranges $n \geq 4$ \"uber dem
 kommutativen K\"orper $K$. Dann ist $G_2(V)$ Schnitt von
 $$ \sum_{i:=0}^{n-4} {i + 2 \choose i}(n - 3 - i) $$
 Quadriken von $\La(\bigwedge_K^2(V))$.}
 \smallskip
       Beweis. Es sei $I$ eine $(n - 3)$-Teilmenge von $\{1, \dots, n\}$ und
 es sei $J = \{k\}$ mit $k \in I$. Ist $1 \not\in I$, so setzen wir
 $I' := I \bt \{1,k\}$ und $J' := I \bt \{1,k\}$. Dann ist $J' = \{1\}$ und
 $|I'| = n - 3$ sowie $1 \in I'$. Mit Satz 7.7 folgt, dass $Q_{J,I}$ und
 $Q_{J',I'}$ die gleiche Quadrik darstellen. Daher d\"urfen wir im Folgenden
 stets annehmen, dass $1 \in I$ gilt.
 \par
       Mit $Z_j$ bezeichnen wir die Menge der ersten $j$ nat\"urlichen
 Zahlen. Ferner setzen wir
 $$ \Pi_i := \{I \mid |I| = n - 3, Z_{n-3-i} \subseteq I, n - 2 - i \not\in
       I\} $$
 f\"ur $i := 0$, \dots, $n - 4$.
 \par
 Es ist $\Pi_0 = \{Z_{n-3}\}$. Daher liefern die $I$ aus $\Pi_0$ --- es gibt
 nur eines --- genau $n - 3$ quadratische Formen $Q_{J,I}$. Es sei bereits
 bewiesen, dass die
 $ I \in \bigcup_{j:=0}^i \Pi_j $
 genau
 $ \sum_{j:=0}^i {j + 2 \choose j} (n - 3 - j) $
 verschiedene quadratische Formen liefern.
 \par
       Es sei $I \in \Pi_{i+1}$. Ferner sei $x \in I$. Setze $Z :=
 Z_{n-3-(i+1)}$. Ist $x \not\in Z$, so folgt
 $$ |I - (Z \cup \{x\})| = n - 3 - (n - 3 - i - 1) - 1 = i. $$
 Setze $L := I - (Z \cup \{x\}$ und $I' := Z_{n-3-i} \cup L$, so ist 
 $I' \in \Pi_i$ und es folgt mit 7.7, dass $Q_{\{x\},I}$ und
 $Q_{\{n - 4 - i\},I'}$ die gleiche Quadrik darstellen. Wir erhalten also
 h\"ochstens dann eine neue quadratische Form, wenn $x \in Z$ ist. Ist
 $x \in Z$, so kommt $Q_{\{x\},I}$ unter den bereits konstruierten Formen
 wegen 7.7 auch nicht vor. Schlie\ss lich folgt, dass die neu konstruierten
 Formen auch alle untereinander verschieden sind. Wegen $Z = Z_{n-3-i-1}
 \subseteq I$ und $n - 2 - i - 1 \not\in I$ gibt es f\"ur $I - Z$
 dann noch
 $$ i + 3 \choose i + 1 $$
 M\"oglichkeiten. Damit liefern die $I \in \Pi_{i+1}$ weitere
 $$ {i + 1 + 2 \choose i + 1}(n - 3 - (i + 1)) $$
 Formen. Damit ist der Satz bewiesen.
 \medskip
       F\"ur $n = 4$ erhalten wir, wie schon zuvor, dass $G_2(V)$ eine
 Quadrik ist. F\"ur $n = 5$, 6, 7, 8, 9 ist $G_2(V)$ Schnitt von 5, 15, 35,
 70, 126 Quadriken in $\La(\bigwedge_K^2(V))$.

\markboth{}{7. Gra\ss mannsche Mannigfaltigkeiten II}


\chapter*{Anhang}{}
\addcontentsline{toc}{chapter}{\protect\numberline{Anhang: Der Aufbruch der Geometrie}}
\markboth{Anhang}{Der Aufbruch der Geometrie}

\vspace*{-10pt}

{\bf\Large\noindent
Der Aufbruch der Geometrie um Reinhold Baer
in Frankfurt und seine Protagonisten

\vspace*{3pt}\noindent --- in memoriam Heinz L\"uneburg}

\vspace*{20pt}

\noindent Vortrag von Karl Strambach
beim Baer--Kolloquium\\ in Kaiserslautern am 7.\ November 2009

\vspace*{30pt}

\noindent\qquad\quad Liebe Frau L\"uneburg, lieber Martin, liebe Freunde ! 
\medskip 

\noindent
 Als ich versp\"atet, durch pers\"onliche Umst\"ande bedingt, in Mai 1961 begann, in Frankfurt
 Physik zu studieren, haben mich, obgleich von der Mathematik voll\-kom\-men unbeleckt, die
 Mathematikvorlesungen am meisten angezogen. Es er\"offnete sich mir eine neue faszinierende Welt, 
 die ich n\"aher kennenlernen wollte, obgleich ich mir meines Unverm\"ogens halbwegs bewusst war.
 Es war ein Gl\"ucksfall, dass ich Annemarie Schlette traf, die mir begeistert vom Baerschen Laden
 erz\"ahlte, wo man von allen Zw\"angen frei Mathematik lernen k\"onne, und mir half, bei Helmut
 Salzmann einen Proseminarvortrag zu ergattern. So kam ich 1962 ins Baersche Wunderland, in dem 
 neue S\"atze wie Pilze aus dem Boden sprossen. Ich war fast geblendet von der Vielzahl der Talente,
 die um mich umherschwirrten und fast w\"ochentlich Neues zu berichten hatten. Auch die Spannweite
 der erzielten Resultate war immens: In der Algebra waren es etwa gruppentheoretische Eigenschaften,
 Engelsche Elemente, Faktorisierung von Gruppen, distributive Quasigruppen, Erweiterungen abelscher
 Gruppen und ringtheoretische Radikale, in der Geometrie waren es endliche und topologische Ebenen,
 M\"obius- und Laguerregeometrie und gruppentheoretische Methoden in der Geometrie,
 die die Gem\"uter erhitzten.

 Ich f\"uhlte mich in diesem mathematischen Karpfenteich zur Geometrie hingezogen, weil ich mir
 einbildete, die Geometrie ist etwas Handfestes, Greifbares und bei meinen geringen mathematischen
 Kenntnissen eher Zug\"anglicheres. Dass ich heute vor Ihnen stehen darf, ist nicht mein Verdienst,
 sondern das von Helmut Salzmann, der mich behutsam an die Hand nahm und mir nicht nur die Geometrie,
 ja die Mathematik erschloss.

 Dass die Geometrie in Frankfurt in der Zeit des Baerschen Wirkens sich zu prachtvoller Bl\"ute entfalten
 konnte, war der faszinierenden,  einnehmenden, junge Mathematiker begeisternden und sie ermutigenden
 Pers\"onlichkeit von Reinhold Baer zu verdanken, der von 1940 
 bis 1952 selbst grundlegende Arbeiten \"uber endliche projektive Ebenen, man denke etwa an den von daher
 stammenden Begriff einer Baer-Unterebene, und ein Buch \anff Linear Algebra and Projective Geometry``
 verfasst hat.
 Die Geometrie lag Reinhold Baer lebenslang am Herzen; davon zeugt auch diese Tagungsreihe, die er zusammen
 mit Herrn Pickert begr\"undet hat und die auch heute seinen Namen tr\"agt. Ein weiterer gl\"ucklicher
 Umstand f\"ur die Geometrie war Ruth Moufang, die Kollegin von Reinhold Baer war und ihre Sch\"uler ebenfalls
 f\"ur die Grundlagen der Geometrie und konvexe K\"orper einzunehmen wusste. Die Symbiose zwischen dem
 Baerschen und Moufangschen Kreis war so perfekt, dass f\"ur mich Peter Dembowski ein fester Bestandteil
 der Baerschen Geometriegruppe war, obgleich sein Talent wohl von Frau Moufang entdeckt worden war.

 Die Baersche Geometriegruppe, die sich der endlichen Geometrie, ihren kombinatorischen Aspekten sowie
 ihren Verbindungen zur Gruppentheorie verschrieben hatte, hatte nach meiner damaligen  \"Uberzeugung,
 und daran hat sich bis heute nichts ge\"andert, drei tragende S\"aulen: Peter Dembowski, Heinz L\"uneburg
 und Christoph Hering. Obgleich die zeitlichen Unterschiede beim Einstieg in mathematische Publikationst\"atigkeit
 dieser drei Pilaster
 aus heutiger Sicht vernachl\"assigbar erscheinen, die erste Arbeit von Peter Dembowski erschien
 1958, die  von Heinz L\"uneburg 1960 und die von Christoph Hering 1963, f\"ur mich als Studenten und
 Bewunderer geh\"orten sie drei verschiedenen Ge\-ne\-ra\-ti\-on\-en an. 

 Peter Dembowski war f\"ur mich der abgekl\"arte Wissenschaftler, der aus dem Nichts, d.h.\ ohne
 kompliziertere Theorien benutzend, bleibende markante Resultate hervorzaubern konnte. Als hervorstechendstes
 Beispiel sei hier an seine Habilitationsschrift \"uber M\"obiusebenen gerader Ordnung erinnert, in der
 bewiesen wird, dass jede endliche M\"obiusebene, in der jeder Kreis eine ungerade Anzahl von Punkten
 tr\"agt, ovoidal ist; dies geschieht dadurch, dass er durch raffinierte Abz\"ahlk\"unste den
 dreidimensionalen Raum mit Hilfe der Winternitzschen Axiome erschafft und die M\"obiusebene als ein
 Ovoid hineinlegt. Durch diese Arbeit hat insbesondere die Klassifikation von Ovoiden in projektiven
 R\"aumen \"uber Galoisfeldern gerader Charakteristik an Bedeutung gewonnen und dazu gef\"uhrt, dass
 man heutzutage intensiv sogar den Computer einsetzt, um au\ss{}er den elliptischen und den Tits--Ovoiden
 neue Ovoide zu entdecken. Bis jetzt ist man zu projektiven R\"aumen \"uber Galoisfeldern der Ordnung
 $2^n$ mit $n\leq 5$ vorgedrungen, hat aber nach meiner Kenntnis leider keine neuen Ovoide gefunden.
 Peter Dembowski war und ist f\"ur mich derjenige der drei Frankfurter Sterne der diskreten Geometrie,
 der es liebte, endliche Geometrien in einen kombinatorischen Kontext zu r\"ucken.

 Christoph Hering empfand ich als ein jugendliches Nachwuchsgenie, der Tolles leistet und durch
 seinen augenzwinkernden Humor es verhindert, als Ze\-le\-bri\-t\"at wahrgenommen zu werden. Unvergessen bleiben
 mir seine Vortr\"age  \"uber die Lenz--Barlotti--Klassifikation von M\"obiusebenen, heutzutage
 Hering Klassifikation genannt, bei den Kindertagungen in Oberwolfach. Obgleich er bis Mitte der achtziger
 Jahre im Geiste der Frankfurter endlichen Geometrie gearbeitet und publiziert hat, empfinde ich ihn als
 einen Geometer, der den endlichen Gruppen, gesehen als Automorphismengruppen, sein unbedingtes
 Interesse und seine Phantasie geschenkt hat.

 Heinz L\"uneburg altersm\"a\ss{}ig zwischen den beiden, Dembowski und Hering, stehend war f\"ur mich die
 Achse der Frankfurter Forschung \"uber nicht\-to\-po\-lo\-gi\-sche Geometrie. In seinen Arbeiten \"uber endliche
 Geometrie bewegte er sich \"aquidistant zwischen Kombinatorik und Gruppentheorie und klopfte au\ss{}erdem
 die damaligen Str\"omungen der geometrischen Forschung auf ihre Entwicklungsf\"ahigkeit ab. Davon zeugen
 seine Arbeiten \"uber $\lambda$--R\"aume, Hjelmslev--Ebenen, Blockpl\"ane, insbesondere Steinersche 
 Tripel- und Quadrupelsysteme, Fundamentals\"atze der projektiven Geometrie, die Rolle der
 Zentral- und Axialkollineationen sowie \"uber elliptische Ebenen. 

 Die erste gro\ss e Leidenschaft von Heinz L\"uneburg, der 1956 sein Studium in Frankfurt begann, also in dem
 Jahr, in dem Reinhold Baer aus Amerika zur\"uckkehrend einen Lehrstuhl in Frankfurt akzeptierte, galt der
 kleinen Reidemeisterbedingung. Die Reidemeisterbedingung ist ein Schlie\ss{}ungssatz, der in der Theorie der
 Loops und der zugeh\"origen Netze die Assoziativit\"at garantiert und am bequemsten mit Hilfe von
 Projektivit\"aten in 3-Netzen formulierbar ist. Eine Spezialisierung dieses Schlie\ss{}ungssatzes war f\"ur
 Andrew Gleason der Angelpunkt beim Beweis des Satzes, dass eine endliche Fano--Ebene, d.h.\ eine endliche Ebene,
 in der die Diagonalpunkte eines jeden Vierecks kollinear sind, pappossch sein  muss. Dass man, um dieses Resultat
 zu erreichen, die Kollinearit\"at der Diagonalpunkte eines jeden Vierecks wirklich verlangen muss, sieht man
 bereits in der projektiven Ebene \"uber einem Fastk\"orper der Ordnung 9, denn in dieser existieren gewisse
 Vierecke mit kollinearen Diagonalpunkten. Die Gleasonsche Arbeit muss im Baerschen Seminar heftig studiert
 worden sein und Heinz L\"uneburg zu der Vermutung gef\"uhrt haben, dass jede endliche projektive Ebene,
 in der die kleine Reidemeisterbedingung gilt, bereits desarguessch sein muss. Er hat diese Vermutung in
 drei Arbeiten eindrucksvoll best\"atigt, wobei er in der zweiten Arbeit sein Ziel bis auf eine Ausnahme,
 n\"amlich die der Ebenen der Ordnung 60 erreicht hat. Die Schlachtung dieser Ebenen gelang ihm zusammen
 mit Otto Kegel, dem Meister der Faktorisierung, mit Hilfe des folgenden Satzes: Ist $G =A\cdot B$ Produkt
 zweier echter Untergruppen $A$ und $B$ und sind sowohl $A$ als auch $B$ isomorph zur alternierenden Gruppe
 des Grades 5, so ist $G$ entweder das direkte Produkt $A\times B$ oder die alternierende Gruppe des
 Grades 6. Diese gemeinsame Arbeit von Heinz L\"uneburg mit Otto Kegel gibt ein beredtes Zeugnis davon,
 dass die an verschiedenen Beeten des Baerschen Forschungsgartens T\"atigen miteinander in enger
 Kommunikation standen und einander unterst\"utzten.

 Innerhalb der Gruppentheorie sind es die Suzuki--Gruppen, die sich zu Lieb\-lings\-ob\-jek\-ten von Heinz
 L\"uneburg entwickelten  und ihn jahrelang treu be\-glei\-te\-ten. Ausgangspunkt f\"ur diese Zuneigung waren
 die von Dembowski gepflegten endlichen M\"obiusebenen. Heinz L\"uneburg beweist 1964, dass eine endliche
 M\"obiusebene $M$, deren Automorphismengruppe auf den Punkten von $M$ zwei\-fach transitiv ist, wobei aber nur
 die Identit\"at drei verschiedene Fixpunkte hat, entweder miquelsch oder die Geometrie der ebenen
 Schnitte eines Tits--Ovoids ist. Sp\"ater zeigt er, dass sich die gleiche Folgerung ergibt auch
 f\"ur kreishomogene endliche M\"obiusebenen oder f\"ur endliche  M\"obiusebenen gerader Ordnung, die eine
 transitive Automorphismengruppe gerader Ordnung gestatten. Da die Automorphismengruppen von
 Tits--Ovoiden Suzuki--Gruppen sind, ist die Begegnung von Heinz L\"uneburg mit diesen Gruppen
 unumg\"anglich. Und nachdem er sich mit ihnen eingehender befasst hat, sieht er, dass sie f\"ur ihn
 eine Br\"ucke zu den Translationsebenen schlagen, indem er den folgenden Satz
 beweist: Ist $q = 2^{2r+1}\geq 8$, so gibt es genau eine Translationsebene, 
 die die Ordnung $q^2$ hat und eine zur Suzuki--Gruppe 
 der Ordnung $(q^2 + 1)q^2(q - 1)$ isomorphe Kollineationsgruppe besitzt; diese Ebene ist nicht desarguessch. 

 Die Translationsebenen waren in Frankfurt wohlbekannt, denn T.\ G.\ Ostrom, dessen Lecture
 notes \anff Finite translation planes`` die erste zusammenfassende Darstellung diese Gebietes war,
 besuchte h\"aufig den Baerschen Kreis in Frankfurt und seine Vortr\"age \"uber
 replaceable nets klingen mir noch heute in den Ohren. So war Heinz L\"uneburg bestens \"uber
 den Forschungsstand bez\"uglich endlicher Translationsebenen informiert und konnte diese Ebenen,
 die seinen Namen tragen, in einer gro\ss en Arbeit f\"ur die Hamburger Abhandlungen 
 in einen gr\"o\ss eren Rahmen der endlichen projektiven Ebenen einbetten, deren von den Elationen
 erzeugte Kollineationsgruppe Punktstabilisatoren vorgeschriebener Struktur hat. Eine sch\"one
 Charakterisierung der L\"uneburgebenen hat 1972 Liebler gegeben: Die L\"uneburgebenen sind
 diejenigen affinen Translationsebenen, bei denen eine Gruppe $G$ von Kollineationen als Gruppe
 vom Rang drei auf der Menge der eigentlichen Punkte (d.h.\ der Stabilisator jedes eigentlichen
 Punktes in $G$ hat drei Bahnen) und als Gruppe vom Rang zwei auf der Menge der uneigentlichen Punkte operiert.

 Den Kulminationspunkt L\"uneburgscher Forschungen \"uber endliche Translationsebenen stellt
 seine 1980 erschienene Monographie \anff Translation planes`` dar. Nachdem er in ihr die n\"otigen
 Grundlagen zusammengestellt hatte, diskutiert er  permutationstheoretische  Kriterien, die es
 sicherstellen, dass eine endliche projektive oder affine Ebene, auf der eine Gruppe $G$ von
 Kollineationen operiert, eine Translationsebene ist und $G$ die Translationsgruppe enth\"alt.
 Von den dort behandelten Kriterien von Ascher Wagner, Kallaher, Ostrom sei eine Charakterisierung
 von Heinz L\"uneburg, die f\"ur seine Ebenen zutrifft, besonders erw\"ahnt: Ist $A$ eine endliche affine
 Ebene und $G$ eine Gruppe von Kollineationen von $A$, so dass der Stabilisator jeder Geraden von $A$
 in $G$ auf den Punkten dieser Geraden zweifach transitiv ist, so ist $A$ eine Translationsebene
 und $G$ enth\"alt die Translationsgruppe. 

 Der L\"owenanteil der L\"uneburgschen Monographie ist jedoch speziellen, signifikanten
 Translationsebenen gewidmet, die in keinen gro\ss en Familien leben, sondern eher ein
 Einsiedlerdasein f\"uhren. Von den von Heinz L\"uneburg  entdeckten und zuerst in Geometriae
 Dedicata publizierten Translationsebenen, die er merkw\"urdige Translationsebenen nennt, gibt es
 zwei Typen, die sich in der Struktur gewisser Elationsgruppen unterscheiden; von dem einem Typ
 gibt es 8 nichtisomorphe, von dem andern Typ gibt es 14 nichtisomorphe Ebenen. Auch im Kapitel VII,
 welches der Klassifikation von Translationsebenen der Ordnung $q^2$ gewidmet ist, die den K\"orper
 $\GF(q)$ in ihrem Kern enthalten und eine zu $\SL_2(q)$ isomorphe Kollineationsgruppe gestatten,
 ist das Hauptaugenmerk auf die Hering-- und Sch\"affer--Ebenen gerichtet. 

 Dass sich Heinz L\"uneburg in den achtziger Jahren den endlichen Translationsebenen entfremdet hat,
 kann man verstehen, wenn man in das neueste, 2007 erschienene und 888 Seiten umfassende
 \anff Handbook of finite translation planes`` von Norman Johnson, Vikram Jha und Mauro Biliotti
 hineinschaut, in dem alle bisher bekannten endlichen Translationsebenen gesammelt sind.  
 Um in diesem Karpfenteich auch heutzutage noch erfolgreich zu fischen, braucht man nicht nur
 Ausdauer und Geduld, sondern auch eine Buchhalterveranlagung, die Heinz L\"uneburg wohl abging.

 Neben den Suzuki--Gruppen galt L\"uneburgs Interesse auch anderen extravaganten endlichen Gruppen,
 so etwa den Ree--Gruppen und  den Mathieu--Grup\-pen sowie deren Darstellungen als Automorphismengruppen
 geeigneter Block\-pl\"ane.

 L\"uneburgsche Arbeiten, in denen keine abschlie\ss enden Antworten der behandelten Probleme gegeben
 wurden, enthielten stets eine F\"ulle  von Anregungen und Methoden, die die Leser dieser Arbeiten
 zum Weiterforschen anregten. So haben seine Arbeiten zur Existenz der endlichen projektiven Ebenen
 vom Lenz--Barlotti--Typ I.6 und III.2  Hering und Kantor inspiriert zu beweisen, dass es weder eine
 endliche projektive Ebene vom Lenz--Barlotti--Typ I.6 noch eine Ebene vom Lenz--Typ III gibt.
 Ich glaube, dass bis heute ungekl\"art ist, ob es unendliche projektive Ebenen vom Lenz--Barlotti--Typ
 I.6 gibt. Der einzige andere offene Fall scheinen die endlichen projektiven Ebenen vom
 Lenz--Barlotti--Typ II.1 zu sein. Meine Quelle f\"ur diese Behauptungen ist:
 Gina Ghinelli and Francesca Merola, Lenz--Barlotti--Classification and related open
 problems: an update, Roma 2005, Quaderni Elettronici del Seminario di Geometria Combinatoria (Quaderno 20).

 Obgleich ich mich nie direkt
 im unmittelbaren Umkreis L\"uneburgscher For\-schung\-en aufgehalten habe,
 gibt es zwei Themen, wo sich unsere Interessen ber\"uhrt haben.  Einmal ist es die
 Funktionalgleichung $f(x + y f(x)) = f(x)f(y)$ von Golab und Schinzel, die die  
 Beziehungen zwischen Komplementen eines semidirekten Produkts regelt. Peter Plaumann und
 ich haben diese Funktionalgleichung \"uber den $p$-adischen Zahlen traktiert, Heinz L\"uneburg
 und Peter Plaumann haben sie \"uber den Galoisfeldern behandelt. Das zweite Thema, das 
 Heinz L\"uneburg und mich faszinierte, sind die Gruppen der Projektivit\"aten eines Blocks auf
 sich in verschiedenen Geometrien, also der von Staudtsche Standpunkt.

 Nachdem Carl Georg Christian von Staudt Mitte des neunzehnten Jahrhunderts die Gruppe $\Pi$
 der Projektivit\"aten einer Geraden auf sich in einer pap\-pos\-schen projektiven Ebene $E$ betrachtet
 und bewiesen hatte, dass sie scharf dreifach transitiv ist, und nachdem mit Hessenberg klar war,
 dass $E$ genau dann pappossch ist, wenn $\Pi$ scharf dreifach transitiv und dann isomorph zu $\PGL(2, K)$
 ist, hat die Sch\"onheit dieses Satzes die Geometer so eingelullt, dass die naheliegende Frage, wie
 es um die Gruppe $\Pi$ der Projektivit\"aten einer Ge\-ra\-den auf sich in anderen, insbesondere
 nichtdesarguesschen projektiven Ebenen steht, lange nicht gestellt wurde. Erst 1959 hat sich
 Adriano Barlotti diesem Problem gestellt und bewiesen, dass die Gruppe $\Pi$ der Projektivit\"aten
 in den drei nichtdesarguesschen Ebenen der Ordnung 9 die symmetrische Gruppe des Grades 10 und
 in der Hallebene der Ordnung 16 die alternierende Gruppe des Grades 17 ist. Damit war die Hatz
 auf die Gruppe $\Pi$ der Projektivit\"aten in nichtdesarguesschen endlichen Ebenen er\"offnet.
 Nachdem sie in vielen endlichen projektiven Ebenen bestimmt worden war, hatte sich schnell
 die Vermutung herauskristallisiert, dass die Gruppe $\Pi$ der Projektivit\"aten in endlichen
 nichtdesarguesschen projektiven Ebenen der Ordnung $n$ stets die alternierende Gruppe des
 Grades $n + 1$  ent\-hal\-ten muss. Da jedoch in den sechziger Jahren die Klassifikation der
 endlichen einfachen Gruppen noch ausstand, trotzte diese Vermutung jedem Angriff.

 Obgleich Heinz L\"uneburg in den Beweisen seiner Arbeiten ausgiebig Projektivit\"aten
 zwischen Geraden benutzt hatte, trat er mit der Gruppe der Projektivit\"aten an die \"Offentlichkeit erst
 mit seiner 1967 ver\"offentlichten Arbeit, in der die desarguesschen affinen Ebenen als diejenigen
 affinen  Ebenen gekennzeichnet werden, in denen es zu je drei verschiedenen Punkten $P$, $Q$, $R$ ein
 (axiomatisch definiertes) Teilverh\"altnis $r$ gibt, so dass $P r Q = R$ gilt. Der Satz, in dem
 die Gruppe der Projektivit\"aten dabei mit voller Wucht auftritt, ist das folgende Nebenergebnis
 der Arbeit: Eine endliche projektive Ebene ungerader Ordnung $q$ ist genau dann desarguessch,
 wenn der Stabilisator eines Punktes in der Gruppe der Projektivit\"aten einer Geraden auf sich
 einen Normalteiler der Ordnung $q$ enth\"alt. Au\ss erdem wird einem bei der Durchsicht dieser
 Arbeit klar, dass der Verfasser wusste, man m\"usse mit gleichem Nachdruck wie die Gruppe $\Pi$
 der Projektivit\"aten einer projektiven Ebene in affinen Ebenen auch die Gruppe der affinen
 Projektivit\"aten einer affinen Geraden auf sich studieren, deren Elemente Produkte affiner
 Parallelperspektivit\"aten sind. 

 Die Besch\"aftigung von Heinz L\"uneburg mit den Gruppen von Projektivit\"aten erreicht
 ihren H\"ohepunkt in dem Bericht \anff Some new results on groups of projectivities``, welcher
 die Frucht seiner Vortr\"age bei der Bad Windsheimer Tagung \anff Geometry -- von Staudt's Point of View`` ist.
 In diesem Beitrag wird auch 
 \"uber die Ergebnisse seines Sch\"ulers Theo Grundh\"ofer, insbesondere bez\"uglich der Gruppe
 der affinen Projektivit\"aten, ausf\"uhrlich berichtet.

 Da inzwischen eine Klassifikation der endlichen einfachen Gruppen vorgelegen hatte, war Theo
 Grundh\"ofer in der Lage, die Vermutung \"uber die  Gruppen $\Pi$ der Projektivit\"aten in
 endlichen nichtdesarguesschen projektiven Ebenen fast vollst\"andig zu beweisen:  Die Gruppe $\Pi$
 der Projektivit\"aten einer endlichen nichtdesarguesschen projektiven Ebene der Ordnung $n$ enth\"alt
 entweder die alternierende Gruppe des Grades $n + 1$ oder es ist $n = 23$ und $\Pi$ ist die
 gr\"o\ss te Mathieugruppe $M_{24}$. Dass die zweite M\"oglichkeit ausgeschlossen ist, haben
 neuerdings Peter M\"uller und Gabor Nagy gezeigt. Eine analoge Situation liegt nach Theo
 Grundh\"ofer, Peter M\"uller und Gabor Nagy auch f\"ur die Gruppe $\Sigma$ der affinen Projektivit\"aten
 einer Geraden auf sich in einer affinen Ebene der Ordnung $n$ vor, welche keine Translationsebene ist:
 $\Sigma$
 enth\"alt entweder die alternierende Gruppe der Ordnung $n$ oder, was nat\"urlich unwahrscheinlich
 ist, es ist $n=24$ und die Gruppe $\Sigma$ ist die Mathieu--Gruppe des Grades
 24. F\"ur Translationsebenen der Ordnung $q^d$ mit dem Kern $\GF(q)$ hat
 bereits in den achtziger Jahren Theo Grundh\"ofer bewiesen, dass die Gruppe der affinen
 Projektivit\"aten eine Gruppe von affinen Abbildungen ist, so dass der Stabilisator eines
 Punktes zwischen $\SL(d, q)$ und $\GL(d, q)$ liegt.

 Heinz L\"uneburg war f\"ur mich ein leidenschaftlicher Mathematiker, dessen Vortr\"age und
 Vorlesungen so sorgf\"altig vorbereitet waren, dass sie  sofort zum Druck gehen konnten.
 Auf diese Weise sind die folgenden seiner B\"ucher entstanden: \anff Kombinatorik``, erschienen 1971,
 \anff Einf\"uhrung in die Algebra``, erschienen 1973, \anff Vorlesungen \"uber Zahlentheorie``,
 erschienen  1978,
 \anff Galoisfelder, Kreisteilungsk\"orper und Schieberegisterfolgen``, erschienen 1979,
 \anff Vorlesungen \"uber Analysis``, erschienen 1981, \anff On the rational normal form of
 endomorphisms. A primer to 
 constructive algebra`` und \anff Kleine Fibel der Arithmetik``, beide erschienen 1987,
 \anff Tools and fundamental constructions of combinatorial mathematics``, erschienen 1989,
 \anff Vorlesungen \"uber lineare Algebra.  Versehen mit der zu ihrem Verst\"andnis n\"otigen
 Algebra sowie einigen Bemerkungen zu ihrer Didaktik``, erschienen 1993,
 \anff Die euklidische Ebene und ihre Verwandten`` und \anff Gruppen, Ringe, K\"orper.
 Die grundlegenden Strukturen der Algebra``, beide erschienen 1999 sowie
 \anff Rekursive Funktionen``, erschienen 2002. Diese f\"ur die Studierenden bestimmten Werke
 wollen keine Enzyklop\"adien \"uber die in ihnen  behandelten Gebiete sein, sondern spiegeln
 die pers\"onliche Sicht des Verfassers wider, welche Sachverhalte man mehr gewichten soll
 und was vernachl\"assigbar ist, weil es in anderen B\"uchern schon hundertmal gesagt worden ist. 
 In ihnen erweist sich Heinz L\"uneburg als ein begnadeter Didaktiker, der nicht durch Beiwerk,
 sondern durch die Klarheit des Gedankens und durch strenge Be\-weis\-f\"uh\-rung \"uberzeugt. Dies
 konnte man \"ubrigens bei jedem seiner Vortr\"age bewundern. Obwohl voll innerer Begeisterung
 und in der \"Uberzeugung Wertvolles beizutragen, stand sein schn\"orkelloser direkter Vortragsstil
 im Gegensatz zu dem derjenigen Mathematiker, die in ihren Vortr\"agen an einer Wagneroper stricken.

 Heinz L\"uneburg war f\"ur mich ein Freund, der seinen Prinzipien ein Leben lang treu blieb und
 dessen erfolgreiche Karriere ihn weder verbogen noch vom t\"aglichen Leben separiert hat. Wenn er
 einmal nach reifer \"Uberlegung zu einem Urteil gekommen ist, war er davon durch 
 Opportunit\"atsgr\"unde nicht abzubringen. Dies hat ihm, dem Bildungswertkonservativen,
 in den Fakult\"atssitzungen, in denen h\"aufig nach    
 Anpassung an Umbr\"uche verlangende politische Vorlagen gehechelt wird, sicherlich nicht
 nur Beifall beschert.

 Heinz L\"uneburg hatte ein sensibles Ohr f\"ur Ver\"anderungen in der Mathematik. Er suchte
 nach zeitloser Mathematik, nahm jedoch wahr, dass seit den achtziger Jahren durch den steigenden
 Publikationszwang bei den Mathematikern die Lust abnahm, einzelne Texte zu studieren, und wie im
 Sport der Wett\-be\-werbs\-charakter in den Vordergrund r\"uckte. Heutzutage z\"ahlt nicht nur der Gehalt
 einer mathematischen Arbeit, sondern auch im gro\ss en Ma\ss e in welch aggressivem Umfeld sie entstanden
 ist und ob sie ein  Problem l\"ost, an welchem sich schon einige mehr oder minder prominente
 Mathematiker die Z\"ahne ausgebissen haben. Der Kampf der Teilgebiete  der Mathematik um
 die Oberhoheit in Zeitschriften mit hohem impact factor, der eine scheinbare Objektivierung
 mathematischer Leistung erlauben soll, ist voll entbrannt.

 Heinz L\"uneburg war gegen\"uber den modernen Entwicklungen der Mathematik aufgeschlossen, wollte
 aber nicht in Abstraktionen der Schemata, Moduli und Motive einsteigen, deren Wolkenhaftigkeit
 man \"uberwinden muss, wenn man sich bei ihnen wohlf\"uhlen will. Als Geometer ist man letztlich
 Greifbares gewohnt. Dar\"uber hinaus hat ihn die Entwicklung elektronischer Rechenanlagen
 fasziniert und so hat er sich den Algorithmen zugewandt, mit denen man effektiv elementare
 zahlentheoretische Funktionen berechnen kann, wie etwa den gr\"o\ss ten gemeinsamen Teiler
 bzw. das kleinste gemeinsame Vielfache zweier ganzer Zahlen, die Einheiten im Ring der
 ganzen Zahlen modulo $n$, die Primitivwurzeln modulo einer Primzahl, oder die gr\"o\ss te
 nat\"urliche Zahl, die die Quadratwurzel einer gegebenen nat\"urlichen Zahl nicht \"ubertrifft.
 Er zeigte, wie man auch f\"ur inneralgebraische Fragestellungen, wie etwa die effektive
 Konstruktion der algebraischen Erweiterungen der Primk\"orper $\GF(p)$ sowie der vollst\"andigen
 Kreisteilungsk\"orper, Algorithmen einsetzen kann.

 Heinz L\"uneburg war der festen \"Uberzeugung, dass der Computer ein Diener der Mathematiker
 bleiben muss (so etwa bei Schleifeninvarianten, von denen er behauptet, sie seien ihm fast
 zu einer fixen Idee geworden) und sich nicht als Herr 
 etablieren darf.
 Da er wohl mit der Zeit merkte, dass diese Tendenz dem Zeitgeist immanent war, hat er sich, anstatt
 unn\"utze Scharm\"utzel zu f\"uhren, aufgemacht, nach den Wurzeln der Mathematik zu suchen.
 Daf\"ur war er bestens pr\"apariert: Als Sch\"uler des bekannten Frankfurter humanistischen
 Lessing--Gymnasiums, des Altgriechischen wie des Lateinischen 
 m\"achtig,
 stand er fest und mit voller \"Uberzeugung zu den Grunds\"atzen des Abendlandes.  Nicht nur im
 Deutschen, sondern auch im Eng\-li\-schen, Franz\"osischen und Italienischen zu Hause, war er ein
 Verfechter einer abendl\"andischen Idee, welche nicht nur aus der anglo-s\"achsischen Monokultur
 ihre Kraft bezieht. 

 Der gro\ss artige Einstieg von Heinz L\"uneburg bei seiner Suche nach den ge\-schicht\-lich\-en Wurzeln
 der Mathematik war das Ergebnis eines akribischen Studiums des Liber Abacci von Leonardo Pisano.
 In dem daraus entstandenen Buch \anff Leonardi Pisani Liber Abacci oder Lesevergn\"ugen eines
 Ma\-the\-ma\-ti\-kers``, das zwei Auflagen erlebte, wird nicht nur eine Interpretation des Textes von
 Leonardo Pisano in der Sprache der uns gel\"aufigen Mathematik gegeben, sondern auch ein Bild
 der Zeit ausgebreitet, in der Leonardo Pisano lebte und wirkte. Es ist keine kritische Edition,
 aber Heinz L\"uneburg findet hier seine in allen sp\"ateren historischen Arbeiten angewandte
 Methode, den Lesern die Ergebnisse aus der Geschichte der Mathematik in modernem Gewand nahezubringen.
 So l\"asst er etwa Euler, Gau\ss{} und Lagrange von K\"orpern reden, obwohl das Konzept eines K\"orpers
 erst von Steinitz 1910 entwickelt wurde. Ebenso zeigt sich  bereits hier, auch durch seine
 etymologischen Kenntnisse bef\"ordert, dass Heinz L\"uneburg zum Urgrund der Mathematik, zur Zahl,
 zur Ziffer und und deren Herkunft  vordringen will; nat\"urlich spielt dabei die Null,
 urspr\"unglich as-sifr, latinisiert ciffra genannt, eine besondere Rolle.

\goodbreak

 Mir hat die Lekt\"ure des ersten Kapitels von L\"uneburgs Leonardi Pisani Liber Abacci
 neue \"uberraschende Einsichten gebracht. Von der Geschichte der Mathematik unbeleckt, war mir
 zwar klar, dass die Griechen, indem sie die Nabelschnur zwischen Wirklichkeit und insbesondere
 den Anwendungen durchgeschnitten hatten, die Geburtshelfer der Mathematik im heutigen Sinn
 sind. Doch von den R\"omern dachte ich, dass die Eroberung der Welt sie so in Anspruch nahm, dass
 keine Zeit zum mathematischen Spintisieren \"ubrig blieb. Und von der katholischen Kirche des
 Mittelalters hatte ich den Eindruck der Gleich\-g\"ul\-tig\-keit gegen\"uber der Mathematik.
 Meine Meinung war, dass die Mathematik nach den Griechen ihr Dasein in arabischem und j\"udischem
 Umfeld, etwa in Spanien, gefristet hatte, ehe sie von der Renaissance wiedergeboren wurde.
 Heinz L\"uneburg hat mich belehrt, dass diese naive Vorstellung falsch ist. Das r\"omische Reich
 hatte ein gut ausgebautes Bildungssystem, in den Kl\"ostern wurden griechische Schriften studiert
 und es gab hohe kirchliche W\"urdentr\"ager, die sich um den Erhalt des naturwissenschaftlichen
 Wissens der Griechen sorgten. Und das Buch Liber Abacci selbst gibt Zeugnis davon, dass der
 Fluss der Mathematik in keiner Periode des Mittelalters ausgetrocknet ist. 

 Die zweite \"Uberraschung, die mir Heinz L\"uneburg mit seinem Buch bescherte, betrifft Palermo,
 wo ich jetzt jedes Jahr einige Zeit zubringe.  Mir erschien die Eroberung von Palermo durch
 die Normannen im positiven Licht, wobei ich diese Ansicht aus der Solidit\"at des Palazzo dei
 Normanni, der Sch\"onheit der Capella Palatina und der imponierenden Pracht des Doms von Monreale bezog.
 Doch bei Heinz L\"uneburg kann man nachlesen, dass die Pisaner, die die Normannen bei deren Kampf
 gegen die Sarazener unterst\"utzten, Palermo nach dessen Eroberung so t\"uchtig pl\"underten, dass
 sie aus dem Erl\"os der Beute in Pisa den von den  Touristen bestaunten Dom auf der Piazza dei miracoli
 erbauen konnten. Meine Meinung \"uber Kaiser Friedrich den Zweiten, an dessen Grab ich in der
 Palermitaner Kathedrale ab und zu vorbeischlendre, hat Heinz L\"uneburg vollauf
 best\"atigt: Friedrich der Zweite war nicht nur ein ausgewiesener Experte der Falknerei,
 sondern auch ein entschiedener F\"orderer der Wissenschaften und Begr\"under der Universit\"at Neapel.  

 Das letzte von Heinz L\"uneburg erschienene Werk sind die beiden B\"ande betitelt \anff Von Zahlen und
 Gr\"o\ss en. Dritthalbtausend Jahre Theorie und Praxis``, die der Geschichte der algebraischen
 Gleichungen und der K\"orper gewidmet sind. Sie sind, wie bei L\"uneburg \"ublich, in der Sprache
 der heutigen Mathematik verfasst, enthalten aber unz\"ahlige historische und etymologische Ausfl\"uge,
 nat\"urlich auch in Latein und Italienisch, und einige Bemerkungen, die zeigen, dass der Autor
 der verlorenen klassischen Bildung der Jugend und der gerade jetzt vor unseren Augen sterbenden
 Humboldtschen Universit\"at nachtrauert. 

 Heinz L\"uneburg hat sich selbst nicht als Mathematikhistoriker betrachtet. Er sagt, dass er bei
 einem Gegenstand mehr in die Tiefe, w\"ahrend der Historiker mehr in die Breite geht. Ich m\"ochte
 hier aus den Vorlesungen zur Geschichte der Mathematik von Hans Wussing, dem auch heute noch
 ber\"uhmten Ma\-the\-ma\-tik\-his\-to\-ri\-ker aus Leipzig, zitieren, um Ihnen die Entscheidung
 zu er\-leich\-tern,
 ob es unbedingt eine Ehre sein muss, zur professionellen Zunft der Wissenschaftshistoriker zu
 geh\"oren. Wussing f\"uhrt auf Seite 18 aus: Die marxistische Historiographie der Mathematik
 beruht methodologisch auf dem historischen und dialektischen Materialismus. Danach ist jede
 Wissenschaft eine gesellschaftliche Erscheinung. Auch die Mathematik ist eine spezifische Form
 des gesellschaftlichen Bewusstseins. Sie ist mehr als das Ergebnis von Kenntnissen und Erkenntnissen,
 von Theorien und Methoden; sie ist zugleich geformt von den materiellen und ideellen Interessen der
 jeweils herrschenden Klassen. 

 Heinz L\"uneburg wurde j\"ah aus dem Leben gerissen, er war bis zum letz\-ten Tag voller mathematischer
 Ideen und beabsichtigter Buchprojekte. Martin L\"uneburg und Theo Grundh\"ofer
 fanden denn auch im Nachlass drei Werke, die es lohnten publiziert zu werden.
 Zwei von ihnen waren von Heinz L\"uneburg so vollst\"andig bearbeitet, dass sie, versehen mit marginalen
 Bemerkungen der Herausgeber, zum Druck im Oldenbourg-Verlag angenommen wurden und im n\"achsten Jahr
 erscheinen werden. Das eine dieser zwei B\"ucher tr\"agt den Titel \anff Zahlentheorie``, hat zweihundert
 Seiten und zeugt davon, wie nachhaltig sich Heinz L\"uneburg mit der Entwicklung der Zahlentheorie
 besch\"aftigt hat. Mit den drei B\"uchern von Euklid startend, die die zahlentheoretischen hei\ss{}en,
 gelangt er zu den Ringen der ganzen algebraischen Zahlen und dem fermatschen Zwei--Quadrate--Satz.
 Im zweiten 220 Seiten umfassenden Buch mit dem Titel \anff Gr\"o\ss{}en und Zahlen. Ein Aufbau des
 Zahlensystems auf der Grundlage der eudoxischen Proportionenlehre`` l\"asst er von seinem
 fr\"uheren Vervollst\"andigen der rationalen Zahlen durch Cauchyfolgen ab und wendet sich den
 dedekindschen Schnitten auf der Menge der positiven rationalen Zahlen zu, wobei er anschlie\ss{}end
 die negativen reellen Zahlen auf die gleiche Weise gewinnt, wie man die negativen ganzen Zahlen
 aus den nat\"urlichen bekommt. So erh\"alt er gleichzeitig mit den reellen Zahlen alle
 Logarithmus- sowie alle Exponentialfunktionen.

 Das dritte hinterlassene Werk ist ein 108-seitiges Manuskript \anff Streifz\"uge durch die Geschichte der
 Mathematik``, 
 das besonders reizvoll ist, weil Heinz L\"uneburg darin in lockerer Form die im Laufe seines
 Lebens gewonnenen ge\-schicht\-lich\-en Erkenntnisse \"uber Zahlen versammelt, Erkenntnisse  
 von den na\-t\"ur\-lichen bis zu den rationalen Zahlen. Ich pers\"onlich bin am meisten von dem
 in dem Buch vorgesehen Vokabular gefesselt, in welchem Heinz L\"uneburg vorhatte, die Herkunft
 der von den Mathematikern meistbenutzten W\"orter blo\ss{}zulegen, denn das meiste war f\"ur
 mich \"uberraschend und neu und zeigte mir meinen Mangel an klassischer Bildung. Die Lekt\"ure
 dieses Vokabulars ist sehr vergn\"uglich, was hier an dem Wort Korollar erl\"autert sei: Das
 lateinische Corollarium ist das Kr\"anzchen, das der Gastgeber seinen Zechkumpanen aufs Haupt
 dr\"uckt, wenn sie zum Symposion (griechisch f\"ur Gelage) kommen. Leider ist dieses Manuskript
 nicht ganz vollst\"andig.

 Heinz L\"uneburg war ein Freund, auf den man sich verlassen konnte. Er hat mit seinen Meinungen
 nie hinterm Berg gehalten und sich seine Gradlinigkeit ein Leben lang bewahrt. Seine Mathematik
 bleibt in uns gut aufgehoben, so lange wir leben. Die L\"uneburgebenen bleiben aktuell, so lange
 man sich f\"ur projektive  Ebenen und Suzuki--Gruppen interessiert. Dass er sich aus dem Kraal
 der Ma\-the\-ma\-tik\-er hinausgewagt hat, wird sich f\"ur sein Gedenken in der Zukunft lohnen. W\"ahrend
 die wenigen kompetenten B\"ucher f\"urs brei\-te Publikum, die die Mathematik aus der Geschichte
 herauswachsen lassen, naive, unmittelbare Begeisterung entfachen  und sich jahrzehntelang
 eines Zu\-spruchs erfreuen k\"onnen, muss ein mathematischer Aufsatz au\ss{}erhalb des engsten
 Spezialistenkreises um die Aufmerksamkeit der in Anbetracht der Flut der Ver\"offentlichungen
 abgebr\"uhten Mathematikergilde hart und manchmal ohne Erfolg k\"ampfen.

 \markboth{}{Der Aufbruch der Geometrie}

\backmatter


\chapter*{Literatur}
\addcontentsline{toc}{chapter}{\protect\numberline{Literatur}}
\markboth{Literatur}{Literatur}

\frenchspacing

\parindent=0pt

\def\abs{\par\medskip\noindent\hangindent=12pt\hangafter=1}
\def\absabs{\par\noindent\hspace{0pt}\hangindent=12pt\hangafter=1}

\vspace*{13pt}
 
 \abs
       Johannes Andr\'e,
  \emph{\"Uber nicht-Desarguessche Ebenen mit transitiver Trans\-la\-tions\-grup\-pe.} Math.
 Z. 60, 156--186, 1954

 \abs
       Emil Artin,
 \emph{Geometric Algebra.} New York, London 1957
 
 \abs
       Reinhold Baer,
 \emph{Homogeneity of projective planes.} Amer. J. Math. 64, 137--152, 1942
 \absabs
 --- \emph{Polarities of finite projective planes.} Bull. Amer. Math. Soc. 53, 77--93,
 1946a
 \absabs
 --- \emph{Projectivities with fixed points on every line of the plane.} Bull. Amer. Math.
 Soc. 52, 273-- 286, 1946b
 \absabs
 --- \emph{Linear Algebra and Projective Geometry.} New York 1952 (Strukturs\"atze)

 \abs
       Richard Brauer,
  \emph{A Characterization of Null Systems in Projective Space.} Bull. Amer. Math.
 Soc. 42, 247--254, 1936 (Beweist u. A., dass jede Kollineation von einer
 semilinearen Abbildung induziert wird; vorausgesetzt wird jedoch, dass der
 K\"orper kommutativ ist. Benutzt Doppelverh\"altnisse. Der Satz sei bekannt,
 sein Beweis jedoch k\"urzer.)
 \absabs
 --- \emph{On the connection between the ordinary and the modular characters of groups
 of finite order.} Ann. Math. 42, 926--935, 1941 (Satz IIII.3.4)

 \abs
       Egbert Brieskorn und Horst Kn\"orrer,
  \emph{Algebraische Kurven.} Math. Institut Bonn, o. J. (Birkh\"auser 1981)
 
 \abs
       Richard H. Bruck and Herbert J. Ryser,
 \emph{The non-existence of certain finite projective planes.} Canadian J. Math. 1,
 88--93, 1949

 \abs
       W. L. Chow,
 \emph{On the Geometry of Algebraic Homogeneous Spaces.} Ann. Math. 50, 32--67,
 1949 (Satz von Chow f\"ur beliebige irreduzible projektive R\"aume, aber
 $r = s$.)
 
 \abs
       Paul M. Cohn,
 \emph{Skew field constructions.} Cambridge 1977\\(K\"orper mit $(K^*)' = K^*$.)
 
 \abs
       Richard Dedekind,
 \emph{\"Uber Zerlegungen von Zahlen durch ihre gr\"o\ss ten gemeinsamen Teiler.}
 Festschrift der TH Braunschweig. 1--40, 1897. Werke, Band 2, 103--147, 1931
 
 \abs
       Peter Dembowski,
 \emph{Verallgemeinerungen von Transitivit\"atsklassen endlicher projektiver
 Ebenen.} Math. Z. 69, 59--89, 1958
 \absabs
 --- \emph{Finite Geometries.} Ergebnisse der Mathematik und ihrer Grenzgebiete,
 Band 44. Berlin, etc. 1968
 
 \abs
       Ulrich Dempwolff,
 \emph{\"Uber die Determinante.} Math. Semesterberichte 40, 193--197, 1993
 
 \abs
       Max Deuring,
 \emph{Algebren.} Ergebnisse der Mathematik und ihrer Grenzgebiete, Band 41. 2.\ Auflage.
 Berlin, Heidelberg, New York 1968
 
 \abs
       Leonard Eugene Dickson,
 \emph{Linear Groups with an Exposition of the Galois Field Theory.} Nachdruck der
 ersten Auflage. New York 1958. In dieser Auflage ist das Vorwort zur ersten
 Auflage mit November 1900 datiert.
 
 \abs
       Jean Dieudonn\'e,
 \emph{Les d\'eterminants sur un corps non commutatif.} Bull. Soc. Math. France 71,
 27--45, 1943 (Zitiert nach Dieudonn\'e 1955)
 \absabs
 --- \emph{La g\'eom\'etrie des groupes classiques.} Ergebnisse der Mathematik und ihrer
 Grenzgebiete. Neue Folge Band 5. Berlin, etc. 1955
 
 \abs
       R. A. Fisher,
 \emph{An examination of the different possible solutions of a problem in
 incomplete blocks.} Ann. of Eugenics 10, 52--75, 1940 (Zitiert nach Dembowski
 1968. Fisher'sche Ungleichung.)
 
 \abs     G. Frobenius,
 \emph{\"Uber lineare Substitutionen und bilineare Formen.} J. reine angew. Math.
 84, 1--63, 1877
 

 \abs
       Carl Friedrich Gau\ss,
 \emph{Demonstration nova altera theorematis omnem functionem algebraicam
 rationalem integram unius variabilis in factores reales primi vel secundi
 gradus resolvi posse.} Comm. soc. reg. sci. Gottingensis. rec. III. G\"ottingen
 1816. Werke, Band 3, 31--56, 1876

 \abs
       William Rowan Hamilton,
 \emph{Memorandum respecting a new system of roots of unity.} Phil. Mag. (4),
 12, 446, 1856 (Satz III.4.5: Pr\"asentation von $A_5$.)

 \abs
       I. N. Herstein,
 \emph{Noncommutative Rings.} The Carus Mathematical Monographs 15. Published by
 The Mathematical Association of America. Kein Ort. 1968
 
 \abs
       Armin Herzer,
 \emph{Dualit\"aten mit zwei Geraden aus absoluten Punkten in projektiven Ebenen.}
 Math. Z. 129, 235--257, 1972

 \abs
       Otto Hesse,
 \emph{\"Uber die Elimination der Variabeln aus drei algebraischen Gleichungen vom
 zweiten Grade mit zwei Variabeln.} J. reine angew. Mathematik 28, 68--96, 1844
 \absabs
 --- \emph{Algebraische Aufl\"osung derjenigen Gleichungen 9ten Grades, deren Wurzeln
 die Eigenschaft haben, dass eine gegebene rationale und symmetrische Function
 $\theta (x_\lambda,x_\mu)$ je zweier Wurzeln $x_\lambda$, $x_\mu$ eine dritte
 Wurzel $x_\kappa$ giebt, so dass gleichzeitig: $x_\kappa = \theta(x_\lambda,
 x_\mu)$, $x_\lambda = \theta(x_\mu,x_\kappa)$, $x_\mu =
 \theta(x_\kappa,x_\lambda)$.} J. reine angew. Mathematik 34, 193--208, 1847

 \abs
       Gerhard Hessenberg,
 \emph{Beweis des Desarguesschen Satzes aus dem Pascalschen.} Math. Ann. 61,
 161--172, 1905

 \abs
       David Hilbert,
 \emph{Grundlagen der Geometrie.} Leipzig 1899

 \abs
       Otto H\"older,
 \emph{Die einfachen Gruppen im ersten und zweiten Hundert der Ordnungszahlen.}
 Math. Ann. 40, 55--88, 1892

 \abs
       Daniel R. Hughes und Fred C. Piper,
 \emph{Projective Planes.} New York, etc. 1973

 \abs
       Bertram Huppert,
 \emph{Endliche Gruppen I.} Berlin, etc. 1973

 \abs
       K. Iwasawa,
 \emph{\"Uber die Einfachheit der speziellen projektiven Gruppe.} Proc. Imp. Acad.
 Tokyo 17, 57--59, 1962 (Zitiert nach Huppert 1973)

 \abs
       John Jackson,
 \emph{Rational Amusements for Winter Evenings.} London 1821

 \abs
       Camille Jordan,
 \emph{Trait\'e des substitutions et des \'equations alg\'ebriques.} Nachdruck der
 Ausgabe Paris 1870. Sceaux 1989

 \abs
       Helmut Karzel,
 \emph{Zweiseitige Inzidenzgruppen.} Abh. Math. Sem. Univ. Hamburg 29, 118--136, 1965

 \abs
       Heinz L\"uneburg,
 \emph{Ein neuer Beweis eines Hauptsatzes der projektiven Geometrie.} Math. Z.
 87, 32--36. 1965
 \absabs
 --- \emph{\"Uber die Strukturs\"atze der projektiven Geometrie.} Arch. Math. 17,
 206--209, 1966
 \absabs
 --- \emph{Einige methodische Bemerkungen zur Theorie der elliptischen
 Bewegungsgruppen.} Abh. Math. Sem. Univ. Hamburg 34, 59--72, 1969a
 \absabs
 --- \emph{Lectures on projective planes.} Notes prepared by Mark Pankin and William
 Patton. Lectures given during 1968/69 at the University of Illinois at
 Chicago Circle. Dept. of Mathematics, University of Illinois at Chicago Circle
 1969b
 \absabs
 --- \emph{Einf\"uhrung in die Algebra.} Berlin, etc. 1973
 \absabs
 --- \emph{Translation Planes.} Berlin, Heidelberg, New York 1980
 \absabs
 --- \emph{Tools and Fundamental Constructions of Combinatorial Mathematics.}
 Mann\-heim 1989

 \abs
       Fumitomo Maeda,
 \emph{Kontinuierliche Geometrien.} Berlin, G\"ottingen, Heidelberg 1958

 \abs
       Helmut M\"aurer,
 \emph{Zur Automorphismengruppe der symmetrischen Gruppe.} Mitt. Math. Ges.
 Hamburg XI, 265--266, 1983

 \abs
       Eliakim Hastings Moore,
 \emph{Concerning the abstract Group of order $k!$ and ${1 \over 2}k!$}.
 Proc. London Math. Soc. (1), 28, 357--366, 1897
 \absabs
 --- \emph{Concerning the general equation of the seventh and eighth degrees.} Math.
 Ann. 51, 417--444, 1899

 \abs
       P. J. Morandi, B. A. Sethuraman und J.-P. Tignol, 
 \emph{Division algebras with
 an anti-automorphism but with no involution.} Adv. Geom. 5, 485--495, 2005

 \abs
       Eugen Netto,
 \emph{Substitutionentheorie und ihre Anwendungen auf die Algebra.} Leipzig 1882

 \abs
       Peter Neumann,
 \emph{A Lemma that is not Burnside's.} The Mathematical Scientist 4, 133--141,
 1979 (S\"atze IIII.1.6 und IIII.1.7.)

 \abs
       John von Neumann,
 \emph{Continuous Geometry.} Princeton, London 1960

 \abs
       Theodore Ostrom und Ascher Wagner,
 \emph{On projective and affine planes with transitive collineation groups.} Math.
 Z. 71, 186--199, 1959

 \abs
       G\"unter Pickert,
 \emph{Projektive Ebenen.} Berlin, G\"ottingen, Heidelberg 1955 (Zweite Auflage 1975)

 \abs
       A. Rosenberg,
 \emph{The Structure of the Infinite General Linear Group.} Ann. Math. 68,
 278--297, 1958

 \abs
      G. Salmon,
 \emph{Lettre de Mr. G. Salmon de Dublin \`a l'\'editeur de ce journal.}
 J. reine angew. Math. 39, 365--366, 1850

 \abs
       Issai Schur,
 \emph{Neue Begr\"undung der Theorie der Gruppencharaktere.} Sitzungsberichte der
 Preussischen Akademie der Wissenschaften 1905, Phys.-Math. Klasse,
 406--432. Werke Band 1, 143--169

 \abs
       Helga Tecklenburg,
 \emph{A Proof of the Theorem of Pappos in Finite Desarguesian Affine Planes.}
 J. Geom. 30, 172--181, 1987

 \abs
       O. Veblen und J. W. Young,
 \emph{Projective Geometry.} 2. Auflage. Band I, Boston 1916. Band II, Boston 1917

 \abs
       Ascher Wagner,
 \emph{A theorem on doubly transitive permutation groups.} Math. Z. 85,
 451--453, 1964 (S\"atze IIII.8.4, IIII.8.5, IIII.8.6.)
 \absabs
 --- \emph{On finite affine line transitive planes.} Math. Z. 87, 1--11, 1965
 (Geradenhomogene affine Ebenen sind Translationsebenen. Die\-ser Satz wurde in
 Ostrom und Wagner 1959 vermutet. Satz IIII.9.11.)

 \abs
       Helmut Wielandt,
 \emph{Unendliche Permutationsgruppen.} Vorlesungen an der Univ.\ T\"ubingen 1959/60 (Englische
 \"Ubersetzung in: 
 Mathematische Werke vol.~1. De Gruyter, Berlin 1994).

 \abs
       Max Zorn,
 \emph{Theorie der alternativen Ringe.} Abh. Math. Sem. Univ. Hamburg 8,
 123--147, 1931 (Hierin der Satz IIII.6.14 und der Satz von Artin-Zorn. Beide
 S\"atze stammen laut Zorn von E. Artin.)


\newpage
\addcontentsline{toc}{chapter}{\protect\numberline{Index}}
\printindex
\clearpage  
\end{document}